\pdfoutput=1

\documentclass[onecolumn, 12 pt, fullpage, letterpaper]{report}

\usepackage{amssymb}
\usepackage{textcomp}
\usepackage{graphicx}
\usepackage{graphics}
\usepackage{epsfig}
\usepackage{epstopdf}
\usepackage{float}
\usepackage{color}
\usepackage[cmex10]{amsmath}
\usepackage{latexsym,amsfonts}
\usepackage{amsthm}
\usepackage{url}
\usepackage{longtable}
\usepackage[figuresright]{rotating}
\usepackage{listings}
\usepackage{etoolbox}
\usepackage{sectsty}
\usepackage{array}
\usepackage{bibentry}
\makeatletter\let\saved@bibitem\@bibitem\makeatother
\usepackage{hyperref} 
\makeatletter\let\@bibitem\saved@bibitem\makeatother

\chapterfont{\fontsize{14}{15}\selectfont}
\sectionfont{\fontsize{14}{15}\selectfont}
\subsectionfont{\fontsize{14}{15}\selectfont}
\subsubsectionfont{\fontsize{14}{15}\selectfont}

\usepackage{fancyhdr}       
\fancyheadoffset[L]{0.25in}

\makeatletter
\patchcmd{\@makechapterhead}{50\p@}{20pt}{}{}
\patchcmd{\@makeschapterhead}{50\p@}{20pt}{}{}
\makeatother

\fancypagestyle{plain}{
\fancyhf{} 
\fancyhead[C]{\thepage} 

}

\usepackage{setspace}
\usepackage{subfigure}

   \usepackage{ifthen,xkeyval,xfor,amsgen}
   \usepackage[acronym,toc, nogroupskip]{glossaries}
   \newglossary[slg]{symbols}{syi}{sbl}{List of Symbols}
 
   \makeglossaries
   
\usepackage{amsfonts}
\usepackage{amsmath}
\usepackage{amsthm}
\usepackage{amssymb} 

\usepackage{mdframed} 
\usepackage{thmtools}

\usepackage{enumerate}
\usepackage{pifont}
\usepackage{enumitem}

 \usepackage{xcolor}
\usepackage{color}
\usepackage{graphicx}

\usepackage{algpseudocode} 
\usepackage{algorithmicx}
\usepackage{algorithm}

\newcommand{\R}{\mathbb{R}} 
\newcommand{\N}{\mathbb{N}} 

\newcommand{\cA}{{\cal A}}
\newcommand{\cB}{{\cal B}}
\newcommand{\cC}{{\cal C}}
\newcommand{\cD}{{\cal D}}

\newcommand{\cF}{{\cal F}}
\newcommand{\cG}{{\cal G}}

\newcommand{\cI}{{\cal I}}

\newcommand{\cL}{{\cal L}}
\newcommand{\cM}{{\cal M}}
\newcommand{\cN}{{\cal N}}
\newcommand{\cO}{{\cal O}}
\newcommand{\cP}{{\cal P}}

\newcommand{\cR}{{\cal R}}
\newcommand{\cS}{{\cal S}}

\newcommand{\cU}{{\cal U}}

\newcommand{\cX}{{\cal X}}
\newcommand{\cY}{{\cal Y}}
\newcommand{\cW}{{\cal W}}
\newcommand{\cZ}{{\cal Z}}

\newcommand{\mA}{{\bf A}}
\newcommand{\mB}{{\bf B}}
\newcommand{\mC}{{\bf C}}
\newcommand{\mD}{{\bf D}}
\newcommand{\mE}{{\bf E}}

\newcommand{\mG}{{\bf G}}
\newcommand{\mH}{{\bf H}}
\newcommand{\mI}{{\bf I}}
\newcommand{\mJ}{{\bf J}}

\newcommand{\mL}{{\bf L}}
\newcommand{\mM}{{\bf M}}
\newcommand{\mN}{{\bf N}}
\newcommand{\mO}{{\bf O}}
\newcommand{\mP}{{\bf P}}
\newcommand{\mQ}{{\bf Q}}
\newcommand{\mR}{{\bf R}}
\newcommand{\mS}{{\bf S}}

\newcommand{\mU}{{\bf U}}
\newcommand{\mV}{{\bf V}}
\newcommand{\mW}{{\bf W}}
\newcommand{\mX}{{\bf X}}
\newcommand{\mY}{{\bf Y}}
\newcommand{\mZ}{{\bf Z}}
\newcommand{\qwerty}{t}

\newcommand{\cmark}{{ \color{green} \ding{51}}}%
\newcommand{\xmark}{{ \color{red} \ding{55}}}%

\usepackage[colorinlistoftodos,bordercolor=orange,backgroundcolor=orange!20,linecolor=orange,textsize=scriptsize]{todonotes}

\newcommand{\eqdef}{\stackrel{\text{def}}{=}}
\newcommand{\dotprod}[1]{\left \langle #1\right \rangle} 
\newcommand{\norm}[1]{\left\| #1 \right\|}      

 \newcommand{\argmin}{\arg \min}
 \newcommand{\argmax}{\arg \max}

\DeclareMathOperator{\dom}{dom}         

\DeclareMathOperator{\prox}{prox}       

\DeclareMathOperator{\diag}{\mathbf{Diag}}       

\newcommand{\Diag}[1]{\mathbf{Diag}\left( #1\right)}

\providecommand{\range}[1]{{\rm Range}\left( #1\right)}

\newcommand{\Prob}{\mathbb{P}}
\newcommand{\Var}{{\bf V}}
\newcommand{\Exp}[1]{{\bf E}\left[#1\right] }    
\newcommand{\E}[1]{\mathbb{E}\left[#1\right] } 
\newcommand{\ED}[1]{\mathbb{E}\left[#1\right] } 

\newcommand{\EEE}{\mathbb{E} } 
\newcommand{\EE}[2]{{\bf E}_{#1}\left[#2\right] } 
\providecommand{\Range}[1]{\mathbf{Range}\left( #1\right)}
\newcommand{\Probbb}[1]{\mathbb{P} \left(#1\right)}

\newcommand{\A}{\mathcal{A}}

\declaretheorem[within=section]{definition}
\declaretheorem[sibling=definition]{theorem}
\declaretheorem[sibling=definition]{proposition}
\declaretheorem[sibling=definition]{assumption}
\declaretheorem[sibling=definition]{corollary}
\declaretheorem[sibling=definition]{lemma}

\theoremstyle{remark}
\newtheorem{example}{Example} 
\newtheorem{remark}{Remark} 

\newcommand{\compactify}{} 

\providecommand{\Null}[1]{\mathbf{Null}\left( #1\right)}
\providecommand{\Rank}[1]{\mathbf{Rank}\left( #1\right)}

\providecommand{\Range}[1]{\mathbf{Range}\left( #1\right)}

\providecommand{\range}[1]{\mathbf{Range}\left( #1\right)}
\providecommand{\Tr}[1]{\mathbf{Tr}\left( #1\right)}
\providecommand{\trace}[1]{\mathbf{Tr}\left( #1\right)}
\providecommand{\diag}[1]{\mathbf{Diag}\left( #1\right)}
\providecommand{\E}[1]{\mathbf{E}\left[ #1\right]}

  \providecommand{\dotprod}[1]{\langle #1\rangle} 

\providecommand{\norm}[1]{\lVert#1\rVert}

\newcommand{\tracee}{{ \mathbf{Tr}}}

\newcommand{\bigZ}{{ \mathbf{Z'}}}

\newcommand{\gS}{ {\mS_0}}

\providecommand{\Vect}[1]{\mathbf{Vec}\left( #1\right)}

\newcommand{\Proj}{{\bf \Pi}} 
\newcommand{\Jac}{{ \bf \nabla F}} 
\newcommand{\ones}{e}

\newcommand{\LL}{\mathcal{L}}

\newcommand{\Lgen}{{\Phi}}
\newcommand{\Lnacc}{{\Psi}}
\newcommand{\Lacc}{{\Upsilon}}

\newcommand{\mVdiag}{{\bf \hat{V}}}
\newcommand{\Probmat}{{\mP}}
\newcommand{\mPdiag}{{\bf \hat{P}}}
\newcommand{\TD}{{\rm \eta}}

\algnewcommand{\IIf}[1]{\State\algorithmicif\ #1\ \algorithmicthen}
\algnewcommand{\EndIIf}{\unskip\ \algorithmicend\ \algorithmicif}

\newcommand{\sumin}{\sum_{i=1}^n}

\newcommand{\avein}{\frac{1}{n}\sum_{i=1}^n}
\newcommand{\avejn}{\frac{1}{n}\sum_{j=1}^n}
\newcommand{\proxR}{\prox_{\alpha \psi}}

\newcommand{\RR}{\mathbb{R}}

\newcommand{\smo}{f'}
\newcommand{\seganu}{{\omega}}
 \newcommand{\quadb}{{o}}

\def\<#1,#2>{\left\langle #1,#2\right\rangle}

\newcolumntype{M}[1]{>{\centering\arraybackslash}m{#1}}
\newcolumntype{N}{@{}m{0pt}@{}}

\newcommand{\eL}{{\color{red}e}} 
\newcommand{\eR}{{\color{blue}e}} %
\newcommand{\eRR}{{\color{blue}e_R}} 
\newcommand{\pL}{{\color{red}p}} 
\newcommand{\pR}{{\color{blue}p}} 
\newcommand{\eLi}{{\color{red}e_i}} 
\newcommand{\eRj}{{\color{blue}e_j}} 
\newcommand{\pLi}{{\color{red}p_i}} 
\newcommand{\pRj}{{\color{blue}p_j}} 
\newcommand{\pLL}{{\color{red}p_L}} 
\newcommand{\pRR}{{\color{blue}p_R}} 
\newcommand{\dL}{{\color{red}d}} 
\newcommand{\nR}{{\color{blue}n}} 

\newcommand{\NRt}{{\color{blue}N_t}} 

\newcommand{\tR}{{t}} 
\newcommand{\TR}{{T}} 

\newcommand{\cDR}{{\color{blue}\cD}} 
\newcommand{\cDL}{{\color{red} \cD}} 

\newcommand{\PR}{{\color{blue} \bf P}} 
\newcommand{\qR}{{\color{blue}q}} 

\newcommand{\qRj}{{\color{blue}q_j}} 

\newcommand{\ptRj}{{\color{blue} p^t_j}} 
\newcommand{\ptLi}{{\color{red}p}_{\color{red} i}^{\color{blue}t}} 

\newcommand{\ptR}{{\color{blue} p^t}} 
\newcommand{\ptL}{{\color{red}p^{\color{blue} t}}} 

\newcommand{\PtR}{{\color{blue} {\bf P}^t}} 

\newcommand{\qtR}{{\color{blue} q^t}} 
\newcommand{\qtRj}{{\color{blue} q^t_j}} 

\newcommand{\probx}{{\color{black} \rho}}
\newcommand{\proby}{{\color{cyan} \delta }} 

\newcommand{\ugly}{{\eta}}

\newcommand{\NORMG}[1]{\left \| #1\right  \|^2}

\newcommand{\piop}{{\Gamma}}

\newcommand{\Popt}{{\color{black}\mW}} 

\newcommand{\Ind}[1]{I_{#1} } 
\newcommand{\Lift}[1]{U \left(#1\right) }

\newcommand{\pp}{{\color{red}\tilde{p}}} 
\newcommand{\vv}{{\color{red}\tilde{v}}}
\newcommand{\ff}{{\color{red} \tilde{f}} }

\newcommand{\LLL}{{\color{red} \tilde{L}} }
\newcommand{\mmu}{{\color{red} \tilde{\mu}} }
\newcommand{\dd}{{\color{red} \tilde{d}} }
\newcommand{\sS}{{\color{red} \tilde{S}} }
\newcommand{\ppp}{{\pp}} 
\newcommand{\ppsi}{{\color{red} \tilde{\psi}} }
\newcommand{\aalpha}{{\color{red} \tilde{\alpha}} }
\newcommand{\eeta}{{\color{red} \tilde{\eta}} }
\newcommand{\ggamma}{{\color{red} \tilde{\gamma}}} 

\newcommand{\xx}{{\color{red} \tilde{x}} }
\newcommand{\yy}{{\color{red} \tilde{y}} }
\newcommand{\zzz}{{\color{red}\tilde{ z}} }

\newcommand{\ww}{{\color{red} \tilde{w}} }
\newcommand{\ggggg}{{\color{red} \tilde{g}} }
\newcommand{\ueWMC}{{\color{red}  {\bf \tilde{M}}}}
\newcommand{\ee}{{\color{red}  \tilde{e}}}
\newcommand{\PpP}{{\color{red} \tilde{ F}}}

\newcommand{\ccL}{{\cal L}}
\newcommand{\Odd}{{O_{dd} }}

\newcommand{\Lcac}{{\ccL'}}
\newcommand{\LcLc}{{\color{red} \tilde{ \ccL}}} 

\newcommand{\BD}{{\bf D_B}}
\newcommand{\blockdiag}{\text{BlockDiag}}

\newcommand{\Corrloc}[1]{{\text{Corr} \left[ #1\right]}} 

\newcommand{\flocc}{{\color{red}f}} 

\newcommand{\nlocc}{N} 
\newcommand{\tRloc}{{i}} 
\newcommand{\TRloc}{{n}} 

\newcommand{\NRtloc}{m_i} 

\newcommand{\pagg}{{p}} 

\newcommand{\onestloc}{{\bf 1}^{(\NRt)}}  
\newcommand{\onesmloc}{{\bf 1}}

\newcommand{\pRloc}{{\color{blue}p}_i} 
\newcommand{\pRlocjloc}{{\color{blue}p}_{i,j}} 
\newcommand{\ptRloc}{{\color{blue}p}^i} 
\newcommand{\Lloc}{{\color{red} \tilde{L}}} 

\newcommand{\maploc}{\Omega} 
\newcommand{\xxloc}{{\color{blue} \tilde{x}}} 
\newcommand{\xbloc}{{\color{blue} x}} 

\newcommand{\mmM}{{\bf Q}}
\newcommand{\mJf}{{\bf J}}
\newcommand{\mJpsi}{{\bf \Psi}}

\newcommand{\nloc}{{\color{green} n}}
\newcommand{\floc}{{\color{green} f}}
\newcommand{\Philoc}{{\color{green} \psi}}

\newcommand{\CS}{{\color{purple}S}} 
\newcommand{\ptg}{{\color{purple}p}_i} 
\newcommand{\pg}{{\color{purple}p}}

\newcommand{\proxop}{\mathop{\mathrm{prox}}\nolimits}

\newcommand{\proxt}{\proxop_{\alpha \psi_i}}
\newcommand{\psvrg}{{\color{red} p}} 
\newcommand{\vg}{{\color{green}v}} 
\newcommand{\mug}{{\color{green}\mu}} 
\newcommand{\Ug}{{\color{green}\Upsilon}}

			\newcommand{\level}{\chi^0 }

\def\ba{\begin{array}}
\def\ea{\end{array}}
\def\beann{\begin{eqnarray*}}
\def\eeann{\end{eqnarray*}}
\def\bea{\begin{eqnarray}}
\def\eea{\end{eqnarray}}

\def\BMP{\begin{minipage}{9.5cm}}
\def\EMP{\end{minipage}}
\def\MPT{\begin{minipage}{11.5cm}}
\def\EPT{\end{minipage}}

\def\Def{\stackrel{\mathrm{def}}{=}}

\def\dom{{\rm dom \,}}
\def\beq{\begin{equation}}
\def\eeq{\end{equation}}

\newcommand{\refLE}[1]{\ensuremath{\stackrel{(\ref{#1})}{\leq}}}

\def\la{\langle}
\def\ra{\rangle}

\newacronym{SEGA}{{\tt SEGA}}{SkEtched Gradeint Algorithm}
\newacronym{ACD}{{\tt ACD}}{Accelerated {\tt CD}}
\newacronym{AGD}{{\tt AGD}}{Accelerated {\tt GD}}
\newacronym{PGD}{{\tt PGD}}{Proximal {\tt GD}}
\newacronym{CD}{{{\tt CD}}}{Coordinate Descent}
\newacronym{GJS}{{\tt GJS}}{Generalized Jacobian Sketching}
\newacronym{SGD}{{\tt SGD}}{Stochastic {\tt GD}}
\newacronym{ISEGA}{{\tt ISEGA}}{Independent {\tt SEGA}}
\newacronym{SSCN}{{\tt SSCN}}{Stochastic Subspace Cubic Newton}
\newacronym{SVRCD}{{\tt SVRCD}}{Stochastic Variance Reduced {\tt CD}}
\newacronym{ASVRCD}{{\tt ASVRCD}}{Accelerated  {\tt SVRCD}}
\newacronym{AMI}{{\tt AMI}}{Accelerated Matrix Inversion}
\newacronym{VR}{VR}{Variance Reduction}
\newacronym{LGD}{{\tt LGD}}{Local {\tt GD}}
\newacronym{ESO}{ESO}{Expected Separable Overapproximation}
\newacronym{ERM}{ERM}{Empirical Risk Minimizatiom}
\newacronym{L2GD}{{\tt L2GD}}{Loopless {\tt LGD}}
\newacronym{BFGS}{{\tt BFGS}}{Broyden-Fletcher-Goldfarb-Shanno}
\newacronym{FL}{FL}{Federated Learning}
\newacronym{LSGD}{{\tt LSGD}}{Local {\tt SGD}}
\newacronym{L2SGD}{{\tt L2SGD}}{{\tt LSGD}}
\newacronym{GD}{{\tt GD}}{Gradient Descent}
\newacronym{ASEGA}{{\tt ASEGA}}{Accelerated {\tt SEGA}}

\newacronym{kaust}{KAUST}{King Abdullah University of Science and Technology}
\newacronym{IBCD}{{\tt IBCD}}{Independent Block Coordinate Descent}

\usepackage[lmargin=1.5in, rmargin=1in, vmargin=1in, headsep=0.083334in]{geometry}
\newcommand{\mathsym}[1]{{}}
\newcommand{\unicode}[1]{{}}

\renewcommand\bibname{\centering BIBLIOGRAPHY}

\begin{document}


\vspace{2pt}
\thispagestyle{empty}
\addvspace{10mm}

\begin{center}

{\textbf{{\large Optimization for Supervised Machine Learning: \\ Randomized Algorithms for Data and Parameters}}}\vfill 
{Dissertation by}\\
{ Filip Hanzely}\vfill

{ In Partial Fulfillment of the Requirements}\\[12pt]
{ For the Degree of}\\[12pt]
{Doctor of Philosophy} \vfill
{King Abdullah University of Science and Technology }\\
{Thuwal, Kingdom of Saudi Arabia}
\vfill
{August, 2020}

\end{center}

\newpage

\begin{center}

\end{center}

\begin{center}

{ \textbf{{\large EXAMINATION COMMITTEE PAGE}}}\\\vspace{1cm}

\end{center}
\noindent{The dissertation of Filip Hanzely is approved by the examination committee}
\addcontentsline{toc}{chapter}{Examination Committee Page}

\vspace{4\baselineskip}

\noindent{Committee Chairperson: Peter Richt\'arik}\\
Committee Members: Stephen Wright, Tong Zhang, Ra\'ul Fidel Tempone, Bernard Ghanem \vfill


\newpage
\addcontentsline{toc}{chapter}{Copyright}
\vspace*{\fill}
\begin{center}
{ \copyright \; August, 2020}\\
{Filip Hanzely}\\
{All Rights Reserved}
\end{center}

\singlespacing

\begin{center}

\end{center}

\begin{center}
{{\bf\fontsize{14pt}{14.5pt}\selectfont \uppercase{ABSTRACT}}}
\end{center}

\addcontentsline{toc}{chapter}{Abstract}

\begin{center}
{{\fontsize{14pt}{14.5pt}\selectfont {Optimization for Supervised Machine Learning:\\ Randomized Algorithms for Data and Parameters\\
\phantom{xx}
 Filip Hanzely}}}
\end{center}

Many key problems in machine learning and data science are routinely modeled as optimization problems and solved via optimization algorithms. With the increase of the volume of data and the size and complexity of the statistical models used to formulate these often ill-conditioned optimization tasks, there is a need for new efficient algorithms able to cope with these challenges. 

In this thesis, we deal with each of these sources of difficulty in a different way. To efficiently address the big data issue, we develop new methods which in each iteration examine a small random subset of the training data only.  To handle the big model issue, we develop methods which in each iteration  update a random subset of the model parameters only. Finally, to deal with ill-conditioned problems, we devise methods that incorporate either higher-order information or Nesterov's acceleration/momentum. In all cases, randomness is viewed as a powerful algorithmic tool that we tune, both in theory and in experiments, to achieve the best results.

Our algorithms have their primary application in training supervised machine learning models via regularized empirical risk minimization, which is the dominant paradigm for training such models. However, due to their generality, our methods can be applied in many other fields, including but not limited to data science, engineering, scientific computing, and statistics.



\begin{center}

\end{center}

\begin{center}

{\bf\fontsize{14pt}{14.5pt}\selectfont \uppercase{Acknowledgements}}\\\vspace{1cm}
\end{center}

\addcontentsline{toc}{chapter}{Acknowledgements}

I owe my deepest gratitude to my supervisor Peter Richt\'arik. Thank you very much for your guidance; it allowed me to get the best out of myself. Thanks a lot for the extraordinary support, career advice, and tons of encouragement. You showed me each aspect of being a complete researcher and always guided me in that direction. 

Next, I would like to thank all members of our research group for countless stimulating discussions, namely: Konstantin Mishchenko, Samuel Horv\'ath, Slavom\'ir Hanzely, Robert Gower, Aritra Dutta, Nicolas Loizou, Alibek Sailanbayev, Jakub Kone\v{c}n\'y, Dominik Csiba, Elnur Gasanov, Eduard Gorbunov, Dmitry Kovalev, Adil Salim, Yazeed Basyoni, Mher Safaryan, El Houcine Bergou, Xun Qian, Zhize Li, and Egor Shulgin.

I am very grateful to all the great researchers I had a chance to collaborate with, especially Lin Xiao, Yurii Nesterov, Sebastian Stich, Jingwei Liang, and Nikita Doikov. I would also like to thank Michael Mahoney, Martin Jaggi, Alex D'Aspremont, Adrien Taylor, Praneeth Karimireddy, and Haihao Lu for multiple fruitful discussions. Further, I owe a big thanks to my internship hosts Rodolphe Jenatton and Sashank Reddi at Amazon and Google respectively as well as to other people I had a chance to interact with, namely Mathias Seeger, Srinadh Bhojanapalli, C\'edric Archambeau and Sanjiv Kumar. I learned a lot from all of you!

I appreciate a lot all the support I received both from KAUST and from the Visual Computing Center at KAUST; I feel extremely lucky for all the opportunities I had. I am also very grateful to my defense committee, namely Stephen J Wright, Tong Zhang, Ra\'{u}l F Tempone, and Bernard Ghanem.

 I would like to thank all my friends that made my stay at KAUST pleasant. Last but not least, I am eminently grateful to my family for their love and support.




\addcontentsline{toc}{chapter}{Table of Contents}
\renewcommand{\contentsname}{\centerline{\textbf{{\large TABLE OF CONTENTS}}}}
\tableofcontents
\cleardoublepage

\printglossary[type=\acronymtype,style=long3col, title=\centerline{LIST OF ABBREVIATIONS}, toctitle=List of Abbreviations, nonumberlist=true] 


\cleardoublepage
\addcontentsline{toc}{chapter}{\listfigurename} 
\renewcommand*\listfigurename{\centerline{LIST OF FIGURES}} 
\listoffigures

\cleardoublepage
\addcontentsline{toc}{chapter}{\listtablename}
\renewcommand*\listtablename{\centerline{LIST OF TABLES}} 
\listoftables



%

\nobibliography*

\singlespacing 
\chapter{Introduction}

Over the past several decades, optimization has become a key tool in the toolbox of modern technology, enabling a multitude of areas of engineering, computer science, physics, economics, finance, chemistry, computational biology, as well as many other fields of human endeavor. 

In this thesis, we predominantly focus on continuous optimization problems arising in the training of supervised machine learning models\footnote{Our results are applicable beyond supervised machine learning (training of regression/classification models), as we shall see. Supervised machine learning is, however, the primary application we have in mind.}. Informally, the training of such models  can be described as the search for the parameters characterizing the model that best fits the observed data. In particular, the dominant  paradigm for solving supervised machine learning problems is to cast them as {\em regularized empirical risk minimization} (ERM)  problems, often also called {\em finite-sum} optimization problems, which take the form
 \begin{equation} \label{eq:finitesum}
\min_{x\in \R^d} \left\{  F(x) \eqdef \underbrace{\frac1n \sum_{i=1}^n f_i(x)}_{\eqdef f(x)} + \psi(x) \right\}.
 \end{equation}
In the above problem,  the vector $x\in \R^d$ represents the parameters describing the  model we wish to train (e.g., support vector machine, logistic regression or a neural network), the function $f_i$ measures the misfit of model $x$ with respect to the $i$th data point, and function $\psi:\R^d \to \R\cup \{+\infty\}$ is a regularizer whose role is to incorporate prior  information  or impart  desirable properties onto the model. The objective function $F$ measures the (regularized) empirical loss of model $x$.
 
 The training of machine learning models carries a multitude of  challenges, with the two most pronounced being  the size of the training dataset (i.e., big $n$) and the size of the model (i.e., big $d$). {\em Big data} and {\em big model} scenarios render standard deterministic optimization methods, such as gradient descent and Newton's method, inefficient at solving \eqref{eq:finitesum}. In the past decade, this  led to a ``Cambrian explosion'' of  new iterative algorithms utilizing {\em randomness} in various ingenious ways aimed at addressing the big data and big model problems. In order to identify a model of suitable (optimization or generalization) properties, these new randomized methods typically rely on significantly cheaper iterations than their deterministic counterparts at the cost of requiring many more iterations. However, the benefits  of such an approach often vastly outweigh the costs, both in theory and in practice, which makes them the methods of choice in the big data or big model regime. The per-iteration savings are due to 
the inclusion of suitable randomization strategies such as {\em subsampling the data}, i.e., working with a small subset of the functions $f_i$ in each iteration only, or  {\em subsampling the parameters}, i.e., updating a small subset of the parameters in each iteration only.

Informally speaking, the main goal of this thesis is to develop, under appropriate assumptions on the properties of the regularized empirical loss function $F$, through its constituents $\{f_i\}_{i=1}^n$ and $\psi$, {\em new  state-of-the-art randomized optimization algorithms} for solving the ERM problem~\eqref{eq:finitesum}, both in theory (by establishing improved convergence and complexity results) and in practice (by extensive experimental testing on synthetic and real data). While the structure of $F$  varies slightly among the individual chapters of this thesis, we mostly assume that $f$ is differentiable and convex, while $\psi$ is convex, possibly non-smooth, but assumed to be {\em proximable}\footnote{The well-known notion of ``proximability'' is formally introduced in Section~\ref{sec:proxgrad}.}.

\section{Technical preliminaries and basic algorithms}

In this section, we introduce typical assumptions that we impose on the functions $\{f_i\}$ and $\psi$ appearing in \eqref{eq:finitesum} throughout the individual chapters, as well as introduce standard tricks and results in optimization which we build upon in this work. We describe gradient descent---the cornerstone of first order optimization---followed by three standard tricks from the literature that gradient descent can be furnished with: Nesterov's acceleration~\cite{nesterov83}, proximal operator~\cite{beck2009fista} and randomness~\cite{robbins}.

We shall first equip $\R^d$ with an inner  product and a norm. The standard Euclidean inner product of vectors $x,y\in \R^d$ is  $\langle x,y\rangle  \eqdef \sum_{i=1}^d  x_i\cdot y_i $ and the (induced) Euclidean norm is $\| x\| \eqdef \langle x,x\rangle ^{1/2}$. For the reader's convenience, we present a table of frequently used notation in Appendix~\ref{sec:table}.

\subsection{Smoothness and convexity}

We now introduce two key concepts  which will be used in various places throughout this text: {\em convexity} and {\em smoothness}.  We will often assume that the objective $F$ (or some part of $F$) is convex and smooth. The exact assumptions used differ from chapter to chapter, and are described therein. Let us first start with (strong) convexity.
\begin{definition}[Strong convexity and convexity]
Let $\mu \geq 0$. Function $h:\R^d\to \R$ is $\mu$-strongly convex if for all $x,y \in \R^d$, and all $t\in [0,1]$:
\[
h(tx + (1-t)y)\leq th(y) +(1-t) h(x) - \frac{\mu t(1-t)}{2} \|x-y \|^2.
\]
In the special case where $\mu=0$, we say that $h$ is convex. 
\end{definition} 

The following standard result states  that for a sufficiently smooth function $h$,  strong convexity provides a global quadratic (or linear in the $\mu=0$ case) lower bound on $h$  and a uniform lower bound on the eigenvalues of its Hessian. 
 
 \begin{proposition}[Nesterov~\cite{nesterov2018lectures}]
Let $h: \R^d \rightarrow \R$ be differentiable. Then, $h$ is $\mu$-strongly convex if and only if for all $x,y\in \R^d$:
\[
h(x)\geq h(y) + \<\nabla h(y),x-y > + \frac{\mu}{2}\| x-y\|^2.
\]
If $h$ is further twice differentiable, it is $\mu$-strongly convex if and only if for all $x\in \R^d$ we have $\nabla^2 h(x) \succeq \mu \mI$, where $\mI \in \R^{d\times d}$ is the identity matrix. 
 \end{proposition}

Next, we introduce a typical {\em smoothness} assumption we make throughout the thesis. 

\begin{definition}[$L$-smoothness]
Differentiable function $h:\R^d\to \R$ is  $L$-smooth if it has $L$-Lipschitz gradient, namely for all $x,y \in \R^d$:
\[
\| \nabla h(x)- \nabla h(y)\| \leq L \| x-y\|.
\]
\end{definition} 

Analogously to strong convexity, smoothness provides us with both an upper bound on the function value as well as with an upper bound on the Hessian at each point in the domain.

\begin{proposition}[Nesterov~\cite{nesterov2018lectures}]
A differentiable function $h: \R^d \rightarrow \R$ is $L$-smooth if and only if for all $x,y\in \R^d$:
\begin{equation}\label{eq:intro_smoothness}
h(y) + \<\nabla h(y),x-y > - \frac{L}{2}\| x-y\|^2 \leq h(x)\leq h(y) + \<\nabla h(y),x-y > + \frac{L}{2}\| x-y\|^2.
\end{equation}
If $h$ is further twice differentiable, it is $L$-smooth if and only if for all $x\in \R^d$ we have $\nabla^2 h(x) \preceq L \mI$, where $\preceq$ designates the L\"{o}wner ordering of matrices.
 \end{proposition}

We are now ready to present the backbone of the world of  first-order optimization algorithms---gradient descent---along with a few basic and well known extensions.

\subsection{Gradient descent}
For the sake of expositional simplicity,  consider optimization problem \eqref{eq:finitesum} in its simplest form: $\psi \equiv 0$ and $n=1$. That is, we consider the unregularized case and  ignore the finite-sum structure of $f$. In this case, $F=f$. 

Note that if $f$  is $L$-smooth, the second inequality in~\eqref{eq:intro_smoothness} provides us with a global convex quadratic upper bound on $f$ using zero and first-order information about $f$ at arbitrary point $y$: 
\begin{equation}\label{eq:nbig98gidd}
f(x) \leq f(y) + \left \langle \nabla f(y), x-y \right \rangle + \frac{L}{2}\norm{x-y}^2. 
\end{equation}
Minimizing this upper bound in the variable $x$ gives 
\[x = y - \frac{1}{L}\nabla f(y).\]
 Doing this iteratively, we arrive at the famous {\em gradient descent} method (Algorithm~\ref{alg:gd}), which is a trivial baseline we build on throughout this thesis.

\begin{algorithm}[!h]
\begin{algorithmic}[1]
\State \textbf{Input:} Starting point $x^0\in \R^d$, smoothness constant $L>0$ 
\For {$k= 0,1, 2, \dots $} 
\State $x^{k+1}=x^k - \frac1L \nabla f(x^k)$
 \EndFor
\end{algorithmic}
\caption{Gradient descent ({\tt GD})}
\label{alg:gd}
\end{algorithm} 

Convergence properties of gradient descent are described in Proposition~\ref{prop:GDrate}. This standard result posits sublinear convergence for the class of smooth and convex functions and linear convergence for the class of smooth and strongly convex functions. 

\begin{proposition}[Nesterov~\cite{nesterov2018lectures}]\label{prop:GDrate}
Let $f^*\eqdef \min_{x\in \R^d} f(x)$ and $x^*\eqdef \argmin_{x\in \R^d} F(x) =  \argmin_{x\in \R^d} f(x)  $. Suppose that the sequence of iterates $\{x^k\}_{k=0}^\infty$ is generated by Algorithm~\ref{alg:gd}. If $f$ is $L$-smooth and convex, then
\begin{equation*}
f\left(x^{k}\right)-f^{*} \leq \frac{2 L\left\|x^{0}-x^{*}\right\|^{2}}{k+4}.
\end{equation*}
If we additionally assume that $f$ is $\mu$-strongly convex\footnote{This implies that, necessarily, $\mu \leq L$.}, then
\[
f\left(x^{k}\right)-f^{*} \leq \left(\frac{L-\mu}{L+\mu}\right)^{2 k} \frac{L \left\|x^{0}-x^{*}\right\|^{2}}{2}.
\]
\end{proposition}

Next, we introduce a handful of tricks that can be incorporated on top of gradient descent: {\em Nesterov's acceleration}, {\em proximal operator}, and {\em stochasticity}. As we shall see, these tricks are mutually ``orthogonal'',  which means that, generally speaking, they can be built on top of each other for a more pronounced additive benefit. 

\subsection{Nesterov's acceleration \label{sec:intro_acceleration}}

Notice that gradient descent is a {\em greedy} method. Indeed,  the next iterate is constructed to find  a point with the smallest guaranteed function value given the information we have about $f$: zero and first-order information about $f$ at the current iterate, and the smoothness parameter $L$. As a byproduct, gradient descent forgets all the past information gathered throughout the optimization process. It turns out that in this case, greediness as an algorithmic design tool is suboptimal since appropriate use of history can yield to a significant improvement in iteration complexity. Nesterov's {\em accelerated gradient descent} method (stated as Algorithm~\ref{alg:acd_intro}) is an algorithm that achieves this. 

\begin{algorithm}[!h]
\begin{algorithmic}[1]
\State \textbf{Input:} Starting point $x^0=y^0 \in \R^d$, smoothness constant $L>0$,  strong convexity $\mu \geq0$
\For {$k= 0,1, 2, \dots $} 
\State $x^{k+1}=y^k - \frac1L \nabla F(y^k)$
\If{$\mu = 0$}
\State $y^{k+1}=x^{k+1}+\frac{k}{k+3}\left(x^{k+1}-x^{k}\right)$
\Else
\State $y^{k+1}=x^{k+1}+\frac{\sqrt{L}-\sqrt{\mu}}{\sqrt{L}+\sqrt{\mu}}\left(x^{k+1}-x^{k}\right)$
\EndIf
 \EndFor
\end{algorithmic}
\caption{Nesterov's accelerated gradient descent ({\acrshort{AGD}})}
\label{alg:acd_intro}
\end{algorithm} 

The following proposition describes the convergence rate of Nesterov's accelerated gradient descent method. 

\begin{proposition}[Nesterov~\cite{nesterov2018lectures, nesterov83}]
Suppose that sequence $\{x^k\}_{k=0}^\infty$ was generated by Algorithm~\ref{alg:acd_intro}. If $f$ is $L$-smooth and $\mu$-strongly convex\footnote{We allow for $\mu=0$.}, then 
\[
f\left(x^{k}\right)-f^{*} \leq \min \left\{\left(1-\sqrt{\frac{\mu}{L}}\right)^{k}, \frac{4 }{\left(2+k\right)^{2}}\right\} \left(f\left(x^{0}\right)-f^{*}+\frac{L}{2}\left\|x^{0}-x^{*}\right\|^{2}\right).
\]

\end{proposition}

Up to a constant factor, the method only requires a square root of the number of iterations needed by gradient descent.

It is important to mention that Algorithm~\ref{alg:acd_intro} is {\em optimal} in terms of oracle complexity for both smooth convex and smooth strongly convex problems as it (up to a constant) matches the corresponding lower bound~\cite{nesterov2018lectures}.

\subsection{Proximal operator and proximal gradient descent \label{sec:proxgrad}}
In their simplest form, neither gradient descent nor Nesterov's accelerated gradient descent are applicable in the presence of a non-smooth regularizer $\psi$. In this section, we offer a brief overview of the proximal gradient descent method, which is capable of solving~\eqref{eq:finitesum} for any convex closed\footnote{We say that a convex function $h:\R^d\rightarrow \R \cup \{ +\infty\}$ is closed if for any $c\in \R$ the sublevel set $\{x 
;|\; h(x) \leq c \}$ is a closed set.} provided that $\psi$  is \emph{proximable}, which means that the \emph{proximal operator} of $\psi$, defined as
\begin{eqnarray}\label{eq:intro_proxdef}
    \prox_{\alpha \psi} (x) \eqdef \argmin_{y\in \R^d} \left\{\psi(y) + \frac{1}{2\alpha}\|y - x\|^2  \right\},
\end{eqnarray}
where $\alpha>0$,  is easily computable (e.g., in closed form). 
 
 \begin{example}In the following two examples we give formulas for the proximal operators of two commonly used regularizers. 
\begin{itemize}
\item In some applications, $\psi$ is used to represent a hard constraint on model $x$. In particular, let $Q \subseteq \R^d$ be any nonempty closed convex set. It is easy to see that the optimization problem \[ \min_{x\in Q} f(x) \] can be equivalently written in the form~\eqref{eq:finitesum} by setting $\psi$ to be the ``indicator'' function of $Q$: 
\[\min_{x\in \R^d} f(x) + \psi(x), \qquad \psi(x)  \eqdef \begin{cases}
0 & \text{if} \;\; x\in Q \\
\infty & \text{if} \;\; x \not \in Q
\end{cases}.
\]
Consequently, the proximal operator of $\psi$ becomes the projection operator onto $Q$, i.e.,
\[
\forall \alpha > 0: \quad   \prox_{\alpha \psi} (x) = \argmin_{y\in Q} \| x-y\|^2.
\] 
\item In applications where one prefers a sparse solution $x^*$, one can set $\psi$ to be the sparsity-inducing $\ell_1$ norm: $\psi(x) =  \| x\|_1 \eqdef \sum_{i=1}^d |x_i|$. In a such case, the proximal operator of $\psi$ is equivalent to applying elementwise soft-thresholding; i.e., 
\[\forall \, \alpha\geq 0: \quad  \left(\prox_{\alpha \psi} (x)\right)_i  = \begin{cases}
0 & \text{if} \;\; |x_i|\leq \alpha \\
\text{sign}(x_i) (|x_i| -\alpha )& \text{if} \;\;  |x_i| > \alpha
\end{cases}
\]
  for all $i \in \{ 1,2,\dots, d\}$. 
\end{itemize}
\end{example}

Clearly, in both these examples, its is not possible to apply gradient descent to minimize $F=f+\psi$ since $F$ is not differentiable due to the presence of $\psi$. While it is possible to replace gradient with subgradient\footnote{Informally speaking, subgradient is a generalization of the gradient for convex, possibly non-differentiable functions.}---resulting in the  \emph{subgradient method}---such an approach suffers from inferior convergence guarantees~\cite{nesterov2018lectures, Grimmer2017}.

 Alternatively, we might take advantage of  the proximability of $\psi$ and incorporate the proximal operator into the optimization procedure. The most natural approach is to alternate the gradient step with the proximal step, which results in the {\em proximal gradient descent} method ({\acrshort{PGD}})~\cite{beck2009fista, beck17book}. If the regularizers considered in the above example are used, the method is often alternatively known under the name  {\em projected gradient descent} and {\tt ISTA}, respectively.

\begin{algorithm}[!h]
\begin{algorithmic}[1]
\State \textbf{Input:} Starting point $x^0 \in \R^d$, smoothness constant $L>0$
\For {$k= 0,1, 2, \dots $} 
\State $x^{k+1}=\prox_{\frac1L \psi} \left( x^k - \frac1L \nabla f(x^k)\right)$
 \EndFor
\end{algorithmic}
\caption{Proximal gradient descent ({\tt PGD})}
\label{alg:pgd_intro}
\end{algorithm} 

The following proposition describes the convergence rate of proximal gradient descent (Algorithm~\ref{alg:pgd_intro}). The method is, up to a small constant factor, as fast as gradient descent. Consequently, incorporating the regularizer $\psi$ into the optimization method does not hurt the convergence rate. In some cases, the presence of $\psi$ might make optimization easier; we will elaborate on this soon.

\begin{proposition}[Beck~\cite{beck17book}]
Let $F^*\eqdef \min_{x\in \R^d} F(x)$ and $x^*\eqdef \argmin_{x\in \R^d} F(x) $. Suppose that the sequence $\{x^k\}_{k=0}^\infty$ is generated by Algorithm~\ref{alg:pgd_intro}. If $f$ is $L$-smooth and convex, we have
\[
F(x^k)-F^*\leq \frac{L \left\|x^{0}-x^{*}\right\|^{2} }{2k}.
\]
If further $f$ is $\mu$-strongly convex (with $\mu>0$), we have 
\[
\| x^k -x^*\|^2 \leq \left(1-\frac{\mu}{L}\right)^{k}\left\|x^{0}-x^{*}\right\|^{2}.
\]
\end{proposition}

\subsection{Incorporating randomness \label{sec:intro_randomized}}
All of the optimization algorithms introduced so far are agnostic to the finite-sum structure of $f$. Consequently, if $f_i$ measures the misfit of the current model $x$ at the $i$th datapoint, both {\tt GD} and {\tt AGD} are passing through the entire dataset every iteration. The larger the number of datapoints $n$ is, the more expensive it is to perform an iteration of these methods, which makes them impractical. 

How can one effectively deal with big $n$ then? The most natural approach is simply to replace the expensive computation of the (full) gradient of $f$,
\[\nabla f(x) = \frac{1}{n} \sum_{i=1}^n \nabla f_i(x),\]
via a cheap stochastic approximation thereof, resulting in the celebrated {\em stochastic gradient descent} ({\tt SGD}) method~\cite{robbins}, which we state as Algorithm~\ref{alg:sgd_intro}.

For simplicity of exposition, we have once again adopted the assumption that $\psi \equiv 0$. However, as we shall show later, the same result holds in the regularized case by incorporating the proximal operator into the algorithm.

\begin{algorithm}[!h]
\begin{algorithmic}[1]
\State \textbf{Input:} Starting point $x^0 \in \R^d$, stepsize $\alpha>0$
\For {$k= 0,1, 2, \dots $} 
\State Sample $i \in \{1,2\dots, n\} $ uniformly at random
\State $x^{k+1}= x^k - \alpha \nabla f_i(x^k)$
 \EndFor
\end{algorithmic}
\caption{Stochastic gradient descent ({\tt SGD})}
\label{alg:sgd_intro}
\end{algorithm} 

Convergence rate of {\tt SGD} is presented in Proposition~\ref{prop:sgdrate}. In particular, under smoothness and strong convexity assumptions, {\tt SGD} enjoys a fast, linear rate to a specific neighborhood of the optimum.

\begin{proposition}[Nguyen et al~\cite{nguyen2018sgd}, Gower et al~\cite{pmlr-v97-qian19b}]\label{prop:sgdrate}
Suppose that function $f_i$ is $L$-smooth and convex for all $i$, while function $f$ is $\mu$-strongly convex with $\mu>0$. Then, for any $\alpha \leq \frac{1}{2L}$ we have
\[
\E{\| x^k - x^*\|^2} \leq \left( 1-\alpha \mu \right)^k \| x^0 - x^*\|^2 + \frac{2\alpha \sum_{i=1}^n \|\nabla f_i(x^*) \|^2}{n\mu}.
\]
\end{proposition}

Note that if the gradients $\nabla f_i(x^*)$ are all zero, which typically happens for over-parameterized models, the above result posits a linear convergence rate to the optimal solution $x^*$. In general,  the right hand side in the complexity guarantee can be made arbitrarily small by choosing the stepsize $\alpha$ sufficiently small and $k$ sufficiently large. Alternatively, this can be achieved by choosing a suitable decreasing stepsize schedule. However, such adjustments  will lead to a worse convergence rate: we get a $\tilde{O}(1/k)$ rate towards the true optimum. It is possible to  preserve the linear rate even if the gradients at the optimum are not zero, but for this to happen, one needs to adjust SGD to employ one of the many {\em variance-reduction} techniques proposed in the literature.

\section{From finite sum to coordinate descent and back \label{sec:history}} 

Let us consider a very specific form of the objective~\eqref{eq:finitesum}: assume that for all $i$, function $f_i$ corresponds to a loss of a linear model,\footnote{Our results go beyond linear models. The assumption is made here in order to provide a simple motivation for coordinate descent methods.} i.e., $f_i(x) =  \phi_i( \la a_i, x \ra )$ for some convex $\phi: \R \rightarrow \R$ and $a_i \in \R^d$,  while $\psi \equiv 0$. The considered objective thus becomes: 
\[
\min_{x\in \R^d} \left \{ F(x) = f(x) =  \frac1n\sum_{i=1}^n \phi_i( \la a_i, x \ra ) \right\}.
\]
 Until very recently, such models were predominantly optimized using standard deterministic methods such as gradient descent, accelerated gradient/{\tt FISTA} or Newton's method. As already mentioned, these classical methods require an evaluation of $\nabla f(x^k)$ at every iteration.\footnote{Newton's method requires the computation of $\nabla^2 f(x^k)$ on top of that.} Consequently, the deterministic methods have to evaluate the dot product $\langle a_i ,x \rangle$ for all $i\in \{1,2,\dots, n \}$ in each step, and thus are either very expensive and may even be  infeasible in the big data setting (i.e., when $n,d$ are large).\footnote{In particular, the cost of performing a single iteration is $\cO(nd)$ for (accelerated) gradient descent and $\cO(nd^2 + d^3)$ for Newton's method.} The demand for solving such big data problems   resulted in the development of algorithms  working with a small random subset of the training data in each iteration only. 

The most natural approach is to use {\tt SGD} (i.e., subsample the finite sum) as described in Section~\ref{sec:intro_randomized}. The main idea of {\tt SGD} is to in each iteration pick a random index $j, \, 1\leq j \leq n$, and move the current iterate $x$ in the direction of the stochastic gradient 
\begin{equation}\label{eq:intro_stochgrad}
\nabla f_j(x) = a_j  \phi_j'( \la a_j, x \ra ).
\end{equation} 
While the cost of performing a simple {\tt SGD} iteration is often  $\cO(d)$ only, {\tt SGD} is slow in terms of how many iterations are required to get to an $\epsilon$-neighborhood to the optimum (where $\epsilon$ is relatively small). In particular, the stochastic gradient estimator~\eqref{eq:intro_stochgrad} has a (non-zero) variance at the optimum, causing {\tt SGD} to be gradually slower over time. Consequently, {\tt SGD} either converges linearly to a neighborhood of the solution only (Proposition~\ref{prop:sgdrate}), or converges sublinearly to the true optimum using a decreasing stepsize policy.

Fortunately, the issue of sublinear convergence of {\tt SGD} has been resolved using a more sophisticated stochastic gradient estimator whose variance progressively diminishes as $x^k\rightarrow x^*$. Methods based on such sophisticated estimators are commonly known as {\em variance reduced algorithms}; the most famous among them are {\tt SAG}~\cite{sag}, {\tt SAGA}~\cite{saga}, {\tt SDCA}~\cite{sdca}, {\tt SVRG}~\cite{svrg}, {\tt S2GD}~\cite{konevcny2013semi}, Finito~\cite{defazio2014finito}, {\tt MISO}~\cite{mairal2015incremental}, {\tt QUARTZ}~\cite{quartz} and {\tt SARAH}~\cite{nguyen2017sarah}. {\tt SAGA} and {\tt SVRG} achieve the variance reduction property by incorporating {\em control variates}~\cite{hickernell2005control} into the stochastic gradient -- we will exploit this idea multiple times throughout this text.

\subsection{From finite sum to coordinate descent}
An orthogonal approach to subsampling the finite sum is to subsample the domain (parameter space) and use (Randomized) Coordinate Descent ({\acrshort{CD}})~\cite{rcdm}. In its most basic form, {\tt CD} samples a random index $i$ (where $ 1\leq i \leq d$) and updates the $i$th coordinate of the current iterate  $x$ in the direction of $\nabla_if(x)e_i$, where $\nabla_if(x)$ is the $i$th partial derivative of $f$ at $x$ and $e_i$ is $i$th standard basis vector. Unlike {\tt SGD}, {\tt CD} does not suffer from the intrinsic variance at the optimum. At the same time, a single iteration of {\tt CD} can\footnote{The trick lies in the memorization of the dot products $\la a_i,x \ra$, see~\cite{rcdm} for details.} be implemented in time $\cO(n)$ and consequently, {\tt CD} is a serious competition to variance reduced {\tt SGD} algorithms. To decide which approach is superior to solve a given problem is rather complex~\cite{face-off}. However, the general rule of thumb suggests to use {\tt CD} if $d>n$ and variance reduced {\tt SGD} if $n>d$.

The above described, most straightforward version of {\tt CD}, is still fairly inefficient. Firstly, it is suboptimal in terms of iteration complexity;\footnote{Total number of iteration to reach $\varepsilon$-solution.} one shall combine it with Nesterov's acceleration as per~\cite{qu2016coordinate1, allen2016even, nesterov:2017}. Secondly, currently used hardware often allows evaluating a subset of partial derivatives in parallel almost as fast as a single partial derivative. This leads to the need to develop a tight theory of {\tt CD} methods under arbitrary sampling of the subsets\footnote{I.e., we wish to give as tight rate as possible for any given probability distribution over all subsets of $\{1, 2, \dots, d\}$ and corresponding sampling strategy for {\tt CD}.} to allow the user to tune {\tt CD} for his/her own specific hardware~\cite{qu2016coordinate1}. However, those two {\tt CD} adjustments were never combined before and this is where the story of this thesis starts. In particular, in a part of Chapter~\ref{acd}, we propose an accelerated {\tt CD} method with arbitrary sampling ({\acrshort{ACD}}).

One of the main disadvantages of {\tt CD} methods over {\acrshort{SGD}} algorithms is that they do not allow for a proximable regularizer $\psi$ that is non-separable.\footnote{We say that a function $h$ is separable if it can be written as $h(x)=\sum_{i=1}^d h_i(x_i)$.} In particular, non-separable $\psi$ prevents attainment of a linear convergence rate for {\tt CD} as the corresponding stochastic gradient estimator suffers from the inherent (non-zero) variance at the optimum, which very much resembles the story of {\tt SGD}. Since the mechanism of variance reduction has already successfully ``fixed'' the issue for {\tt SGD}, one might ask whether it is possible to incorporate an analogous trick into {\tt CD} methods. Fortunately, we were successful: in Chapter~\ref{sega}, we propose a new randomized algorithm---{\acrshort{SEGA}}---which accesses only a block of partial derivatives of $f$ each iteration and still converges linearly to the solution despite the presence of a \emph{non-separable} regularizer $\psi$.  This is the first variance-reduced {\tt CD} method in the literature.

\subsection{From coordinate descent to finite sum: three approaches}

The development of {\tt SEGA} provided us with many insights and ideas for future research. It brought us back to the finite-sum minimization in three somewhat independent ways, which we describe next. 

\subsubsection{Distributed optimization and random sparsification}
In many applications, the scale of the problem we are solving is so large that the dataset does not fit into the memory of a single machine. Consequently, multiple machines need to be employed to both store the data and train the model. In this thesis, we consider a specific, centralized case of distributed optimization/learning, where the machines are not allowed to communicate directly among themselves, but instead are allowed to communicate with a central server/master, also known as \emph{parameter server}.

Note that the optimization problem~\eqref{eq:finitesum} provides convenient notation for the mentioned scenario: function $f_i$ might represent a loss of the model on data owned by $i$th machine. In such a case, the value of $n$ corresponds to the number of machines/workers instead of the size of the dataset. 

Distributed optimization brings up several new challenges that are not present in standard optimization. Specifically, the communication between the workers and the parameter server/master takes a non-trivial time, often much more than the computation itself. There are several different ways to reduce communication complexity of gradient-type methods, one of which is {\em gradient sparsification}. Specifically, in order to communicate some non-sparse gradient $\nabla f_i(x)\in \R^d$, one should  send $d$ real numbers (often this is $32d$ or $64d$ bits). In contrast, to communicate a randomly sparsified gradient $ \nabla_{j_i} f_i(x) e_{j_i}$, where $1\leq j_i \leq d$ is selected uniformly at random, we only need to send a single real number along with its position, which is at least $d/2$ times cheaper in practice.

The major drawback of random sparsification is that the estimator $g(x)$ of $\nabla f(x)$ constructed as a naive aggregation of sparsified gradients from the workers
\[
g(x) \eqdef \frac1n \sum_{i=1}^n  \nabla_{j_i} f_i(x) e_{j_i}
\]
is very noisy, and its variance does not diminish as the method progresses through its iterations. Indeed, $g(x)$ has a non-zero variance at the optimum. In Chapter~\ref{99}, we incorporate control variates (similarly to {\tt SEGA}) on top of the sparsified gradient, which enables us to eliminate the adverse effect of the variance at the optimum on the convergence rate. Consequently, we show that our method can reduce worker$\rightarrow$server communication by as much as the factor of $n$ without hurting the convergence rate by more than a small constant. To illustrate the scale of this effect, consider a setup with 100 workers. In this case, we prove that  only $1\%$ of the usual worker$\rightarrow$server communication is needed to preserve the fast convergence rate.

\subsubsection{Unification of algorithms}
Since the variance reduced methods in three different setups (classical finite sum, Chapter~\ref{sega}, and Chapter~\ref{99}) share certain intrinsic similarities, one may wonder whether it is possible to unify them in a single algorithm, admitting a single analysis, so that one would not have to keep developing novel variance reduced algorithms along with their analyses from scratch. In Chapter~\ref{jacsketch}, we propose a general method---{\acrshort{GJS}} (Generalized Jacobian Sketching)---which constructs a gradient estimator given that a randomized linear transformation (a {\em sketch}) of the Jacobian matrix $$\mG(x) \eqdef \left[ \nabla f_1(x), \nabla f_2(x),\dots,  \nabla f_n(x) \right]\in \R^{d\times n}$$ is observed in each iteration. The  sketch is allowed to follow an arbitrary fixed distribution, and in special cases includes right matrix multiplication (in such a case we can recover {\tt SAGA} or {\tt SVRG}), and left matrix multiplication (in such a case we can recover {\tt SEGA}). many more sketches are possible, which gives rise to novel method not considered in literature before. This work is the first unification of stochastic optimization algorithms which subsample the finite sum, such as {\tt SAGA}, and algorithms which subsample the parameters, such as {\tt SEGA}. Our theory gives the currently best-known  convergence rate in each special case, and also allows for the development of importance sampling rates that exploit the smoothness structure of the objective.

We did not stop here, the story of this thesis unfolds further. 

Our findings made us realize that we can go one step further in terms of generality. In particular, the analysis of variance reduced {\tt SGD} algorithms and non-variance reduced {\tt SGD} shares a number of similar steps that can be abstracted to a {\em unified analysis framework}, which is what we do in Chapter~\ref{sigmak}. We provide a convergence rate for {\tt SGD} given that the unbiased stochastic gradient $g^k$ at iteration $k$ satisfies the novel general parametric bound
\[
        \EEE\left[\norm{g^k -\nabla f(x^*)}^2\mid x^k, \sigma^2_k \right] \le 2A (f(x^k) - f(x^*)) + B\sigma_k^2 + D_1,
\]
where $A, B, D_1 \in \R$ are some nonegative constants, while the sequence of nonegative random variables numbers $\sigma_k$ satisfies
\[
        \EEE\left[\sigma_{k+1}^2 \, \mid \, x^k, \sigma^2_k\right] \le (1-\rho) \sigma_k^2 + 2CD_f(x^k,x^*)  + D_2
        \]
 for some nonegative constants $\rho\leq 1, C, D_2$. Remarkably,  the above inequalities enable us to analyze {\tt SGD}, variance reduced methods for both finite sum and subspace gradients (i.e., {\tt SAGA} and {\tt SEGA}), quantized methods~\cite{mishchenko2019distributed}, and to develop and analyze several new algorithms of intriguing properties. Specifically, we introduce quantized methods with arbitrary sampling, partially variance reduced algorithms and an efficient, with-replacement importance sampling for minibatch {\tt SGD}.

Both of the above-mentioned frameworks have many different applications besides recovering well-known algorithms. In particular, we have noticed an application in \emph{federated learning}, which we describe in Chapter~\ref{local}: Local {\tt SGD} method ({\tt LSGD}) with imperfect aggregation can be seen as (non-uniform) {\tt SGD} applied to a carefully constructed 2-sum objective that we introduce. The corresponding variance reduced algorithm we propose (a special case of {\tt GJS}) achieves a linear rate which does not rely on the assumption of data homogeneity, and is  favorable to classical variants of local {\tt SGD} in terms of the convergence speed and communication complexity. Besides the importance of the newly proposed objective from the modeling perspective, our results suggest that the celebrated {\tt LSGD} method should better be seen as minimizing our objective than the classical finite sum, which explains the difficulties in the standard analysis of {\tt LSGD}, and reveals that the method implicitly aims to find {\em personalized} models.

\subsubsection{Product space objective}

Having previously discovered variance reduced variants of {\tt CD} methods, and their subspace generalizations, we realized that there is a {\em new deep connection} between these methods and modern variance reduced methods for finite sum minimization.  Specifically, we found that subspace VR algorithms are {\em more general} than finite sum VR algorithms: {\it applying subspace VR methods ({\tt SEGA}) to minimize a particular product space (i.e., in the domain of $\R^{nd}$) objective is equivalent to applying {\tt SAGA} to minimize arbitrary finite sum objective.} In order to obtain the best-known convergence rate of {\tt SAGA} from {\tt SEGA}, we had to tighten {\tt SEGA} theory to take advantage of the structure of the non-smooth function $\psi$. As a by-product, we have improved upon the rate of {\tt GJS} as well. More details are provided in Chapter~\ref{asvrcd}.

\subsection{Towards better stochastic condition numbers}

The iteration complexity of each proposed algorithm in this work is determined by the so-called {\em stochastic condition number}, which is itself a function of the objective smoothness, strong convexity, and randomness of the algorithm. Thus a natural question arises: what is the best possible stochastic condition number, assuming that the source of stochasticity is fixed (we want to have the freedom to develop arbitrary stochastic algorithm)? Intuitively speaking, the stochastic condition number is non-decreasing in the smoothness and non-increasing in the strong convexity parameters, and  is minimized  if these parameters are equal. Such a setting corresponds to a quadratic objective, where both the smoothness and strong convexity are measured with respect to the same Euclidean norm, given via the Hessian of the objective. 

Therefore, minimizing a general convex objective should not be simpler than minimizing the corresponding quadratic. We can now ask the reverse question: is there an algorithm which can minimize a non-quadratic convex objective with the same rate as if the function was in fact  quadratic, with its Hessian being the Hessian of the non-quadratic function at the optimum?  In Chapter~\ref{sscn}, we provide an affirmative answer: we develop a second-order\footnote{I.e., the method is allowed to access  second derivatives of the objective.} subspace descent method---{\tt SSCN} (Stochastic Subspace Cubic Newton)---capable of achieving so. In particular, the local convergence rate of {\tt SSCN} matches the rate of stochastic subspace descent applied to the problem of minimizing the quadratic function $$x \mapsto \frac12 (x-x^*)^\top \nabla^2f(x^*)(x-x^*),$$ where $x^*$ is the minimizer of $f$.

However, {\tt SSCN} does not achieve the optimal stochastic convergence rate as it does not incorporate Nesterov's momentum or another acceleration mechanism. In Chapter~\ref{ami}, we introduce an accelerated sketch-and-project\footnote{Sketch-and-project is a general stochastic method to minimize quadratic objective that recovers subspace descent in a special case.} method with a superior rate to its non-accelerated counterpart developed by Gower and Richt\'arik~\cite{gower2015randomized}. In particular, besides direct applications we elaborate on in the text, the fast rate from Chapter~\ref{ami} may also serve as an ambitious goal for the local rates of stochastic higher-order methods.

\section{Relationship among the chapters \label{sec:relation}}
Section~\ref{sec:history} describes how the chapters of this thesis were developed historically, outlining the chain of thought that led from one project to another.\footnote{With one exception -- Chapter~\ref{ami} was developed before everything else.} In this section, we elaborate on some non-historical connections among the chapters.

\subsubsection{High-level picture: A step towards the optimization utopia}

In the utopian optimization universe, a complexity\footnote{A complexity in a broader sense, for example the number of gradient evaluations, number of communication rounds, number of flops, or any other value of the interest.} would be known for any algorithm applied to solve any optimization problem. Such a knowledge would enable the practitioners to always apply an ideal algorithm given the problem to be solved and the complexity of the interest. This thesis presents a multiple steps towards the optimization utopia:

\begin{itemize}
\item {\bf We fill the missing gaps} in the current literature in terms of the tightening best-known theory of well-established algorithms (Chapters~\ref{jacsketch},~\ref{asvrcd},~\ref{local}), generalizing/extending the well-established algorithms (Chapters~\ref{acd},~\ref{jacsketch},~\ref{sigmak}) and proposing a brand-new methods (Chapters~\ref{sega},~\ref{99},~\ref{jacsketch},~\ref{sigmak},~\ref{asvrcd},~\ref{local},~\ref{sscn},~\ref{ami}). 
\item {\bf We establish novel and often surprising connections} between various algorithms, providing a better understanding of the optimization field (Chapters~\ref{jacsketch},~\ref{sigmak},~\ref{asvrcd},~\ref{local}). 
\item {\bf We unify and generalize} both the known and the newly introduced algorithms, allowing to tailor the randomized optimization strategy for a broad range of different applications (Chapters~\ref{jacsketch},~\ref{sigmak}). 
\end{itemize} 

Next, we describe specific topics that the thesis chapters focus on.

\subsubsection{Self-variance reduced methods, sublinear rates and control variates}
The algorithms  proposed in this thesis can be  categorized based on their relation to control variates into three different classes:
\begin{itemize}
\item {\bf Fast stochastic algorithms that do not require the aid of control variates.} This category includes {\tt ACD} (Chapter~\ref{acd}), {\tt SSCN} (Chapter~\ref{sscn}), accelerated sketch-and-project (Chapter~\ref{ami}), and over-parameterized {\tt SGD} (i.e., {\tt SGD} applied to a finite-sum problem where $\nabla f_i(x^*)=0$ for all $i$; see Chapter~\ref{sigmak} for the general method and rate or Chapter~\ref{99} for an application to distributed optimization). 

\item  {\bf Stochastic algorithms that do not use control variates despite the fact that control variates would improve the rate.} Such methods converge sublinearly (or converge linearly to a certain neighborhood of the optimal solution) due to the inherent variance of the gradient estimator at the optimum. This category includes local {\tt SGD} (Chapter~\ref{local}), some variants of sparsified parallel algorithms (Chapter~\ref{99}), and a number of other {\tt SGD} variants that can be analyzed using the framework of Chapter~\ref{sigmak}.

\item  {\bf Linearly converging stochastic algorithms aided by control variates.} Those include {\tt SEGA} (Chapter~\ref{sega}), sparsified VR algorithms {\tt ISEGA}, {\tt ISAGA} (Chapter~\ref{99}), as well as local {\tt SGD} with variance reduction (Chapter~\ref{local}). All of these algorithms can also be obtained as a special case of the {\tt GJS} framework (Chapter~\ref{jacsketch}) -- {\tt GJS} tightens the rate of {\tt SEGA} and extends both {\tt ISEGA} and {\tt ISAGA} (and allows for their combination). The rate of {\tt GJS} is further improved in Chapter~\ref{asvrcd}, which allows for exploiting the specific structure of the regularizer $\psi$. 
\end{itemize}

\subsubsection{Randomization over the data or parameters}

As already mentioned, there are two different ways in which randomization can enter an optimization procedure -- either subsampling the domain (parameters) or subsampling the finite sum (data).

\begin{itemize}
\item  {\bf Subsampling the space.} Generally speaking, methods in this category in each iteration compute the gradient over a randomly chosen subspace only. This corresponds to a subset of partial derivatives in the special case when the subspace is spanned by a subset of the standard unit basis vectors. While some algorithms update the current iterate along the selected random subspace only ({\tt ACD} from Chapter~\ref{acd} and {\tt SSCN} from Chapter~\ref{sscn}), the others perform a full dimensional update due to the presence of control variates ({\tt SEGA} from Chapter~\ref{sega}, {\tt ISEGA} from Chapter~\ref{99} or {\tt SVRCD} from Chapters~\ref{jacsketch},~\ref{asvrcd}). We shall also mention that the methods aided by the control variates are usually somewhat slower than the methods moving along the subspace only.

\item {\bf Subsampling the data.} Various chapters of this thesis propose or improve upon known methods that subsample the finite sum~\eqref{eq:finitesum}. As a special case of the {\tt GJS} framework (Chapter~\ref{jacsketch}), we were able to introduce Loopless {\tt SVRG} ({\tt LSVRG})~\cite{hofmann2015variance, kovalev2019don} with arbitrary sampling and proximal step (thus making it significantly faster). Next, we introduce a linearly convergent variance reduced local {\tt SGD} method in Chapter~\ref{local}; which is by an order of magnitude faster than other variants of local SGD in the literature. Lastly, the unified {\tt SGD} analysis we provide in Chapter~\ref{sigmak} allowed us to both analyze a new, with-replacement minibatch {\tt SGD} method with importance sampling which is cheaper to implement than the without-replacement variant, and improve upon several quantized {\tt SGD} algorithms (for example, we propose the first quantized {\tt SGD} method with arbitrary sampling).

\item {\bf Subsampling both the domain and the space at the same time. }Two chapters of this work consider random linear measurements of the Jacobian as an oracle model: {\tt GJS} (Chapter~\ref{jacsketch}) is a variance reduced algorithm for minimizing a general finite-sum objective, while accelerated sketch-and-project (Chapter~\ref{ami}) is an algorithm for minimizing quadratics. Our oracle model allows for sampling from both the space and the finite sum at the same time. In a special case, this reduces to the gradient sparsification approach we propose in Chapter~\ref{99}, and thus recovers the {\tt ISEGA}, {\tt ISAGA} or {\tt ISAEGA} algorithms we which proposed previously. Needless to say, the unified {\tt SGD} analysis from Chapter~\ref{sigmak} captures this level of generality as well.

\end{itemize}

To conclude this section, we shall mention that domain-subsampling algorithms are often capable of performing finite sum subsampling, either through the product space objective, which we introduce in Chapter~\ref{asvrcd}, or via the duality trick from~\cite{sdca}.

\subsubsection{Distributed optimization}
Two chapters of this work consider predominantly distributed optimization, where the bottleneck of the optimization system is communication. Chapter~\ref{99} and Chapter~\ref{local} present two orthogonal approaches in two different distributed setups. Specifically, Chapter~\ref{99} introduces a new method based on random sparsification of the gradient, which {\em provably reduces the worker$\rightarrow$server communication by order of the number of  the workers at essentially no cost}.\footnote{In some distributed computation systems, communication from the workers to the server, is 10-20 times more expensive than the communication from the server to workers~\cite{mishchenko2019distributed}.} On the other hand, Chapter~\ref{local} focuses on the {\em federated learning} paradigm. In it, we introduce a novel personalization-encouraging objective which we argue is more natural to be optimized by local gradient methods, and for the first time prove communication complexity benefits of local gradient decent methods. We shall note that all variance reduced algorithms introduced in these chapters are a special case of {\tt GJS} (Chapter~\ref{jacsketch}), and at the same time, Chapter~\ref{jacsketch} extends the results of Chapter~\ref{99} allowing for both subsampling the local objective and gradient sparsification while keeping linear rate. Further, all convergence rates of Chapters~\ref{99},~\ref{local}, as well as the rates of other quantized algorithms for distributed optimization~\cite{mishchenko2019distributed, horvath2019stochastic} can be obtained as a special case of the framework of Chapter~\ref{sigmak}. 

\subsubsection{Importance sampling for minibatches} 
While minibatch variants of \texttt{{\tt CD}}  methods are very popular in practice, until now, there was no importance sampling for \texttt{{\tt CD}} that outperforms the standard uniform minibatch sampling in terms of worst-case guarantees. In Chapter~\ref{acd} we design new importance sampling for minibatch \texttt{{\tt CD}} and minibatch \texttt{{\tt ACD}} which significantly outperforms previous state-of-the-art minibatch \texttt{{\tt ACD}} in practice. Surprisingly, the sampling strategy applies to stochastic minibatch methods that subsample the finite sum objective -- it can improve upon {\tt SGD}, {\tt SAGA}, {\tt SVRG} and others.\footnote{It applies to all special cases covered by Chapters~\ref{jacsketch} and~\ref{sigmak}; see the corresponding appendices.} Further, Chapter~\ref{sigmak} presents a with-replacement variant of {\tt SGD}, where the importance minibatch sampling is particularly cheap to implement.

\subsubsection{Proximal methods}
Most of the algorithms proposed in this work support arbitrary proximable regularizer $\psi$ which is proper, closed, convex,  and possibly non-smooth. This includes Chapters~\ref{sega},~\ref{jacsketch},~\ref{sigmak},~\ref{asvrcd} and ~\ref{local} and a part of Chapter~\ref{99} (further generalized in Chapter~\ref{jacsketch}). Next, Chapter~\ref{sscn} requires $\psi$ to be separable as it proposes a subspace descent method without control variates. 

While the standard analysis of proximal methods provides a rate identical to the corresponding non-proximal variants, in Chapter~\ref{asvrcd} we show that the presence of $\psi$ with a specific structure might significantly simplify the problem and thus enable faster optimization. As a consequence of this observation, we show that fast rates of variance reduced algorithms that subsample the finite sum can be obtained from variance reduced methods that subsample the space.  This establishes a new and deep link between two strands of optimization methods.

\subsubsection{Accelerated algorithms}
Many of the algorithms proposed in this work incorporate some form of Nesterov's acceleration~\cite{nesterov83}. In some chapters, the acceleration is the or one of the key contributions (i.e., {\tt ACD} in Chapter~\ref{acd},  {\tt ASVRCD} in Chapter~\ref{asvrcd}, accelerated sketch-and-project in Chapter~\ref{ami}), while some other chapters merely demonstrate that acceleration can be incorporated into the loop (i.e., {\tt ASEGA} in Chapter~\ref{sega} or {\tt  IASGD} in Chapter~\ref{99}).

\subsubsection{Second order methods}
While this thesis focuses predominantly on first-order optimization, Chapters~\ref{sscn} and~\ref{ami} study second-order algorithms as well. Specifically, Chapter~\ref{sscn} introduces the Stochastic Subspace Cubic Newton method ({\tt SSCN}) -- a new globally convergent second-order subspace descent method. On the other hand, Chapter~\ref{ami} introduces accelerated sketch-and-project method for solving linear systems in Euclidean spaces, which can be seen as a first-order and second-order method at the same time due to the quadratic nature of the objective.

\subsubsection{Summary of the links among the chapters}
To conclude this section, we summarize both what the chapters are about, as well outline several links among them.

\begin{figure}[!h]
\centering
\includegraphics[width =  \textwidth]{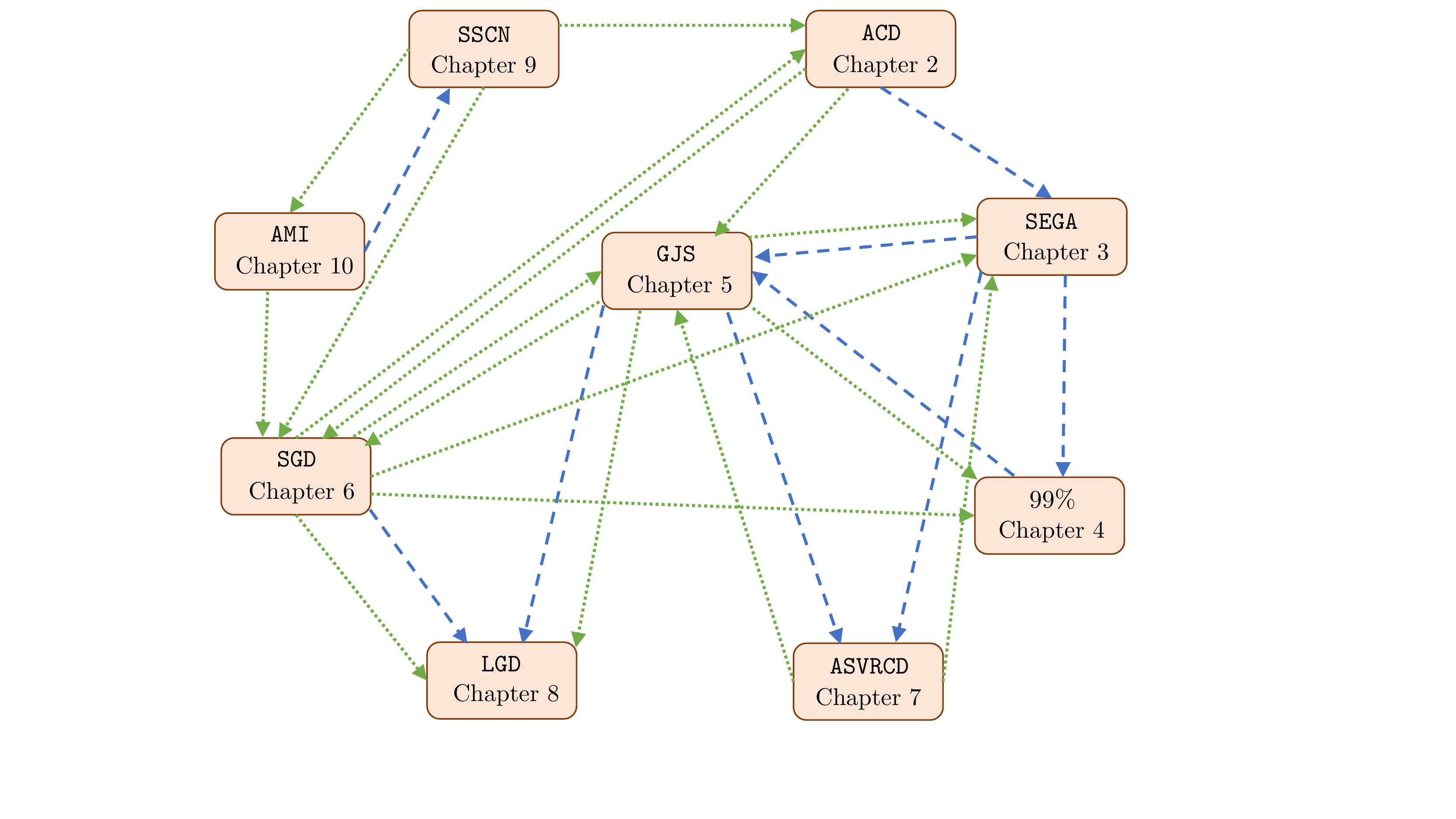}
\caption{Graph depicting the relationships among the chapters of this thesis. Blue dashed arrow indicates motivation among chapters, while green dotted arrow indicates a significant insight that chapters shed on each other. As an example, let us explain the edges of Chapter~\ref{sega} ({\tt SEGA}): the development of {\tt SEGA} was enabled by our results on {\tt CD} (Chapter~\ref{acd}) and motivated us to develop the results contained in Chapters~\ref{99},~\ref{jacsketch} and~\ref{asvrcd}.  Further, Chapter~\ref{jacsketch} recovers/improves upon the convergence rate of {\tt SEGA}, Chapter~\ref{sigmak} enables a partial variance reduction in {\tt SEGA} and lastly, Chapter~\ref{jacsketch} shows that {\tt SAGA} is a special case of {\tt SEGA}.} \label{fig:relation}
\end{figure}

First, Table~\ref{tbl:algorithms} presents a representative algorithm for each chapter of this thesis, as well as the covered topics. Next, Table~\ref{tbl:algorithms_all} highlights which algorithms presented in this thesis are novel and which are not. Lastly, Figure~\ref{fig:relation} summarizes the essential connections among the chapters of this thesis that were outlined above.

\begin{table}[!t]
{
\small
\begin{center}
\begin{tabular}{|c|c|c|c|c|c|c|c|c|}
\hline
{\bf Chapter} & {\bf Ref} & {\bf Alg}  &  {\bf {\acrshort{VR}} } &  {\bf Accel } & {\bf Subsp} & {\bf Prox} & {\bf Distrib}  & {\bf Note} \\
\hline
\hline
 \ref{acd} & \cite{hanzely2018accelerated}&{\tt ACD} & \xmark & \cmark & \cmark${}^\star$  & \xmark  &\xmark &  \begin{tabular}{c} Minibatch  \\ sampling \end{tabular} \\  \hline
\ref{sega}  &  \cite{sega} &{\tt SEGA} & \cmark & \cmark & \cmark${}^\dagger$  &\cmark   &  \xmark &   \begin{tabular}{c} Non-separable  \\ regularizer \end{tabular}   \\  \hline
 \ref{99}&\cite{mishchenko201999} & {\tt ISEGA}  & \cmark   &  \cmark  & \cmark${}^\dagger$ &\cmark   & \cmark &  \begin{tabular}{c} Reduced  \\ communication \end{tabular}   \\  \hline
\ref{jacsketch} &  \cite{hanzely2019one}&{\tt GJS}  & \cmark &  \xmark & \cmark${}^\dagger$ &\cmark   &  \cmark &   \begin{tabular}{c} General VR  \\ framework \end{tabular}  \\  \hline
\ref{sigmak} & \cite{sigma_k} &{\tt SGD}  &\cmark   &  \xmark & \cmark${}^\dagger$&\cmark & \cmark  &  \begin{tabular}{c}  General {\tt SGD} \\ analysis  \end{tabular}  \\  \hline
 \ref{asvrcd} & \cite{asvrcd} &  {\tt ASVRCD}& \cmark & \cmark  & \cmark${}^\dagger$ &\cmark  & \xmark &   \begin{tabular}{c} Subspace $\geq$  \\ finite sum \end{tabular} \\ \hline
\ref{local}& \cite{local} &  {\acrshort{LGD} }& \cmark & \xmark  & \xmark & \cmark  & \cmark &   \begin{tabular}{c} Local {\tt SGD} \\ with VR \end{tabular} \\ \hline
\ref{sscn} & \cite{sscn}&{\tt SSCN}  & \xmark & \xmark  & \cmark &\cmark${}^\ddagger$  & \xmark &  \begin{tabular}{c} $2^\mathrm{nd}$  order + \\ cubic regularizer  \end{tabular}   \\ \hline
 \ref{ami} & \cite{gower2018accelerated} &{\acrshort{AMI}}  & \xmark & \cmark & \cmark & \xmark & \xmark &  \begin{tabular}{c} Quadratic \\ objective \end{tabular} \\ \hline
\end{tabular}
\end{center}
\caption{Summary of representative algorithms proposed in each chapter and topics covered in each chapter. Columns (chapter topics): VR = variance reduced method, Accel = Nesterov's acceleration, Subsp = subspace descent, Prox = proximal setup, Distrib = distributed setup. Further clarifications: ${}^\star$ {\tt ACD} allows for subspaces spanned by standard basis vectors only;  ${}^\dagger$ these methods consider a general subspace oracle, but perform full dimensional updates; ${}^\ddagger$ {\tt SSCN} requires  the regularizer $\psi$ to be separable. }
\label{tbl:algorithms}
}
\end{table}

\begin{table}[!t]
\centering
\footnotesize
\begin{tabular}{|c|c||c|c||c|c||c|c||c|c||c|c||c|c|}
\hline
{\bf  \#} & {\bf New} &   {\bf  \#} & {\bf New}&{\bf  \#} & {\bf New} &   {\bf  \#} & {\bf New} &   
{\bf  \#} & {\bf New} &   {\bf  \#} & {\bf New} &   {\bf  \#} & {\bf New}  \\
\hline
1	&	\xmark	&	11	&	\cmark	&	21	&	\cmark	&	31	&	\cmark	&	41	&	\cmark	&	51	&	\cmark	&	61  	&	\cmark		\\
2	&	\xmark	&	12	&	\cmark	&	22	&	\cmark	&	32	&	\cmark	&	42	&	\cmark	&	52	&	\xmark	&	62	&	\cmark		\\
3	&	\xmark	&	13	&	\xmark	&	23	&	\cmark	&	33	&	 \cmark	&	43	&	\cmark	&	53	&	\xmark	&	63	&	\cmark		\\
4	&	\xmark	&	14	&	\cmark	&	24	&	\cmark	&	34	&	\cmark	&	44	&	\xmark	&	54	&	\xmark	&	64	&	\cmark		\\
5	&	\cmark	&	15	&	\cmark	&	25	&	\cmark	&	35	&	\cmark	&	45	&	\xmark	&	55	&	\xmark	&		&			\\
6	&	\cmark	&	16	&	\cmark	&	26	&	\cmark	&	36	&	\cmark	&	46	&	\cmark	&	56	&	\cmark	&		&			\\
7	&	\cmark	&	17	&	\cmark	&	27	&	\xmark	&	37	&	\cmark	&	47	&	\cmark	&	57	&	\xmark	&		&			\\
8	&	\cmark	&	18	&	\cmark	&	28	&	\xmark	&	38	&	\cmark	&	48	&	\xmark	&	58	&	\xmark	&		&			\\
9	&	\cmark	&	19	&	\cmark	&	29	&	\cmark	&	39	&	\cmark	&	49	&	\cmark	&	59	&	\cmark	&		&			\\
10	&	\cmark	&	20	&	\cmark	&	30	&	\cmark	&	40	&	\cmark	&	50	&	\cmark	&	60	&	\cmark	&		&			\\
\hline
\end{tabular}
\caption{List of all algorithms stated in this work. Marker \cmark\,  indicates that the algorithm is new (i.e., proposed in this work) while marker \xmark\, indicates that the algorithm is known.}
\label{tbl:algorithms_all}
\end{table}

\section{Outline and individual contributions}
Each chapter of this work consists of a single paper; some of them are already published while the others are at various stages of the submission process. Let us now give a brief overview of the contents of each chapter individually.

\subsection{Accelerated coordinate descent with arbitrary sampling and best rates for minibatches (Chapter~\ref{acd})}
Accelerated coordinate descent is a widely popular optimization algorithm due to its efficiency in large-dimensional problems. It achieves state-of-the-art complexity on an important class of empirical risk minimization problems.  In this work, we design and analyze an accelerated coordinate descent (\texttt{{\tt ACD}}) method, which in each iteration updates a random subset of coordinates according to an arbitrary but fixed probability law, which is a parameter of the method.  While minibatch variants of \texttt{{\tt ACD}} are more popular and relevant in practice, there is no importance sampling for \texttt{{\tt ACD}} that outperforms the standard uniform minibatch sampling. Through insights enabled by our general analysis, we design new importance sampling for minibatch \texttt{{\tt ACD}}, which significantly outperforms previous state-of-the-art minibatch \texttt{{\tt ACD}} in practice. We prove a rate that is at most $\cO(\sqrt{\tau})$ times worse than the rate of minibatch \texttt{{\tt ACD}} with uniform sampling, but can be $\cO(d/\tau)$ times better, where $\tau$ is the minibatch size. Since in modern supervised learning training systems, it is standard practice to choose $\tau \ll d$, and often $\tau=\cO(1)$, our method can lead to dramatic speedups. We obtain similar results for minibatch non-accelerated \texttt{{\tt CD}} as well, achieving improvements on previous best rates. Further, the importance sampling for non-accelerated \texttt{{\tt CD}} can be incorporated into stochastic algorithms that decompose finite sums such as {\tt SGD}, {\tt SAGA}, and others. 

The chapter is based on the paper:
\begin{quote}
\cite{hanzely2018accelerated} \bibentry{hanzely2018accelerated}.
\end{quote}

\subsection{{\tt SEGA}: Variance reduction via gradient sketching (Chapter~\ref{sega})}
In Chapter~\ref{sega}, we propose a randomized first-order optimization method---\texttt{{\tt SEGA}} (SkEtched GrAdient)---which progressively throughout its iterations builds a variance-reduced estimate of the gradient from random linear measurements (sketches) of the gradient obtained from an oracle. In each iteration, \texttt{{\tt SEGA}} updates the current estimate of the gradient through a sketch-and-project operation using the information provided by the latest sketch, and this is subsequently used to compute an unbiased estimate of the true gradient through a random relaxation procedure. This unbiased estimate is then used to perform a gradient step. Unlike standard subspace descent methods, such as coordinate descent, \texttt{{\tt SEGA}} can be used for optimization problems with a {\em non-separable} proximal term. We provide a general convergence analysis and prove linear convergence for strongly convex objectives. In the special case of coordinate sketches, \texttt{{\tt SEGA}} can be enhanced with various techniques such as {\em importance sampling}, {\em minibatching}, and {\em acceleration}, and its rate is up to a small constant factor identical to the best-known rate of coordinate descent from Chapter~\ref{acd}. 

The chapter is based on the paper:
\begin{quote}
\cite{sega} \bibentry{sega}.
\end{quote}

\subsection{99\% of Worker-Master Communication in Distributed Optimization is Not Needed (Chapter~\ref{99})}

We improve upon algorithms that fit the following template: a local gradient estimate is computed independently by each worker, then communicated to a master, which subsequently performs averaging. The average is broadcast back to the workers, which uses it to perform a gradient-type step to update the local version of the model. We observe that the above template is fundamentally inefficient in that too much data is unnecessarily communicated from the workers to the server, which slows down the overall system.  We propose a fix based on a new update-sparsification method we develop in this work, which we suggest be used on top of existing methods. Namely, we develop a new variant of parallel block coordinate descent based on independent sparsification of the local gradient estimates before communication. We demonstrate that with only $m/n$ blocks sent by each of $n$ workers, where $m$ is the total number of parameter blocks, the theoretical iteration complexity of the underlying distributed methods is essentially unaffected. As an illustration, this means that when $n=100$ parallel workers are used, the communication of  $99\%$  blocks is redundant, and hence a waste of time. Our theoretical claims are supported through extensive numerical experiments that demonstrate an almost perfect match with our theory on a number of synthetic and real datasets. 

The chapter is based on the paper:
\begin{quote}
\cite{mishchenko201999} \bibentry{mishchenko201999}.
\end{quote}

\subsection{One method to rule them all: Variance reduction for data, parameters and many new methods (Chapter~\ref{jacsketch})}

Next, in Chapter~\ref{sega}, we propose a remarkably general variance-reduced method suitable for solving regularized empirical risk minimization problems with either a large number of training examples, or a large model dimension, or both. In special cases, our method reduces to several known and previously thought to be unrelated methods, such as {\tt {\tt SAGA}}~\cite{saga}, {\tt L{\tt SVRG}}~\cite{hofmann2015variance, kovalev2019don}, {\tt JacSketch}~\cite{jacsketch}, {\tt {\tt SEGA}}~\cite{sega} and {\tt {\acrshort{ISEGA}}}~\cite{mishchenko201999}, and their arbitrary sampling and proximal generalizations. However, we also highlight a large number of new specific algorithms with interesting properties. We provide a single theorem establishing linear convergence of the method under smoothness and quasi strong convexity assumptions. With this theorem, we recover best-known and sometimes improved rates for known methods arising in special cases. As a by-product, we provide the first unified method and theory for stochastic gradient and stochastic coordinate descent type methods. 

The chapter is based on the paper:
\begin{quote}
\cite{hanzely2019one} \bibentry{hanzely2019one}.
\end{quote}

\subsection{A unified theory of {\tt SGD}: Variance reduction, sampling, quantization  and coordinate descent (Chapter~\ref{sigmak})}
We introduce a unified analysis of a large family of variants of proximal stochastic gradient descent ({\tt {\tt SGD}}), which so far have required different intuitions, convergence analyses, have different applications, and which have been developed separately in various communities. We show that our framework includes methods with and without the following tricks, and their combinations: variance reduction, importance sampling, mini-batch sampling, quantization, and coordinate sub-sampling.  As a by-product, we obtain the first unified theory of {\tt {\tt SGD}} and randomized coordinate descent ({\tt CD}) methods,  the first unified theory of variance reduced and non-variance-reduced {\tt {\tt SGD}} methods, and the first unified theory of quantized and non-quantized methods. A key to our approach is a parametric assumption on the iterates and stochastic gradients. In a single theorem, we establish a linear convergence result under this assumption and strong-quasi convexity of the loss function. Whenever we recover an existing method as a special case, our theorem gives the best-known complexity result.
Our approach can be used to motivate the development of new useful methods and offers pre-proved convergence guarantees. To illustrate the strength of our approach, we develop five new variants of {\tt {\tt SGD}}, and through numerical experiments, demonstrate some of their properties.  

The chapter is based on the paper:
\begin{quote}
\cite{sigma_k} \bibentry{sigma_k}.
\end{quote}

\subsection{Variance reduced coordinate descent with acceleration: New method with a surprising application to finite-sum problems (Chapter~\ref{asvrcd})}

Further, in Chapter~\ref{asvrcd}, we propose  {\acrshort{ASVRCD}}: an accelerated version of stochastic variance reduced coordinate descent. As other variance reduced coordinate descent methods such as {\tt SEGA} or {\tt SVRCD}, our method can deal with problems that include a non-separable and non-smooth regularizer, while accessing a random block of partial derivatives in each iteration only. However, {\tt ASVRCD} incorporates Nesterov's momentum, which offers favorable iteration complexity guarantees over both  {\tt SEGA} and {{\acrshort{SVRCD}}. As a by-product of our theory, we show that a variant of Katyusha~\cite{allen2017katyusha} is a specific case of {\tt ASVRCD}, recovering the optimal oracle complexity for the finite sum objective. 

The chapter is based on the paper:
\begin{quote}
\cite{asvrcd} \bibentry{asvrcd}.
\end{quote}

\subsection{Federated learning of a mixture of global and local models (Chapter~\ref{local})}
We propose a new optimization formulation for training federated learning models. The standard formulation has the form of an empirical risk minimization problem constructed to find a single global model trained from the private data stored across all participating devices. In contrast, our formulation seeks an explicit trade-off between this traditional global model and the local models, which can be learned by each device from its own private data without any communication. Further, we develop several efficient variants of {\tt SGD} (with and without partial participation and with and without variance reduction) for solving the new formulation and prove communication complexity guarantees. Notably, our methods are similar but not identical to federated averaging / local {\tt SGD}, thus shedding some light on the essence of the elusive method. In particular, our methods do not perform full averaging steps and instead merely take steps towards averaging. We argue for the benefits of this new paradigm for federated learning. 

The chapter is based on the paper:
\begin{quote}
\cite{local} \bibentry{local}.
\end{quote}

\subsection{Stochastic subspace cubic Newton (Chapter~\ref{sscn}) }
In Chapter~\ref{sscn}, we propose a new randomized second-order optimization algorithm---Stochastic Subspace Cubic Newton ({\acrshort{SSCN}})---for minimizing a high dimensional convex function $f$. Our method can be seen both as a {\em stochastic} extension of the cubically-regularized Newton method~\cite{nesterov2006cubic}, and a {\em second-order} enhancement of stochastic subspace descent~\cite{kozak2019stochastic}. We prove that as we vary the minibatch size, the global convergence rate of {\tt SSCN} interpolates between the rate of stochastic coordinate descent ({\tt CD}) and the rate of cubic regularized Newton, thus giving new insights into the connection between first and second-order methods. Remarkably, the local convergence rate of {\tt SSCN} matches the rate of stochastic subspace descent applied to the problem of minimizing the quadratic function $x \mapsto\frac12 (x-x^*)^\top \nabla^2f(x^*)(x-x^*)$, where $x^*$ is the minimizer of $f$, and hence depends on the properties of $f$ at the optimum only. Our numerical experiments show that {\tt SSCN} outperforms non-accelerated first-order {\tt CD} algorithms while being competitive to their accelerated variants. 

The chapter is based on the paper:
\begin{quote}
\cite{sscn} \bibentry{sscn}.
\end{quote}

\subsection{Accelerated stochastic matrix inversion:  General theory and  speeding up {\tt BFGS} rules for faster second-order optimization (Chapter~\ref{ami})}
In Chapter~\ref{ami}, we present the first accelerated randomized algorithm for solving linear systems in Euclidean spaces. One essential problem of this type is the matrix inversion problem. In particular, our algorithm can be specialized to invert positive definite matrices in such a way that all iterates (approximate solutions) generated by the algorithm are positive definite matrices themselves. This opens the way for many applications in the field of optimization and machine learning.  As an application of our general theory, we develop the {\em first accelerated (deterministic and stochastic) quasi-Newton updates}. Our updates lead to provably more aggressive approximations of the inverse Hessian and lead to speedups over classical non-accelerated rules in numerical experiments. Experiments with empirical risk minimization show that our rules can accelerate the training of machine learning models. 

The chapter is based on the paper:
\begin{quote}
\cite{gower2018accelerated} \bibentry{gower2018accelerated}.
\end{quote}

\subsection{Excluded papers}
I had a chance to co-author four more papers during my studies, which are not included in this work: one about an accelerated mirror descent method for relatively smooth optimization~\cite{abpg}, two about robust principal component analysis~\cite{dutta2019nonconvex, dutta2019best} and the last one about optimal algorithms for personalized federated learning~\cite{fedoptimal}.

\chapter{Accelerated Coordinate Descent with Arbitrary Sampling and Best Rates for Minibatches}
\label{acd}

\graphicspath{{ACD/experiments/images/}}

In this chapter we consider a particular instance of the general optimization problem~\eqref{eq:finitesum} with $\psi\equiv 0$ and $f$ not necessarily having a finite-sum structure, i.e.,
\begin{equation}\label{eq:acd_problem}
\min_{x\in \R^d} f(x).
\end{equation}
Specifically, we assume that  $f$ is a smooth and strongly convex function, and the main difficulty  comes from the dimension $d$ being very large (e.g., millions or billions). 
In this regime, {\em coordinate descent (\texttt{CD})} variants of gradient methods are the state of the art. 

The simplest variant of  \texttt{CD} in each iteration updates a single variable of $x$ by taking a one dimensional gradient step along the direction of the $i$th unit basis vector $e_i\in \R^d$, which leads to the update rule \begin{equation} \label{eq:acd_CD-intro}x^{k+1} = x^k - \alpha_i \nabla_i f(x^k) e_i,\end{equation} where $\nabla_i f(x^k) \eqdef e_i^\top \nabla f(x^k)$ is the $i$th partial derivative and $\alpha_i$ is a suitably chosen stepsize.  The classical smoothness assumption used in the analysis of \texttt{CD} methods~\cite{rcdm} is to require the existence of constants $L_i>0$ such that \begin{equation}\label{eq:acd_98g98gf} f(x+ t e_i) \leq f(x) + t \nabla_i f(x)+ \frac{L_i}{2} t^2\end{equation} holds for all $x\in \R^d$, $t\in \R$ and $i\in [d]\eqdef \{1,2,\dots,d\}$. In this setting, one can choose the stepsizes to be $\alpha_i = 1/L_i$. 

There are several rules studied in the literature for  choosing the coordinate $i$ in iteration $k$, including cyclic rules~\cite{luo1992convergence,tseng2001convergence,saha2013nonasymptotic,wright2015coordinate,gurbuzbalaban2017cyclic}, Gauss-Southwell or other greedy rules~\cite{nutini2015coordinate,you2016asynchronous,stich2017approximate}, random (stationary) rules~\cite{rcdm, richtarik2014iteration, pcdm, shalev2014accelerated, lin2014accelerated, approx} and adaptive random rules~\cite{csiba15,stich2017safe}. In this work we focus on stationary random rules, which are popular by practitioners and well understood in theory.

\paragraph{Updating one coordinate at a time.} The simplest randomized \texttt{CD} method  of the form \eqref{eq:acd_CD-intro} chooses coordinate $i$  in each iteration uniformly at random. If $f$ is $\mu$-strongly convex, then this method converges in $(d \max_i  L_i/\mu) \log (1/\epsilon)$ iterations in expectation. If index $i$ is chosen with probability $p_i \propto L_i$, then the iteration complexity  improves to $(\sum_i L_i/\mu) \log (1/\epsilon)$. The latter result is always better than the former, and can be up to $d$ times better. These results were established in a seminal paper by Nesterov~\cite{rcdm}. The analysis was later generalized to  arbitrary probabilities $p_i>0$ by Richt\'arik and Tak\'a\v{c}~\cite{richtarik2014iteration}, who obtained the complexity 
\begin{equation} \label{eq:acd_089ff} \left(\max_i \frac{L_i}{p_i \mu}\right) \log \frac{1}{\epsilon}.\end{equation}  
Clearly, \eqref{eq:acd_089ff} includes the previous two results as special cases. Note that the importance sampling probabilities given by $p_i\propto L_i$ minimizes the complexity bound \eqref{eq:acd_089ff} and are therefore in this sense optimal.

\paragraph{Minibatching: updating more coordinates at a time.} In many situations it  is advantageous to update a small {\em subset (minibatch)} of coordinates in each iteration, which leads to the 
{\em minibatch \texttt{CD} method} which has  the form
 \begin{equation} \label{eq:acd_Parallel-CD-intro}x^{k+1}_i = \begin{cases} x^k_i - \alpha_i \nabla_i f(x^k) & \quad i\in S^k, \\
 x^k_i & \quad i\notin S^k.
 \end{cases} 
 \end{equation} 

   For instance, it is often equally easy to fetch information about a small batch of coordinates $S^k$ from memory at the same or comparable  time as it is to fetch information about a single coordinate. If this memory access time is the bottleneck as opposed to computing the actual updates to coordinates $i\in S^k$, then it is more efficient to update all coordinates belonging to the minibatch $S^k$. Alternatively, in situations where parallel processing is available, one is able to compute the updates to a small batch of coordinates simultaneously, leading to speedups in wall clock time. With this application in mind, minibatch \texttt{CD} methods are also often called {\em parallel} \texttt{CD} methods \cite{pcdm}.



\section{Arbitrary sampling and  minibatching}

\paragraph{Arbitrary sampling.}  The method \eqref{eq:acd_Parallel-CD-intro} was analyzed in \cite{pcdm}   for {\em uniform samplings} $S^k$, i.e., assuming that $\Prob( i\in S^k) = \Prob(j\in S^k)$  for all $i,j$.  However, the ultimate generalization is captured by the notion of {\em arbitrary sampling} \cite{nsync}. A {\em sampling} refers to a set-valued random mapping $S$ with values being the subsets of $[d]$. The word {\em arbitrary} refers to the fact that no additional assumptions on the sampling, such as uniformity, are made. This result generalizes the results mentioned above.

\paragraph{$\mM$-smoothness.}
For minibatch \texttt{CD} methods it is useful to assume a more general notion of smoothness parameterized by a positive semidefinite matrix $\mM\in \R^{d\times d}$. We say that $f$ is $\mM$-smooth if 
\begin{equation}\label{eq:acd_M-smooth-intro}f(x+h) \leq f(x) + \nabla f(x)^\top h + \frac{1}{2}h^\top \mM h\end{equation}
for all $x,h\in \R^d$. The standard $L$-smoothness condition is obtained in the special case when $\mM = L \mI$, where $\mI$ is the identity matrix in $\R^d$. Note that if $f$ is $\mM$-smooth, then \eqref{eq:acd_98g98gf} holds for $L_i=\mM_{ii}$. Conversely, it is known that if \eqref{eq:acd_98g98gf} holds, then \eqref{eq:acd_M-smooth-intro} holds for $\mM = d \Diag{L_1,L_2,\dots,L_d}$~\cite{rcdm}. If $h$ has at most $\omega$ nonzero entries, then this result can be strengthened and \eqref{eq:acd_M-smooth-intro} holds with $\mM = \omega \Diag{L_1,L_2,\dots,L_d}$~\cite[Theorem 8]{pcdm}. In many situations, $\mM$-smoothness is a very natural assumption. For instance, in the context of empirical risk minimization ({\acrshort{ERM}}), which is a key problem in supervised machine learning, $f$ is of the form
$f(x) = \frac{1}{n}\sum_{i=1}^n \phi_i(\mA_i x) + \frac{\mu}{2} \|x\|^2,$
where $\mA_i\in \R^{m\times d}$ are data matrices, $\phi_i:\R^m\to \R$ are loss functions and $\mu\geq 0$ is a regularization constant.  If $\phi_i$ is convex and  $\gamma_i$-smooth for all $i$, then $f$ is $\mu$-strongly convex and $\mM$-smooth with $\mM = (\frac{1}{n}\sum_i \gamma_i \mA_i^\top \mA_i) + \mu \mI$~\cite{qu2016coordinate2}. In these situations it is useful to design \texttt{CD} algorithms making full use of the information contained in the data as captured in the smoothness matrix $\mM$.

Given a sampling $S$ and $\mM$-smooth function $f$, let $v=(v_1,\dots,v_d)$ be positive constants satisfying the {\acrshort{ESO}} (expected separable overapproximation) inequality
\begin{equation}\label{eq:acd_v_def}
\mP\circ \mM\preceq \Diag{p_1 v_1, \dots, p_d v_d},
\end{equation}
where $\mP$ is the {\em probability matrix} associated with sampling $S$, defined by $\mP_{ij}\eqdef \Prob(i\in S \; \& \; j\in S)$, $p_i \eqdef \mP_{ii}=\Prob(i\in S)$ and $\circ$ denotes the Hadamard (i.e., elementwise) product of matrices. From now on we define the {\em probability vector} as $p\eqdef (p_1,\dots,p_d)\in \R^d$ and let $v=(v_1,\dots,v_d)\in \R^d$ be the vector of ESO parameters. With this notation, \eqref{eq:acd_v_def} can be equivalently written as $\mP\circ \mM\preceq \Diag{p\circ v}$. We say that $S$ is {\em proper} if $p_i>0$ for all $i$. 

It can be show by combining the results of \cite{nsync} and \cite{qu2016coordinate2} that under the above assumptions, the minibatch \texttt{CD} method \eqref{eq:acd_Parallel-CD-intro} with stepsizes $\alpha_i= 1/v_i$ enjoys the iteration complexity \begin{equation}\label{eq:acd_NSync}\left(\max_{i}\frac{ v_i}{ p_i \mu}\right) \log \frac{1}{\epsilon} .\end{equation}
Since in situations when $|S^k|=1$ with probability 1 once can choose $v_i=L_i$,  the complexity result \eqref{eq:acd_NSync} generalizes \eqref{eq:acd_089ff}. Inequality~\eqref{eq:acd_v_def} is  standard in minibatch coordinate descent literature. It was studied extensively in \cite{qu2016coordinate2}, and has been used to analyze parallel \texttt{CD} methods \cite{pcdm, nsync, approx}, distributed \texttt{CD} methods \cite{Hydra, Hydra2}, accelerated \texttt{CD} methods  \cite{approx, Hydra2, qu2016coordinate1, SCP}, and dual methods \cite{quartz, SCP}. 

\paragraph{Importance sampling for minibatches.}  It is easy to see, for instance, that if we do not restrict the class of samplings over which we optimize, then the trivial {\em full sampling} $S^k = [d]$ with probability 1 is optimal. For this sampling, $\mP$ is the matrix of all ones, $p_i=1$ for all $i$, and \eqref{eq:acd_v_def} holds for $v_i=L\eqdef \lambda_{\max}(\mM)$ for all $i$. The minibatch \texttt{CD} method \eqref{eq:acd_Parallel-CD-intro} reduces  to gradient descent, and the complexity estimate \eqref{eq:acd_NSync} becomes $(L/\mu) \log (1/\epsilon)$, which is the standard rate of gradient descent. However, typically we are interested in  finding the best sampling from the class of samplings which use a minibatch of size $\tau$, where $\tau\ll d$. While we have seen that the importance sampling $p_i=L_i/\sum_j L_j$ is optimal for $\tau=1$, in the minibatch case $\tau>1$  the problem of determining a sampling which minimizes the bound \eqref{eq:acd_NSync} is much more difficult. For instance, \cite{nsync} consider a certain parametric family of samplings where the problem of finding the best sampling from this family reduces to a linear program. 

Surprisingly, and in contrast to the situation in the $\tau=1$ case where an optimal sampling is known and is in general non-uniform, there is no minibatch sampling that is guaranteed to outperform $\tau$--nice sampling. We say that $S$ is $\tau$--nice if it samples uniformly from among all subsets of $[d]$ of cardinality $\tau$. The probability matrix of this sampling is given by 
$\mP = \frac{\tau}{d}\left((1-\beta) \mI + \beta \ones \ones^\top \right),$
where $\beta = \frac{\tau-1}{d-1}$ (assume $d>1$) and $\ones\in \R^d$ is the vector of all ones, and $p_i = \frac{\tau}{d}$ \cite{qu2016coordinate2}. It follows that the ESO inequality \eqref{eq:acd_v_def} holds for $v_i = (1-\beta) \mM_{ii} + \beta L.$ By plugging into  \eqref{eq:acd_NSync}, we get the iteration complexity
\begin{equation}\label{eq:acd_tau-nice-rate} \frac{d}{\tau}\left( \frac{ (1-\beta) \max_i \mM_{ii} + \beta L}{\mu}  \right) \log \frac{1}{\epsilon}.\end{equation}

This rate interpolates between the rate of \texttt{CD} with uniform probabilities (for $\tau=1$) and the rate of gradient descent (for $\tau=d$).


\section{Contributions}

\begin{table*}[t]

\centering
\begin{tabular}{|M{2.5cm}|M{5cm}|M{5cm}|N}
\hline
& \texttt{CD}  &   \texttt{ACD} &  \\ [0.1cm]
\hline
\hline
$\tau=1$, $p_i >0 $ & \begin{tabular}{c} $\displaystyle \left(\max_i \frac{L_i}{p_i \mu}\right) \log \frac{1}{\epsilon} $ \\ \cite{richtarik2014iteration} \end{tabular}& \begin{tabular}{c} $\displaystyle \sqrt{ \max_i \frac{L_i}{p_i^2 \mu} }   \log \frac{1}{\epsilon}$\\ {\bf (this work)} \end{tabular}  &\\
\hline
$\tau=1$, best $p_i$ &  \begin{tabular}{c} $\displaystyle \frac{\sum_i L_i}{\mu }\log \frac{1}{\epsilon}$; \quad  $p_i\propto L_i$ \\ \cite{rcdm} \end{tabular} & \begin{tabular}{c} $\displaystyle\frac{\sum_i \sqrt{L_i}}{\sqrt{\mu}}  \log \frac{1}{\epsilon}$;  $p_i\propto \sqrt{L_i}$  \\ \cite{allen2016even}\end{tabular} &\\ 
\hline
arbitrary sampling $S$ & \begin{tabular}{c} $\displaystyle\left(\max_i \frac{v_i}{p_i \mu}\right) \log \frac{1}{\epsilon}$ \\ \cite{nsync}  \end{tabular}& \begin{tabular}{c} $\displaystyle \sqrt{ \max_i \frac{v_i}{p_i^2 \mu} }   \log \frac{1}{\epsilon}$  \\ {\bf (this work)} \end{tabular}&\\ 
\hline
\end{tabular}
\caption{ Complexity results for non-accelerated (\texttt{CD}) and accelerated (\texttt{ACD}) coordinate descent methods for $\mu$-strongly convex functions and arbitrary sampling $S$. The last row corresponds to the setup with arbitrary proper sampling $S$ (i.e., a random subset of $[d]$ with the property that $p_i\eqdef \Prob(i\in S)>0$). We let  $\tau\eqdef \Exp{|S|}$ be the expected mini-batch size.  We assume that $f$ is $\mM$-smooth (see \eqref{eq:acd_M-smooth-intro}).  The positive constants $v_1,v_2,\dots,v_d$ are the ESO parameters (depending on $f$ and $S$), defined in \eqref{eq:acd_v_def}. The first row arises as a special of the third row in the non-minibatch (i.e., $\tau=1$) case. Here we have $v_i=L_i\eqdef \mM_{ii}$. The second row is a special case of the first row for the optimal choice of the probabilities $p_1,p_2,\dots,p_d$.}
\label{tab:main_acd}
\end{table*}

For {\em accelerated coordinate descent (\texttt{ACD})} without minibatching (i.e., when $\tau=1$), the currently best known iteration complexity result, due to~\cite{allen2016even}, is
\begin{equation}\label{eq:acd_ZAZ}\cO\left(\frac{\sum_i \sqrt{L_i}}{\sqrt{\mu}}  \log \frac{1}{\epsilon}\right).\end{equation}
 The probabilities used in the algorithm are proportional to the square roots of the coordinate-wise Lipschitz constants: $p_i \propto \sqrt{L_i}$. This is the first \texttt{CD} method with a complexity guarantee which does not explicitly depend on the dimension $n$, and is an improvement on the now-classical result of \cite{rcdm} giving the complexity
\[\cO\left(\sqrt{ \frac{d \sum_i L_i }{ \mu}}  \log \frac{1}{\epsilon} \right).\]
The rate \eqref{eq:acd_ZAZ} is always better than this, and can be up to $\sqrt{d}$ times better if the distribution of $L_i$ is extremely non-uniform. Unlike in the non-accelerated case described in the previous section, there is no complexity result for \texttt{ACD} with general probabilities such as \eqref{eq:acd_089ff}, or with an arbitrary sampling such as \eqref{eq:acd_NSync}. In fact, an \texttt{ACD} method was not even designed in such settings, despite a significant recent development in accelerated coordinate descent methods~\cite{rcdm,Lee2013,lin2014accelerated,qu2016coordinate1,allen2016even}. 

To summarize, our key contributions are:

\begin{itemize}
\item {\bf \texttt{ACD} with arbitrary sampling.} We design an \texttt{ACD} method which is able to operate with an {\em arbitrary sampling} of subsets of coordinates. We describe our method in Section~\ref{sec:acd_ACD}.
 
\item {\bf  Iteration complexity.} We prove (see Theorem~\ref{th:acd}) that the iteration complexity of \texttt{ACD} is 
\begin{equation}\label{eq:acd_ug98sg98s}\cO\left( \sqrt{ \max_i \frac{v_i}{p_i^2 \mu} }   \log \frac{1}{\epsilon}\right),\end{equation}
where $v_i$ are ESO parameters given by \eqref{eq:acd_v_def} and $p_i>0$ is the probability that coordinate $i$ belongs to the sampled set $S^k$: $p_i\eqdef \Prob(i\in S^k)$. The result of Allen-Zhu et al.\ \eqref{eq:acd_ZAZ} (\texttt{NUACDM}) can be recovered as a special case of  \eqref{eq:acd_ug98sg98s} by focusing on samplings defined by $S^k=\{i\}$ with probability $p_i \propto \sqrt{L_i} $ (recall that in this case $v_i=L_i$). When $S^k=[d]$ with probability 1, then our method reduces to accelerated gradient descent (Algorithm~\ref{alg:acd_intro}, \texttt{AGD}), and since $p_i=1$  and $v_i=L$ (the Lipschitz constant of $\nabla f$) for all $i$,  \eqref{eq:acd_ug98sg98s} reduces to the standard complexity of \texttt{AGD}:
$\cO(\sqrt{L/\mu} \log (1/\epsilon)).$

\item {\bf  Weighted strong convexity.} We prove a slightly more general result than \eqref{eq:acd_ug98sg98s} in which we allow the strong convexity of $f$ to be measured in a weighted Euclidean norm with weights $v_i/p_i^2$. In situations when $f$ is naturally strongly convex with respect to a weighted norm, this more general result will typically lead to a better complexity result than \eqref{eq:acd_ug98sg98s}, which is fine-tuned for standard strong convexity.
There are applications when $f$ is naturally a strongly convex with respect to some weighted norm~\cite{allen2016even}. 

\item {\bf  Minibatch methods.} We design several {\em new} importance samplings for minibatches, calculate the associated complexity results, and show through experiments that they significantly outperform the standard uniform samplings used in practice and constitute the state of the art. Our importance sampling leads to rates which are provably within a small factor from the best known rates, but can lead to an improvement by a factor of $\cO(d)$. We are the first to establish such a result, both for \texttt{CD} (Appendix~\ref{sec:acd_cd_imp}) and \texttt{ACD} (Section~\ref{sec:acd_import}). Further, the importance sampling we design for \texttt{CD} can be applied beyond coordinate descent algorithms: Chapters~\ref{jacsketch} and~\ref{sigmak} discuss an application in stochastic algorithms that subsample the finite sum.

\end{itemize}

The key complexity results obtained in this chapter are summarized and compared to prior results in Table~\ref{tab:main_acd}.

\section{The {\tt ACD} algorithm} \label{sec:acd_ACD}

The accelerated coordinate descent method (\texttt{ACD}) we propose is formalized as Algorithm~\ref{alg:acd_acd}. If we removed \eqref{eq:acd_x_update_acd} and \eqref{eq:acd_z_update} from the method, and replaced $y^{k+1}$ in \eqref{eq:acd_y_update} by $x^{k+1}$, we would recover the   \texttt{CD} method. Acceleration is obtained by the inclusion of the extrapolation steps \eqref{eq:acd_x_update_acd} and \eqref{eq:acd_z_update}.  
As mentioned before, we will analyze our method under a more general strong convexity assumption.

\begin{assumption}\label{ass:acd_sc}
Function $f$ is $\mu_w$-strongly convex with respect to the  $\|\cdot \|_{w}$ norm. That is, 
\begin{equation}\label{eq:acd_sc}
f(x+h)\geq f(x) +\langle \nabla f(x),h\rangle +\frac{\mu_w}{2}\|h\|_w^2,
\end{equation}
for all $x,h\in \R^d$, where $\mu_w>0$.
\end{assumption}

Note that if $f$ is $\mu$-strongly convex in the standard sense (i.e., for $w=(1,\dots,1)$), then  $f$ is $\mu_w$-strongly convex for any $w>0$ with $\mu_w = \min_i \frac{\mu}{w_i}.$ Considering a general $\mu_w$-strong convexity allows us to get a tighter convergence rate in some cases~\cite{allen2016even}.

\begin{algorithm}[!h]
\begin{algorithmic}[1]
\State \textbf{Parameters:} i.i.d.\ proper samplings $S^k\sim \cD$; $v,w\in \R^d_{++}$; $\mu_w>0$; stepsize parameters $\eta,\theta>0$. 
\State  Initial iterate $y^0=z^0\in \R^d $ 
\For {$k= 0,1,2, \dots $} 
\State $x^{k+1}=(1-\theta)y^{k}+\theta z^{k}$ \label{eq:acd_x_update_acd}
\State Get $S^k \sim \cD$
 \State $ y^{k+1}=x^{k+1}-\sum_{i\in S^k} \frac{1}{v_i} \nabla_if(x^{k+1}) e_i$\label{eq:acd_y_update}
  \State $z^{k+1}=\frac{1}{1+\eta\mu_w}\left(z^k+\eta\mu_w x^{k+1}-\sum_{i\in S^k}\frac{\eta}{p_i w_i } \nabla_i f(x^{k+1}) e_i \right)$ \label{eq:acd_z_update}
 \EndFor
\end{algorithmic}
\caption{\texttt{ACD} (Accelerated coordinate descent with arbitrary sampling)}
\label{alg:acd_acd}
\end{algorithm}

Using the tricks developed in \cite{Lee2013, approx, lin2014accelerated}, Algorithm~\ref{alg:acd_acd} can be implemented so that only $|S^k|$ coordinates are updated in each iteration. We are now ready derive a convergence rate of \texttt{ACD}.

\begin{theorem}[Convergence of \texttt{ACD}]\label{th:acd} Let $S^k$ be i.i.d.\ proper (but otherwise arbitrary) samplings. Let $\mP$ be the associated probability matrix and  $p_i\eqdef \Prob(i\in S^k)$.   Assume $f$ is $\mM$-smooth (see \eqref{eq:acd_M-smooth-intro}) and let $v$ be ESO parameters satisfying \eqref{eq:acd_v_def}.  Further, assume that $f$ is $\mu_w$- strong convex (with $\mu_w>0$) for
 \begin{equation}
w_i\eqdef \frac{v_i}{p_i^2}, \qquad i=1,2,\dots,d,\label{eq:acd_w_def}
\end{equation}
with respect to the weighted Euclidean norm $\|\cdot\|_w$ (i.e., we enforce Assumption~\ref{ass:acd_sc}).  Then  \begin{equation} \label{eq:acd_998dgff}\mu_w \leq \frac{\mM_{ii} p_i^2}{v_i} \leq p_i^2 
\leq 1, \qquad i=1,2,\dots d.\end{equation}
In particular, if $f$ is $\mu$-strongly convex with respect to the standard Euclidean norm, then  we can choose \begin{equation}\label{eq:acd_8h8hs8s}\mu_w = \min_i \frac{p_i^2 \mu}{v_i}.\end{equation} 
Finally, if we choose
\begin{eqnarray}
\theta&\eqdef&\frac{\sqrt{\mu_w^2+4\mu_w}-\mu_w}{2}=\frac{2\mu_w}{\sqrt{\mu_w^2+4\mu_w}+\mu_w} \nonumber \geq  0.618 \sqrt{\mu_\omega}
\label{eq:acd_tau_def_acd}
\end{eqnarray}
and 
$
\eta\eqdef \frac{1}{\theta},
$
then the random iterates of \texttt{ACD} satisfy
\begin{equation}\label{eq:acd_recur}
\E{P^{k}}
\leq
(1-\theta)^kP^0,
\end{equation}
where $P^k\eqdef \frac{1}{\theta^2}\left( f(y^k)-f(x^*)\right)+\frac{1}{2(1-\theta)}\|z^k-x^* \|_w^2$ and $x^*$ is the optimal solution of \eqref{eq:acd_problem}. 
\end{theorem}

 Noting that $1/0.618\leq 1.619$, as an immediate consequence of \eqref{eq:acd_recur} and  \eqref{eq:acd_tau_def_acd} we get bound
\begin{equation} \label{eq:acd_ius98g9skkk} k \geq \frac{1.619}{\sqrt{\mu_w}} \log \frac{1}{\epsilon} \quad \Rightarrow \quad \E{P^k} \leq \epsilon P^0.\end{equation}
 If $f$ is $\mu$-strongly convex, then by plugging \eqref{eq:acd_8h8hs8s} into \eqref{eq:acd_ius98g9skkk} we obtain the iteration complexity bound
\begin{equation}\label{eq:acd_98g98df} 1.619 \cdot \sqrt{\max_i \frac{v_i}{p_i^2 \mu}} \log \frac{1}{\epsilon}.\end{equation}
Complexity \eqref{eq:acd_98g98df} is our key result (also  mentioned in \eqref{eq:acd_ug98sg98s} and Table~\ref{tab:main_acd}).

\section{Importance sampling for minibatches \label{sec:acd_import}}

Let $\tau\eqdef \E{|S^k|}$ be the expected minibatch size. The next theorem provides an insightful lower bound for the complexity of \texttt{ACD} we established, one independent of  $p$ and $v$. 

\begin{theorem}[Limits of minibatch performance] \label{thm:acd_LB} Let the assumptions of Theorem~\ref{th:acd}  be satisfied and let $f$ be $\mu$-strongly convex. Then the dominant term in the rate~\eqref{eq:acd_98g98df} of \texttt{ACD} admits the lower bound
\begin{equation}\label{eq:acd_ineq-minibatch-speedup} \sqrt{\max_i \frac{v_i}{p_i^2 \mu}}  \geq \frac{\sum_{i} \sqrt{\mM_{ii}}}{\tau \sqrt{\mu}}. \end{equation}
\end{theorem}
Note that for $\tau=1$ we have $\mM_{ii}=v_i=L_i$, and the lower bound is achieved by using the importance sampling $p_i\propto \sqrt{L_i}$.
Hence, this bound gives a  limit on how much speedup, compared to the best known complexity in the  $\tau=1$ case, we can hope for  as we increase $\tau$. The bound says we can not hope for better than linear speedup in the minibatch size.  An analogous result (obtained by removing all the squares and square roots in \eqref{eq:acd_ineq-minibatch-speedup}) was established in \cite{nsync} for \texttt{CD}.


In what follows, it will be useful to write the complexity result \eqref{eq:acd_98g98df} in a new form by considering a specific choice of the ESO vector $v$.

\begin{lemma}\label{thm:acd_special-ESO-result} Choose any proper sampling $S$ and let $\mP$ be its probability matrix and $p$ its probability vector. Let $c(S,\mM) \eqdef \lambda_{\max}(\mP'\circ \mM'),$ where $\mP' \eqdef \mD^{-1/2} \mP \mD^{-1/2}$, $\mM' \eqdef \mD^{-1}\mM \mD^{-1}$ and $\mD \eqdef \Diag{p}$. Then the vector $v$ defined by $v_i \eqdef c(S,\mM) p_i^2$ satisfies  the ESO inequality \eqref{eq:acd_v_def} and the total complexity~\eqref{eq:acd_98g98df} becomes
\begin{equation}\label{eq:acd_98g98df2} 
1.619 \cdot \frac{\sqrt{c(S,\mM)}}{\sqrt{\mu}} \log \frac{1}{\epsilon}.
\end{equation}
\end{lemma}
Let $\trace{\cdot}$ be a trace function. Since $\frac{1}{d}\trace{\mP'\circ \mM'}\leq c(S,\mM) \leq \trace{\mP'\circ \mM'}$ and \[\trace{\mP'\circ \mM'}=\sum_i \mP'_{ii} \mM'_{ii} = \sum_i  \mM'_{ii} =\sum_i \mM_{ii}/p_i^2,\] we get the bounds:
\begin{eqnarray}
 \sqrt{\frac{1}{d}\sum_i \frac{\mM_{ii}}{p_i^2\mu}} \log \frac{1}{\epsilon} \leq   \sqrt{\frac{c(S,\mM)}{\mu}} \log \frac{1}{\epsilon}  \leq  \sqrt{\sum_i \frac{\mM_{ii}}{p_i^2\mu}} \log \frac{1}{\epsilon}.
   \label{eq:acd_nbisg8dd}
\end{eqnarray}

\begin{table*}[t]
{
\footnotesize
\centering
\begin{tabular}{|M{1.9cm}|M{4.2cm}|M{3.2cm}|M{4.4cm}|N}
\hline
 Lower bound & $\displaystyle S_1: p_i = \frac{\tau}{d}$  & $\displaystyle S_2: \frac{p_i^2}{\mM_{ii}}\propto 1$ & $\displaystyle S_3: \frac{p_i^2}{\mM_{ii}} \propto 1-p_i$ &\\
\hline
\phantom{xxx} &&&& \\
 $\displaystyle \frac{\sum_i \sqrt{\mM_{ii}}}{\tau \sqrt{\mu}}$  & $ \displaystyle  \frac{d\sqrt{(1-\beta) \max_i \mM_{ii} + \beta L}}{ \tau\sqrt{\mu}}   $ &  $\displaystyle \frac{\gamma\sum_i \sqrt{\mM_{ii}}}{\tau \sqrt{\mu}}$ & $\displaystyle \omega  \frac{d\sqrt{(1-\beta) \max_i \mM_{ii} + \beta L}}{ \tau\sqrt{\mu}} $ & \\ 
\phantom{xxx} &&&& \\ 
\hline
\phantom{xxx} &&&& \\
\eqref{eq:acd_ineq-minibatch-speedup} & 
\begin{tabular}{c} = uniform \texttt{ACD} for $\tau=1$ \\
= \texttt{AGD}  for $\tau=d$ 
\end{tabular} 
& 
\begin{tabular}{c} $\leq \sqrt{d} \times $  lower bound \\ $\bullet$ $\displaystyle \tau\leq \frac{\sum_j \sqrt{\mM_{jj}}}{\max_i \mM_{ii}}$ \end{tabular} & 
\begin{tabular}{c} $\bullet$ fastest in practice \\ $\bullet$ any $\tau$ allowed \end{tabular} &
 \\ 
 \phantom{xxx} &&&& \\
\hline 
\end{tabular}
}
\caption{ New complexity results for  \texttt{ACD}  with minibatch size $\tau=\E{|S^k|}$ and various samplings (we suppress $\log (1/\epsilon)$ factors in all expressions).  Constants: $\mu=$ strong convexity constant of $f$, $L=\lambda_{\max}(\mM)$, $\beta = (\tau-1)/(d-1)$, $1\leq \gamma \leq \sqrt{d}$,  and $\omega \leq \cO(\sqrt{\tau})$ ($\omega$ can be as small as $\cO(\tau/d)$).
} \label{tab:main2_acd}
\end{table*}

\subsection{Sampling 1: standard uniform minibatch sampling \label{sec:acd_sam1}} Let $S_1$ be the $\tau$-nice sampling. It can be shown (see Lemma~\ref{lem:acd_tau-nice-2nd-derivation}) that $c(S_1,\mM) \leq\frac{d^2}{\tau^2} ((1-\beta)\max_i \mM_{ii} + \beta L),$ and  hence the iteration complexity \eqref{eq:acd_98g98df} becomes \begin{equation}\label{eq:acd_ug98sg98s-tau-nice}\cO\left( \frac{d}{\tau}\sqrt{ \frac{(1-\beta) \max_i \mM_{ii} + \beta L}{ \mu} }   \log \frac{1}{\epsilon}\right).\end{equation}
This result interpolates between \texttt{ACD} with uniform probabilities (for $\tau=1$) and accelerated gradient descent (for $\tau=d$). Note that the rate \eqref{eq:acd_ug98sg98s-tau-nice} is a strict improvement on the \texttt{CD} rate  \eqref{eq:acd_tau-nice-rate}.

\subsection{Sampling 2: importance sampling for minibatches} \label{sec:acd_sam2}
Consider now the sampling $S_2$ which  includes every $i \in [d]$ in $S_2$, independently, with probability $p_i = \tau \frac{\sqrt{\mM_{ii}}}{\sum_j \sqrt{\mM_{jj}}}.$ This sampling was not considered in the literature before. Note that $\E{|S_2|}=\sum_i p_i = \tau$. For this sampling, bounds \eqref{eq:acd_nbisg8dd}  become:
\begin{eqnarray}
\frac{\sum_{i} \sqrt{\mM_{ii}}}{\tau \sqrt{\mu}}\log \frac{1}{\epsilon} \leq    \sqrt{\frac{c(S,\mM)}{\mu}} \log \frac{1}{\epsilon} \leq  \frac{\sqrt{d}\sum_{i} \sqrt{\mM_{ii}}}{\tau \sqrt{\mu}} \log \frac{1}{\epsilon}.
 \label{eq:acd_LB--x} 
\end{eqnarray}

Clearly, with this sampling we obtain an \texttt{ACD} method with complexity within a $\sqrt{d}$ factor from the lower bound established in Theorem~\ref{thm:acd_LB}. For $\tau=1$ we have  $\mP'=\mI$ and hence \begin{eqnarray*}c(S,\mM)&=&\lambda_{\max}(\mI\circ \mM') = \lambda_{\max}(\Diag{\mM'}) \\&=& \max_i \mM_{ii}/p_i^2 = \left(\sum_j \sqrt{\mM_{jj}}\right)^2.\end{eqnarray*} Thus, the rate of \texttt{ACD} achieves the lower bound in \eqref{eq:acd_LB--x} (see also \eqref{eq:acd_ZAZ}) and we recover the best current rate of \texttt{ACD} in the $\tau=1$ case, established by Allen-Zhu et. al. \cite{allen2016even}. However, the sampling has an important limitation: it can be used for $\tau\leq \sum_j \sqrt{\mM_{jj}}/\max_i \mM_{ii}$ only  as otherwise the probabilities $p_i$ exceed 1.

\subsection{Sampling 3: another importance sampling for minibatches} \label{sec:acd_sam3}
Now consider sampling $S_3$ which includes each coordinate $i$ within $S_3$ independently, with probability $p_i$ satisfying the relation $p_i^2/\mM_{ii} \propto 1-p_i$. This is equivalent to setting
\begin{equation}\label{eq:acd_p_imp_def_acc1}
p_{i}\eqdef\frac{2 \mM_{ii}}{\sqrt{\mM_{ii}^2+2\mM_{ii}\delta}+\mM_{ii}},
\end{equation}
where $\delta$ is a scalar for which $\sum_i p_i=\tau$. This sampling was not considered in the literature before. Probability vector $p$ was chosen as~\eqref{eq:acd_p_imp_def_acc1} for two reasons: i) $p_i\leq 1$ for all $i$, and therefore the sampling can be used for all $\tau$ in contrast to $S_1$, and ii) we can prove Theorem~\ref{thm:acd_comparison}. 


Let $c_1\eqdef c(S_1,\mM)$ and $c_3\eqdef c(S_3,\mM)$. In light of \eqref{eq:acd_98g98df2}, Theorem~\ref{thm:acd_comparison} compares $S_1$ and $S_3$ and says that \texttt{ACD} with $S_3$ has at most $\cO(\sqrt{\tau})$ times worse rate compared to \texttt{ACD} with $S_1$, but has the capacity to be $\cO(d/\tau)$ times better. We prove in Appendix~\ref{sec:acd_cd_imp} a similar theorem for \texttt{CD}. We stress that, despite some advances in the development of importance samplings for minibatch methods~\cite{nsync, csiba2018importance}, $S_1$ was until now the state-of-the-art in theory for \texttt{CD}. We are the first to give a provably better rate in the sense of Theorem~\ref{thm:acd_nonacc_comp}. The numerical experiments show that $S_3$ consistently outperforms $S_1$, and often dramatically so.

\begin{theorem}\label{thm:acd_comparison} The leading complexity terms $c_1$ and $c_3$ of \texttt{ACD} (Algorithm~\eqref{eq:acd_Parallel-CD-intro}) with samplings $S_1$, and $S_3$, respectively, defined in Lemma~\ref{thm:acd_special-ESO-result}, compare as follows:
 \begin{equation}\label{eq:acd_thm_comparison}
  c_3 \leq 2\frac{(2d-\tau)(d\tau+d-\tau)}{(d-\tau)^2}   c_1 =\cO(\tau) c_1.
 \end{equation}
Moreover, there exists $\mM$ where $c_3\leq \cO\left(\frac{\tau^2}{d^2}\right)c_1$. 
\end{theorem}

In real world applications,  minibatch size $\tau$ is limited by  hardware and in typical situations, one has  $\tau \ll d$, oftentimes $\tau=\cO(1)$. The importance of Theorem~\ref{thm:acd_comparison} is best understood from this  perspective.

\section{Experiments} \label{sec:acd_exp}

We perform extensive numerical experiments to justify that minibatch \texttt{ACD} with importance sampling works well in practice.

We first present some synthetic examples in Section~\ref{sec:acd_artif} in order to have better understanding of both acceleration and importance sampling, and to see how it performs on what type of data. We also study how minibatch size influences the convergence rate. 

Then, in Section~\ref{exp:logreg}, we work with  logistic regression problem on LibSVM~\cite{chang2011libsvm} data.  For small datasets, we choose the parameters of \texttt{ACD} as theory suggests and for large ones, we estimate them, as we describe in the main body of the chapter. Lastly, we tackle dual of SVM problem with squared hinge loss, which we present in Section~\ref{sec:acd_SVM}.\footnote{Coordinate descent methods which allow for separable proximal operator were proven to be efficient to solve ERM problem, when applied on dual~\cite{shalev2011stochastic,sdca,shalev2014accelerated, iprox-sdca}.
Although we do not develop proximal methods in this chapter, we empirically demonstrate that \texttt{ACD} allows for this extension as well. As a specific problem to solve, we choose dual of SVM with hinge loss. The results and a detailed description of the experiment are presented in Section~\ref{sec:acd_SVM}, and are indeed in favour of \texttt{ACD} with importance sampling. Therefore, \texttt{ACD} is not only suitable for \emph{big dimensional} problems, it can handle the \emph{big data} setting as well.}

In most of plots we compare of both accelerated and non-accelerated \texttt{CD} with all samplings $S_1,S_2,S_3$ introduced in Sections~\ref{sec:acd_sam1}, \ref{sec:acd_sam2} and \ref{sec:acd_sam3} respectively. We refer to \texttt{ACD} with sampling $S_3$ as \texttt{AN} (\texttt{A}ccelerated \texttt{N}onuniform), \texttt{ACD} with sampling $S_1$ as \texttt{AU}, \texttt{ACD} with sampling $S_2$ as \texttt{AN2}, \texttt{CD} with sampling $S_3$ as \texttt{NN}, \texttt{CD} with sampling $S_1$as \texttt{NU} and \texttt{CD} with sampling $S_2$ as \texttt{NN2}. As for Sampling 2, it might happen that probabilities become larger than one if $\tau$ is large (see Section~\ref{sec:acd_sam2}), we set those probabilities to 1 while keeping the rest as it is.

 All the experimental results clearly show that  acceleration, importance sampling and minibatching have a significant impact on practical performance of \texttt{CD} methods. Moreover, the difference in the performance of samplings $S_2$ and $S_3$ is negligible, and therefore we recommend using $S_3$, as it is not limited by the bound on expected minibatch size $\tau$.

\subsection{Synthetic quadratics~\label{sec:acd_artif}}
As we mentioned, the goal of this section is to provide a better understanding of both acceleration and importance sampling. For this purpose we consider as simple setting as possible -- minimizing quadratic
\begin{equation}\label{eq:acd_quadratic}
f(x)=\frac12 x^\top \mM x-b^\top x,
\end{equation}
where $b\sim N(0,I)$ and $\mM$ is chosen as one of the 5 types, as Table~\ref{tab:types} suggests. 

\begin{table}
\centering \small
\begin{tabular}{|c|c|}
 \hline
 {\bf Type} & $\mM$\\
 \hline
 \hline
1 & $\mA^\top \mA + \mI$ for $\mA^{\frac{d}{2}\times d}$;  have independent entries from $N(0,1)$
\\  \hline
2 & 
$\mA^\top \mA + \mI$ for $\mA^{2d\times d}$;  have independent entries from $N(0,1)$
\\ \hline
3 & 
$\Diag{1,2,\dots,d}$
\\ \hline
4 & 
$\mA + \mI$, $\mA_{d,d}=d$, $\mA_{1:(d-1),1:(d-1)}=1$, $\mA_{1:(d-1),d}= \mA_{d,1:(d-1)}=0$
\\ \hline
5 & 
$\mA^\top \mD \mA + \mI$ for $\mA^{\frac{d}{2}\times d}$;  have independent entries from $N(0,1)$, \\ & $\mD=\frac{1}{\sqrt{d}}\Diag{1,2,\dots,d}$ 
\\ 
 \hline
\end{tabular}
\caption{Problem types for testing \texttt{ACD}.}\label{tab:types}
\end{table}
In the first example we perform (Figure~\ref{fig:acd_artif_1}), we compare the performance of both accelerated and non-accelerated algorithm with both nonuniform and $\tau$ nice sampling on problems as per Table~\ref{tab:types}. In all experiments, we set $d=1 000$ and we plot a various choices of $\tau$. 

\subsubsection{Comparison of methods on synthetic data \label{sec:acd_marginal}}

Figure~\ref{fig:acd_artif_1} presents the numerical performance of \texttt{ACD} for various types of synthetic problems given by~\eqref{eq:acd_quadratic} and Table~\ref{tab:types}. It suggests what our theory shows: accelerated algorithm is always faster than its non-accelerated counterpart, and on top of that, performance of $\tau$--nice sampling ($S_1$) can be negligibly faster than  importance sampling ($S_2,S_3$), but is usually significantly slower. A significance of the importance sampling is mainly demonstrated on problem type 4, which roughly coincides with Examples~\ref{ex:nonacc_imp_diff} and~\ref{ex:acc_imp_diff}. Figure~\ref{fig:acd_artif_1} presents Sampling 2 only for the cases when the bound on $\tau$ form Section~\ref{sec:acd_sam2} is satisfied. 
\begin{figure}[H]
\centering
\begin{minipage}{0.25\textwidth}
  \centering
\includegraphics[width =  \textwidth ]{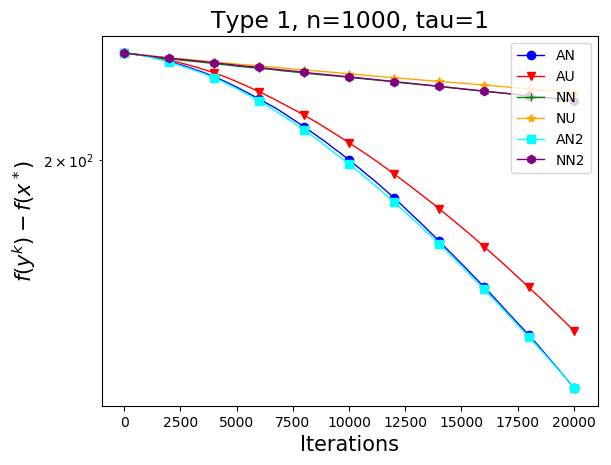}
\end{minipage}%
\begin{minipage}{0.25\textwidth}
  \centering
\includegraphics[width =  \textwidth ]{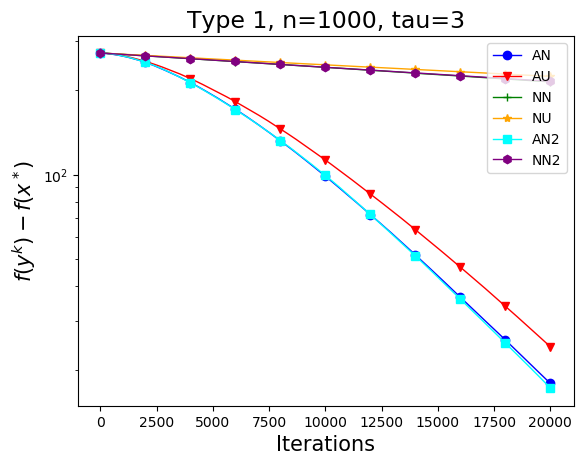}
\end{minipage}%
\begin{minipage}{0.25\textwidth}
  \centering
\includegraphics[width =  \textwidth ]{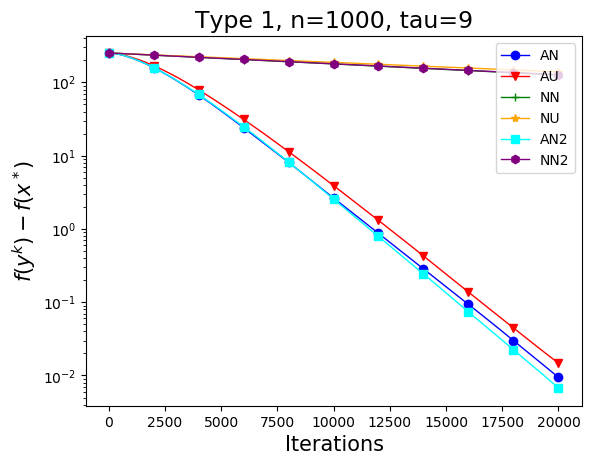}
\end{minipage}%
\begin{minipage}{0.25\textwidth}
  \centering
\includegraphics[width =  \textwidth ]{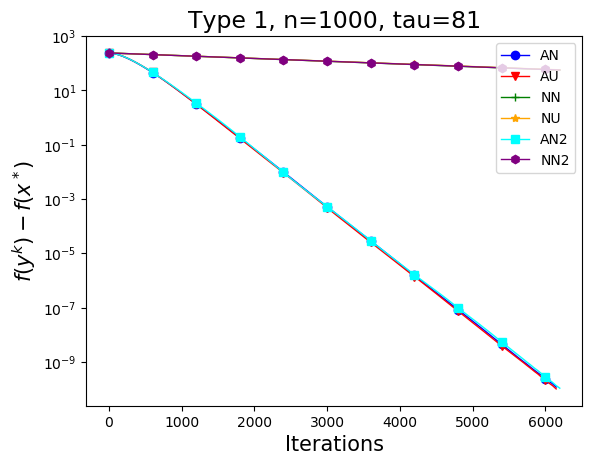}
\end{minipage}%
\\
\begin{minipage}{0.25\textwidth}
  \centering
\includegraphics[width =  \textwidth ]{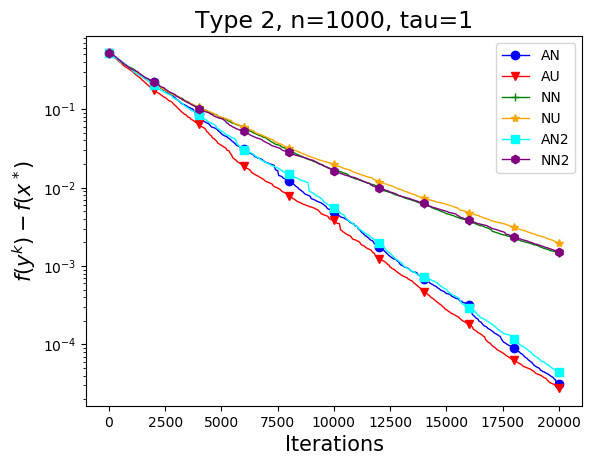}
\end{minipage}%
\begin{minipage}{0.25\textwidth}
  \centering
\includegraphics[width =  \textwidth ]{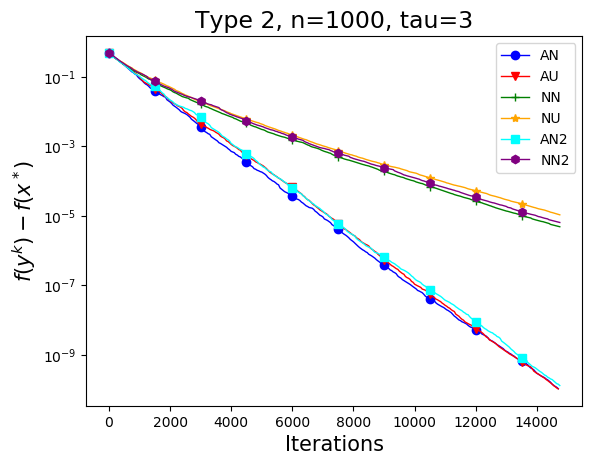}
\end{minipage}%
\begin{minipage}{0.25\textwidth}
  \centering
\includegraphics[width =  \textwidth ]{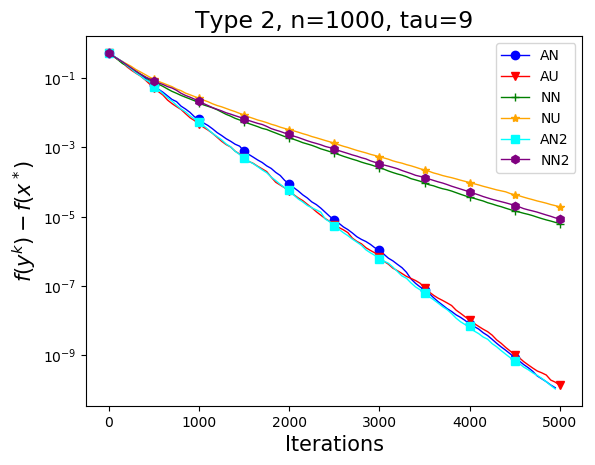}
\end{minipage}%
\begin{minipage}{0.25\textwidth}
  \centering
\includegraphics[width =  \textwidth ]{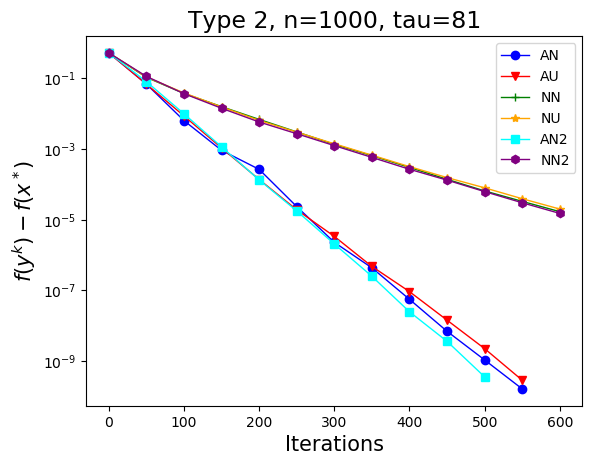}
\end{minipage}%
\\
\begin{minipage}{0.25\textwidth}
  \centering
\includegraphics[width =  \textwidth ]{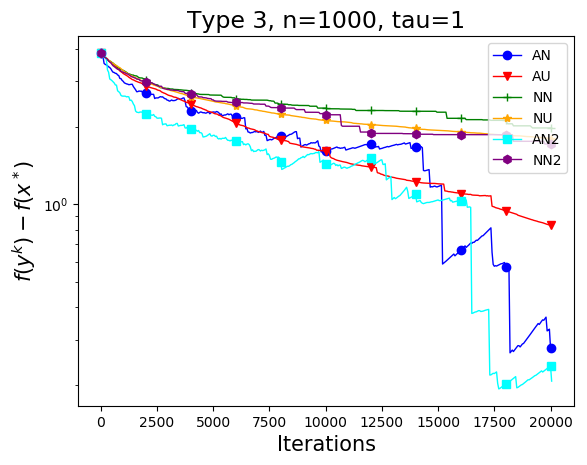}
\end{minipage}%
\begin{minipage}{0.25\textwidth}
  \centering
\includegraphics[width =  \textwidth ]{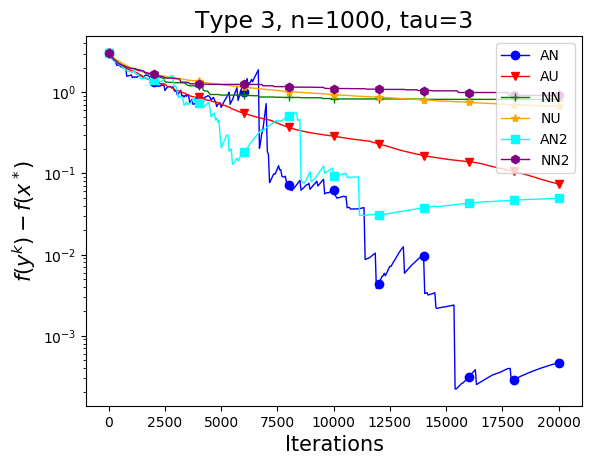}
\end{minipage}%
\begin{minipage}{0.25\textwidth}
  \centering
\includegraphics[width =  \textwidth ]{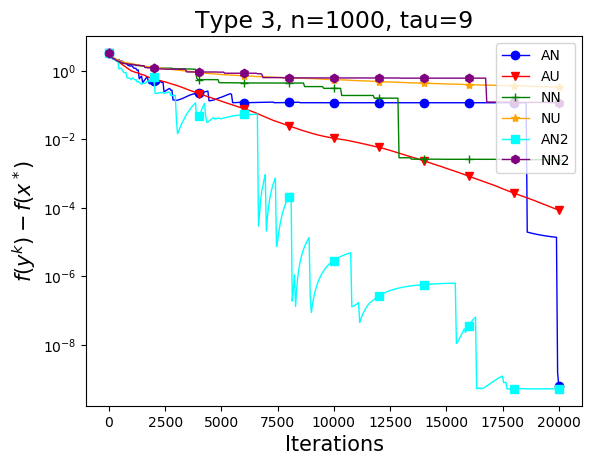}
\end{minipage}%
\begin{minipage}{0.25\textwidth}
  \centering
\includegraphics[width =  \textwidth ]{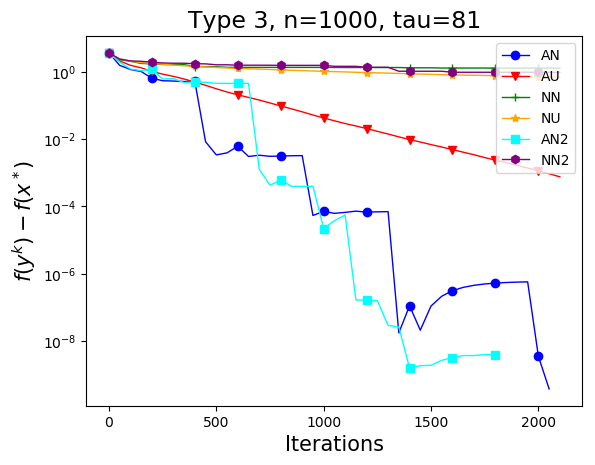}
\end{minipage}%
\\
\begin{minipage}{0.25\textwidth}
  \centering
\includegraphics[width =  \textwidth ]{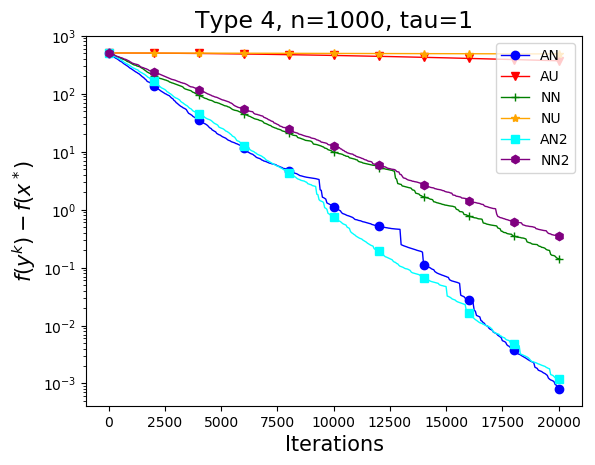}
\end{minipage}%
\begin{minipage}{0.25\textwidth}
  \centering
\includegraphics[width =  \textwidth ]{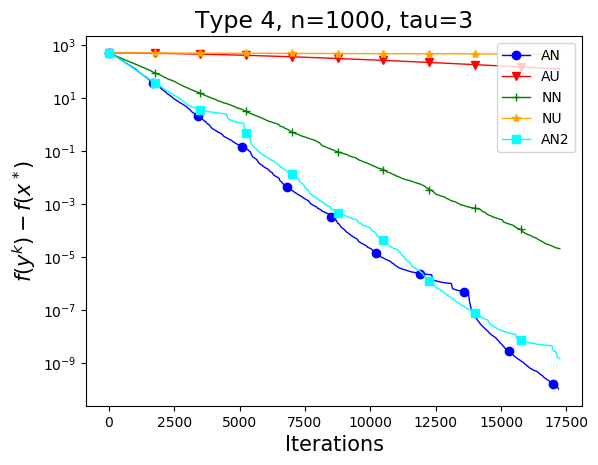}
\end{minipage}%
\begin{minipage}{0.25\textwidth}
  \centering
\includegraphics[width =  \textwidth ]{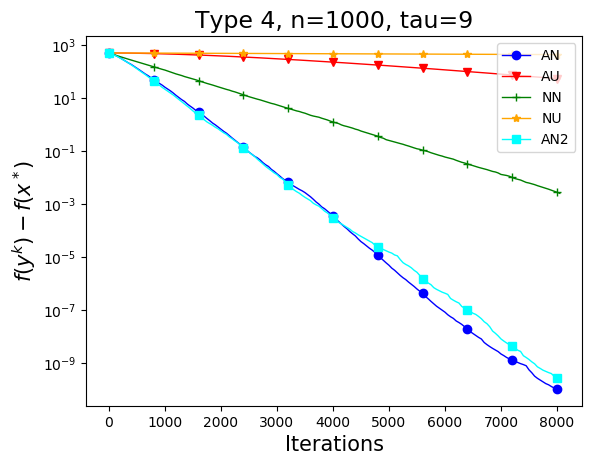}
\end{minipage}%
\begin{minipage}{0.25\textwidth}
  \centering
\includegraphics[width =  \textwidth ]{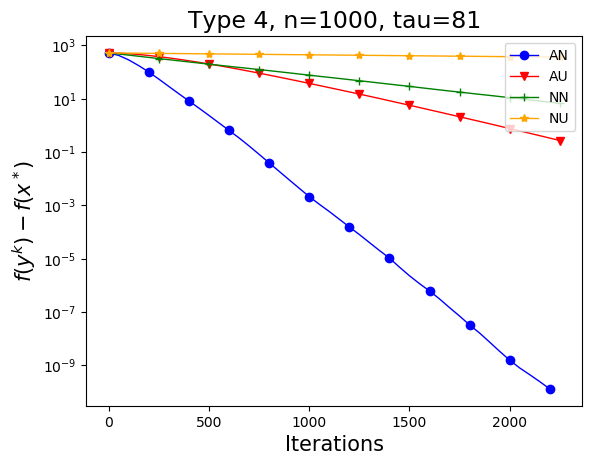}
\end{minipage}%
\\
\begin{minipage}{0.25\textwidth}
  \centering
\includegraphics[width =  \textwidth ]{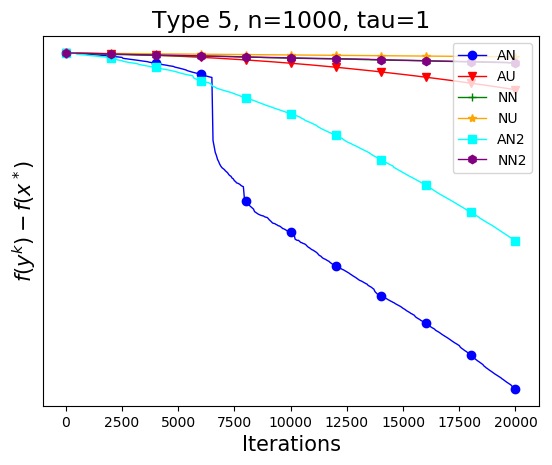}
\end{minipage}%
\begin{minipage}{0.25\textwidth}
  \centering
\includegraphics[width =  \textwidth ]{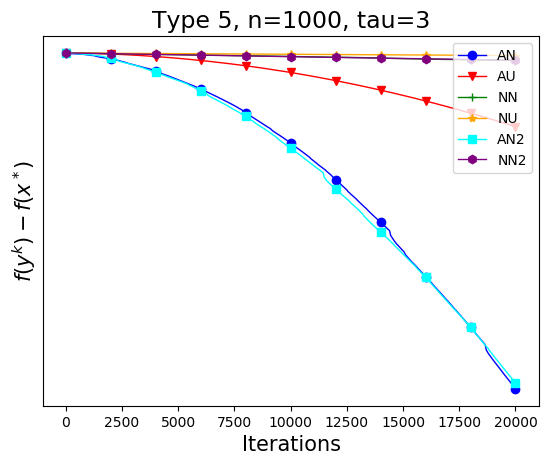}
\end{minipage}%
\begin{minipage}{0.25\textwidth}
  \centering
\includegraphics[width =  \textwidth ]{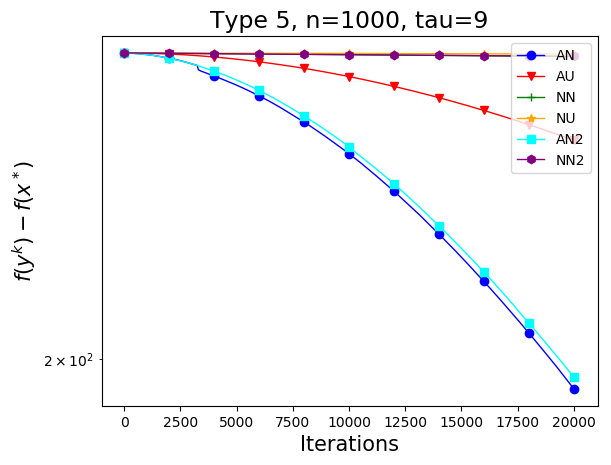}
\end{minipage}%
\begin{minipage}{0.25\textwidth}
  \centering
\includegraphics[width =  \textwidth ]{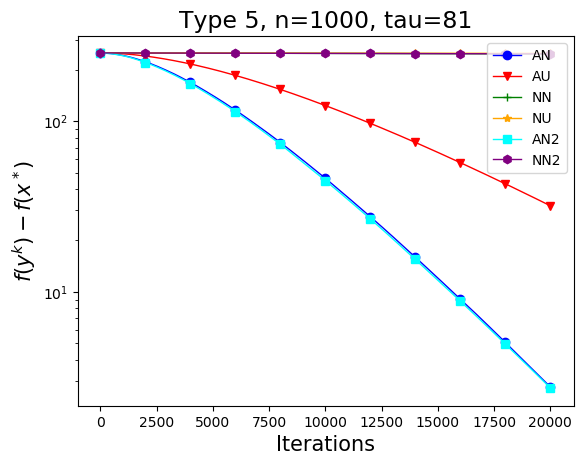}
\end{minipage}%
\\
\caption{Coordinate descent. Comparison of accelerated, nonaccelerated algorithm with both importance and $\tau$ nice sampling for a various quadratic problems.}\label{fig:acd_artif_1}
\end{figure}

\subsubsection{Speedup in $\tau$ \label{sec:acd_tau}}
The next experiment shows an empirical speedup for the coordinate descent algorithms for a various types of problems. For simplicity, we do not include Sampling 2. Figure~\ref{fig:acd_tau_comp} provides the results. Oftentimes, the empirical speedup (in terms of the number of iteration) in $\tau$ is close to linear, which demonstrates the power and significance of minibatching.

\begin{figure}[H]
\centering
\begin{minipage}{0.25\textwidth}
  \centering
\includegraphics[width =  \textwidth ]{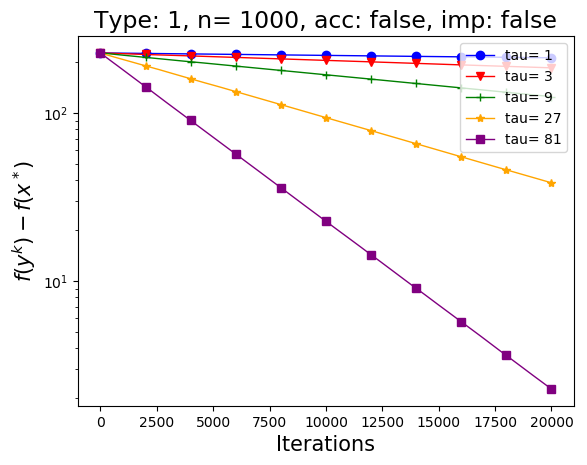}
\end{minipage}%
\begin{minipage}{0.25\textwidth}
  \centering
\includegraphics[width =  \textwidth ]{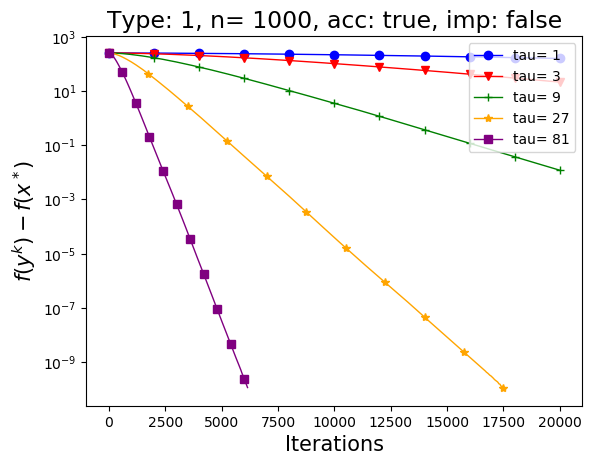}
\end{minipage}%
\begin{minipage}{0.25\textwidth}
  \centering
\includegraphics[width =  \textwidth ]{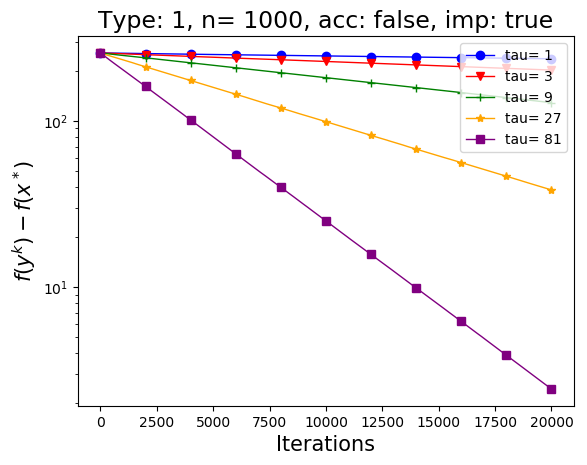}
\end{minipage}%
\begin{minipage}{0.25\textwidth}
  \centering
\includegraphics[width =  \textwidth ]{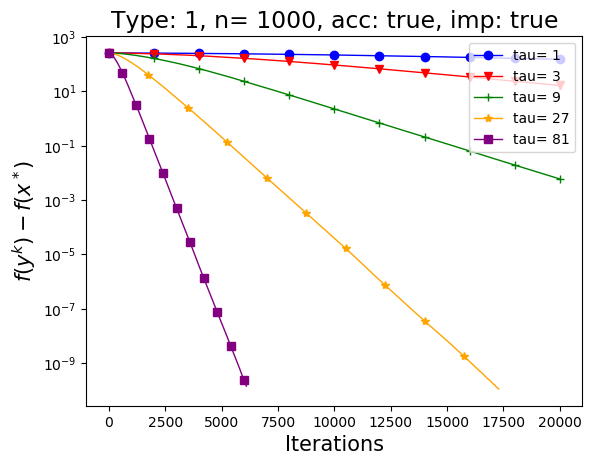}
\end{minipage}%
\\
\begin{minipage}{0.25\textwidth}
  \centering
\includegraphics[width =  \textwidth ]{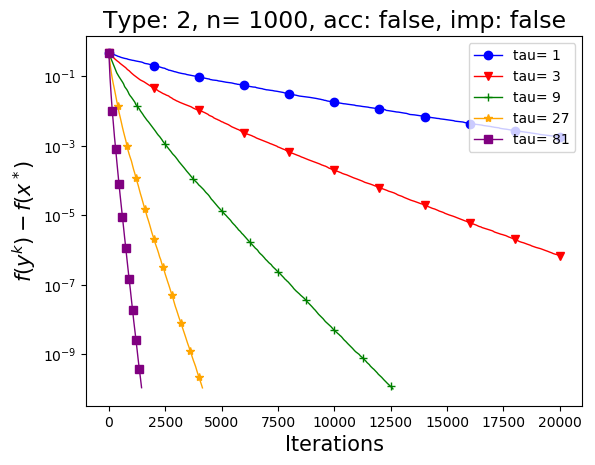}
\end{minipage}%
\begin{minipage}{0.25\textwidth}
  \centering
\includegraphics[width =  \textwidth ]{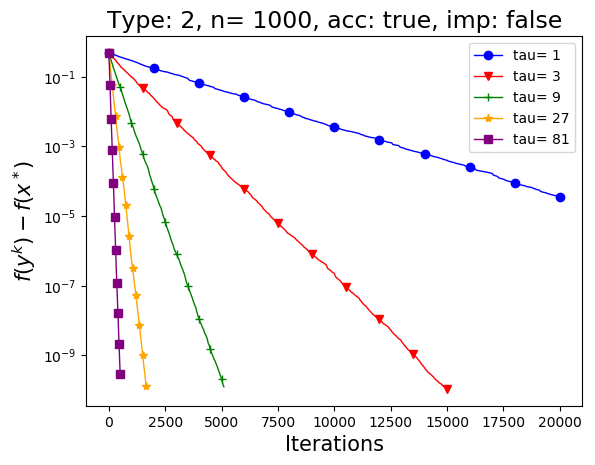}
\end{minipage}%
\begin{minipage}{0.25\textwidth}
  \centering
\includegraphics[width =  \textwidth ]{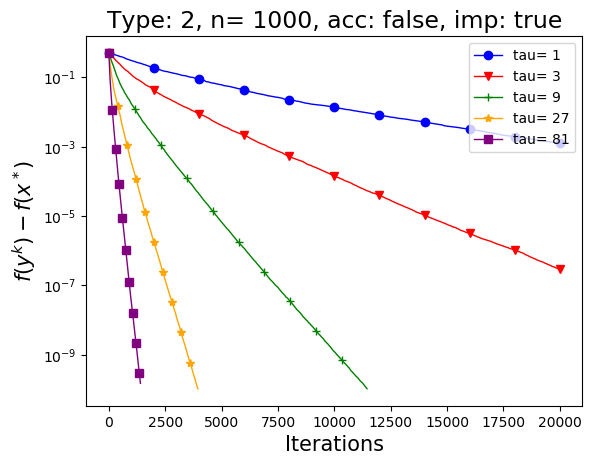}
\end{minipage}%
\begin{minipage}{0.25\textwidth}
  \centering
\includegraphics[width =  \textwidth ]{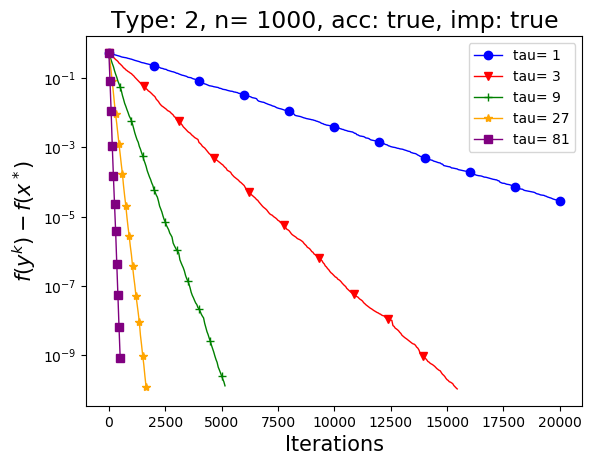}
\end{minipage}%
\\
\begin{minipage}{0.25\textwidth}
  \centering
\includegraphics[width =  \textwidth ]{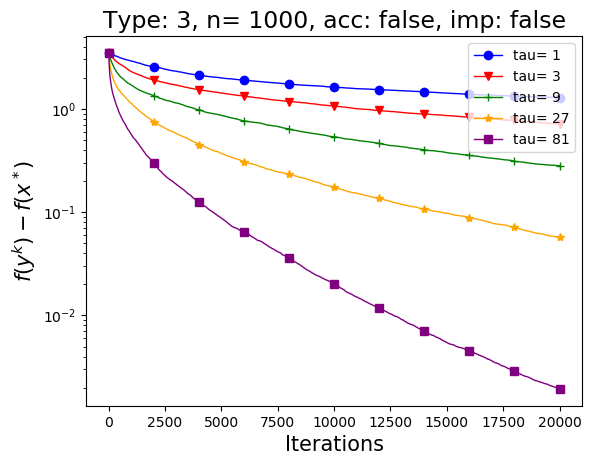}
\end{minipage}%
\begin{minipage}{0.25\textwidth}
  \centering
\includegraphics[width =  \textwidth ]{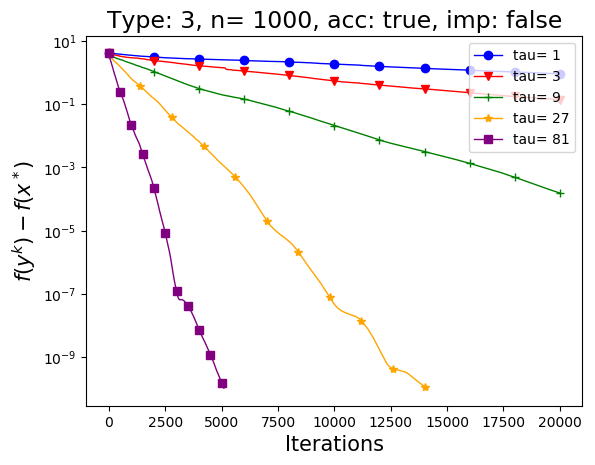}
\end{minipage}%
\begin{minipage}{0.25\textwidth}
  \centering
\includegraphics[width =  \textwidth ]{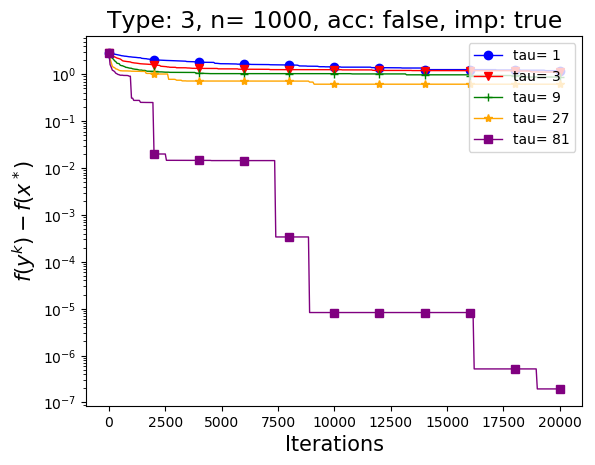}
\end{minipage}%
\begin{minipage}{0.25\textwidth}
  \centering
\includegraphics[width =  \textwidth ]{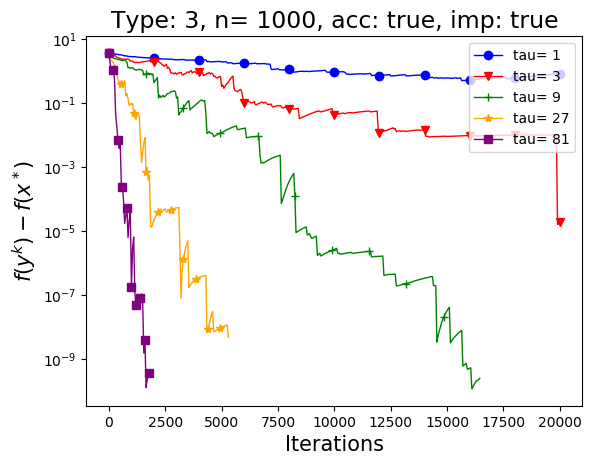}
\end{minipage}%
\\
\begin{minipage}{0.25\textwidth}
  \centering
\includegraphics[width =  \textwidth ]{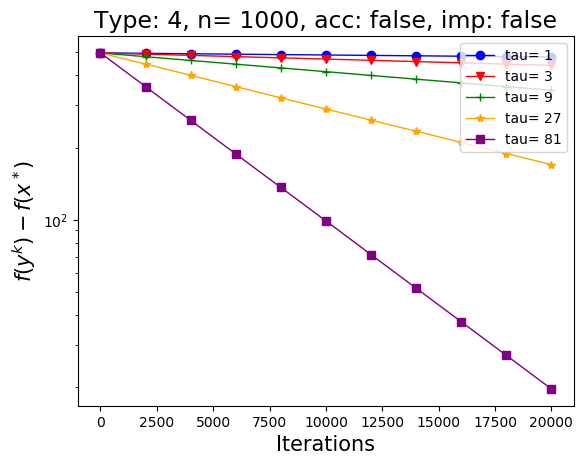}
\end{minipage}%
\begin{minipage}{0.25\textwidth}
  \centering
\includegraphics[width =  \textwidth ]{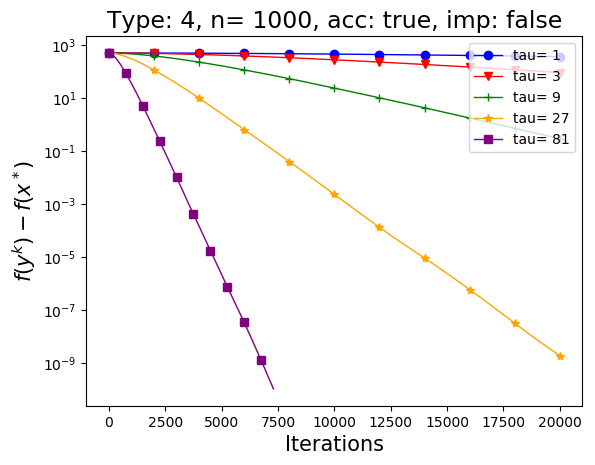}
\end{minipage}%
\begin{minipage}{0.25\textwidth}
  \centering
\includegraphics[width =  \textwidth ]{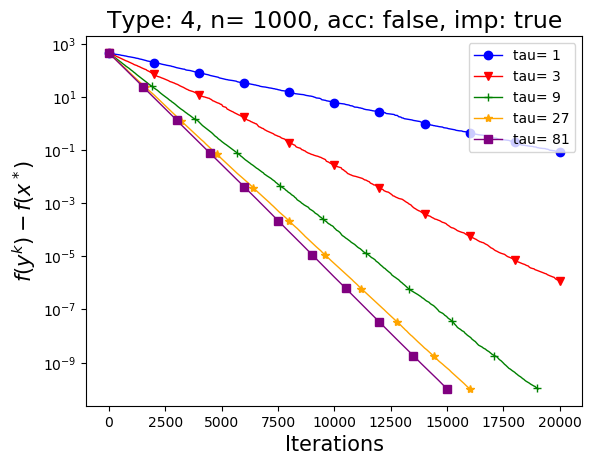}
\end{minipage}%
\begin{minipage}{0.25\textwidth}
  \centering
\includegraphics[width =  \textwidth ]{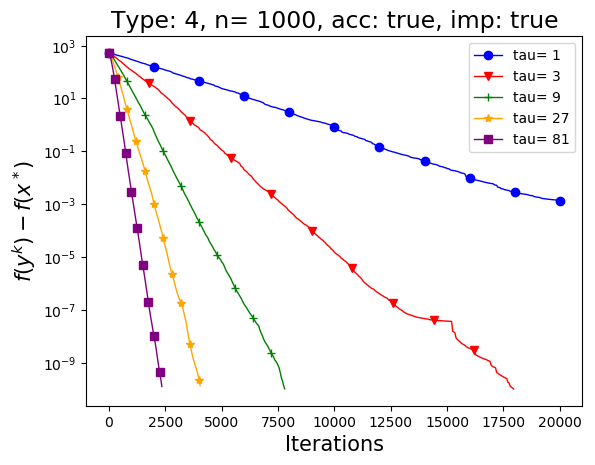}
\end{minipage}%
\\
\begin{minipage}{0.25\textwidth}
  \centering
\includegraphics[width =  \textwidth ]{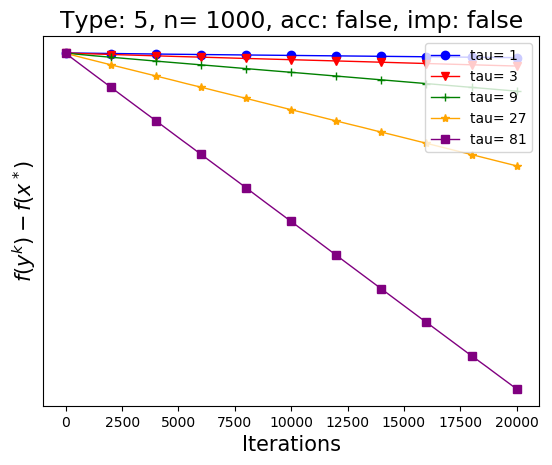}
\end{minipage}%
\begin{minipage}{0.25\textwidth}
  \centering
\includegraphics[width =  \textwidth ]{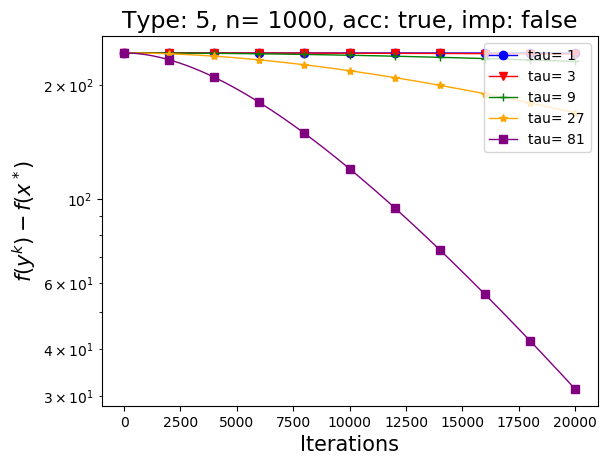}
\end{minipage}%
\begin{minipage}{0.25\textwidth}
  \centering
\includegraphics[width =  \textwidth ]{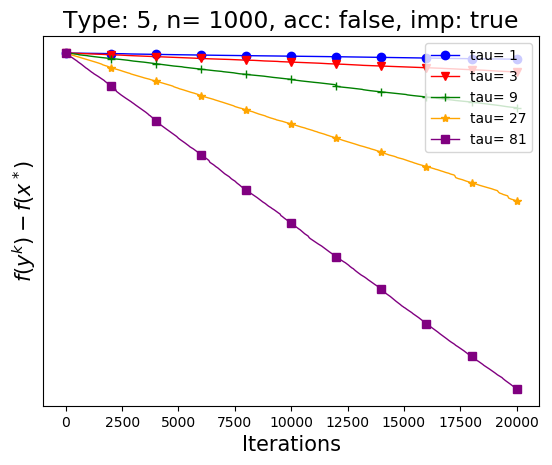}
\end{minipage}%
\begin{minipage}{0.25\textwidth}
  \centering
\includegraphics[width =  \textwidth ]{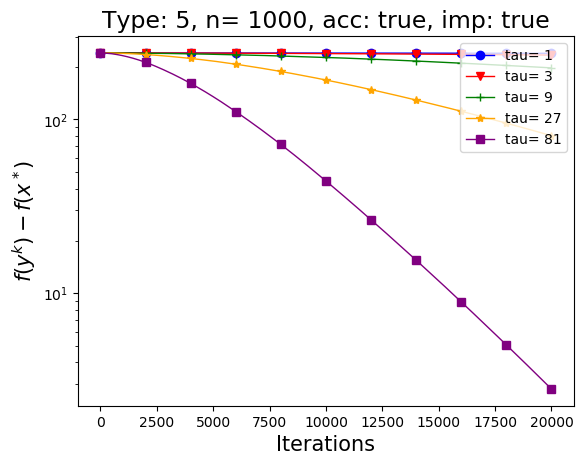}
\end{minipage}%
\\
\caption{Coordinate descent. Comparison of speedup gained by both $\tau$-nice sampling and importance sampling with and without acceleration on various quadratic problems.} \label{fig:acd_tau_comp}
\end{figure}

\subsection{Logistic regression \label{exp:logreg}}

In this section we apply \texttt{ACD} on the regularized logistic regression problem, i.e.
\[
f(x)= \frac1n \sum_{i=1}^n \log \left(1+\exp\left(\mA_{i,:}x\cdot  b\right) \right)+\frac{\lambda}{2} \| x\|^2,
\]
for $b\in \{-1,1 \}$ and data matrix $\mA$ comes from LibSVM. In each experiment in this section, we have chosen regularization parameter $\lambda$ to be the average diagonal element of the smoothness matrix. We first apply the methods with the optimal parameters as our theory suggests on smaller datasets. On larger ones (Section~\ref{sec:acd_practical}), we set them in a cheaper way, which is not guaranteed to work by theory we provide. 

In our first experiment, we apply \texttt{ACD} on LibSVM data directly for various minibatch sizes $\tau$. Figure~\ref{fig:acd_logreg_noncor} shows the results. As expected, \texttt{ACD} is always better to \texttt{CD}, and importance sampling is always better to uniform one. 
 
\begin{figure}[H]
\centering
\begin{minipage}{0.25\textwidth}
  \centering
\includegraphics[width =  \textwidth ]{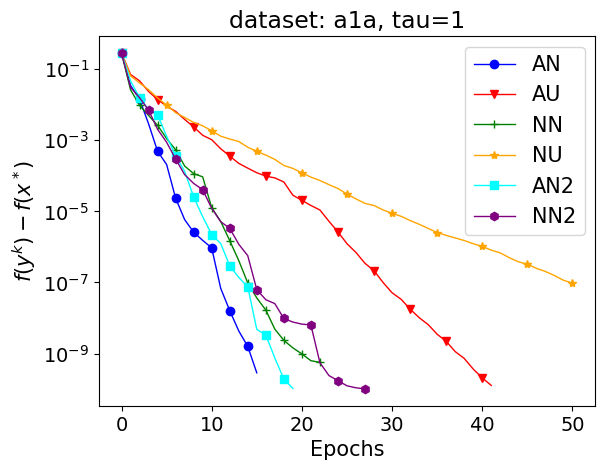}
\end{minipage}%
\begin{minipage}{0.25\textwidth}
  \centering
\includegraphics[width =  \textwidth ]{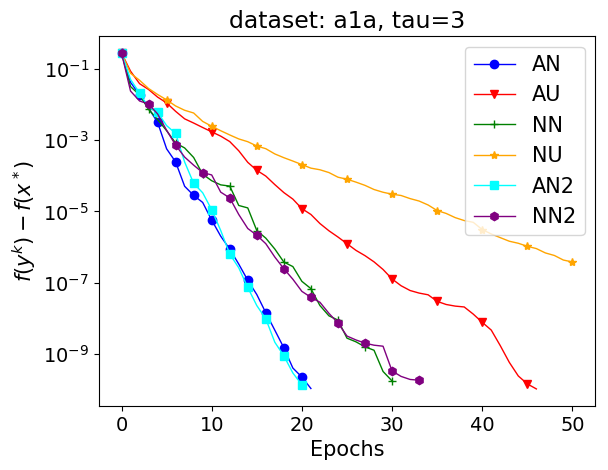}
\end{minipage}%
\begin{minipage}{0.25\textwidth}
  \centering
\includegraphics[width =  \textwidth ]{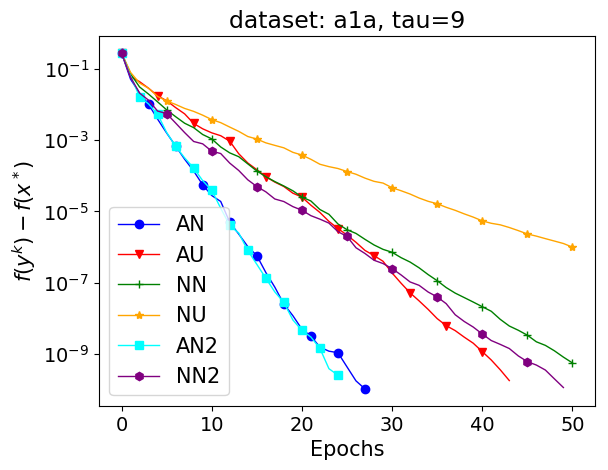}
\end{minipage}%
\begin{minipage}{0.25\textwidth}
  \centering
\includegraphics[width =  \textwidth ]{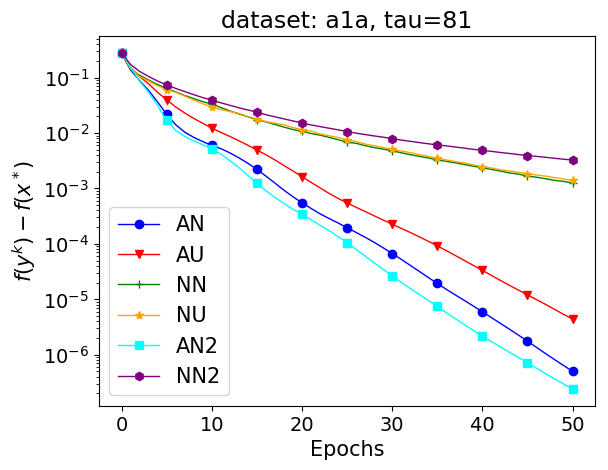}
\end{minipage}%
\\
\begin{minipage}{0.25\textwidth}
  \centering
\includegraphics[width =  \textwidth ]{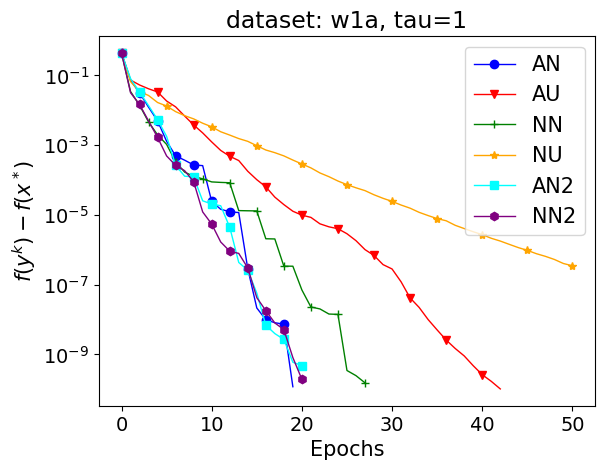}
\end{minipage}%
\begin{minipage}{0.25\textwidth}
  \centering
\includegraphics[width =  \textwidth ]{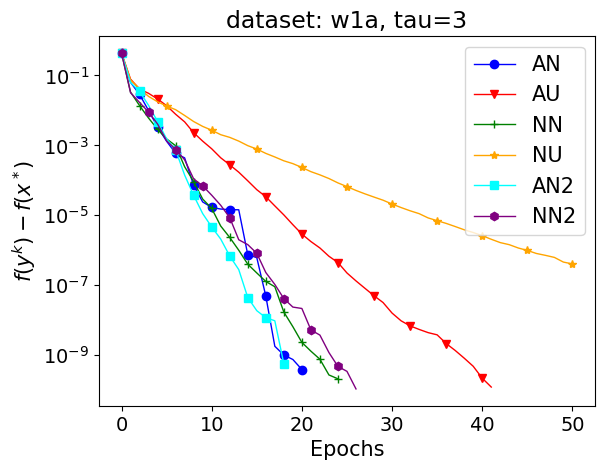}
\end{minipage}%
\begin{minipage}{0.25\textwidth}
  \centering
\includegraphics[width =  \textwidth ]{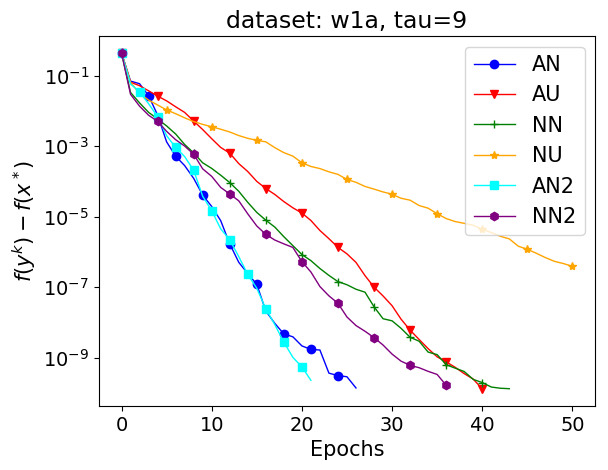}
\end{minipage}%
\begin{minipage}{0.25\textwidth}
  \centering
\includegraphics[width =  \textwidth ]{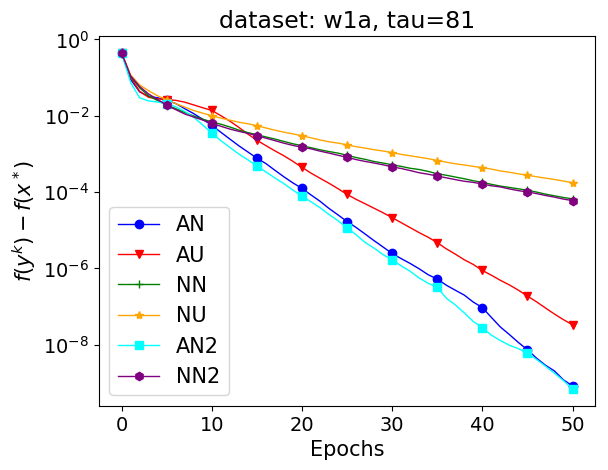}
\end{minipage}%
\\
\begin{minipage}{0.25\textwidth}
  \centering
\includegraphics[width =  \textwidth ]{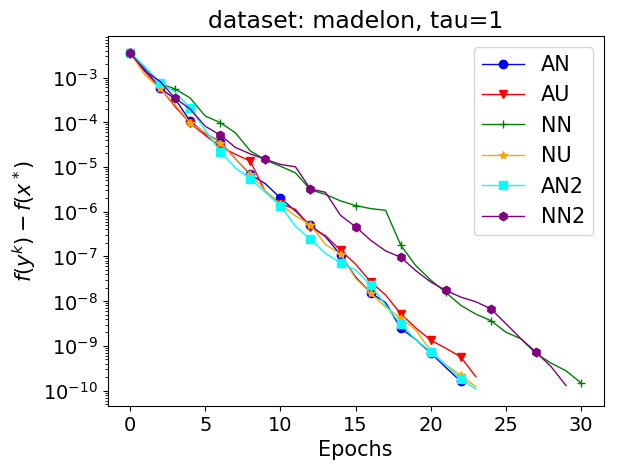}
\end{minipage}%
\begin{minipage}{0.25\textwidth}
  \centering
\includegraphics[width =  \textwidth ]{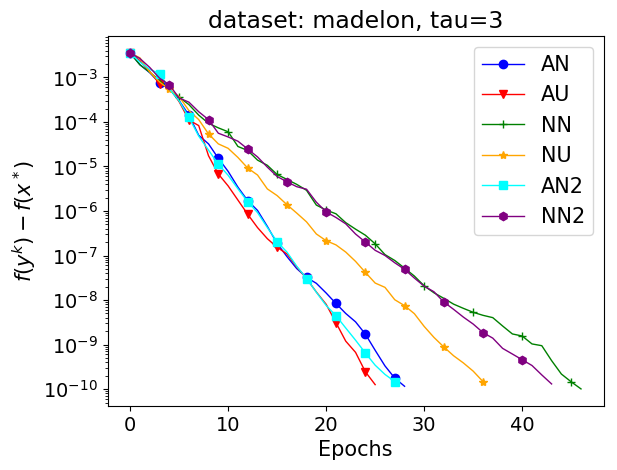}
\end{minipage}%
\begin{minipage}{0.25\textwidth}
  \centering
\includegraphics[width =  \textwidth ]{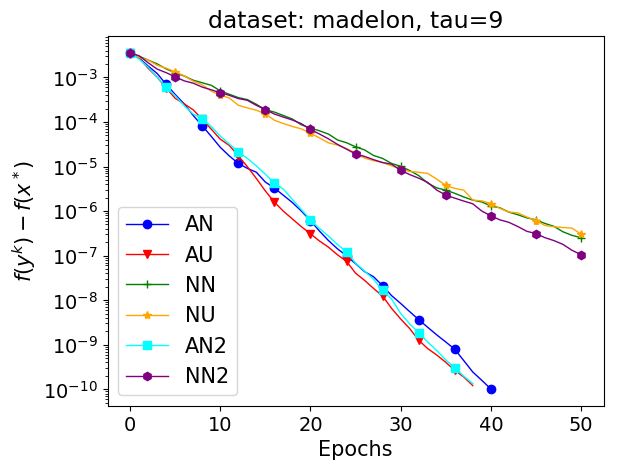}
\end{minipage}%
\begin{minipage}{0.25\textwidth}
  \centering
\includegraphics[width =  \textwidth ]{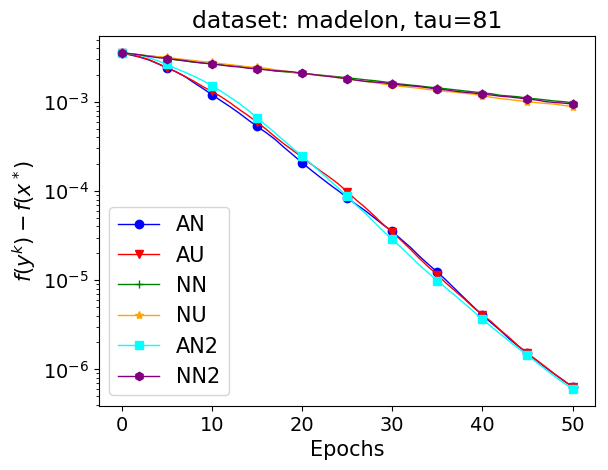}
\end{minipage}%
\\
\centering
\begin{minipage}{0.25\textwidth}
  \centering
\includegraphics[width =  \textwidth ]{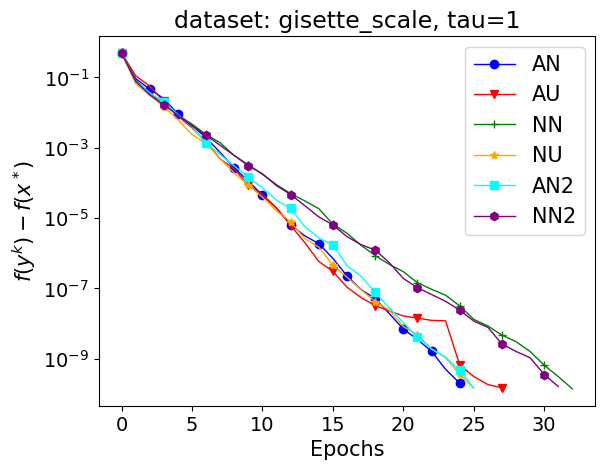}
\end{minipage}%
\begin{minipage}{0.25\textwidth}
  \centering
\includegraphics[width =  \textwidth ]{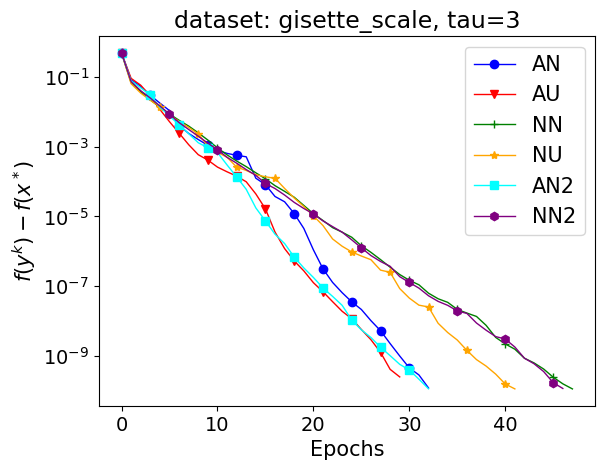}
\end{minipage}%
\begin{minipage}{0.25\textwidth}
  \centering
\includegraphics[width =  \textwidth ]{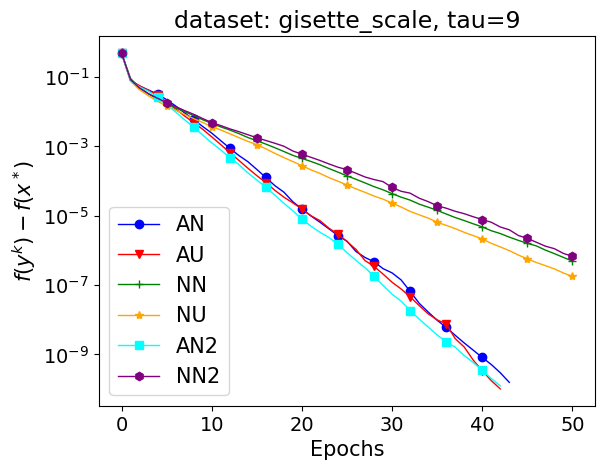}
\end{minipage}%
\begin{minipage}{0.25\textwidth}
  \centering
\includegraphics[width =  \textwidth ]{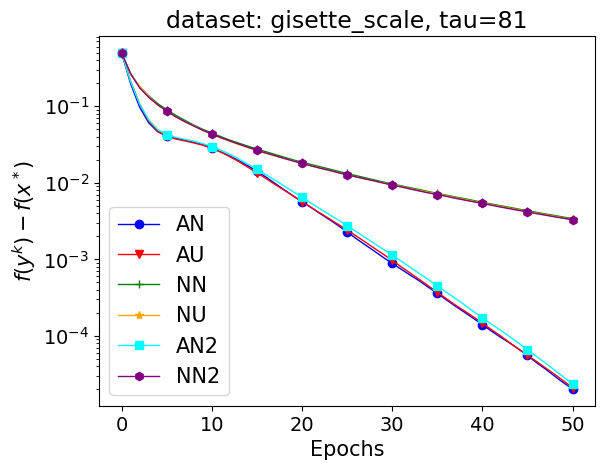}
\end{minipage}%
\\
\begin{minipage}{0.25\textwidth}
  \centering
\includegraphics[width =  \textwidth ]{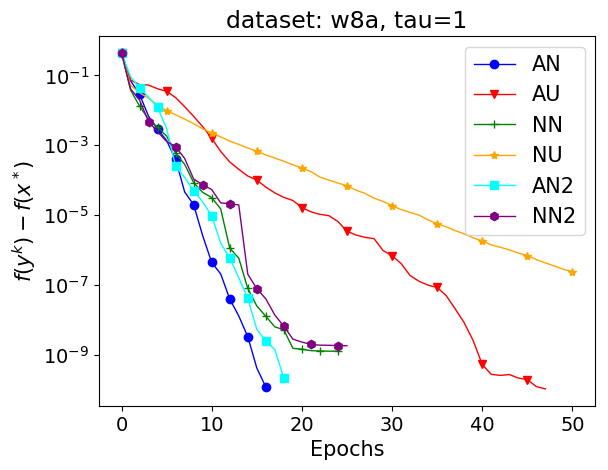}
\end{minipage}%
\begin{minipage}{0.25\textwidth}
  \centering
\includegraphics[width =  \textwidth ]{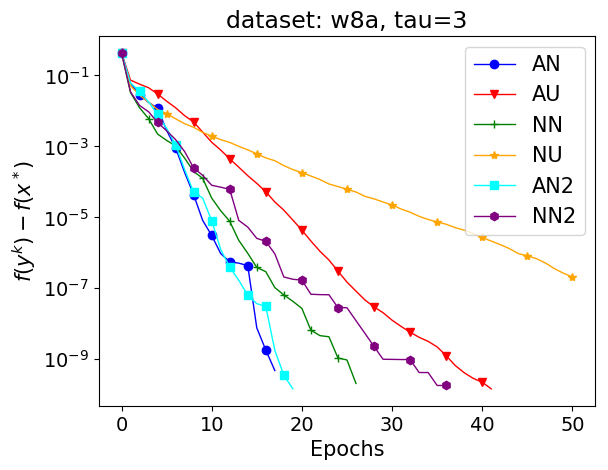}
\end{minipage}%
\begin{minipage}{0.25\textwidth}
  \centering
\includegraphics[width =  \textwidth ]{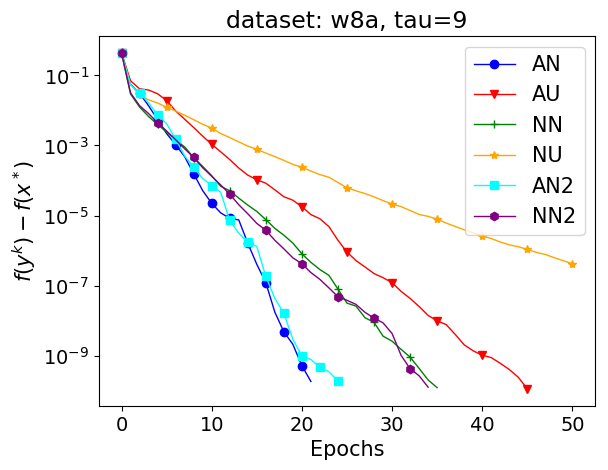}
\end{minipage}%
\begin{minipage}{0.25\textwidth}
  \centering
\includegraphics[width =  \textwidth ]{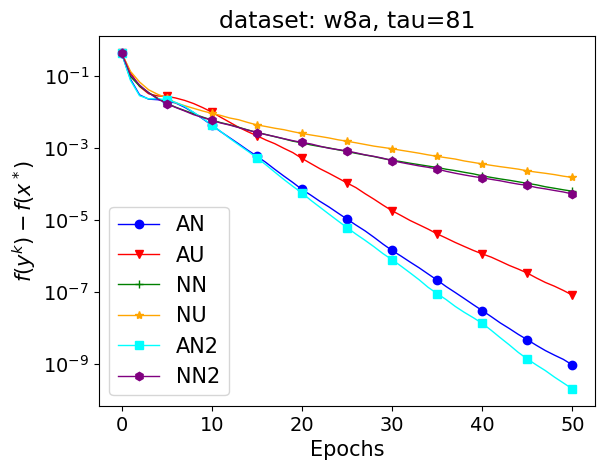}
\end{minipage}%
\caption{Accelerated coordinate desent applied on the logistic regression problem, for various LibSVM datasets and minibatch sizes $\tau$}\label{fig:acd_logreg_noncor}
\end{figure}

Note that, for some datasets and especially bigger minibatch sizes, the effect of importance sampling is sometimes negligible. To demonstrate the power of importance sampling, in the next experiment, we first corrupt the data -- we multiply each row and column of the data matrix $\mA$ by random number from uniform distribution over $[0,1]$. The results can be seen in Figure~\ref{fig:acd_logreg_cor}. As expected, the effect of importance sampling becomes more significant.

\begin{figure}[H]
\centering
\begin{minipage}{0.25\textwidth}
  \centering
\includegraphics[width =  \textwidth ]{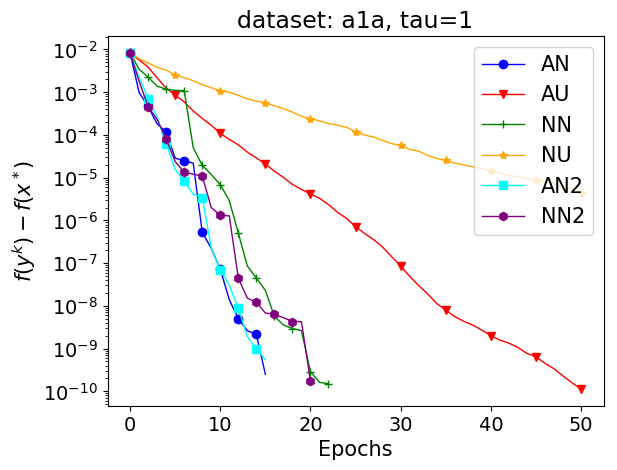}
\end{minipage}%
\begin{minipage}{0.25\textwidth}
  \centering
\includegraphics[width =  \textwidth ]{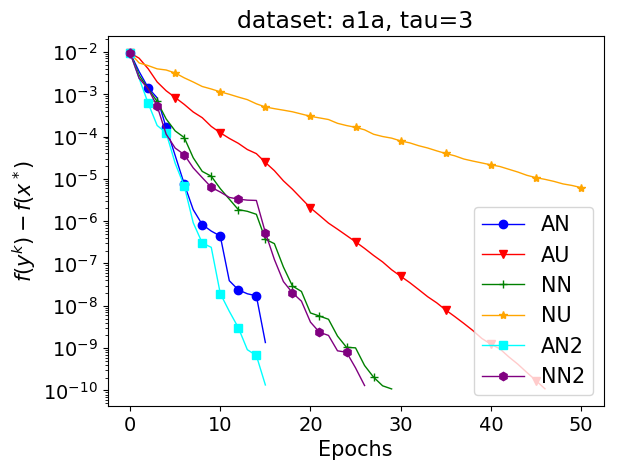}
\end{minipage}%
\begin{minipage}{0.25\textwidth}
  \centering
\includegraphics[width =  \textwidth ]{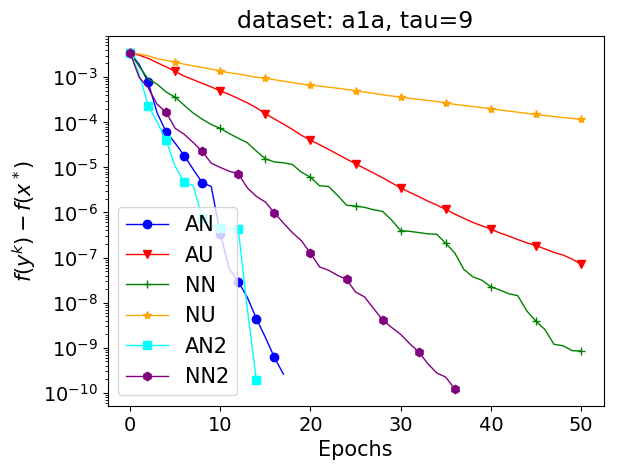}
\end{minipage}%
\begin{minipage}{0.25\textwidth}
  \centering
\includegraphics[width =  \textwidth ]{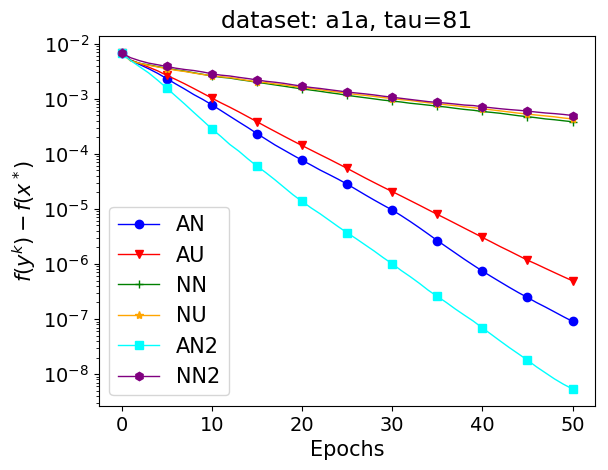}
\end{minipage}%
\\
\begin{minipage}{0.25\textwidth}
  \centering
\includegraphics[width =  \textwidth ]{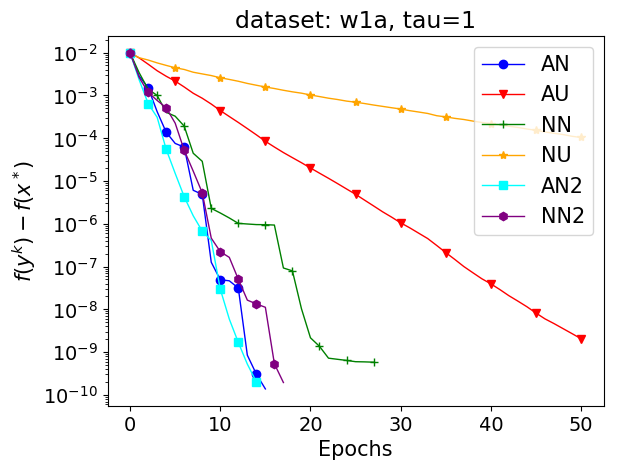}
\end{minipage}%
\begin{minipage}{0.25\textwidth}
  \centering
\includegraphics[width =  \textwidth ]{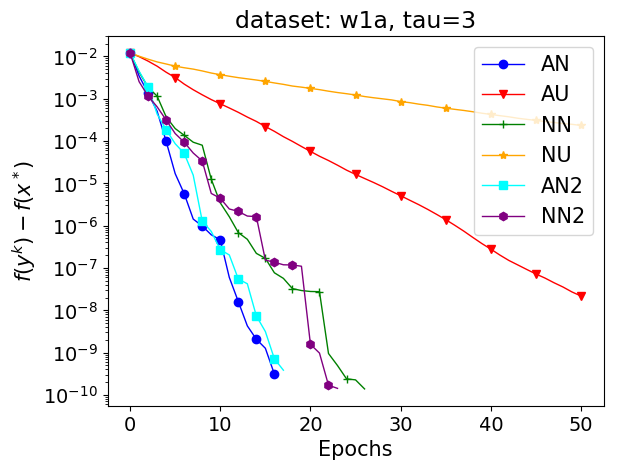}
\end{minipage}%
\begin{minipage}{0.25\textwidth}
  \centering
\includegraphics[width =  \textwidth ]{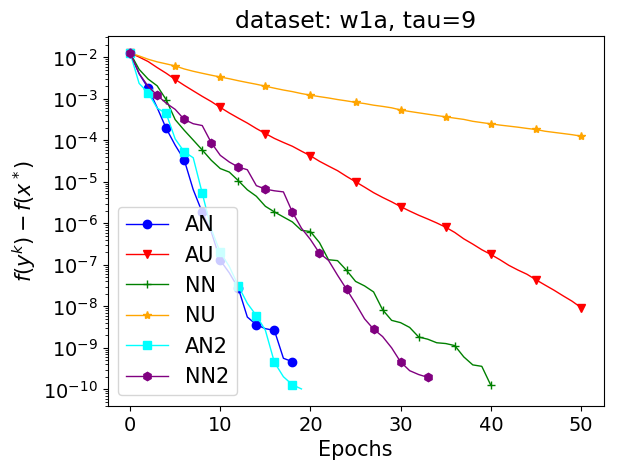}
\end{minipage}%
\begin{minipage}{0.25\textwidth}
  \centering
\includegraphics[width =  \textwidth ]{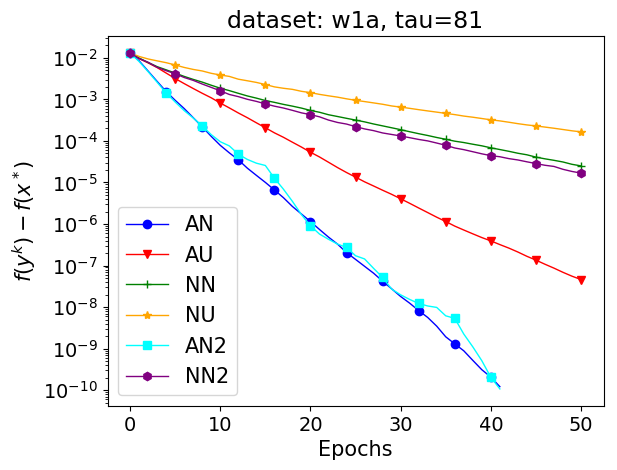}
\end{minipage}%
\\
\begin{minipage}{0.25\textwidth}
  \centering
\includegraphics[width =  \textwidth ]{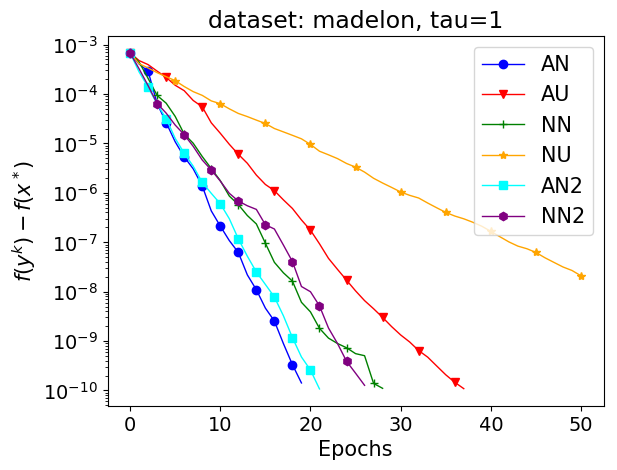}
\end{minipage}%
\begin{minipage}{0.25\textwidth}
  \centering
\includegraphics[width =  \textwidth ]{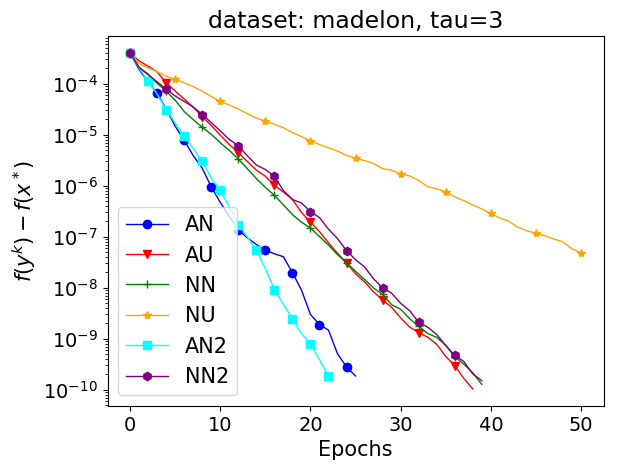}
\end{minipage}%
\begin{minipage}{0.25\textwidth}
  \centering
\includegraphics[width =  \textwidth ]{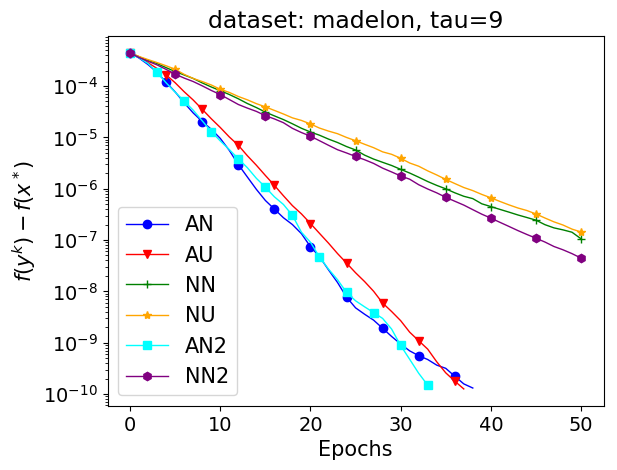}
\end{minipage}%
\begin{minipage}{0.25\textwidth}
  \centering
\includegraphics[width =  \textwidth ]{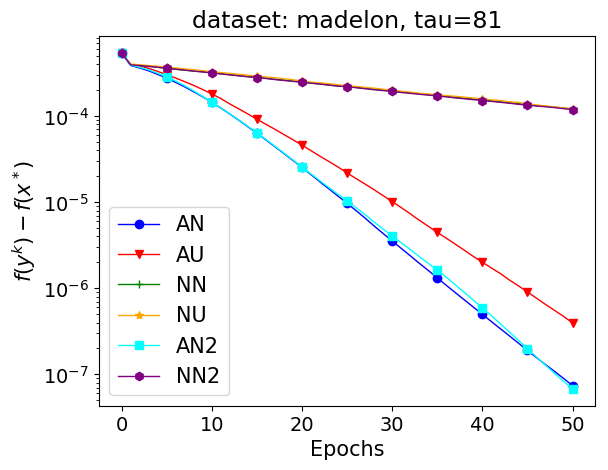}
\end{minipage}%
\\
\centering
\begin{minipage}{0.25\textwidth}
  \centering
\includegraphics[width =  \textwidth ]{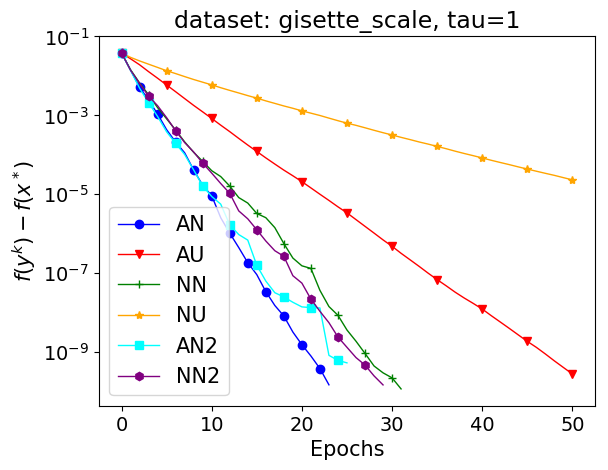}
\end{minipage}%
\begin{minipage}{0.25\textwidth}
  \centering
\includegraphics[width =  \textwidth ]{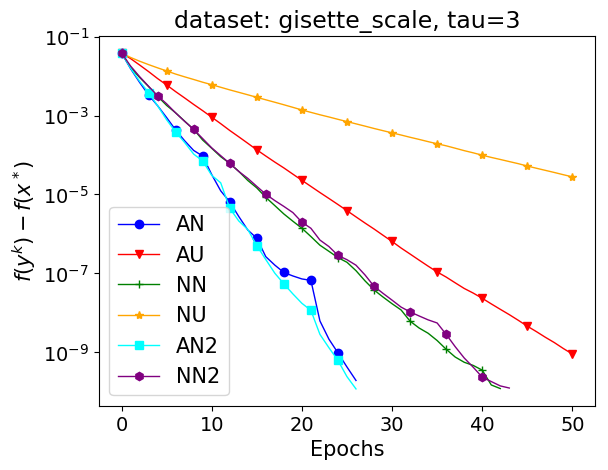}
\end{minipage}%
\begin{minipage}{0.25\textwidth}
  \centering
\includegraphics[width =  \textwidth ]{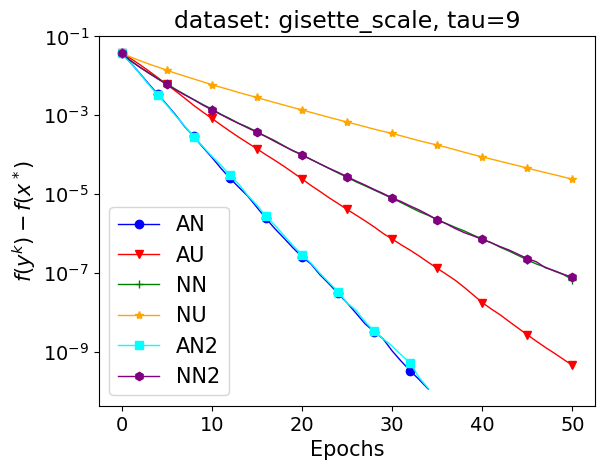}
\end{minipage}%
\begin{minipage}{0.25\textwidth}
  \centering
\includegraphics[width =  \textwidth ]{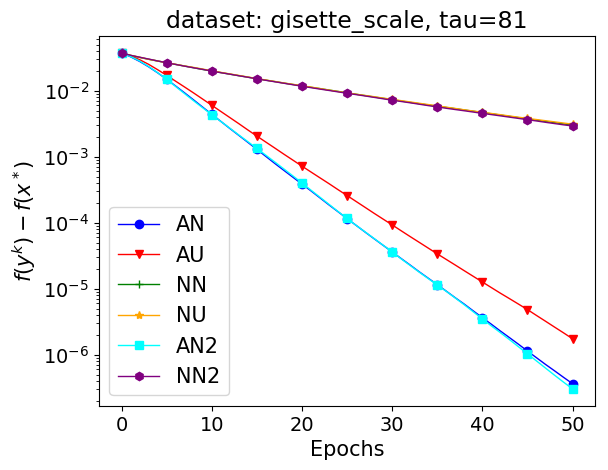}
\end{minipage}%
\\
\begin{minipage}{0.25\textwidth}
  \centering
\includegraphics[width =  \textwidth ]{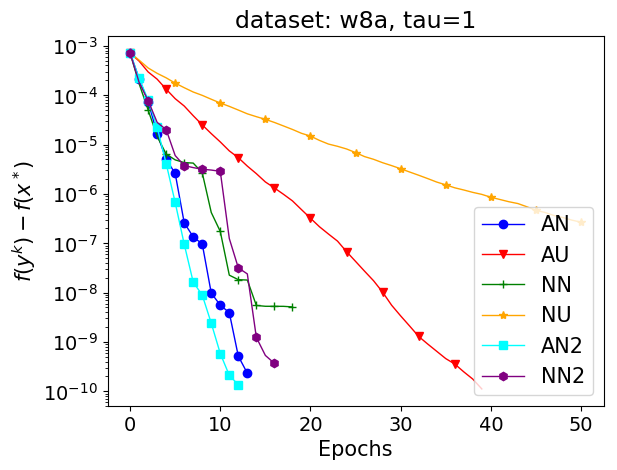}
\end{minipage}%
\begin{minipage}{0.25\textwidth}
  \centering
\includegraphics[width =  \textwidth ]{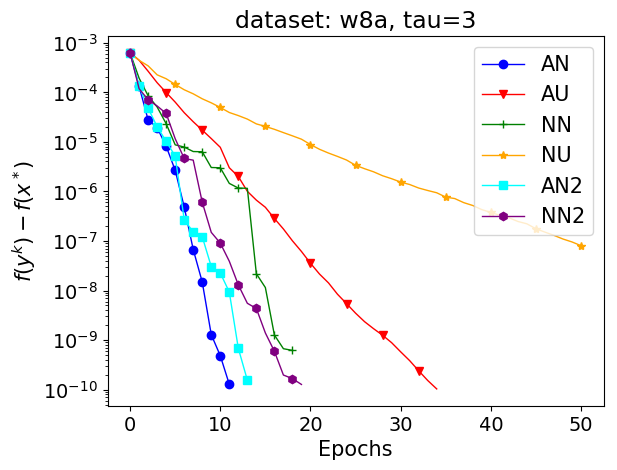}
\end{minipage}%
\begin{minipage}{0.25\textwidth}
  \centering
\includegraphics[width =  \textwidth ]{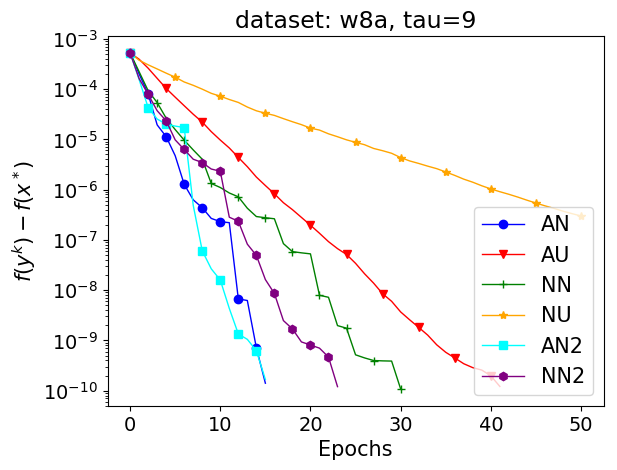}
\end{minipage}%
\begin{minipage}{0.25\textwidth}
  \centering
\includegraphics[width =  \textwidth ]{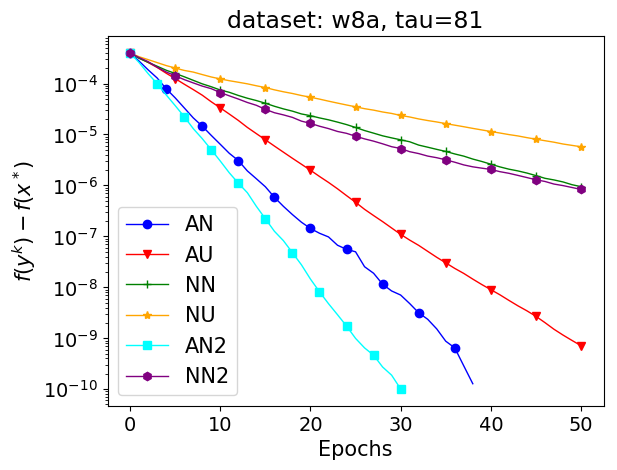}
\end{minipage}%
\caption{{\tt ACD} applied on the logistic regression problem, for various rescaled LibSVM datasets and minibatch sizes $\tau$.}\label{fig:acd_logreg_cor}
\end{figure}

\subsubsection{Practical method on larger dataset \label{sec:acd_practical}}

In Figure~\ref{fig:acd_logreg_big2}, we report on a logistic regression problem with a few selected LibSVM~\cite{chang2011libsvm} datasets.  For larger datasets, pre-computing both strong convexity parameter $\mu$ and  $v$ may be expensive (however, recall that for $v$ we need to tune only one scalar). Therefore, we choose ESO parameters $v$ from Lemma~\ref{thm:acd_special-ESO-result}, while estimating the smoothness matrix as $10\times$ its diagonal. An estimate of the strong convexity $\mu$ for acceleration was chosen to be the minimal diagonal element of the smoothness matrix. We provide a formal formulation of the logistic regression problem, along with  more experiments applied to further datasets in Appendix~\ref{exp:logreg}, where we choose $v$ and $\mu$ in full accord with the theory.

We have chosen regularization parameter $\lambda$ to be the average diagonal element of the smoothness matrix and estimated $v,\mu$ as described in Section~\ref{sec:acd_exp}. 

\begin{figure}[H]
\centering
\begin{minipage}{0.25\textwidth}
  \centering
\includegraphics[width =  \textwidth ]{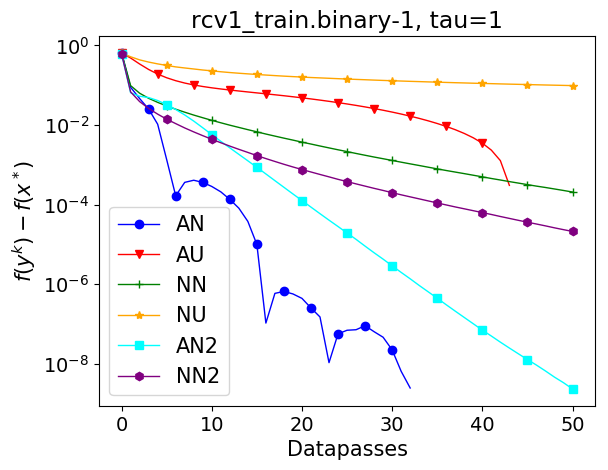}
\end{minipage}%
\begin{minipage}{0.25\textwidth}
  \centering
\includegraphics[width =  \textwidth ]{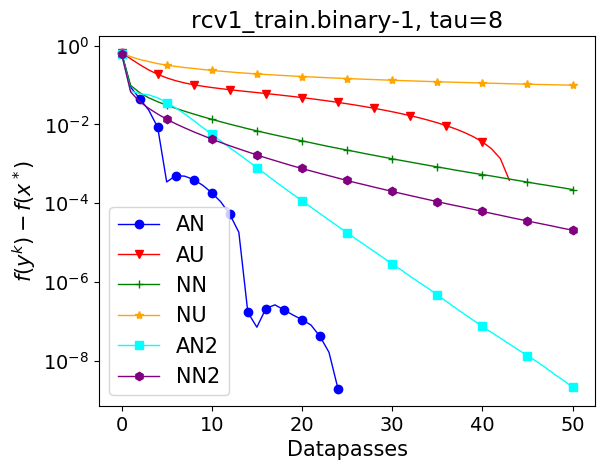}
\end{minipage}%
\begin{minipage}{0.25\textwidth}
  \centering
\includegraphics[width =  \textwidth ]{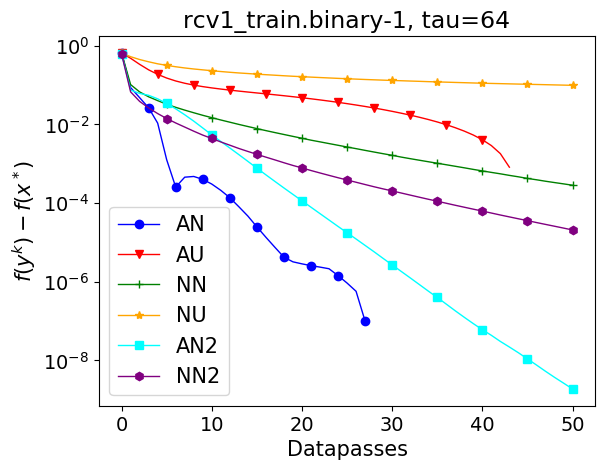}
\end{minipage}%
\begin{minipage}{0.25\textwidth}
  \centering
\includegraphics[width =  \textwidth ]{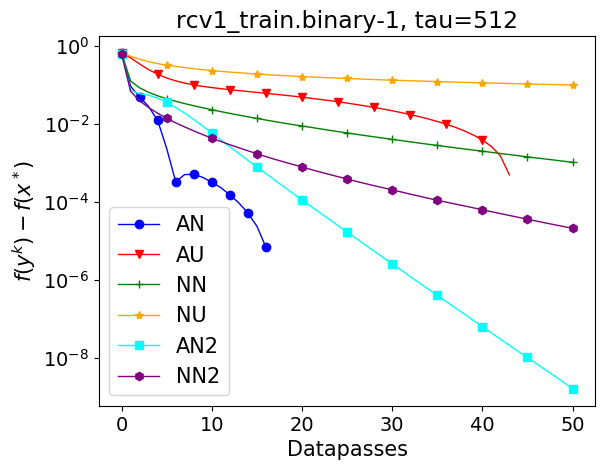}
\end{minipage}%
\\
\begin{minipage}{0.25\textwidth}
  \centering
\includegraphics[width =  \textwidth ]{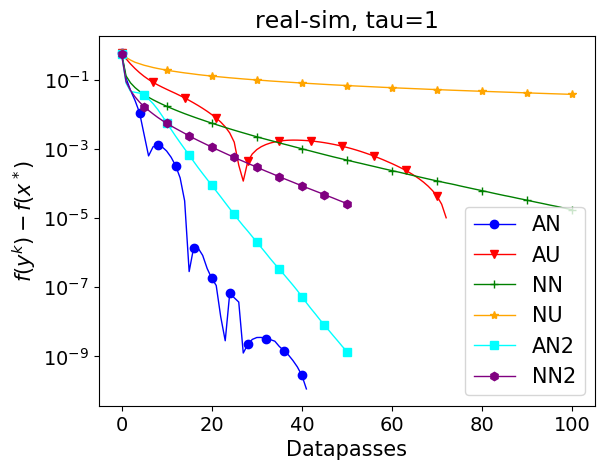}
\end{minipage}%
\begin{minipage}{0.25\textwidth}
  \centering
\includegraphics[width =  \textwidth ]{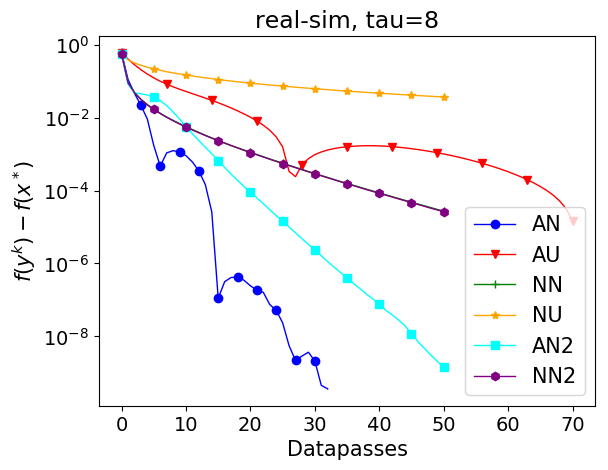}
\end{minipage}%
\begin{minipage}{0.25\textwidth}
  \centering
\includegraphics[width =  \textwidth ]{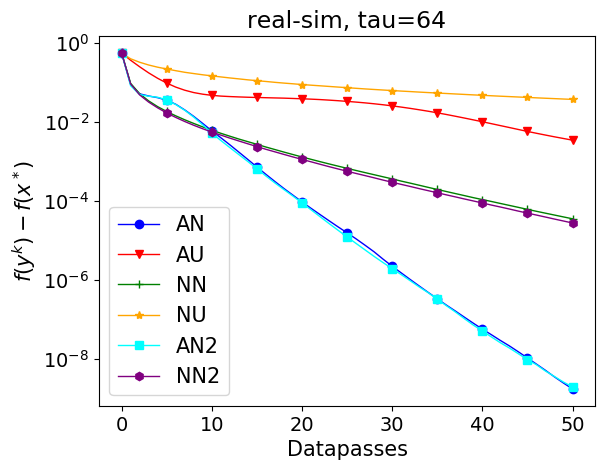}
\end{minipage}%
\begin{minipage}{0.25\textwidth}
  \centering
\includegraphics[width =  \textwidth ]{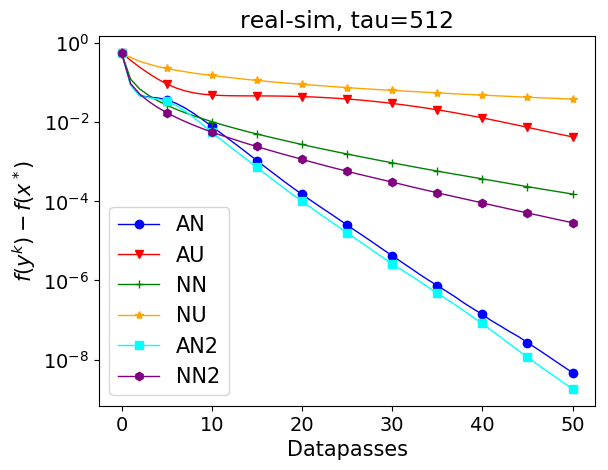}
\end{minipage}%
\\
\begin{minipage}{0.25\textwidth}
  \centering
\includegraphics[width =  \textwidth ]{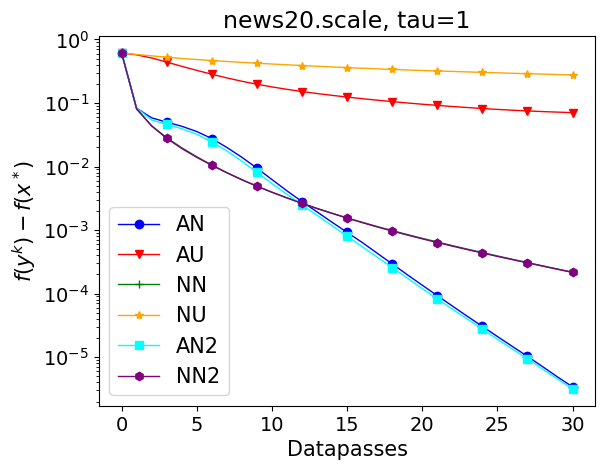}
\end{minipage}%
\begin{minipage}{0.25\textwidth}
  \centering
\includegraphics[width =  \textwidth ]{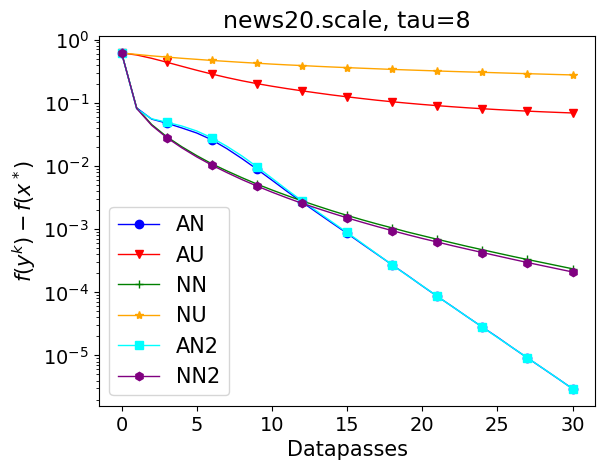}
\end{minipage}%
\begin{minipage}{0.25\textwidth}
  \centering
\includegraphics[width =  \textwidth ]{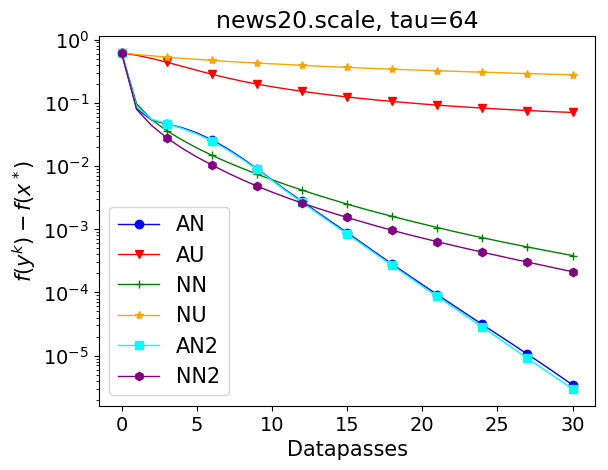}
\end{minipage}%
\begin{minipage}{0.25\textwidth}
  \centering
\includegraphics[width =  \textwidth ]{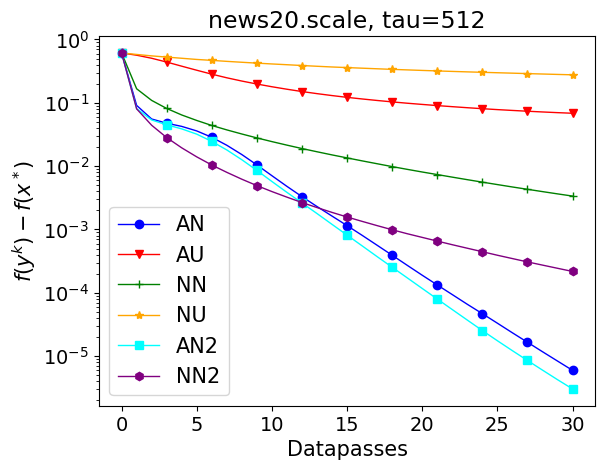}
\end{minipage}%
\caption{Six variants of coordinate descent (\texttt{AN}, \texttt{AU}, \texttt{NN}, \texttt{NU}, \texttt{AN2} and \texttt{AU2})  applied to a logistic regression problem, with minibatch sizes $\tau=1, 8, 64$ and $512$.}\label{fig:acd_logreg_big2}
\end{figure}

\subsection{Support vector machines \label{sec:acd_SVM}}
In this section we apply \texttt{ACD} on the dual of  SVM problem with squared hinge loss, i.e.,
\[
f(x)= \frac{1}{\lambda d^2} \sum_{j=1}^n\left( \sum_{i=1}^d b_i \mA_{ji} x_i\right)^2-\frac1d \sum_{i=1}^d x_i +\frac{1}{4d}\sum_{i=1}^d x_i^2 + \cI_{[0,\infty]}(x),
\]
where $\cI_{[0,\infty]}$ stands for indicator function of set $[0,\infty]$, i.e. $\cI_{[0,\infty]}(x)=0$ if $x\in \R^d_+$, otherwise $\cI_{[0,\infty]}(x)=\infty$. As for the data, we have rescaled each row and each column of the data matrix coming frol LibSVM by random scalar generated from uniform distribution over $[0,1]$. We have chosen regularization parameter $\lambda$ to be maximal diagonal element of the smoothness matrix divided by 10 in each experiment below. We deal with nonsmooth indicator function using proximal operator, which happens to be a projection in this case. We choose ESO parameters $v$ from Lemma~\ref{thm:acd_special-ESO-result}, while estimating the smoothness matrix as $\sqrt{d}$--times multiple of its diagonal. An estimate of the strong convexity $\mu$ for acceleration was chosen to be minimal diagonal element of the smoothness matrix, therefore we adapt a similar approach as in Section~\ref{sec:acd_practical}.

Recall that we did not provide a theory for the proximal steps. However, we make the experiment to demonstrate that \texttt{ACD} can solve big data problems on top of large dimensional problems. Although the results are presented in the main body, we restate them here again (Figure~\ref{fig:acd_SVM}) for the sake of readibility.

\begin{figure}[H]
\centering
\begin{minipage}{0.25\textwidth}
  \centering
\includegraphics[width =  \textwidth ]{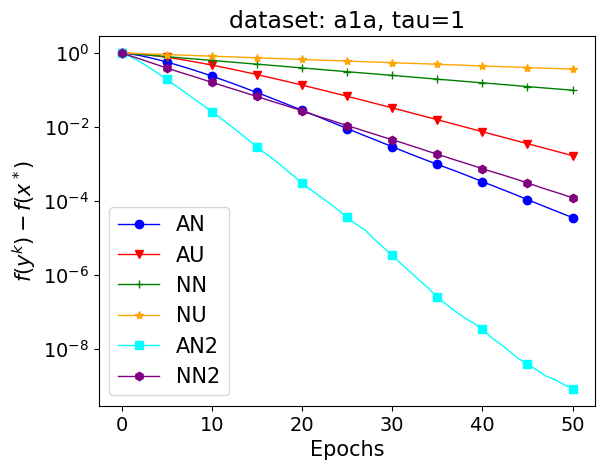}
\end{minipage}%
\begin{minipage}{0.25\textwidth}
  \centering
\includegraphics[width =  \textwidth ]{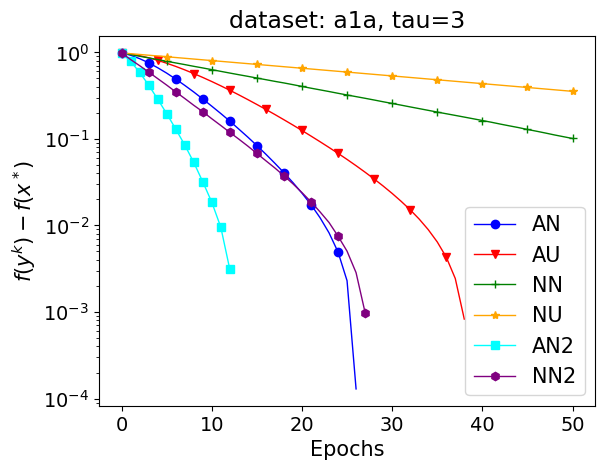}
\end{minipage}%
\begin{minipage}{0.25\textwidth}
  \centering
\includegraphics[width =  \textwidth ]{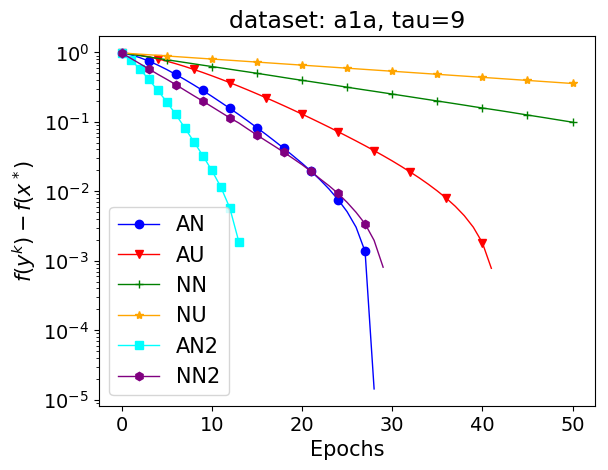}
\end{minipage}%
\begin{minipage}{0.25\textwidth}
  \centering
\includegraphics[width =  \textwidth ]{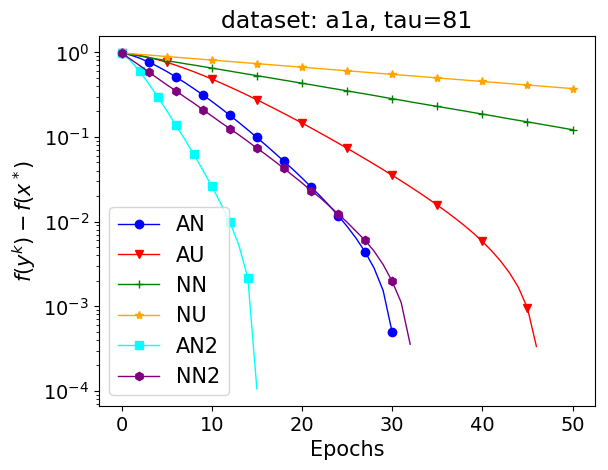}
\end{minipage}%
\\
\begin{minipage}{0.25\textwidth}
  \centering
\includegraphics[width =  \textwidth ]{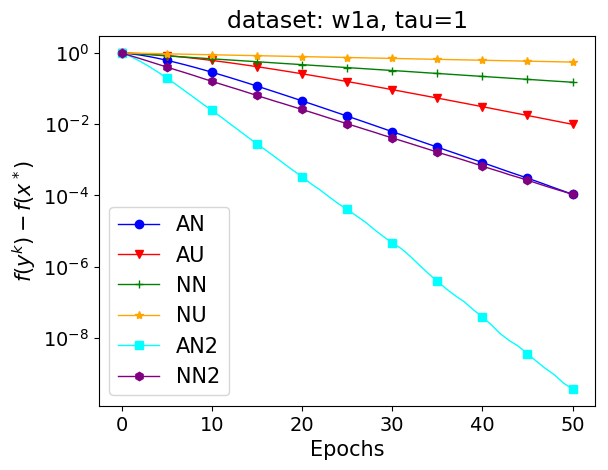}
\end{minipage}%
\begin{minipage}{0.25\textwidth}
  \centering
\includegraphics[width =  \textwidth ]{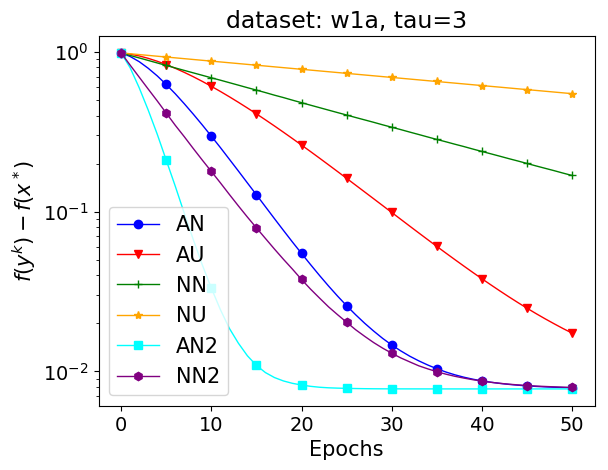}
\end{minipage}%
\begin{minipage}{0.25\textwidth}
  \centering
\includegraphics[width =  \textwidth ]{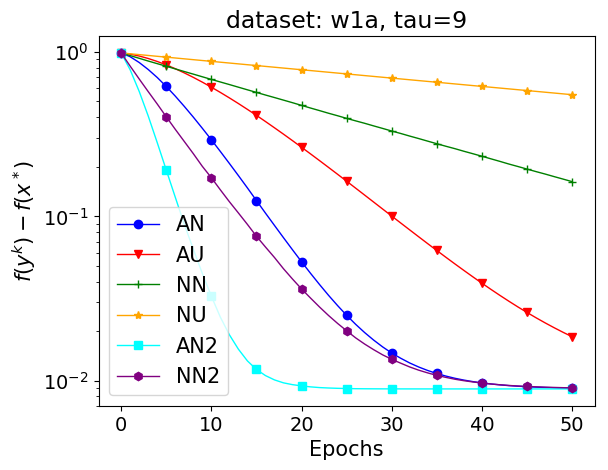}
\end{minipage}%
\begin{minipage}{0.25\textwidth}
  \centering
\includegraphics[width =  \textwidth ]{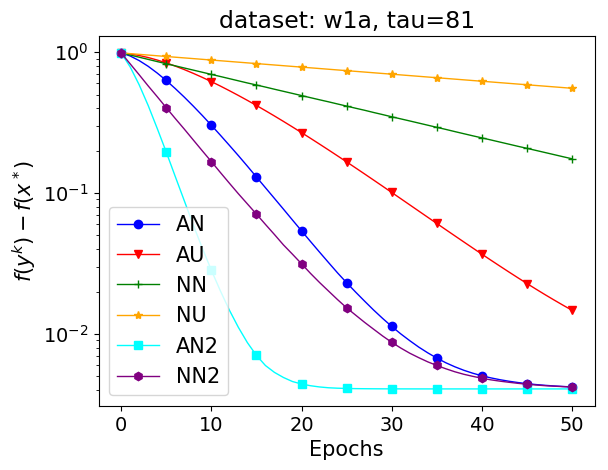}
\end{minipage}%
\\
\begin{minipage}{0.25\textwidth}
  \centering
\includegraphics[width =  \textwidth ]{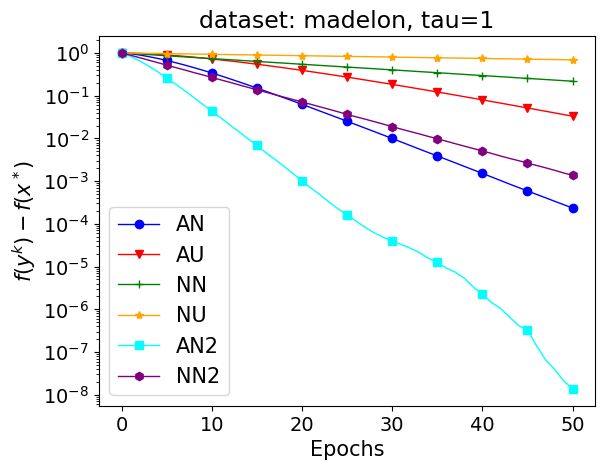}
\end{minipage}%
\begin{minipage}{0.25\textwidth}
  \centering
\includegraphics[width =  \textwidth ]{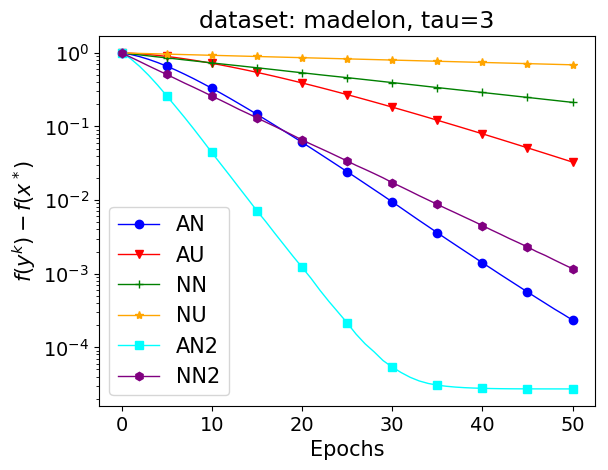}
\end{minipage}%
\begin{minipage}{0.25\textwidth}
  \centering
\includegraphics[width =  \textwidth ]{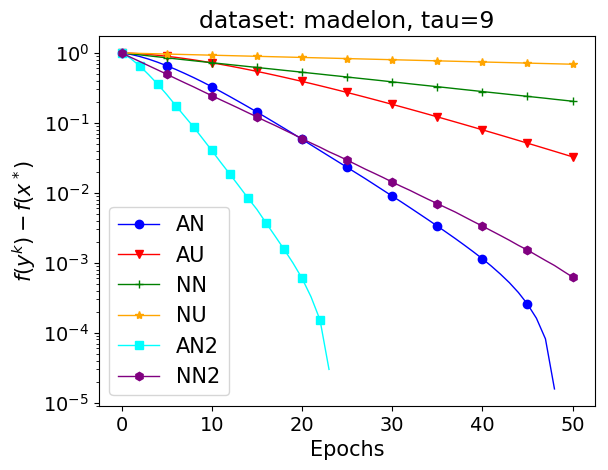}
\end{minipage}%
\begin{minipage}{0.25\textwidth}
  \centering
\includegraphics[width =  \textwidth ]{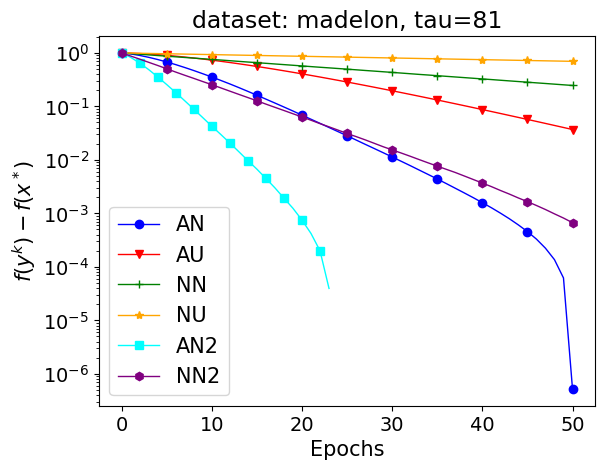}
\end{minipage}%
\\
\centering
\begin{minipage}{0.25\textwidth}
  \centering
\includegraphics[width =  \textwidth ]{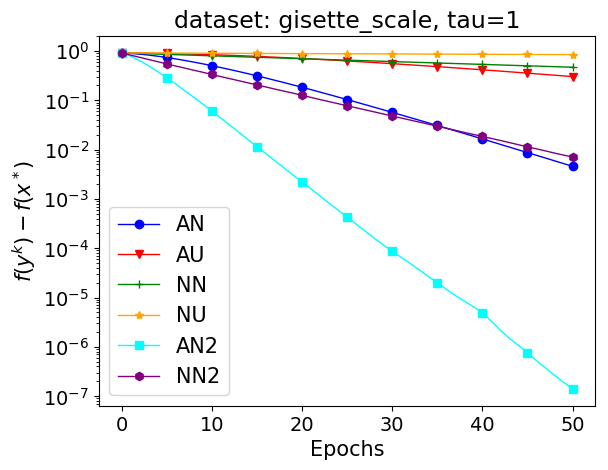}
\end{minipage}%
\begin{minipage}{0.25\textwidth}
  \centering
\includegraphics[width =  \textwidth ]{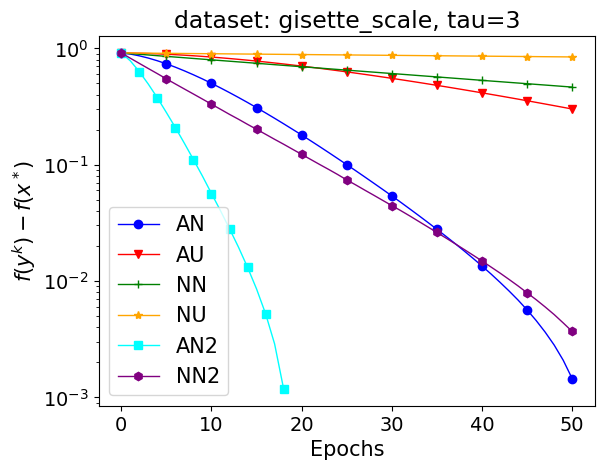}
\end{minipage}%
\begin{minipage}{0.25\textwidth}
  \centering
\includegraphics[width =  \textwidth ]{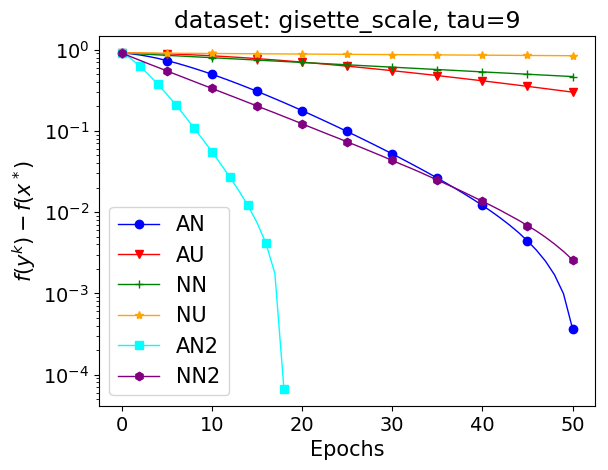}
\end{minipage}%
\begin{minipage}{0.25\textwidth}
  \centering
\includegraphics[width =  \textwidth ]{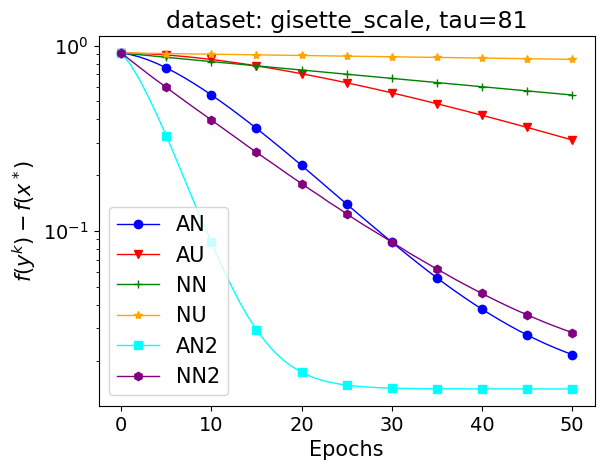}
\end{minipage}%
\caption{Accelerated coordinate desent applied on the dual of of SVM with squared hinge loss, for various LibSVM datasets. } \label{fig:acd_SVM}
\end{figure}

\section{Conclusion}
In this chapter we have presented an minibatch version of accelerated coordinate descent and provided best rates for arbitrary sampling. We have introduced the importance sampling for minibatches, which can be arbitrarily better to uniform sampling, but can be at most constant times worse to uniform sampling. This is the first result of the kind for minibatch coordinate descent samplings. 

As mentioned throughout the chapter, setting of Algorithm~\ref{alg:acd_acd} has a limitation -- it does not allow a minimization with non-separable regularizer using proximal operator. In particular,  objective with non-separable proximal regularizer is not expected to have zero gradient at optimum; and therefore coordinate descent methods can not be expected to converge, unless a decreasing step size is used which leads to significantly slower method. The next chapter solves the issue using variance reduction technique called {\tt SEGA}. 

\chapter{{\tt SEGA}: Variance Reduction via Gradient Sketching}
\label{sega}

\graphicspath{{SEGA/experiments/}}


In  this chapter, we again consider a specific instance of the optimization problem~\eqref{eq:finitesum}. In particular, $f$ is not necessarily assumed to have a  finite-sum structure. However, we allow the presence of a closed convex regularizer $\psi:\R^d\to \R\cup \{+\infty\}$, of which a proximal operator~\eqref{eq:intro_proxdef} is available.  In summary, we aim to solve the following optimization task:
\begin{equation}\label{eq:sega_main_sega} \min_{x\in \R^d} \left\{ F(x)\eqdef f(x) + \psi(x)\right\}.
\end{equation}

 \section{Gradient sketching}
 
The main goal of this chapter is to design provably fast proximal gradient-type methods for solving \eqref{eq:sega_main_sega} {\em without assuming that the true gradient  of $f$ is available.}  Instead, we assume that an oracle provides a {\em random linear transformation (i.e., a sketch) of the gradient}, which is the information available to drive the iterative process. 
 In particular, given a fixed distribution $\cD$  over matrices $\mS \in \R^{d\times \tau}$ ($b\geq 1$ can but does not need to be fixed), and  a query point $x\in \R^d$, our oracle provides us the random linear transformation  of the gradient given by
\begin{equation}
\label{eq:sega_sketched_grad}
\zeta(\mS, x) \eqdef \mS^\top \nabla f(x) \in \R^{\tau}, \qquad \mS \sim \cD.
\end{equation}

Information of this type is available/used in a variety of scenarios. For instance, randomized coordinate descent (\texttt{CD}) methods use oracle~\eqref{eq:sega_sketched_grad}  with $\cD$ corresponding to a distribution over standard  basis vectors. Minibatch/parallel variants of \texttt{CD} methods utilize  oracle \eqref{eq:sega_sketched_grad} with $\cD$ corresponding to a distribution over random column submatrices of the identity matrix. If one is prepared to use difference of function values to approximate  directional derivatives, then one can apply our oracle model  to   zeroth-order optimization~\cite{conn2009introduction}.  Indeed, the directional derivative of $f$ in a random direction $\mS=s \in \R^{d\times 1}$ can be approximated by $\zeta(s, x)  \approx \frac{1}{\epsilon}(f(x+ \epsilon s) - f(x))$, where $\epsilon>0$ is sufficiently small.

\begin{example}[Sketches]  We now illustrate this concept using two examples.
\begin{itemize}
\item[(i)]  {\bf Coordinate sketch.} Let $\cD$ be the uniform distribution over standard unit basis vectors $e_1,e_2,\dots,e_d$ of $\R^d$. Then
$\zeta(e_i,x) = e_i^\top \nabla f(x)$, i.e., the $i^{\text{th}}$ {\em partial derivative} of $f$ at $x$. 
\item[(ii)]{\bf  Gaussian sketch.} Let $\cD$ be the standard Gaussian distribution in $\R^d$. Then for $s\sim \cD$ we have
$\zeta(s, x) = s^\top \nabla f(x) $, i.e., the {\em directional derivative} of $f$ at $x$ in direction $s$.
\end{itemize}
\end{example}

We describe \texttt{SEGA} in Section~\ref{sec:sega_SEGA}. Convergence results for general sketches are described in Section~\ref{sec:sega_analysis}. Refined results for coordinate sketches are presented in Section~\ref{sec:sega_CD}, where we also describe and analyze an accelerated variant of \texttt{SEGA}. Experimental results can be found in Section~\ref{sec:sega_experiments}. We also include here experiments with a {\em subspace} variant of \texttt{SEGA}, which is described and analyzed in Appendix~\ref{sec:sega_subSEGA}. Conclusions are drawn and potential extensions outlined in Section~\ref{sec:sega_conclusion}. A simplified analysis of \texttt{SEGA} in the case of coordinate sketches and for $\psi\equiv 0$ is developed in Appendix~\ref{sec:sega_simple_SEGA} (under standard assumptions as in the main body). 

We introduce notation when and where needed. For convenience, we provide a table of frequently used notation in Appendix~\ref{sec:table}.

\subsection{Related work}

In the last decade, stochastic gradient-type methods for solving problem~\eqref{eq:sega_main_sega} have received unprecedented attention by theoreticians and practitioners alike. Specific examples of such methods  are stochastic gradient descent (\texttt{SGD})~\cite{robbins}, variance-reduced variants of \texttt{SGD} such as \texttt{SAG}~\cite{sag}, \texttt{SAGA}~\cite{saga}, \texttt{SVRG}~\cite{svrg}, and their accelerated counterparts~\cite{lin2015universal, allen2017katyusha}. While these methods are specifically designed for objectives formulated as an expectation or a finite sum, we do not assume such a structure. Moreover, these methods  utilize a fundamentally different stochastic gradient information: they have access to an unbiased estimator of the gradient. In contrast, we do not assume that \eqref{eq:sega_sketched_grad} is an unbiased estimator of $\nabla f(x)$. In fact,  $\zeta(\mS, x)\in \R^{\tau}$ and $\nabla f(x)\in \R^d$ do not even necessarily  belong to the same space.  Therefore, our algorithms and results should be seen as complementary to the above line of research.

While the gradient sketch $\zeta(\mS, x)$  does not immediatey lead to an unbiased estimator of the gradient,  \texttt{SEGA} uses  the information provided in the sketch  to {\em construct} an unbiased estimator of the gradient  via a {\em sketch-and-project} process. Sketch-and-project iterations were introduced in \cite{gower2015randomized} in the contex of linear feasibility problems. A dual view  uncovering a direct relationship with  stochastic subspace ascent methods was developed in \cite{sda}. The latest and most in-depth treatment of sketch-and-project for linear feasibility is based on the idea of stochastic reformulations  \cite{richtarik2017stochastic}. Sketch-and-project can be combined with Polyak~\cite{SMOMENTUM, SHB-NIPS} and Nesterov momentum~\cite{gower2018accelerated}, extended to convex~ feasibility problems \cite{necoara2019randomized}, matrix inversion \cite{gower:2017, pseudoinverse, gower2018accelerated}, and empirical risk minimization \cite{sbfgs, jacsketch}. Connections to gossip algorithms for average consensus  were made in \cite{new-perspective, agossip}.

The line of work most closely related to our setup is that on randomized coordinate/subspace descent methods~\cite{rcdm, sda}. Indeed, the information available to these methods is compatible with our oracle for specific distributions $\cD$. However, the main disadvantage of these methods is that they are not able to handle non-separable regularizers $\psi$. In contrast, the algorithm we propose---\texttt{SEGA}---works for any regularizer $\psi$.    In particular,  \texttt{SEGA} can handle non-separable constraints even with coordinate sketches, which is out of range of current coordinate descent methods. Hence, our work could be understood as extending the reach of coordinate and subspace descent methods from separable to arbitrary regularizers, which allows for a plethora of new applications.  Our method is able to work with an arbitrary regularizer due to its ability to {\em build an unbiased variance-reduced estimate of the gradient} of $f$ throughout the iterative process from the random linear measurements thereof provided by the oracle.  Moreover, and unlike coordinate descent,  \texttt{SEGA} allows for general sketches from essentially any distribution $\cD$. 

Another stream of work on designing gradient-type methods without assuming perfect access to the gradient is represented by the {\em inexact gradient descent} methods~\cite{d2008smooth,devolder:2011inexact,schmidt2011convergence}. However, these methods deal with deterministic estimates of the gradient and are not based on linear transformations of the gradient. Therefore, this second line of research is also significantly different from what we do here. 


\section{Contributions}
We now list the main contributions of this chapter.
 
\begin{itemize}
\item
 \textbf{Subspace oracle with non-separable regularizer.} \texttt{SEGA} is the first iterative proximal algorithm with a subspace gradient oracle that achieves linear convergence. Unlike coordinate descent,  \texttt{SEGA} does not require the regularizer to be separable and thus has a much broader range of applications. It achieves by constructing control variance to progressively reduce the variance of stochastic gradient estimator.  

\item
 \textbf{Generality and Subspace \texttt{SEGA} .} We provide the convergence rate of \texttt{SEGA} under the full generality -- we allow for arbitrary distribution of sketching matrices $\mS$. In some scenarios, this might lead to a very fast convergence, especially when $\nabla f$ always belongs to a particular subspace.

\item
 \textbf{Fast rates without $\psi$.} Given that $\psi \equiv 0$, we show that \texttt{SEGA} is, up to a small constant, as fast as the state-of-the-art coordinate descent. Specifically, we show that {\tt SEGA} with importance sampling and acceleration converges, up to a constant, as fast as the analogous version of {\tt CD}.
\end{itemize}

\section{The \texttt{SEGA} algorithm} \label{sec:sega_SEGA}

In this section we introduce a learning process for estimating the gradient from the sketched information provided by \eqref{eq:sega_sketched_grad}; this will be used as a subroutine of \texttt{SEGA}.

Let $x^k$ be the current iterate, and let $h^k$ be the current estimate of the gradient of $f$. We then query the oracle, and receive new gradient information in the form of the sketched gradient \eqref{eq:sega_sketched_grad}. At this point, we would like to update $h^k$ based on this new information. We do this using a {\em sketch-and-project} process~\cite{gower2015randomized, sda, richtarik2017stochastic}: we set $h^{k+1}$ to be the closest vector to $h^k$ satisfying~\eqref{eq:sega_sketched_grad}:
\begin{eqnarray}
h^{k+1} &=& \arg \min_{h\in \R^{d}} \| h -  h^k\|^2 \notag \\
&& \text{subject to} \quad \mS_{k}^\top h =  \mS_{k}^\top \nabla f(x^k). \label{eq:sega_sketch-n-project}
\end{eqnarray}

The closed-form solution of \eqref{eq:sega_sketch-n-project} is
\begin{equation} h^{k+1} = h^k - \mZ_{k} (h^k - \nabla f(x^k)) = (\mI-\mZ_{k} ) h^k +\mZ_{k} \nabla f(x^k),\label{eq:sega_h^{k+1}}
\end{equation}
where $\mZ_{k} \eqdef  \mS_{k} \left(\mS_{k}^\top  \mS_{k}\right)^\dagger\mS_{k}^\top$. Notice that $h^{k+1}$ is a \emph{biased} estimator of $\nabla f(x^k)$. In order to obtain an unbiased gradient estimator, we introduce a random variable\footnote{Such a random variable may not exist. Some sufficient conditions are provided later.} $\theta_k=\theta(\mS_{k})$ for which
\begin{equation} \label{eq:sega_unbiased} 
\E{\theta_k  \mZ_{k}} = \mI.
\end{equation}

If $\theta_k$ satisfies~\eqref{eq:sega_unbiased}, it is straightforward to see that the random vector
\begin{equation} \label{eq:sega_g^k} g^k \eqdef (1-\theta_k) h^k + \theta_k h^{k+1} \overset{\eqref{eq:sega_h^{k+1}}}{=} h^k + \theta_k\mZ_{k} (\nabla f(x^k) - h^k)
\end{equation}
is an {\em unbiased estimator} of the gradient:
\begin{eqnarray}
\E{g^k} &\overset{\eqref{eq:sega_unbiased} +\eqref{eq:sega_g^k}}{=}&  \nabla f(x^k). \label{eq:sega_unbiased_estimator}
\end{eqnarray}

Finally, we use $g^k$ instead of the true gradient, and perform a proximal step with respect to $\psi$. This leads to a new randomized optimization method, which we call {\em SkEtched Gradient Algorithm (\texttt{SEGA})}. The method is formally described in Algorithm~\ref{alg:sega_gs}. We stress again that the method does not need the access to the full gradient. 

\begin{algorithm}[!h]
\begin{algorithmic}[1]
\State \textbf{Parameters:} $x^0, h^0\in \R^d$;  distribution $\cD$; stepsize $\alpha>0$ 
\For {$k= 0,1,2,\dots $} 
 \State Sample $\mS_{k} \sim \cD$
   \State$g^{k} = h^k + \theta_k  \mZ_{k} (\nabla f(x^k) - h^k)$ \label{eq:sega_g_update}
  \State$x^{k+1} = \prox)_{\alpha \psi}(x^k - \alpha g^k) $  \label{eq:sega_x_update} 
    \State$h^{k+1} = h^k +\mZ_{k} (\nabla f(x^k) - h^k) $\label{eq:sega_h_update}
    \EndFor
   \end{algorithmic}
\caption{\texttt{SEGA} (SkEtched Gradient Algorithm)}
 \label{alg:sega_gs}
\end{algorithm}

\begin{figure}
\centering
\includegraphics[width = 0.4\textwidth ]{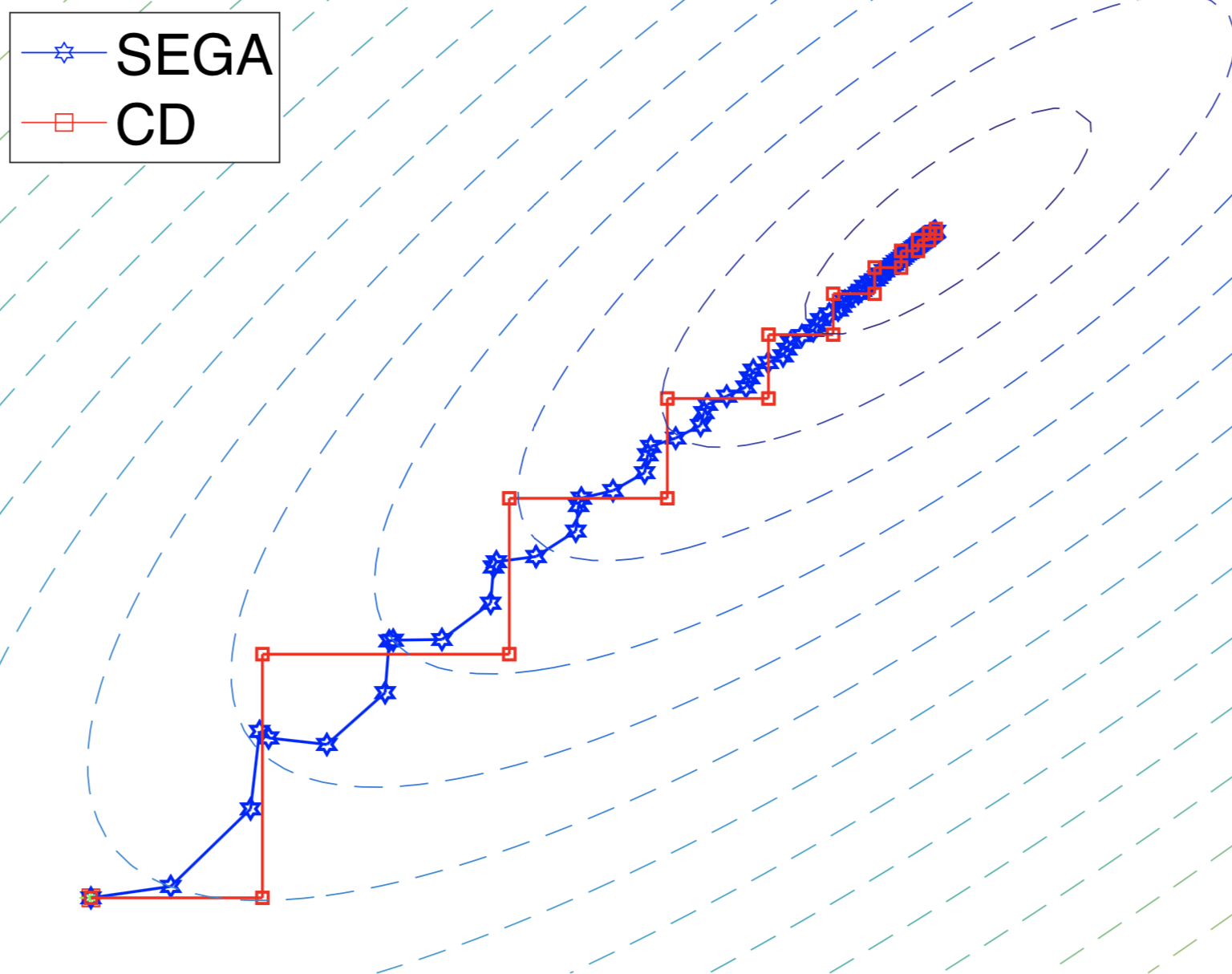}
\caption{Iterates of \texttt{SEGA} and \texttt{CD}}\label{fig:sega_trajectory}
\end{figure}

\subsection{\texttt{SEGA} as a variance-reduced method}

As we shall show, both $h^k$ and $g^k$ are becoming better at approximating $\nabla f(x^k)$ as the iterates $x^k$ approach the optimum. Hence, the variance of $g^k$ as an estimator of the gradient tends to zero, which means that  \texttt{SEGA} is a \emph{variance-reduced} algorithm. The structure of \texttt{SEGA} is inspired by the \texttt{JackSketch} algorithm  introduced in~\cite{jacsketch}. However,  as \texttt{JackSketch} is aimed at solving a finite-sum optimization problem with many components, it does not make much sense to apply it to \eqref{eq:sega_main_sega}. Indeed, when applied to \eqref{eq:sega_main_sega} (with $\psi\equiv0$, since  \texttt{JackSketch} was analyzed for smooth  optimization only), \texttt{JackSketch}  reduces to gradient descent. While \texttt{JackSketch} performs {\em Jacobian} sketching (i.e., multiplying the Jacobian by a random matrix from the right, effectively sampling a subset of the gradients forming the finite sum), \texttt{SEGA} multiplies the Jacobian by a random matrix from the left. In doing so, \texttt{SEGA} becomes oblivious to the finite-sum structure and transforms into the gradient sketching mechanism described in~\eqref{eq:sega_sketched_grad}.



\subsection{\texttt{SEGA} versus coordinate descent}

We now illustrate the above general setup on the simple example when $\cD$ corresponds to a distribution over standard unit basis vectors in $\R^d$. 

\begin{example}\label{ex:coord_setup} Let $\cD$ be defined as follows. We choose $\mS_{k} = e_{i}$ with probability $p_i>0$, where  $e_1,e_2,\dots, e_d$ are the unit basis vectors in $\R^d$. Then
 \begin{equation} \label{eq:sega_988fgf} h^{k+1} \overset{\eqref{eq:sega_h^{k+1}}}{=} h^k + e_{i}^\top (\nabla f(x^k) - h^k) e_{i},\end{equation}
which can equivalently be written as $h^{k+1}_i = e_{i}^\top \nabla f(x^k)$ and $h^{k+1}_j = h^k_j$ for $j\neq i$. If we choose $\theta_k=\theta(\mS_{k}) = 1/p_i$, then 
\[\E{\theta_k \mZ_{k}} = \sum_{i=1}^d p_i \frac{1}{p_i} e_i (e_i^\top e_i)^{-1} e_i^\top = \sum_{i=1}^d e_i e_i^\top = \mI,\]
which means that $\theta_k$ is a bias-correcting random variable. We then get
\begin{equation} \label{eq:sega_8h0h09ffs}g^k \overset{\eqref{eq:sega_g^k} }{=}   h^k + \frac{1}{p_{i}}  e_{i}^\top (\nabla f(x^k) - h^k) e_{i} . \end{equation}
\end{example}

In the setup of Example~\ref{ex:coord_setup}, both \texttt{SEGA} and \texttt{CD}  obtain  new gradient information in the form of a random partial derivative of $f$. However, the two methods process this information differently, and perform a different update: 
\begin{itemize}
\item[(i)]  While \texttt{SEGA} allows for arbitrary proximal term, \texttt{CD} allows for separable proximal term only~\cite{proxsdca,lin2014accelerated,approx}.  
\item[(ii)] While \texttt{SEGA} updates all coordinates in every iteration,  \texttt{CD} updates  a single coordinate only. 
\item[(iii)] If we force $h^k=0$ in  \texttt{SEGA} and use coordinate sketches, the method transforms  into \texttt{CD}.
\end{itemize}

 Based on the above observations, we conclude that \texttt{SEGA} can be applied in more general settings for the price of potentially more expensive iterations\footnote{Forming vector $g$ and computing the prox.}.  For intuition-building illustration of how \texttt{SEGA} works, Figure~\ref{fig:sega_trajectory} shows  the evolution of iterates of both \texttt{SEGA} and \texttt{CD} applied to minimizing a simple quadratic function in 2 dimensions. For more figures of this type, including the composite case where \texttt{CD} does not work, see Appendix~\ref{sec:sega_evolution_extra}.

In Section~\ref{sec:sega_CD} we show that \texttt{SEGA} enjoys the same theoretical iteration complexity rates as \texttt{CD}, up to a small constant factor. This remains true when comparing state-of-the-art variants of \texttt{CD} utilizing importance-sampling, parallelism/mini-batching and acceleration with the appropriate corresponding  variants of \texttt{SEGA}.

\begin{remark}
Nontrivial sketches $\mS$ might, in some applications, bring a substantial speedup against the baseline choices mentioned in Example~\ref{ex:coord_setup}. Appendix~\ref{sec:sega_subSEGA} provides one setting where this can happen: there are problems where the gradient of $f$ always lies in a particular $m$-dimensional subspace of $\R^d$. In such a case, suitable choice of $\mS$ leads to $\cO\left(\frac{d}{m}\right)$--times faster convergence compared to the setup of Example~\ref{ex:coord_setup}. In Section~\ref{sec:sega_exp_aggressive} we numerically verify this claim. 
\end{remark}

\section{Convergence of \texttt{SEGA} for general sketches \label{sec:sega_analysis}}

In this section we state a linear convergence result for \texttt{SEGA} (Algorithm~\ref{alg:sega_gs}) for general sketch distributions $\cD$ under  smoothness and strong convexity assumptions. 

\subsection{Smoothness assumptions}

We will use the following general version of smoothness. 

\begin{assumption}[$\mQ$-smoothness]  \label{ass:sega_M_smooth_inv} Function $f$ is $\mQ$-smooth for some $\mQ\succ 0$, that is, for all $x,y\in \R^d$, the following inequality is satisfied:
    \begin{align}\label{eq:sega_M_smooth_inv}
         f(x) - f(y) - \langle \nabla f(y),x - y \rangle\ge \frac{1}{2}\|\nabla f(x) - \nabla f(y)\|_{\mmM}^2.
    \end{align}
\end{assumption}
Assumption~\ref{ass:sega_M_smooth_inv} is not standard in the literature. However, as Lemma~\ref{lem:sega_relate} states, for twice differentiable $f$ with $\mQ=\mM^{-1}$, Assumption~\ref{ass:sega_M_smooth_inv} is equivalent to $\mM$-smoothness (see~\eqref{eq:acd_M-smooth-intro}), which is a common assumption in modern analysis of \texttt{CD} methods. As discussed in Chapter~\ref{acd}, $\mM$-smoothness appears naturally in various application such as empirical risk minimization with linear predictors and is a baseline in the development of minibatch \texttt{CD} methods~\cite{nsync, qu2016coordinate1, qu2016coordinate2, sdna}. We will adopt this notion in Section~\ref{sec:sega_CD}, when comparing \texttt{SEGA} to coordinate descent. Until then, let us consider the almost equivalent Assumption~\ref{ass:sega_M_smooth_inv}.

\subsection{Main result}

We are now ready to present one of the key theorems of the chapter, which states that the iterates of  \texttt{SEGA} converge linearly to the optimal solution.  

\begin{theorem}\label{thm:sega_main}
    Assume that $f$ is $\mmM$-smooth and $\mu$-strongly convex. Choose stepsize $\alpha>0$ and Lyapunov parameter $\sigma>0$ so that 
\begin{equation}
        \alpha\left(2(\mC - \mI) +\sigma \mu \mI\right) \le \sigma\E{\mZ},\qquad \alpha \mC \le \frac{1}{2}\left(\mmM - \sigma \E{\mZ}\right), \label{eq:sega_general_bound_on_stepsize}        
    \end{equation} where $\mC\eqdef \E{\theta_k^2 \mZ_k}$. Fix $x^0,h^0\in {\rm dom} (F)$ and let $x^k,h^k$ be the
 random iterates produced by  \texttt{SEGA}.     Then
\[
        \E{\Lgen^{k}} \le (1 - \alpha\mu)^k \Lgen^0,
\]
where  $\Lgen^k \eqdef \|x^k - x^*\|^2 + \sigma \alpha \|h^k - \nabla f(x^*)\|^2$ is a Lyapunov function and $x^*$ is the solution of~\eqref{eq:sega_main_sega}.
\end{theorem}

Note that the convergence of the Lyapunov function $\Lgen^k$ implies both $x^k \rightarrow x^*$ and $h^k \rightarrow \nabla f(x^*)$. The latter means that \texttt{SEGA} is {\em variance reduced}, in contrast to \texttt{CD} in the proximal setup with  non-separable $\psi$, which does not converge to the solution.

To clarify on the assumptions, let us mention that if $\sigma$ is small enough so that $\mmM - \sigma \E{\mZ}\succ 0$, one can always choose stepsize $\alpha$ satisfying
\begin{align}\label{eq:sega_alfa_0}
    \alpha \leq \min\left\{
\frac {\lambda_{\text{min}}(\E{\mZ})}{\lambda_{\max} (2\sigma^{-1}(\mC - \mI) + \mu\mI)}, \frac{\lambda_{\min}(\mmM - \sigma \E{\mZ})}{2\lambda_{\max}(\mC)} \right\}
\end{align}
and inequalities~\eqref{eq:sega_general_bound_on_stepsize} will hold. Therefore, we get the next corollary. 
\begin{corollary}\label{cor:sega_general}
    If $\sigma < \frac{\lambda_{\min}(\mmM)}{\lambda_{\max}(\E{\mZ})}$, $\alpha$ satisfies~\eqref{eq:sega_alfa_0} and $k\ge\frac{1}{\alpha \mu}\log \frac{\Lgen^0}{\epsilon}$, then \[\E{\|x^k - x^*\|^2} \le \epsilon.\]
\end{corollary}
As Theorem~\ref{thm:sega_main} is rather general, we also provide a simplified version thereof, complete with a simplified analysis (Theorem~\ref{thm:sega_simple} in Appendix~\ref{sec:sega_simple_SEGA}). In the simplified version we remove the proximal setting (i.e., we set $\psi\equiv 0$), assume $L$-smoothness\footnote{The standard $L$-smoothness assumption is a special case of $\mM$-smoothness for $\mM =L \mI$ and special case of $\mQ$-smoothness for $\mQ =L^{-1} \mI$.}, and only consider coordinate sketches with uniform probabilities. The result is provided as Corollary~\ref{cor:sega_simple}.
\begin{corollary} \label{cor:sega_simple} Let $\cD$ be the uniform distribution over the standard unit basis vectors in $\R^d$. If the stepsize satisfies 
\[0<\alpha \leq \min\left\{ \frac{1-\frac{L\sigma}{d}}{2Ld}, \frac{1}{n\left(\mu + \frac{2(d-1)}{\sigma}\right)} \right\},\]
then
$$\E{\Lgen^{k}} \leq (1-\alpha \mu)^k \Lgen^{0}.$$ Therefore, the iteration complexity is $\tilde{\cO}(dL/\mu)$.
\end{corollary}

\begin{remark}\label{rem:aggressive}
In the fully general setting, one might choose $\alpha$ to be bigger than bound~\eqref{eq:sega_alfa_0}, which depends on eigen properties of matrices $\E{\mZ}, \mC, \mQ$, leading to a better overall complexity according to Corollary~\ref{cor:sega_general}. However, in the simple case with $\mmM=\mI$ and $\mS_k = e_{i_k}$ with uniform probabilities, bound~\eqref{eq:sega_alfa_0} is tight. 
\end{remark}

\section{Convergence of \texttt{SEGA} for coordinate sketches\label{sec:sega_CD}}

In this section we compare \texttt{SEGA} with  coordinate descent. We demonstrate that, specialized to a particular  choice of the distribution $\cD$ (where $\mS$ is a random column submatrix of the identity matrix), which makes \texttt{SEGA} use the same random gradient information as that used in modern state-of-the-art randomized \texttt{CD} methods,  \texttt{SEGA} attains, up to a small constant factor, the same convergence rate as \texttt{CD} methods.

Firstly, in Section~\ref{sec:sega_nonacc} we develop \texttt{SEGA} with arbitrary ``coordinate sketches'' (Theorem~\ref{t:imp_dacc}). Then, in Section~\ref{s:acc} we develop an {\em accelerated  variant of \texttt{SEGA}} in a very general setup known as {\em arbitrary sampling} (see Theorem~\ref{t:imp_acc}) \cite{nsync, quartz, qu2016coordinate1, qu2016coordinate2}. Lastly, Corollary~\ref{cor:sega_imp_dacc} and Corollary~\ref{cor:sega_acc_imp} provide us with \emph{importance sampling} for both nonaccelerated and accelerated method, which matches up to a constant factor cutting-edge coordinate descent rates \cite{nsync,allen2016even} under the same oracle and assumptions\footnote{There was recently introduced a notion of importance minibatch sampling for coordinate descent~\cite{hanzely2018accelerated}. We state, without a proof, that \texttt{SEGA} with block coordinate sketches allows for the same importance sampling as developed in the mentioned chapter. 
}. Table~\ref{tab:CDcmp} summarizes the results of this section. We provide a dedicated analysis for the methods from this section in Appendix~\ref{sec:sega_proofs_CD}.

\begin{table}[t]
\centering
\begin{tabular}{|c|c|c|c|c|}
\hline
& \texttt{CD}  &   \texttt{SEGA}  \\
\hline
\hline
\begin{tabular}{c}Nonaccelerated method\\ importance sampling, $b=1$ \end{tabular} & 
$\frac{\tracee(\mM)}{ \mu} \log \frac{1}{\epsilon}$ \cite{rcdm}
& 
$8.55\cdot \frac{\tracee(\mM)}{ \mu} \log \frac{1}{\epsilon}$ 
\\
\hline
\begin{tabular}{c}Nonaccelerated method\\ arbitrary sampling  \end{tabular} & 
$ \left(\max_i \frac{v_i}{p_i \mu}\right) \log \frac{1}{\epsilon}$ \hfill \cite{nsync} 
& 
$8.55 \cdot \left(\max_i \frac{v_i}{p_i \mu}\right) \log \frac{1}{\epsilon}$
\\
\hline
\begin{tabular}{c}Accelerated method\\  importance sampling, $b=1$  \end{tabular} & 
$1.62\cdot\frac{\sum_i \sqrt{\mM_{ii}}}{\sqrt{\mu}}  \log \frac{1}{\epsilon} $ \cite{allen2016even} 
&
$9.8 \cdot \frac{\sum_i \sqrt{\mM_{ii}}}{\sqrt{\mu}}  \log \frac{1}{\epsilon}$ 
\\
\hline
\begin{tabular}{c}Accelerated method\\ arbitrary sampling  \end{tabular} & 
$1.62 \cdot\sqrt{ \max_i \frac{v_i}{p_i^2 \mu} }   \log \frac{1}{\epsilon} $  \cite{hanzely2018accelerated}
&
$9.8 \cdot \sqrt{ \max_i \frac{v_i}{p_i^2 \mu} }   \log \frac{1}{\epsilon}$
\\
\hline
\end{tabular}
\caption{Complexity results for coordinate descent (\texttt{CD}) and our sketched gradient method (\texttt{SEGA}), specialized to coordinate sketching, for $\mM$-smooth and $\mu$-strongly convex functions. }
\label{tab:CDcmp}
\end{table}

We now describe the setup and technical assumptions for this section. In order to facilitate a direct comparison with \texttt{CD} (which 
does not work with non-separable regularizer $\psi$), for simplicity we consider problem~\eqref{eq:sega_main_sega} in the simplified setting with $\psi\equiv 0$.  Further, function $f$ is assumed to be $\mM$-smooth~\eqref{eq:acd_M-smooth-intro} and $\mu$-strongly convex. 


\subsection{Defining $\cD$: samplings}

In order to draw a direct comparison with general variants of \texttt{CD} methods (i.e., with those analyzed in the {\em arbitrary sampling} paradigm), we consider sketches in~\eqref{eq:sega_sketch-n-project} that are column submatrices of the identity matrix:
$\mS  = \mI_S,$
where $S$ is a random subset (aka {\em sampling}) of $[d]\eqdef \{1,2,\dots,d\}$.  Note that the columns of $\mI_S$  are the standard basis vectors $e_i$ for $i\in S$ and hence \[\Range{\mS} = \Range{e_i\;:\; i\in S}.\] So, distribution $\cD$ from which we draw matrices  is uniquely determined  by the distribution of sampling $S$. Given a sampling $S$, define $p = (p_1,\dots,p_d)\in \R^d$ to be the vector satisfying $p_i=\Probbb{e_i\in \Range{\mS}} = \Probbb{i\in S}$, and $\Probmat$ to be the matrix for which 
\[\Probmat_{ij}=\Probbb{\{i,j\}\subseteq S}.\]

Note that $p$ and $\Probmat$ are the {\em probability vector} and {\em probability matrix} of sampling $S$, respectively~\cite{qu2016coordinate2}. We assume throughout the chapter that $S$ is proper, i.e., we assume that $p_i>0$ for all $i$. State-of-the-art minibatch \texttt{CD} methods (including the ones we compare against~\cite{nsync, hanzely2018accelerated}) utilize large stepsizes related to the so-called ESO \emph{Expected Separable Overapproximation (ESO)}~\cite{qu2016coordinate2} parameters $v=(v_1,\dots,v_d)$. ESO parameters play a key role in \texttt{SEGA} as well, and are defined next. 
\begin{assumption}[ESO]\label{ass_ESO}
There exists a vector $v$ satisfying the following inequality
\begin{equation}\label{eq:sega_ESO}
\Probmat \circ \mM \preceq \diag(p) \diag(v) ,
\end{equation}
where $\circ$ denotes the Hadamard (i.e., element-wise) product of matrices.
\end{assumption}
In case of single coordinate sketches, parameters $v$ are equal to coordinate-wise smoothness constants of $f$. An extensive study on how to choose them in general was performed in~\cite{qu2016coordinate2}. 
For notational brevity, let us set $\mPdiag \eqdef\diag(p)$ and $\mVdiag\eqdef \diag(v)$ throughout this section.

%

\subsection{Non-accelerated method \label{sec:sega_nonacc}}

We now state the convergence rate of (non-accelerated) \texttt{SEGA} for coordinate sketches with {\em arbitrary sampling} of subsets of coordinates. The corresponding \texttt{CD} method was developed in~\cite{nsync}.

\begin{theorem}\label{t:imp_dacc}
Assume that $f$ is $\mM$-smooth and $\mu$-strongly convex. Denote $\Lnacc^{k} \eqdef f(x^{k})-f(x^*)+ \sigma \|h^{k} \|^2_{\mPdiag^{-1}}$. Choose $\alpha, \sigma>0$ such that 
\begin{equation}\label{eq:sega_assumption}
\sigma \mI-\alpha^2(\mVdiag\mPdiag^{-1}-\mM) \succeq \gamma \mu \sigma \mPdiag^{-1},
\end{equation}
where $\gamma \eqdef \alpha - \alpha^2\max_{i}\{\frac{v_{i}}{p_{i}}\}-\sigma$. Then the iterates of \texttt{SEGA} satisfy
$$
\E{\Lnacc^{k}}\leq (1-\gamma\mu )^k \Lnacc^0.
$$
\end{theorem}

We now give an importance sampling result for a coordinate version of  \texttt{SEGA}. We recover, up to a constant factor, the same convergence rate as standard \texttt{CD}~\cite{rcdm}. The probabilities we chose are optimal in our analysis and are proportional to the diagonal elements of matrix $\mM$. 

\begin{corollary}\label{cor:sega_imp_dacc}
Assume that $f$ is $\mM$-smooth and $\mu$-strongly convex. Suppose that $\cD$ is such that at each iteration standard unit basis vector $e_i$ is sampled with probability $p_i\propto \mM_{ii}$.
If we choose $ \alpha=\frac{0.232}{\tracee(\mM)}, \sigma=\frac{0.061}{ \tracee(\mM)}$, then
$$
\E{\Lnacc^{k}}\leq \left(1-\frac{0.117 \mu}{\tracee(\mM)} \right)^k \Lnacc^0.
$$
\end{corollary}

The iteration complexities provided in Theorem~\ref{t:imp_dacc} and Corollary~\ref{cor:sega_imp_dacc} are summarized in Table~\ref{tab:CDcmp}. We also state that $\sigma, \alpha$ can be chosen so that~\eqref{eq:sega_assumption} holds, and the rate from Theorem~\ref{t:imp_dacc} coincides with the rate from Table~\ref{tab:CDcmp}.

 \begin{remark}
 Theorem~\ref{t:imp_dacc} and Corollary~\ref{cor:sega_imp_dacc} hold even under a non-convex relaxation of strong convexity -- Polyak-\L{}ojasiewicz inequality: $\mu (f(x)- f(x^*))\leq \frac{1}{2}\|\nabla f(x) \|_2^2$. Therefore, \texttt{SEGA} also converges for a certain class of non-convex problems. For an overview on different relaxations of strong convexity, see~\cite{karimi2016linear}.
 \end{remark}

\subsection{Accelerated method \label{s:acc}}
 
 In this section, we propose an accelerated (in the sense of Nesterov's method~\cite{nesterov83,nesterov2018lectures}) version of \texttt{SEGA}, which we call {\acrshort{ASEGA}}. The analogous accelerated \texttt{CD} method, in which a single coordinate is sampled in every iteration, was developed and analyzed in~\cite{allen2016even}. The general variant utilizing arbitrary sampling was developed and analyzed in~\cite{hanzely2018accelerated}.

\begin{algorithm}[!h]
\begin{algorithmic}[1]
\State \textbf{Parameters:}  $x^0=y^0=z^0\in \R^d$; $h^0\in \R^d$; 	 $S$; parameters $\alpha, \beta,\eta, \mu>0$
\For {$k= 0,1,2,\dots $} 
    \State $x^{k}=(1-\eta)y^{k-1}+\eta z^{k-1}$ \label{eq:sega_x_update}
\State Sample $\mS_k = \mI_{S_k}$, where $S_k\sim S$, and compute $g^{k},h^{k+1}$ according to~\eqref{eq:sega_h^{k+1}},~\eqref{eq:sega_g^k}
 \State$y^{k}=x^{k}-\alpha\mPdiag^{-1} g^k $\label{eq:sega_y_update}  
   \State $z^{k}=\frac{1}{1+\beta \mu}(z^k+\beta \mu x^{k}-\beta g^k) $\label{eq:sega_z_update}
   \EndFor
      \end{algorithmic}
\caption{\texttt{ASEGA}: Accelerated \texttt{SEGA}}
   \label{alg:sega_acc}
\end{algorithm}

The method and analysis is inspired by~\cite{allen2014linear}. Due to space limitations and technicality of the content, we state the main theorem of this section in Appendix~\ref{sec:sega_acc_thm}. Here, we provide  Corollary~\ref{cor:sega_acc_imp}, which shows that Algorithm~\ref{alg:sega_acc} with single coordinate sampling enjoys, up to a constant factor, the same convergence rate as state-of-the-art accelerated coordinate descent method \texttt{NUACDM} of Allen-Zhu et al.~\cite{allen2016even}. 

 \begin{corollary}\label{cor:sega_acc_imp}
Let the sampling be defined as follows: $S = \{i\}$ with probability $p_i \propto \sqrt{\mM_{ii}}$, for $i \in [d]$. Then there exist acceleration parameters and a Lyapunov function $\Lacc^k$ such that $f(y^k) - f(x^*) \le \Lacc^k$ and
 \[
\E{\Lacc^{k}}\leq (1-\eta)^k\Lacc^0=\left(1-\cO\left(\frac{ \sqrt{\mu}}{\sum_{i} \sqrt{\mM_{ii}}} \right) \right)^k\Lacc^0.
\]
 \end{corollary}

The iteration complexity guarantees provided by Theorem~\ref{t:imp_acc} and Corollary~\ref{cor:sega_acc_imp} are summarized in Table~\ref{tab:CDcmp}.

 \section{Experiments} \label{sec:sega_experiments}
  In this section we perform numerical experiments to illustrate the potential of \texttt{SEGA}. Firstly, in Section~\ref{sec:sega_exp_pgd}, we compare it  to projected gradient descent (\texttt{PGD}) algorithm. Then in Section~\ref{sec:sega_exp_zero}, we study the performance of zeroth-order \texttt{SEGA} (when sketched gradients are being estimated through function value evaluations) and compare it to the analogous zeroth-order method. Next, in Section~\ref{sec:sega_exp_aggressive} we verify the claim from Remark~\ref{rem:aggressive} that in some applications, particular sketches might lead to a significantly faster convergence. Lastly, Section~\ref{sec:sega_cd_exp} demonstrates that {\tt SEGA} is competitive to {\tt CD} methods when $\psi \equiv 0$ as the results from Section~\ref{sec:sega_CD} predict.

  In the all experiments where theory-supported stepsizes were used -- we obtained them by precomputing strong convexity and smoothness measures.

\subsection{Comparison to projected gradient descent\label{sec:sega_exp_pgd}}

In this experiment, we illustrate the potential superiority of our method to \texttt{PGD}. We consider the $\ell_2$ ball constrained problem ($\psi$ is the indicator function of the unit ball) with the oracle providing the sketched gradient in the random Gaussian direction. As we mentioned in the introduction, a method moving in the gradient direction (analogue of \texttt{CD}), will not converge due to the proximal nature of the problem. Therefore, we can only compare against the projected gradient. However, in order to obtain the full gradient, one needs to gather $n$ sketched gradients and solve a linear system to recover the gradient. To illustrate this, we choose 4 different quadratic problems of the form \[f(x)\eqdef \frac12 x^\top \mM x-b^\top x,\] where $b$ is a random vector with independent entries from $\cN(0,1)$ and $\mM\eqdef\mU\Sigma \mU^\top$ according to Table~\ref{tab:problem} for $\mU$ obtained from QR decomposition of random matrix with independent entries from $\cN(0,1)$. For each problem, the starting point was chosen to be a vector with independent entries from $\cN(0,1)$.  
\begin{table}[!h]
\centering
\small
\begin{tabular}{|c|c|}
\hline
Type  &   $\Sigma$  \\
\hline
\hline
1 &  Diagonal matrix with first  $n/2$ components equal to 1, the rest equal to $n$\\ \hline
2&  Diagonal matrix with first  $n-1$ components equal to 1, the last one equal to $n$\\ \hline
3&  Diagonal matrix with $i$th component equal to $i$\\ \hline
4&   Diagonal matrix with components coming from uniform distribution over $[0,1]$\\ \hline
\end{tabular}
\caption{Spectrum of $\mM$.}
\label{tab:problem}
\end{table}
We stress that these are synthetic problems generated for the purpose of illustrating the potential of our method against a natural baseline.  Figure~\ref{fig:sega_pgd_comp} compares \texttt{SEGA} and \texttt{PGD} under various relative cost scenarios of solving the linear system compared to the cost of the oracle calls. The results show that \texttt{SEGA}  significantly outperforms \texttt{PGD} as soon as solving the linear system is expensive, and is as fast as \texttt{PGD} even if solving the linear system comes for free.

\begin{figure}
\centering
\begin{minipage}{0.35\textwidth}
  \centering
\includegraphics[width =  \textwidth ]{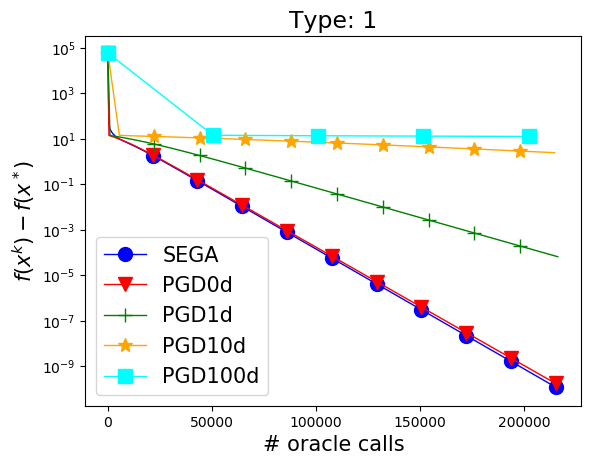}
\end{minipage}%
\begin{minipage}{0.35\textwidth}
  \centering
\includegraphics[width =  \textwidth ]{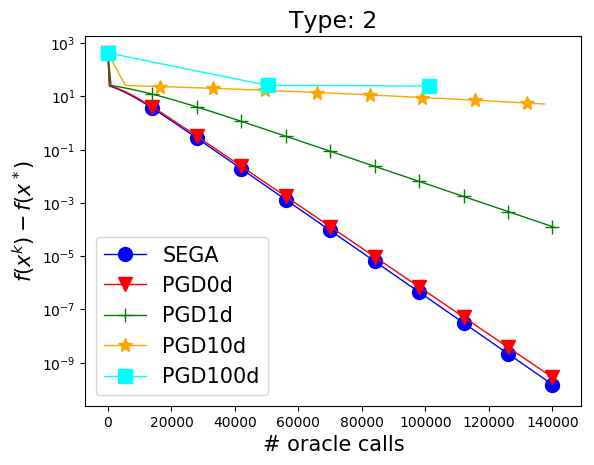}
\end{minipage}%
\\
\begin{minipage}{0.35\textwidth}
  \centering
\includegraphics[width =  \textwidth ]{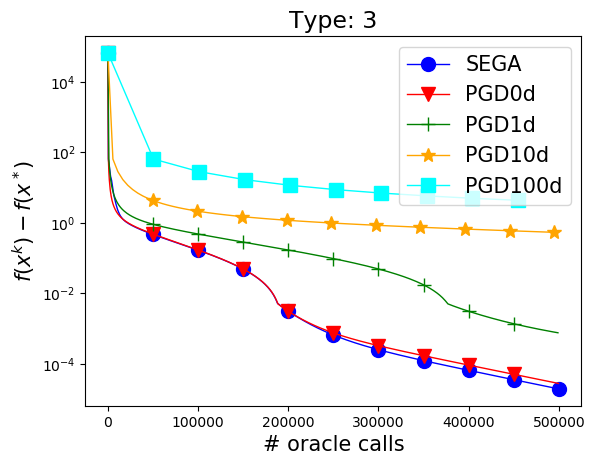}
\end{minipage}%
\begin{minipage}{0.35\textwidth}
  \centering
\includegraphics[width =  \textwidth ]{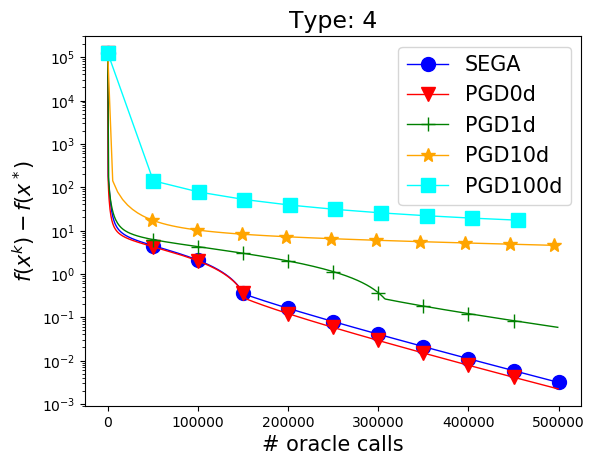}
\end{minipage}%
\caption{ 
Convergence of \texttt{SEGA} and \texttt{PGD} on synthetic problems with  $d=500$. The indicator ``Xd'' in the label indicates  the setting where the cost of solving linear system is $Xd$ times higher comparing to the cost ov=f evaluating a single directional derivative. Recall that a linear system is solved after each $d$ oracle calls. Stepsizes $1/\lambda_{\max}(\mM)$ and $1/(d\lambda_{\max}(\mM))$ were used for \texttt{PGD} and \texttt{SEGA}, respectively. }\label{fig:sega_pgd_comp}
\end{figure}

\subsection{Comparison to zeroth-order optimization methods\label{sec:sega_exp_zero}}
In this section, we compare \texttt{SEGA} to the {\em random direct search} (\texttt{RDS}) method~\cite{RDS} under a zeroth-order oracle for unconstrained optimization. For \texttt{SEGA}, we estimate the sketched gradient using finite differences. Note that \texttt{RDS} is a randomized version of the classical direct search method~\cite{hooke1961direct, kolda2003optimization, konevcny2014simple}. At iteration $k$, \texttt{RDS} moves to \[\argmin \left(f(x^k+\alpha^k s^k),f(x^k-\alpha^k s^k),f(x^k)\right)\] for a random direction $s^k\sim \cD$ and a suitable stepszie $\alpha^k$. For illustration, we choose $f$ to be a quadratic problem based on Table~\ref{tab:problem} and compare both Gaussian and coordinate directions. Figure~\ref{fig:sega_DFO} shows that \texttt{SEGA} outperforms \texttt{RDS}. 

\begin{figure}
\centering
\begin{minipage}{0.35\textwidth}
  \centering
\includegraphics[width =  \textwidth ]{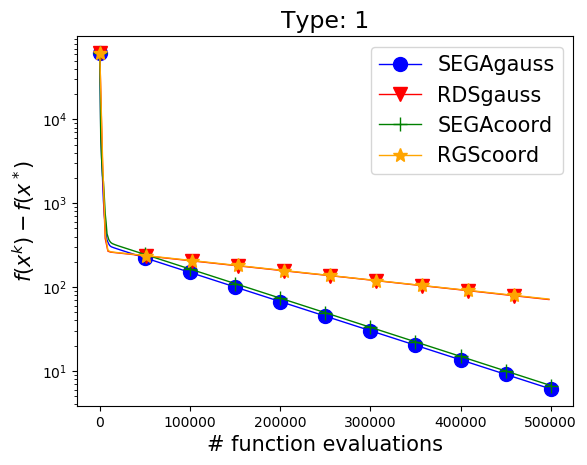}
\end{minipage}%
\begin{minipage}{0.35\textwidth}
  \centering
\includegraphics[width =  \textwidth ]{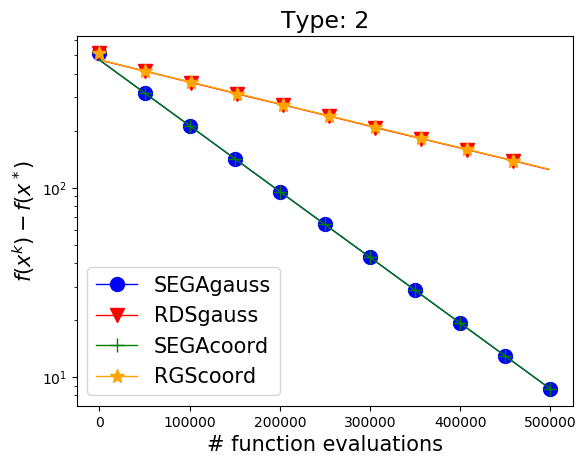}
\end{minipage}%
\\
\begin{minipage}{0.35\textwidth}
  \centering
\includegraphics[width =  \textwidth ]{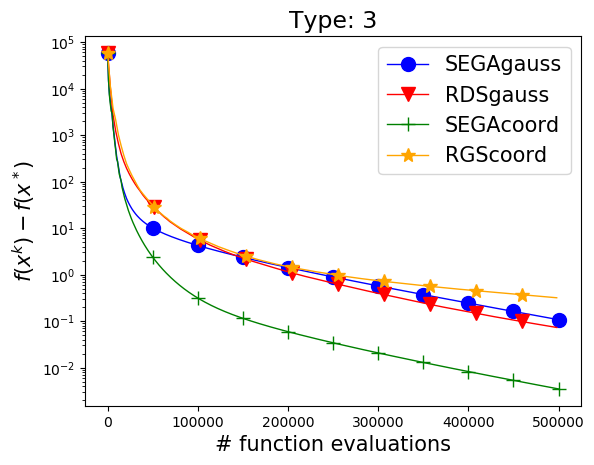}
\end{minipage}%
\begin{minipage}{0.35\textwidth}
  \centering
\includegraphics[width =  \textwidth ]{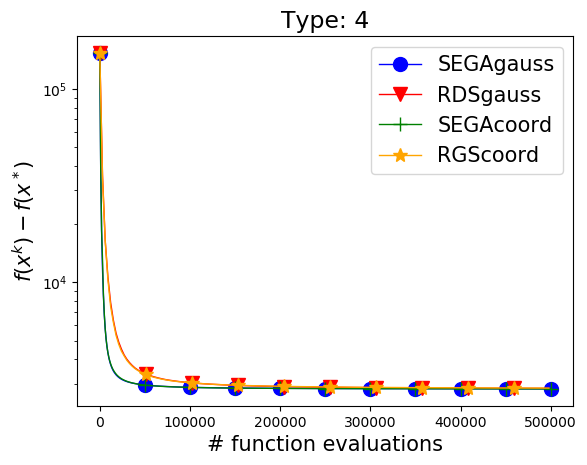}
\end{minipage}%
\caption{ Comparison of \texttt{SEGA} and randomized direct search for various problems. Theory supported stepsizes were chosen for both methods. 500 dimensional problem. }\label{fig:sega_DFO}
\end{figure}

\subsection{Subspace \texttt{SEGA} \label{sec:sega_exp_aggressive}}
As mentioned in Remark~\ref{rem:aggressive}, well designed sketches are capable of exploiting  structure of $f$ and lead to a better rate. We address this in detail Appendix~\ref{sec:sega_subSEGA} where we develop and analyze a subspace variant of \texttt{SEGA}. 

To illustrate this phenomenon in a simple setting, we perform experiments for problem \eqref{eq:sega_main_sega} with $f(x)=\| \mA x -b\|^2,$ where $b\in \R^{m}$ and $\mA\in \R^{m\times d}$ has orthogonal rows, and with $\psi$ being the indicator function of the unit ball in $\R^d$.  That is, we solve the problem
\[\min_{\|x\|_2 \leq 1} \| \mA x -b\|^2.\]
We assume that $d\gg m$.  We compare two methods: \texttt{naiveSEGA}, which uses coordinate sketches, and \texttt{subspaceSEGA}, where sketches are chosen as rows of $\mA$. Figure~\ref{fig:sega_aggressive} indicates that \texttt{subspaceSEGA} outperforms  \texttt{naiveSEGA} roughly by the factor $\frac{d}{m}$, as claimed in Appendix~\ref{sec:sega_subSEGA}.

\begin{figure}
\centering
\begin{minipage}{0.35\textwidth}
  \centering
\includegraphics[width =  \textwidth ]{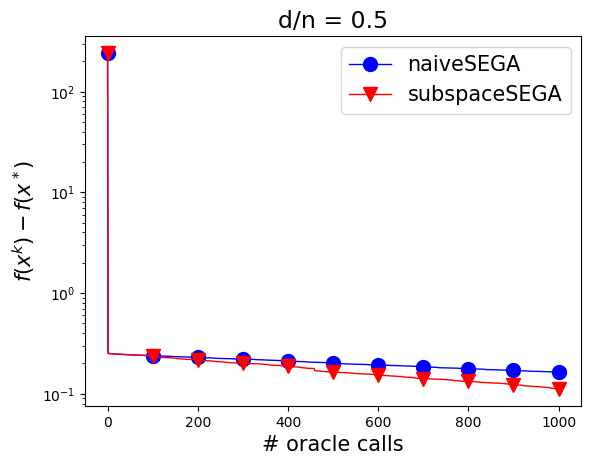}
\end{minipage}%
\begin{minipage}{0.35\textwidth}
  \centering
\includegraphics[width =  \textwidth ]{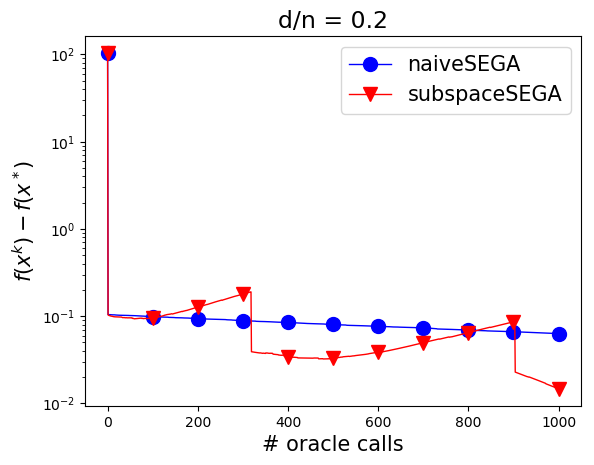}
\end{minipage}%
\\
\begin{minipage}{0.35\textwidth}
  \centering
\includegraphics[width =  \textwidth ]{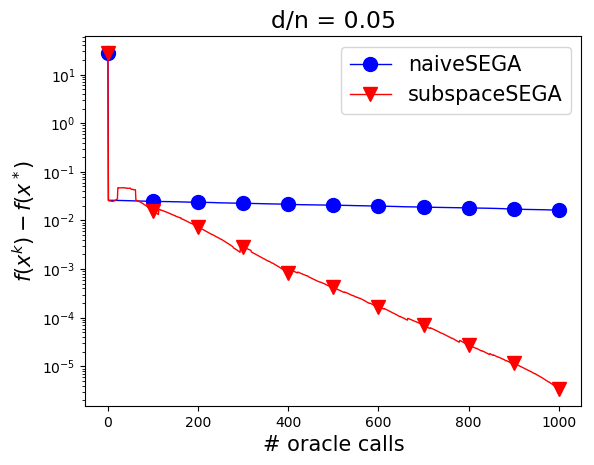}
\end{minipage}%
\begin{minipage}{0.35\textwidth}
  \centering
\includegraphics[width =  \textwidth ]{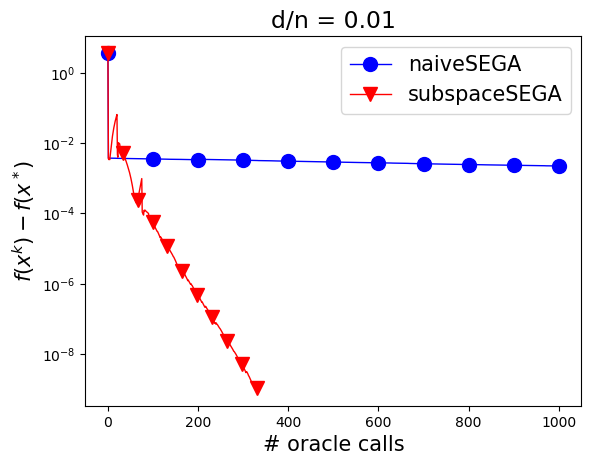}
\end{minipage}%
\caption{ Comparison of \texttt{SEGA} with sketches from a correct subspace versus coordinate sketches \texttt{naiveSEGA}. Stepsize chosen according to theory. 1000 dimensional problem. }\label{fig:sega_aggressive}
\end{figure}

\subsection{Comparison to randomized coordinate descent \label{sec:sega_cd_exp}}
In this section we numerically compare the results from Section~\ref{sec:sega_CD} to analogous results for coordinate descent (as indicated in Table~\ref{tab:CDcmp}).  We consider the ridge regression problem on LibSVM~\cite{chang2011libsvm} data, for both primal and dual formulation. For all methods, we have chosen parameters as suggested from theory Figure~\ref{fig:sega_cd_cmp} shows the results. We can see that in all cases, \texttt{SEGA} is slower to the corresponding coordinate descent method, but still is competitive. We however observe only constant times difference in terms of the speed, as suggested by Table~\ref{tab:CDcmp}.

\begin{figure}[!h]
\centering
\begin{minipage}{0.25\textwidth}
  \centering
\includegraphics[width =  \textwidth ]{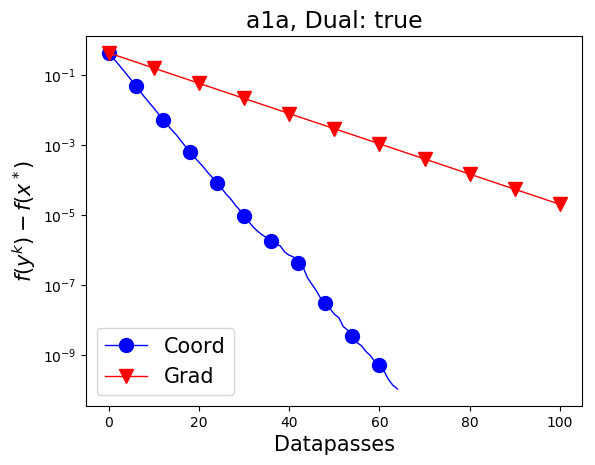}
\end{minipage}%
\begin{minipage}{0.25\textwidth}
  \centering
\includegraphics[width =  \textwidth ]{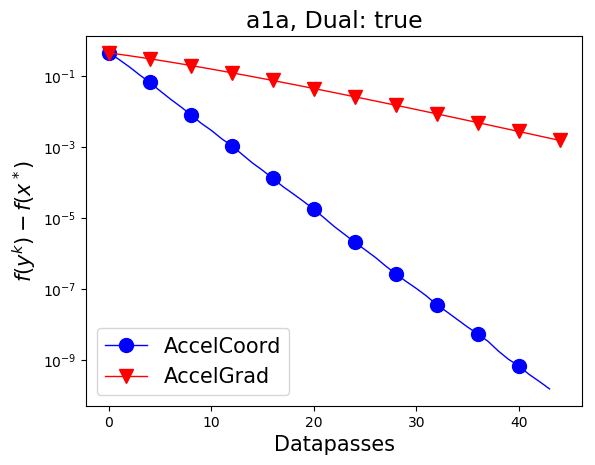}
\end{minipage}%
\begin{minipage}{0.25\textwidth}
  \centering
\includegraphics[width =  \textwidth ]{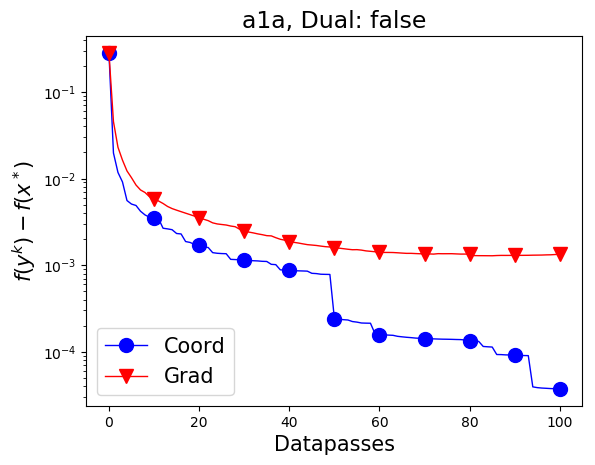}
\end{minipage}%
\begin{minipage}{0.25\textwidth}
  \centering
\includegraphics[width =  \textwidth ]{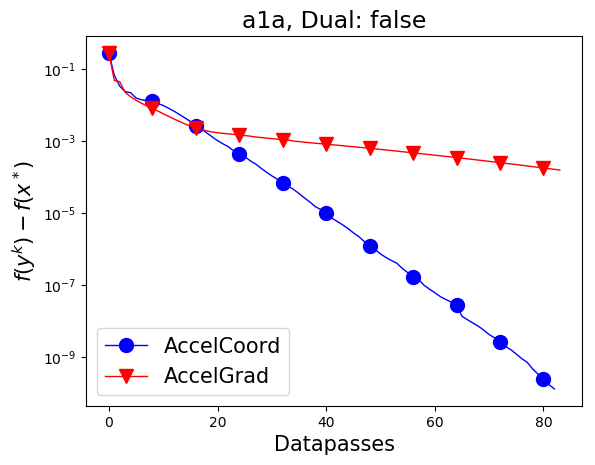}
\end{minipage}%
\\
\begin{minipage}{0.25\textwidth}
  \centering
\includegraphics[width =  \textwidth ]{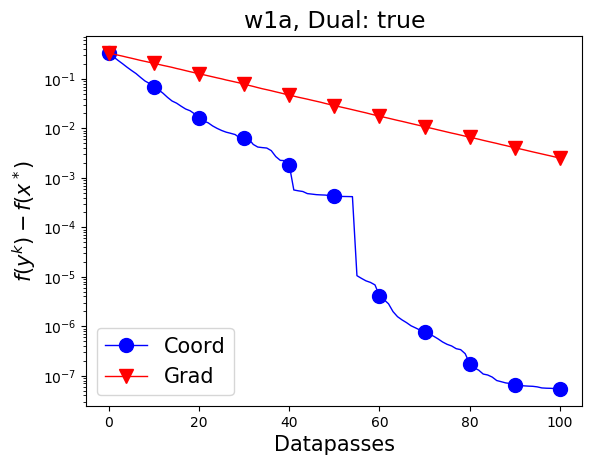}
\end{minipage}%
\begin{minipage}{0.25\textwidth}
  \centering
\includegraphics[width =  \textwidth ]{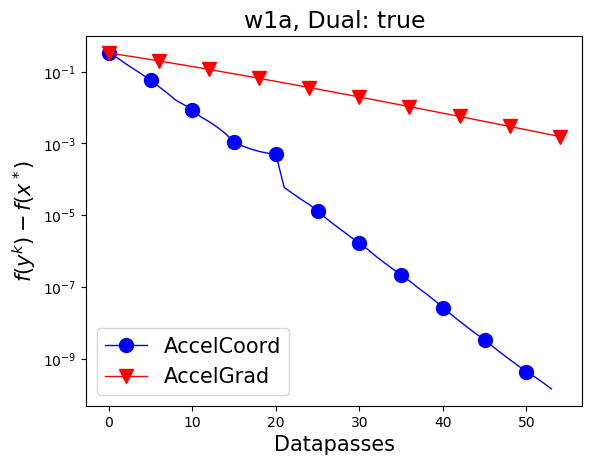}
\end{minipage}%
\begin{minipage}{0.25\textwidth}
  \centering
\includegraphics[width =  \textwidth ]{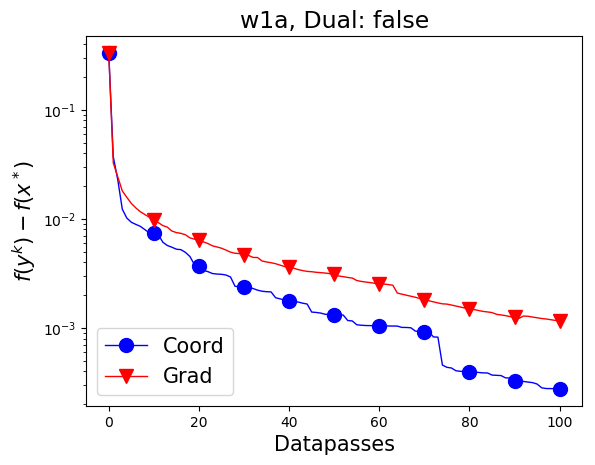}
\end{minipage}%
\begin{minipage}{0.25\textwidth}
  \centering
\includegraphics[width =  \textwidth ]{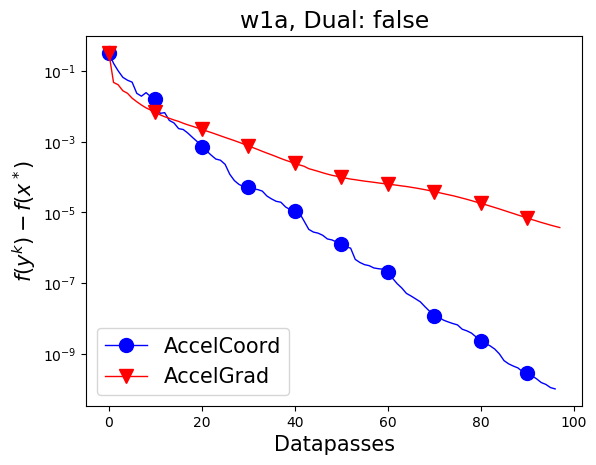}
\end{minipage}%
\\
\begin{minipage}{0.25\textwidth}
  \centering
\includegraphics[width =  \textwidth ]{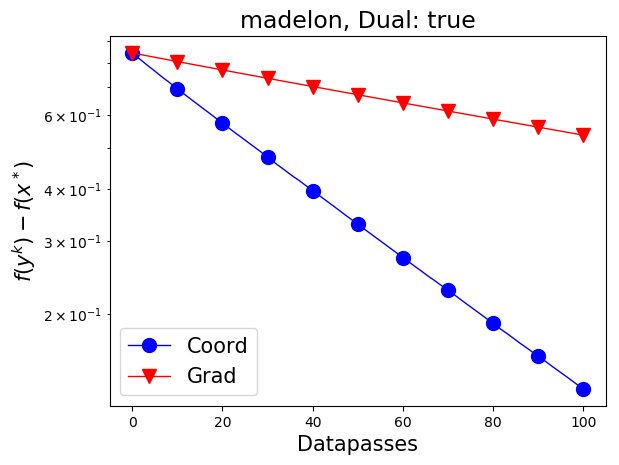}
\end{minipage}%
\begin{minipage}{0.25\textwidth}
  \centering
\includegraphics[width =  \textwidth ]{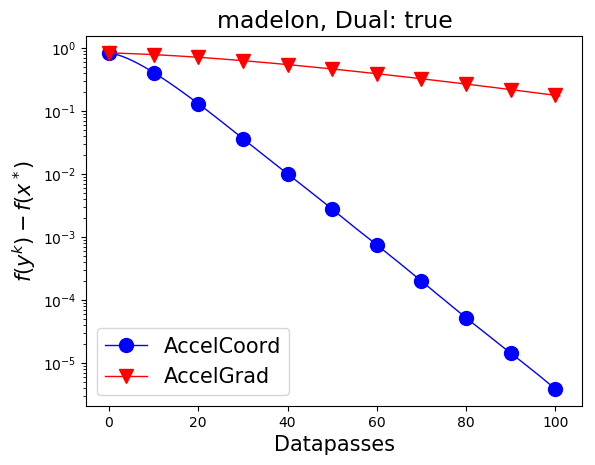}
\end{minipage}%
\begin{minipage}{0.25\textwidth}
  \centering
\includegraphics[width =  \textwidth ]{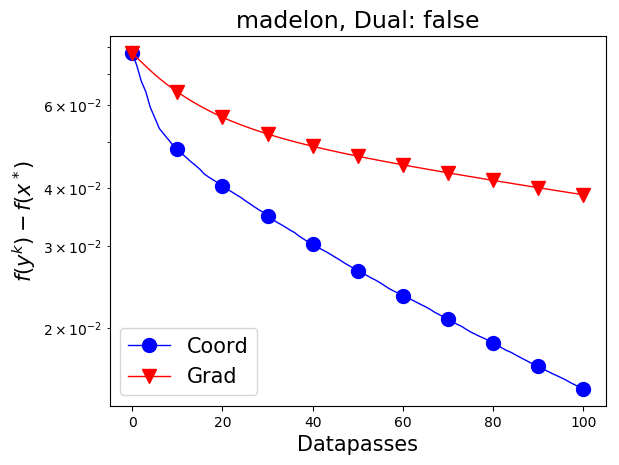}
\end{minipage}%
\begin{minipage}{0.25\textwidth}
  \centering
\includegraphics[width =  \textwidth ]{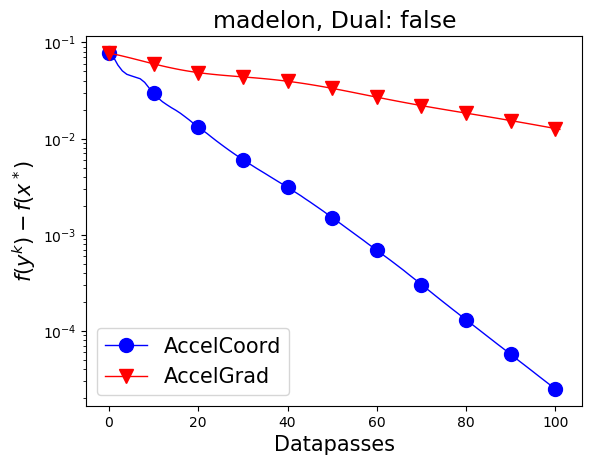}
\end{minipage}%
\\
\begin{minipage}{0.25\textwidth}
  \centering
\includegraphics[width =  \textwidth ]{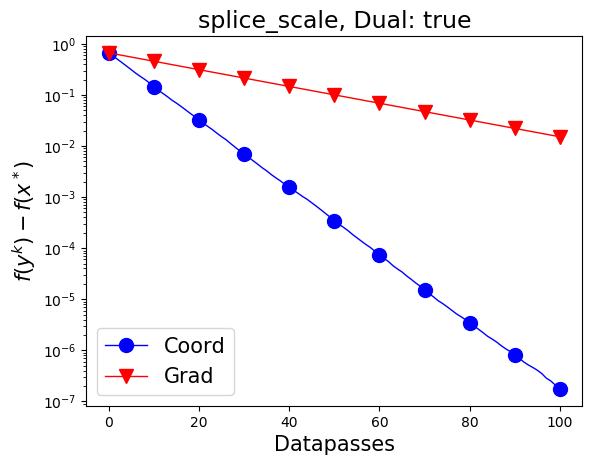}
\end{minipage}%
\begin{minipage}{0.25\textwidth}
  \centering
\includegraphics[width =  \textwidth ]{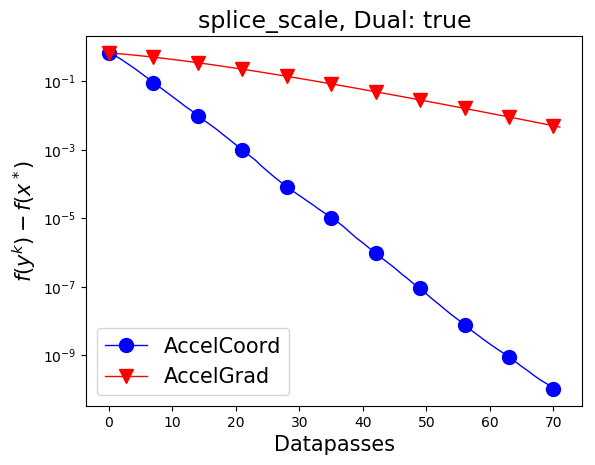}
\end{minipage}%
\begin{minipage}{0.25\textwidth}
  \centering
\includegraphics[width =  \textwidth ]{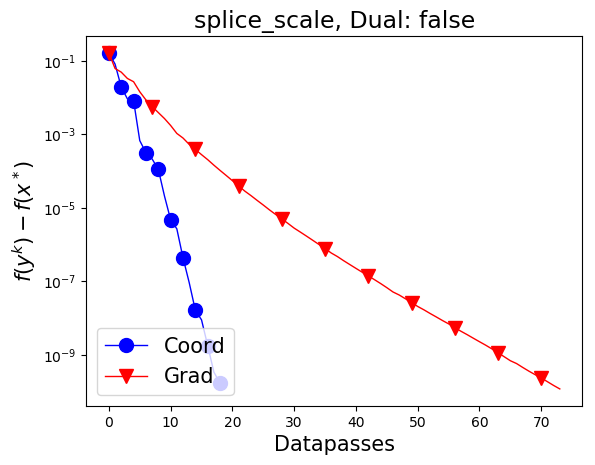}
\end{minipage}%
\begin{minipage}{0.25\textwidth}
  \centering
\includegraphics[width =  \textwidth ]{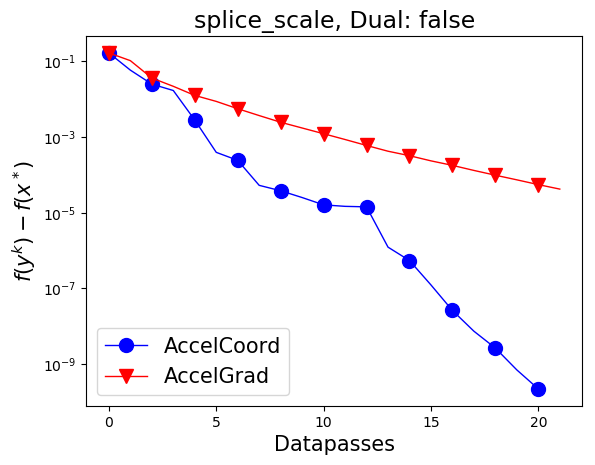}
\end{minipage}%
\\
\begin{minipage}{0.25\textwidth}
  \centering
\includegraphics[width =  \textwidth ]{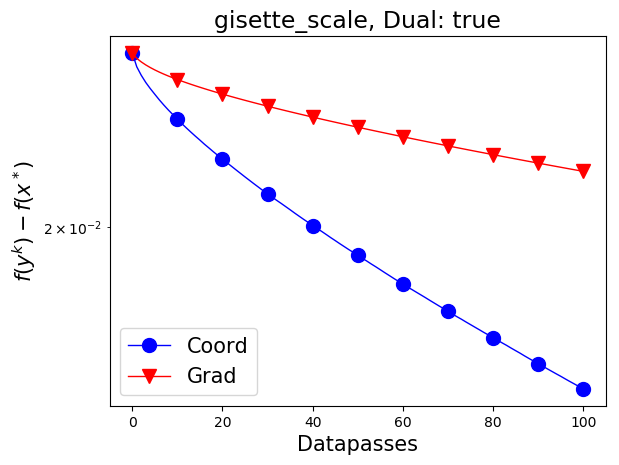}
\end{minipage}%
\begin{minipage}{0.25\textwidth}
  \centering
\includegraphics[width =  \textwidth ]{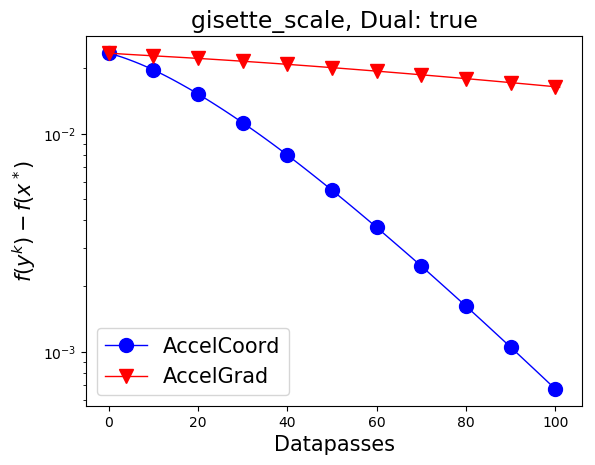}
\end{minipage}%
\begin{minipage}{0.25\textwidth}
  \centering
\includegraphics[width =  \textwidth ]{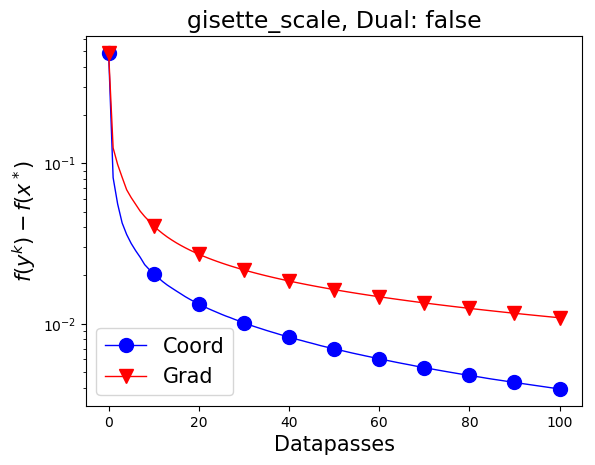}
\end{minipage}%
\begin{minipage}{0.25\textwidth}
  \centering
\includegraphics[width =  \textwidth ]{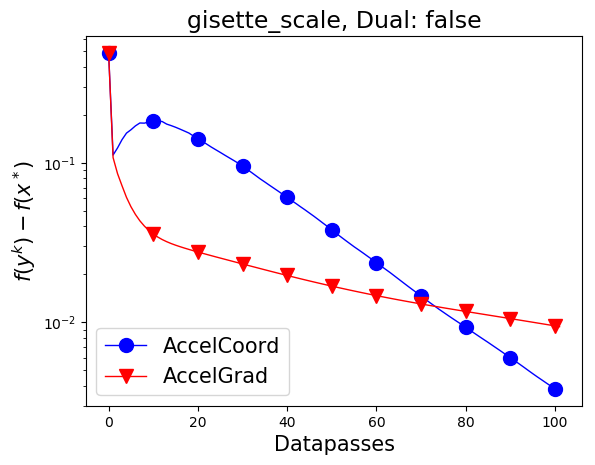}
\end{minipage}%
\caption{Comparison of \texttt{SEGA} and \texttt{ASEGA} with corresponding coordinate descent methods for $\psi\equiv 0$.}\label{fig:sega_cd_cmp}
\end{figure}

\subsection{Evolution of iterates: Extra plots\label{sec:sega_evolution_extra}}

Here we show some additional plots similar to Figure~\ref{fig:sega_trajectory}, which we believe help to build intuition about how the iterates of  \texttt{SEGA} behave. We also include plots for \texttt{biasSEGA}, which uses biased estimators of the gradient instead. We found that the iterates of \texttt{biasSEGA} often behave in a more stable way, as could be expected given the fact that they enjoy lower variance. However, we do not have any theory supporting the convergence of \texttt{biasSEGA}; this is left for future research.

\begin{figure}[!h]
\centering
\begin{minipage}{0.35\textwidth}
\centering
\includegraphics[width = \textwidth ]{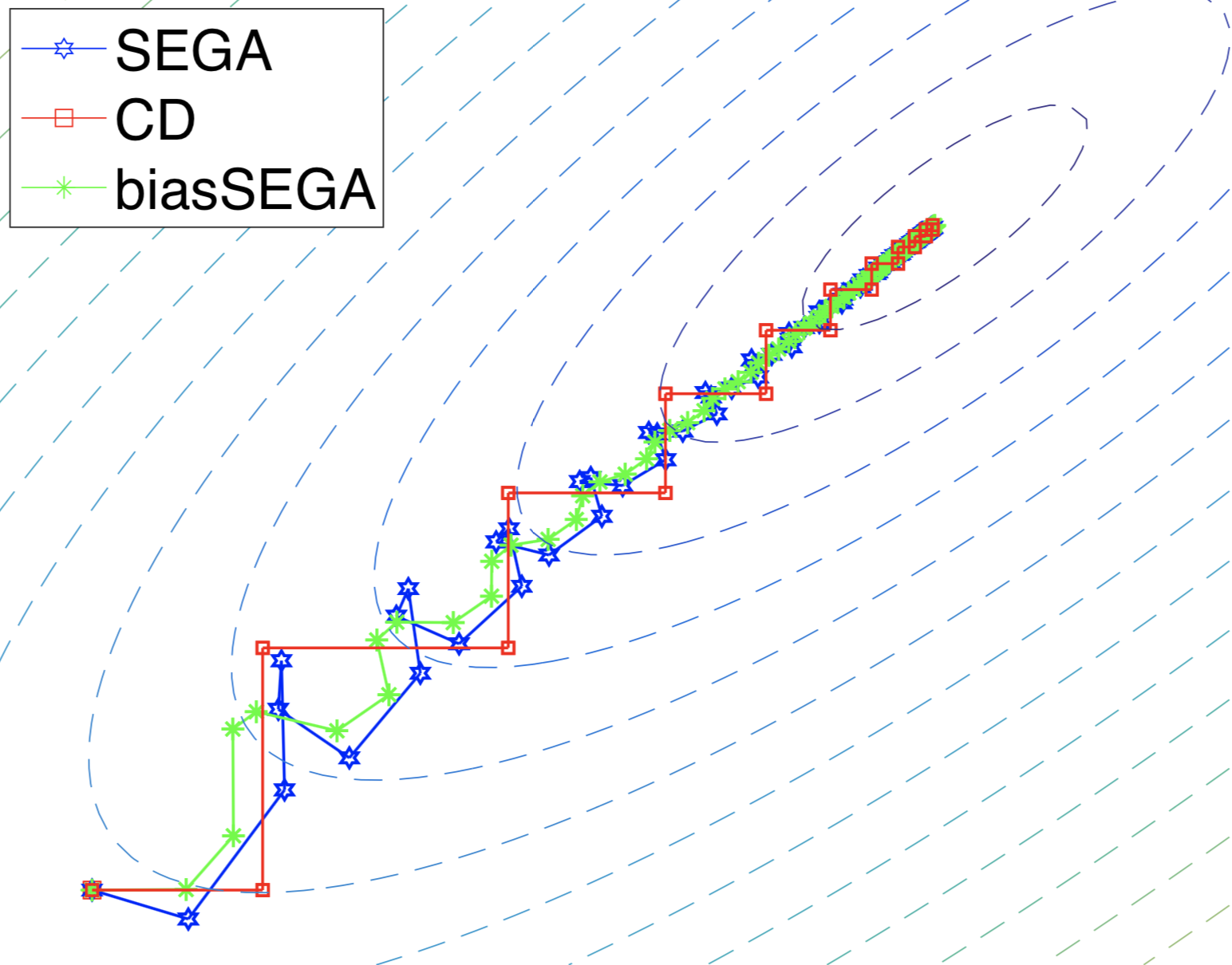}
\caption{Evolution of iterates of \texttt{SEGA}, \texttt{CD} and \texttt{biasSEGA} (updates made via $h^{k+1}$ instead of $g^k$).}
\end{minipage}
\hskip 1cm
\begin{minipage}{0.35\textwidth}
\centering
\includegraphics[width = \textwidth ]{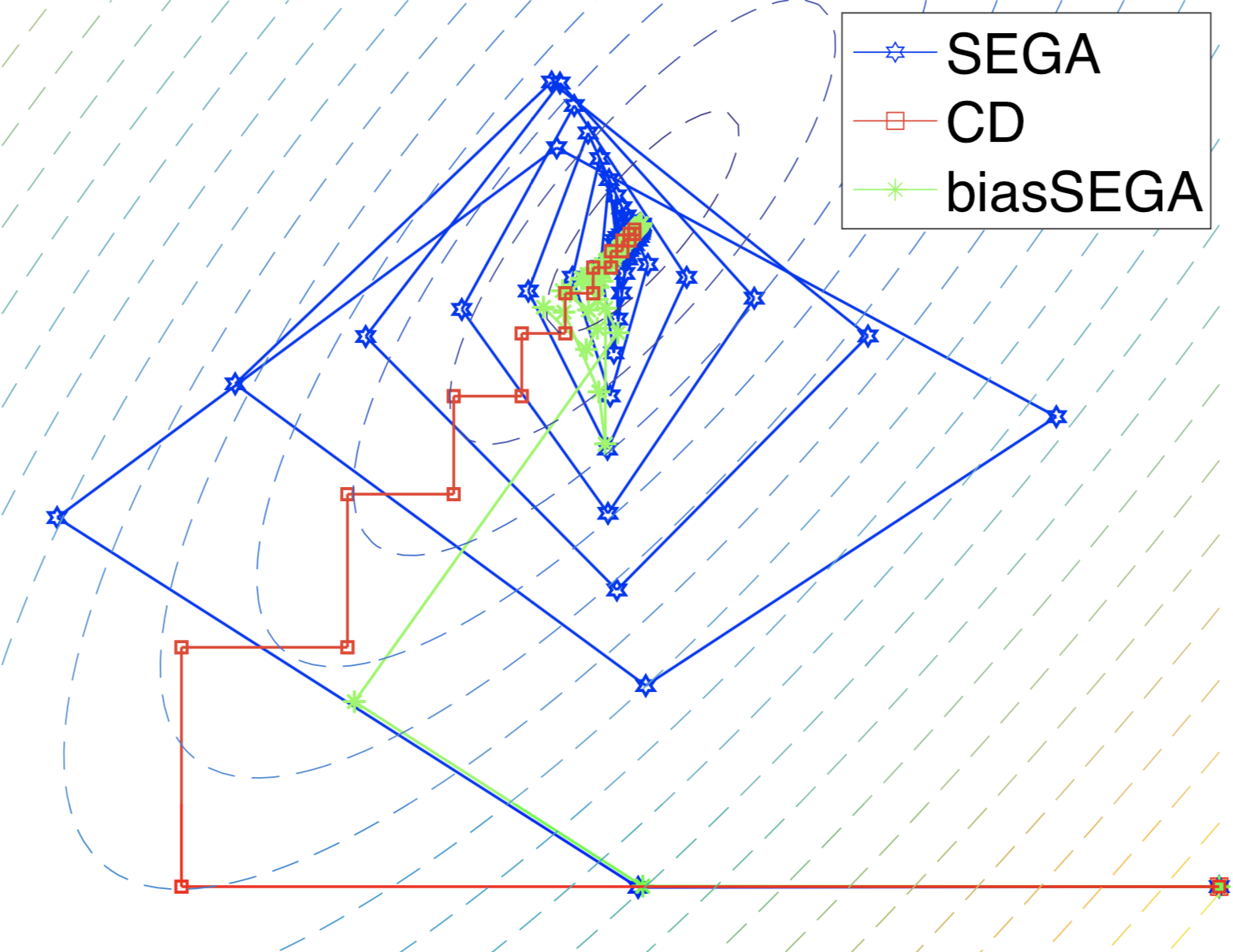}
\caption{Iterates of \texttt{SEGA}, \texttt{CD} and \texttt{biasSEGA} (updates made via $h^{k+1}$ instead of $g^k$). Different starting point.}
\end{minipage}
\\
\hskip 1cm
\begin{minipage}{0.35\textwidth}
\centering
\includegraphics[width = \textwidth ]{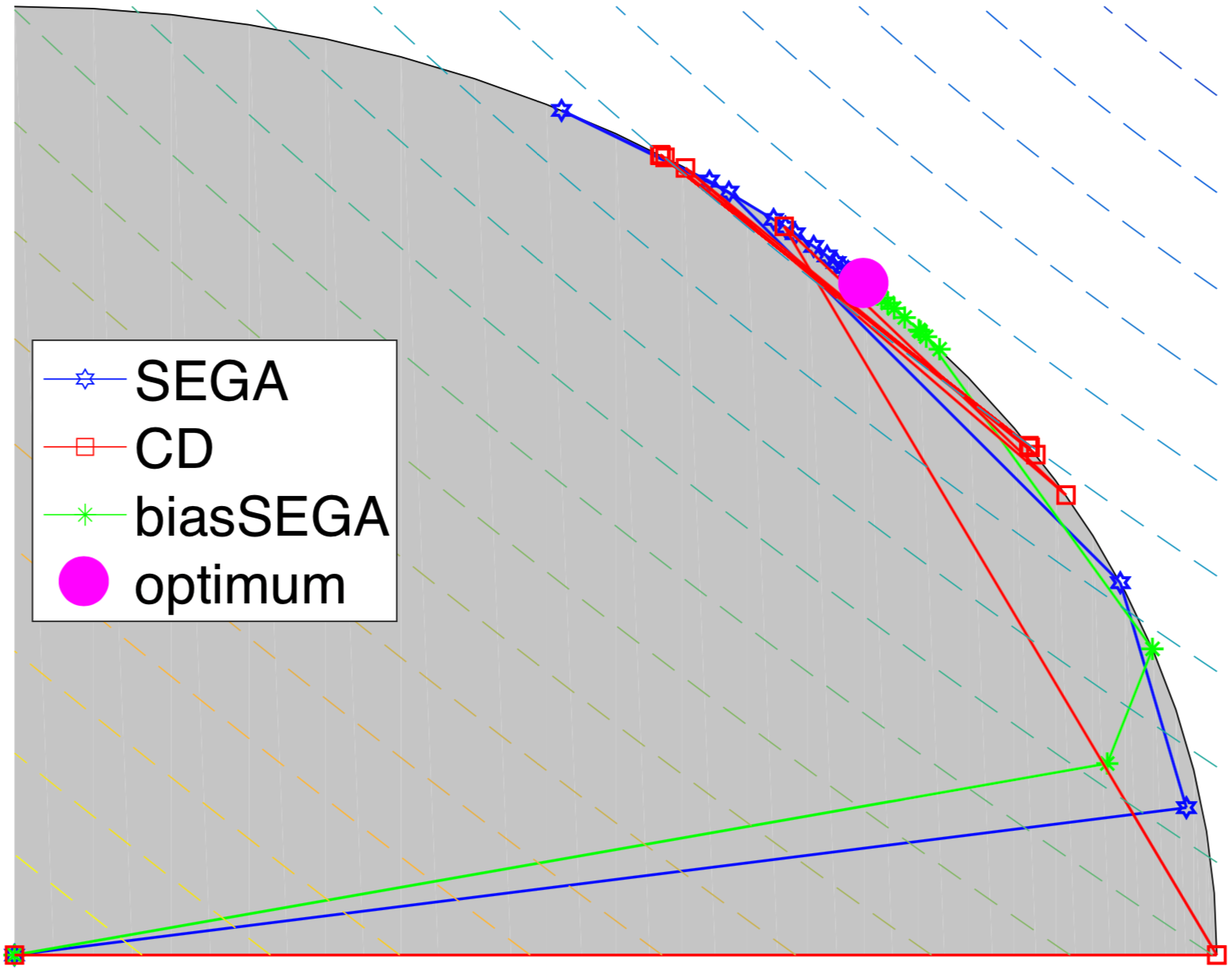}
\caption{Iterates of projected \texttt{SEGA}, projected \texttt{CD} (which do not converge) and projected \texttt{biasSEGA} (updates made via $h^{k+1}$ instead of $g^k$). The constraint set is represented by the shaded region.}
\end{minipage}
\end{figure}

\section{Conclusion} \label{sec:sega_conclusion}

We proposed \texttt{SEGA}, a  method for solving composite optimization problems under a novel stochastic linear first-order oracle. \texttt{SEGA} is variance-reduced, and this is achieved via  sketch-and-project updates of gradient estimates. We provided an analysis for smooth and strongly convex functions and general sketches, and  a refined analysis for coordinate sketches. For coordinate sketches we  also proposed an accelerated variant of \texttt{SEGA}, and  our theory matches that of state-of-the-art \texttt{CD} methods. However, in contrast to  \texttt{CD}, \texttt{SEGA} can be used for optimization problems with a {\em non-separable} proximal term. We develop a more aggressive subspace variant of the method---\texttt{subspaceSEGA}---which leads to improvements in the $d\gg m$ regime. In the Appendix we give several further results, including simplified and alternative analyses of \texttt{SEGA}  in the coordinate setup from Example~\ref{ex:coord_setup}.  Our experiments are encouraging and substantiate our theoretical predictions.

Next, we point to several potential  extensions of our work.

\paragraph{Speeding up the general method.}  We believe that it should be possible to extend \texttt{ASEGA} to the general setup from Theorem~\ref{thm:sega_main}. In such a case, it might be possible a distribution of sketches $\cD$ so as to outperform accelerated proximal gradient methods~\cite{Nesterov05:smooth, beck2009fista}. 

\paragraph{Biased gradient estimator.} Recall that \texttt{SEGA} uses unbiased gradient estimator $g^k$ for updating the iterates $x^k$ in a similar way \texttt{JacSketch}~\cite{jacsketch} or \texttt{SAGA}~\cite{saga} do this for the stochastic finite sum optimization. Recently, a  stochastic method for finite sum  optimization using biased gradient estimators was proven to be more efficient~\cite{nguyen2017sarah}. Therefore, it might be possible to establish better properties for a biased variant of \texttt{SEGA}. To demonstrate the potential of this approach, in Appendix~\ref{sec:sega_evolution_extra} we plot the evolution of iterates for the very simple biased method which uses $h^k$ as an update for line~\ref{eq:sega_x_update} in Algorithm~\ref{alg:sega_gs}.

\paragraph{Applications.} We believe that \texttt{SEGA} might work well in applications where a zeroth-order approach is inevitable, such as reinforcement learning. We therefore believe that \texttt{SEGA} might be an efficient proximal method in some reinforcement learning applications. We also believe that communication-efficient variants of \texttt{SEGA} can be used for  distributed training of machine learning models. This is because \texttt{SEGA} can be adapted to communicate sparse model updates only.

In the next chapter we introduce a different scenario where \texttt{SEGA} can be superior to \texttt{CD} even for problems without non-separable regularizer. The setups goes as follows: instead of minimizing a single function, we aim to minimize a finite sum. The oracle provides us with mutually independent random set of partial derivatives of each function from the sum. In such case, the gradient in the optimum does not have to be zero for each function, and thus \texttt{SEGA} trick might be necessary to keep fast convergence. However, the motivation for the mentioned setup does not come from \texttt{SEGA}, but rather that the independent sampling of coordinates yields surprisingly fast convergence.


\chapter{99\% of Worker-Master Communication in Distributed Optimization is Not Needed}
\label{99}

\renewcommand{\EE}{\mathbb{E}}

\graphicspath{{99/experiments/}}

In this work we are concerned with parallel/distributed  algorithms for solving finite sum minimization problems 
\begin{align}
  \min_{x\in \RR^d}  \left \{ f(x) \eqdef  \frac{1}{n}\sum \limits_{i=1}^n f_i(x) \right\},\label{eq:99_problem}
\end{align}
where each $f_i$ is convex and smooth. In particular, we are interested in methods which employ $n$ parallel units/workers/nodes/processors, each of which has access to a single function $f_i$ and its gradients (or unbiased estimators thereof). Let $x^*$ be an optimal solution of~\eqref{eq:99_problem}. In many practical scenarios, $f_i$ is often of the form
\begin{equation}\label{eq:99_stoch-f_i}      f_i(x) = \EE_{\xi} \phi_i(x; \xi),\end{equation}
where the expectation is with respect to a distribution of training examples stored locally at machine $i$. More typically, however, each machine contains a very large but finite number of examples (for simplicity, say there are $l$ examples on each machine), and $f_i$ is of the form \begin{equation}\label{eq:99_problem_saga_dist}
    f_i(x) = \frac{1}{l}\sum_{j=1}^l f_{ij}(x).
\end{equation}

 In the rest of this section we provide some basic motivation and intuitions in support of our approach. To this purpose, assume, for simplicity of exposition, that $f_i$ is of the finite-sum form \eqref{eq:99_problem_saga_dist}. In typical modern machine learning workloads, the number of machines $n$  is much smaller than the number of data points on each machine $l$. In a large scale regime (i.e., when the model size $d$,  the number of data points $nl$, or both are large),  problem \eqref{eq:99_problem} needs to be solved by a combination of efficient methods and modern hardware. In recent years there has been a lot of progress in designing new algorithms for solving this problem using techniques such as stochastic approximation~\cite{robbins}, variance reduction~\cite{sag, svrg, saga}, coordinate descent~\cite{rcdm,richtarik2014iteration, wright2015coordinate} and acceleration~\cite{nesterov83}, resulting in excellent theoretical and practical performance.

The computational power of the hardware is increasing as well. In recent years, a very significant amount of such increase is due to parallelism. Since many methods, such as minibatch Stochastic Gradient Descent ({\tt SGD}), are embarrassingly parallel, it is very simple to use them in big data applications. However, it has been observed in practice that adding more resources beyond a certain limit does not improve iteration complexity significantly. Moreover, having more parallel units makes their synchronization harder due to so-called communication bottleneck.  Minibatch versions of most variance reduced methods\footnote{We shall mention that there are already a few variance reduced methods that scale, up to some level, linearly in a parallel setup: Quartz for sparse data~\cite{quartz}, Katyusha~\cite{allen2017katyusha},  or {\tt SAGA}/{\tt SVRG}/SARAH with importance sampling for non-convex problems~\cite{horvath2018nonconvex}. } such as {\tt SAGA}~\cite{saga} or {\tt SVRG}~\cite{svrg} scale even worse in parallel setting -- they do not guarantee, in the worst case,  any speedup from using more than one function at a time. Unfortunately, numerical experiments show that this is not a proof flaw, but rather a real property of these methods~\cite{jacsketch}. A similar observation was made for {\tt SVRG} by~\cite{zhao2014accelerated}, where it was shown that only a small number of partial derivatives are needed at each iteration.

Since there are too many possible situations, we choose to focus on black-box optimization, although we admit that much can be achieved by assuming the sparsity structure. In fact, for any method there exists a toy situation where the method would scale perfectly -- one simply needs to assume that each function $f_i$ depends on its own subset of coordinates and minimize each $f_i$ independently. This can be generalized assuming sparsity patterns~\cite{leblond2017asaga, leblond2018improved} to get almost linear scaling if any coordinate appears in a small number of functions. Our interest, however, is in explaining situations as in~\cite{jacsketch} where the models almost do not scale.

In this chapter, we demonstrate that a simple trick of {\em independent} block sampling can remedy the problem of scaling,  to a substantial but limited extent. To illustrate one of the key insights on a simple example, in what follows consider a thought experiment in which {\tt GD} is a baseline method we would want to improve on.

\section{From gradient descent to block coordinate descent and back}

A simple benchmark in the distributed setting is a parallel implementation of gradient descent ({\tt GD}). {\tt GD} arises as a special case of the more general class of block coordinate descent methods ({\tt BCD})~\cite{rcdm}.  The conventional way to run  {\tt BCD} for problem \eqref{eq:99_problem} is to update a single or several blocks\footnote{Assume the entries of $x$ are partitioned into several non-overlapping blocks.} of $x$, chosen at random, on all  $n$ machines~\cite{rcdm, approx}, followed by an update aggregation step. Such updates on each worker typically involve a gradient step on a subspace corresponding to the selected blocks.  Importantly, and this is a key structural property of {\tt BCD} methods, {\em the same set of blocks is updated on each machine}. If  communication is expensive, it often makes sense to do more work on each machine,  which in the context of {\tt BCD} means updating more blocks.
A particular special case is to update {\em all} blocks, which leads to parallel implementation of {\tt GD} for problem \eqref{eq:99_problem}, as mentioned above. Moreover, it is known that the theoretical iteration complexity of {\tt BCD} improves as the number of blocks updated increases~\cite{rcdm, qu2016coordinate1, qu2016coordinate2}. For these and similar reasons, {\tt GD} (or one of its variants, such as {\tt GD} with momentum), is often the preferable method to {\tt BCD} (in terms of iteration complexity). Having said that, we did not choose to describe {\tt BCD} only to discard it at this point; we shall soon return to it, albeit with a twist.

\subsection{From gradient descent to independent block coordinate descent}

Because of what we have just said, iteration complexity of {\tt GD} will not improve by any variant running {\tt BCD}; it can only get worse.  Despite this, {\em we propose to run  {\tt BCD}, but a new variant which allows each worker to sample an independent subset of blocks} instead. This variant of {\tt BCD} for \eqref{eq:99_problem} was not considered before. As we shall show,  our {\em independent sampling} approach leads to a better-behaved aggregated gradient estimator when compared to that of {\tt BCD}, which in turn leads to better overall iteration complexity. We call our method {\em independent block coordinate descent} ({\tt IBCD}). 

We provide a unified analysis of our method, allowing for a random subset of $\tau m$ out of a total of $m$ blocks to be sampled on each machine, independently from other machines. {\tt GD}  arises as a special case of this method by setting $\tau=1$.  However, as we show (see Corollary~\ref{cor:99_0893y83}), {\em the same iteration complexity guarantee can be obtained by choosing $\tau$ as low as $\tau=\frac{1}{n}$.} The immediate consequence of this result is that {\em it is suboptimal to run {\tt GD} in terms of communication complexity.} Indeed, {\tt GD} needs to communicate all $m$ blocks per machine, while {\tt IBCD} achieves the same rate with $\frac{m}{n}$ blocks per machine only. Coming back to the abstract, consider an example with $n=100$ machines. In this case, when compared to {\tt GD}, {\tt IBCD} only communicates $1\%$ of the data. Because the iteration complexities of the two methods are the same, and if communication cost is dominant, this means that the problem can be solved in just $1\%$ of the time. In contrast, and when compared to the potential of {\tt IBCD}, parallel implementation of {\tt GD} inevitably wastes 99\% of the time.

The intuition behind why our approach works lies in the law of large numbers. By averaging independent noise we reduce the total variance of the resulting estimator by the factor of $n$. If, however, the noise is already tiny, as, in non-accelerated variance reduced methods, there is no improvement. On the other hand, (uniform)  block coordinate descent ({\tt CD}) has variance proportional to $\frac{1}{\tau}$~\cite{wangni2018gradient}, where $\tau < 1$ is the ratio of used blocks. Therefore, after the averaging step the variance is $\frac{1}{\tau n}$, which illustrates why setting any $\tau > \frac{1}{n}$ should not yield a significant speedup when compared to the choice $\tau = \frac1n$. It also indicates that it should be possible to throw away a $(1-\frac1n)$ fraction of blocks while keeping the same convergence rate.

\section{Contributions}

The goal of the above discussion was to introduce one of the ideas of this chapter in a gentle way. However, our independent sampling idea has immense consequences beyond the realm of {\tt GD}, as we show in the rest of the chapter. Let us summarize the contributions here:

\begin{itemize}
\item We show that the independent sampling idea can be coupled with variance reduction/{\tt SAGA} (see Section~\ref{sec:99_saga}),  {\tt SGD} for problem  \eqref{eq:99_problem}+\eqref{eq:99_stoch-f_i} (see Section~\ref{sec:99_sgd}),  acceleration (under mild assumption on stochastic gradients; see Section~\ref{sec:99_ABCDE}) and regularization/{\tt SEGA} (see Section~\ref{sec:99_sega}). We call the new methods {\tt ISAGA}, {\tt ISGD}, {\tt IASGD} and {\tt ISEGA}, respectively. We also develop {\tt ISGD} variant for asynchronous distributed optimization -- {\tt IASGD} (Section~\ref{sec:99_asynch}). 

\item We present two versions of  {\tt SAGA} coupled with {\tt IBCD}. The first one is for a distributed setting, where each machine owns a subset of data and runs a {\tt SAGA} iteration with block sampling locally, followed by aggregation. The second version is in a shared data setting, where each machine has access to all functions. This allows for linear convergence even if $\nabla f_i(x^*)\neq 0$. 
 
\item We show that when combined with {\tt IBCD}, the {\tt SEGA} trick (Chapter~\ref{sega}) leads to a method that enjoys a   linear rate for problems where $\nabla f_i(x^*)\neq 0$ and allows for more general objectives which may include a non-separable non-smooth regularizer. 
\end{itemize}

A comprehensive summary of all algorithms proposed  in this chapter is given in Table~\ref{tbl:99_algs}.

 \begin{table*}[ht]
 \footnotesize
\begin{center}
\small
\begin{tabular}{|c|c|c|c|c|c|c|c|}
\hline
{\bf \#}& {\bf Name}  & {\bf Origin}& { \bf  \begin{tabular}{c} $\nabla f_i(x^*)$ \\ $\neq0 $ \end{tabular}}&  {  \bf \begin{tabular}{c} Linear \\ rate \end{tabular}}  &{  \bf \begin{tabular}{c} Stochastic \\ gradient \end{tabular}} &{\bf Note} \\
 \hline
  \hline
\ref{alg:cd}   & {\tt IBCD} &  I+ {\tt CD}~\cite{rcdm} & \xmark & \cmark &\xmark & Simplest \\
  \hline
\ref{alg:sega}   & {\tt ISEGA}   & I + {\tt SEGA}~\cite{rcdm}& \cmark  & \cmark&\xmark & Allows prox \\
  \hline
\ref{alg:ibd}   & {\tt IBGD} &  I + {\tt GD}  & \xmark & \cmark &\xmark & Bernoulli \\
  \hline
\ref{alg:saga}   & {\tt ISAGA}  &  + {\tt SAGA}~\cite{saga}& \cmark &  \cmark& \cmark &Shared memory \\
  \hline
\ref{alg:saga_dist}   & {\tt ISAGA}    &I + {\tt SAGA}~\cite{saga}& \xmark & \cmark&\cmark & \\
  \hline
\ref{alg:sgd}   & {\tt ISGD}  & I + {\tt SGD}~\cite{robbins}&\cmark  &  \xmark &\cmark & + Non-convex\\
  \hline
\ref{alg:acc}   & {\tt IASGD}  &I + {\tt ASGD}~\cite{vaswani2019-overparam}& \cmark  &  \xmark&\cmark & Accelerated \\
 \hline
\ref{alg:asynch_sgd}   & {\tt IASGD}   & I + {\tt ASGD}~\cite{recht2011hogwild}& \cmark  & \xmark&\cmark & Asynchronous \\
\hline
\end{tabular}
\end{center}
\caption{Summary of all algorithms proposed in the chapter.}
\label{tbl:99_algs}
\end{table*}

\section{Practical implications and limitations \label{sec:99_practical}}
In this section, we outline some further limitations and practical implications of our framework.
\subsection{Main limitation}

The main limitation of this work is that  independent sampling does not generally result in a sparse aggregated update. Indeed, since each machine might sample a different subset of blocks, all these updates add up to a dense one, and this problem gets worse as $n$ increases, other things equal. For instance, if every parallel unit updates a single unique block\footnote{Assume $x$ is partitioned into several ``blocks'' of variables.}, the total number of updated blocks is equal $n$. In contrast, standard {\tt BCD}, one that samples the {\em same} block on each worker, would update a single block only. For simple linear problems, such as logistic regression, sparse updates allow for a fast implementation of {\tt BCD} via memorization of the residuals. However, this limitation is not crucial in common settings where broadcast is much faster than reduce.

\subsection{Practical implications}

The main body of this work focuses on theoretical analysis and on verifying our claims via experiments. However, there are several straightforward and important applications of our technique.

\paragraph{Distributed synchronous learning.} A common way to run a distributed optimization method is to perform a local update, communicate the result to a parameter server using a {\tt reduce} operation, and inform all workers using {\tt broadcast}. Typically, if the number of workers is significantly large, the bottleneck of such a system is communication. In particular, the {\tt reduce} operation takes much more time than {\tt broadcast} as it requires to add up different vectors computed locally, while {\tt broadcast} informs the workers about \textit{the same} data (see~\cite{mishchenko2019distributed} for a numerical validation that {\tt broadcast} is 10-20 times faster across a wide range of dimensions).
 Nevertheless, if every worker can instead send to the parameter server only $\tau = \frac{1}{n}$ fraction of the $d$-dimensional update, essentially the server node will receive just one full $d$-dimensional vector, and thus our approach can compete against methods like {\tt QSGD}~\cite{alistarh2017qsgd}, {\tt signSGD}~\cite{bernstein2018signsgd}, TernGrad~\cite{wen2017terngrad}, {\tt DGC}~\cite{lin2017deep} or {\tt ATOMO}~\cite{wang2018atomo}. In fact, our approach may completely remove the communication bottleneck. 

\paragraph{Distributed asynchronous learning.} The main difference with the synchronous case is that only one-to-one communications will be used instead of highly efficient {\tt reduce} and {\tt broadcast}. Clearly, the communication to the server will be much faster with $\tau=\frac{1}{n}$, so the main question is how to make the communication back fast as well. Hopefully, the parameter server can copy the current vector and send it using non-blocking communication, such as \textit{isend()} in {\tt MPI4PY}~\cite{dalcin2011parallel}. Then, the communication back will not prevent the server from receiving the new updates. We combine the {\tt IBCD} approach with asynchronous updates, which leads to a new method:  {\tt IASGD} (Algorithm~\ref{alg:asynch_sgd}).

\paragraph{Distributed sparse learning.} Large datasets, such as binary classification data from LibSVM, often have sparse gradients. In this case, the {\tt reduce} operation is not efficient and one needs to communicate data by sending positions of nonzeros and their values. Moreover, as we prove later, one can use independent sampling with $\ell_1$-penalty, which makes the problem solution sparse. In that case, only communication from a worker to the parameter server is slow, so both synchronous and asynchronous methods gain in performance.

\paragraph{Methods with local subproblems.} One can also try to extend our analysis to methods with exact block-coordinate minimization or primal-dual and proximal methods such as Point-{\tt SAGA}~\cite{defazio2016simple}, {\tt PDHG}~\cite{chambolle2011first}, {\tt DANE}~\cite{dane}, etc. There, by restricting ourselves to a subset of coordinates, we may obtain a subproblem that is easier to solve by orders of magnitude.

\paragraph{Block-separable problems within machines.} 
Given that the local problem on each machine is block coordinate-wise separable, partial derivative blocks can be evaluated $\frac{1}{\tau}$ times cheaper than the gradients. Thus, independent sampling improves scalability at no cost. Such problems can be obtained considering the dual problem, as is done in~\cite{COCOA+}, for example.

For a comprehensive list of frequently used notation that is specific to this chapter, see Table~\ref{tbl:notation_99} in the supplementary material.



\section{Independent block coordinate descent \label{sec:99_basic}}
Before presenting the algorithm, we shall assume smoothness and convexity of the objective.

\begin{assumption}\label{as:99_smooth_sc}
For every $i$, function $f_i$ is convex, $L$-smooth while function $f$ is $\mu$-strongly convex. 
\end{assumption}

 Let $\RR^d$  be partitioned into $m$ blocks, $u_1, \dotsc, u_m$, of arbitrary sizes, so that the parameter space is $\RR^{|u_1|}\times\dotsb \RR^{|u_m|}$. For any vector $x\in \RR^d$ and a set of blocks $U$ we denote by $x_U$ the vector that has the same coordinate as $x$ in the set of blocks $U$ and zeros elsewhere.

\subsection{The {\tt IBCD} algorithm}

In order to provide a quick taste of our results, we first present the {\tt IBCD} method described in the introduction and formalized as Algorithm~\ref{alg:cd}. 

\begin{algorithm}[h]
  \caption{Independent Block Coordinate Descent ({\tt IBCD})}
  \label{alg:cd}
\begin{algorithmic}[1]
\State {\bfseries Input: } {$x^0\in\RR^d$, partition of $\RR^d$ into $m$ blocks $u_1,\dotsc, u_m$, ratio of blocks to be sampled $\tau$, stepsize $\alpha$, \# of parallel units $n$}
  \For{$k=0,1,2,\dotsc$}
    \For{$i=1,\dotsc,n$ in parallel}
        \State Sample independently and uniformly a subset of $\tau m$ blocks $U_i^k \subseteq \{u_1, \dotsc, u_m\}$
        \State $x_i^{k+1} = x^k - \alpha  (\nabla f_i(x^k))_{U_i^k}$
    \EndFor
    \State $x^{k+1} = \frac{1}{n}\sumin x_i^{k+1}$
  \EndFor
\end{algorithmic}
\end{algorithm}

A key parameter of the method is $\frac{1}{m} \leq \tau \leq 1$ (chosen so that $\tau m$ is an integer), representing a fraction of blocks to be sampled by each worker. At iteration $k$, each machine independently samples a subset of $\tau m$ blocks $U_i^k \subseteq \{u_1,\dots,u_m\}$, uniformly at random. The $i$th worker then performs a subspace gradient step of the form $x_i^{k+1} = x^k - \alpha (\nabla f_i(x^k))_{U_i^k},$ where $\alpha>0$ is a stepsize.  Note that only coordinates of $x^k$ belonging to $U_i^k$ get updated. This is then followed by aggregating all $n$ gradient updates: $x^{k+1} = \frac{1}{n}\sum_i x_i^{k+1}$.

\subsection{Convergence of {\tt IBCD}}

Theorem~\ref{th:cd} provides a convergence rate for Algorithm~\ref{alg:cd}. Admittedly, the assumptions of Theorem~\ref{th:cd} are somewhat restrictive; in particular, we  require $\nabla f_i(x^*)=0$ for all $i$. However, this is necessary. Indeed, in general one can not expect to have $\sum_{i=1}^n(\nabla f_i(x^*))_{U_i}=0$ (which would be required for the method to converge to $x^*$) for independently sampled sets of blocks $U_i$ unless $\nabla f_i(x^*)=0$ for all $i$. As mentioned, the issue is resolved in Section~\ref{sec:99_sega} using the {\tt SEGA} trick from Chapter~\ref{sega}.
\begin{theorem}\label{th:cd} 
Suppose that Assumptions~\ref{as:99_smooth_sc} holds and $\nabla f_i(x^*) = 0$ for all $i$.\footnote{ The requirement of $\nabla f_i(x^*) = 0$ is only necessary for the plainest results; which we present to better explain the main idea of the chapter; and there are ways to go around it. In particular, in Section~\ref{sec:99_saga} we show that it can be dropped once the memory is shared among the machines.
Further, in Section~\ref{sec:99_sega} we show that $\nabla f_i(x^*) = 0$ can be dropped even in the fully distributed setup using the {\tt SEGA} trick. Lastly, $\nabla f_i(x^*) = 0$ is naturally satisfied in many applications. For example, in least squares setting $\min \|Ax-b\|^2$, it is equivalent to existence of $x^*$ such that $Ax^*=b$. On the other hand, current state-of-the-art deep learning models are often overparameterized so that they allow zero training loss, which is again equivalent to $\nabla f_i(x^*)=0$ for all $i$  (however, such problems are typically non-convex).
} For Algorithm~\ref{alg:cd} with $\alpha = \frac{n}{\tau n + 2(1 - \tau)}\frac{1}{2L}$ we have
\[
\EE\left[ \|x^{k} -x^*\|^2 \right] \leq \left(1-\frac{\mu}{2L} \frac{\tau n}{\tau n + 2(1-\tau)}\right)^k\|x^{0} -x^*\|^2.
\]
\end{theorem}
As a consequence of Theorem~\ref{th:cd}, we can choose $\tau$ as small as $\frac{1}{n}$ and get, up to a constant factor, the same convergence rate as gradient descent, as described next.
\begin{corollary}\label{cor:99_0893y83}
If $\tau=\frac1n$, the iteration complexity\footnote{Number of iterations to reach $\epsilon$ accurate solution.} of Algorithm~\ref{alg:cd} is ${\cal O} \left(\frac{L}{\mu} \log\frac{1}{\epsilon}\right)$. 
\end{corollary}
\subsection{Optimal block sizes}
If we naively use coordinates as blocks, i.e.\ all blocks have size equal 1, the update will be very sparse and the efficient way to send it is by providing positions of nonzeros and the corresponding values. If, however, we partition $\RR^d$ into blocks of size approximately equal $d/n$, then on average only one block will be updated by each worker. This means that it will be just enough for each worker to communicate the block number and its entries, which is twice less data sent than when using coordinates as blocks.

\section{Variance reduction \label{sec:99_saga}}
As the first extension of {\tt IBCD}, we inject independent coordinate sampling into {\tt SAGA}\footnote{Independent coordinate sampling is not limited to {\tt SAGA} and can be similarly applied to other variance reduction techniques.}~\cite{saga}, resulting in a new method we call {\tt ISAGA}. We consider two different settings for {\tt ISAGA}. The first one is standard distributed setup~\eqref{eq:99_problem}, where each $f_i$ is of the fine-sum form \eqref{eq:99_problem_saga_dist}. The idea is to run {\tt SAGA} with independent coordinate sampling locally on each worker, followed by aggregating the updates. However, as for {\tt IBCD}, we require $\nabla f_i(x^*) = 0$ for all $i$. The second setting is a {\em shared data/memory} setup; i.e., we assume that all workers have access to all functions from the finite sum.

\subsection{Shared data {\tt ISAGA}}
We now present a different setup for {\tt ISAGA} in which the requirement $\nabla f_i(x^*)=0$ is not needed. Instead of \eqref{eq:99_problem}, we rather solve the problem
\begin{align} \label{eq:99_problem_saga_sm}
      \min_{x\in \RR^d} \left\{  f(x) \eqdef \frac{1}{N}\sum_{j=1}^N \smo_{j}(x) \right \}
\end{align}
with $n$ workers all of which have  access to all  data describing $f$. Therefore, all workers can evaluate $\nabla \smo_{j}(x)$ for any $1\leq j\leq N$. Similarly to plain {\tt SAGA}, we remember the freshest gradient information in table $\mJ$, which we update as follows: 
\begin{equation}
\mJ_{j_i^k}^{k+1} = \mJ_{j_i^k}^k+ (\nabla \smo_{j_i^k}(x^k)- \mJ_{j_i^k}^k )_{U_i^k}, \quad \mJ_{j'}^{k+1} = \mJ_{j'}^k, \label{eq:99_saga_alpha_sm}
\end{equation}
where  $j_i^k$ is the  index  sampled at iteration $k$ by machine $i$, and $j'$ refers to all indices that were not sampled at iteration $k$ by any machine. The iterate updates within each machine are taken only on a sampled set of coordinates, i.e.,  
$
x_i^{k+1} = x^k - \alpha  (\nabla \smo_{j_i^k}(x^k) - \mJ_{j_i^k}^k + \overline \mJ^k)_{U_i^k}.
$
where $\overline \mJ^k$ stands for the average of all $\mJ$, and thus it is a delayed estimate of $\nabla f(x^k)$. Lastly, we set the next iterate as the average of proposed iterates by each machine  $x^{k+1} = \frac1n \sumin x_i^{k+1} $. The formal statement of the algorithm is given in the supplementary as Algorithm~\ref{alg:saga}.
\begin{algorithm}[h]
  \caption{{\tt ISAGA} with shared data}
  \label{alg:saga}
\begin{algorithmic}[1]
\State{\bfseries Input: }{$x^0\in\RR^d$, $\mJ_1^0,\dotsc, \mJ_N^0$ partition of $\RR^d$ into $m$ blocks $u_1,\dotsc, u_m$, ratio of blocks to be sampled $\tau$, stepsize $\alpha$, \# parallel units $n$}
  \State Set $\overline \mJ^0 \eqdef \frac{1}{N}\sumin \mJ_i^0$
  \For{$k=0,1,2,\dotsc$}
          \State Sample uniformly set of indices $\{j_{1}^k, \dots, j_n^k\} \subseteq  \{1,\dots, N\}$ without replacement
    \For{$i=1,\dotsc,n$ in parallel}
        \State Sample independently and uniformly a subset of $\tau m$ blocks $U_t^i$
        \State $x_i^{k+1} = x^k - \alpha  (\nabla \smo_{j_i^k}(x^k) - \mJ_{j_i^k}^k + \overline \mJ^k)_{U_i^k}$
        \State $(\mJ_{j_i^k}^{k+1})_{U_i^k} = \mJ_{j_i^k}^k+ (\nabla \smo_{j_i^k}(x^k)- \mJ_{j_i^k}^k )_{U_i^k}$
    \EndFor
    \State For $j\not\in \{j_{1}^k, \dots, j_n^k\} $ set $(\mJ_{j}^{k+1})= \mJ_{i}^k$
    \State $x^{k+1} = \frac{1}{n}\sumin x_i^{k+1}$
    \State $\overline \mJ^{k+1} = \frac{1}{n}\sum_{j=1}^N \mJ_j^{k+1}$
  \EndFor
\end{algorithmic}
\end{algorithm}

\begin{theorem}\label{th:saga_shared}
 Suppose that function $f$ is $\mu$-strongly convex and each $\smo_i$ is $L$ smooth and convex. If $\alpha\le \frac{1}{L\left(\frac{3}{n} + \tau\right)}$, then for iterates of Algorithm~\ref{alg:saga} we have
    \begin{align*}
        \E{ \|x^k - x^*\|^2} \le (1 - \vartheta)^k\left(\|x^0 - x^*\|^2 + c \alpha^2 \Psi^0\right),
    \end{align*}
    where $\Psi^0\eqdef \sum_j \|\mJ_{j}^0 - \nabla \smo_{j}(x^*)\|^2$, $\vartheta\eqdef \tau\min\left\{\alpha\mu, \frac{ n}{N} - \frac{2}{nNc} \right\}\ge 0$ and $c\eqdef \frac{1}{n}( \frac{1}{\alpha L} - \frac{1}{n} - \tau) > 0$.
\end{theorem}

As in Section~\ref{sec:99_saga_dist}, the choice $\tau = \frac{1}{n}$ yields a convergence rate which is, up to a constant factor, the same as the convergence rate of {\tt SAGA}. Therefore, Algorithm~\ref{alg:saga} enjoys the desired parallel linear scaling without extra assumptions. Corollary~\ref{cor:99_saga_sm} formalizes the claim.

\begin{corollary} \label{cor:99_saga_sm}
  Consider the setting from Theorem~\ref{th:saga_shared}.  Set $\tau = \frac{1}{n}$ and $\alpha = \frac{n}{5L}$. Then $c=\frac{3}{n^2}$, $\rho = \min \left\{ \frac{\mu}{5L}, \frac{1}{3N}\right\}$ and the complexity of Algorithm~\ref{alg:saga} is $O\left(\max\left\{\frac{L}{\mu}, N \right\}\log\frac{1}{\varepsilon} \right)$.
\end{corollary}

\subsection{Distributed {\tt ISAGA} \label{sec:99_saga_dist}}
In this section we consider problem \eqref{eq:99_problem} with $f_i$ of the finite-sum structure \eqref{eq:99_problem_saga_dist}.  Just like {\tt SAGA}, every machine remembers the freshest gradient information of all local functions (stored in arrays $\mJ_{ij}$), and updates them once a new gradient information is observed. Given that index $j_i^k$ is sampled on $i$th machine at iteration $k$, the iterate update step within each machine is taken only on a sampled set of coordinates: 
\[
x_i^{k+1} = x^k - \alpha  (\nabla f_{ij_i^k}(x^k) - \mJ_{ij_i^k}^k + \overline \mJ_i^k)_{U_i^k} .
\]
Above, $\overline \mJ_i^k$ stands for the average of $\mJ$ variables on $i$th machine, i.e.\ it is a delayed estimate of $\nabla f_i(x^k)$. Since the new gradient information is a set of partial derivatives of $\nabla f_{ij_i^k}(x^k)$, we shall update
\begin{equation}
\mJ_{ij}^{k+1} = 
\left\{
                \begin{array}{ll}
                    \mJ_{ij}^k+ (\nabla f_{ij}(x^k)- \mJ_{ij}^k )_{U_i^k}   & j=j_i^k\\
                       \mJ_{ij}^k & j\neq j_i^k
                \end{array}
     \right.\label{eq:99_saga_alpha_dist}
\end{equation}
Lastly, the local results are aggregated. See Algorithm~\ref{alg:saga_dist} for details. 

\begin{algorithm}[h]
  \caption{Distributed {\tt ISAGA}}
  \label{alg:saga_dist}
\begin{algorithmic}[1]
\State{\bfseries Input: }{$x^0\in\RR^d$, \# parallel units $n$, $i$th unit owns $l$ functions $f_{i1} ,\dots, f_{il}$, partition of $\RR^d$ into $m$ blocks $u_1,\dotsc, u_m$, ratio of blocks to be sampled $\tau$, stepsize $\alpha$, initial vectors $\mJ_{ij}^0 \in \RR^d$ for $1\leq i \leq n, 1\leq j\leq l$  }
  \State Set $\overline \mJ^0 \eqdef \frac{1}{N}\sumin \mJ_i^0$
  \For{$k=0,1,2,\dotsc$}
    \For{$i=1,\dotsc,n$ in parallel}
         \State Sample independently \& uniformly $ j_i^k\in [l]$ 
        \State Sample independently \& uniformly a subset of $\tau m$ blocks $U_t^i$
        \State $x_i^{k+1} = x^k - \alpha  (\nabla f_{ij_i^k}(x^k) - \mJ_{ij_i^k}^k + \overline \mJ_i^k)_{U_i^k}$
        \State $\mJ_{ij^k}^{k+1} = \mJ_{ij_i^k}^k+ (\nabla f_{ij_i^k}(x^k)- \mJ_{ij_i^k}^k )_{U_i^k}$
        \State For any $j\neq j_i^k$ set $\mJ_{ij}^{k+1} = \mJ_{ij}^k$
        \State $\overline \mJ^{k+1} = \frac{1}{l}\sum_{j=1}^l \mJ_{ij}^{k+1}$
    \EndFor
    \State $x^{k+1} = \frac{1}{n}\sumin x_i^{k+1}$
  \EndFor
\end{algorithmic}
\end{algorithm}

The next result provides a convergence rate of distributed {\tt ISAGA}.
\begin{theorem}\label{th:saga_dist}
    Suppose that Assumption~\ref{as:99_smooth_sc} holds and $\nabla f_i(x^*)  = 0$ for all $i$. If $\alpha\le \frac{1}{L\left(\frac{3}{n} + \tau\right)}$, for iterates of distributed {\tt ISAGA} we have
    \begin{align*}
        \EE \|x^k - x^*\|^2 \le (1 - \vartheta )^k\left(\|x^0 - x^*\|^2 + c \alpha^2 \Psi^0\right),
    \end{align*}
    where $\Psi^0\eqdef\sum_{i=1}^n\sum_{j=1}^l \|\mJ_{ij}^k - \nabla f_{ij}(x^*)\|^2$, $\vartheta\eqdef \tau\min\left\{\alpha\mu, \frac{ 1}{l} - \frac{2}{n^2lc} \right\}\ge 0$ and $c\eqdef \frac{1}{n} ( \frac{1}{\alpha L} - \frac{1}{n} - \tau) > 0$.
\end{theorem}

The choice $\tau = n^{-1}$ yields a convergence rate which is, up to a constant factor, the same as convergence rate of original {\tt SAGA}. Thus, distributed {\tt ISAGA} enjoys the desired parallel linear scaling. Corollary~\ref{cor:99_saga_dist} formalizes this claim. 

\begin{corollary}\label{cor:99_saga_dist}
 Consider the setting from Theorem~\ref{th:saga_dist}. Set $\tau = \frac{1}{n}$ and $\alpha = \frac{n}{5L}$. Then $c=\frac{3}{n^2}$, $\rho = \min \left\{ \frac{\mu}{5L}, \frac{1}{3nl}\right\}$ and the complexity of distributed {\tt ISAGA} is \[{\cal O}\left(\max \left\{\frac{L}{\mu}, nl \right\}\log\frac{1}{\varepsilon} \right).\]
\end{corollary}

\section{{\tt SGD} \label{sec:99_sgd}}

In this section, we apply independent sampling in a setup with a stochastic objective. In particular, we consider problem~\eqref{eq:99_problem} where $f_i$ is given as an expectation; see \eqref{eq:99_stoch-f_i}.
We assume we have access to a stochastic gradient oracle which, when queried at $x^k$,  outputs a random vector $g_i^k$ whose mean is $\nabla f_i(x^k)$:  $\EE g_i^k = \nabla f_i(x^k)$. 

Our proposed algorithm---{\tt ISGD}---evaluates a subset of stochastic partial derivatives for the local objective and takes a step in the given direction for each machine. Next, the results are averaged and followed by the next iteration. We stress that the coordinate blocks have to be sampled independently within each machine.

\begin{algorithm}[h]
  \caption{{\tt ISGD}}
  \label{alg:sgd}
\begin{algorithmic}[1]
\State{\bfseries Input: }{$x^0\in\RR^d$, partition of $\RR^d$ into $m$ blocks $u_1,\dotsc, u_m$, ratio of blocks to be sampled $\tau$, stepsize sequence $\{\alpha^k\}_{k=1}^\infty$, \# parallel units $n$}
  \For{$k=0,1,2,\dotsc$}
    \For{$i=1,\dotsc,n$ in parallel}
        \State Sample independently and uniformly a subset of $\tau m$ blocks $U_i^k \subseteq \{u_1, \dotsc, u_m\}$
        \State Sample blocks of stochastic gradient $(g_i^k)_{U_i^k}$  such that $\EE [g_i^k \, |\, x^k] = \nabla f_i(x^k)$
        \State $x_i^{k+1} = x^k - \alpha^k  (g_i^k)_{U_i^k}$
    \EndFor
    \State $x^{k+1} = \frac{1}{n}\sumin x_i^{k+1}$
  \EndFor
\end{algorithmic}
\end{algorithm}

In order to establish a convergence rate of {\tt ISGD}, we shall assume boundedness of stochastic gradients for each worker.

\begin{assumption}\label{as:99_bounded_noise}
   Consider a sequence of iterates $\{x^k \}_{k=0}^\infty$ of Algorithm~\ref{alg:sgd}. Assume that $g_i^k$ is an unbiased estimator of $\nabla f_i(x^k)$ satisfying
   $
        \EE \|g_i^k - \nabla f_i(x^k)\|^2 \le \sigma^2.
$
    
\end{assumption}
\begin{assumption}\label{as:99_bounded_noise_at_opt}
   Stochastic gradients of function $f_i$ have bounded variance at the optimum of $f$:
   $
        \EE \|g_i - \nabla f_i(x^*)\|^2 \le \sigma^2,
$
    where $g_i$ is a random vector such that $\EE g_i = \nabla f_i(x^*)$.
\end{assumption}

Next, we present the convergence rate of Algorithm~\ref{alg:sgd}. Since {\tt SGD} is not a variance reduced algorithm, it does not enjoy a linear convergence rate and one shall use decreasing step sizes. As a consequence, it is not required to assume that $\nabla f_i(x^*)= 0$ for all $i$ since there is no variance reduction property to be broken. 

\begin{theorem} \label{th:sgd}
Let Assumptions~\ref{as:99_smooth_sc} and~\ref{as:99_bounded_noise} hold. If $\alpha^k = \frac{1}{a + ck}$, where $a= 2\left(\tau + \frac{2(1 - \tau)}{n} \right)L$, $c= \frac14 \mu\tau$, then for Algorithm~\ref{alg:sgd} we can upper bound $ \EE [f(\hat x^k) - f(x^*) ]$ by
    \begin{align*}
 \frac{a^2\left(1 - \frac{\tau\mu}{a}\right)\|x^0 - x^*\|^2}{\tau(k+1)a+ \frac{c\tau}{2}k(k+1)}  +  \frac{\sigma^2 + (1 - \tau)\frac{2}{n}\sumin\|\nabla f_i(x^*)\|^2 }{n\left(1+\frac1k\right)a+ \frac{nc}{2}(k+1)},
    \end{align*}
where $\hat x^k\eqdef \frac{1}{(k+1)a+ \frac{c}{2}k(k+1)}\sum_{t=0}^k (\alpha^\qwerty)^{-1} x^\qwerty$.
\end{theorem}

Note that the residuals decrease as ${\cal O}(k^{-1})$, which is a behavior one expects from standard {\tt SGD}. Moreover, the leading complexity term scales linearly: if the number of workers $n$ is doubled, one can afford to halve $\tau$ to keep the same complexity. 
\begin{corollary}\label{cor:99_sgd}
Consider the setting from Theorem~\ref{th:sgd}. Then, iteration complexity of Algorithm~\ref{alg:sgd} is 
\[
{\cal O} \left( 
\frac{\sigma^2 + \frac{1}{n}\sumin\|\nabla f_i(x^*)\|^2 }{n\tau \mu \epsilon}
\right).
\]
\end{corollary}

Although problem~\eqref{eq:99_problem} explicitly assumes convex $f_i$, we also consider a non-convex extension, where smoothness of each individual $f_i$ is not required either. Theorem~\ref{th:sgd_ncvx} provides the result. 

\begin{theorem}[Non-convex rate]\label{th:sgd_ncvx}
    Assume $f$ is $L$ smooth, Assumption~\ref{as:99_bounded_noise} holds and for all $x\in\RR^d$ the difference between gradients of $f$ and $f_i$'s is bounded: $\frac{1}{n}\sumin \|\nabla f(x) - \nabla f_i(x)\|^2\le \nu^2$ for some constant $\nu\ge 0$. If $\hat x^k$ is sampled uniformly from $\{x^0, \dotsc, x^k\}$, then for Algorithm~\ref{alg:sgd} we have
    \begin{align*}
        \EE \|\nabla f(\hat x^k)\|^2
        \le \frac{\frac{f(x^0) - f^*}{k\tau\alpha} + \alpha  L\frac{\left(1 - \tau\right)\nu^2+\frac12\sigma^2}{n}}{1 - \frac{\alpha \tau L}{2} - \alpha L\left(1 - \tau\right)\frac{1}{n}}.
    \end{align*}
    \end{theorem}

Again, the convergence rate from Theorem~\ref{th:sgd_ncvx} scales almost linearly with $\tau$: with doubling the number of workers one can afford to halve $\tau$ to keep essentially the same guarantees. Note that if $n$ is sufficiently large, increasing $\tau$ beyond a certain threshold does not improve convergence. This is a slightly weaker conclusion to the rest of our results where increasing $\tau$ beyond $n^{-1}$ might still offer speedup. The main reason behind this is  the fact that {\tt SGD} may be  noisy enough on its own  to still benefit from the averaging step.

    \begin{corollary}\label{cor:99_ncvx}
 Consider the setting from Theorem~\ref{th:sgd_ncvx}. i) Choose $\tau\ge \frac{1}{n}$ and $\alpha = \frac{\sqrt{n}}{L\sqrt{\tau k}} \le \frac{1}{2L\left(\tau/2 + (1 - \tau)/n \right)}$. Then \[\EE \|\nabla f(\hat x^k)\|^2 \le \frac{2}{\sqrt{k\tau n}}\left(\frac{f(x^0) - f^*}{L} + (1 - \tau)\nu^2\right) = O\left(\frac{1}{\sqrt{k}}\right).\] ii) For any $\tau$ there is sufficiently large $n$ such that choosing $\alpha = {\cal O}\left( \frac{\epsilon}{\tau L^2}\right)$ yields complexity ${\cal O} \left( \frac{L^2}{\epsilon^2}\right)$. The complexity does not improve significantly when $\tau$ is increased. 

    \end{corollary}

  \section{Acceleration~\label{sec:99_ABCDE}}
Here we describe an accelerated variant of {\tt IBCD} in the sense of~\cite{nesterov83}. In fact, we will do something more general and accelerate {\tt ISGD}, obtaining the {\tt IASGD} algorithm. We again assume that machine $i$ owns $f_i$, which is itself a stochastic objective as in~\eqref{eq:99_stoch-f_i} with an access to an unbiased stochastic gradient $g^k$ every iteration: $\EE g_i^k = \nabla f_i(x^k)$. A key assumption for the accelerated {\tt SGD} used to derive the best known rates~\cite{vaswani2019-overparam} is so the called strong growth of the unbiased gradient estimator.
\begin{definition}
Function $\phi(x)=\EE_\zeta \phi(x,\zeta)$ satisfies the strong growth condition with parameters $\rho, \sigma^2$, if for all $x$ we have
\[
\EE_{\zeta} \| \nabla \phi(x,\zeta)\|^2\leq \rho \| \nabla \phi(x)\|^2 +\sigma^2.
\]  
\end{definition} 

In order to derive a strong growth property of the gradient estimator coming from the independent block coordinate sampling, we require a strong growth condition on $f$ with respect to $f_1, \dots, f_n$ and also a variance bound on stochastic gradients of each individual $f_i$.

\begin{assumption} \label{as:99_strong_growth}
Function $f$ satisfies the strong growth condition with respect to $f_1, \dots, f_n$ : 
\begin{equation}\label{eq:99_acc_sg_f}
\frac{1}{n}\sum_{i=1}^n \|\nabla f_i(x)\|^2 \leq  \tilde{\rho} \|\nabla f(x) \|^2+  \tilde{\sigma}^2.
\end{equation}
Similarly, given that $g_i = g_i(x)$ provides an unbiased estimator of $\nabla f_i(x)$, i.e.\ $\EE g_i =\nabla f_i(x)$, variance of $g_i$ is bounded as follows for all $i$:
\begin{equation}\label{eq:99_acc_sg_fi}
 \Var\left[  g_i\right] \leq  \bar{\rho} \|\nabla f_i(x) \|^2+  \bar{\sigma}^2.
\end{equation}
\end{assumption}

Note that the variance bound \eqref{eq:99_acc_sg_fi} is weaker than the strong growth property as we always have $ \Var\left[  g_i\right]  \leq  \EE\left[ \| g_i \|^2\right] $.

Given that Assumption~\ref{as:99_strong_growth} is satisfied, we derive a strong growth property for the unbiased gradient estimator $q \eqdef \frac{1}{n\tau}\sum_{i=1}^n(\nabla g_i)_{U_i}$ in Lemma~\ref{lem:99_stronggrowth}. Next, {\tt IASGD} is nothing but the scheme from~\cite{vaswani2019-overparam} applied to stochastic gradients $q$. For completeness, we state {\tt IASGD} as Algorithm~\ref{alg:acc}.

\begin{algorithm}[h]
  \caption{{\tt IASGD}}
  \label{alg:acc}
\begin{algorithmic}[1]
\State{\bfseries Input: } {Starting point $y^0=v^0\in\RR^d$, partition of $\RR^d$ into $m$ blocks $u_1,\dotsc, u_m$, ratio of blocks to be sampled $\tau$, stepsize $\alpha$, number of parallel units $n$, acceleration parameter sequences $\{a,b,\eta\}_{k=0}^\infty$}
  \For{$k=0,1,2,\dotsc$}
         \State $x^k  = a^k v^k + (1 - a^k)y^k $ 
    \For{$i=1,\dotsc,n$ in parallel}
        \State Sample independently and uniformly a subset of $\tau m$ blocks $U_i^k \subset \{u_1, \dotsc, u_m\}$
        \State Sample blocks of stochastic gradient $(g_i^k)_{U_i^k}$  such that $\EE [g_i^k \, |\, x^k] = \nabla f_i(x^k)$
       \EndFor
        \State  $q^k  =\frac{1}{n\tau}\sumin (g_i^k)_{U_i^k}$ 
    	   \State  $y^{k+1}  = x^k -\alpha  q^k$ 
		\State $v^{k+1} = b^k v^k + (1 - b^k)x^k  - \eta^k \gamma q^k$.
  \EndFor
\end{algorithmic}
\end{algorithm}

\begin{lemma}\label{lem:99_stronggrowth}
Suppose that Assumption~\ref{as:99_strong_growth} is satisfied. Then, we have
$
\EE\left[\|q\|^2\right] \leq \hat{\rho} \|\nabla f(x) \|^2 + \hat{\sigma}^2
$
for
\begin{eqnarray}
\label{eq:99_acc_rho}
\hat{\rho} &\eqdef& \left(1+ \frac{\tilde{\rho}}{n}   \left(\frac1\tau-1+\frac{\bar{\rho}}{\tau} \right) \right), 
\\
\hat{\sigma}^2 &\eqdef& \frac{\bar{\sigma}^2}{n\tau} + \frac{\tilde{\sigma}^2}{n}\left(\frac1\tau-1+\frac{\bar{\rho}}{\tau} \right). 
\label{eq:99_acc_sigma}
\end{eqnarray}
\end{lemma}

It remains to use the stochastic gradient $q$ (with the strong growth bound from Lemma~\ref{lem:99_stronggrowth}) as a gradient estimate in~\cite{vaswani2019-overparam}[Theorem 6], which we restate as Theorem~\ref{th:accelerated} for completeness. 

\begin{theorem}\label{th:accelerated}
Suppose that $f$ is $L$ smooth, $\mu$ strongly convex and Assumption~\ref{as:99_strong_growth} holds. Then, for a specific choice of parameter sequences $\{a,b, \eta\}_{k=0}^\infty$ (See~\cite{vaswani2019-overparam}[Theorem 6] for details), iterates of {\tt IASGD} admit an upper bound on $ \EE \left[f(x^{k+1}) \right]-f(x^*)$ of the form
\begin{eqnarray*}
 \left(1-\sqrt{\frac{\mu}{L\hat{\rho}^2}} \right)^k \left(f(x^0)-f(x^*) + \frac{\mu}{2}\|x^0-x^* \|^2\right) + \frac{\hat{\sigma}^2}{\hat{\rho} \sqrt{L\mu}}.
\end{eqnarray*}

\end{theorem}

The next corollary provides a complexity of Algorithm~\ref{alg:acc} in a simplified setting where $\bar{\sigma}^2=\tilde{\sigma}^2=0$. Note that $\tilde{\sigma}^2=0$ implies $\nabla f_i(x^*) = 0$ for all $i$. It again shows a desired linear scaling: given that we double the number of workers, we can halve the number of  blocks to be evaluated on each machine and still keep the same convergence guarantees. It also shows that increasing $\tau$ beyond $\frac{\tilde{\rho}\bar{\rho}}{n}$ does not improve the convergence significantly. 
\begin{corollary}
Suppose that $\bar{\sigma}^2=\tilde{\sigma}^2=0$. Then, complexity of {\tt IASGD} is \[{\cal O}\left(\frac{1}{\hat{\rho}}\sqrt{\frac{\mu}{L}}\log\frac1\epsilon \right) ={\cal O}\left(\frac{1}{1+\frac{\tilde{\rho}}{\tau n}(1+\bar{\rho})}\sqrt{\frac{\mu}{L}}\log\frac1\epsilon \right)  .\]
\end{corollary}

Theorem~\ref{th:accelerated} shows an accelerated rate for strongly convex functions applying~\cite[Theorem 6] {vaswani2019-overparam} to the bound. A non-strongly convex rate can be obtained analogously from~\cite[Theorem 7]{vaswani2019-overparam}.

\section{Beyond interpolation without shared data and regularization \label{sec:99_sega}}
For this section only, let us consider a regularized objective of the form
\begin{align} \label{eq:99_problem_sega}
      \min_{x\in \RR^d} \left \{  f(x) \eqdef \frac{1}{n}\sum_{i=1}^n f_i(x) + \psi(x) \right \}, 
\end{align}
where $\psi$ is a closed convex regularizer such that its proximal operator,
$$   \prox_{\alpha\psi}(x) \eqdef \argmin_{y} \left\{\psi(y) + \frac{1}{2\alpha}\|y - x\|^2  \right\},
$$ 
is computable. In this section we propose {\tt ISEGA}:  an independent sampling variant of {\tt SEGA}. We do this in order to both i) avoid assuming $\nabla f_i(x^*)=0$ (while keeping linear convergence) and ii) allow for $R$. Original {\tt SEGA} learns gradients $\nabla f(x^k)$ from sketched gradient information via the so called sketch-and-project process~\cite{gower2015randomized}, constructing a vector sequence $h^k$. In {\tt ISEGA} on each machine $i$ we iteratively construct a sequence of vectors $h^k_i$ which play the role of estimates of $\nabla f_i(x^k)$. This is done via the following rule: 
\begin{equation} \label{eq:99_sega_h}
h_i^{k+1} = h_i^k + (\nabla f_i(x^k)- h_i^k)_{U_i^k}.
\end{equation}
The key idea is again that these vectors are created from random blocks independently sampled on each machine. Next, using $h^k$, {\tt SEGA} builds an unbiased gradient estimator $g_i^k$ of $\nabla f_i(x^k)$ as follows:
\begin{equation} \label{eq:99_sega_g}
   g_i^k = h_i^k + \frac{1}{\tau} (\nabla f_i(x^k) - h_i^k)_{U_i^k}. 
\end{equation}
Then, we average the vectors $g_i^k$ and take a proximal step.

Unlike coordinate descent, {\tt SEGA} (or {\tt ISEGA}) is not limited to separable proximal operators since, as follows from our analysis, $h_i^k\to \nabla f_i(x^*)$. Therefore, {\tt ISEGA} can be seen as a variance reduced version of {\tt IBCD} for problems with non-separable regularizers. 

In order to be consistent with the rest of the chapter, we only develop a simple variant of  {\tt ISEGA} (Algorithm~\ref{alg:sega}) in which we consider block coordinate sketches with uniform probabilities. While is possible to develop the theory in full generality (done in Chapter~\ref{jacsketch}) we avoid this for the sake of simplicity. 

\begin{algorithm}[h]
  \caption{{\tt ISEGA}}\label{alg:sega}
\begin{algorithmic}[1]
\State{\bfseries Input: }{$x^0\in \RR^d$, initial gradient estimates $h_1^0, \dotsc, h_n^0\in \RR^d$, partition of $\RR^d$ into $m$ blocks $u_1,\dotsc, i_m$, ratio of blocks to be sampled $\tau$, stepsize $\alpha$,  \# parallel units $n$}
  \For{$k=0,1,2,\dotsc$}
    \For{$i=1,\dotsc,n$ in parallel}
        \State Sample independently and uniformly a subset of $\tau m$ blocks $U_i^k$
        \State $g_i^k = h_i^k + \frac{1}{\tau} (\nabla f_i(x^k) - h_i^k)_{S_i^k}$ 
        \State  $h_i^{k+1} = h_i^k + \tau (g_i^k - h^k)$  
    \EndFor
    \State $x^{k+1} = \proxR\left( x^k - \alpha \frac{1}{n}\sumin g_i^{t} \right)$  
  \EndFor
\end{algorithmic}
\end{algorithm}

We next present the convergence rate of {\tt ISEGA} (Algorithm~\ref{alg:sega}). 

\begin{theorem}\label{thm:99_sega}
Suppose Assumption~\ref{as:99_smooth_sc} holds and choose stepsize
$$ \alpha = \min \left\{ \frac{1}{4L\left( 1+\frac{1}{n\tau }\right)}, \frac{1}{\frac{\mu}{\tau}+ \frac{4L}{n\tau}} \right\}
.$$  Then Algorithm~\ref{alg:sega} 
satisfies
\[
 \EE[\|x^k - x^*\|^2] \leq (1-\alpha\mu)^k\Lgen^0,
\]
where the Lyapunov function is given by $\Lgen^0 \eqdef \|x^0 - x^*\|^2 + \frac{ \alpha}{2L\tau n} \sum \limits_{i=1}^n\|h^0 - \nabla f(x^*)\|^2$.
\end{theorem}
Note that if the condition number of the problem is not too small so that $n= {\cal O}\left(L/\mu\right)$ (which is usually the case in practice), {\tt ISEGA} scales linearly in the parallel setting. In particular, when doubling the number of workers, each worker can afford to evaluate only half of  the block partial derivatives while keeping the same convergence speed. Moreover, setting $\tau = \frac1n$, the rate corresponds, up to a constant factor, to the rate of gradient descent. Corollary~\ref{cor:99_sega} states the result.

\begin{corollary}\label{cor:99_sega}
Consider the setting from Theorem~\ref{thm:99_sega}. Suppose that $\frac{L}{\mu}\geq n$ and choose $\tau = \frac1n$. The complexity of Algorithm~\ref{alg:sega} is ${\cal O} \left(\frac{L}{\mu}\log\frac1\epsilon \right)$.
\end{corollary}

\begin{remark}
Parallel implementation Algorithm~\ref{alg:sega} would be to always send $(\nabla f_i(x^k))_{U_i^k}$ to the server; which keeps updating vector $h^k$ and takes the prox step. 
\end{remark}

\section{Experiments}
In this section, we numerically verify our theoretical claims. Recall that there are various settings where it is possible to make practical experiments (see Section~\ref{sec:99_practical}), however, we do not restrain ourselves to any of them in order to deliver as clear a message as possible.  

We present exhaustive numerical experiments to verify the theoretical claims of the chapter. The experiments are performed in a simulated environment instead of the honestly distributed setup, as we only aim to verify the iteration complexity of proposed methods.

First, in Section~\ref{sec:99_exp_quad} provides the simplest setting in order to gain the best possible insight -- Algorithm~\ref{alg:cd} is tested on the artificial quadratic minimization problem. We compare Algorithm~\ref{alg:cd} against both gradient descent (GD) and standard CD (in our setting: when each machine samples the same subset of coordinates). We also study the effect of changing $\tau$ on the convergence speed. 

In the remaining parts, we consider a logistic regression problem on LibSVM data~\cite{chang2011libsvm}. Recall that logistic regression problem is given as

\begin{equation}
\label{eq:99_logreg}
f(x)\eqdef \frac1N \sum_{j=1}^N \left( \log \left(1+\exp\left(\mA_{j,:}x\cdot  b_j\right) \right)+\frac{\lambda}{2} \| x\|^2\right),
\end{equation}
where $\mA$ is data matrix and $b$ is vector of data labels: $b_j\in \{-1,1 \}$\footnote{
The datapoints (rows of $\mA$) have been normalized so that each is of norm $1$. Therefore, each $f_i$ is $\frac14$ smooth in all cases. We set regularization parameter as $\lambda = 0.00025$ in all cases. 
}. In the distributed scenario (everything except of Algorithm~\ref{alg:saga}), we imitate that the data is evenly distributed to $n$ workers (i.e.\ each worker owns a subset of rows of $\mA$ and corresponding labels, all subsets have almost the same size). 

As our experiments are not aimed to be practical at this point (we aim to properly prove the conceptual idea), we consider multiple of rather smaller datasets: \texttt{a1a} ($d=123, n =1605$), \texttt{mushrooms} ($d=112, n =8124$), \texttt{phishing} ($d=68, n = 11055$), \texttt{w1a} ($d=300, n =2477$). The experiments are essentially of 2 types: one shows that setting $n\tau=1$ does not significantly violate the convergence of the original method. In the second type of experiments we study the behavior for varying $\tau$, and show that beyond certain threshold, increasing $\tau$ does not significantly improve the convergence. The threshold is smaller as $n$ increases, as predicted by theory.

\subsection{Simple, well understood experiment \label{sec:99_exp_quad}}

In this section we study the simplest possible setting -- we test the behavior of Algorithm~\ref{alg:cd} on a quadratic minimization problem with artificial data. The considered quadratic objective is set as 
\begin{equation} f_i(x) \eqdef \frac12 x^\top \mM_ix,  \quad  \mM_i \eqdef vv^\top +  \frac{ \left(\mI- vv^\top\right) \mA_i\mA_i^\top \left(\mI- vv^\top\right) }{ \lambda_{\max} \left( \mA_i\mA_i^\top \right) }, \quad  v=\frac{v'}{\|v'\|}, \label{eq:99_quadratic}
\end{equation}
where entries of $v'\in \RR^d$ and $\mA_i\in \RR^{d\times \quadb}$ are sampled independently from standard normal distribution.

In the first experiment (Figure~\ref{fig:99_artif_1}), we compare Algorithm~\ref{alg:cd} with $n\tau=1$ against gradient descent (GD) and two versions of coordinate descent - a default version with stepsize $\frac1L$, and a coordinate descent with importance sampling (sample proportionally to coordinate-wise smoothness constants) and optimal step sizes (inverse of coordinate-wise smoothness constants). In all experiments, gradient descent enjoys twice better iteration complexity than Algorithm~\ref{alg:cd} which is caused by twice larger stepsize. However, in each case, Algorithm~\ref{alg:cd} requires fewer iterations to CD with importance sampling, which is itself significantly faster to plain CD.

\begin{figure}[H]
\centering
\begin{minipage}{0.3\textwidth}
  \centering
\includegraphics[width =  \textwidth ]{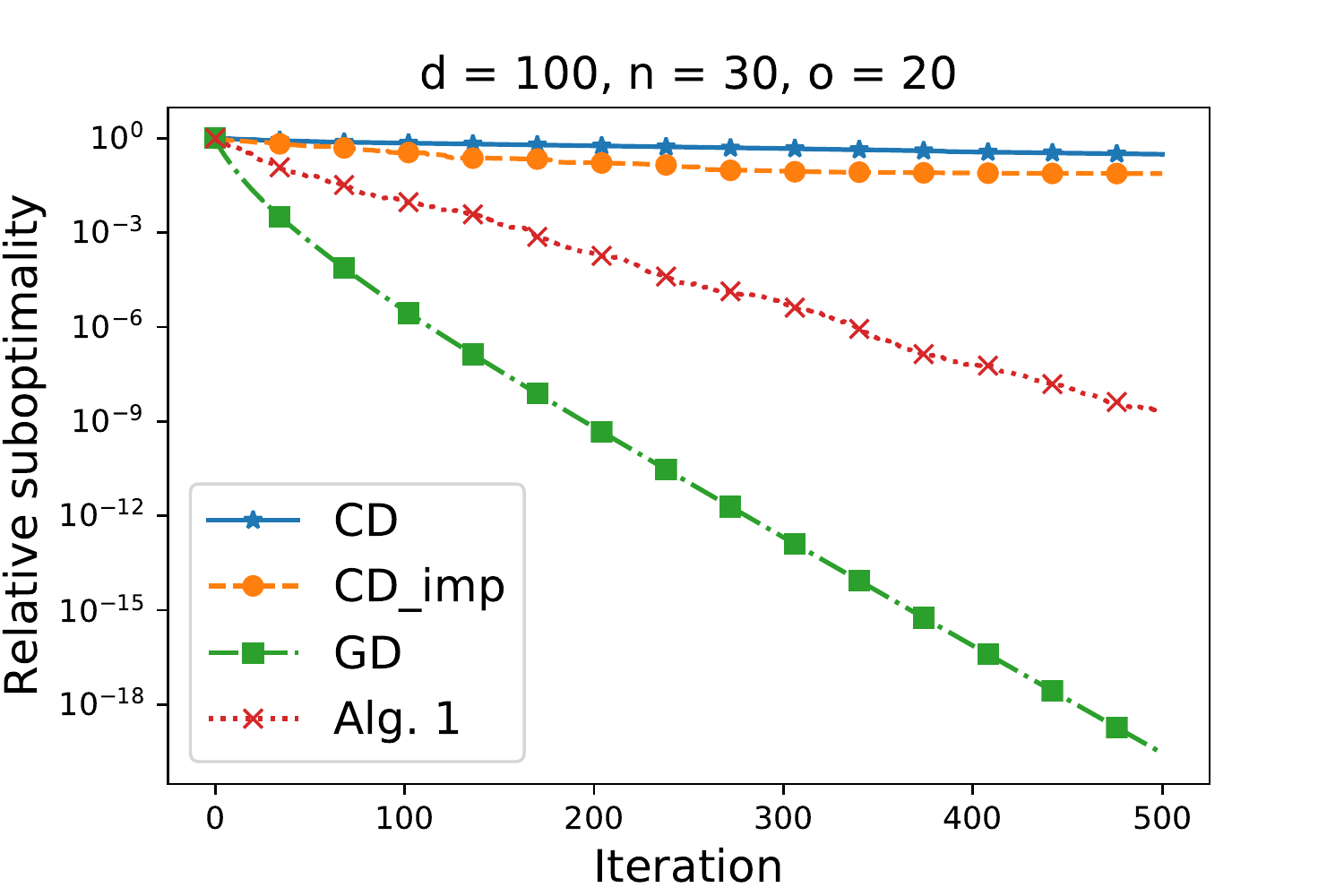}
\end{minipage}%
\begin{minipage}{0.3\textwidth}
  \centering
\includegraphics[width =  \textwidth ]{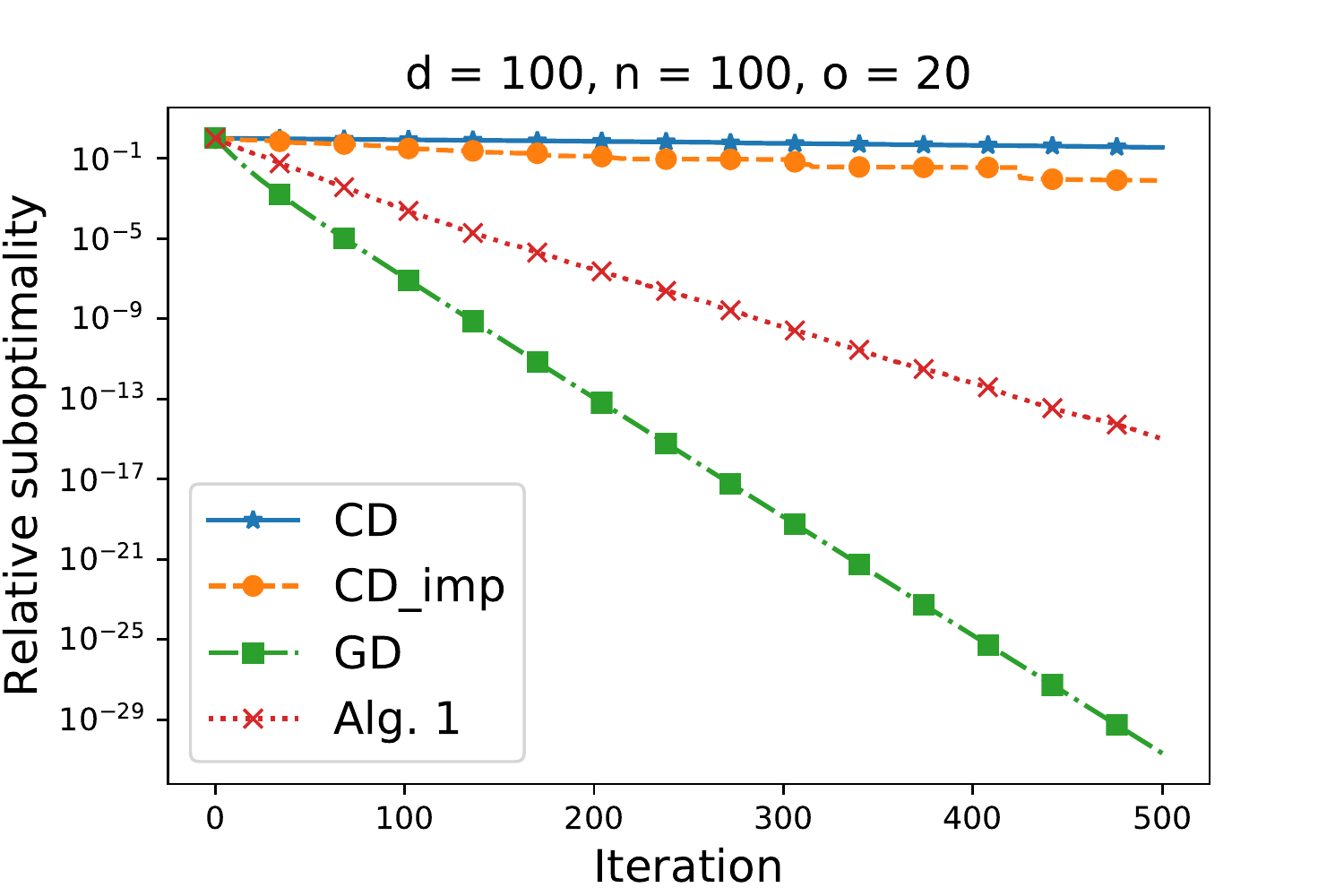}
\end{minipage}%
\begin{minipage}{0.3\textwidth}
  \centering
\includegraphics[width =  \textwidth ]{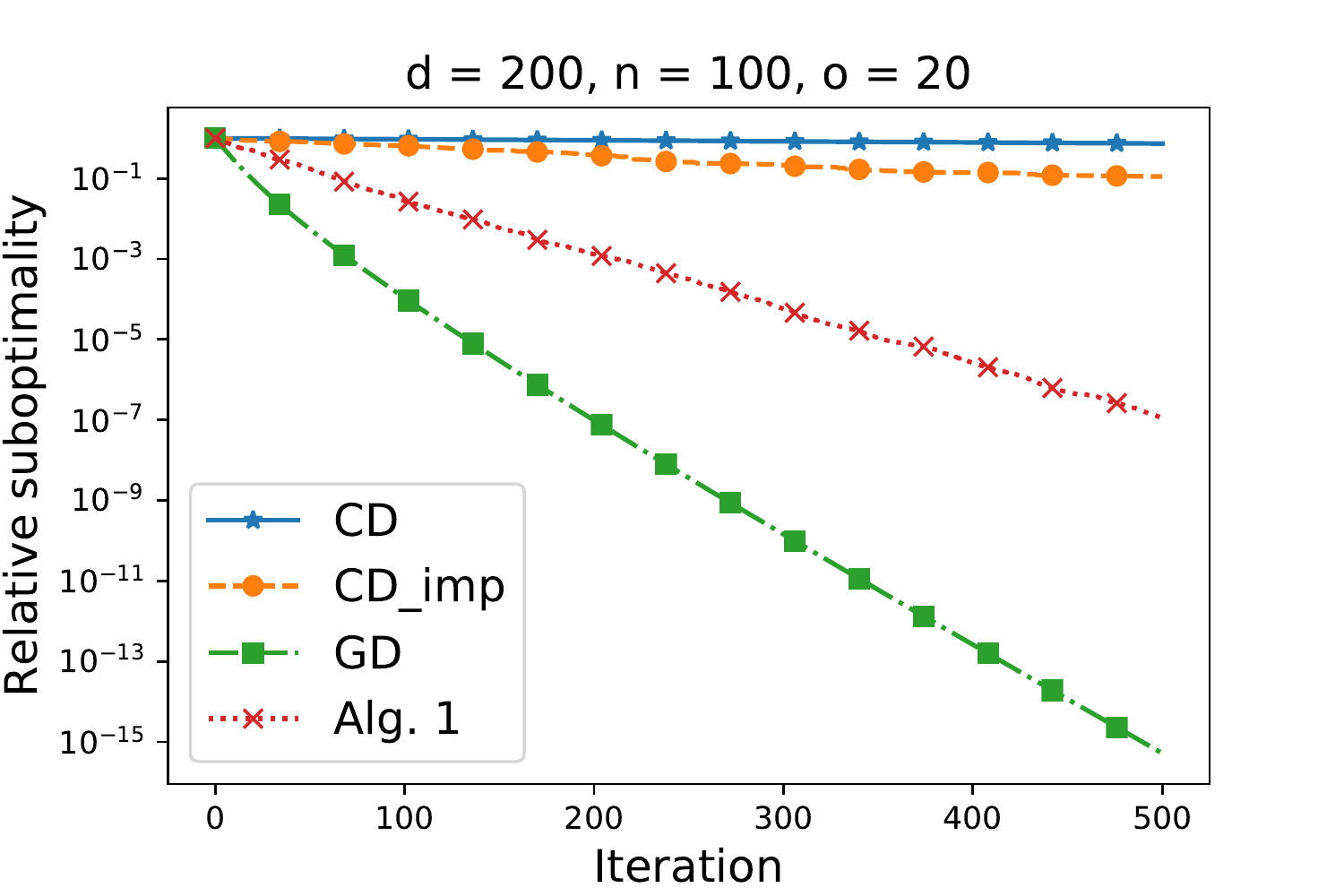}
\end{minipage}%
\\
\begin{minipage}{0.3\textwidth}
  \centering
\includegraphics[width =  \textwidth ]{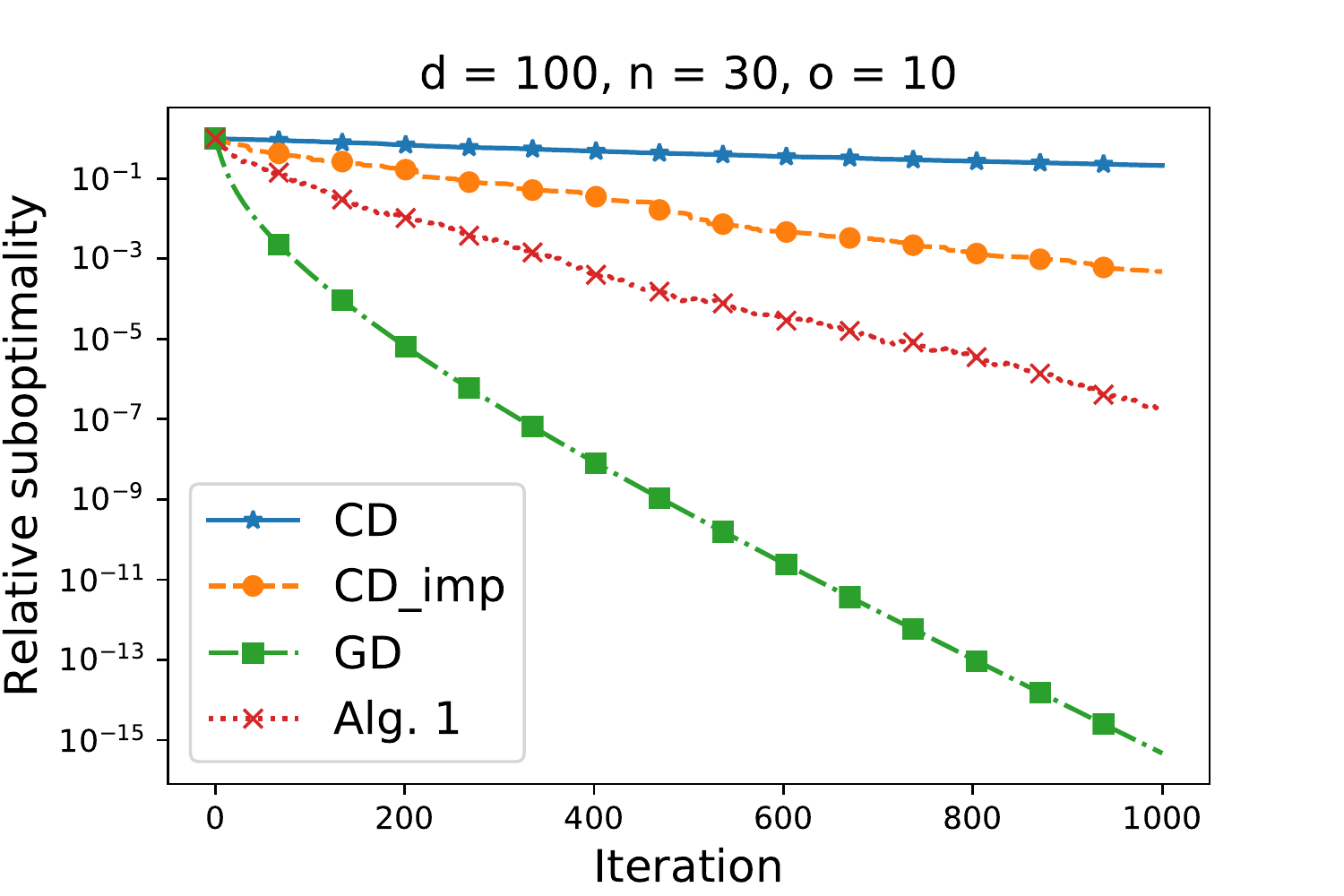}
\end{minipage}%
\begin{minipage}{0.3\textwidth}
  \centering
\includegraphics[width =  \textwidth ]{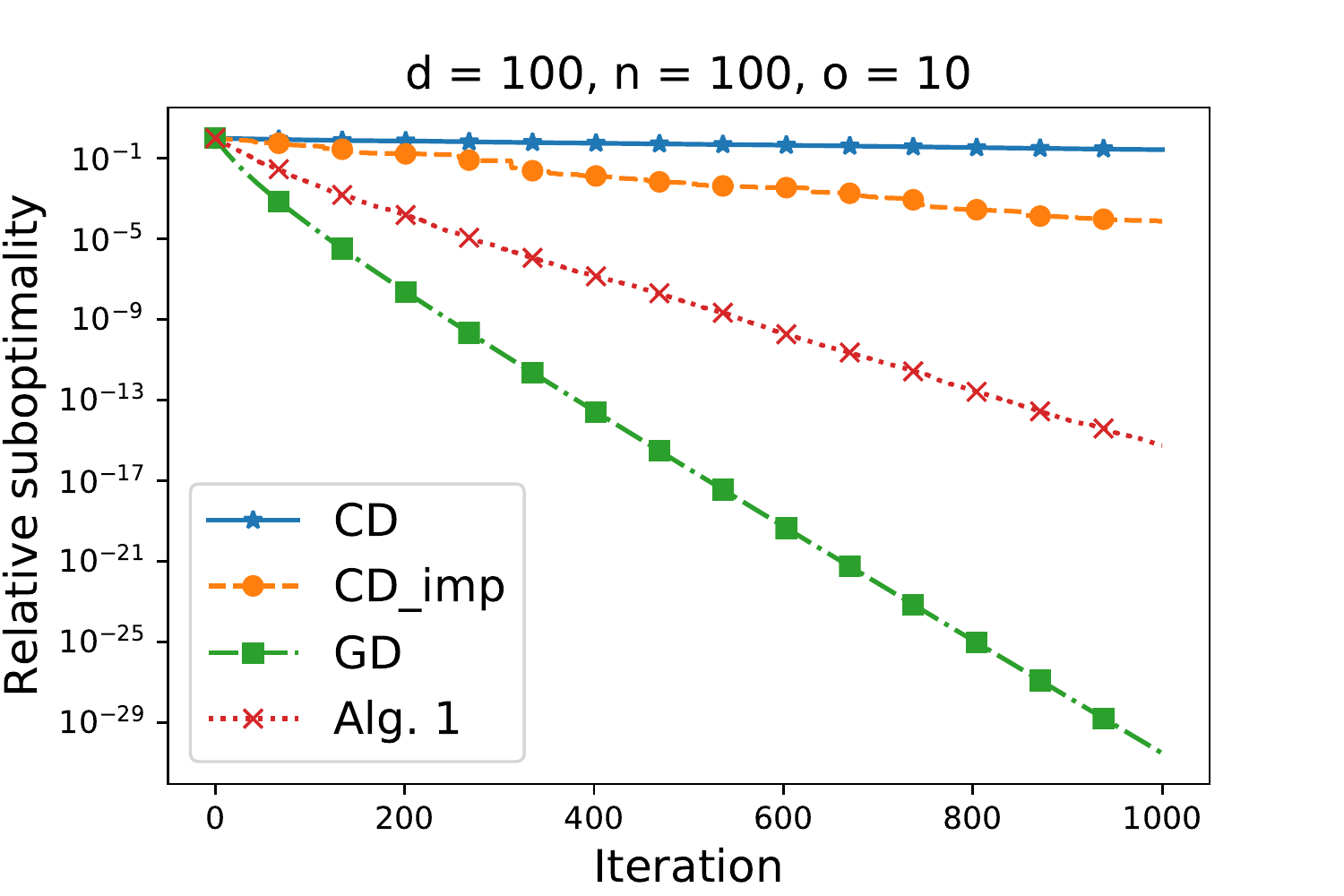}
\end{minipage}%
\begin{minipage}{0.3\textwidth}
  \centering
\includegraphics[width =  \textwidth ]{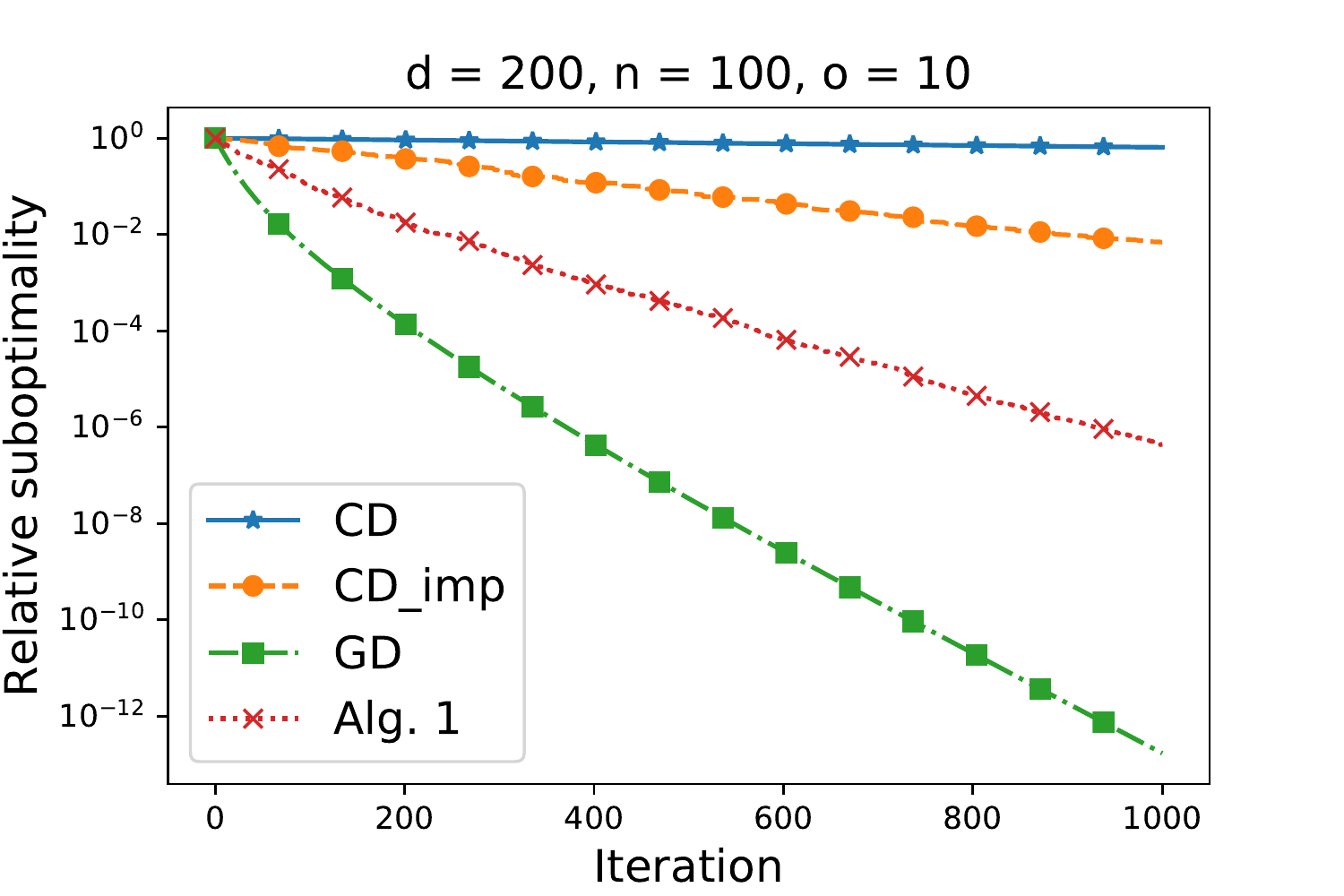}
\end{minipage}%
\\
\begin{minipage}{0.3\textwidth}
  \centering
\includegraphics[width =  \textwidth ]{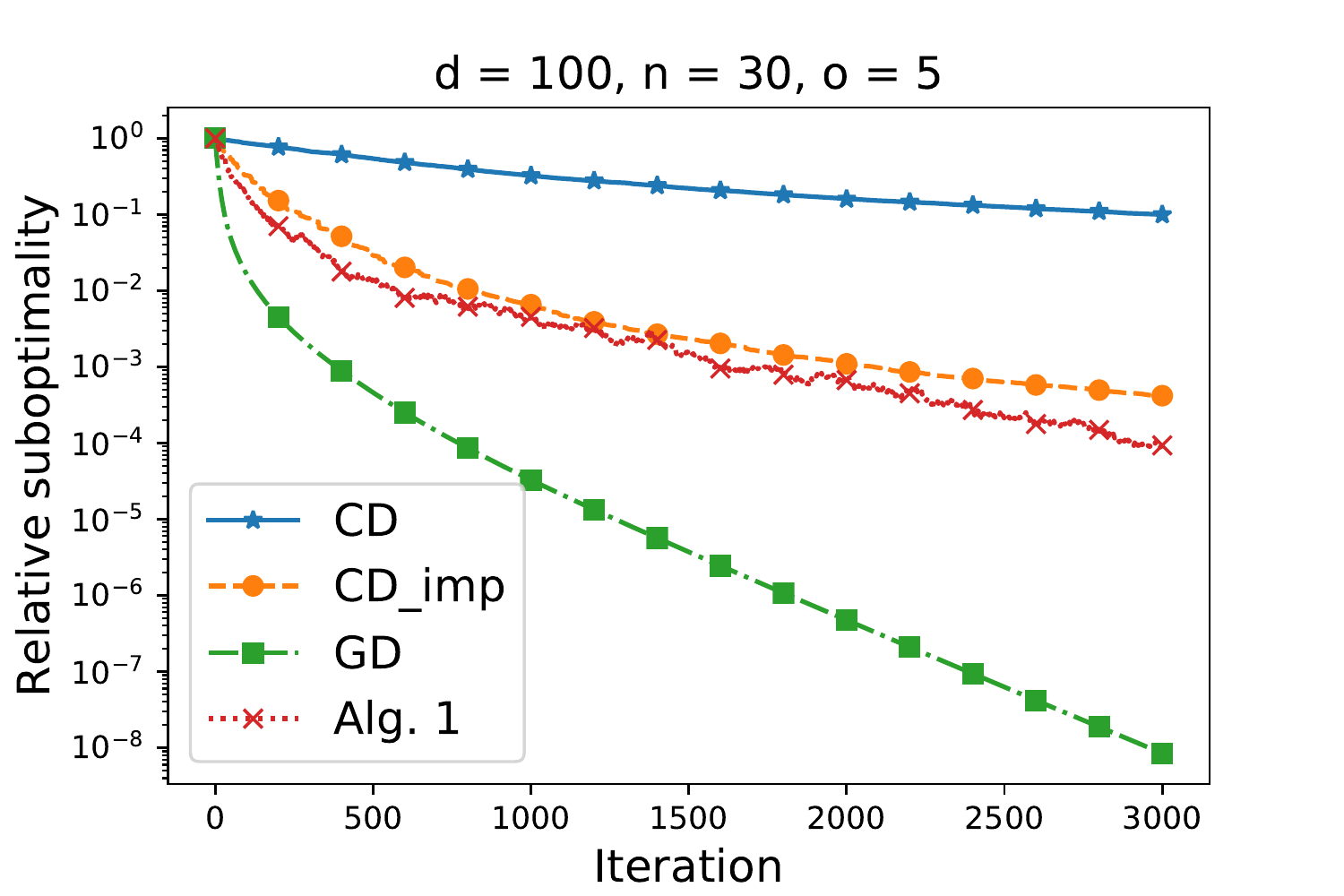}
\end{minipage}%
\begin{minipage}{0.3\textwidth}
  \centering
\includegraphics[width =  \textwidth ]{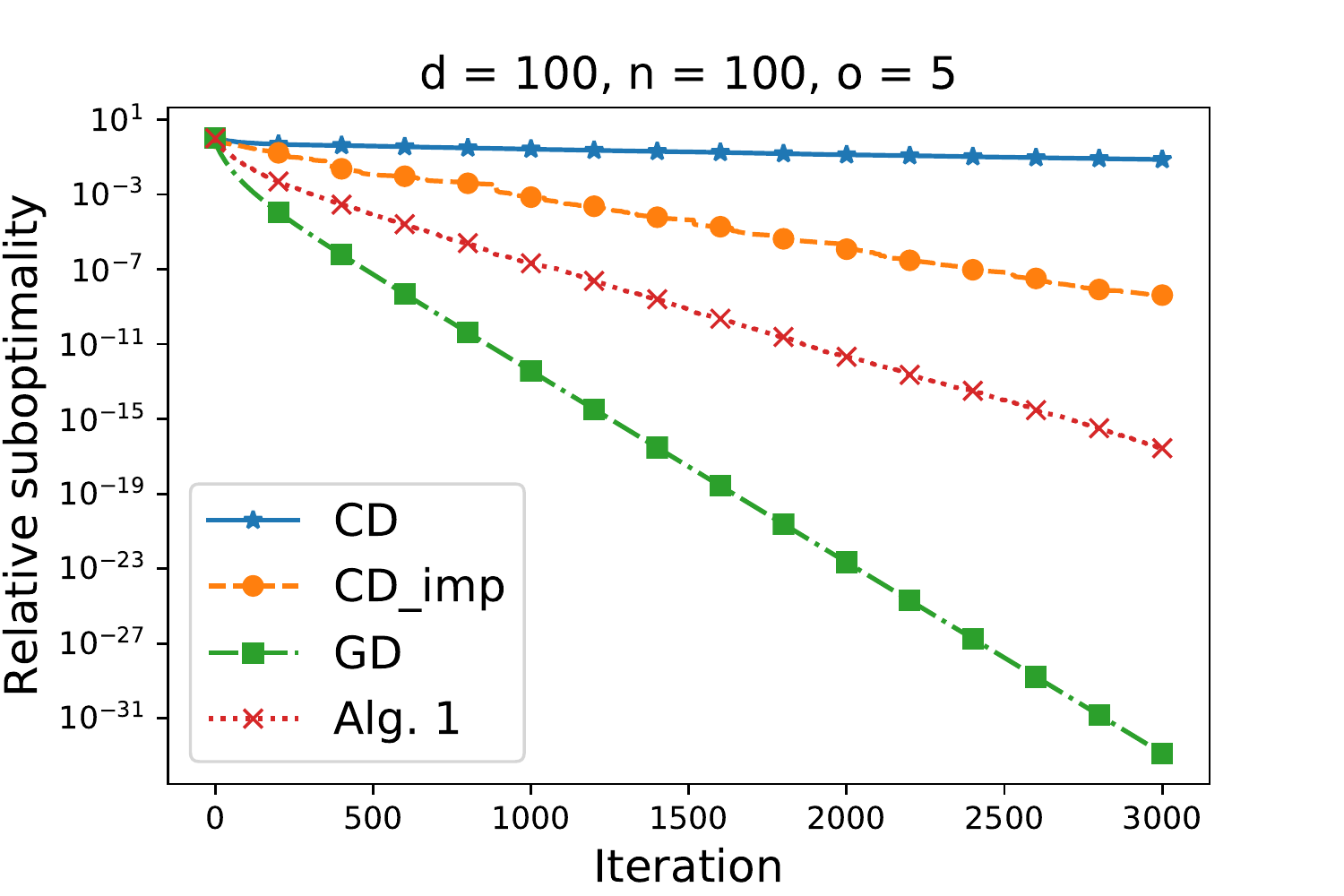}
\end{minipage}%
\begin{minipage}{0.3\textwidth}
  \centering
\includegraphics[width =  \textwidth ]{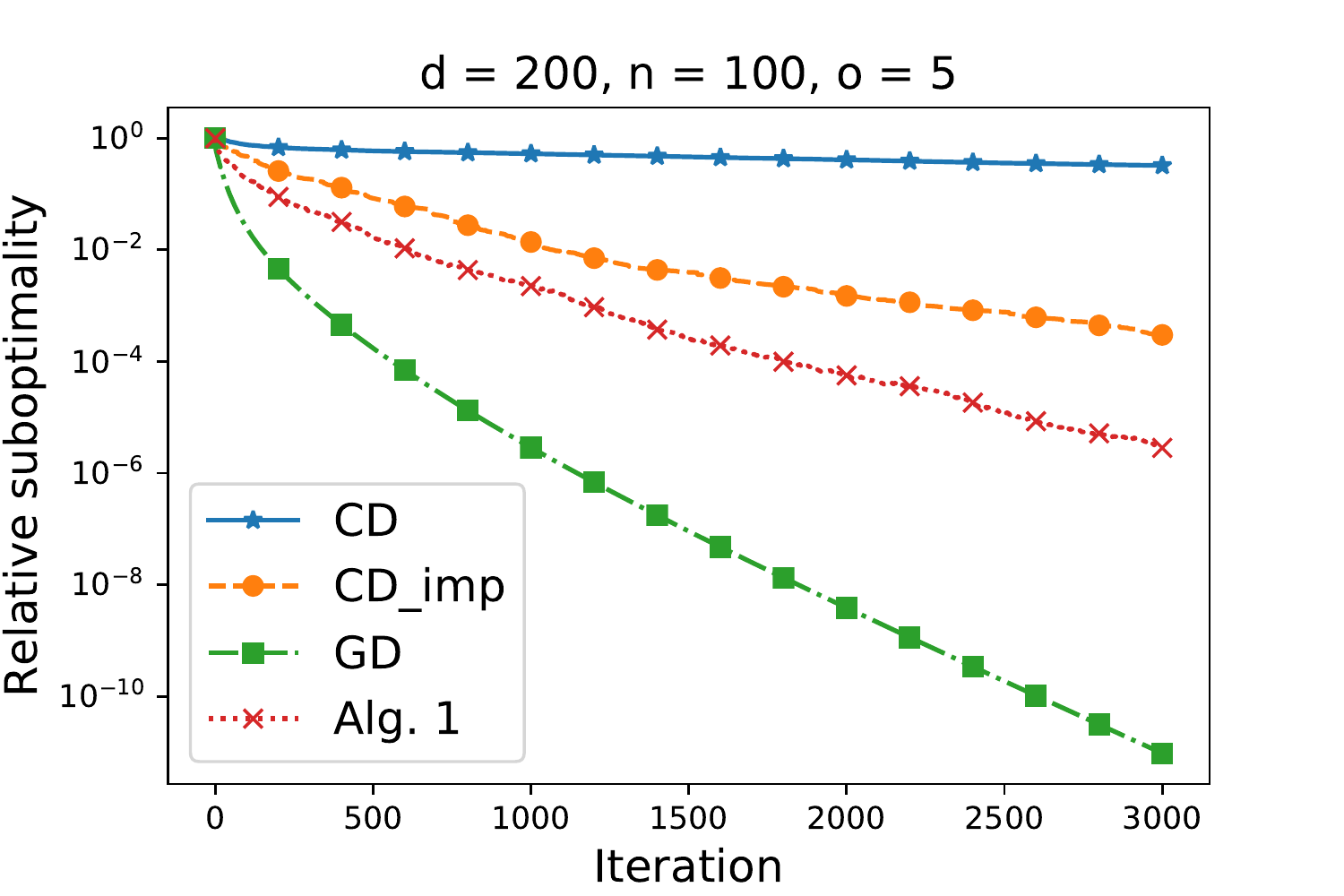}
\end{minipage}%
\\
\begin{minipage}{0.3\textwidth}
  \centering
\includegraphics[width =  \textwidth ]{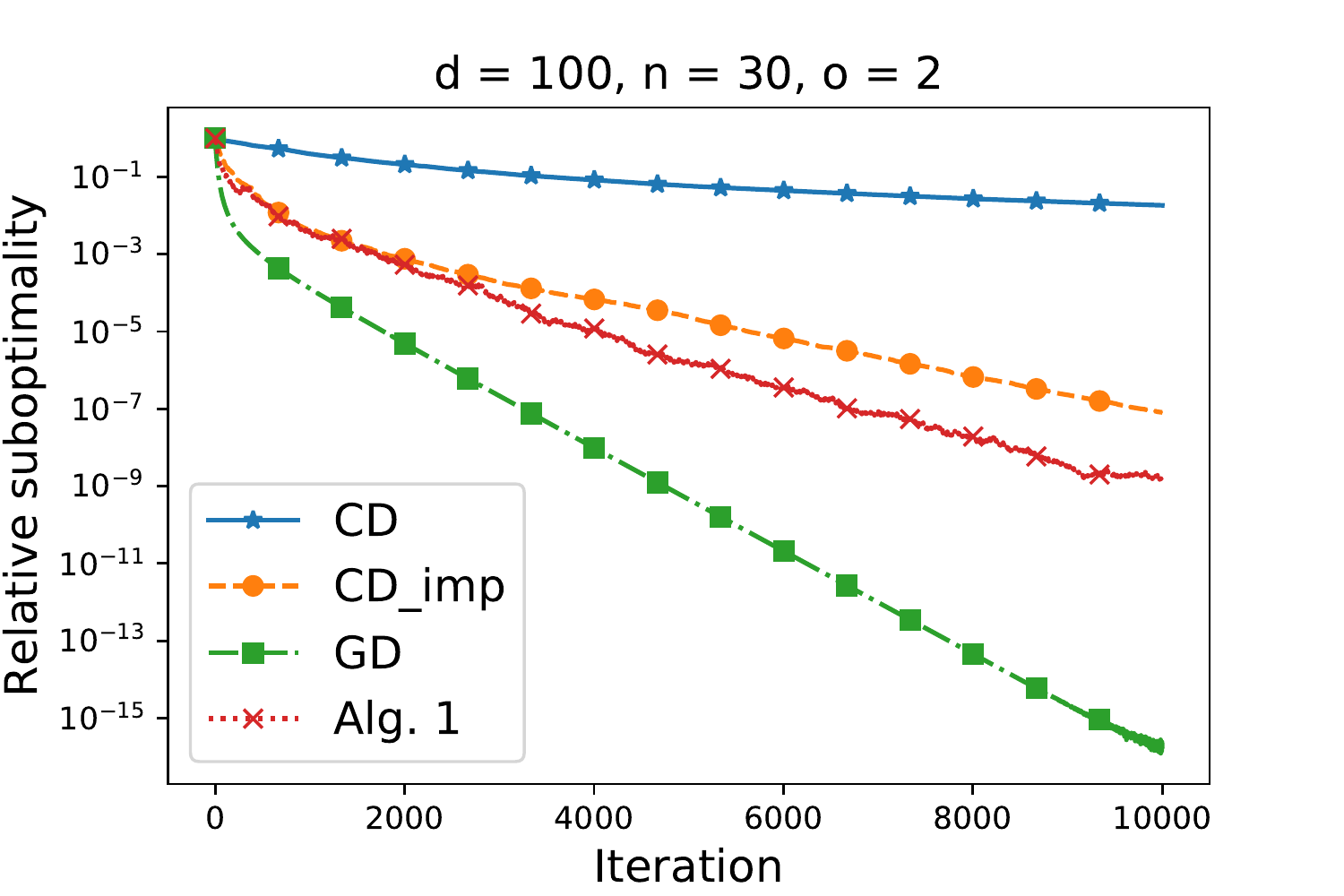}
\end{minipage}%
\begin{minipage}{0.3\textwidth}
  \centering
\includegraphics[width =  \textwidth ]{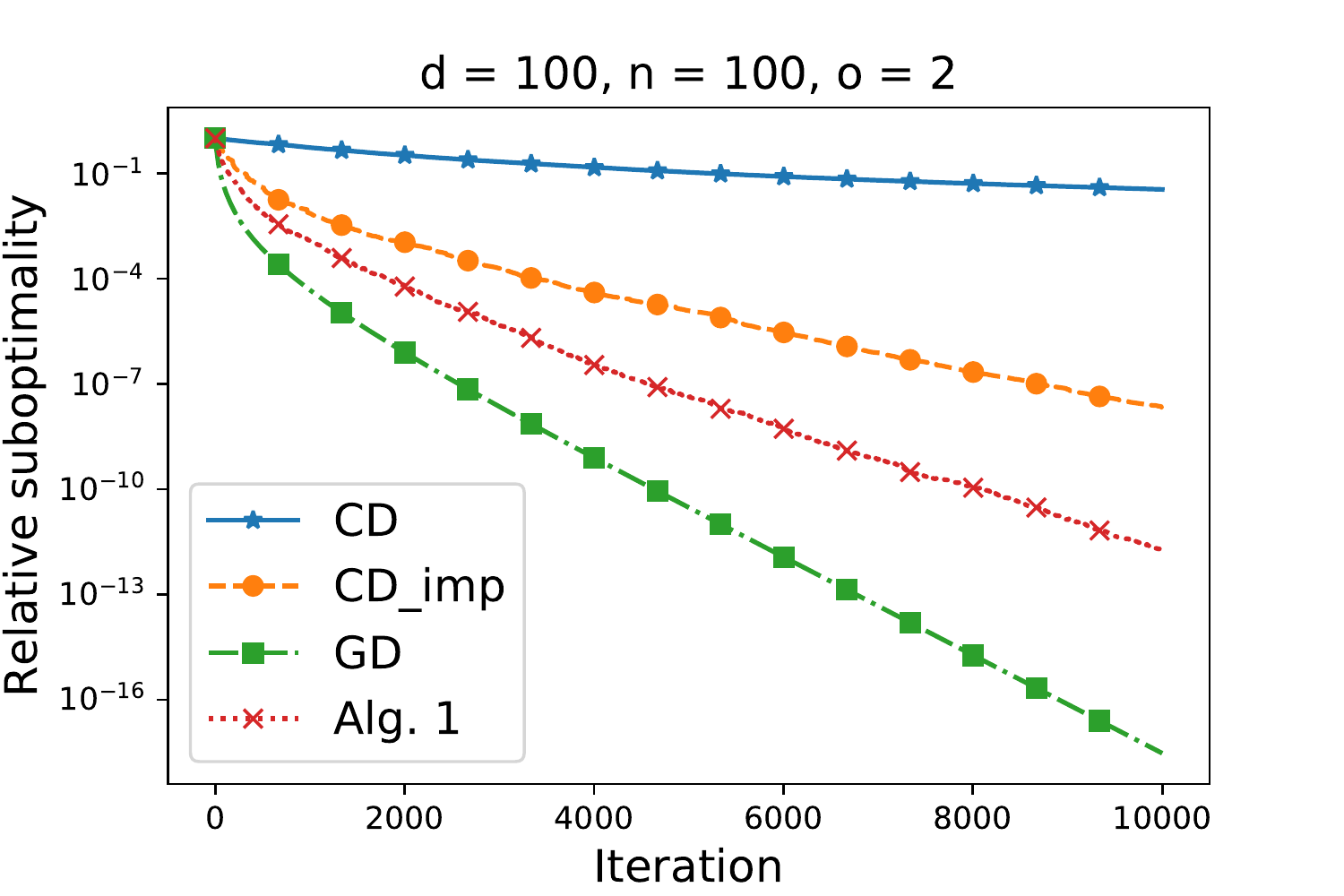}
\end{minipage}%
\begin{minipage}{0.3\textwidth}
  \centering
\includegraphics[width =  \textwidth ]{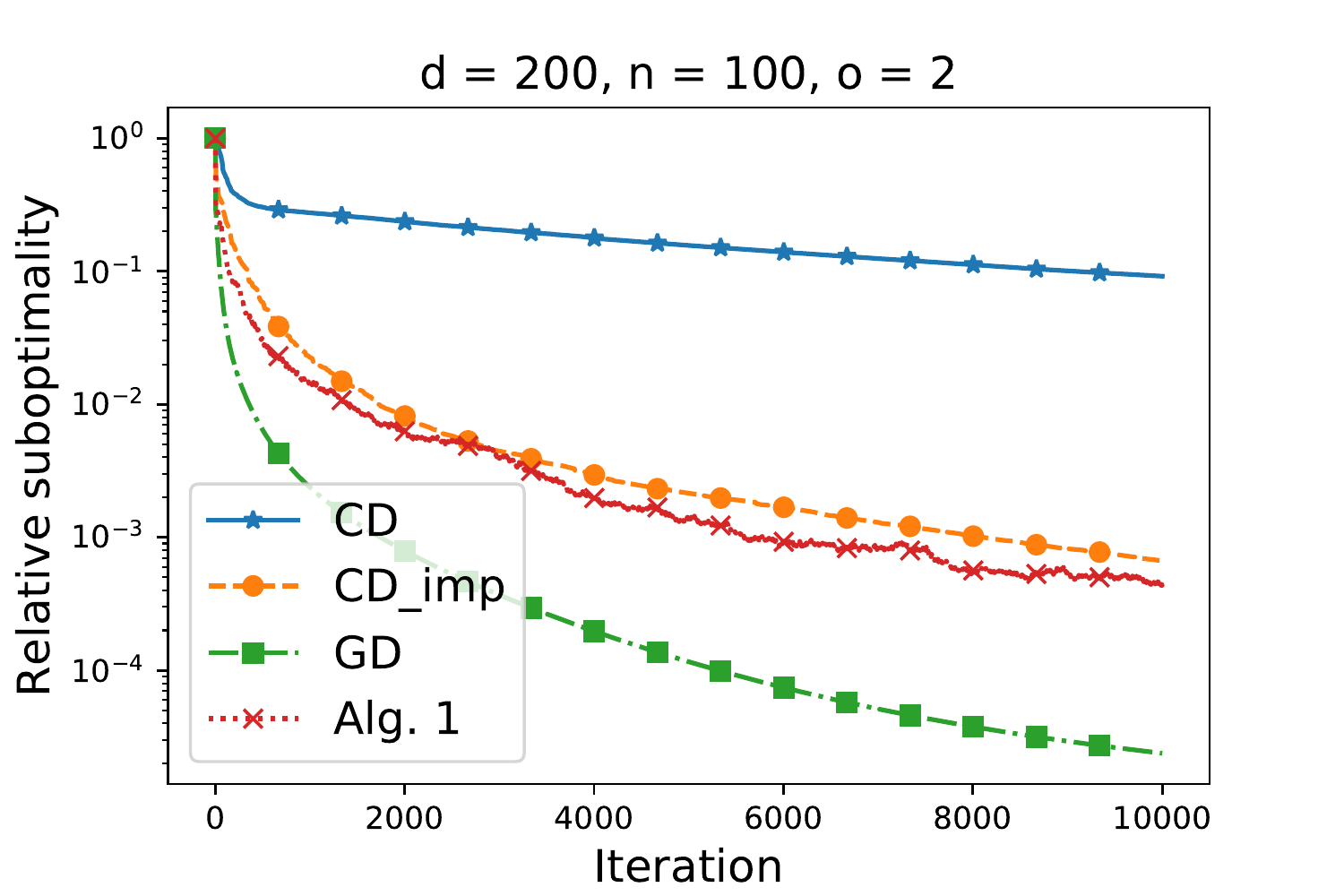}
\end{minipage}%
\\
\caption{Comparison of gradient descent, (standard) coordinate descent, (standard) coordinate descent with importance sampling  and Algorithm~\ref{alg:cd} on artificial quadratic problem~\eqref{eq:99_quadratic}.}\label{fig:99_artif_1}
\end{figure}

Next, we study the effect of changing $\tau$ on the iteration complexity of Algorithm~\ref{alg:cd}. Figure~\ref{fig:99_artif_2} provides the result. The behavior predicted from theory is observed --  increasing  $\tau$ over $n^{-1}$ does not significantly improve the convergence speed, while decreasing it below $n^{-1}$ slows the algorithm notably.

\begin{figure}[H]
\centering
\begin{minipage}{0.3\textwidth}
  \centering
\includegraphics[width =  \textwidth ]{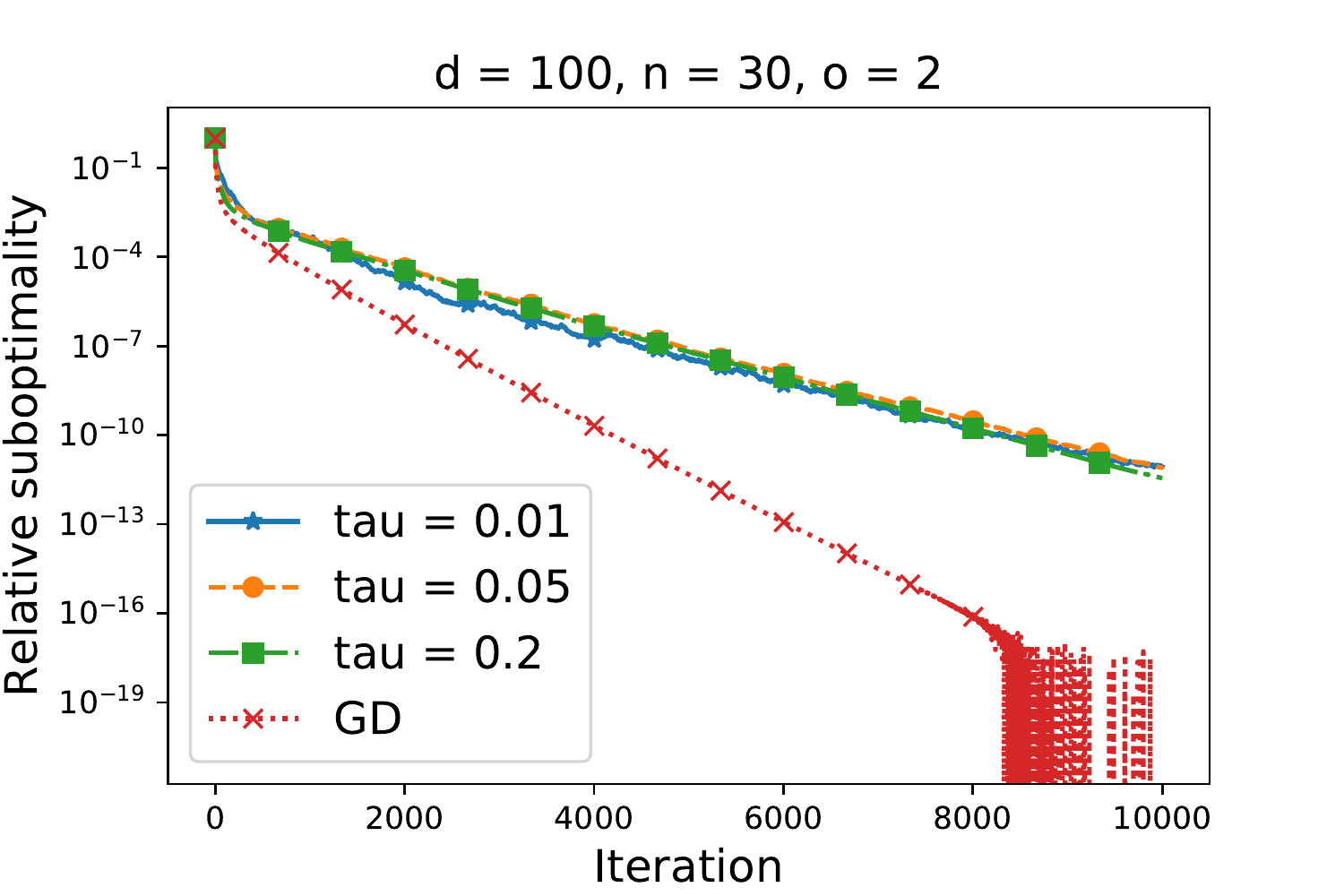}
\end{minipage}%
\begin{minipage}{0.3\textwidth}
  \centering
\includegraphics[width =  \textwidth ]{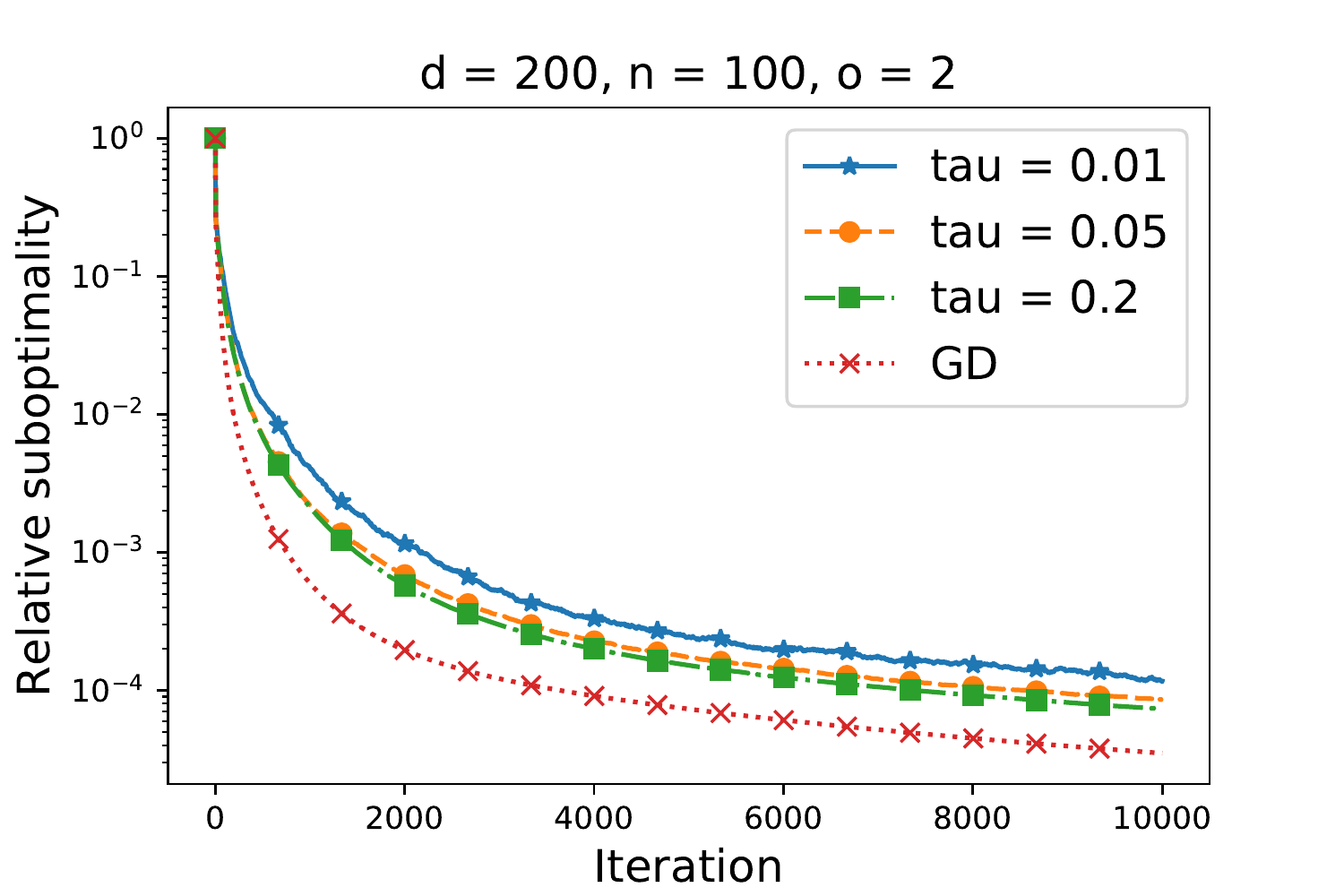}
\end{minipage}%
\begin{minipage}{0.3\textwidth}
  \centering
\includegraphics[width =  \textwidth ]{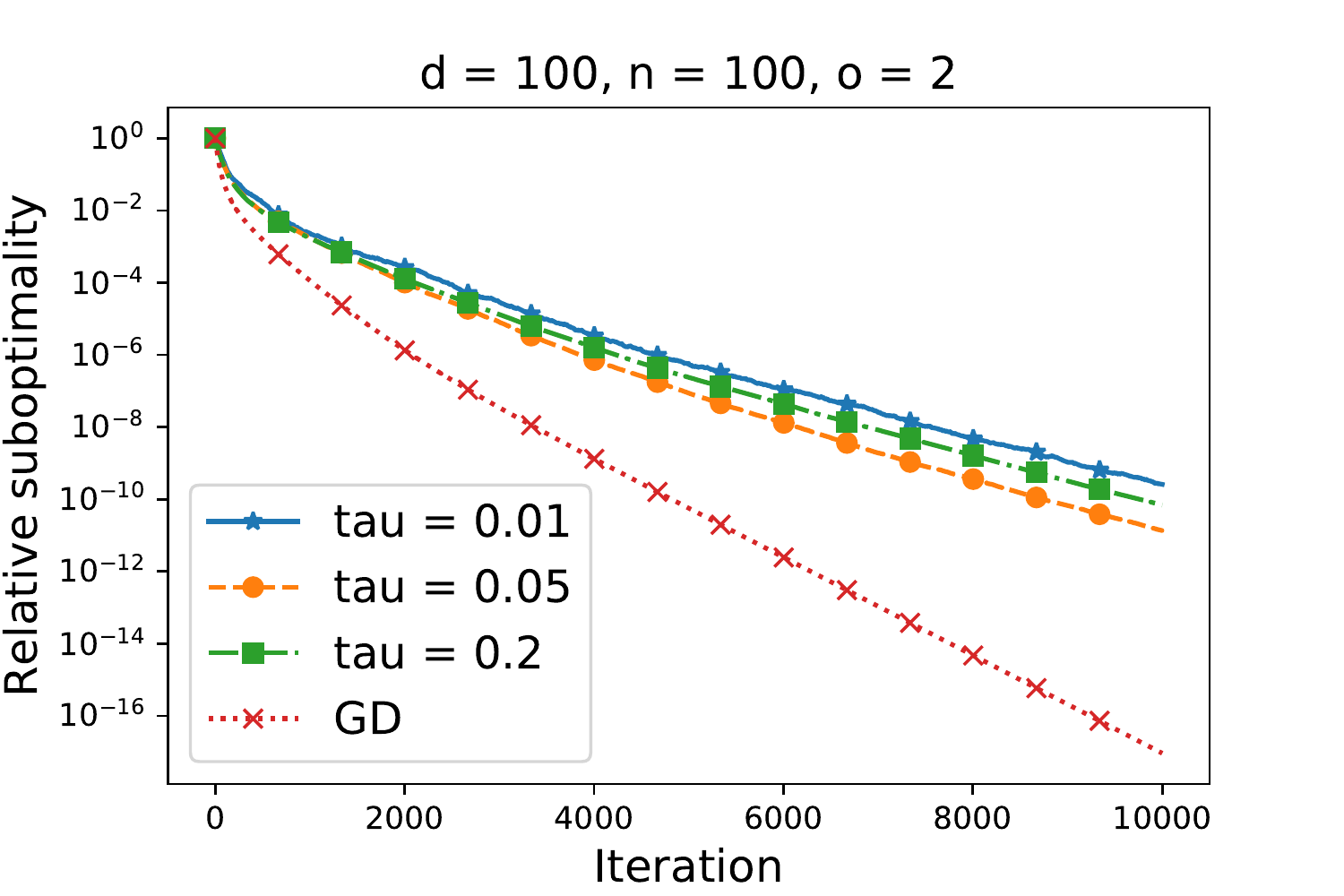}
\end{minipage}%
\\
\begin{minipage}{0.3\textwidth}
  \centering
\includegraphics[width =  \textwidth ]{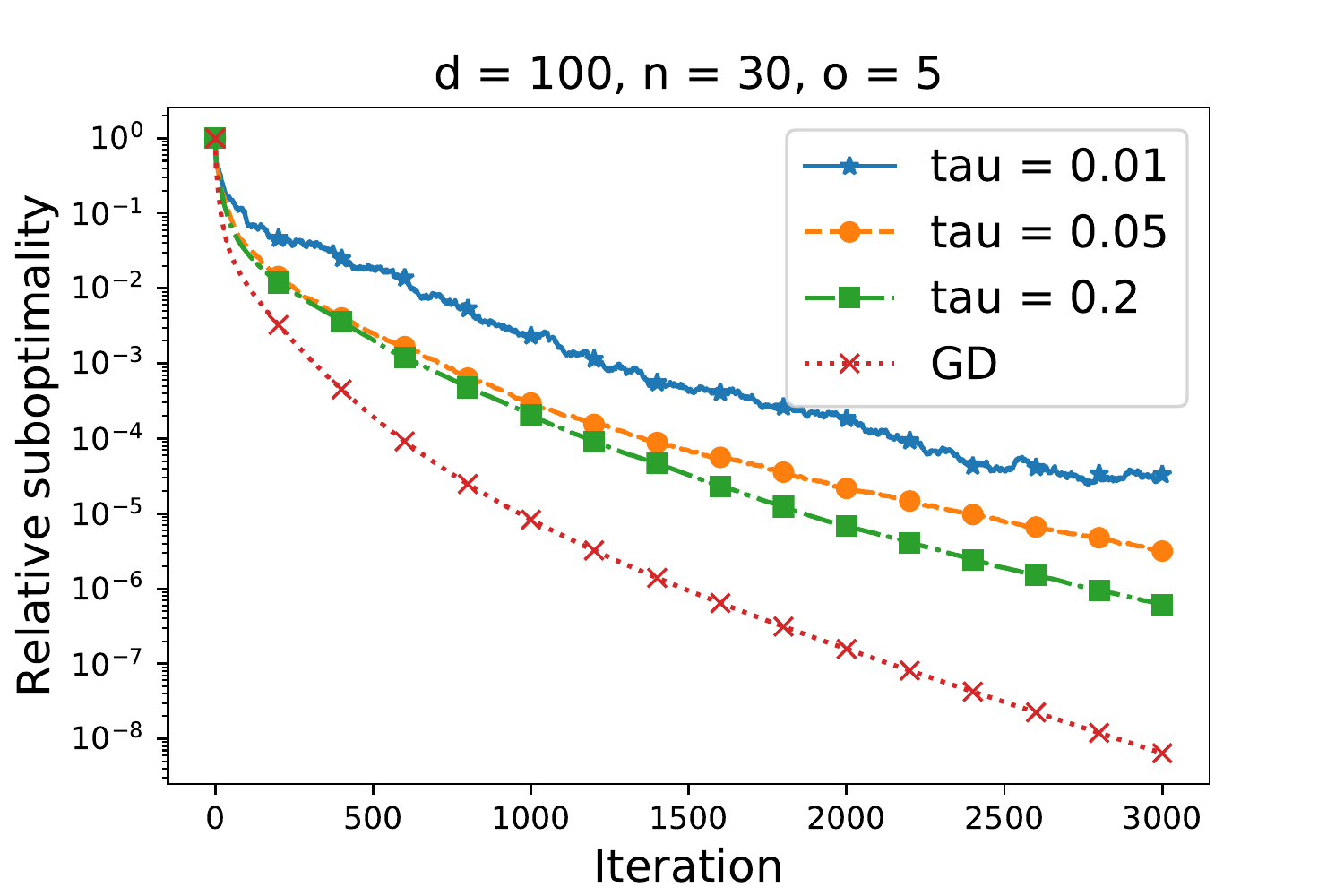}
\end{minipage}%
\begin{minipage}{0.3\textwidth}
  \centering
\includegraphics[width =  \textwidth ]{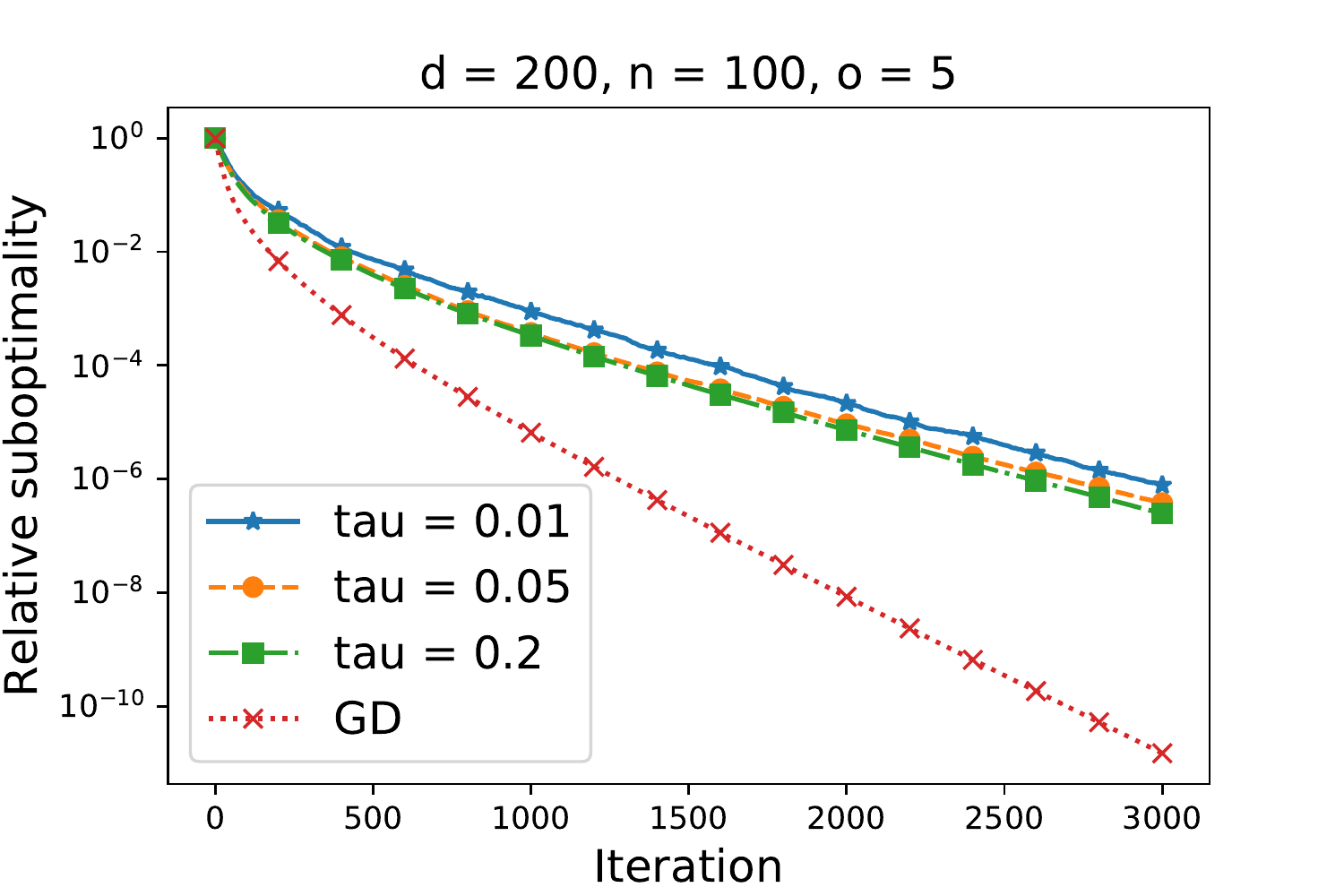}
\end{minipage}%
\begin{minipage}{0.3\textwidth}
  \centering
\includegraphics[width =  \textwidth ]{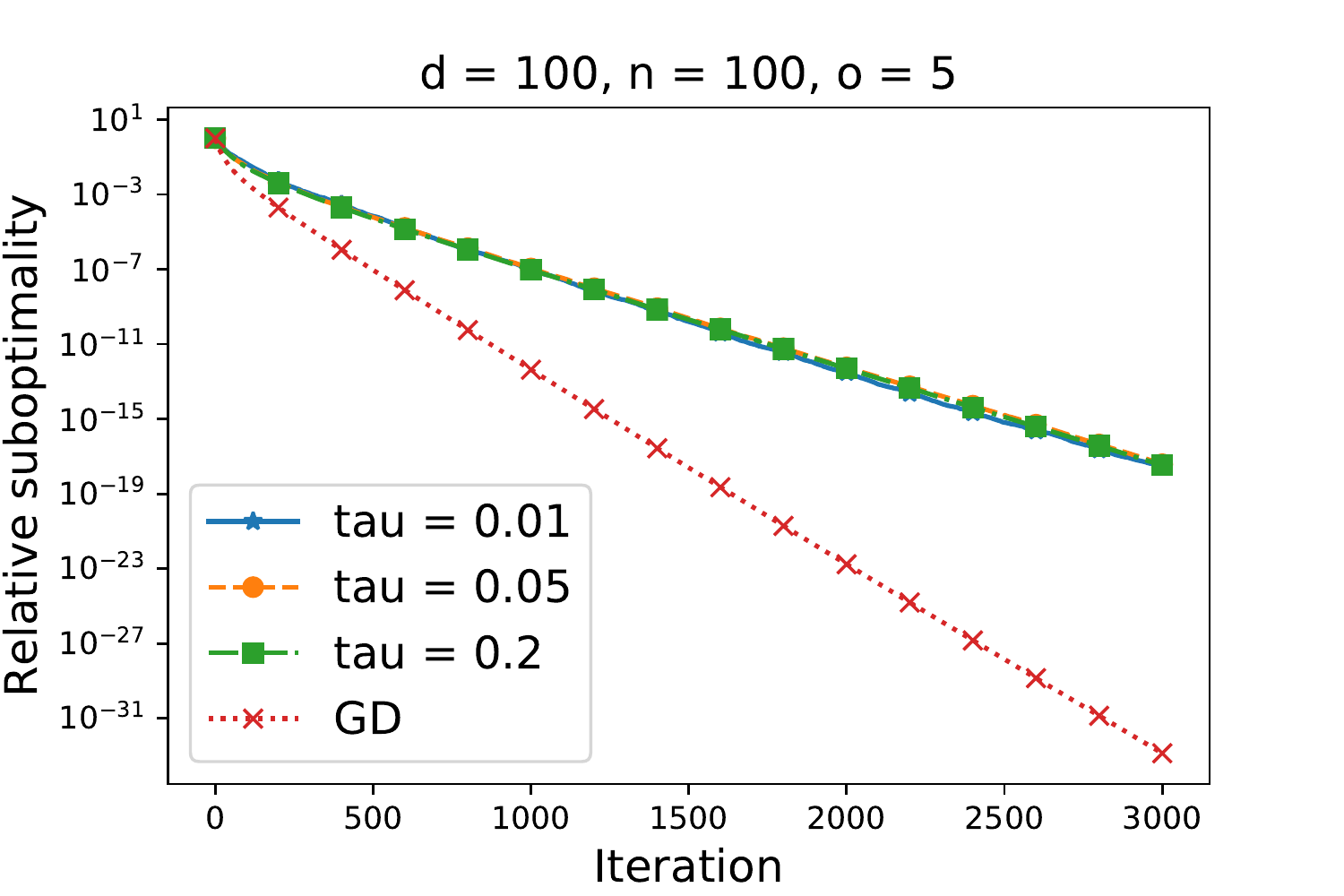}
\end{minipage}%
\\
\begin{minipage}{0.3\textwidth}
  \centering
\includegraphics[width =  \textwidth ]{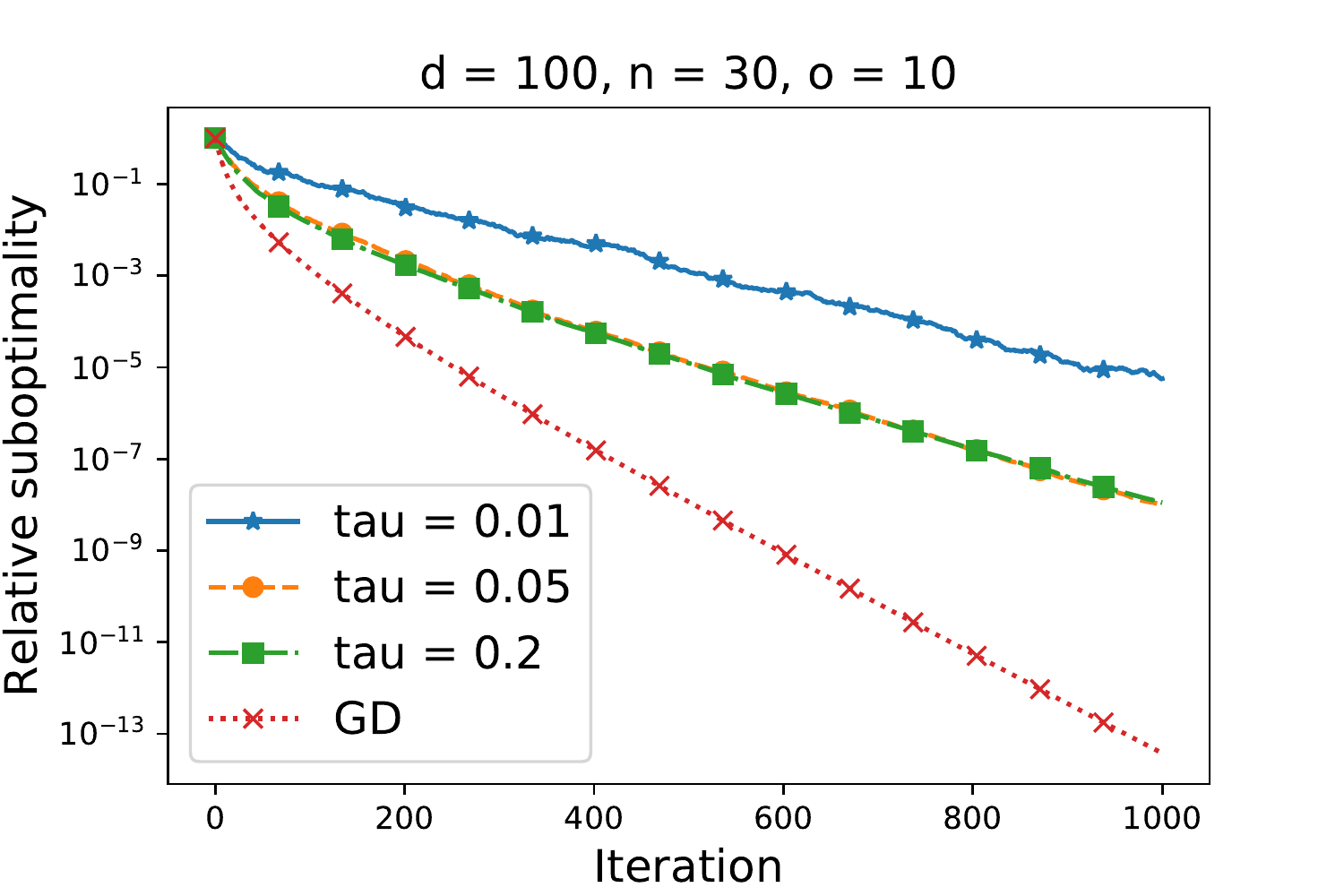}
\end{minipage}%
\begin{minipage}{0.3\textwidth}
  \centering
\includegraphics[width =  \textwidth ]{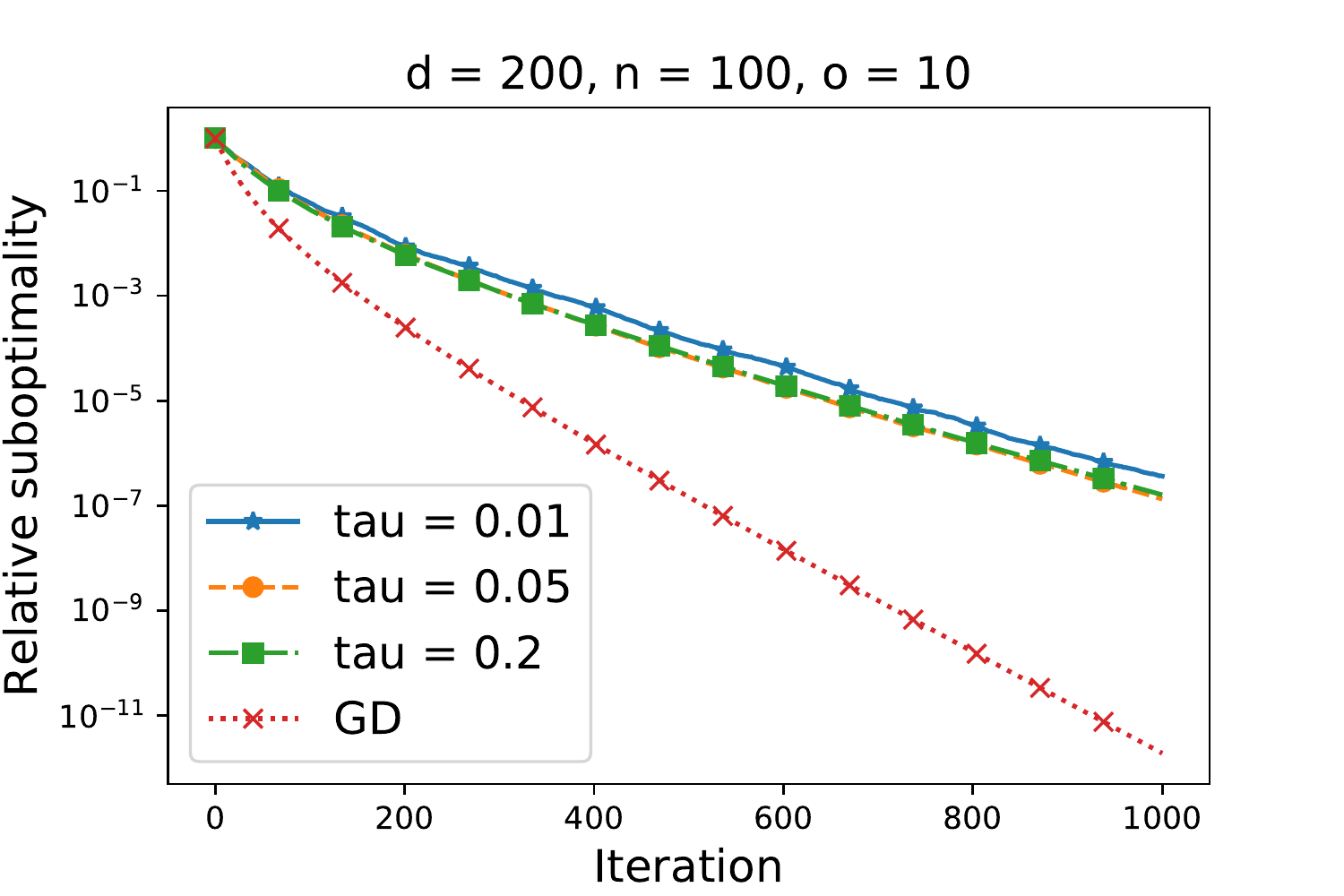}
\end{minipage}%
\begin{minipage}{0.3\textwidth}
  \centering
\includegraphics[width =  \textwidth ]{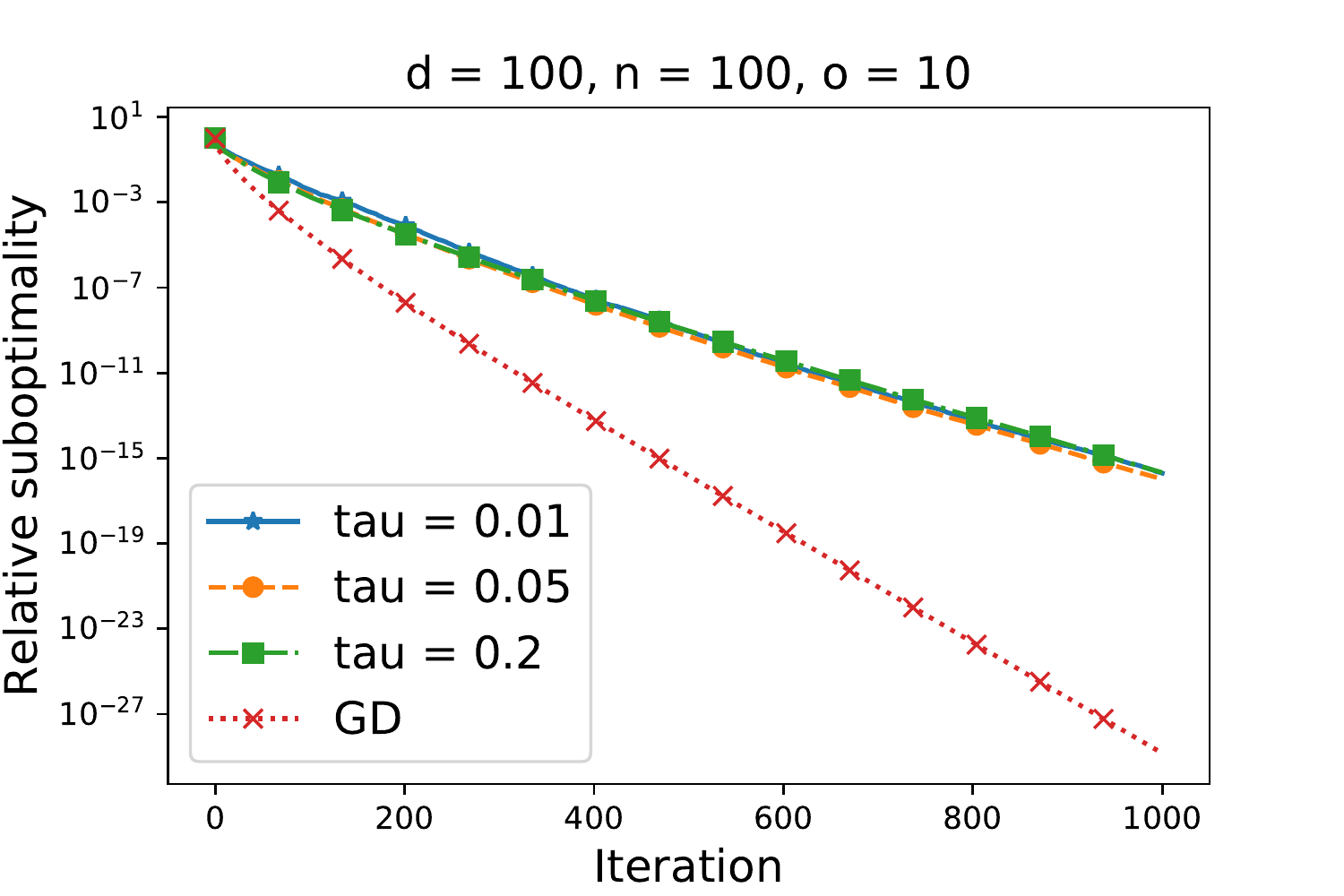}
\end{minipage}%
\\
\begin{minipage}{0.3\textwidth}
  \centering
\includegraphics[width =  \textwidth ]{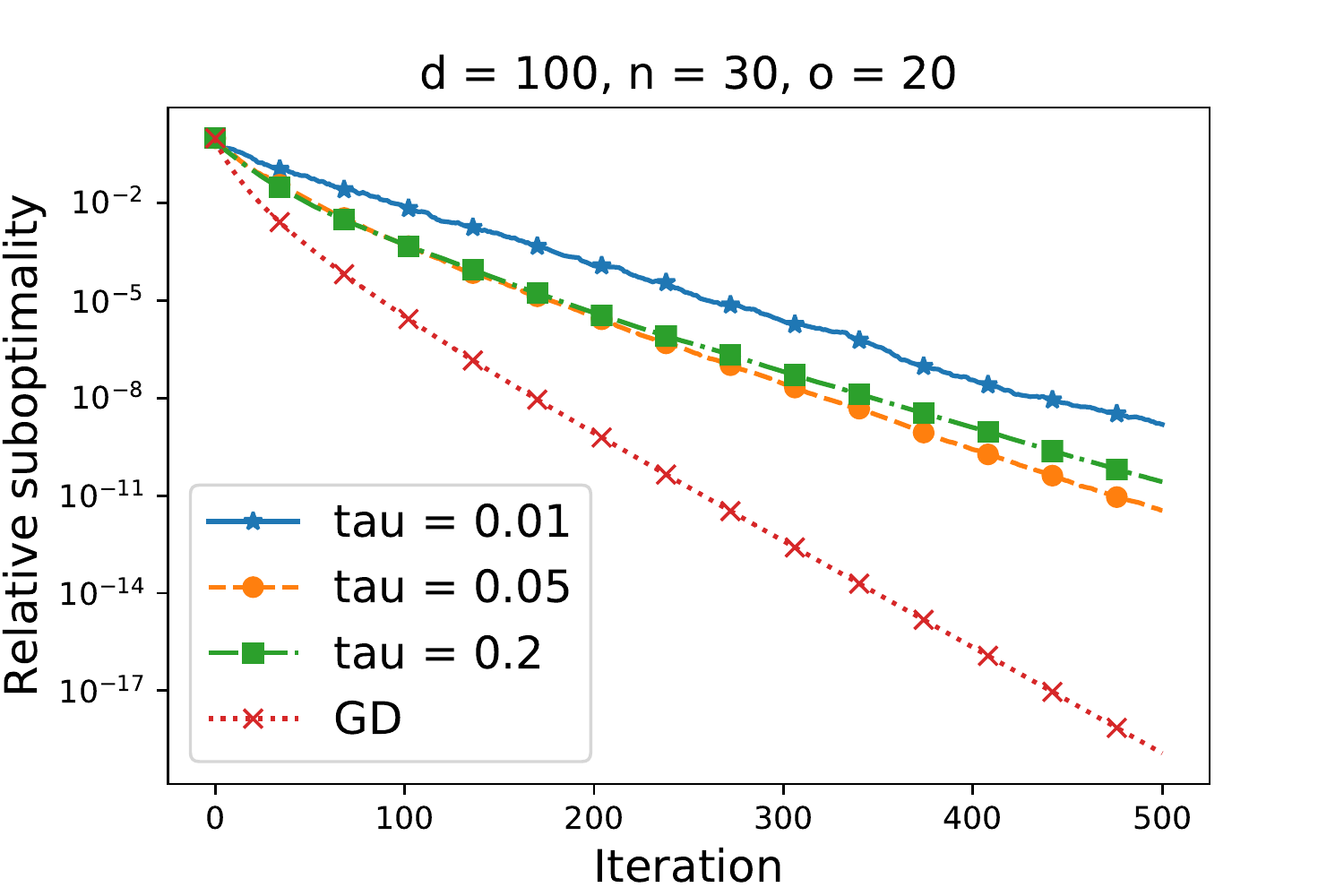}
\end{minipage}%
\begin{minipage}{0.3\textwidth}
  \centering
\includegraphics[width =  \textwidth ]{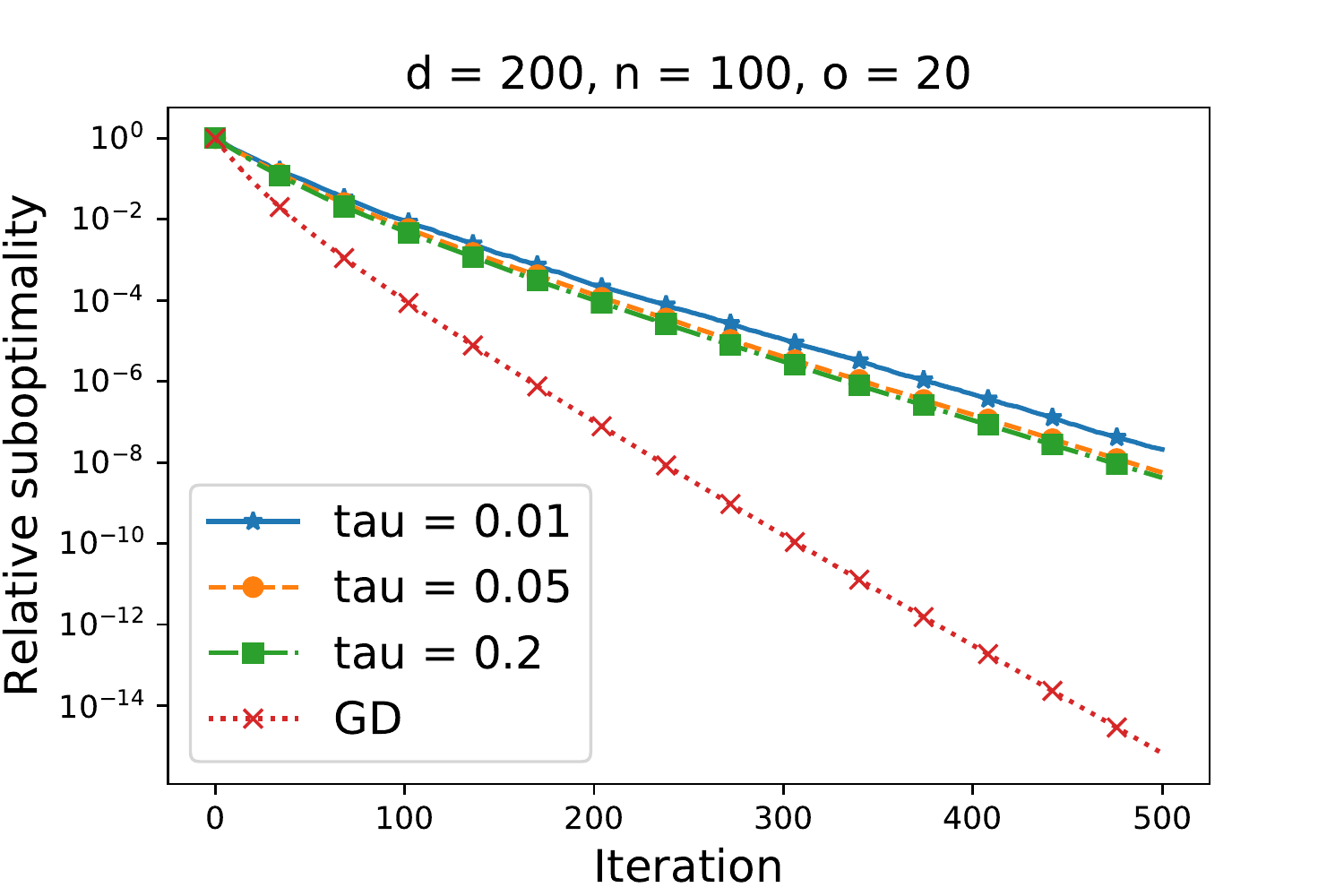}
\end{minipage}%
\begin{minipage}{0.3\textwidth}
  \centering
\includegraphics[width =  \textwidth ]{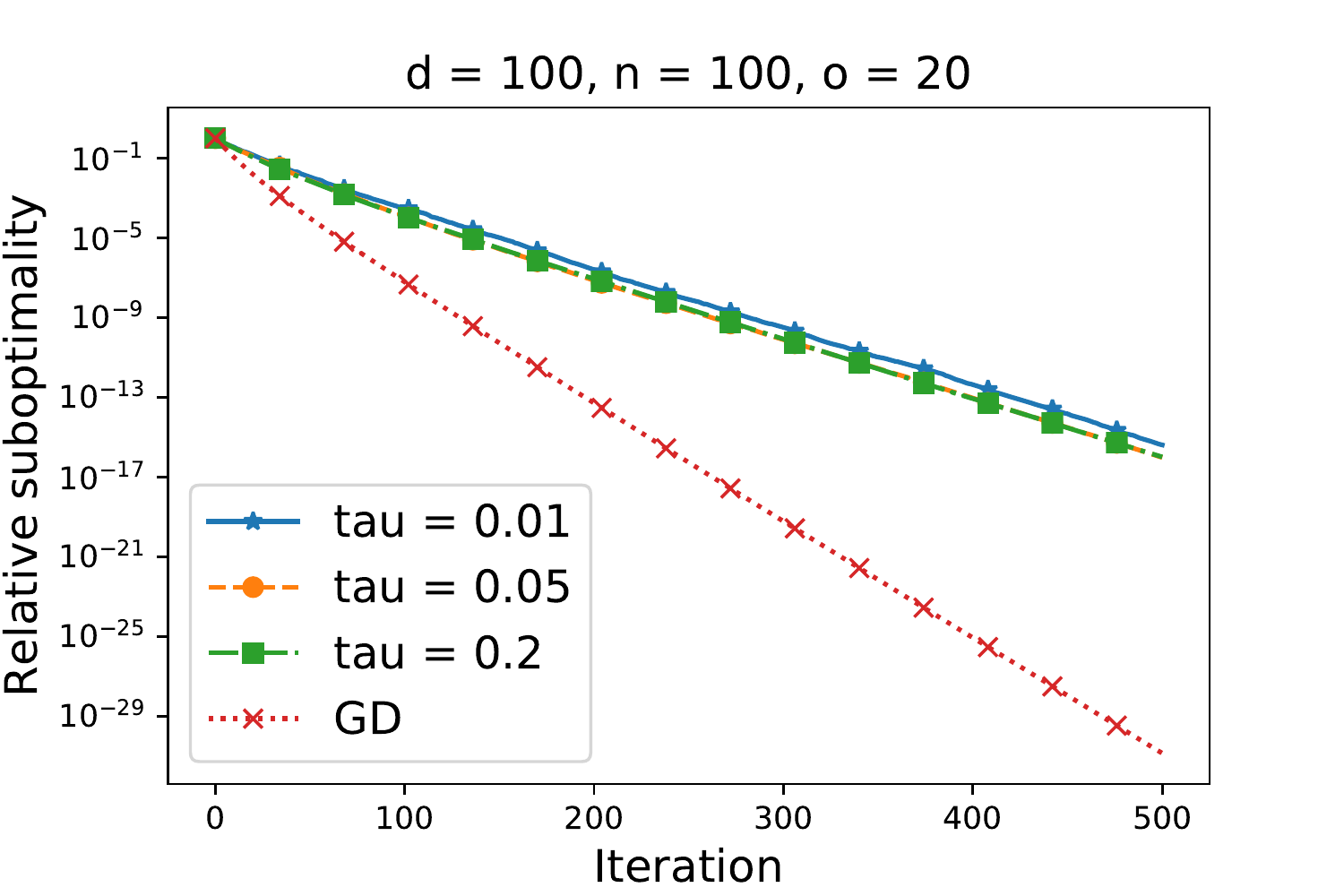}
\end{minipage}%
\\
\caption{Behavior of Algorithm~\ref{alg:cd} for different $\tau$ on a simple artificial quadratic problem~\eqref{eq:99_quadratic}.}\label{fig:99_artif_2}
\end{figure}

\subsection{{\tt ISGD} \label{sec:99_exp_sgd}}

In this section we numerically test Algorithm~\ref{alg:sgd} for logistic regression problem. As mentioned, $f_i$ consists of set of (uniformly distributed) rows of $\mA$ from~\eqref{eq:99_logreg}. We consider the most natural unbiased stochastic oracle for the $\nabla f_i$: the gradient computed on a subset of the data points from $f_i$. 

In all experiments of this section, we consider constant step sizes in order to keep the setting as simple as possible and gain as much insight from the experiments as possible. Therefore, one can not expect convergence to the exact optimum. 

In the first experiment, we compare standard {\tt SGD} (stochastic gradient is computed on single, randomly chosen datapoint every iteration) against Algorithm~\ref{alg:sgd} varying $n$ and choosing $\tau=\frac{1}{n}$ for each $n$. The results are presented by Figure~\ref{fig:99_sgd1}. We see that, as our theory suggests, {\tt SGD} and Algorithm~\ref{alg:sgd} have always very similar performance.

\begin{figure}[H]
\centering
\begin{minipage}{0.33\textwidth}
  \centering
\includegraphics[width =  \textwidth ]{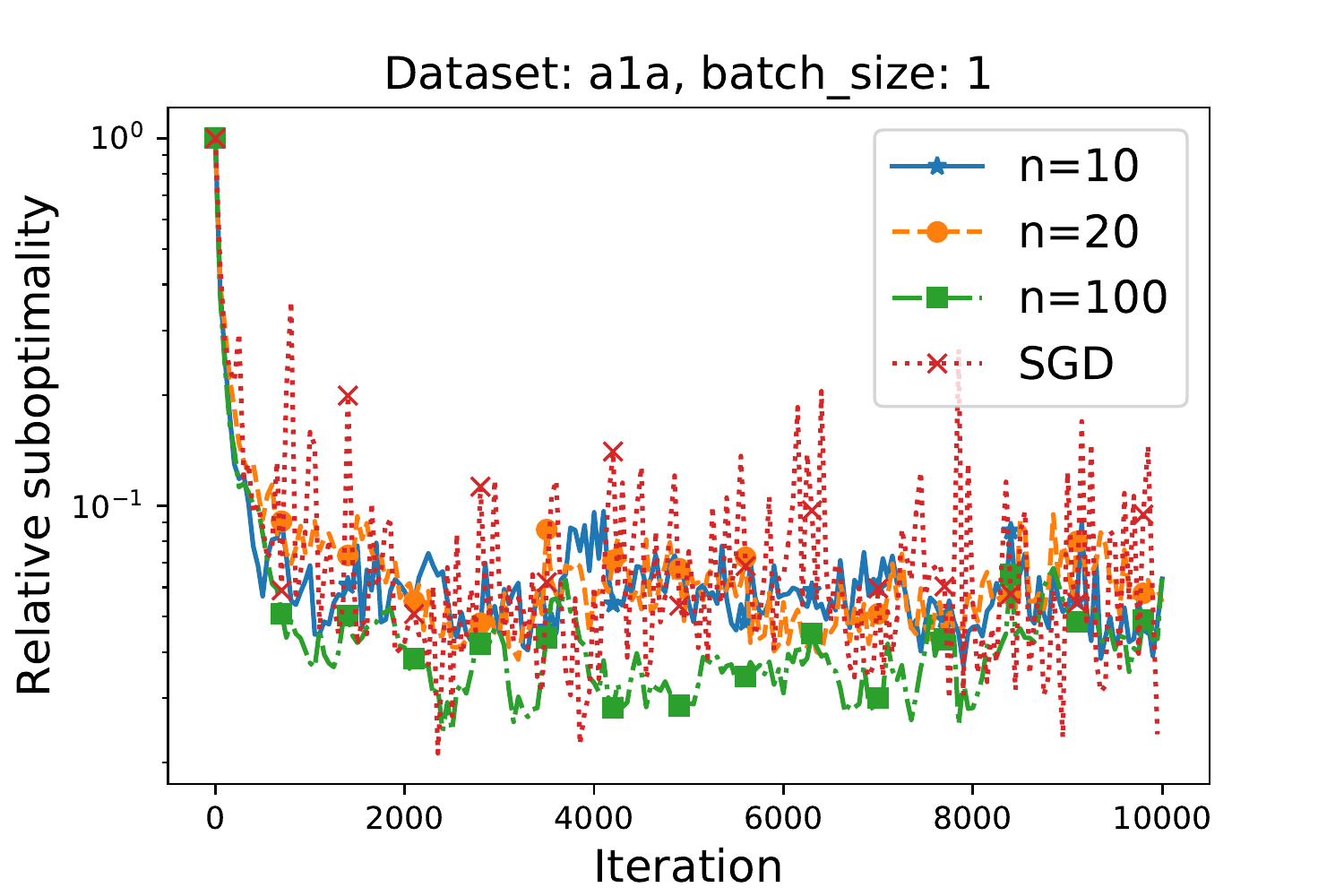}
\end{minipage}%
\begin{minipage}{0.33\textwidth}
  \centering
\includegraphics[width =  \textwidth ]{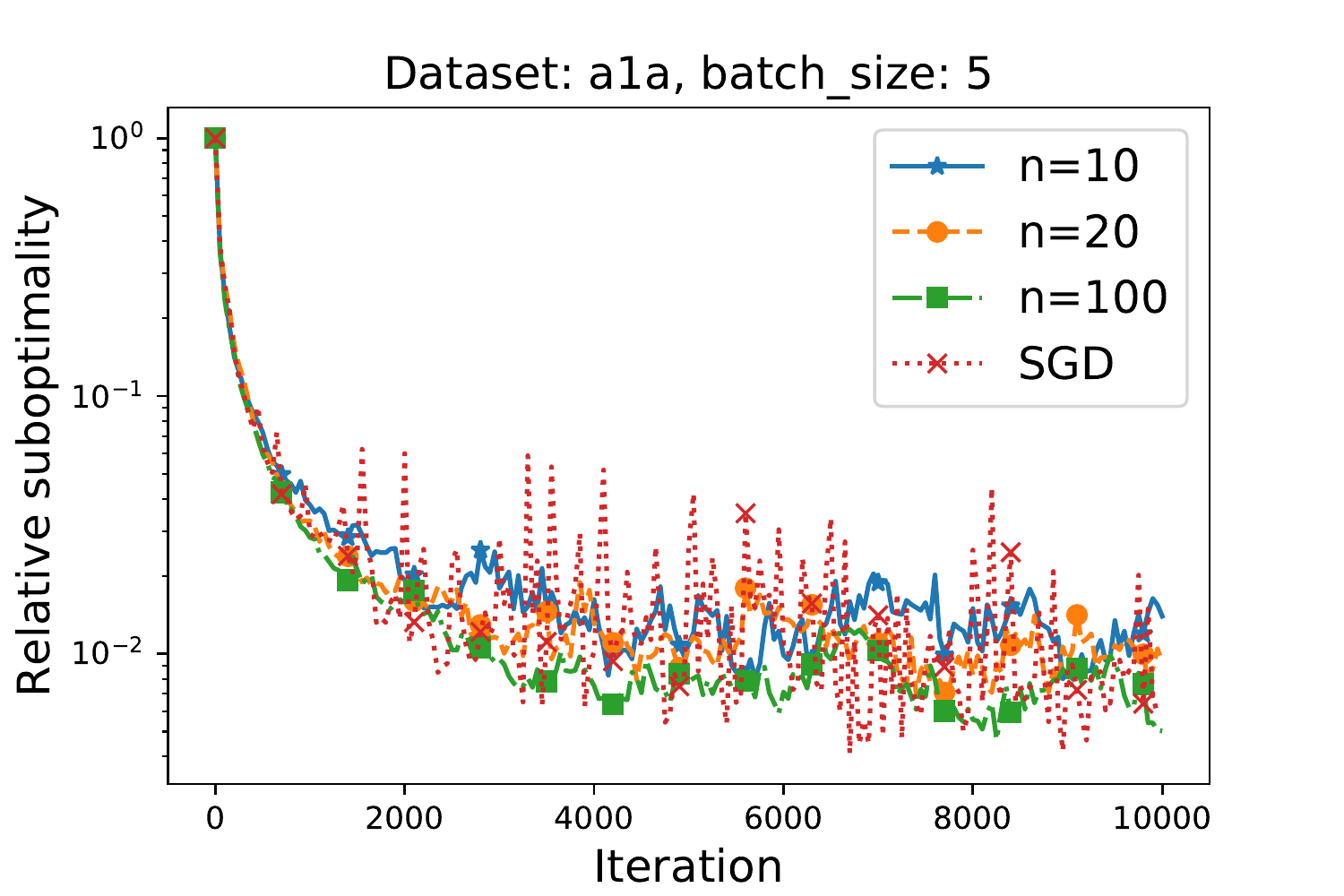}
\end{minipage}%
\begin{minipage}{0.33\textwidth}
  \centering
\includegraphics[width =  \textwidth ]{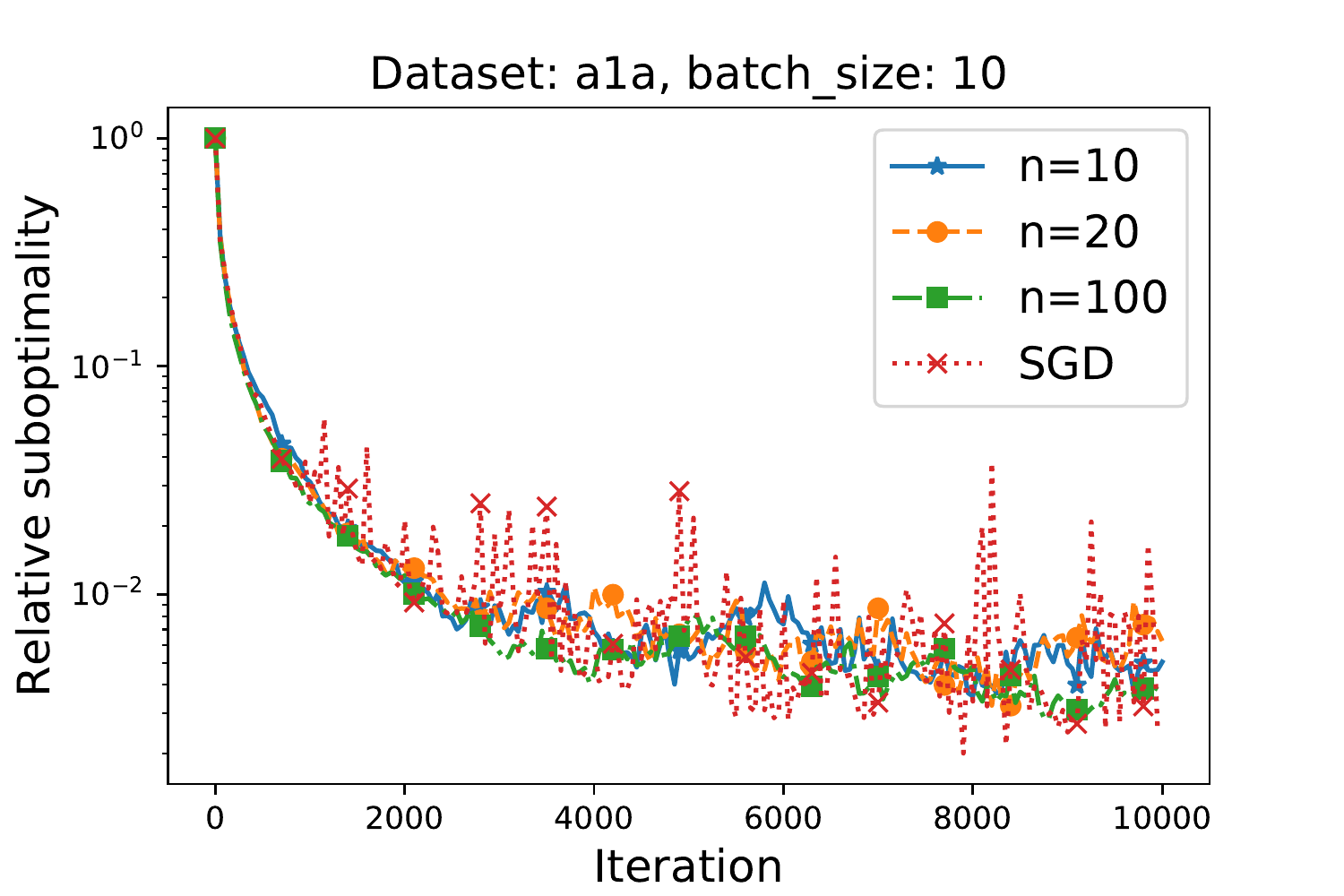}
\end{minipage}%
\\
\begin{minipage}{0.33\textwidth}
  \centering
\includegraphics[width =  \textwidth ]{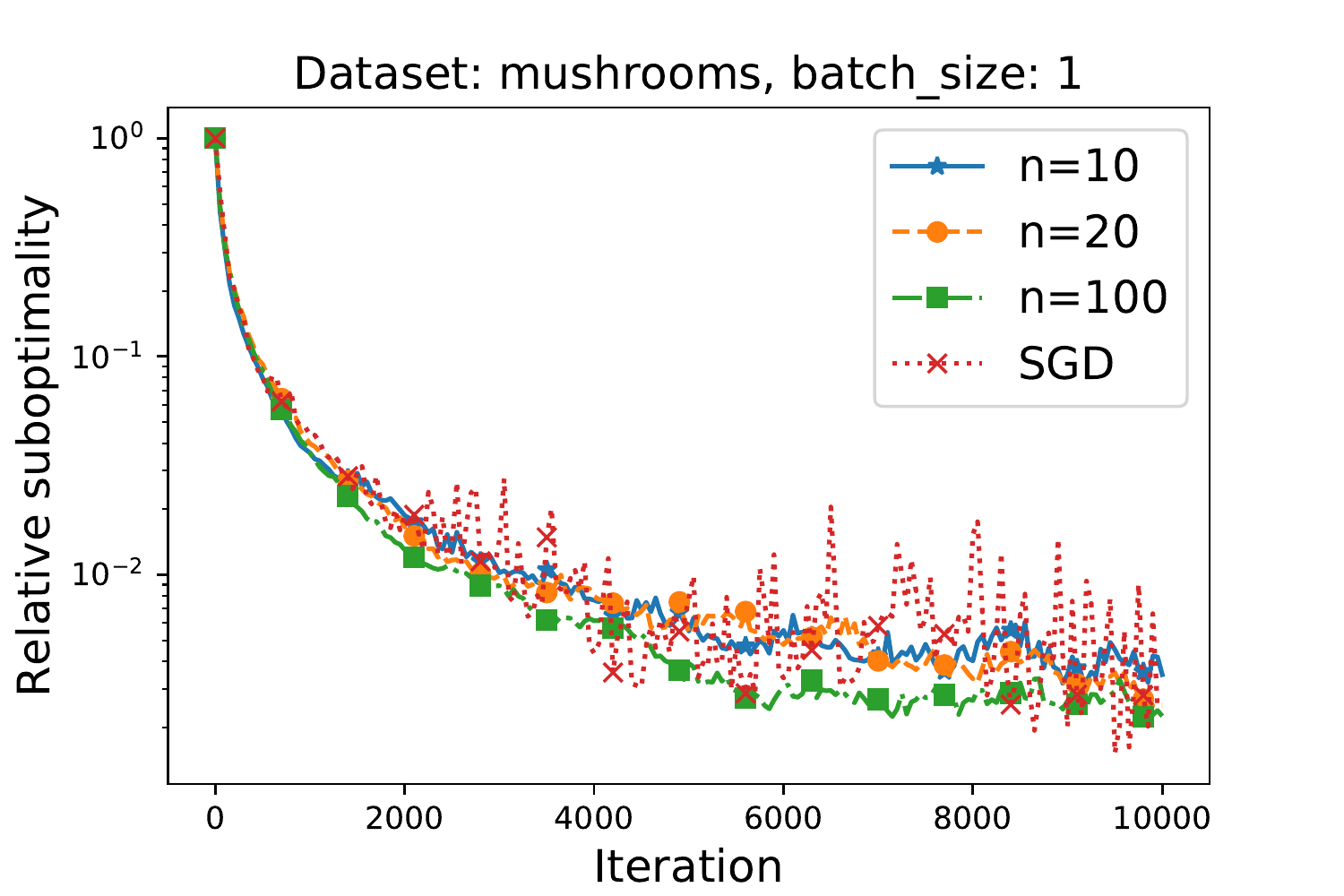}
\end{minipage}%
\begin{minipage}{0.33\textwidth}
  \centering
\includegraphics[width =  \textwidth ]{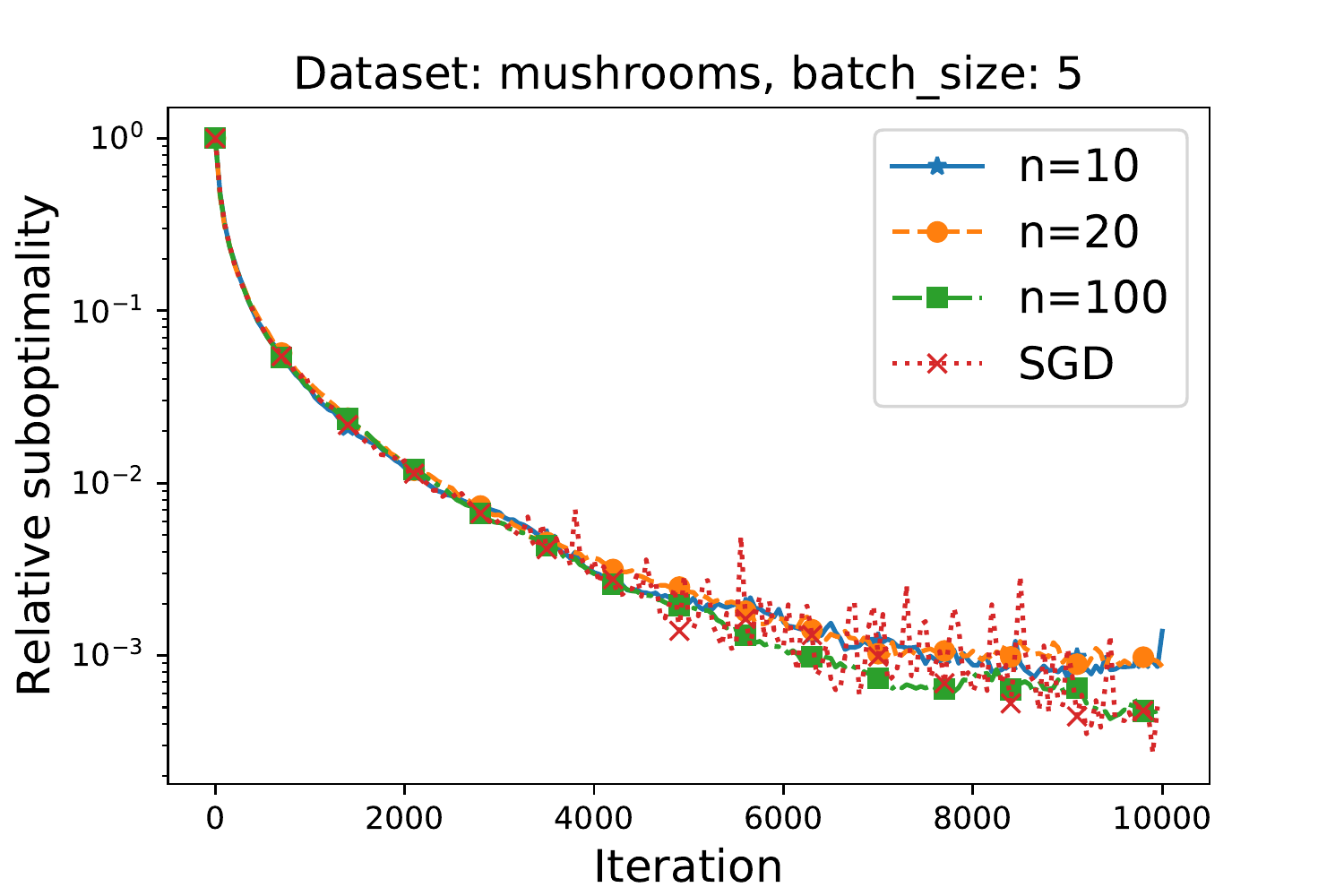}
\end{minipage}%
\begin{minipage}{0.33\textwidth}
  \centering
\includegraphics[width =  \textwidth ]{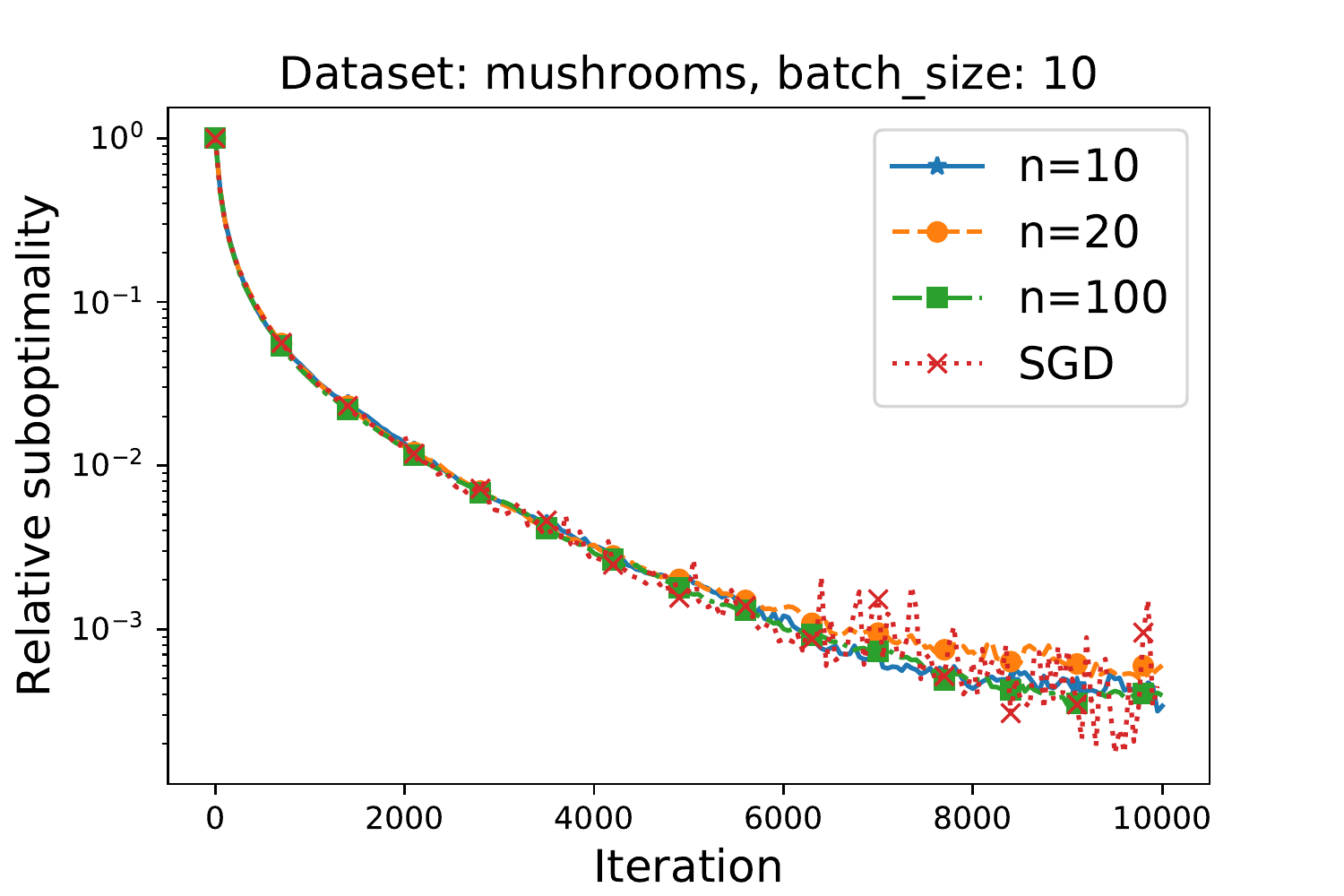}
\end{minipage}%
\\
\begin{minipage}{0.33\textwidth}
  \centering
\includegraphics[width =  \textwidth ]{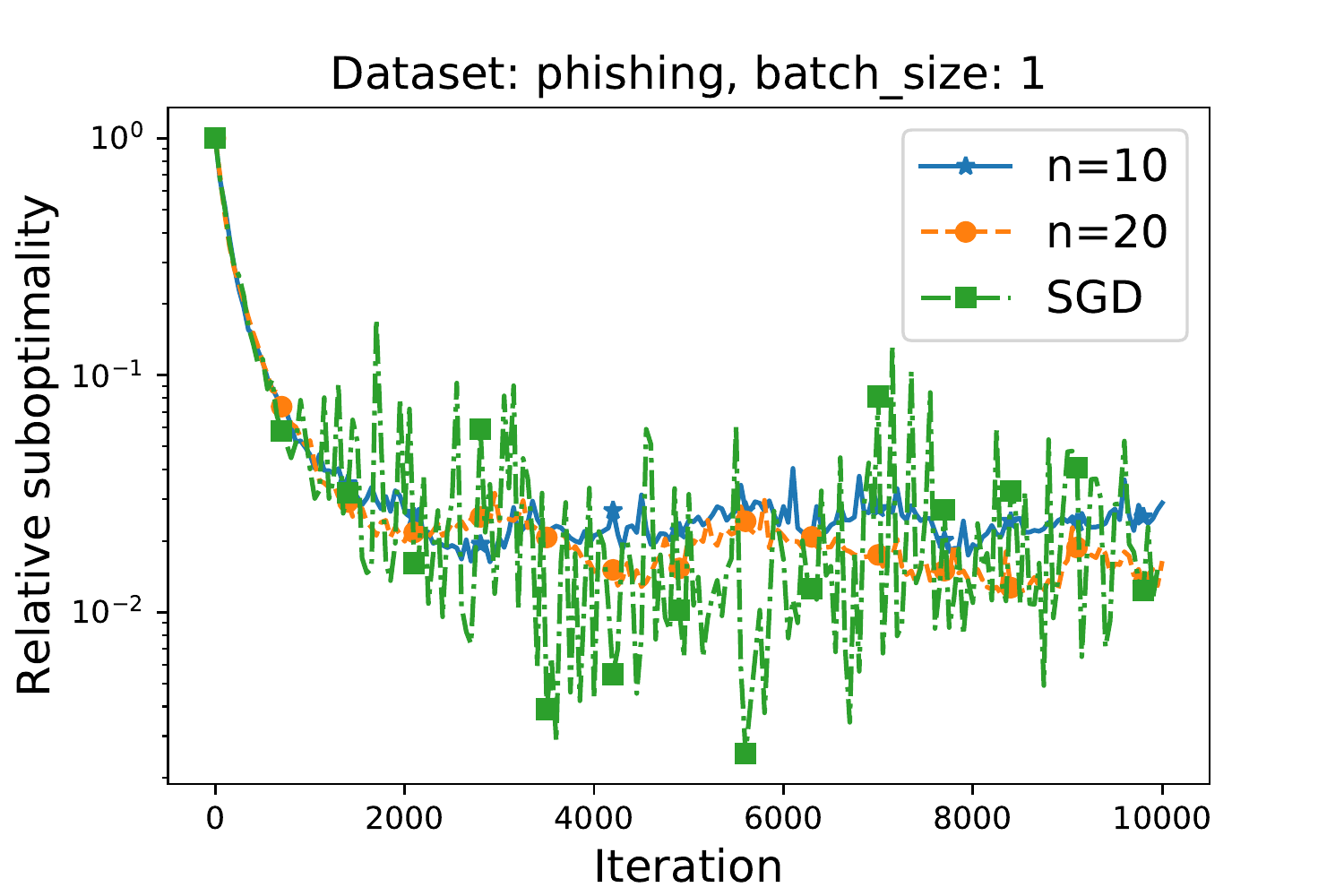}
\end{minipage}%
\begin{minipage}{0.33\textwidth}
  \centering
\includegraphics[width =  \textwidth ]{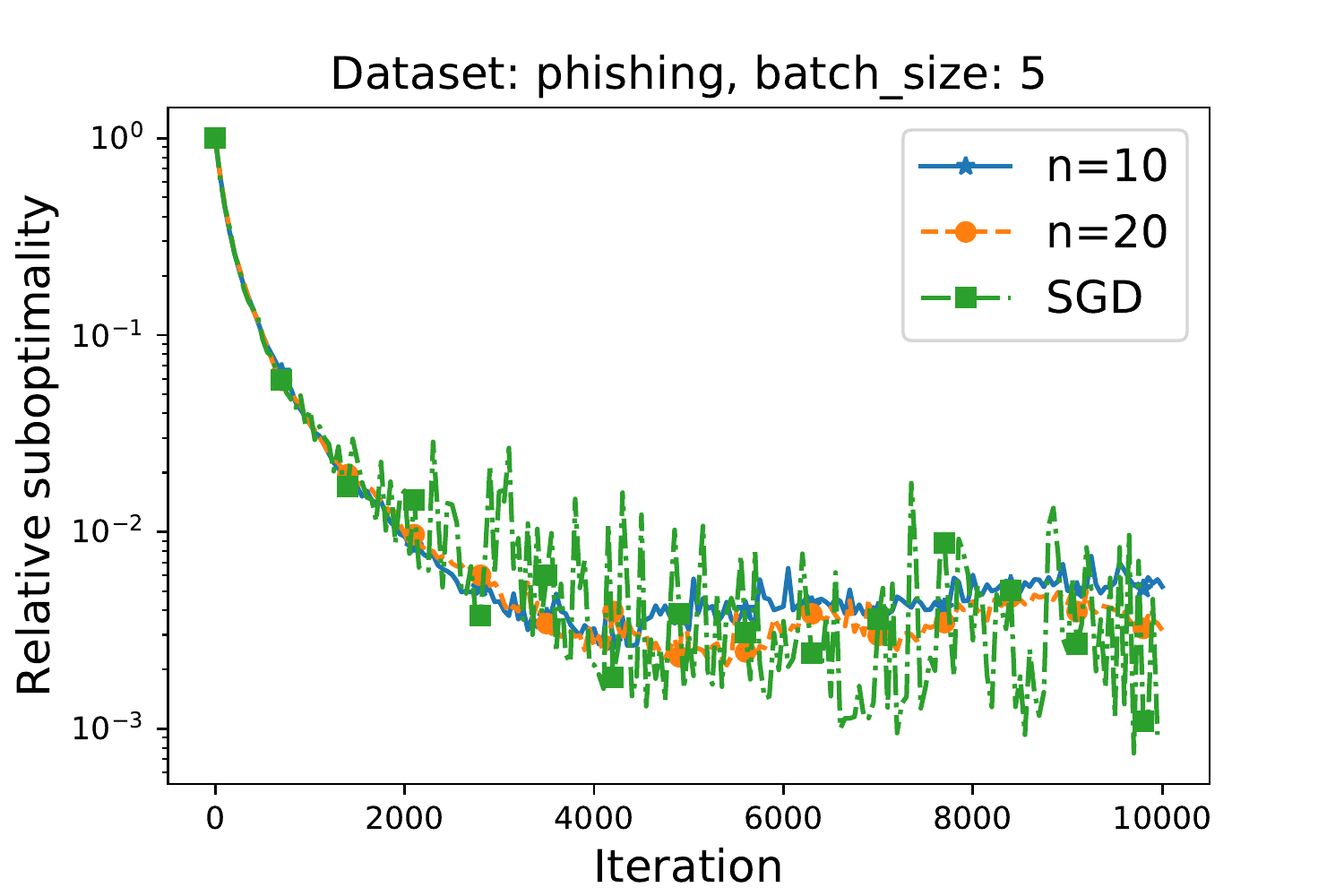}
\end{minipage}%
\begin{minipage}{0.33\textwidth}
  \centering
\includegraphics[width =  \textwidth ]{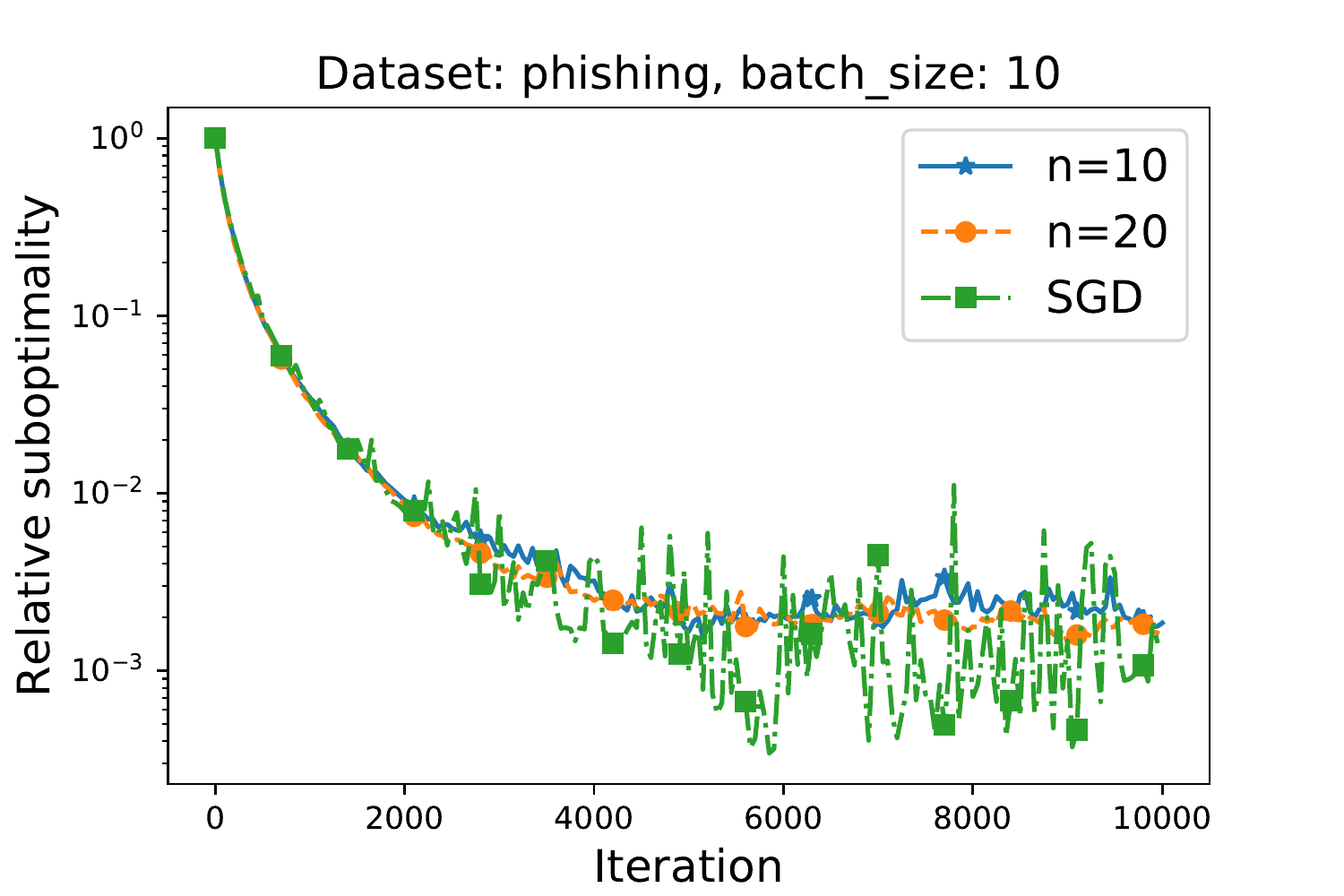}
\end{minipage}%
\\
\begin{minipage}{0.33\textwidth}
  \centering
\includegraphics[width =  \textwidth ]{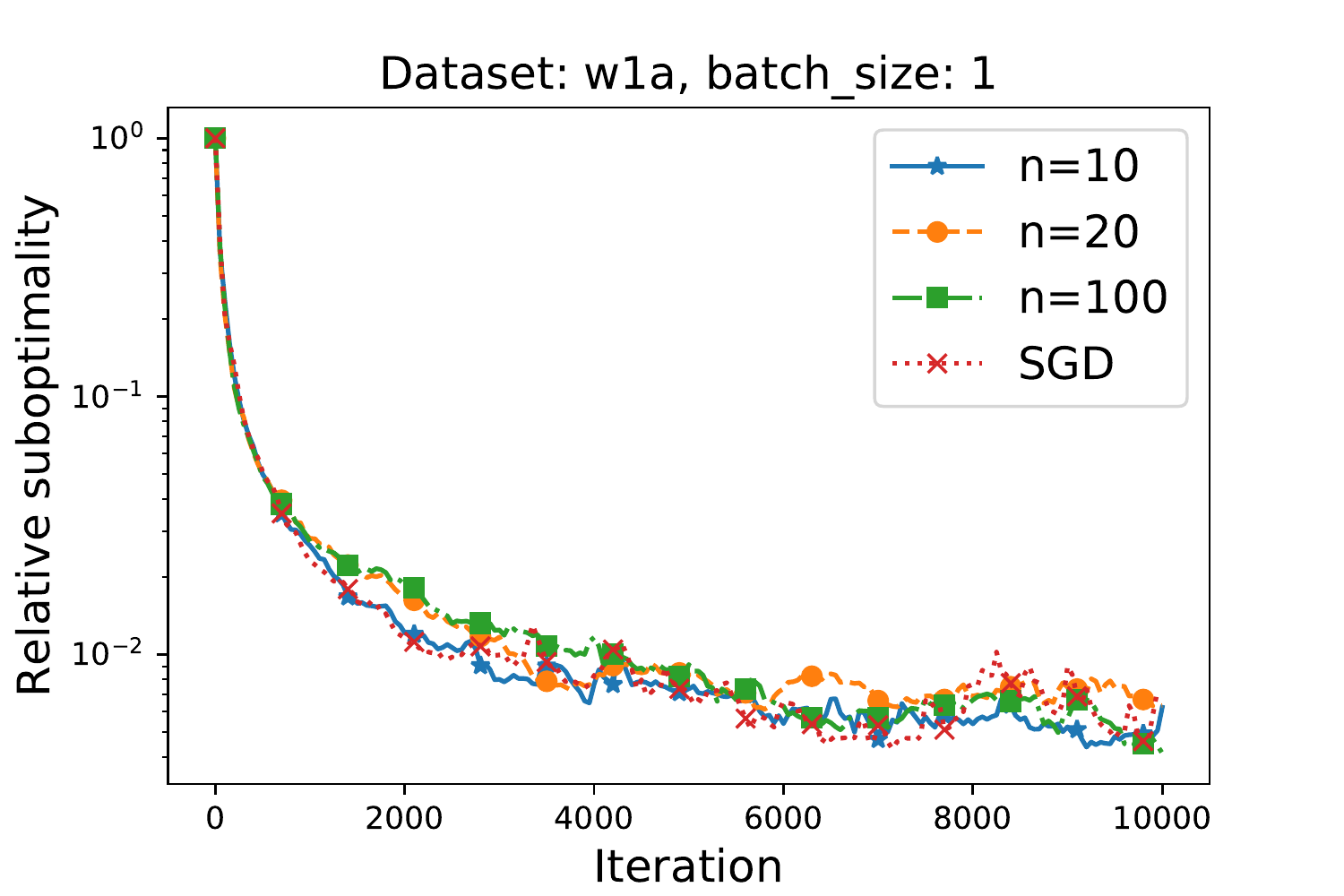}
\end{minipage}%
\begin{minipage}{0.33\textwidth}
  \centering
\includegraphics[width =  \textwidth ]{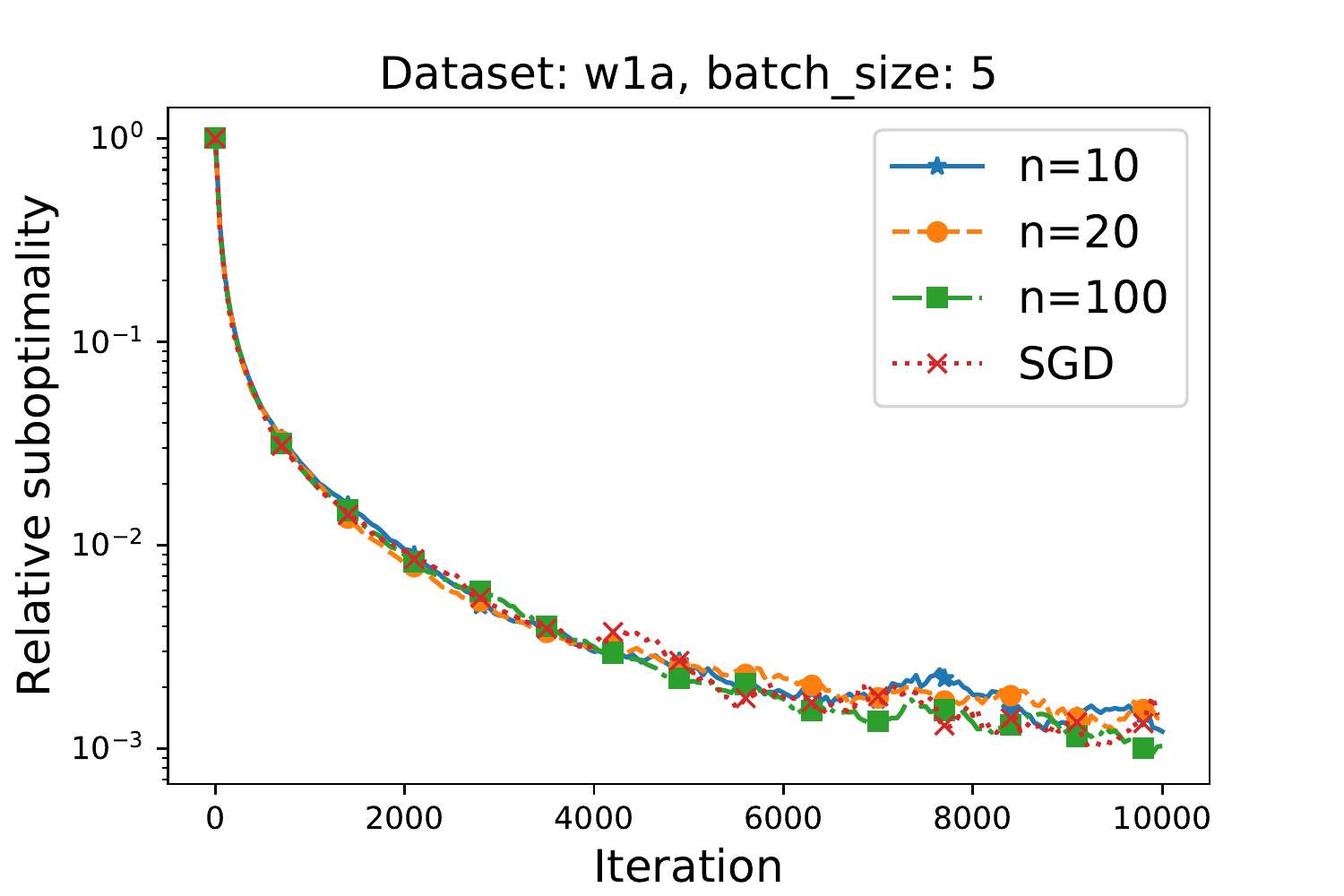}
\end{minipage}%
\begin{minipage}{0.33\textwidth}
  \centering
\includegraphics[width =  \textwidth ]{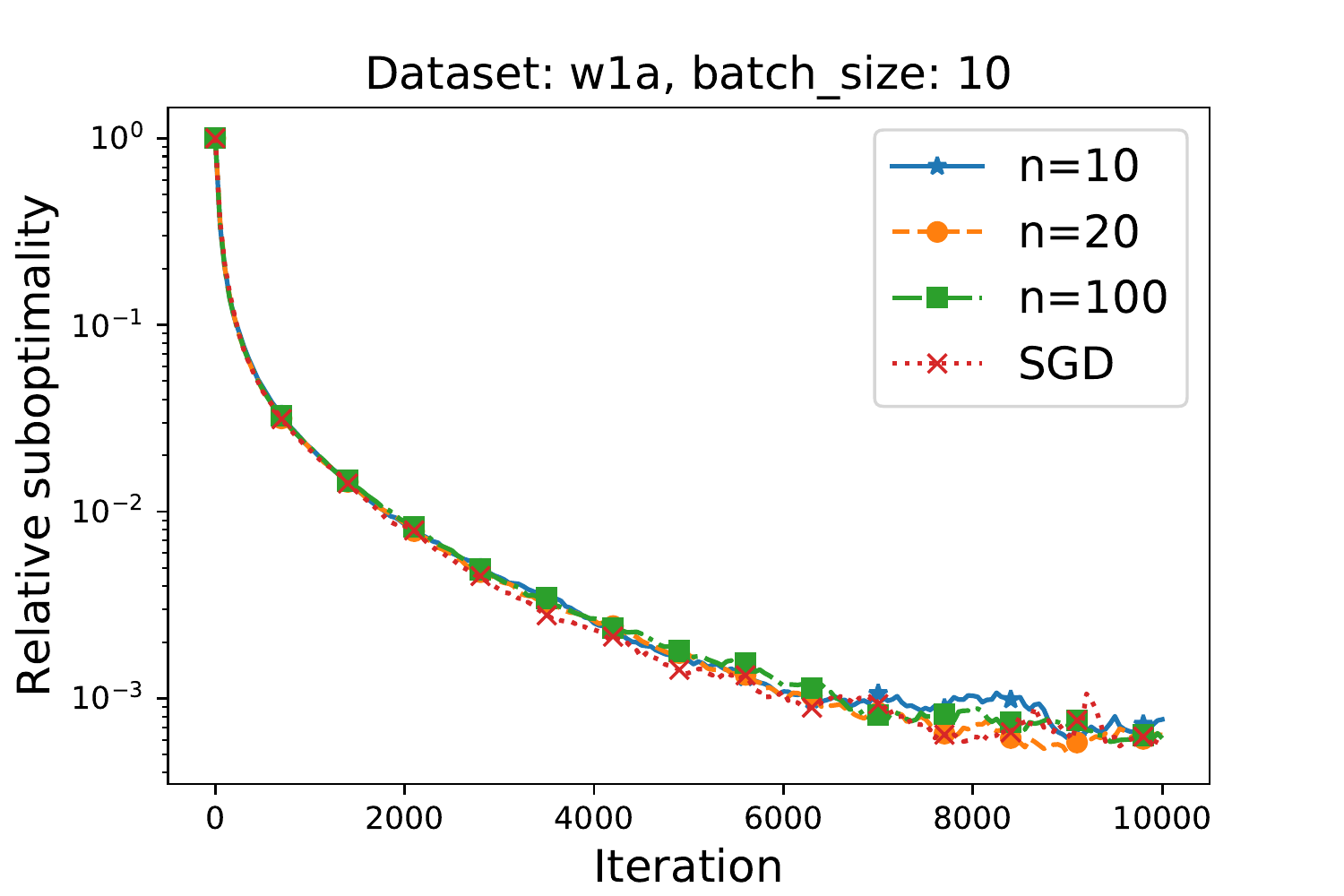}
\end{minipage}%
\\
\caption{Comparison of {\tt SGD} (gradient evaluated on a single datapoint) and Algorithm~\ref{alg:sgd} with $n\tau=1$. Constant $\alpha  = \frac{1}{5L}$ was used for each algorithm. Label ``batch\_size'' indicates how big minibatch was chosen for stochastic gradient of each worker's objective.} \label{fig:99_sgd1}
\end{figure}

Next, we study the dependence of the convergence speed on $\tau$ for various values of $n$. Figure~\ref{fig:99_sgd2} presents the results. In each case, $\tau$ influences the convergence rate (or the region where the iterates oscillate) significantly, however, the effect is much weaker for larger $n$. This is in correspondence with Corollary~\ref{cor:99_sgd}.

\begin{figure}[H]
\centering
\begin{minipage}{0.33\textwidth}
  \centering
\includegraphics[width =  \textwidth ]{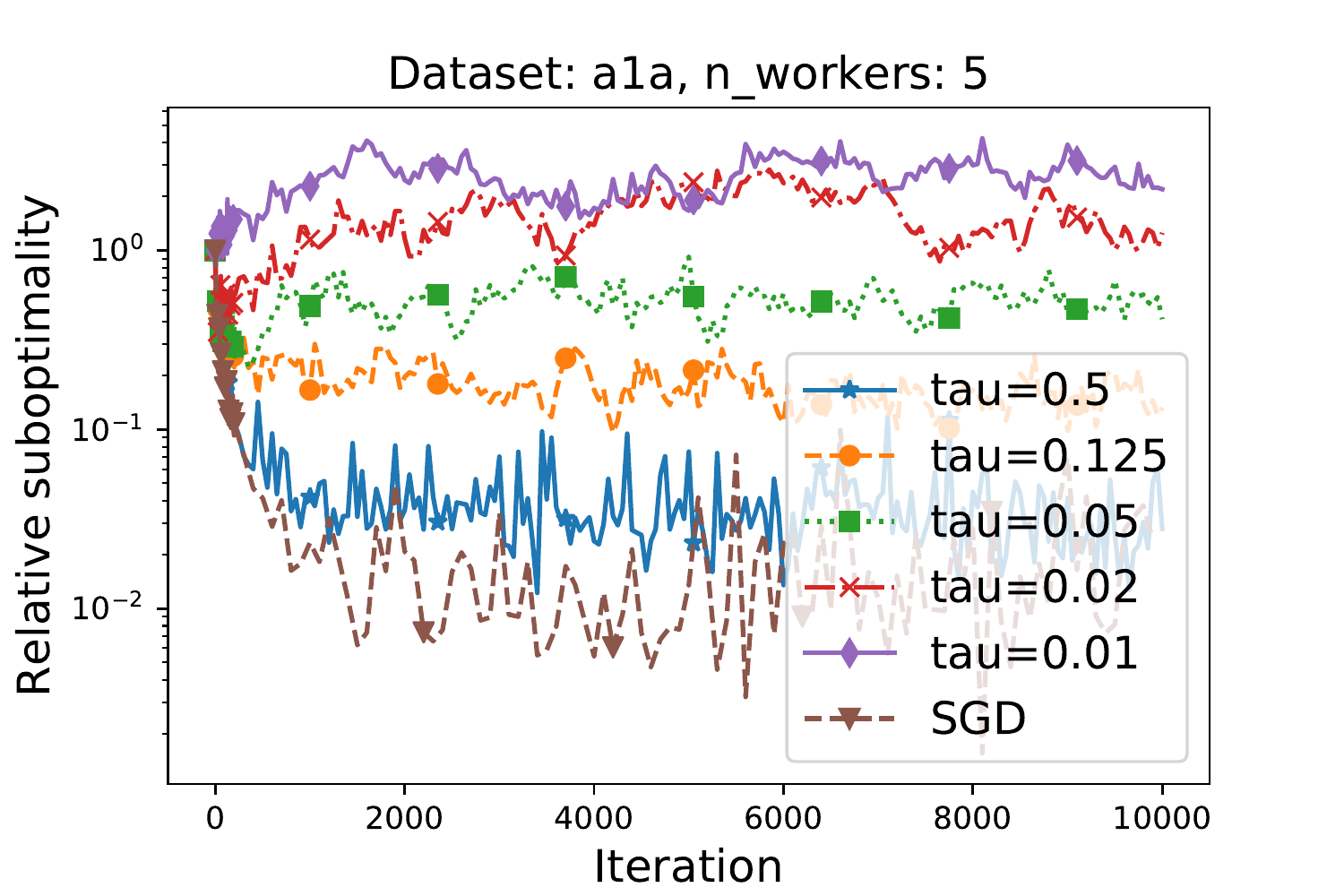}
\end{minipage}%
\begin{minipage}{0.33\textwidth}
  \centering
\includegraphics[width =  \textwidth ]{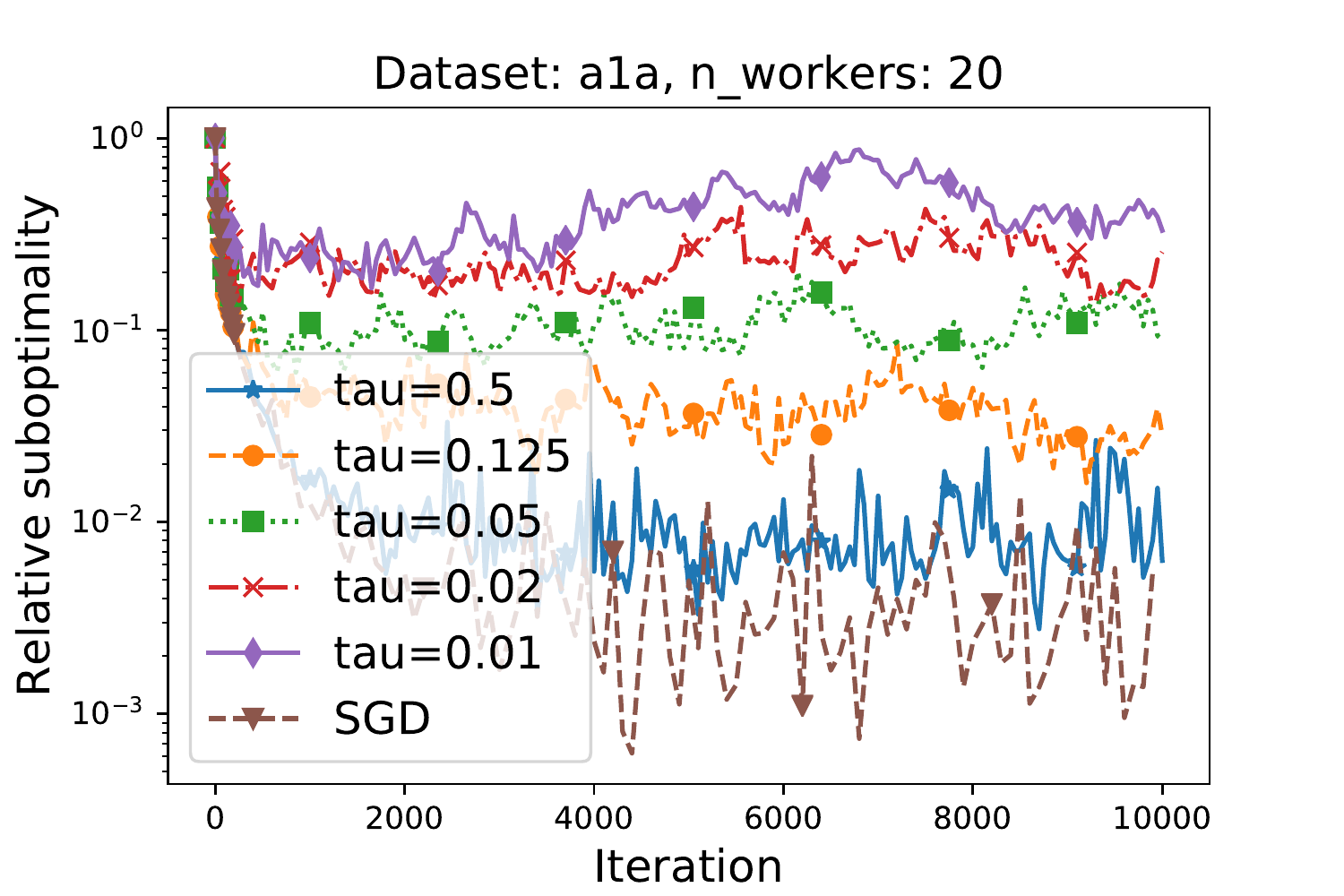}
\end{minipage}%
\begin{minipage}{0.33\textwidth}
  \centering
\includegraphics[width =  \textwidth ]{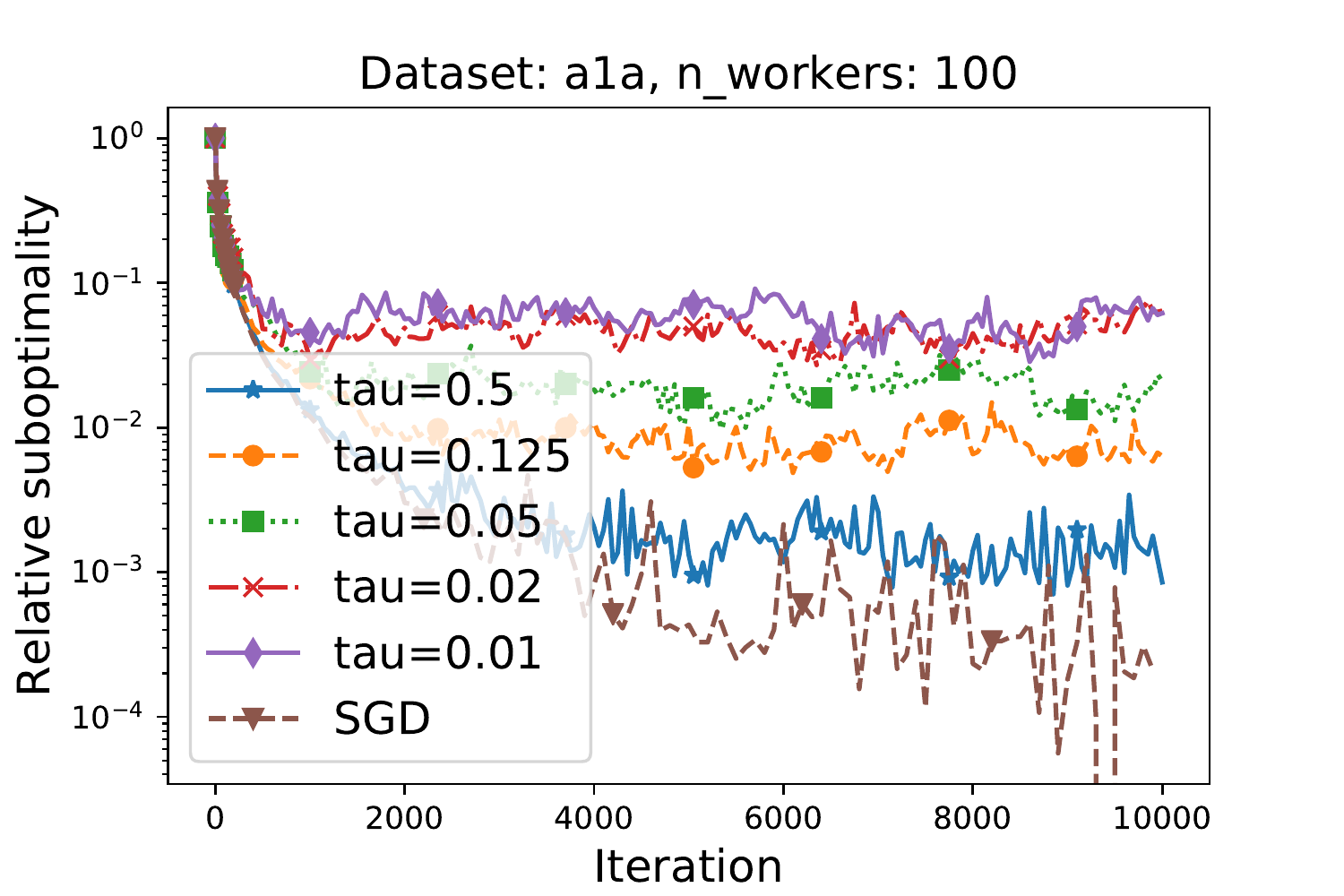}
\end{minipage}%
\\
\begin{minipage}{0.33\textwidth}
  \centering
\includegraphics[width =  \textwidth ]{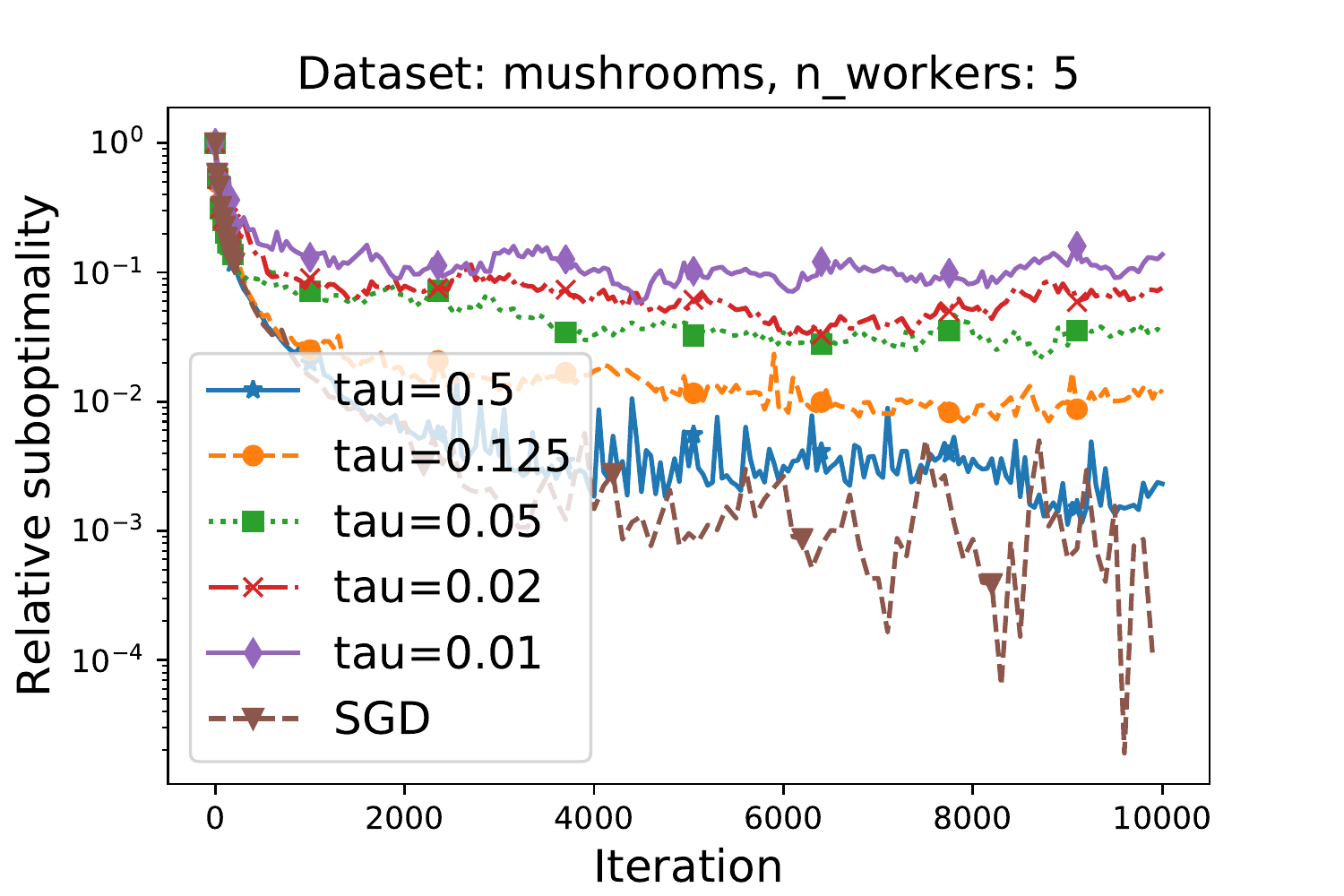}
\end{minipage}%
\begin{minipage}{0.33\textwidth}
  \centering
\includegraphics[width =  \textwidth ]{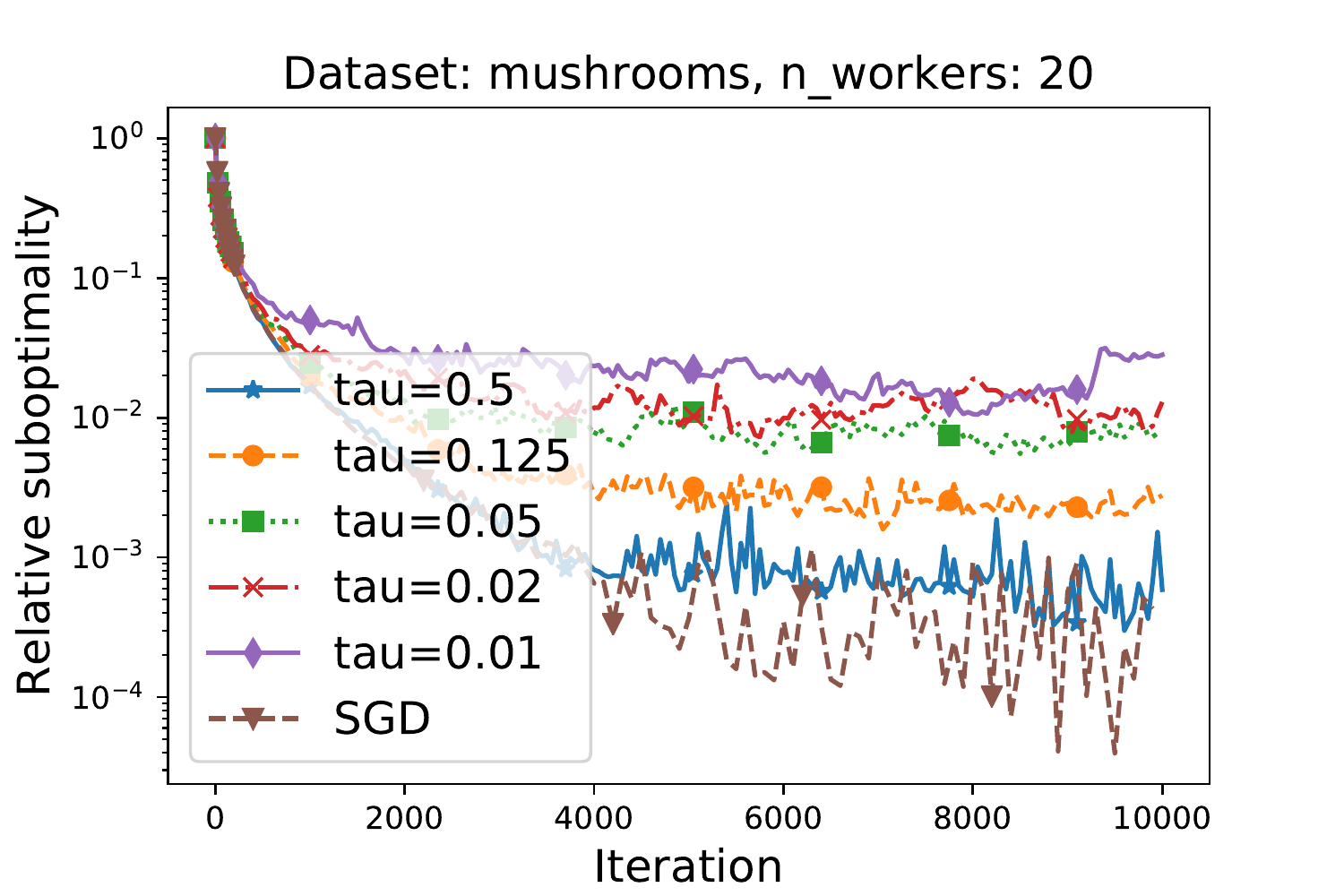}
\end{minipage}%
\begin{minipage}{0.33\textwidth}
  \centering
\includegraphics[width =  \textwidth ]{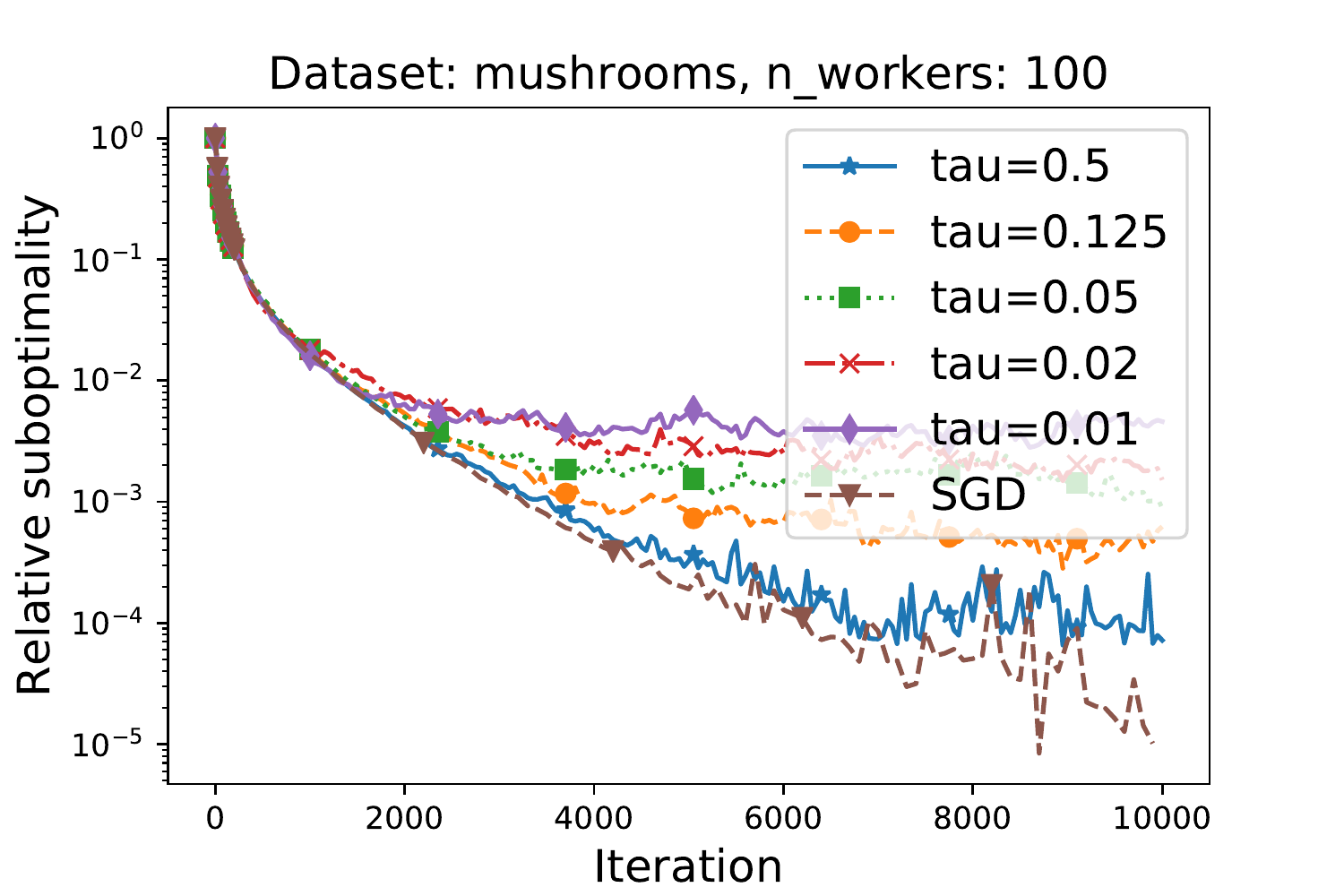}
\end{minipage}%
\\
\begin{minipage}{0.33\textwidth}
  \centering
\includegraphics[width =  \textwidth ]{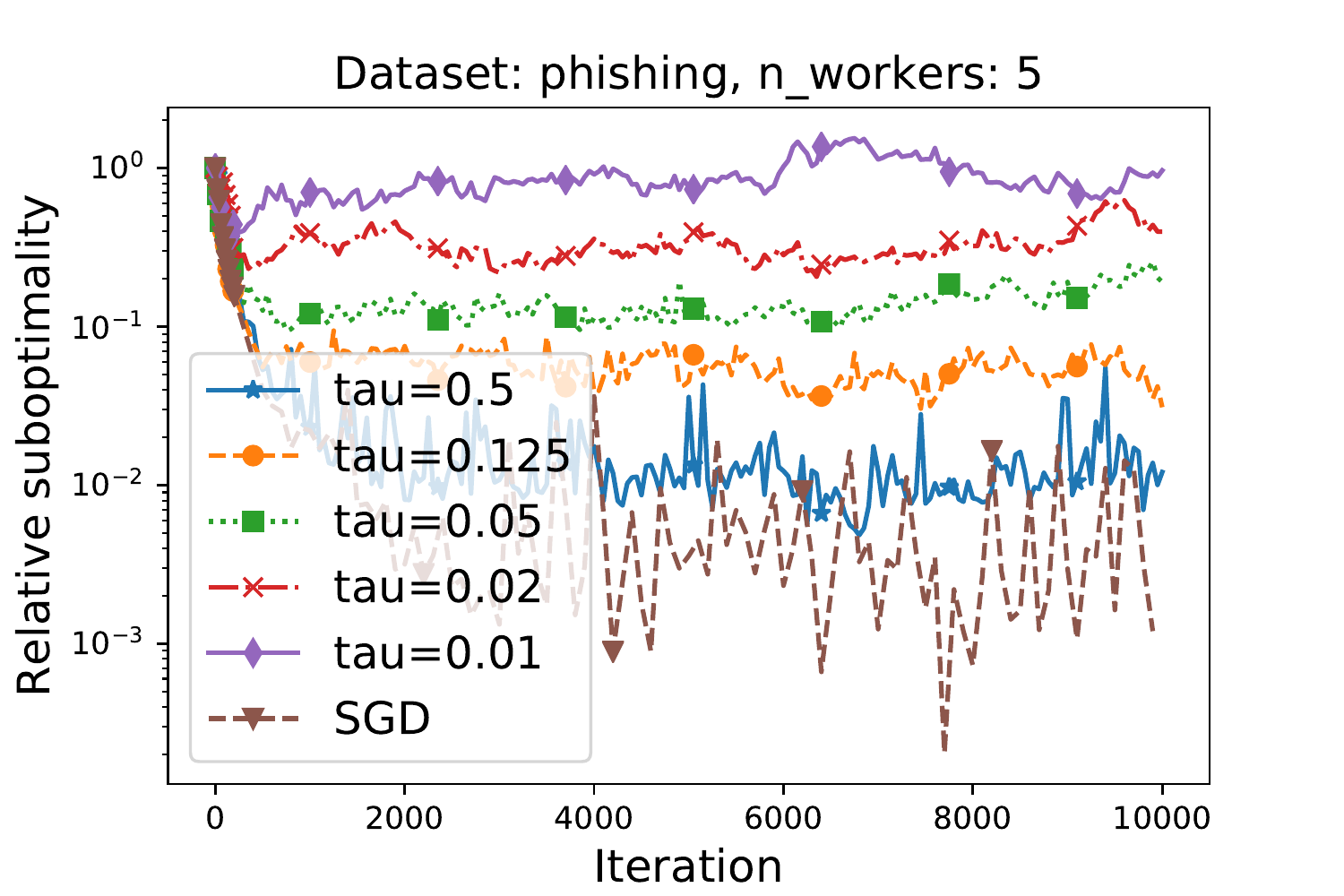}
\end{minipage}%
\begin{minipage}{0.33\textwidth}
  \centering
\includegraphics[width =  \textwidth ]{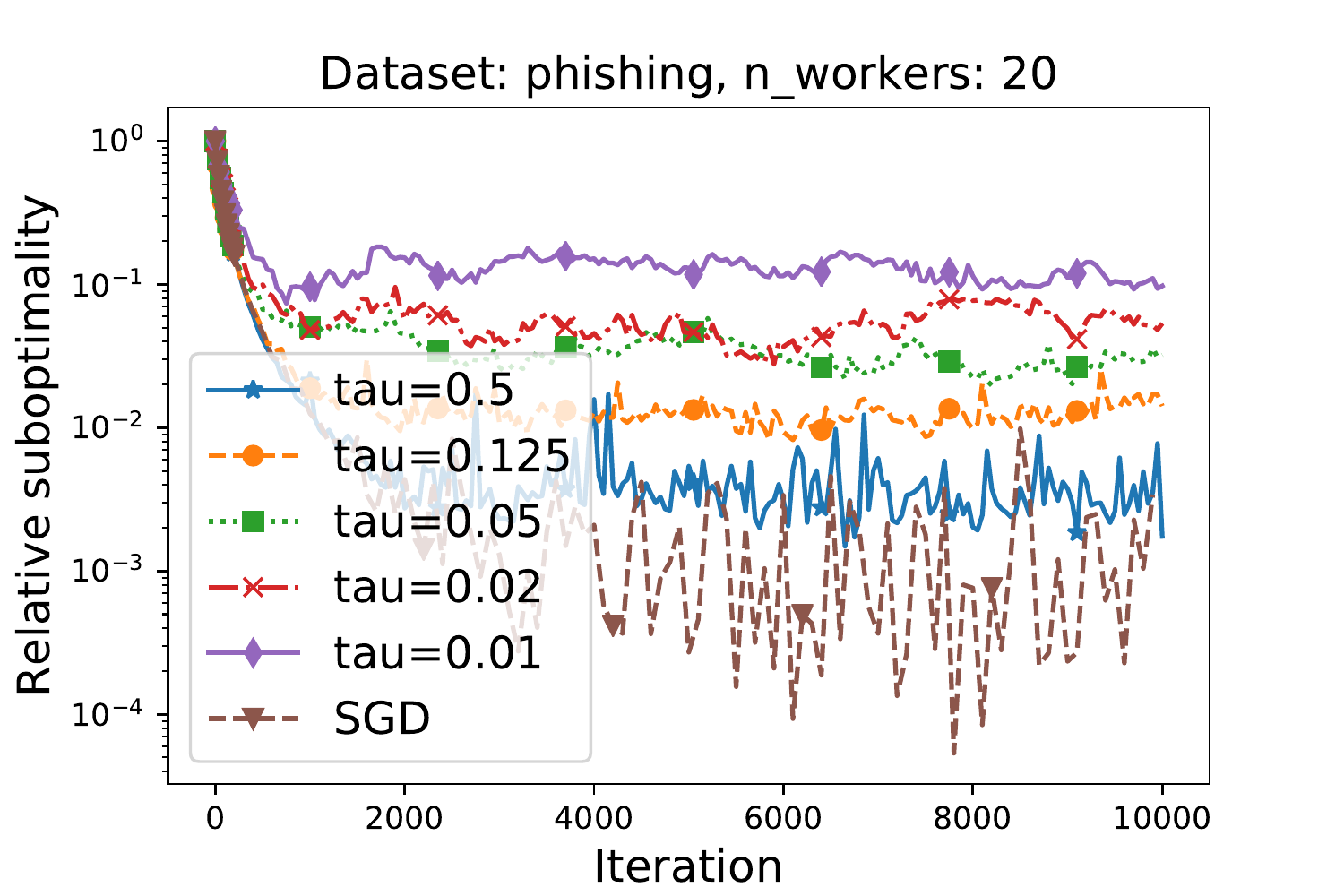}
\end{minipage}%
\\
\begin{minipage}{0.33\textwidth}
  \centering
\includegraphics[width =  \textwidth ]{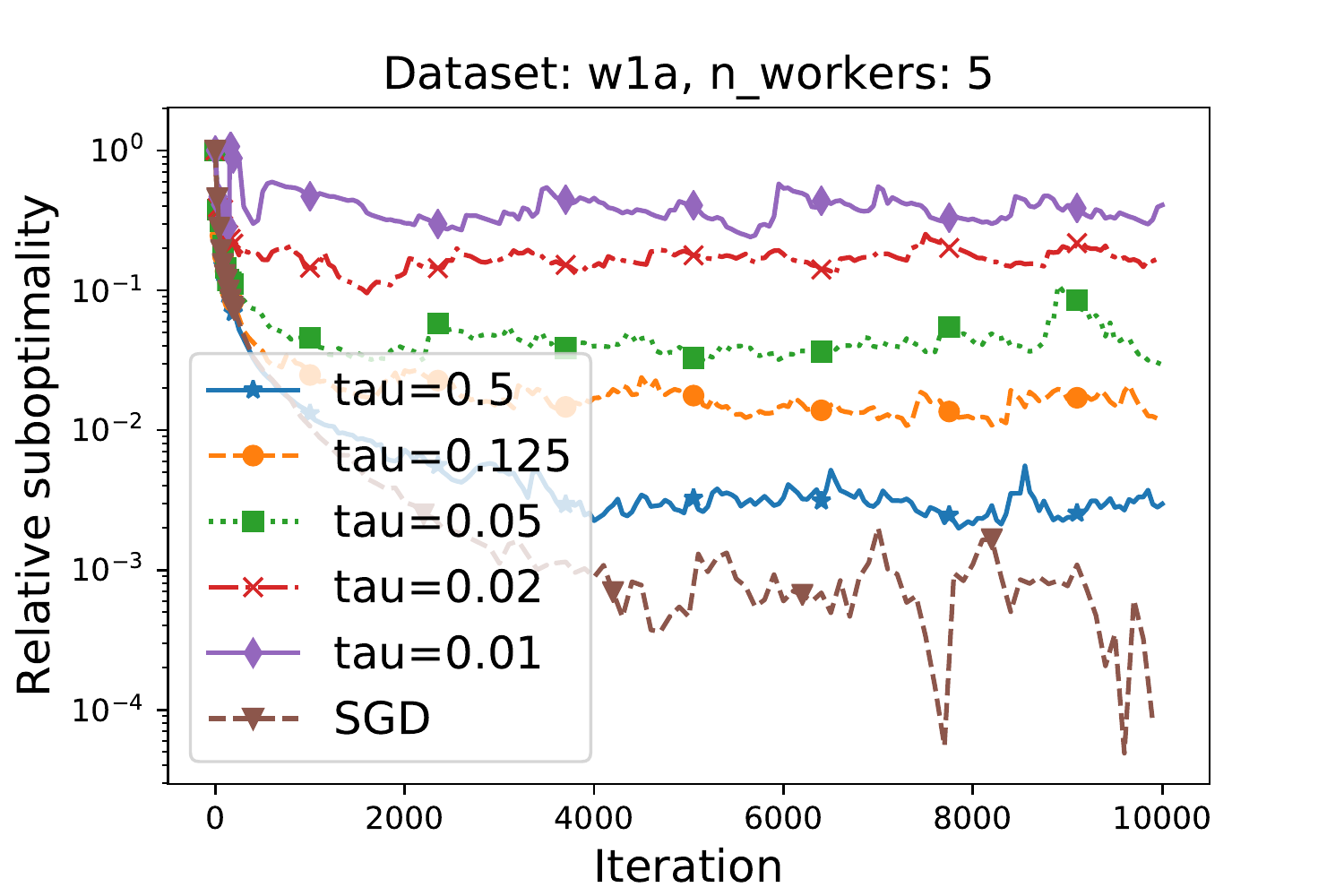}
\end{minipage}%
\begin{minipage}{0.33\textwidth}
  \centering
\includegraphics[width =  \textwidth ]{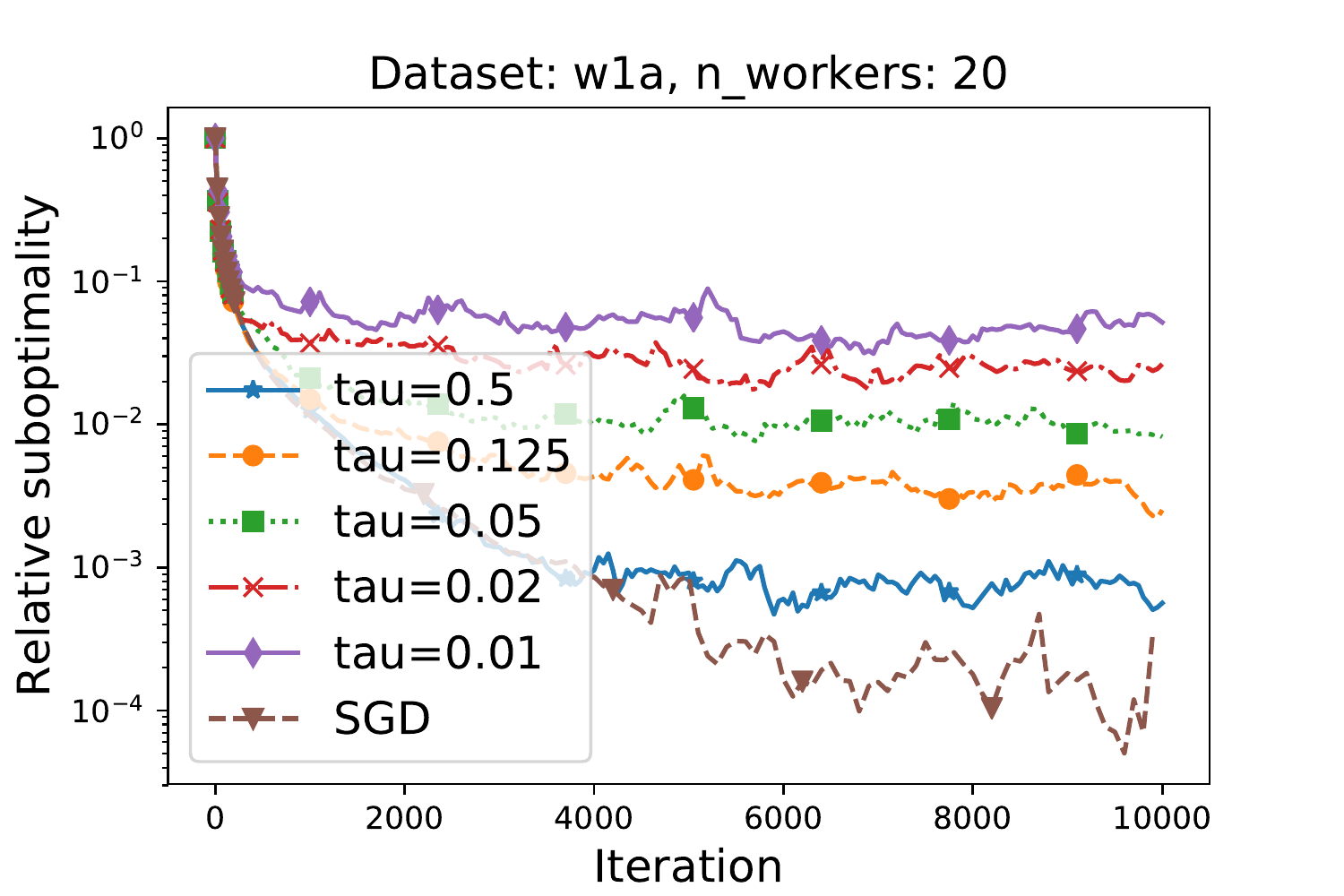}
\end{minipage}%
\begin{minipage}{0.33\textwidth}
  \centering
\includegraphics[width =  \textwidth ]{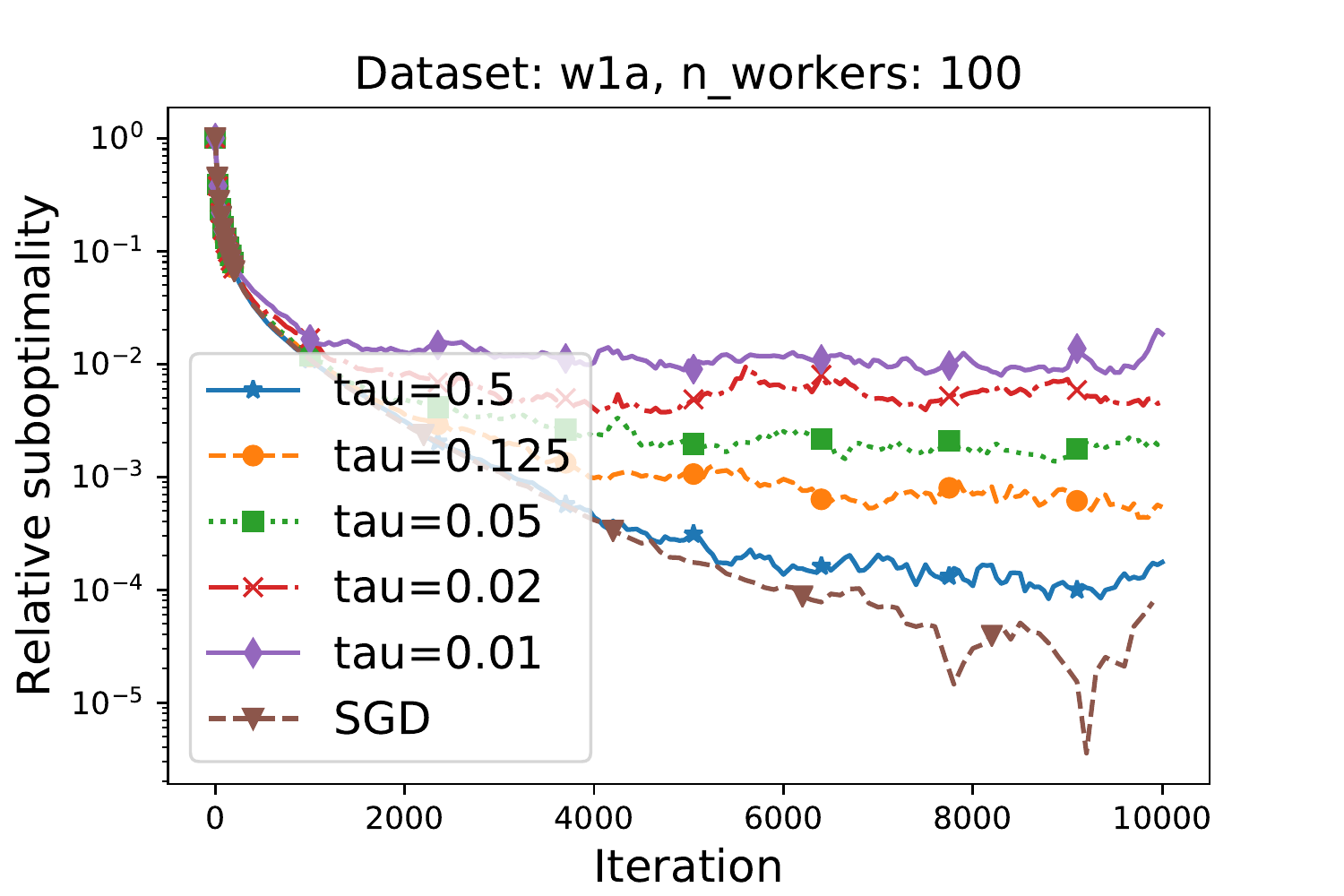}
\end{minipage}%
\\
\caption{Behavior of Algorithm~\ref{alg:sgd} while varying $\tau$. Label ``{\tt SGD}'' corresponds to the choice $n=1, \tau = 1$. Stepsize $\alpha = \frac1{3L}$ was used in every case.}\label{fig:99_sgd2}
\end{figure}

\subsection{{\tt IASGD} \label{sec:99_exp_asgd}}

In this section we numerically test Algorithm~\ref{alg:acc} for logistic regression problem. As in the last section, $f_i$ consists of set of (uniformly distributed) rows of $\mA$ from~\eqref{eq:99_logreg}. The stochastic gradient is taken as a gradient on a subset data points from each $f_i$. Note that Algorithm~\ref{alg:acc} depends on a priori unknown strong growth parameter $\hat{\rho}$ of unbiased stochastic gradient $q$\footnote{Formulas to obtain parameters of Algorithm~\ref{alg:acc} are given in~\cite{vaswani2019-overparam}. }.  Therefore, we first find empirically optimal $\hat{\rho}$ for each algorithm run by grid search and report only the best performance for each algorithm. 

The first experiment (Figure~\ref{fig:99_acc1}) verifies the linearity claim -- we vary $(n,\tau)$ such that $n\tau=1$. As predicted by theory, the behavior of presented algorithms is almost indistinguishable. 

\begin{figure}[H]
\centering
\begin{minipage}{0.33\textwidth}
  \centering
\includegraphics[width =  \textwidth ]{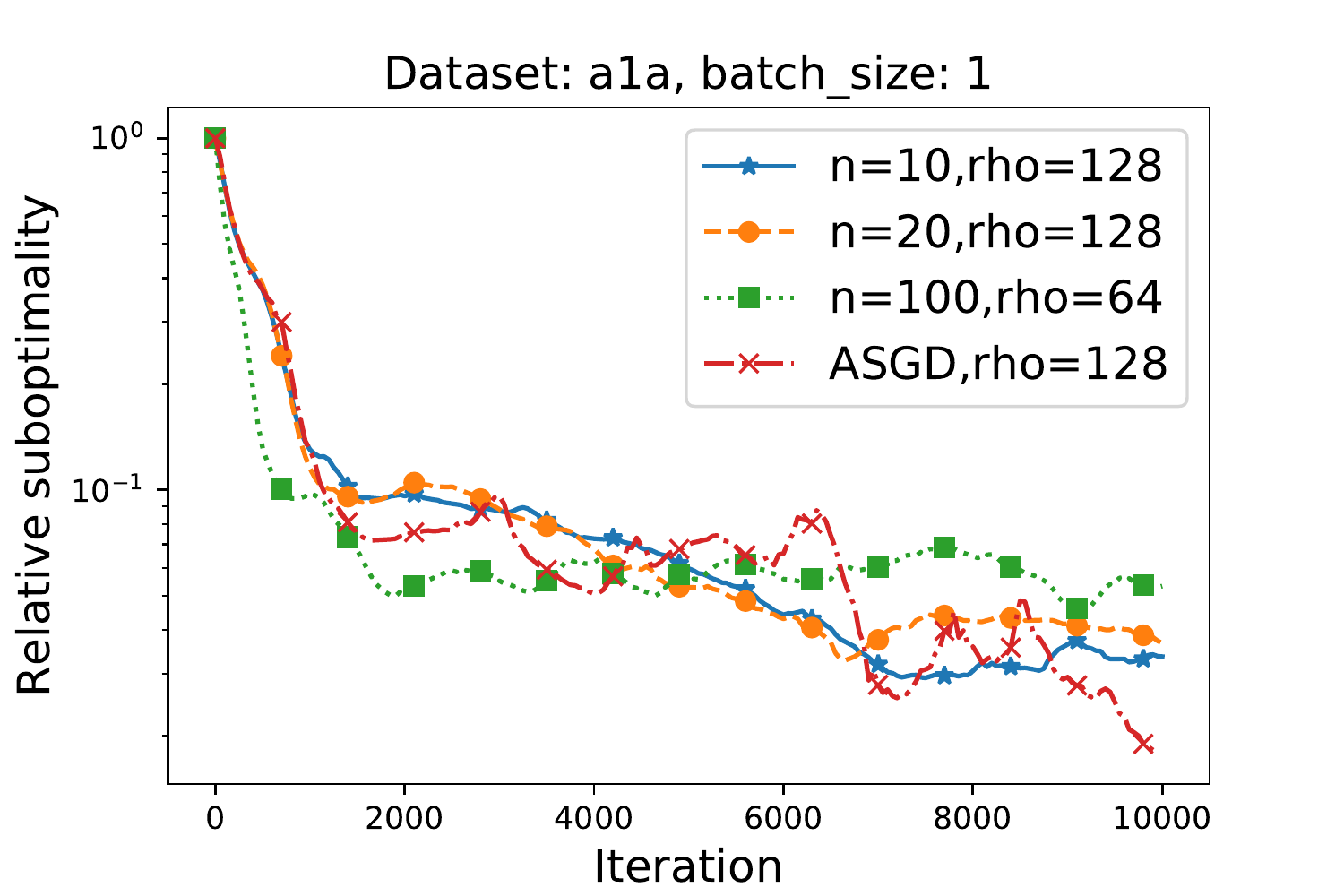}
\end{minipage}%
\begin{minipage}{0.33\textwidth}
  \centering
\includegraphics[width =  \textwidth ]{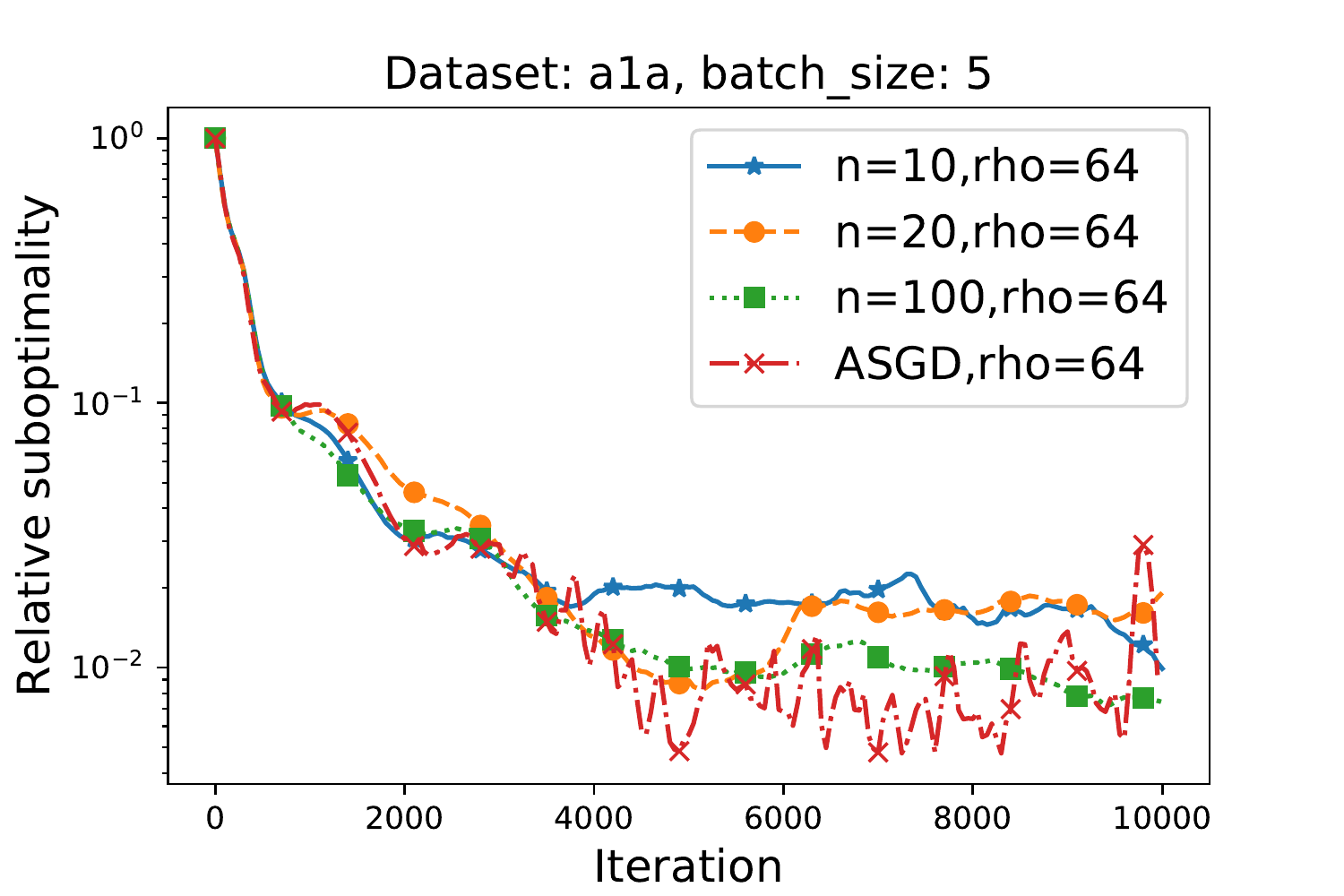}
\end{minipage}%
\begin{minipage}{0.33\textwidth}
  \centering
\includegraphics[width =  \textwidth ]{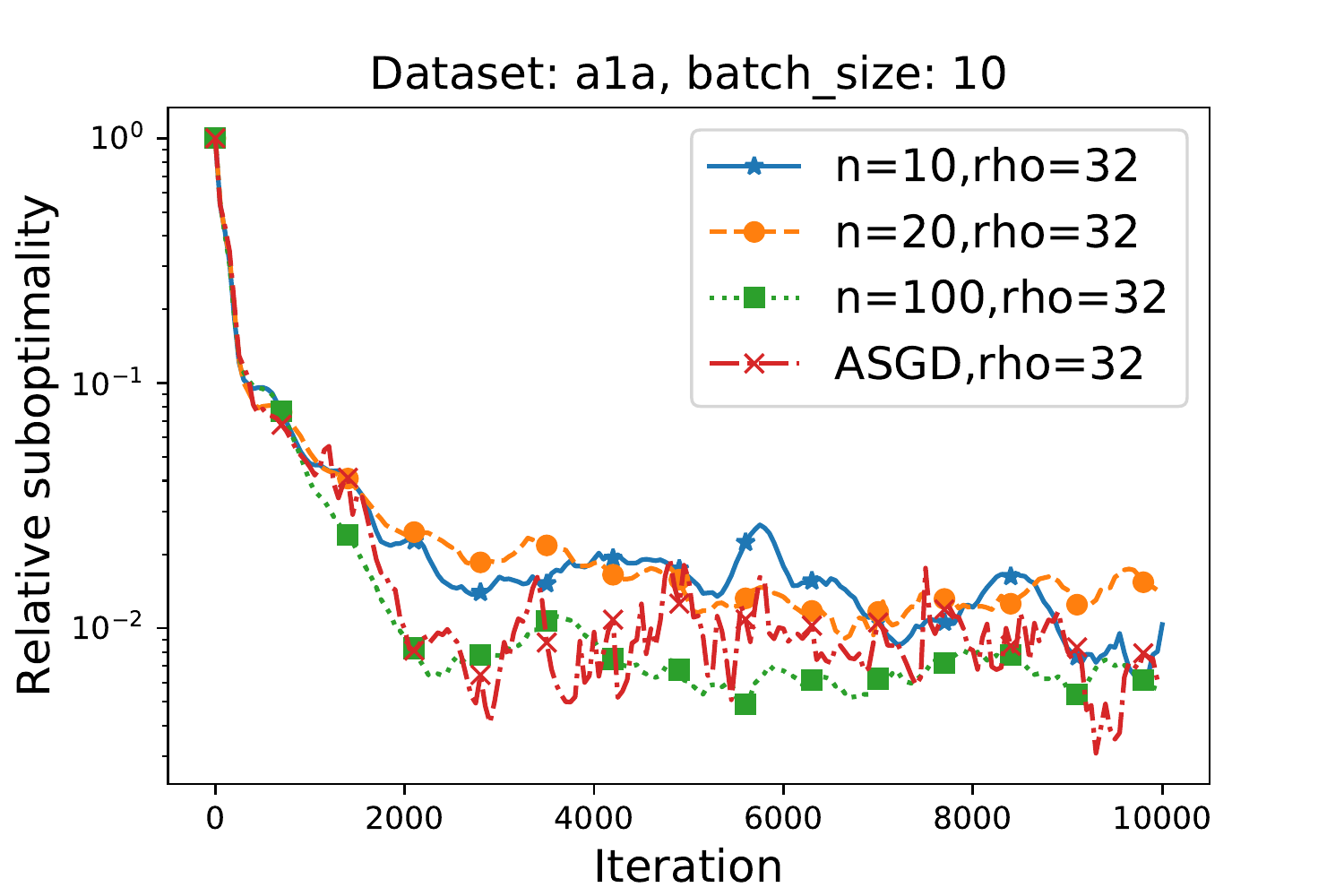}
\end{minipage}%
\\
\begin{minipage}{0.33\textwidth}
  \centering
\includegraphics[width =  \textwidth ]{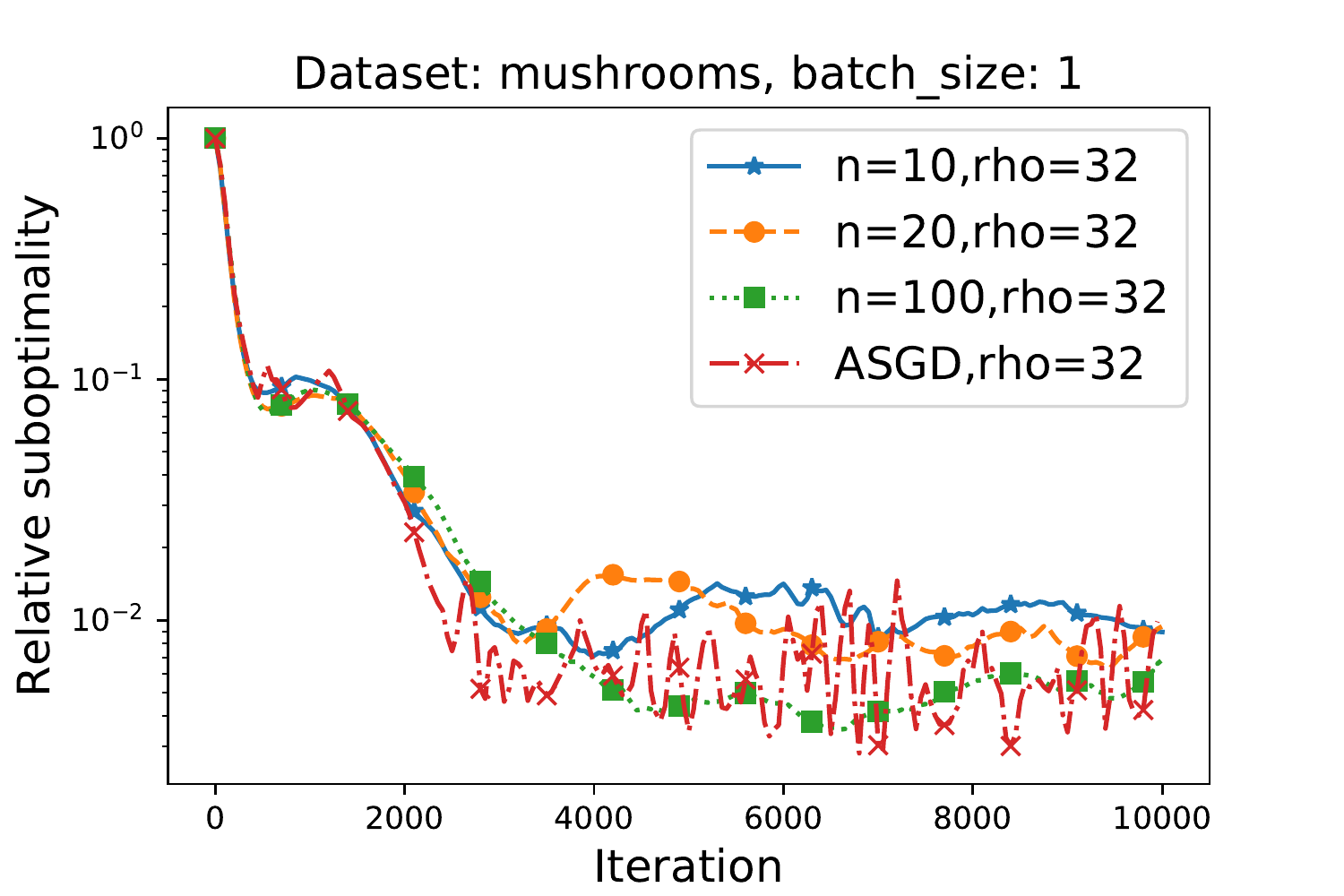}
\end{minipage}%
\begin{minipage}{0.33\textwidth}
  \centering
\includegraphics[width =  \textwidth ]{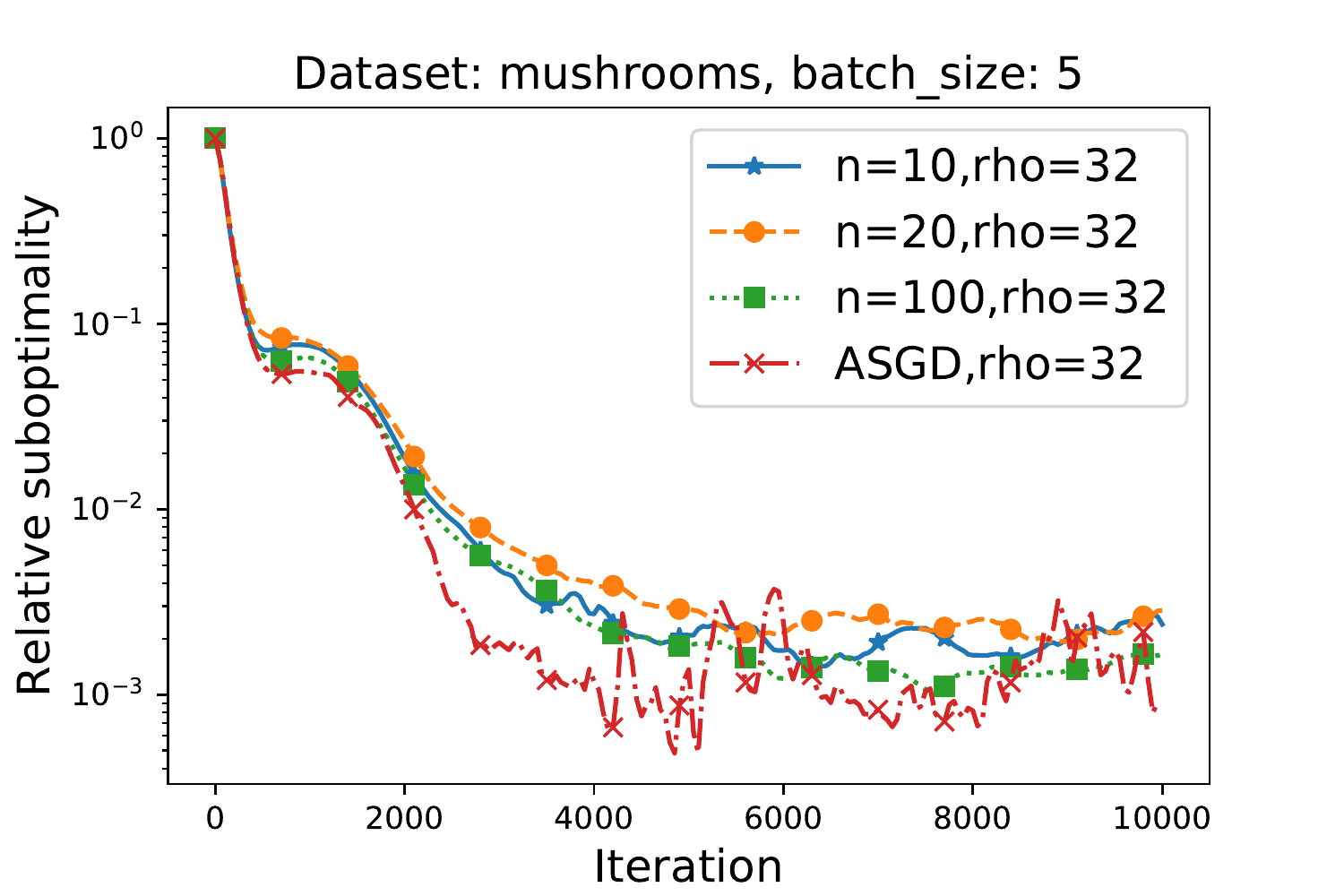}
\end{minipage}%
\begin{minipage}{0.33\textwidth}
  \centering
\includegraphics[width =  \textwidth ]{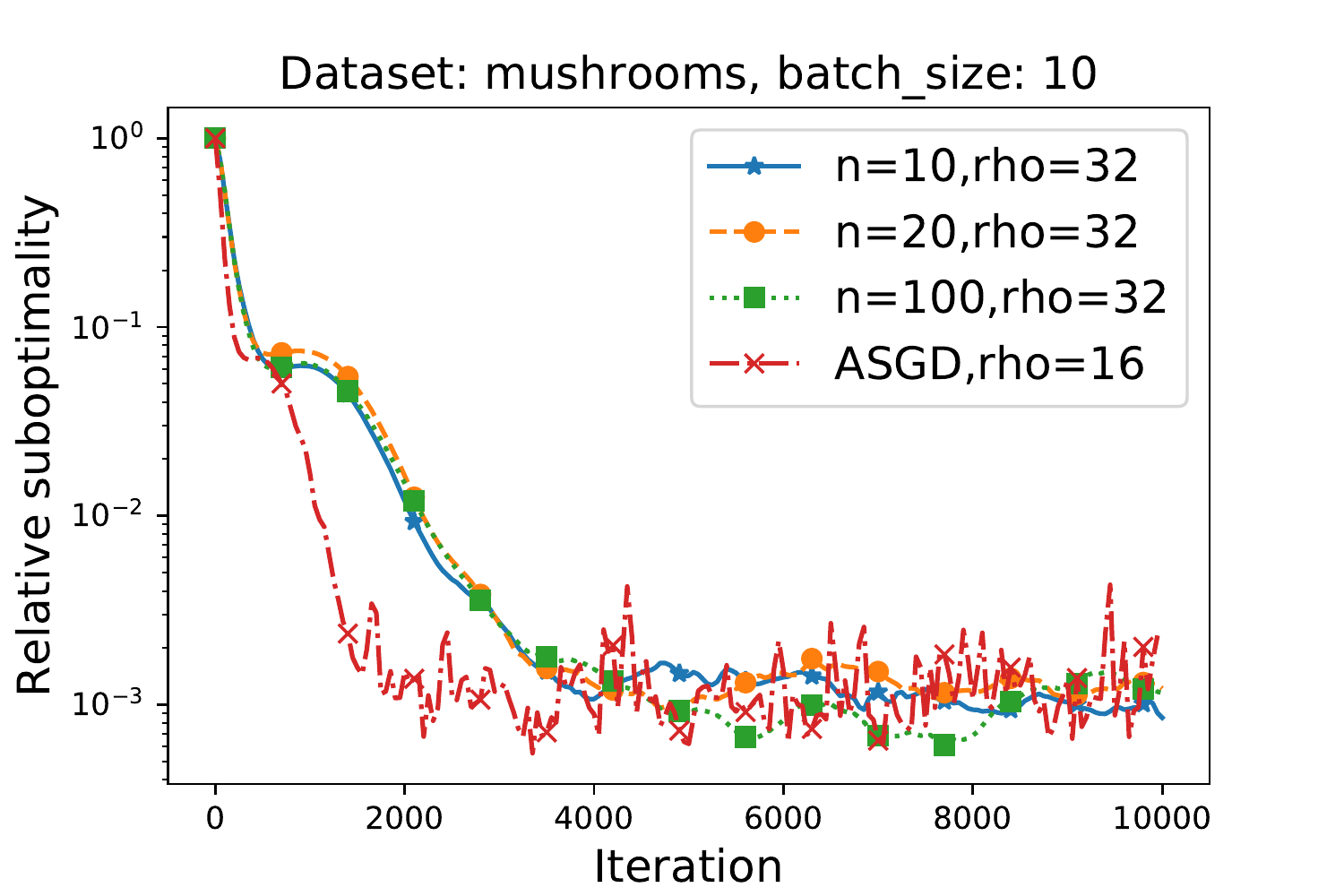}
\end{minipage}%
\\
\begin{minipage}{0.33\textwidth}
  \centering
\includegraphics[width =  \textwidth ]{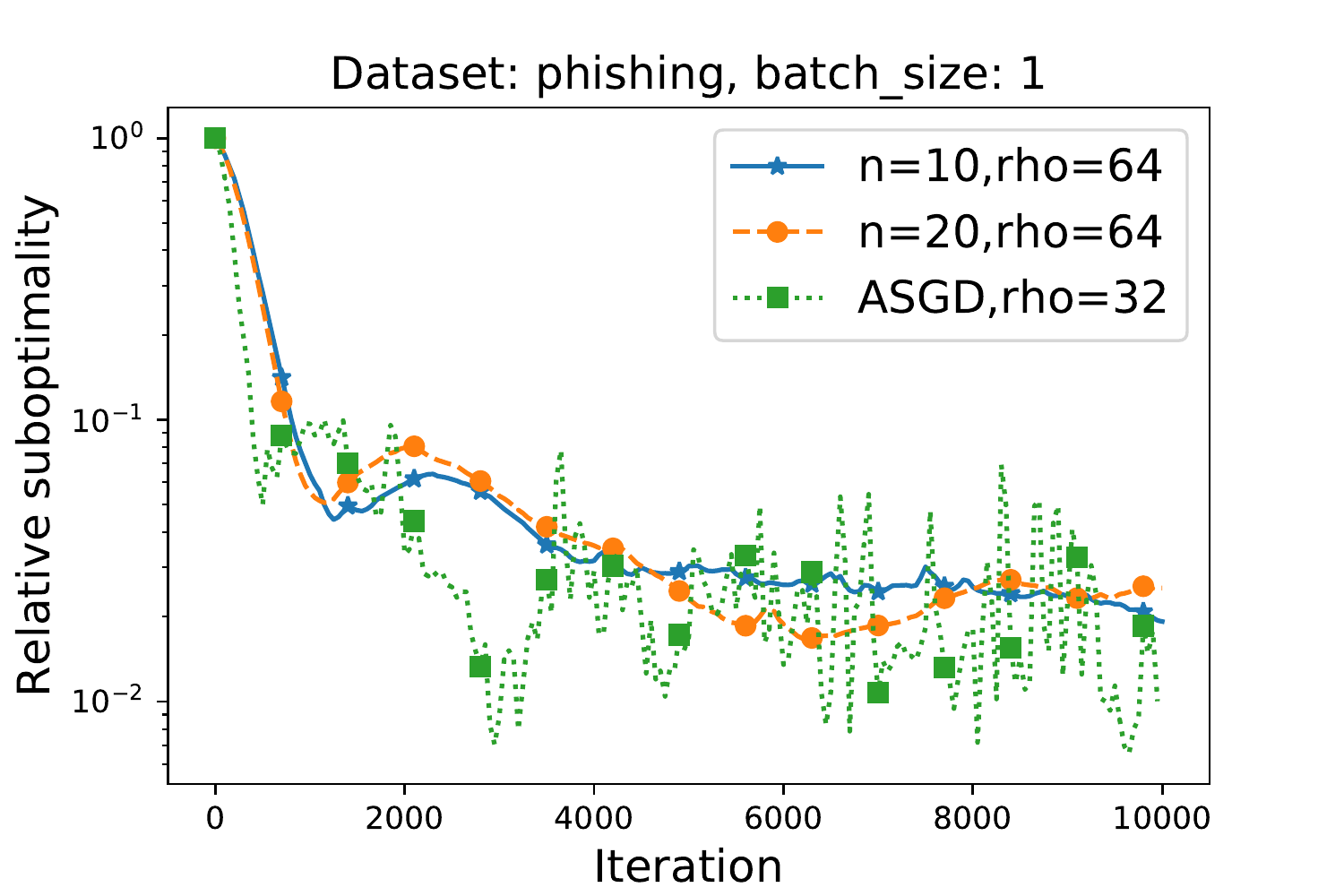}
\end{minipage}%
\begin{minipage}{0.33\textwidth}
  \centering
\includegraphics[width =  \textwidth ]{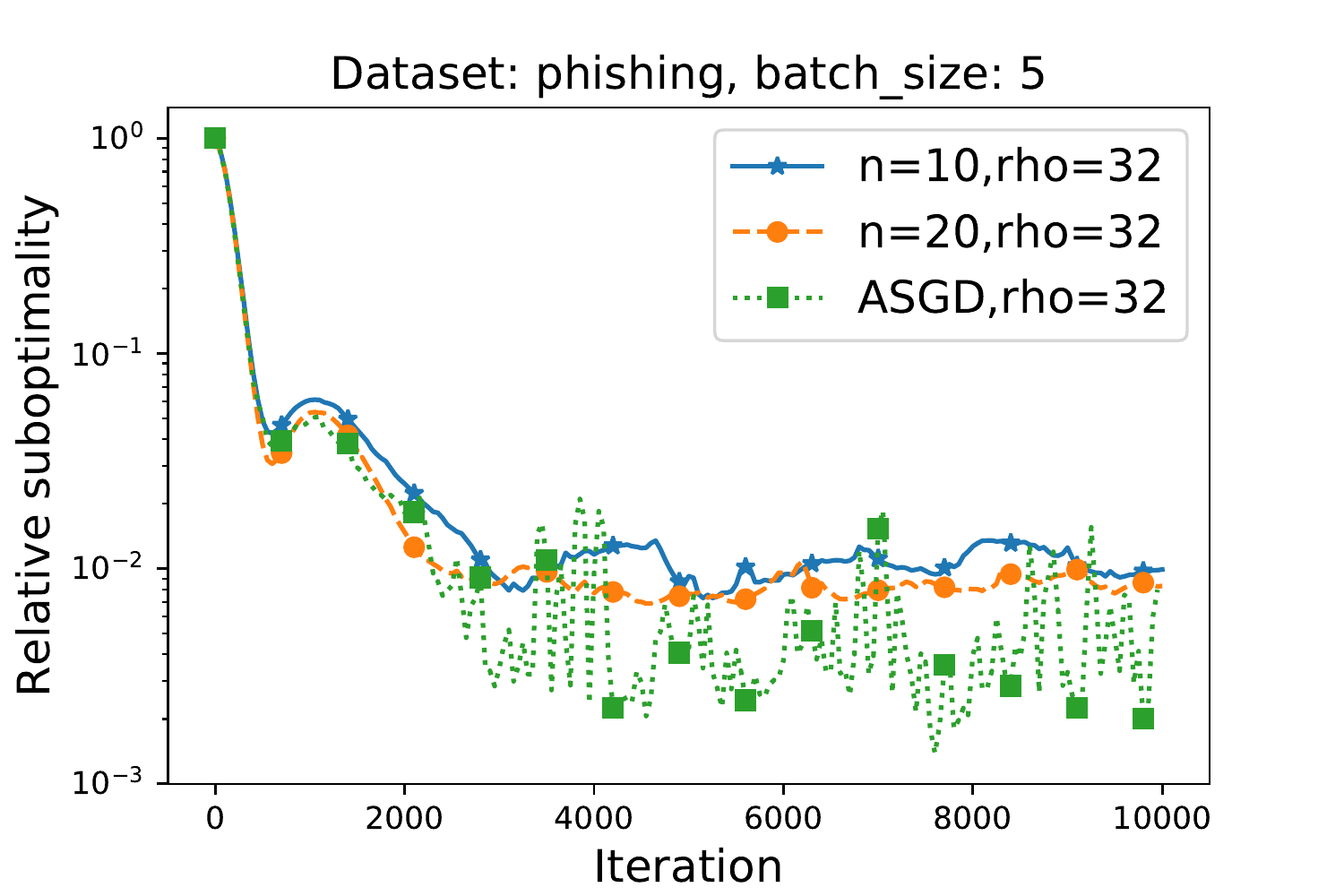}
\end{minipage}%
\begin{minipage}{0.33\textwidth}
  \centering
\includegraphics[width =  \textwidth ]{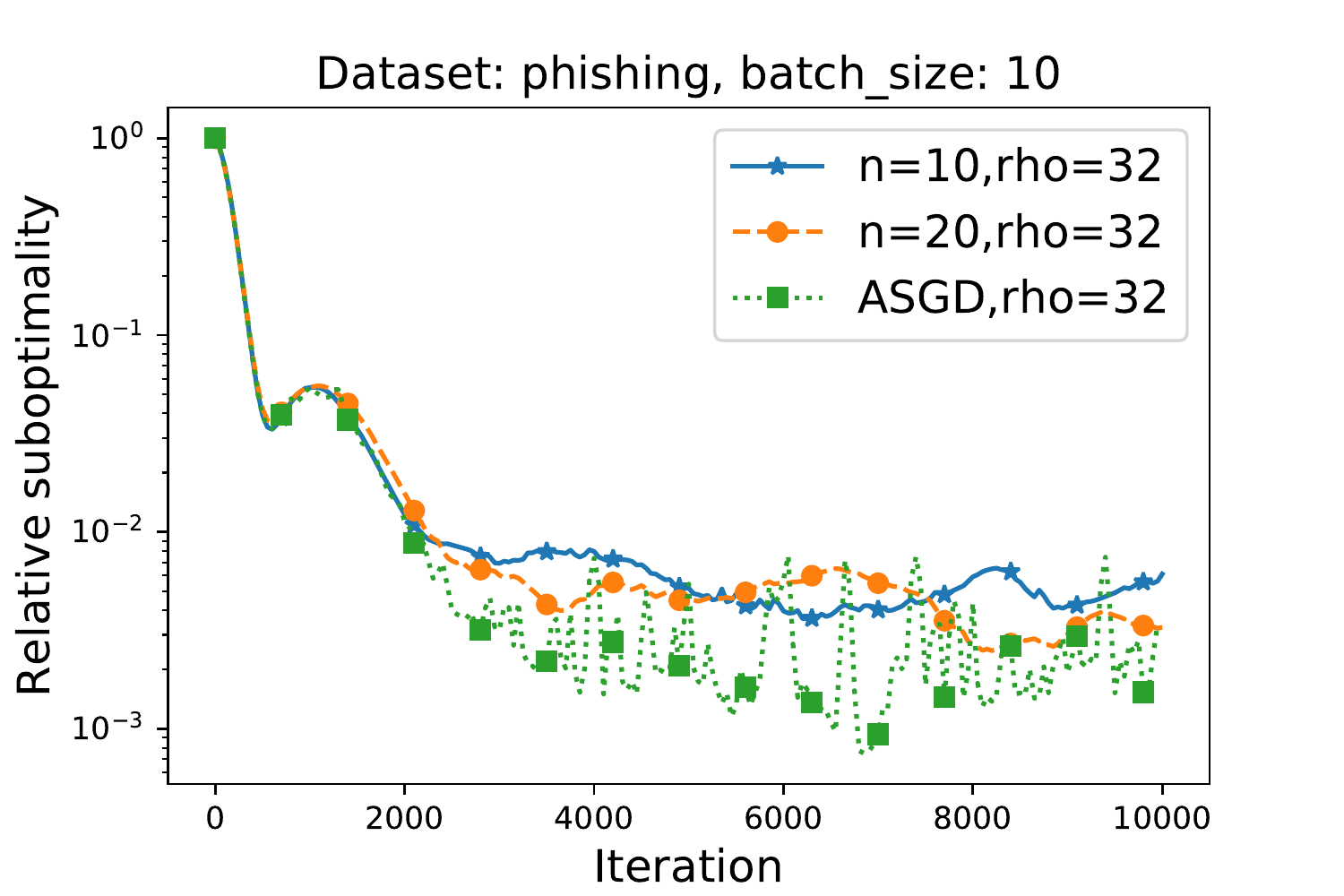}
\end{minipage}%
\\
\begin{minipage}{0.33\textwidth}
  \centering
\includegraphics[width =  \textwidth ]{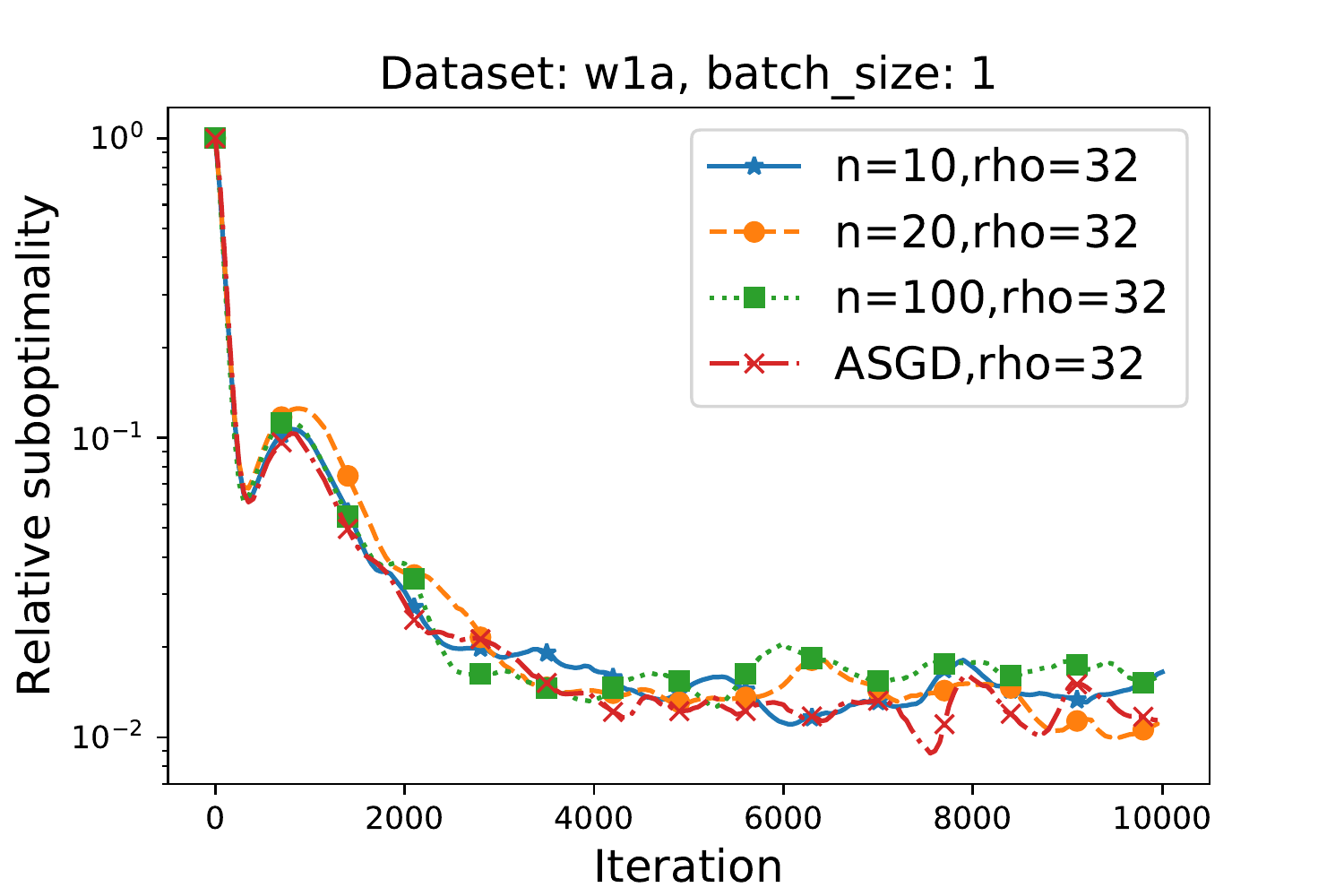}
\end{minipage}%
\begin{minipage}{0.33\textwidth}
  \centering
\includegraphics[width =  \textwidth ]{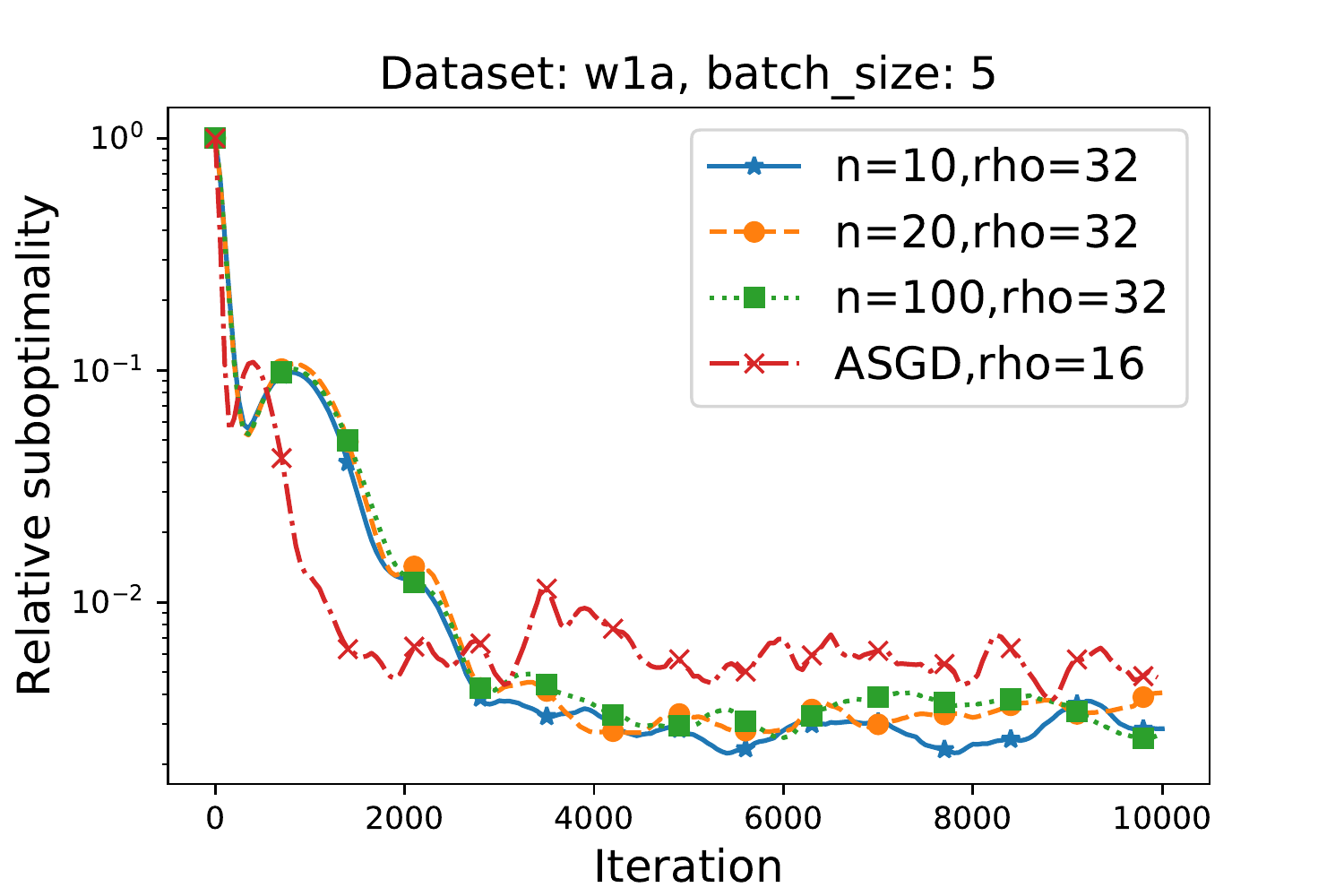}
\end{minipage}%
\begin{minipage}{0.33\textwidth}
  \centering
\includegraphics[width =  \textwidth ]{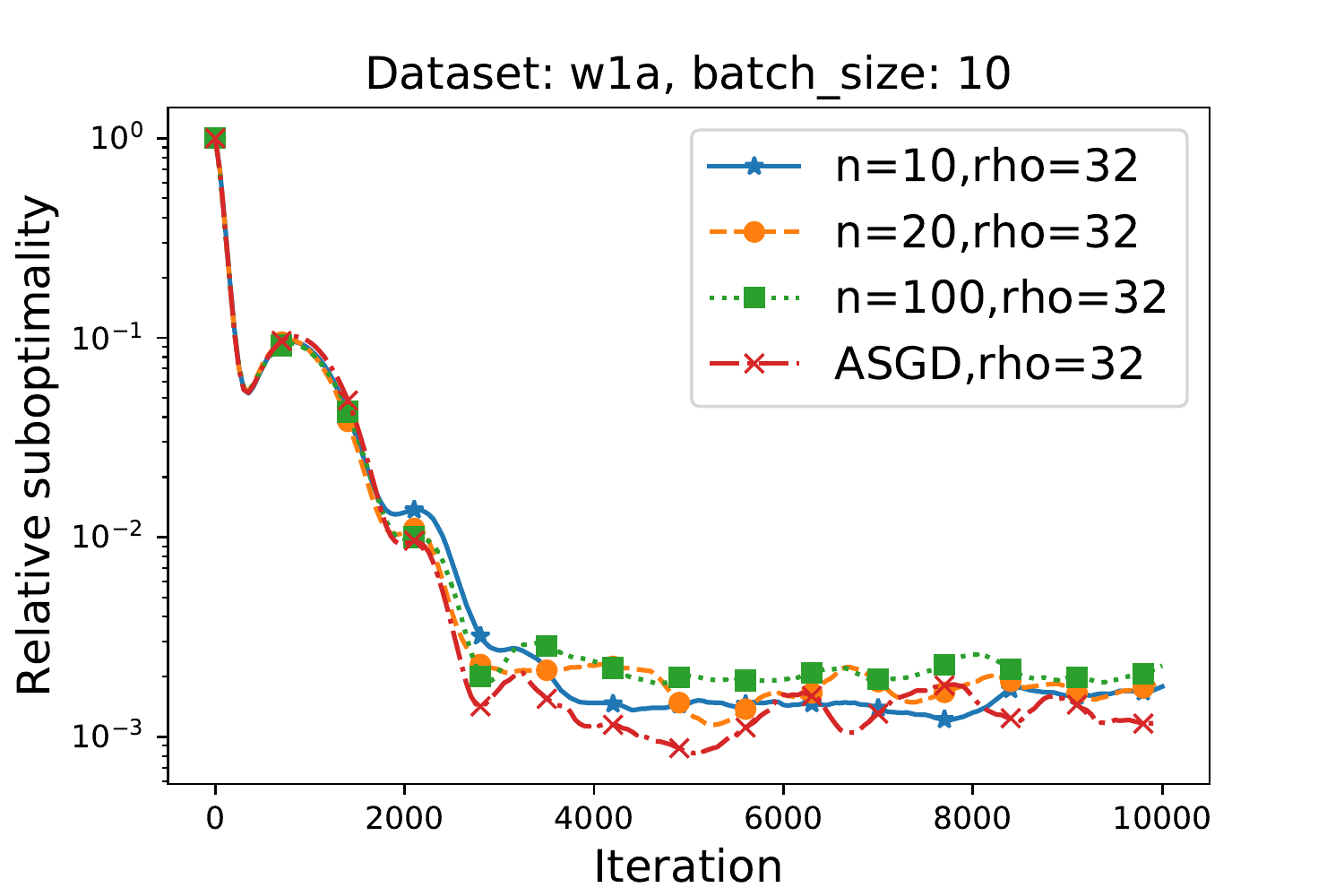}
\end{minipage}%
\\
\caption{Comparison of Algorithm~\ref{alg:acc} for various $(n,\tau)$ such that $n\tau=1$. Label ``ASGD'' corresponds to the choice $n=1, \tau = 1$. Label ``batch\_size'' indicates how big minibatch was chosen for stochastic gradient of each worker's objective. Parameter $\rho$ was chosen by grid search.  } \label{fig:99_acc1}
\end{figure}

Now, we once again check how different values of $\tau$ affect the convergence speed for several values of $n$. Figure~\ref{fig:99_acc2} presents the results. In every case, $\tau$ slightly influences the convergence rate (or the region where the iterates oscillate), although the effect is weaker for larger $n$. Note that theory predicts diminishing effect of $\tau$ only above $\frac{\bar{\rho}\tilde{\rho}}{n}$, in contrast to other sections, where the limit is $\frac1n$.

\begin{figure}[H]
\centering
\begin{minipage}{0.33\textwidth}
  \centering
\includegraphics[width =  \textwidth ]{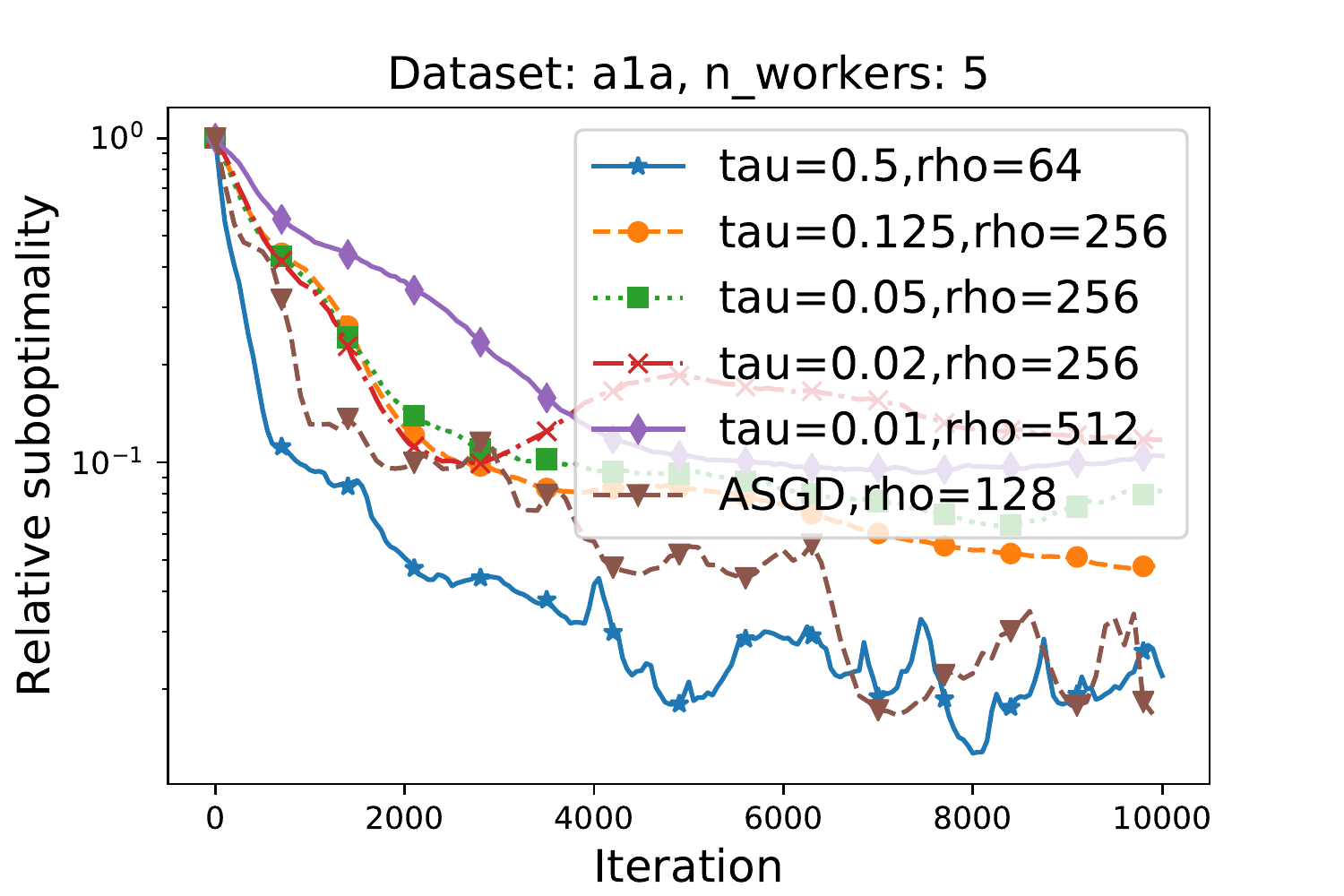}
\end{minipage}%
\begin{minipage}{0.33\textwidth}
  \centering
\includegraphics[width =  \textwidth ]{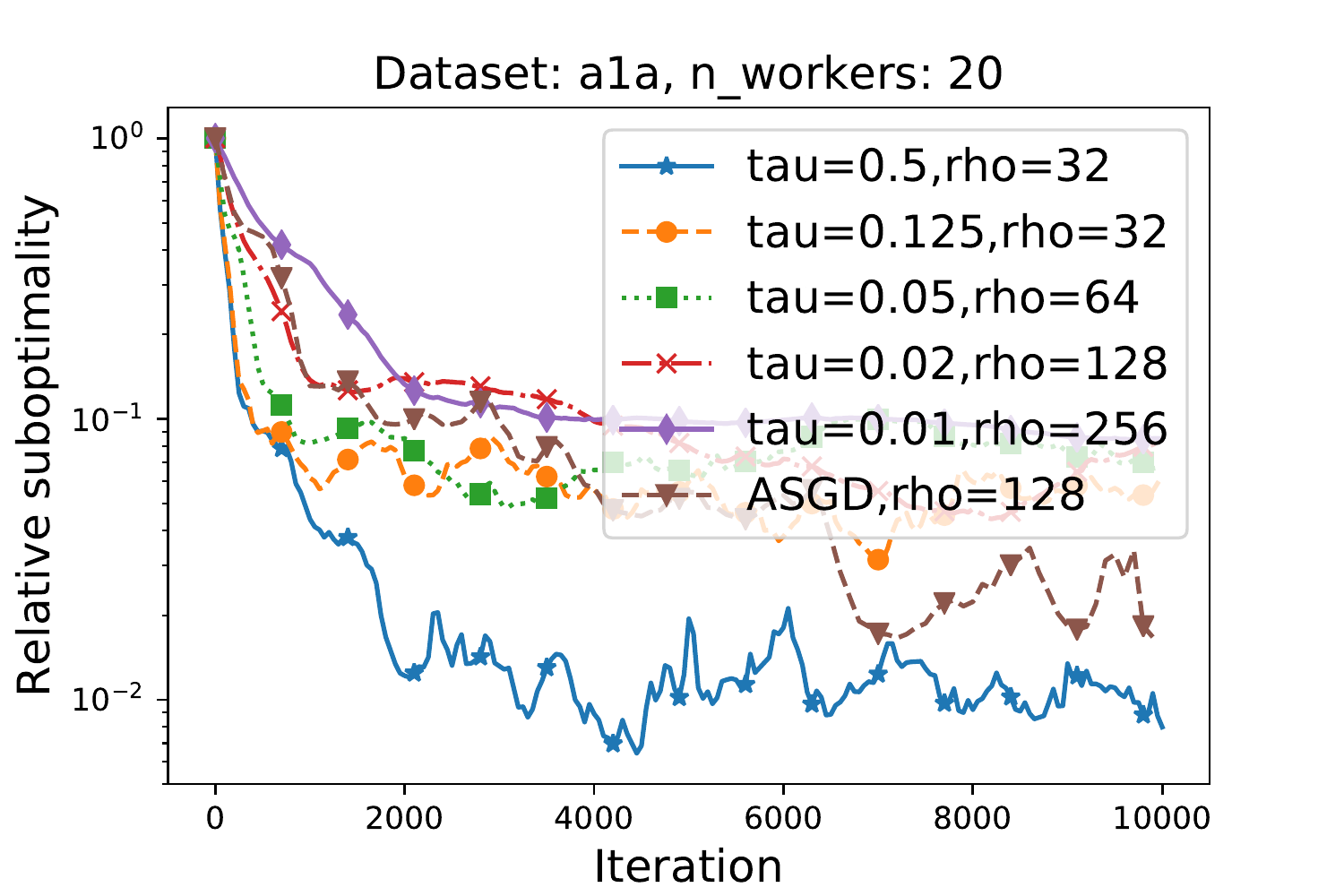}
\end{minipage}%
\begin{minipage}{0.33\textwidth}
  \centering
\includegraphics[width =  \textwidth ]{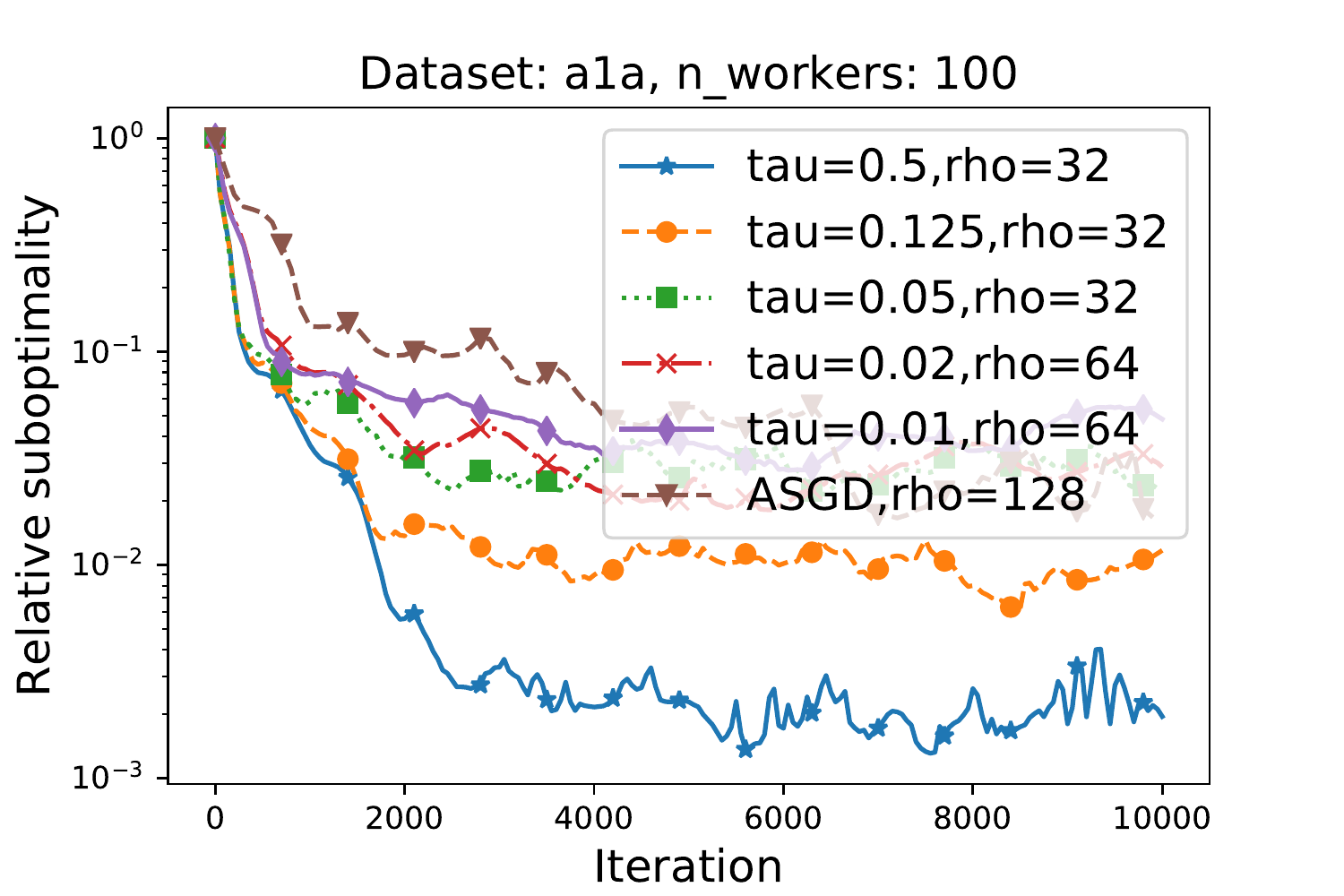}
\end{minipage}%
\\
\begin{minipage}{0.33\textwidth}
  \centering
\includegraphics[width =  \textwidth ]{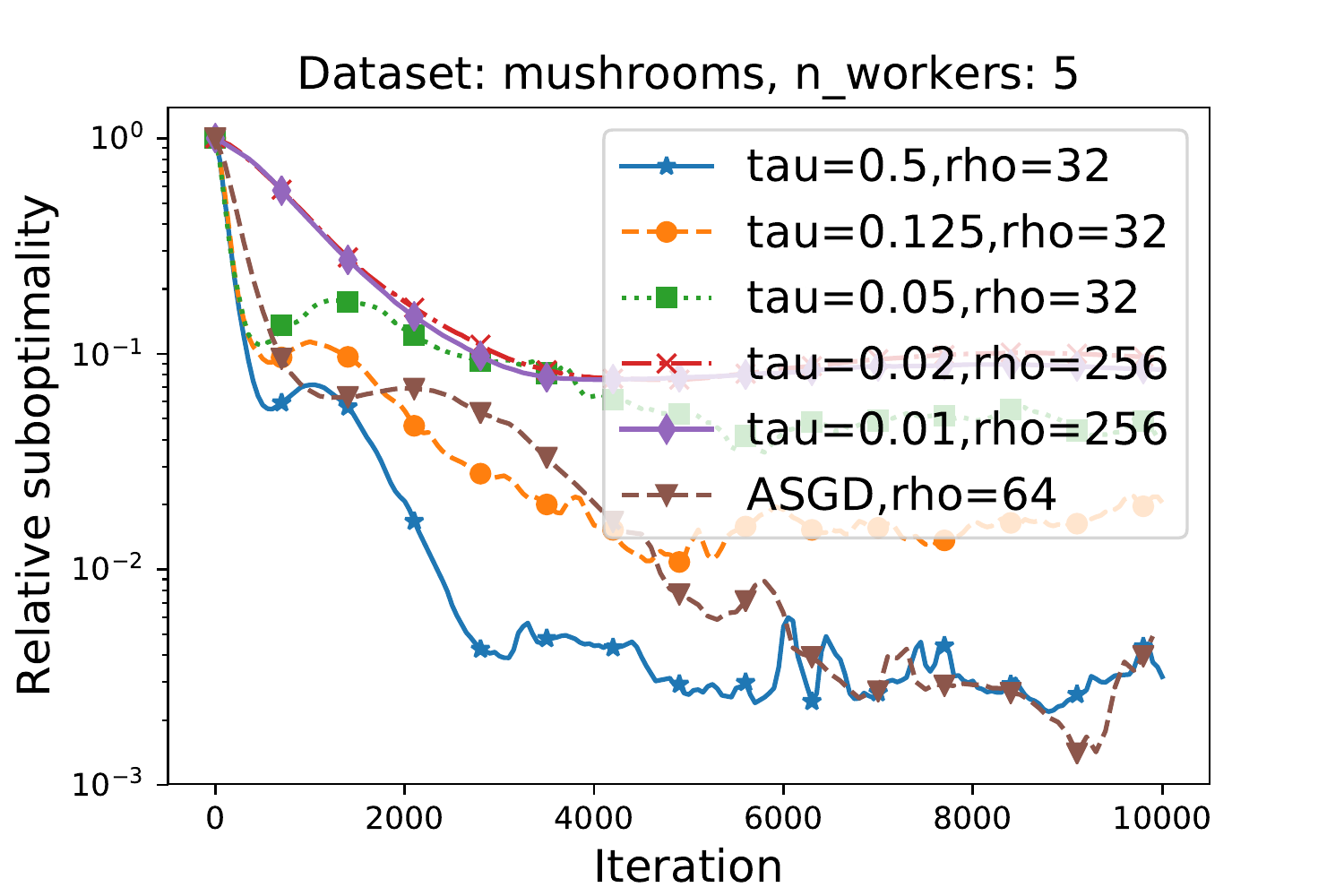}
\end{minipage}%
\begin{minipage}{0.33\textwidth}
  \centering
\includegraphics[width =  \textwidth ]{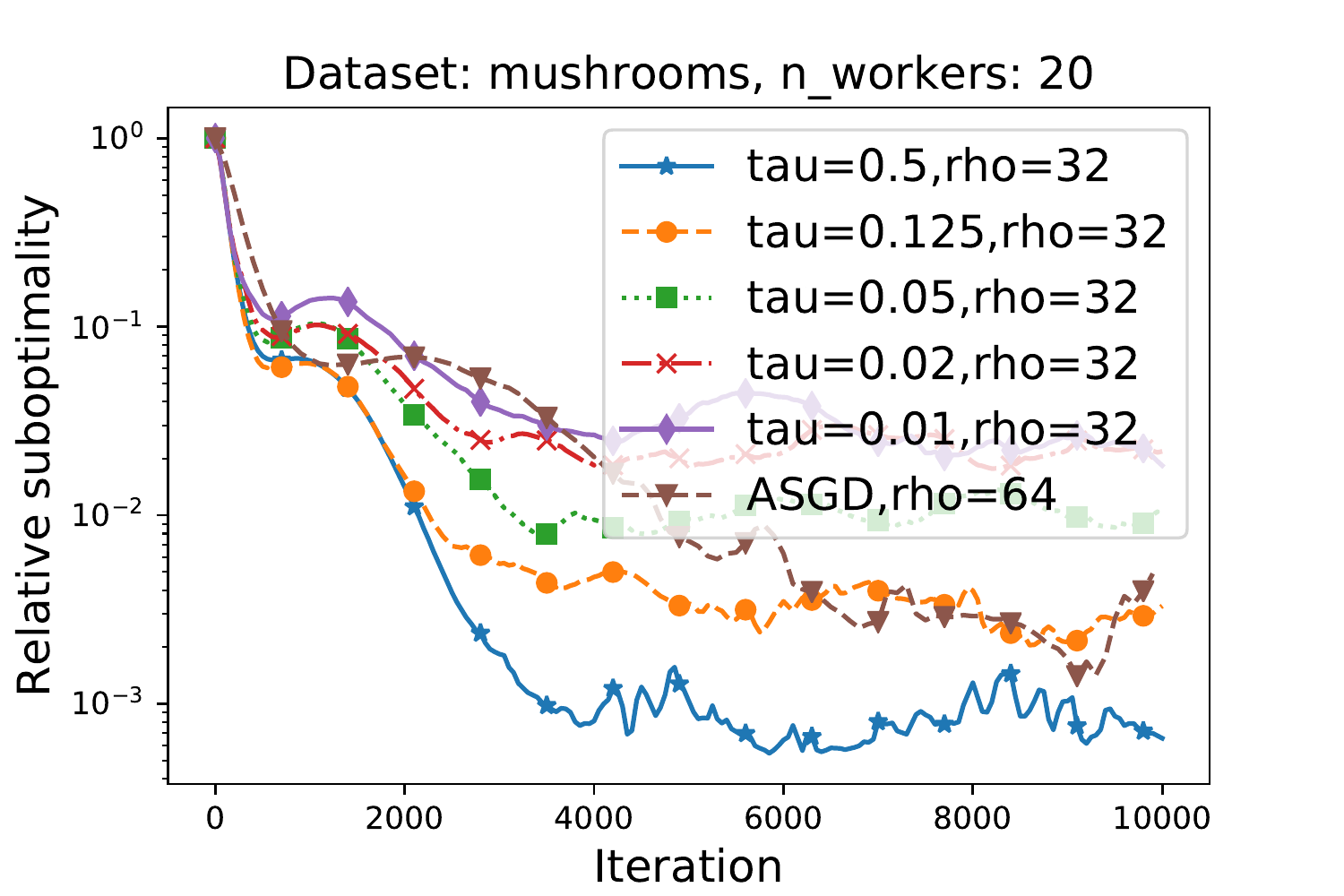}
\end{minipage}%
\begin{minipage}{0.33\textwidth}
  \centering
\includegraphics[width =  \textwidth ]{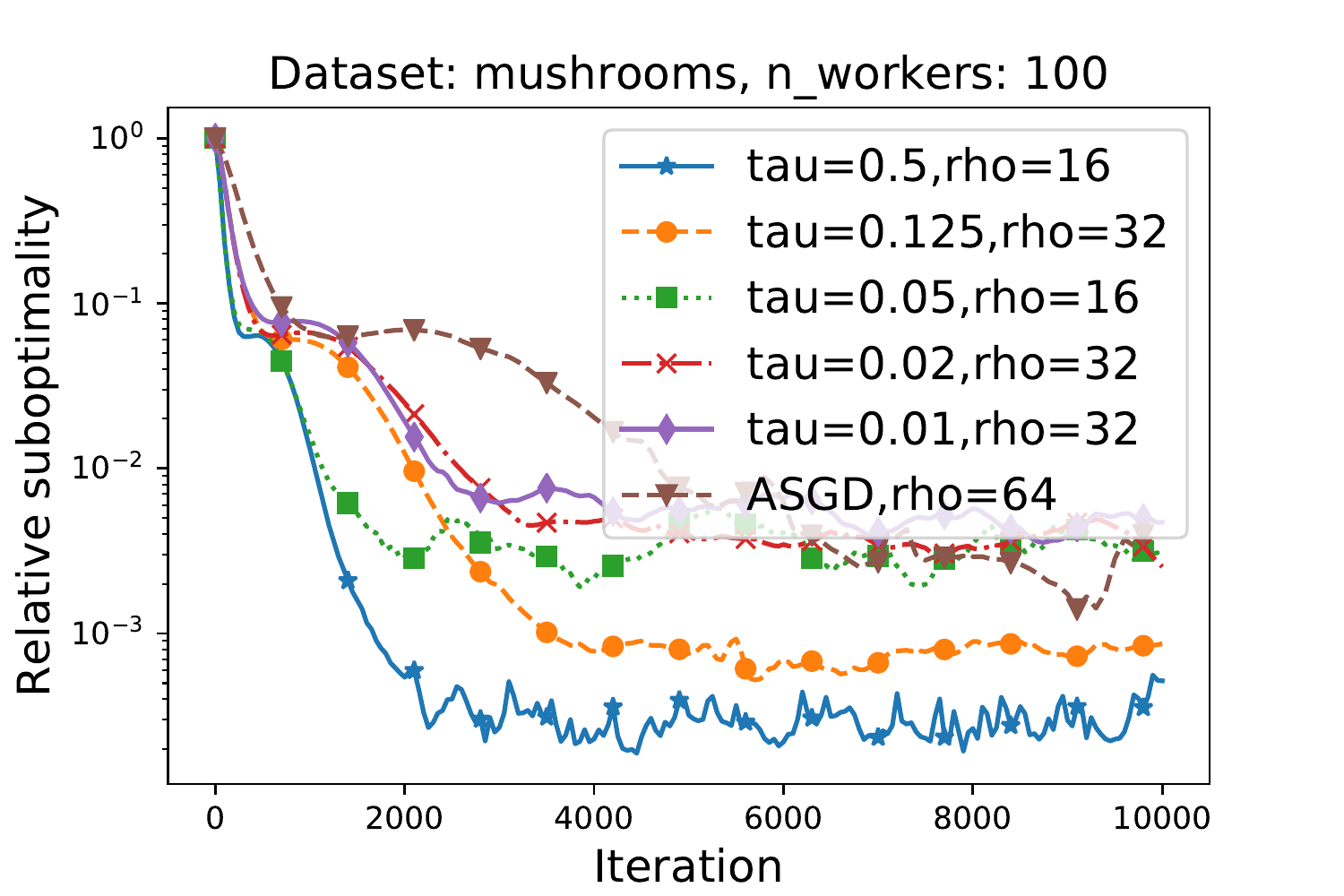}
\end{minipage}%
\\
\begin{minipage}{0.33\textwidth}
  \centering
\includegraphics[width =  \textwidth ]{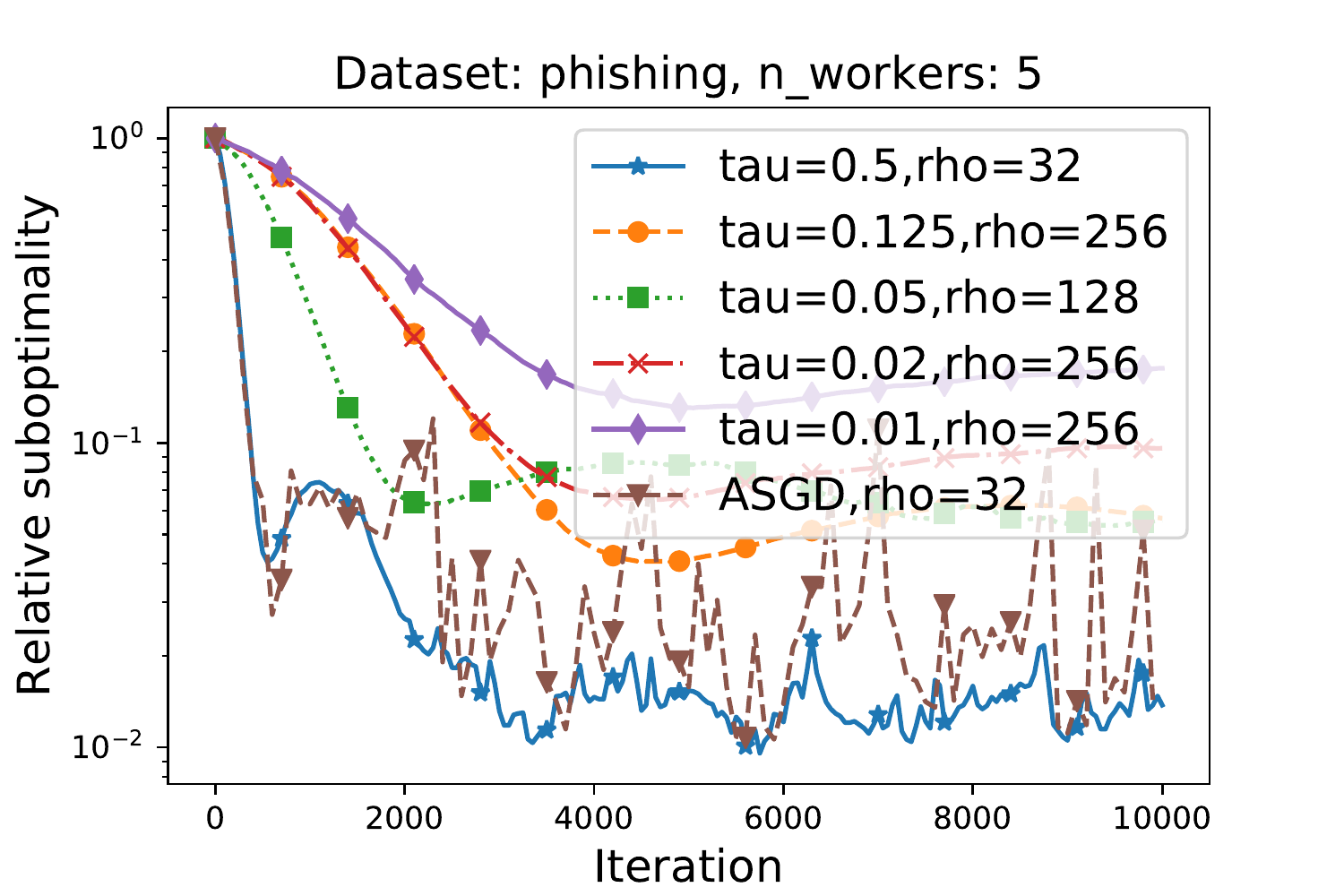}
\end{minipage}%
\begin{minipage}{0.33\textwidth}
  \centering
\includegraphics[width =  \textwidth ]{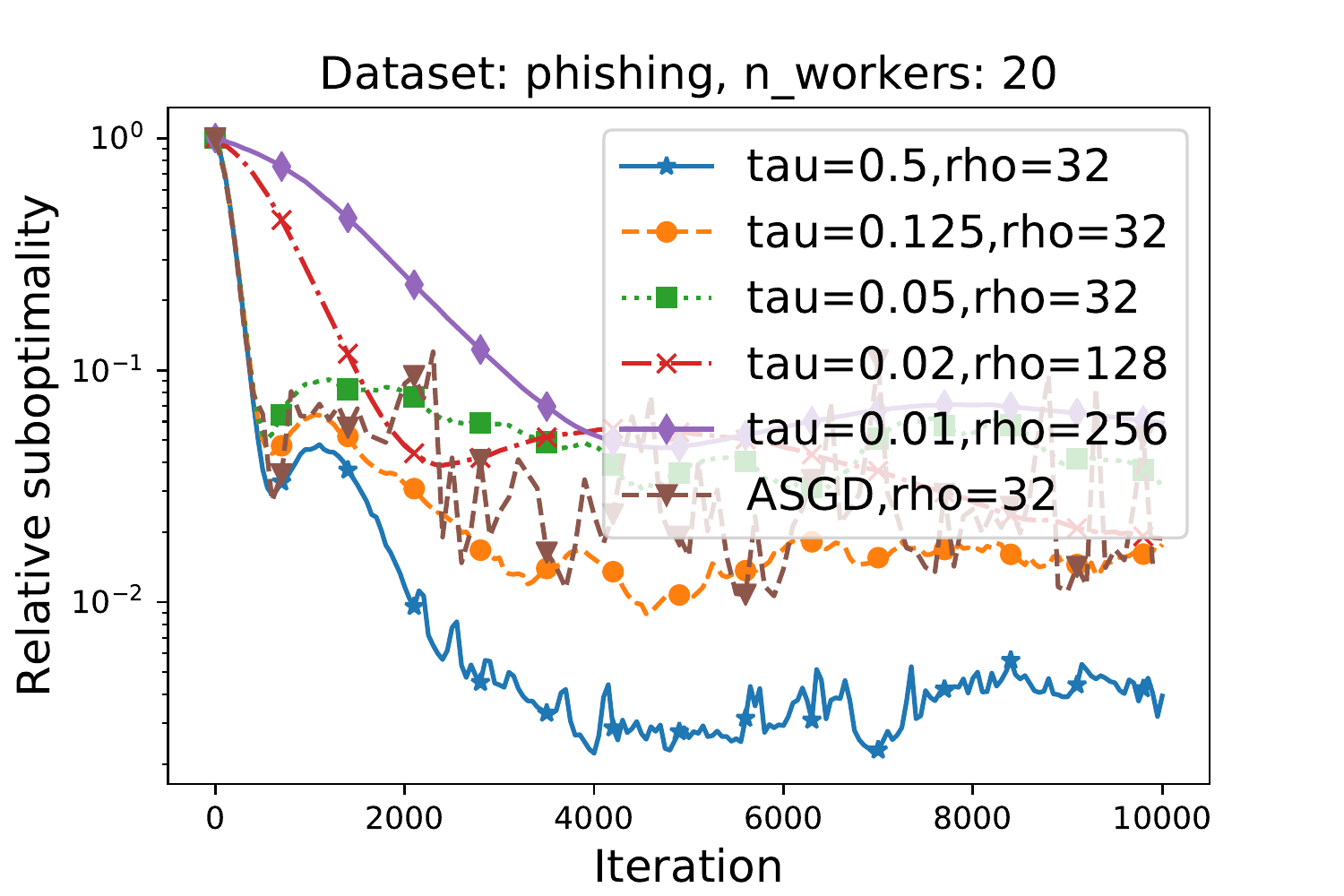}
\end{minipage}%
\\
\begin{minipage}{0.33\textwidth}
  \centering
\includegraphics[width =  \textwidth ]{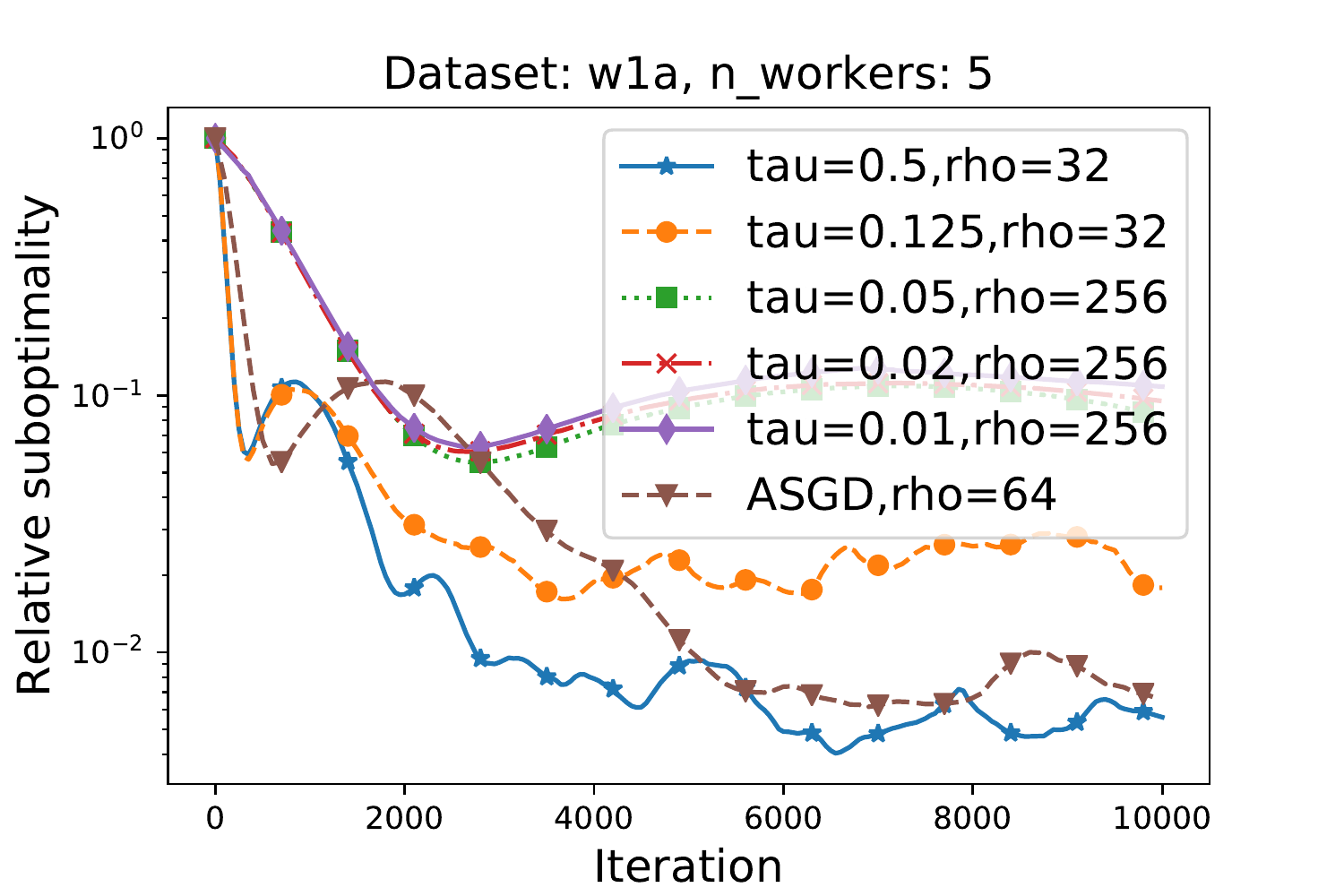}
\end{minipage}%
\begin{minipage}{0.33\textwidth}
  \centering
\includegraphics[width =  \textwidth ]{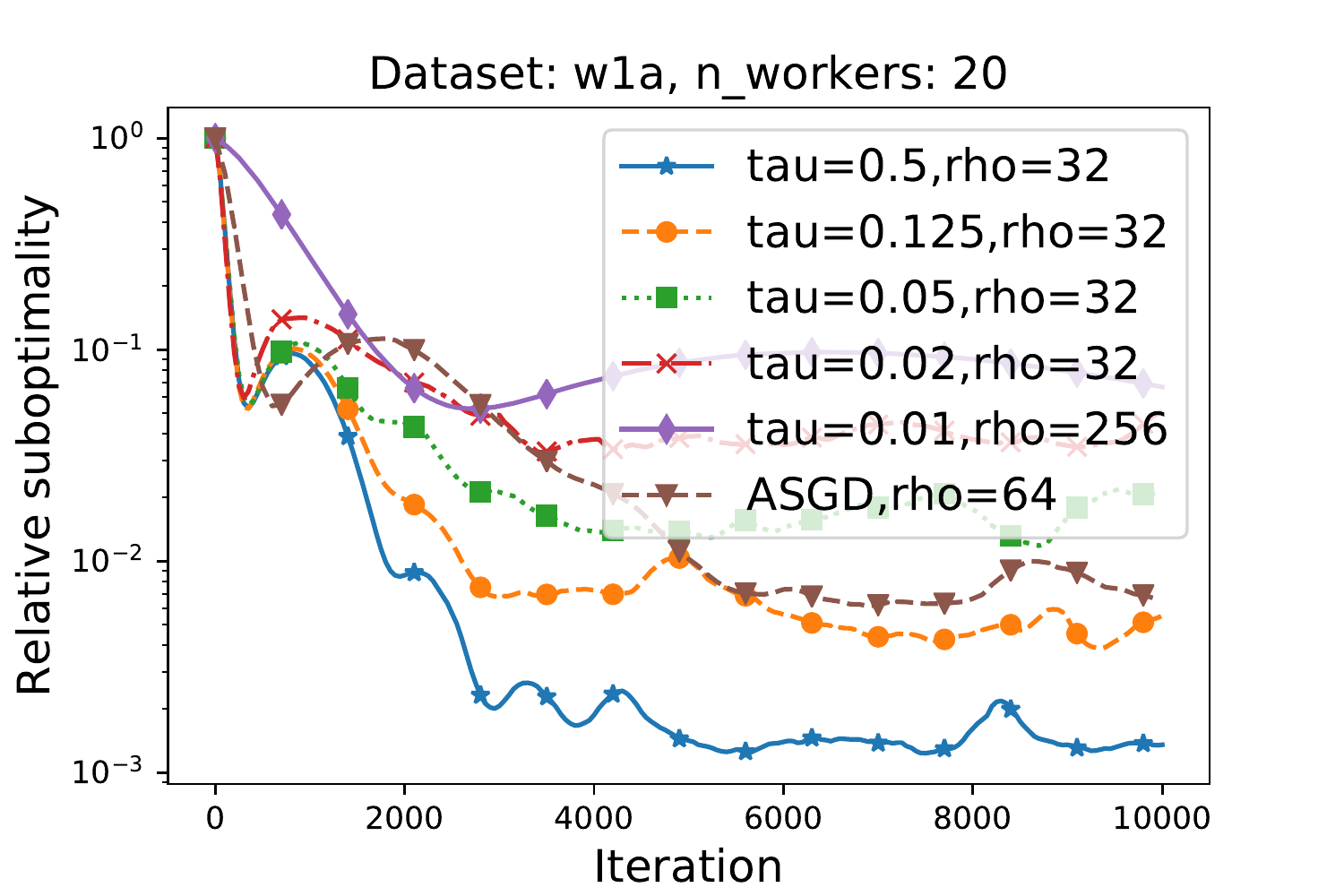}
\end{minipage}%
\begin{minipage}{0.33\textwidth}
  \centering
\includegraphics[width =  \textwidth ]{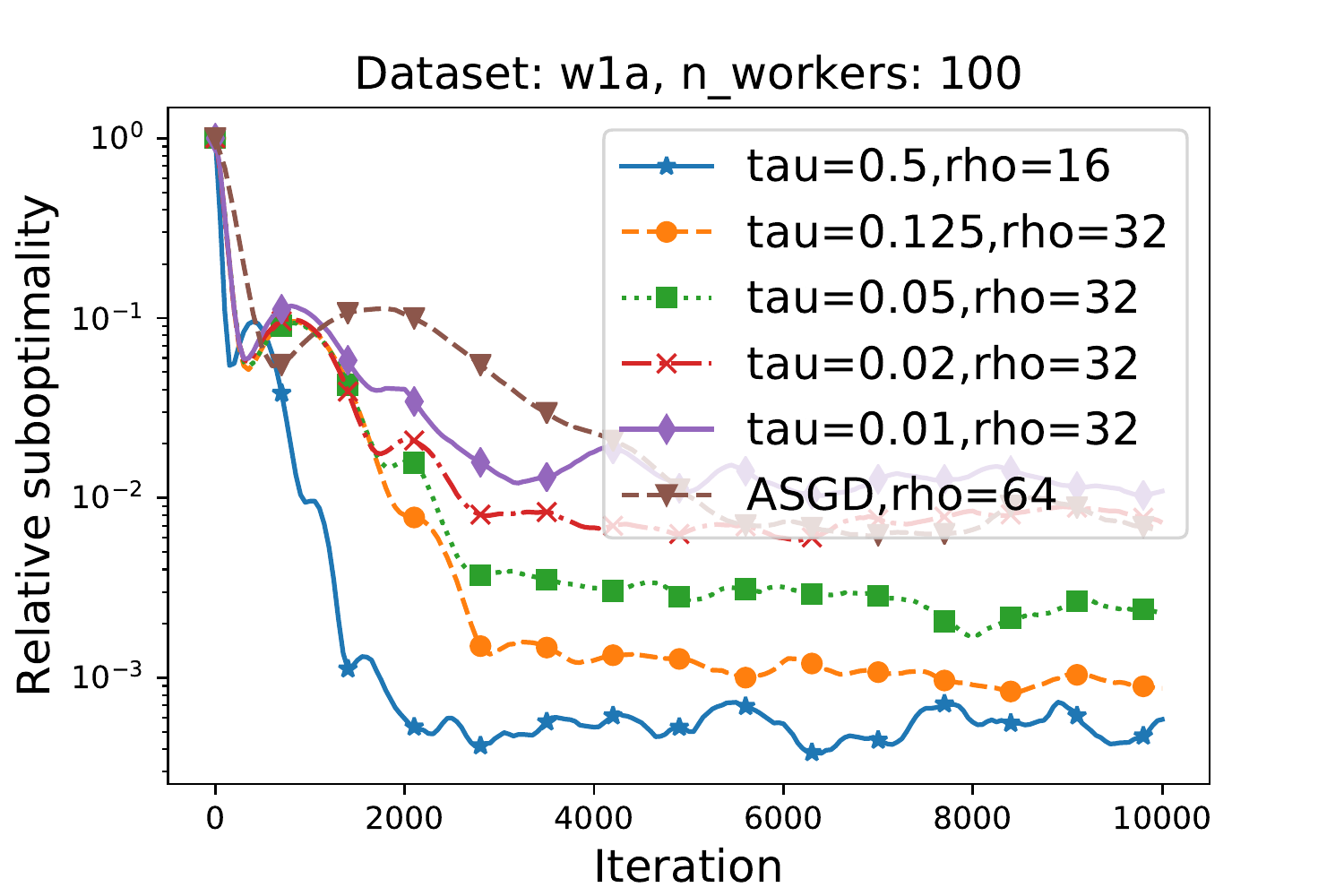}
\end{minipage}%
\\
\caption{Behavior of Algorithm~\ref{alg:acc} while varying $\tau$. Label ``ASGD'' corresponds to the choice $n=1, \tau = 1$. Parameter $\rho$ was chosen by grid search.  } \label{fig:99_acc2}
\end{figure}

\subsection{{\tt ISAGA} \label{sec:99_exp_saga}}
In the next experiment, we compare {\tt SAGA} against {\tt ISAGA} in a shared data setup (Algorithm~\ref{alg:saga}) for various values of $n$ with $\tau = \frac{1}{n}$ in order to demonstrate linear scaling. We consider logistic regression problem on LibSVM data~\cite{chang2011libsvm}. The results (Figure~\ref{fig:99_saga1}) corroborate our theory: indeed, setting $n\tau=1$ does not lead to a decrease in the convergence rate when compared to the original {\tt SAGA}.

\begin{figure}[H]
\centering
\begin{minipage}{0.35\textwidth}
  \centering
\includegraphics[width =  \textwidth ]{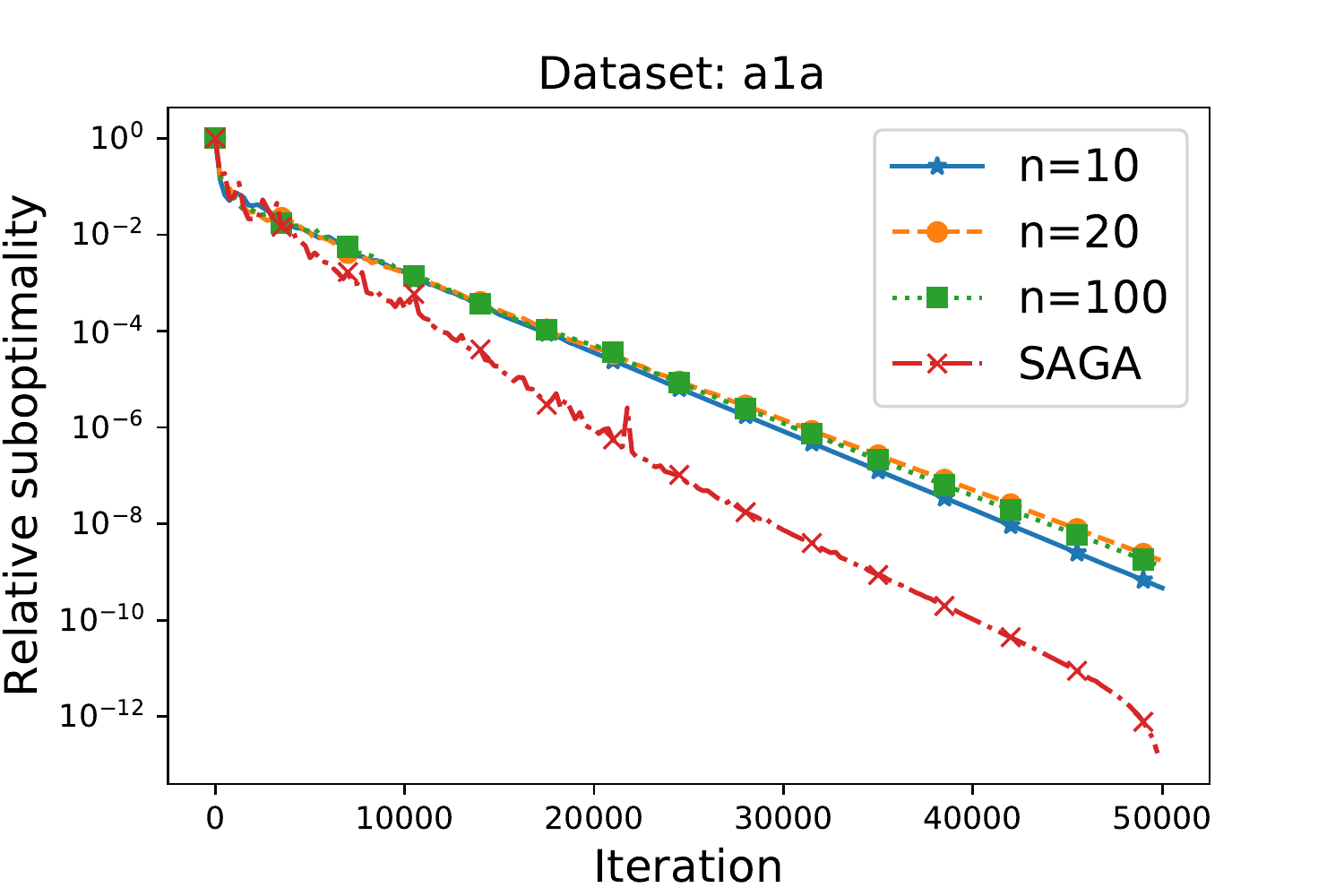}
\end{minipage}%
\begin{minipage}{0.35\textwidth}
  \centering
\includegraphics[width =  \textwidth ]{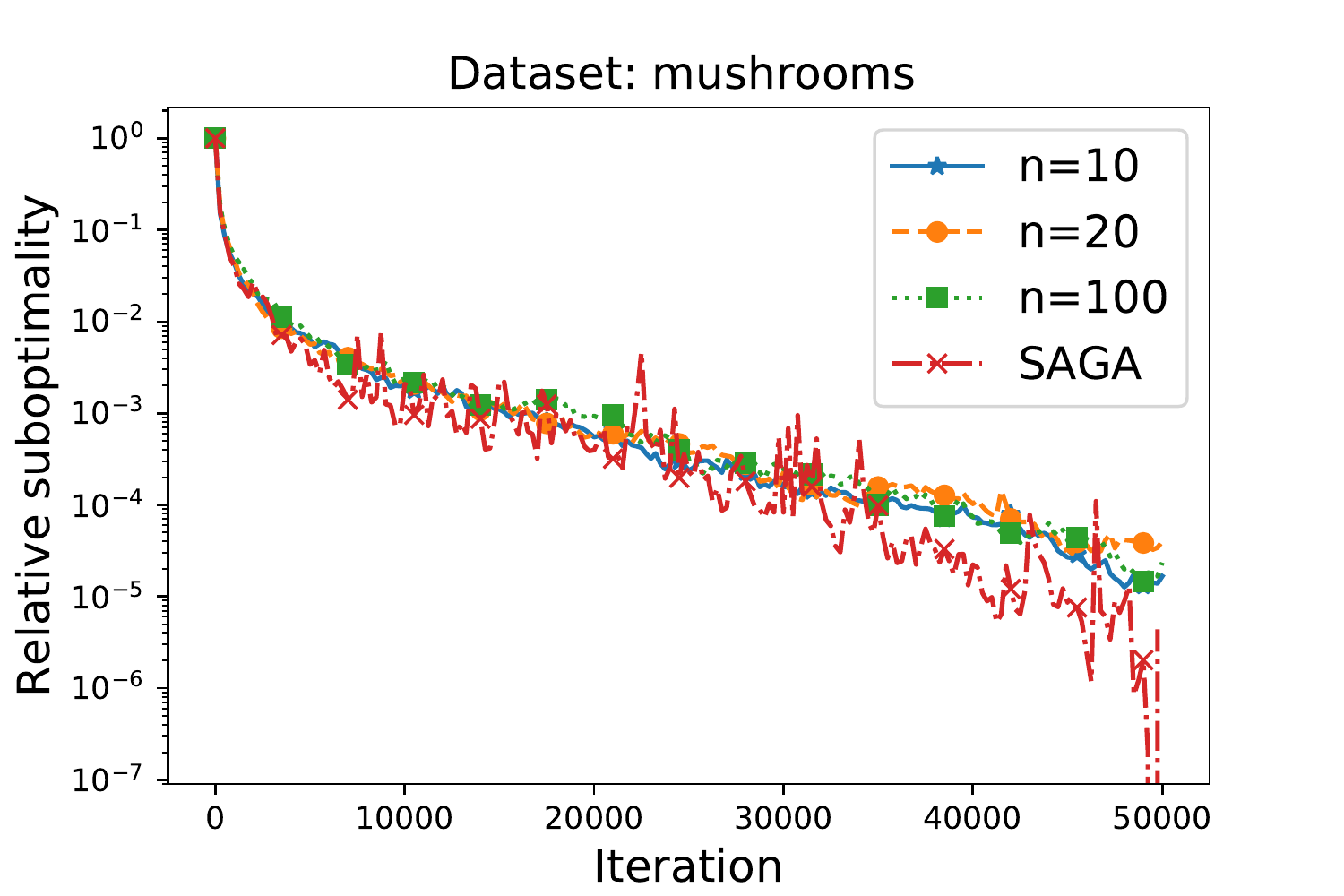}
\end{minipage}%
\\
\begin{minipage}{0.35\textwidth}
  \centering
\includegraphics[width =  \textwidth ]{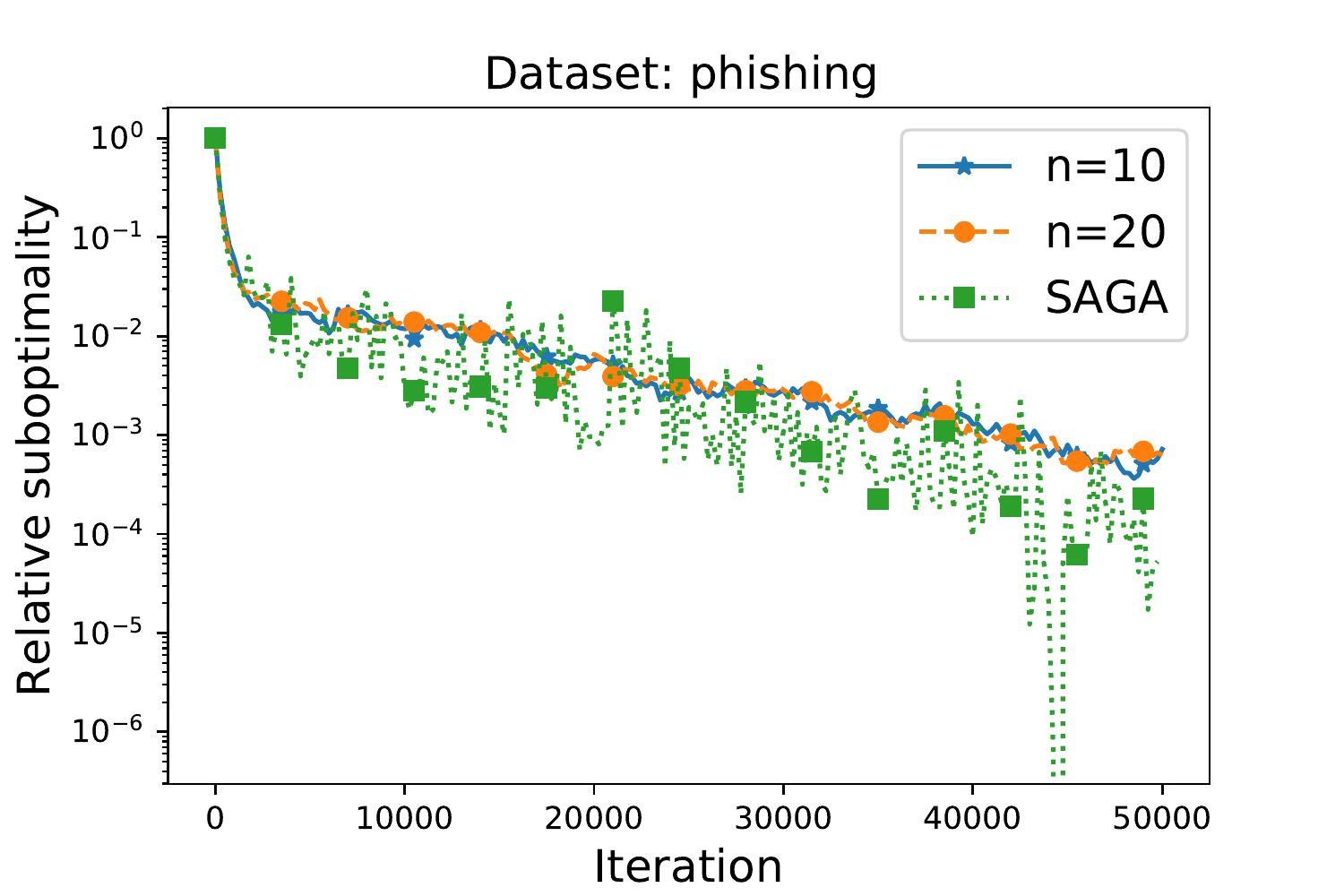}
\end{minipage}%
\begin{minipage}{0.35\textwidth}
  \centering
\includegraphics[width =  \textwidth ]{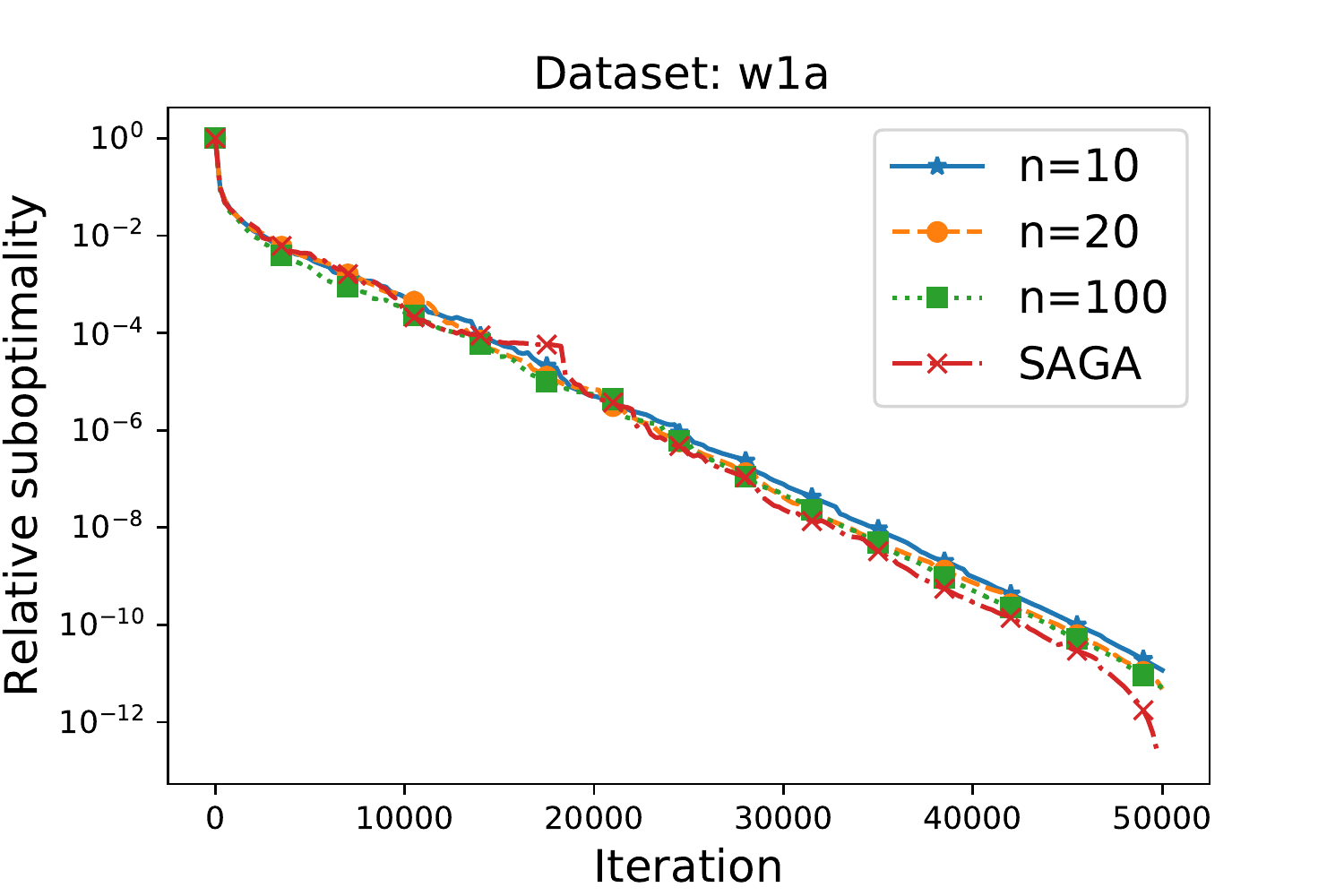}
\end{minipage}%
\\
\caption{Comparison of {\tt SAGA} and Algorithm~\ref{alg:saga} for various values $n$ and $\tau=n^{-1}$. Stepsize $\alpha = \frac{1}{L(3n^{-1}+\tau)}$ is chosen in each case. } \label{fig:99_saga1}
\end{figure}

The second experiment of this section shows the convergence behavior for varying $\tau$ of Algorithm~\ref{alg:saga}. The results (Figure~\ref{fig:99_saga2}) show that, for small $n$, the ratio of coordinates $\tau$ affects the speed heavily. However, as $n$ increases, the effect of $\tau$ is diminishing. 

\begin{figure}[H]
\centering
\begin{minipage}{0.33\textwidth}
  \centering
\includegraphics[width =  \textwidth ]{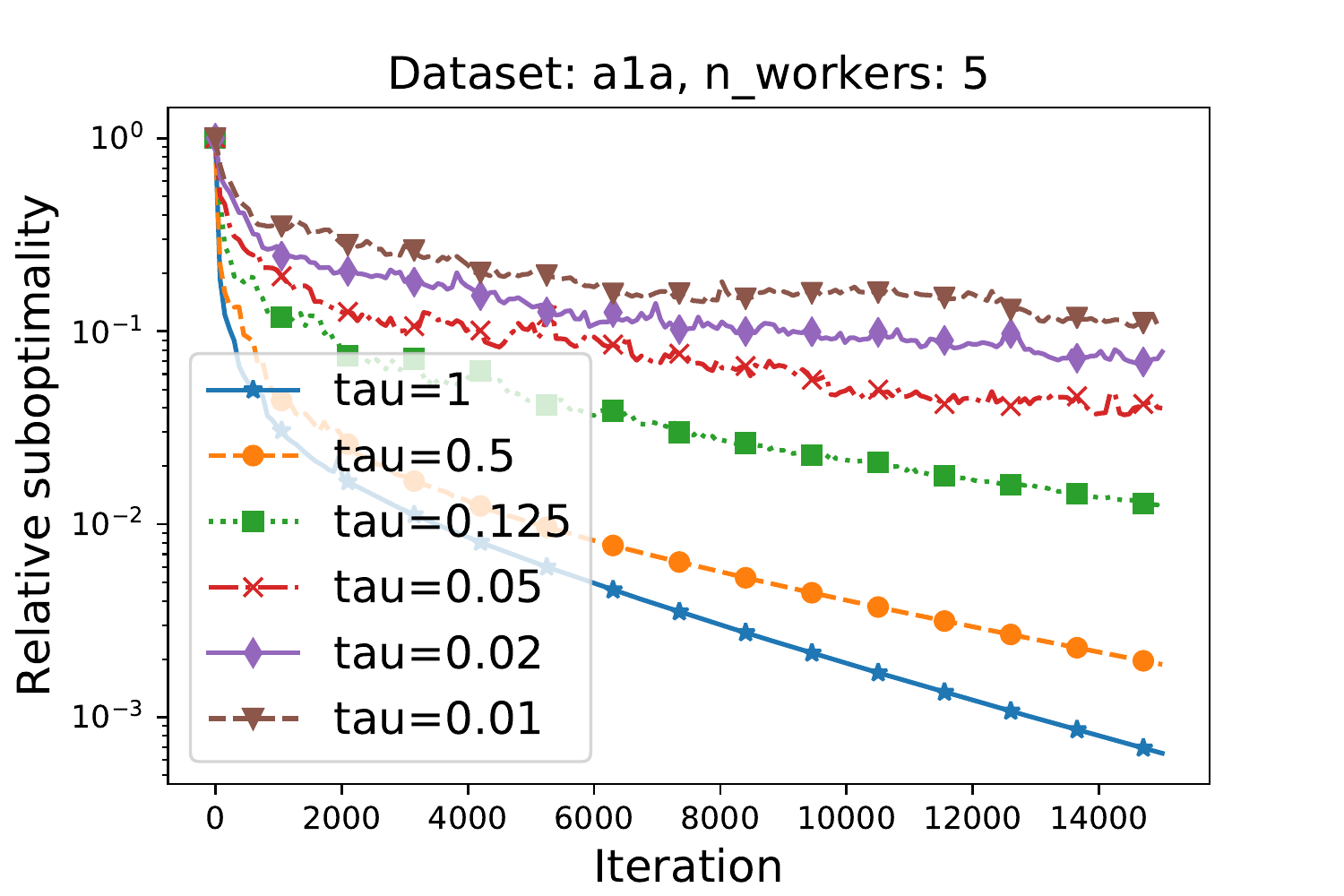}
\end{minipage}%
\begin{minipage}{0.33\textwidth}
  \centering
\includegraphics[width =  \textwidth ]{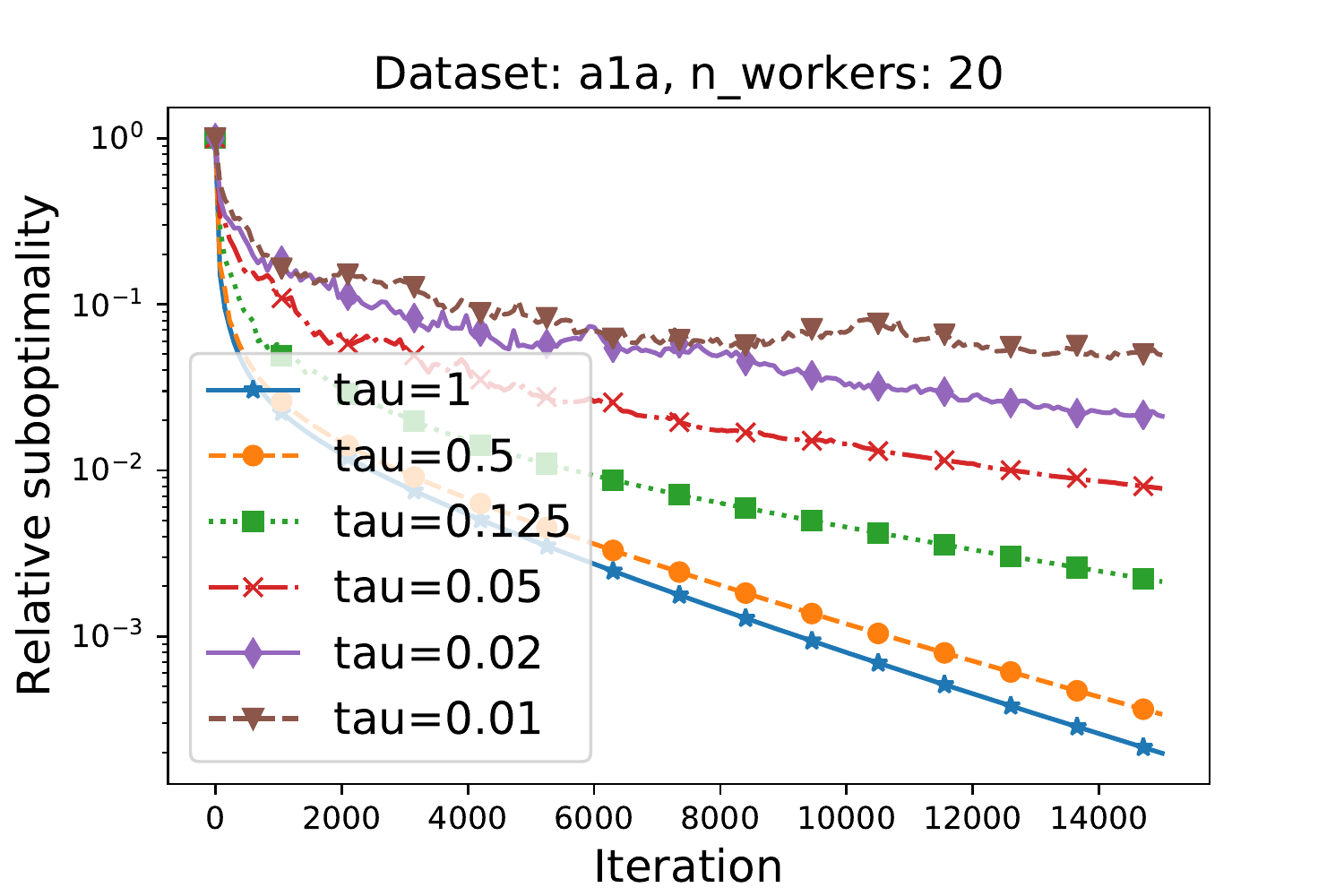}
\end{minipage}%
\begin{minipage}{0.33\textwidth}
  \centering
\includegraphics[width =  \textwidth ]{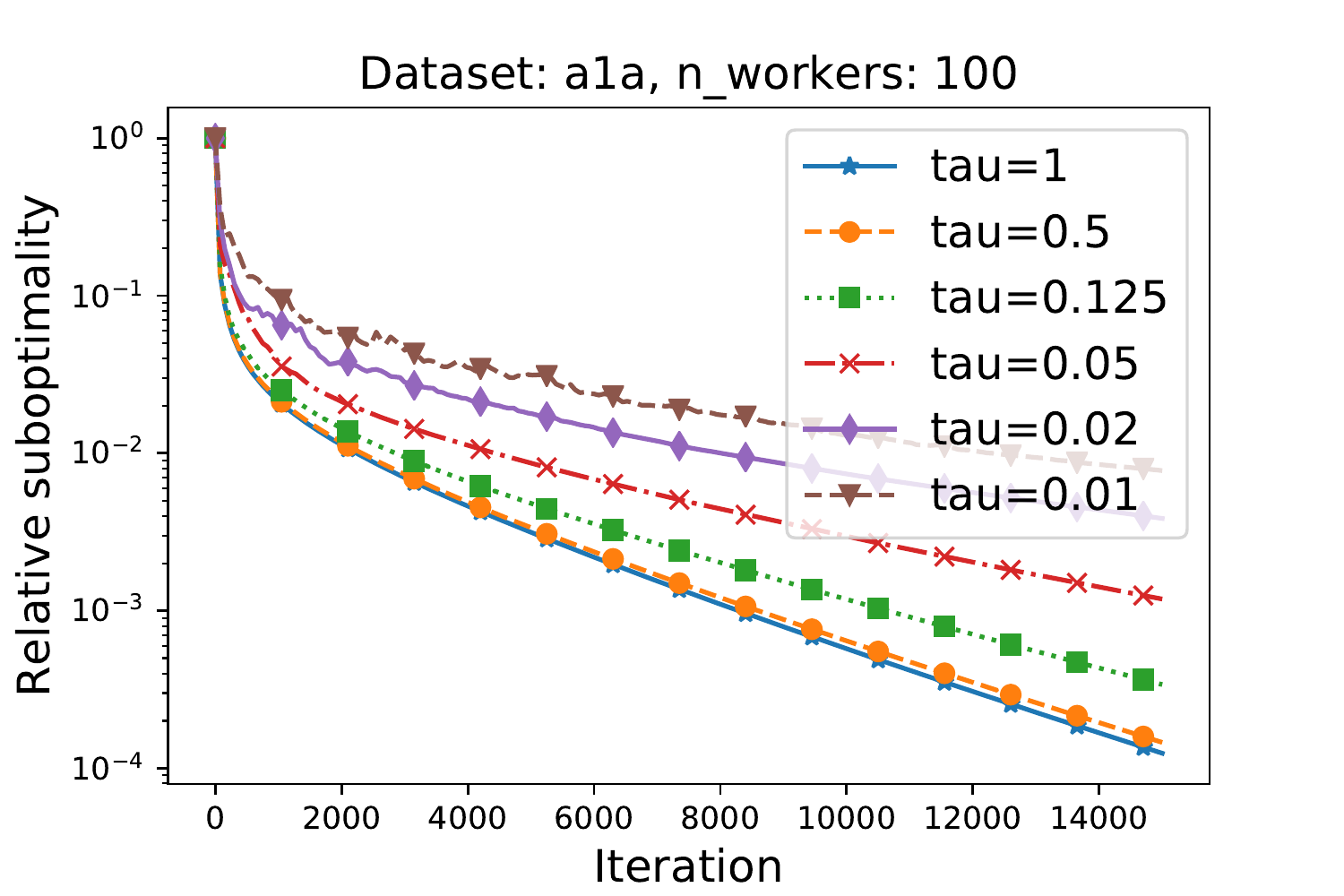}
\end{minipage}%
\\
\begin{minipage}{0.33\textwidth}
  \centering
\includegraphics[width =  \textwidth ]{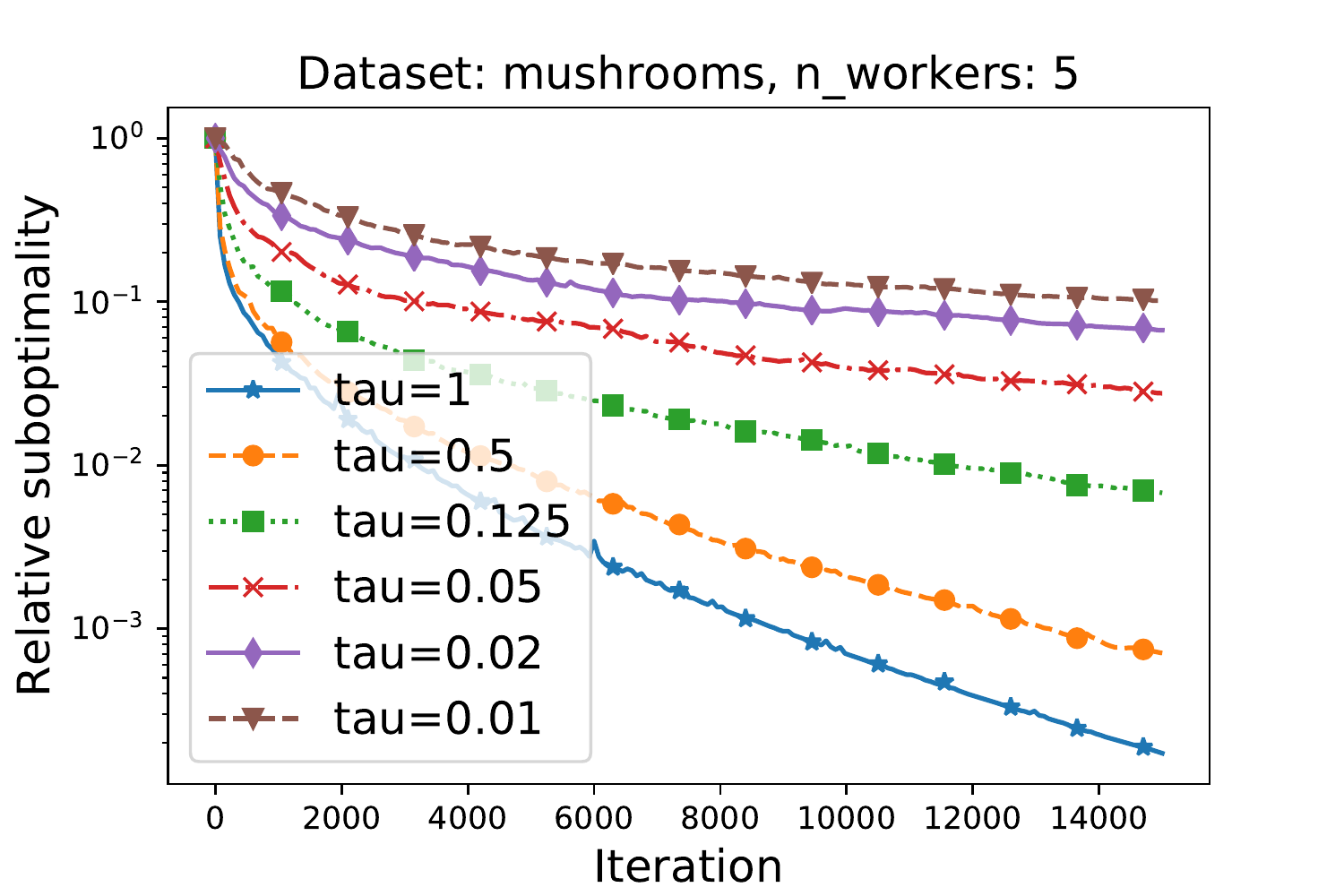}
\end{minipage}%
\begin{minipage}{0.33\textwidth}
  \centering
\includegraphics[width =  \textwidth ]{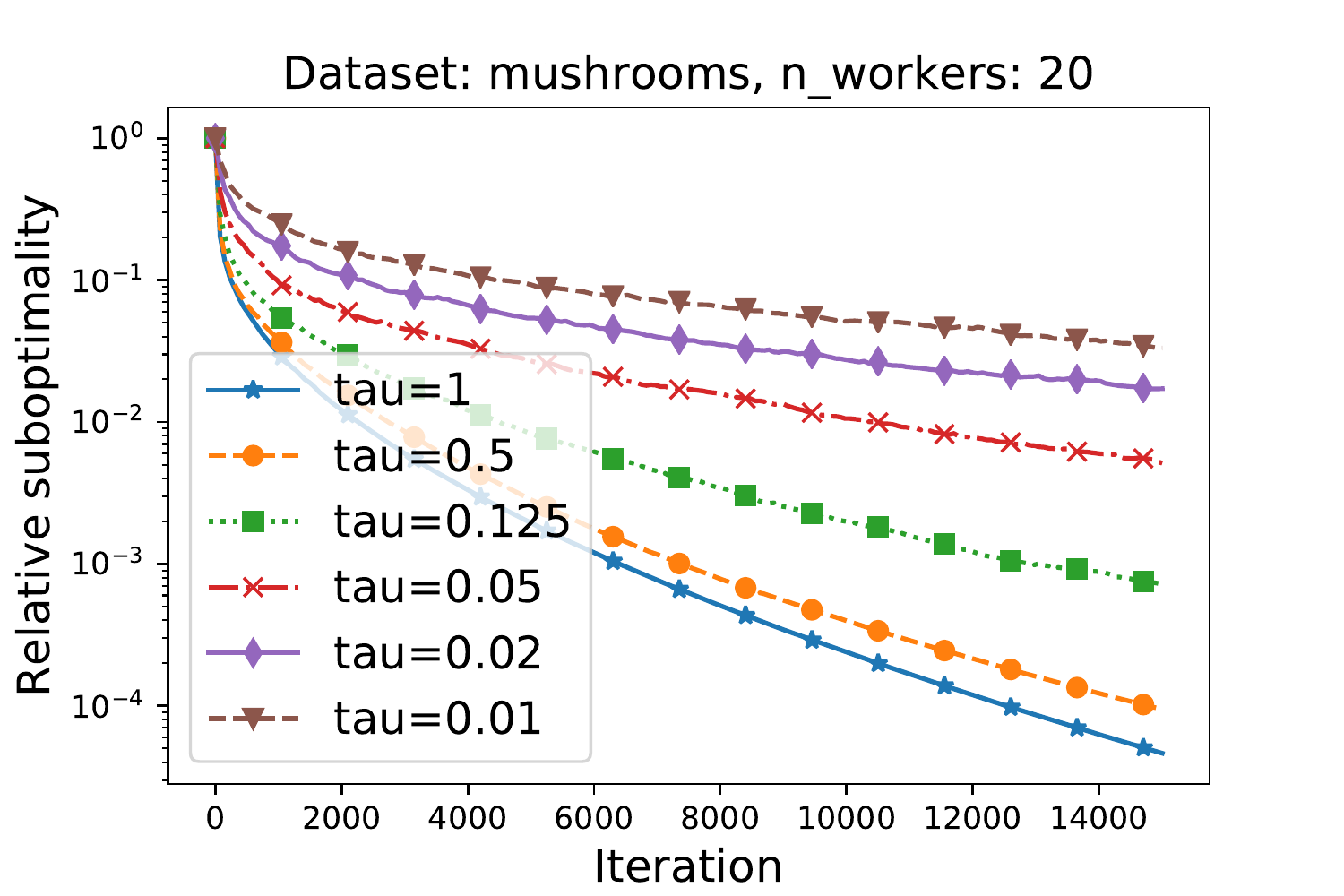}
\end{minipage}%
\begin{minipage}{0.33\textwidth}
  \centering
\includegraphics[width =  \textwidth ]{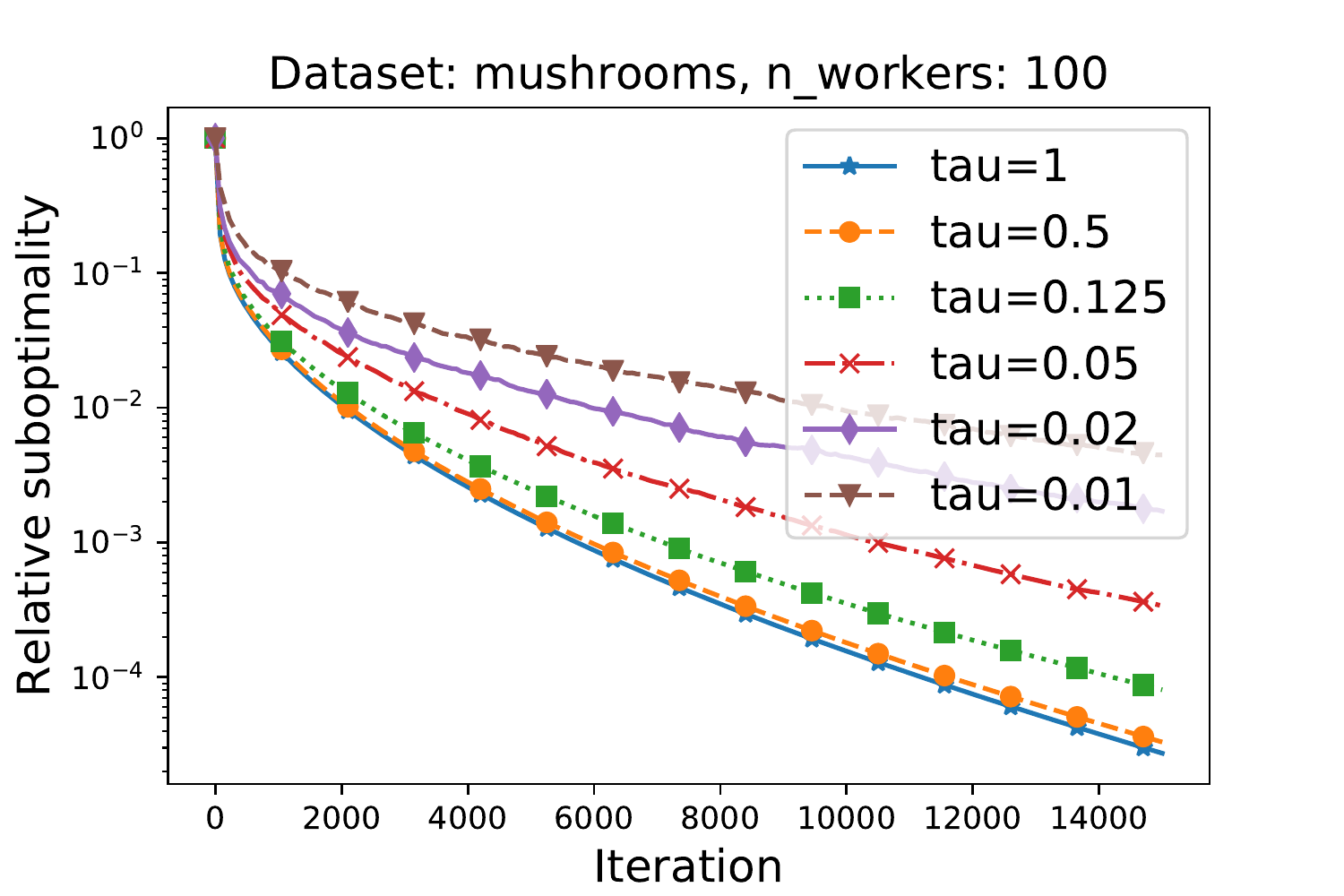}
\end{minipage}%
\\
\begin{minipage}{0.33\textwidth}
  \centering
\includegraphics[width =  \textwidth ]{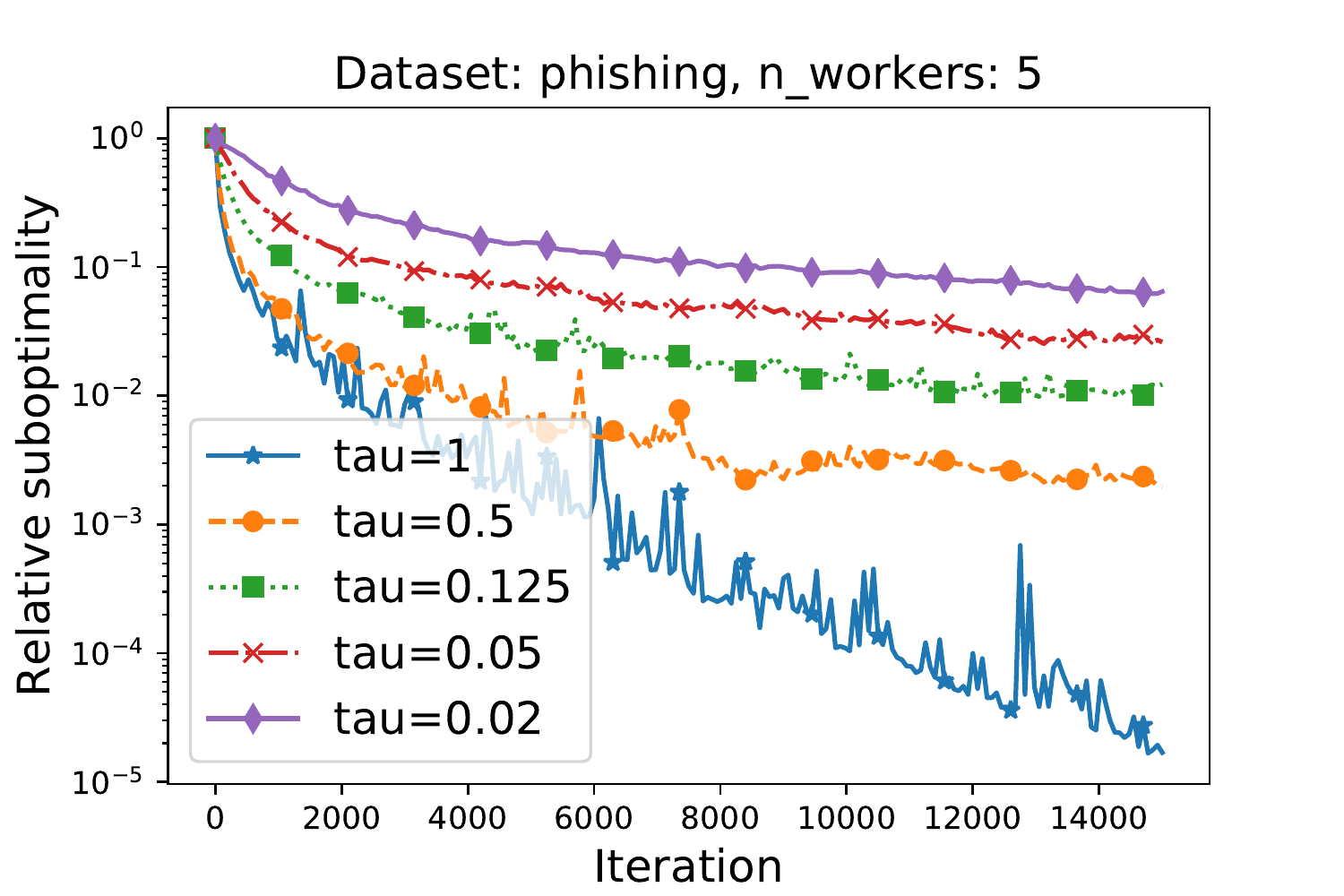}
\end{minipage}%
\begin{minipage}{0.33\textwidth}
  \centering
\includegraphics[width =  \textwidth ]{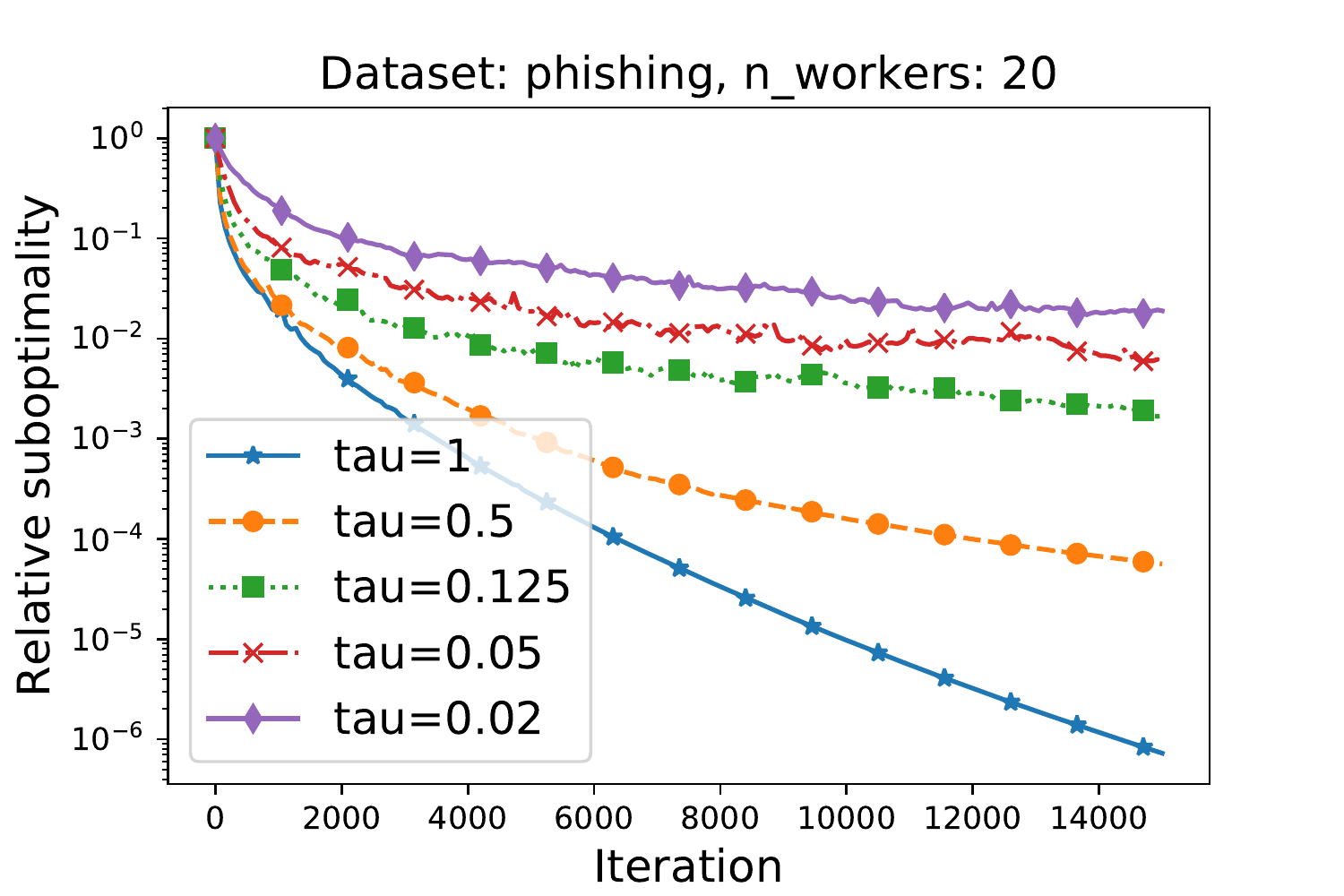}
\end{minipage}%
\begin{minipage}{0.33\textwidth}
  \centering
\includegraphics[width =  \textwidth ]{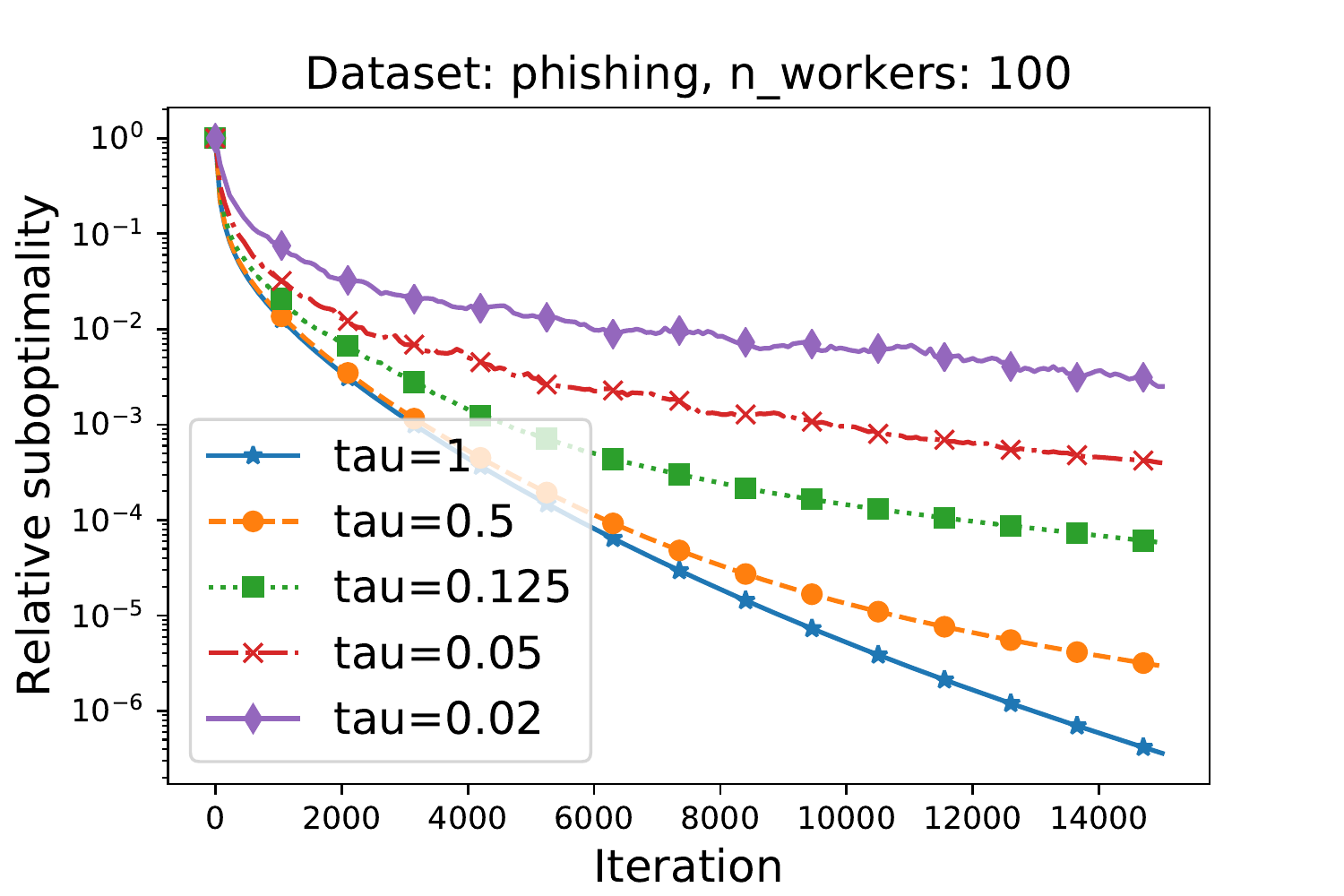}
\end{minipage}%
\\
\begin{minipage}{0.33\textwidth}
  \centering
\includegraphics[width =  \textwidth ]{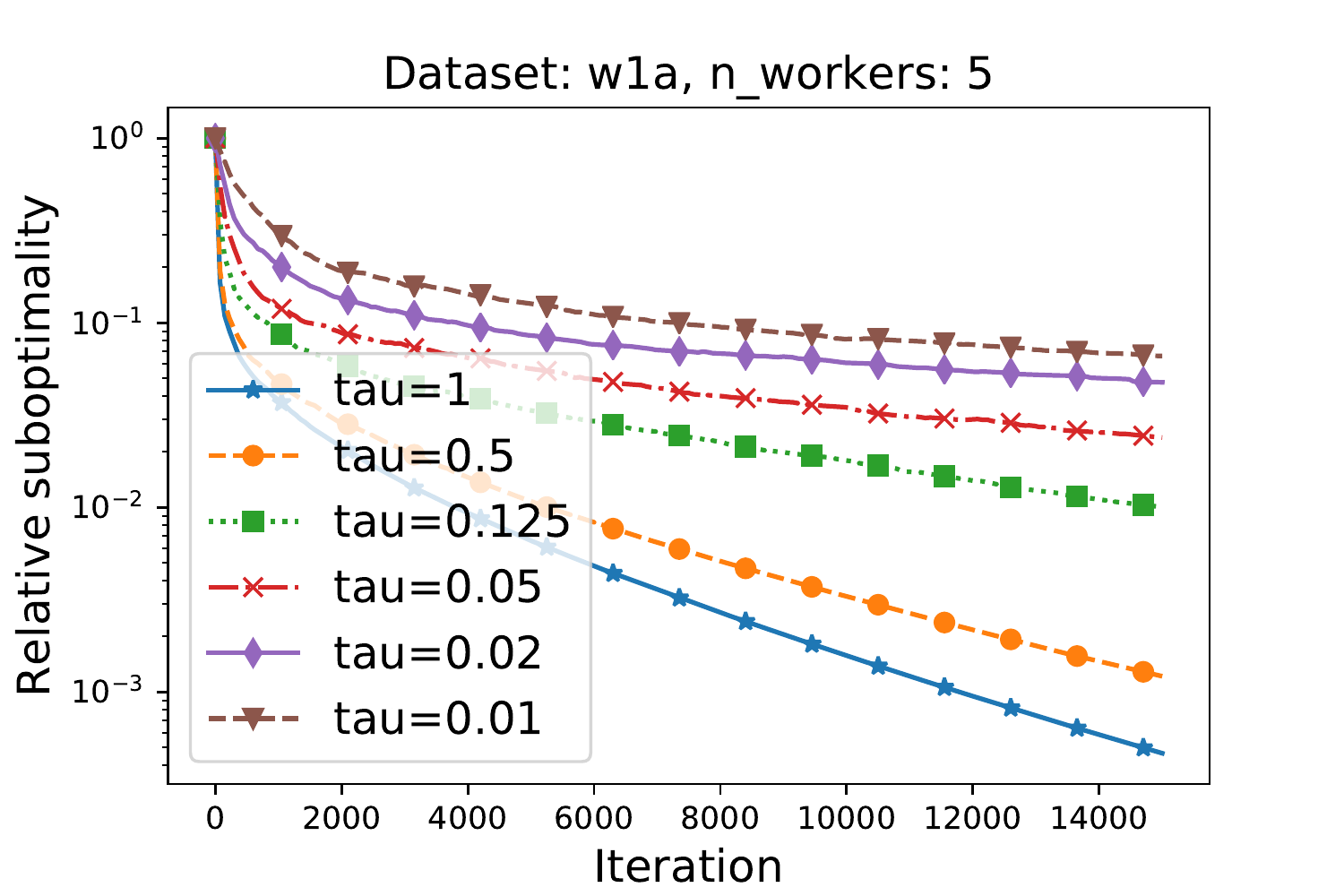}
\end{minipage}%
\begin{minipage}{0.33\textwidth}
  \centering
\includegraphics[width =  \textwidth ]{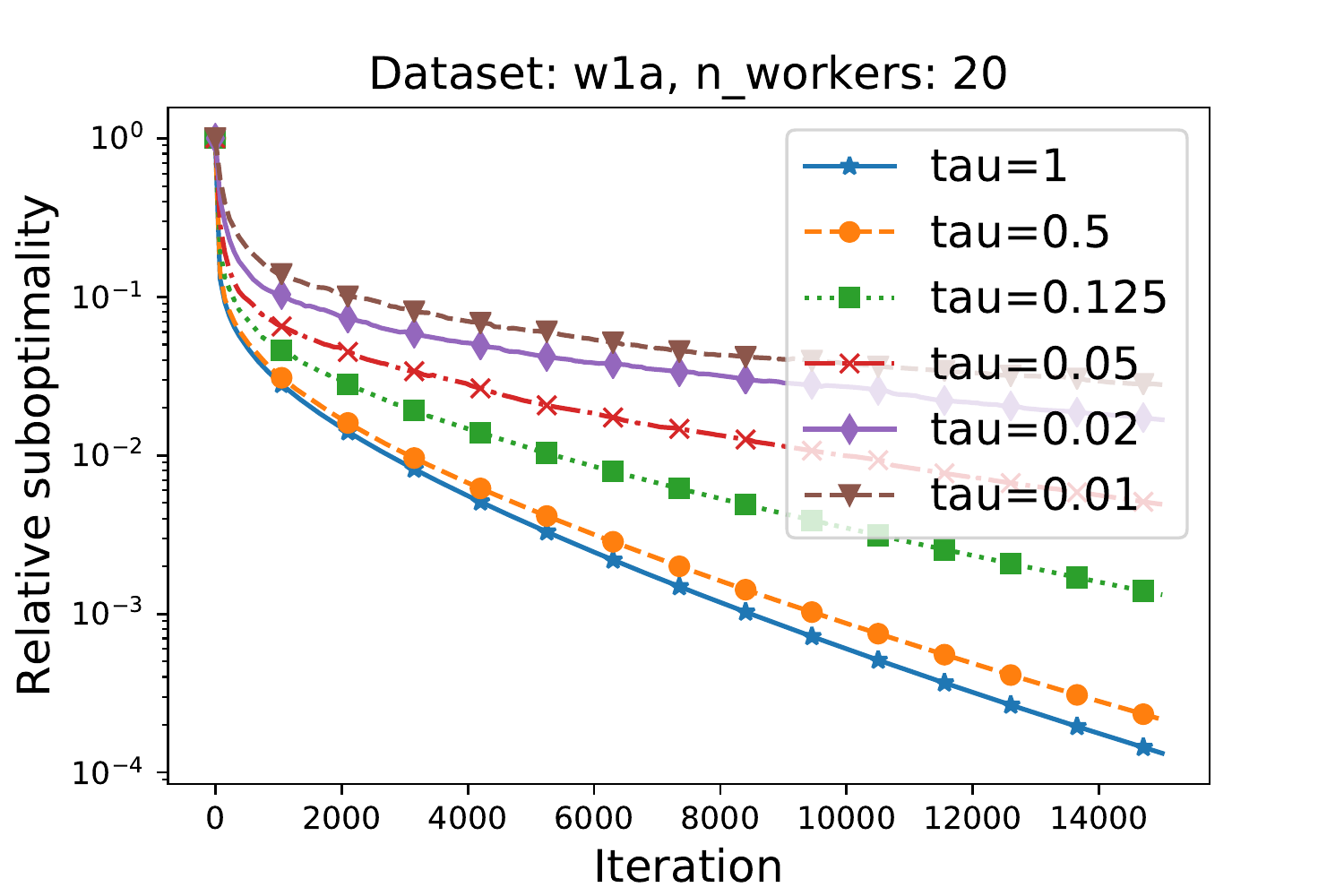}
\end{minipage}%
\begin{minipage}{0.33\textwidth}
  \centering
\includegraphics[width =  \textwidth ]{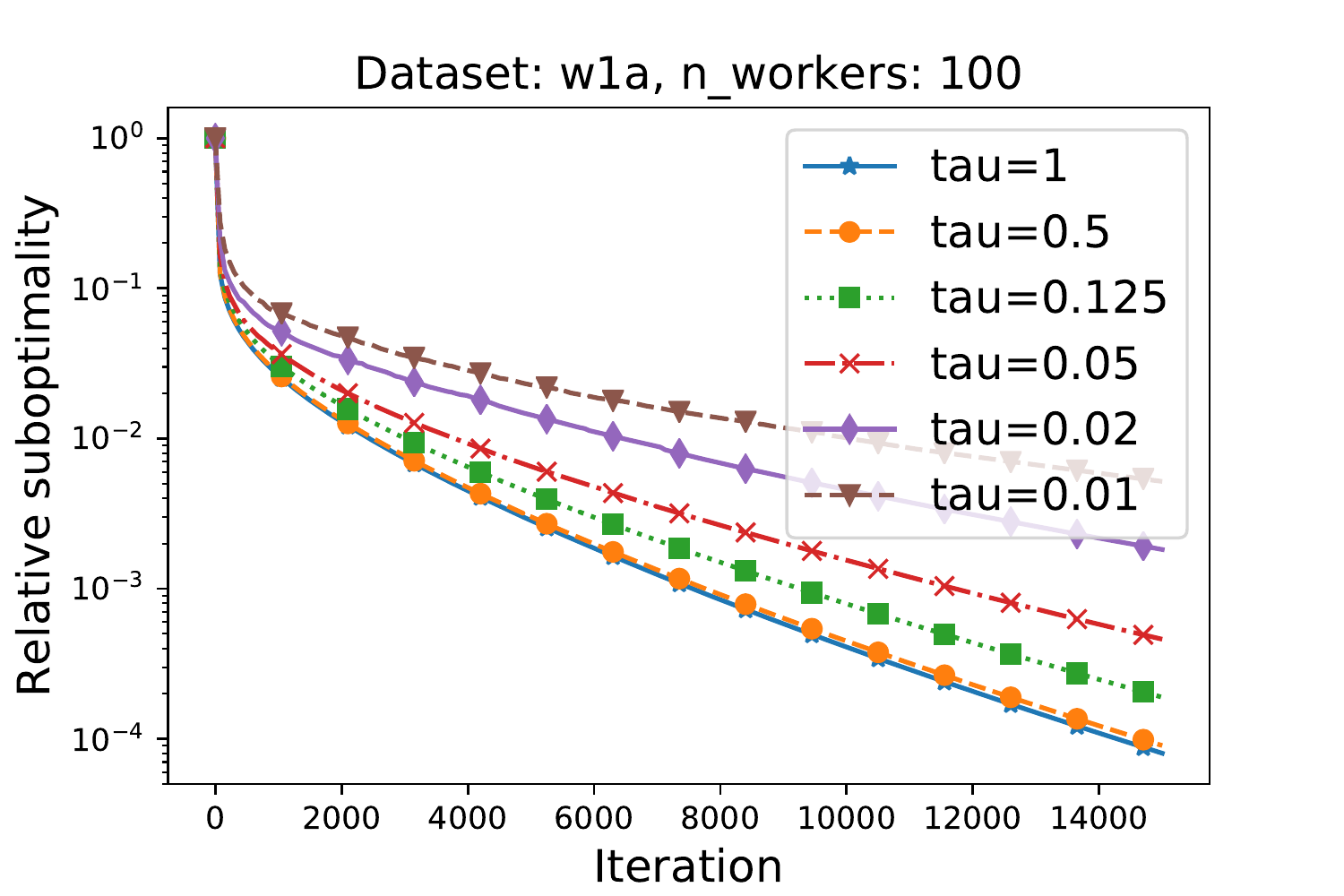}
\end{minipage}%
\\
\caption{Comparison of Algorithm~\ref{alg:saga} for different values of $\tau$. Stepsize $\alpha = \frac{1}{L(3n^{-1}+\tau)}$ is chosen in each case. For this experiment, we choose smaller regularization; $\lambda  = 0.000025 $. }\label{fig:99_saga2}
\end{figure}

\subsection{{\tt ISEGA} \label{sec:99_exp_sega}}
Lastly, we numerically test Algorithm~\ref{alg:sega}, and its linear convergence. For simplicity, we consider $\psi\equiv0$ in~\eqref{eq:99_problem_sega}. 

In the first experiment (Figure~\ref{fig:99_sega1}), we compare Algorithm~\ref{alg:sega} for various $(n,\tau)$ such that $n\tau=1$. For illustration, we also plot convergence of gradient descent with the analogous stepsize. As theory predicts, the method has almost same convergence speed.\footnote{We have chosen stepsize $\alpha = \frac{1}{2L}$ for GD, as this is the baseline to Algorithm~\ref{alg:sega} with zero variance. One can in fact set $\alpha = \frac{1}{L}$ for GD and get 2 times faster convergence. However, this is still only a constant factor. }

\begin{figure}[H]
\centering
\begin{minipage}{0.35\textwidth}
  \centering
\includegraphics[width =  \textwidth ]{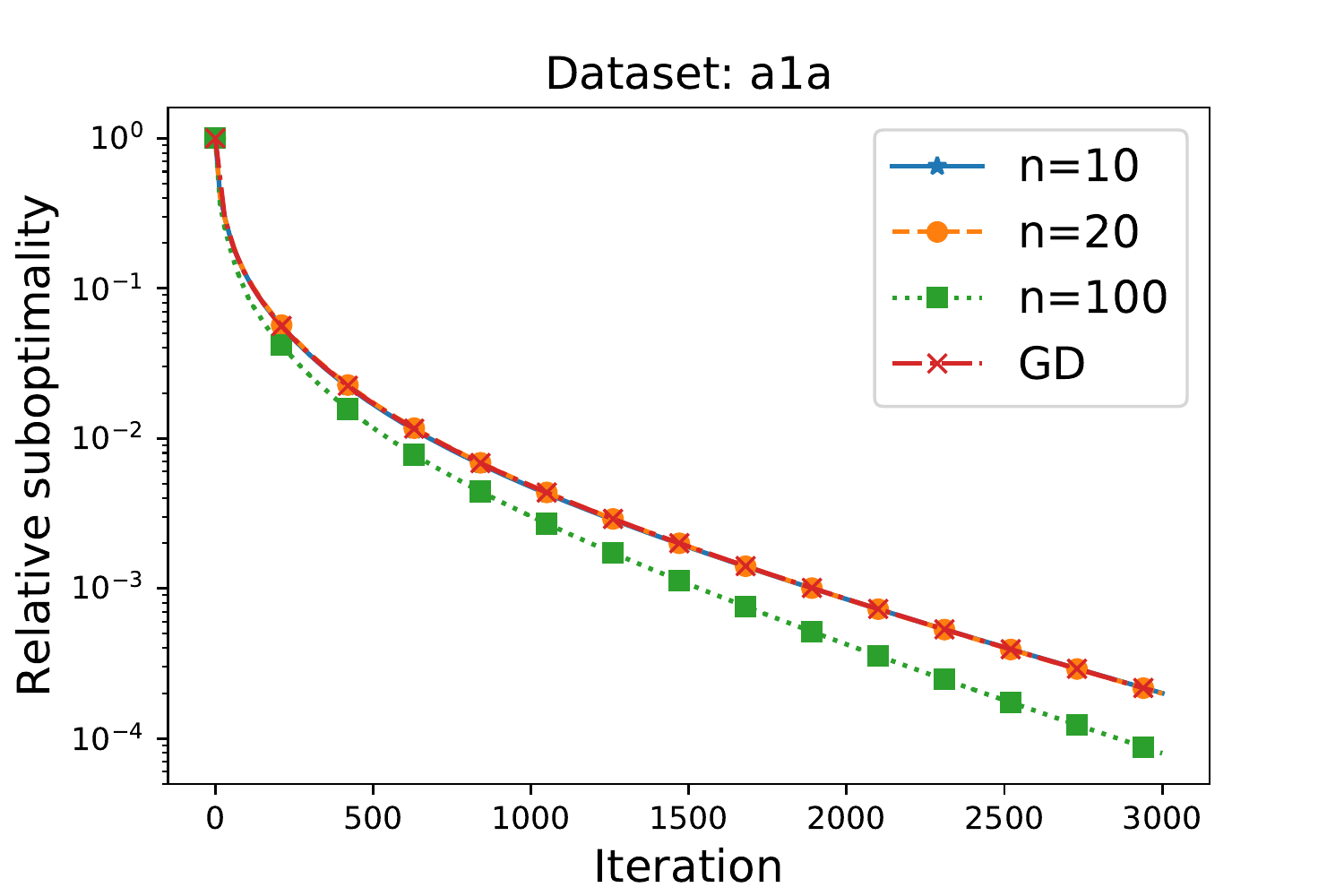}
\end{minipage}%
\begin{minipage}{0.35\textwidth}
  \centering
\includegraphics[width =  \textwidth ]{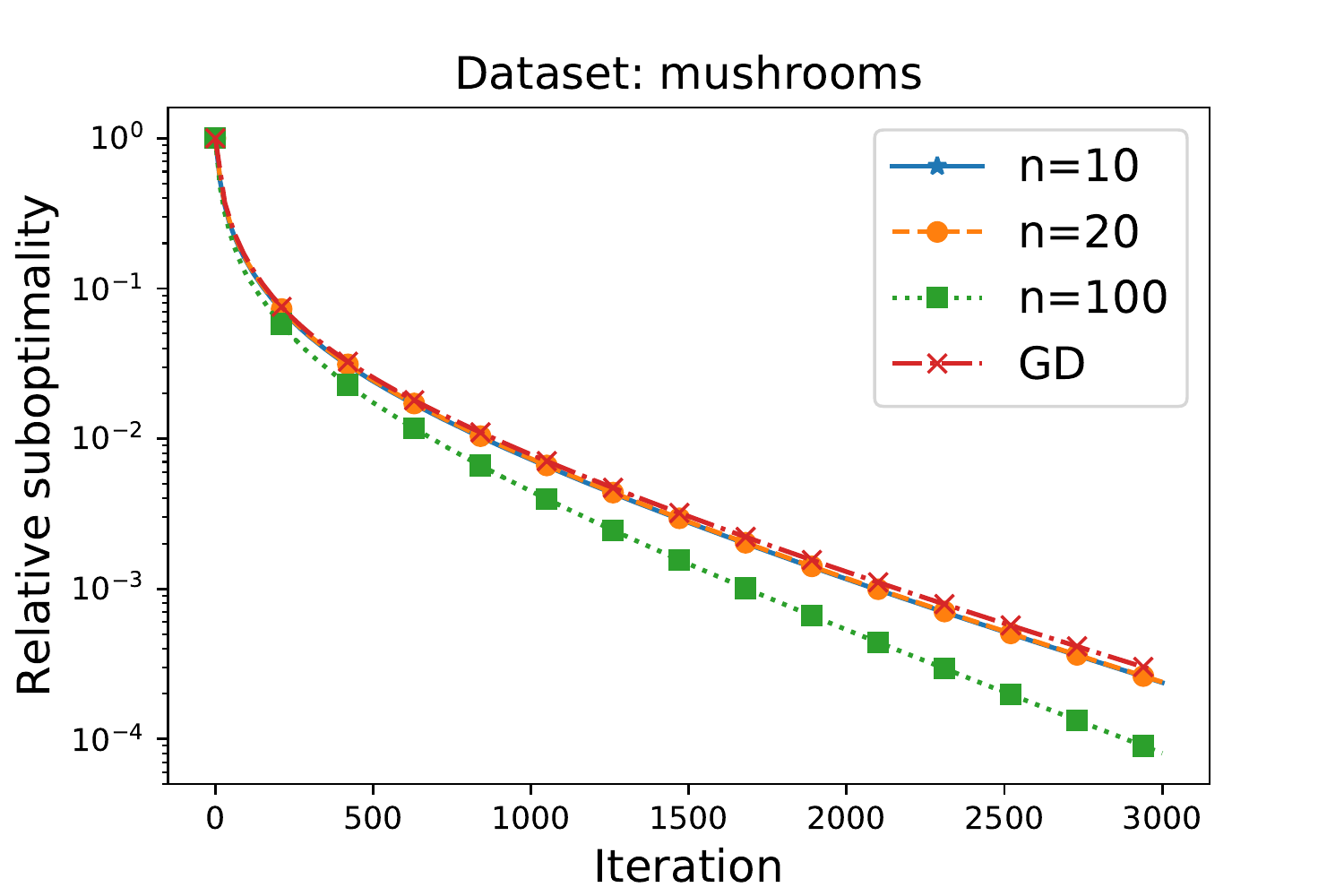}
\end{minipage}%
\\
\begin{minipage}{0.35\textwidth}
  \centering
\includegraphics[width =  \textwidth ]{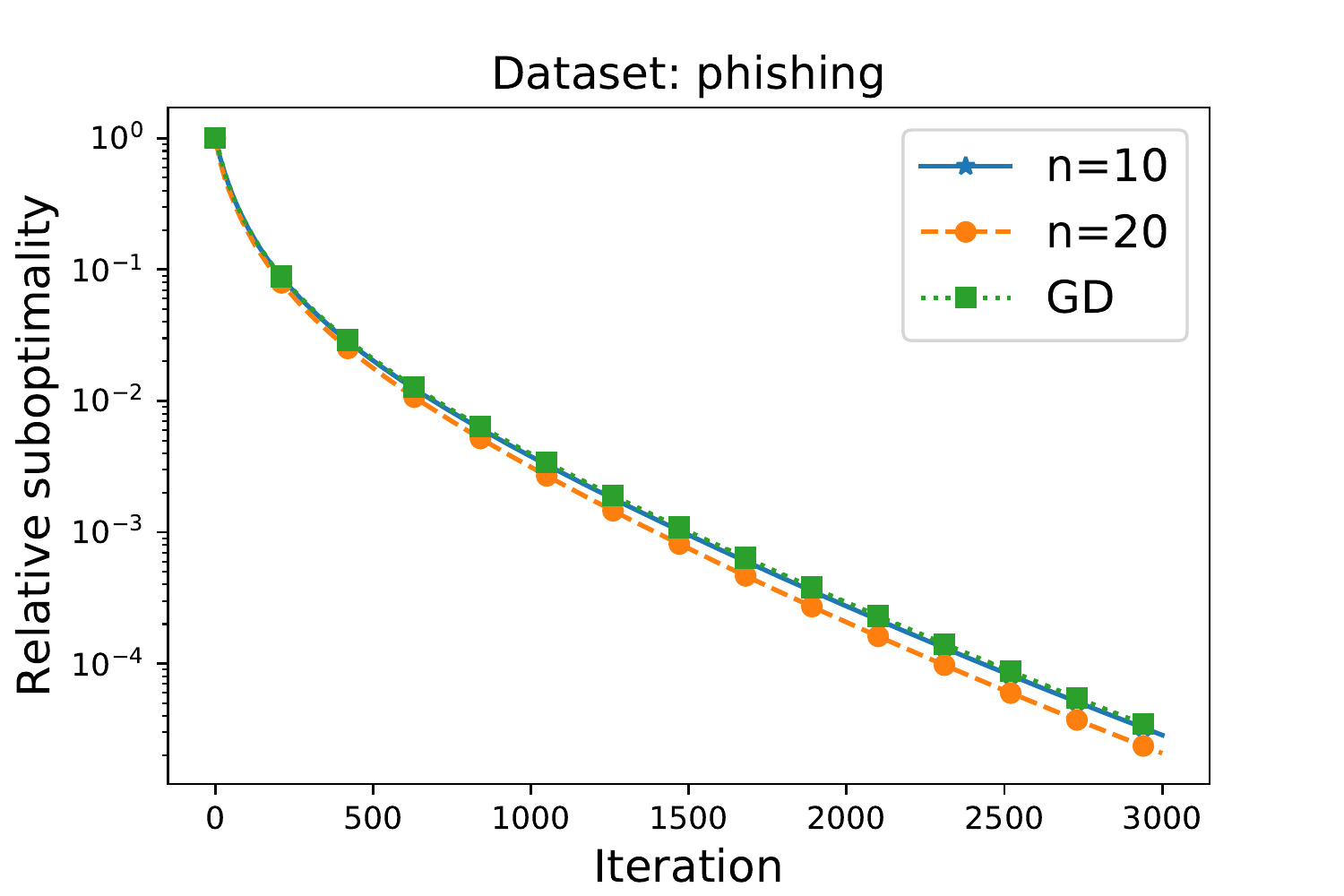}
\end{minipage}%
\begin{minipage}{0.35\textwidth}
  \centering
\includegraphics[width =  \textwidth ]{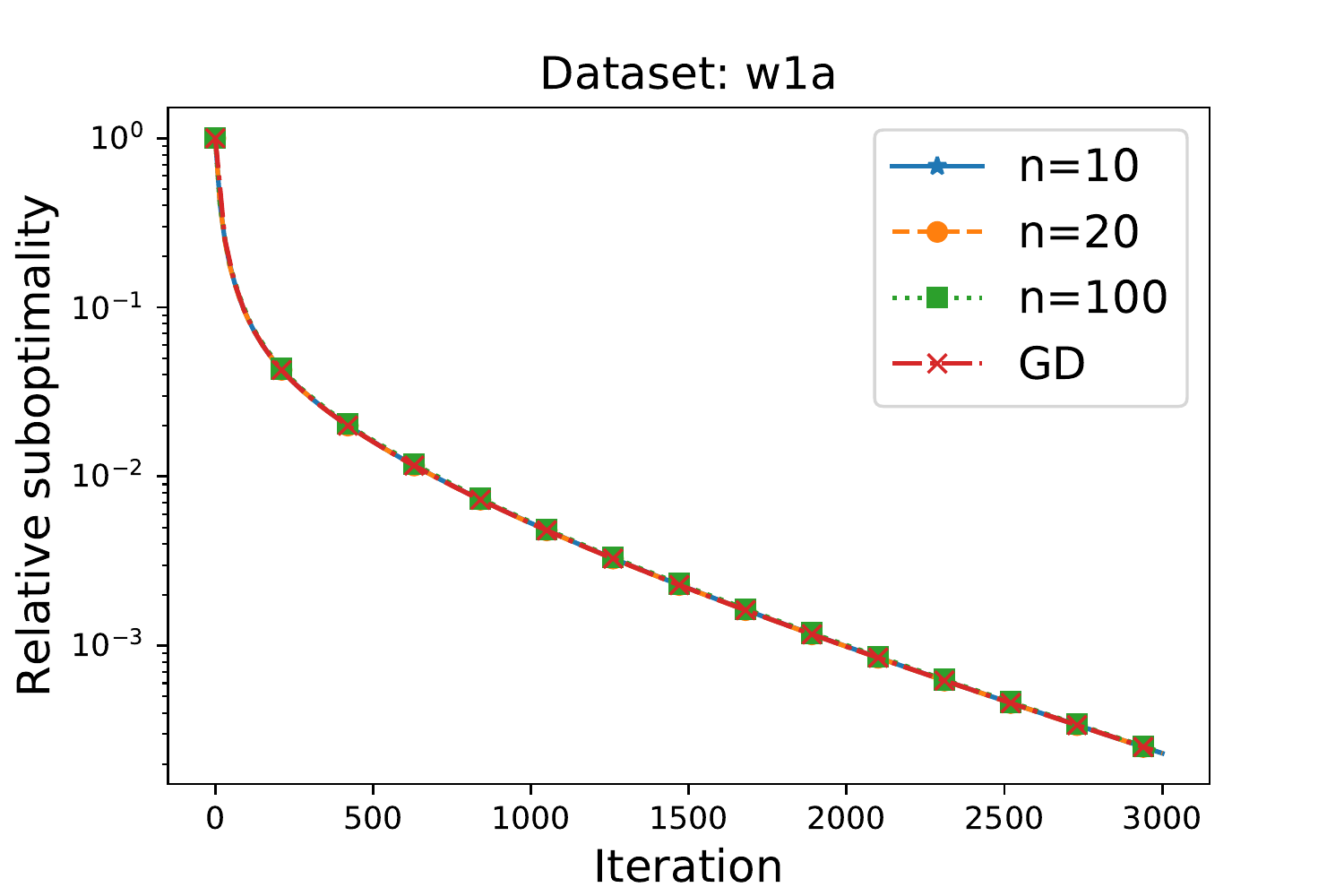}
\end{minipage}%
\\
\caption{Comparison of Algorithm~\ref{alg:sega} for various $(n,\tau)$ such that $n\tau=1$ and GD. Stepsize $\frac{1}{L\left(1+\frac{1}{n\tau}\right)} $ was chosen for Algorithm~\ref{alg:sega} and $\frac1{2L}$ for GD.}\label{fig:99_sega1}
\end{figure}

The second experiment of this section shows the convergence behavior for varying $\tau$ of Algorithm~\ref{alg:sega}. Again, the results (Figure~\ref{fig:99_sega2}) indicate that $\tau$ has a heavy impact on the convergence speed for small $n$. However, as $n$ increases, the effect of $\tau$ is diminishing. In particular, for increasing $\tau$ beyond $n^{-1}$ does not yield a significant speedup.

\begin{figure}[H]
\centering
\begin{minipage}{0.33\textwidth}
  \centering
\includegraphics[width =  \textwidth ]{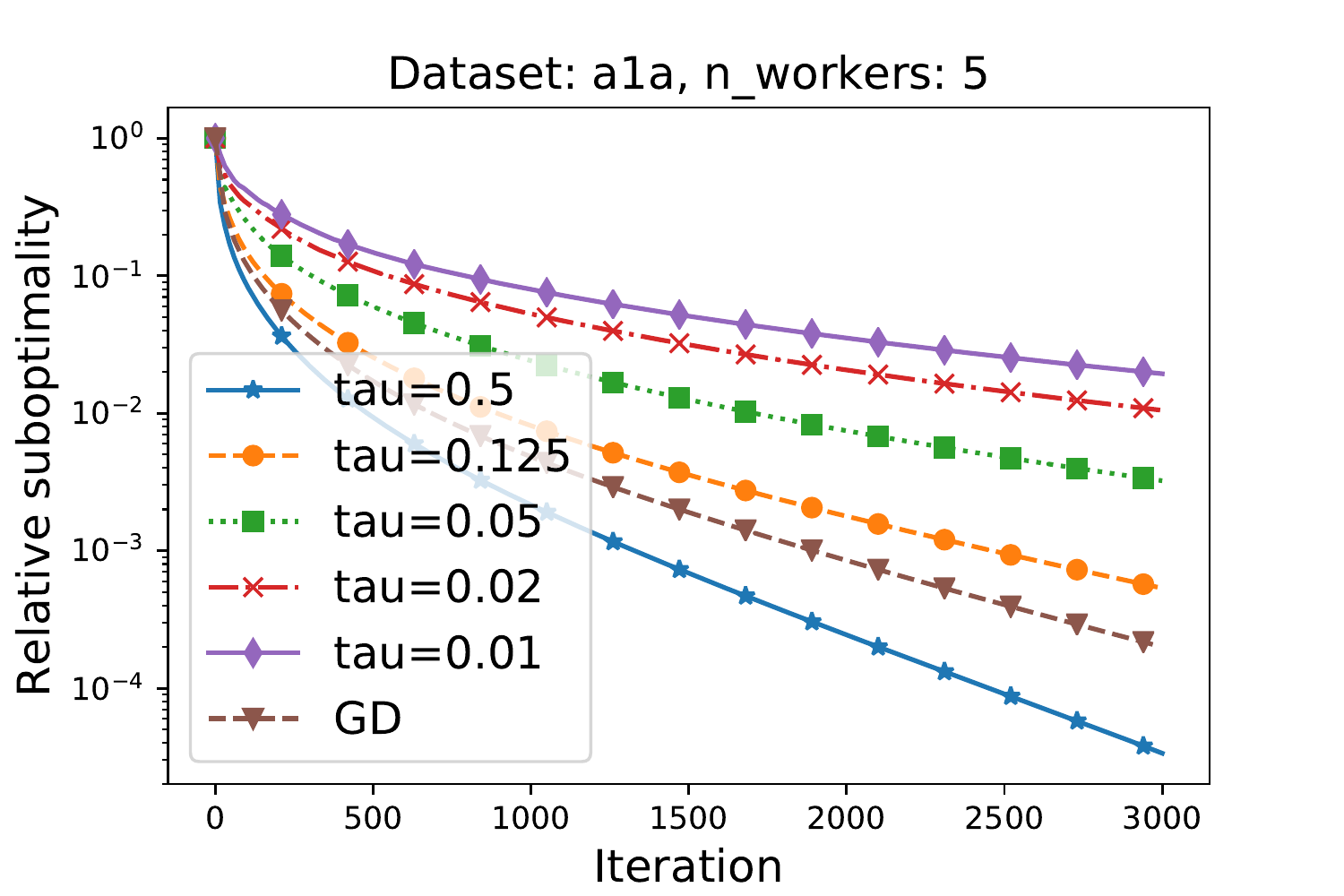}
\end{minipage}%
\begin{minipage}{0.33\textwidth}
  \centering
\includegraphics[width =  \textwidth ]{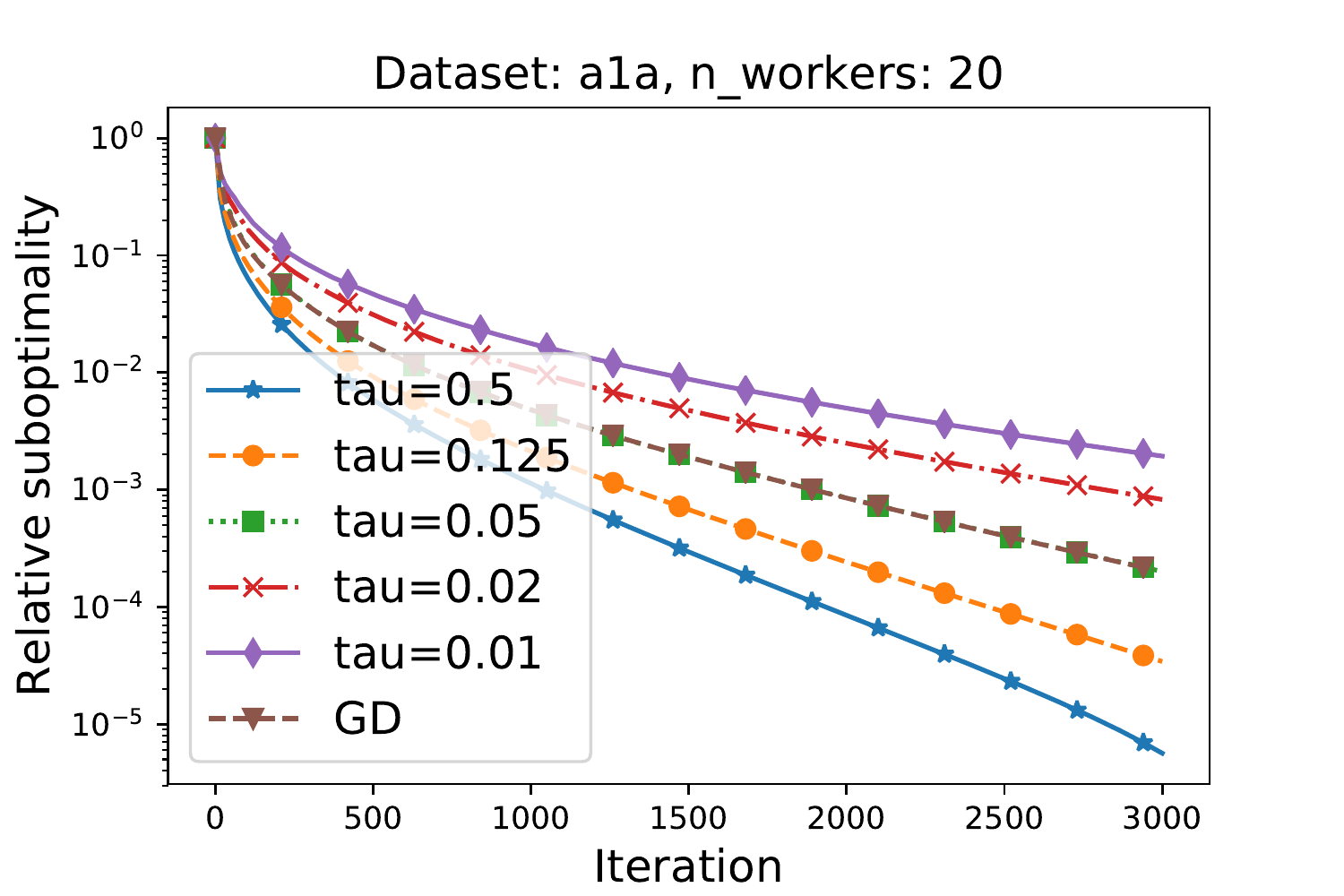}
\end{minipage}%
\begin{minipage}{0.33\textwidth}
  \centering
\includegraphics[width =  \textwidth ]{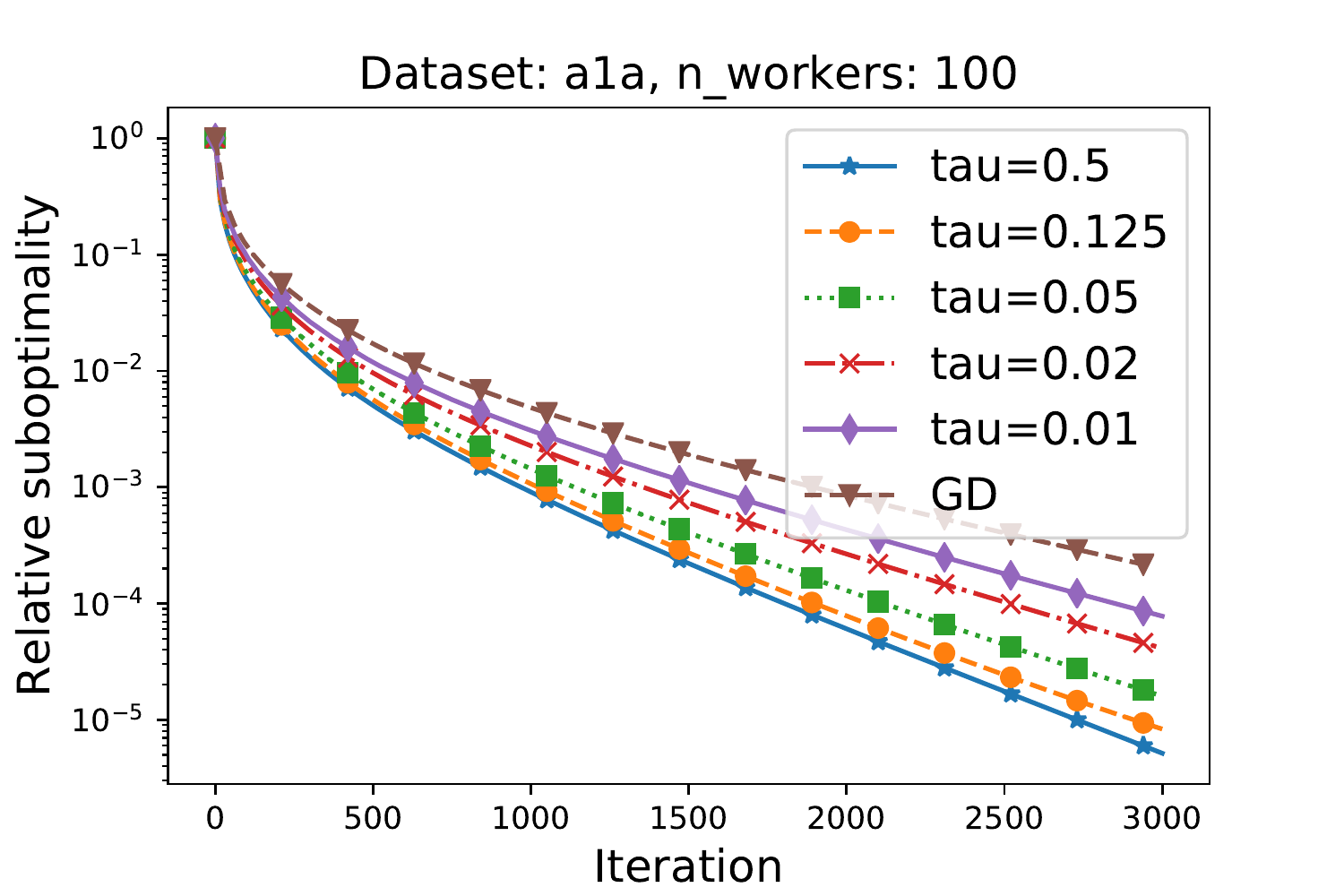}
\end{minipage}%
\\
\begin{minipage}{0.33\textwidth}
  \centering
\includegraphics[width =  \textwidth ]{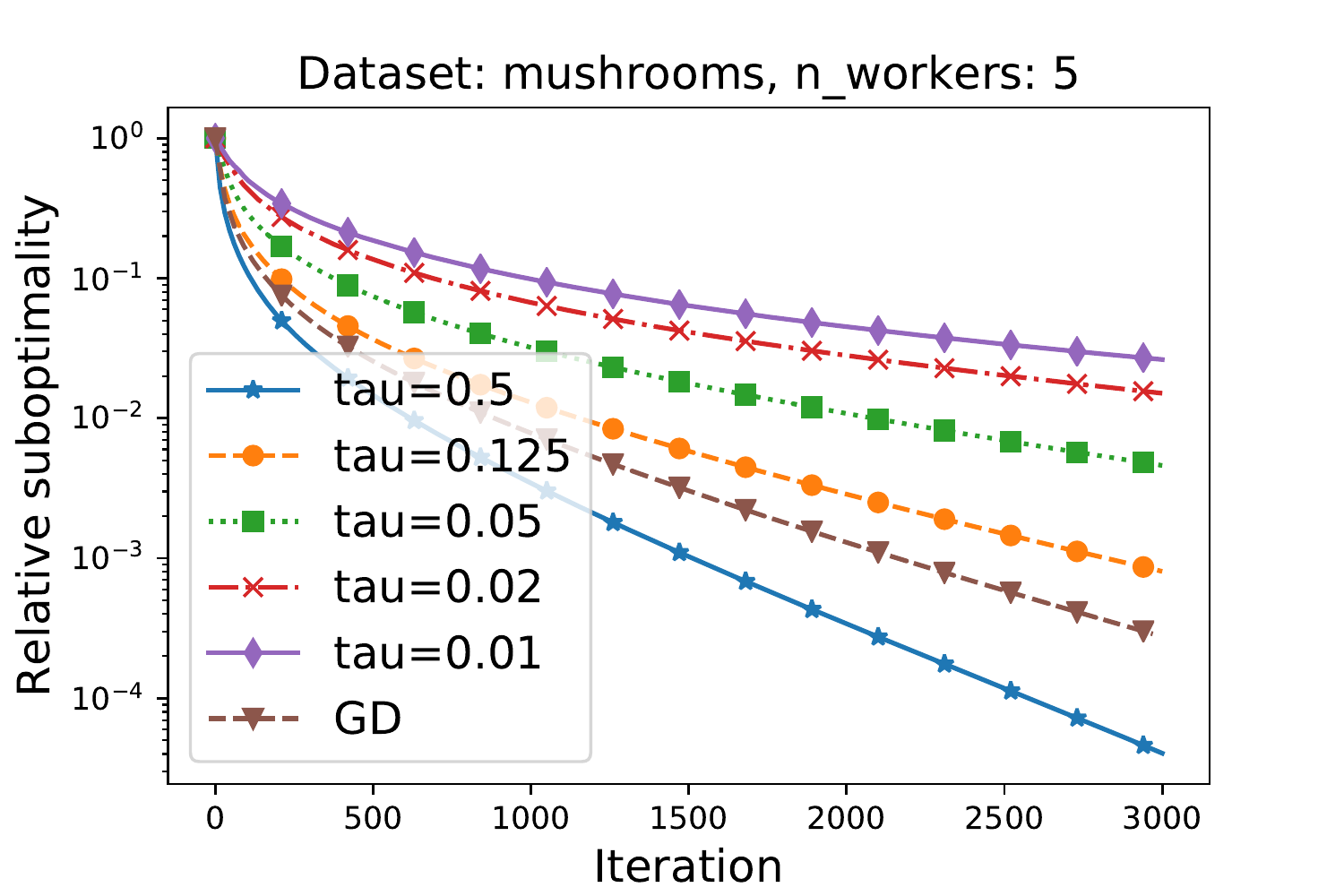}
\end{minipage}%
\begin{minipage}{0.33\textwidth}
  \centering
\includegraphics[width =  \textwidth ]{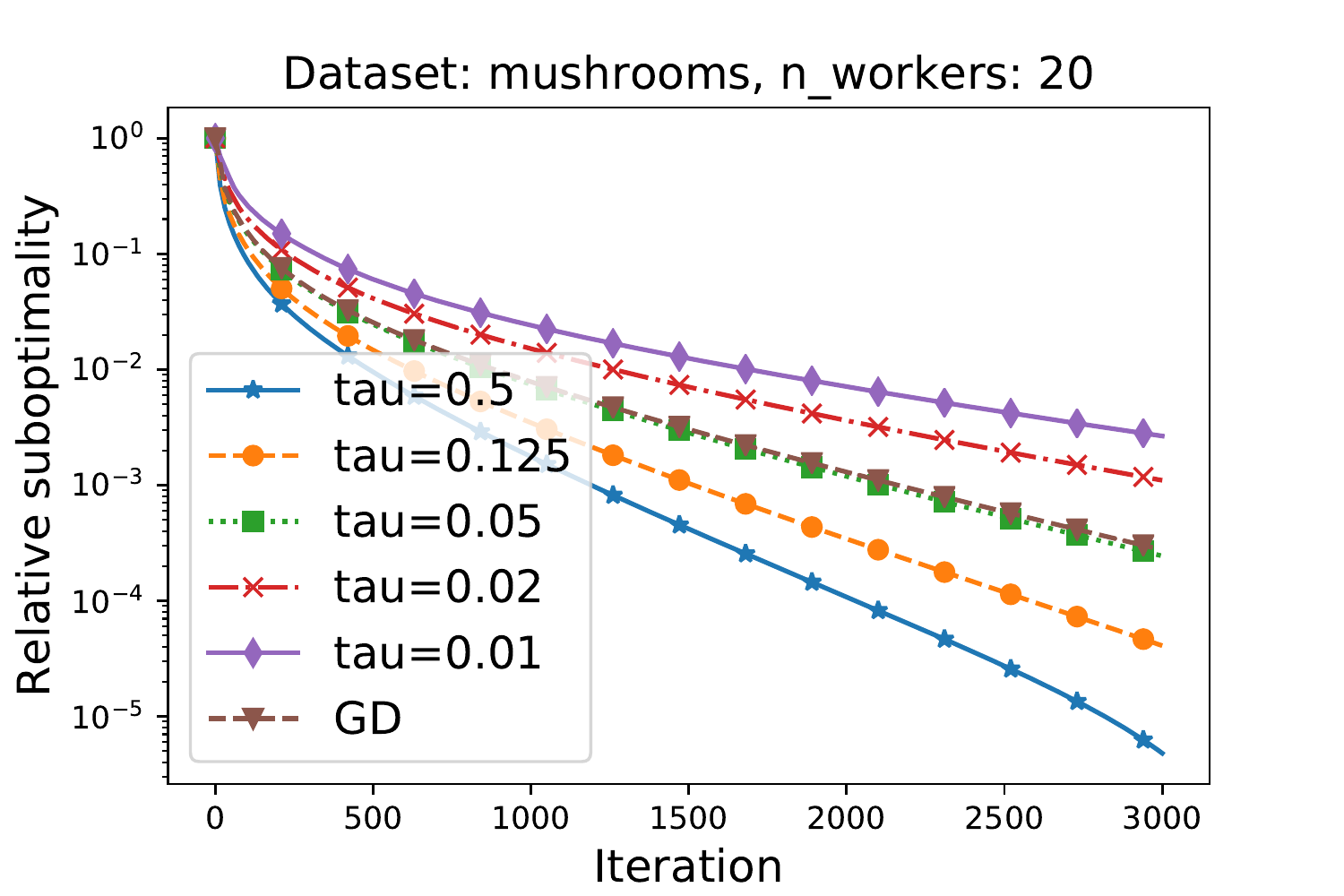}
\end{minipage}%
\begin{minipage}{0.33\textwidth}
  \centering
\includegraphics[width =  \textwidth ]{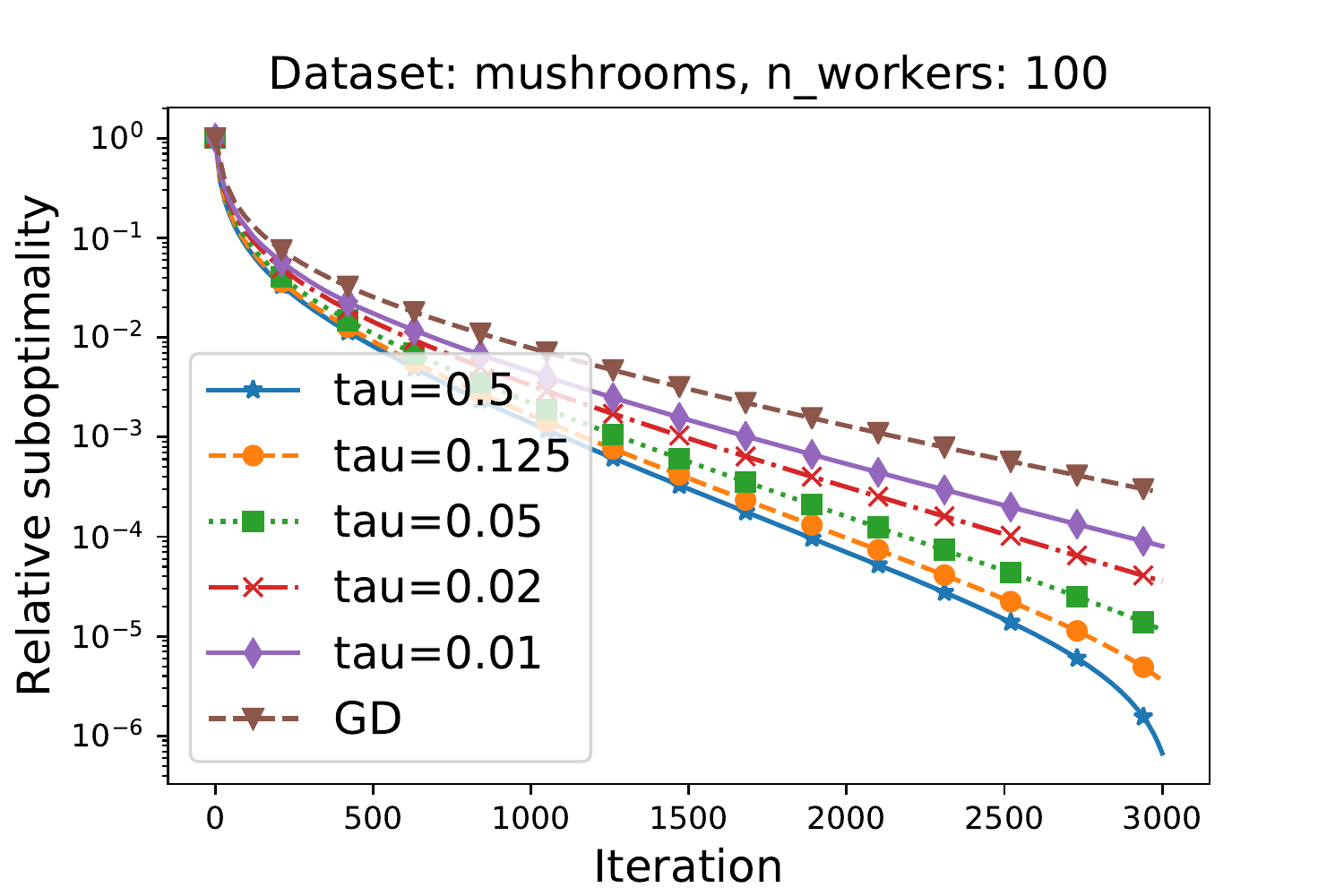}
\end{minipage}%
\\
\begin{minipage}{0.33\textwidth}
  \centering
\includegraphics[width =  \textwidth ]{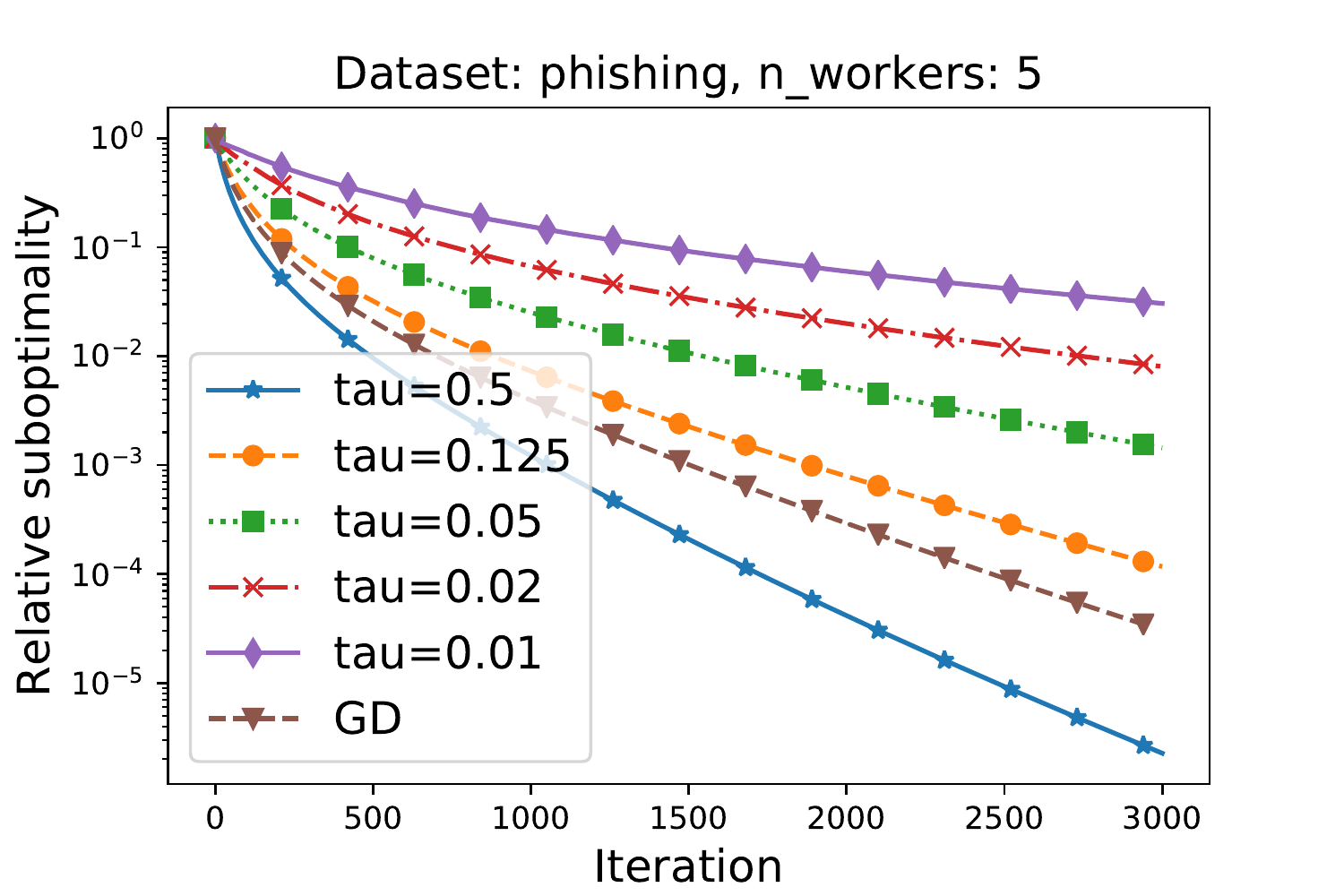}
\end{minipage}%
\begin{minipage}{0.33\textwidth}
  \centering
\includegraphics[width =  \textwidth ]{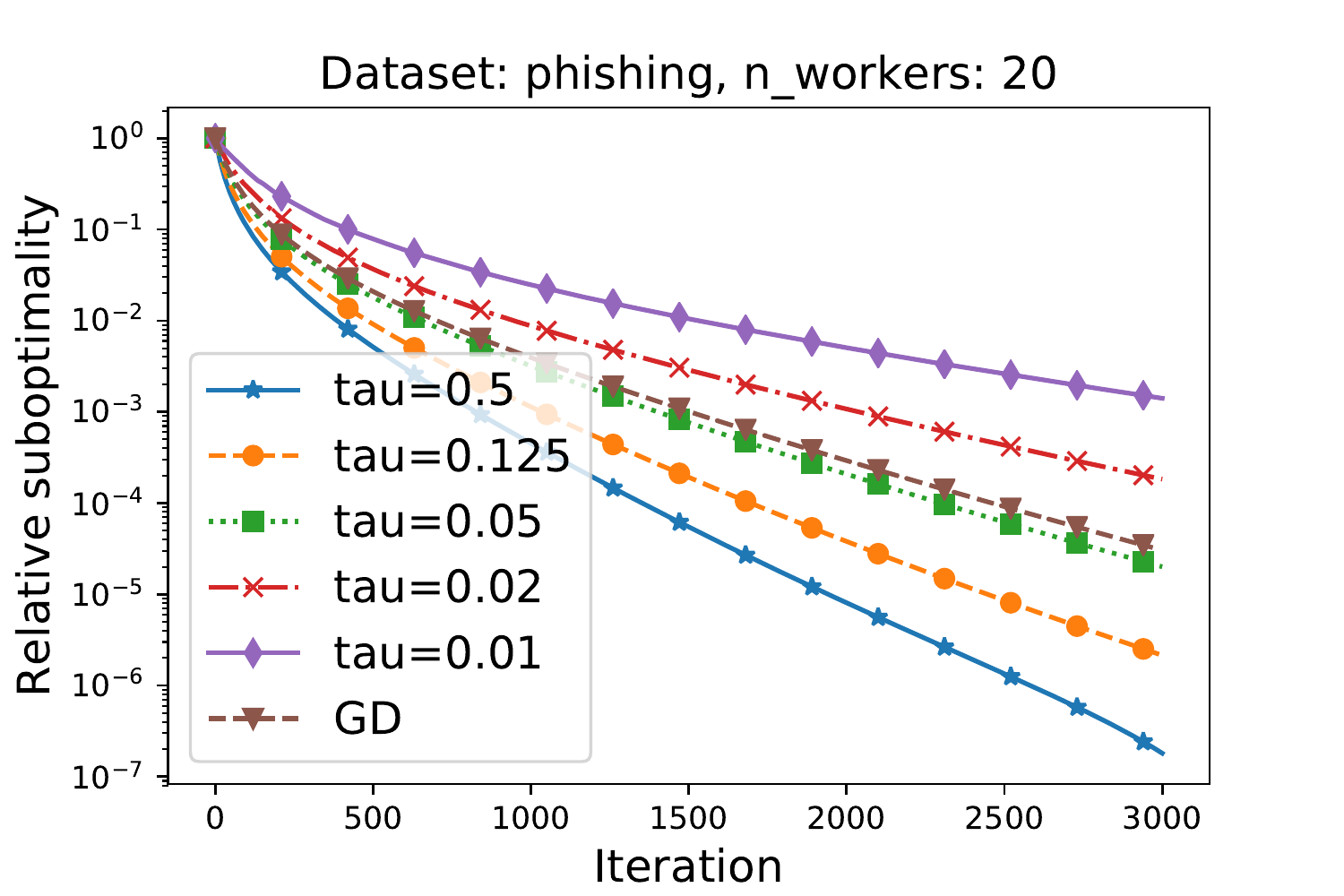}
\end{minipage}%
\\
\begin{minipage}{0.33\textwidth}
  \centering
\includegraphics[width =  \textwidth ]{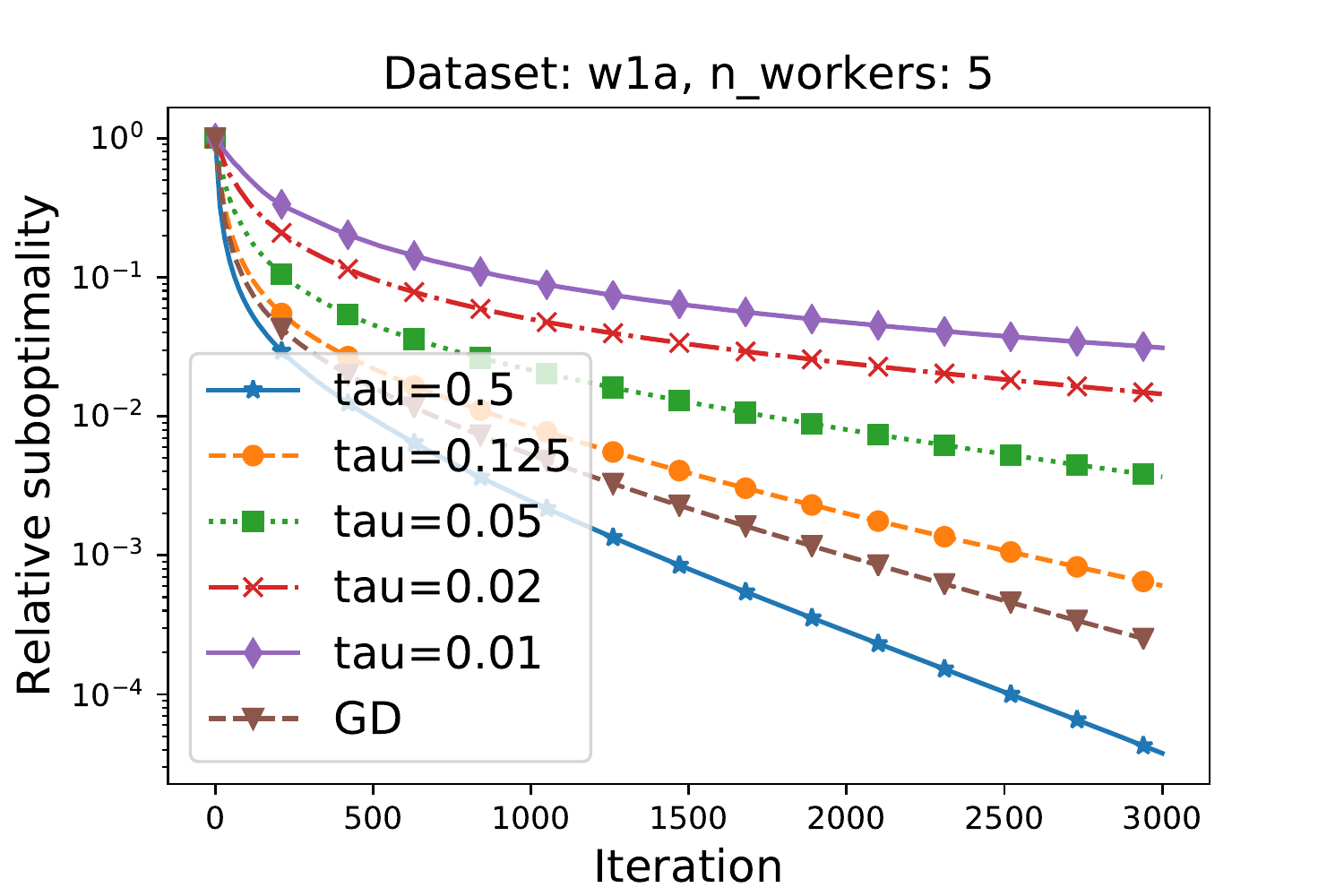}
\end{minipage}%
\begin{minipage}{0.33\textwidth}
  \centering
\includegraphics[width =  \textwidth ]{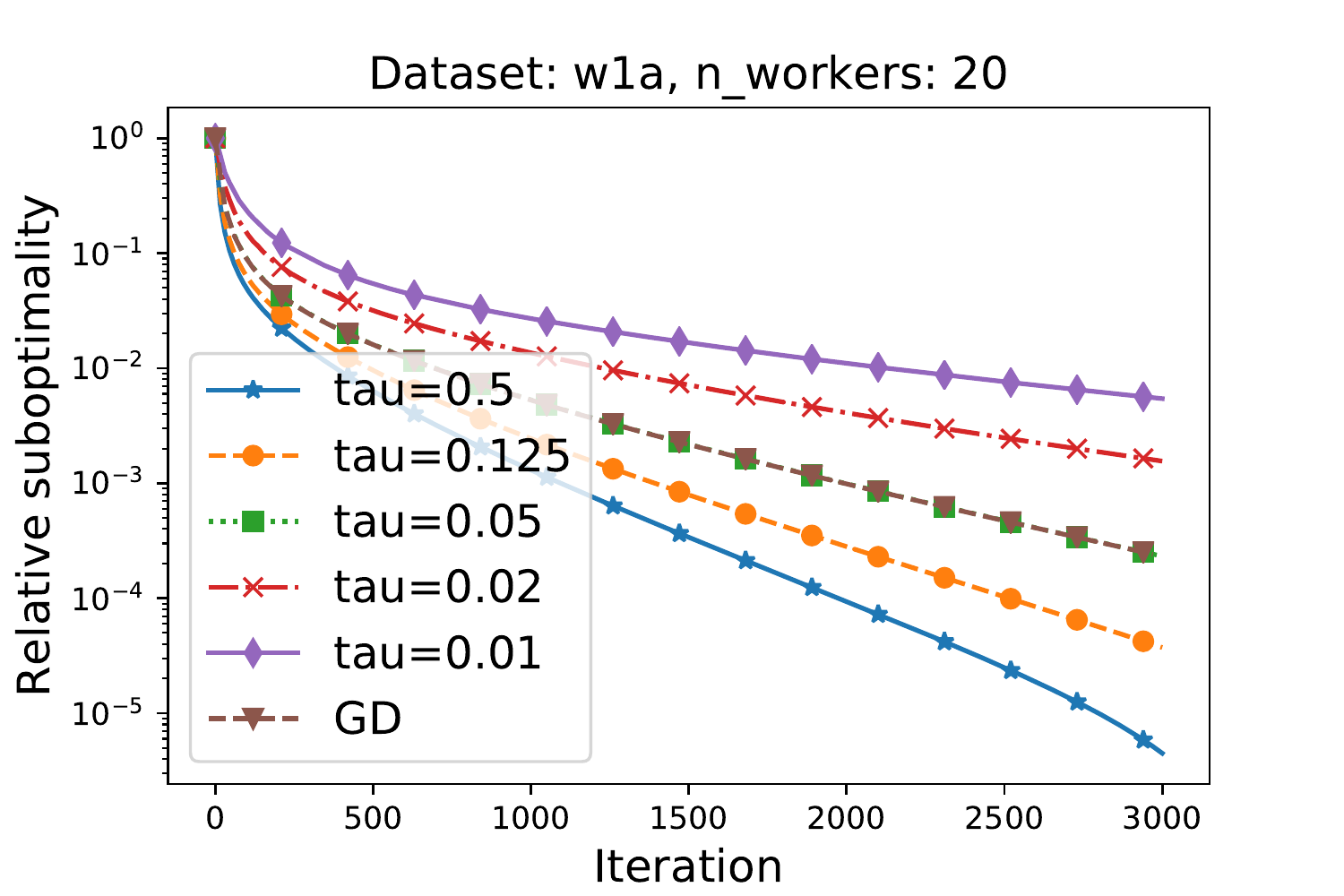}
\end{minipage}%
\begin{minipage}{0.33\textwidth}
  \centering
\includegraphics[width =  \textwidth ]{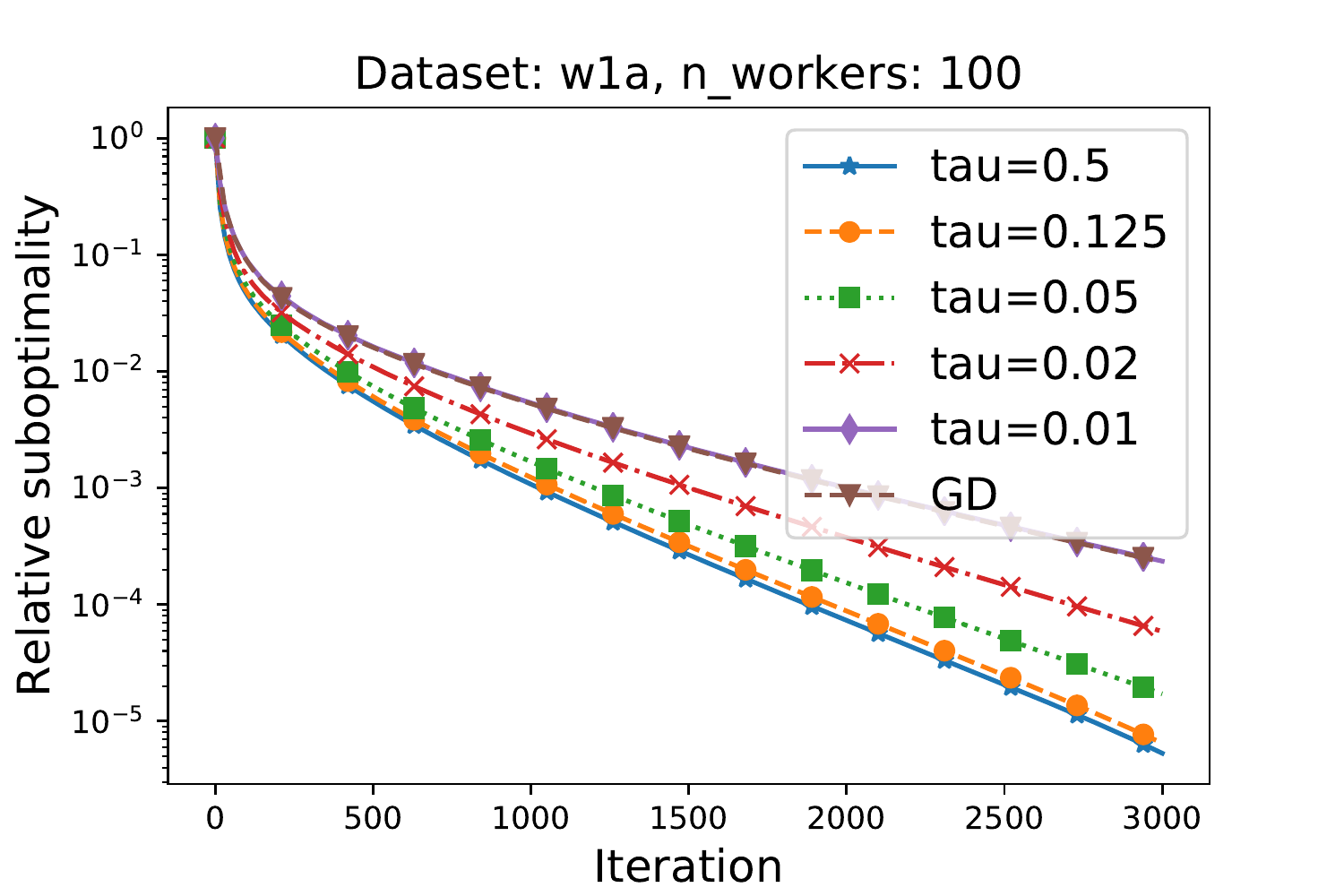}
\end{minipage}%
\caption{Comparison of Algorithm~\ref{alg:sega} for different values of $\tau$. Stepsize $\alpha = \frac{1}{L\left(1+ \frac{1}{n\tau}\right)}$ is chosen in each case.}\label{fig:99_sega2}
\end{figure}

\section{Conclusion}
In this chapter, we have proposed a strategy for reducing the worker$\rightarrow$server communication by  $\cO(\frac{n-1}{n}) \times 100  \%$, where $n$ is the number of workers. The algorithms we introduced are merely act as demonstrations of what can be achieved using our main insight, and many further extensions are possible. Specifically, in  the next chapter we propose {\tt GJS}: a new algorithm that obtains several further extensions of the methods developed in this  chapter in special cases:
\begin{itemize}
\item Distributed {\tt ISAGA} requires $\nabla f_i(x)=0$. {\tt GJS} allows to develop {\tt SEGA} approach on top of it in order to drop this requirement.

\item  Standard coordinate descent is able to exploit a complex smoothness structure of objective in order to sample coordinates non-uniformly~\cite{qu2016coordinate2, csiba2018importance}. As a special case of {\tt GJS}, we obtain importance sampling variants of multiple algorithms proposed here.

\end{itemize}

\chapter{One Method to Rule Them All: Variance Reduction for Data, Parameters and Many New Methods}
\label{jacsketch}

\graphicspath{{jacsketch/images/}}

In this chapter we finally consider problem \eqref{eq:finitesum} in its fully general form; i.e., we aim to solve the problem
\begin{equation}\label{eq:gjs_problem_gjs}
\compactify \min_{x\in \R^d}   \frac1n \sum \limits_{j=1}^n f_j(x) + \psi(x).
\end{equation}
We assume  that the functions $f_j:\R^d \to \R$ are smooth and convex, and $\psi:\R^d\to \R\cup \{+\infty\}$ is a proper, closed and convex regularizer, admitting a cheap proximal operator. As usual, we write $f\eqdef \frac{1}{n}\sum_j f_j$. 

\paragraph{Proximal gradient descent.} A baseline method for solving problem \eqref{eq:gjs_problem_gjs} is  {\em (proximal) gradient descent}, described in detail in Section~\ref{sec:proxgrad} of the introduction. For the sake of simplcity, let us call it {\tt PGD} throughout this section.  As already stressed, {\tt PGD} performs well when both $n$ and $d$ are not too large. However, in the big data (large $n$) and/or big parameter (large $d$) case,  the formation of the gradient becomes overly expensive, rendering {\tt PGD} inefficient in both theory and practice. A typical remedy is to replace the gradient by a cheap-to-compute random approximation. Typically, one replaces $\nabla f(x^k)$ with a random vector $g^k$ whose mean is the gradient: $\E{g^k} = \nabla f(x^k)$, i.e., with a stochastic gradient. This results in the {\em (proximal) stochastic gradient descent} ({\tt SGD}) method:
\begin{equation}\label{eq:gjs_PSGD} x^{k+1} = \prox_{\alpha \psi}( x^k - \alpha g^k).\end{equation}

Below we comment on the typical approaches to constructing $g^k$ in the big $n$ and big $d$ regimes (this was, to some extend, mentioned in the introduction already).

\paragraph{Proximal {\tt SGD}.} In the big $n$ regime, the simplest choice is to set \begin{equation} \label{eq:gjs_nbu9gff}g^k = \nabla f_j(x^k)\end{equation} for an index $j\in [n]\eqdef \{1,2,\dots,n\}$ chosen uniformly at random. By construction, it is $n$ times cheaper to compute this estimator than the gradient, which is a key driving force behind the efficiency of this variant of {\tt SGD}. However, there is an infinite array of other possibilities of constructing an unbiased estimator~\cite{needell2016batched, pmlr-v97-qian19b}. Depending on how $g^k$ is formed, \eqref{eq:gjs_PSGD} specializes to one of the many existing variants of proximal {\tt SGD}, each with different convergence properties and proofs.

\paragraph{Proximal {\tt RCD}.} In the big $d$ regime (this is interesting even if $n=1$), the simplest choice is to set \begin{equation} \label{eq:gjs_b98gf98f} g^k = d  \langle \nabla f(x^k), \eLi \rangle \eLi,\end{equation} where $\langle x,y\rangle =\sum_i x_i y_i$ is the standard Euclidean inner product, $\eLi$ is the $i$th standard unit basis vector in $\R^d$, and $i$ is chosen uniformly at random from $[d]\eqdef \{1,2,\dots,d\}$. \footnote{The algorithm proposed in this chapter subsabples both the finite sum and the domain. For the notational simplicity, we distinguish the two different spaces using color where necessary.} With this estimator, \eqref{eq:gjs_PSGD} specializes to  (proximal) randomized coordinate decent ({\tt RCD}). There are situations where it is $d$ times cheaper to compute the partial derivative $\nabla_i f(x^k)\eqdef \langle \nabla f(x^k), \eLi \rangle$ than the gradient, which is a key driving force behind the efficiency of {\tt RCD} \cite{rcdm}.  However, there is an infinite array of other possibilities for constructing an unbiased estimator of the gradient in a similar way~\cite{nsync,rcdm,qu2016coordinate2}. 

\paragraph{Issues.}  For the sake of argument in the rest of this section, assume that $f$ is a $\mu$-strongly convex function, and let $x^*$ be the (necessarily) unique solution of \eqref{eq:gjs_problem_gjs}.  It is well known that in this case, method \eqref{eq:gjs_PSGD}  with estimator $g^k$ defined as in \eqref{eq:gjs_nbu9gff} does {\em not} in general converge to $x^*$. Instead, {\tt SGD}  converges linearly to a neighborhood of $x^*$ of size proportional to the stepsize $\alpha$, noise $\sigma^2\eqdef \frac{1}{n}\sum_j \norm{\nabla f_j(x^*)}^2$, and inversely proportional to $\mu$~\cite{moulines2011non, needellward2015}.  In the generic regime with $\sigma^2>0$, the neighbourhood is nonzero, causing issues with convergence. This situation does not change even when tricks such as {\em mini-batching} or {\em importance sampling} (or a combination of both) are applied~\cite{needellward2015,needell2016batched, pmlr-v97-qian19b}. While these tricks affect both the (linear convergence) rate and the size of the neighbourhood, they are incapable\footnote{Unless, of course, in the special case when one  uses the full batch approximation $g^k = \nabla f(x^k)$.} of ensuring convergence to the solution. 

However, a remedy does exist: the situation with non-convergence can be  resolved by using one of the many {\em variance-reduction}  strategies for constructing $g^k$ developed over the last several years~\cite{sag, saga, svrg, mairal2013optimization, sdca}.  

Further, while it is well known that method \eqref{eq:gjs_PSGD}  with estimator $g^k$  defined as in \eqref{eq:gjs_b98gf98f} (i.e., randomized coordinate descent) converges to $x^*$ for $\psi\equiv 0$~\cite{rcdm, richtarik2014iteration, nsync}, it is also known that it does {\em not}  generally converge to $x^*$ unless the regularizer $\psi$ is separable (e.g., $\psi(x) = \norm{x}_1$ or $\psi(x)=c_1 \|x\|_1 + c_2 \|x\|_2^2$). In~\cite{sega}, an alternative estimator (known as {\tt SEGA}) was constructed from the same (random) partial derivative information $\nabla f_i(x^k)$, one that does not suffer from this incompatibility with general regularizers $\psi$. This work  resolved a long standing open problem in the theory of {\tt RCD} methods.


\paragraph{\bf Notation.} Let $\eR$ (resp.\ $\eL$) be the vector of all ones in $\R^n$ (resp.\ $\R^d$), and  $\eRj$ (resp.\ $\eLi$) be the $j$th (resp.\ $i$th) unit basis vector in $\R^n$ (resp.\ $\R^d$).  By $\|\cdot \|$ we denote the standard Euclidean norm in $\R^d$ and $\R^n$.  Matrices are denoted by upper-case bold letters. Given $\mX,\mY\in \R^{d\times n}$, let $\langle \mX, \mY\rangle \eqdef \Tr{\mX^\top \mY}$ and  $\|\mX\|\eqdef \langle \mX, \mX \rangle^{1/2}$ be the Frobenius norm. By $\mX_{:j}$ (resp.\ $\mX_{i:}$) we denote the $j$th column (resp.\ $i$th row) of matrix $\mX$.  By $\mI_n$ (resp.\ $\mI_d$) we denote the $n\times n$ (resp.\ $d\times d$) identity matrices. Upper-case calligraphic letters, such as $\cS,\cU, \cI, \cM, \cR$, are used to denote (deterministic or random) linear operators mapping $\R^{d\times n}$ to $\R^{d\times n}$. Most used notation is summarized in Table~\ref{tbl:notation_jacsketch}  in Appendix~\ref{sec:table}.

\section{Contributions} 

Having experienced a ``Cambrian explosion'' in the last 10 years, the world of efficient {\tt SGD} methods is remarkably complex. There is a large and growing set of rules for constructing the gradient estimators $g^k$, with differing levels of sophistication and varying theoretical and practical properties. It includes the classical estimator \eqref{eq:gjs_nbu9gff}, as well as an infinite array of mini-batch \cite{li2014efficient} and importance sampling \cite{needellward2015, iprox-sdca} variants, and a growing list of variance-reduced variants~\cite{saga}. Furthermore, there are estimators of the coordinate descent variety, including the simplest one based on \eqref{eq:gjs_b98gf98f}~\cite{rcdm}, more elaborate variants utilizing the arbitrary sampling paradigm \cite{qu2016coordinate1}, and variance reduced methods capable of handling general non-separable regularizers~\cite{sega}.

\begin{itemize}

\item {\bf New general method and a single convergence theorem}. In this chapter we propose a {\em general method}---which we call {\tt GJS}---which reduces to many of the aforementioned classical and several recently developed {\tt SGD} type methods in special cases.  We provide a {\em single convergence theorem}, establishing a linear convergence rate for {\tt GJC}, assuming $f$ to be smooth and quasi strongly convex. In particular, we obtain  the following methods in special cases, or their generalizations, always recovering the best-known convergence guarantees or improving upon them: {\tt SAGA}~\cite{saga, qian2019saga, gazagnadou2019optimal}, {\tt JacSketch}~\cite{jacsketch},  {\tt LSVRG}~\cite{hofmann2015variance, kovalev2019don}, {\tt SEGA}~\cite{sega},  and {\tt ISEGA}~\cite{mishchenko201999}  (see Table~\ref{tbl:gjs_all_special_cases}, in which we list 17 special cases).  This is the first time such a direct connection is made between many of these methods, which previously required different intuitions and dedicated analyses. Our general method, and hence also all special cases we consider, can work with a regularizer. This provides novel (although not hard) results for some methods, such as {\tt LSVRG}.

\item {\bf Unification of {\tt SGD} and {\tt RCD}.} As a by-product of the generality of {\tt GJS}, we obtain the {\em unification of variance-reduced {\tt SGD} and variance reduced {\tt RCD} methods.} To the best of our knowledge, there is no algorithm besides {\tt GJS}, one whose complexity is captured by a single theorem, which specializes to {\tt SGD} and {\tt RCD} type methods at the same time and recovers best known rates in both cases.\footnote{A single theorem (not a single algorithm) to obtain rates for both variance-reduced {\tt SGD} and variance reduced {\tt RCD} methods was done in the concurrent work~\cite{sigma_k}. However,~\cite{sigma_k} focuses in orthogonal direction instead -- it is a tool to analyze stochastic gradient algorithms which includes non-variance reduced methods as well.}

\item {\bf Generalizations to arbitrary sampling.}  Many specialized methods we develop are cast in a very general {\em arbitrary sampling} paradigm~\cite{nsync, quartz, qu2016coordinate1}, which allows for the estimator $g^k$ to be formed through information contained in a random subset $R^k\subseteq [n]$ (by computing $\nabla f_j(x^k)$ for $j\in \R^k$) or a random subset $L^k\subseteq [d]$ (by computing $\nabla_i f(x^k)$ for $i\in L^k$), where these subsets are allowed to follow an arbitrary distribution. In particular, we  extend {\tt SEGA}~\cite{sega}, {\tt LSVRG}~\cite{hofmann2015variance, kovalev2019don} or {\tt ISEGA}~\cite{mishchenko201999} to this setup.  Likewise, {\tt GJS} specializes to an arbitrary sampling extension of the {\tt SGD}-type method {\tt SAGA}~\cite{saga, qian2019saga}, obtaining state-of-the-art rates. As a special case of the arbitrary sampling paradigm, we obtain \emph{importance sampling} versions of all mentioned methods.

\item {\bf New methods.} {\tt GJS} can be specialized to many new specific methods. To illustrate this, we construct 10 specific {\em new} methods in special cases, some with intriguing structure and properties (see Section~\ref{sec:gjs_Special_Cases}; Table~\ref{tbl:gjs_all_special_cases}; and Table~\ref{tbl:gjs_all_special_cases_theory} for a summary of the rates). 


\item {\bf Relation to {\tt JacSketch}.} Our method can be seen as a vast generalization of the recently proposed Jacobian sketching method {\tt JacSketch}~\cite{jacsketch} in several directions, notably by enabling {\em arbitrary randomized linear} (i.e., sketching) operators, allowing different linear operators to learning Jacobian and constructing control variates, extending the analysis to the proximal case, and replacing strong convexity assumption by quasi strong convexity or strong growth (see Appendix~\ref{sec:gjs_sg}). In particular, from all methods we recover, only variants of {\tt SAGA} can be obtained from {\tt JacSketch}~\cite{jacsketch} (even in that case, rates obtained from~\cite{jacsketch} are suboptimal).

\item {\bf Limitations.} 
We  focus on developing methods capable of enjoying a linear convergence rate with a fixed stepsize $\alpha$ and do not consider the non-convex setting.  Although there exist several {\em accelerated} variance reduced algorithms~\cite{lan2018optimal, allen2017katyusha, zhou2018simple, zhou2018direct, kovalev2019don, kulunchakov2019estimate}, we do not consider such methods here. 

\end{itemize}

\section{Sketching}

A key object in this chapter is the Jacobian matrix $\mG(x) = [\nabla f_1(x),\dots, \nabla f_n(x)] \in \R^{d\times n}.$
Note that  \begin{equation}\label{eq:gjs_nbifg98dz}
\compactify \nabla f(x) = \frac{1}{n}\mG(x) \eR.\end{equation} 
Extending the insights from \cite{jacsketch}, one of the key observations of this work is that {\em random linear transformations} (sketches) of $\mG$ can be used to {\em construct} unbiased estimators of the gradient of $f$.  For instance,  $\mG(x^k) \eRj$  leads to the simple {\tt SGD} estimator \eqref{eq:gjs_nbu9gff},  and $\frac{d}{n} \eLi \eLi^\top \mG(x^k) \eR $ gives the simple {\tt RCD} estimator \eqref{eq:gjs_b98gf98f}. We will consider more elaborate examples later on.  It will be useful to embed these estimators into $\R^{d\times n}$. For instance, instead of $\mG(x^k)\eRj$ we consider the matrix  $\mG(x^k) \eRj \eRj^\top$. Note that all columns of this matrix are zero, except for the $j$th column, which is equal to  $\mG(x^k)\eRj$. Similarly, instead of $\frac{d}{n} \eLi \eLi^\top \mG(x^k) \eR $ we will consider the matrix $\frac{d}{n} \eLi \eLi^\top \mG(x^k)$. All rows of this matrix are zero, except for the $i$th row, which consists of the $i$th partial derivatives of functions $f_j(x^k)$ for $j\in [n]$, scaled by $\frac{d}{n}$.

\paragraph{Random projections.}
 Generalizing from these examples, we  consider a random linear operator (``sketch'') $\cA:\R^{d\times n}\to \R^{d\times n}$. By $\cA^\ast$ we denote the adjoint of $\cA$, i.e., linear operator  satisfying $\langle \cA \mX, \mY \rangle = \langle \mX, \cA^\ast \mY \rangle$ for all $\mX,\mY\in \R^{d\times n}$.  Given $\cA$, we let $\cP_{\cA}$ be the (random) projection operator onto $\Range{\cA^\ast}$. That is, \[\cP_{\cA}(\mX) = \arg\min_{\mY \in \Range{\cA^\ast}} \norm{\mX - \mY } = \cA^\ast (\cA \cA^\ast)^\dagger \cA \mX,\]
where ${}^\dagger$ is the Moore-Penrose pseudoinverse. The identity operator is denoted by $\cI$. We say that $\cA$ is {\em identity in expectation}, or {\em unbiased} when $\E{\cA} = \cI$; i.e., when if $\E{\cA  \mX} =\mX$ for all $\mX\in \R^{d\times n}$.

\begin{definition} We will often consider the following\footnote{The algorithm we develop is, however, not limited to such sketches.} sketching operators $\cA$:
\begin{itemize}
\item [(i)] {\bf Right sketch.} Let $\mR\in \R^{n\times n}$ be a random matrix. Define $\cA$ by $\cA \mX = \mX \mR$ (``R-sketch''). Notice that $\cA^\ast \mX = \mX \mR^\top$. In particular, if $R$ is random subset of $[n]$, we can define $\mR = \sum_{j \in R} \eRj \eRj^\top$. The resulting operator $\cA$ (``R-sampling'') satisfies: $\cA = \cA^\ast = \cA^2 = \cP_{\cA}$.  If we let $\pRj \eqdef \Probbb{j \in R}$, and instead define $\mR = \sum_{j \in R} \frac{1}{\pRj} \eRj \eRj^\top$, then $\E{\mR} = \mI_n$ and hence  $\cA$ is unbiased.\\
\item [(ii)] {\bf Left sketch.} Let $\mL\in \R^{d\times d}$ be a random matrix. Define $\cA$ by $\cA \mX = \mL \mX $ (``L-sketch'').  Notice that $\cA^\ast \mX = \mL^\top \mX $. In particular, if $L$ is random subset of $[d]$, we can define $\mL = \sum_{i\in L} \eLi \eLi^\top$. The resulting operator $\cA$ (``L-sampling'') satisfies: $\cA = \cA^\ast = \cA^2 = \cP_{\cA}$. If we let $\pLi\eqdef \Probbb{i\in L}$, and instead define $\mL = \sum_{i\in L} \frac{1}{\pLi} \eLi \eLi^\top$, then $\E{\mL} = \mI_d$  and hence $\cA$  us unbiased.\\
\item [(iii)] {\bf Scaling/Bernoulli.} Let $\xi$ be a Bernoulli random variable, i.e., $\xi = 1$ with probability $\probx $ and $\xi = 0$ with probability $1-\probx $, where $\probx \in [0,1]$.  Define $\cA$ by $\cA \mX = \xi \mX$ (``scaling''). Then $\cA = \cA^\ast = \cA^2 = \cP_{\cA}$. If we instead define $\cA \mX = \frac{1}{\probx } \xi \mX$, then  $\cA$ is unbiased.\\
\item [(iv)]   {\bf LR sketch.}  All the above operators can be combined. In particular, we can define $\cA \mX = \xi \mL \mX \mR$. All of the above arise as special cases of this: (i) arises for $\xi\equiv 1$ and $\mL\equiv \mI_d$, (ii) for  $\xi\equiv 1$ and $\mR\equiv \mI_n$, and (iii) for $\mL\equiv \mI_d$ and $\mR\equiv \mI_n$.
\end{itemize}
\end{definition}


\section{The {\tt GJS} algorithm}
We are now ready to describe our method (formalized as Algorithm~\ref{alg:gjs_SketchJac}). 
\begin{algorithm}[!h]
\begin{algorithmic}[1]
\State \textbf{Parameters:} Stepsize $\alpha>0$, random projector $\cS$ and unbiased sketch $\cU$
\State \textbf{Initialization:} Choose  solution estimate $x^0 \in \R^d$ and Jacobian estimate $ \mJ^0\in \R^{d\times n}$ 
\For{$k =  0, 1,2, \dots$}
\State Sample realizations of $\cS$ and $\cU$, and perform sketches $\cS\mG(x^k)$ and $\cU\mG(x^k)$
\State  $\mJ^{k+1} = \mJ^k - \cS(\mJ^k - \mG(x^k))$ \quad \hfill update the Jacobian estimate via  \eqref{eq:gjs_nio9h8fbds79kjh}
\State $g^k = \frac1n \mJ^k \eR + \frac1n \cU \left(\mG(x^k) -\mJ^k\right)\eR$  \hfill construct the gradient estimator via \eqref{eq:gjs_ni98hffs}
    \State $x^{k+1} = \prox_{\alpha \psi} (x^k - \alpha g^k)$ \label{eq:gjs_alg_update} \hfill perform the proximal {\tt SGD} step \eqref{eq:gjs_PSGD}
\EndFor
\end{algorithmic}
\caption{Generalized {\tt JacSketch} ({\tt GJS}) }
\label{alg:gjs_SketchJac}
\end{algorithm}
Let $\cS$ be a random linear operator (e.g., right sketch, left sketch, or scaling)  such that $\cS = \cP_\cS$ and let $\cU$ be an unbiased operator. We propose to construct the gradient estimator as
\begin{equation}\label{eq:gjs_ni98hffs}
\compactify g^k = \frac{1}{n} \mJ^k \eR + \frac{1}{n}\cU (\mG(x^k) - \mJ^k) \eR,
\end{equation}
where the matrices $\mJ^k\in \R^{d\times n}$ are constructed iteratively. Note that, taking expectation in $\cU$, we get 
\begin{equation}\label{eq:gjs_unbiased_xx}
\compactify \E{g^k} \overset{\eqref{eq:gjs_ni98hffs}}{=}\frac{1}{n} \mJ^k \eR +  \frac{1}{n}(\mG(x^k) - \mJ^k) \eR = \frac{1}{n} \mG(x^k)e \overset{\eqref{eq:gjs_nbifg98dz}}{=} \nabla f(x^k),\end{equation}
and hence $g^k$ is indeed unbiased. We will construct $\mJ^k$ so that  $\mJ^k\to \mG(x^*)$. By doing so, the variance of $g^k$ decreases throughout the iterations, completely vanishing at $x^*$. The sequence $\{\mJ^k\}$ is updated as follows:
\begin{equation}\label{eq:gjs_nio9h8fbds79kjh}\mJ^{k+1} = \arg\min_{\mJ} \left\{ \|\mJ - \mJ^k\| \;:\; \cS \mJ = \cS \mG(x^k) \right\} = \mJ^k - \cS(\mJ^k - \mG(x^k)).\end{equation}
That is, we sketch the Jacobian $\mG(x^k)$, obtaining the sketch $\cS \mG(x^k)$, and seek to use this information to construct a new matrix $\mJ^{k+1}$ which is consistent with this sketch, and as close to $\mJ^k$ as possible. The intuition here is as follows: if we repeated the sketch-and-project process \eqref{eq:gjs_nio9h8fbds79kjh} for fixed $x^k$, the matrices $\mJ^k$ would converge to $\mG(x^k)$, at a linear rate~\cite{gower2015randomized, gower:2017}. This process can be seen as {\tt SGD} applied to a certain quadratic stochastic optimization problem~\cite{richtarik2017stochastic, jacsketch}. Instead, we take just one step of this iterative process, change $x^k$, and repeat. Note that the unbiased sketch $\cU$ in \eqref{eq:gjs_ni98hffs} also claims access to $\mG(x^k)$. Specific variants of {\tt GJS} are obtained by choosing specific operators $\cS$ and $\cU$ (see Section~\ref{sec:gjs_Special_Cases}).



\section{Theory}

We now describe the main result of this chapter, which depends on a relaxed strong convexity assumption and a more precise smoothness assumption on $f$.

\begin{assumption} \label{as:gjs_smooth_strongly_convex}
Problem \eqref{eq:gjs_problem_gjs} has a unique minimizer $x^*$, and $f$ is $\mu$-quasi strongly convex, i.e., 
\begin{equation} 
\compactify f(x^*) \geq   f(y) + \< \nabla f(y) ,x^*- y> + \frac{\mu}{2} \norm{y - x^*}^2, \quad \forall y\in \R^d, \label{eq:gjs_strconv3}
\end{equation}
Functions $f_j$ are convex and $\mM_j$-smooth for some $\mM_j\succeq 0$, i.e.,
\begin{equation} 
\compactify f_j(y) + \< \nabla f_j(y) ,x- y> \leq f_j(x)\leq    f_j(y) + \< \nabla f_j(y) ,x- y> +\frac{1}{2} \norm{y - x}_{{\mM_j}}^2, \quad \forall x,y\in \R^d.  \label{eq:gjs_smooth_ass}
\end{equation}
\end{assumption}

Assumption~\ref{eq:gjs_smooth_ass} generalizes classical $L$-smoothness, which is obtained in the special case $\mM_j=L\mI_d$. The usefulness of this assumption comes from i) the fact that ERM problems typically satisfy \eqref{eq:gjs_smooth_ass} in a non-trivial way~\cite{qu2016coordinate2, pmlr-v97-qian19b}, ii) our method is able to utilize the full information contained in these matrices for further acceleration (via increased stepsizes). Given matrices $\{\mM_{j}\}$ from Assumption~\ref{as:gjs_smooth_strongly_convex}, let $\cM$ be the linear operator defined via
$\left(\cM \mX\right)_{:j} = \mM_j \mX_{:j}$ for $j\in [n]$. It is easy to check that this operator is self-adjoint and positive semi-definite, and that its square root is given by $$\left(\cM^{\frac{1}{2}}  \mX\right)_{:j} = \mM_j^{\frac{1}{2}} \mX_{:j}.$$ The pseudoinverse $\cM^\dagger$ of this operator  plays an important role in our  main result.

\begin{theorem}\label{thm:gjs_main} Let Assumption~\ref{as:gjs_smooth_strongly_convex} hold. Let $\cB$ be any linear operator commuting with $\cS$, and assume ${\cM^\dagger}^{\frac{1}{2}}$ commutes with $\cS$. Let $\cR$ be any linear operator for which $\cR(\mJ^k) = \cR(\mG(x^*))$ for every $k\geq 0$. 
Define the Lyapunov function 
\begin{eqnarray}\label{eq:gjs_Lyapunov}
\Psi^k & \eqdef & \norm{ x^k - x^* }^2 + \alpha \NORMG{ \cB {\cM^\dagger}^{\frac12} \left( \mJ^k - \mG(x^*)\right)},
\end{eqnarray}
where $\{x^k\}$ and $\{\mJ^k\}$ are the random iterates produced by Algorithm~\ref{alg:gjs_SketchJac} with stepsize $\alpha>0$. Suppose that $\alpha$ and $\cB$ are chosen so that 
\begin{eqnarray}\label{eq:gjs_small_step}
\compactify \frac{2\alpha}{n^2}  \E{ \norm{ \cU  \mX \eR }^2 }   +   \NORMG{  \left(\cI - \E{\cS} \right)^{\frac12}\cB  {\cM^\dagger}^{\frac12} \mX }  &\leq & \compactify (1-\alpha \mu) \NORMG{ \cB {\cM^\dagger}^{\frac12}\mX }
\end{eqnarray}
whenever  $ \mX\in \Range{\cR}^\perp$ and
\begin{eqnarray}
\compactify \frac{2\alpha}{n^2} \E{  \norm{ \cU  \mX  \eR }^2  } 
+   \NORMG{\left(\E{\cS}\right)^{\frac12}  \cB  {\cM^\dagger}^{\frac12}\mX  }   &  \leq & \compactify \frac{1}{n} \norm{{\cM^\dagger}^{\frac12}\mX}^2.\label{eq:gjs_small_step2}
\end{eqnarray}
for all $\mX\in \R^{d\times n}$.  Then  for all $k\geq 0$, we have
$\E{\Psi^{k}}\leq \left( 1-\alpha\mu\right)^k \Psi^{0}.$
\end{theorem}

The above theorem is very general as it applies to essentially arbitrary random linear operators $\cS$ and $\cU$. It postulates a linear convergence rate of a Lyapunov function composed of two terms: distance of $x^k$ from $x^*$, and weighted distance of the Jacobian $\mJ^k$ from $\mG(x^*)$. Hence, we obtain convergence of both the iterates and the Jacobian to $x^*$ and $\mG(x^*)$, respectively. 
Inequalities \eqref{eq:gjs_small_step} and \eqref{eq:gjs_small_step2} are mainly assumptions one stepsize $\alpha$, and are used to define suitable weight operator $\cB$. See Lemma~\ref{lem:gjs_existence} for a general statement on when these inequalities are satisfied. However, we give concrete and simple answers in all special cases of {\tt GJS} in the appendix. For a summary of how the operator $\cB$ is chosen in special cases, and the particular complexity results derived from this theorem, we refer to Table~\ref{tbl:gjs_all_special_cases_theory}.

\begin{remark}We use the trivial choice $\cR\equiv 0$ in almost all special cases. With this choice of $\cR$,  the condition $\cR(\mJ^k) = \cR(\mG(x^*))$ is automatically satisfied, and inequality \eqref{eq:gjs_small_step2} is requested to hold for {\em all} matrices $\mX\in \R^{d\times n}$. However, a non-trivial choice of $\cR$ is sometimes useful;  e.g., in the analysis of a subspace variant of {\tt SEGA}~\cite{sega}.  Further, the results of Theorem~\ref{thm:gjs_main} can be generalized from a quasi strong convexity to a strong growth condition \cite{karimi2016linear}  on $f$~(see Appendix~\ref{sec:gjs_sg}). While interesting, these are not the key results of this work and we therefore suppress them to the appendix. 
\end{remark}


\section{Special cases} \label{sec:gjs_Special_Cases}

As outlined in the introduction, {\tt GJS} (Algorithm~\ref{alg:gjs_SketchJac}) is a surprisingly versatile method. In Table~\ref{tbl:gjs_all_special_cases} we list {\em 7 existing methods }(in some cases, generalizations of existing methods), and construct also {\em 10 new variance reduced methods.} We also provide a summary of all specialized iteration complexity results, and a guide to the corollaries which state them (see Table~\ref{tbl:gjs_all_special_cases_theory} in the appendix).

\begin{table}[!t]
\begin{center}
\tiny
\begin{tabular}{|c|c|c|c|c|c|}
\hline
\multicolumn{2}{|c|}{\bf  Choice of random operators $\cS$ and $\cU$ defining Algorithm~\ref{alg:gjs_SketchJac}} &  \multicolumn{4}{|c|}{\bf Algorithm} \\
\hline
$\cS \mX $ & $ \cU \mX $  &  \#  &  Name  & Comment & Sec.  \\
\hline
\hline
$\mX \eRj \eRj^\top$ \text{w.p.}\; $\pRj = \frac1n$ & $\mX n \eRj \eRj^\top $ \text{w.p.}\; $\pRj = \frac1n$ & \ref{alg:gjs_SAGA} &  {\tt SAGA} & basic variant of {\tt SAGA} \cite{saga} & \ref{sec:gjs_saga_basic}  \\
\hline
 $\mX \sum \limits_{j\in R} \eRj \eRj^\top$ \text{w.p.}\; $\pRR$ & $\mX \sum \limits_{j\in R} \frac{1}{\pRj} \eRj \eRj^\top$ \text{w.p.}\; $\pRR$ &  \ref{alg:gjs_SAGA_AS_ESO} & {\tt SAGA} & 
 {\tt SAGA} with AS \cite{qian2019saga}   &  \ref{sec:gjs_saga_as} \\
\hline
  $\eLi \eLi^\top \mX $ \text{w.p.}\; $\pLi = \frac{1}{d}$
 & $d \eLi \eLi^\top  \mX $ \text{w.p.}\; $\pLi = \frac{1}{d}$ &   \ref{alg:gjs_SEGA} & {\tt SEGA} & basic variant  of  {\tt SEGA} \cite{sega} & \ref{sec:gjs_sega}  \\
\hline
 $\sum \limits_{i\in L}\eLi \eLi^\top \mX $ \text{w.p.}\; $\pLL$
 & $\sum \limits_{i\in L} \frac{1}{\pLi} \eLi \eLi^\top  \mX $ \text{w.p.}\; $\pLL$ &   \ref{alg:gjs_SEGAAS} & {\tt SEGA} &   {\tt SEGA}  \cite{sega} with AS and prox & \ref{sec:gjs_sega_is_v1}  \\
\hline
 $
=\begin{cases}
    0  & \text{w.p.}\; 1- \probx \\
    \mX              & \text{w.p.}\;  \probx
\end{cases} $
 & $\sum \limits_{i\in L} \frac{1}{\pLi} \eLi \eLi^\top \mX $ \text{w.p.}\; $\pLL$ & \ref{alg:gjs_SVRCD} & {\tt SVRCD} & {\bf NEW}  & \ref{sec:gjs_svrcd_is2} \\
\hline
 0 & $\mX \sum \limits_{j\in R} \frac{1}{\pRj} \eRj \eRj^\top$ \text{w.p.}\; $\pRR$ &  \ref{alg:gjs_SGD_AS} & {\tt SGD-star} & {\tt SGD-star}~\cite{sigma_k} with AS & \ref{sec:gjs_SGD-AS-star}  \\
\hline
 $
=\begin{cases}
    0  & \text{w.p.}\; 1- \probx \\
    \mX              & \text{w.p.}\;  \probx
\end{cases} $
 & $\mX \sum \limits_{j\in R} \frac{1}{\pRj} \eRj \eRj^\top$ \text{w.p.}\; $\pRR$ &   \ref{alg:gjs_LSVRG-AS} & {\tt LSVRG}  & {\tt LSVRG} \cite{kovalev2019don} with AS and prox & \ref{sec:gjs_LSVRG-AS}  \\
\hline
$
=\begin{cases}
    0  & \text{w.p.}\; 1- \probx \\
    \mX              & \text{w.p.}\;  \probx
\end{cases} $
 &   
 $=\begin{cases}
    0  & \text{w.p.}\; 1- \proby \\
   \frac{1}{\proby} \mX              & \text{w.p.}\;  \proby
\end{cases} $ & \ref{alg:gjs_B2} & {\tt B2} & {\bf NEW} & \ref{sec:gjs_B2}  \\
\hline
$\mX \sum \limits_{j\in R} \eRj \eRj^\top  $ \text{w.p.}\; $\pRR$
 &   
 $=\begin{cases}
    0  & \text{w.p.}\; 1- \proby \\
   \frac{1}{\proby} \mX              & \text{w.p.}\;  \proby
\end{cases} $ &   \ref{alg:gjs_invsvrg}  & {\tt LSVRG-inv} &  {\bf NEW} & \ref{sec:gjs_SVRG-1}  \\
\hline
$\sum \limits_{i\in L} \eLi \eLi^\top \mX $ \text{w.p.}\; $\pLL$
 &   
 $=\begin{cases}
    0  & \text{w.p.}\; 1- \proby \\
  \frac{1}{\proby} \mX              & \text{w.p.}\;  \proby
\end{cases} $ & \ref{alg:gjs_B_sega} & {\tt SVRCD-inv} & {\bf NEW} &  \ref{sec:gjs_SVRCD-inv}\\
\hline
$\mX \sum \limits_{j\in R} \eRj \eRj^\top$ \text{w.p.}\;  $\pRR
$ 
& 
$\sum \limits_{i\in L} \frac{1}{\pLi} \eLi  \eLi^\top \mX$  \text{w.p.}\;  $\pLL
$ 
& \ref{alg:gjs_RL} & {\tt RL} & {\bf NEW} &  \ref{sec:gjs_RL} \\
\hline
$\sum \limits_{i\in L} \eLi  \eLi^\top \mX$ \text{w.p.}\;  $\pLL
$ 
& 
 $\mX \sum \limits_{j\in R} \frac{1}{\pRj} \eRj \eRj^\top$ \text{w.p.}\;  $\pRR
$ 
& \ref{alg:gjs_LR} & {\tt LR} & {\bf NEW} &  \ref{sec:gjs_LR} \\
\hline
$ \mI_{L:} \mX \mI_{:R}$ \text{w.p.}\; $\pLL \pRR$
 &   
 $  \mI_{L:}\left( \left(    \pL^{-1} \left(\pR^{-1}\right)^\top  \right) \circ \mX\right)\mI_{:R}$ \text{w.p.}\; $\pLL \pRR$& \ref{alg:gjs_saega} & {\tt SAEGA} & {\bf NEW} &  \ref{sec:gjs_SAEGA} \\
 \hline
$=\begin{cases}
0 & \text{w.p.}\; 1-\probx  \\
\mX & \text{w.p.}\; \probx \\
\end{cases}$
 &   
 $  \mI_{L:}\left( \left( \pL^{-1} \left(\pR^{-1}\right)^\top\right) \circ \mX\right)\mI_{:R}$ \text{w.p.}\; $\pLL \pRR$&  \ref{alg:gjs_svrcdg} & {\tt SVRCDG} & {\bf NEW} & \ref{sec:gjs_SVRCDG}  \\
 \hline
  $\sum \limits_{\tR=1}^\TR  \mI_{L_{\tR}:}\mX_{:\NRt} \mI_{:R_{\tR}} $
  &   
$\sum \limits_{\tR=1}^\TR   \left( {(\ptL)^{-1}} {(\ptR)^{-1}}^\top\right) \circ \left(\mI_{L_{\tR}:}\mX_{:\NRt}  \mI_{:R_{\tR}}\right)$ &  \ref{alg:gjs_isaega} & {\tt ISAEGA} & {\bf NEW} (reminiscent of \cite{mishchenko201999}) &  \ref{sec:gjs_ISAEGA}  \\
\hline
  $\sum \limits_{\tR=1}^\TR  \mI_{L_{\tR}:}\mX_{:\NRt} $ 
  &   
$\sum \limits_{\tR=1}^\TR   \left( {(\ptL)^{-1}} {\eR}^\top\right) \circ \left(\mI_{L_{\tR}:}\mX_{:\NRt}  \right)$ &  \ref{alg:gjs_isega} & {\tt ISEGA} &   {\tt ISEGA} \cite{mishchenko201999} with AS&  \ref{sec:gjs_ISAEGA}  \\
\hline
  $\mX \mR $  & $\mX \mR \E{\mR}^{-1} $ &   \ref{alg:gjs_jacsketch} & {\tt JS}  & {\tt JacSketch} \cite{jacsketch}  with AS and prox & \ref{sec:gjs_jacsketch}  \\
\hline
\end{tabular}
\end{center}
\caption{Selected special cases of {\tt GJS} (Algorithm~\ref{alg:gjs_SketchJac}) arising by choosing operators $\cS$ and $\cU$ in particular ways. $R$ is a random subset of $[n]$, $L$ is a random subset of $[d]$, $\pLi = \Probbb{i\in L}$, $\pRj = \Probbb{j\in R}$.}
\label{tbl:gjs_all_special_cases}
\end{table}

\begin{itemize}

\item  {\tt SGD-star.} In order to illustrate why variance reduction is needed in the first place, let us start by describing one of the methods---{\tt SGD-star}  (Algorithm~\ref{alg:gjs_SGD_AS})---which happens to be particularly suitable to shed light on this issue. In {\tt SGD-star} we assume that the Jacobian at optimum, $\mG(x^*)$, is known. While this is clearly an unrealistic assumption, let us see where it leads us.
If this is the case, we can choose $\mJ^0 = \mG(x^*)$, and  let $\cS \equiv 0$. This implies that $\mJ^k=\mJ^0$ for all $k$.  We then choose $\cU$ to be the right unbiased sampling operator, i.e.,  $\cU \mX =  \mX  \sum_{j\in R} \frac{1}{\pRj} \eRj \eRj^\top$, which gives
\[\compactify g^k = \frac{1}{n}\sum \limits_{j=1}^n \nabla f_j(x^*)  +\sum \limits_{j\in R^k} \frac{1}{n \pRj} \left( \nabla f_j(x^k) -\nabla f_j(x^*) \right).\]
This method does not need to learn the Jacobian at $x^*$ as it is known, and instead moves in a direction of average  gradient at the optimum, perturbed by a random estimator of the direction $\nabla f(x^k) - \nabla f(x^*)$ formed via sub-sampling $j\in R^k\subseteq [n]$ . What is special about this perturbation? As the method converges, $x^k\to x^*$ and the perturbations converge to zero, for {\em any} realization of the random set $R^k\subseteq [n]$. So,  gradient estimation stabilizes, we get $g^k\to \nabla f(x^*)$, and hence the variance of $g^k$ converges to zero. In view of Corollary~\ref{cor:gjs_sgd} of our main result (Theorem~\ref{thm:gjs_main}), the iteration complexity of {\tt SGD-star} is $\max_j  \frac{v_j}{\mu  n \pRj} \log \frac{1}{\epsilon}$, where $\mu$ is the quasi strong convexity parameter of $f$, and the smoothness constants $v_j$ are defined in Appendix~\ref{sec:gjs_SGD-AS-star}.

Since knowing  $\mG(x^*)$ is unrealistic,  {\tt GJS} is instead  {\em learning} these perturbations on the fly. Different variants of {\tt GJS} do this differently, but ultimately all attempt to learn the gradients $\nabla f_j(x^*)$ and use this information to stabilize the gradient estimation. Due to space restrictions, we do not describe all remaining 9 new methods in the main body of the chapter, let alone the all 17 methods. We will briefly outline 2 more (not necessarily the most interesting) new methods:

\item
 {\tt SVRCD.} This method belongs to the {\tt RCD} variety, and constructs the gradient estimator via the rule
\[\compactify g^k = h^k+\sum \limits_{i\in L^k}  \frac{1}{\pLi}(\nabla_i f(x^k) - h_i^k) \eLi\; ,\]
where $L^k \subseteq [d]$  is sampled afresh in each iteration.
The auxiliary vector $h^k$ is updated using a simple biased coin flip: $h^{k+1} =   h^k$ with probability $1- \probx$, and   $h^{k+1} =  \nabla f(x^k) $        with probability $\probx$. So, a full pass over all coordinates is made in each iteration with probability $\probx$, and a partial derivatives $\nabla_i f(x^k)$ for $i\in L^k$  are computed in each iteration. This method has  a similar structure to {\tt LSVRG}, which instead sub-sampling coordinates sub-samples functions $f_j$ for $j\in R^k$ (see Table~\ref{tbl:gjs_all_special_cases}). The iteration complexity of this method is $\left(\frac{1}{\probx} + \max_i \frac{1}{\pLi}\frac{4 m_i}{\mu}\right)\log \frac{1}{\epsilon} $, where $m_i$ is a smoothness parameter of $f$ associated with coordinate $i$ (see Table~\ref{tbl:gjs_all_special_cases_theory} and Corollary~\ref{cor:gjs_svrcd}).

\item  {\tt ISAEGA.} In Chapter~\ref{99}, a strategy of running {\tt RCD} on top of a parallel implementation  of optimization algorithms such as {\tt PGD}, {\tt SGD} or {\tt SAGA} was proposed. Surprisingly, it was shown that the runtime of the overall algorithm is unaffected whether one computes and communicates {\em all entries} of the stochastic gradient on each worker, or only a {\em fraction} of all entries of size  inversely proportional to the number of all workers. However, {\tt ISAGA}~\cite{mishchenko201999} (distributed {\tt SAGA} with {\tt RCD} on top of it), as proposed, requires the gradients with respect to the data owned by a given machine to be zero at the optimum. On the other hand, {\tt ISEGA}~\cite{mishchenko201999} does not have the issue, but it requires a computation of the exact partial derivatives on each machine and thus is expensive. As a special case of {\tt GJS} we propose {\tt ISAEGA} -- a method which cherry-picks the best properties from both {\tt ISAGA} (allowing for stochastic partial derivatives) and {\tt ISEGA} (not requiring zero gradients at the optimum). Further, we present the method in the arbitrary sampling paradigm. See Appendix~\ref{sec:gjs_ISAEGA} for more details.

\end{itemize}

\section{Experiments}
We perform extensive numerical testing for various special cases of Algorithm~\ref{alg:gjs_SketchJac}. We first start with perfectly understood example -- minimizing artificial quadratics. After that, we present experiments on logistic regression with real-world data.

\subsection{{\tt SEGA} and {\tt SVRCD} with importance sampling \label{sec:gjs_sega_exp}}
In Sections~\ref{sec:gjs_sega_is_v1} and~\ref{sec:gjs_svrcd_is2} we develop an arbitrary (and thus importance in special case) sampling for {\tt SEGA}, as well as new method {\tt SVRCD} with arbitrary sampling.  In this experiment, we compare them to its natural competitors -- basic {\tt SEGA} from~\cite{sega} and proximal gradient descent. 

Consider artificial quadratic minimization with regularizer $\psi$ being an indicator of the unit ball\footnote{In such case, proixmal operator of $\psi$ becomes a projection onto the unit ball.}:
\[ f(x) = x^\top \mM x  - b^\top x, \quad \psi(x) = \begin{cases} x& 0\leq 1 \\ \infty & \| x\| > 1\end{cases} .
\]
 Specific choices of $\mM,b$ are given by by Table~\ref{tbl:gjs_quadratics}. As both {\tt SEGA} and {\tt SVRCD} (from Section~\ref{sec:gjs_sega_is_v1} and~\ref{sec:gjs_svrcd_is2}) require a diagonal smoothness matrix, we shall further consider vector $m$ such that the upped bound $\mM\preceq \diag(m)$ holds. As the choice of $m$ is not unique, we shall choose the one which minimizes $\sum_{i=1}^d m_i$ for importance sampling and $m = \lambda_{\max}(\mM) \eL$  for uniform. Further, stepsize $\gamma = \frac{1}{4\sum_{i=1}^d m_i}$ was chosen in each case.
 Figure~\ref{fig:gjs_sega_cmp} shows the results of this experiment. As theory suggests, importance sampling for both {\tt SEGA} and {\tt SVRCD} outperform both plain {\tt SEGA} and proximal gradient always. The performance difference depends on the data; the closer $\mM$ is to a diagonal matrix with non-uniform elements, the larger stronger is the effect of importance sampling.

 \begin{table}[!h]
\begin{center}
\begin{tabular}{|c|c|c|}
\hline
Type & $\mM $ & $b$ \\
 \hline
 1   & $\diag\left(1.3^{[d]}\right)$ & $\gamma u$  \\
 \hline
  2   & $\diag((d,1,1,\dots, 1))$ & $\gamma u$ \\
\hline
3   & $\diag\left(1.1^{[d]}\right)+\mN \mN^\top \frac{1.1^d}{1000 d}$, $\mN~\sim N(0,\mI)$ & $\gamma  u$\\
\hline
4   & $\mN \mN^\top$, $\mN~\sim N(0, \mI)$ & $\gamma u$ \\
\hline
\end{tabular}
\end{center}
\caption{Four types of quadratic problems. We choose $u \sim N(0,\mI_d)$, and  $\gamma$ to be such that $\|\gamma \mM^{-1} u\| = \frac32$. Notation $c^{[d]}$ stands for a vector $(c, c^2, \dots c^d)$.}
\label{tbl:gjs_quadratics}
\end{table}

\begin{figure}[!h]
\centering
\begin{minipage}{0.25\textwidth}
  \centering
\includegraphics[width =  \textwidth ]{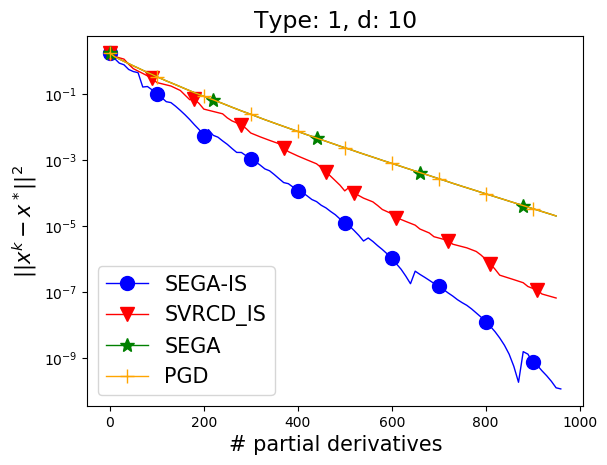}
\end{minipage}%
\begin{minipage}{0.25\textwidth}
  \centering
\includegraphics[width =  \textwidth ]{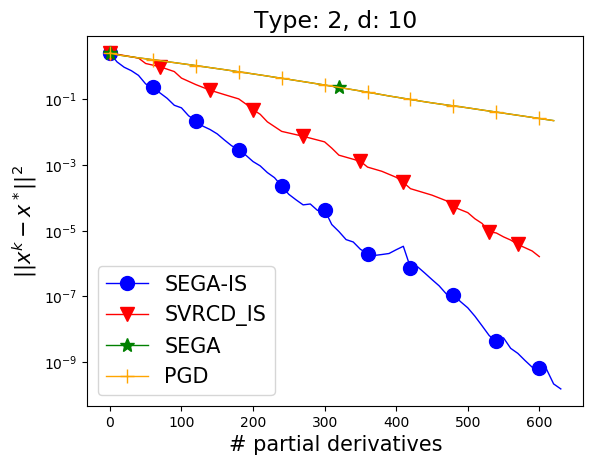}
\end{minipage}%
\begin{minipage}{0.25\textwidth}
  \centering
\includegraphics[width =  \textwidth ]{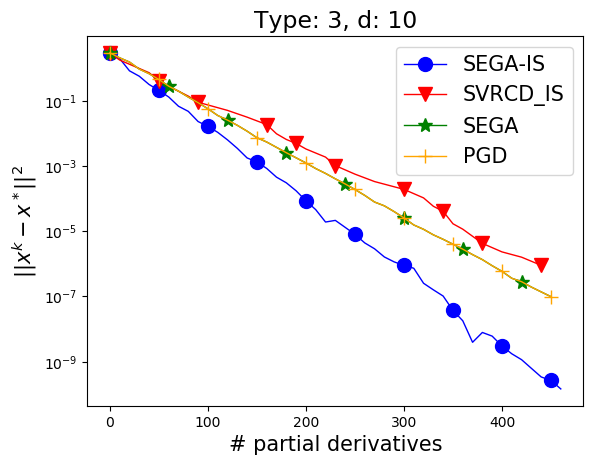}
\end{minipage}%
\begin{minipage}{0.25\textwidth}
  \centering
\includegraphics[width =  \textwidth ]{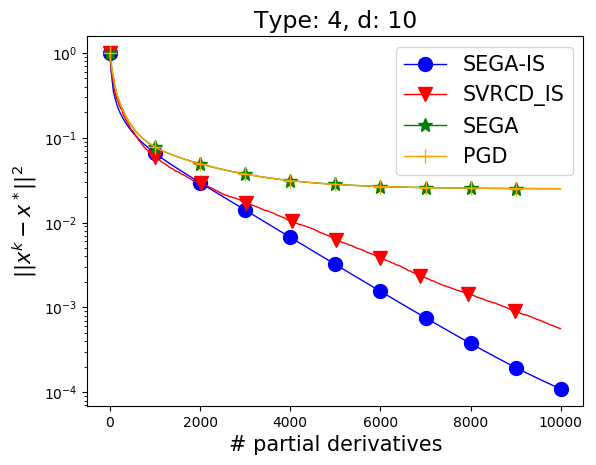}
\end{minipage}%
\\
\begin{minipage}{0.25\textwidth}
  \centering
\includegraphics[width =  \textwidth ]{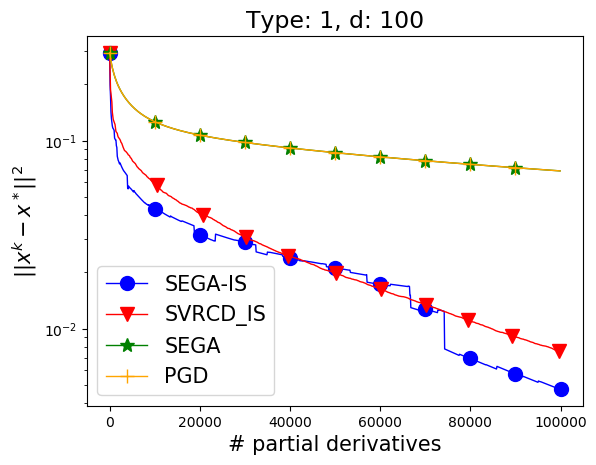}
\end{minipage}%
\begin{minipage}{0.25\textwidth}
  \centering
\includegraphics[width =  \textwidth ]{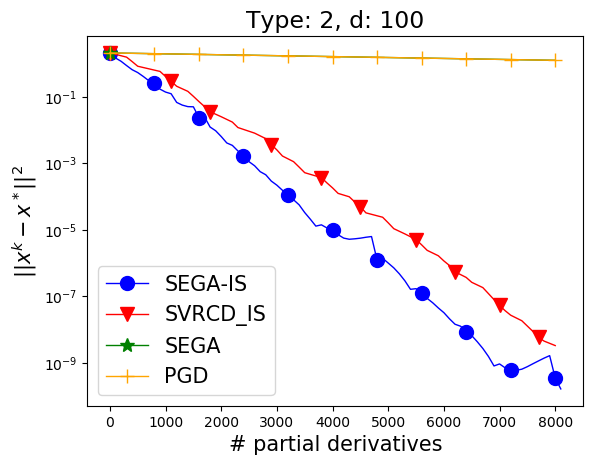}
\end{minipage}%
\begin{minipage}{0.25\textwidth}
  \centering
\includegraphics[width =  \textwidth ]{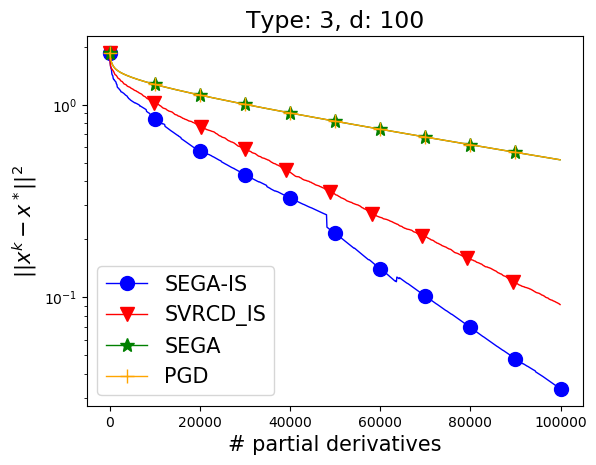}
\end{minipage}%
\begin{minipage}{0.25\textwidth}
  \centering
\includegraphics[width =  \textwidth ]{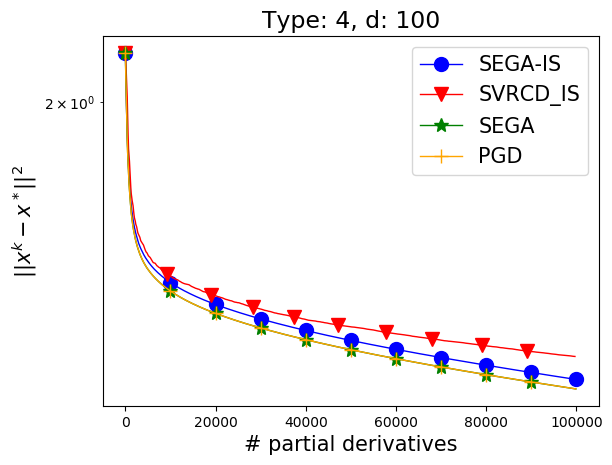}
\end{minipage}%
\\
\begin{minipage}{0.25\textwidth}
  \centering
\includegraphics[width =  \textwidth ]{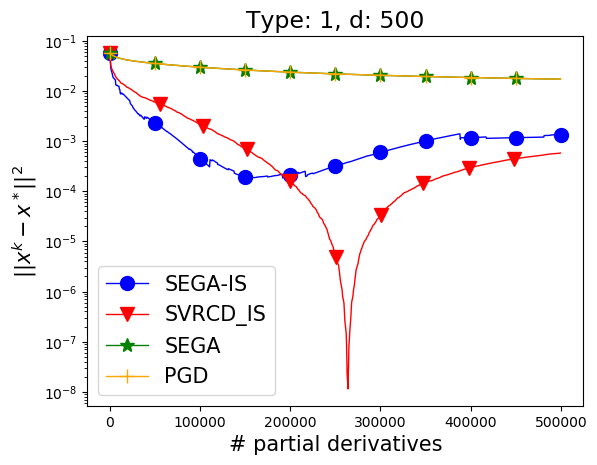}
\end{minipage}%
\begin{minipage}{0.25\textwidth}
  \centering
\includegraphics[width =  \textwidth ]{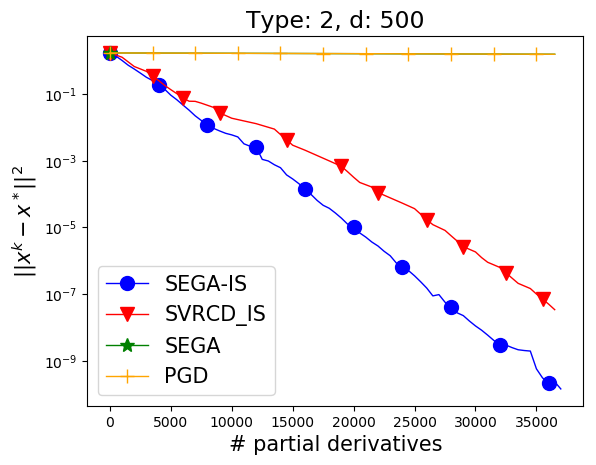}
\end{minipage}%
\begin{minipage}{0.25\textwidth}
  \centering
\includegraphics[width =  \textwidth ]{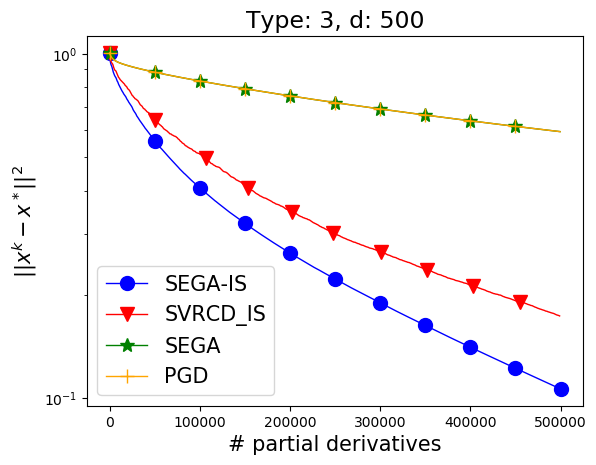}
\end{minipage}%
\begin{minipage}{0.25\textwidth}
  \centering
\includegraphics[width =  \textwidth ]{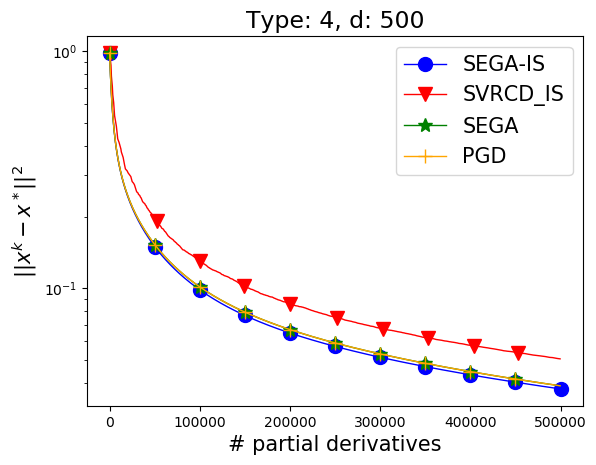}
\end{minipage}%
\caption{Comparison of {\tt SEGA-AS}, {\tt SVRCD-AS}, {\tt SEGA} and proximal gradient on 4 quadratic problems given by Table~\ref{tbl:gjs_quadratics}. {\tt SEGA-AS}, {\tt SVRCD-AS} and {\tt SEGA} compute single partial derivative each iteration ({\tt SVRCD} computes all of them with probability $\probx$), {\tt SEGA-AS}, {\tt SVRCD-AS} with probabilities proportional to diagonal of $\mM$.  }
\label{fig:gjs_sega_cmp}
\end{figure}

\subsection{{\tt SVRCD}: effect of $\probx$ \label{sec:gjs_svrcd_exp}}
In this experiment we demonstrate very broad range of $\probx$ can be chosen to still attain almost best possible rate for {\tt SVRCD} for problems from Table~\ref{tbl:gjs_quadratics} and $m,\gamma$ as described in Section~\ref{sec:gjs_sega_exp} Results can be found in Figure~\ref{fig:gjs_svrcd}. They indeed show that in many cases, varying $\probx$ from $\frac1n$ down to $\frac{2\lambda_{\min}(\mM)}{\sum_{i=1}^d m_i}$ does not influences the complexity significantly. However, too small $\probx$ leads to significantly slower convergence. Note that those findings are in accord with Corollary~\ref{cor:gjs_svrcd}. Similar results were shown in~\cite{kovalev2019don} for {\tt LSVRG}. 

\begin{figure}[!h]
\centering
\begin{minipage}{0.25\textwidth}
  \centering
\includegraphics[width =  \textwidth ]{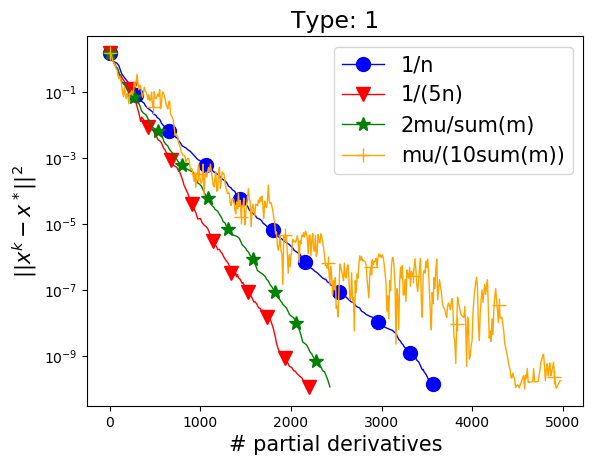}
\end{minipage}%
\begin{minipage}{0.25\textwidth}
  \centering
\includegraphics[width =  \textwidth ]{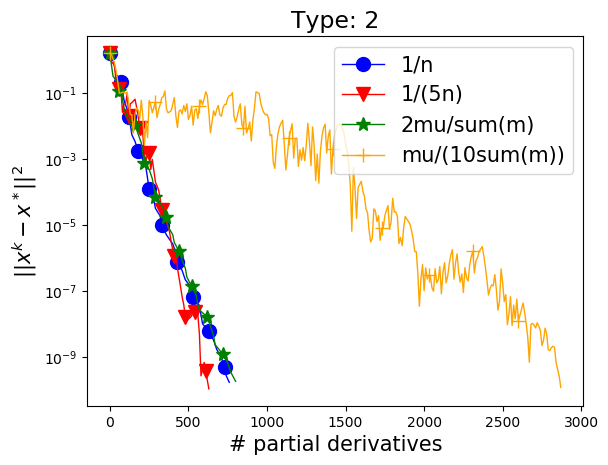}
\end{minipage}%
\begin{minipage}{0.25\textwidth}
  \centering
\includegraphics[width =  \textwidth ]{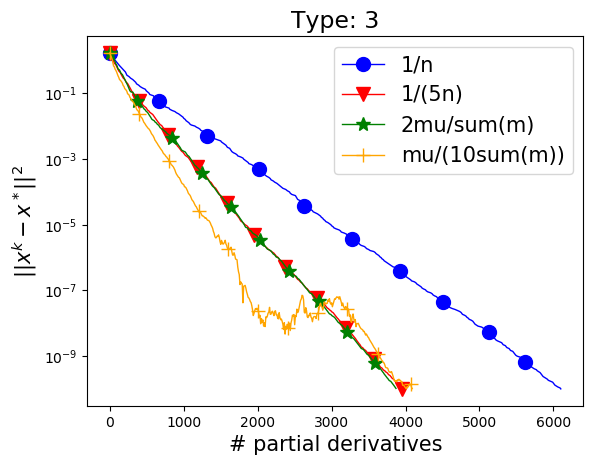}
\end{minipage}%
\begin{minipage}{0.25\textwidth}
  \centering
\includegraphics[width =  \textwidth ]{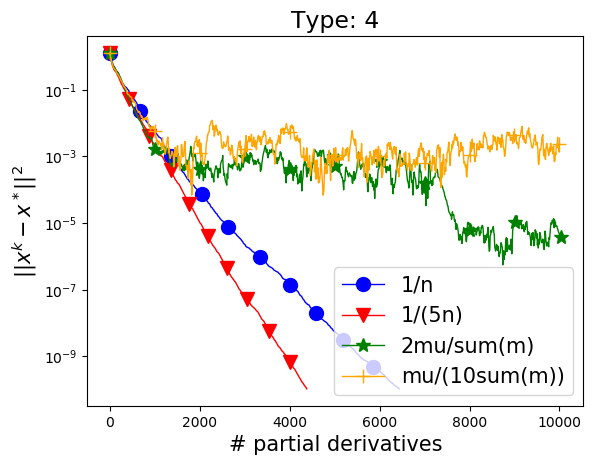}
\end{minipage}%
\\
\begin{minipage}{0.25\textwidth}
  \centering
\includegraphics[width =  \textwidth ]{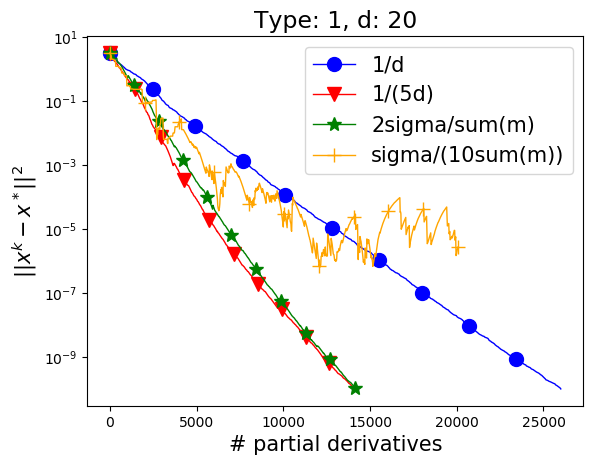}
\end{minipage}%
\begin{minipage}{0.25\textwidth}
  \centering
\includegraphics[width =  \textwidth ]{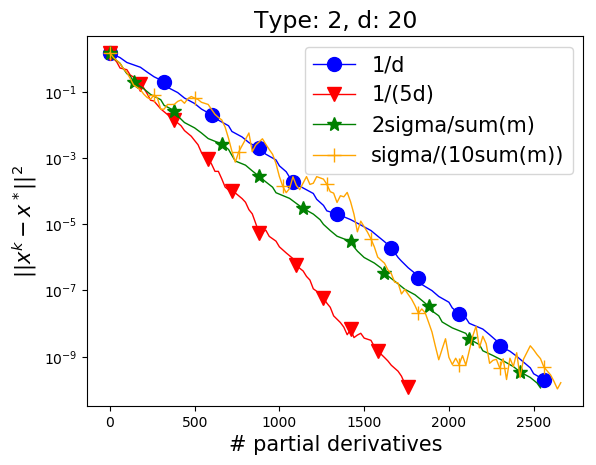}
\end{minipage}%
\begin{minipage}{0.25\textwidth}
  \centering
\includegraphics[width =  \textwidth ]{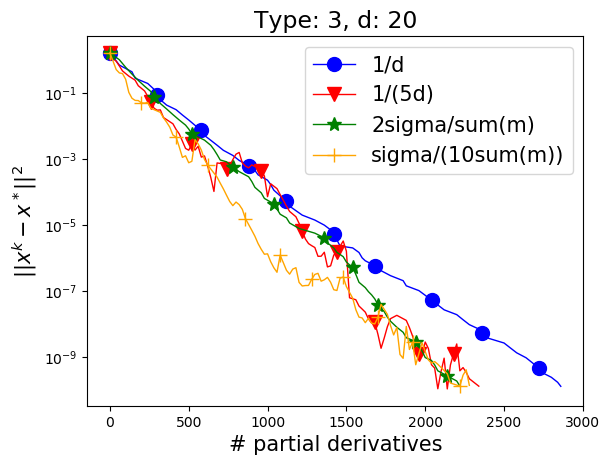}
\end{minipage}%
\begin{minipage}{0.25\textwidth}
  \centering
\includegraphics[width =  \textwidth ]{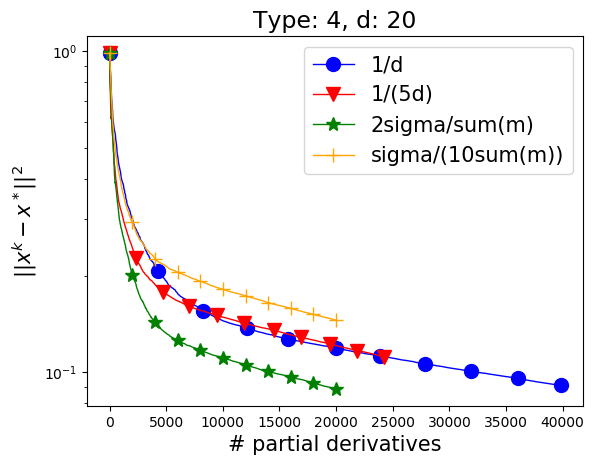}
\end{minipage}%
\\
\begin{minipage}{0.25\textwidth}
  \centering
\includegraphics[width =  \textwidth ]{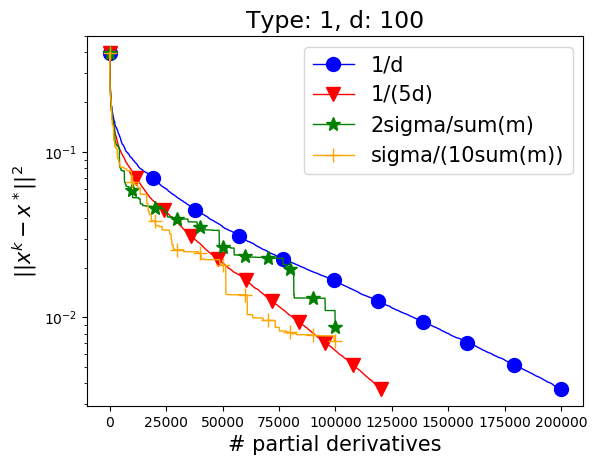}
\end{minipage}%
\begin{minipage}{0.25\textwidth}
  \centering
\includegraphics[width =  \textwidth ]{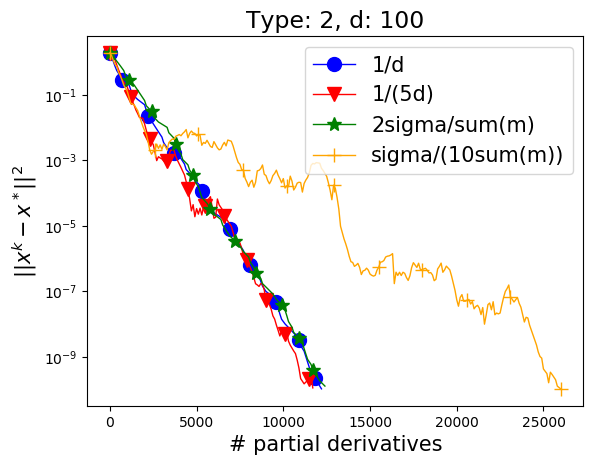}
\end{minipage}%
\begin{minipage}{0.25\textwidth}
  \centering
\includegraphics[width =  \textwidth ]{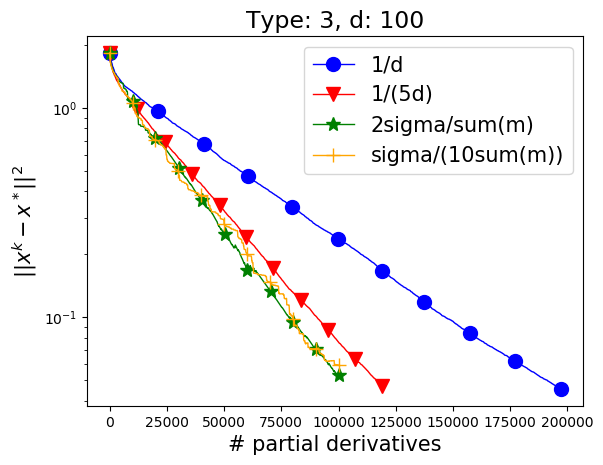}
\end{minipage}%
\begin{minipage}{0.25\textwidth}
  \centering
\includegraphics[width =  \textwidth ]{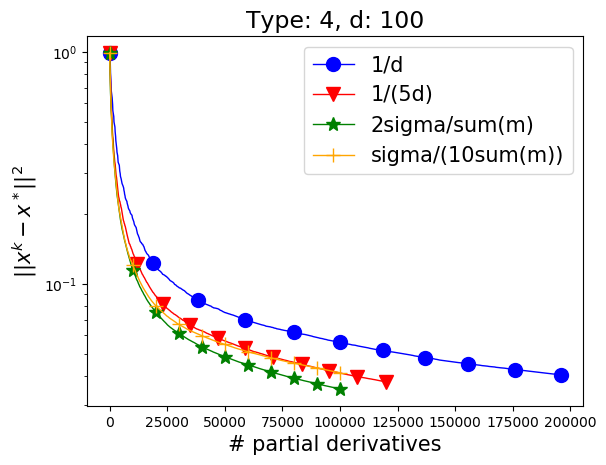}
\end{minipage}%
\\
\caption{The effect of $\probx$ on convergence rate of {\tt SVRCD} on quadratic problems from Table~\ref{tbl:gjs_quadratics}. In every case, probabilities were chosen proportionally to the diagonal of $\mM$ and only a single partial derivative is evaluated in $\cS$.}
\label{fig:gjs_svrcd}
\end{figure}

\subsection{{\tt ISAEGA} \label{sec:gjs_ISAEGA_exp}}

In this section we test a simple version of {\tt ISAEGA} (Algorithm~\ref{alg:gjs_isaega})\footnote{The full description of {\tt ISAEGA}, together with convergence guarantees are provided in Section~\ref{sec:gjs_ISAEGA}}. As mentioned, {\tt ISAEGA} is an algorithm for distributed optimization which, at each iteration, computes a subset of partial derivatives of stochastic gradient on each machine, and constructs corresponding Jacobian estimate and stochastic gradient.

For simplicity, we consider only the simple version which assumes $\mM_j = m \mI_d$ for all $j$ (i.e. we do not do importance sampling), and we suppose that $|R_\tR |=1$ always for all $\tR$ (i.e. each machine always looks at a single function from the local finite sum). Further, we consider $\psi(x)=0$. Corollary~\ref{cor:gjs_isaega} shows that, if the condition number of the problem is not too small, {\tt ISAEGA} with $|L_{\tR} | \approx \frac1\TR$ (where $\TR$ is a number of parallel units) enjoys, up to small constant factor, same rate as {\tt SAGA} (which is, under a convenient smoothness, the same rate as the convergence rate of gradient descent). Thus, {\tt ISAEGA} scales linearly in terms of partial derivative complexity in parallel setup. In other words, given that we have twice more workers, each of them can afford to evaluate twice less partial derivatives\footnote{Practical implications of the method are further explained in~\cite{mishchenko201999}.}. The experiments we propose aim to verify this claim. 

We consider $\ell_2$ regularized logistic regression (for the binary classification). In particular, 
\[
\forall j: \quad f_j(x) \eqdef  \log \left(1+\exp\left(\mA_{j,:}x\cdot  y_i\right) \right)+\frac{\lambda}{2} \| x\|^2,
\]
where $\mA\in \R^{n\times d}$ is a data matrix, $y\in \{-1,1\}^{n}$ is a vector of labels and $\lambda\geq 0$ is the regularization parameter.  
Both $\mA,y$ are provided from LibSVM~\cite{chang2011libsvm} datasets: {\tt a1a}, {\tt a9a}, w1a, {\tt w8a}, {\tt gisette}, {\tt madelon}, {\tt phishing} and {\tt mushrooms}. Further,  $\mA$ was normalized such that $\| \mA_{j,:}\|^2=1$. 
Next, it is known that $f_j$ is $(\frac14+ \lambda)$-smooth, convex, while $f$ is $\lambda$-strongly convex. Therefore, as a stepsize for all versions of ${\tt ISAEGA}$, we set $\gamma = \frac{1}{6\lambda + \frac{3}{2}}$ (this is an approximation of theoretical stepsize). 

In each experiment, we compare 4 different setups for {\tt ISEAGA} -- given by 4 different values of $\TR$. Given a value of $\TR$, we set $|L_\tR|=\frac1\TR$ for all $\tR$. Further, we always sample $L_\tR$ uniformly. The results are presented in Figure~\ref{fig:gjs_ISAEGA}. Indeed, we observe the almost perfect parallel linear scaling.

For completeness, we provide dataset sized in Table~\ref{tbl:gjs_libsvm}. 

 \begin{table}[!h]
\begin{center}
\begin{tabular}{|c|c|c|}
\hline
Name & $n $ & $d$ \\
 \hline
  \hline
{\tt a1a}   & $1605$ & $123$  \\
  \hline
{\tt a9a}   & $32561$ & $123$  \\
  \hline
{\tt w1a}   & $2477$ & $300$  \\
  \hline
{\tt w8a}   & $49749$ & $300$  \\
  \hline
{\tt gisette}   & $6000$ & $5000$  \\
  \hline
{\tt madelon}   & $2000$ & $500$  \\
  \hline
{\tt phishing}   & $11055$ & $68$  \\
  \hline
{\tt mushrooms}   & $8124$ & $112$  \\
\hline
\end{tabular}
\end{center}
\caption{Table of LibSVM data used for our experiments. }
\label{tbl:gjs_libsvm}
\end{table}

\begin{figure}[!h]
\centering
\begin{minipage}{0.35\textwidth}
  \centering
\includegraphics[width =  \textwidth ]{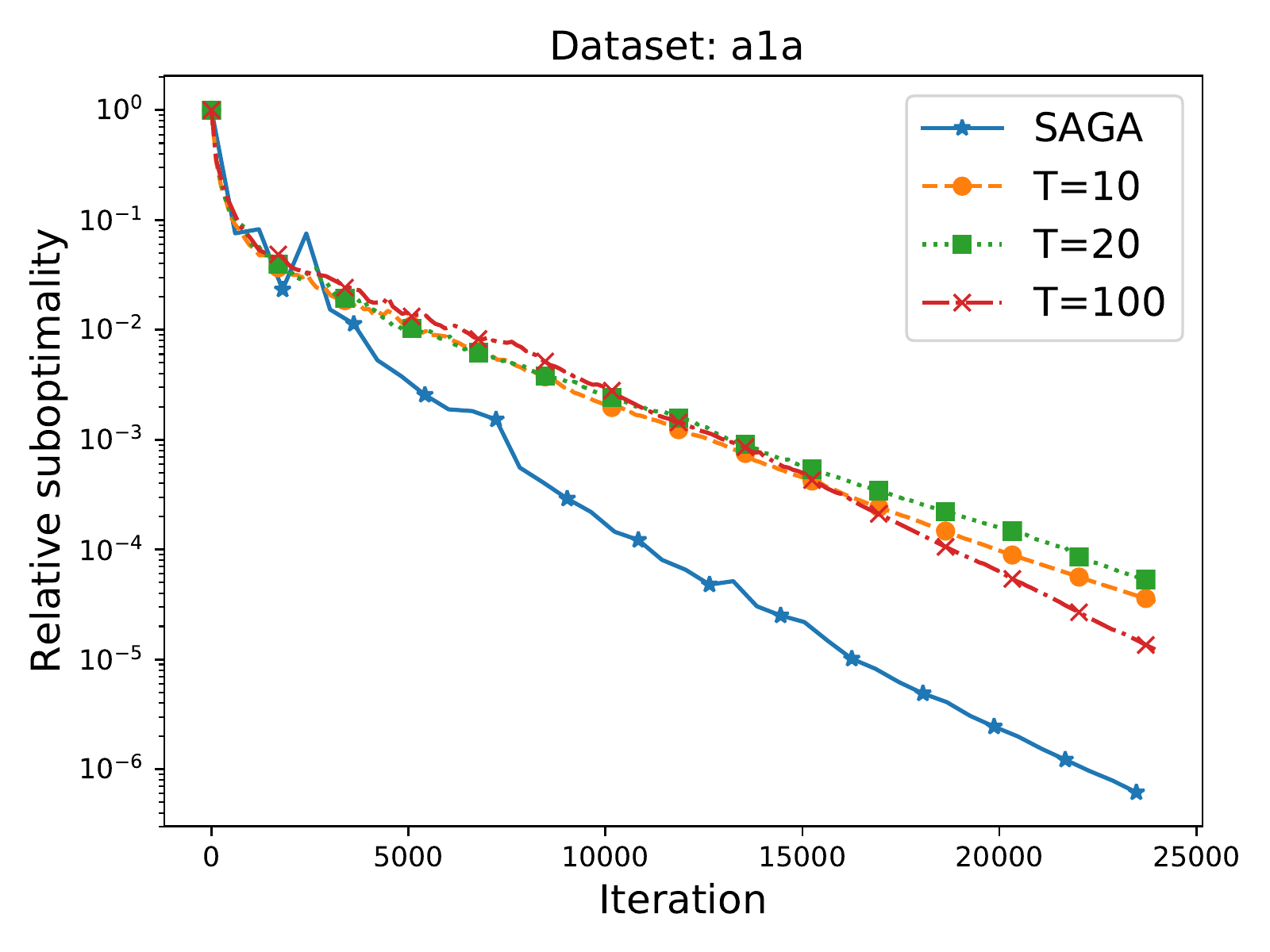}
\end{minipage}%
\begin{minipage}{0.35\textwidth}
  \centering
\includegraphics[width =  \textwidth ]{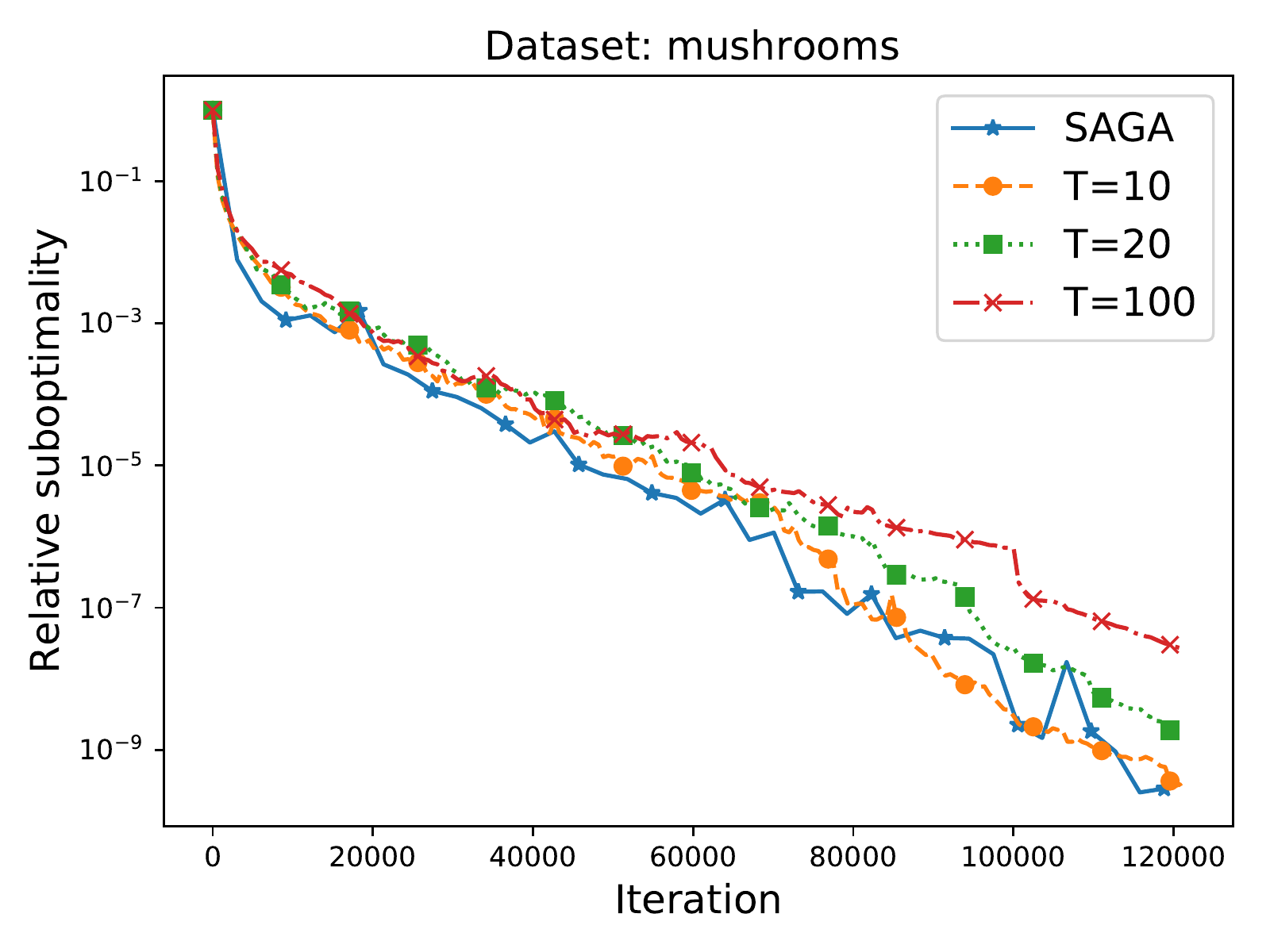}
\end{minipage}%
\\
\begin{minipage}{0.35\textwidth}
  \centering
\includegraphics[width =  \textwidth ]{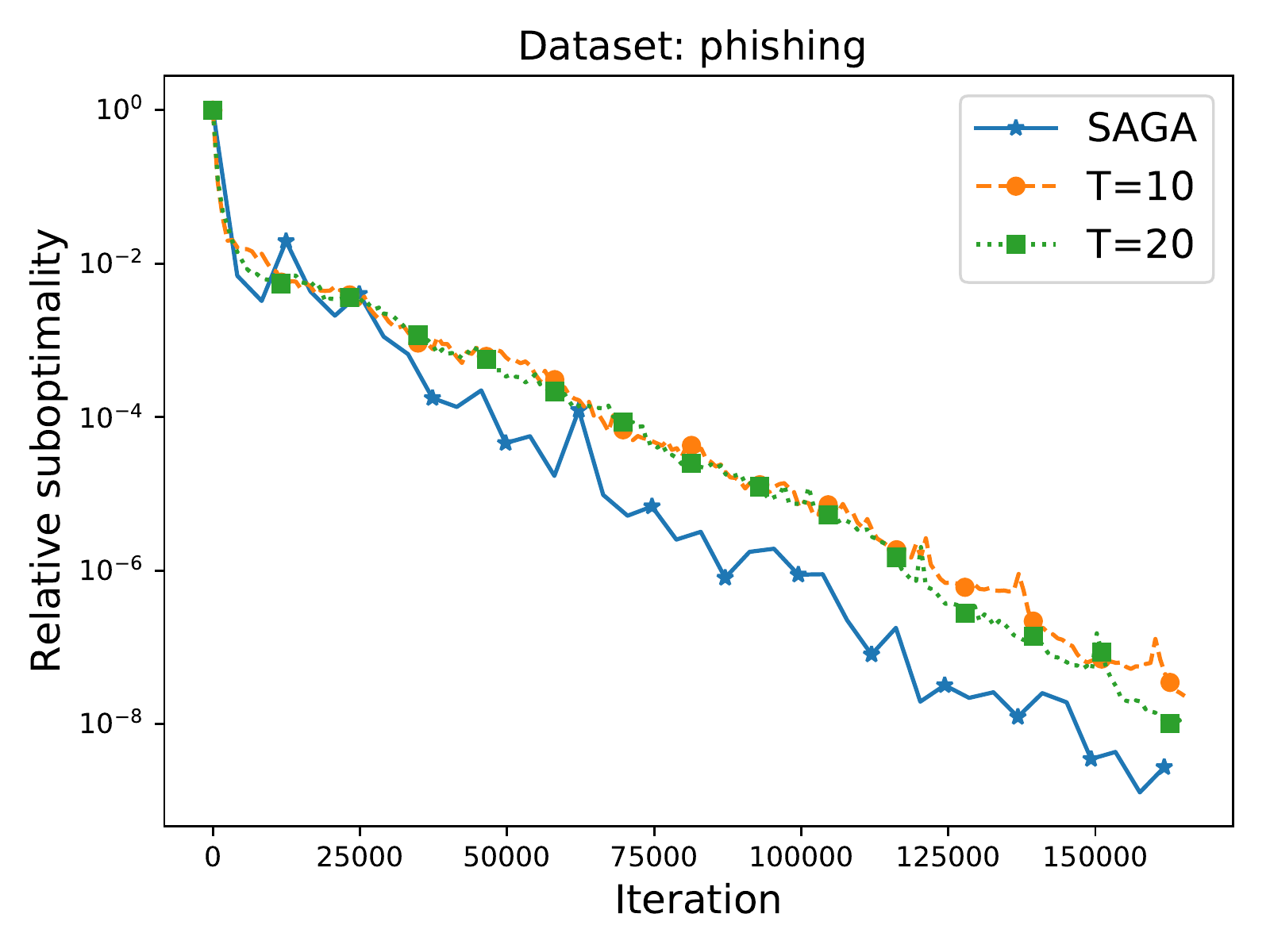}
\end{minipage}%
\begin{minipage}{0.35\textwidth}
  \centering
\includegraphics[width =  \textwidth ]{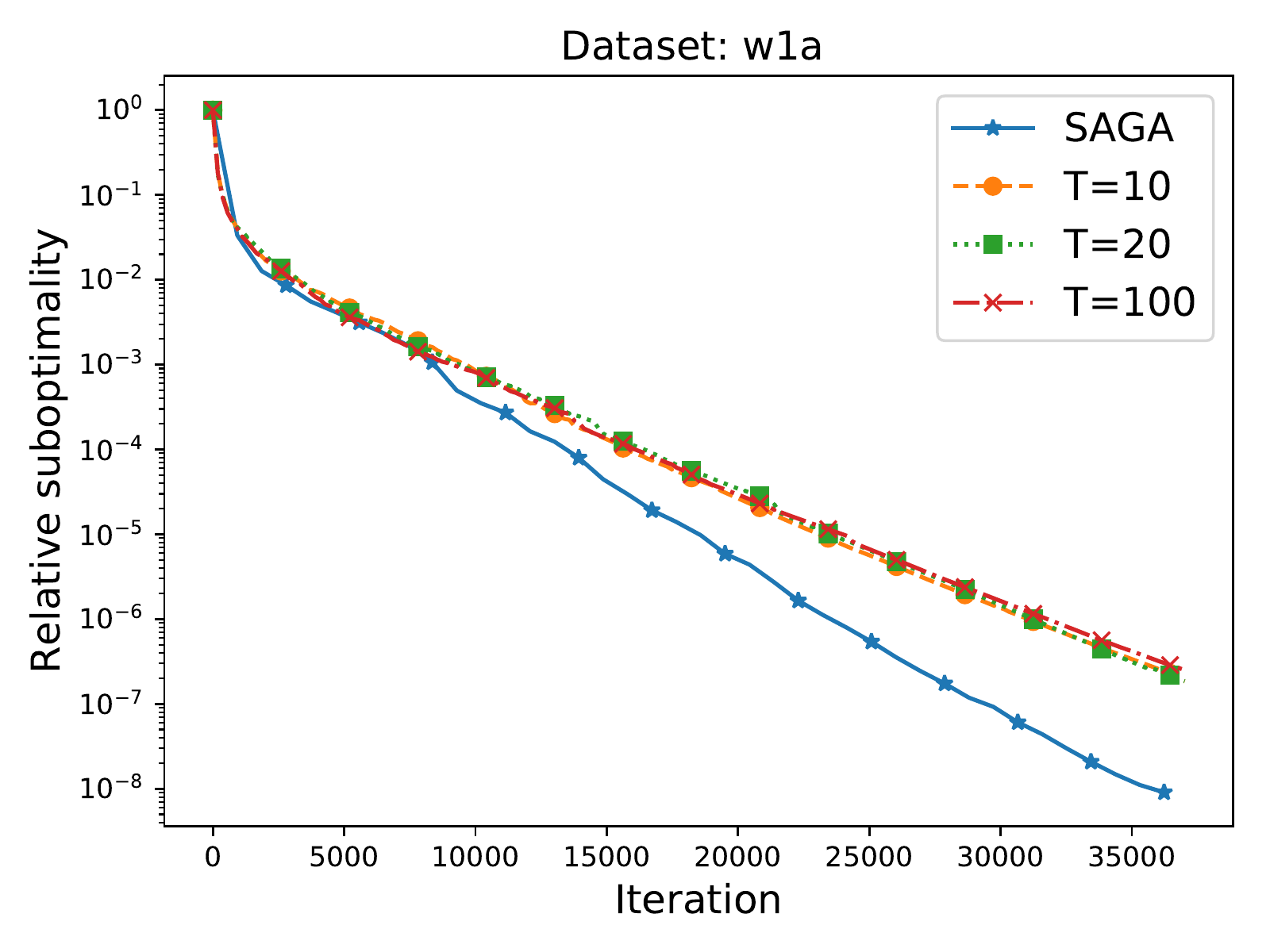}
\end{minipage}%
\\
\begin{minipage}{0.35\textwidth}
  \centering
\includegraphics[width =  \textwidth ]{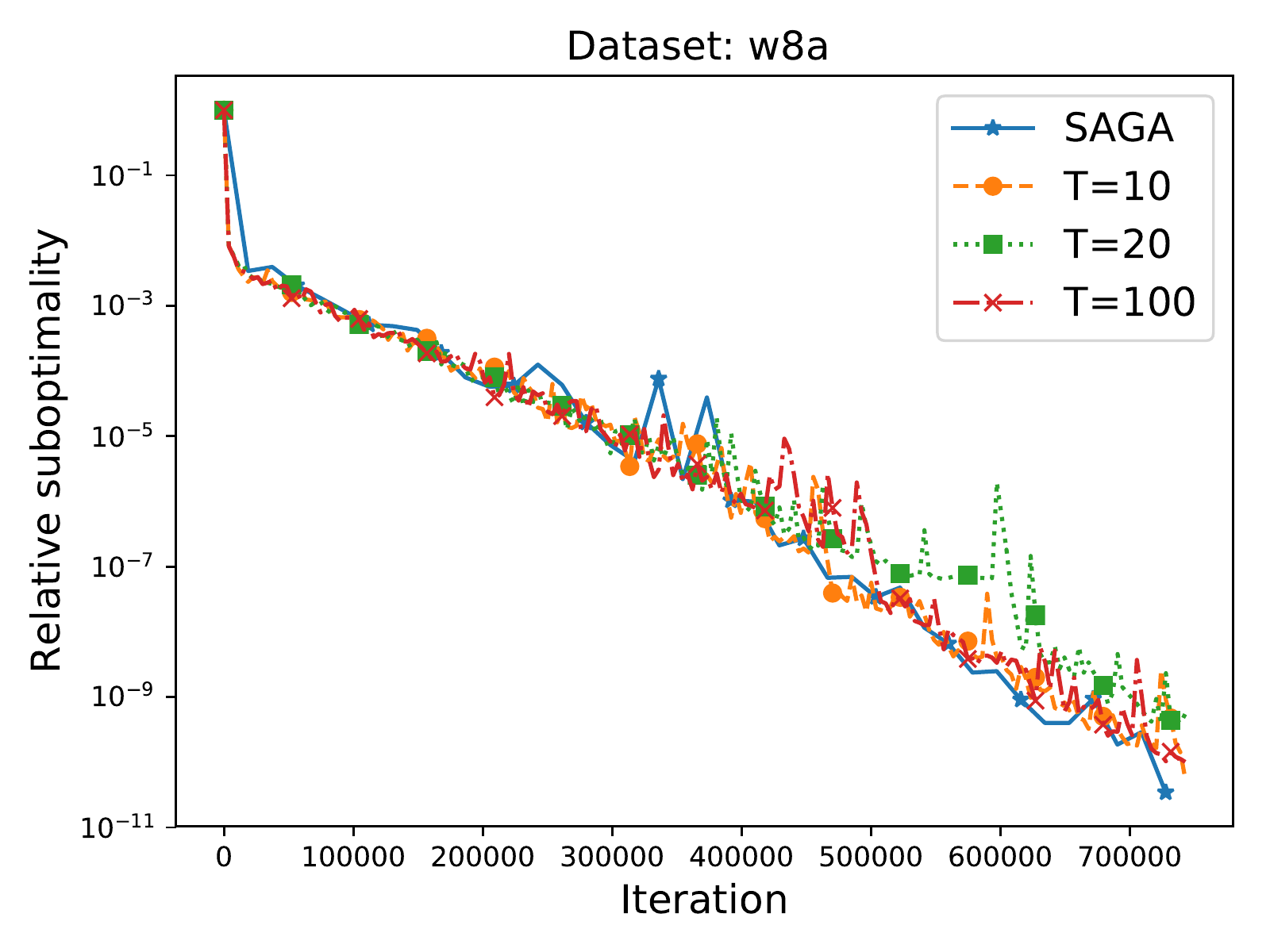}
\end{minipage}%
\begin{minipage}{0.35\textwidth}
  \centering
\includegraphics[width =  \textwidth ]{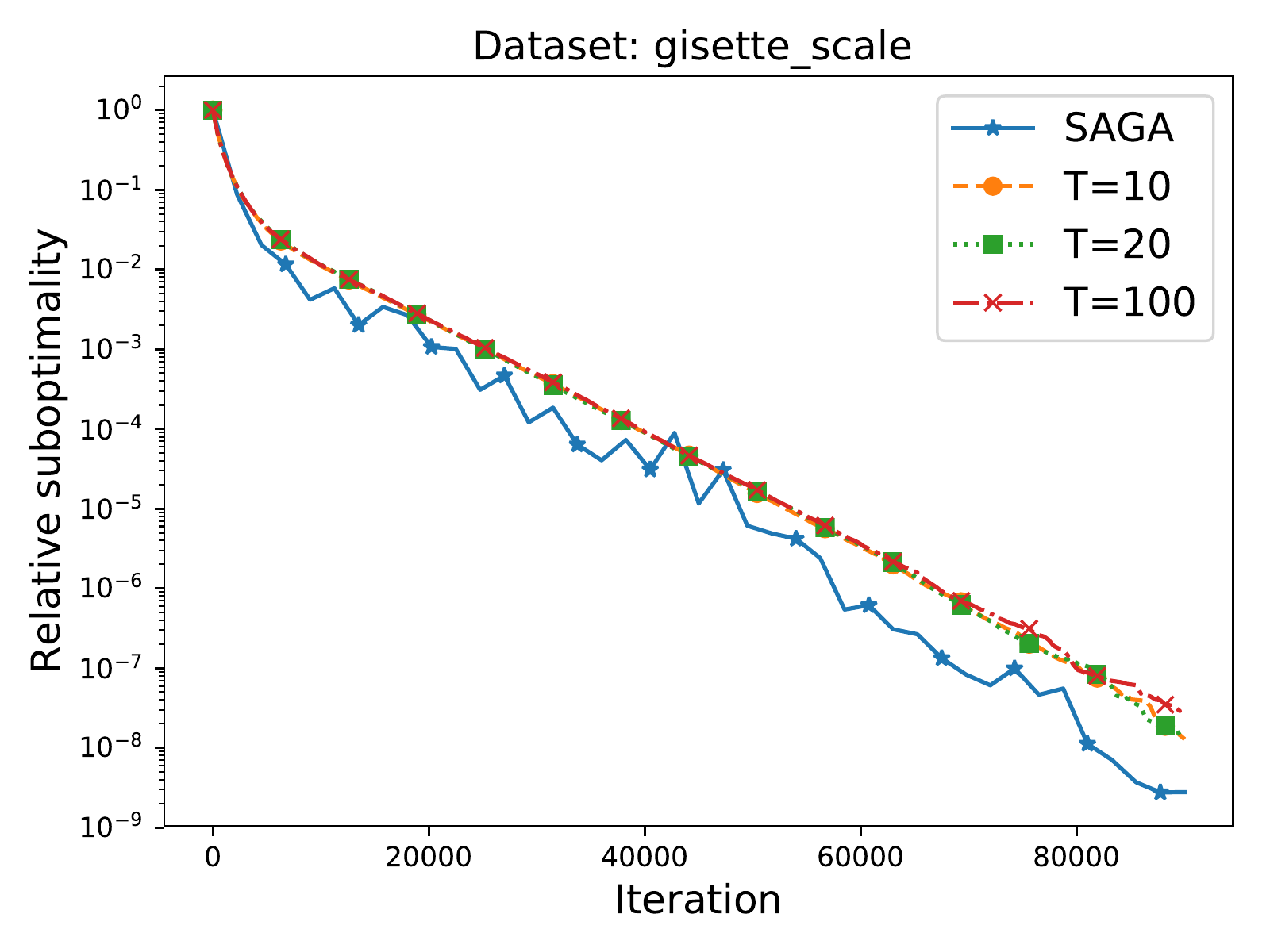}
\end{minipage}%
\\
\begin{minipage}{0.35\textwidth}
  \centering
\includegraphics[width =  \textwidth ]{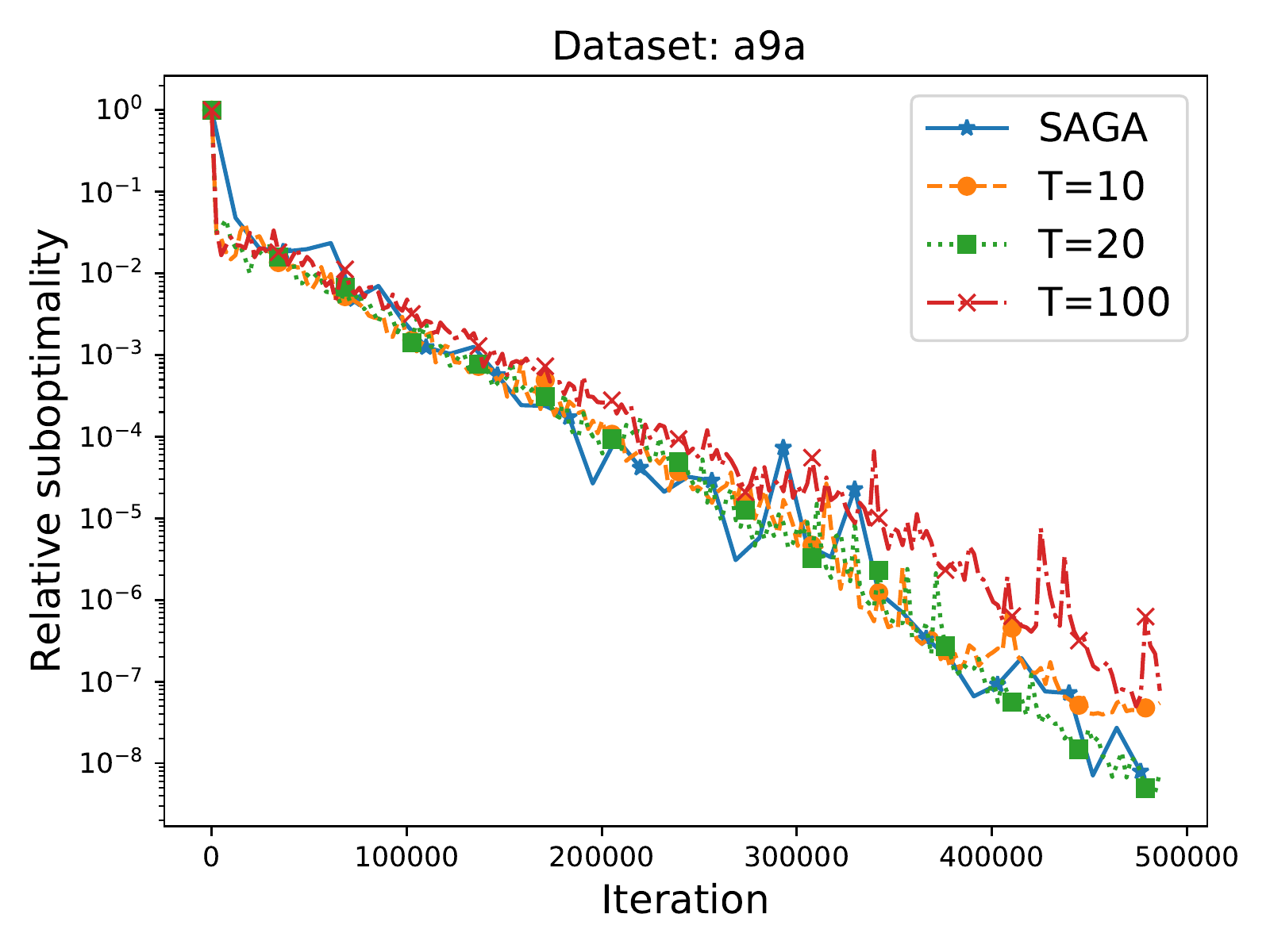}
\end{minipage}%
\begin{minipage}{0.35\textwidth}
  \centering
\includegraphics[width =  \textwidth ]{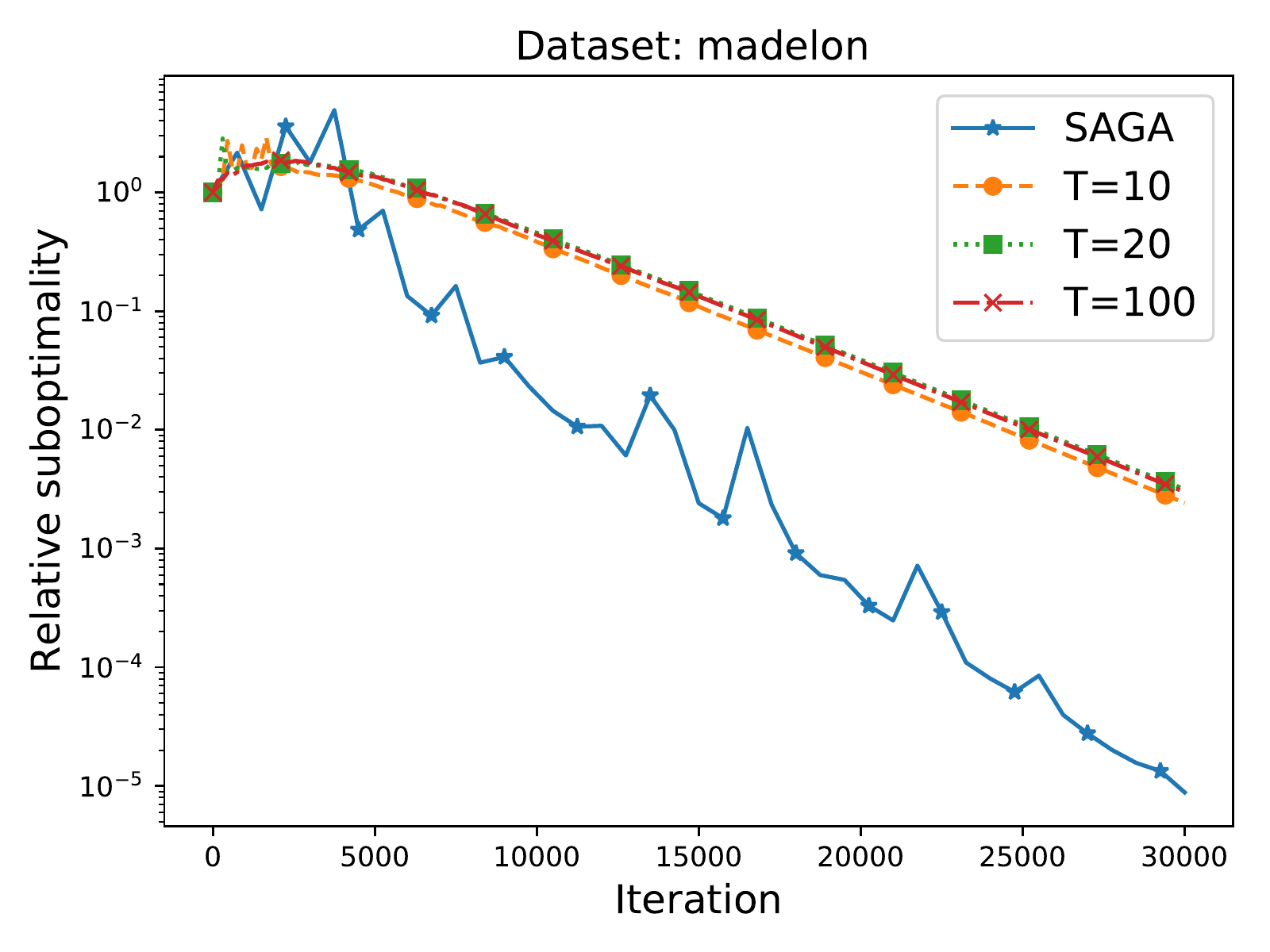}
\end{minipage}%
\caption{{\tt ISAEGA} applied on LIBSVM~\cite{chang2011libsvm} datasts with $\lambda = 4\cdot 10^{-5}$. Axis $y$ stands for relative suboptimality, i.e. $\frac{f(x^k)-f(x^*)}{f(x^k)-f(x^0)}$.}
\label{fig:gjs_ISAEGA}
\end{figure}

\subsection{{\tt LSVRG} with importance sampling \label{sec:gjs_extra_lsvrg}}
As mentioned, one of the contributions of this work is {\tt LSVRG} with arbitrary sampling. In this section, we demonstrate that designing a good sampling can yield a significant speedup in practice. We consider logistic regression problem on LibSVM~\cite{chang2011libsvm} data, as described in Section~\ref{sec:gjs_ISAEGA_exp}. However, since LibSVM data are normalized, we pre-multiply each row of the data matrix by a random scaling factor. In particular, the scaling factors are proportional to $l^2$ where $l$ is sampled uniformly from $[1000]$ such that the Frobenius norm of the data matrix is $n$. For the sake of simplicity, consider case $\lambda=0$.

\paragraph{Choice vector $v$.} 
Note that since $\mM_j = \mA_{j:}^\top \mA_{j:}$, the following claim must hold: \emph{ Consider fixed $v$. Then if~\eqref{eq:gjs_ESO_saga} holds for any set of vector $\{h_j\}_{j=1}^n$ such that $h_j$ is parallel to $\mA_{j:}$, then ~\eqref{eq:gjs_ESO_saga} holds for any set of vector $\{h_j\}_{j=1}^n$.} Thus, we can set $h_j = c_j\mA_{j:}^\top/\| \mA_{j:}\|$ without loss of generality. Thus, $\mM_j^{\frac12} h_j= c_j \mA_{j:}^\top$, and~\eqref{eq:gjs_ESO_saga} becomes equivalent to $\PR \circ \left(\mA^\top \mA\right)\preceq \diag(p\circ v)$ where $\PR_{jj'} = \Probbb{j\in R, j' \in R}$. Note that this is exactly \emph{expected separable overapproximation (ESO)} for coordinate descent~\cite{qu2016coordinate2}. Thus we choose vector $v$ to be proportional to $\pR$ such that $\PR \circ \left(\mA^\top \mA\right)\preceq \diag(p\circ v)$ holds (as proposed in Chapter~\ref{acd}). In order to compute the scaling constant, one needs to evaluate maximum eigenvalue of PSD $n\times n$ matrix, which is of $\cO(n^2)$ cost. We do so in the experiments. Note that there is a suboptimal, but cheeaper way to obtain $v$ described in~\cite{qian2019saga}. Lastly, if $\lambda>0$, we set $v$ such that $\PR \circ \left(\mA^\top \mA +\lambda \mI \right)\preceq \diag(p\circ v)$.

\paragraph{Choice of probabilities.}
In order to be fair, we only compare methods where $\E{|R|}=\tau$. For the case $\tau =1$, we consider a sampling such that $|R|=1$ according to a given probability vector $\pR$. For uniform sampling, we have $\pR=n^{-1}\eR$, while for importance sampling, we set $\pR_j= \frac{\lambda_{\max}(\mM_j)}{\sum_{j'=1}^n \lambda_{\max}(\mM_{j'})}$. In the case $\tau >1$, we consider independent sampling from Chapter~\ref{acd}. In particular, $\Probbb{j\in R} = \pRj$ with $\sum \pRj = \tau$ and binary random variables $(j\in R)$ are jointly independent. For uniform sampling we have $\pR = \tau n^{-1} \eR$. For importance sampling, probability vector $\pR$ is chosen such that $p_j =  \frac{\lambda_{\max}(\mM_j)}{\varrho+ \lambda_{\max}(\mM_{j}}$, where $\varrho$ is such that $\sum \pRj = \tau$. The mentioned sampling was proven to be superior over uniform minibatching in Chapter~\ref{acd}. Next, stepsize $\gamma = \frac{1}{6}\min_j\frac{n\pRj}{v_j}$ was chosen for all methods.

Lastly, $\probx = \frac{1}{2n}$ was chosen for {\tt LSVRG}. The results are presented in Figures~\ref{fig:gjs_LSVRG1} and~\ref{fig:gjs_LSVRG2}.

In all cases, {\tt LSVRG} with importance sampling was the fastest method. As provided theory suggests, it outperformed methods with importance sampling especially significantly for small $\tau$; and the larger $\tau$, the smaller the effect of importance sampling is. However, our experiments indicate the superiority of {\tt LSVRG} to  {\tt SAGA} in the importance sampling setup. In particular, stepsize $\gamma = \frac{1}{6}\min_j\frac{n\pRj}{v_j}$ is often too large for {\tt SAGA}. Note that both optimal stepsize and optimal probabilities require the prior knowledge of the quasi strong convexity constant $\mu$\footnote{Or more generally, strong growth constant, see Appendix~\ref{sec:gjs_sg}} which is, in our case unknown (see the importance serial sampling proposed in~\cite{jacsketch}, and {\tt SAGA} is more sensitive to that choice. One can still estimate it as $\lambda$, however, this would yield suboptimal performance as well.

\begin{figure}[!h]
\centering
\begin{minipage}{0.3\textwidth}
  \centering
\includegraphics[width =  \textwidth ]{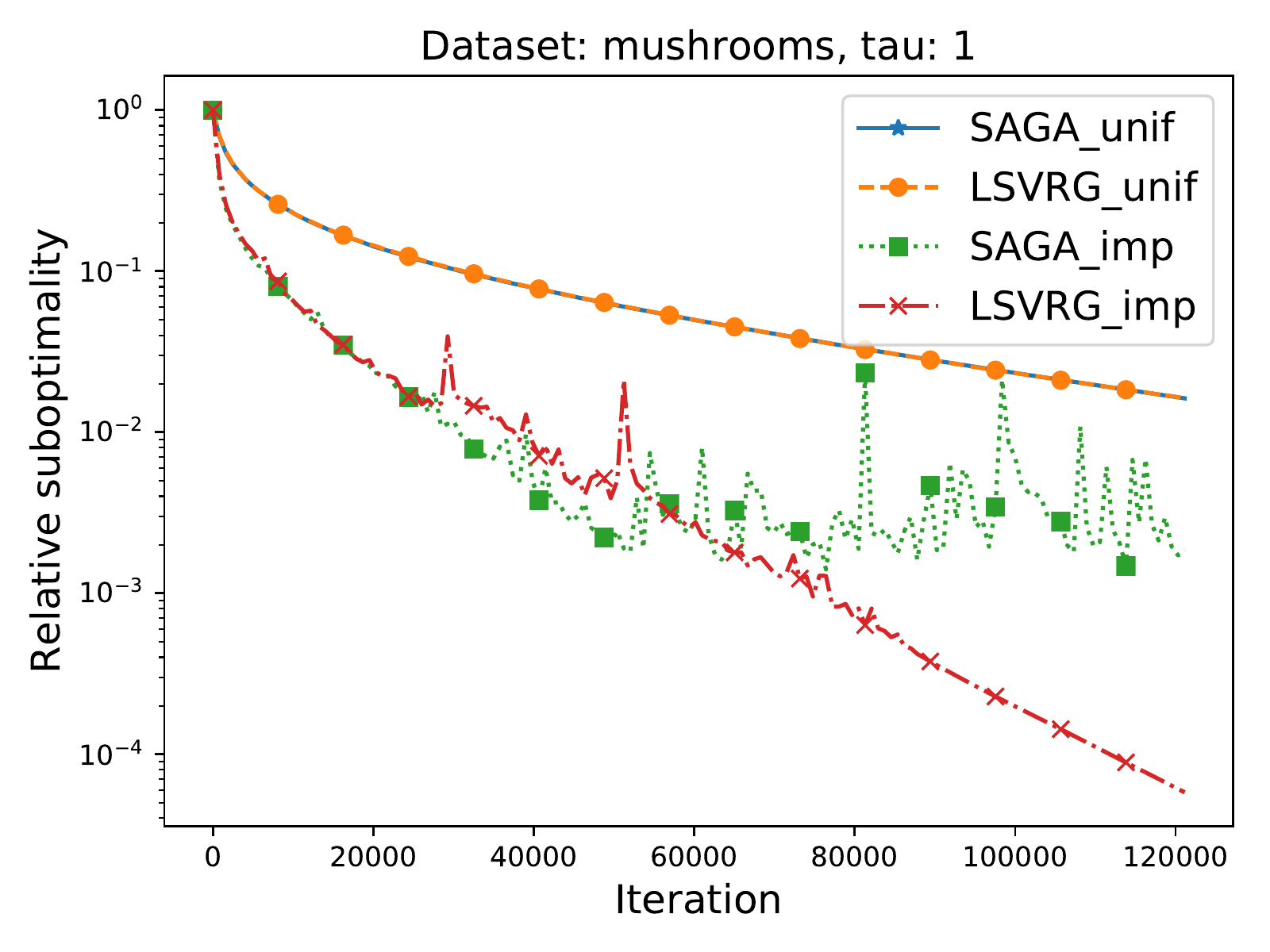}
\end{minipage}%
\begin{minipage}{0.3\textwidth}
  \centering
\includegraphics[width =  \textwidth ]{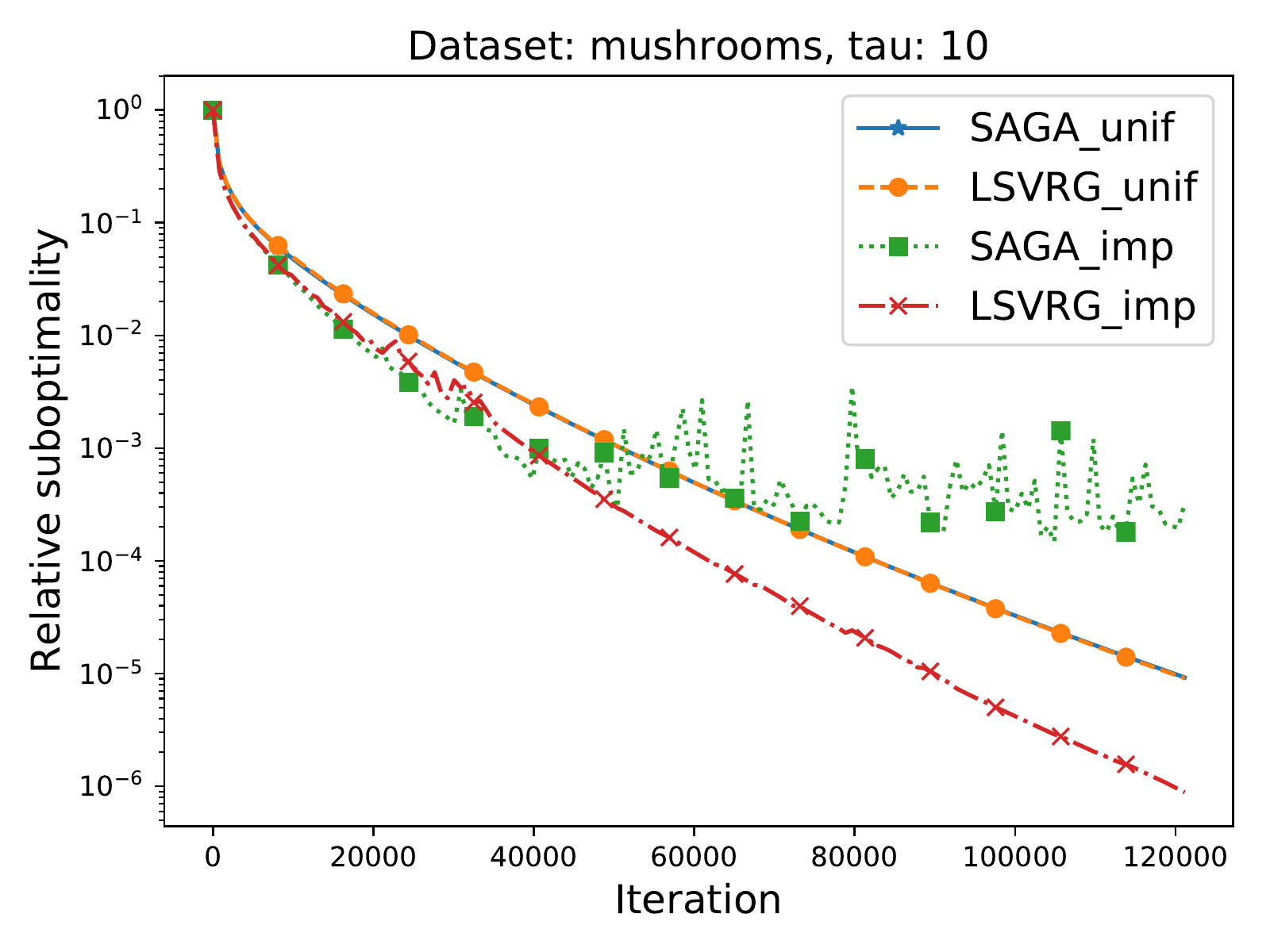}
\end{minipage}%
\begin{minipage}{0.3\textwidth}
  \centering
\includegraphics[width =  \textwidth ]{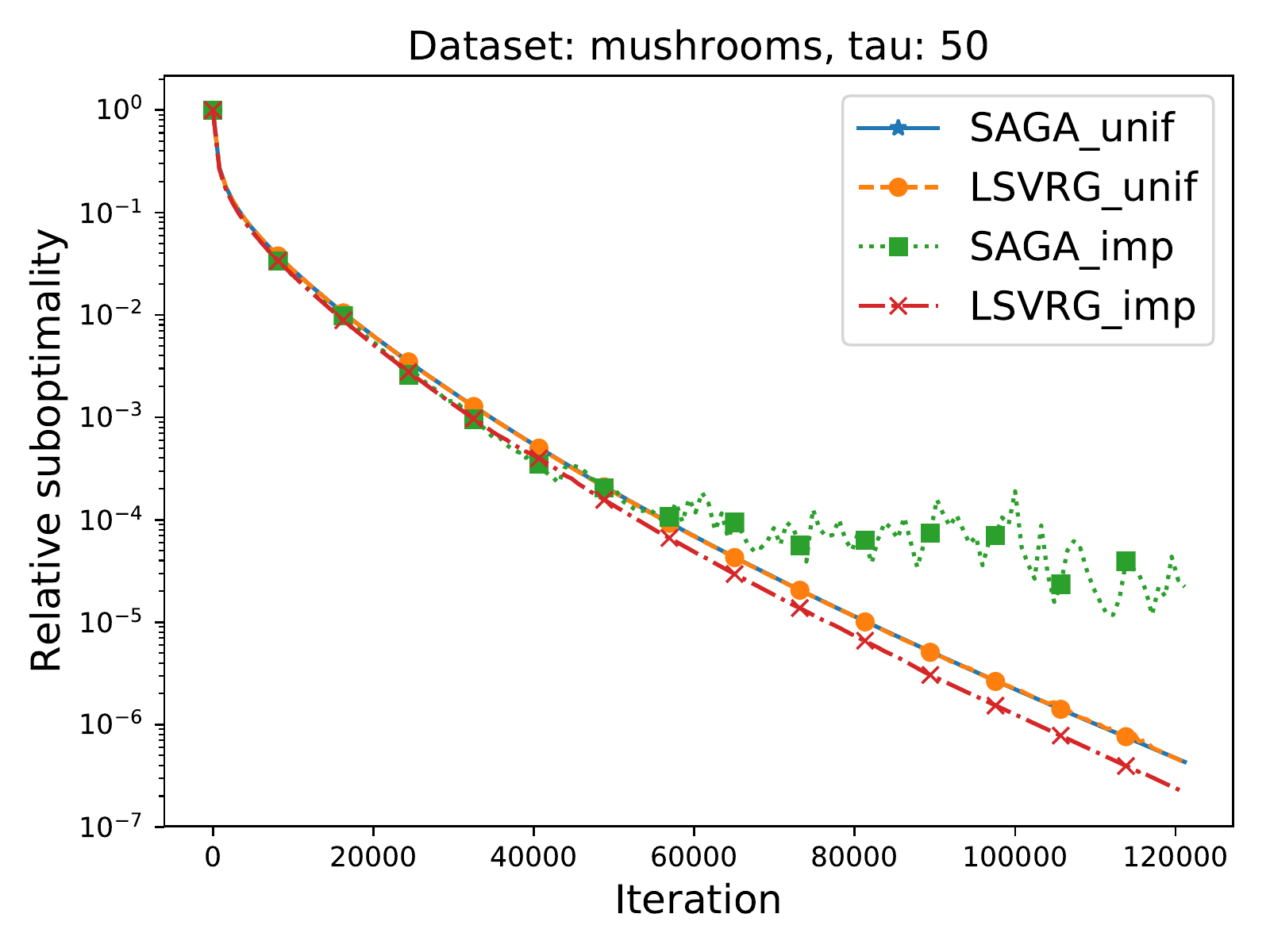}
\end{minipage}%
\\
\begin{minipage}{0.3\textwidth}
  \centering
\includegraphics[width =  \textwidth ]{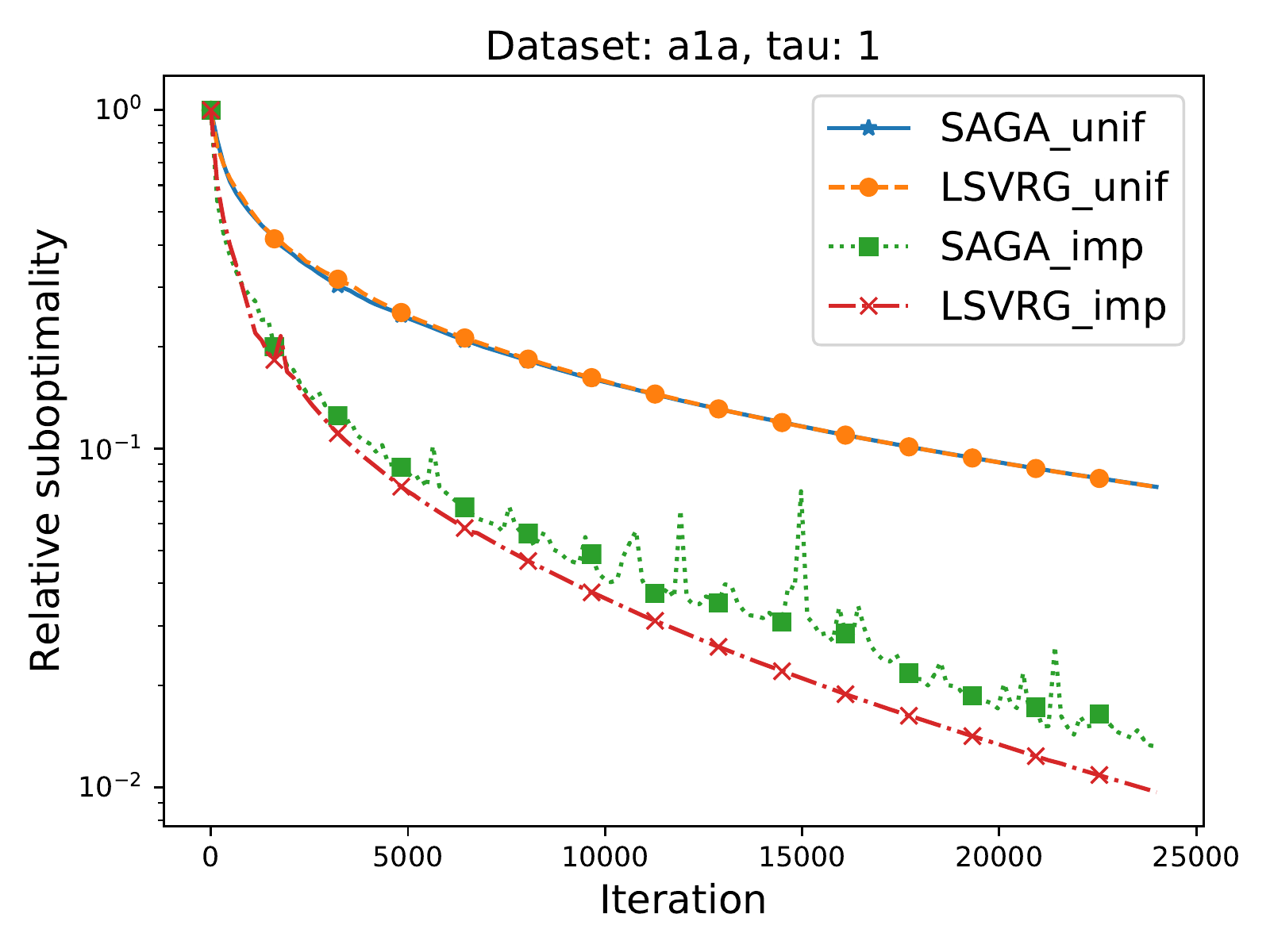}
\end{minipage}%
\begin{minipage}{0.3\textwidth}
  \centering
\includegraphics[width =  \textwidth ]{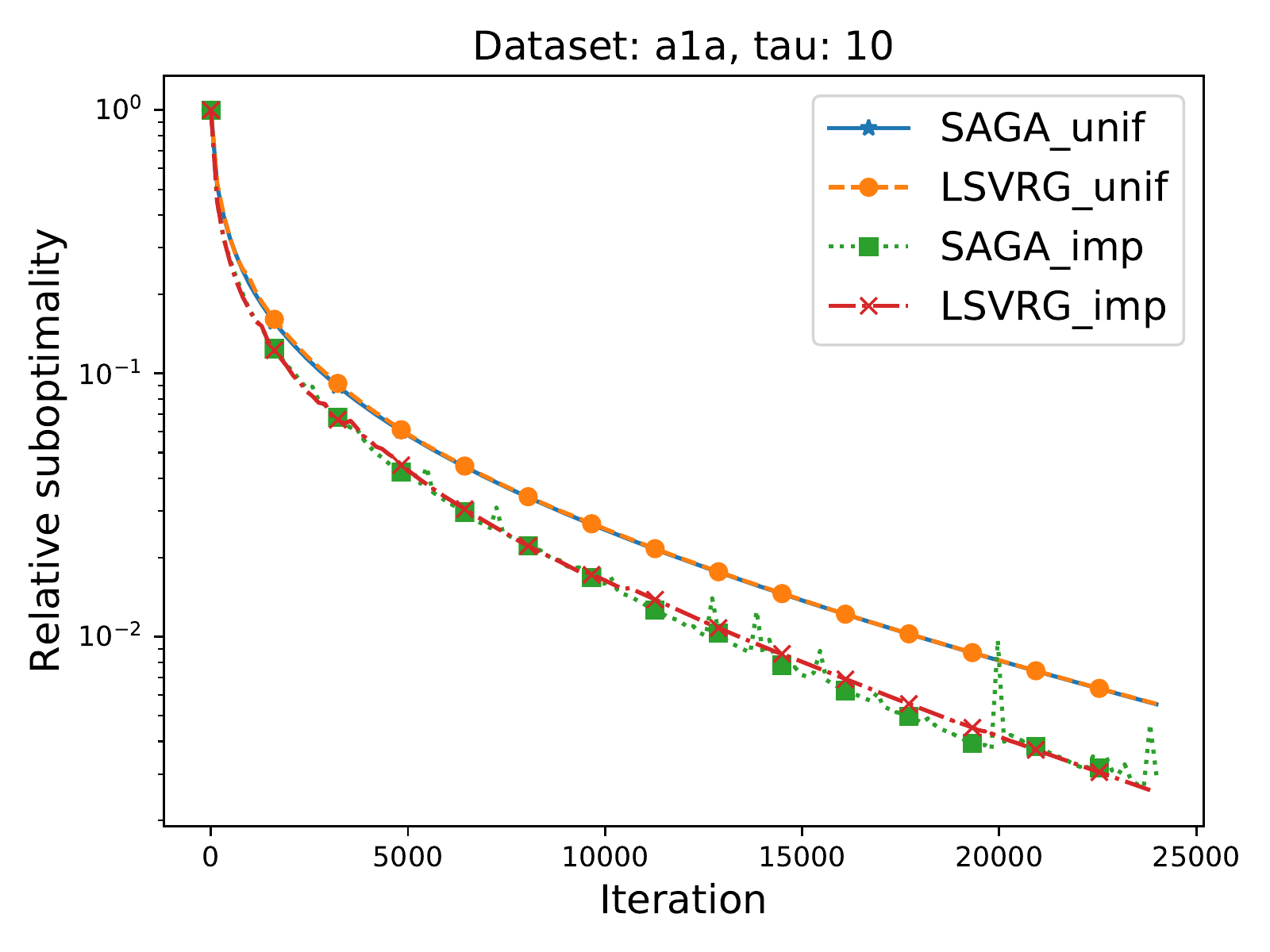}
\end{minipage}%
\begin{minipage}{0.3\textwidth}
  \centering
\includegraphics[width =  \textwidth ]{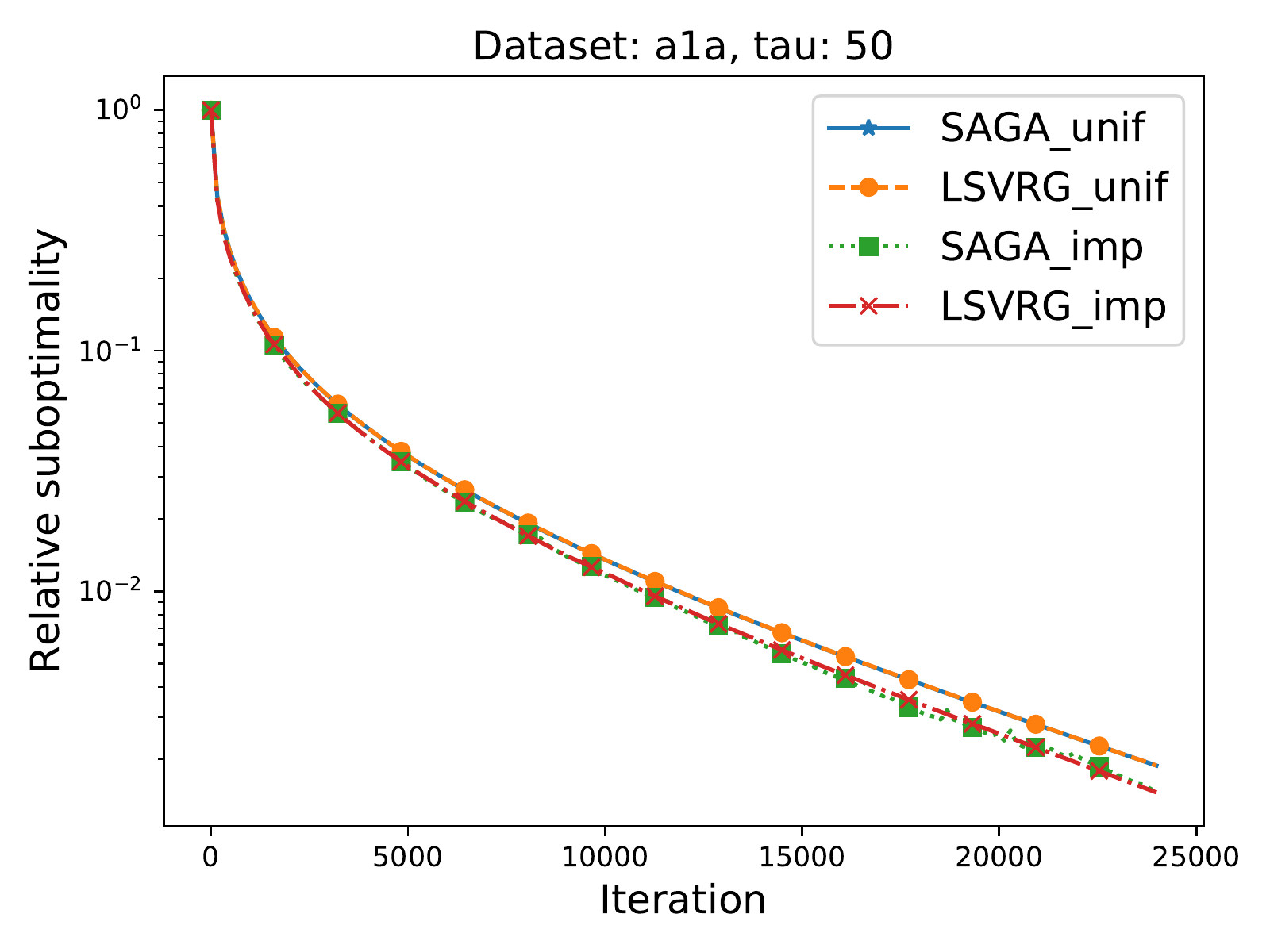}
\end{minipage}%
\\
\begin{minipage}{0.3\textwidth}
  \centering
\includegraphics[width =  \textwidth ]{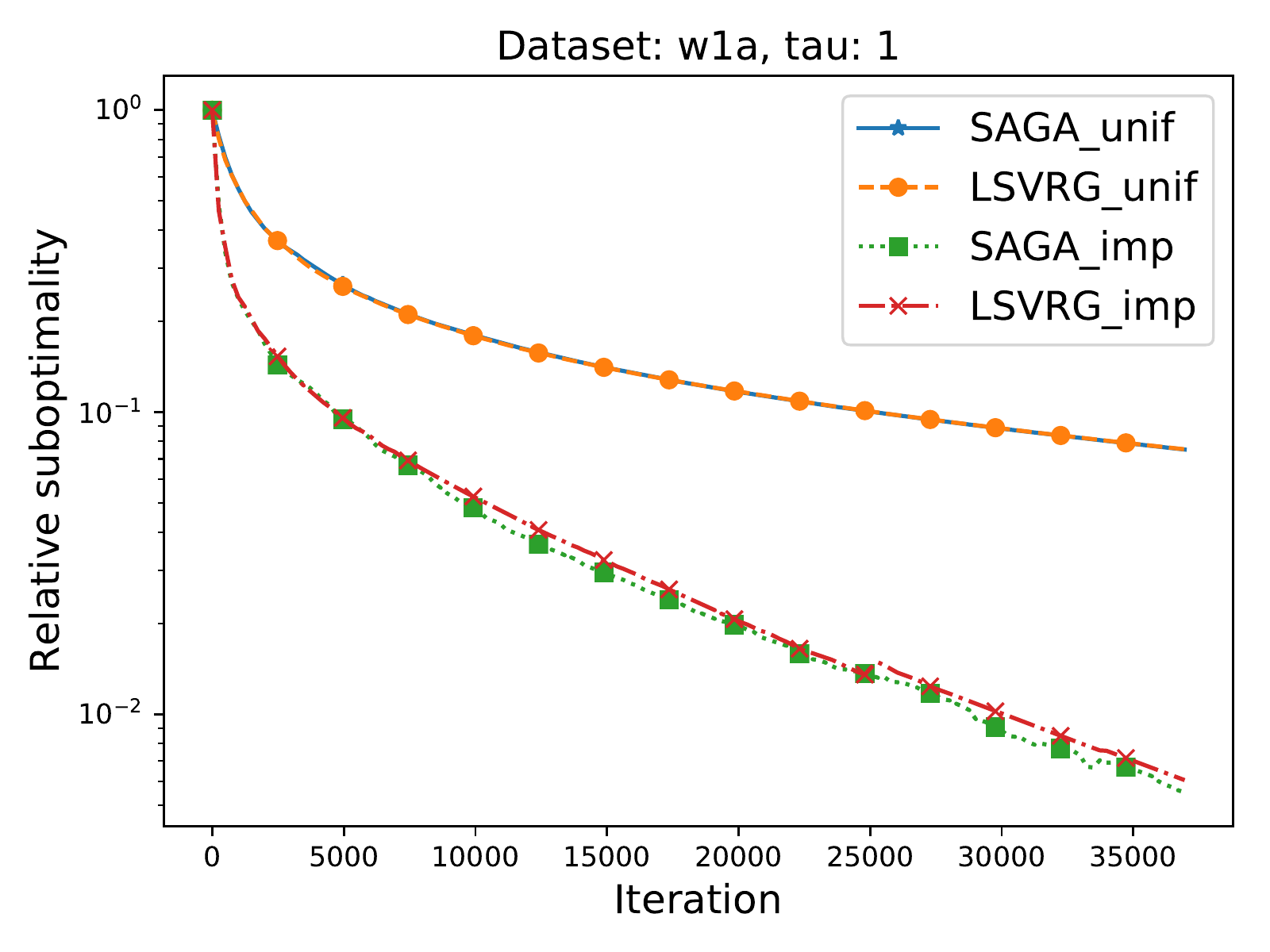}
\end{minipage}%
\begin{minipage}{0.3\textwidth}
  \centering
\includegraphics[width =  \textwidth ]{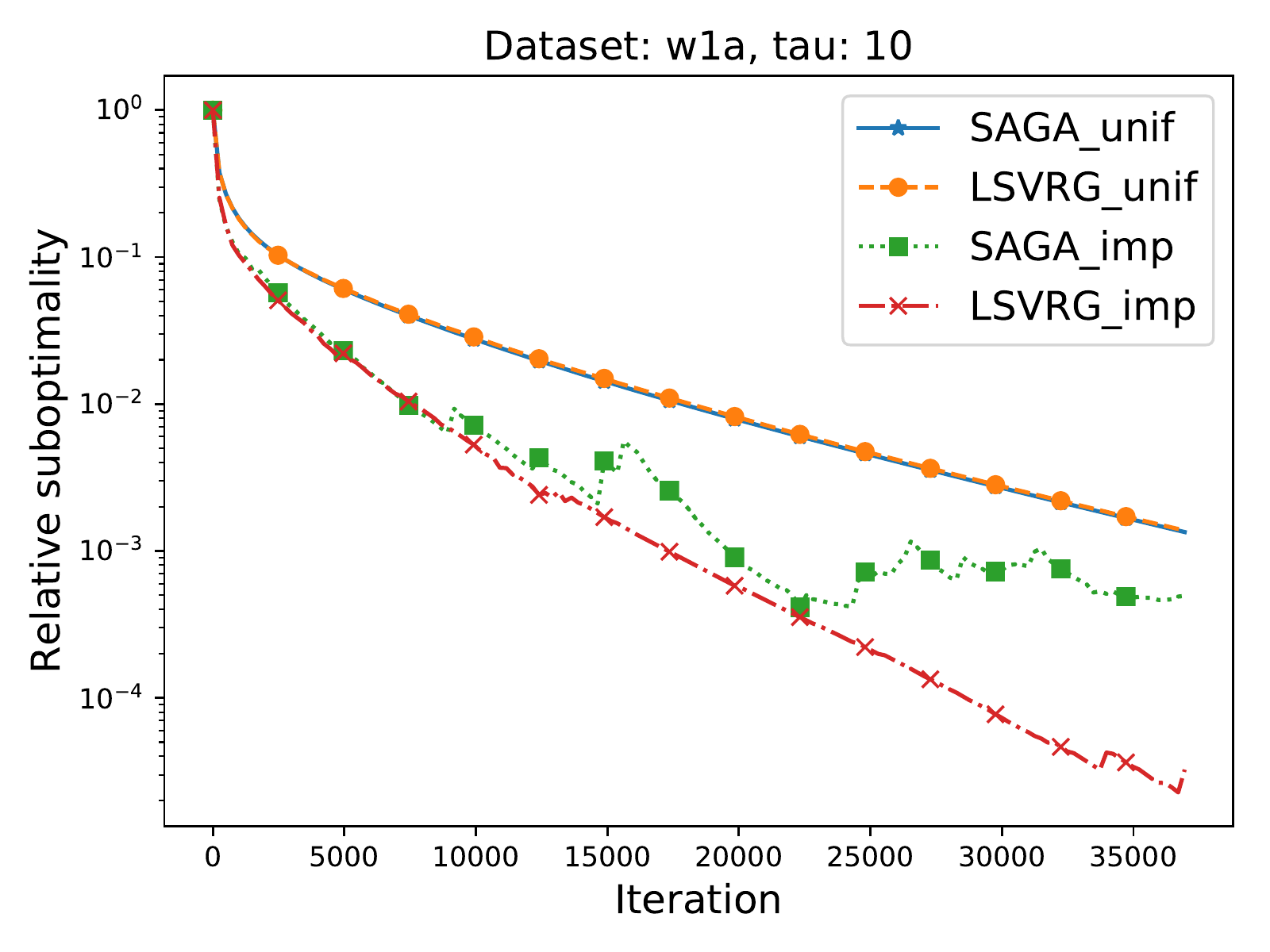}
\end{minipage}%
\begin{minipage}{0.3\textwidth}
  \centering
\includegraphics[width =  \textwidth ]{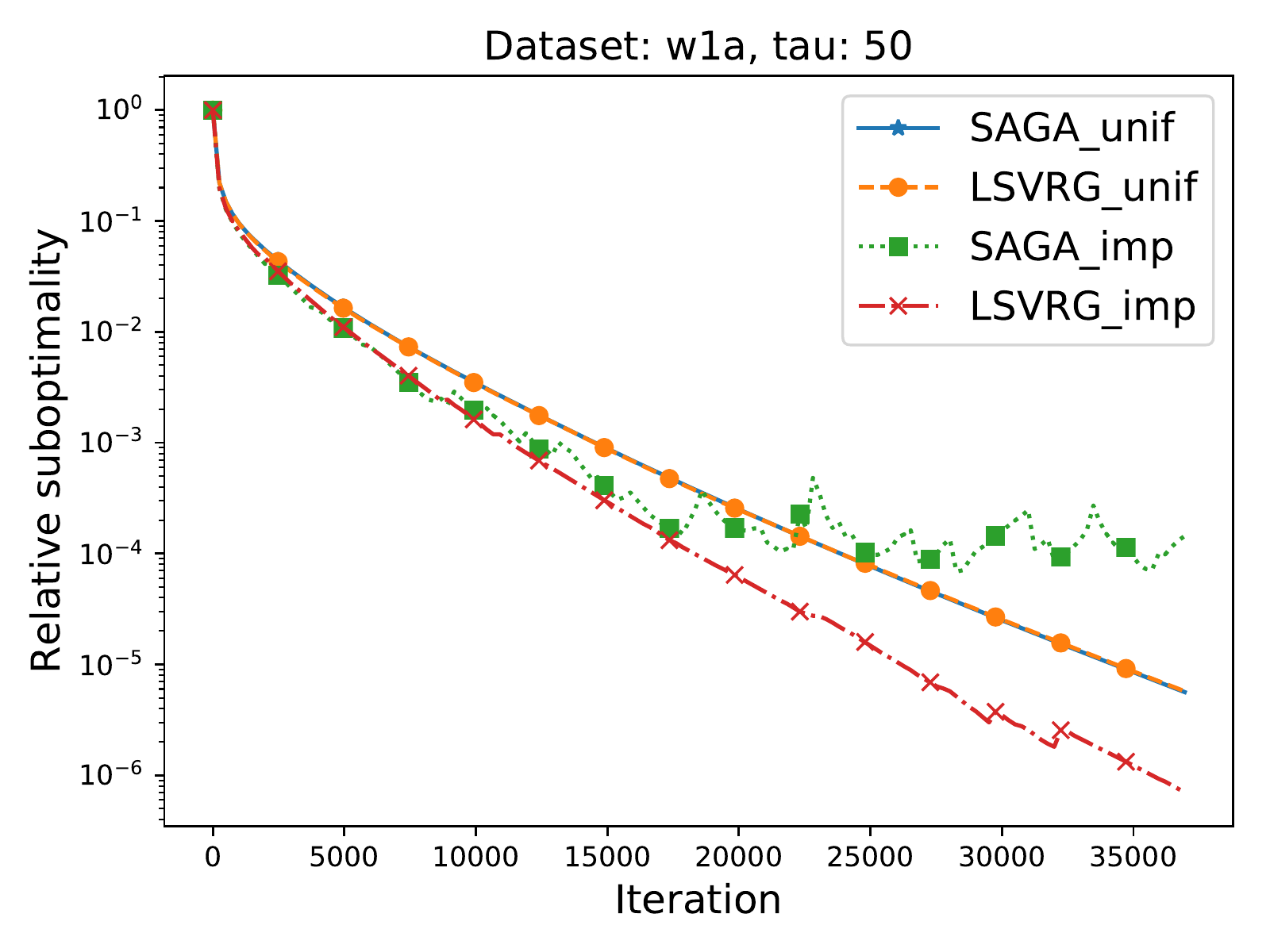}
\end{minipage}%
\\
\begin{minipage}{0.3\textwidth}
  \centering
\includegraphics[width =  \textwidth ]{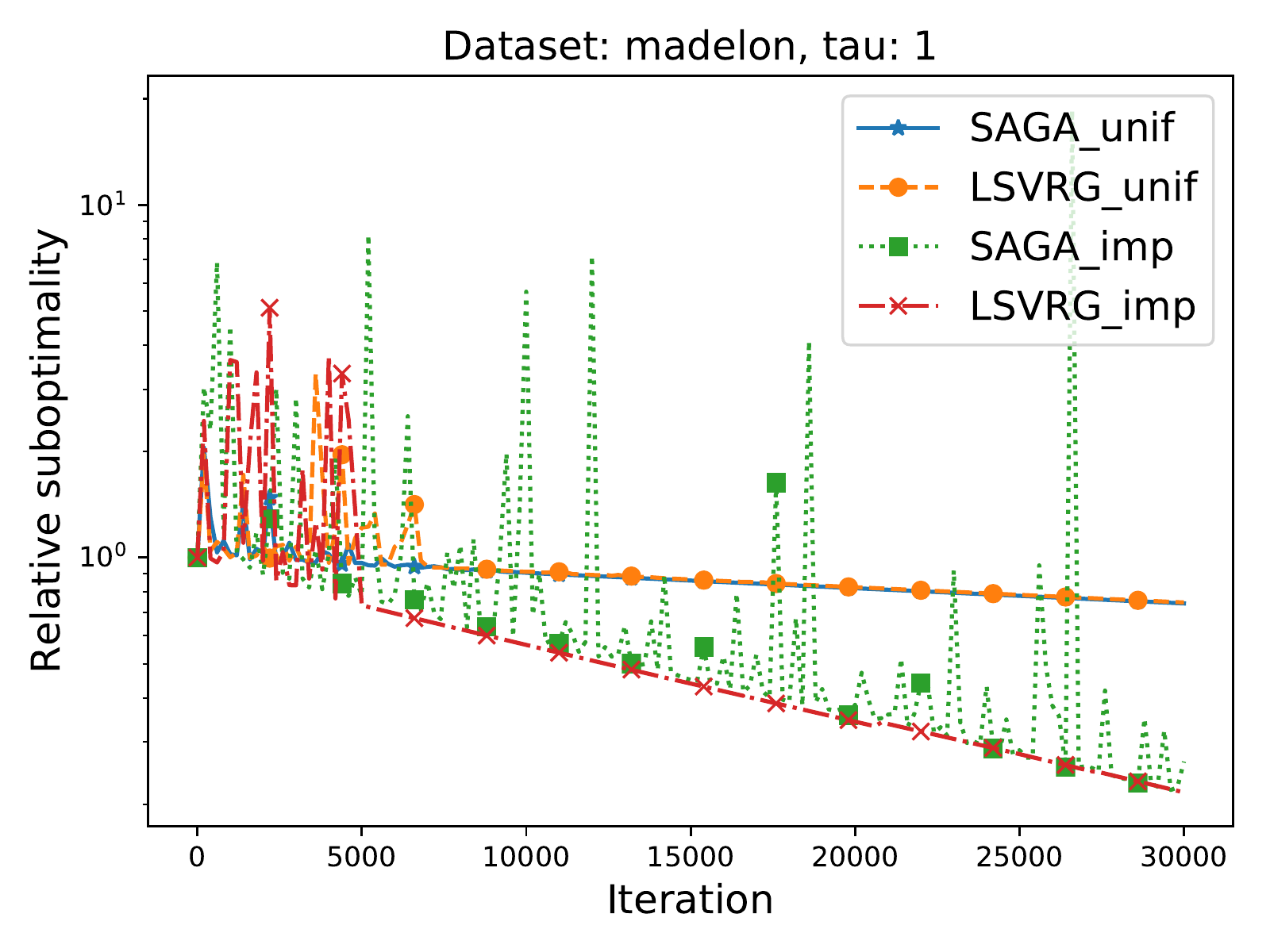}
\end{minipage}%
\begin{minipage}{0.3\textwidth}
  \centering
\includegraphics[width =  \textwidth ]{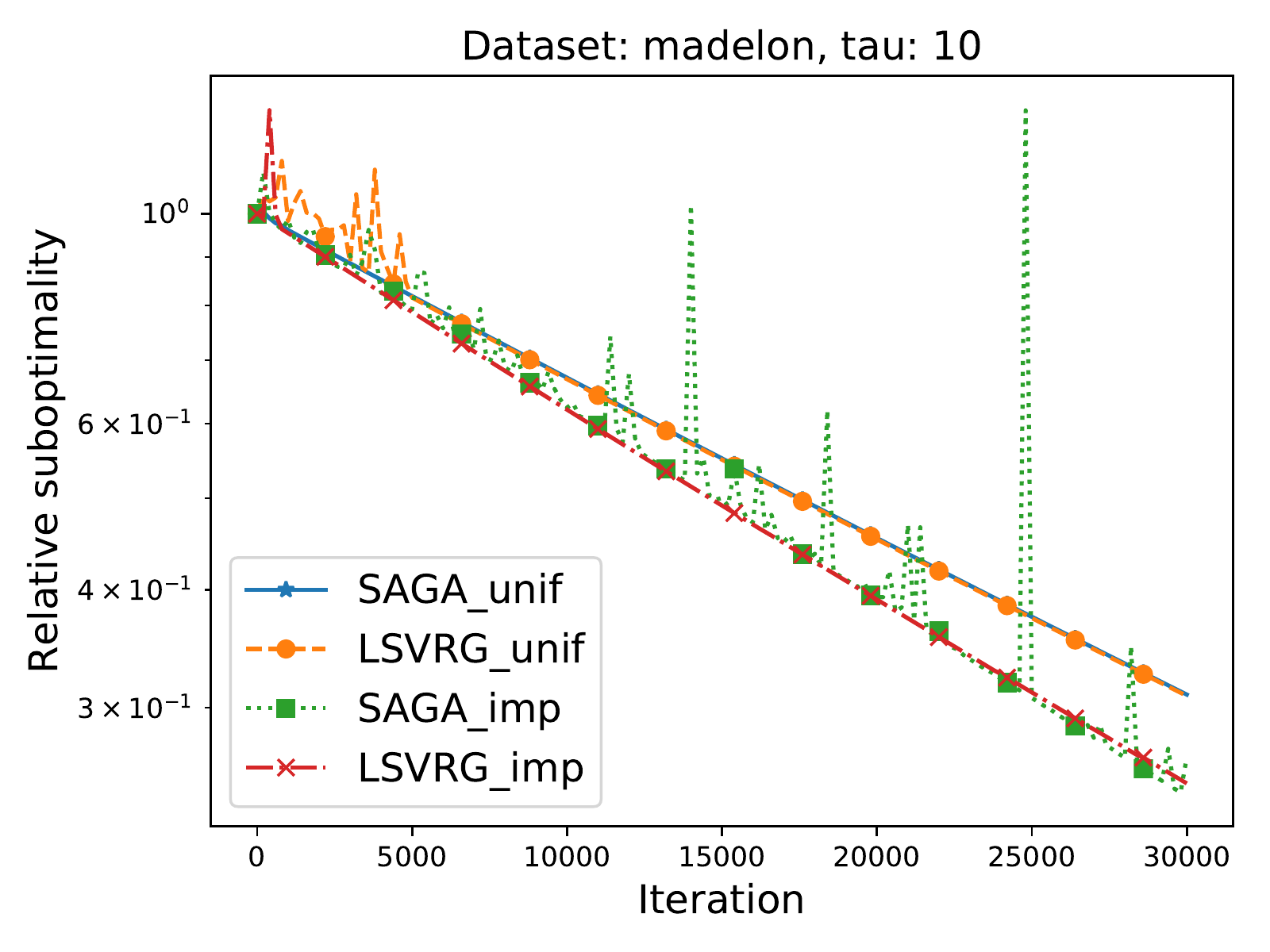}
\end{minipage}%
\begin{minipage}{0.3\textwidth}
  \centering
\includegraphics[width =  \textwidth ]{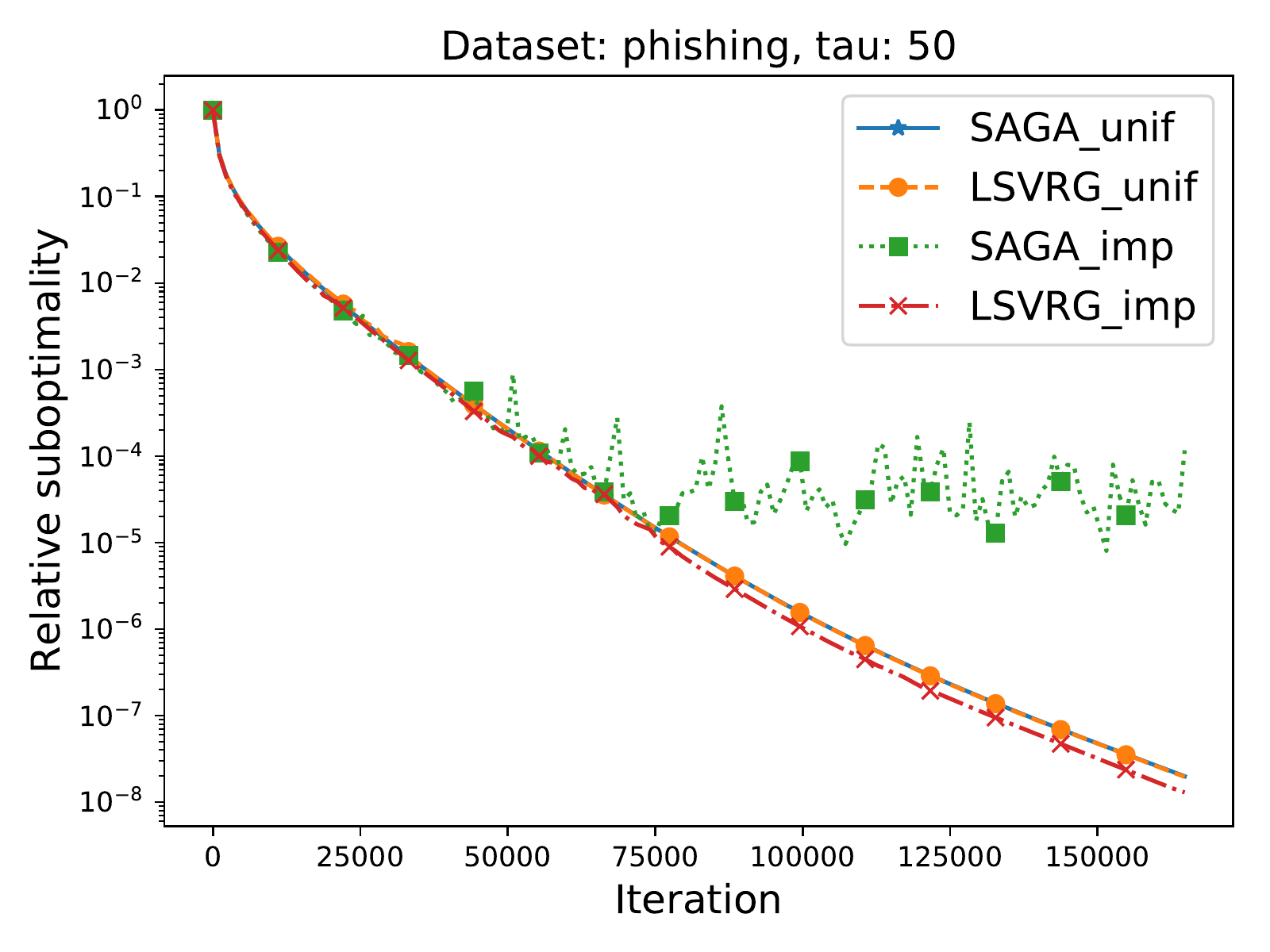}
\end{minipage}%
\\
\begin{minipage}{0.3\textwidth}
  \centering
\includegraphics[width =  \textwidth ]{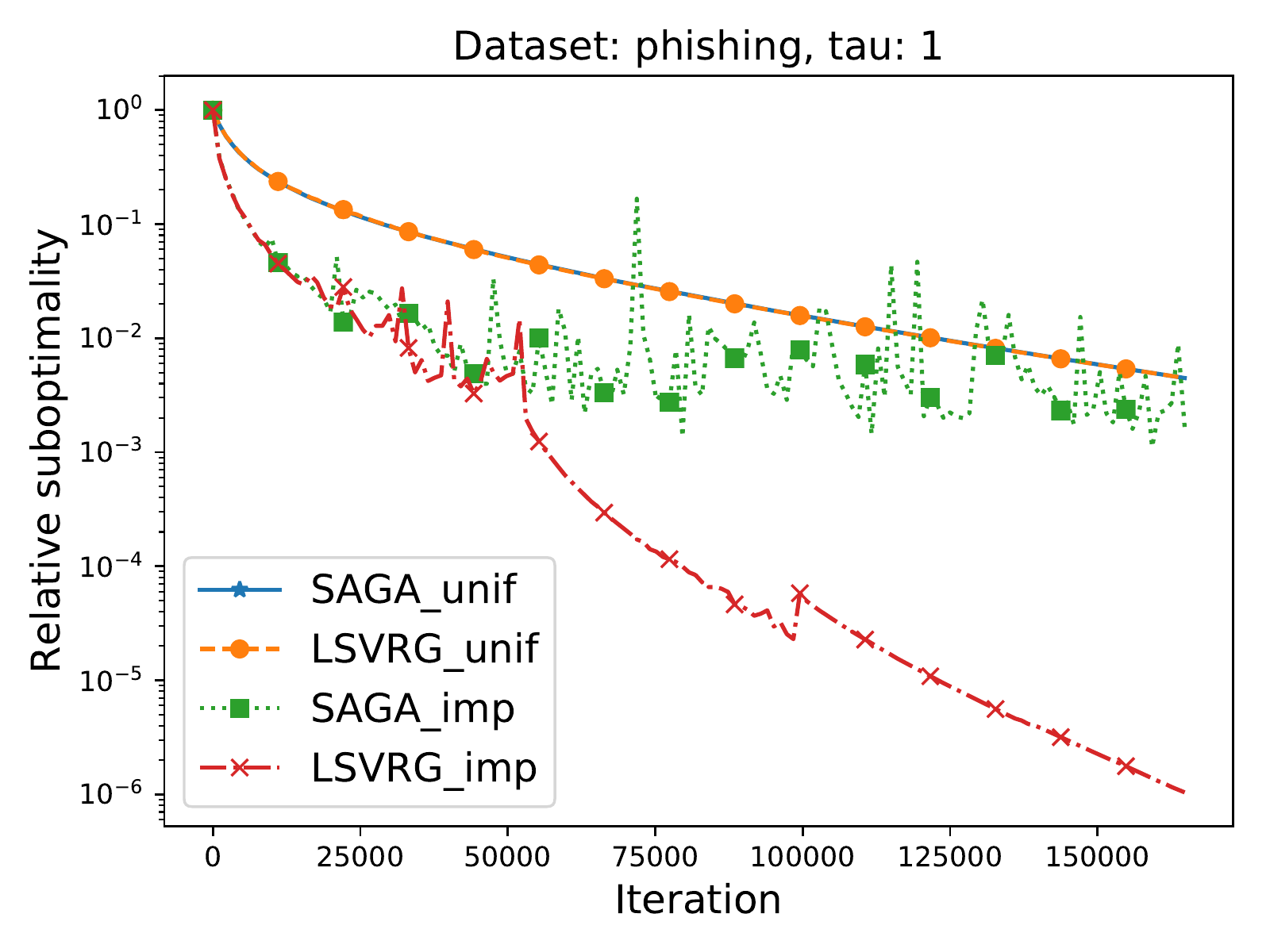}
\end{minipage}%
\begin{minipage}{0.3\textwidth}
  \centering
\includegraphics[width =  \textwidth ]{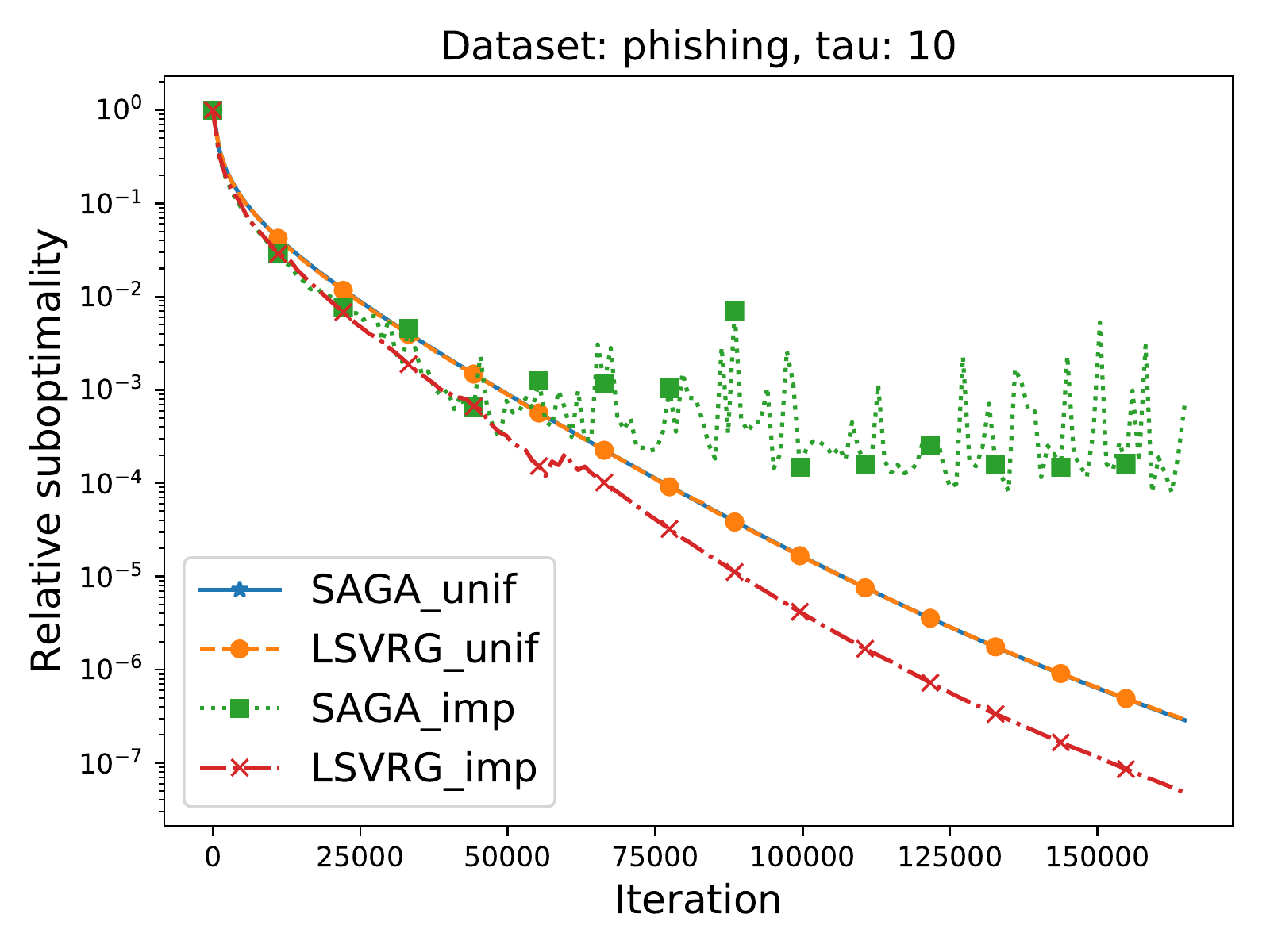}
\end{minipage}%
\begin{minipage}{0.3\textwidth}
  \centering
\includegraphics[width =  \textwidth ]{LSVRGDatasetphishingtau50meth3.pdf}
\end{minipage}%
\\
\begin{minipage}{0.3\textwidth}
  \centering
\includegraphics[width =  \textwidth ]{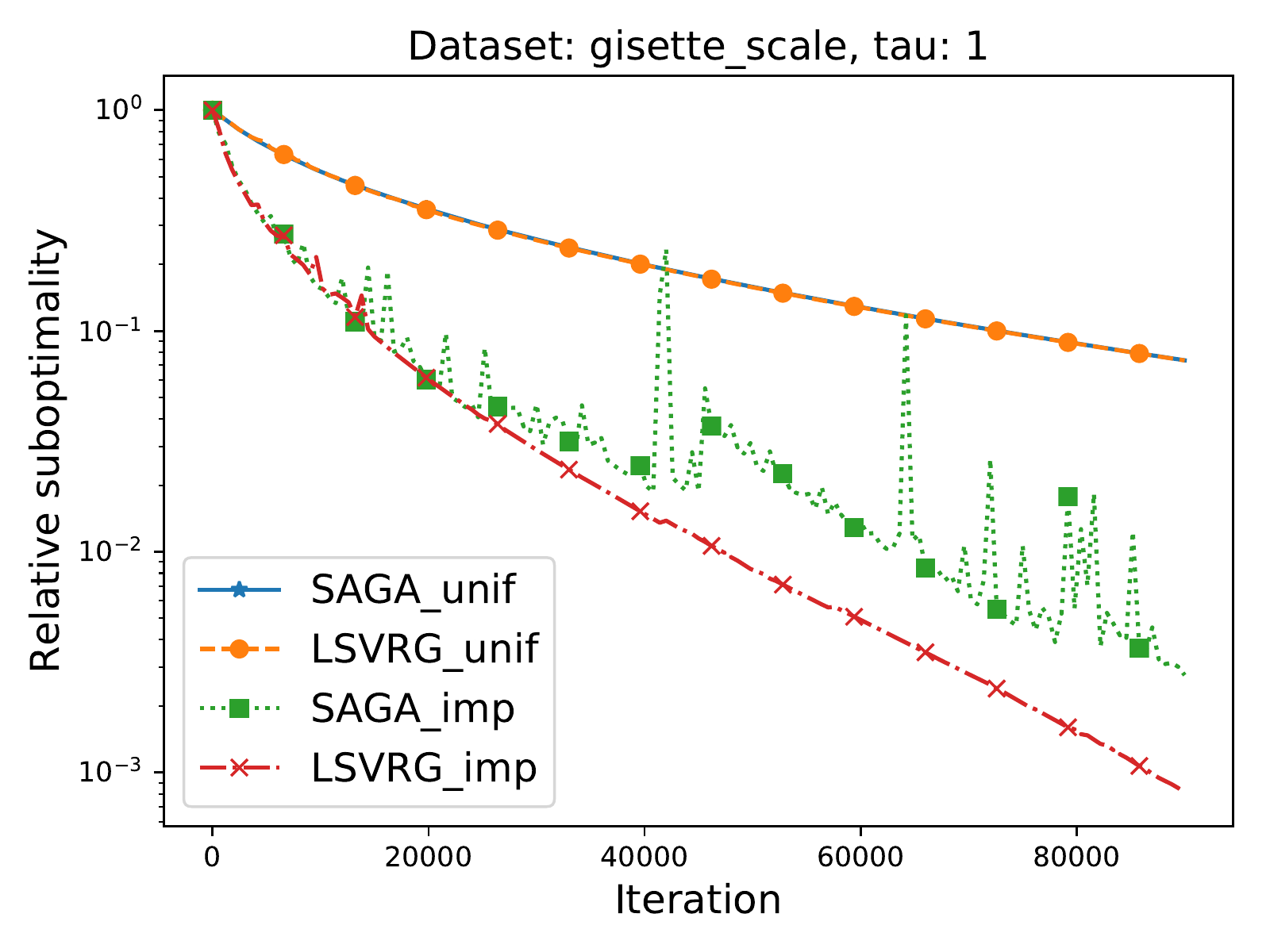}
\end{minipage}%
\begin{minipage}{0.3\textwidth}
  \centering
\includegraphics[width =  \textwidth ]{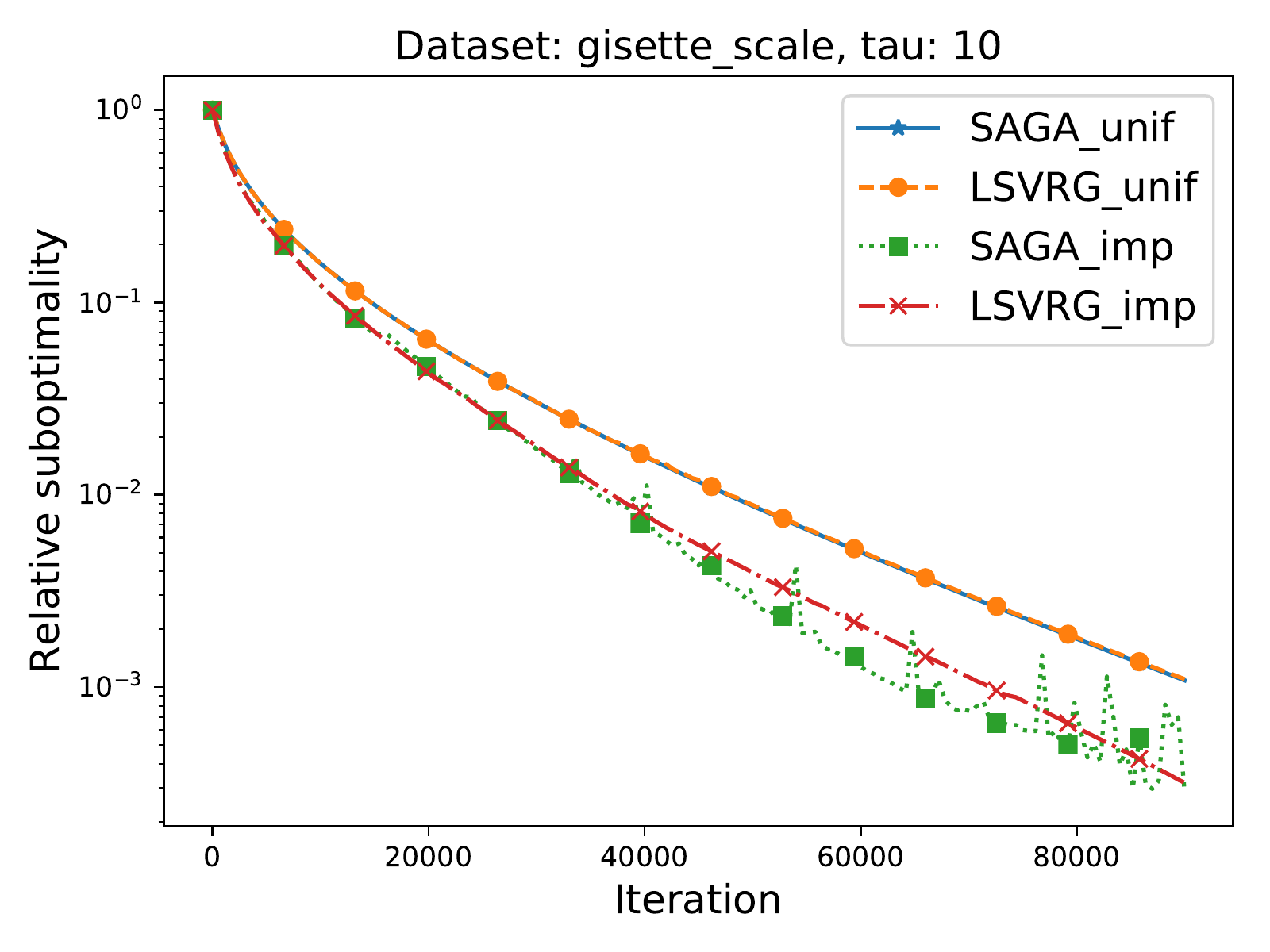}
\end{minipage}%
\begin{minipage}{0.3\textwidth}
  \centering
\includegraphics[width =  \textwidth ]{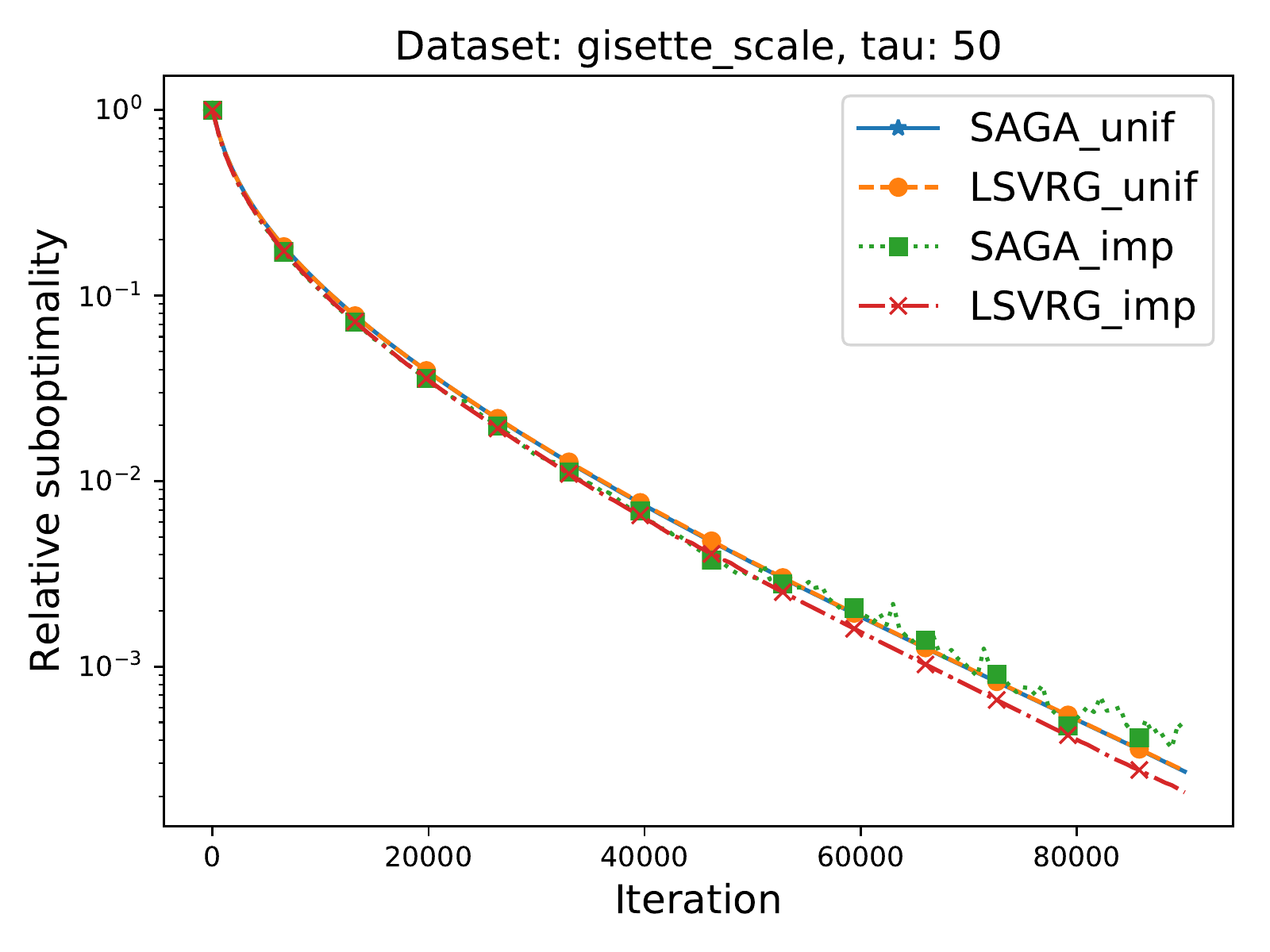}
\end{minipage}
\caption{{\tt LSVRG} applied on LIBSVM~\cite{chang2011libsvm} datasets with $\lambda = 10^{-5}$. Axis $y$ stands for relative suboptimality, i.e. $\frac{f(x^k)-f(x^*)}{f(x^k)-f(x^0)}$.}
\label{fig:gjs_LSVRG1}
\end{figure}

\begin{figure}[!h]
\centering
\begin{minipage}{0.3\textwidth}
  \centering
\includegraphics[width =  \textwidth ]{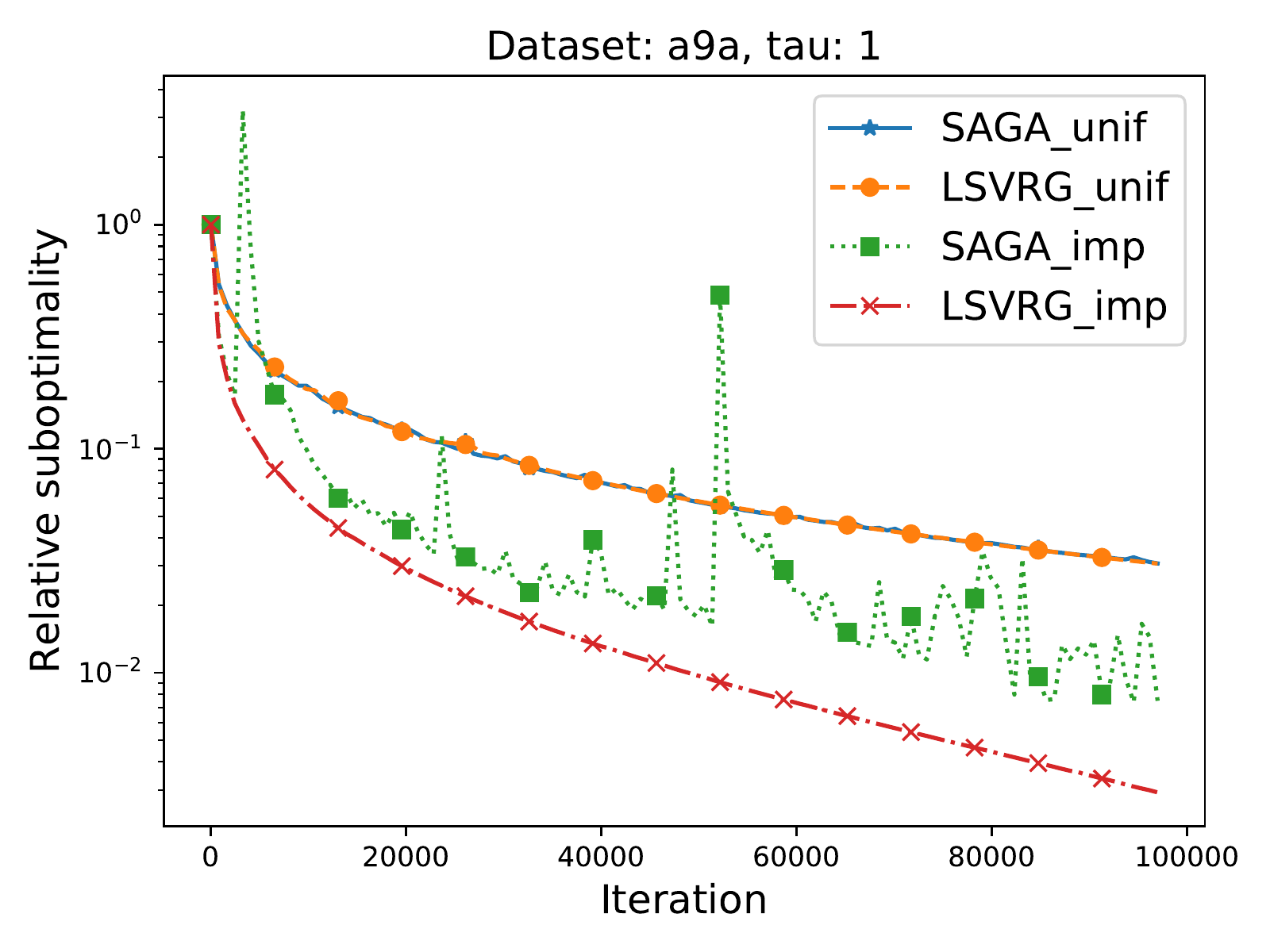}
\end{minipage}%
\begin{minipage}{0.3\textwidth}
  \centering
\includegraphics[width =  \textwidth ]{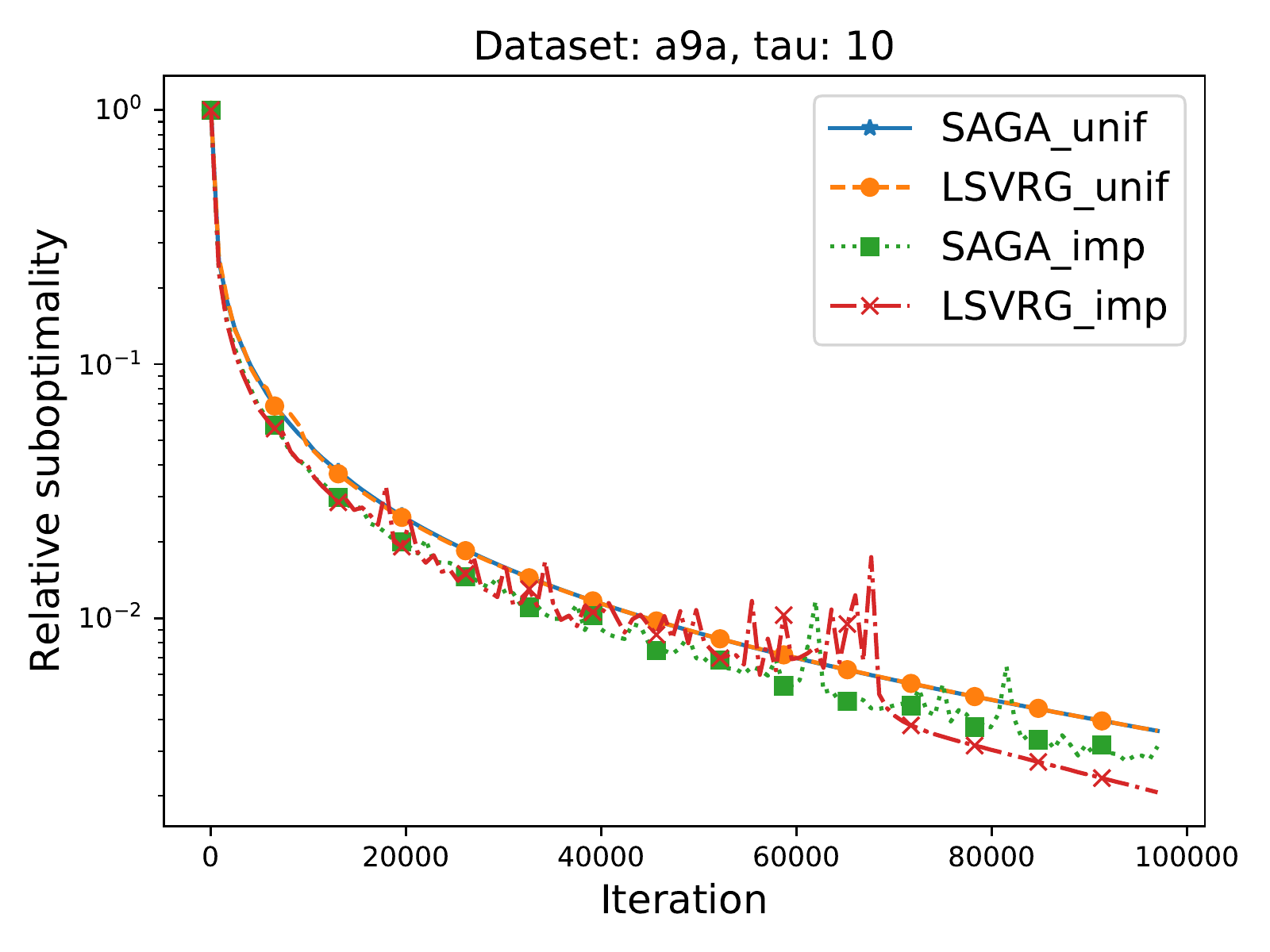}
\end{minipage}%
\begin{minipage}{0.3\textwidth}
  \centering
\includegraphics[width =  \textwidth ]{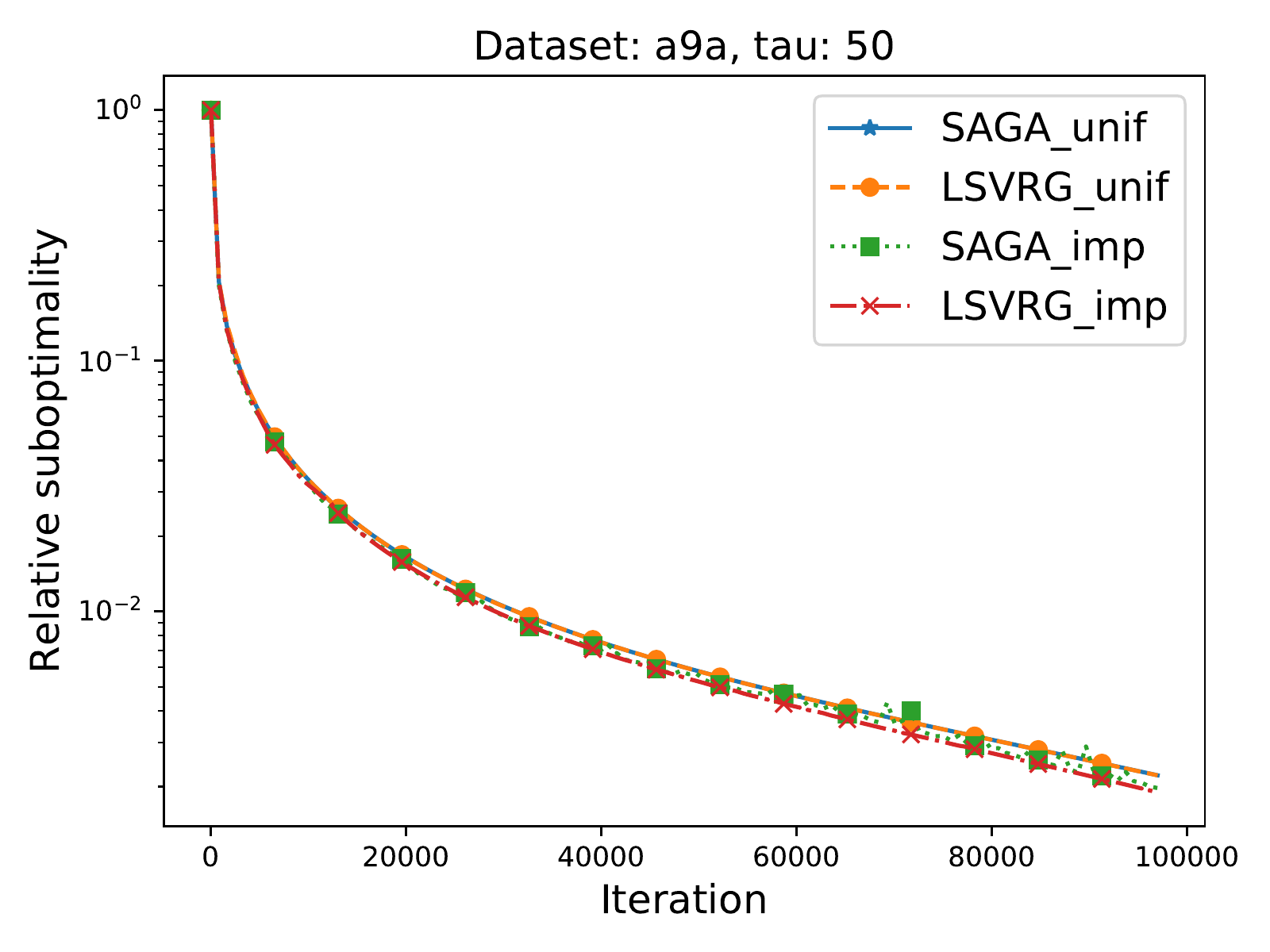}
\end{minipage}
\\
\begin{minipage}{0.3\textwidth}
  \centering
\includegraphics[width =  \textwidth ]{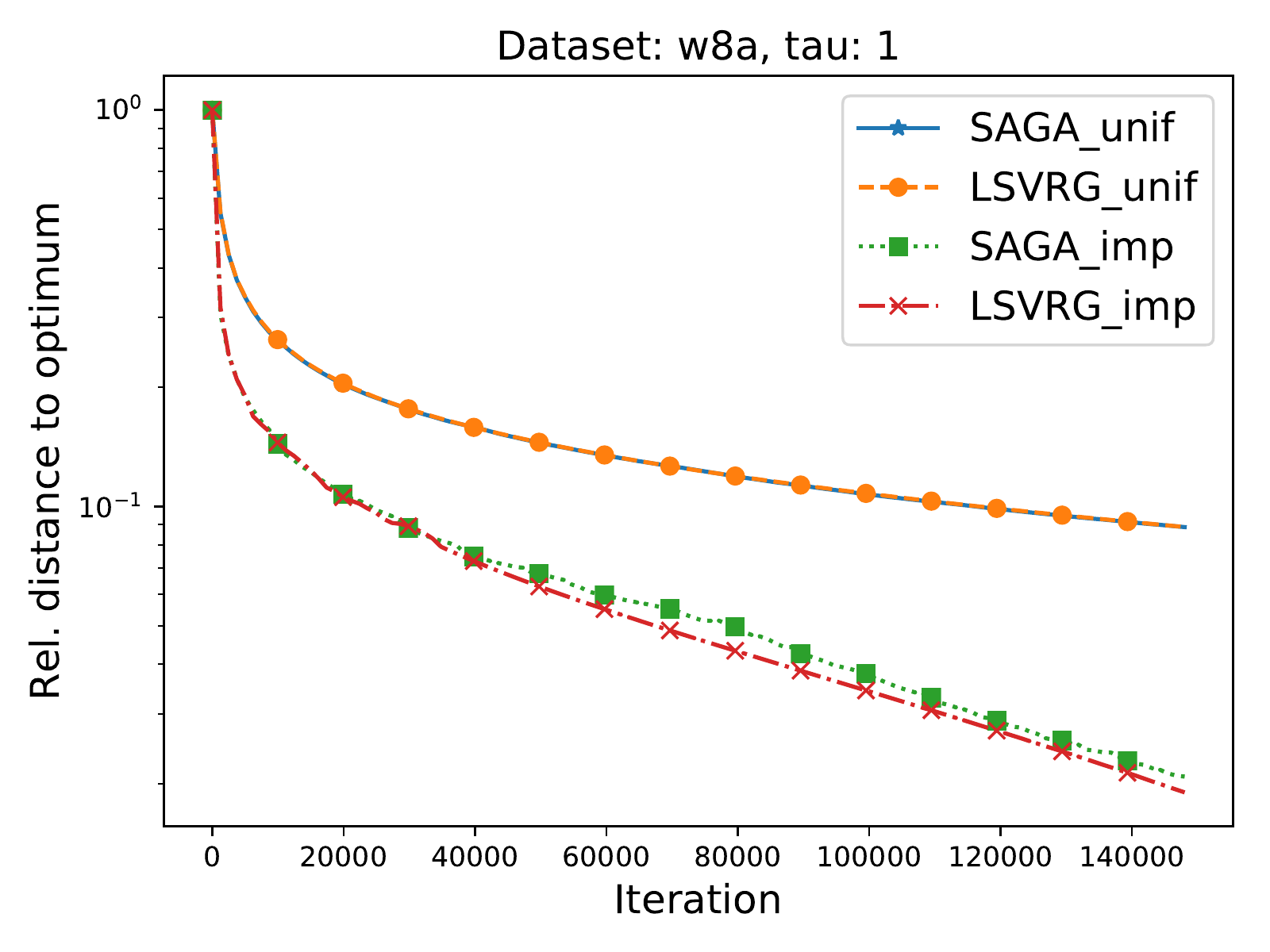}
\end{minipage}%
\begin{minipage}{0.3\textwidth}
  \centering
\includegraphics[width =  \textwidth ]{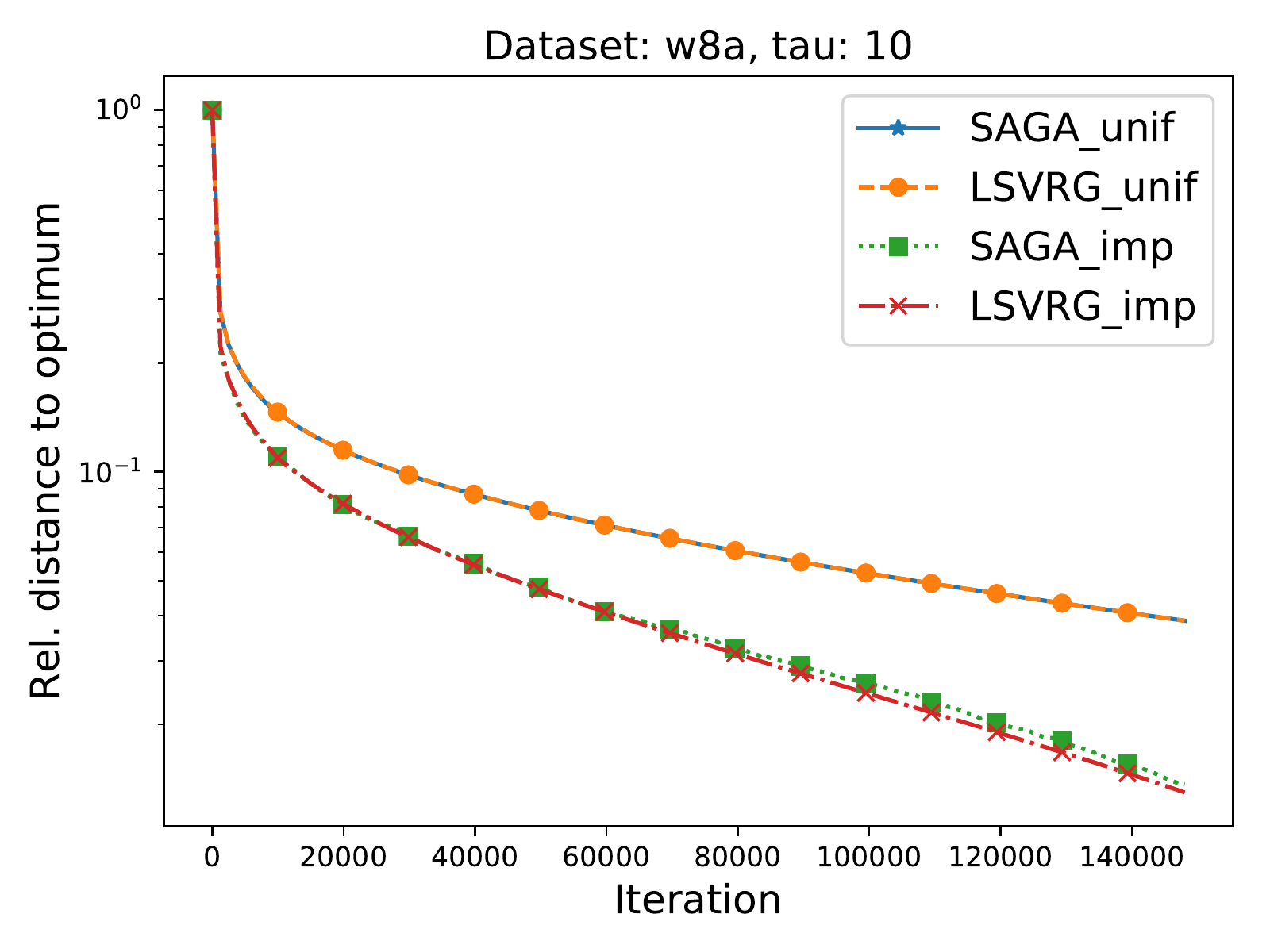}
\end{minipage}%
\begin{minipage}{0.3\textwidth}
  \centering
\includegraphics[width =  \textwidth ]{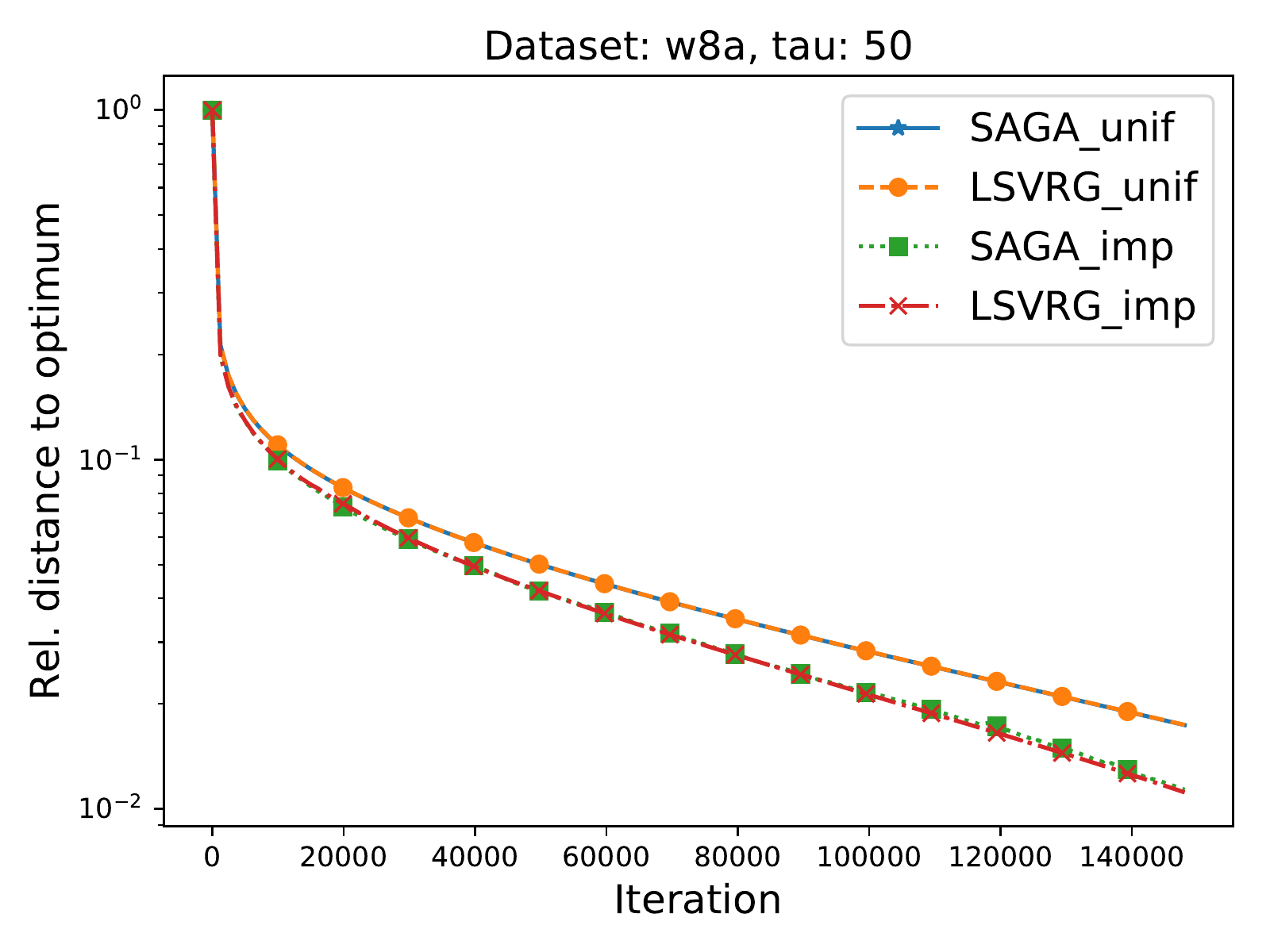}
\end{minipage}
\caption{{\tt LSVRG} applied on LIBSVM~\cite{chang2011libsvm} datasets. For {\tt a9a},  $\lambda = 0$ and $\probx = \frac1n$ was chosen; for {\tt w8a}, $\lambda = 10^{-8}$ and $\probx = \frac3n$ was chosen. Axis $y$ stands for relative suboptimality, i.e. $\frac{f(x^k)-f(x^*)}{f(x^k)-f(x^0)}$.
}
\label{fig:gjs_LSVRG2}
\end{figure}

\section{Conclusion}

In this chapter we proposed a fairly general algorithm---{\tt GJS}---capable of inserting the variance reduction mechanism under arbitrary random first-order oracle. Each special case either recovers a known algorithm with its tight rate, or improves a known algorithm or is a new algorithm. In the next chapter we go even further: we introduce a general technique to analyze unbiased stochastic gradient algorithms that are not necessarily variance reduced.

\chapter{A Unified Theory of {\tt SGD}: Variance Reduction, Sampling, Quantization  and Coordinate Descent}
\label{sigmak}

\graphicspath{{sigmak/images/}}

In Chapter~\ref{jacsketch}, we have proposed a general variance-reduced algorithm applicable in many different scenarios. In this chapter, we go a step further. In particular, we propose a new generic analysis technique capable of providing complexity bounds for a significantly broader class of stochastic gradient algorithms. 

\paragraph{Stochastic optimization.} 
In this chapter we are primarily concerned with regularized stochastic optimization problems of the form
    \begin{equation}\label{eq:sk_problem_gen}
        \min_{x\in\R^d} \EEE_{\xi \sim \cD} \left[ f_\xi(x)\right] + \psi(x),
    \end{equation}
    and let    \begin{equation} \label{eq:sk_f_exp}
    f(x)=\EEE_{\xi \sim \cD} \left[ f_\xi(x)\right] .\end{equation} 
    
  As usual, function $f$ is assumed to be convex, differentiable with Lipschitz gradient, and $\psi:\R^d\to \R\cup \{+\infty\}$ is a proximable (proper closed convex) regularizer. Specifically for this section, we assume that $\xi$ is a random variable, and  $f_\xi:\R^d\to \R$  is smooth function for all $\xi$. 
  
  Stochastic optimization problems are of key importance in statistical supervised learning theory. In this setup, $x$ represents a machine learning model described by $d$ parameters (e.g., logistic regression or a deep neural network), $\cD$  is an unknown distribution of labelled examples,  $f_\xi(x)$ represents the loss of model $x$ on datapoint $\xi$, and $f$ is the generalization error. Problem \eqref{eq:sk_problem_gen} seeks to find the model $x$ minimizing the generalization error. In statistical learning theory one assumes that while $\cD$ is not known, samples $\xi\sim \cD$ are available. In such a case, $\nabla f(x)$ is not computable, while $\nabla f_\xi(x)$, which is  an unbiased estimator of the gradient of $f$ at $x$, is easily computable. 

\paragraph{Finite-sum problems.} 
Another prominent example, one of special interest in this work, are functions $f$ which arise as averages of a very large number of smooth functions: \begin{equation}\label{eq:sk_f_sum} 
\compactify f(x)=\frac{1}{n}\sum \limits_{i=1}^n f_i(x).\end{equation} This problem often arises by approximation of the stochastic optimization loss function \eqref{eq:sk_f_exp} via Monte Carlo integration, and is in this context known as the empirical risk minimization (ERM) problem.  ERM is currently the dominant paradigm for solving supervised learning problems \cite{shai_book}.  If index $i$ is chosen uniformly at random from $[n]\eqdef \{1,2,\dots,n\}$, $\nabla f_i(x)$ is an unbiased estimator of $\nabla f(x)$. Typically, $\nabla f(x)$ is about $n$ times more expensive to compute than $\nabla f_i(x)$. 

\paragraph{Distributed optimization.} 
Lastly, in some applications, especially in distributed training of supervised models, one considers problem \eqref{eq:sk_f_sum}, with $n$ being the number of machines, and each $f_i$ also having a finite sum structure, i.e.,
\begin{equation}\label{eq:sk_f_i_sum}
\compactify    f_i(x) = \frac{1}{m}\sum\limits_{j=1}^m  f_{ij}(x),
\end{equation}
where $m$ corresponds to the number of training examples stored on machine $i$.

\section{The many faces of stochastic gradient descent} \label{sec:sk_many_faces_of_SGD}

 Stochastic gradient descent ({\tt SGD}) \cite{robbins, nemirovski2009robust, vaswani2019-overparam} is a state-of-the-art algorithmic paradigm for solving optimization problems \eqref{eq:sk_problem_gen} in situations when $f$ is either of structure 
 \eqref{eq:sk_f_exp} or \eqref{eq:sk_f_sum}. In its generic form, (proximal) {\tt SGD} defines the new iterate by subtracting a multiple of a stochastic gradient from the current iterate, and subsequently applying the proximal operator of $\psi$:
\begin{equation} \label{eq:sk_SGD} x^{k+1} = \proxR(x^k - \alpha g^k).\end{equation} 
Here, $g^k$ is an unbiased estimator of the gradient (i.e., a stochastic gradient), \begin{equation}\label{eq:sk_stoch_grad}\EEE \left[g^k \;|\; x^k\right] = \nabla f(x^k),\end{equation}
and $\proxR(x) = \argmin_y \{\alpha \psi(y) + \frac{1}{2}\norm{y-x}^2\}$. However, and this is the starting point of our journey in this chapter, there are {\em infinitely many} ways of obtaining a random vector $g^k$ satisfying \eqref{eq:sk_stoch_grad}. On the one hand, this gives algorithm designers the flexibility to {\em construct}  stochastic gradients in various ways in order to target desirable properties such as convergence speed, iteration cost, parallelizability and generalization.  On the other hand, this poses considerable challenges in terms of convergence analysis. Indeed, if one aims to, as one should, obtain the sharpest bounds possible, dedicated analyses are needed to handle each of the particular variants of {\tt SGD}.

\paragraph{Vanilla {\tt SGD}.} The flexibility in the design of efficient strategies for constructing $g^k$ has led to a creative renaissance in the optimization and machine learning communities, yielding a large number of immensely powerful new variants\footnote{In this chapter, by {\em vanilla} {\tt SGD} we refer to {\tt SGD} variants with or without importance sampling and mini-batching, but {\em excluding} variance-reduced variants, such as {\tt SAGA} \cite{saga} and {\tt SVRG} \cite{svrg}.} of {\tt SGD}, such as those employing  {\em importance sampling} \cite{iprox-sdca, needellward2015}, and {\em mini-batching} \cite{mS2GD}. These efforts are subsumed by the recently developed and remarkably sharp analysis of {\tt SGD} under  {\em arbitrary sampling} paradigm \cite{pmlr-v97-qian19b}, first introduced in the study of randomized coordinate descent methods by \cite{nsync}. The arbitrary sampling paradigm covers virtually all stationary mini-batch and importance sampling strategies in a unified way, thus making headway towards theoretical unification of two separate strategies for constructing stochastic gradients. For strongly convex $f$, the {\tt SGD} methods analyzed in \cite{pmlr-v97-qian19b} converge linearly to a neighbourhood of the solution $x^* = \argmin_x f(x)$ for a fixed stepsize $\alpha^k=\alpha$. The size of the neighbourhood is proportional to the second moment of the stochastic gradient at the optimum ($\sigma^2 \eqdef \frac{1}{n}\sum_{i=1}^n \norm{\nabla f_i(x^*)}^2$), to the stepsize ($\alpha$), and inversely proportional to the modulus of strong convexity. The effect of various sampling strategies, such as importance sampling and mini-batching, is twofold: i) improvement of the linear convergence rate by enabling larger stepsizes, and ii) modification of $\sigma^2$. However, none of these strategies\footnote{Except for the full batch strategy, which is prohibitively expensive.} is able to completely eliminate the adverse effect of $\sigma^2$. That is,  {\tt SGD} with a fixed stepsize does not reach the optimum, unless one happens to be in the overparameterized case characterized by the identity $\sigma^2=0$.

\paragraph{Variance reduced {\tt SGD}.} While sampling strategies such as importance sampling and mini-batching reduce the variance of the stochastic gradient, in the finite-sum case \eqref{eq:sk_f_sum} a new type of {\em variance reduction} strategies has been developed over the last few years \cite{sag, saga, svrg, sdca, quartz,nguyen2017sarah, kovalev2019don, horvath2020adaptivity} (see also~Chapter~\ref{jacsketch}). These variance-reduced {\tt SGD} methods differ from the sampling strategies discussed before in a significant way: they can iteratively {\em learn} the stochastic gradients at the optimum,  and in so doing are able to eliminate the adverse effect of the gradient noise $\sigma^2>0$ which, as mentioned above, prevents the iterates of vanilla {\tt SGD} from converging to the optimum. As a result, for strongly convex $f$, these new variance-reduced {\tt SGD} methods converge linearly to $x^*$, with a fixed stepsize. At the moment, these variance-reduced variants require a markedly different convergence theory from the vanilla variants of {\tt SGD}. An exception to this is the situation when $\sigma^2=0$ as then variance reduction is not needed; indeed, vanilla {\tt SGD} already converges to the optimum, and with a fixed stepsize. We end the discussion here by remarking that this {\em hints} at a possible existence of a more unified theory, one that would include both vanilla and variance-reduced {\tt SGD}.

\paragraph{Distributed {\tt SGD}, quantization and variance reduction.} When {\tt SGD} is implemented in a distributed fashion,  the problem is often expressed in the form \eqref{eq:sk_f_sum}, where $n$ is the number of workers/nodes, and  $f_i$ corresponds to the loss based on data stored on node $i$. Depending on the number of data points stored on each node, it may or may not be efficient to compute the gradient of $f_i$ in each iteration. In general, {\tt SGD} is implemented in this way: each  node $i$ first computes a stochastic gradient $g_i^k$ of $f_i$ at the current point $x^k$ (maintained individually by each node). These gradients are then aggregated by a master node \cite{dane, RDME}, in-network by a switch~\cite{switchML}, or a different technique best suited to the architecture used. To alleviate the communication bottleneck, various lossy update compression strategies such as quantization~\cite{1bit, Gupta:2015limited, zipml}, sparsification~\cite{RDME, alistarh2018convergence, wangni2018gradient} and dithering~\cite{alistarh2017qsgd} were proposed. The basic idea is for each worker to apply a randomized transformation $Q:\R^d\to \R^d$ to $g_i^k$, resulting in a vector which is still an unbiased estimator of the gradient, but one that can be communicated with fewer bits. Mathematically, this amounts to injecting additional noise into the already noisy stochastic gradient $g_i^k$. The field of quantized {\tt SGD} is still young, and even some basic questions remained open until recently. For instance, there was no distributed quantized {\tt SGD} capable of provably solving  \eqref{eq:sk_problem_gen} until the {\tt DIANA} algorithm \cite{mishchenko2019distributed} was introduced. {\tt DIANA} applies quantization to {\em gradient differences}, and in so doing is able to learn the gradients at the optimum, which makes it able to work for any regularizer $\psi$. {\tt DIANA} has some structural similarities with {\tt SEGA}~\cite{sega}---the first coordinate descent type method which works for non-separable regularizers---but a more precise relationship remains elusive.  When the functions of $f_i$ are of a finite-sum structure as in \eqref{eq:sk_f_i_sum}, one can apply variance reduction to reduce the variance of the stochastic gradients $g_i^k$ together with quantization, resulting in the {\tt VR-DIANA} method~\cite{horvath2019stochastic}. This is the first distributed quantized {\tt SGD} method which provably converges to the solution of \eqref{eq:sk_problem_gen}+\eqref{eq:sk_f_i_sum} with a fixed stepsize.


\paragraph{Randomized coordinate descent ({\tt RCD}).} Lastly, in a distinctly separate strain, there are {\tt SGD} methods for the coordinate/subspace descent variety~\cite{rcdm}. While it is possible to see  {\em some} {\tt RCD} methods as special cases of \eqref{eq:sk_SGD}+\eqref{eq:sk_stoch_grad}, most of them do not follow this algorithmic template. First, standard {\tt RCD} methods use different stepsizes for updating different coordinates~\cite{qu2016coordinate1}, and this seems to be crucial to their success. Second,  until the recent discovery of the {\tt SEGA} method, {\tt RCD} methods were not able to converge with non-separable regularizers. Third, {\tt RCD} methods are naturally variance-reduced in the $\psi\equiv0$ case as partial derivatives at the optimum are all zero. As a consequence, attempts at creating variance-reduced {\tt RCD} methods seem to be futile. Lastly, {\tt RCD} methods are typically analyzed using different techniques. While there are deep links between standard {\tt SGD} and {\tt RCD} methods, these are often indirect and rely on duality \cite{sdca, face-off, sda}. 


\section{Contributions}

As outlined in the previous section, the world of {\tt SGD} is vast and beautiful. It is formed by many largely disconnected islands populated by elegant and efficient methods, with their own applications, intuitions, and convergence analysis techniques. While some links already exist (e.g., the unification of importance sampling and mini-batching variants under the arbitrary sampling umbrella), there is no comprehensive general theory. 
It is becoming increasingly difficult for the community to understand the relationships between these variants, both in theory and practice. New variants are yet to be discovered, but it is not clear what tangible principles one should adopt beyond intuition to aid the discovery. This situation is exacerbated by the fact that a number of different assumptions on the stochastic gradient, of various levels of strength, is being used in the literature.


The main contributions of this work include: 

\begin{itemize}
\item {\bf  Unified analysis.} In this work we propose a {\em unifying theoretical framework}
 which covers all of the  variants of {\tt SGD} outlined in Section~\ref{sec:sk_many_faces_of_SGD}. As a by-product, we obtain the {\em first unified analysis} of vanilla and variance-reduced {\tt SGD} methods.  For instance, our analysis covers as special cases vanilla {\tt SGD} methods from~\cite{nguyen2018sgd} and~\cite{pmlr-v97-qian19b}, variance-reduced {\tt SGD} methods such as {\tt SAGA}~\cite{saga}, {\tt LSVRG}~\cite{hofmann2015variance, kovalev2019don} and {\tt JacSketch}~\cite{jacsketch}. Another by-product  is {\em the unified analysis of {\tt SGD} methods which include {\tt RCD}.} For instance, our theory covers the subspace descent method {\tt SEGA}~\cite{sega} as a special case. Lastly, our framework is general enough to capture the phenomenon of {\em quantization}. For instance, we obtain the {\tt DIANA} and {\tt VR-DIANA} methods in special cases. 


\item {\bf Generalization of existing methods.} An important yet {\em relatively} minor contribution of our work is that it enables  {\em generalization} of knowns methods.  For instance, some particular methods we consider, such as {\tt LSVRG} (Algorithm~\ref{alg:sk_L-SVRG})~\cite{kovalev2019don}, were not analyzed in the proximal ($\psi\neq 0$) case before. To illustrate how this can be done within our framework, we do it here for {\tt LSVRG}. Further, most of the methods we analyze can be extended to the {\em arbitrary sampling} paradigm. 

\item {\bf Sharp rates.} In all known special cases, the rates obtained from our general theorem (Theorem~\ref{thm:sk_main_gsgm}) are the {\em best known rates} for these methods.

\item {\bf New methods.}   Our general analysis provides estimates for a possibly infinite array of new and yet-to-be-developed variants of {\tt SGD}. One only needs to verify that Assumption~\ref{as:sk_general_stoch_gradient} holds, and a complexity estimate is readily furnished by Theorem~\ref{thm:sk_main_gsgm}. Selected existing and new methods that fit our framework are summarized in Table~\ref{tbl:sk_special_cases2}. This list is for illustration only, we believe that future work by us and others will lead to its rapid expansion.

\item {\bf Experiments.} We show through extensive experimentation that some of the {\em new} and {\em generalized} methods proposed here and analyzed via our framework have some intriguing practical properties when compared against appropriately selected existing methods.

 \end{itemize}

\section{Main result \label{sec:sk_main_res}}

We first introduce the key assumption on the stochastic gradients $g^k$ enabling our general analysis (Assumption~\ref{as:sk_general_stoch_gradient}), then state our assumptions on $f$ (Assumption~\ref{as:sk_mu_strongly_quasi_convex}), and finally  state and comment on our  unified convergence result (Theorem~\ref{thm:sk_main_gsgm}). 

\paragraph{\bf Notation.} Consistently with the rest of the thesis, we use the following notation: $\langle x, y\rangle \eqdef \sum_i x_i y_i$ is the standard Euclidean inner product, and $\norm{x}\eqdef \langle x, x\rangle ^{1/2} $ is the induced $\ell_2$ norm. For simplicity we assume that \eqref{eq:sk_problem_gen} has a unique minimizer, which we denote $x^*$. Let $D_f(x,y)$ denote the \textit{Bregman divergence} associated with $f$: $D_f(x,y) \eqdef f(x) - f(y) - \<\nabla f(y), x-y>$. We  often write $[n]\eqdef \{1,2,\dots,n\}$.

\subsection{Key assumption}
Our first assumption is of key importance. It is mainly an assumption on the sequence of stochastic gradients $\{g^k\}$ generated by an arbitrary randomized algorithm. Besides unbiasedness (see \eqref{eq:sk_general_stoch_grad_unbias}), we require two recursions to hold for the iterates $x^k$ and the stochastic gradients $g^k$ of a randomized method. We allow for flexibility by casting these inequalities in a parametric manner.

\begin{assumption}\label{as:sk_general_stoch_gradient} Let $\{x^k\}$ be the random iterates produced by proximal {\tt SGD} (Algorithm in Eq~\eqref{eq:sk_SGD}).
We first assume that the stochastic gradients $g^k$ are unbiased
    \begin{equation}\label{eq:sk_general_stoch_grad_unbias}
        \EEE\left[ g^k\mid x^k\right] = \nabla f(x^k),
    \end{equation}
    for all $k\geq 0$. Further, we assume that there exist
non-negative constants $A, B, C, D_1, D_2, \rho$ and a (possibly) random sequence $\{\sigma_k^2\}_{k\ge 0}$ such that the following two relations hold\footnote{For convex and $L$-smooth $f$, one can show that $    \norm{\nabla f(x) - \nabla f(y)}^2 \le 2LD_{f}(x,y).$ Hence, $D_f$ can be used as a measure of proximity for the gradients.}
    \begin{equation}\label{eq:sk_general_stoch_grad_second_moment}
        \EEE\left[\norm{g^k -\nabla f(x^*)}^2\mid x^k, \sigma^2_k\right] \le 2AD_f(x^k,x^*) + B\sigma_k^2 + D_1,
    \end{equation}
    \begin{equation}\label{eq:sk_gsg_sigma}
        \EEE\left[\sigma_{k+1}^2 \, \mid \, x^k, \sigma^2_k\right] \le (1-\rho) \sigma_k^2 + 2CD_f(x^k,x^*)  + D_2,
    \end{equation} 
The expectation above is with respect to the randomness of the algorithm. 

\end{assumption}

The unbiasedness assumption \eqref{eq:sk_general_stoch_grad_unbias} is standard.  The key innovation we bring is inequality \eqref{eq:sk_general_stoch_grad_second_moment} coupled with \eqref{eq:sk_gsg_sigma}. We argue, and justify this statement by furnishing many examples in Section~\ref{sec:sk_main-paper-special-cases}, that these inequalities capture the essence of a wide array of existing and some new {\tt SGD} methods, including vanilla, variance reduced, arbitrary sampling, quantized and coordinate descent variants. Note that in the case when $\nabla f(x^*) = 0$ (e.g., when $\psi\equiv0$),  the inequalities in Assumption~\ref{as:sk_general_stoch_gradient}  reduce to
\begin{equation}\label{eq:sk_general_stoch_grad_second_moment_special}
        \EEE\left[\norm{g^k}^2\mid x^k, \sigma^2_k\right] \le 2A(f(x^k) - f(x^*)) + B\sigma_k^2 + D_1,
\end{equation}
\begin{equation}\label{eq:sk_gsg_sigma_special}
    \EEE\left[\sigma_{k+1}^2 \, \mid \,x^k, \sigma^2_k\right] \le (1-\rho) \sigma_k^2 + 2C(f(x^k) - f(x^*)) + D_2.
\end{equation}
Similar inequalities can be found in the analysis of stochastic first-order methods. However, this is the first time that such inequalities are generalized, equipped with parameters, and elevated to the status of an assumption that can be used on its own, independently from any other details defining the underlying method that generated them.

To give a further intuition about inequalities~\eqref{eq:sk_general_stoch_grad_second_moment} and~\eqref{eq:sk_gsg_sigma},    we shall note that sequence $\sigma_k$ usually represents the portion of noise that can gradually decrease over the course of optimization while constants $D_1, D_2$ represent a static noise. On the other hand, constants $A, C$ are usually related to some measure of smoothness of the objective. For instance, the parameters for (deterministic) gradient descent can be chosen as $A = L, B = C = D_1 =D_2 = \sigma_k^2 = \rho =0$. For an overview of parameter choices for specific instances of~\eqref{eq:sk_SGD}, see Table~\ref{tbl:sk_special_cases-parameters}. Note also that the choice of parameters of~\eqref{eq:sk_general_stoch_grad_second_moment} and~\eqref{eq:sk_gsg_sigma} is not unique, however this has no impact on convergence rates we provide.

\subsection{Main theorem}

 For simplicity, we shall assume throughout that $f$ is $\mu$-strongly quasi-convex, which is a generalization of $\mu$-strong convexity. We leave an analysis under different assumptions on $f$ to future work.

\begin{assumption}[$\mu$-strong quasi-convexity]\label{as:sk_mu_strongly_quasi_convex}
There exists $\mu>0$ such that $f:\R^d\to \R$ is \textit{$\mu$-strongly quasi-convex}.  That is, the following inequality holds for all $ x\in\R^d$:
    \begin{equation}\label{eq:sk_mu_strongly_quasi_convex}
    \compactify
        f(x^*) \ge f(x) + \langle \nabla f(x), x^* - x\rangle + \frac{\mu}{2}\norm{x^* - x}^2.    
    \end{equation}
\end{assumption}

We are now ready to present the key lemma of this chapter which states per iteration recurrence to analyze~\eqref{eq:sk_SGD}. Due to space limitations, we present the proof in Section~\ref{sec:sk_main_res} of the Appendix. 
 
\begin{lemma}  \label{lem:sk_iter_dec}
Let Assumptions~\ref{as:sk_general_stoch_gradient}~and~\ref{as:sk_mu_strongly_quasi_convex}  be satisfied. Then the following inequality holds for all $k\geq 0$:
\begin{align*}
& \EEE\left[\norm{x^{k+1}-x^*}^2\right] + M\alpha^2\EEE\left[\sigma_{k+1}^2\right]  + 2\alpha\left(1-\alpha(A+CM)\right)\EEE\left[D_f(x^k,x^*)\right] \\ 
    &   \qquad   \le
        (1-\alpha\mu)\EEE \left[\norm{x^k - x^*}^2 \right] + \left(1 - \rho\right)M\alpha^2\EEE\left[\sigma_k^2\right] + B \alpha^2\EEE\left[\sigma_k^2\right]   + (D_1+MD_2)\alpha^2.
\end{align*}
\end{lemma}

Using recursively Lemma~\ref{lem:sk_iter_dec}, we obtain the convergence rate of proximal SGD, which we state as Theorem~\ref{thm:sk_main_gsgm}.

 \begin{theorem}\label{thm:sk_main_gsgm}
Let Assumptions~\ref{as:sk_general_stoch_gradient}~and~\ref{as:sk_mu_strongly_quasi_convex} be satisfied. Choose constant $M$ such that $M > \frac{B}{\rho}$. Choose a stepsize  satisfying
    \begin{equation}\label{eq:sk_gamma_condition_gsgm}
    \compactify
        0 < \alpha \le \min\left\{\frac{1}{\mu}, \frac{1}{A+CM}\right\}.
    \end{equation}
    Then the iterates $\{x^k\}_{k\geq 0}$ of proximal {\tt SGD} (Algorithm~\eqref{eq:sk_SGD}) satisfy
    \begin{align}
     \EEE\left[V^k\right] \le & \max\left\{(1-\alpha\mu)^k,\left(1+\frac{B}{M}-\rho\right)^k\right\} V^0   + \frac{(D_1+MD_2)\alpha^2 }{\min\left\{\alpha\mu, \rho - \frac{B}{M}\right\}},\label{eq:sk_main_result_gsgm}
    \end{align}
    where the Lyapunov function $V^k$ is defined by $V^k \eqdef \norm{x^k - x^*}^2 + M\alpha^2\sigma_k^2$.
\end{theorem}

This theorem establishes a linear rate for a wide range of proximal {\tt SGD} methods up to a certain oscillation radius, controlled by the additive term in \eqref{eq:sk_main_result_gsgm}, and namely, by parameters $D_1$ and $D_2$. As we shall see  in Section~\ref{sec:sk_special_cases} (refer to Table~\ref{tbl:sk_special_cases-parameters}), the main difference between the vanilla and variance-reduced {\tt SGD} methods is that while the former  satisfy inequality \eqref{eq:sk_gsg_sigma} with $D_1>0$ or $D_2>0$, which in view of \eqref{eq:sk_main_result_gsgm} prevents them from reaching the optimum $x^*$ (using a fixed stepsize), the latter methods satisfy inequality \eqref{eq:sk_gsg_sigma} with $D_1=D_2=0$, which in view of \eqref{eq:sk_main_result_gsgm} enables them to reach the optimum.

\section{The classic, the recent and the brand new} \label{sec:sk_main-paper-special-cases}

In this section we deliver on the promise from the introduction and show how many existing and some new variants of {\tt SGD} fit our general framework (see Table~\ref{tbl:sk_special_cases2}). 

\paragraph{\bf An overview.} As claimed, our framework is powerful enough to include vanilla methods (\xmark\; in the ``VR'' column) as well as variance-reduced methods (\cmark\; in the ``VR'' column), methods which generalize to arbitrary sampling (\cmark\; in the ``AS'' column), methods supporting gradient quantization (\cmark\; in the ``Quant'' column) and finally, also {\tt RCD} type methods (\cmark\; in the ``RCD'' column). 

\begin{table*}[!t]
\begin{center}
\footnotesize
\begin{tabular}{|c|c|c|c|c|c|c|c|c|c|}
\hline
Problem & Method &  Alg & Citation &   VR & AS & Quant &  RCD & Sec  & Cor \\
\hline
\eqref{eq:sk_problem_gen}+\eqref{eq:sk_f_exp} & {\tt SGD}  &  \ref{alg:sk_sgd_prox} & \cite{nguyen2018sgd}  & \xmark &  \xmark & \xmark &  \xmark & \ref{sec:sk_SGD} & 
\ref{cor:sk_recover_sgd_rate} \\
\eqref{eq:sk_problem_gen}+\eqref{eq:sk_f_sum}  & {\tt SGD-SR} &  \ref{alg:sk_sgdas} & \cite{pmlr-v97-qian19b} & \xmark &  \cmark & \xmark & \xmark & \ref{SGD-AS} &   \ref{cor:sk_recover_sgd-as_rate} \\
\eqref{eq:sk_problem_gen}+\eqref{eq:sk_f_sum} &  {\tt SGD-MB} & \ref{alg:sk_SGD-MB} & {\bf NEW} & \xmark & \xmark & \xmark & \xmark & \ref{sec:sk_SGD-MB} & \ref{cor:sk_mb}   \\
\eqref{eq:sk_problem_gen}+\eqref{eq:sk_f_sum} &  {\tt SGD-star} &  \ref{alg:sk_SGD-star} & {\bf NEW} & \cmark & \cmark  & \xmark & \xmark & \ref{sec:sk_SGD-star} & \ref{cor:sk_SGD-star} \\
\eqref{eq:sk_problem_gen}+\eqref{eq:sk_f_sum}  & {\tt SAGA} & \ref{alg:sk_SAGA} & \cite{saga} & \cmark & \xmark  & \xmark & \xmark & \ref{sec:sk_saga} & \ref{thm:sk_recover_saga_rate} \\
\eqref{eq:sk_problem_gen}+\eqref{eq:sk_f_sum}  & {\tt N-SAGA} & \ref{alg:sk_N-SAGA} & {\bf NEW} & \xmark &  \xmark & \xmark & \xmark & \ref{N-SAGA} & \ref{cor:sk_N-SAGA} \\
\eqref{eq:sk_problem_gen} & {\tt SEGA}  & \ref{alg:sk_SEGA} &  \cite{sega}  & \cmark &   \xmark & \xmark & \cmark & \ref{sec:sk_sega} & \ref{cor:sk_sega} \\
\eqref{eq:sk_problem_gen} & {\tt N-SEGA}  & \ref{alg:sk_N-SEGA} &  {\bf NEW}  & \xmark & \xmark  & \xmark  & \cmark & \ref{N-SEGA} &
\ref{cor:sk_N-SEGA}\\
\eqref{eq:sk_problem_gen}+\eqref{eq:sk_f_sum}  & {\tt SVRG}${}^{a}$  & \ref{alg:sk_SVRG} & \cite{svrg} & \cmark & \xmark & \xmark & \xmark & \ref{sec:sk_svrg} & \ref{cor:sk_svrg} \\
\eqref{eq:sk_problem_gen}+\eqref{eq:sk_f_sum}  & {\tt LSVRG} & \ref{alg:sk_L-SVRG} & \cite{hofmann2015variance} & \cmark &  \xmark & \xmark & \xmark & \ref{sec:sk_L-SVRG} & \ref{cor:sk_recover_l-svrg_rate}  \\
\eqref{eq:sk_problem_gen}+\eqref{eq:sk_f_sum}  & {\tt DIANA} & \ref{alg:sk_diana} & 
\cite{mishchenko2019distributed} & \xmark & \xmark & \cmark &  \xmark & \ref{sec:sk_diana} & \ref{cor:sk_main_diana}  \\
\eqref{eq:sk_problem_gen}+\eqref{eq:sk_f_sum}  & {\tt DIANA}${}^b$  & \ref{alg:sk_diana_case}  & \cite{mishchenko2019distributed} & \cmark &  \xmark & \cmark  &  \xmark &\ref{sec:sk_diana}  & \ref{cor:sk_main_diana_special_case}\\
\eqref{eq:sk_problem_gen}+\eqref{eq:sk_f_sum}  & {\tt Q-SGD-SR} &  \ref{alg:sk_qsgdas} & {\bf NEW} & \xmark &  \cmark & \cmark & \xmark & \ref{Q-SGD-AS} &   \ref{cor:sk_recover_q-sgd-as_rate} \\ 
\eqref{eq:sk_problem_gen}+\eqref{eq:sk_f_sum}+\eqref{eq:sk_f_i_sum}  & {\tt VR-DIANA} & \ref{alg:sk_vr-diana} & \cite{horvath2019stochastic}& \cmark & \xmark & \cmark & \xmark & \ref{sec:sk_VR-DIANA} & \ref{cor:sk_main_vr_diana} \\
\eqref{eq:sk_problem_gen}+\eqref{eq:sk_f_sum}& {\tt JacSketch} & \ref{alg:sk_jacsketch} & \cite{jacsketch} & \cmark &  \cmark \xmark & \xmark & \xmark & \ref{sec:sk_JacSketch} & \ref{thm:sk_main_jacsketch}\\
\hline
\end{tabular}
\end{center}
\caption{List of specific existing (in some cases generalized) and new methods which fit our general analysis framework. VR = variance reduced method, AS = arbitrary sampling, Quant = supports gradient quantization, RCD = randomized coordinate descent type method. ${}^{a}$ Special case of {\tt SVRG} with  1 outer loop only;  ${}^{b}$ Special case of {\tt DIANA} with $1$ node and quantization of exact gradient. }
\label{tbl:sk_special_cases2}
\end{table*}

For existing methods we provide a citation; new methods developed in this chapter are marked accordingly.  Due to space restrictions, all algorithms are described (in detail) in the Appendix; we provide a link to the appropriate section for easy navigation. While these details are important, the main message of this chapter, i.e., the generality of our approach, is captured by Table~\ref{tbl:sk_special_cases2}. The ``Result'' column of Table~\ref{tbl:sk_special_cases2} points to a corollary of Theorem~\ref{thm:sk_main_gsgm}; these corollaries state in detail the convergence  statements for the various methods. In all cases where known methods are recovered, these corollaries of Theorem~\ref{thm:sk_main_gsgm} recover the best known rates.

\paragraph{\bf Parameters.} From the point of view of Assumption~\ref{as:sk_general_stoch_gradient}, the methods listed in Table~\ref{tbl:sk_special_cases2} exhibit certain patterns.  To shed some light on this, in Table~\ref{tbl:sk_special_cases-parameters} we summarize the values of these parameters. 

\begin{table*}[!t]
\caption{The parameters for which the methods from Table~\ref{tbl:sk_special_cases2} (special cases of~\eqref{eq:sk_SGD}) satisfy Assumption~\ref{as:sk_general_stoch_gradient}. The meaning of the expressions appearing in the table, as well as their justification is defined in detail in the Appendix (Section~\ref{sec:sk_special_cases}). }
\label{tbl:sk_special_cases-parameters}
\begin{center}
\begin{tabular}{|c|c|c|c|c|c|c|}
\hline
 Method &   $A$ & $B$ & $\rho$ & $C$ & $D_1$ & $D_2$\\
\hline
 {\tt SGD}   &  $2L $ & $0$ & $1$ & $0$ & $2\sigma^2$ & $0$ \\
 {\tt SGD-SR} &  $2 \cL$ & $0$ & $1$ & $0$ & $2\sigma^2$ & $0$ \\
  {\tt SGD-MB}  &  $\frac{A' + L(\tau-1)}{\tau}$ & 0 & $1$ & $0$ & $ \frac{D'}{\tau}$ & $0$ \\
  {\tt SGD-star}  &  $2 \cL$ & $0$ & $1$ & $0$ & $0$ & $0$ \\
 {\tt SAGA}   &  $2L$ & $2$ & $1/n$ & $L/n$ & $0$ & $0$\\
{\tt N-SAGA}   &  $2L$ & $2$ & $1/n$ & $L/n$ & $2\sigma^2$ & $\frac{\sigma^2}{n}$\\
{\tt SEGA}   &   $2dL$ & $2d$ & $1/d$ & $L/d$ & $0$ & $0$\\
{\tt N-SEGA}  &   $2dL$ & $2d$ & $1/d$ & $L/d$ & $2d\sigma^2$ & $\frac{\sigma^2}{d}$\\
 {\tt SVRG}${}^{a}$   & $2L$  & $2$ & $0$ & $0$& $0$ & $0$ \\
 {\tt LSVRG}   &  $2 L$ & $2$ & $p$ & $Lp$ & $0$ & $0$ \\
 {\tt DIANA}   & $\left(1+\frac{2\omega}{n}\right)L$ & $\frac{2\omega}{n}$ & $\gamma$ & $L\gamma$ & $\frac{(1+\omega)\sigma^2}{n}$ & $\gamma\sigma^2$  \\
 {\tt DIANA}${}^b$  &  $\left(1+2\omega\right)L$ & $2\omega$ & $\gamma$ & $L\gamma$ & $0$ & $0$ \\
{\tt Q-SGD-SR} &  $2(1+\omega)\cL$ & $0$ & $1$ & $0$ & $2(1+\omega)\sigma^2$ & $0$ \\
 {\tt VR-DIANA} &  $\left(1+\frac{4\omega + 2}{n}\right)L$ & $\frac{2(\omega+1)}{n}$ & $\gamma$ & $\left(\frac{1}{m}+4\gamma\right)L$ & $0$ & $0$ \\
 {\tt JacSketch}   &  $2\cL_1$ & $\frac{2\lambda_{\max}}{n}$ & $\lambda_{\min}$ & $\frac{\cL_2}{n}$ & $0$ & $0$ \\
\hline
\end{tabular}
\end{center}
\end{table*}

Note, for example, that for all methods the parameter $A$ is non-zero. Typically, this a multiple of an appropriately defined smoothness parameter (e.g., $L$ is the Lipschitz constant of the gradient of $f$, $\cL$ and $\cL_1$ in {\tt SGD-SR}\footnote{{\tt SGD-SR} is first {\tt SGD} method analyzed in the {\em arbitrary sampling} paradigm. It was developed using the  {\em stochastic reformulation} approach (whence the ``SR'') pioneered in \cite{richtarik2017stochastic} in a numerical linear algebra setting, and later extended to develop the {\tt JacSketch} variance-reduction technique for finite-sum optimization  \cite{jacsketch}.}, {\tt SGD-star} and {\tt JacSketch} are {\em expected smoothness} parameters). In the three variants of the {\tt DIANA} method, $\omega$ captures the variance of the  quantization operator $Q$. That is, one assumes that $\E{Q(x)}=x$ and $\E{\norm{Q(x)-x}^2} \leq \omega \norm{x}^2$ for all $x\in \R^d$. In view of \eqref{eq:sk_gamma_condition_gsgm}, large $A$ means a smaller stepsize, which slows down the rate. Likewise, the variance $\omega$ also affects the parameter $B$, which in view of \eqref{eq:sk_main_result_gsgm} also has an adverse effect on the rate.
Further, as predicted by Theorem~\ref{thm:sk_main_gsgm}, whenever either $D_1>0$ or $D_2>0$, the corresponding method converges to an oscillation region only. These methods are not variance-reduced. All symbols used in Table~\ref{tbl:sk_special_cases-parameters} are defined in the appendix, in the same place where the methods are described and analyzed.

\paragraph{\bf Five new methods.} To illustrate the usefulness of our general framework, we develop  {\em 5 new variants} of {\tt SGD} never explicitly considered in the literature before (see Table~\ref{tbl:sk_special_cases2}). Here we briefly motivate them; details can be found in the Appendix. 

\begin{itemize}
\item {\tt SGD-MB} (Algorithm~\ref{alg:sk_SGD-MB}).
This method is specifically designed for functions of the finite-sum structure  \eqref{eq:sk_f_i_sum}. As we show through experiments, this is a powerful mini-batch {\tt SGD} method, with mini-batches formed with replacement as follows: in each iteration, we repeatedly ($\tau$ times) and independently pick $i\in [n]$ with probability $p_i>0$. Stochastic gradient $g^k$ is then formed by averaging the stochastic gradients $\nabla f_i(x^k)$ for all selected indices $i$ (including each $i$ as many times as this index was selected). This allows for a more practical importance mini-batch sampling implementation than what was until now possible (see Remark~\ref{rem:sgdmb} in the Appendix for more details and experiment in Figure~\ref{fig:sk_SGDMB_full}).

\item {\tt SGD-star} (Algorithm~\ref{alg:sk_SGD-star}). 
This new method forms a bridge between vanilla and variance-reduced {\tt SGD} methods. While  not practical, it sheds light on the role of variance reduction. Again, we consider functions of the finite-sum form~\eqref{eq:sk_f_i_sum}. This methods answers the following question:  assuming that the gradients $\nabla f_i(x^*)$, $i\in [n]$ are {\em known}, can they be used to design a more powerful {\tt SGD} variant? The answer is {\em yes}, and {\tt SGD-star} is the method. In its most basic form, {\tt SGD-star}  constructs the stochastic gradient via $g^k=\nabla f_i(x^k) - \nabla f_i(x^*)+\nabla f(x^*)$, where $i\in [n]$ is chosen uniformly at random. 
Inferring from Table~\ref{tbl:sk_special_cases-parameters}, where $D_1=D_2=0$, this method converges to $x^*$, and not merely to some oscillation region.  Variance-reduced methods essentially work by iteratively constructing increasingly more accurate {\em estimates} of $\nabla f_i(x^*)$. Typically, the term $\sigma_k^2$ in the Lyapunov function of variance reduced methods will contain a term of the form $\sum_i \norm{h_i^k - \nabla f_i(x^*)}^2$, with $h_i^k$ being the estimators maintained by the method. Remarkably, {\tt SGD-star} was never explicitly considered in the literature before.

\item {\tt N-SAGA} (Algorithm~\ref{alg:sk_N-SAGA}). This is a novel variant of {\tt SAGA}~\cite{saga}, one in which one does not have access to the gradients of $f_i$, but instead only has access to {\em noisy} stochastic estimators thereof (with noise $\sigma^2$). Like {\tt SAGA}, {\tt N-SAGA} is able to reduce the variance inherent in the finite sum structure~\eqref{eq:sk_f_i_sum} of the problem. However, it necessarily pays the price of noisy estimates of $\nabla f_i $, and hence, just like vanilla {\tt SGD} methods, is ultimately unable to converge to $x^*$. The oscillation region is governed by the noise level $\sigma^2$ (refer to $D_1$ and $D_2$ in Table~\ref{tbl:sk_special_cases-parameters}). This method will be of practical importance for problems where each $f_i$ is of the form~\eqref{eq:sk_f_exp}, i.e., for problems of the ``average of expectations'' structure. Batch versions of {\tt N-SAGA} would be well suited for distributed optimization, where each $f_i$ is owned by a different worker, as in such a case one wants the workers to work in parallel.

\item {\tt N-SEGA} (Algorithm~\ref{alg:sk_N-SEGA}). This is a {\em noisy} extension of the {\tt RCD}-type method {\tt SEGA}, in complete analogy with the relationship between {\tt SAGA} and {\tt N-SAGA}. Here we assume that we only have noisy estimates of partial derivatives (with noise $\sigma^2$). This situation is common in derivative-free optimization, where such a noisy estimate can be obtained by taking (a random) finite difference approximation~\cite{nesterov2017randomDFO}. Unlike {\tt SEGA}, {\tt N-SEGA} only converges to an oscillation region the size of which is governed by $\sigma^2$.

\item {\tt Q-SGD-SR} (Algorithm~\ref{alg:sk_qsgdas}). This is a quantized version of {\tt SGD-SR}, which is the first {\tt SGD} method analyzed in the  arbitrary sampling paradigm. As such, {\tt Q-SGD-SR}  is a vast generalization of the celebrated {\tt QSGD} method \cite{alistarh2017qsgd}.

\end{itemize}

\section{Experiments \label{sec:sk_exp}}

In this section we numerically verify the claims from the chapter. We perform three differnent experiments: we verify the usefulness of {\tt SGD-MB} alongside with testing both {\tt SGD-star} and {\tt N-SEGA}.

\subsection{{\tt SGD-MB}: remaining experiments and exact problem setup.}

In Section~\ref{sec:sk_SGD-MB}, we describe in detail the {\tt SGD-MB} method already outlined before. The main advantage of {\tt SGD-MB} is that the sampling procedure it employs can be implemented in just $\cO(\tau \log n)$ time. In contrast, even the simplest without-replacement sampling which selects each function into the minibatch with a prescribed probability independently (we will refer to it as independent {\tt SGD}) requires $n$ calls of a uniform random generator.  We demonstrate numerically that {\tt SGD-MB} has essentially identical iteration complexity to  independent {\tt SGD} in practice. We consider logistic regression with Tikhonov regularization of order $\lambda$:
\begin{equation}\label{eq:sk_logreg}
\frac1n \sum_{i=1}^n   \log \left(1+\exp\left(a_i^\top x\cdot  b_i\right) \right)+\frac{\lambda}{2} \norm{ x}^2,
\end{equation}
where $a_i\in \R^{n}$, $b_i \in \{-1,1\}$ is $i$th data-label pair is a vector of labels and $\lambda\geq 0$ is the regularization parameter. The data and labels were obtained from LibSVM datasets {\tt a1a}, {\tt a9a}, w1a, {\tt w8a}, {\tt gisette}, {\tt madelon}, {\tt phishing} and {\tt mushrooms}. Further, the data were rescaled by a random variable $cu_i^2$ where $u_i$ is random integer from $1,2,\dots, 1000$ and $c$ is such that the mean norm of $a_i$ is $1$.  Note that we have now an infinite array of possibilities on how to write~\eqref{eq:sk_logreg} as~\eqref{eq:sk_f_sum}. For simplicity, distribute $l2$ term evenly among the finite sum.

For a fixed  expected sampling size $\tau$, we consider two options for the probability of sampling the $i$th function:
\begin{enumerate}[label=(\roman*)]
\item \label{item:unif} $\frac{\tau}{n}$, or
\item \label{item:imp} $\frac{\norm{a_i}^2+\lambda}{\delta+\norm{a_i}^2+\lambda}$, where $\delta$ is such that\footnote{An {\tt RCD} version of this sampling was proposed in~\cite{hanzely2018accelerated}; it was shown to be superior to uniform sampling both in theory and practice.} $\sum_{i=1}^n\frac{\norm{a_i}^2+\lambda}{\delta+\norm{a_i}^2+\lambda}=1$. 
\end{enumerate}
The results can be found in Figure~\ref{fig:sk_SGDMB_full}, where we also report the choice of stepsize $\alpha$ and the choice of $\tau$ in the legend and title of the plot, respectively.

\begin{figure}[!h]
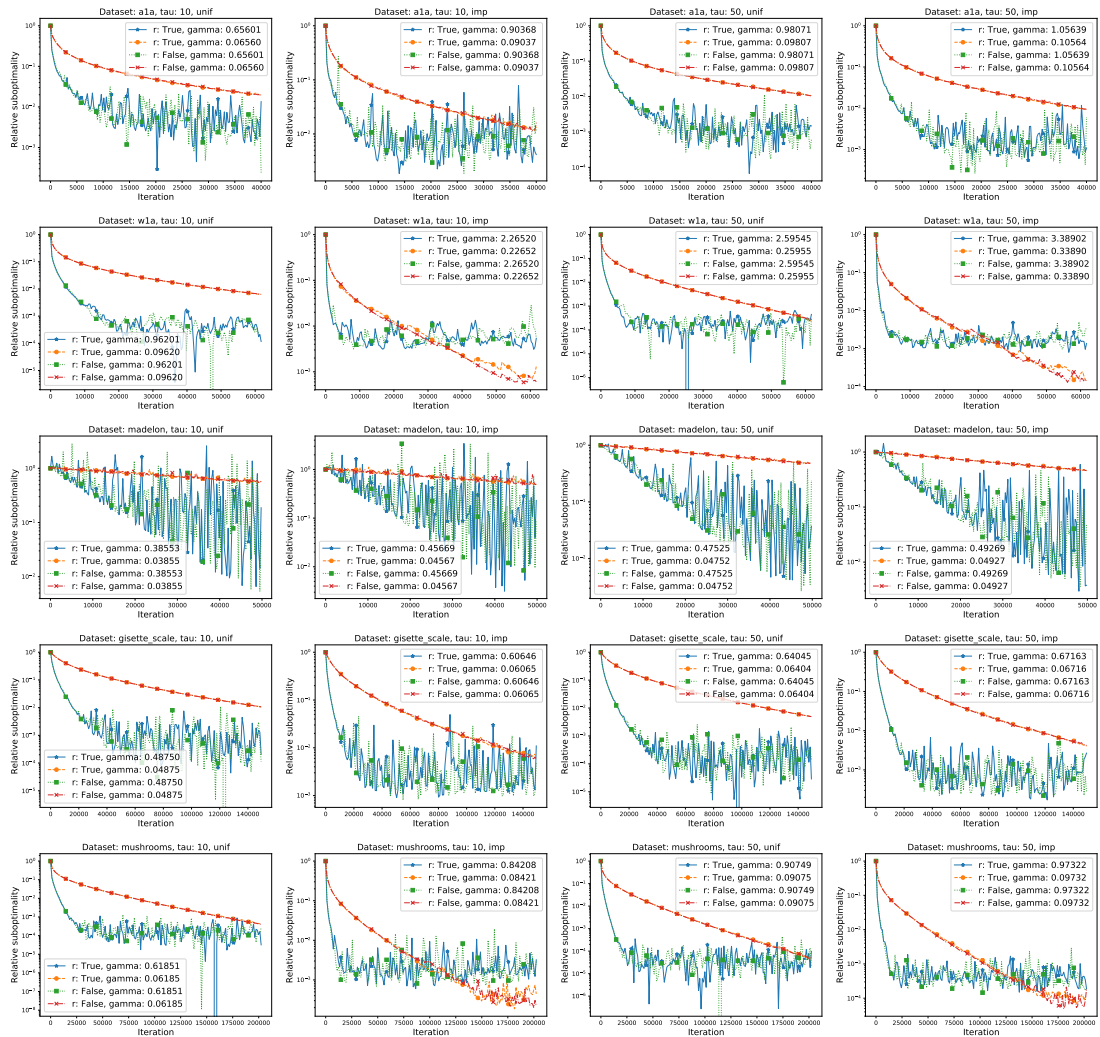

\centering
\begin{minipage}{0.24\textwidth}
  \centering
\includegraphics[width =  \textwidth ]{SGD_a1a_tau_10_meth_3_unif.pdf}
\end{minipage}%
\begin{minipage}{0.24\textwidth}
  \centering
\includegraphics[width =  \textwidth ]{SGD_a1a_tau_10_meth_3_imp.pdf}
\end{minipage}%
\begin{minipage}{0.24\textwidth}
  \centering
\includegraphics[width =  \textwidth ]{SGD_a1a_tau_50_meth_3_unif.pdf}
\end{minipage}%
\begin{minipage}{0.24\textwidth}
  \centering
\includegraphics[width =  \textwidth ]{SGD_a1a_tau_50_meth_3_imp.pdf}
\end{minipage}%
\\
\begin{minipage}{0.24\textwidth}
  \centering
\includegraphics[width =  \textwidth ]{SGD_w1a_tau_10_meth_3_unif.pdf}
\end{minipage}%
\begin{minipage}{0.24\textwidth}
  \centering
\includegraphics[width =  \textwidth ]{SGD_w1a_tau_10_meth_3_imp.pdf}
\end{minipage}%
\begin{minipage}{0.24\textwidth}
  \centering
\includegraphics[width =  \textwidth ]{SGD_w1a_tau_50_meth_3_unif.pdf}
\end{minipage}%
\begin{minipage}{0.24\textwidth}
  \centering
\includegraphics[width =  \textwidth ]{SGD_w1a_tau_50_meth_3_imp.pdf}
\end{minipage}%
\\
\begin{minipage}{0.24\textwidth}
  \centering
\includegraphics[width =  \textwidth ]{SGD_madelon_tau_10_meth_3_unif.pdf}
\end{minipage}%
\begin{minipage}{0.24\textwidth}
  \centering
\includegraphics[width =  \textwidth ]{SGD_madelon_tau_10_meth_3_imp.pdf}
\end{minipage}%
\begin{minipage}{0.24\textwidth}
  \centering
\includegraphics[width =  \textwidth ]{SGD_madelon_tau_50_meth_3_unif.pdf}
\end{minipage}%
\begin{minipage}{0.24\textwidth}
  \centering
\includegraphics[width =  \textwidth ]{SGD_madelon_tau_50_meth_3_imp.pdf}
\end{minipage}%
\\
\begin{minipage}{0.24\textwidth}
  \centering
\includegraphics[width =  \textwidth ]{SGD_gisette_scale_tau_10_meth_3_unif.pdf}
\end{minipage}%
\begin{minipage}{0.24\textwidth}
  \centering
\includegraphics[width =  \textwidth ]{SGD_gisette_scale_tau_10_meth_3_imp.pdf}
\end{minipage}%
\begin{minipage}{0.24\textwidth}
  \centering
\includegraphics[width =  \textwidth ]{SGD_gisette_scale_tau_50_meth_3_unif.pdf}
\end{minipage}%
\begin{minipage}{0.24\textwidth}
  \centering
\includegraphics[width =  \textwidth ]{SGD_gisette_scale_tau_50_meth_3_imp.pdf}
\end{minipage}%
\\
\begin{minipage}{0.24\textwidth}
  \centering
\includegraphics[width =  \textwidth ]{SGD_mushrooms_tau_10_meth_3_unif.pdf}
\end{minipage}%
\begin{minipage}{0.24\textwidth}
  \centering
\includegraphics[width =  \textwidth ]{SGD_mushrooms_tau_10_meth_3_imp.pdf}
\end{minipage}%
\begin{minipage}{0.24\textwidth}
  \centering
\includegraphics[width =  \textwidth ]{SGD_mushrooms_tau_50_meth_3_unif.pdf}
\end{minipage}%
\begin{minipage}{0.24\textwidth}
  \centering
\includegraphics[width =  \textwidth ]{SGD_mushrooms_tau_50_meth_3_imp.pdf}
\end{minipage}%
\caption{{\tt SGD-MB} and independent {\tt SGD} applied on LIBSVM~\cite{chang2011libsvm} datasets with regularization parameter $\lambda = 10^{-5}$. Axis $y$ stands for relative suboptimality, i.e. $\frac{f(x^k)-f(x^*)}{f(x^k)-f(x^0)}$. Title label ``unif'' corresponds to probabilities chosen by~\ref{item:unif} while label ``imp'' corresponds to probabilities chosen by~\ref{item:imp}. Lastly, legend label ``r'' corresponds to ``replacement'' with value ``True'' for {\tt SGD-MB} and value ``False'' for independent {\tt SGD}.}
\label{fig:sk_SGDMB_full}
\end{figure}

Indeed, iteration complexity of {\tt SGD-MB} and independent {\tt SGD} is almost identical. Since the cost of each iteration of {\tt SGD-MB} is cheaper\footnote{The relative difference between iteration costs of {\tt SGD-MB} and independent {\tt SGD} can be arbitrary, especially for the case when cost of evaluating $\nabla f_i(x)$ is cheap, $n$ is huge and $n\gg \tau$. In such case, cost of one iteration of {\tt SGD-MB} is $\tau \text{Cost}(\nabla f_i) +\tau\log(n)$ while the cost of one iteration of independent {\tt SGD} is $\tau \text{Cost}(\nabla f_i) + n$.}, we conclude superiority of {\tt SGD-MB} to independent {\tt SGD}.

\subsection{Experiments on {\tt SGD-star} \label{sec:sk_exp_star}}
In this section, we study {\tt SGD-star} and numerically verify claims from Section~\ref{sec:sk_SGD-star}. In particular, Corollary~\ref{cor:sk_SGD-star} shows that {\tt SGD-star} enjoys linear convergence rate which is constant times better to the rate of {\tt SAGA} (given that problem condition number is high enough). We compare 3 methods -- {\tt SGD-star}, {\tt SGD} and {\tt SAGA}. We consider simple and well-understood least squares problem $\min_x \frac12 \| \mA x-b\|^2$ where elements of $\mA,b$ were generated (independently) from standard normal distribution. Further, rows of $\mA$ were normalized so that $\|\mA_{i:}\|=1$. Thus, denoting $f_i(x) = \frac12 (\mA_{i:}^\top x-b_i )^2$, $f_i$ is 1-smooth. For simplicity, we consider {\tt SGD-star} with uniform serial sampling, i.e. $\cL=1$.

 Next, for both {\tt SGD-star} and {\tt SGD} we use stepsize $\alpha = \frac{1}{2}$ (theory supported stepsize for {\tt SGD-star}), while for {\tt SAGA} we set $\alpha = \frac{1}{5}$ (almost theory supported stepsize). Figure~\ref{fig:sk_star} shows the results.

\begin{figure}[!h]
\centering
\begin{minipage}{0.3\textwidth}
  \centering
\includegraphics[width =  \textwidth ]{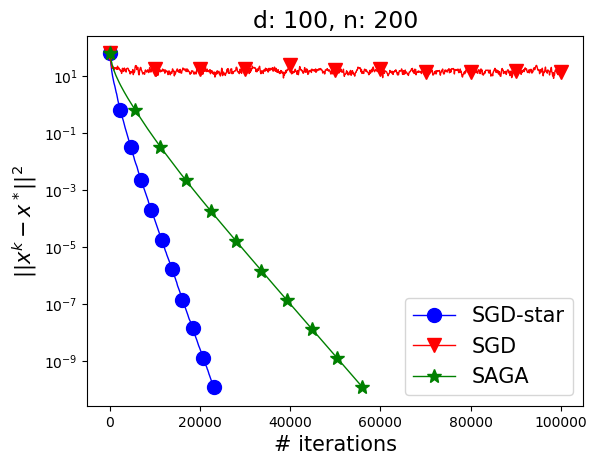}
\end{minipage}%
\begin{minipage}{0.3\textwidth}
  \centering
\includegraphics[width =  \textwidth ]{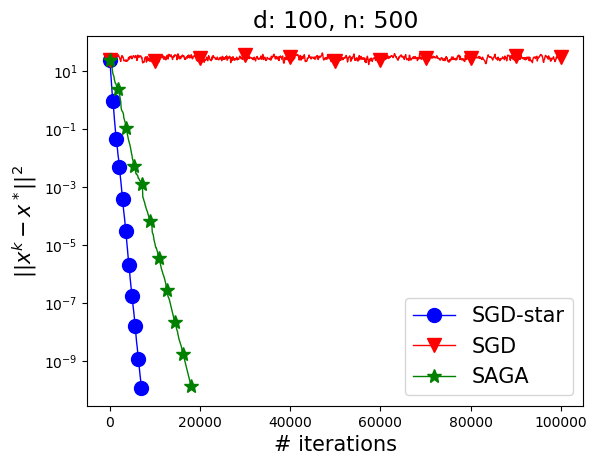}
\end{minipage}%
\begin{minipage}{0.3\textwidth}
  \centering
\includegraphics[width =  \textwidth ]{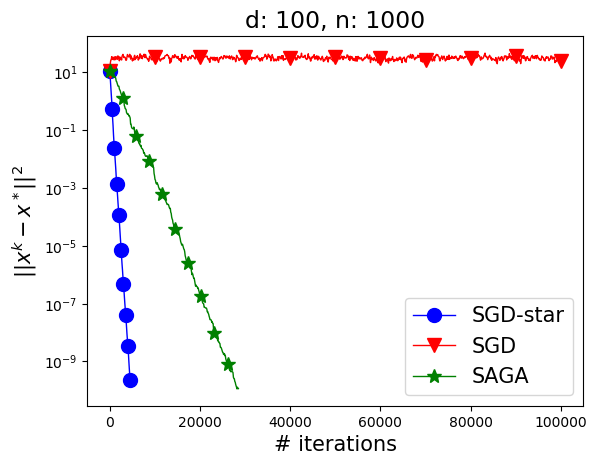}
\end{minipage}%
\\
\begin{minipage}{0.3\textwidth}
  \centering
\includegraphics[width =  \textwidth ]{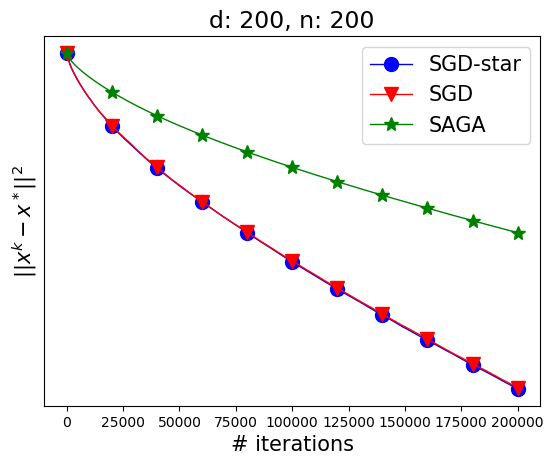}
\end{minipage}%
\begin{minipage}{0.3\textwidth}
  \centering
\includegraphics[width =  \textwidth ]{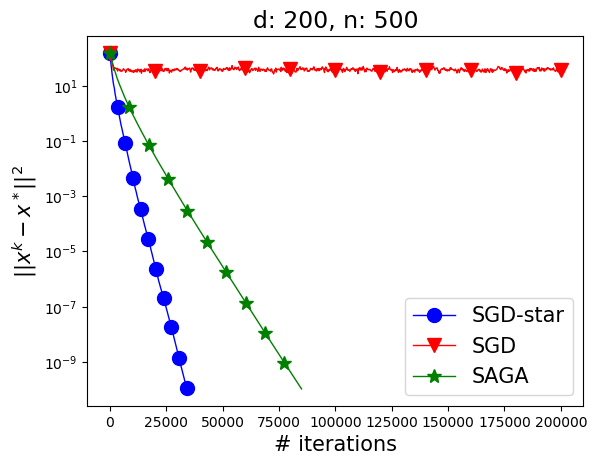}
\end{minipage}%
\begin{minipage}{0.3\textwidth}
  \centering
\includegraphics[width =  \textwidth ]{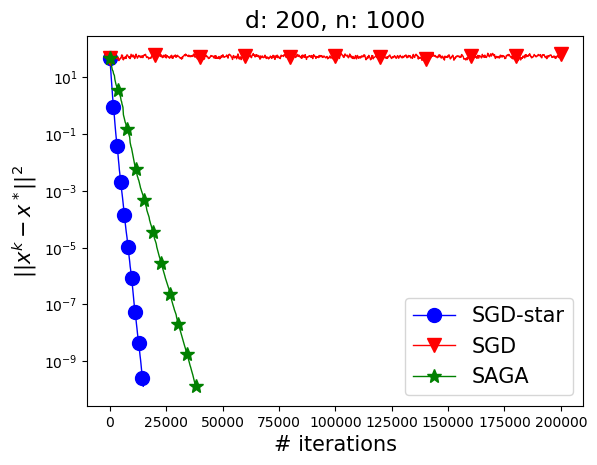}
\end{minipage}%
\\
\begin{minipage}{0.3\textwidth}
  \centering
\includegraphics[width =  \textwidth ]{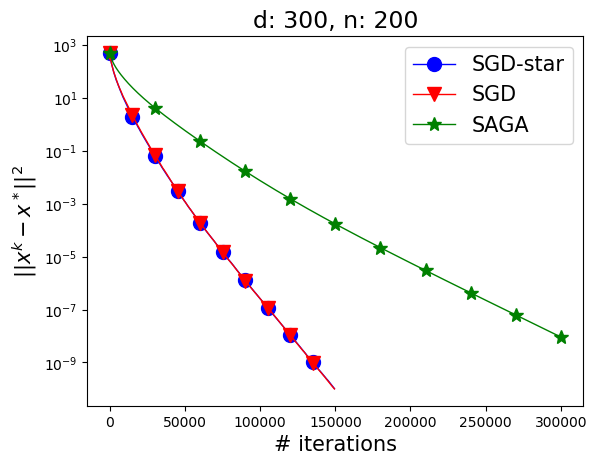}
\end{minipage}%
\begin{minipage}{0.3\textwidth}
  \centering
\includegraphics[width =  \textwidth ]{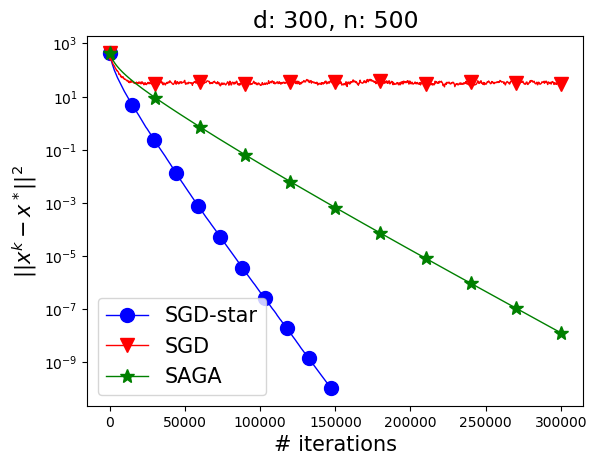}
\end{minipage}%
\begin{minipage}{0.3\textwidth}
  \centering
\includegraphics[width =  \textwidth ]{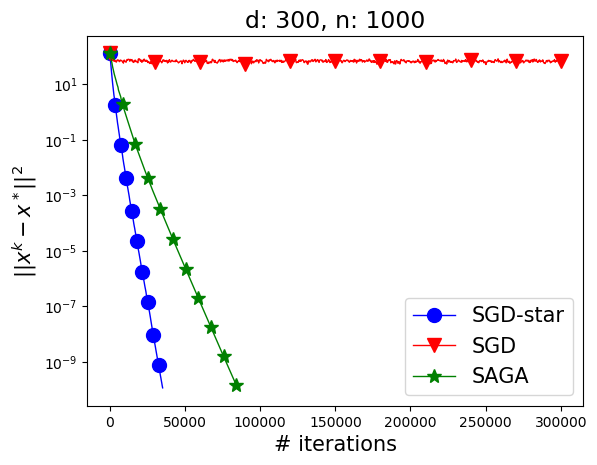}
\end{minipage}%
\caption{Comparison of {\tt SGD-star}, {\tt SGD} and {\tt SAGA} on least squares problem.}
\label{fig:sk_star}
\end{figure}

Note that, as theory predicts, {\tt SGD-star} is always faster to {\tt SAGA}, although only constant times. Further, in the cases where $d\geq n$, performance of {\tt SGD} seems identical to the performance of {\tt SGD-shift}. This is due to a simple reason: if $d\geq n$, we must have $\nabla f_i(x^*) = 0$ for all $i$, and thus {\tt SGD} and {\tt SGD-shift} are in fact identical algorithms. 

\subsection{Experiments on {\tt N-SEGA} \label{sec:sk_exp_nsega}}
In this experiment we study the effect of noise on {\tt N-SEGA}. We consider unit ball constrained least squares problem: $\min_{\|x\|\leq 1} f(x)$ where $f(x)=\|\mA x-b\|^2$. and we suppose that there is an oracle providing us with 
noised partial derivative $g_i(x,\zeta) = \nabla_i f(x) +\zeta$, where $\zeta \sim N(0,\sigma^2)$. For each problem instance (i.e. pair $\mA, b$), we compare performance of  {\tt N-SEGA} under various noise magnitudes $\sigma^2$.

The specific problem instances are presented in Table~\ref{tbl:sk_leastsquares}. Figure~\ref{fig:sk_nsega} shows the results. 

\begin{table}[!h]
\begin{center}
\begin{tabular}{|c|c|c|}
\hline
Type & $\mA $ & $b$ \\
 \hline
 \hline
 1   & $\mA_{ij}\sim N(0,1)$ (independently)  & vector of ones  \\
 \hline
  2   & Same as 1, but scaled so that $\lambda_{\max}(A^\top A)=1$ & vector of ones \\
\hline
 3   & $\mA_{ij} = \varrho_{ij}\varpi_{j}$ $\forall i,j:\varrho_{ij},\varpi_{j} \sim N(0,1)$ (independently)  & vector of ones  \\
 \hline
  4   & Same as 3, but scaled so that $\lambda_{\max}(A^\top A)=1$ & vector of ones \\
\hline
\end{tabular}
\end{center}
\caption{Four types of least squares. }
\label{tbl:sk_leastsquares}
\end{table}

We shall mention that this experiment serves to support and give a better intuition about the results from Section~\ref{N-SEGA} and is by no means practical. The results show, as predicted by theory, linear convergence to a specific neighborhood of the objective. The effect of the noise varies, however, as a general rule, the larger strong convexity $\mu$ is (i.e. problems 1,3 where scaling was not applied), the smaller the effect of noise is.

\begin{figure}[!h]
\centering
\begin{minipage}{0.35\textwidth}
  \centering
\includegraphics[width =  \textwidth ]{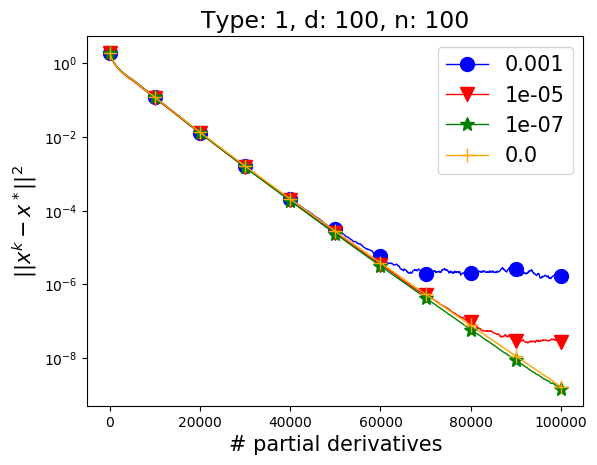}
\end{minipage}%
\begin{minipage}{0.35\textwidth}
  \centering
\includegraphics[width =  \textwidth ]{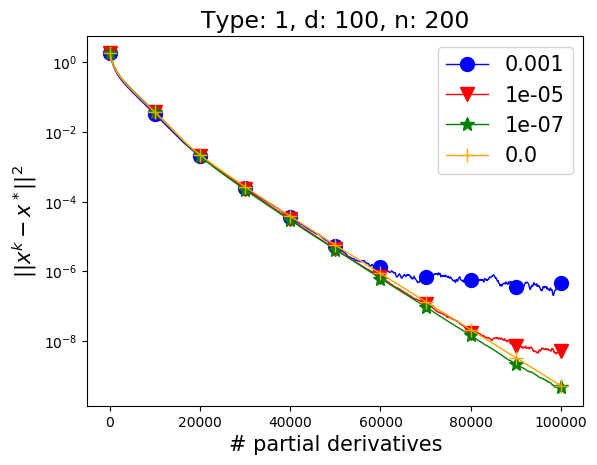}
\end{minipage}%
\\
\begin{minipage}{0.35\textwidth}
  \centering
\includegraphics[width =  \textwidth ]{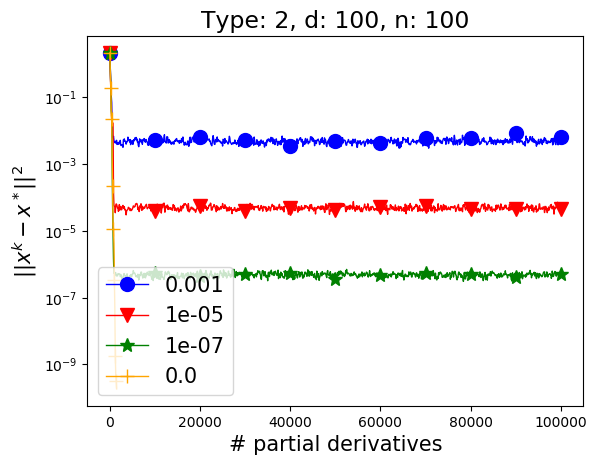}
\end{minipage}%
\begin{minipage}{0.35\textwidth}
  \centering
\includegraphics[width =  \textwidth ]{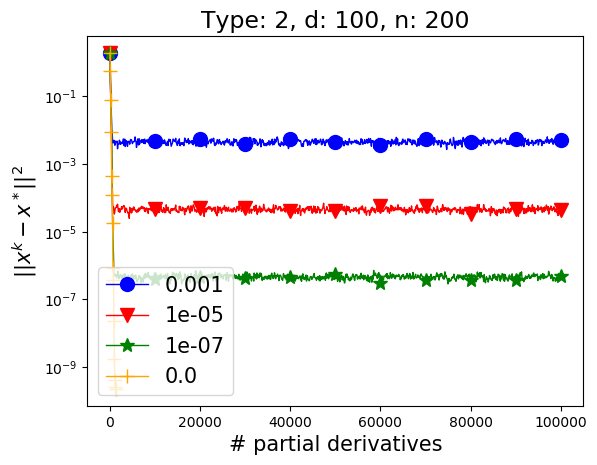}
\end{minipage}%
\\
\begin{minipage}{0.35\textwidth}
  \centering
\includegraphics[width =  \textwidth ]{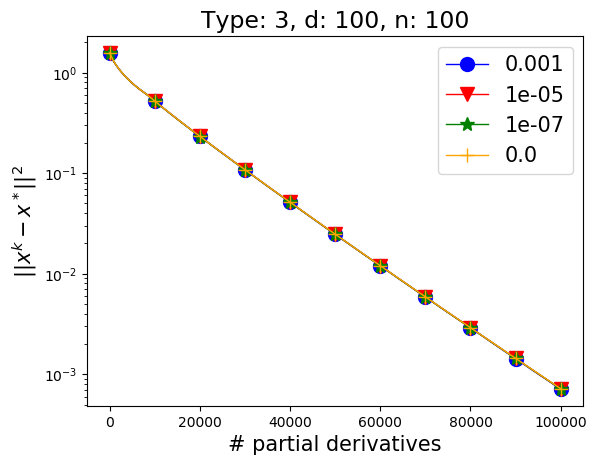}
\end{minipage}%
\begin{minipage}{0.35\textwidth}
  \centering
\includegraphics[width =  \textwidth ]{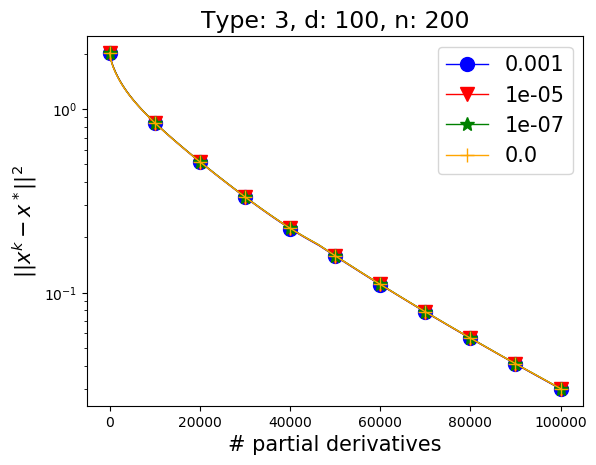}
\end{minipage}%
\\
\begin{minipage}{0.35\textwidth}
  \centering
\includegraphics[width =  \textwidth ]{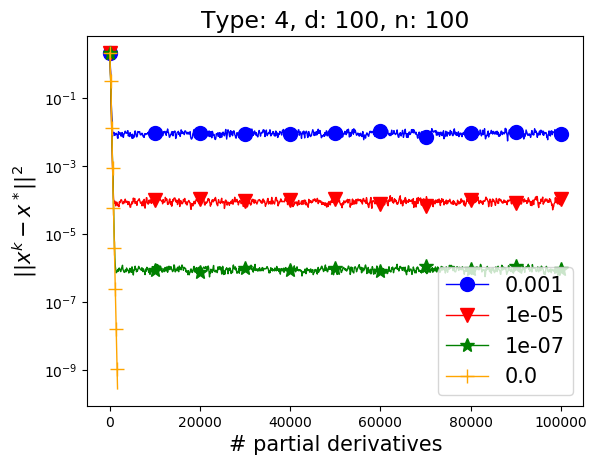}
\end{minipage}%
\begin{minipage}{0.35\textwidth}
  \centering
\includegraphics[width =  \textwidth ]{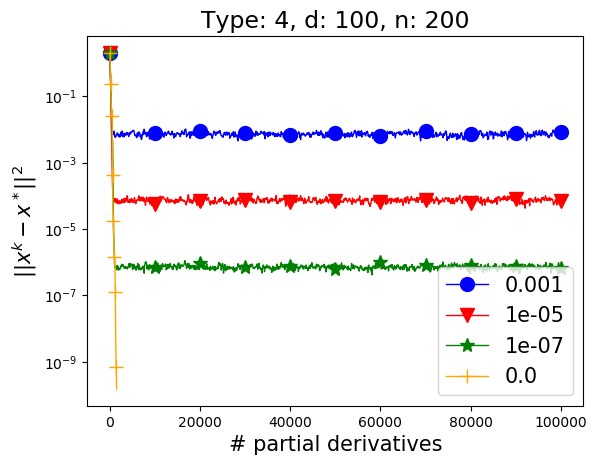}
\end{minipage}%
\caption{ {\tt N-SEGA} applied on constrained least squares problem with noised partial derivative oracle. Legend labels stand for the magnitude $\sigma^2$ of the oracle noise. }
\label{fig:sk_nsega}
\end{figure}

\section{Conclusion}
In this chapter we have introduced a general scheme to analyze to analyze stochastic gradient algorithms with many different applications. Although the presented approach is rather general, we still see several possible directions for  future extensions, including: 
\begin{itemize}
\item  We believe our results can be extended to {\em weakly convex} functions. However, producing a comparable result in the {\em nonconvex} case remains a major open problem. 

\item It would be further interesting to unify our theory with {\em biased} gradient estimators. If this was possible, one could recover methods as {\tt SAG}~\cite{sag} in special cases, or obtain rates for the zero-order optimization. We have some preliminary results in this direction already.

\item Although our theory allows for non-uniform stochasticity, it does not recover the best known rates for {\tt RCD} type methods with {\em importance sampling}. It would be thus interesting to provide a more refined  analysis capable of capturing importance sampling phenomena more accurately.

\item An extension of Assumption~\ref{as:sk_general_stoch_gradient} to {\em iteration dependent} parameters $A,B,C,D_1, D_2, \rho$ would enable an array of new methods, such as {\tt SGD} with decreasing stepsizes. Such an extension is rather very straightforward. 

\item It would be interesting to provide a unified analysis of stochastic methods with {\em acceleration} and {\em momentum}. In fact,~\cite{kulunchakov2019estimate} provide (separately) a unification of some methods with and without variance reduction. The next chapter provides another step towards the unified accelerated analysis -- we introduce an accelerated {\tt SVRCD} algorithm.

\end{itemize}

\chapter{Variance Reduced Coordinate Descent with Acceleration: New Method With a Surprising Application to Finite-Sum Problems}
\label{asvrcd}

\graphicspath{{asvrcd/images/}}

In this chapter, we aim to solve the regularized optimization problem
\begin{equation}\label{eq:asvrcd_problem}
\compactify \min_{x\in \R^d}  \left\{  F(x) = f(x) + \psi(x) \right \},
\end{equation}
where function $f$ is convex and differentiable (not necessarily of a finite-sum structure), while the regularizer $\psi$ is convex and non-smooth. Furthermore, we assume that the dimensionality $d$ is large. 

The most standard approach to deal with the huge $d$ is to decompose the space, i.e., use coordinate descent, or, more generally, subspace descent methods~\cite{rcdm, wright2015coordinate, kozak2019stochastic}. Those methods are especially popular as they achieve a linear convergence rate on strongly convex problems while enjoying a relatively cheap cost of performing each iteration.

However, coordinate descent methods are only feasible if the regularizer $\psi$ is separable~\cite{richtarik2014iteration}. In contrast, if $\psi$ is not separable, the corresponding stochastic gradient estimator has an inherent (non-zero) variance at the optimum, and thus the linear convergence rate is not achievable. 

This phenomenon is, to some extent, similar when applying Stochastic Gradient Descent ({\tt SGD})~\cite{robbins, nemirovski2009robust} on finite sum objective -- the corresponding stochastic gradient estimator has a (non-zero) variance at the optimum, which prevents {\tt SGD} from converging linearly. Recently, the issue of sublinear convergence of {\tt SGD} has been resolved using the idea of control variates~\cite{hickernell2005control}, resulting in famous variance reduced methods such as {\tt SVRG}~\cite{svrg} and {\tt SAGA}~\cite{saga}. 

Motivated by the massive success of variance reduced methods for finite sums, control variates have been proposed to ``fix'' coordinate descent methods to minimize problem~\eqref{eq:asvrcd_problem} with non-separable $\psi$. To best of our knowledge, there are two such algorithms in the literature---{\tt SEGA} (proposed in Chapter~\ref{sega}) and {\tt SVRCD} (proposed in Chapter~\ref{jacsketch})---which we now quickly describe.\footnote{VRSSD~\cite{kozak2019stochastic} is yet another stochastic subspace descent algorithm aided by control variates; however, it was proposed to minimize $f$ only (i.e., considers $\psi\equiv0$). }

Let $x^k$ be the current iterate of {\tt SEGA} (or {\tt SVRCD}) and suppose that the oracle reveals $\nabla_i f(x^k)$ (for $i$ chosen uniformly at random). The simplest unbiased gradient estimator of $\nabla f(x^k)$ can be constructed as $\tilde{g}^k= d\nabla_i f(x^k)e_i$ (where $e_i\in \R^d$ is the $i$th standard basis vector). The idea behind these methods is to enrich $\tilde{g}^k$ using a control variate $h^k\in \R^d$, resulting in a new (still unbiased) gradient estimator $g^k$:
\[
g^k  = d\nabla_i f(x^k)e_i - dh_i^ke_i + h^k.
\]

\emph{How to choose the sequence of control variates $\{h^k\}$?} Intuitively, we wish for both sequences $\{h^k\}$ and $\{\nabla f(x^k\})$ to have an identical limit point. In such case, we have $\lim_{k \rightarrow \infty }\mathrm{Var}(g^k) =0$, and thus one shall expect faster convergence. There is no unique way of setting $\{h^k\}$ to have the mentioned property satisfied -- this is where {\tt SEGA} and {\tt SVRCD} differ. See Algorithm~\ref{alg:asvrcd_SEGAAS} for details.

\begin{algorithm}[h]
    \caption{{\tt SEGA} and {\tt SVRCD}}
    \label{alg:asvrcd_SEGAAS}
    \begin{algorithmic}
        \Require Stepsize $\alpha>0$, starting point $x^0\in\R^d$, probability vector $p$: $p_i\eqdef \Probbb{i\in S} $
        \State Set $h^0 = 0 \in \R^d$
        \For{ $k=0,1,2,\ldots$ }
        \State{Sample random  $S\subseteq \{1,2,\dots,d\}$}
        \State{$g^k = \sum \limits_{i\in S}  \frac{1}{p_i}(\nabla_i f(x^k) - h_i^k)e_i + h^k$}
        \State{$x^{k+1} = \prox_{\alpha \psi}(x^k - \alpha g^k)$}
        \State{$h^{k+1} = \begin{cases} h^{k} +\sum \limits_{i\in S}( \nabla_i f(x^k) - h^{k}_i)e_i & \text{for {\tt SEGA}} \\  \begin{cases} \nabla f(x^k) & \text{with probability } \probx \\ h^k & \text{with probability } 1-\probx\end{cases}&  \text{for {\tt SVRCD}} \end{cases} $}
        \EndFor
    \end{algorithmic}
\end{algorithm}

In this work, we continue the above  research  along the lines of variance reduced coordinate descent algorithms, with surprising consequences.

\section{Contributions\label{sec:asvrcd_contrib}}
Here we list the main contributions of this chapter. 

\begin{itemize}
\item {\bf Exploiting prox in {\tt SEGA}/{\tt SVRCD}.} Assume that the regularizer $\psi$ includes an indicator function of some affine subspace of $\R^d$. We show that both {\tt SEGA} and {\tt SVRCD} might exploit this fact, resulting in a faster convergence rate. As a byproduct, we establish the same result in the more general {\tt GJS} framework from Chapter~\ref{jacsketch} (presented in the appendix).

\item {\bf Accelerated {\tt SVRCD}.} We propose an accelerated version of {\tt SVRCD} - {\tt ASVRCD}. {\tt ASVRCD} is the first accelerated variance reduced coordinate descent to minimize objectives with non-separable, proximable regularizer.\footnote{We shall note that an accelerated version of {\tt SEGA} was already proposed in~\cite{sega} for $\psi = 0$ -- this was rather an impractical result demonstrating that {\tt SEGA} can match state-of-the art convergence rate of accelerated coordinate descent from~\cite{allen2016even, nesterov2017efficiency, hanzely2018accelerated}. In contrast, our results cover any convex $\psi$.}

\item {\bf {\tt SEGA}/{\tt SVRCD}/{\tt ASCRVD} generalizes {\tt SAGA}/{\tt LSVRG}/L-Katyusha.} We show a surprising link between {\tt SEGA} and {\tt SAGA}. In particular, {\tt SAGA} is a special case of {\tt SEGA}; and the new rate we obtain for {\tt SEGA} recovers the tight complexity of {\tt SAGA}~\cite{qian2019saga, gazagnadou2019optimal}. Similarly, we recover loopless {\tt SVRG} ({\tt LSVRG})~\cite{hofmann2015variance, kovalev2019don} along with its best-known rate~\cite{hanzely2019one, l-svrg-as} using a result for {\tt SVRCD}. Lastly, as a particular case of {\tt ASVRCD}, we recover an algorithm which is marginally preferable to loopless Katyusha (L-Katyusha)~\cite{l-svrg-as}: while we recover their iteration complexity result, our proof is more straightforward, and at the same time, the stepsize for the proximal operator is smaller.\footnote{This is preferable especially if the proximal operator has to be estimated numerically.}
\end{itemize}

\section{Preliminaries}

As mentioned in Section~\ref{sec:asvrcd_contrib}, the new results we provide i are particularly interesting if the regularizer $\psi$ contains an indicator function of some affine subspace of $\R^d$. 
 
 \begin{assumption}\label{ass:asvrcd_indicator}

 Assume that $\Popt$ is a projection matrix such that 
 \begin{equation}\label{eq:asvrcd_danjadsou}
\psi (x) = \begin{cases}
 \psi'(x) & {\text if}\quad x \in \{ x^0 + \Range{ \Popt} \} \\
\infty & {\text if} \quad  x \not \in \{ x^0 + \Range{ \Popt} \} \end{cases}  \end{equation}
for some convex function $\psi'(x)$. Furthermore, suppose that the proximal operator of $\psi$ is cheap to compute. 
 \end{assumption}

\begin{remark}
If $\psi$ is convex, there is always some $\Popt$ such that~\eqref{eq:asvrcd_danjadsou} holds as one might choose $\Popt=\mI$.
\end{remark}

Next, we require smoothness of the objective, as well as the strong convexity over the affine subspace given by $\Range{\Popt}$.

\begin{assumption}\label{as:asvrcd_smooth_strongly_convex}
Function $f$ is $\mM$-smooth, i.e.,  for all $ x,y\in \R^d$:\footnote{We define $\| x\|^2 \eqdef \langle x,x\rangle$ and $ \|x \|^2_{\mM} \eqdef \langle \mM x, x\rangle$.}
\[
f(x)\leq f(y) + \langle \nabla f(y),x-y\rangle + \frac{1}{2}\| x-y\|^2_{\mM}.
\]
Function $f$ is $\mu$-strongly convex over $ \{ x^0 + \Range{ \Popt} \}$, i.e.,  for all $ x,y\in\{ x^0 + \Range{ \Popt} \}$:
\begin{align}
 f(x)&\geq f(y) + \langle \nabla f(y),x-y\rangle + \frac{\mu}{2}\| x-y\|^2.\label{eq:asvrcd_sc}
\end{align}
\end{assumption}

\begin{remark}
Smoothness with respect to matrix $\mM$ arises naturally in various applications. For example, if $f(x) = f'(\mA x)$, where $f'$ is $L'$-smooth (for scalar $L'>0$), we can derive that $f$ is $\mM = L'\mA^\top \mA$-smooth.
\end{remark}

In order to stress the distinction between the finite sum setup and the setup from the rest of the chapter, we are denoting the finite-sum variables that differ from the non-finite sum case in ${\color{red} \text{red}}$. We thus, {\emph{recommend printing this chapter in color}}. 

\section{Better rates for {\tt SEGA} and {\tt SVRCD}  \label{sec:asvrcd_sega_is2}}
In this section, we show that a specific structure of nonsmooth function $\psi$ might lead to faster convergence of {\tt SEGA} and {\tt SVRCD}.

The next lemma is a direct consequence of Assumption~\ref{ass:asvrcd_indicator} -- it shows that proximal operator of $\psi$ is contractive under $\Popt$-norm.
  \begin{lemma}
Let $\{x^k\}_{k\geq0}$ be a sequence of iterates of Algorithm~\ref{alg:asvrcd_SEGAAS} and let $x^*$ be optimal solution of~\eqref{eq:asvrcd_problem}. Then 
\begin{equation}
x^k  \in \{x^0 + \Range{\Popt}\},\;  x^* \in\{x^0 + \Range{\Popt)}\} . \label{eq:asvrcd_q_identity}
\end{equation}
for all $k$. Furthermore, for any $x,y\in \R^{d}$ and $\alpha>0$ we have
\begin{equation}\label{eq:asvrcd_prox_cont}
\| \prox_{\alpha \psi}(x ) - \prox_{\alpha \psi}(y) \|^2\leq \| x-y\|^2_{\Popt}.
\end{equation}

  \end{lemma}

Next, we state the convergence rate of both {\tt SEGA} and {\tt SVRCD} under Assumption~\ref{ass:asvrcd_indicator} as Theorem~\ref{thm:asvrcd_sega_as}. We also generalize the main theorem from Chapter~\ref{jacsketch} (fairly general algorithm which covers {\tt SAGA}, {\tt SVRG}, {\tt SEGA}, {\tt SVRCD}, and more as a special case; see Section~\ref{sec:asvrcd_analysis2} of the appendix); from which the convergence rate of {\tt SEGA}/{\tt SVRCD} follows as a special case. 

\begin{theorem}\label{thm:asvrcd_sega_as}

Let Assumptions~\ref{ass:asvrcd_indicator},~\ref{as:asvrcd_smooth_strongly_convex} hold and denote $p_i\eqdef \Probbb{i\in S}$. Consider vector $v = \sum_{i=1}^d e_i v_i, v_i\geq 0$ such that
\begin{equation} \label{eq:asvrcd_ESO_sega_good}
\mM^{\frac12} \E{  \sum_{i\in S} \frac{1}{p_i} e_i e_i^\top\Popt \sum_{i\in S} \frac{1}{p_i} e_i e_i^\top} \mM^{\frac12}  \preceq \diag(p^{-1}\circ v),
\end{equation}
where $\diag(\cdot)$ is a diagonal operator.\footnote{Returns matrix with the input on the diagonal, zeros everywhere else.}
Then, iteration complexity of {\tt SEGA} with $\alpha  = \min_i \frac{p_i}{4v_i+ \mu}$ is $  \max_i\left(\frac{4v_i + \mu}{p_i \mu} \right)\log\frac1\epsilon$. At the same time, iteration complexity of {\tt SVRCD} with $\alpha  = \min_i \frac{1}{4v_ip_i^{-1}+  \mu\probx^{-1}}$ is $  \left(\frac{4\max_i(v_ip_i^{-1} )+ \mu \probx^{-1}}{ \mu} \right)\log\frac1\epsilon$. 
\end{theorem}

Let us look closer to convergence rate of {\tt SVRCD} from Theorem~\ref{thm:asvrcd_sega_as}. The optimal vector $v$ is a solution to the following optimization problem
\[
\min_{v\in \R^d}   \;\; \left(\frac{4\max_i\{v_ip_i^{-1} \}+ \mu \probx^{-1}}{ \mu} \right)\log\frac1\epsilon  \quad \text{such that}  \quad \eqref{eq:asvrcd_ESO_sega_good}\; \text{holds}. 
\]
Clearly, there exists a solution of the form $v\propto p$; let us thus choose $v \eqdef \ccL p$ with $\ccL>0$. In this case, to satisfy~\eqref{eq:asvrcd_ESO_sega_good} we must have 
\begin{equation}\label{eq:asvrcd_ccLdef}
\ccL = \lambda_{\max}\left( \mM^{\frac12} \E{  \sum_{i\in S} \frac{1}{p_i} e_i e_i^\top\Popt \sum_{i\in S} \frac{1}{p_i} e_i e_i^\top} \mM^{\frac12} \right)
\end{equation}
and the iteration complexity of {\tt SVRCD} becomes $\left(\frac{4\ccL+ \mu \probx^{-1}}{ \mu} \right)\log\frac1\epsilon$.\footnote{We decided to not present this, simplified rate in Theorem~\ref{thm:asvrcd_sega_as} for the following two reasons: 1) it would yields a slightly subpotimal rate of {\tt SEGA} and 2) the connection of to the convergence rate of {\tt SAGA} from~\cite{qian2019saga} is more direct via~\eqref{eq:asvrcd_ESO_sega_good}.}

\emph{How does $\Popt$ influence the rate?} As mentioned, one can always consider $\Popt=\mI$. In such a case, we recover the convergence rate of {\tt SEGA} and {\tt SVRCD} from Chapter~\ref{jacsketch}. However, the smaller rank of $\Popt$ is, the faster rate is Theorem~\ref{thm:asvrcd_sega_as} providing. To see this, it suffices to realize that if $\ccL$ is increasing in $\Popt$ (in terms of Loewner ordering). 

\begin{example}
Let $\mM= \mI$ and $S=\{i\}$ with probability $d^{-1}$ for all $1\leq i\leq d$. Given that $\Popt=\mI$, it is easy to see that $\ccL = d$. In such case, the iteration complexity of {\tt SVRCD} is $\left(\frac{4d + \mu \probx^{-1}}{ \mu} \right)\log\frac1\epsilon$. In the other extreme, if $\Popt=\frac1d ee^\top$, we have $\ccL=1$, which yields complexity (of {\tt SVRCD}) $\left(\frac{4 + \mu \probx^{-1}}{ \mu} \right)\log\frac1\epsilon$. Therefore, given that $\mu = \cO(\probx)$, the low rank of $\Popt$ caused the speedup of order $\Theta(d)$.
\end{example}

We shall also note that the tight rate of {\tt SAGA} and {\tt LSVRG} might be recovered from Theorem~\ref{thm:asvrcd_sega_as} only using a non-trivial $\Popt$ (see Section~\ref{sec:asvrcd_sagasega}), while the original theory of {\tt SEGA} and {\tt SVRCD} only yield a suboptimal rate for both {\tt SAGA} and {\tt LSVRG}.

\paragraph{Connection with Subspace {\tt SEGA} (from Section~\ref{sec:sega_subSEGA}).} Assume that function $f$ is of structure $f(x)=h(\mA x)$. As a consequence, we have $\nabla f(x) = \mA^\top \nabla h(\mA x)$ and thus $\nabla f(x)\in \Range{\mA^\top}$. This fact was exploited by Subspace {\tt SEGA} in order to achieve a faster convergence rate. Our results can mimic Subspace {\tt SEGA} by setting $\psi$ to be an indicator function of $x^0 + \Range{\mA^\top}$, given that there is no extra non-smooth term in the objective.

\begin{remark}
Throughout all proofs of this section, we have used a weaker conditions than Assumption~\ref{as:asvrcd_smooth_strongly_convex}. In particular, instead of-$\mM$-smoothness, it is sufficient to have\footnote{By $D_f(x,y)$ we denote Bregman distance between $x,y$, i.e., $D_f(x,y)\eqdef f(x)-f(y)-\langle \nabla f(x)$ } 
$
D_{f}(x,x^*) \geq \frac12 \norm{ \nabla f(x)-\nabla f(x^*) }^2_{\mM^{-1}}
$
for all $x\in \R^d$ (Lemma~\ref{lem:asvrcd_smooth_consequence} shows that it is indeed a consequence of $\mM$ smoothness and convexity). At the same time, instead of $\mu$-strong convexity, it is sufficient to have $\mu$-quasi strong convexity, i.e.,  for all $x\in \{x^0 + \Range{\Popt} \}$:
$
f(x^*) \geq f(x) + \langle \nabla f(x), x^*-x\rangle + \frac{\mu}{2}\|x-x^* \|^2. 
$
However, the accelerated method (presented in Section~\ref{sec:asvrcd_acc}) requires the fully general version of Assumption~\ref{as:asvrcd_smooth_strongly_convex}.
\end{remark}

\section{Connection between {\tt SEGA} ({\tt SVRCD}) and {\tt SAGA} ({\tt LSVRG})  \label{sec:asvrcd_sagasega}}

In this section, we show that {\tt SAGA} and {\tt LSVRG} are special cases of {\tt SEGA} and {\tt SVRCD}, respectively. At the same time, the previously tightest convergence rate of {\tt SAGA}~\cite{gazagnadou2019optimal, qian2019saga} and {\tt LSVRG}~\cite{hanzely2019one, l-svrg-as} follow from Theorem~\ref{thm:asvrcd_sega_as} (convergence rate of {\tt SEGA} and {\tt SVRCD}). 

\subsection{Convergence rate of {\tt SAGA} and {\tt LSVRG} \label{sec:asvrcd_saga}}
We quickly state the best-known convergence rate for both {\tt SAGA} and {\tt LSVRG} to minimize the following objective:
\begin{equation}\label{eq:asvrcd_problem_finitesum}
\compactify \min_{\xx\in \R^\dd}   \left \{ \PpP(\xx) \eqdef \underbrace{\frac1n \sum \limits_{j=1}^n \ff_j(\xx)}_{\eqdef \ff(\xx)} + \ppsi(\xx) \right \}.
\end{equation}

\begin{assumption} \label{as:asvrcd_finitesum}
 Each $\ff_j$ is convex, $\ueWMC_j$-smooth and $\ff $ is $\mmu$-strongly convex. 
\end{assumption}

Assuming the oracle access to $\nabla \ff_i(\xx^k)$ for $i\in \sS$ (where $\sS$ is a random subset of $\{1,\dots,n \}$), the minibatch {\tt SGD}~\cite{pmlr-v97-qian19b} uses moves in the direction of the ``plain'' unbiased stochastic gradient $\frac1n \sum \limits_{i\in \sS}  \frac{1}{\ppp_i}\nabla \ff_i(\xx^k) $ (where $\pp_i \eqdef \Probbb{i\in \sS} $).

In contrast, variance reduced methods such as {\tt SAGA} and {\tt LSVRG} enrich the ``plain'' unbiased stochastic gradient with control variates: 
\begin{equation}
\ggggg^k = \frac1n\sum \limits_{i\in \sS}  \frac{1}{\ppp_i}\left(\nabla \ff_i(\xx^k) - \mJ^k_{:,i} \right)  +  \frac1n \mJ^k \ee.
\end{equation} 
where $\mJ^k \in \R^{\dd\times n}$ is the control matrix and $\ee \in\R^n$ is vector of ones.
The difference between {\tt SAGA} and {\tt LSVRG} lies in the procedure to update $\mJ^k$; {\tt SAGA} uses the freshest gradient information to replace corresponding columns in $\mJ^k$; i.e.
\begin{equation}\label{eq:asvrcd_saga_update}
\mJ^{k+1}_{:,i} = \begin{cases}
\nabla \ff_i(\xx^k)  & \text{if } i\in \sS \\
\mJ^{k}_{:,i}   & \text{if } i \not \in \sS.
\end{cases}
\end{equation}
On the other hand, {\tt LSVRG} sets $\mJ^k$ to the true Jacobian of $f$ upon a successful, unfair coin toss:
 \begin{equation}\label{eq:asvrcd_lsvrg_update}
\mJ^{k+1} = \begin{cases}
\left[ \nabla \ff_1(\xx^k), \dots,  \nabla \ff_n(\xx^k)\right] & \text{w. p. } \probx \\
\mJ^{k}   & \text{w. p. } 1- \probx.
\end{cases}
\end{equation}

The formal statement of {\tt SAGA} and {\tt LSVRG} is provided in as Algorithm~\ref{alg:asvrcd_saga}, while Proposition~\ref{prop:sagarate} states their convergence rate.

\begin{algorithm}[!h]
	\caption{{\tt SAGA}/{\tt LSVRG}}
	\label{alg:asvrcd_saga}
	\begin{algorithmic}
		\Require $\alpha > 0$, $\probx \in (0,1)$
		\State $\xx^0 \in \R^\dd, \mJ^{0} = 0 \in \R^{\dd\times n}$
		\For{$k=0,1,2,\ldots$}
		\State Sample random $\sS \subseteq \{1,\dots,n \}$
        \State{$\ggggg^k = \frac1n \mJ^k \ee+ \frac1n\sum \limits_{i\in \sS}  \frac{1}{\ppp_i}(\nabla \ff_i(\xx^k) - \mJ^k_{:,i}) $}
			\State $\xx^{k+1} = \prox_{\aalpha \psi} (\xx^k - \aalpha \ggggg^k)$
        \State{Update $\mJ^{k+1}$ according to~\eqref{eq:asvrcd_saga_update} or~\eqref{eq:asvrcd_lsvrg_update}}
		\EndFor
	\end{algorithmic}
\end{algorithm}

\begin{proposition}\label{prop:sagarate}(\cite{hanzely2019one})
Suppose that Assumption~\ref{as:asvrcd_finitesum} holds and let $\vv$ be a nonegative vector such that for all $h_1,\dots,h_n \in \R^{\dd}$ we have
\begin{equation} \label{eq:asvrcd_ESO_saga}
\E{\left\|\sum_{j \in \sS} \ueWMC^{\frac12}_{j} h_{j}\right\|^{2}} \leq \sum_{j=1}^{n} \pp_j \vv_{j}\left\|h_{j}\right\|^{2} .
\end{equation}
Then the iteration complexity of {\tt SAGA} with $ \aalpha  = \min_j \frac{n \pp_j }{4\vv_j + n\mmu}$ is $\max_j \left(  \frac{4\vv_j  + n \mmu }{n  \mmu \pp_j }\right) \log\frac1\epsilon$. At the same time, iteration complexity of {\tt LSVRG} with $ \aalpha  = \min_j \frac{n}{4 \frac{\vv_j}{ \pp_j } + \frac{  \mmu n}{\probx}}$is $\max_j \left(  4 \frac{\vv_j}{n  \mmu \pp_j  }  + \frac{1}{ \probx } \right) \log\frac1\epsilon$.
\end{proposition}

\subsection{{\tt SAGA} is a special case of {\tt SEGA}}
Consider setup from Section~\ref{sec:asvrcd_saga}; i.e.,  problem~\eqref{eq:asvrcd_problem_finitesum} along with Assumption~\ref{as:asvrcd_finitesum} and $\vv$ defined according to~\eqref{eq:asvrcd_ESO_saga}. We will construct an instance of~\eqref{eq:asvrcd_problem} (i.e.,  specific $f$, $\psi$), which is equivalent to~\eqref{eq:asvrcd_problem_finitesum}, such that applying {\tt SEGA} on~\eqref{eq:asvrcd_problem} is equivalent applying {\tt SAGA} on~\eqref{eq:asvrcd_problem_finitesum}.

Let $d \eqdef \dd n$. 

For convenience, define $R_j \eqdef \{ \dd(j-1)+1,  \dd(j-1)+1, \dots,  \dd j  \}$ (i.e., $|R_j|= \dd$) and lifting operator $\Lift{\cdot}: \R^\dd \rightarrow \R^{d}$ defined as 
$\Lift{\xx} \eqdef \left [\underbrace{\xx^\top, \dots ,\xx^\top}_{n\, \mathrm{times}}\right]^\top$.

\paragraph{Construction of $f$, $\psi$.}  Let $\Ind{}$ be indicator function of the set\footnote{Indicator function of a set returns 0 for each point inside of the set and $\infty$ for each point outside of the set.} $x_{R_1}=\dots=x_{R_n}$ and choose
\begin{eqnarray}
f(x) \eqdef  \frac1n\sum_{j=1}^n \ff_j(x_{R_j}),  \;\; 
 \psi(x) \eqdef  \Ind{} (x) + \ppsi(x_{R_1})\label{eq:asvrcd_equivalent}
\end{eqnarray}
Now, it is easy to see that problem~\eqref{eq:asvrcd_problem_finitesum} and problem~\eqref{eq:asvrcd_problem} with the choice~\eqref{eq:asvrcd_equivalent} are equivalent; each $x\in \R^d$ such that $F(x)<\infty$ must be of the form $x = \Lift{\xx}$ for some $\xx\in \R^\dd$. In such case, we have $F(x) = \PpP(\xx)$. The next lemma goes further, and derives the values $\mM, \mu, \Popt$ and $v$ based on $\ueWMC_i$ ($1\leq i\leq n$), $\mmu, \vv$.

\begin{lemma}\label{lem:asvrcd_equivalent_objectives}
Consider $f, \psi$ defined by~\eqref{eq:asvrcd_equivalent}. Function $f$ satisfies Assumption~\ref{as:asvrcd_smooth_strongly_convex} with $\mu \eqdef \frac{ \mmu}{n}$ and 
$\mM \eqdef \frac1n \blockdiag(\ueWMC_{1}, \dots, \ueWMC_n)$. Function $\psi$ and $x^0 = \Lift{\xx^0}$ satisfy Assumption with $\Popt \eqdef \frac1n \ee\ee^\top \otimes \mI$. At the same time, given that $\vv$ satisfies~\eqref{eq:asvrcd_ESO_saga}, inequality~\eqref{eq:asvrcd_ESO_sega_good} holds with $v = \vv n^{-1}$. 
\end{lemma} 

Next, we show that running Algorithm~\ref{alg:asvrcd_SEGAAS} in this particular setup is equivalent to running Algorithm~\ref{alg:asvrcd_saga} for the finite sum objective.

\begin{lemma}\label{lem:asvrcd_saga_from_sega}
Consider $f, \psi$ from~\eqref{eq:asvrcd_equivalent}, $S$ as described in the last paragraph and 
$x^0 =\Lift{\xx^0}$. Running {\tt SEGA} ({\tt SVRCD}) on~\eqref{eq:asvrcd_problem} with $S \eqdef \cup_{j\in \sS} R_j$ and $\alpha \eqdef n\aalpha$ is equivalent to running {\tt SAGA} ({\tt LSVRG}) on~\eqref{eq:asvrcd_problem_finitesum}; i.e.,  we have for all $k$  
\begin{equation}x^k = \Lift{\xx^k}.\label{eq:asvrcd_iterates_equivalence}\end{equation}
\end{lemma} 

As a consequence of Lemmas~\ref{lem:asvrcd_equivalent_objectives} and~\ref{lem:asvrcd_saga_from_sega}, we get the next result.

\begin{corollary}\label{cor:asvrcd_saga_as2++}  Let $f, \psi, S$ be as described above. Convergence rate of {\tt SAGA} ({\tt LSVRG}) given by Proposition~\ref{prop:sagarate} to solve~\eqref{eq:asvrcd_problem} is identical to convergence rate of {\tt SEGA} ({\tt SVRCD}) given by Theorem~\ref{thm:asvrcd_sega_as}.
\end{corollary}

\section{The {\tt ASVRCD} algorithm\label{sec:asvrcd_acc}}
In this section we present {\tt SVRCD} with Nesterov's momentum~\cite{nesterov83} -- {\tt ASVRCD}. 
The development of {\tt ASVRCD} along with the theory (Theorem~\ref{thm:asvrcd_acc}) was motivated by Katyusha~\cite{allen2017katyusha}, {\tt ASVRG}~\cite{asvrg} and their loopless variants~\cite{kovalev2019don, l-svrg-as}. In Section~\ref{sec:asvrcd_kat_special}, we show that a variant of L-Katyusha (Algorithm~\ref{alg:asvrcd_katyusha2}) is a special case of {\tt ASVRCD}, and argue that it is slightly superior to the methods mentioned above.

The main component of {\tt ASVRCD} is the gradient estimator $g^k$ constructed analogously to {\tt SVRCD}. In particular, $g^k$ is an unbiased estimator of $\nabla f(x^k)$ controlled by $\nabla f(w^k)$:\footnote{This is efficient to implement as sequence of iterates $\{w^k\}$ is updated rarely.}
\begin{equation}
g^k = \nabla f(w^k)+\sum \limits_{i\in S}  \frac{1}{p_i}(\nabla_i f(x^k) - \nabla_i f(w^k))\ones_i.
\end{equation}

Next, {\tt ASVRCD} requires two more sequences of iterates $\{ y^k\}_{k\geq0}, \{ z^k\}_{k\geq 0}$ in order to incorporate Nesterov's momentum. The update rules of those sequences consist of subtracting $g^k$ alongside with convex combinations or interpolations of the iterates. See Algorithm~\ref{alg:asvrcd_acc} for specific formulas.  

\begin{algorithm}[h]
	\caption{Accelerated {\tt SVRCD} ({\tt ASVRCD})}
	\label{alg:asvrcd_acc}
	\begin{algorithmic}
		\Require $0< \theta_1, \theta_2 <1$, $\eta, \beta , \gamma > 0$, $\probx \in (0,1)$, $y^0 = z^0 = x^0 \in \R^d$
		\For{$k=0,1,2,\ldots$}
			\State $x^k = \theta_1 z^k + \theta_2 w^k + ( 1 -\theta_1 -\theta_2) y^k$
        \State{Sample random  $S\subseteq \{1,2,\dots,d\}$}
        \State{$g^k = \nabla f(w^k)+\sum \limits_{i\in S}  \frac{1}{p_i}(\nabla_i f(x^k) - \nabla_i f(w^k))\ones_i$}
			\State $y^{k+1} = \prox_{\eta \psi} (x^k - \eta g^k)$
			\State $z^{k+1} = \beta z^k + (1-\beta)x^k + \frac{\gamma}{\eta}(y^{k+1} - x^k)$
			\State $w^{k+1} = \begin{cases}
				y^k, &\text{ with probability } \probx\\
				w^k, &\text{ with probability } 1-\probx\\
			\end{cases}$
		\EndFor
	\end{algorithmic}
\end{algorithm}

We are now ready to present {\tt ASVRCD} along with its convergence guarantees.

\begin{theorem} \label{thm:asvrcd_acc}
Let Assumption~\ref{ass:asvrcd_indicator},~\ref{as:asvrcd_smooth_strongly_convex} hold and denote $L \eqdef \lambda_{\max} \left( \mM^{\frac12} \Popt \mM^{\frac12}\right)$.
Further, let $\Lcac$ be such that for all $k$ we have
	\begin{equation}\label{eq:asvrcd_exp:smooth}
		\E{\norm{g^k - \nabla f(x^k)}^2_\Popt} \leq 2\Lcac D_f(w^k,x^k).
	\end{equation}
 Define the following Lyapunov function:
	\begin{eqnarray*}
	\Psi^k &\eqdef&  \norm{z^k - x^*}^2 + \frac{2\gamma\beta}{\theta_1}\left[F(y^k) - F(x^*)\right] + \frac{(2\theta_2 + \theta_1)\gamma\beta}{\theta_1\probx}\left[F(w^k) - F(x^*)\right],
	\end{eqnarray*}
and let
\begin{eqnarray*}
\eta &=&  \frac14 \max\{\Lcac, L\}^{-1}, \\
\theta_2 &=& \frac{\Lcac}{2\max\{L, \Lcac\}}, \\
\gamma &=& \frac{1}{\max\{2\mu, 4\theta_1/\eta\}},\\
 \beta &=& 1 - \gamma\mu \; \mathrm{and} \\
 \theta_1 &=& \min\left\{\frac{1}{2},\sqrt{\eta\mu \max\left\{\frac{1}{2}, \frac{\theta_2}{\rho}\right\}}\right\} .
\end{eqnarray*}
	Then the following inequality holds:
	\begin{equation*}
		\E{\Psi^{k+1}} \leq
		\left[1 -  \frac{1}{4}\min\left\{\probx, \sqrt{\frac{\mu}{2\max\left\{L, \frac{\Lcac}{\rho}\right\}}} \right\} \right]\Psi^0.
	\end{equation*}

	As a consequence, iteration complexity of Algorithm~\ref{alg:asvrcd_acc} is 
	$\cO\left( \left( \frac{1}{\probx} +  \sqrt{\frac{L}{\mu}}  + \sqrt{\frac{\Lcac}{\probx\mu}}  \right)\log\frac1\epsilon\right)$.
\end{theorem}

Convergence rate of {\tt ASVRCD} depends on constant $\Lcac$ such that~\eqref{eq:asvrcd_exp:smooth} holds. The next lemma shows that $\Lcac$ can be obtained indirectly from $\mM$-smoothness (via $\ccL$), in which case the convergence rate provided by Theorem~\ref{thm:asvrcd_acc} significantly simplifies. 

\begin{lemma}\label{lem:asvrcd_exp_smooth_vs_ESO}
Inequality~\ref{eq:asvrcd_exp:smooth} holds for $\Lcac = \ccL$ (defined in~\eqref{eq:asvrcd_ccLdef}). Further, we have $L\leq \ccL$. Therefore, setting $ \probx \geq \sqrt{\frac{\mu}{\ccL}}$ yields the following complexity of {\tt ASVRCD}: 
\begin{equation}\label{eq:asvrcd_rate_simple}
\cO\left( \sqrt{\frac{\ccL}{\probx\mu}} \log\frac1\epsilon\right).
\end{equation}
\end{lemma}

Setting $\Lcac = \ccL$ might be, however, loose in some cases. In particular, inequality~\eqref{eq:asvrcd_exp:smooth} is slightly weaker than~\eqref{eq:asvrcd_ESO_sega_good} and consequently, the bound bound from Theorem~\ref{thm:asvrcd_acc} is slightly better than~\eqref{eq:asvrcd_rate_simple}. To see this, notice that the proof of Lemma~\ref{lem:asvrcd_exp_smooth_vs_ESO} bounds variance of $g^k + \nabla f(w^k)$ by its second moment. Admittedly, this bound might not worsen the rate by more than a constant factor when $\frac{\E{|S|}}{d}$ is not close to 1. Therefore, bound~\eqref{eq:asvrcd_rate_simple} is good in essentially all practical cases. The next reason why we keep inequality~\eqref{eq:asvrcd_exp:smooth} is that an analogous assumption was required for the analysis of L-Katyusha in~\cite{l-svrg-as} (see Section~\ref{sec:asvrcd_katyusha}) -- and so we can now recover L-Katyusha results directly.

Let us give a quick taste how the rate of {\tt ASVRCD} behaves depending on $\Popt$. In particular, Lemma~\ref{lem:asvrcd_acc_example} shows that nontrivial $\Popt$ might lead to speedup of order $\Theta(\sqrt{d})$ for {\tt ASVRCD}. 

\begin{lemma} \label{lem:asvrcd_acc_example}
Let $S = i$ for each $1\leq i\leq d$ with probability $\frac1d$ and $\probx = \frac1d$. Then, if $\Popt = \mI$, iteration complexity of {\tt ASVRCD} is $\cO\left(d\sqrt{\frac{\lambda_{\max} \mM}{\mu}} \log \frac1\epsilon\right)$. If, however, $\Popt = \frac1d ee^\top$, iteration complexity of {\tt ASVRCD} is $\cO\left(\sqrt{\frac{d\lambda_{\max} \mM}{\mu}} \log \frac1\epsilon\right)$.
\end{lemma}

\section{Connection between {\tt ASVRCD} and L-Katyusha}
Next, we show that L-Katyusha can be seen as a particular case of {\tt ASVRCD}. 

\subsection{Convergence rate of L-Katyusha\label{sec:asvrcd_katyusha}}
In this section, we quickly introduce the loopless Katyusha (L-Katyusha) from~\cite{l-svrg-as} along with its convergence guarantees. In the next section, we show that an improved version of L-Katyusha can be seen as a special case of {\tt ASVRCD}, and at the same time, the tight convergence guarantees from~\cite{l-svrg-as} can be obtained as a special case of Theorem~\ref{thm:asvrcd_acc}.

Consider problem~\eqref{eq:asvrcd_problem_finitesum} and suppose that  $\ff $ is $\LLL$-smooth and $\mmu$-strongly convex. Let $\sS$ be a random subset of $\{1, \dots ,n \}$ (sampled from arbitrary fixed distribution) such that $\ppp_i \eqdef \Probbb{i\in \sS}$. For each $k$ let $\ggggg^k$ be the following unbiased, variance reduced estimator of $\nabla \ff(x^k)$: 
\[
\ggggg^k = \frac1n \left( \sum_{i\in \sS} \ppp_i^{-1} \left(\nabla \ff_i(\xx^k) - \nabla \ff_i(\ww^k) \right) \right) + \nabla \ff(\ww^k).
\]

Next, L-Katyusha requires the variance of $\ggggg^k$ to be bounded by Bregman distance between $\ww^k$ and $\xx^k$ with constant $\LcLc$, as the next assumption states.

\begin{assumption}\label{as:asvrcd_Katyusha} For all $k$ we have
\begin{equation}
\E{\| \ggggg^k - \nabla\ff(\xx^k) \|^2} \leq 2 \LcLc D_f(\ww^k,\xx^k).
\end{equation}
\end{assumption}
 Proposition~\ref{prop:katyusha} provides a convergence rate of L-Katyusha. 
\begin{proposition}(\cite{l-svrg-as})\label{prop:katyusha} Let $\ff $ be $\LLL$-smooth and $\mmu$-strongly convex while Assumption~\ref{as:asvrcd_Katyusha} holds. 
Iteration complexity of L-Katyusha is $\cO\left( \left(\frac{1}{\pp} + \sqrt{\frac{\LLL}{\mmu}} + \sqrt{\frac{\LcLc}{\mmu \pp}} \right)\log \frac1\epsilon \right)$.
\end{proposition}

\subsection{L-Katyusha is a special case of {\tt ASVRCD}\label{sec:asvrcd_kat_special}}
In this section, we show that a modified version of L-Katyusha (Algorithm~\ref{alg:asvrcd_katyusha2}) is a special case of {\tt ASVRCD}. Furthermore, we show that the tight convergence rate of L-Katyusha~\cite{l-svrg-as} follows from Theorem~\ref{thm:asvrcd_acc} (convergence rate of {\tt ASVRCD}).

Consider again $f,\psi$ chosen according to~\eqref{eq:asvrcd_equivalent}. With this choice, problem~\eqref{eq:asvrcd_problem} and~\eqref{eq:asvrcd_problem_finitesum} are equivalent. At the same time, Lemma~\ref{lem:asvrcd_saga_from_sega} establishes that $f$ satisfies Assumption~\ref{as:asvrcd_smooth_strongly_convex} with $\mu =\frac{ \mmu}{n}$ and $\mM = \frac1n \blockdiag(\ueWMC_{1}, \dots, \ueWMC_n)$ while $\psi$ and $x^0$ satisfy Assumption with $\Popt \eqdef \frac1n \ee\ee^\top \otimes \mI$.

 Note that the update rule of sequences $x^k, z^k, w^k$ are identical for both algorithms; we shall thus verify that the update rule on $y^k$ is identical as well. The last remaining thing is to relate $\Lcac$ and $\LcLc$. The next lemma establishes both results.

\begin{lemma}\label{lem:asvrcd_katyusha_from_asvrcd}
Running {\tt ASVRCD} on~\eqref{eq:asvrcd_problem} with $S \eqdef \cup_{j\in \sS} R_j$ and $\eta \eqdef n\eeta$, $\gamma \eqdef n\ggamma$ is equivalent to running Algorithm~\ref{alg:asvrcd_katyusha2} on~\eqref{eq:asvrcd_problem_finitesum}. At the same time, inequality~\ref{eq:asvrcd_exp:smooth} holds with $\Lcac = n^{-1}\LcLc$, while we have $L = n^{-1}\LLL$.
\end{lemma} 

As a direct consequence of Lemma~\ref{lem:asvrcd_katyusha_from_asvrcd} and Theorem~\ref{thm:asvrcd_acc}, we obtain the next corollary.

\begin{corollary}\label{cor:asvrcd_saga_as2++}  Let $f, \psi, S$ be as described above. Iteration complexity of Algorithm~\ref{alg:asvrcd_katyusha2} is \[\cO\left( \left(\frac{1}{\pp} + \sqrt{\frac{\LLL}{\mmu}} + \sqrt{\frac{\LcLc}{\mmu \pp}} \right)\log \frac1\epsilon \right).\]
\end{corollary}

As promised, the convergence rate of Algorithm~\ref{alg:asvrcd_katyusha2} matches the convergence rate of L-Katyusha from Proposition~\ref{prop:katyusha} and thus matches the lower bound for finite sum minimization by~\cite{woodworth2016tight}. Let us now argue that Algorithm~\ref{alg:asvrcd_katyusha2} is slightly superior to other accelerated {\tt SVRG} variants. 

First, Algorithm~\ref{alg:asvrcd_katyusha2} is loopless; thus has a simpler analysis and slightly better properties (as shown by~\cite{kovalev2019don}) over Katyusha~\cite{allen2017katyusha} and {\tt ASVRG}~\cite{asvrg}. Next, the analysis is simpler than~\cite{l-svrg-as} (i.e., we do not require one page of going through special cases). At the same time, Algorithm~\ref{alg:asvrcd_katyusha2} uses a smaller stepsize for the proximal operator than L-Katyusha, which is useful if the proximal operator does is estimated numerically. However, Algorithm~\ref{alg:asvrcd_katyusha2} is almost indistinguishable from L-Katyusha if $\ppsi=0$.

\begin{remark}
The convergence rate of L-Katyusha from~\cite{l-svrg-as} allows exploiting the strong convexity of regularizer $\psi$ (given that it is strongly convex). While such a result is possible to obtain in our case, we have omitted it for simplicity. 
\end{remark}

\begin{algorithm}[h]
	\caption{Variant of L-Katyusha (special case of Algorithm~\ref{alg:asvrcd_acc})}
	\label{alg:asvrcd_katyusha2}
	\begin{algorithmic}
		\Require $0< \theta_1, \theta_2 <1$, $\eeta, \beta , \ggamma > 0$, $\probx \in (0,1)$
		\State $\yy^0 = \zzz^0 = \xx^0 \in \R^\dd$
		\For{$k=0,1,2,\ldots$}
			\State $\xx^k = \theta_1 \zzz^k + \theta_2 \ww^k + ( 1 -\theta_1 -\theta_2) \yy^k$
        \State{Sample random  $\sS\subseteq \{1,2,\dots,n\}$}
        \State{$g^k = \nabla \ff(\ww^k)+\sum \limits_{i\in \sS}  \frac{1}{\ppp_i}(\nabla \ff_i(\xx^k) - \nabla \ff_i(\ww^k)) $}
			\State $\yy^{k+1} = \prox_{\eeta \psi} (\xx^k - \eeta \ggggg^k)$
			\State $\zzz^{k+1} = \beta \zzz^k + (1-\beta)\xx^k + \frac{\ggamma}{\eeta}(\yy^{k+1} - \xx^k)$
			\State $\ww^{k+1} = \begin{cases}
				\yy^k, &\text{ with probability } \probx\\
				\ww^k, &\text{ with probability } 1-\probx\\
			\end{cases}$
		\EndFor
	\end{algorithmic}
\end{algorithm}

\section{Experiments \label{sec:asvrcd_experiments}}

In this section, we numerically verify the performance of {\tt ASVRCD}, as well as the improved performance of {\tt SVRCD} under Assumption~\ref{ass:asvrcd_indicator}. In order to better understand and control the experimental setup, we consider a quadratic minimization (four different types) over the unit ball intersected with a linear subspace.\footnote{Note that the practicality of {\tt ASVRCD} immediately follows as it recovers Algorithm~\ref{alg:asvrcd_katyusha2} as a special case, which is (especially for $\psi\equiv0$) almost indistinguishable to L-Katyusha -- state-of-the-art method for smooth finite sum minimization. For this reason, we decided to focus on less practical, but better-understood experiments.} 

In all experiments, we have chosen $$f(x) = \frac12 x^\top \mM x - b^\top x,$$ where $x\in \R^{1000}$, while $\psi$ is an indicator function of the unit ball intersected with $\range{\Popt}$. First, matrix $\mM$ was chosen according to Table~\ref{tbl:asvrcd_M}. Next, vector $b$ was chosen as follows: first we generate $\tilde{x}\in \R^d$ with independent normal entries, then compute $\tilde{b} = \mM^{-1} \tilde{x}$ and set $b = \frac{3}{2 \| \tilde{b} \|}  \tilde{b}$. Lastly, for Figure~\ref{fig:asvrcd_variable}, the projection matrix $\Popt$ of rank $r$ was chosen as a block diagonal matrix with $r$ blocks, each of them being the matrix of ones multiplied by $\frac{r}{d}$.

\begin{table}[!h]
\caption{Choice of $\mM$. $\Odd$ is set of all odd positive integers smaller than $d+1$, while matrix $\mU$ was set as random orthonormal matrix (generated by QR decomposition from a matrix with independent standard normal entries).}
\label{tbl:asvrcd_M}
\begin{center}
\small
\begin{tabular}{|c|c|c|c|}
\hline
Type & $\mM$  & Fig.~\ref{fig:asvrcd_identity}: $L$ &  Fig.~\ref{fig:asvrcd_variable}: $L$ \\
\hline
\hline
1 & $\mU \left( \mI + \mI_{:, \Odd}  \diag\left( ((L-1)^{\frac{1}{500}})^{(1:500)} \right)  \mI_{\Odd, :} \right) \mU^\top $  & 100 & 1000 \\
\hline
2 & $\mU \left( \mI +\sum_{i=1}^{100} (L-1) e_ie_i^\top \right) \mU^\top $ & 100 & 1000 \\
\hline
3&$\mU \left( \kappa \mI -\sum_{i=1}^{100} (L-1) e_ie_i^\top \right) \mU^\top $ &100 & 1000 \\
\hline
4 & $\left( \mI + \frac{L}{500}\mI_{:, \Odd}  \diag\left(1:500 \right)  \mI_{\Odd, :} \right) $& 100 & 1000 \\
\hline
\end{tabular}
\end{center}
\end{table}

\subsection{The effect of acceleration and importance sampling}

In the first experiment we demonstrate the superiority of {\tt ASVRCD} to {\tt SVRCD} for problems with $\Popt = \mI$. We consider four different methods -- {\tt ASVRCD} and {\tt SVRCD}, both with uniform and importance sampling such that $|S|=1$ with probability 1. The importance sampling is the same as one from Chapter~\ref{jacsketch}. In short, the goal is to have $\ccL$ from~\eqref{eq:asvrcd_ccLdef} as small as possible. Using $\Popt = \mI$, it is easy to see that $\ccL = \lambda_{\max} \left( \diag(p)^{-\frac12}\mM  \diag(p)^{-\frac12}\right)$. While the optimal $p$ is still hard to find, we set $p_i\propto \mM_{i,i}$ (i.e., the effect of importance sampling is the same as the effect of Jacobi preconditioner). Figure~\ref{fig:asvrcd_identity} shows the result. As expected, accelerated {\tt SVRCD} always outperforms non-accelerated variant, while at the same time, the importance sampling improves the performance too.

\begin{figure}[!h]
\centering
\begin{minipage}{0.33\textwidth}
  \centering
\includegraphics[width =  \textwidth ]{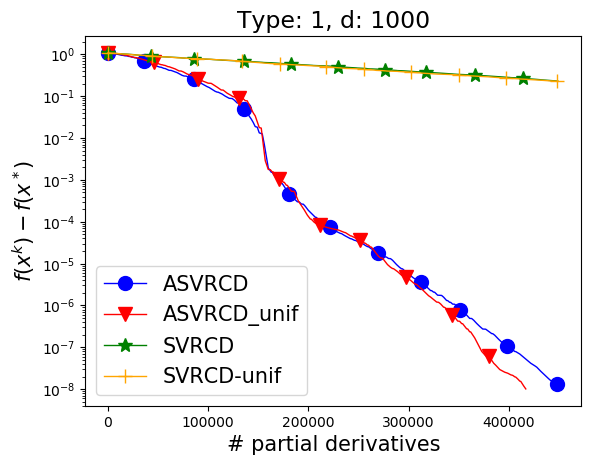}
\end{minipage}%
\begin{minipage}{0.33\textwidth}
  \centering
\includegraphics[width =  \textwidth ]{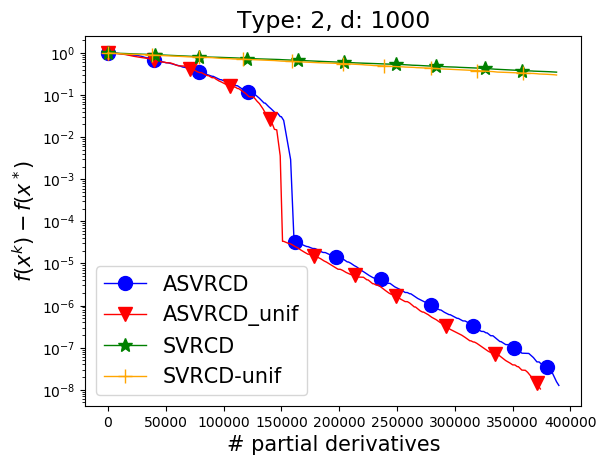}
\end{minipage}%
\\
\begin{minipage}{0.33\textwidth}
  \centering
\includegraphics[width =  \textwidth ]{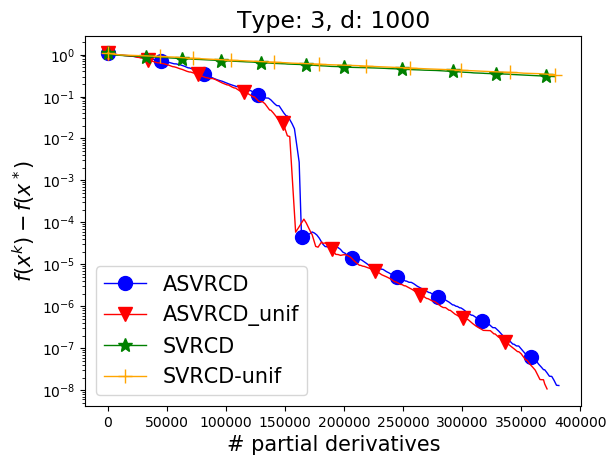}
\end{minipage}%
\begin{minipage}{0.33\textwidth}
  \centering
\includegraphics[width =  \textwidth ]{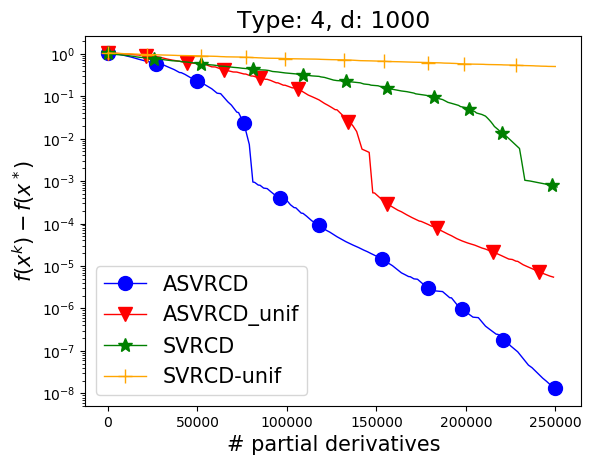}
\end{minipage}%
\caption{Comparison of both {\tt ASVRCD} and {\tt SVRCD} with importance and uniform sampling.} 
\label{fig:asvrcd_identity}
\end{figure}

\subsection{The effect of $\mW$}

The second experiment compares the performance of both {\tt ASVRCD} and {\tt SVRCD} for various $\mW$. We only consider methods with importance sampling ($p_i\propto \mM_{i,i} \mW_{i,i}$) and theory supported stepsize. Figure~\ref{fig:asvrcd_variable} presents the result.  We see that the smaller $\Range{\mW}$ is, the faster the convergence is. This observation is well-aligned with our theory: $\ccL$ is increasing as a function of $\mW$ (in terms of Loewner ordering).

\begin{figure}[!h]
\centering
\begin{minipage}{0.33\textwidth}
  \centering
\includegraphics[width =  \textwidth ]{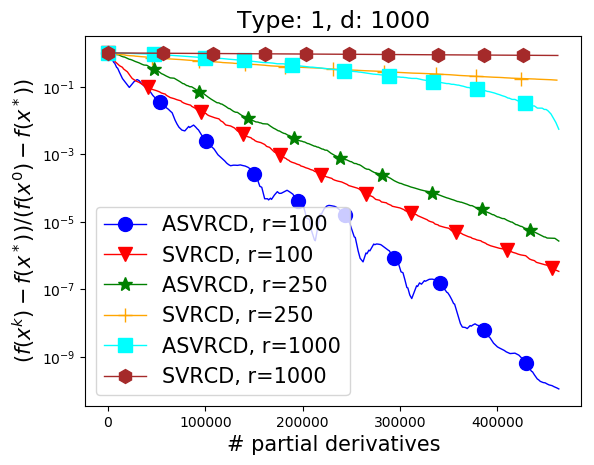}
\end{minipage}%
\begin{minipage}{0.33\textwidth}
  \centering
\includegraphics[width =  \textwidth ]{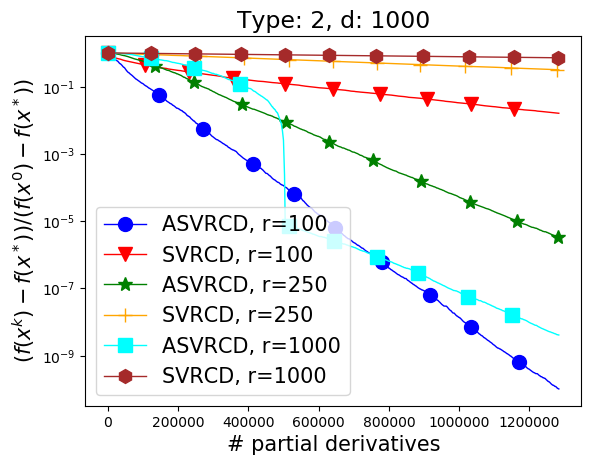}
\end{minipage}%
\\
\begin{minipage}{0.33\textwidth}
  \centering
\includegraphics[width =  \textwidth ]{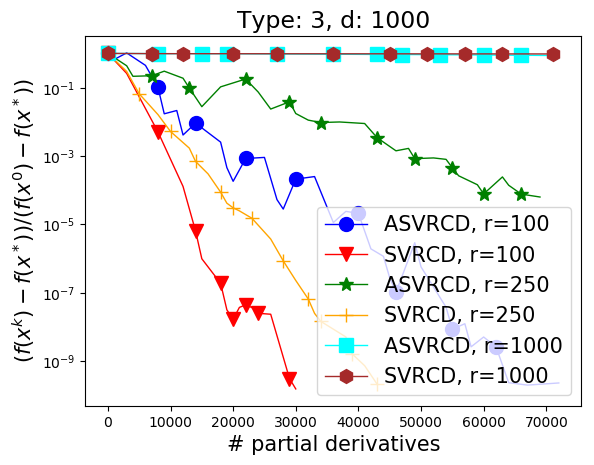}
\end{minipage}%
\begin{minipage}{0.33\textwidth}
  \centering
\includegraphics[width =  \textwidth ]{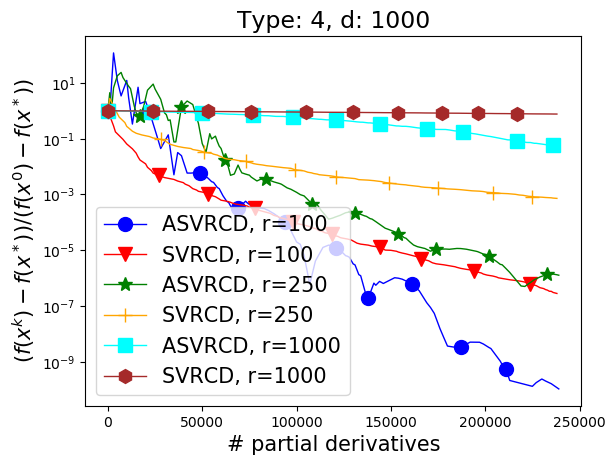}
\end{minipage}%
\caption{Comparison of {\tt ASVRCD} and {\tt SVRCD} for various $\Popt$. Label 'r' indicates the dimension of $\Range{\Popt}$.} 
\label{fig:asvrcd_variable}
\end{figure}

\section{Conclusion}
In this chapter we have introduced {\tt ASVRCD} -- an accelerated {\tt SVRCD} algorithm. Besides that, we have shown that {\tt SAGA}/{tt L-Katyusha} are a special case of {\tt SEGA}/{\tt ASVRCD}, while their convergence guarantees can be recovered. This rationale can be further generalized: it is possible to show that essentially any finite-sum stochastic algorithm is a special case of analogous method with partial derivative oracle (those are yet to be discovered/analyzed) in a given setting (i.e.,  strongly convex, convex, non-convex). Those include, but are not limited to {\tt SGD}~\cite{robbins, nemirovski2009robust}, over-parametrized {\tt SGD}~\cite{vaswani2019-overparam}, {\tt SAG}~\cite{sag}, {\tt SVRG}~\cite{svrg}, {\tt S2GD}~\cite{konevcny2013semi}, {\tt SARAH}~\cite{nguyen2017sarah}, incremental methods such as Finito~\cite{defazio2014finito}, {\tt MISO}~\cite{mairal2015incremental} or accelerated algorithms such as point-{\tt SAGA}~\cite{defazio2016simple}, Katyusha~\cite{allen2017katyusha}, {\tt MiG}~\cite{zhou2018simple}, {\tt SAGA-SSNM}~\cite{zhou2018direct}, Catalyst~\cite{lin2015universal, kulunchakov2019generic}, non-convex variance reduced algorithms~\cite{reddi2016stochastic, allen2016variance, fang2018spider} and others. In particular, {\tt SGD} can be seen as a special case of block coordinate descent, while {\tt SAG} is a special case of bias-{\tt SEGA} from~\cite{sega} (neither of {\tt CD} with non-separable prox, nor bias-{\tt SEGA} were analyzed yet).

\chapter{Federated Learning of a Mixture of Global and Local Models}

\label{local}

\graphicspath{{local/plots/}}

With the proliferation of mobile phones, wearable devices, tablets, and smart home devices comes an increase in the volume of data captured and stored on them. This data contains a wealth of potentially useful information to the owners of these devices, and more so if appropriate machine learning models could be trained on the heterogeneous data stored across the network of such devices.  The traditional approach involves moving the relevant data to a  data center where centralized machine learning techniques can be efficiently applied~\cite{Dean2012, AIDE}.  However, this approach is not without issues.  First, many device users are increasingly sensitive to privacy concerns and prefer their data to never leave their devices. Second, moving data from their place of origin to a centralized location is very inefficient in terms of energy and time. 

\section{Federated learning}

{\em Federated learning (\acrshort{FL})} \cite{FedAvg2016, FEDLEARN,   FEDOPT, FL2017-AISTATS} has emerged as an interdisciplinary  field focused on addressing these issues by training machine learning models directly on  edge devices.  The currently prevalent paradigm~\cite{FL_survey_2019, FL-big} casts supervised FL as an empirical risk minimization problem of the form
\begin{equation}\label{eq:local_ERM-traditional}  
\min \limits_{x\in \R^d} \frac{1}{n}\sum \limits_{i=1}^n f_i(x),
\end{equation}
where $n$ is the number of devices participating in training, $x\in \R^d$ encodes the $d$ parameters of  a {\em global model}  (e.g., weights of a neural network) and $f_i(x)\eqdef \EE_{\xi\sim \cD_i}{f(x,\xi)}$ represents the aggregate loss of model $x$ on the local data represented by distribution $\cD_i$ stored on device $i$. One of the defining characteristics of FL is that the data distributions $\cD_i$ may possess very different properties across the devices. Hence, any potential FL method  is explicitly required to be able to work under the {\em heterogeneous} data setting.

The most popular method for solving \eqref{eq:local_ERM-traditional} in the context of FL is the FedAvg algorithm~\cite{FedAvg2016}. In its most simple form, when one does not employ partial participation, model compression, or stochastic approximation, FedAvg reduces to Local Gradient Descent ({\tt LGD})~\cite{localGD, localSGD-AISTATS2020}, which is an extension of {\tt GD} performing more than a single gradient step on each device before aggregation.  FedAvg has been shown to work well empirically, particularly for non-convex problems, but comes without convergence guarantees and can diverge in practical settings when data are heterogeneous.

\subsection{Some issues with current approaches to FL}

%

The first motivation for our research comes from the appreciation that data heterogeneity does not merely present challenges to the design of new provably efficient training  methods for solving  \eqref{eq:local_ERM-traditional}, but {\em also inevitably raises questions about the utility of such a global solution to individual users.} Indeed, a global model trained across all the data from all devices might be so removed from the typical data and usage patterns experienced by an individual user as to render it virtually useless. This issue has been observed before, and various approaches have been proposed to address it. For instance, the MOCHA \cite{NIPS2017_FL-multitask} framework uses a multi-task learning approach to allow for personalization. A generic online algorithm for gradient-based parameter-transfer meta-learning \cite{Khodak2019FLmeta} was demonstrated to improve practical performance over FedAvg~\cite{FL2017-AISTATS}. Approaches based on variational inference~\cite{FL2019variational}, cyclic patterns in practical FL data sampling  \cite{Eichner2019semi-cyclicSGD_for_FL} and transfer learning \cite{Zhao2018FL-transfer-learn} have been proposed.

The second motivation for our work is the realization that even very simple variants of FedAvg, such as {\tt LGD}, which should be easier to analyze, fail to provide theoretical improvements in communication complexity over their non-local cousins, in this case, {\tt GD} \cite{localGD, localSGD-AISTATS2020}.  This observation is at odds with the practical success of local methods in FL. This leads us to ask the question: 
\begin{quote}{\em If {\tt LGD} does not theoretically improve upon {\tt GD} as a solver for the traditional global problem~\eqref{eq:local_ERM-traditional}, perhaps  {\tt LGD} should not be seen as a method for solving~\eqref{eq:local_ERM-traditional} at all. In such a case, what problem does {\tt LGD} solve?}
\end{quote}

 A good answer to this question would shed light on the workings of {\tt LGD}, and by analogy, on the role local steps play in more elaborate  FL methods such as local {\tt SGD}~\cite{localSGD-Stich, localSGD-AISTATS2020} and FedAvg.

\section{Contributions \label{sec:local_contrib}}

In our work we argue that the two motivations mentioned in the introduction  point in the same direction, i.e., we show that {\em a single solution can be devised addressing both problems at the same time.} 

Our main contributions are:

\begin{itemize}
\item  {\bf New formulation of FL which seeks a mixture of  global and local models.} We propose a {\em new optimization formulation of FL.}  Instead of learning a single global model by solving \eqref{eq:local_ERM-traditional}, we propose to learn a mixture of the global model and the purely local models which can be trained by each device $i$ on its own, using its data $\cD_i$ only.  Our formulation (see \eqref{eq:local_main}  in Section~\ref{sec:local_new_formulation}) lifts the problem from $\R^d$ to $\R^{nd}$, allowing each device $i$ to learn a personalized model $x_i\in \R^d$. However, these personalized models are explicitly encouraged to not depart too much from their mean by the inclusion of a quadratic penalty $\Phi$ multiplied by a penalty parameter $\lambda\geq 0$.\footnote{The idea of softly-enforced similarity of the local models was already introduced in the domain of decentralized optimization~\cite{lan2018communication, gorbunov2019optimal}. However, their motivation is vastly different to ours (besides not considering FL or local algorithms) -- the mentioned methods still aim to find the global model by having the penalty parameter inversely proportional to the target accuracy $\varepsilon$.}

\item  {\bf Theoretical properties of the new formulation.} We study the properties of the optimal solution of our formulation, thus developing an algorithmic-free theory. When the penalty parameter is set to zero, then obviously, each device is allowed to train their own model without any dependence on the data stored on other devices. Such purely local models are rarely useful. We prove that the optimal local models converge to the traditional global model characterized by \eqref{eq:local_ERM-traditional} at the rate $\cO(1/\lambda)$. We also show that the total loss evaluated at the local models is always not higher than the total loss evaluated at the global model (see Theorem~\ref{thm:local_penalty}).  Moreover, we prove an insightful structural result for the optimal local models: the optimal model learned by device $i$ arises by subtracting the gradient of the loss function stored on that device evaluated at the same point (i.e., a local model) from the average of the optimal local models (see Theorem~\ref{thm:local_characterization}).  As a byproduct, this theoretical result sheds new light on the key update step in the model agnostic meta-learning (MAML) method \cite{MAML2017}, which has a similar but subtly different structure. The subtle difference is that the MAML update obtains the local model by subtracting the gradient evaluated at the {\em global} model. While MAML is a heuristic, we provide rigorous theoretical guarantees.

\item  {\bf Loopless {\tt LGD}: non-uniform {\tt SGD} applied to our formulation.} We then propose a randomized gradient-based method---{\em Loopless Local Gradient Descent ({\acrshort{L2GD}})}---for solving our new formulation (Algorithm~\ref{alg:local_L2GD}). This method is, in fact, a non-standard application of {\tt SGD} to our problem, and can be seen as an instance of {\tt SGD} with non-uniform sampling applied to the problem of minimizing the sum of two convex functions~\cite{iprox-sdca, pmlr-v97-qian19b}: the average loss, and the penalty. When the loss function is selected by the randomness in our {\tt SGD} method, the resultant stochastic gradient step can be {\em interpreted} as the execution of a single local {\tt GD} step on each device. Since we set the probability of the loss being sampled to be high, this step is typically repeated multiple times, and this has the effect of taking multiple local {\tt GD} steps. In contrast to standard {\tt LGD}, the number of local steps is not fixed, but random, and follows a geometric distribution. This mechanism is similar in spirit to how the recently proposed loopless variants of {\tt SVRG} \cite{hofmann2015variance, kovalev2019don} work in comparison with the original {\tt SVRG}~\cite{svrg, proxsvrg}. Once the penalty is sampled by our method, the resultant {\tt SGD} step can be interpreted as the execution of an aggregation step. In contrast with standard aggregation, which performs full averaging of the local models, our method is more sophisticated and merely takes a {\em step towards averaging}. However, the step is relatively large. This suggests that perhaps full averaging in modern FL methods such as FedAvg or {\tt LGD} and {\acrshort{LSGD}} is too aggressive, and should be re-examined.

\item {\bf Convergence theory.} By adapting the general theory from \cite{pmlr-v97-qian19b} to our setting, we obtain theoretical convergence guarantees assuming that each $f_i$ is $L$-smooth and $\mu$-strongly convex (see Theorem~\ref{thm:local_local_gd}). Interestingly, by optimizing  the sampling probability (we get $p^* = \frac{\lambda}{\lambda + L}$) which is an  indirect way of fixing the {\em expected number of local steps} to $1+\frac{L}{\lambda}$, we prove the communication complexity result (i.e., bound on the expected number of communication rounds; see Corollary~\ref{cor:local_optimalp})  \[\frac{2 \lambda }{\lambda + L} \frac{L}{\mu}\log \frac{1}{\varepsilon}.\] We believe that this is remarkable in several ways.  By choosing $\lambda$ small, we tilt our goal towards pure local models, and the number of communication rounds is very small, tending to 0 as $\lambda\to 0$. If $\lambda\to \infty$, our the solution to our formulation converges to the optimal global model,  and {\tt L2GD} obtains the communication bound  $\cO\left(\frac{L}{\mu} \log \frac{1}{\varepsilon}\right)$, which matches the efficiency of {\tt GD}.   Our results can be extended to convex and non-convex regimes, but we do not explore such generalizations here.

\item  {\bf Generalizations: partial participation, local {\tt SGD} and variance reduction.} We further generalize and improve our method and convergence results by allowing for 
\begin{itemize}
\item[(i)] stochastic {\em partial participation} of devices in each communication round,
\item[(ii)]  {\em subsampling} on each device which means we can perform local {\tt SGD} steps instead of local {\tt GD} steps, and  
\item[(iii)]  {\em total variance reduction mechanism} to tackle the variance  coming from three sources: locality of the updates induced by non-uniform sampling (already present in {\tt L2GD}), partial participation and subsampling from local data.
\end{itemize}

Due to its level of generality, this method, which we call {\tt L2SGD++}, is presented in the Appendix only, alongside the associated complexity results.  In the main body of this chapter, we instead present a simplified version thereof, one that does not include partial participation. We call this method {\tt L2SGD+} (Algorithm~\ref{alg:local_L2SGD}). The convergence theory for it is presented in 
Theorem~\ref{thm:local_saga_simple} and Corollary~\ref{cor:local_lsd_optimal_p}.

\item  {\bf Allowing for heterogeneous data.} All our methods and convergence results allow for fully heterogeneous data and do not depend on any assumptions on data similarity across the devices.

\item  {\bf Superior empirical performance.}
We show through ample numerical experiments that our theoretical predictions can be observed in practice.

\end{itemize}

\section{New formulation of FL} \label{sec:local_new_formulation}



 We now introduce our new formulation for training supervised FL models:
\begin{equation}\label{eq:local_main}
\begin{aligned}
\min_{x_1,\dots,x_n\in \R^d} &\left\{ F(x) \eqdef f(x) + \lambda \Phi(x) \right\}\\
  f(x) \eqdef  \frac{1}{n}\sum \limits_{i=1}^n  f_i(&x_i), \quad   \Phi(x) \eqdef \frac{1}{2 n}\sum \limits_{i=1}^n \norm{x_i-\bar{x}}^2,
\end{aligned}
\end{equation}
where $\lambda \geq 0 $ is a  penalty parameter,  $x_1,\dots,x_n \in \R^d$ are local models, $x\eqdef (x_1,x_2,\dots,x_n) \in \R^{nd}$ and  $\bar{x}\eqdef \frac{1}{n}\sum_{i=1}^n x_i$ is the average of the local models. 

Due to the assumptions on $f_i$ we will make in Section~\ref{bd798gsd}, $F$ is strongly convex and hence \eqref{eq:local_main} has a unique solution, which we denote $$x(\lambda) \eqdef (x_1(\lambda),\dots,x_n(\lambda)) \in \R^{nd}.$$ We further let $$\bar{x}(\lambda) \eqdef \frac{1}{n}\sum_{i=1}^n x_i(\lambda).$$ 
We now comment on the rationale behind the new formulation.

\paragraph{\bf Local models ($\lambda=0$).} Note that  for each $i$, $x_i(0)$ solves the {\em local problem}
\[ \min_{x_i \in \R^d} f_i(x_i).\]
That is,  $x_i(0)$ is the local model based on data $\cD_i$ stored on device $i$ only. This model can be computed by device $i$ without any communication whatsoever. Typically, $\cD_i$ is not rich enough for this local model to be useful. In order to learn a better model, one has to take into account the date from other clients as well. This, however, requires communication.

\paragraph{\bf Mixed models ($\lambda \in (0,\infty)$).} As $\lambda$ increases, the penalty $\lambda \Phi(x)$ has an increasingly more substantial effect, and communication is needed to ensure that the models are not too dissimilar, as otherwise $\Phi$ would be too large. 

\paragraph{\bf Global model ($\lambda = \infty$).} Let us now look at the limit case $\lambda \to \infty$. Intuitively, this limit case should force the optimal local models to be mutually identical, while minimizing the loss $f$. In particular, this limit case will solve\footnote{If $\lambda=\infty$ and $x_1=x_2=\dots=x_n$ does not hold, we have $F(x)=\infty$. Therefore, we can restrict ourselves on set $x_1=x_2=\dots=x_n$ without loss of generality. }
\[ \min_{x_1,\dots,x_n\in \R^{d}} \left\{ f(x)  \;:\; x_1=x_2=\dots=x_n \right\},\]
which is equivalent to the global formulation \eqref{eq:local_main}.  Because of this, let us defined $x_i(\infty)$ for each $i$ to be the optimal global solution of \eqref{eq:local_ERM-traditional}, and let $x(\infty) \eqdef (x_1(\infty),\dots,x_n(\infty))$. 

\subsection{Technical preliminaries}\label{bd798gsd}
Similarly to the rest of the thesis, we make the following assumption on the functions $f_i$:
\begin{assumption}\label{as:local_smooth_sc_main}
For each $i$, the function $f_i:\R^d\to \R$ is $L$-smooth and $\mu$-strongly convex.
\end{assumption}

Note that the separable structure of $f$ implies that $\left(\nabla f(x)\right)_i = \frac{1}{n}\nabla f_i (x_i)$, i.e., 
\begin{equation}\label{eq:local_nabla f}\nabla f(x) = \frac{1}{n} ( \nabla f_1 (x_1), \nabla f_2 (x_2), \dots,   \nabla f_n (x_n)).\end{equation}   Hence,
 the norm of $\nabla f(x)\in \R^{nd}$ decomposes as
$ \norm{\nabla f(x)}^2 = \frac1n \sum \limits_{i=1}^n \norm{\nabla f_i (x_i)}^2.$

Note that Assumption~\ref{as:local_smooth_sc_main} implies that $f$ is $L_f$-smooth with $L_f \eqdef \frac{L}{n}$ and $\mu_f$-strongly convex with $\mu_f \eqdef \frac{\mu}{n}$. Clearly, $\Phi$ is convex by construction. It can be shown that $\Phi$ is $L_{\Phi}$-smooth with $L_{\Phi}=\frac1n$ (see Appendix). 
 We can also easily see that 
\begin{equation} \label{eq:local_nhid8gf8f}
\left(\nabla \Phi(x)\right)_i = \frac1n(x_i - \bar{x})
\end{equation} 
(see Appendix), which implies
$$  \Phi(x)   \overset{\eqref{eq:local_main}+\eqref{eq:local_nhid8gf8f}}{=} \frac{n}{2}\sum \limits_{i=1}^n \norm{ (\nabla \Phi(x))_i}^2 = \frac{n}{2} \norm{\nabla \Phi(x)}^2.$$

\subsection{Characterization of optimal solutions}

Our first result describes the behavior of $f(x(\lambda))$ and $\Phi(x(\lambda))$ as a function of $\lambda$.

\begin{theorem} \label{thm:local_penalty}
The function $\lambda \to \Phi(x(\lambda))$ is non-increasing, and
for all $\lambda>0$ we have
\begin{equation} \label{eq:local_dissimilarity} \Phi(x(\lambda)) \leq \frac{f(x(\infty)) - f(x(0))}{\lambda}. \end{equation}
Moreover, the function $\lambda \to f(x(\lambda))$ is non-decreasing, and for all $\lambda \geq 0$ we have
\begin{equation} \label{eq:local_bifg9dd8} f(x(\lambda)) \leq f(x(\infty)).\end{equation}
\end{theorem}

Inequality~\eqref{eq:local_dissimilarity} says that the penalty decreases to zero as $\lambda$ grows, and hence the optimal local models $x_i(\lambda)$ are increasingly similar as $\lambda$ grows. The second statement suggest that the loss $f(x(\lambda))$ increases with $\lambda$, but never exceeds the optimal global loss $f(x(\infty))$ of the standard FL formulation \eqref{eq:local_ERM-traditional}.

We now characterize the optimal local models which connect our model to the MAML framework~\cite{MAML2017}, as mentioned in the introduction.
\begin{theorem} \label{thm:local_characterization}
For each $\lambda >0$ and $1\leq i \leq n$ we have
\begin{equation}\label{eq:local_step_from_mean}
x_i(\lambda) =  \overline{x}(\lambda) - \frac{1}{\lambda }\nabla f_i(x_i(\lambda)).
\end{equation}
Further, we have $ \sum \limits_{i=1}^n \nabla f_i(x_i(\lambda)) = 0$
and $ \Phi(x(\lambda)) = \frac{1}{2 \lambda^2  } \norm{\nabla f(x(\lambda))}^2$.
\end{theorem}
The optimal local models $\eqref{eq:local_step_from_mean}$ are obtained from the average model by subtracting a multiple of the local gradient. Moreover,  observe that the local gradients always sum up to zero at optimality. This is obviously true for $\lambda=\infty$, but it is a bit less obvious that this holds for any $\lambda>0$.

\section{The {\tt L2GD} algorithm} \label{sec:local_L2GD}

In this section we  describe a new randomized gradient-type method for solving our new formulation \eqref{eq:local_main}.  Our method is a non-uniform {\tt SGD} for \eqref{eq:local_main} seen as a 2-sum problem, sampling either  $\nabla f$ or $\nabla \Phi$ to estimate $\nabla F$. Letting $0<p<1$,  we define a stochastic gradient of $F$ at  $x \in \R^{nd}$ as follows
\begin{equation}\label{eq:local_g(x)}g(x) \eqdef \begin{cases} \frac{\nabla f(x)}{1-p} & \text{with probability} \quad 1-p\\
 \frac{\lambda \nabla \Phi(x)}{p} & \text{with probability} \quad p
\end{cases}.\end{equation}

Since $$\E{g(x)} =(1-p) \frac{\nabla f(x)}{1-p} +  p \frac{\lambda \nabla \Phi(x)}{p} = \nabla F(x),$$ the vector $g(x)$ is an unbiased estimator of $\nabla F(x)$. This leads to the following method for minimizing $F$, which we call {\tt L2GD}:
\begin{equation} \label{eq:local_SGD}x^{k+1} = x^k -\alpha G(x^k).\end{equation}
Plugging formulas \eqref{eq:local_nabla f} and \eqref{eq:local_nhid8gf8f} for $\nabla f(x)$ and $\nabla \Phi(x)$ into \eqref{eq:local_g(x)} and subsequently into \eqref{eq:local_SGD}, and writing the resulting method in a distributed manner, we arrive at Algorithm~\ref{alg:local_L2GD}. In each iteration, a coin $\xi$ is tossed and lands $1$  with probability $p$ and $0$ with probability $1-p$. If $\xi=0$, all {\color{red}Devices} perform one local {\tt GD} step \eqref{eq:local_damkpdaosdaon}, and if $\xi=1$, {\color{blue}Master} shifts each local model towards the average via \eqref{eq:local_buiufg98fb}.  As we shall see in Section~\ref{sec:local_analysis}, our convergence theory  limits the value of the stepsize $\alpha$, which has the effect that the ratio $\frac{\alpha \lambda}{n p}$ cannot exceed $\frac{1}{2}$. Hence, \eqref{eq:local_buiufg98fb} is a convex combination of $x_i^k$ and $\bar{x}^k$, which justifies the statement we have made above: $x_i^{k+1}$ shifts towards $\bar{x}^k$ along the line joining these two points.

\begin{algorithm}[h]
  \caption{{\tt L2GD}: Looples Local Gradient Descent}
  \label{alg:local_L2GD}
\begin{algorithmic}
\State{\bfseries Input: }{$x^0_1 = \dots = x_{n}^0\in\R^{d}$, stepsize $\alpha$, probability $p$ }
  \For{$k=0,1,2,\dotsc$}
  \State $\xi = 1$ with probability $p$ and $0$ with probability $1-p$ 
  \If {$\xi=0$}
     \State All {\color{red}Devices} $i=1,\dots,n$ perform a local {\tt GD} step:

        \begin{equation} \label{eq:local_damkpdaosdaon}
        x^{k+1}_i  =x^k_i - \frac{\alpha}{n(1-p)} \nabla f_i(x_i^k)
\end{equation}   
  \Else
  \State {\color{blue}Master} computes the average $\bar{x}^k = \frac{1}{n}\sum_{i=1}^n x_i^k$
  \State {\color{blue}Master}  for each $i$ computes step towards aggregation
\State 
 \begin{flalign}
 		x^{k+1}_i  = \left(1-\frac{\alpha \lambda}{np}\right) x^k_i + \frac{\alpha \lambda}{np} \bar{x}^k  \label{eq:local_buiufg98fb}
 		\end{flalign}     
  \EndIf
  \EndFor
\end{algorithmic}
\end{algorithm}

\subsection{Understanding communication}

\begin{example} In order to better understand when communication takes place  in Algorithm~\ref{alg:local_L2GD}, consider the following possible sequence of coin tosses: 
$0, 0, 1, 0, 1, 1, 1, 0.$
The first two coin tosses lead to two local {\acrshort{GD}} steps \eqref{eq:local_damkpdaosdaon} on all devices. The third coin toss lands $1$, at which point all local models $x_i^k$ are communicated to the master, averaged to form $\bar{x}^k$, and the step \eqref{eq:local_buiufg98fb} towards averaging is taken. The fourth coin toss is $0$, and at this point, the master communicates the updated local models back to the devices, which subsequently perform a single local {\tt GD} step \eqref{eq:local_damkpdaosdaon}. Then come three consecutive coin tosses landing $1$, which means that the local models are again communicated to the master, which performs three averaging steps  \eqref{eq:local_buiufg98fb}. Finally, the eight coin toss lands $0$, which makes the master send the updated local models back to the devices, which subsequently perform a single local {\tt GD} step.
\end{example}

This example illustrates that communication needs to take place whenever two consecutive coin tosses land a different value. If $0$ is followed by a $1$, all devices communicate to the master, and if $1$ is followed by a $0$, the master communicates back to the devices. It is standard to count each pair of communications, {\color{red}Device}$\to${\color{blue}Master} and the subsequent {\color{blue}Master}$\to${\color{red}Device}, as a single communication round. 

\begin{lemma} \label{lem:local_exp_no_comm}The expected number of communication rounds in $k$ iterations of {\tt L2GD} is $p(1-p)k$.
\end{lemma}

\subsection{The dynamics of local {\tt GD} and averaging steps}

Further, notice that  {\em the average of the local models does not change} during an aggregation step. Indeed, $\bar{x}^{k+1}$ is equal to
\[  \frac{1}{n} \sum \limits_{i=1}^n x_i^{k+1}   \overset{\eqref{eq:local_buiufg98fb}}{=}   \frac{1}{n} \sum \limits_{i=1}^n  \left[\left(1-\frac{\alpha \lambda}{np}\right) x^k_i + \frac{\alpha \lambda}{np} \bar{x}^k \right] = \bar{x}^k.\]

If several averaging steps take place in a sequence, the point $a=\bar{x}^k$ in \eqref{eq:local_buiufg98fb} remains unchanged, and each local model $x_i^k$ merely moves along the line joining the initial value of the local model at the start of the sequence and  $a$, with each step pushing $x_i^k$ closer to the average $a$. 

{\em In summary, the more local {\tt GD} steps are taken, the closer the local models get to the pure local models, and the more averaging steps are taken, the closer the local models get to their average value. The relative number of local {\tt GD} vs. averaging steps is controlled by the parameter $p$: the expected number of local {\tt GD} steps is $\frac{1}{p}$, and the expected number of consecutive aggregation steps is $\frac{1}{1-p}$.}


\subsection{Convergence theory }\label{sec:local_analysis}

We fist show that our gradient estimator $g(x)$ satisfies the expected smoothness property \cite{jacsketch, pmlr-v97-qian19b}.

\begin{lemma} \label{lem:local_exp_smoothnes-basic}
Let $\cL\eqdef \frac1n \max\left\{\frac{L}{1-p}, \frac{\lambda }{p}\right\}$ and \[\sigma^2 \eqdef  \frac{1}{n^2}  \sum_{i=1}^n \left( \frac{1}{1-p} \| \nabla f_i(x_i(\lambda))\|^2 + \frac{\lambda^2}{p} \|x_i(\lambda) - \overline{x}(\lambda) \|^2 \right).\] Then for all $x\in \R^d$ we have the inequalities
$$\E{ \norm{g(x) - G(x(\lambda))}^2}  \leq 2\cL \left(F(x)-F(x(\lambda))\right)$$ and
$$\E{\norm{g(x)}^2}  \leq 4 \cL  (F(x)-F(x(\lambda))) + 2 \sigma^2.$$
\end{lemma}

We now present our convergence result for {\tt L2GD}.

\begin{theorem}\label{thm:local_local_gd} Let Assumption~\ref{as:local_smooth_sc_main} hold.
If $\alpha \leq \frac{1}{2\cL}$, then
\[\E{\norm{x^k-x(\lambda)}^2} \leq \left(1- \frac{\alpha \mu}{n}\right)^k \norm{x^0-x(\lambda)}^2 + \frac{2n\alpha \sigma^2}{\mu}.\]
If we choose $\alpha = \frac{1}{2\cL}$, then 
$\frac{\alpha \mu}{n} =  \frac{\mu}{2 \max\left\{\frac{L}{1-p}, \frac{\lambda }{p}\right\} }$
and
\[\frac{2n\alpha \sigma^2}{\mu} = \frac{ \sum\limits_{i=1}^n \left( \frac{1}{1-p} \| \nabla f_i(x_i(\lambda))\|^2 +  \frac{\lambda^2}{p} \|x_i(\lambda) - \overline{x}(\lambda) \|^2 \right)}{ \max\left\{\frac{L}{1-p}, \frac{\lambda }{p}\right\} \mu}.\]
\end{theorem}

\begin{remark}[Full averaging  not supported]
Is a setup such that conditions of Theorem~\ref{thm:local_local_gd} are satisfied and the aggregation update~\eqref{eq:local_buiufg98fb} is identical to full averaging? This is equivalent requiring $0<p<1$ such that $\alpha \lambda =np$. However, we have $\alpha \lambda   \leq \frac{\lambda  }{2\cL}  \leq np$, which means that full averaging is not supported by our theory.
\end{remark}

\subsection{Optimizing the rate and communication} \label{sec:local_rate}

Let us find the parameters $p$ and $\alpha$ which lead to the fastest rate, in terms of either iterations or communication rounds, to push the error within $\varepsilon$ of the neighborhood\footnote{In Section~\ref{sec:local_lsgd} we propose a variance reduced algorithm which is able to get rid of the neighborhood in the convergence result completely. In that setting, our goal will be to achieve $\E{\norm{x^k-x(\lambda)}^2} \leq \varepsilon \norm{x^0-x(\lambda)}^2$.  } from Theorem~\ref{thm:local_local_gd}, i.e., to achieve
\begin{equation}\label{eq:local_bi7g97ibudd} \E{\norm{x^k-x(\lambda)}^2} \leq \varepsilon \norm{x^0-x(\lambda)}^2 + \frac{2n\alpha \sigma^2}{\mu}. \end{equation}

\begin{corollary}\label{cor:local_optimalp}
The value $p^{*} = \frac{\lambda}{L + \lambda}$ minimizes both the number of iterations and the expected number of communications for achieving \eqref{eq:local_bi7g97ibudd}. In particular, the optimal number of iterations is $2\frac{L+\lambda}{\mu} \log \frac{1}{\varepsilon}$, and the optimal expected number of communications is $\frac{2\lambda}{\lambda+ L} \frac{L}{\mu} \log \frac{1}{\varepsilon}$.
\end{corollary}

If we choose $p=p^*$, then $\frac{\alpha \lambda }{np} = \frac{1}{2}$, and the aggregation  rule~\eqref{eq:local_buiufg98fb} in Algorithm~\ref{alg:local_L2GD} becomes
\begin{equation}\label{eq:local_bui87g9f}
x^{k+1}_i = \frac12 \left( x^k_i +\bar{x}^k\right) 
\end{equation}
while the local {\tt GD} step~\eqref{eq:local_damkpdaosdaon} becomes
$
x^{k+1}_i  =x^k_i - \frac{1}{2L} \nabla f_i(x_i^k).
$
Notice that while our method does not support full averaging as that is too unstable,  \eqref{eq:local_bui87g9f} suggests that one should take a large step {\em towards} averaging.  

\begin{figure}[t]
\centering
\begin{minipage}{0.5\textwidth}
  \centering
\includegraphics[width =  \textwidth ]{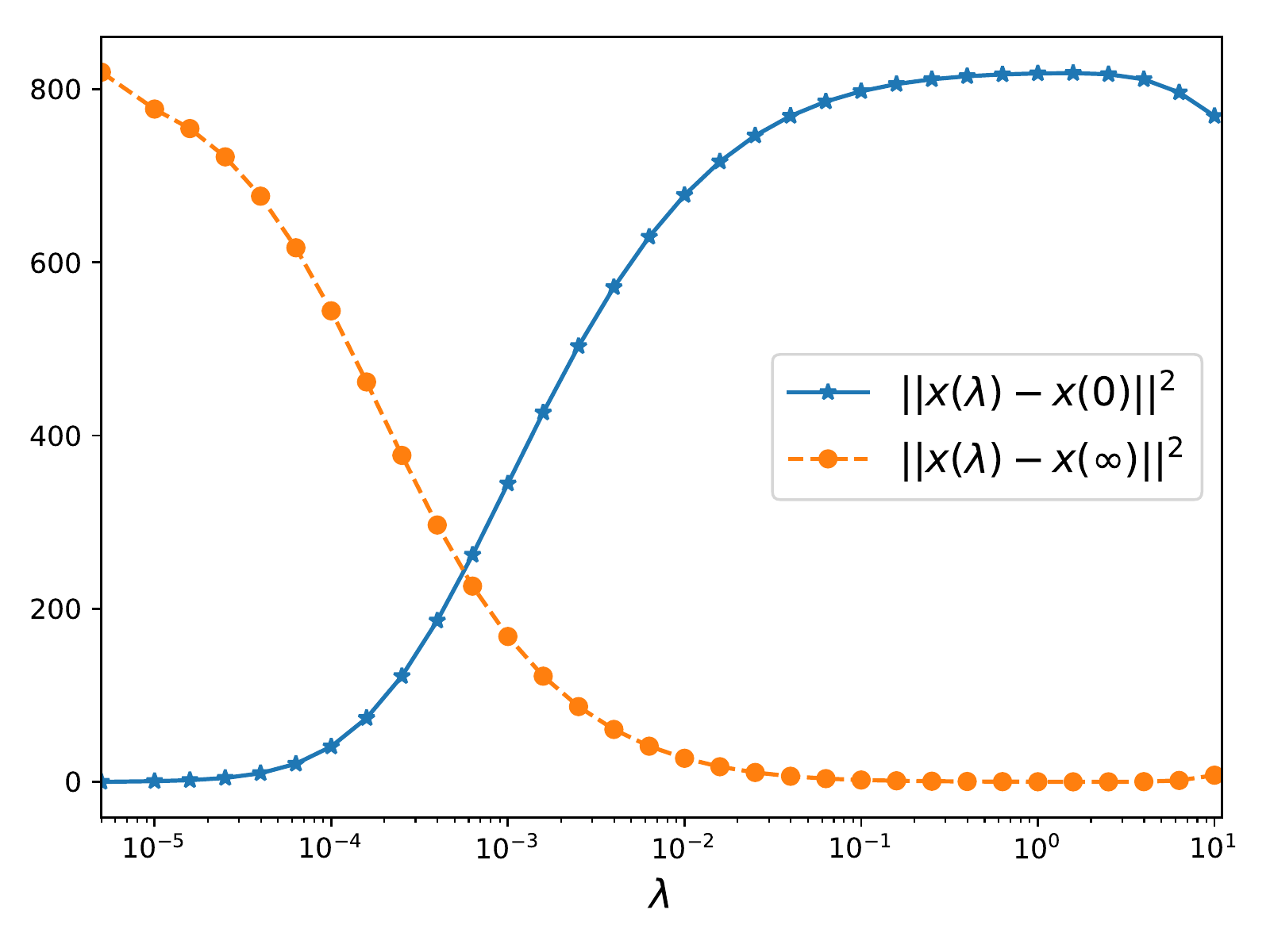}
\end{minipage}%
\vskip -0.3cm
\caption{Distance of solution $x(\lambda)$ of~\eqref{eq:local_main} to pure local solution $x(0) $ and global solution $x(\infty)$ as a function of $\lambda$. Logistic regression on a1a dataset. See Appendix for experimental setup.} 
\label{fig:local_lambda}
\end{figure}

As $\lambda $ get smaller, the solution to the optimization problem \eqref{eq:local_main} will increasingly favour pure local models, i.e., $x_i(\lambda) \to x_i(0) \eqdef \arg\min f_i$ for all $i$ as $\lambda \to 0$. Pure local models can be computed without any communication whatsoever and Corollary~\ref{cor:local_optimalp} confirms this intuition: the optimal number of communication round decreases to zero as $\lambda \to 0$. On the other hand, as $\lambda \to \infty$, the optimal number of communication rounds converges to $2\frac{L}{\mu}\log \frac{1}{\varepsilon}$, which recovers the performance of {\tt GD} for finding the globally optimal model (see Figure~\ref{fig:local_lambda}).

{\em In summary, we recover the communication efficiency of {\tt GD} for finding the globally optimal model as $\lambda \to \infty$. However, for other values of $\lambda$, the communication complexity of {\tt L2GD} is better and decreases to $0$ as $\lambda \to 0$. Hence, our communication complexity result interpolates between the communication complexity of {\tt GD} for finding the global model and the zero communication complexity for finding the pure local models. }

\section{The {\tt L2SGD+} algorithm} \label{sec:local_lsgd}

As we have seen in Section~\ref{sec:local_analysis}, {\tt L2GD} is a specific instance of {\tt SGD}, thus only converges linearly to the neighborhood of the optimum. In this section, we resolve the mentioned issue by incorporating control variates to the stochastic gradient~\cite{svrg, saga}.

We also go further. We assume that each local objective has a finite-sum structure and propose an algorithm---{\tt L2SGD+}---which takes \emph{local stochastic} gradient steps, while maintaining (global) linear convergence rate. As a consequence, {\tt L2SGD+} is the first local {\tt SGD} with linear convergence.\footnote{We are aware that a linearly converging local {\tt SGD} (with $\lambda = \infty$) might be obtained as a particular instance of the decoupling method from~\cite{mishchenko2019stochastic}, although this was not stated in the mentioned paper. Other variance reduced local {\tt SGD} algorithms~\cite{liang2019variance, karimireddy2019scaffold, wu2019federated} are not capable of achieving linear convergence.} For the reader's convenience, we present a variance reduced local gradient descent (i.e., no subsampling) in the Appendix. 

\subsection{Setup}
Consider $f_\tRloc(x_{\tRloc}) = \frac1m \sum_{j =1}^m \flocc_{i,j}(x_{\tRloc})$. Therefore, the objective function~\eqref{eq:local_main} becomes
\[
  F(x)  =\frac1n \sum \limits_{\tRloc=1}^{\TRloc} \underbrace{  \left( \frac1m\sum_{j=1}^m\flocc_{i,j}(x_{\tRloc})\right)}_{=  f_{\tRloc}(x) } + \lambda\underbrace{\frac{1}{2\TRloc} \sum_{\tRloc=1}^{\TRloc} \| x_{\tRloc} - \bar{x}\|^2.}_{= \Phi(x)}
\]

\begin{assumption}\label{as:local_smooth_sc_simple}
Function $\flocc_{i,j}$ is convex, $\Lloc$ smooth while $f_i$ is $\mu$-strongly convex (for each $1\leq j \leq m, 1\leq \tRloc \leq \TRloc$).
\end{assumption}
Denote $\onesmloc\in \R^m$ to be vector of ones. We are now ready to state {\tt L2SGD+} as Algorithm~\ref{alg:local_L2SGD}.

\begin{algorithm}[h]
  \caption{{\tt L2SGD+}: Loopless Local  {\tt SGD} with Variance Reduction}
  \label{alg:local_L2SGD}
\begin{algorithmic}
\State{\bfseries Input: }{$x^0_1 = \dots = x_{n}^0\in\R^{d}$, stepsize $\alpha$, probability $p$ }
\State $\mJf^0_{\tRloc}  = 0 \in \R^{d\times m}, \mJpsi^0_{\tRloc}  = 0 \in \R^{d}$ (for $\tRloc = 1,\dots, \TRloc$)
  \For{$k=0,1,2,\dotsc$}
  \State $\xi = 1$ with probability $p$ and $0$ with probability $1-p$ 
  \If {$\xi=0$}
 \State All {\color{red}Devices} $i=1,\dots,n$:      
          \State \hskip .3cm Sample $j \in \{1,\dots, m\}$ (uniformly at random)
       \State \hskip .3cm $g^k_\tRloc = \frac{1}{\TRloc(1-\pagg)} \left(\nabla \flocc_{i,j}(x^k_{\tRloc})- \left(\mJf^k_{\tRloc}\right)_{:,j}\right) + \frac{ \mJf^k_{\tRloc} \onesmloc}{\TRloc m}  +   \frac{\mJpsi^k_{\tRloc}}{\TRloc}  $ 
         \State \hskip .3cm $x^{k+1}_\tRloc = x^{k}_\tRloc - \alpha g^k_\tRloc$
        \State  \hskip .3cm Set $(\mJ^{k+1}_i)_{:,j} = \nabla \flocc_{i,j}(x^k_\tRloc)$, $\mJpsi^{k+1}_{\tRloc}= \mJpsi^k_{\tRloc}$,
        \State \hskip 0.83cm  $(\mJ^{k+1}_i)_{:,l} = (\mJ^{k+1}_i)_{:,l} $ for all $l\neq j$
         		\vspace*{-0.1cm}
 \Else 
  \State {\color{blue}Master} computes the average $\bar{x}^k = \frac{1}{n}\sum_{i=1}^n x_i^k$
\State {\color{blue}Master} does for all $i=1,\dots,n$:       
   \State \hskip .3cm $g^k_\tRloc=   \frac{\lambda }{\TRloc \pagg} (x^k_{\tRloc} - \bar{x}^k) - \frac{\pagg^{-1} -1}{ \TRloc}\mJpsi^k_{\tRloc} + \frac{1}{\TRloc m} \mJf_{\tRloc}^k \onesmloc$
 		\State  \hskip .3cm Set $x_{\tRloc}^{k+1} = x_{\tRloc}^k  - \alpha g^k_\tRloc $
 		\State  \hskip .3cm Set $\mJpsi^{k+1}_{\tRloc} =  \lambda (x_{\tRloc}^k - \bar{x}^k )$, $\mJf^{k+1}_{\tRloc} = \mJf^{k}_{\tRloc} $
    \EndIf
      \EndFor
\end{algorithmic}
\end{algorithm}

\begin{remark}
{\tt L2SGD+} is the simplest local {\tt SGD} method with variance reduction. In the Appendix, we present general {\tt L2SGD++} which allows for 1) an arbitrary number of data points per client and arbitrary local subsampling strategy, 2) partial participation of clients, and 3) local {\tt SVRG}-like updates of control variates (thus potentially better memory). Lastly, {\tt L2SGD++} is able exploit the smoothness structure of the local objectives, resulting in tighter rates.
\end{remark}

{\tt L2SGD+} only communicates when a two consecutive coin tosses land a different value, thus, on average $p(1-p)k$ times per $k$ iterations. However, {\tt L2SGD+} requires communication of control variates $\mJf_i \onesmloc,  \mJpsi_i$ as well -- each communication round is thus three times more expensive. In the Appendix, we provide an implementation of {\tt L2SGD+} that does not require the communication of $\mJf_i \onesmloc,  \mJpsi_i$.

\subsection{Theory}

We are now ready to present a convergence rate of {\tt L2SGD+}.

\begin{theorem} \label{thm:local_saga_simple}

Let Assumption~\ref{as:local_smooth_sc_simple} hold and choose
$  \alpha =\TRloc \min \left \{  \frac{(1-\pagg) }{ 4  \Lloc + \mu m} , \frac{\pagg}{4\lambda+ \mu} \right\}.$
Then the iteration complexity of Algorithm~\ref{alg:local_L2SGD} is 
$  \max \left\{   \frac{4\Lloc  + \mu m}{ (1-\pagg)\mu}, \frac{4\lambda +\mu}{\pagg\mu}\right \}\log\frac1\varepsilon.$
\end{theorem}

Next, we find a probability $p$ that yields both the best iteration and communication complexity. 

\begin{corollary}\label{cor:local_lsd_optimal_p}
Both communication and iteration complexity of {\tt L2SGD+} are minimized for $\pagg  = \frac{4\lambda +\mu}{4\lambda +4\Lloc +(m+1)\mu}$. The resulting iteration complexity is $\left(4 \frac{\lambda }{\mu} +4 \frac{\Lloc}{\mu} +m+1\right) \log \frac1\varepsilon$, while the communication complexity is $\frac{4 \lambda + \mu }{4\Lloc + 4 \lambda + (m+1)\mu} \left( 4\frac{\Lloc}{\mu}+m\right) \log \frac{1}{\varepsilon}$. 
\end{corollary}

Note that with $\lambda \rightarrow \infty$, the communication complexity of {\tt L2SGD+} tends to $ \left( 4\frac{\Lloc}{\mu}+m\right) \log \frac{1}{\varepsilon}$, which is communication complexity of minibatch {\tt SAGA} to find the globally optimal model (see Chapter~\ref{jacsketch}). On the other hand, in the pure local setting ($\lambda=0$), the communication complexity becomes $\log \frac1\epsilon$ -- this is because the Lyapunov function involves a term that measures the distance of local models, which requires communication to be estimated.

\begin{figure}[th]
\centering
\begin{minipage}{0.5\textwidth}
  \centering
\includegraphics[width =  \textwidth ]{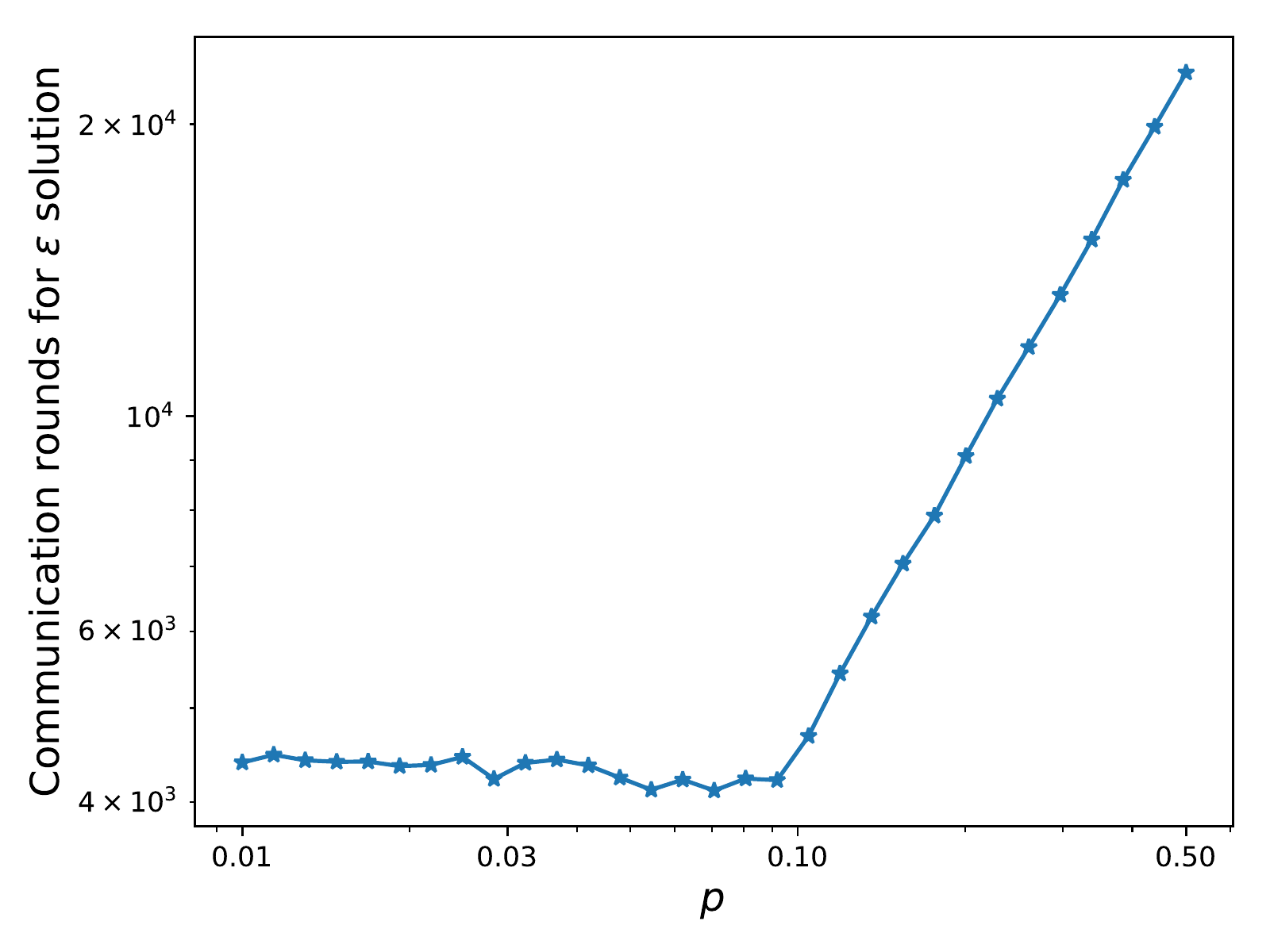}
\end{minipage}%
\caption{ Communication rounds to get $\frac{F(x^k)-F(x^*)}{F(x^0)-F(x^*)}\leq 10^{-5}$ as a function of $p$ with $p^* \approx 0.09$ (for {\tt L2SGD+}). Logistic regression on a1a dataset with $\lambda = 0.1$; details in the Appendix. } 
\label{fig:local_comm_p}
\end{figure}

\section{Experiments \label{sec:local_exp}}
In this section, we numerically verify the theoretical claims from this chapter.  In all experiments in this chapter, we consider a simple binary classification model -- logistic regression. In particular, suppose that device $i$ owns data matrix $\mA_{i}\in \R^{m \times d}$ along with corresponding labels $b_{i}\in \{-1, 1\}^{m}$. The local objective for client $i$ is then given as follows
\[
f_{i}(x) \eqdef  \frac1m \sum_{j=1}^m \flocc_{i, j}(x) +\frac{\mu}{2} \| x\|^2, \quad \text{where} \quad \flocc_{i m+ j}(x) \eqdef \log \left(1+\exp\left((\mA_{i})_{j,:}x\cdot  b_{i}\right) \right).
\]
The rows of data matrix $\mA$ were normalized to have length 4 so that each $\flocc_{i,j}$ is $1$-smooth for each $j$. At the same time, the local objective on each device is $10^{-4}$-strongly convex. Next, datasets are from LibSVM~\cite{chang2011libsvm}. 

In each case, we consider the simplest locally stochastic algorithm. In particular, each dataset is evenly split among the clients, while the local stochastic method samples a single data point each iteration.

We have chosen a different number of clients for each dataset -- so that we cover different possible scenarios. See Table~\ref{tbl:local_data} for details (it also includes sizes of the datasets). Lastly, the stepsize was always chosen according to Theorem~\ref{thm:local_saga_simple}.

\begin{table}[!h]
\caption{Setup for the experiments.}
\label{tbl:local_data}
\begin{center}
\small
\begin{tabular}{|c|c|c|c|c|c|c|c|c|c|}
\hline
Dataset &  \begin{tabular}{c} $\nlocc$ \\ $=\TRloc m$ \end{tabular}   & $d$ &  $\TRloc$ & $m$   & $\mu$&$L$& \begin{tabular}{c} $\pagg$ \\ \ref{sec:local_exp_comp} \end{tabular} & \begin{tabular}{c} $\lambda$ \\ \ref{sec:local_exp:pagg} \end{tabular} & \begin{tabular}{c} $\pagg$ \\ \ref{sec:local_exp:omega} \end{tabular}\\
\hline
\hline
\texttt{a1a} & 1 605	& 123 & 5 & 321 & $10^{-4}$ & $1$&0.1 &0.1 &0.1\\ \hline
\texttt{mushrooms} & 8 124 & 112 & 12 &677 & $10^{-4}$ & $1$ &0.1 & 0.05 & 0.3 \\ \hline
\texttt{phishing} & 11 055	&	68 & 11 & 1 005 & $10^{-4}$ & $1$ &0.1 &0.1 & 0.001\\ \hline
\texttt{madelon} & 2 000& 500& 50 &40 & $10^{-4}$& $1$ &0.1  &0.02&0.05\\ \hline
\texttt{duke} &44 & 7 129& 4 & 11& $10^{-4}$ & $1$ & 0.1 & 0.4&0.1\\ \hline
\texttt{gisette\_scale} & 6 000& 5 000& 100 & 60 & $10^{-4}$ & $1$&0.1 &0.2&0.003  \\ \hline
\texttt{a8a} &22 696 & 123& 8 & 109& $10^{-4}$ & $1$ &0.1&0.1  &0.1\\ \hline
\end{tabular}
\end{center}
\end{table}

\subsection{Comparison of the methods \label{sec:local_exp_comp}}

In our first experiment, we verify two phenomena:
\begin{itemize}
\item Effect of variance reduction on the convergence speed of local methods. We compare 3 different methods: local {\tt SGD} with full variance reduction (Algorithm~\ref{alg:local_L2SGD}), shifted local {\tt SGD} (Algorithm~\ref{alg:local_lsgd_partial}) and local {\tt SGD} (Algorithm~\ref{alg:local_lsgd_none}). Our theory predicts that a fully variance reduced algorithm converges to the global optimum linearly, while both shifted local {\tt SGD} and local {\tt SGD} converge to a neighborhood of the optimum. At the same time, the neighborhood should be smaller for shifted local {\tt SGD}.

\item The claim that heterogeneity of the data does not influence the convergence rate. We consider two splits of the data heterogeneous and homogenous. For the homogenous split, we first randomly reshuffle the data and then construct the local objectives according to the current order (i.e., the first client owns the first $m$ indices, etc.). For heterogenous split, we first sort the data based on the labels and then construct the local objectives accordingly (thus achieving the worst-case heterogeneity). Note that the overall objective to solve is different in homogenous and heterogenous case -- we thus plot relative suboptimality of the objective (i.e., $\frac{F(x^k)-F(x^{*})}{F(x^0)-F(x^{*})}$) to directly compare the convergence speed. 
\end{itemize} 

In all cases, we choose $\pagg=0.1$ and $\lambda = \frac19$ -- such choice mean that $\pagg$ is very close to optimal. The other parameters (i.e. number of clients) are provided in Table~\ref{tbl:local_data}. Figure~\ref{fig:local_libsvm_methods} presents the result.

\begin{figure}[!h]
\centering
\begin{minipage}{0.3\textwidth}
  \centering
\includegraphics[width =  \textwidth ]{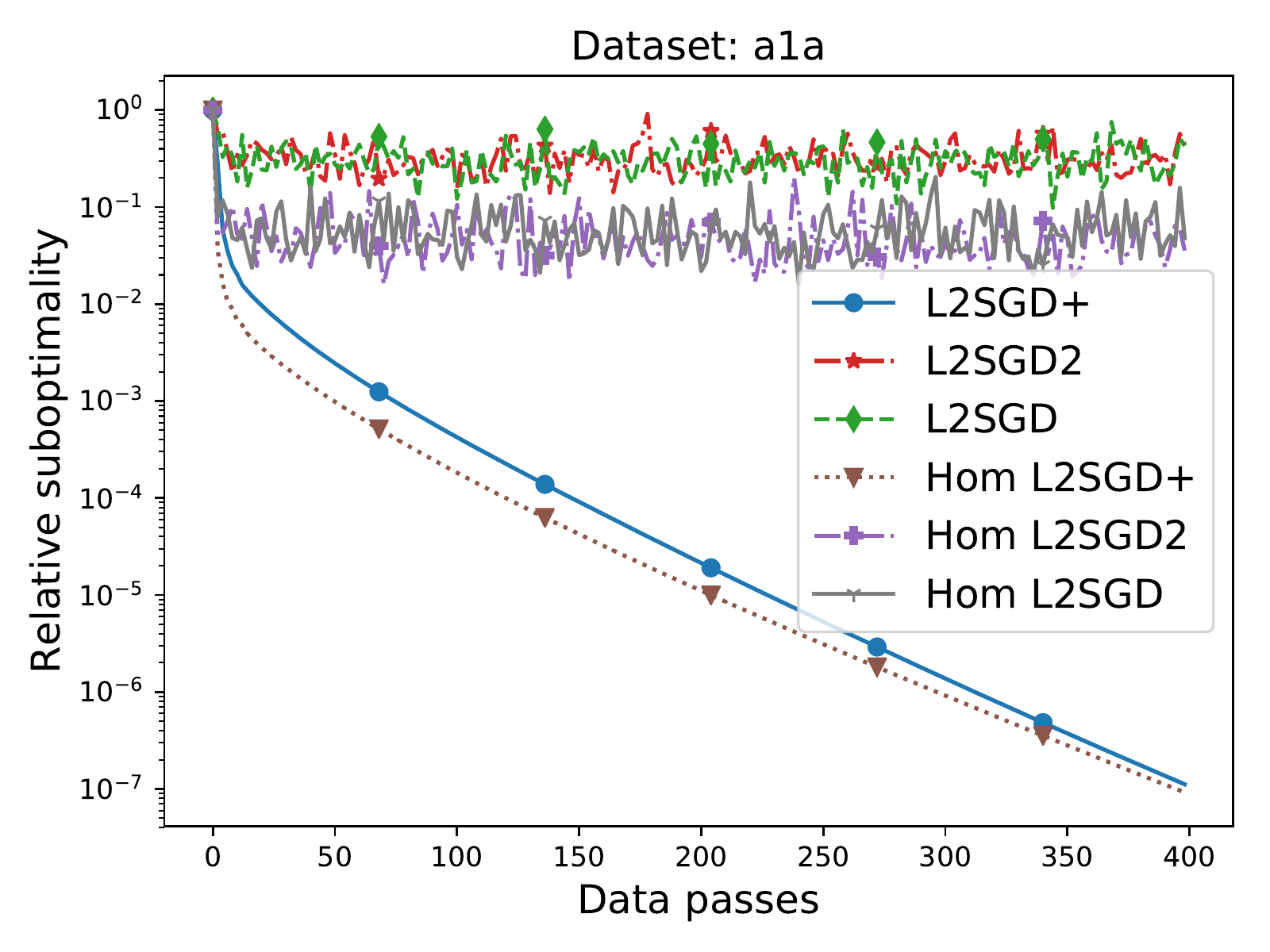}
\end{minipage}%
\begin{minipage}{0.3\textwidth}
  \centering
\includegraphics[width =  \textwidth ]{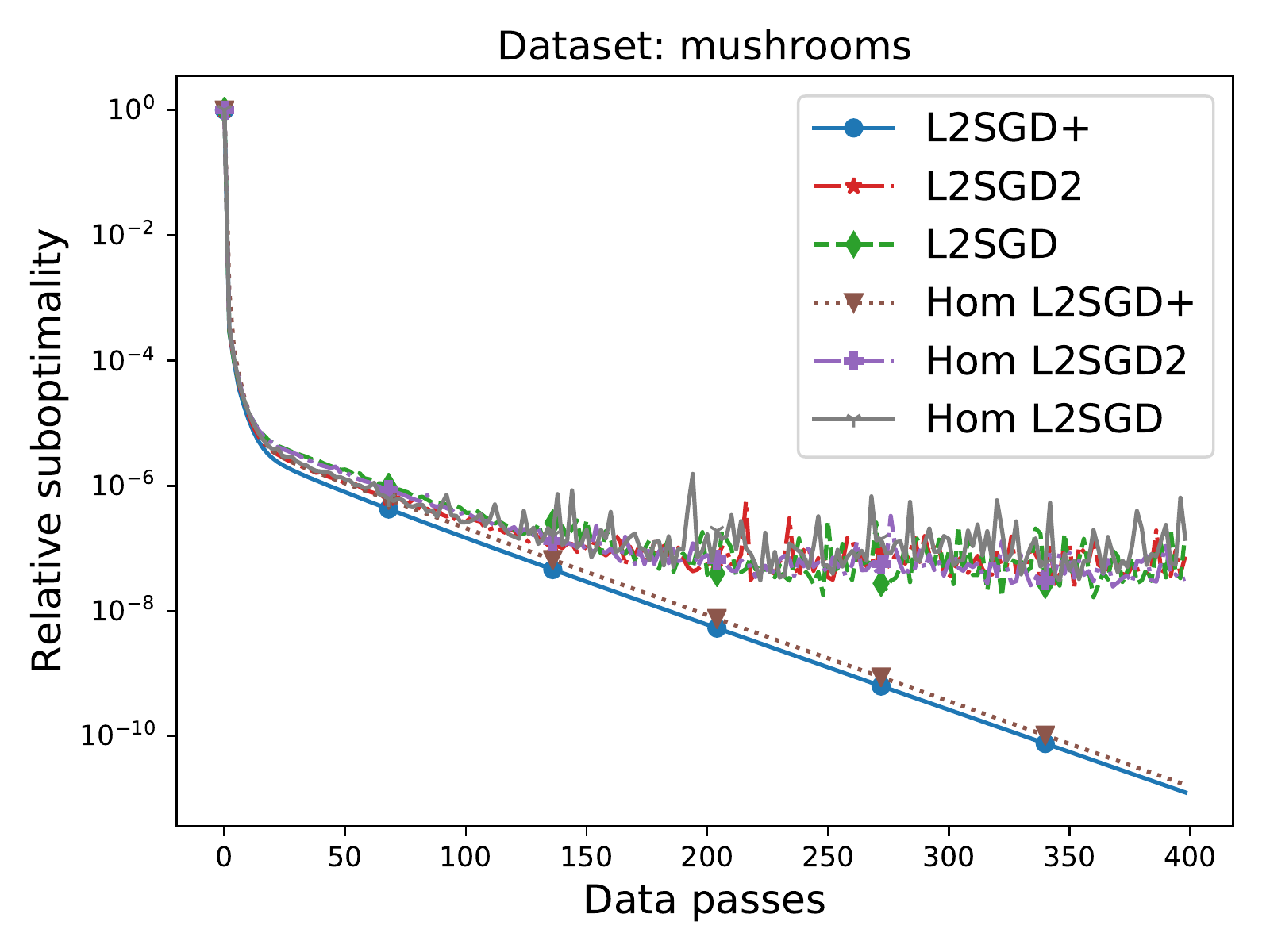}
\end{minipage}%
\begin{minipage}{0.3\textwidth}
  \centering
\includegraphics[width =  \textwidth ]{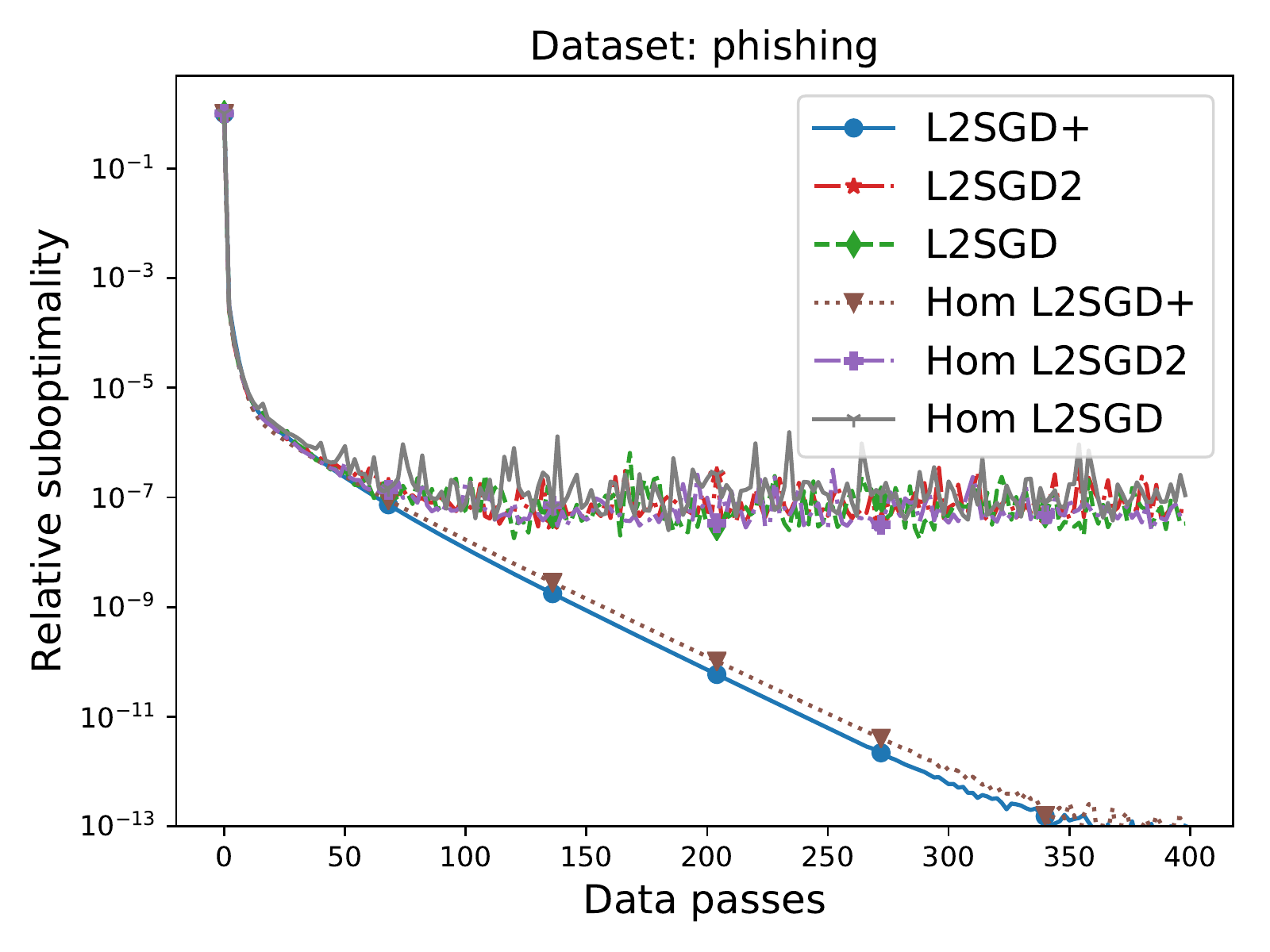}
\end{minipage}%
\\
\begin{minipage}{0.3\textwidth}
  \centering
\includegraphics[width =  \textwidth ]{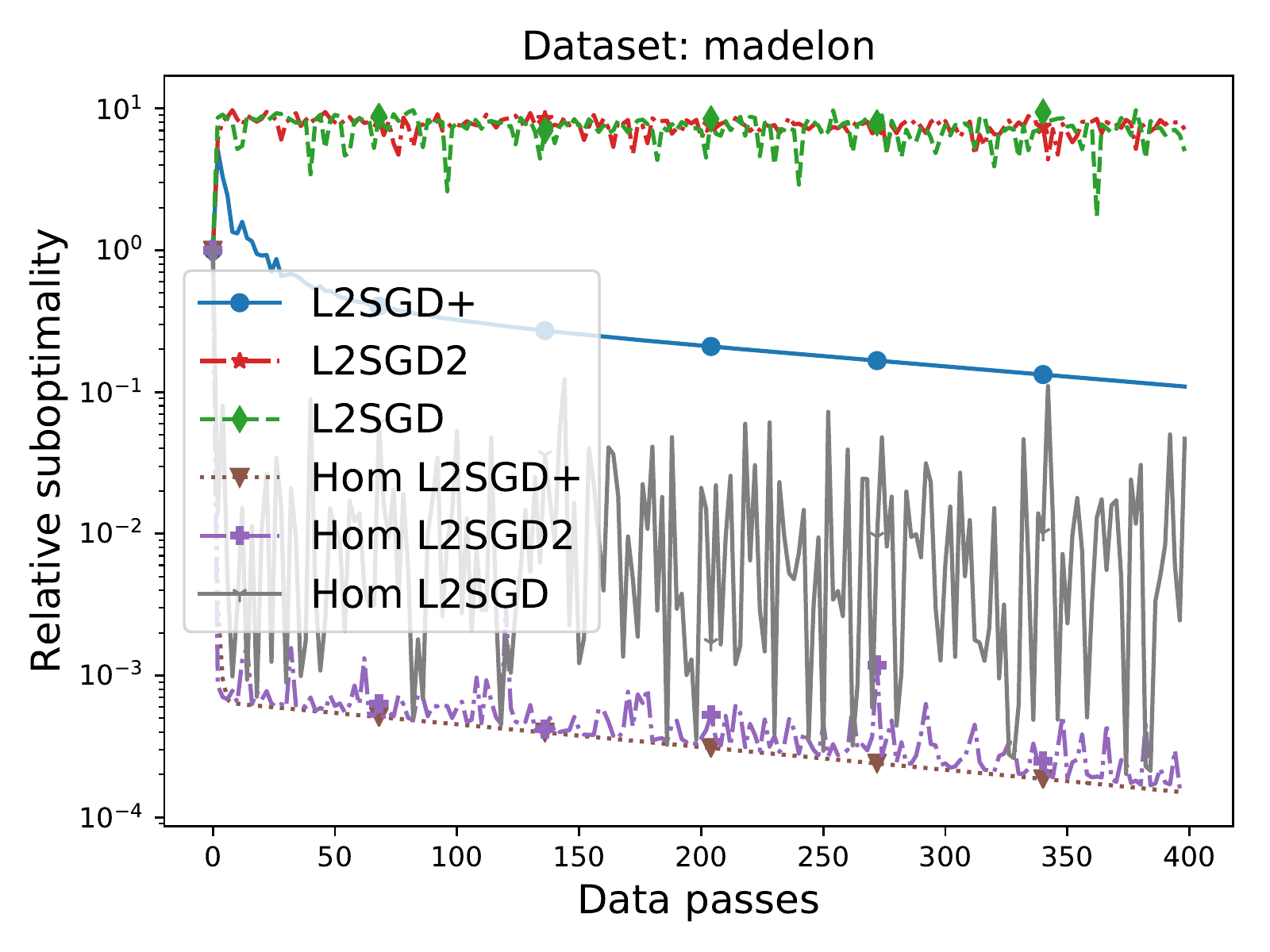}
\end{minipage}%
\begin{minipage}{0.3\textwidth}
  \centering
\includegraphics[width =  \textwidth ]{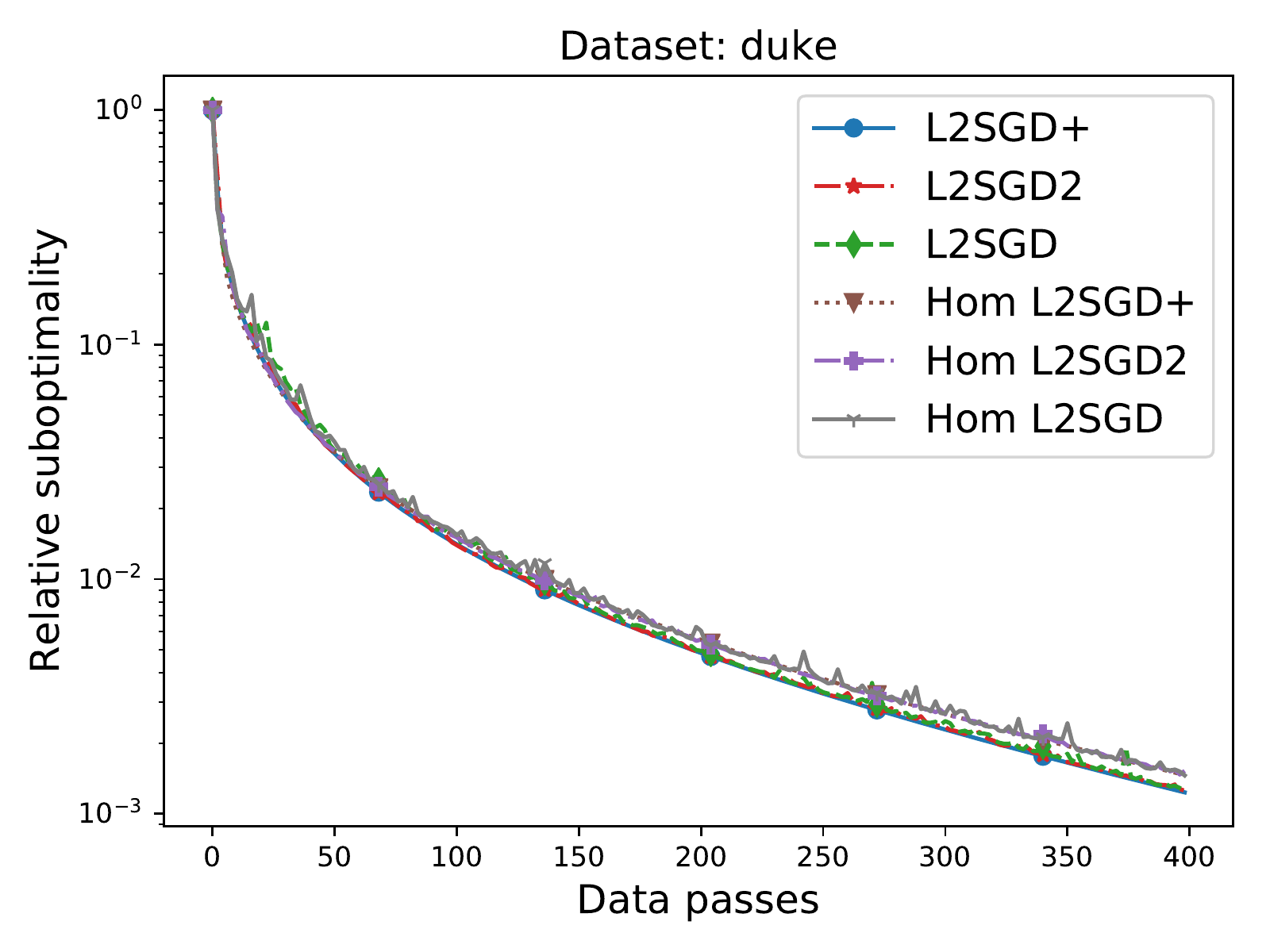}
\end{minipage}%
\begin{minipage}{0.3\textwidth}
  \centering
\includegraphics[width =  \textwidth ]{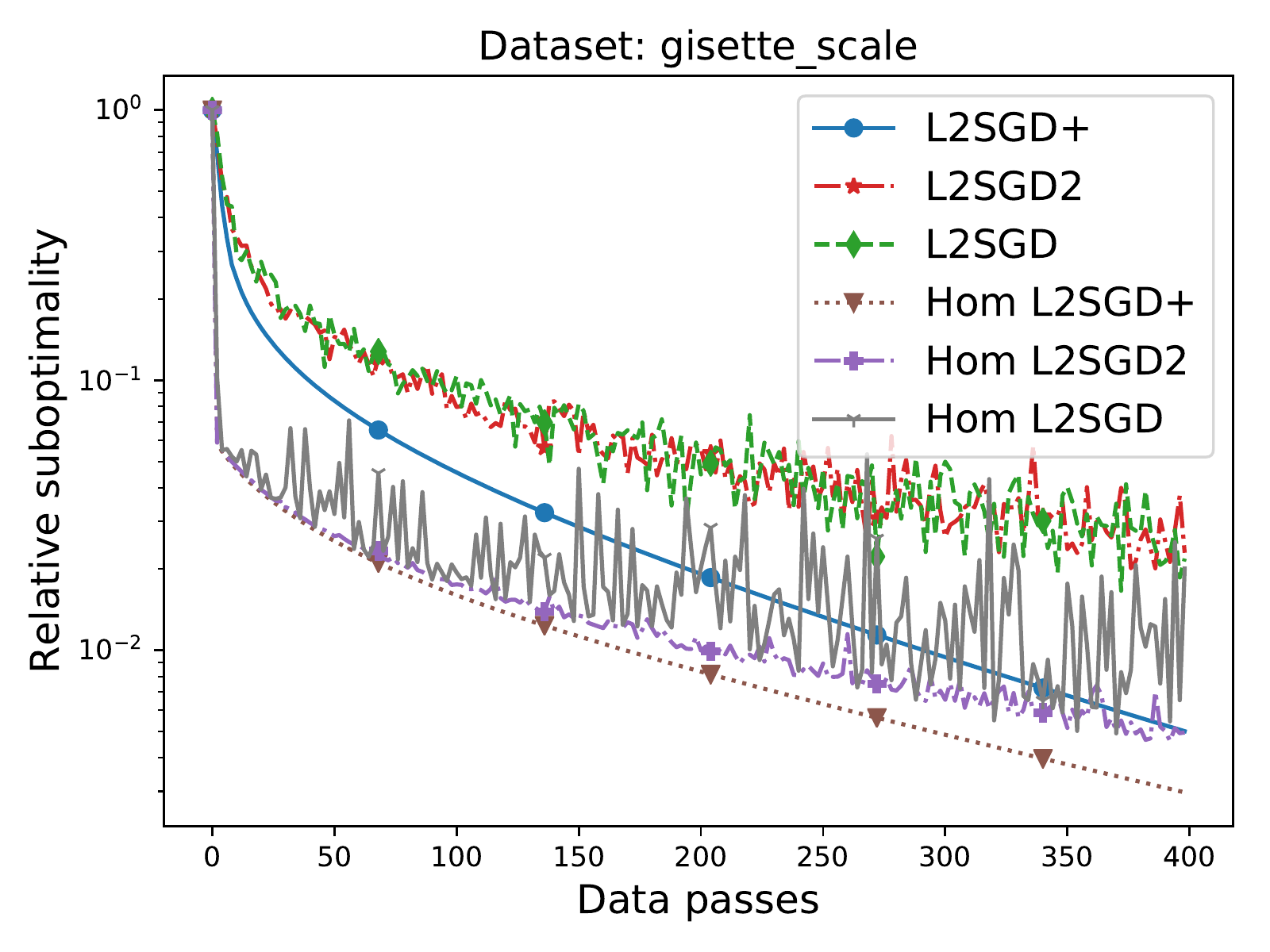}
\end{minipage}%
\\
\begin{minipage}{0.3\textwidth}
  \centering
\includegraphics[width =  \textwidth ]{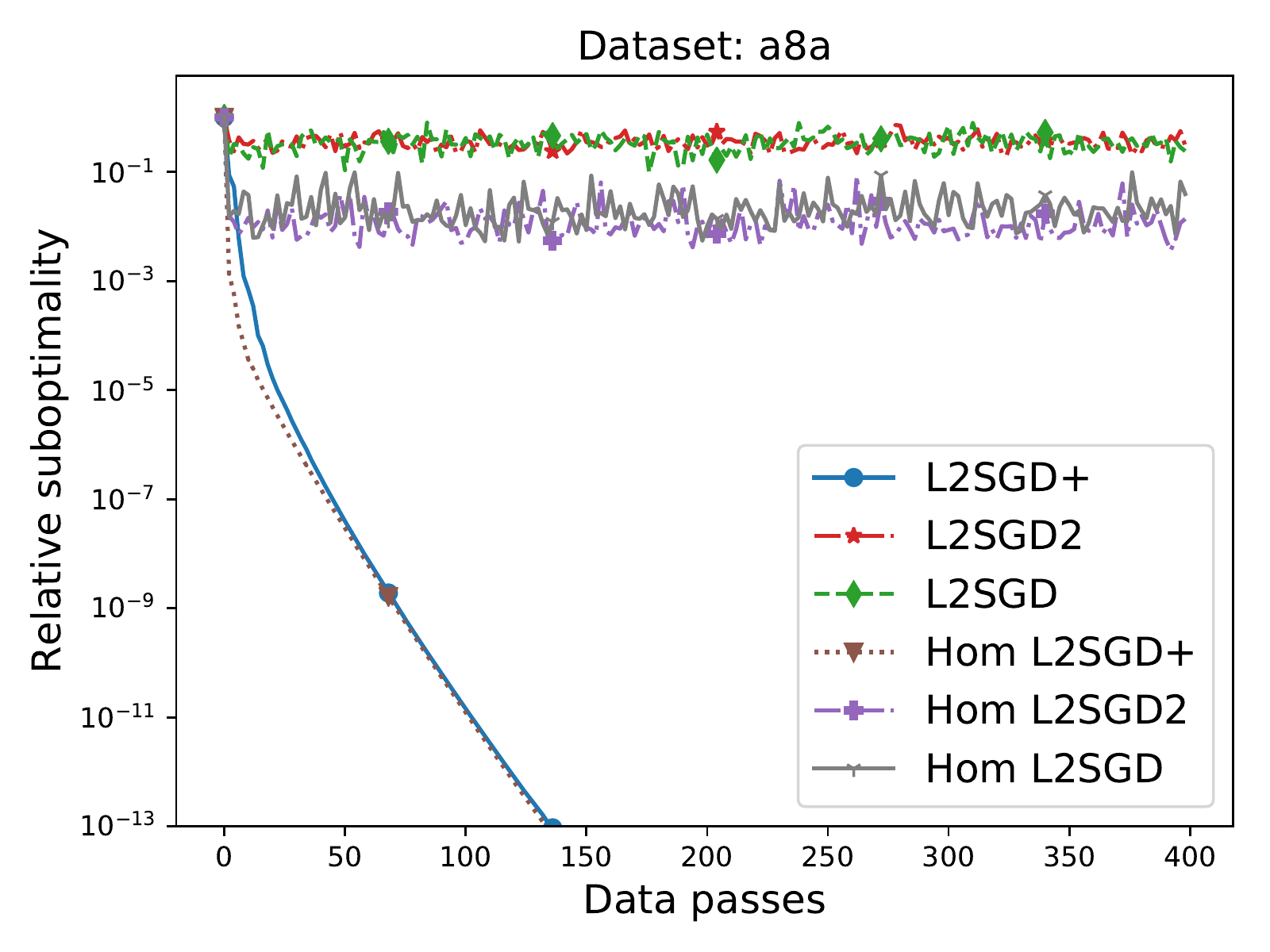}
\end{minipage}%
\caption{Variance reduced local {\tt SGD} (Algorithm~\ref{alg:local_L2SGD}), shifted local {\tt SGD} (Algorithm~\ref{alg:local_lsgd_partial}) and local {\tt SGD} (Algorothm~\ref{alg:local_lsgd_none}) applied on LibSVM problems for both homogenous split of data and Heterogenous split of the data. Stepsize for non-variance reduced method was chosen the same as for the analogous variance reduced method.} 
\label{fig:local_libsvm_methods}
\end{figure}

As expected, Figure~\ref{fig:local_libsvm_methods} clearly demonstrates the following: 
\begin{itemize}
\item Full variance reduction always converges to the global optima, methods with partial variance reduction only converge to a neighborhood of the optimum.
\item Partial variance reduction (i.e., shifting the local {\tt SGD}) is better than not using control variates at all. Although the improvement in the performance is rather negligible.
\item Data heterogeneity does not affect the convergence speed of the proposed methods. Therefore, unlike standard local {\tt SGD}, mixing the local and global models does not suffer the problems with heterogeneity.
\end{itemize}

\subsection{Effect of $\pagg$\label{sec:local_exp:pagg}}

In the second experiment, we study the effect of $\pagg$ on the convergence rate of variance reduced local {\tt SGD}. Note that $\pagg$ immediately influences the number of communication rounds -- on average, the clients take $(\pagg^{-1}-1)$ local steps in between two consecutive rounds of communication (aggregation). 

In Section~\ref{sec:local_lsgd}, we argue that, it is optimal (in terms of the convergence rate) to choose $\pagg$ of order $\pagg^{*}\eqdef \frac{\lambda}{\Lloc+\lambda}$. Figure~\ref{fig:local_libsvm_pagg} compares $\pagg=  \pagg^{*}$ against other values of $\pagg$ and confirms its optimality (in terms of optimizing the convergence rate). 

While the slower convergence of Algorithm~\ref{alg:local_L2SGD} with $\pagg<\pagg^{*}$ is expected (i.e., communicating more frequently yields a faster convergence), slower convergence for $\pagg>\pagg^{*}$ is rather surprising; in fact, it means that communicating less frequently yields faster convergence. This effect takes place due to the specific structure of problem~\eqref{eq:local_main}; it would be lost when enforcing $x_1=\dots = x_{\TRloc}$ (corresponding to $\lambda=\infty$).

\begin{figure}[!h]
\centering
\begin{minipage}{0.3\textwidth}
  \centering
\includegraphics[width =  \textwidth ]{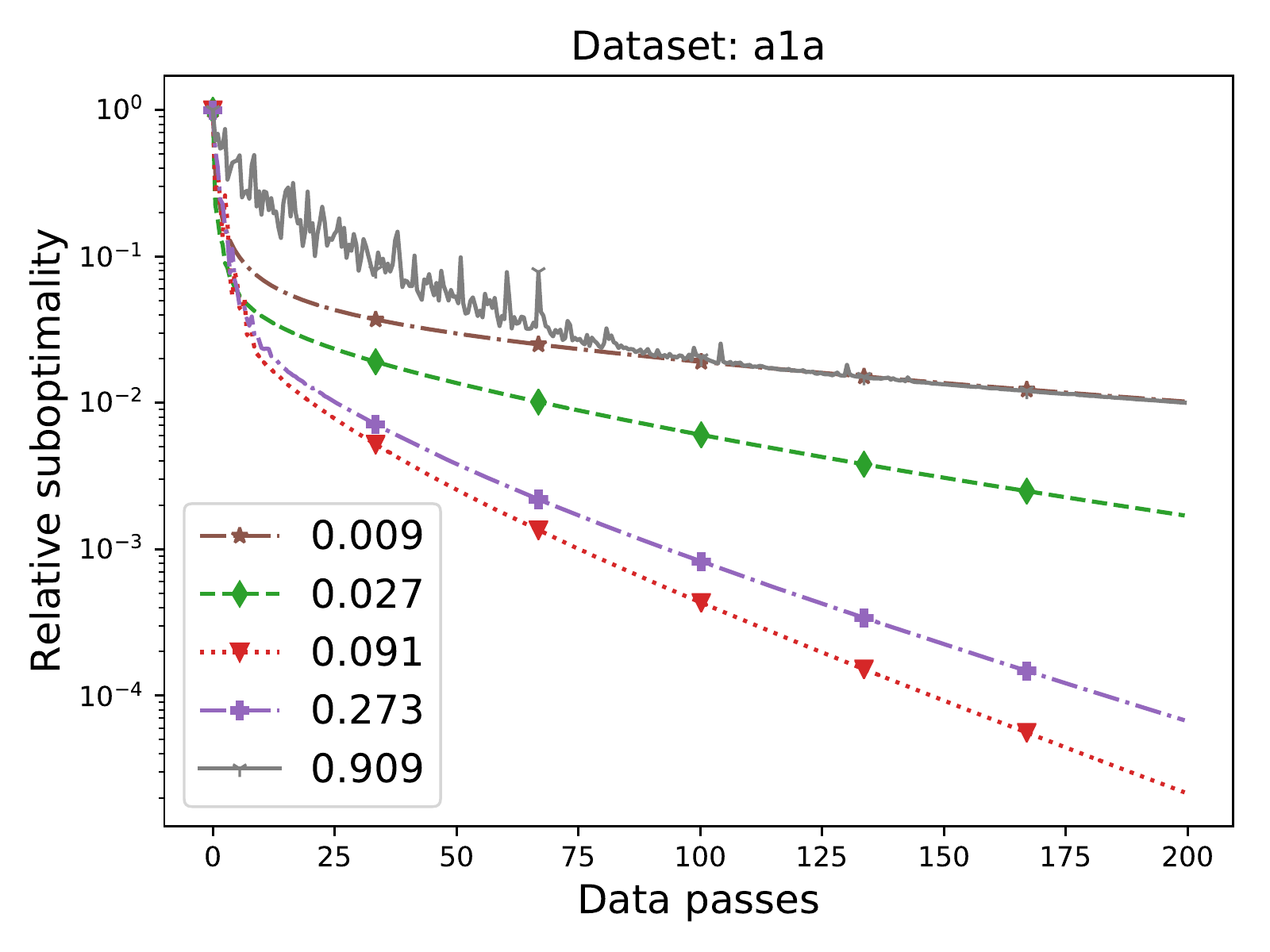}
\end{minipage}%
\begin{minipage}{0.3\textwidth}
  \centering
\includegraphics[width =  \textwidth ]{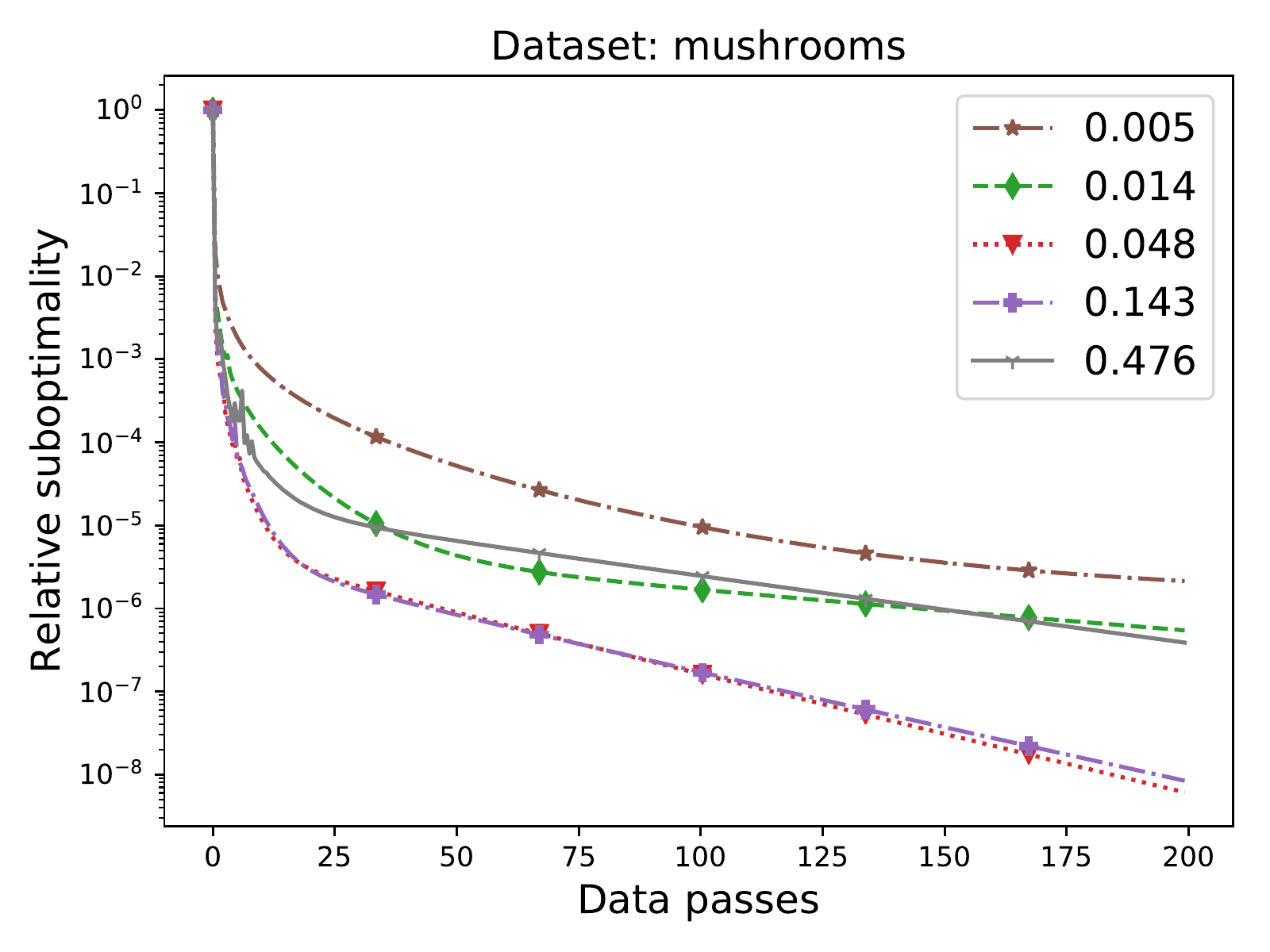}
\end{minipage}%
\begin{minipage}{0.3\textwidth}
  \centering
\includegraphics[width =  \textwidth ]{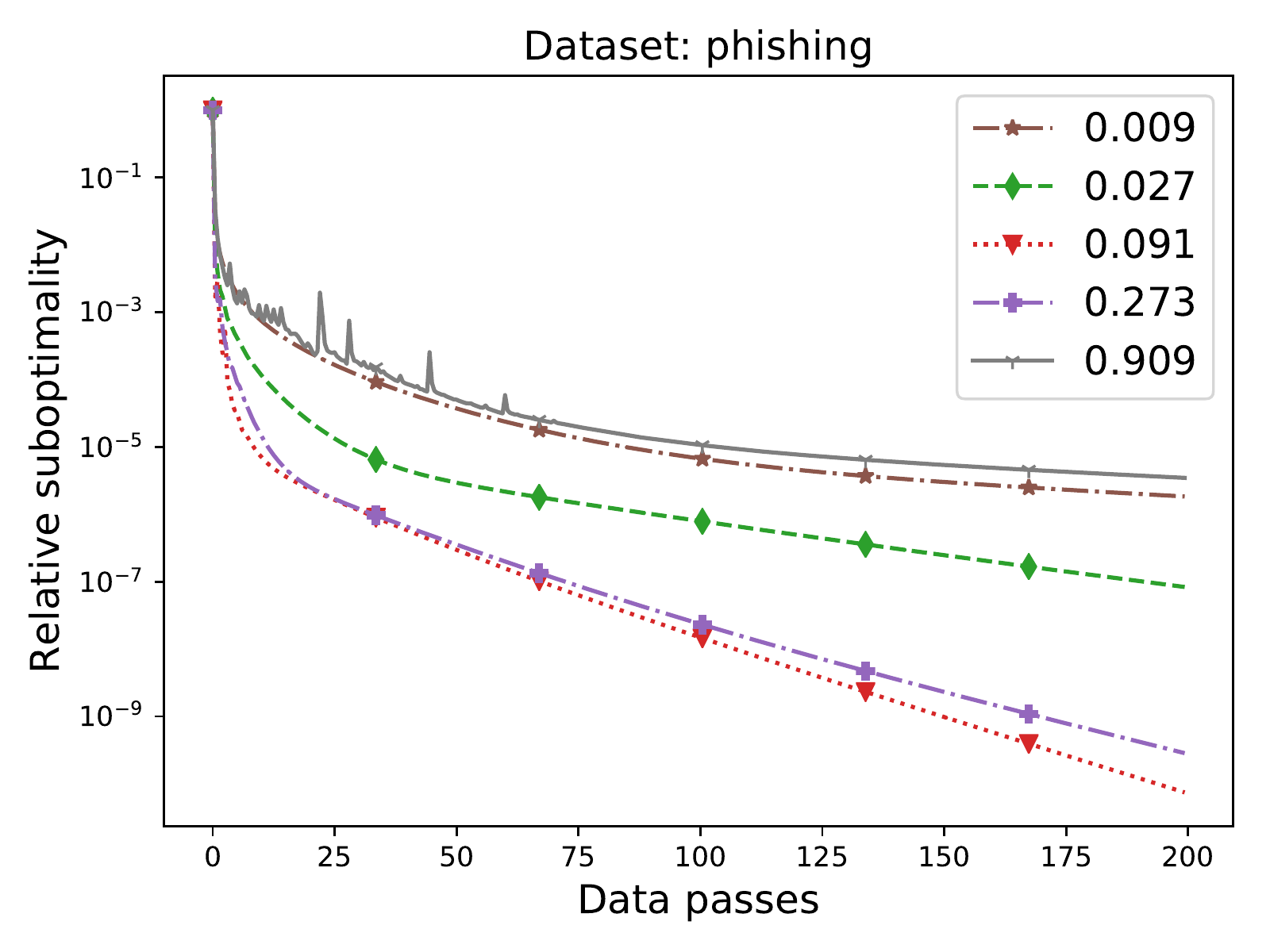}
\end{minipage}%
\\
\begin{minipage}{0.3\textwidth}
  \centering
\includegraphics[width =  \textwidth ]{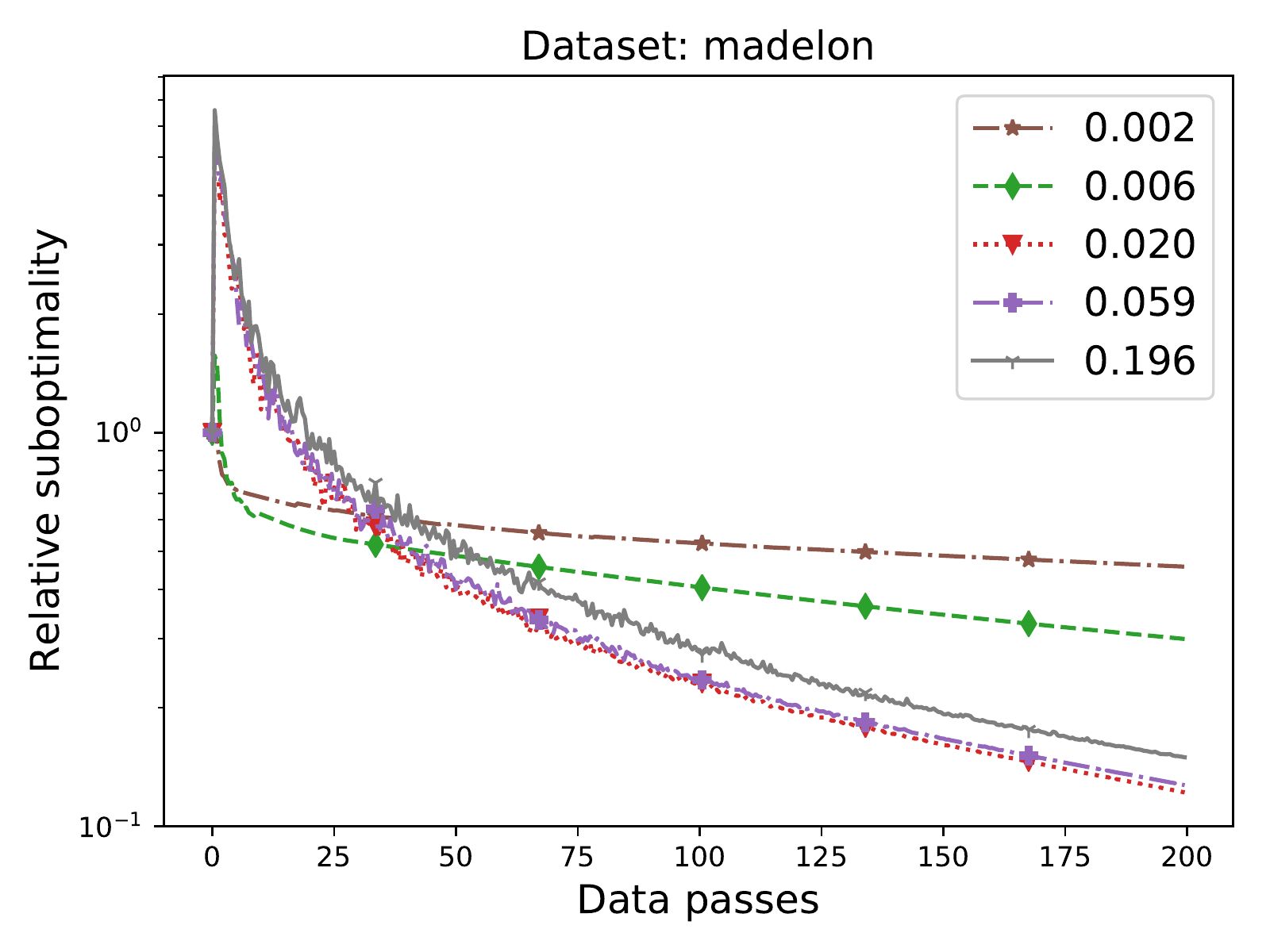}
\end{minipage}%
\begin{minipage}{0.3\textwidth}
  \centering
\includegraphics[width =  \textwidth ]{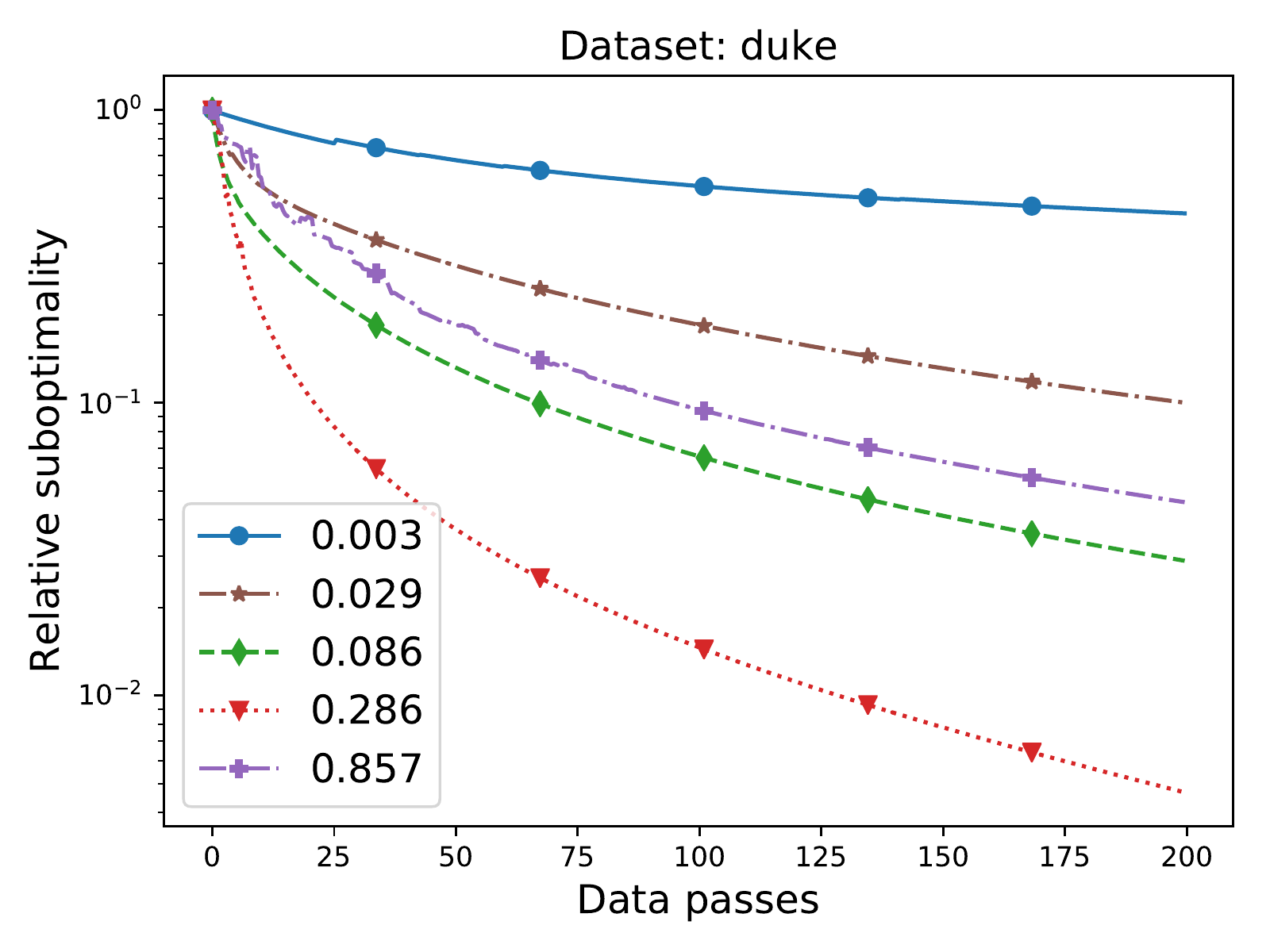}
\end{minipage}%
\begin{minipage}{0.3\textwidth}
  \centering
\includegraphics[width =  \textwidth ]{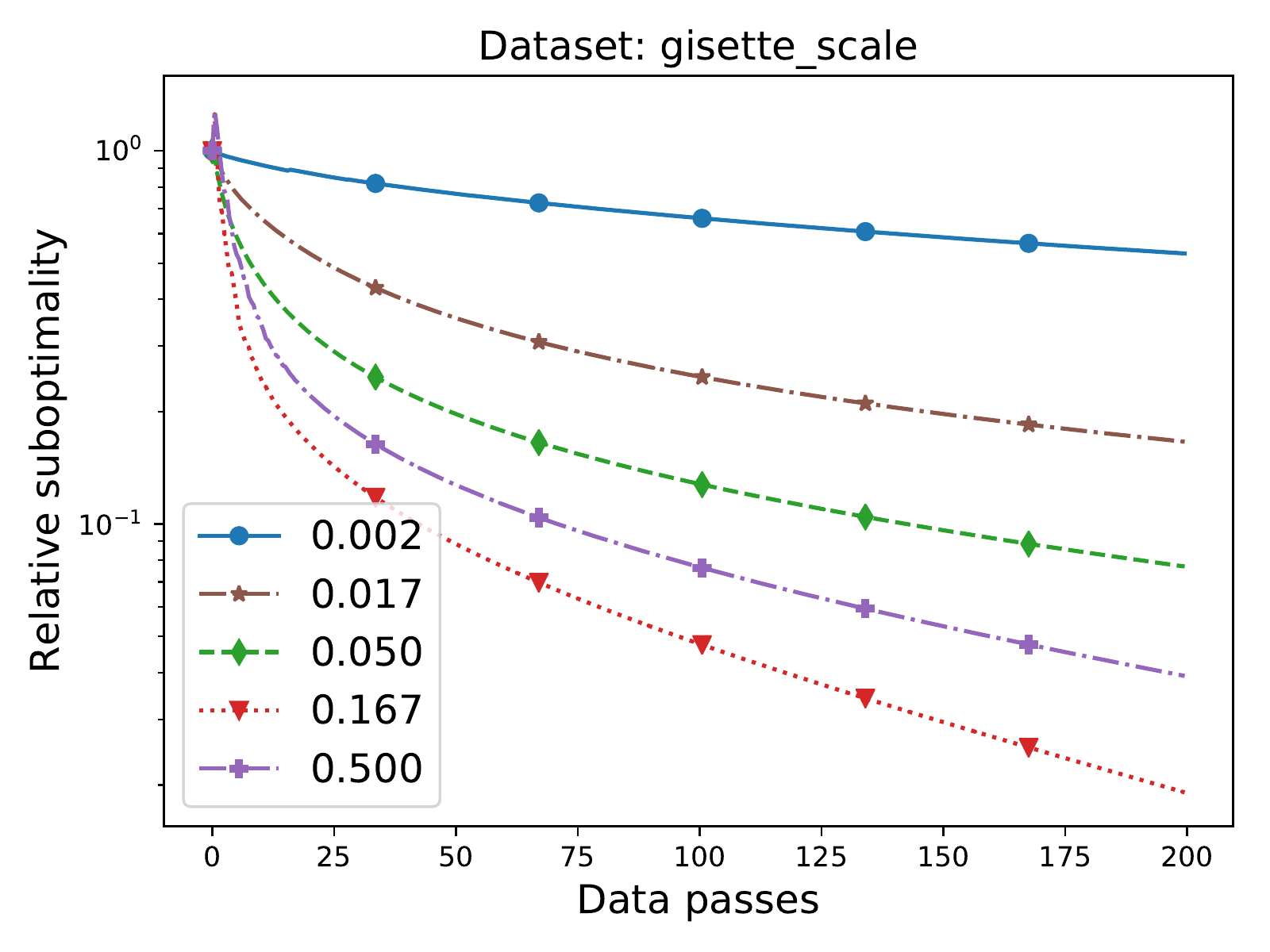}
\end{minipage}%
\\
\begin{minipage}{0.3\textwidth}
  \centering
\includegraphics[width =  \textwidth ]{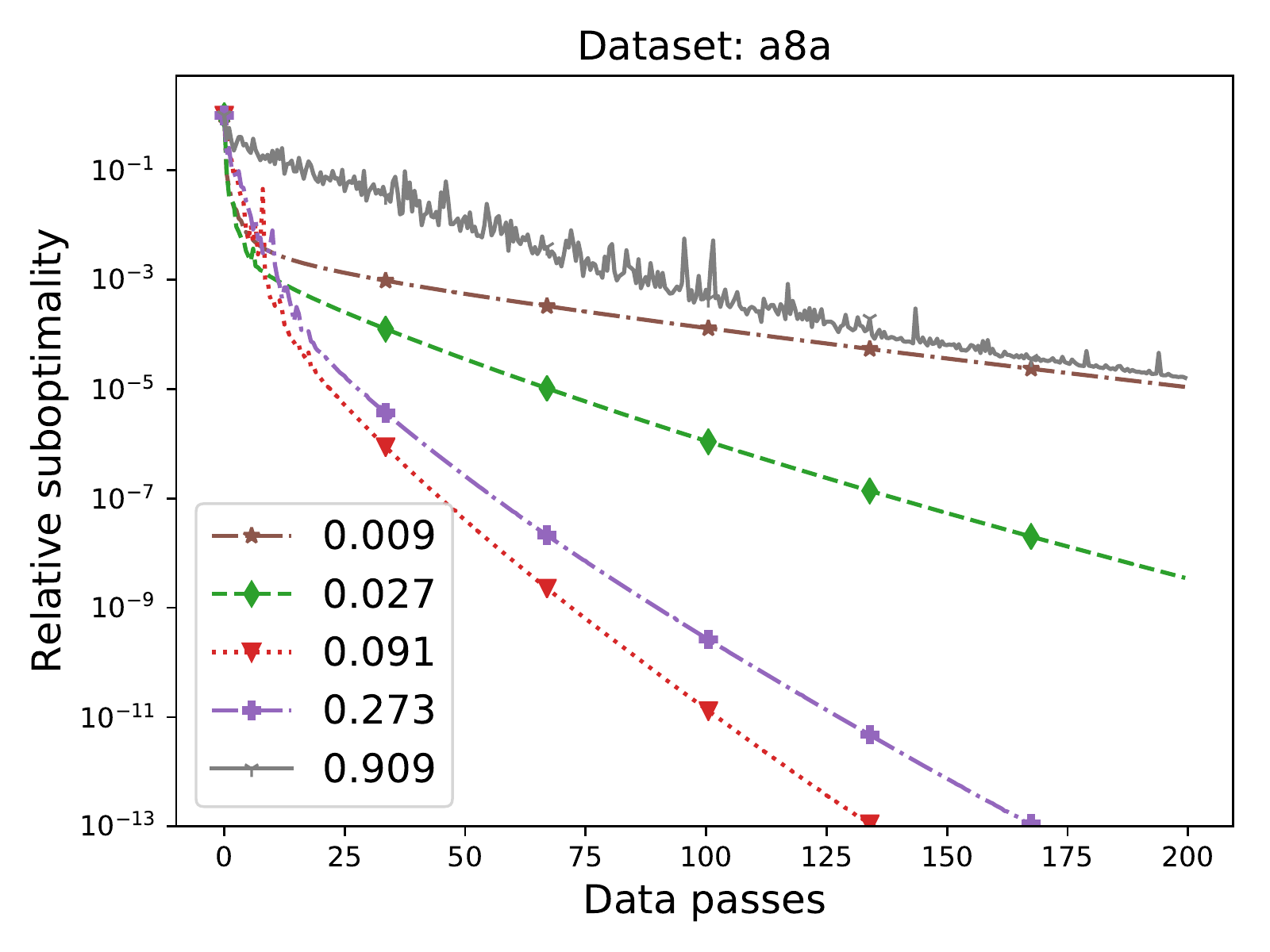}
\end{minipage}%
\caption{Effect of the aggregation probability $\pagg$ (legend of the plots) on the convergence rate of Algorithm~\ref{alg:local_L2SGD}. Choice $\pagg = \pagg^{*}$ corresponds to red dotted line with triangle marker. Parameter $\lambda$ was chosen in each case as Table~\ref{tbl:local_data} indicates. } 
\label{fig:local_libsvm_pagg}
\end{figure}

\subsection{Effect of $\lambda$ \label{sec:local_exp:omega}}
In this experiment we study how different values of $\lambda$ influence the convergence rate of Algorithm~\ref{alg:local_L2SGD}, given that everything else (i.e. $\pagg$) is fixed. Note that for each value of $\lambda$ we get a different instance of problem~\eqref{eq:local_main}; thus the optimal solution is different as well. Therefore, in order to make a fair comparison between convergence speeds, we plot the relative suboptimality (i.e. $\frac{F(x^k)-F(x^{*})}{F(x^0)-F(x^{*})}$) against the data passes. Figure~\ref{fig:local_libsvm_omega} presents the results.

The complexity of Algorithm~\ref{alg:local_L2SGD} is\footnote{Given that $\mu$ is small.} $\cO\left(\frac{\Lloc}{(1-\pagg)\mu} \right)\log\frac1\varepsilon$ as soon as $\lambda< \lambda^{*} \eqdef \frac{L\pagg}{(1-\pagg)}$; otherwise the complexity is $\cO\left(\frac{\lambda}{\pagg\mu} \right)\log\frac1\varepsilon$. This perfectly consistent with what  Figure~\ref{fig:local_libsvm_omega} shows -- the choice $\lambda<\lambda^{*}$ resulted in comparable convergence speed than $\lambda=\lambda^{*}$; while the choice $\lambda>\lambda^{*}$ yields noticeably worse rate than $\lambda=\lambda^{*}$.

\begin{figure}[!h]
\centering
\begin{minipage}{0.3\textwidth}
  \centering
\includegraphics[width =  \textwidth ]{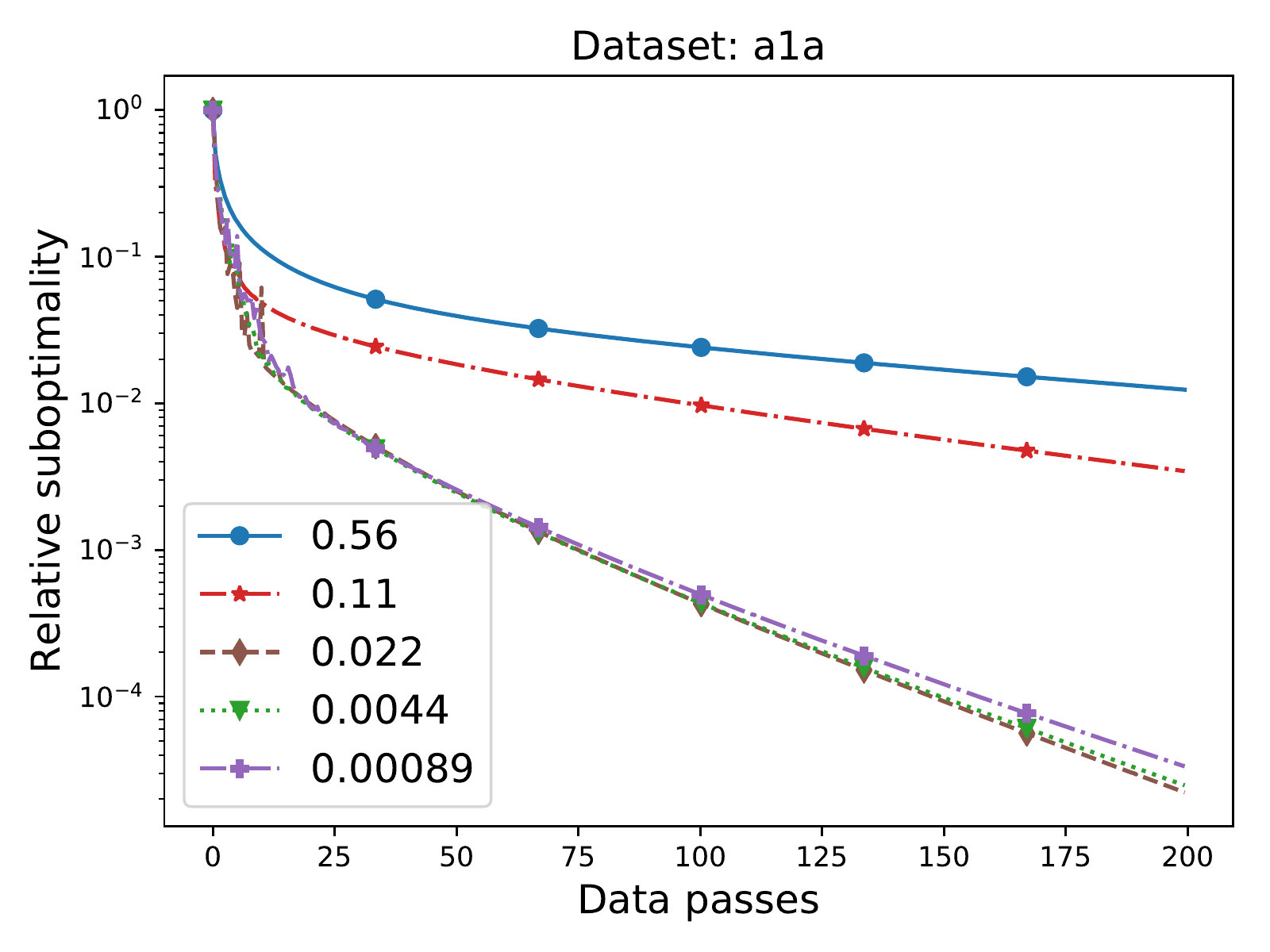}
\end{minipage}%
\begin{minipage}{0.3\textwidth}
  \centering
\includegraphics[width =  \textwidth ]{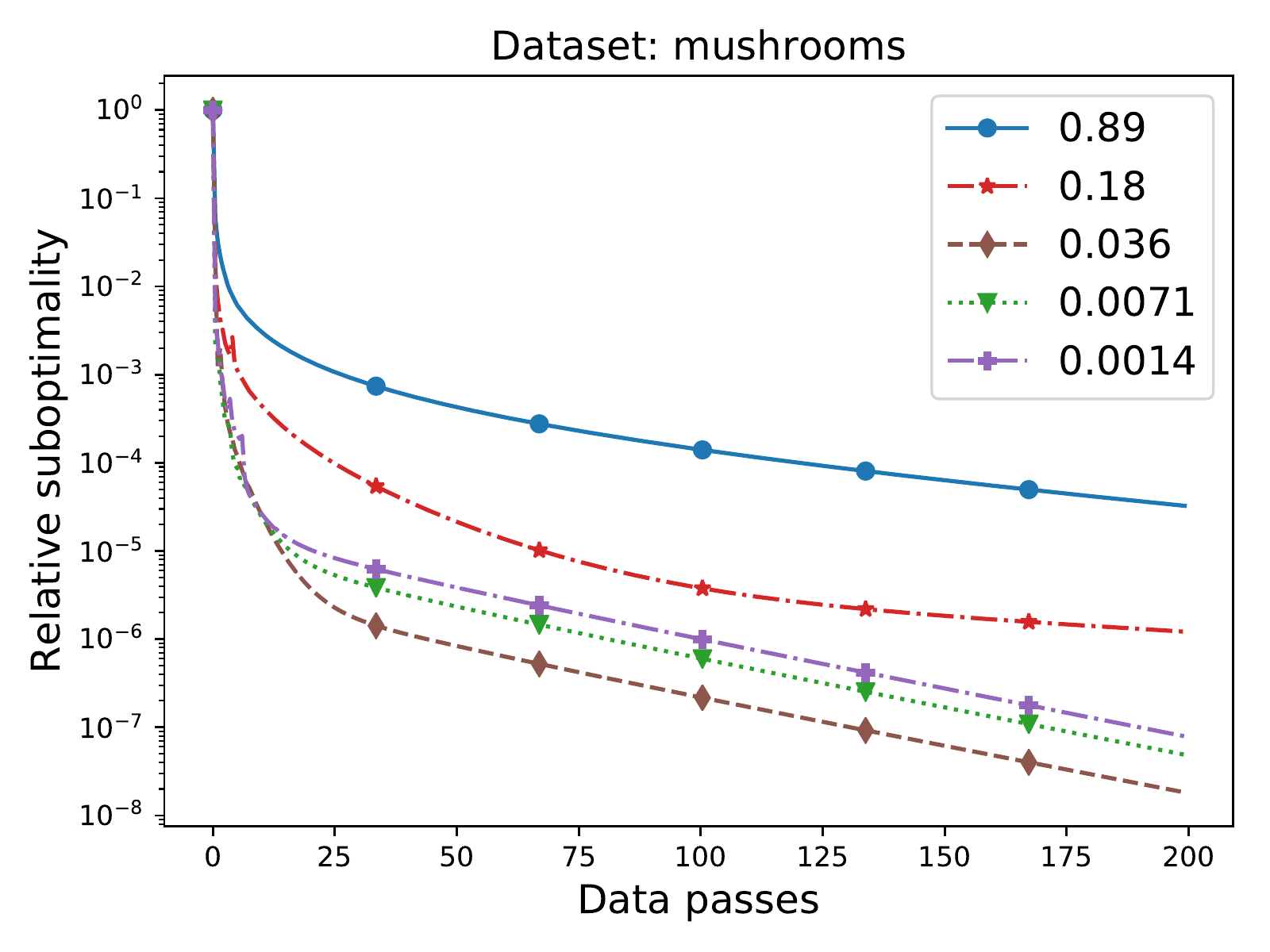}
\end{minipage}%
\begin{minipage}{0.3\textwidth}
  \centering
\includegraphics[width =  \textwidth ]{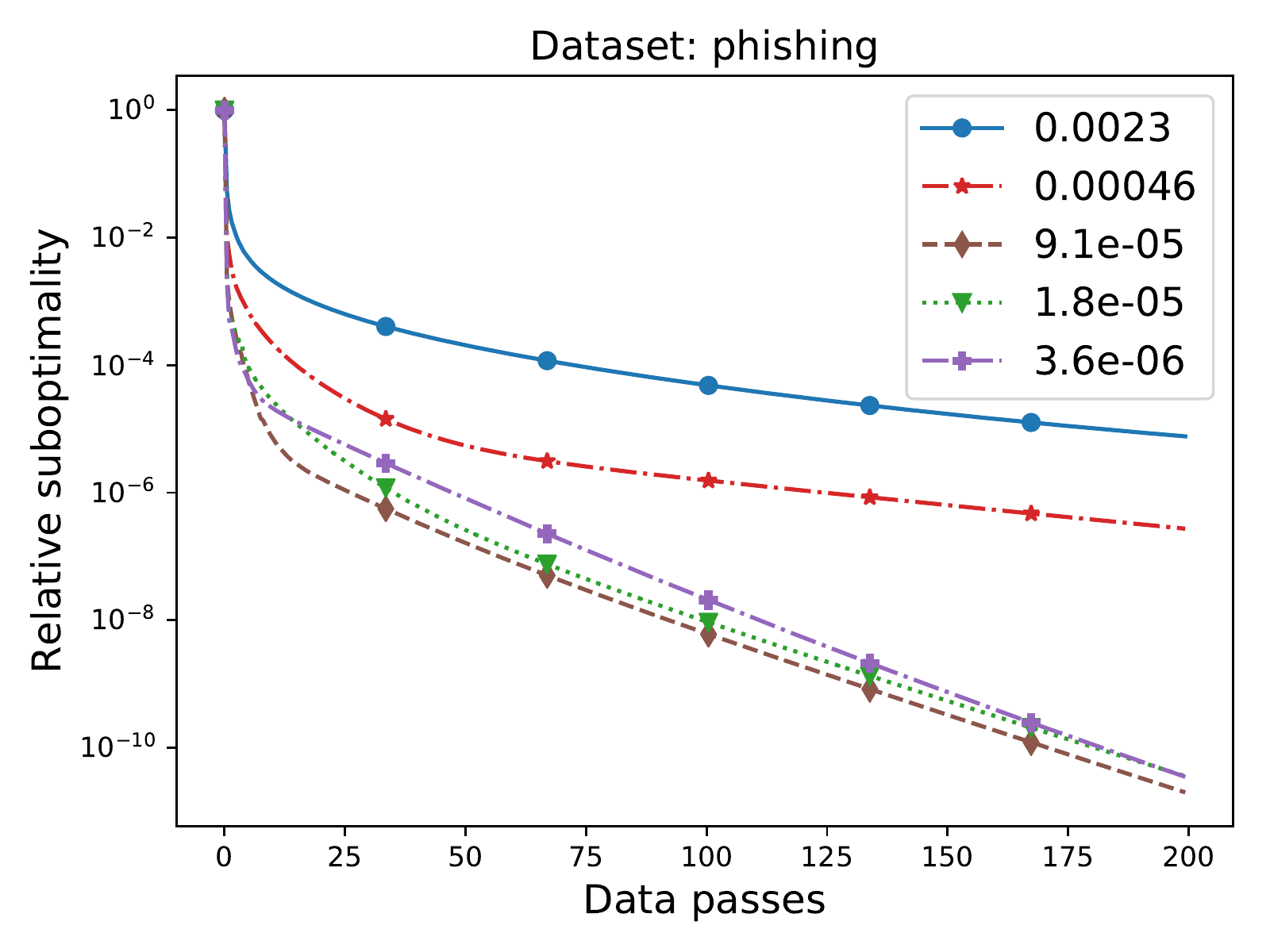}
\end{minipage}%
\\
\begin{minipage}{0.3\textwidth}
  \centering
\includegraphics[width =  \textwidth ]{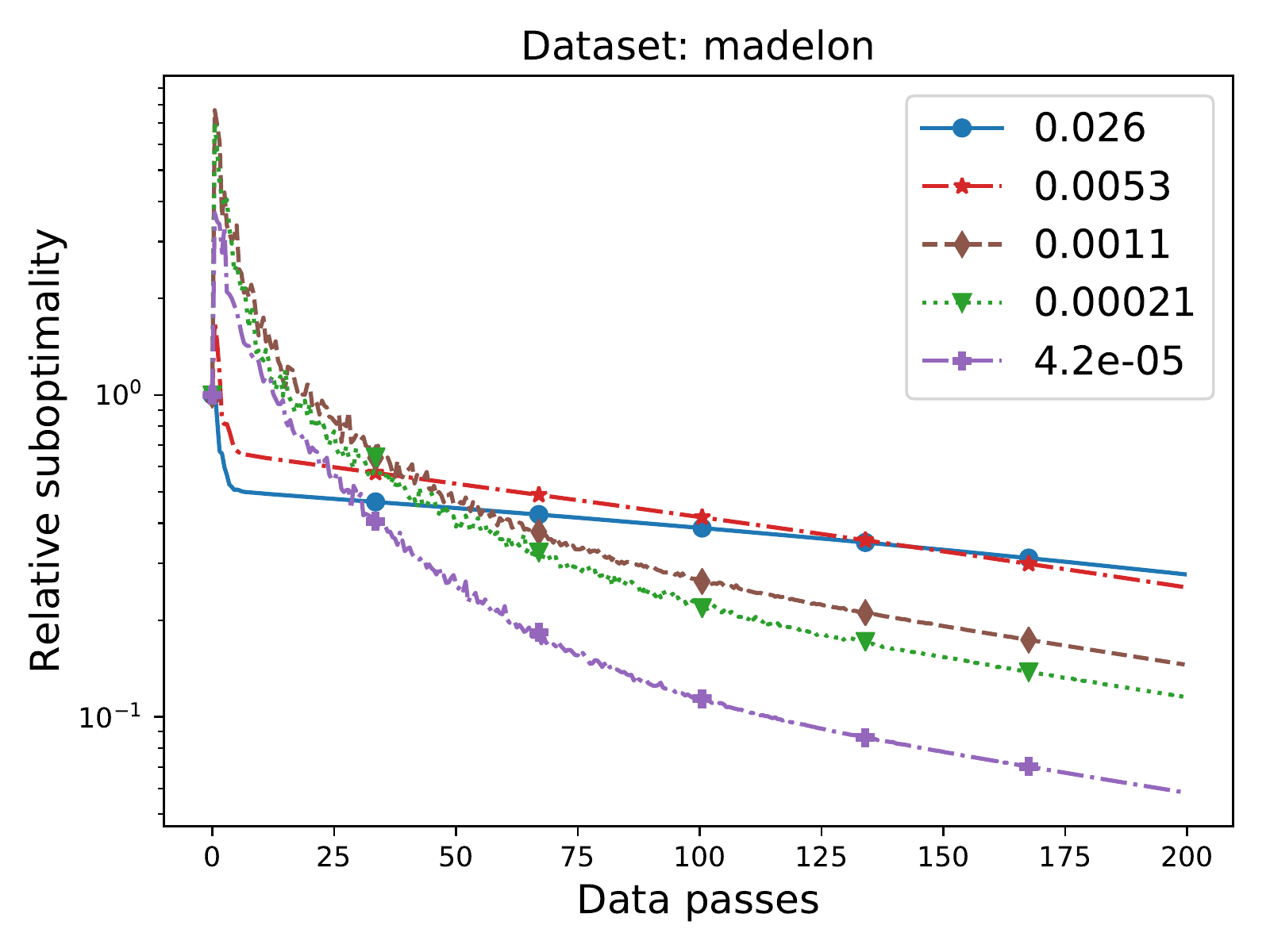}
\end{minipage}%
\begin{minipage}{0.3\textwidth}
  \centering
\includegraphics[width =  \textwidth ]{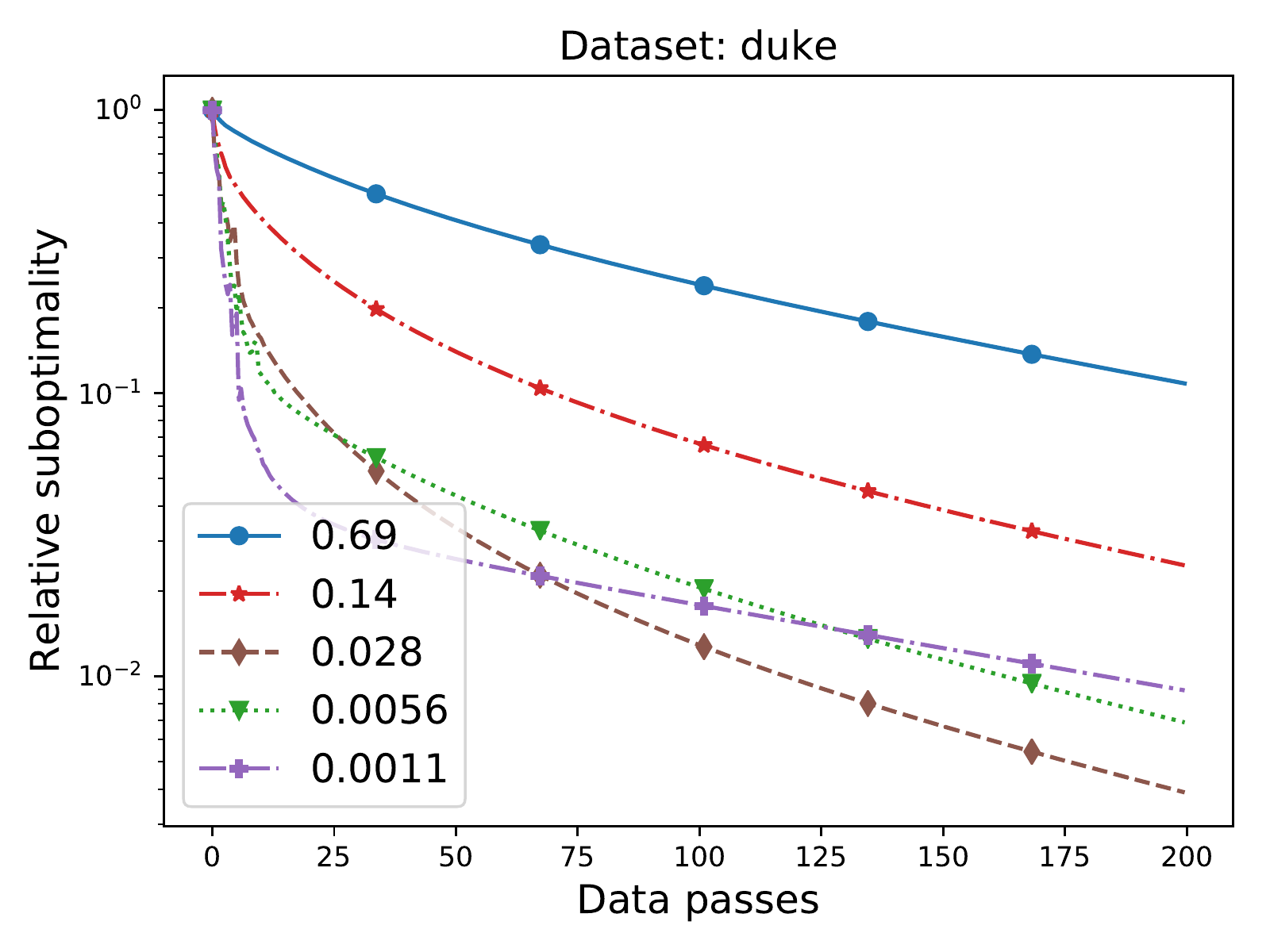}
\end{minipage}%
\begin{minipage}{0.3\textwidth}
  \centering
\includegraphics[width =  \textwidth ]{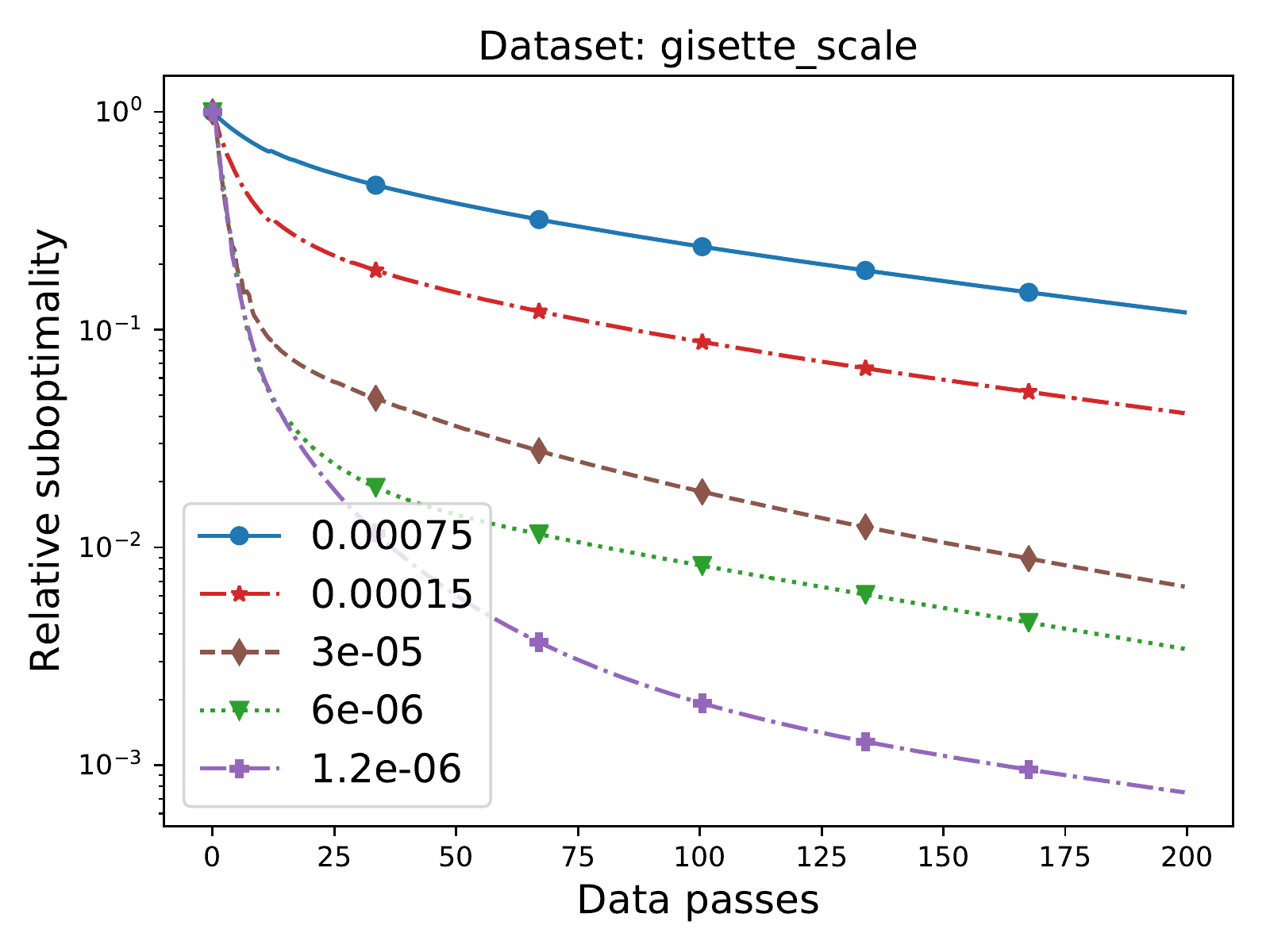}
\end{minipage}%
\\
\begin{minipage}{0.3\textwidth}
  \centering
\includegraphics[width =  \textwidth ]{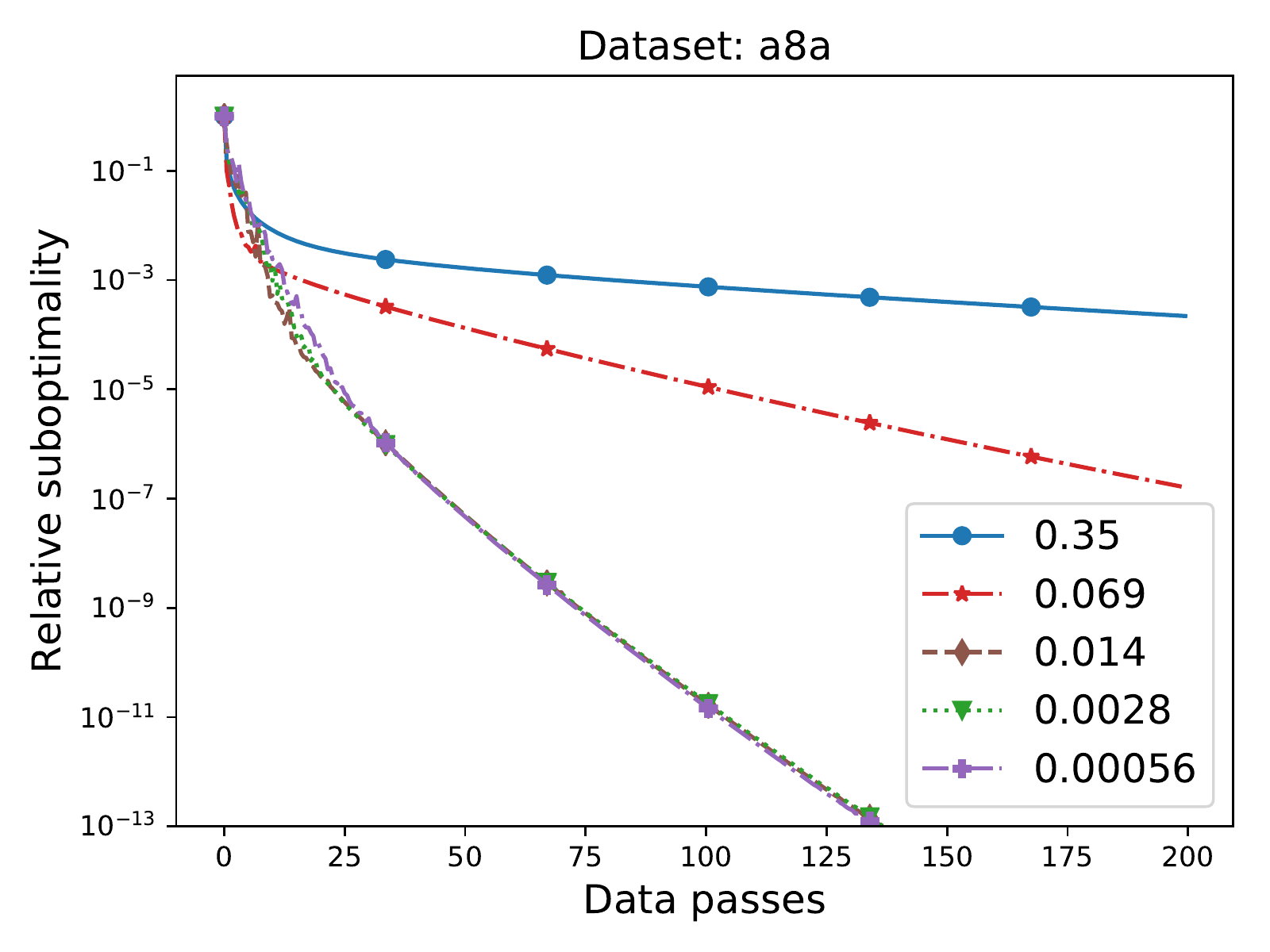}
\end{minipage}%
\caption{Effect of parameter $\lambda$ (legend of the plot) on the convergence rate of Algorithm~\ref{alg:local_L2SGD}. The choice $\lambda=\lambda^{*}$ corresponds to borwn dash-dotted line with diamond marker (the third one from the legend). Aggregation probability $\pagg$ was chosen in each case as Table~\ref{tbl:local_data} indicates. } 
\label{fig:local_libsvm_omega}
\end{figure}

\section{Conclusion}

In this chapter we have proposed a new optimization formulation for federated learning. The algorithms (i.e., {\tt L2GD}) we propose to solve the new formulation are similar the classical local {\tt SGD}, however, the rates we have provided are superior to classical local {\tt SGD} analysis. 

Our analysis of {\tt L2GD} can be extended to cover smooth convex and non-convex loss functions $f_i$ (we do not explore these directions). Further, our methods can be extended to a decentralized regime where the devices correspond to devices of a connected network, and communication is allowed along the edges of the graph only. This can be achieved by introducing an additional randomization over the penalty $\Phi$. Further, our approach can be accelerated in the sense of Nesterov~\cite{nesterov2018lectures} by adapting the results from \cite{allen2017katyusha, l-svrg-as} to our setting, thus further reducing the number of communication rounds.

\chapter{Stochastic Subspace Cubic Newton Method}
\label{sscn}

\graphicspath{{sscn/plots/}}

In this chapter we consider a regularized not necessarily finite-sum optimization problem 
\begin{equation}\label{eq:sscn_problem}
    \min_{x \in \R^d}  \, \left \{ F(x) \eqdef f(x) + \psi(x)\right \},
\end{equation}
where $f:\R^d\to \R$ is convex and twice differentiable and $\psi:\R^d\to \R\cup \{+\infty\}$ is a proximable convex function. We are interested in the regime where the dimension $d$ is very large, which arises in many contexts, such as  the training of modern over-parameterized machine learning models. In this regime, coordinate descent ({\tt CD}) methods, or more generally subspace descent methods, are the methods of choice. 

\section{Subspace descent methods}

Subspace descent methods rely on  update rules of the form
\begin{equation} \label{eq:sscn_update_general}
x^+ = x + \mS h,\qquad \mS \in \R^{d\times \tau(\mS)}, \qquad h\in \R^{\tau(\mS)},
\end{equation}
where $\mS$ is a thin matrix, typically with a negligible number of columns compared to the dimension (i.e., $\tau(\mS)\ll d$). That is, they move from $x$ to $x^+$ along the subspace spanned by the columns of $\mS$.

In these methods,  the subspace matrix $\mS$ is typically chosen first, followed by the determination of the parameters $h$ which define the linear combination of the columns determining the update direction. Several different rules have been proposed in the literature for choosing the matrix $\mS$, including greedy,  cyclic and randomized rules. In this work we consider a \textit{randomized} rule. In particular, we assume that $\mS$ is sampled from an arbitrary but fixed distribution $\cD$ restricted to requiring that $\mS$ be of full column rank\footnote{It is rather simple to extend our results to matrices $\mS$ which are column-rank deficient. However, this would introduce a rather heavy notation burden which we decided to avoid for the sake of clarity and readability.} with probability one. 

Once $\mS\sim \cD$ is sampled, a rule for deciding  the stepsize $h$ varies from algorithm to algorithm, but is mostly determined by the underlying \emph{oracle model} for information access to function $f$. For instance, first-order methods require  access to the subspace gradient $\nabla_{\mS} f(x) \eqdef  \mS^\top \nabla f(x)$, and are relatively well studied~\cite{rcdm, stich2013optimization, richtarik2014iteration, wright2015coordinate, kozak2019stochastic}. At the other extreme are  variants performing a full subspace minimization, i.e.,  $f$ is minimized over the affine subspace given by $\{ x+\mS h \, | \, h\in \R^{\tau(\mS)}\}$~\cite{chang2008coordinate}.  In particular, in this chapter we are interested in the \emph{second-order} oracle model; i.e. we claim access both to the subspace gradient $\nabla_{\mS} f(x)$ and the subspace Hessian $\nabla^2_{\mS} f(x)\eqdef \mS^\top \nabla^2f(x)\mS$.

\section{Contributions}

We now summarize our contributions:

\begin{itemize}
\item {\bf New 2nd order subspace method.} We propose a new stochastic subspace method---Stochastic Subspace Cubic Newton ({\tt SSCN})---constructed by minimizing an oracle-consistent global upper bound on the objective $f$ in each iteration (Section~\ref{sec:sscn_Alg}). This bound is formed using both the subspace gradient and the subspace Hessian at the current iterate and relies on Lipschitzness of the subspace Hessian. 

\item {\bf Interpolating global rate.} We prove (Section~\ref{sec:sscn_global}) that {\tt SSCN} enjoys a global convergence rate that interpolates between the rate of stochastic {\tt CD} and the rate of cubic regularized Newton as one varies the expected dimension of the subspace, $\E{\tau(\mS)}$. 

\item {\bf Fast local rate.} Remarkably, we establish a local convergence bound for {\tt SSCN}  (Section~\ref{sec:sscn_local}) that matches the rate of stochastic subspace descent ({\tt SSD})~\cite{gower2015randomized} applied to solving the problem \begin{equation} \label{eq:sscn_quadratic_optimum}
 \min_{x \in \R^d}\; \frac12 (x-x^*)^\top \nabla^2f(x^*)(x-x^*),
\end{equation} where $x^*$ is the solution of~\eqref{eq:sscn_problem}. Thus, {\tt SSCN} behaves as {\em if} it had access to a perfect second-order model of $f$ at the optimum, and was given the (intuitively much simpler) task of minimizing this model instead. 
Furthermore, note that {\tt SSD}~\cite{gower2015randomized} applied to minimize a convex quadratic can be interpreted as doing an exact subspace search in each iteration, i.e., it minimizes the objective exactly along the active subspace~\cite{richtarik2017stochastic}. Therefore, the local rate of {\tt SSCN} matches the rate of the greediest strategy for choosing $h$ in the active subspace, and as such, this rate is the best one can hope for a method that does not incorporate some form of acceleration.

\item {\bf Special cases.} We discuss in Section~\ref{sec:sscn_special} how {\tt SSCN} reduces to several existing stochastic second-order methods in special cases, either recovering the best known rates, or improving upon them. This includes  {\tt  SDSA}~\cite{sda}, {\tt CN}~\cite{griewank1981modification, nesterov2006cubic} and {\tt RBCN}~\cite{doikov2018randomized}. However, our method is more general and hence allows for more applications.

\end{itemize}

We discuss more remotely related literature in Section~\ref{sec:sscn_related_literature}. We now give a simple example of our setting.

\begin{example}[Coordinate subspace setup]
Let $\mI^d\in \R^{d\times d}$ be the identity and let $S$ be a random subset of $\{1, 2, \dots, d \}$. Given that $\mS = \mI^d_{(:,S)}$ with probability 1, the oracle model reveals $(\nabla f(x))_S$ and $(\nabla^2 f(x))_{(S,S)}$. Therefore, we have access to a random block of partial derivatives of $f$ and a block submatrix of its Hessian, both corresponding to the subset of indices $S$. Furthermore, the rule~\eqref{eq:sscn_update_general} updates a subset $S$ of coordinates only. In this setting, our method is a new \textit{second-order coordinate subspace descent} method.
\end{example}

\section{Preliminaries \label{sec:sscn_preliminaries}}

Throughout the chapter, we assume that  $f$ is convex, twice differentiable, and sufficiently smooth and that $\psi$ is convex, albeit possibly non-differentiable, as the next assumption states.\footnote{We will also require separability of $\psi$; see Section~\ref{sec:sscn_setup}.}

\begin{assumption}\label{as:sscn_conv_lipschitz}
Function $f: \R^d \to \R$ is convex and twice differentiable with $M$-Lipschitz continuous Hessian. Function $\psi: \R^d \to \R \cup \{ +\infty \}$ is  proper closed and convex.
\end{assumption}

We always assume that a minimum of $F$  exists and by  $x^*$ denote any of its minimizers. We let $F^{*} \eqdef F(x^{*})$.

Since our method always takes steps along random subspaces spanned by the columns of $\mS \in \R^{d\times \tau(\mS)}$,  it is reasonable to define the Lipschitzness of the Hessian over the range of $\mS$:\footnote{By  $\|x\| \eqdef \la x, x \ra^{1/2}$we denote the standard Euclidean norm.}
\begin{equation}\label{eq:sscn_MS_def}
M_{\mS}  \eqdef \max_{x\in \R^d} \max_{\substack{  \;\; h \in \R^{\tau(\mS)}, \\ h \not= 0 } } \frac{ |\nabla^3 f(x)[\mS h]^3| }{\|\mS h\|^3}\, .
\end{equation}

As the next lemma shows,  the maximal value of $M_{\mS}$ for any $\mS$ of width $\tau$ can be up to $(\frac{d}{\tau})^{\frac32}$ times smaller than $M$ and this will lead to  a tighter approximation of the objective. 
\begin{lemma} \label{lem:sscn_sharpness}
We have $
M \geq \max_{\tau(\mS) = \tau }M_{\mS}.
$
Moreover, there is a problem where $ \max_{\tau(\mS) = \tau}M_{\mS} =\left(\frac{\tau}{d}\right)^{\frac32} M$. Lastly, if $\Range{\mS} = \Range{\mS'}$, then  $M_{\mS} = M_{\mS'}$.
\end{lemma}

The next lemma provides a direct motivation for our algorithm. It gives a global upper bound on the objective over a random subspace, given the first and second-order information at the current point. 

\begin{lemma}\label{lem:sscn_ub}
Let $x\in \R^d$, $\mS\in \R^{d\times \tau(\mS)}$, $h\in \R^{\tau(\mS)}$ and $x^+$ be as in \eqref{eq:sscn_update_general}. Then
\begin{equation}
\left | f(x^+) -  f(x)  - \langle \nabla_{\mS} f(x), h \rangle - \frac12 \la \nabla^2_{\mS}f(x) h, h\ra \right |   \leq  \frac{M_{\mS}}{6} \| \mS h\|^3. \label{eq:sscn_coordinate_ub}
\end{equation}
As a consequence, we have
\begin{eqnarray} 
F(x^+) \leq  F(x) + T_{\mS}(x,h),
\label{eq:sscn_coordinate_ub_full}
\end{eqnarray}
where
$T_{\mS}(x,h) \eqdef \langle \nabla_{\mS} f(x), h \rangle + \frac12 \la  \nabla^2_{\mS}  f(x) h, h \ra
+  \frac{M_{\mS}}{6} \| \mS h\|^3 + \psi(x+\mS h).$
\end{lemma}

We shall also note that for function $\psi$
we require \textit{separability} with respect to the sampling distribution
(see Definition~\ref{def:separable} and the corresponding
Assumption~\ref{as:sscn_separable} in Section~\ref{sec:sscn_setup}).

For better orientation throughout the chapter, we provide a table of frequently used notation in the Appendix.

\section{The {\tt SSCN} algorithm} \label{sec:sscn_Alg}

For a given $\mS$ and current iterate $x^k$, it is a natural idea to choose $h$ as a minimizer of the upper bound~\eqref{eq:sscn_coordinate_ub_full} in $h$ for $x=x^k$, and subsequently set $x^{k+1} = x^+$ via~\eqref{eq:sscn_update_general}.  Note that we are choosing $\mS$ randomly according to a fixed distribution $\cD$ (with a possibly random number of columns). We have just described {\tt SSCN}---Stochastic Subspace Cubic Newton---formally stated as  Algorithm~\ref{alg:sscn_crcd}.
 
\begin{algorithm}[!h]
\begin{algorithmic}[1]
\State \textbf{Initialization:} $x^0$, distribution $\cD$ of random matrices with $d$ rows and full column rank
\For{$k =  0, 1, 2,\dots$}
\State Sample $\mS$ from distribution $\cD$
\State $h^k =  \argmin_{h \in \R^{\tau(\mS)}}  T_{\mS}(x^k,h)$
\State Set $x^{k+1}=x^k + \mS h^k$ \label{eq:sscn_x_update_CRCD}
\EndFor
\end{algorithmic}
\caption{{\tt SSCN}: Stochastic Subspace Cubic Newton}
\label{alg:sscn_crcd}
\end{algorithm}

\begin{remark}\label{rem:monotonic}
Inequality~\eqref{eq:sscn_coordinate_ub_full} becomes an equality with $h=0$. As a consequence, we must have $F(x^{k+1}) \leq F(x^k)$, and thus the sequence $\{ F(x^k) \}_{k \geq 0}$ is non-increasing.
\end{remark}

\subsection{Solving the subproblem \label{sec:sscn_solving}}
Algorithm~\ref{alg:sscn_crcd} requires $T_{\mS}$ to be minimized in $h$ each iteration. As this operation does not have a closed-form solution in general, it requires an optimization subroutine itself of a possibly non-trivial complexity, which we discuss here. 

\paragraph{The subproblem without $\psi$.}
Let us now consider the case when  $\psi(x) \equiv 0$ in which our problem~\eqref{eq:sscn_problem} does not contain any nondifferentiable components. Various techniques for minimizing regularized quadratic functions were developed  during the development of
Trust-region methods (see~\cite{conn2000trust}), and applied
to Cubic regularization in~\cite{nesterov2006cubic}. The classical approach consists in performing some diagonalization of the matrix $\nabla^2_{\mS} f(x)$ first, by computing the \textit{eigenvalue} or \textit{tridiagonal} decomposition, which costs $\cO(\tau(\mS)^3)$ arithmetical operations.
Then, to find the minimizer, it merely remains to solve a one-dimensional nonlinear equation (this part can be done by $\tilde{\cO}(1)$ iterations of the one-dimensional Newton method, with a linear cost per step). More details and analysis of this procedure can be found
in~\cite{gould2010solving}.

The next example gives a setting in which  an explicit formula for the minimizer of $T_{\mS}$ can be deduced.
\begin{example}\label{ex:explicit} 
Let $e_i$ be the $i$th unit basis vector in $\R^d$. If $\mS\in \{e_1, \dots, e_d\}$ with probability 1 and $\psi(x) = 0$, the update rule can be written as
$x^{k+1} = x^k - \alpha_i^k e_i,$ with
\[
\alpha_i^k = \frac{2\nabla_i f(x^{k})}{\nabla^2_i f(x^{k}) + \sqrt{\left(\nabla^2_{ii} f(x^{k})  \right)^2 + 2M_{e_i}|\nabla_i f(x^{k}) |} },
\]
thus the cost of solving the subproblem is $\cO(1)$.
\end{example}

\paragraph{Subproblem with simple $\psi$.} In some scenarios, minimization of $T_{\mS}$  can be done using a simple algorithm if $\psi$ is simple enough. We now give an example of this.

\begin{example} 
	If $\mS\in \{e_1, \dots, e_d\}$ with probability 1, the subproblem can be solved using a binary search given that the evaluation of $\psi$ is cheap. In particular, if we can evaluate $\psi(x^k+ \mS h)- \psi(x^k)$ in $\tilde{\cO}(1)$, the cost of solving the subproblem will be $\tilde{\cO}(1)$.
\end{example}

\paragraph{The subproblem with general $\psi$.} In the case of general regularizers, recent line of work \cite{carmon2019gradient} explores to the use
of \textit{first-order} optimization methods (Gradient Methods) for
computing an approximate minimizer of $T_{\mS}$.
We note that the backbone of such Gradient Methods is an implementation of the following operation 
(for a any given vector $b \in \R^{\tau(\mS)}$, and positive scalars $\alpha, \beta$):
\[
 \argmin \limits_{h \in \R^{\tau(\mS)}} \la b, h \ra + \frac{\alpha}{2}\|\mS h\|^2 + \frac{\beta}{3}\|\mS h\|^3 + \psi(x^k + \mS h).
\]

To the best of our knowledge, the most efficient gradient method is the Fast Gradient Method ({\tt FGM})\cite{nesterov2019inexact}, achieving an $\cO(1/k^6)$ convergence rate. However, {\tt FGM} can deal with any $\psi$ as long as the above subproblem is cheap to solve. We shall also note that gradient methods do not require a storage of $\nabla^2_{\mS} f(x)$; but rather iteratively access partial Hessian-vector products $\nabla^2_{\mS} f(x) h$.

\paragraph{Line search.}

Note that in Algorithm~\ref{alg:sscn_crcd} we use the Lipschitz constants $M_{\mS}$ 
of the subspace Hessian (see Definition~\eqref{eq:sscn_MS_def}) as the regularization parameters. In many application, $M_{\mS}$ can be estimated cheaply (see Section~\ref{sec:sscn_applications}). In general, however, $M_\mS$ might be unknown or hard to estimate. 
In such a case, one might use a simple one-dimensional search on each iteration: multiply the estimate of $M_\mS$ by the factor of two until the bound~\eqref{eq:sscn_coordinate_ub_full} is satisfied, and divide it by two at the start of each iteration. Note that the average number of such line search steps per iteration can be bounded by two (see~\cite{grapiglia2017regularized} for the details).

\subsection{Special cases} \label{sec:sscn_special}
There are several scenarios where {\tt SSCN} becomes an already known algorithm.  We list them below.

\paragraph{Quadratic minimization.}  If $M=0$ and $\psi = 0$, {\tt SSCN} reduces to the stochastic dual subspace ascent ({\tt SDSA}) method~\cite{sda}, first analyzed in an equivalent primal form as a \emph{sketch-and-project} method in~\cite{gower2015randomized}. In such a case, {\tt SSCN} performs both first-order, second-order updates, and exact minimization over a subspace at the same time due to the quadratic structure of the objective~\cite{richtarik2017stochastic}. The convergence rate we provide in Section~\ref{sec:sscn_local} exactly matches the rate of sketch-and-project as well. As a consequence, we recover a subclass of matrix inversion algorithms~\cite{pseudoinverse} together with stochastic spectral (coordinate) descent~\cite{kovalev2018stochastic} along with their convergence theory.

\paragraph{Full-space method.} If $\mS=\mI^d$ with probability 1, {\tt SSCN} reduces to cubically regularized Newton ({\tt CN})~\cite{griewank1981modification, nesterov2006cubic}. In this case, we recover both existing global convergence rates and superlinear local convergence rates.  

\paragraph{Separable non-quadratic part of $f$.} The {\tt RBCN} method~\cite{doikov2018randomized} aims to minimize~\eqref{eq:sscn_problem} with $f(x) = g(x) + \phi(x)$, where $g,\phi$ are both convex, and $\phi$ is separable.\footnote{Separability is defined in Section~\ref{sec:sscn_setup}.}
They assume that $\nabla^2 g(x) \preceq \mA \in \R^{d\times d}, \forall x \in \R^d$, while $\phi$ has Lipschitz continuous Hessian. In each iteration, {\tt RBCN} constructs an upper bound on the objective using first-order information from $g$ only. This is unlike {\tt SSCN}, which uses second-order information from $g$. In a special case when $\nabla^2 g(x) = \mA$ for all $x$, {\tt SSCN} and {\tt RBCN} are identical algorithms. However, {\tt RBCN} is less general: it requires separable $\phi$, and thus does not cover some of our applications, and takes directions along coordinates only. Further, the rates we provide are better even in the setting where the two methods coincide ($\nabla^2 g(x)=\mA$). The simplest way to see that is by looking at local convergence -- {\tt RBCN} does not achieve the local convergence rate of block {\tt CD} to minimize~\eqref{eq:sscn_quadratic_optimum}, which is the best one might hope for.

Besides these particular cases, for a general twice-differentiable $f$, {\tt SSCN} is a new second-order method.
 
\section{Related literature} \label{sec:sscn_related_literature}

Several methods in the literature are related to {\tt SSCN}. We briefly review them below.

\begin{itemize}
\item \emph{Cubic regularization of Newton method} was proposed first in~\cite{griewank1981modification},
and received substantial attention after the work of Nesterov~\cite{nesterov2006cubic}, where its global complexity guarantees were established.
During the last decade, there was a steady increase of research in second-order methods, 
discovering Accelerated~\cite{nesterov2008accelerating,monteiro2013accelerated},
Adaptive~\cite{cartis2011adaptive1,cartis2011adaptive2}, and
Universal~\cite{grapiglia2017regularized,grapiglia2019accelerated,doikov2019minimizing} 
schemes (the latter ones are adjusting automatically  to the smoothness properties of the objective).

\item There is a vast literature on \emph{first-order coordinate descent ({\tt CD})} methods. While {\tt CD} with $\tau=1$ is consistently the same method within the literature~\cite{rcdm, richtarik2014iteration, wright2015coordinate}, there are several ways to deal with $\tau>1$. The first approach constructs a separable upper bound on the objective (in expectation) in the direction of a random subset of coordinates~\cite{qu2016coordinate1, qu2016coordinate2}, which is minimized to obtain the next iterate. The second approach---{\tt SDNA}~\cite{sdna}---works with a tighter non-separable upper bound. {\tt SDNA} is, therefore, more costly to implement but requires a smaller number of iterations to converge. The literature on first-order subspace descent algorithms is slightly less rich, the notable examples are random pursuit~\cite{stich2013optimization} or stochastic subspace descent~\cite{kozak2019stochastic}.

\item \emph{Randomized subspace Newton} ({\tt RSN})~\cite{gower2019rsn} is a method of the form $$x^{k+1} = x^k - \hat{L}^{-1} \mS \left(\nabla^2_\mS f(x^k)\right)^{-1} \nabla_{\mS} f(x^k)$$ for some specific fixed $\hat{L}$. In particular, it can be seen as a method minimizing the following upper bound on the function, which follows from their assumption:
\begin{align*}
 & h^k = \argmin_h\, \langle \nabla_{\mS} f(x^k), h  \rangle + \frac{\hat{L}}{2}  \la \nabla^2_{\mS}  f(x^k) h, h \ra. 
\end{align*}
This is followed by an update over the subspace: $ x^{k+1} = x^k + \mS h^k$. 
Since both {\tt RSN} and {\tt SSCN} are analyzed under different assumptions, the global linear rates are not directly comparable. However, the local rate of {\tt SSCN} is superior to {\tt RSN}. We shall also note that {\tt RSN} is a stochastic subspace version of a method from~\cite{karimireddy2018global}.

\item \emph{Subsampled Newton} (SN) methods~\cite{byrd2011use, erdogdu2015convergence, xu2016sub, roosta2019sub} and \emph{subsampled cubic regularized Newton methods}~\cite{kohler2017sub, xu2017newton, wang2018stochastic} and \emph{stochastic (cubic regularized) Newton methods}~\cite{tripuraneni2018stochastic,cartis2018global, kovalev2019stochastic} are  stochastic second-order algorithms to tackle finite sum minimization. Their major disadvantage is a requirement of an immense sample size, which makes them often impractical if used as theory prescribes. A notable exception that does not require a large sample size was recently proposed in~\cite{kovalev2019stochastic}. However, none of these methods are directly comparable to {\tt SSCN} as they are not subspace descent methods, but rather randomize over data points (or sketch the Hessian from ``inside''~\cite{pilanci2017newton}). 

\end{itemize}

\section{Global complexity bounds}
\label{sec:sscn_global}

We first start presenting the global complexity results of {\tt SSCN}.

\subsection{Setup~\label{sec:sscn_setup}}

Throughout this section, we require some kind of uniformity of the distribution $\cD$ over subspaces given by $\mS$. In particular, we require $\mZ= \mZ(\mS) \eqdef \mS \left(\mS^\top \mS \right)^{-1} \mS^\top$, the projection matrix onto the range of $\mS$, to be a scalar multiple of identity matrix in expectation. 
\begin{assumption}\label{as:sscn_uniform}
$\exists \tau>0$ such that distribution  $\cD$ satisfies
\begin{equation} \label{eq:sscn_uniform_sampling}
\E{\mZ} = \frac{\tau}{d} \mI^d.
\end{equation}
\end{assumption}
A direct consequence of Assumption~\ref{as:sscn_uniform} is that $\tau$ is an expected width of $\mS$, as the next lemma states.
\begin{lemma}\label{lem:sscn_exp_size}
If Assumption~\ref{as:sscn_uniform} holds, then $\E{\tau(\mS)} = \tau$.
\end{lemma}

As mentioned before, the global complexity results are interpolating between convergence rate of (first-order) {\tt CD} and (global) convergence rate of Cubic Newton. However, first-order {\tt CD} requires Lipschitzness of gradients, and thus we will require it as well.
\begin{assumption}\label{as:sscn_smooth}
Function $f$ has $L$-Lipschitz continuous gradients, i.e. $\nabla^2 f(x) \preceq L \mI^d$ for all $x\in \R^d$.
\end{assumption}

We will also need an extra assumption on $\psi$. It is well known that proximal (first-order) {\tt CD} with fixed step size does not converge if $\psi$ is not separable -- in such case, even if $f(x^k) = f(x^*)$ we might have $f(x^{k+1})> f(x^*)$. Therefore, we might not hope that {\tt SSCN} will converge without additional assumptions on $\psi$. Informally speaking, separability of $\psi$ with respect to directions given by columns of $\mS$ is required. 
To define it formally, let us introduce first the notion of a separable set.
\begin{definition}\label{def:separable_set}
Set $Q \subseteq \R^d$ is called $D$-separable, if $\forall x, y \in Q, \mS \in D$:
\[
\mZ x + (\mI^d- \mZ)y  \in  Q.
\]
\end{definition}
Let $e\in \R^d$ be the vector of all ones. Then, for arbitrary functions, we have
\begin{definition}\label{def:separable}
Function $\phi: \R^d \rightarrow \R \cup \{ +\infty \}$ is $D$-separable if
$\dom \phi$ is $D$-separable, and there is map $\phi': \dom \phi \rightarrow \R^d$ such that
\begin{enumerate}
\item  $\forall x \in \dom\phi:\; \phi(x) = \langle \phi'(x), e \rangle$,
\item $\forall x,y \in \dom\phi, \mS\in D: \; \phi'(\mZ x + (\mI^d- \mZ)y) = \mZ\phi'(x) + (\mI^d-\mZ) \phi'(y)$.		
\end{enumerate}
\end{definition}

\begin{example}
If $D$ is a set of matrices whose columns are standard basis vectors, $D$-separability reduces to classical (coordinate-wise) separability.
\end{example}

\begin{example}
If $D$ is set of matrices which are column-wise submatrices of orthogonal $\mU$, $D$-separability of $\phi$ reduces to classical coordinate-wise separability of $\phi(\mU^\top x)$.
\end{example}

\begin{example}
$\phi(x) = \frac12\| x\|^2$ is $D$-separable for any~$D$.
\end{example}

\begin{assumption}\label{as:sscn_separable}
Function $\psi$ is $\Range{\cD}$-separable.
\end{assumption}

We are now ready to present the convergence rate of {\tt SSCN}. 

\subsection{Theory}

First, let us introduce the critical lemma from which the main global complexity results are derived. Our first lemma gives a bound on the expected progress after a single step of {\tt SSCN}.

\begin{lemma}\label{lem:sscn_keylemma}
Let Assumptions~\ref{as:sscn_conv_lipschitz},~\ref{as:sscn_uniform},~\ref{as:sscn_smooth} and~\ref{as:sscn_separable} hold.
Then, for every $k \geq 0$ and $y \in \R^d$ we have
\begin{equation} \label{GlobalUpper}
 \E{F(x^{k + 1}) \, | \, x^k } 
 \leq 
\left( 1 - \frac{\tau}{d} \right) F(x^k) + \frac{\tau}{d} F(y)   + \frac{\tau}{d} \left(  \frac{d - \tau}{d } \frac{L}{2}\|y - x^k\|^2 + \frac{M}{3}\|y - x^k\|^3
\right).
\end{equation}
\end{lemma}

Now we are ready to present global complexity results for the general class of convex functions. The convergence rate is obtained by summing ~\eqref{GlobalUpper} over the different iterations $k$, and with a specific choice of $y$. 

\begin{theorem}\label{thm:sscn_global_weakly}
	Let Assumptions~\ref{as:sscn_conv_lipschitz},~\ref{as:sscn_uniform},~\ref{as:sscn_smooth} and~\ref{as:sscn_separable} hold. Denote
	\beq \label{Rdef}
	R  \Def  \sup\limits_{x \in \R^d} \Bigl\{  \|x - x^{*} \| \; : \; F(x) \leq F(x^{0})  \Bigr\},
	\eeq
	and suppose that $R < +\infty$.
	Then, for every $k \geq 1$ we have
	\beq \label{GlobalConv}
	  \E{ F(x^k) } - F^{*} \leq 
	\frac{d - \tau}{\tau} \cdot \frac{4.5 L R^2}{ k } + \left(\frac{d}{\tau}\right)^2 \cdot \frac{9 M R^3}{k^2}
	+ \frac{F(x^0) - F^{*}}{1 + \frac{1}{4}\left( \frac{\tau}{d} k \right)^{3}}.
	\eeq
	
\end{theorem}

Note that convergence rate of the minibatch version\footnote{Sampling $\tau$ coordinates at a time for objectives with $L$-Lipschitz gradients.} of first-order {\tt CD} is $\cO\bigl( \frac{d}{\tau} \frac{LR^2}{k}\bigr)$. At the same time, (global) convergence rate of cubically regularized Newton method is $\cO\bigl( \frac{MR^3}{k^2}\bigr)$. Therefore, Theorem~\ref{thm:sscn_global_weakly} shows that the global rate of {\tt SSCN} well interpolates between the two extremes,
depending on the sample size $\tau$ we choose.

\begin{remark}
According to estimate~\eqref{GlobalConv}, in order to have $\E{F(x^k)} - F^{*} \leq \varepsilon$,
it is enough to perform
$$
k = \cO\left(  \frac{d - \tau}{\tau} \frac{LR^2}{\varepsilon} + \frac{d}{\tau} \sqrt{\frac{MR^3}{\varepsilon}}
+ \frac{d}{\tau}\left( \frac{F(x^0) - F^{*}}{\varepsilon} \right)^{1/3}  \right)
$$
iterations of {\tt SSCN}.
\end{remark}

Next, we move to the strongly convex case.  

\begin{assumption}\label{as:sscn_sc}
Function $f$ is $\mu$-strongly convex, i.e. $\nabla^2 f(x) \succeq \mu \mI^d$ for all $x\in \R^d$.
\end{assumption}

\begin{remark}
Strong convexity of the objective (assumed for Theorem~\ref{thm:sscn_global_strongly} later) implies:
$ R < +\infty$. Furthermore, due to monotonicity of the sequence $\{ F(x_k) \}_{k \geq 0}$ (see Remark~\ref{rem:monotonic}), we have
$\|x^k - x^{*}\| \leq R$ for all $k$. Therefore, it is sufficient to require Lipschitzness of gradients over the sublevel set, which holds with $L = \lambda_{\max}(\nabla^2 f(x^{*})) + MR$.
\end{remark}

As both extremes cubic regularized Newton (where $\mS = \mI^d$ always) and (first-order) {\tt CD} ($\mS = e_i$ for randomly chosen $i$) enjoy (global) linear rate under strong convexity, linear convergence of {\tt SSCN} is expected as well. At the same time, the leading complexity term should be in between the two extremes. Such a result is established as Theorem~\ref{thm:sscn_global_strongly}.

\begin{theorem}\label{thm:sscn_global_strongly}
Let Assumptions~\ref{as:sscn_conv_lipschitz},~\ref{as:sscn_uniform},~\ref{as:sscn_separable}
and~\ref{as:sscn_sc} hold. 
Then, 
$ \E{F(x^k)} - F^{*} \leq  \varepsilon$,
as long as the number of iterations of {\tt SSCN} is
$$
 k =   \cO\left( \left(  \frac{d - \tau}{\tau} \frac{L}{\mu} + \frac{d}{\tau}\sqrt{\frac{M R}{\mu}} + \frac{d}{\tau}     
\right) 
\cdot \log \frac{F(x^0) - F^{*}}{\varepsilon}  \right) .
$$
\end{theorem}

Indeed, if $\mS =\mI^d$ with probability 1 and $MR\geq \mu$, the leading complexity term becomes $\sqrt{\frac{MR}{\mu}} \log\frac{1}{\varepsilon}$ which corresponds to the global complexity of cubically regularized Newton for
minimizing strongly convex functions~\cite{nesterov2006cubic}.
On the other side of the spectrum if $\mS = e_i$ with probability $\frac1d$, the leading complexity term becomes $\frac{dL}{\mu} \log\frac{1}{\varepsilon}$, which again corresponds to convergence rate of {\tt CD}~\cite{rcdm}. 
Lastly, if $1<\tau<d$, the global linear rate interpolates the rates mentioned above.

\begin{remark}
Proof of Theorem~\ref{thm:sscn_global_strongly} only uses the following consequence of strong convexity:
\beq \label{eq:sscn_SConvex}
\frac{\mu}{2}\|x - x^{*}\|^2  \leq F(x) - F^{*}, \qquad x \in \R^d
\eeq
and thus the conditions of Theorem~\ref{thm:sscn_global_strongly} might be slightly relaxed.\footnote{However, this relaxation is not sufficient to obtain the local convergence results.} For detailed comparison of various relaxations of strong convexity, see~\cite{karimi2016linear}.
\end{remark}

\section{Local convergence \label{sec:sscn_local}}

Throughout this section, assume that $\psi = 0$. We first present the key descent lemma, which will be used to obtain local rates. Let $\mH_{\mS}(x) \eqdef  \nabla^2_{\mS} f(x) +\sqrt{ \frac{M_{\mS}}{2}} \| \nabla_{\mS} f(x)\|^{\frac12} \mI^{\tau(\mS)} $.

\begin{lemma} \label{lem:sscn_decrease}
\begin{equation}\label{eq:sscn_decrease}
f(x^k)-f(x^{k+1})  \geq   \frac12 \| \nabla_{\mS} f(x^k)\|^2_{  \mH^{-1}(x^k)  }.
\end{equation}
\end{lemma}

Before stating the convergence theorem, it will be suitable to define the stochastic condition number of  $\mH_*\eqdef \nabla^2 f(x^*)$:
\begin{equation}\label{eq:sscn_sc_generalized}
\zeta  \eqdef \lambda_{\min} \left( \mH_*^{\frac12}  \E{{\mS} \left(\mS^\top \mH_* \mS \right)^{-1}{\mS}^\top}  \mH_*^{\frac12} \right),
\end{equation}
as it will drive the local convergence rate of {\tt SSCN}.

\begin{theorem}[Local Convergence]\label{thm:sscn_local}
Let Assumptions~\ref{as:sscn_conv_lipschitz},~\ref{as:sscn_sc} hold, and suppose that $\psi = 0$. 
 For any $\varepsilon>0$ there exists $\delta>0$ such that if $F(x^0) - F^* \leq \delta$, we have
\begin{equation}\label{eq:sscn_local_rate}
\E{F(x^k) - F^*} \leq \left( 1- (1-\varepsilon)\zeta  \right)^k\left( F(x^0) - F^*\right)
\end{equation}
and therefore the local complexity of {\tt SSCN} is $\cO\left( \zeta^{-1} \log\frac1\varepsilon\right)$. If further $M=0$ (i.e. $f$ is quadratic), then $\varepsilon =0$ and $\delta = \infty $, and thus the rate is global.

\end{theorem}

The proof of Theorem~\ref{thm:sscn_local} along with the exact formulas for $\varepsilon,\delta$ can be found in Section~\ref{sec:sscn_local_proofs} of the Appendix.

Theorem~\ref{thm:sscn_local} provides a local linear convergence rate of {\tt SSCN}. While one might expect a superlinear rate to be achievable, this is not the case, and we argue that the rate from Theorem~\ref{thm:sscn_local} is the best one can hope for. 

In particular, if $M = 0$, Algorithm~\ref{alg:sscn_crcd} becomes subspace descent for minimizing positive definite quadratic which is a specific instance of sketch-and-project~\cite{gower2015randomized}. However, sketch-and-project only converges linearly -- the iteration complexity of sketch-and-project to minimize $(x-x^*)^\top \mA (x-x^*)$ with $\mA\succ 0$ is 
\[ \cO\left( \left( \mA^{\frac12}  \E{{\mS}\left( \mS^\top \mA \mS \right)^{-1}{\mS}^\top}  \mA^{\frac12} \right)^{-1}\log\frac1\varepsilon\right). \]
 Notice that this rate is matched by Theorem~\ref{thm:sscn_local} in this case.

Next, we compare the local rate of {\tt SSCN} to the rate of {\tt SDNA}~\cite{sdna}. To best of our knowledge, {\tt SDNA} requires the least oracle calls to minimize $f$ among all first-order non-accelerated methods.

\begin{remark}
{\tt SDNA} is a first-order analogue to Algorithm~\ref{alg:sscn_crcd} with $\mS=\mI^d_{(:,S)}$. In particular, given matrix $\mL$ such that $\mL\succeq \nabla^2 f(x) \succ 0 $ for all $x$, the update rule of {\tt SDNA} is $$x^+ = x - \mS \left( \mS^\top \mL \mS \right)^{-1} \nabla_{\mS} f(x),$$ where $\mS = \mI^d_{(:,S)}$ for a random subset of columns $S$.  {\tt SDNA} enjoys linear convergence rate with leading complexity term $\left(\mu \lambda_{\min} \left(\E{\mS(\mS^\top \mL \mS)^{-1}\mS^\top} \right)\right)^{-1}$. The leading complexity term of {\tt SSCN} is $\zeta^{-1}$, and we can bound
\begin{eqnarray*}
\zeta &\geq & \lambda_{\min}\left( \mH_* \right) \lambda_{\min} \left(\E{{\mS} \left(\mS^\top \mH_* \mS\right)^{-1}{\mS}^\top}   \right)
\\
& \geq & \mu \lambda_{\min} \left(\E{\mS \left( \mS^\top \mL \mS \right)^{-1}\mS^\top}\right).
\end{eqnarray*}
Hence, the local rate of {\tt SSCN} is no worse than the  rate of {\tt SDNA}. Furthermore, both of the above inequalities might be very loose in some cases (i.e., there are examples where $\frac{\zeta}{ \mu \lambda_{\min}\E{\mS(\mL_\mS)^{-1}\mS^\top} }$ can be arbitrarily high). Therefore, local convergence rate of {\tt SSCN} might be arbitrarily better than the convergence rate of {\tt SDNA}. As a consequence, the local convergence of {\tt SSCN} is better than convergence rate of any non-accelerated first-order method.\footnote{The rate of {\tt SSCN} and rate of accelerated subspace descent methods are not directly comparable -- while the (local) rate of {\tt SSCN} might be better than rate of ACD, the reverse might happen as well. However, both ACD and {\tt SSCN} are faster than non-accelerated subspace descent.}.
\end{remark}

Lastly, the local convergence rate provided by Theorem~\ref{thm:sscn_local} recovers the superlinear rate of cubic regularized Newton's method, as the next remark states.

 \begin{remark}
 If $\mS = \mI^d$ with probability 1, Algorithm~\ref{alg:sscn_crcd} becomes cubic regularized Newton method~\cite{griewank1981modification, nesterov2006cubic}. For $\mH_*\eqdef \nabla^2 f(x^*)$ we have 
 \begin{align*}
 \zeta = \lambda_{\min} \left( \mH_*^{\frac12} \mH_*^{-1} \mH_*^{\frac12} \right) =  \lambda_{\min} (\mI^d) = 1.
 \end{align*}
As a consequence of Theorem~\ref{thm:sscn_local}, for any $\varepsilon>0$ there exists $\delta>0$ such that if $F(x)-F(x^*)\leq \delta$, we have
 \[
 F(x^+)- F(x^*)\leq \varepsilon(F(x)-F(x^*)).
 \]
 Therefore, we obtain a superlinear convergence rate. 
 \end{remark}

\section{Applications}
\label{sec:sscn_applications}

\subsection{Linear models\label{sec:sscn_linear}}
Consider only $\mS = \mI^d_{(:,S)}$ for simplicity.
 Let 
 \begin{equation} \label{eq:sscn_linear_model}
  F(x) \eqdef \frac1n \sum \limits_{i=1}^n \phi_i( \la a_i, x \ra )+ \psi(x),
 \end{equation}
 and $f(x) \eqdef \frac1n \sum_{i=1}^n \phi_i( \la a_i, x \ra )$ 
  and suppose that $|\nabla^3 \phi_i (y)| \leq c$. Then clearly, for any $h \in \R^d$, we have $$\nabla^3 f(x)[h]^3 =\frac1n \sum_{i=1}^n \nabla^3\phi_i( \la a_i, x\ra ) \la a_i, h \ra^3 .$$ While evaluating \[E\eqdef \max_{\|h\|=1, x} \nabla^3 f(x)[h]^3 \] is infeasible, we might bound it instead via
\begin{align}
 E &\leq \max \limits_{\|h\|=1}  \frac{c}{n}\sum \limits_{i=1}^n| \la a_i, h \ra |^3  \leq    \frac{c}{n}\sum \limits_{i=1}^n\|a_i\|^3, \label{eq:sscn_ndanjdajn}
\end{align}
which means that $M=\frac{c}{n} \sum_{i=1}^n \| a_i\|^3$ is a feasible choice. On the other hand, for $S = \{ j \}$ we have
 \[
\max \limits_{\|h_j\|=1, x} \nabla^3 f(x)[h_j]^3 = \max \limits_{x} \nabla^3 f(x)[e_j]^3 \leq \frac{c}{n} \sum \limits_{i=1}^n|a_{ij}|^3 
 \]
 and thus we might set $M_j =\frac{c}{n} \sum_{i=1}^n | a_{ij}|^3$. The next lemma compares the above choices of $M$ and $M_j$.

\begin{lemma}
We have $M\geq \max_j M_j$. At the same time, there exist vectors $a_i$ that $\max_j M_j  = \frac{M}{d^{\frac32}}$.
\end{lemma}
\begin{proof}
The first part is trivial. For the second part, consider $a_{i,j}\in \{-1,1\}$. 
\end{proof}

\begin{remark}
One might avoid the last inequality from~\eqref{eq:sscn_ndanjdajn} using polynomial optimization; however, this might be more expensive than solving the original optimization problem and thus is not preferable. Another strategy would be to use a line search, see Section~\ref{sec:sscn_solving}.
\end{remark}

Both the formula for $M$ and the formula for $M_j$ require the prior knowledge of $c\geq 0$ such that $|\nabla^3 \phi_i (y)| \leq c$ for all $i$. The next Lemma shows how to compute such $c$ for the logistic regression (binary classification model).

\begin{lemma}
Let $\phi_i(y) = \log(1+e^{-b_iy}),  b_i \in \{-1,1\}$. Then $c = \frac{1}{6\sqrt{3}}$.
\end{lemma}
\begin{proof}
$\nabla^3 \phi_i(y) = -\frac{e^x(e^x-1)}{(1+e^x)^3}$ $\Rightarrow$ $\left| \nabla^3 \phi_i(y) \right|\leq \frac{1}{6\sqrt{3}}$.
\end{proof}

\paragraph{Cost of performing a single iteration}
For the sake of simplicity, let $\tau(\mS)=1$, $\psi\equiv0$. Any {\tt CD} method (i.e. method with update rule~\eqref{eq:sscn_update_general} with $\mS\in \{e_1, \dots, e_d\}$) can be efficiently implemented by memorizing the residuals $\la a_i, x^k \ra$, which is cheap to track since $x^{k+1} - x^k$ is a sparse vector. The overall cost of updating the residuals is $\cO(n)$ while the cost of computing $\nabla_i f(x)$ and $\nabla^2_{i,i} f(x)$ (given the residuals are stored) is $\cO(n)$. Therefore the overall cost of performing a single iteration is $\cO(n)$. Generalizing to $\tau(\mS)=\tau\geq 1$, the overall cost of single iteration of {\tt SSCN} can be estimated as $\cO(n\tau^2 + \tau^3)$,
where $\cO(n\tau^2)$ comes from evaluating subspace gradient and Hessian, while $\cO(\tau^3)$ comes from solving the cubic subproblem.

\subsection{Dual of linear models \label{sec:sscn_dual}}
So far, all results and applications for {\tt SSCN} we mentioned were problems with large model size $d$. In this section we describe how {\tt SSCN} can be efficient to tackle big data problems in some settings. Let $\mA\in \R^{n\times d}$ is data matrix and consider a specific instance of~\eqref{eq:sscn_linear_model} where
\begin{equation}\label{eq:sscn_to_be_dualized}
\ \min_{x\in \R^d} \left\{ F_P(x) \eqdef  \frac1d \sum \limits_{i=1}^n \rho_i(\mA_{(:,i)}x) + \frac{\lambda}{2} \|x \|^2 \right\}.
\end{equation}
where $\rho_i$ is convex for all $i$. One can now formulate a dual problem of~\eqref{eq:sscn_to_be_dualized} as follows:
\begin{equation}\label{eq:sscn_dual_linear}
 \max_{y\in \R^n} \left\{ F_D(y) \eqdef - \frac{1}{2\lambda n^2} \left\| \mA^\top y \right\|^2 -\frac{1}{n} \sum \limits_{i=1}^n\rho_i^*(e_i^\top x) \right\}.
\end{equation}
Note that~\eqref{eq:sscn_dual_linear} is of form~\eqref{eq:sscn_linear_model}, and therefore if $\rho^*_i$ has Lipschitz Hessian, we can apply {\tt SSCN} to efficiently solve it (same as Section~\ref{sec:sscn_linear}). Given the solution of~\eqref{eq:sscn_dual_linear}, we can recover the solution of~\eqref{eq:sscn_to_be_dualized} (duality theory). Thus, {\tt SSCN} can be used as a data-stochastic method to solve finite-sum optimization problems.

The trick described in this section is rather well known. It was  first used in~\cite{sdca}, where  {\tt CD} applied to the problem~\eqref{eq:sscn_dual_linear} ({\tt SDCA}) was shown to be competitive with the variance reduced methods like {\tt  SAG}~\cite{sag}, {\tt SVRG}~\cite{svrg} or {\tt SAGA}~\cite{saga}.


\section{Experiments}
We now numerically verify our theoretical claims. We consider two different objectives: logistic regression (Section~\ref{sec:logreg}) and log-sum-exp (Section~\ref{sec:logsumexp}).

\subsection{Logistic regression\label{sec:logreg}}
Regularized logistic regression is a machine learning model for binary classification. Given data matrix $\mA\in \R^{n\times d}$, labels $b\in \{-1, 1\}^{n}$ and regularization parameter $\lambda \in \R_{+}$, the training corresponds to solving the following optimization problem
\[
f(x)= \frac1n \sum_{i=1}^n \log \left(1+\exp\left(\mA_{i,:}x\cdot  b_i\right) \right)+\frac{\lambda}{2} \| x\|^2.
\]

In the first experiment, we compare {\tt SSCN} to first-order coordinate descent ({\tt CD}) on LIBSVM~\cite{chang2011libsvm}. We consider three different instances of {\tt CD}: {\tt CD} with uniform sampling, {\tt CD} with importance sampling~\cite{rcdm}, and accelerated {\tt CD} with importance sampling~\cite{allen2016even, nesterov2017efficiency}. 

In order to be comparable with the mentioned first-order methods, we consider $\mS\in \{e_1, \dots, e_d\}$ with probability 1  -- the complexity of performing each iteration is about the same for each algorithm now. At the same time, computing $M_{e_i}$ for all $1\leq i \leq d$ is of cost $\cO(nd)$ -- the same cost as computing coordinate-wise smoothness constants for (accelerated) {\tt CD} (see Section~\ref{sec:sscn_linear} for the details). Figure~\ref{fig:sscn_libsvm} shows the result for non-normalized data, while Figure~\ref{fig:sscn_libsvm_normalzied} shows the results for normalized data (thus importance sampling is identical to uniform). 

In all examples, {\tt SSCN} outperformed {\tt CD} with uniform sampling. Moreover, the performance of {\tt SSCN} was always either about the same or significantly better to {\tt CD} with importance sampling. Furthermore, {\tt SSCN} was also competitive to accelerated {\tt CD} with importance sampling (in about half of the cases, {\tt SSCN} was better, while in the other half, ACD was better).

\begin{figure}[!h]
\centering
\begin{minipage}{0.3\textwidth}
  \centering
\includegraphics[width =  \textwidth ]{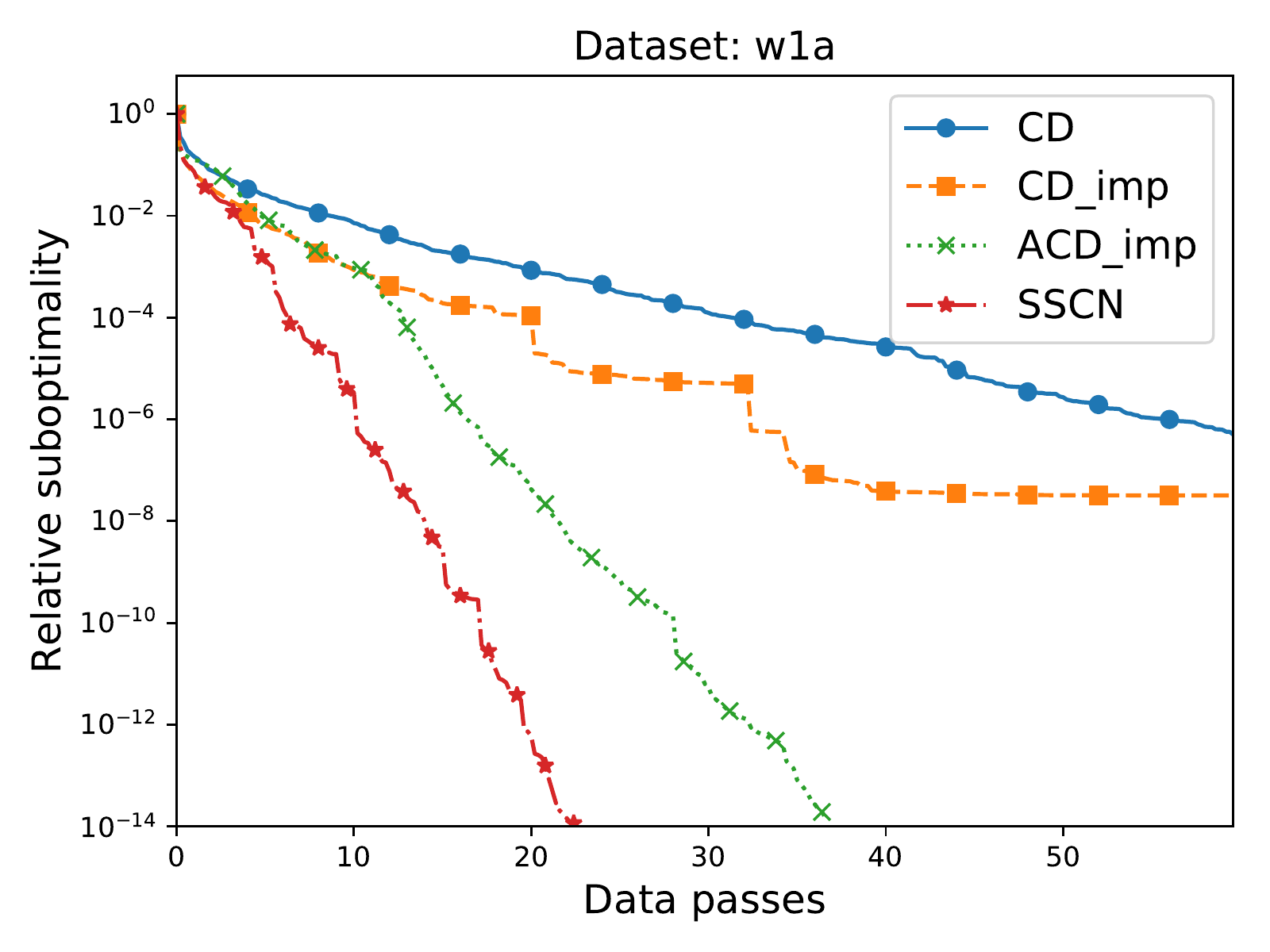}
\end{minipage}%
\begin{minipage}{0.3\textwidth}
  \centering
\includegraphics[width =  \textwidth ]{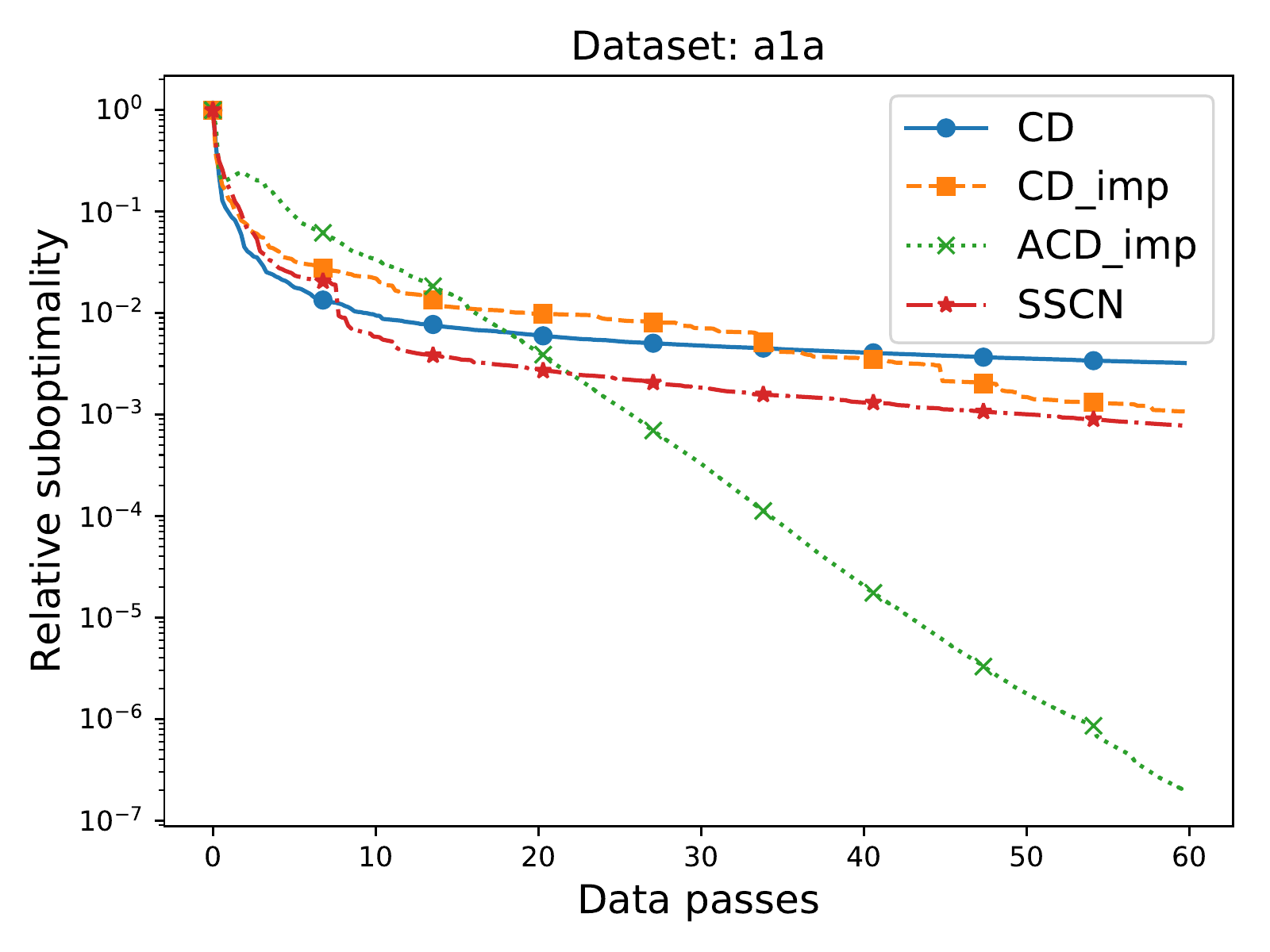}
\end{minipage}%
\begin{minipage}{0.3\textwidth}
  \centering
\includegraphics[width =  \textwidth ]{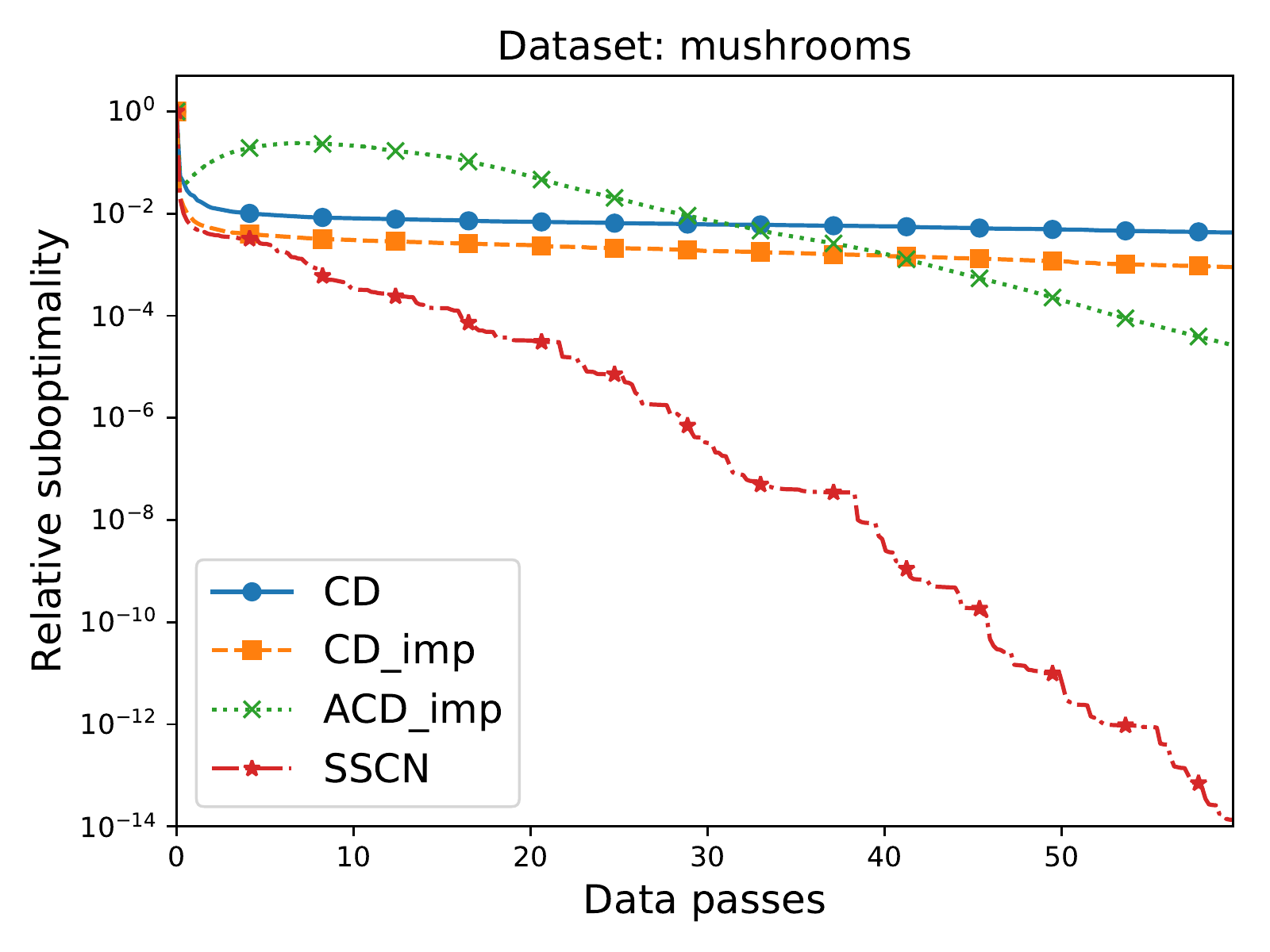}
\end{minipage}%
\\
\begin{minipage}{0.3\textwidth}
  \centering
\includegraphics[width =  \textwidth ]{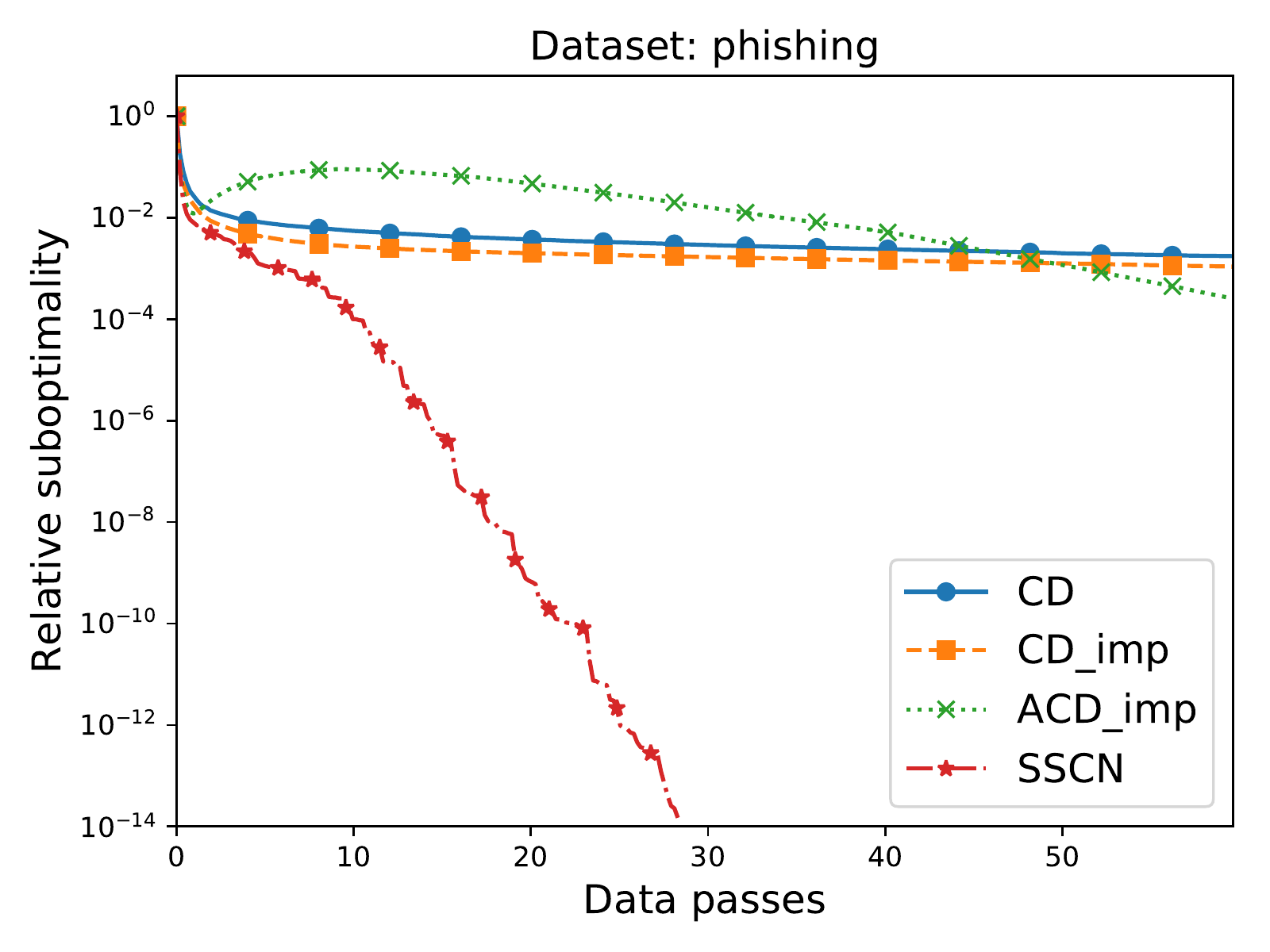}
\end{minipage}%
\begin{minipage}{0.3\textwidth}
  \centering
\includegraphics[width =  \textwidth ]{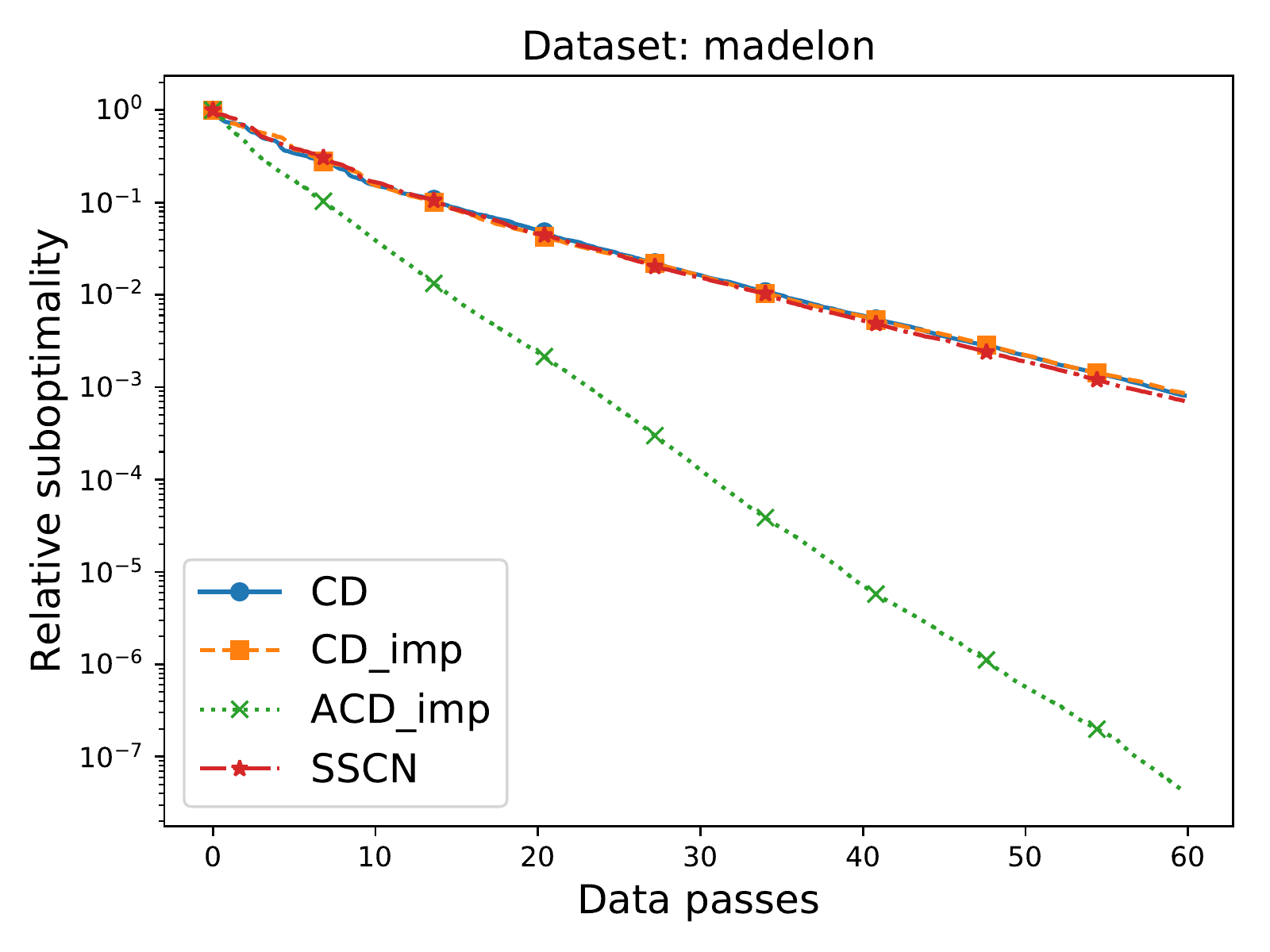}
\end{minipage}%
\begin{minipage}{0.3\textwidth}
  \centering
\includegraphics[width =  \textwidth ]{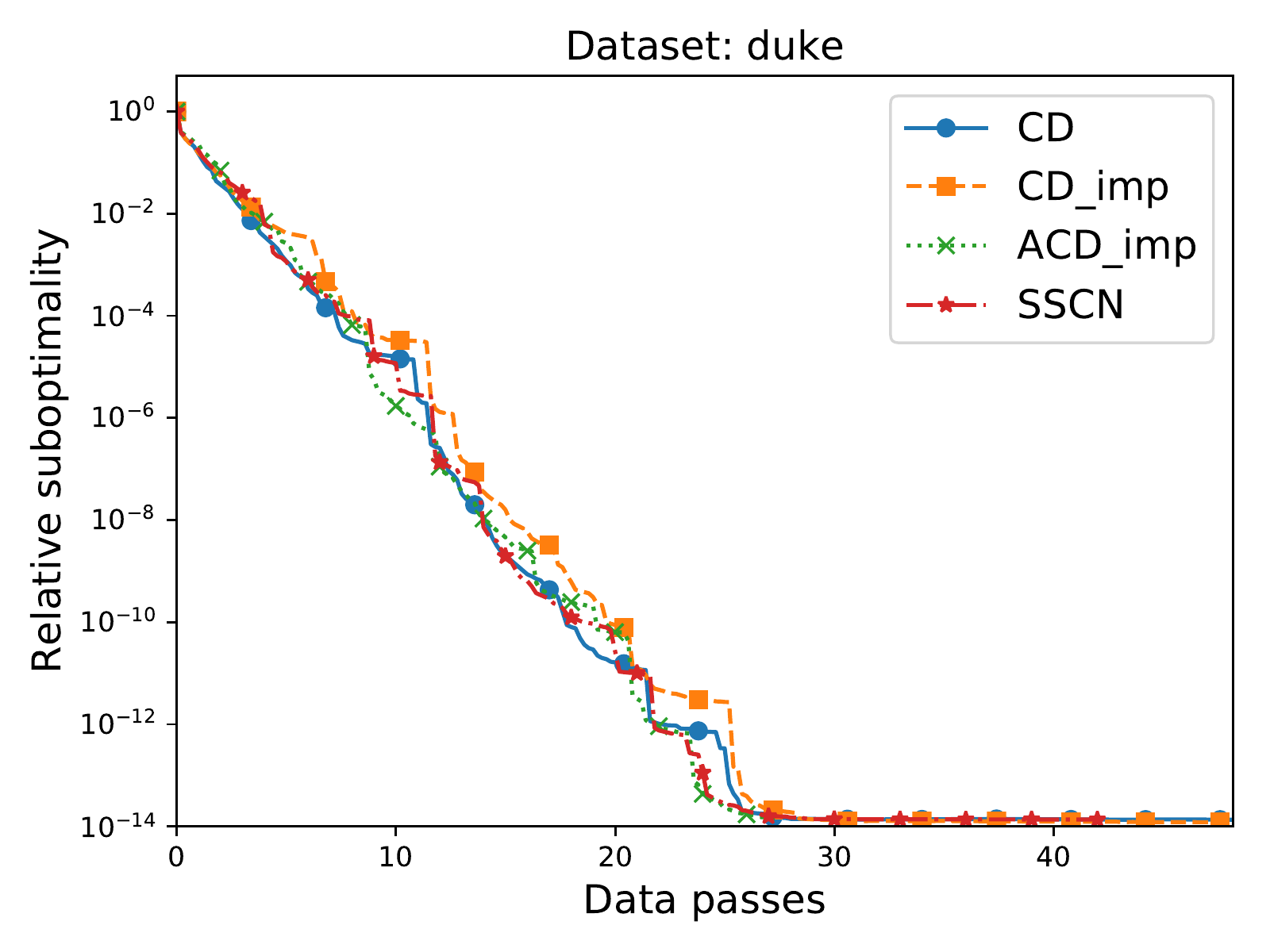}
\end{minipage}%
\\
\begin{minipage}{0.3\textwidth}
  \centering
\includegraphics[width =  \textwidth ]{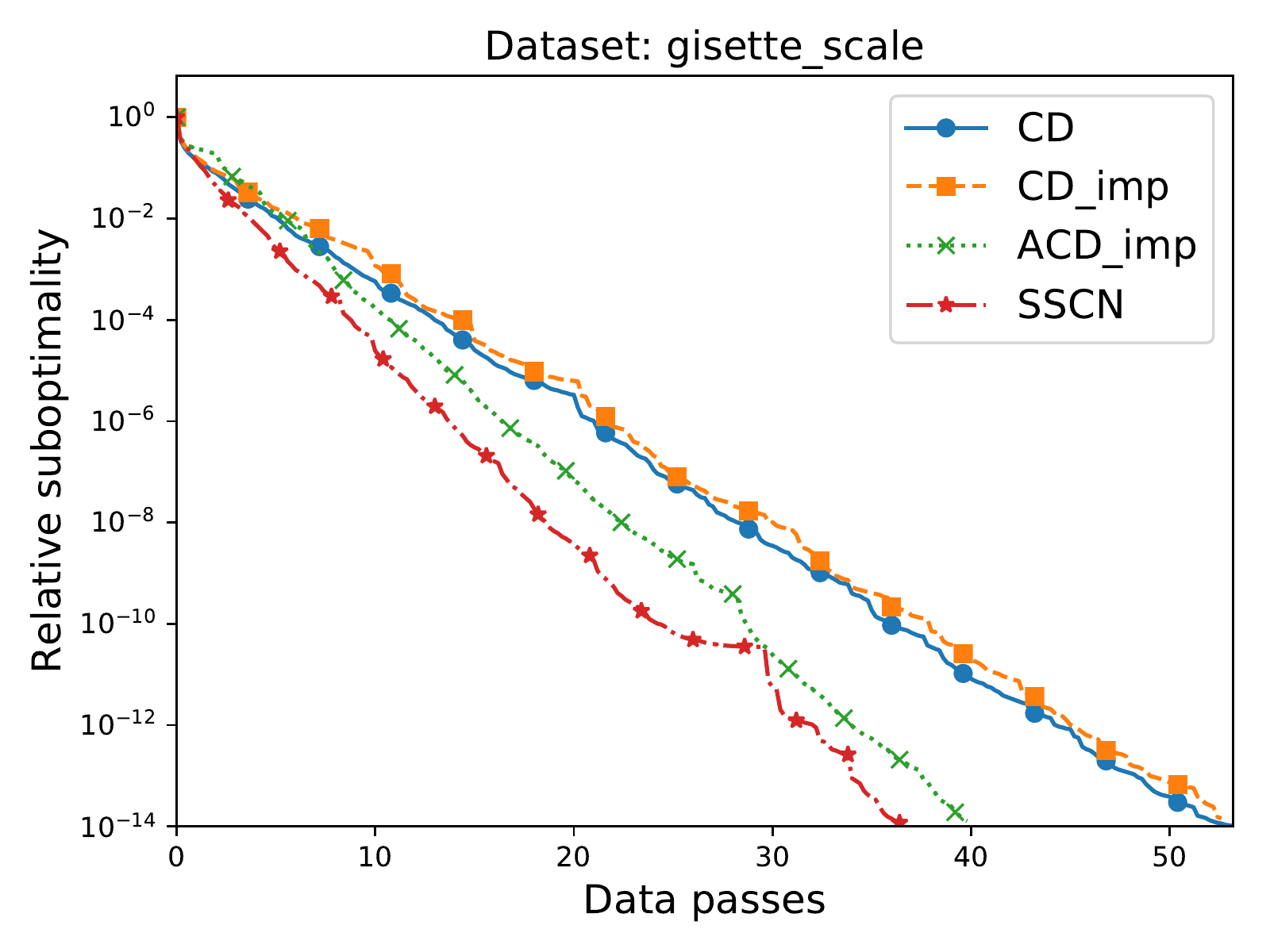}
\end{minipage}%
\begin{minipage}{0.3\textwidth}
  \centering
\includegraphics[width =  \textwidth ]{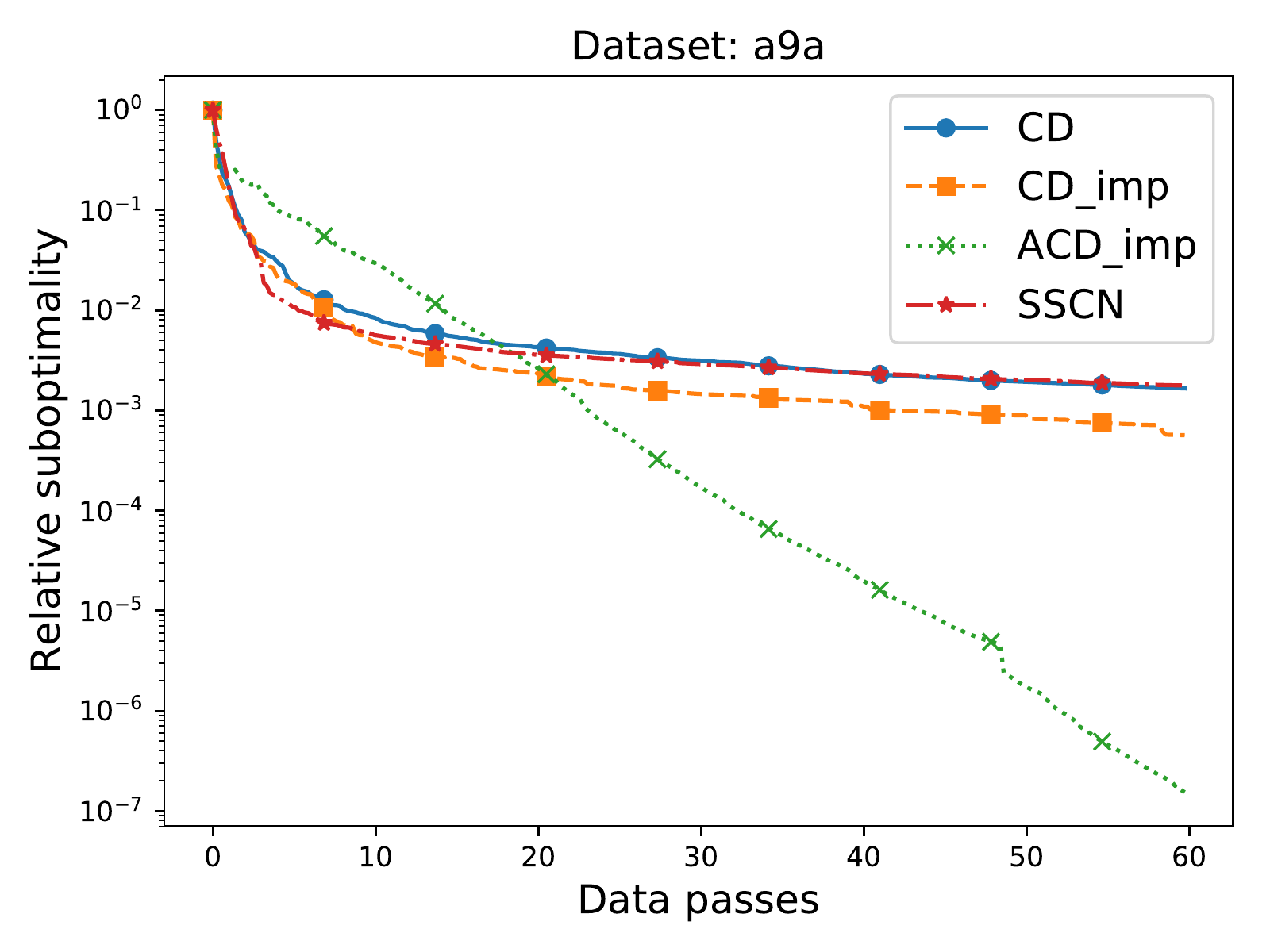}
\end{minipage}%
\begin{minipage}{0.3\textwidth}
  \centering
\includegraphics[width =  \textwidth ]{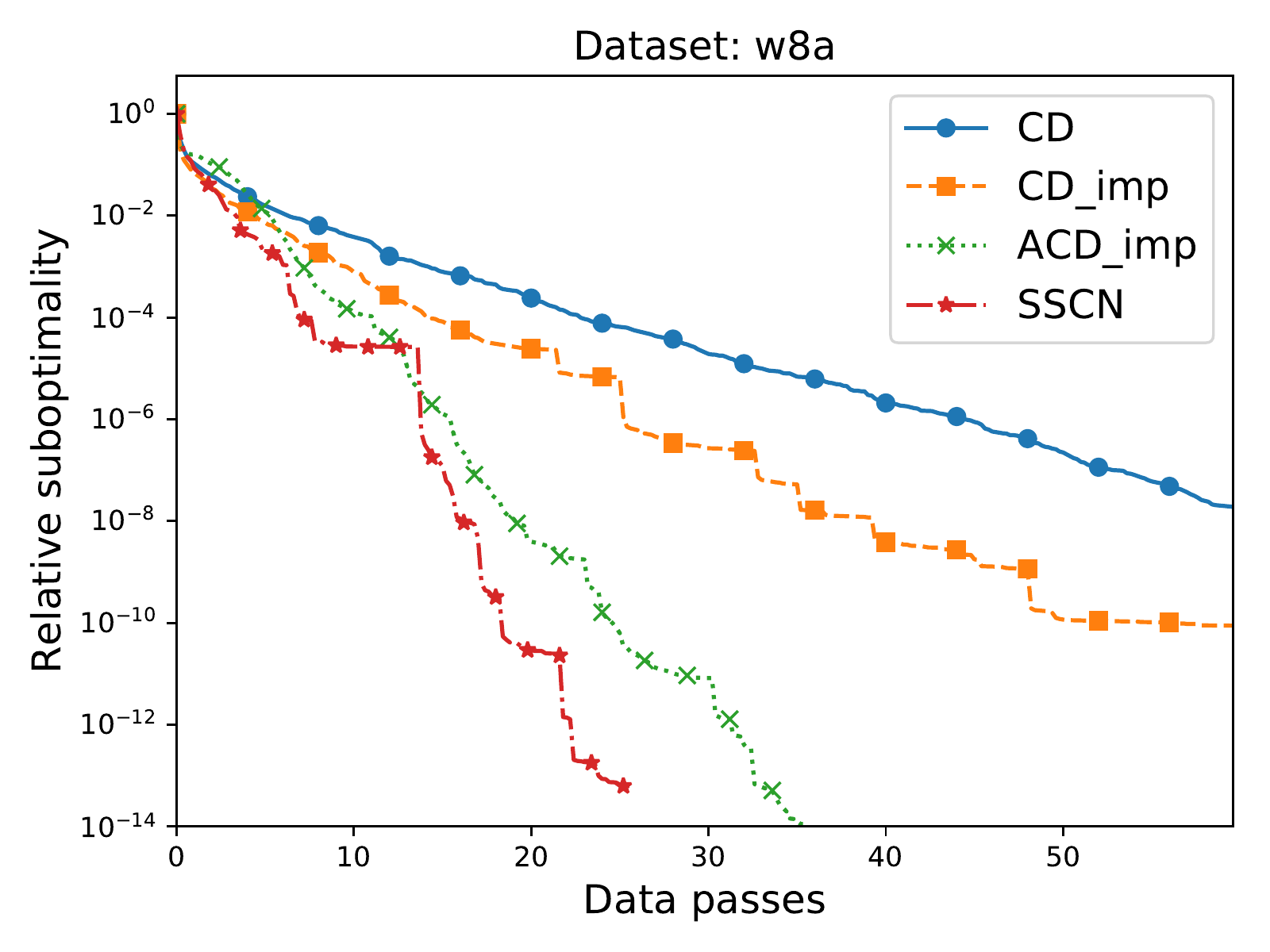}
\end{minipage}%
\caption{Comparison of {\tt CD} with uniform sampling, {\tt CD} with importance sampling, accelerated {\tt CD} with importance sampling and {\tt SSCN} (Algorithm~\ref{alg:sscn_crcd}) with uniform sampling on LibSVM datasets.} 
\label{fig:sscn_libsvm}
\end{figure}

\begin{figure}[!h]
\centering
\begin{minipage}{0.3\textwidth}
  \centering
\includegraphics[width =  \textwidth ]{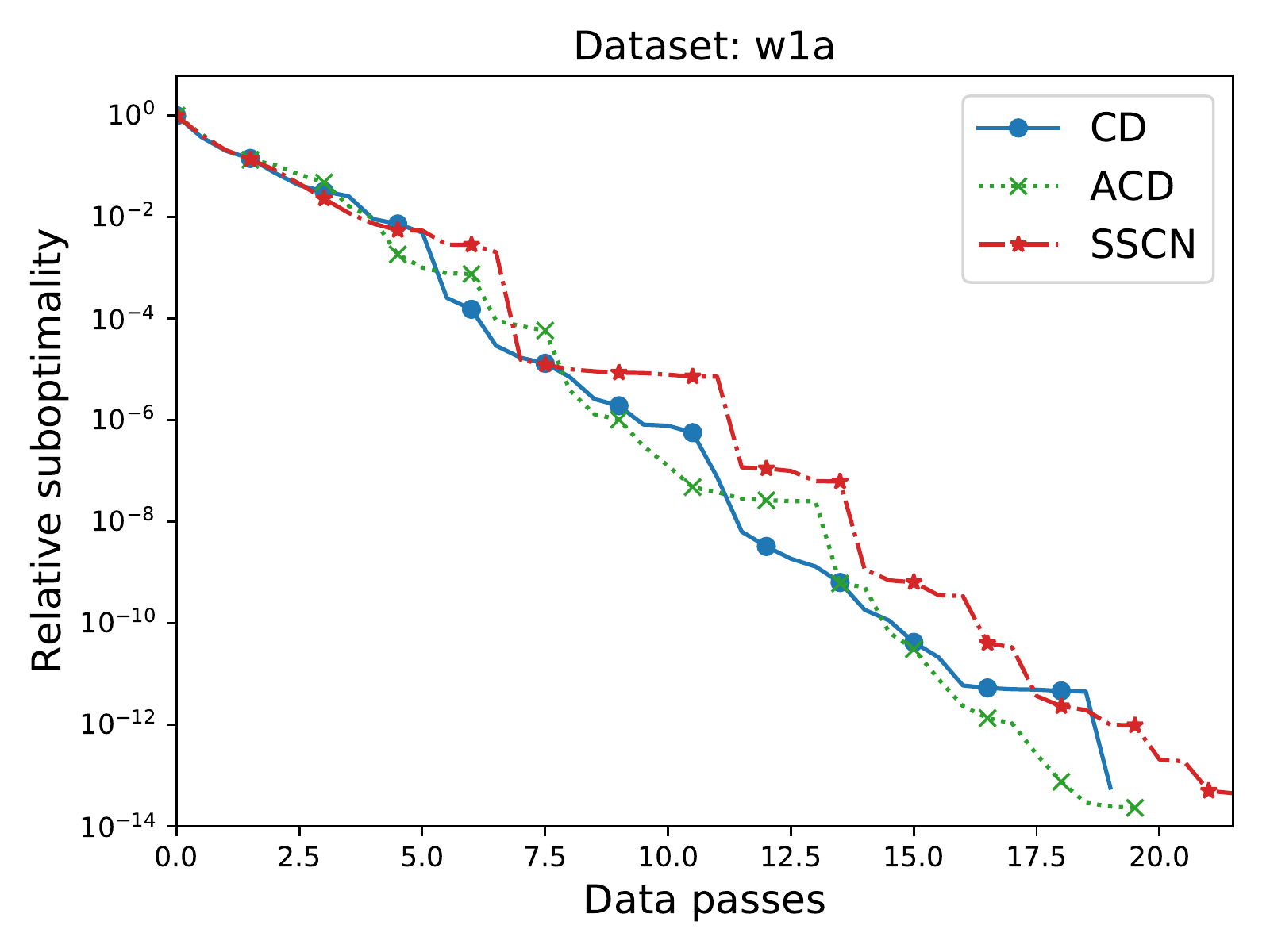}
\end{minipage}%
\begin{minipage}{0.3\textwidth}
  \centering
\includegraphics[width =  \textwidth ]{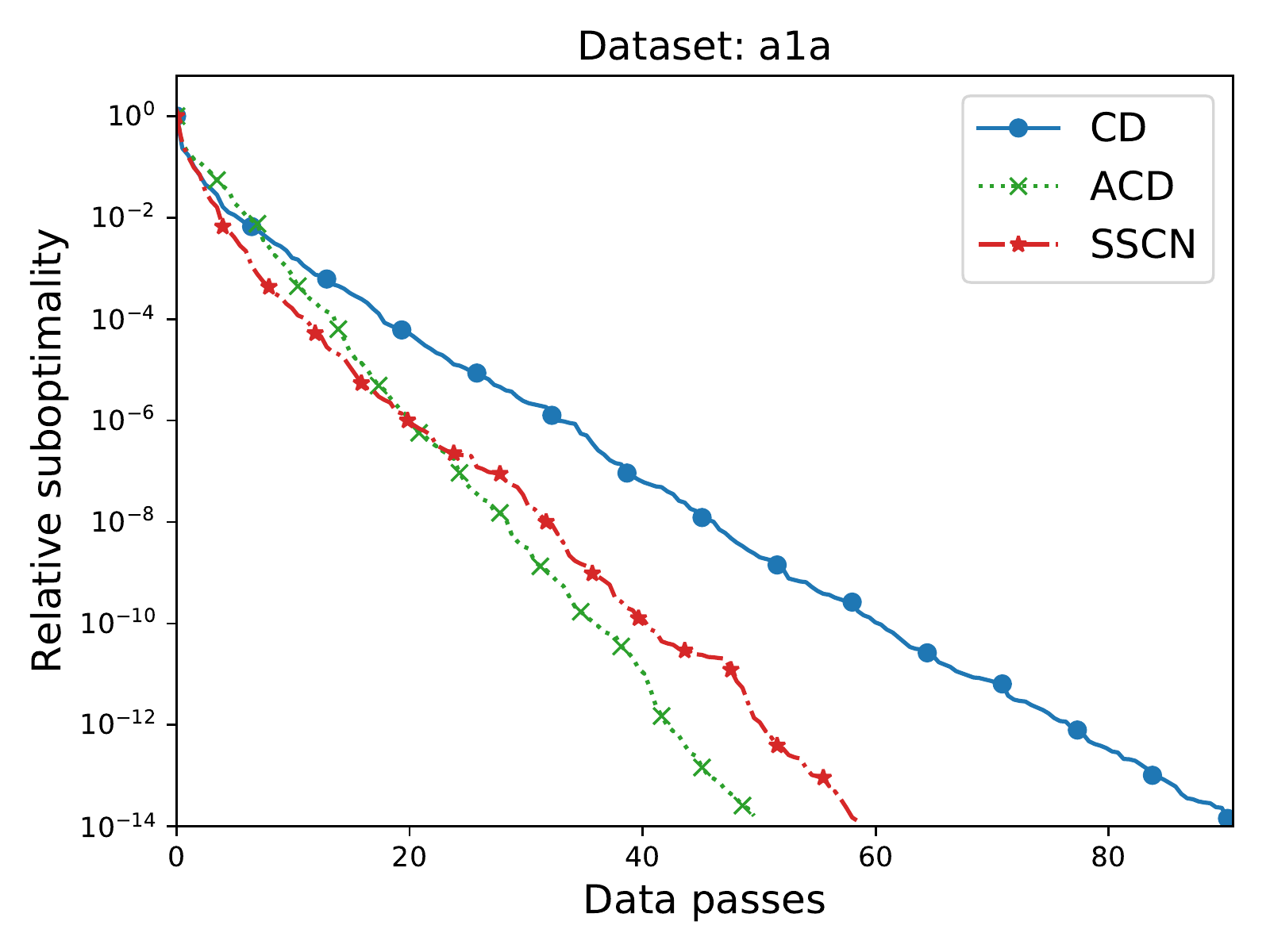}
\end{minipage}%
\begin{minipage}{0.3\textwidth}
  \centering
\includegraphics[width =  \textwidth ]{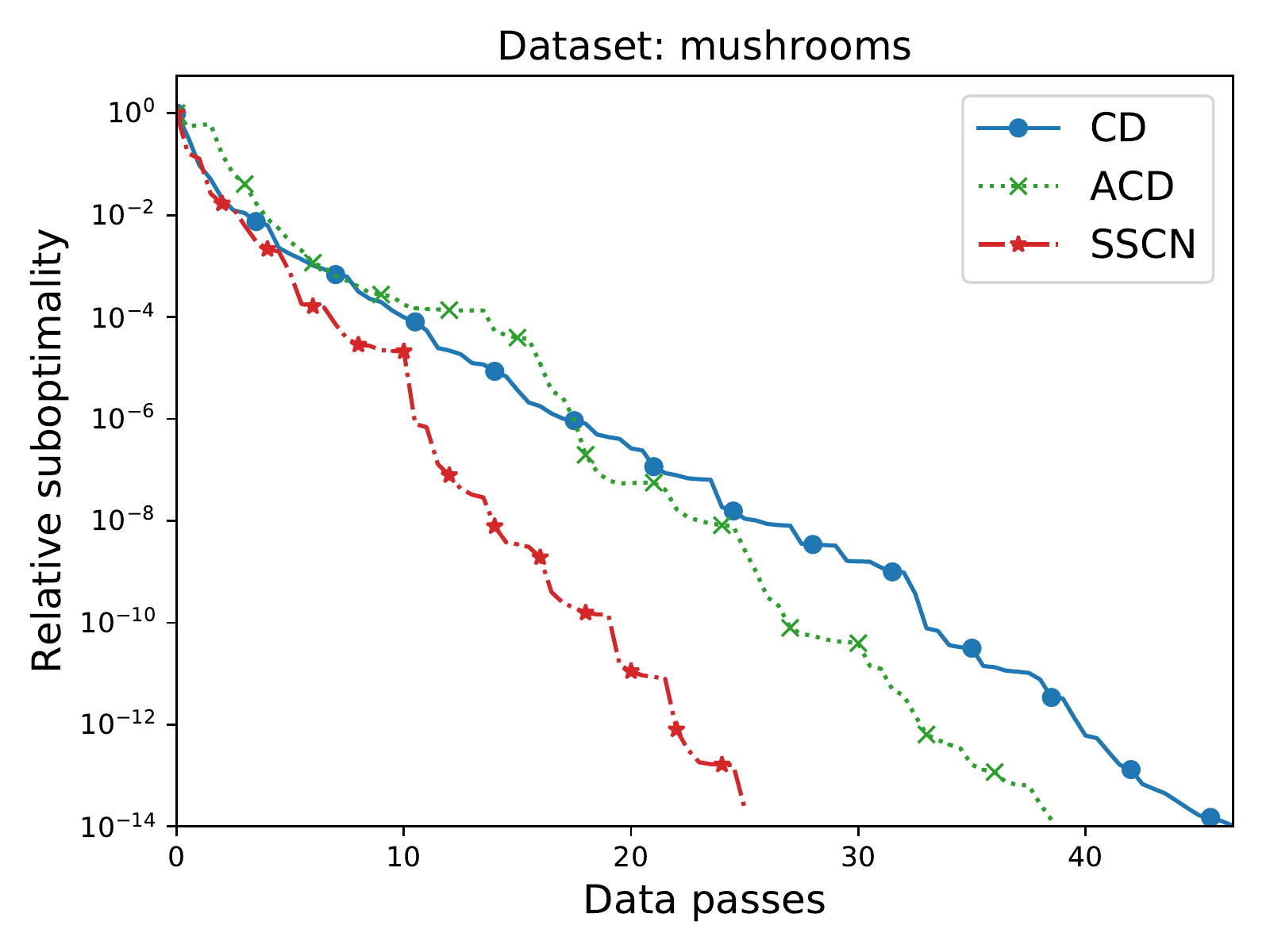}
\end{minipage}%
\\
\begin{minipage}{0.3\textwidth}
  \centering
\includegraphics[width =  \textwidth ]{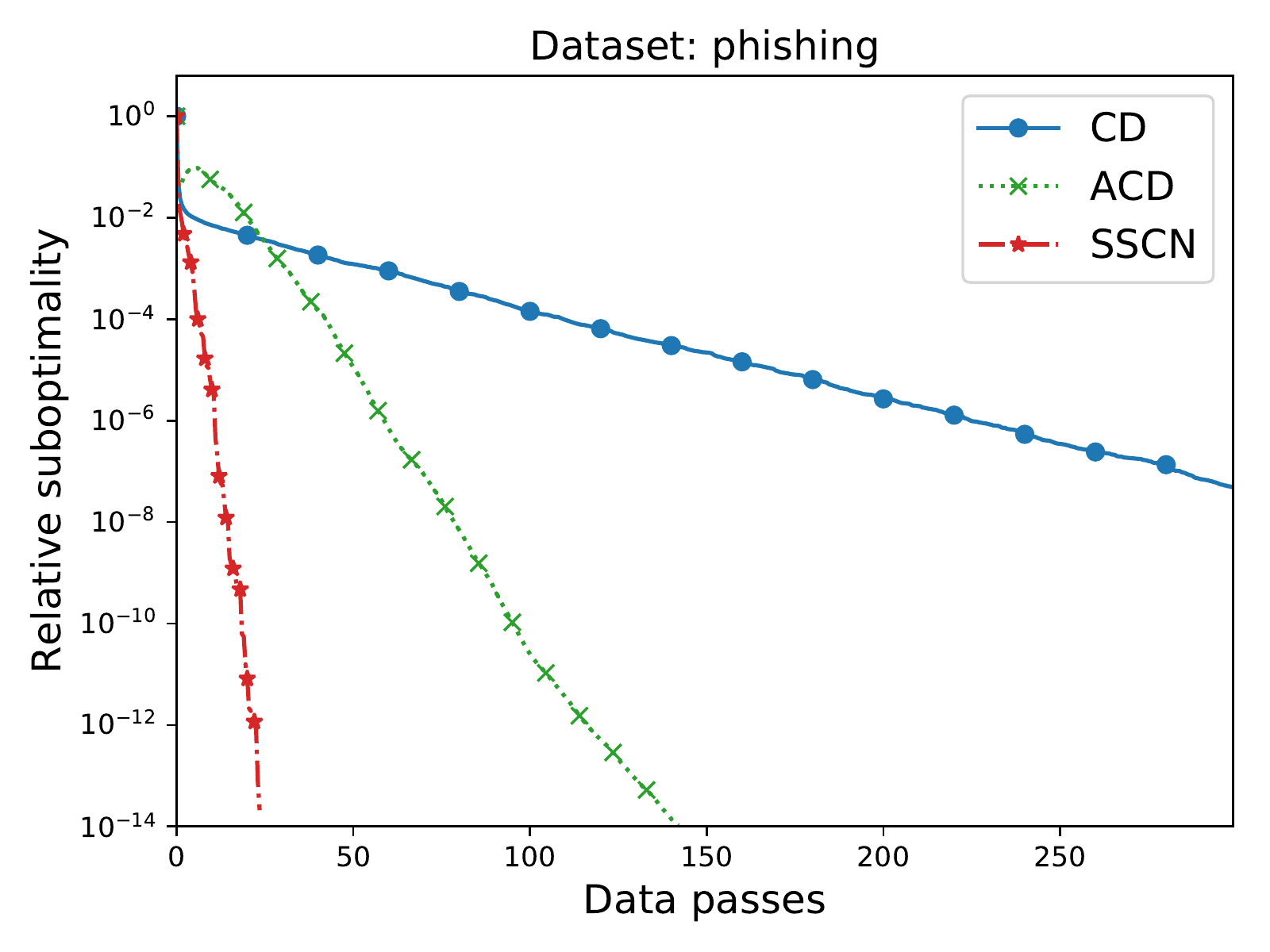}
\end{minipage}%
\begin{minipage}{0.3\textwidth}
  \centering
\includegraphics[width =  \textwidth ]{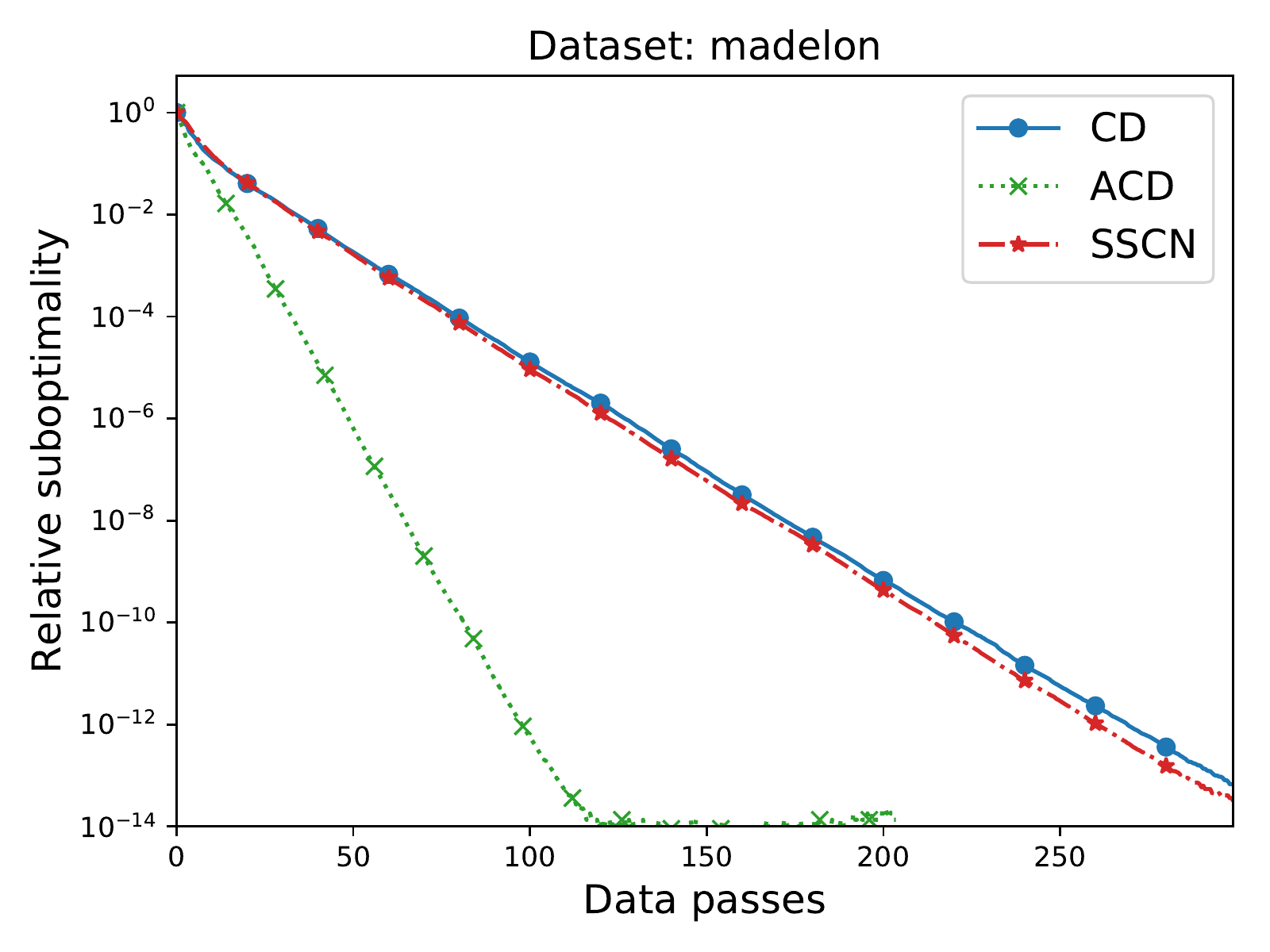}
\end{minipage}%
\begin{minipage}{0.3\textwidth}
  \centering
\includegraphics[width =  \textwidth ]{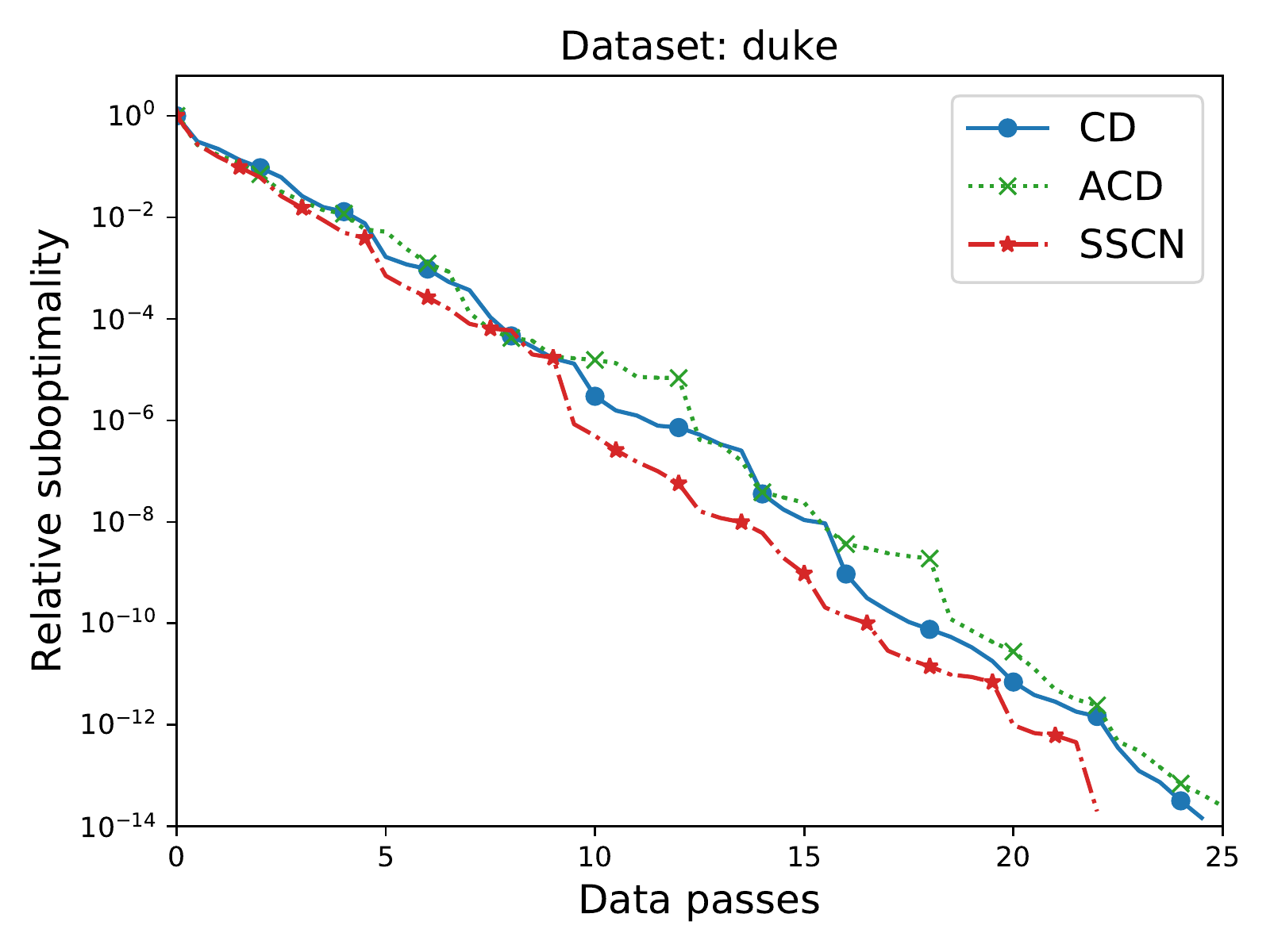}
\end{minipage}%
\\
\begin{minipage}{0.3\textwidth}
  \centering
\includegraphics[width =  \textwidth ]{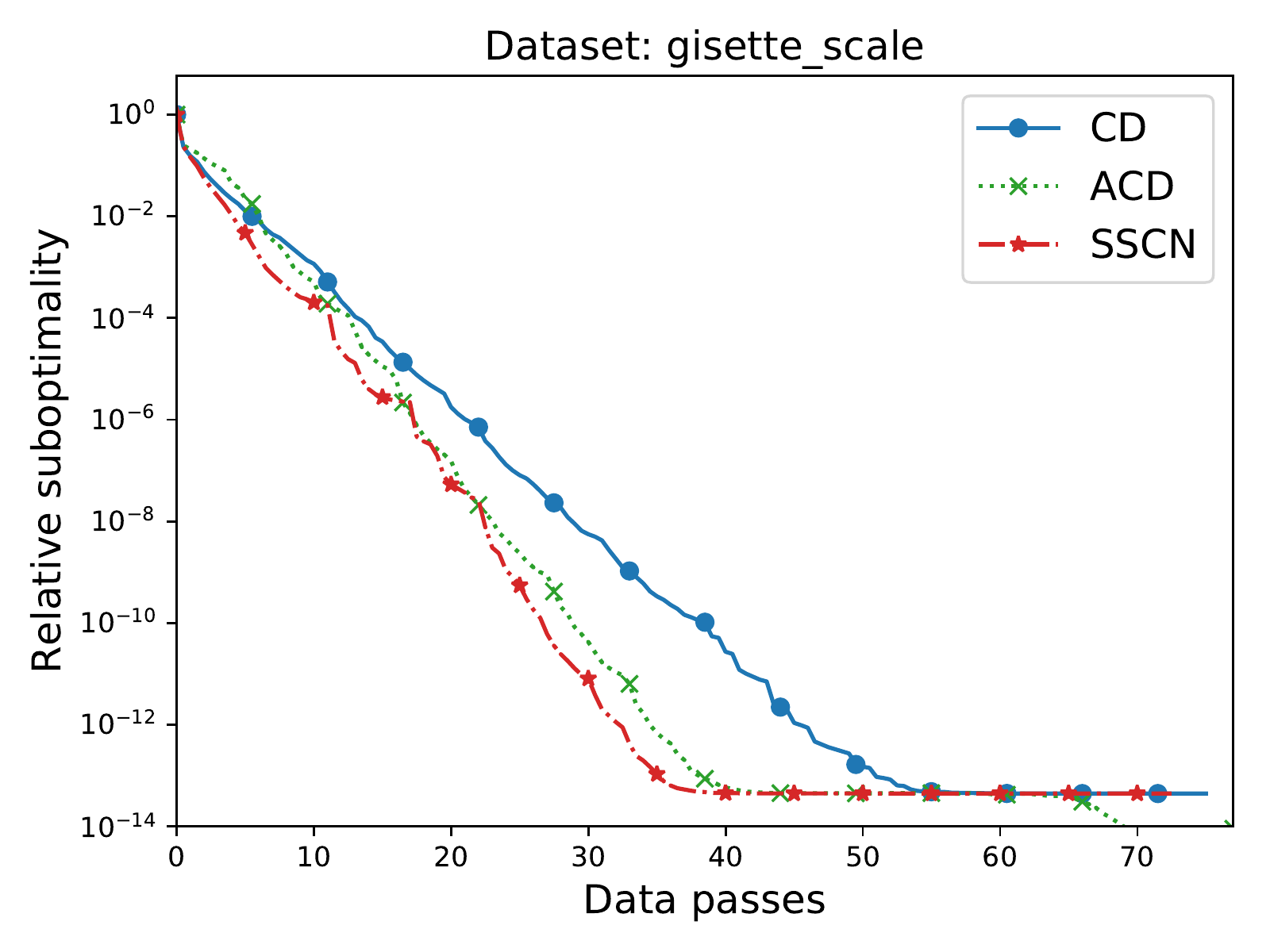}
\end{minipage}%
\begin{minipage}{0.3\textwidth}
  \centering
\includegraphics[width =  \textwidth ]{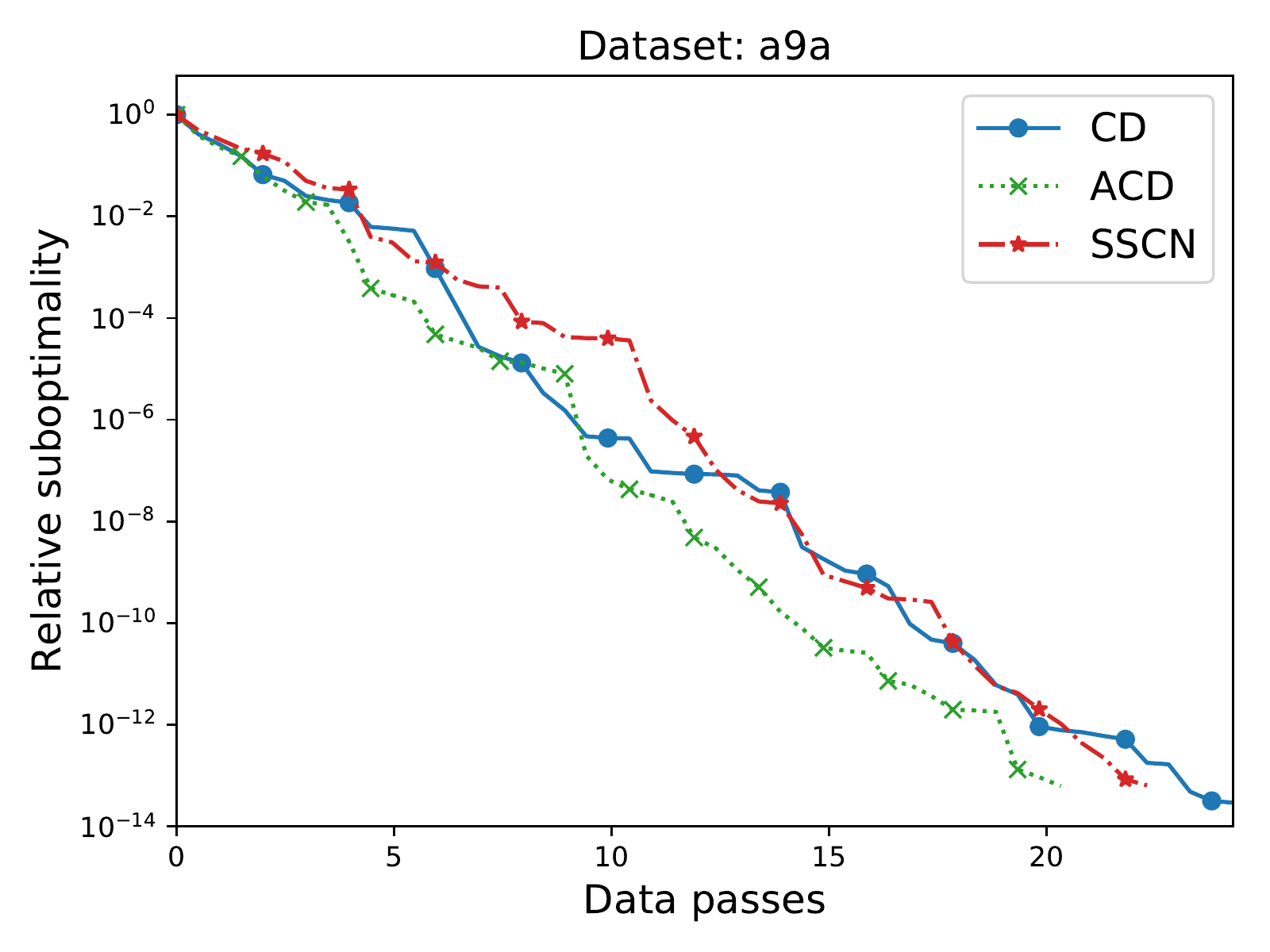}
\end{minipage}%
\begin{minipage}{0.3\textwidth}
  \centering
\includegraphics[width =  \textwidth ]{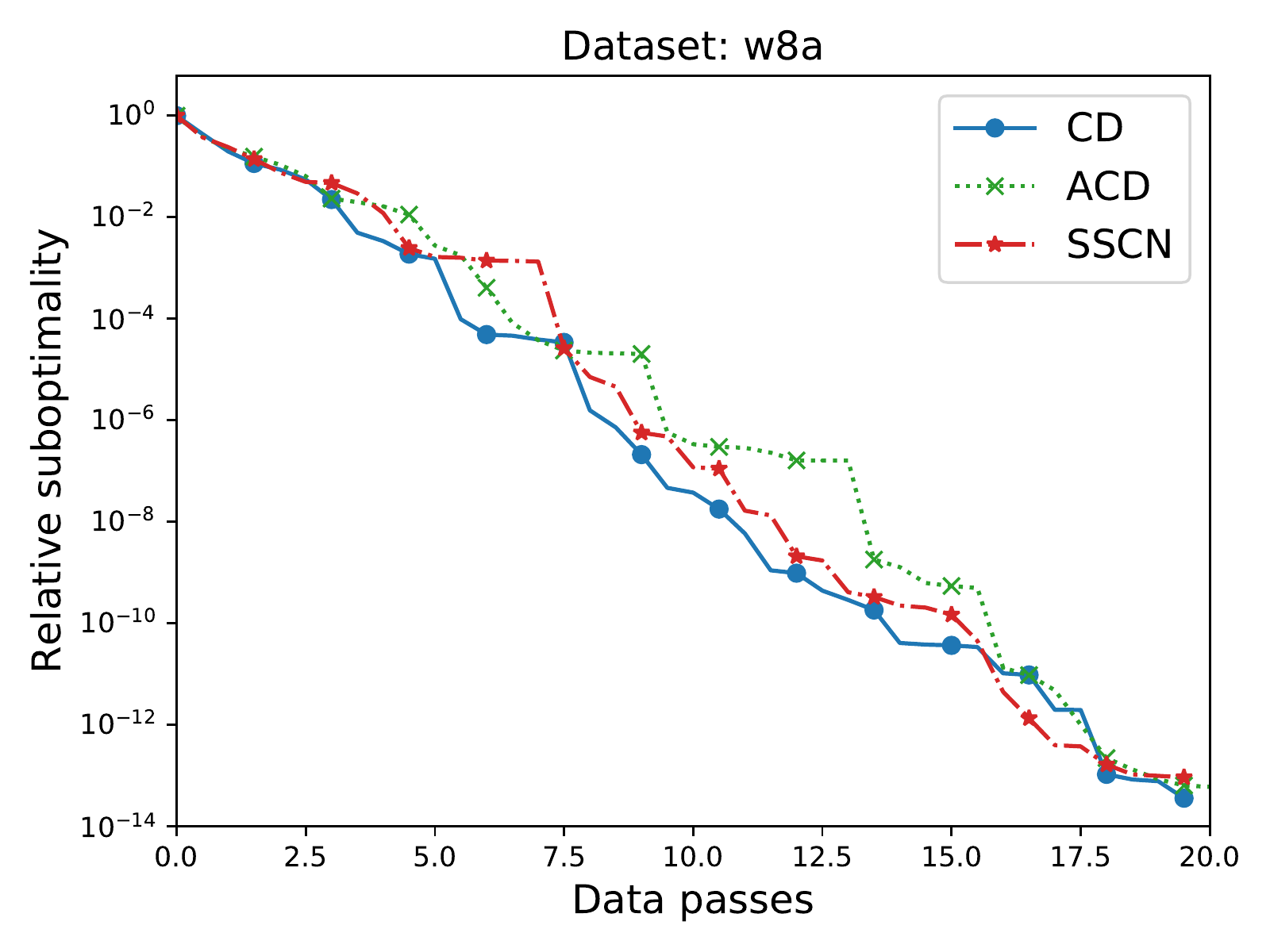}
\end{minipage}%
\caption{Comparison of coordinate descent, accelerated coordinate descent and {\tt SSCN} (all with uniform sampling) on LibSVM datasets. In each case we have normalized the data matrix to have identical norms of all columns.} 
\label{fig:sscn_libsvm_normalzied}
\end{figure}

In the second experiment, we compare methods with $\tau>1$: {\tt SSCN} and {\tt SDNA}~\cite{sdna} (analogous first-order method). Again, we consider the logistic regression problem on LIBSVM data. We consider $\tau \in \{1,5,25\}$. In all cases, we sample uniformly --  every subset of size $\tau$ have equal chance to be chosen at every iteration (independent of the past). 

There is, however, one tricky part in terms of implementation. While we can evaluate and store $M_{e_i}$ ($i\leq d$) cheaply for linear models, this is not the case for evaluating/storing $M_S$ (at least we do not know how to do it efficiently). Therefore, we use $M_S = M$ for $|S|>1$ for {\tt SSCN}. Figure~\ref{fig:sscn_libsvm_sdna} shows the result. As expected, {\tt SSCN} has outperformed {\tt SDNA}.

\begin{figure}[!h]
\centering
\begin{minipage}{0.3\textwidth}
  \centering
\includegraphics[width =  \textwidth ]{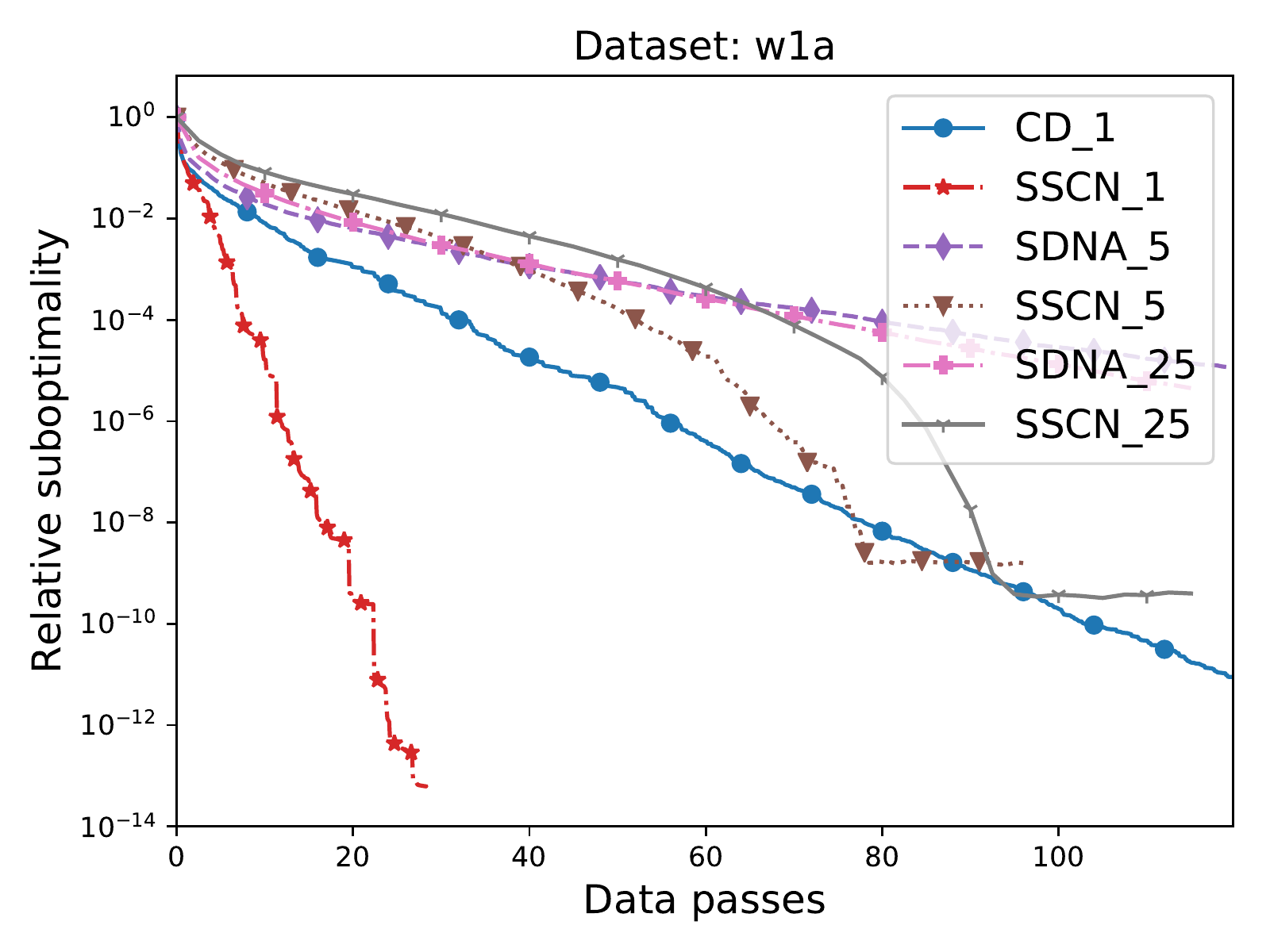}
\end{minipage}%
\begin{minipage}{0.3\textwidth}
  \centering
\includegraphics[width =  \textwidth ]{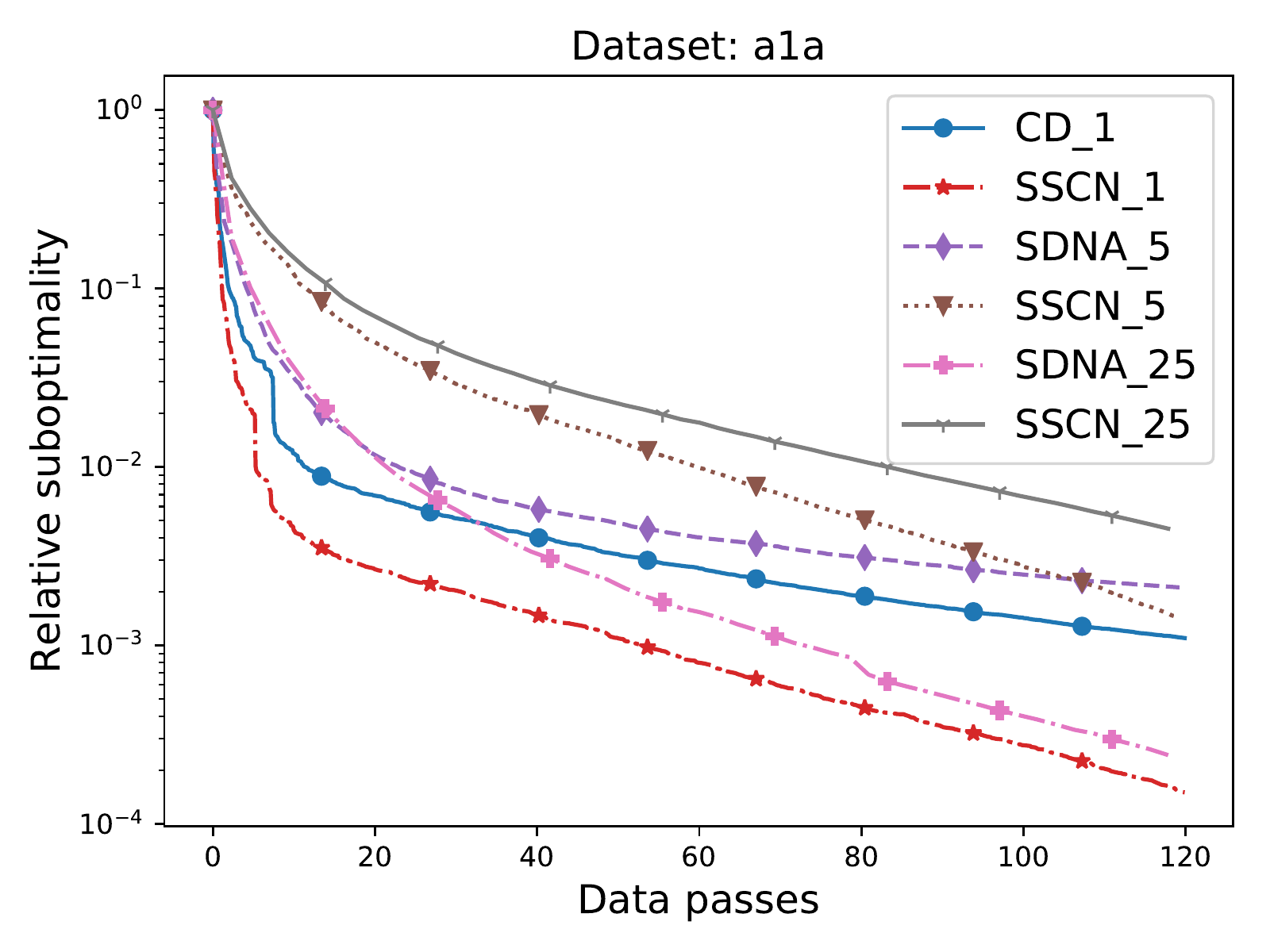}
\end{minipage}%
\begin{minipage}{0.3\textwidth}
  \centering
\includegraphics[width =  \textwidth ]{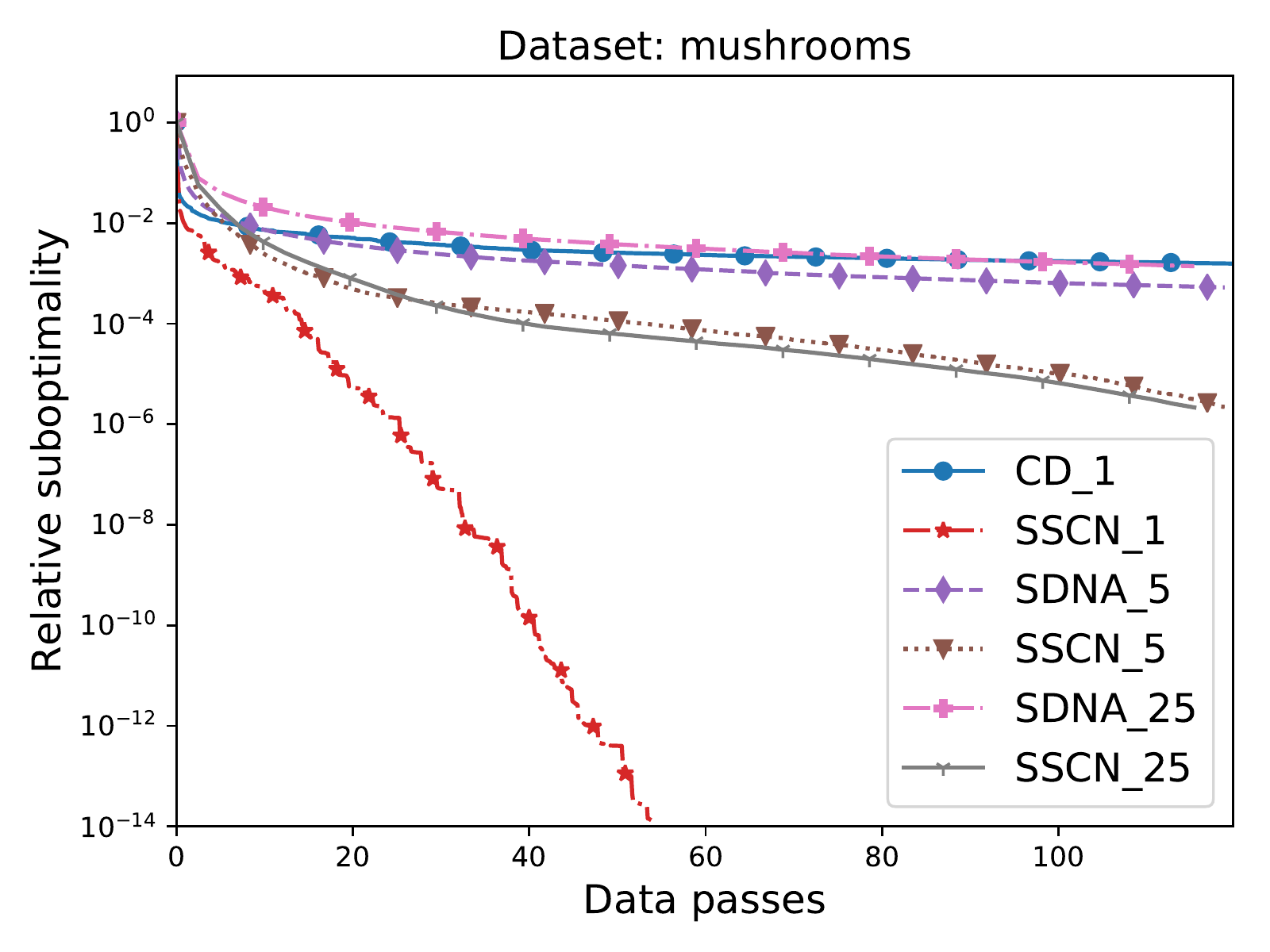}
\end{minipage}%
\\
\begin{minipage}{0.3\textwidth}
  \centering
\includegraphics[width =  \textwidth ]{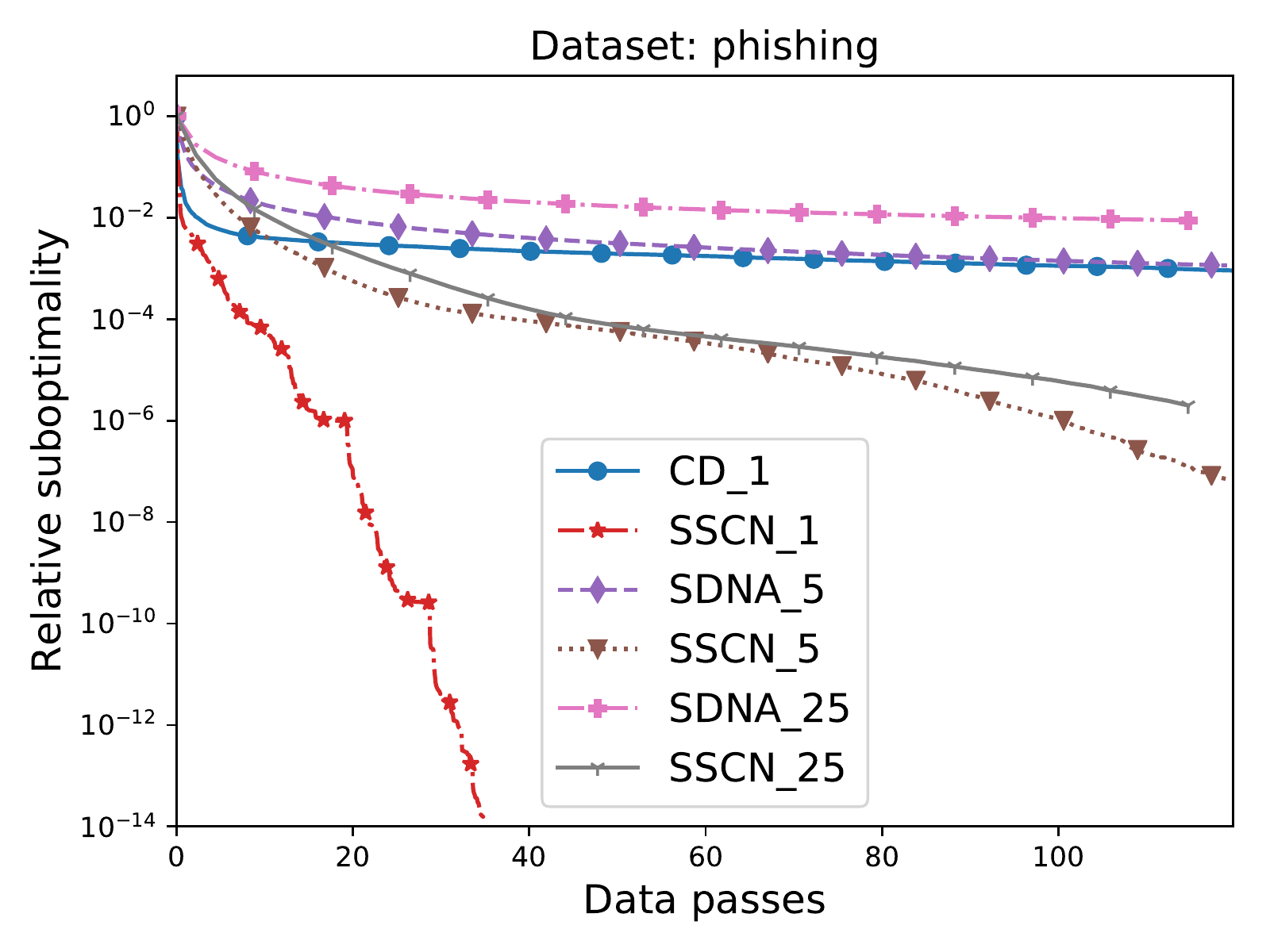}
\end{minipage}%
\begin{minipage}{0.3\textwidth}
  \centering
\includegraphics[width =  \textwidth ]{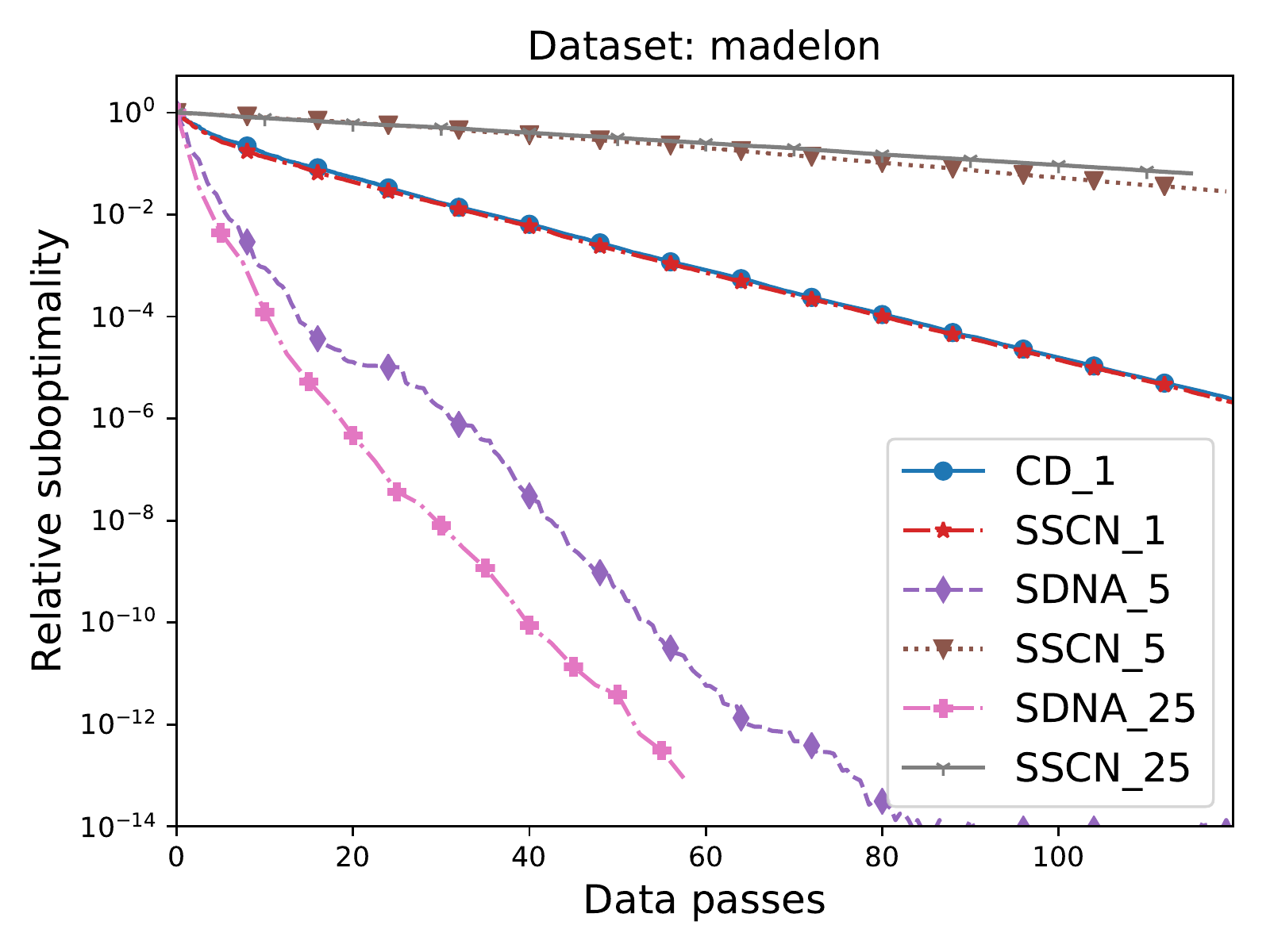}
\end{minipage}%
\begin{minipage}{0.3\textwidth}
  \centering
\includegraphics[width =  \textwidth ]{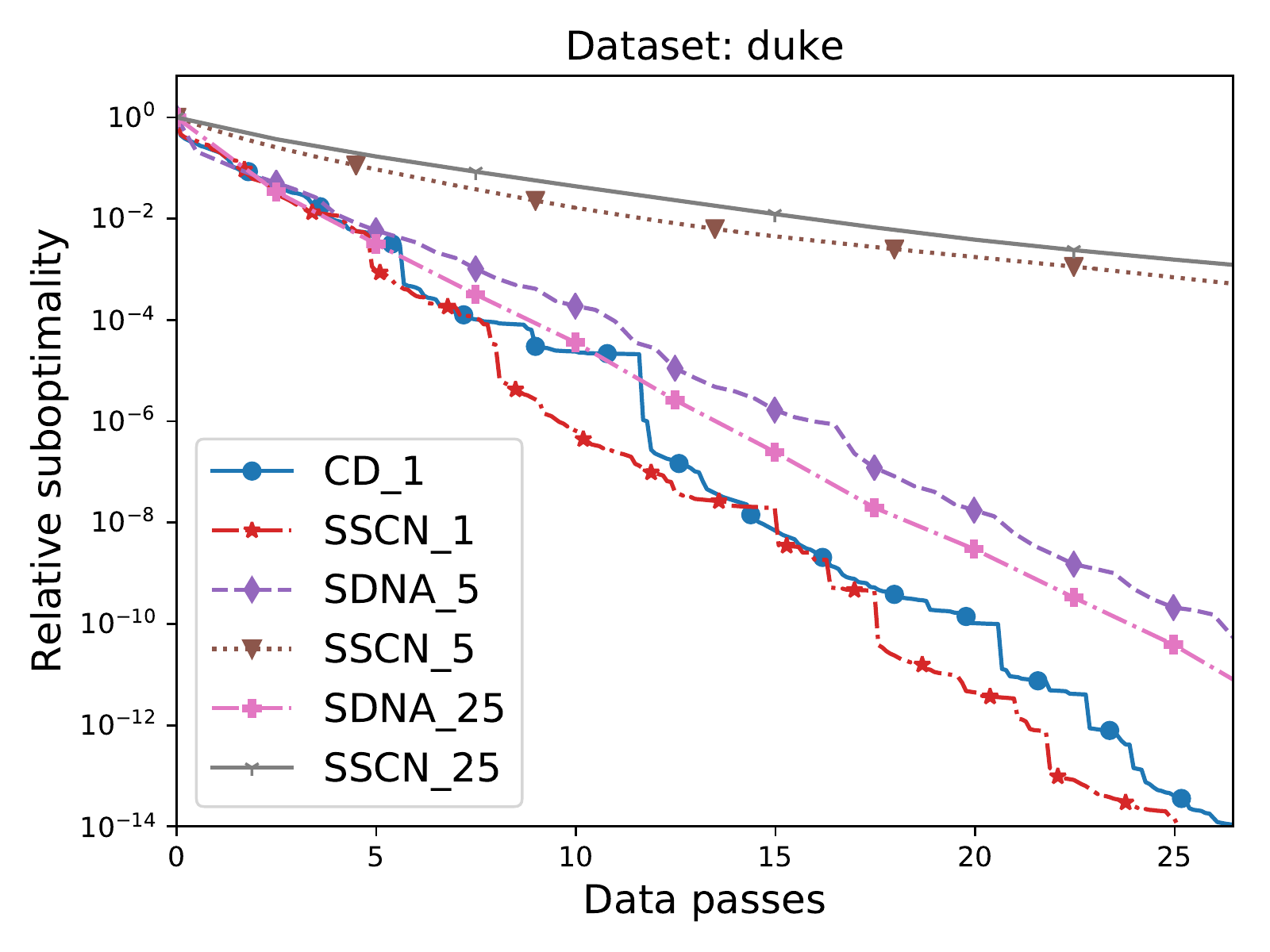}
\end{minipage}%
\\
\begin{minipage}{0.3\textwidth}
  \centering
\includegraphics[width =  \textwidth ]{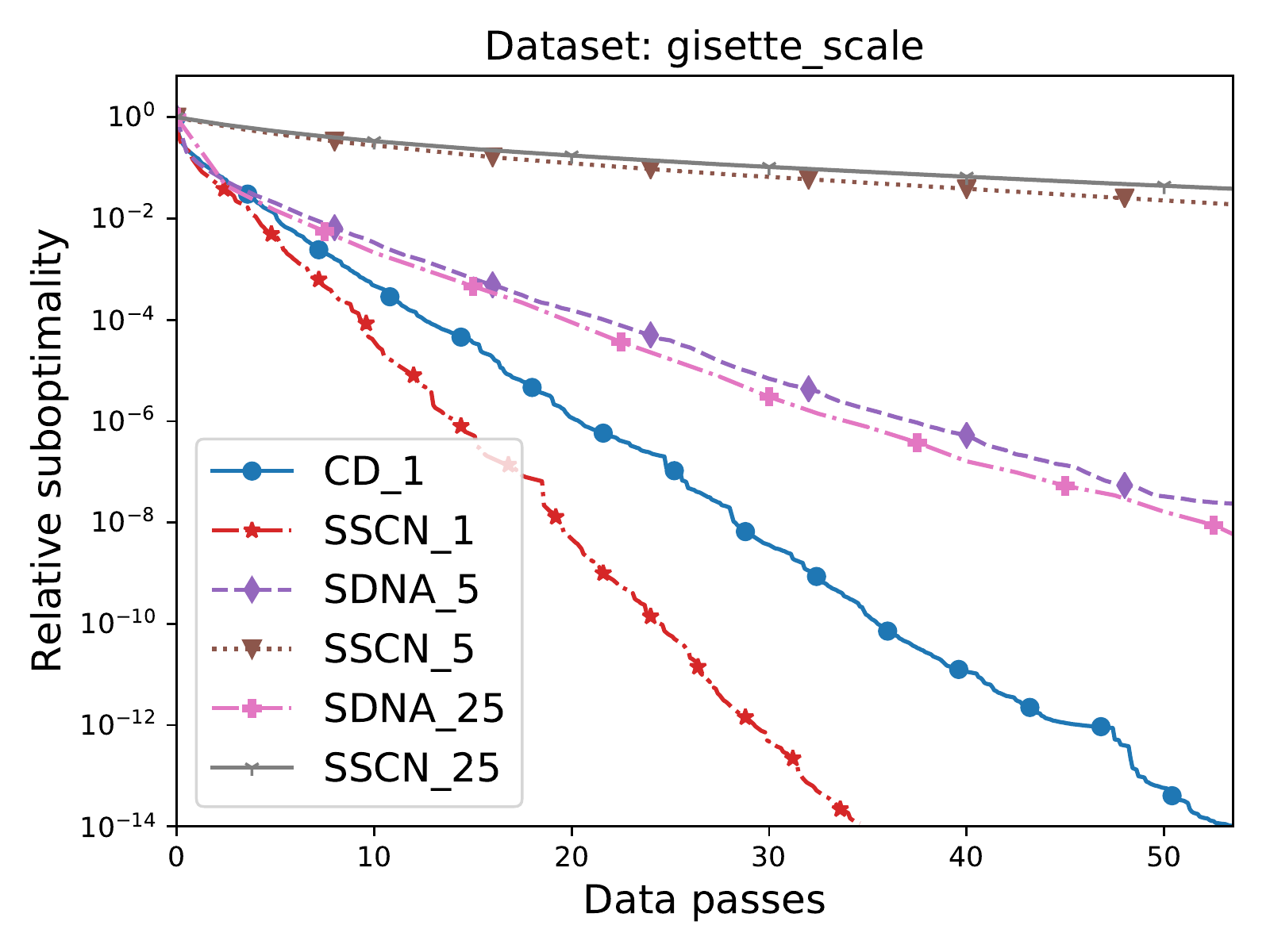}
\end{minipage}%
\begin{minipage}{0.3\textwidth}
  \centering
\includegraphics[width =  \textwidth ]{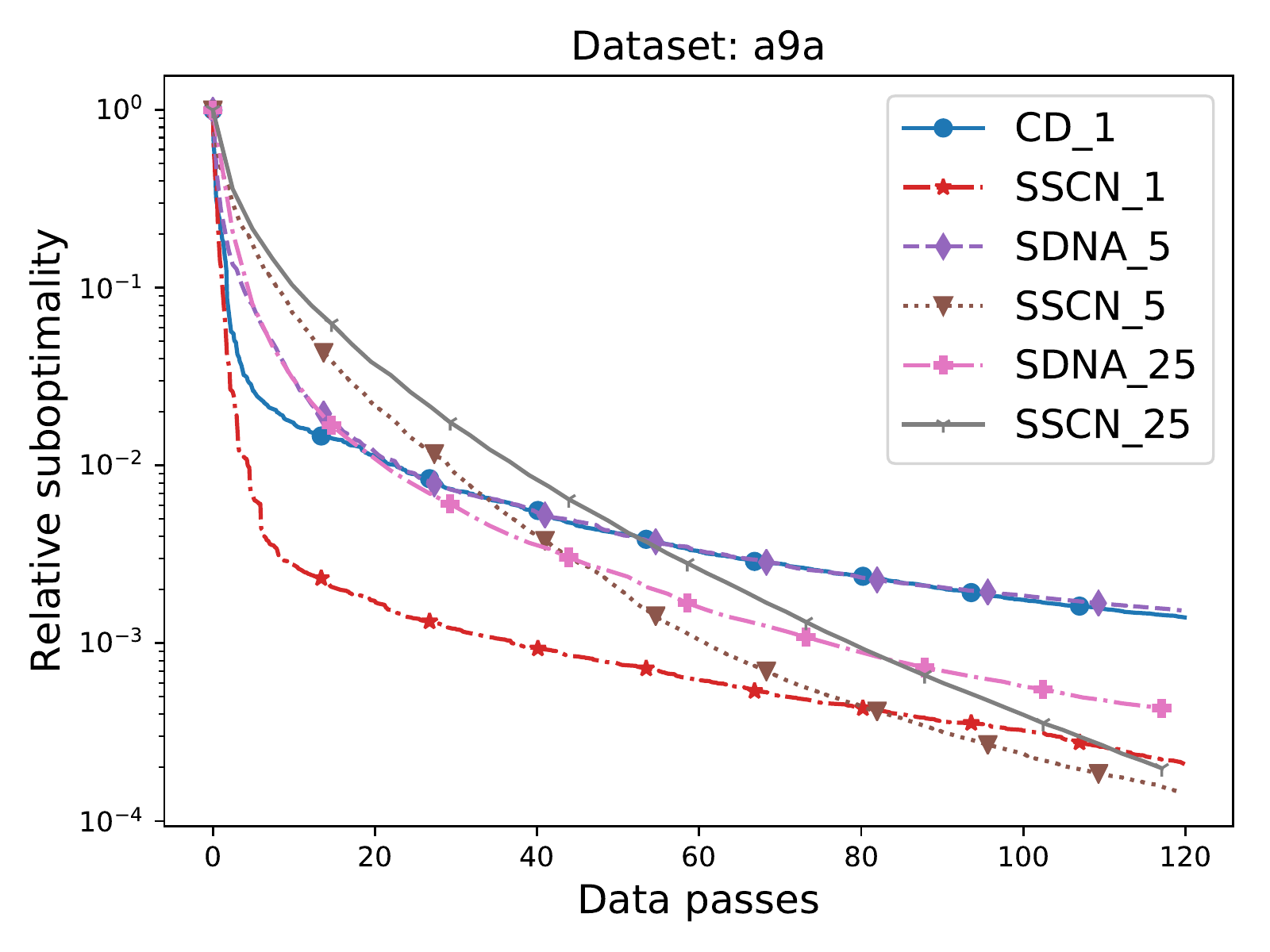}
\end{minipage}%
\begin{minipage}{0.3\textwidth}
  \centering
\includegraphics[width =  \textwidth ]{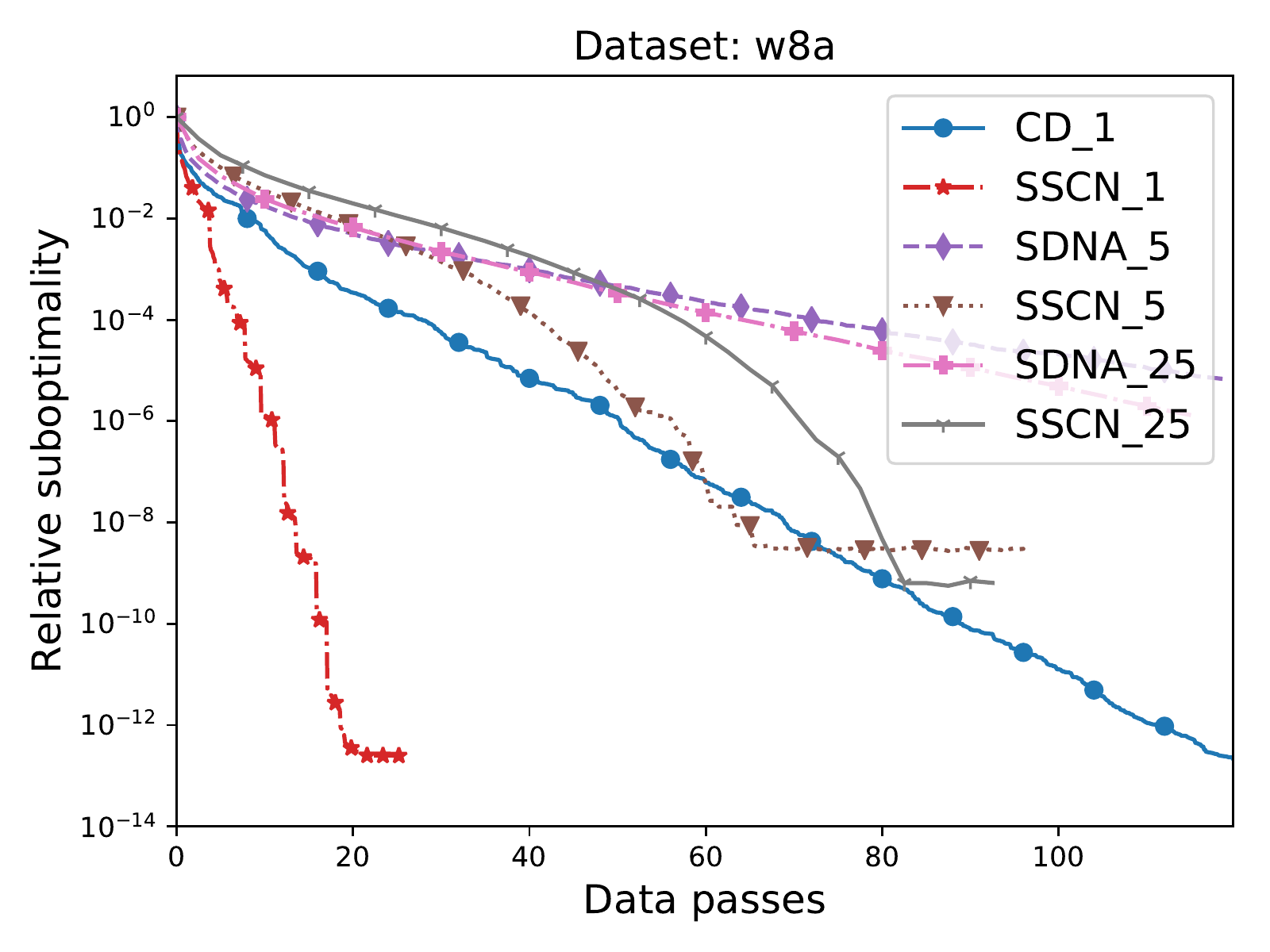}
\end{minipage}%
\caption{{\tt SSCN} vs. {\tt SDNA} on LibSVM datasets. All algorithms with uniform sampling.} 
\label{fig:sscn_libsvm_sdna}
\end{figure}

\subsection{Log-sum-exp \label{sec:logsumexp}}

In this section, let us consider unconstrained minimization of the following Log-sum-exp function
$$
f(x)  =  \sigma \log\left( \sum\limits_{i = 1}^m \exp \left( \frac{\la a_i, x \ra - b_i}{\sigma} \right) \right),
\quad x \in \R^d,
$$
where $\sigma > 0$ is a \textit{smoothing} parameter, while
$a_i \in \R^n, 1 \leq i \leq m$ and $b \in \R^m$ are given data.
This function has both Lipschitz continuous gradient and Lipschitz continuous Hessian
(see Example~1 in~\cite{doikov2019minimizing}). 

In our experiments, we first generate randomly elements of $\{ \tilde{a}_i \}_{i = 1}^m$ and $b$ from
uniform distribution on $[-1, 1]$. Then, we form an auxiliary function 
\[\tilde{f}(x) \eqdef \sigma \log\left( \sum\limits_{i = 1}^m \exp \bigl(  \frac{\la \tilde{a}_i, x \ra - b_i}{\sigma} \bigr) \right),\]
using these parameters, and set
$$
a_i   \eqdef  \tilde{a}_i - \nabla \tilde{f}(0), \quad 1 \leq i \leq m.
$$
Thus, we essentially obtain the optimum $x^{*}$ of $f$ in the origin, since $\nabla f(0) = 0$.
We use $x_0 \eqdef e$ (vector of all ones) as a starting point,
and always set $m  \eqdef 6d$.

For this problem, we compare the performance of {\tt SSCN} with the first-order Coordinate Descent ({\tt CD}), using uniform samples of coordinates $S \subseteq [d]$ of a fixed size $\tau = |S|$. 

Note, that keeping scalar products $\{ \la a_i, x_k \ra \}_{i = 1}^m$
precomputed for a current point $x_k$, we are able to compute the partial gradient
$\nabla_{\mS} f(x^k)$ in time $O(\tau m)$ and the partial Hessian $\nabla^2_{\mS} f(x^k)$ in time $O(\tau^2 m)$.
To find the next direction $h^k$ of {\tt SSCN} (solving the Cubic subproblem), we call Nonlinear Conjugate Gradient method,
and use the following condition as a stopping criterion:
$$
\| \nabla_h T_{\mS}(x^k; h^k) \|  \leq  10^{-4},
$$
where $T_{\mS}(x^k; h)   \eqdef \la \nabla_{\mS} f(x^k), h \ra + \frac{1}{2}\la \nabla_{\mS}^2 f(x^k)h, h \ra + \frac{M_k}{6}\|\mS h\|^3$ is the Cubic model, and $M_k \geq 0$ is a regularization constant.

For both methods, we use one-dimensional search at every iteration, to fit the corresponding parameter:
\begin{enumerate}
	\item For the Coordinate Descent, we find $L_k$ such that $f(x^{k}) - f(x^{k + 1}) \geq \frac{1}{2L_k}\| \nabla_{\mS} f(x^k) \|^2$, 
	where $x^{k + 1}$ is the next point of the method: $x^{k + 1} = x^k + \frac{1}{L_k}\mS \nabla_{\mS} f(x^k)$.
	
	\item For {\tt SSCN}, we find $M_k$ such that~\eqref{eq:sscn_coordinate_ub_full}
	is satisfied, i.e. $f(x^{k}) - f(x^{k + 1}) \geq -T_{\mS}(x^k, h^k)$.
\end{enumerate}
Therefore,  we need to evaluate the function value inside the procedure, which is not very expensive.

The results are shown on Figures~\ref{fig:sscn_log_sum_exp_500},\ref{fig:sscn_log_sum_exp_1000},
for $d = 500$ and $1000$ respectively\footnote{Clock time was evaluated using the machine with Intel Core i7-8700 CPU, 3.20GHz; 16 GB RAM.}.
We see, that {\tt SSCN} outperforms {\tt CD} significantly in terms of the iteration rate.
For {\tt SSCN} with a medium batchsize $\tau$, we may obtain the best performance 
in terms of the total computational time.

\begin{figure}[h!]
	\centering
	\begin{minipage}{0.3\textwidth}
		\centering
		\includegraphics[width =  \textwidth ]{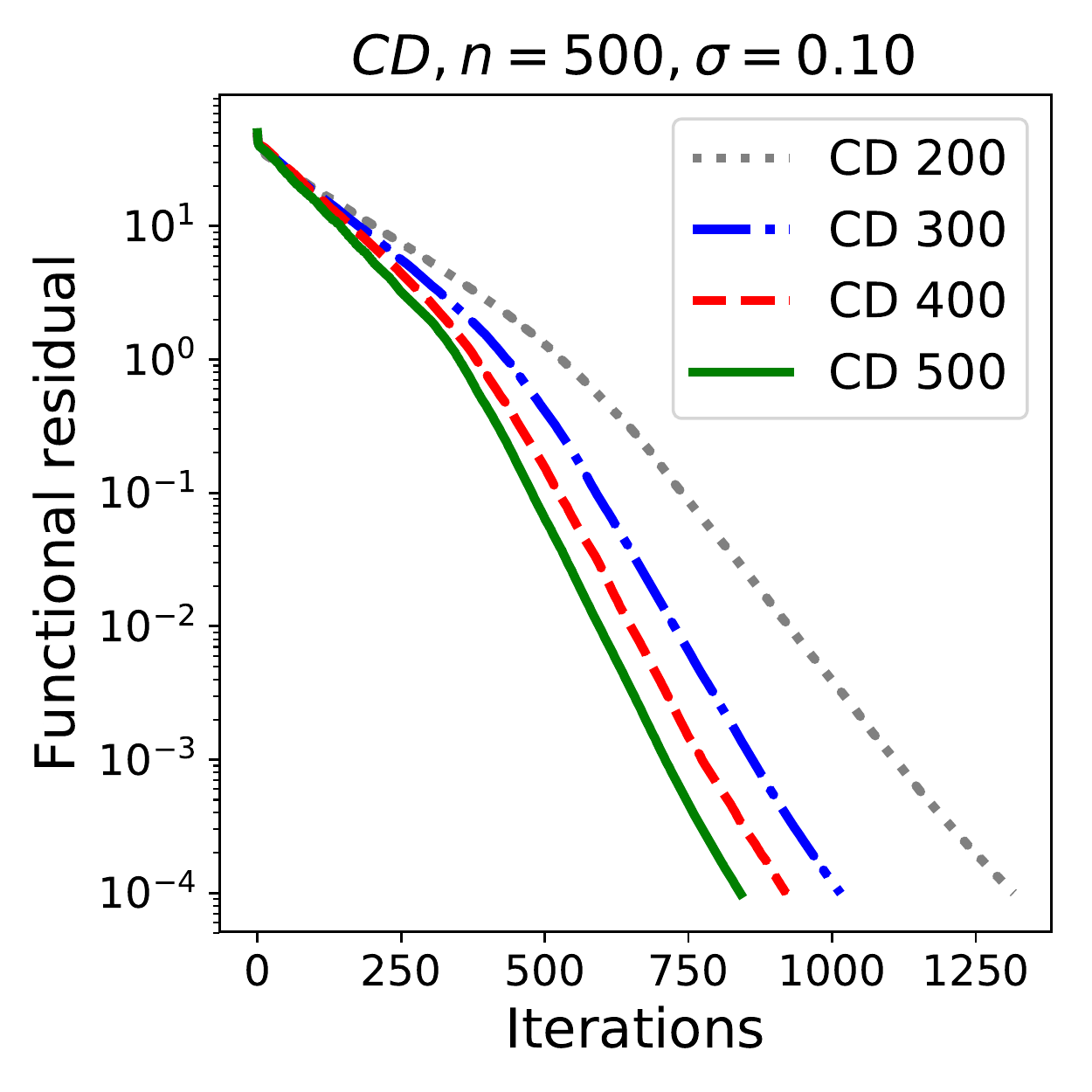}
	\end{minipage}
	\begin{minipage}{0.3\textwidth}
		\centering
		\includegraphics[width =  \textwidth ]{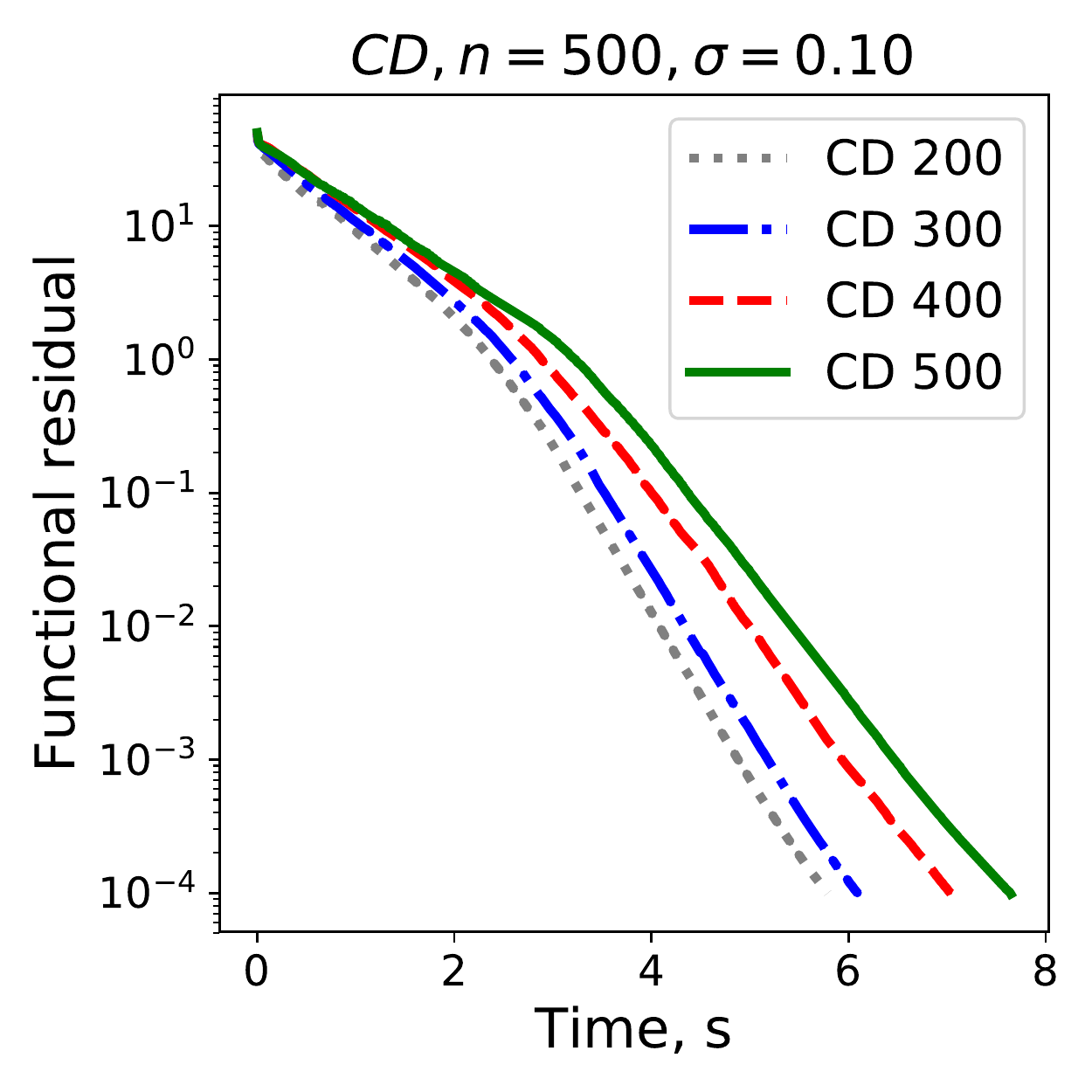}
	\end{minipage}
	\begin{minipage}{0.3\textwidth}
		\centering
		\includegraphics[width =  \textwidth ]{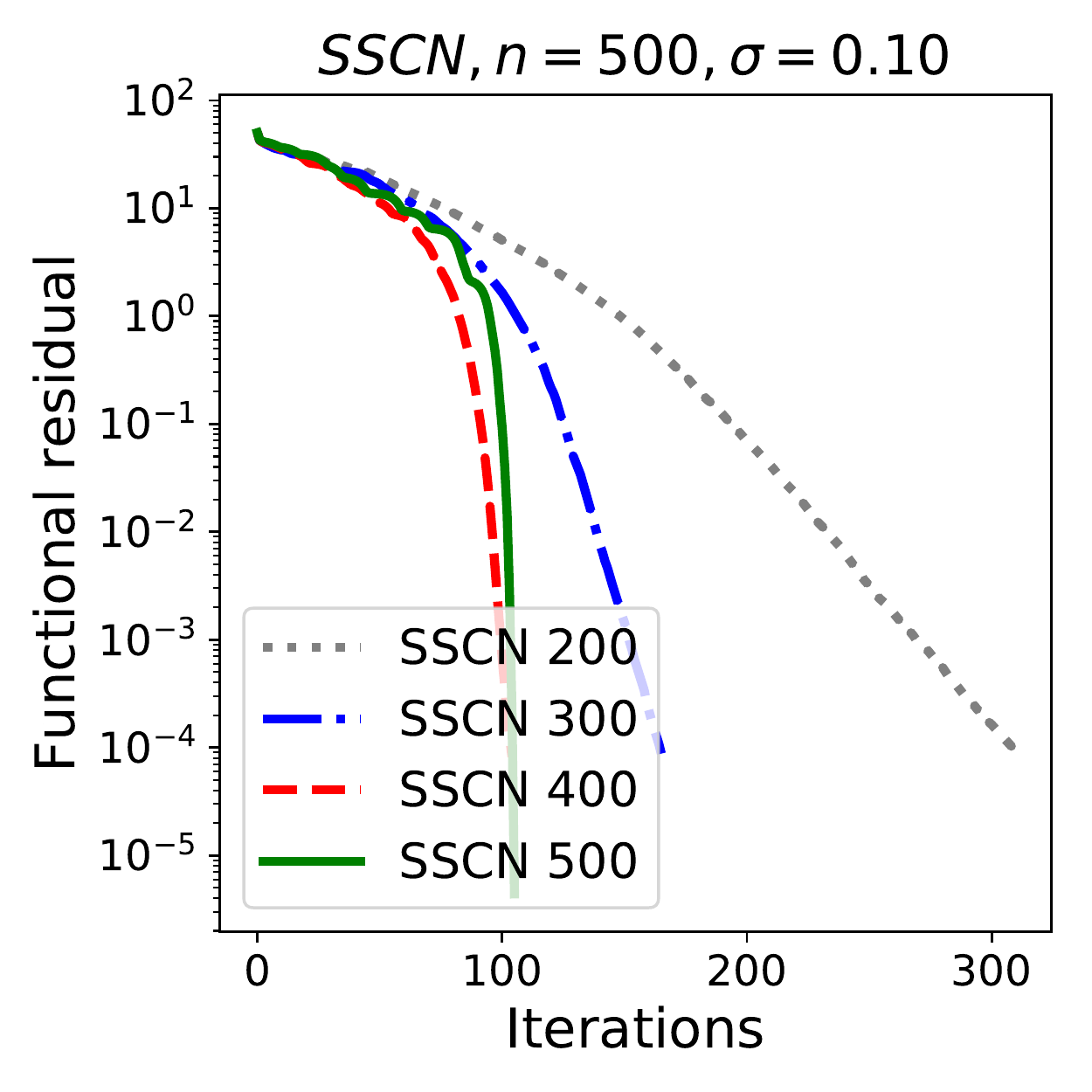}
	\end{minipage}
	\\
	\begin{minipage}{0.3\textwidth}
		\centering
		\includegraphics[width =  \textwidth ]{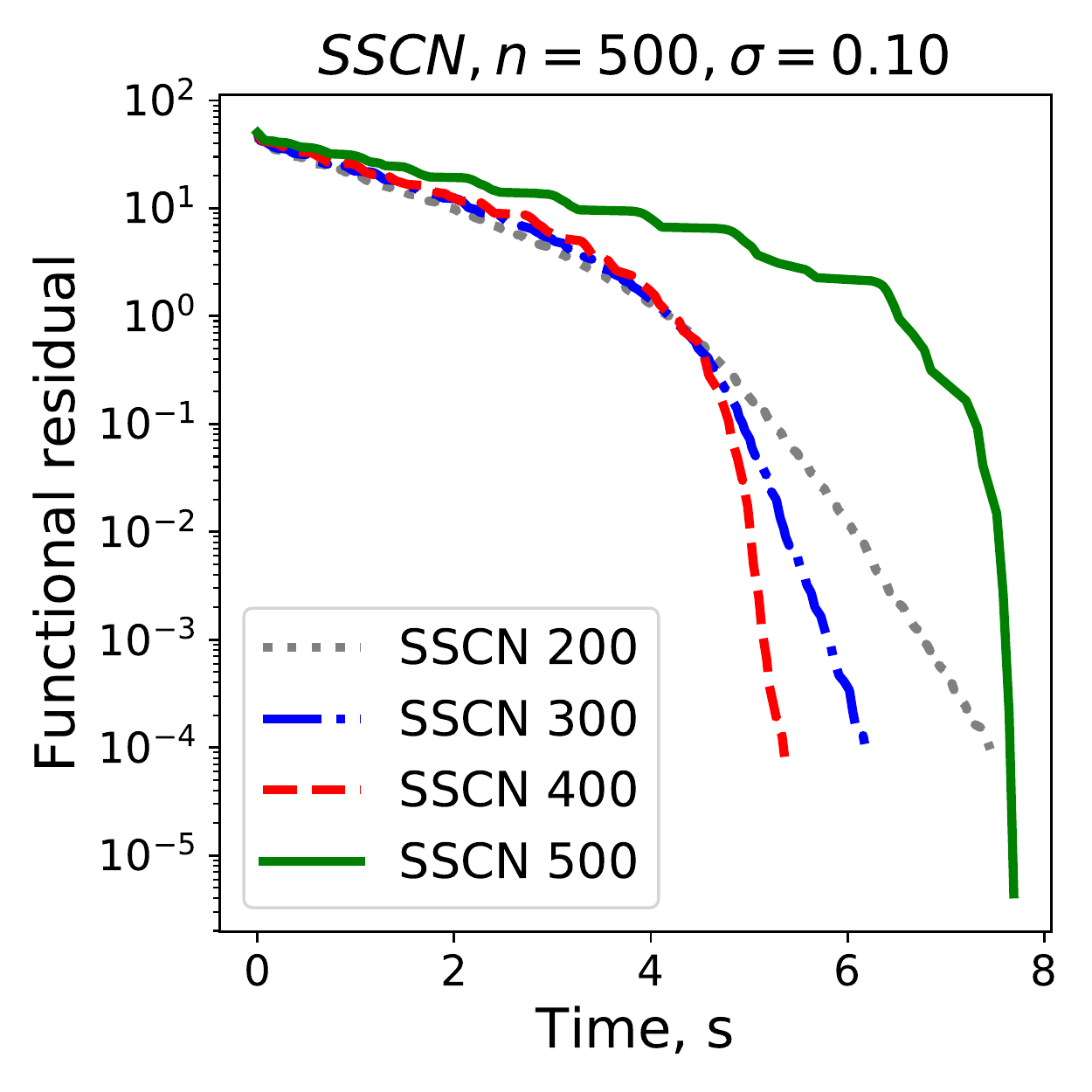}
	\end{minipage}
		\begin{minipage}{0.3\textwidth}
		\centering
		\includegraphics[width =  \textwidth ]{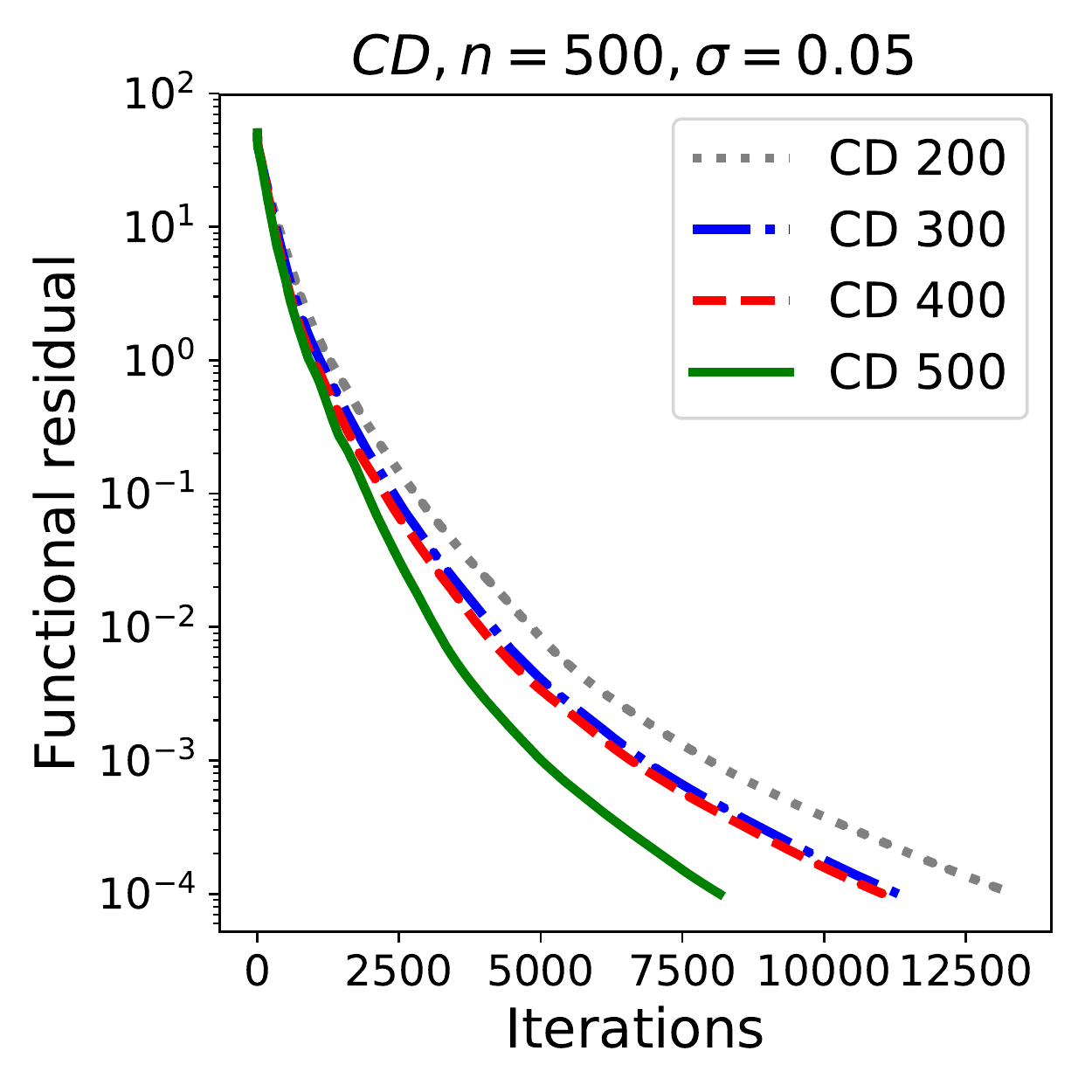}
	\end{minipage}
	\begin{minipage}{0.3\textwidth}
		\centering
		\includegraphics[width =  \textwidth ]{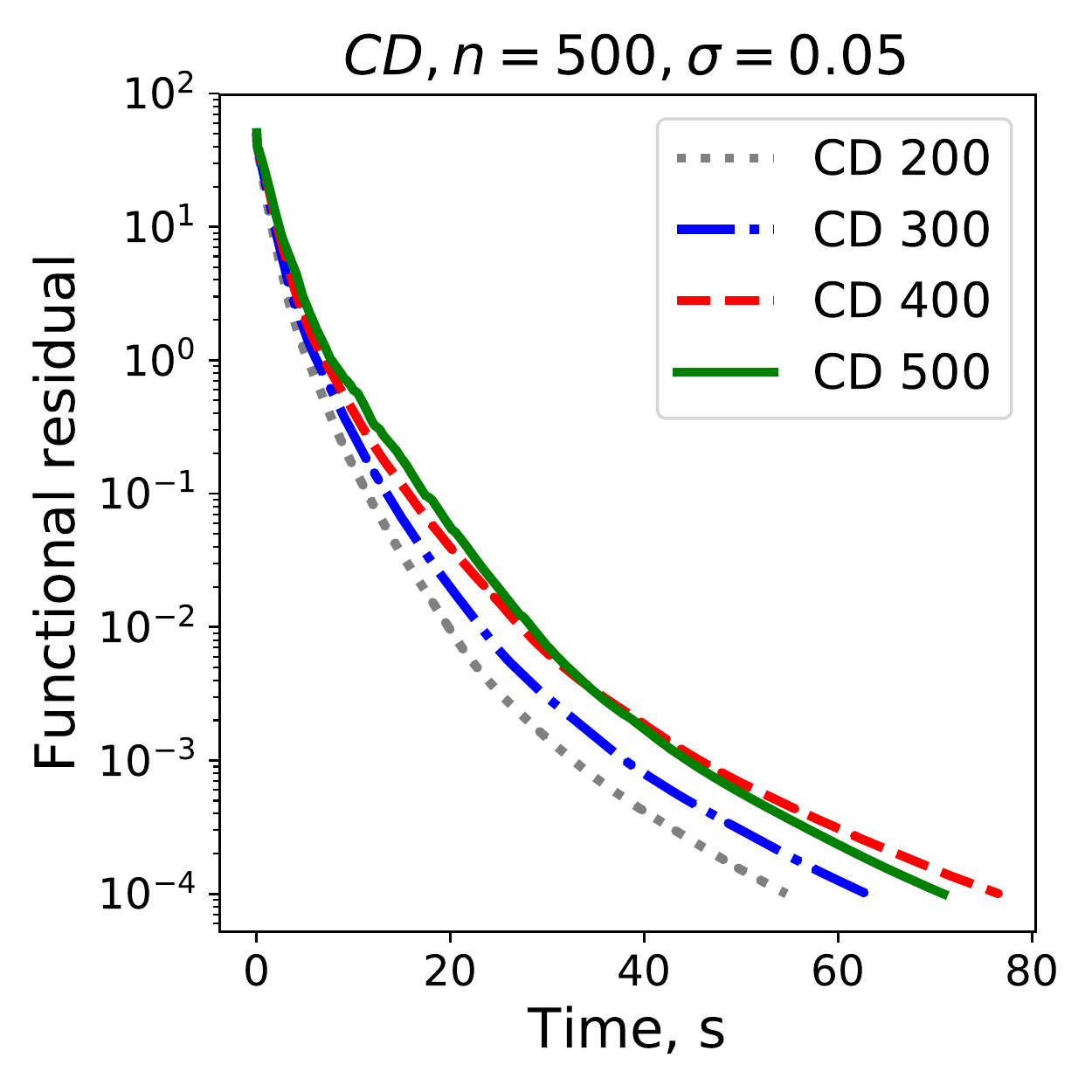}
	\end{minipage}
	\\
	\begin{minipage}{0.3\textwidth}
		\centering
		\includegraphics[width =  \textwidth ]{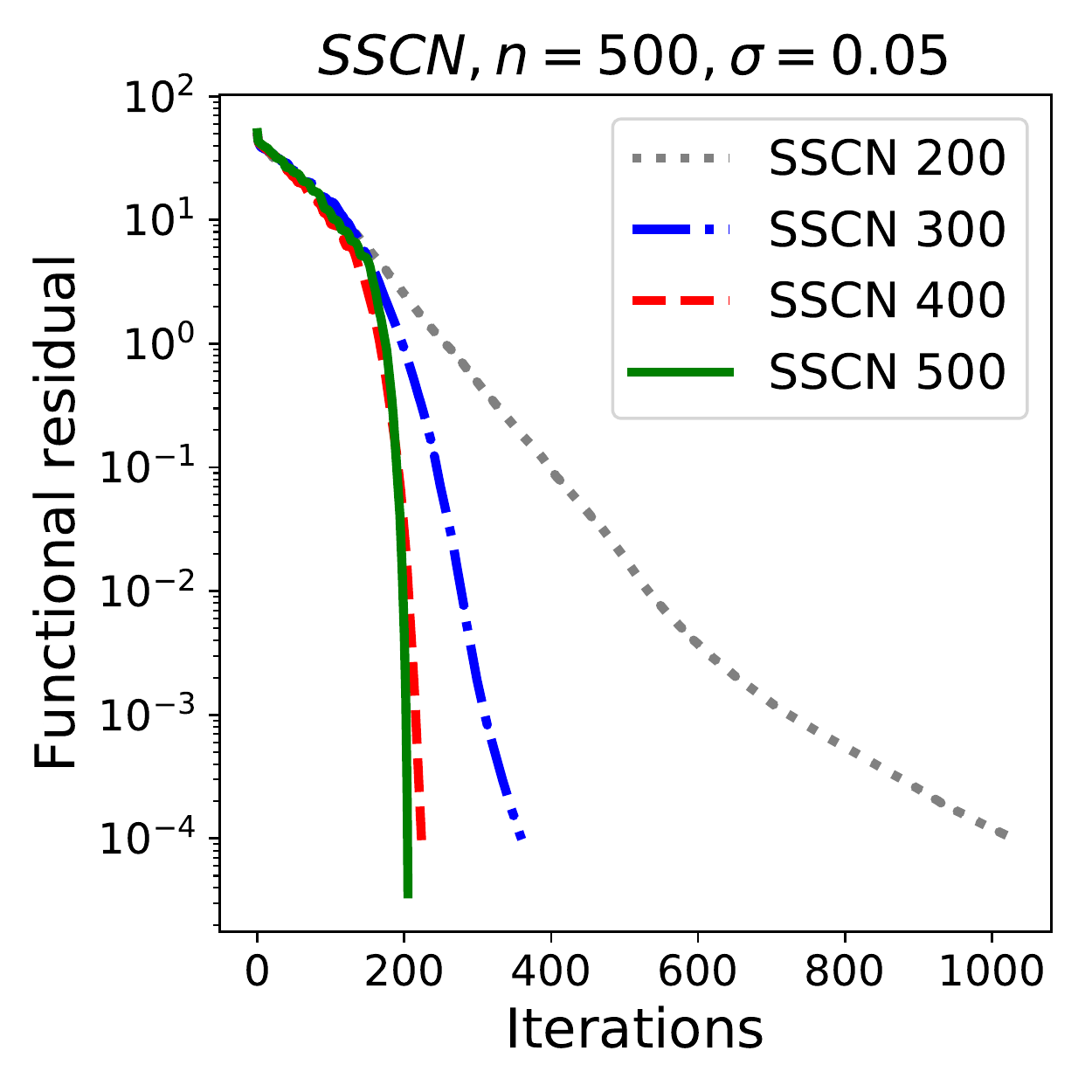}
	\end{minipage}
	\begin{minipage}{0.3\textwidth}
		\centering
		\includegraphics[width =  \textwidth ]{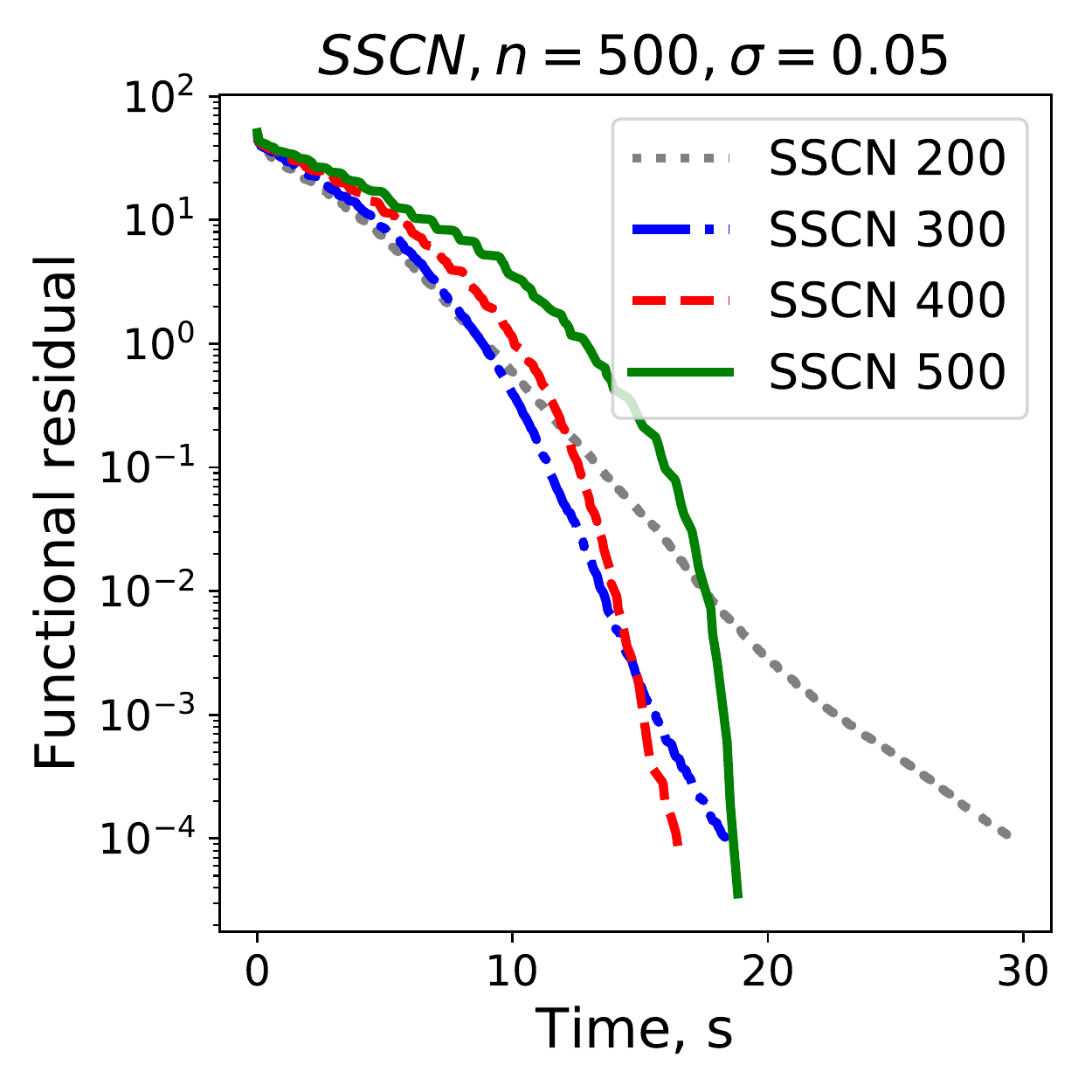}
	\end{minipage}

	\caption{{\tt SSCN} and Coordinate Descent ({\tt CD}) methods, minimizing Log-Sum-Exp function, $d = 500$.} 
	\label{fig:sscn_log_sum_exp_500}
\end{figure}

\begin{figure}[h!]
	\centering
	\begin{minipage}{0.3\textwidth}
		\centering
		\includegraphics[width =  \textwidth ]{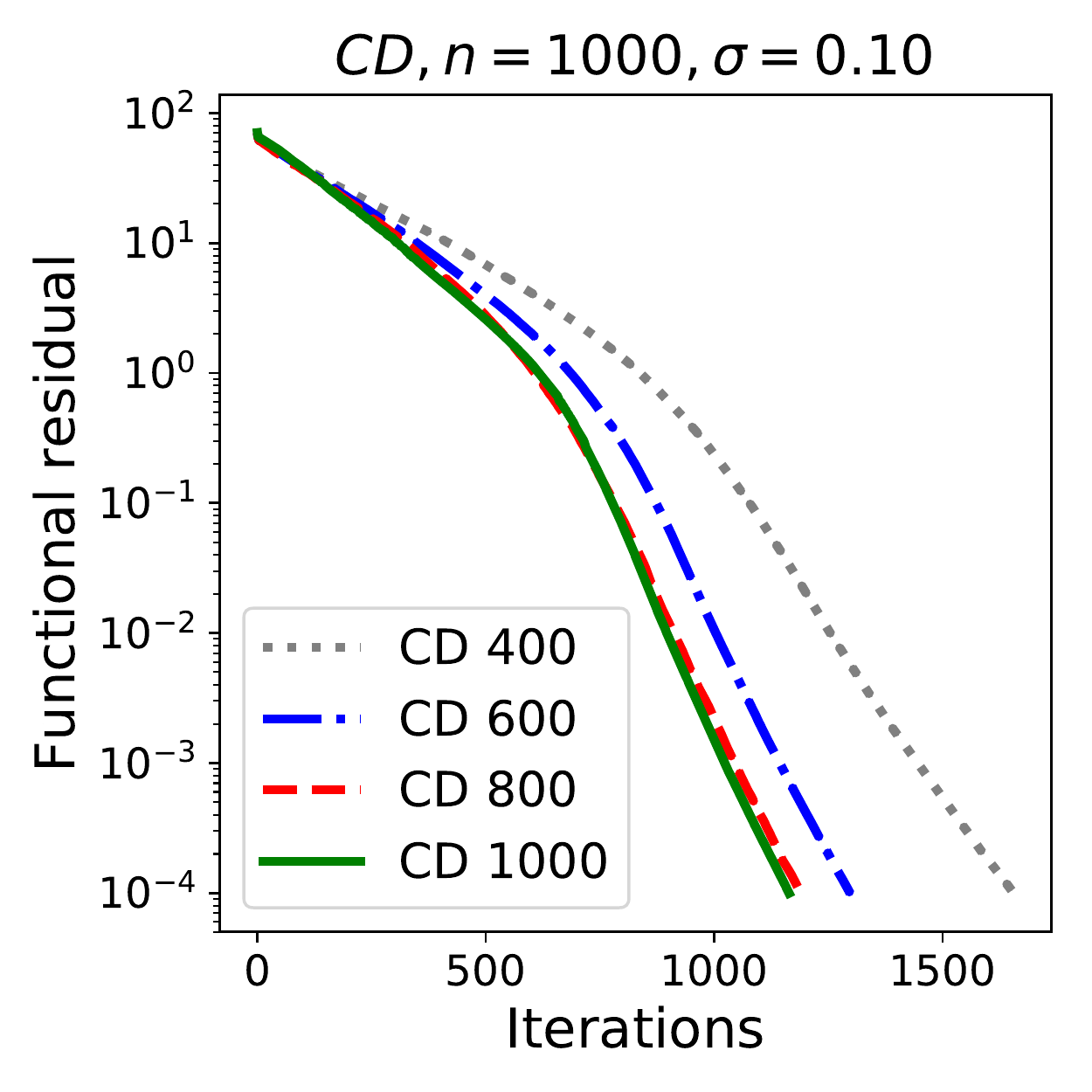}
	\end{minipage}
	\begin{minipage}{0.3\textwidth}
		\centering
		\includegraphics[width =  \textwidth ]{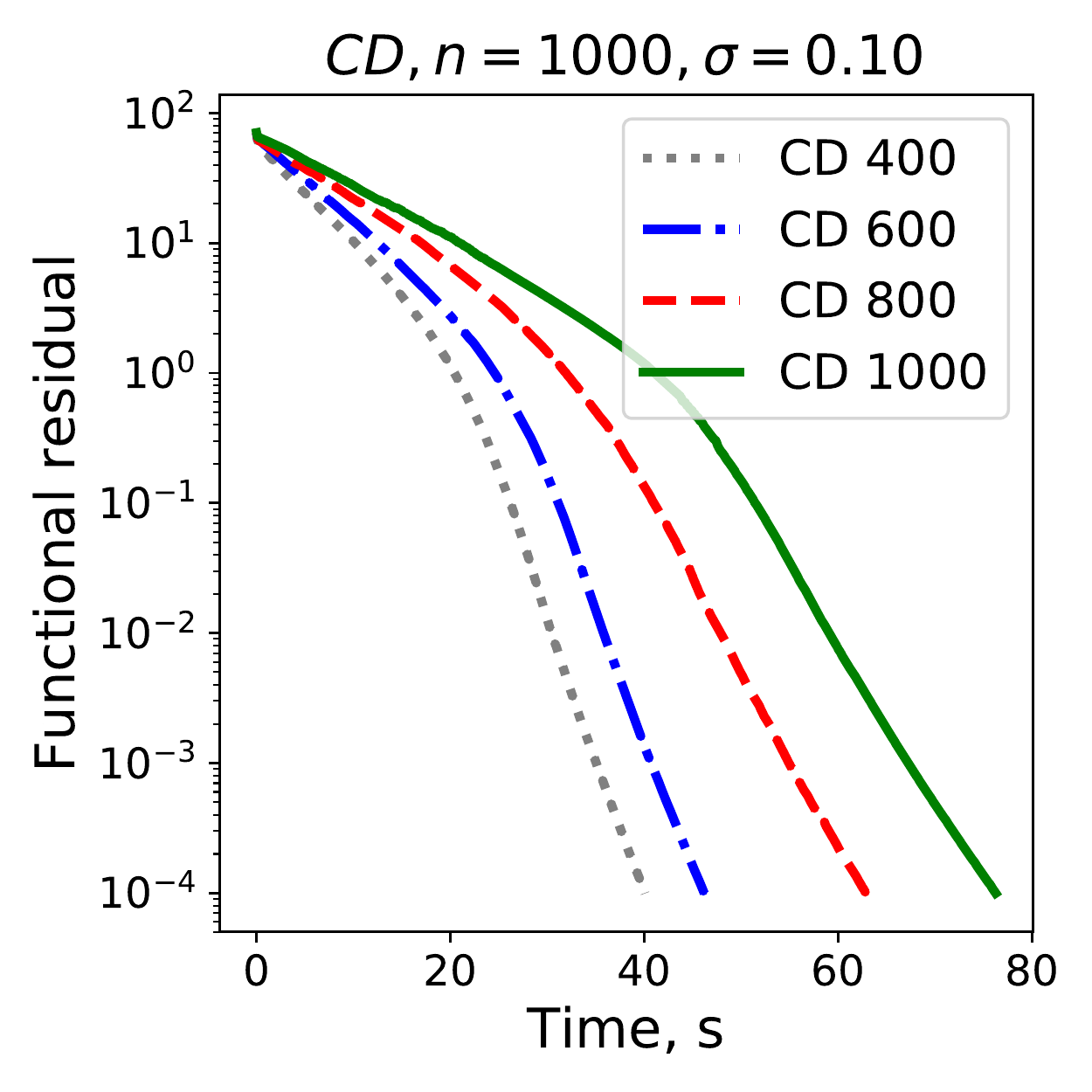}
	\end{minipage}
	\begin{minipage}{0.3\textwidth}
		\centering
		\includegraphics[width =  \textwidth ]{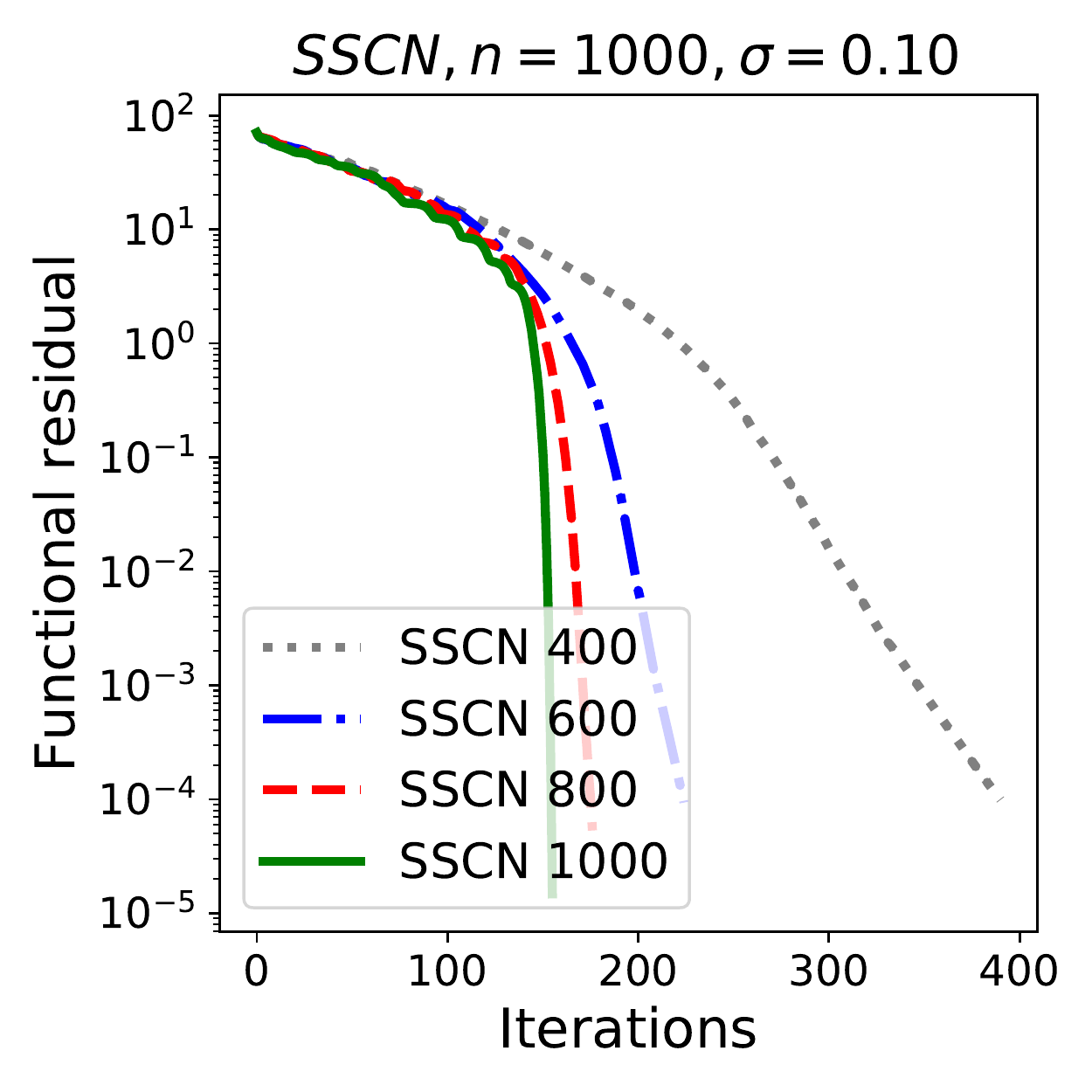}
	\end{minipage}
	\\
	\begin{minipage}{0.3\textwidth}
		\centering
		\includegraphics[width =  \textwidth ]{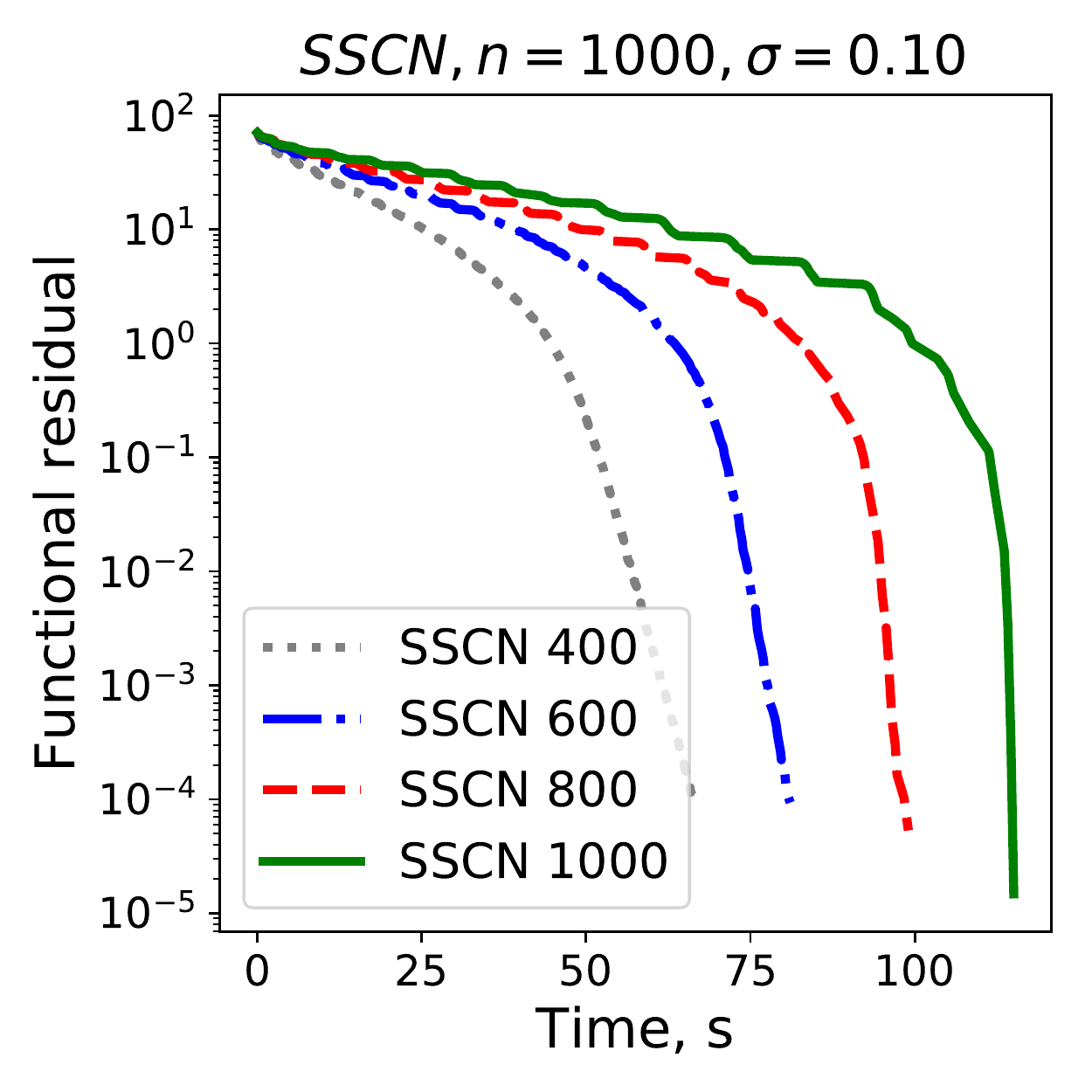}
	\end{minipage}
	\begin{minipage}{0.3\textwidth}
		\centering
		\includegraphics[width =  \textwidth ]{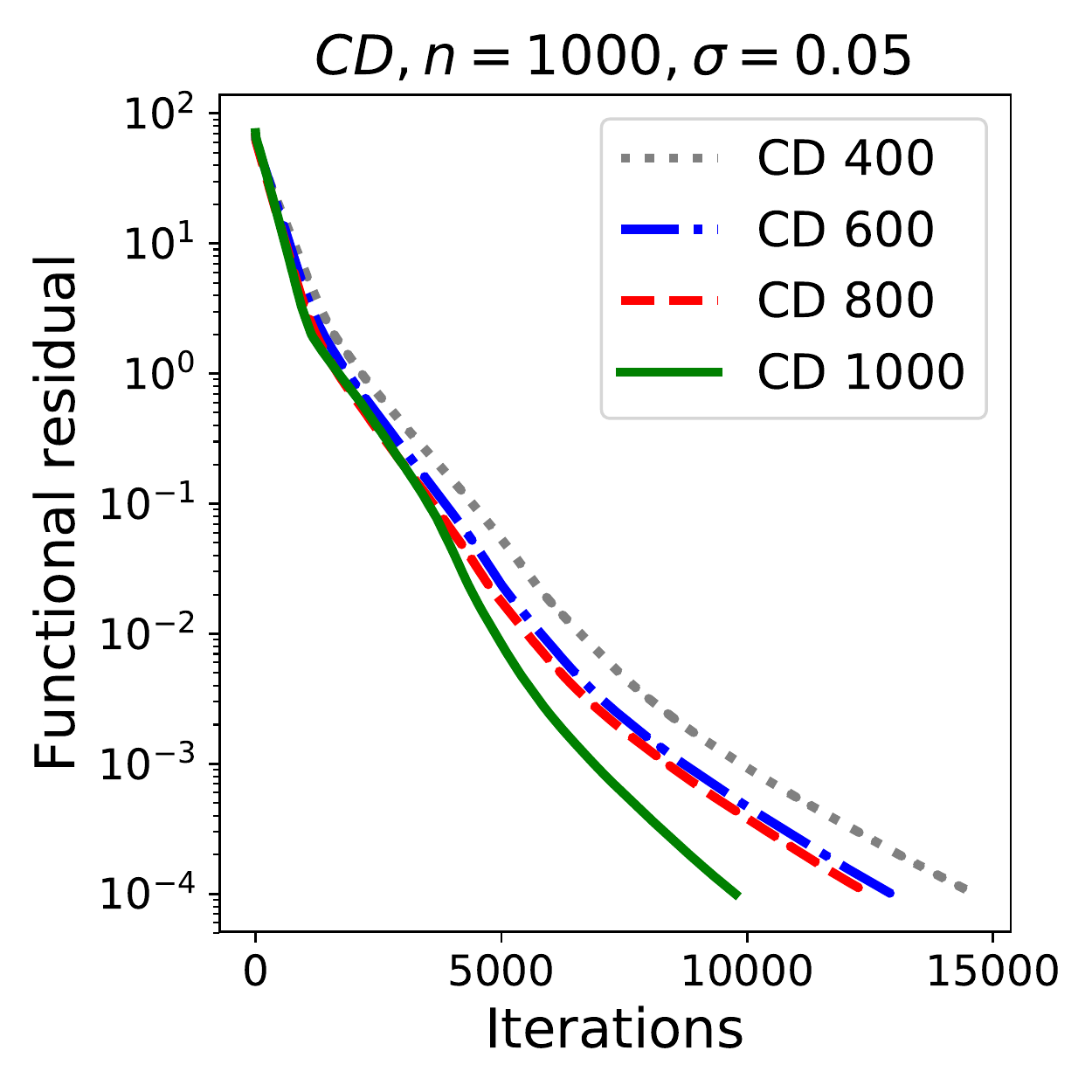}
	\end{minipage}
	\begin{minipage}{0.3\textwidth}
		\centering
		\includegraphics[width =  \textwidth ]{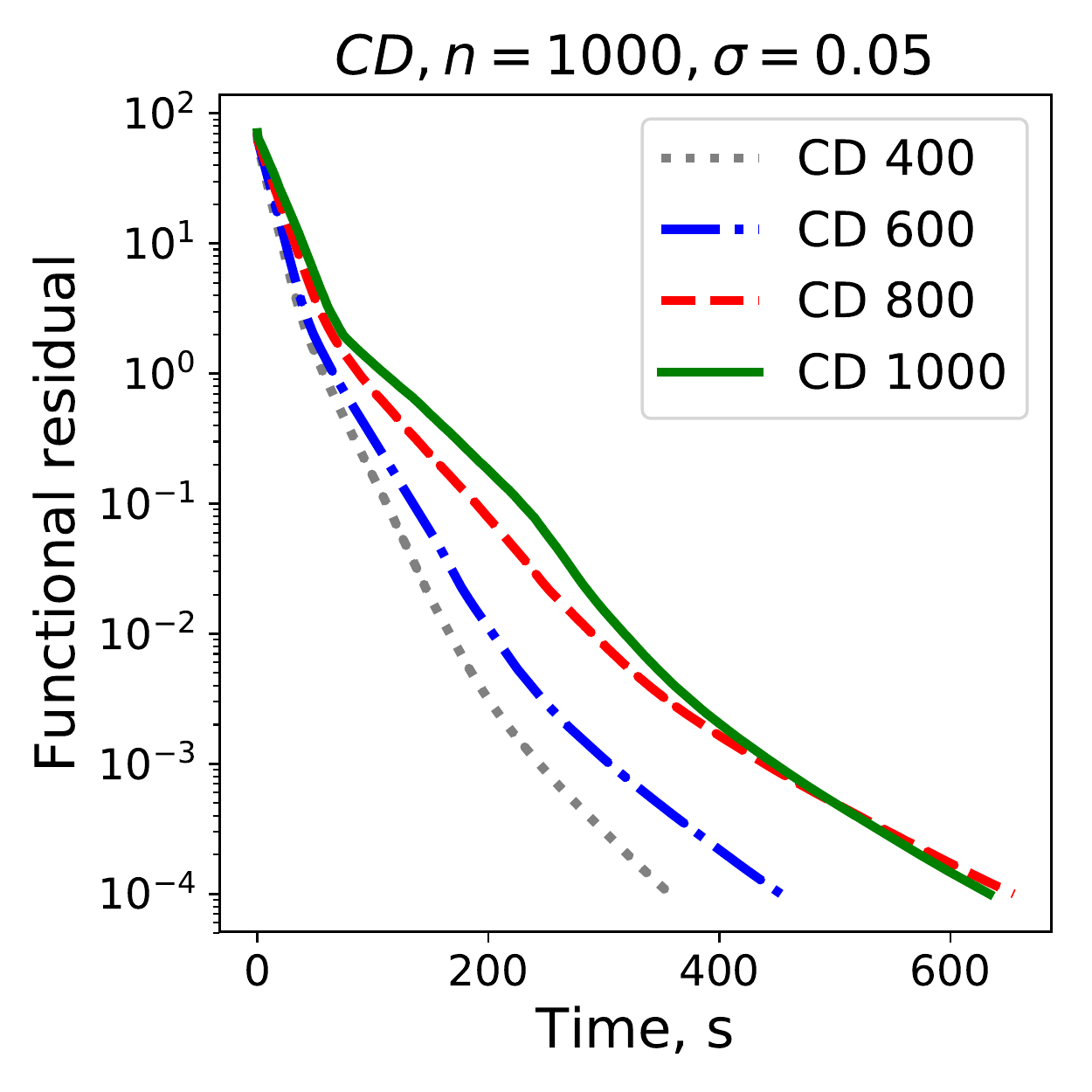}
	\end{minipage}
	\\
	\begin{minipage}{0.3\textwidth}
		\centering
		\includegraphics[width =  \textwidth ]{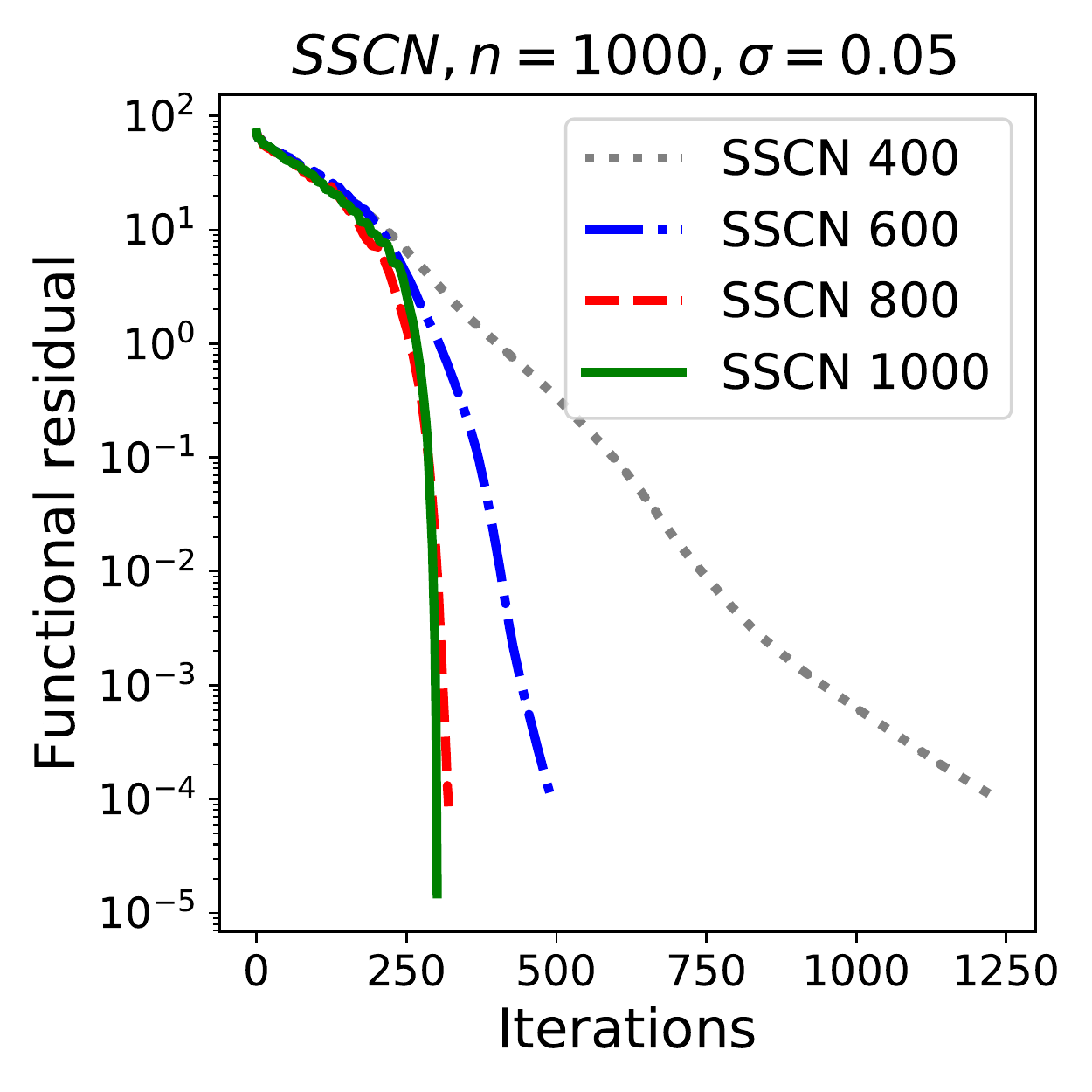}
	\end{minipage}
	\begin{minipage}{0.3\textwidth}
		\centering
		\includegraphics[width =  \textwidth ]{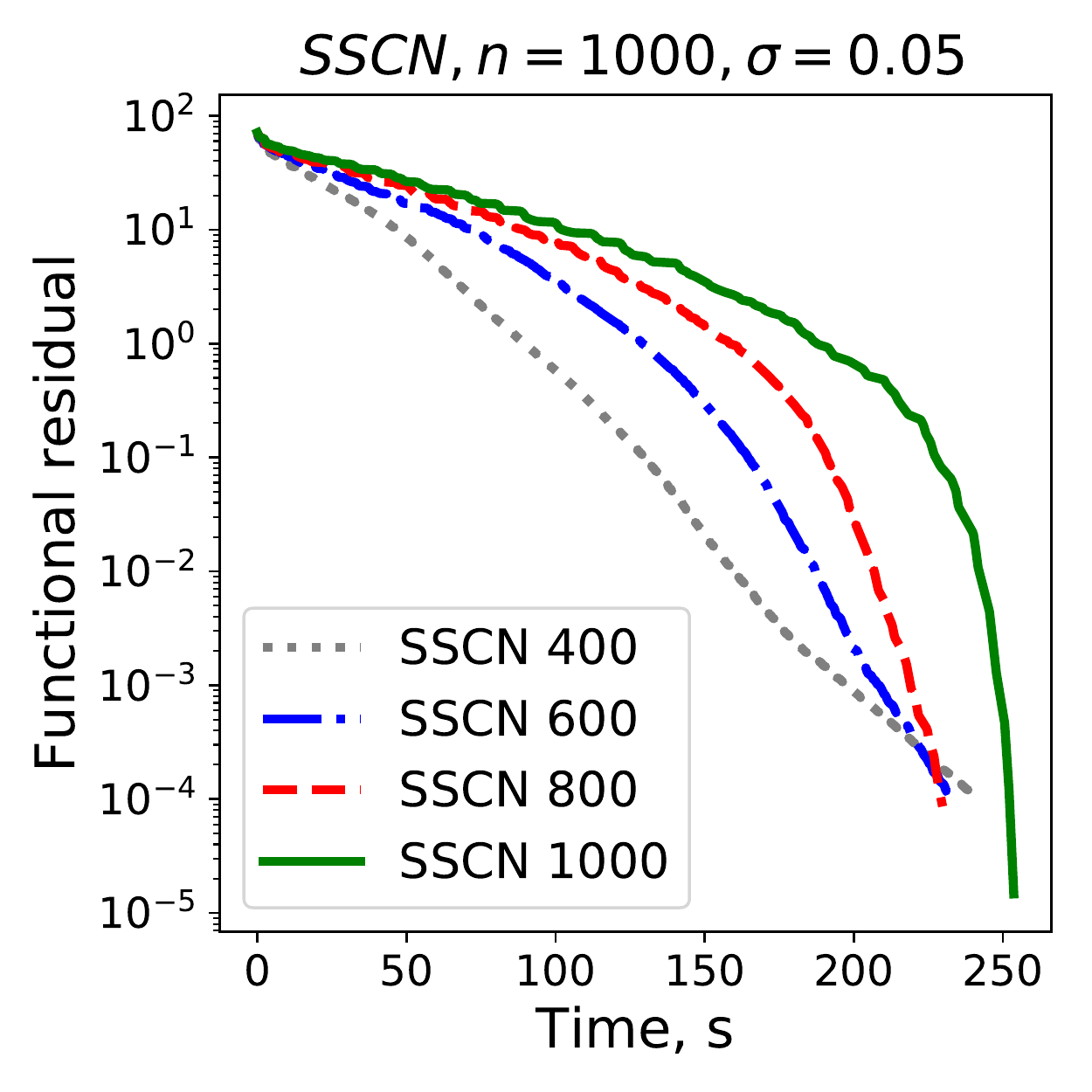}
	\end{minipage}

	\caption{{\tt SSCN} and Coordinate Descent ({\tt CD}) methods, minimizing Log-Sum-Exp function, $d = 1000$.} 
	\label{fig:sscn_log_sum_exp_1000}
\end{figure}

\section{Conclusion}

In this chapter, we have introduced {\acrshort{SSCN}}, which is both a subspace version cubically-regularized Newton method~\cite{nesterov2006cubic}, and a second-order enhancement of stochastic subspace descent~\cite{kozak2019stochastic}. The algorithm enjoys the global convergence to the optimum along with the fast local rates. We believe our method opens up several new avenues for the future research which we list next.

\paragraph{Acceleration.} We believe it would be valuable to incorporate Nesterov's momentum into Algorithm~\ref{alg:sscn_crcd}. Ideally, one would like to get the global rate in between convergence rate of accelerated cubic regularized Newton~\cite{nesterov2008accelerating} and accelerated {\tt CD}~\cite{allen2016even, nesterov2017efficiency}. On the other hand, the local rate (for strongly convex objectives) should recover accelerated sketch-and-project~\cite{tu2017breaking, gower2018accelerated}. If accelerated sketch-and-project is optimal (this is yet to be established), then accelerated {\tt SSCN} (again, given that it recovers accelerated sketch-and-project) would be a locally optimal algorithm as well.

\paragraph{Non-separable $\psi$.} As mentioned in Section~\ref{sec:sscn_setup}, one should not hope for linear convergence of {\tt SSCN} if $\psi$ is not separable, as the iterates can ``jump'' away from the optimum in such case. This issue has been resolved for first-order methods using control variates in Chapter~\ref{sega} via {\tt SEGA} algorithm. Therefore, the development of second-order {\tt SEGA} remains an interesting open problem.

\paragraph{Inexact method.} {\tt SSCN} is applicable in the setup, where function $f$ is accessible via zeroth-order oracle only. In such a case, for any $\mS \in \R^{\tau\times d}$ we can estimate $\nabla_{\mS} f(x)$ and $ \nabla_\mS^2 f(x)$ using $\cO(\tau^2)$ function value evaluations. However, since both $\nabla_{\mS} f(x)$ and $ \nabla_\mS^2 f(x)$ are only evaluated inexactly, a slight modification of our theory is required.

\paragraph{Non-uniform sampling.} Note that our local theory allows for arbitrary non-uniform distribution of $\mS$, which might be potentially exploited. While developing optimal and implementable importance sampling for the local convergence is beyond the scope of this work,\footnote{As this is still an open problem even for sketch-and-project~\cite{gower2015randomized}. }
we sketch several possible sampling strategies that might yield faster convergence.\footnote{This only applies to the local results as the global convergence requires some uniformity; see Assumption~\ref{as:sscn_uniform}.}

\begin{itemize}
\item Let $\Prob(\mS \in \{e_1, e_2, \dots, e_d \})=1$. If we evaluate the diagonal of the Hessian close to optimum (cost $\cO(nd)$ for linear models) and sample proportionally to it, we obtain local linear rate with leading complexity term $\frac{\Tr{\nabla^2 f(x^*)}}{\lambda_{\min}\nabla^2 f(x^*)}$. 

\item It is unclear how to design an efficient importance sampling for minibatch (i.e. $1<\E{\tau(\mS)}<d$) methods. Determinantal point processes (DPP)~\cite{rodomanov2019randomized, mutny2019convergence} were proposed to speed up {\tt SDNA} from~\cite{sdna} (i.e., analogous {\tt CD} with static matrix upper bound) -- we thus believe they might be applicable on our setting too. However, in such a case, one would need to evaluate the whole Hessian close to optimum, which is infeasible for applications where $d$ is large.

\item It is known that {\tt SDNA} (see related literature) is faster than minibatch {\tt CD} under the ESO assumption~\cite{qu2016coordinate1, qu2016coordinate2}. Therefore, we might instead apply minibatch importance sampling for ESO assumption from~\cite{hanzely2018accelerated} (which corresponds to optimizing the upper bound on iteration complexity). Using the mentioned sampling, we only require evaluating the diagonal of Hessian at some point close to optimum, which is of the same cost as computing the full gradient for linear models -- thus is feasible.

\item It is a natural question to ask whether one can speed up the convergence using greedy rule instead of random one. For standard {\tt CD}, greedy rule was shown to have a superior iteration complexity to any randomized rule~\cite{nutini2015coordinate, karimireddy2018efficient}. For simplicity, consider case where $\Prob(\mS \in \{e_1, e_2, \dots, e_d \})=1$. Far from the optimum, (approximate) greedy rule at iteration $k$ chooses index $i = \argmax_{j} | \nabla_j f(x^k)|^{\frac32}M_{e_j}^{-\frac12}$. Close to optimum, if a diagonal of a Hessian was evaluated, (approximate) greedy index would be $\argmax_{j} | \nabla_j f(x^k)|^{2}\nabla_{j,j} f(x)^{-1}$. For linear models, both of the mentioned cases are implementable using the efficient neirest neighbour search~\cite{dhillon2011nearest} with sublinear complexity in terms of $d$.

\end{itemize}

\chapter{Accelerated Stochastic Matrix Inversion:  General Theory and  Speeding up {\tt BFGS} Rules for Faster Second-Order Optimization}
\label{ami}

\graphicspath{{AMI/experiments/}}

A new wave of second-order stochastic methods are being developed nowadays with the aim of solving large scale optimization problems. In particular, many of these new methods are often based on stochastic { \acrshort{BFGS}} updates~\cite{Schraudolph2007,wang2017stochastic,Mokhtari2014, moritz2016linearly,
Byrd2015, curtis2016self, berahas2016multi}.  Another approach to scaling up second-order methods is to use randomized~\emph{sketching} to reduce the dimension, and hence the complexity of the Hessian and the updates involving the Hessian~ \cite{pilanci2017newton,xu2016sub}, or \emph{subsampled} 
  Hessian matrices  when the objective function is a sum of many loss functions~\cite{byrd2011use, BerahasBN17,agarwal2017second, xu2017newton}.

In this chapter we develop a new stochastic accelerated {\tt BFGS} update that can form the backbone of new stochastic quasi-Newton methods. Since the {\tt BFGS} update mechanism which we improve upon is as an optimization routine on its own, this chapter tackles two different objectives in two different domains at the same time. For this reason, the notation will be slightly inconsistent with respect to the rest of the thesis. Specifically, our high-level goal is to minimize smooth function $f$ in variable $w$:
\begin{equation}\label{eq:ami_opt_main}\min_{w\in \R^d} f(w),\end{equation}
while the mentioned {\tt BFGS} subroutine is (as we shall see) a quadratic objective in matrix variable $\mX$ (or $x$ in the vectorized form). Given the (admittedly inconsitent) notation is explained, let us properly motivate our work.

The starting point for developing second-order methods is arguably Newton's method, which performs the iterative process
\begin{align}
 w_{k+1} = w_k - (\nabla^2 f(w_k))^{-1} \nabla f(w_k),
\end{align}
where $\nabla^2 f(w_k)$ and $\nabla f(w_k)$ are the Hessian and gradient of $f$, respectively. However, it is inefficient for solving large scale problems as it requires the computation of the Hessian and then solving a linear system at each iteration. Several methods have been developed to address this issue, based on the idea of approximating the exact update.

\emph{Quasi-Newton} methods, in particular {\tt BFGS}~\cite{broyden1967quasi,fletcher1970new,goldfarb1970family,shanno1970conditioning}, have been the leading optimization algorithm in various fields since the late 60's until the rise of big data, which brought a need for simpler first-order algorithms. It is well known that Nesterov's acceleration \cite{nesterov83} is a reliable way to speed up first-order methods. However until now, acceleration techniques have been applied exclusively to speeding up gradient updates. In this chapter we present an accelerated {\tt BFGS} algorithm, opening up new applications for acceleration. The acceleration in fact comes from an accelerated algorithm for inverting the Hessian matrix.

To be more specific, recall that quasi-Newton rules aim to maintain an estimate of the inverse Hessian $\mX_k$, adjusting it every iteration so that the inverse Hessian acts appropriately in a particular direction, while enforcing symmetry:
\begin{equation} \label{eq:ami_bfgs_const}
\mX_k(\nabla f(w_{k})- \nabla f(w_{k-1})) =w_{k}-w_{k-1}, \qquad \mX_k  =\mX^\top_k.
\end{equation}



A notable research direction is the development of stochastic quasi-Newton methods~\cite{gower:2017}, where the estimated inverse is equal to the true inverse over a subspace:
\begin{equation}
\mX_k\nabla^2f(w_k) \mS_k=\mS_k, \qquad \mX_k=\mX^\top_k,
\label{eq:ami_sbfgs}
\end{equation}
where $\mS_k \in \R^{d \times \tau}$ is a randomly generated matrix.

In fact, \eqref{eq:ami_sbfgs} can be seen as the so called sketch-and-project iteration for inverting $\nabla^2f(w_k)$. In this chapter we first develop the accelerated algorithm for inverting positive definite matrices. As a direct application, our algorithm can be used as a primitive in quasi-Newton methods
 which results in a novel accelerated (stochastic) quasi-Newton method of the type \eqref{eq:ami_sbfgs}. In addition, our acceleration technique can also be incorporated in the classical (non stochastic) {\tt BFGS} method. This results in the accelerated {\tt BFGS} method. Whereas the matrix inversion contribution is accompanied by strong theoretical justifications, this does not apply to the latter. Rather, we verify the effectiveness of this new  accelerated {\tt BFGS} method through numerical experiments.

%

\section{Sketch-and-project for linear systems}
Our accelerated algorithm can be applied to more general tasks than only inverting matrices. In its most general form, it can be seen as an accelerated version of a  \emph{sketch-and-project} method in Euclidean spaces which we present now.
Consider a linear system $\mA x=b$ such that $b\in \Range{\mA}$. One step of the 
sketch-and-project algorithm reads as:
\begin{equation} \label{eq:ami_sap}
x_{k+1}=\argmin_{x} \; \norm{x_k-x}_{\mB}^2 \quad \text{subject to} \quad \mS_k^\top \mA x=\mS_k^\top b,
\end{equation}
where $\norm{x}^2_\mB=\dotprod{\mB x,x}$ for some $\mB\succ 0$ and $\mS_k$ is a random sketching matrix sampled i.i.d at each iteration from a fixed distribution.

Randomized Kaczmarz~\cite{K-1937,strohmer2009randomized} was the first algorithm of this type.
In~\cite{gower2015randomized}, this sketch-and-project algorithm was analyzed in its full generality.
Note that the dual problem of~\eqref{eq:ami_sap} takes the form of a quadratic minimization problem \cite{sda}, and randomized methods such as coordinate descent~\cite{rcdm,wright2015coordinate}, random pursuit~\cite{Stich14,Stich2016} or stochastic dual ascent~\cite{sda} can thus also be captured as special instances of this method. Richt\'{a}rik and Tak\'{a}\v{c} \cite{richtarik2017stochastic} adopt a new point of view through a theory of stochastic reformulations of linear systems. In addition, they consider the addition of a relaxation parameter, as well as  mini-batch and accelerated variants. Acceleration was only achieved for the expected iterates, and not in the L2 sense as we do here. We refer to Richt\'{a}rik and Tak\'{a}\v{c}  \cite{richtarik2017stochastic} for interpretation of sketch-and-project as stochastic gradient descent, stochastic Newton, stochastic proximal point method, and stochastic fixed point method.

Gower \cite{gower:2017} observed that the procedure~\eqref{eq:ami_sap} can also be applied to find the inverse of a matrix. Assume the optimization variable itself is a matrix, $x=\mX$, $b= \mI$, the identity matrix, then sketch-and-project  converges (under mild assumptions) to a solution of $\mA\mX=\mI$. Even the symmetry constraint $\mX =\mX^\top$ can be incorporated into the sketch-and-project framework since it is a linear constraint.

There has been recent development in speeding up the sketch-and-project method using the idea of Nesterov's acceleration \cite{nesterov83}. In~\cite{Liu:2016} an accelerated Kaczmarz algorithm was presented for special sketches of rank one. Arbitrary sketches of rank one where considered in~\cite{Stich14}, block sketches in~\cite{nesterov:2017} and recently, Tu and coathors \cite{tu2017breaking} developed acceleration for special sketching matrices, assuming the matrix $\mA$ is square. This assumption, along with any assumptions on $\mA$, was later dropped in~\cite{martinrichtarikaccell}. Another notable way to accelerate the sketch-and-project algorithm is by using momentum or stochastic momentum~\cite{SMOMENTUM}.

We build on recent work of Richt\'{a}rik and Tak\'{a}\v{c}  \cite{martinrichtarikaccell} and further extend their analysis by studying accelerated sketch-and-project in general Euclidean spaces. This allows us to deduce the result for matrix inversion as a special case. However, there is one additional caveat that has to be considered for the intended application in quasi-Newton methods: ideally, all iterates of the algorithm should be symmetric positive definite matrices. This is not the case in general, but we address this problem by constructing special sketch operators that preserve symmetry and positive definiteness.

Our accelerated sketch-and-project algorithm for solving linear systems in Euclidean spaces is developed and analyzed in Section~\ref{sec:ami_ami_paper}, and is used later in Section~\ref{sec:ami_asqn_pap} to analyze an accelerated sketch-and-project algorithm for matrix inversion. The accelerated sketch-and-project algorithm for matrix inversion is then used to accelerate the {\tt BFGS} update, which in turn leads to the development of an accelerated {\tt BFGS} optimization method. Lastly in Section~\ref{sec:ami_num_pap}, we perform numerical experiments to gain  different insights into the newly developed methods.  Proofs of all results and additional insights can be found in the appendix.

\section{Contributions} \label{sec:ami_contrib}
We now present our main contributions.

\begin{itemize}

\item \textbf{Accelerated Sketch and Project in Euclidean Spaces.} We generalize the analysis of an accelerated version of the sketch-and-project algorithm~\cite{martinrichtarikaccell} to linear operator systems in Euclidean spaces. We provide a self-contained convergence analysis, recovering the original results in a more general setting. 

\item \textbf{Faster Algorithms for Matrix Inversion.} We develop an accelerated algorithm for inverting positive definite matrices. This algorithm can be seen as a special case of the accelerated sketch-and-project in Euclidean space, thus its convergence follows from the main theorem. However, we also provide a different formulation of the proof that is specialized to this setting. Similarly to~\cite{tu2017breaking}, the performance of the algorithm depends on two parameters $\theta$ and $\nu$ that capture spectral properties of the input matrix and the sketches that are used. 
Whilst for the non-accelerated sketch-and-project algorithm for matrix inversion~\cite{gower:2017} the knowledge of these parameters is not necessary, they need to be given as input to the accelerated scheme. When employed with the correct choice of parameters, the accelerated algorithm is always faster than the non-accelerated one. We also provide a theoretical rate for sub-optimal parameters $\theta, \nu$, and we perform numerical experiments to argue the choice of $\theta, \nu$ in practice.

\item \textbf{Randomized Accelerated Quasi-Newton.} 
The proposed iterative algorithm for matrix inversion is designed in such a way that each iterate is a symmetric matrix. This means, we can use the generated approximate solutions as estimators for the inverse Hessian in quasi-Newton methods, which is a direct extension of stochastic quasi-Newton methods. To the best of our knowledge, this yields the first accelerated (stochastic) quasi-Newton method.

\item \textbf{Accelerated Quasi-Newton.}
In the standard {\tt BFGS} method the updates to the Hessian estimate are not chosen randomly, but deterministically. Based on the intuition gained from the accelerated random method, we propose an accelerated scheme for {\tt BFGS}. The main idea is that we replace the random sketching of the Hessian with a deterministic update. The theoretical convergence rates do not transfer to this scheme, but we demonstrate by numerical experiments that it is possible to choose a parameter combination which yields a slightly faster convergence. We believe that the novel idea of accelerating {\tt BFGS} update is extremely valuable, as until now, acceleration techniques were only considered to improve gradient updates.

\end{itemize}


\section{Accelerated stochastic algorithm for matrix inversion \label{sec:ami_ami_paper}}
In this section we propose an accelerated randomized algorithm to solve linear systems in Euclidean spaces. 
This is a very general problem class which comprises the matrix inversion problem as well.
Thus, we will use the result of this section later to analyze our newly proposed matrix inversion algorithm, which we then use to estimate the inverse of the Hessian within a quasi-Newton method.\footnote{Quasi-Newton methods do not compute an exact matrix inverse, rather, they only compute an incremental update. Thus, it suffices to apply \emph{one step} of our proposed scheme per iteration. This will be detailed in Section~\ref{sec:ami_asqn_pap}.}

Let $\cX$ and $\cY$ be finite dimensional Euclidean spaces and 
let $\A:\cX \mapsto \cY$ be  a linear operator. 
Let $L(\cX,\cY)$ denote the space of linear operators that map from $\cX$ to $\cY.$ Consider the linear system  
\begin{equation}\label{eq:ami_system} \A x = b,\end{equation}
where $x\in \cX$ and $b\in \Range{\A}.$
Consequently there exists a solution to the equation~\eqref{eq:ami_system}. In particular,  we aim to find the solution closest to a given initial point $x_0 \in \cX$:
\begin{equation} \label{eq:ami_primal}
x^* \eqdef \arg\min_{x \in \cX} \frac{1}{2}\norm{x-x_0}^2 \quad \mbox{subject to} \quad \A x = b.
\end{equation}
Using the pseudoinverse and Lemma~\ref{lem:ami_pseudo} item~\emph{\ref{it:pseudoleastnorm}}, the solution to~\eqref{eq:ami_primal} is given by
\begin{equation}\label{eq:ami_xsol}
x^* = x_0- \A^{\dagger}(\A x_0 -b) \in x_0+ \Range{\A^*},
\end{equation}
where $\cA^{\dagger}$ and $\cA^*$ denote the pseudoinverse and the adjoint of $\cA,$ respectively.

\subsection{The algorithm}

Let $\cW$ be a  Euclidean space and consider a random linear operator $\cS_k \in L(\cY,\cW)$ chosen from some distribution $\cD$ over $L(\cY,\cW)$ at iteration $k$.
Our method is given in Algorithm~\ref{alg:ami_SketchJac}, where $\cZ_k \in L(\cX)$ is a random linear operator given by the following compositions
\begin{equation}\label{eq:ami_Z}
\cZ_k =\cZ(\cS_k) \eqdef \A^*\cS_k^*(\cS_k\A\A^*\cS_k^*)^{\dagger}\cS_k \A.
\end{equation}
The updates of variables $g_k$ and $x_{k+1}$ on lines~8 and~9, respectively, correspond to what is known as the \emph{sketch-and-project} update:
\begin{equation}  \label{eq:ami_sk}
 x_{k+1} =   \arg\min_{x \in \cX} \frac{1}{2}\norm{x-y_k}^2 \quad \text{subject to} \quad \cS_k \A x = \cS_k b, 
\end{equation}
which can also be written as the following operation 
\begin{equation}\label{eq:ami_IZprojres}
x_{k+1} - x_*  = (\cI- \cZ_k)(y_k -x_*), 
\end{equation}
where $\cI$ is the identity operator. This follows from the fact that $b \in \Range{\A}$, together with item~\ref{it:pseudoTTdagT} of Lemma~\ref{lem:ami_pseudo}. 
Furthermore, note that the adjoint $\A^*$ and the pseudoinverse in Algorithm~\ref{alg:ami_SketchJac} are taken with respect to the norm  in~\eqref{eq:ami_primal}. 


\begin{algorithm}[!h]
\begin{algorithmic}[1]
\State \textbf{Parameters:} $\theta, \nu >0$, ${\cal D}$ = distribution over random linear operators.
\State Choose  $x_0\in \cX$ and set $v_0 = x_0$, $\beta  =1 - \sqrt{\frac{\theta}{\nu}},$ $\gamma = \sqrt{\frac{1}{\theta \nu}},$ $\eta = \frac{1}{1+\gamma\nu}.$  
\For{$k =  0, 1, 2,\dots$}
\State $y_k = \eta v_k + (1-\eta) x_k$ 
	\State Sample an independent copy $\cS_k\sim {\cal D}$
	\State 
	$g_k = \A^*\cS_k^*(\cS_k\A\A^*\cS_k^*)^{\dagger}\cS_k(\A y_k -b)=\cZ_k(y_k -x_*)$  \label{ln:stochgrad}	
\State $x_{k+1} =y_k -g_k$\label{ln:xupdate}
	\State $v_{k+1} = \beta v_k +(1-\beta)y_k - \gamma g_k$ \label{ln:vupdate}
\EndFor
\end{algorithmic}
\caption{Accelerated  Sketch-and-Project  for solving \eqref{eq:ami_sk}  \cite{martinrichtarikaccell}}
\label{alg:ami_SketchJac}
\end{algorithm}

Algorithm~\ref{alg:ami_SketchJac} was first proposed and analyzed by Richt\'{a}rik and Tak\'{a}\v{c}  \cite{martinrichtarikaccell} for the special case when $\cX = \R^d$ and $\cY = \R^m$. Our contribution here is in extending the algorithm and analysis to the more abstract setting of Euclidean spaces.  In addition, we provide some further extensions of this method in Sections~\ref{sec:ami_omega} and~\ref{sec:ami_alpha}, allowing for a non-unit stepsize  and variable $\eta$, respectively.

\subsection{Key assumptions and quantities}
Denote $\cZ=\cZ(\cS)$ for $\cS\sim \cD$. Assume that the \emph{exactness property} holds
\begin{equation}
\Null{\A} = \Null{\E{\cZ}} ; 
 \label{eq:ami_exactness}
\end{equation}
this is also equivalent to $\Range{\A^*} = \Range{\E{\cZ}}$. The exactness assumption is of key importance in the sketch-and-project framework, and indeed it is not very strong. For example, it holds for the matrix inversion problem with every sketching strategy we consider. 
We further assume that $\A \neq 0$ and $\E{\cZ}$ is finite. First we collect a few observation on the $\cZ$ operator

\begin{lemma}\label{lem:ami_Z}
The $\cZ$ operator~\eqref{eq:ami_Z} is a self-adjoint positive projection. Consequently $\E{\cZ}$ is a self-adjoint positive operator.
\end{lemma}

The two parameters that govern the acceleration are
\begin{equation}
 \label{eq:ami_mu+nu}
\theta  \eqdef   \inf_{x \in \Range{\A^*}} \frac{\dotprod{\E{\cZ}x,x}}{\dotprod{x,x}}, \qquad \quad 
\nu  \eqdef   \sup_{x \in \Range{\A^*}} \frac{\dotprod{\E{\cZ\E{\cZ}^\dagger \cZ}x,x}}{\dotprod{\E{\cZ}x,x}}.
\end{equation}

The supremum in the definition of $\nu$ is well defined due to the exactness assumption together with $\A \neq 0.$

\begin{lemma} \label{lem:ami_Z_bounds}
We have 
\begin{equation} \label{eq:ami_nubnds}
1\quad \leq \quad\nu\quad \leq \quad  \frac{1}{\theta} \quad = \quad  \norm{\E{\cZ}^{\dagger}}.
\end{equation}
Moreover, if $\Range{\A^*}=\cX$, we have
\begin{equation} \label{eq:ami_nu_lower}
\frac{\Rank{\A^*} }{\E{\Rank{\cZ}}}\leq \nu.
\end{equation}
\end{lemma}

\subsection{Convergence and change of the norm}
 For a positive self-adjoint $\cG \in L(\cX)$  and $x \in \cX$ let $\norm{x}_\cG \eqdef \sqrt{\dotprod{x,x}_\cG} \eqdef \sqrt{\dotprod{\cG x,x}}$. We now informally state the convergence rate of Algorithm \ref{alg:ami_SketchJac}. Theorem~\ref{theo:conv} generalizes  the main theorem from \cite{martinrichtarikaccell} to linear systems in Euclidean spaces.

\begin{theorem}\label{theo:conv} Let $x_k, v_k$ be the random iterates of Algorithm \ref{alg:ami_SketchJac}. Then
\[
 \E{\norm{v_{k} -x_*}_{\E{\cZ}^\dagger}^2 +\frac{1}{\theta}\norm{x_{k} -x_*}^2 }  \leq \left(1 - \sqrt{\frac{\theta}{\nu}} \right)^k \E{\norm{v_0 -x_*}_{\E{\cZ}^\dagger}^2 + \frac{1}{\theta}\norm{x_0-x_*}^2}.
\]
\end{theorem}

This theorem shows the accelerated Sketch-and-Project algorithm converges linearly with a rate of $ \bigl( 1-\sqrt{\frac{\theta}{\nu}} \bigr),$ which translates to a total of 
$O(\sqrt{\nu/\theta}\log\left( 1/\epsilon\right) )$ iterations to bring the given error in Theorem~\ref{theo:conv} below $\epsilon >0.$ This is in contrast with the non-accelerated 
Sketch-and-Project algorithm which requires $O((1/\theta)\log\left( 1/\epsilon\right) )$ iterations, as shown in~\cite{gower2015randomized} for solving linear systems. From~\eqref{eq:ami_nubnds}, we have the bounds
$1/\sqrt{\theta}  \leq  \sqrt{\nu/\theta}  \leq 1/\theta.$
On one extreme, this  inequality shows that the iteration complexity of the accelerated algorithm is at least as good as its non-accelerated counterpart. On the other extreme, the accelerated algorithm might require as little as the square root of the number of iterations of its non-accelerated counterpart. Since the cost of a single iteration of the accelerated algorithm is of the same order as the non-accelerated algorithm, this theorem shows that acceleration can offer a significant speed-up, which is verified numerically in Section~\ref{sec:ami_num_pap}. 
It is also possible to get the convergence rate of accelerated sketch-and-project where projections are taken with respect to a different weighted norm. For technical details, see Section~\ref{sec:ami_change_norm} of the Appendix. 

\subsection{Coordinate sketches with convenient probabilities \label{sec:ami_convenient_munu}}
Let us consider a simple example in the setting for Algorithm~\ref{alg:ami_SketchJac} where we can understand parameters $\theta, \nu$. 
 In particular, consider a linear system $\mA x=b$ in $\R^d$ where $\mA$ is symmetric positive definite. 
 
\begin{corollary} \label{cor:ami_sss}
 Choose $\mB=\mA$ and $\mS=e_i$ with probability proportional to $\mA_{i,i}$. Then
\begin{equation}
\theta=\frac{\lambda_{\min}(\mA)}{\Tr{\mA}} =: \theta^P \quad \mbox{and} \quad \nu = \frac{\Tr{\mA}}{\min_i \mA_{i,i}} =: \nu^P
\label{eq:ami_munu_conv_paper}
\end{equation}
and therefore the convergence rate given in Theorem~\ref{theo:conv} for the accelerated algorithm is 
\begin{equation}\label{eq:ami_89g9d8g098f}
\biggl( 1- \sqrt{\frac{\theta}{\nu}}\biggr)^k\quad  =\quad   \left( 1-\frac{\sqrt{\lambda_{\min}(\mA) \min_i \mA_{i,i}}}{\Tr{\mA}} \right) ^k.
\end{equation}
\end{corollary}

Rate \eqref{eq:ami_89g9d8g098f} of our accelerated method is to be contrasted with the rate of the non-accelerated method: 
$
(1- \theta )^k  =  ( 1- \lambda_{\min}(\mA)/ \Tr{\mA}) ) ^k.
$
Clearly, we gain from acceleration if the smallest diagonal element of $\mA$ is significantly larger than the smallest eigenvalue.

 In fact, parameters $\theta^P, \nu^P$ above are the correct choice for the matrix inversion algorithm, when symmetry is not enforced, as we shall see later.  Unfortunately, we are not able to estimate the parameters while enforcing symmetry for different sketching strategies. We dedicate a section in numerical experiments to test, if the parameter selection \eqref{eq:ami_munu_conv_paper} performs well under enforced symmetry and different sketching strategies, and also how one might safely choose $\theta,\nu$ in practice.


\section{Accelerated stochastic {\tt BFGS} update \label{sec:ami_asqn_pap}}
The update of the inverse Hessian used in quasi-Newton methods (e.g., in {\tt BFGS}) can be seen as a sketch-and-project update applied to the linear system $\mA\mX=\mI$, while $\mX= \mX^\top$ is enforced, and where $\mA$ denotes and approximation of the Hessian.
In this section, we present an accelerated version of these updates. We provide two different proofs: one based on Theorem~\ref{theo:conv} and one based on vectorization.   By mimicking the updates of the accelerated stochastic {\tt BFGS} method for inverting matrices, we determine a heuristic for accelerating the classic deterministic {\tt BFGS} update. We then incorporate this acceleration into the classic {\tt BFGS} optimization method and show that the resulting  algorithm can offer a speed-up of the standard {\tt BFGS} algorithm.
 

\subsection{The {\tt AMI} algorithm}
Consider the symmetric positive definite matrix $\mA \in \R^{d \times d}$ and the following projection problem
\begin{equation}  \label{eq:ami_primalqN}
 \mA^{-1} =  \arg\min_{\mX} \; \norm{\mX}_{F(\mA)}^2 \quad \text{subject to} \quad \mA\mX= \mI, \quad \mX = \mX^\top,
\end{equation}
where $\norm{\mX}_{F(\mA)}\eqdef \Tr{\mA\mX^\top \mA\mX} = \norm{\mA^{1/2}\mX \mA^{1/2}}_F^2.$
This projection problem can be cast as an instantiation of the general projection problem~\eqref{eq:ami_primal}. Indeed, we need only note that the constraint in~\eqref{eq:ami_primalqN} is linear and equivalent to
$\mathcal{\mA}(\mX) \eqdef \left(\begin{smallmatrix} \mA\mX \\ \mX -\mX^\top \end{smallmatrix} \right)= 
\left(\begin{smallmatrix} \mI \\ 0 \end{smallmatrix}\right). $
The matrix inversion problem can be efficiently solved using sketch-and-project with a symmetric sketch~\cite{gower:2017}. The symmetric sketch is given by
$ \mathcal{S}_k\mathcal{A}(\mX)  = \left(\begin{smallmatrix} \mS_k^\top \mA\mX \\ \mX -\mX^\top \end{smallmatrix}\right), $
where $\mS_k \in \R^{d \times \tau}$ is a random matrix drawn from a distribution $\mathcal{D}$ and $\tau \in \N.$ The resulting sketch-and-project method is as follows
\begin{equation}  \label{eq:ami_primalqNX}
 \mX_{k+1}=   \arg\min_{\mX} \; \norm{\mX - \mX_k}_{F(\mA)}^2 \quad   \text{subject to} \quad \mS_k^\top \mA \mX = \mS_k^\top, \quad \mX = \mX^\top,
\end{equation}
 the closed form solution of which is 
\begin{equation}
\mX_{k+1}  =\mS_k(\mS_k^\top \mA \mS_k)^{-1}\mS_k^\top + \left(\mI-\mS_k(\mS_k^\top \mA\mS_k)^{-1}\mS_k^\top \mA\right) \mX_{k} \left(\mI -\mA\mS_k(\mS_k^\top \mA\mS_k)^{-1}\mS_k^\top  \right).\label{eq:ami_qunac} 
\end{equation}
By observing that~\eqref{eq:ami_qunac} is the sketch-and-project algorithm applied to a linear operator equation, we have constructed an accelerated version in Algorithm~\ref{alg:ami_qn}. We can also apply Theorem~\ref{theo:conv} to prove that 
 Algorithm~\ref{alg:ami_qn} is indeed accelerated.
\begin{theorem}\label{theo:qn}
Let $\mL^k\eqdef \norm{\mV_k -\mA^{-1}}_{\mC}^2 + \frac{1}{\theta}\norm{\mX_k-\mA^{-1}}^2_{F(\mA)}$. The iterates of Algorithm~\ref{alg:ami_qn} satisfy 
\begin{equation}\label{eq:ami_qnaccconv}
 \E{\mL_{k+1}} \leq \left(1 - \sqrt{\frac{\theta}{\nu}} \right) \E{\mL_k},
\end{equation} 
where 
$\norm{\mX}_{\mC}^2 = \Tr{\mA^{1/2}\mX^\top \mA^{1/2} \E{\mZ}^\dagger \mA^{1/2} \mX \mA^{1/2}}.$
Furthermore,
\begin{equation}
\theta  \eqdef   \inf_{\mX \in \R^{d\times d}} \frac{\dotprod{\E{\mZ}\mX,\mX}}{\dotprod{\mX,\mX}}= \lambda_{\min}(\E{\bigZ}), \qquad \nu \eqdef   \sup_{\mX \in \R^{d\times d}} \frac{\dotprod{\E{\mZ\E{\mZ}^\dagger \mZ}\mX,\mX}}{\dotprod{\E{\mZ}\mX,\mX}},\label{eq:ami_nuqn} 
\end{equation}
where 
\begin{equation}
\bigZ \eqdef  \mI\otimes \mI- (\mI-\mP)\otimes(\mI-\mP) , \qquad \mP \eqdef \mA^{1/2}\mS(\mS^\top \mA\mS)^{-1}\mS^\top \mA^{1/2}, \label{eq:ami_bigz}
\end{equation}
and
$\mZ: \mX \in \R^{d\times d} \rightarrow \R^{d\times d}$ is given by
$\mZ(\mX) = \mX - \left(\mI-\mP\right) \mX\left(\mI -\mP \right) = \mX \mP + \mP\mX(\mI-\mP).$
Moreover, $ 2\lambda_{\min}(\E{\mP})\geq \lambda_{\min}(\E{\bigZ})\geq \lambda_{\min}(\E{\mP}).$
\end{theorem}

Notice that preserving symmetry yields $\theta =\lambda_{\min}(\E{\bigZ})$ , which can be up to twice as large as  $\lambda_{\min}(\E{\mP})$, which is the value of the $\theta$ parameter of the method without preserving symmetry. 
This improved rate is new, and was not present in the algorithm's debut publication~\cite{gower:2017}.  In terms of parameter estimation, once symmetry is not preserved, we fall back onto the setting from Section~\ref{sec:ami_convenient_munu}. Unfortunately, we were not able to quantify the effect of enforcing symmetry on the parameter $\nu$.

\begin{algorithm}[!h]
\begin{algorithmic}[1]
\State \textbf{Parameters:} $\theta, \nu >0$, ${\cal D}$ = distribution over random linear operators.
\State Choose  $\mX_0\in \cX$ and set $\mV_0 = \mX_0$, $\beta  =1 - \sqrt{\frac{\theta}{\nu}},$ $\gamma = \sqrt{\frac{1}{\theta \nu}},$ $\eta = \frac{1}{1+\gamma\nu}$  
\For {$k =  0, 1, 2, \dots$}
\State $\mY_k = \eta \mV_k + (1-\eta) \mX_k$ 
	\State Sample an independent copy $S\sim {\cal D}$
	\State $\mX_{k+1} = \mY_k + (\mY_k\mA-\mI)\mS(\mS^\top \mA\mS)^{-1}\mS^\top  
- \mS(\mS^\top \mA\mS)^{-1}\mS^\top \mA\mY_k$ \\
		\qquad \qquad $+ \mS(\mS^\top \mA\mS)^{-1}\mS^\top \mA\mY_k\mA\mS(\mS^\top \mA\mS)^{-1}\mS^\top $  \label{ln:stochgradqn}	
	\State $\mV_{k+1} = \beta \mV_k +(1-\beta)\mY_k - \gamma (\mY_k-\mX_{k+1})$ \label{ln:vupdateqn}
\EndFor
\end{algorithmic}
\caption{{\tt AMI} (Accelerated  {\tt BFGS} Matrix Inversion)}
\label{alg:ami_qn}
\end{algorithm}

\subsection{Vectorizing -- a different insight}
Define ${\bf Vec}: \R^{d\times d}\rightarrow \R^{d^2}$ to be a vectorization operator of column-wise stacking and denote $x\eqdef \Vect{\mX}$. It can be shown that the sketch-and-project operation for matrix inversion \eqref{eq:ami_primalqNX} is equivalent to 
\begin{eqnarray*} 
 x_{k+1} &=&   \arg\min_{x\in \R^{d^2}}\; \norm{x- x_k}_{\mA\otimes \mA}^2  \\
  \text{subject to} \quad  (\mI\otimes \mS_k^\top) (\mI\otimes \mA) x &=& (\mI\otimes \mS_k^\top) \Vect{\mI}, \; \mC x=0,
\end{eqnarray*}
where $\mC$ is defined so that  $\mC x=0$ if and only if $\mX~=~\mX^\top$. The above is a sketch-and-project update for a linear system in $\R^{d^2}$, which allows to obtain an alternative proof of Theorem \ref{theo:qn}, without using our results from Euclidean spaces. The details are provided in Section~\ref{sec:ami_alternate} of the Appendix.

\subsection{Accelerated {\tt BFGS} as an optimization algorithm}
\label{sec:ami_accBFGSmethod}
 As a tweak in the stochastic {\tt BFGS} allows for a faster estimation of Hessian inverse and therefore more accurate steps of the method, one might wonder if a equivalent tweak might speed up the standard, deterministic {\tt BFGS} algorithm for solving~\eqref{eq:ami_opt_main}. The mentioned tweaked version of standard {\tt BFGS} is proposed as Algorithm~\ref{alg:ami_bfgs_opt}. We do not state a convergence theorem for this algorithm---due to the deterministic updates the analysis is currently elusive---nor propose to use it as a default solver, but we rather introduce it as a novel idea for accelerating optimization algorithms. We leave theoretical analysis for the future work. For now, we perform several numerical experiments, in order to understand the potential and limitations of this new method.

\begin{algorithm}[!h]
\begin{algorithmic}[1]
\State \textbf{Parameters:} $\theta, \nu >0$, 
stepsize $\alpha$.
\State Choose $\mX_0\in \cX$, $w_0$ and set $\mV_0 = \mX_0$, $\beta  =1 - \sqrt{\frac{\theta}{\nu}},$ $\gamma = \sqrt{\frac{1}{\theta \nu}},$ $\eta = \frac{1}{1+\gamma\nu}.$  
\For {$k =  0, 1, 2, \dots$}
\State $w_{k+1}=w_{k}-\alpha \mX_{k}\nabla f(w_{k})$
\State $s_k = w_{k+1}-w_k$, \quad $\zeta_k = \nabla f(w_{k+1})- \nabla f(w_k)$
\State $\mY_k = \eta \mV_k + (1-\eta) \mX_k$ 
\State $\mX_{k+1} =   \frac{\delta_k \delta_k^\top}{\delta_k^\top \zeta_k}+ \left(\mI-\frac{\delta_k \zeta_k^\top}{\delta_k^\top \zeta_k} \right) \mY_{k} \left(\mI -\frac{\zeta_k\delta_k^\top}{\delta_k^\top \zeta_k}   \right) $ \label{ln:updateX}
	\State $\mV_{k+1} = \beta \mV_k+(1-\beta)\mY_k-\gamma (\mY_k - \mX_{k+1})$ 
\EndFor
\end{algorithmic}
\caption{{\tt BFGS} method with accelerated {\tt BFGS} update for solving \eqref{eq:ami_opt_main}}
\label{alg:ami_bfgs_opt}
\end{algorithm}

To better understand Algorithm~\ref{alg:ami_bfgs_opt}, recall that the {\tt BFGS} updates an estimate of the inverse Hessian via
\begin{equation}
 \mX_{k+1}=\argmin_{\mX} \; \|\mX-\mX_k \|^2_{F(A)} \quad \text{subject to} \quad \mX\delta_k=\zeta_k ,\, \mX=\mX^\top, 
\end{equation}
where $\delta_k = w_{k+1}-w_k $ and $ \zeta_k = \nabla f(w_{k+1})- \nabla f(w_k)$. The above has the following  closed form solution
$
\mX_{k+1}=\frac{\delta_k \delta_k^\top}{\delta_k^\top \zeta_k}+ \left(\mI-\frac{\delta_k \zeta_k^\top}{\delta_k^\top \zeta_k} \right) \mX_{k} \left(\mI -\frac{\zeta_k\delta_k^\top}{\delta_k^\top \zeta_k}   \right).
$
This update appears on line~\ref{ln:updateX} of Algorithm~\ref{alg:ami_bfgs_opt} with the difference being that it is applied to a matrix $Y_k$.

\section{Experiments \label{sec:ami_num_pap}}

We perform extensive numerical experiments to bring additional insight to both the performance of and to parameter selection for Algorithms~\ref{alg:ami_qn} and~\ref{alg:ami_bfgs_opt}. We first test our  accelerated matrix inversion algorithm, and subsequently perform experiments related to Section~\ref{sec:ami_accBFGSmethod}. 

\subsection{Accelerated matrix inversion} \label{sec:ami_ex-BFGS-opt}

We consider the problem of inverting a symmetric positive matrix $A$. 
We focus on a few particular choices of matrices $A$ (specified when describing each experiment), that differ in their eigenvalue spectra.
Three different sketching strategies are studied: Coordinate sketches with convenient probabilities ($\mS=e_i$ with probability proportional to $\mA_{i,i}$), coordinate sketches with uniform probabilities ($\mS=e_i$ with probability $\frac1n$) and Gaussian sketches ($S\sim \cN(0, \mI)$).
As matrices to be inverted, we use both artificially generated matrices with the access to the spectrum and also Hessians of ridge regression problems from LibSVM.

We compare the speed of the accelerated method with pre-computed estimates of the parameters $\theta, \nu$ to the nonaccelerated method. The pre-computed estimates of $\theta^P, \nu^P$ are set as per \eqref{eq:ami_munu_conv_paper}:
\begin{equation*}
\theta^P=\frac{\lambda_{\min}(\mA)}{\Tr{\mA}}, \qquad \nu^P=\frac{\Tr{\mA}}{\min_i(\mA_{i,i})}, 
\end{equation*}
which is the optimal choice for coordinate sketches with convenient probabilities without enforcing symmetry. In practice we might not have an access to $\lambda_{\min}(\mA)$, thus we cannot compute $\theta^P$ exactly. Therefore we also test sensitivity of the algorithm to the choice of parameters, and we run some experiments where we only guess parameter $\theta^P$. 

Lastly, the tests are performed on both artificial examples and LibSVM~\cite{chang2011libsvm} data. We shall also explain the legend of plots: ``a'' indicates acceleration, ``nsym'' indicates the algorithm without enforcing symmetry and ``h'' indicates the setting when $\nu^P$ is not known, and a naive heuristic choice is casted.

\subsubsection{The first experiment: synthetic and real-world data}
Let us start with a simple experiment (Figure~\ref{fig:ami_artificial_paper}) to give a quick taste of the numerical performance. 
\begin{figure}[!h]
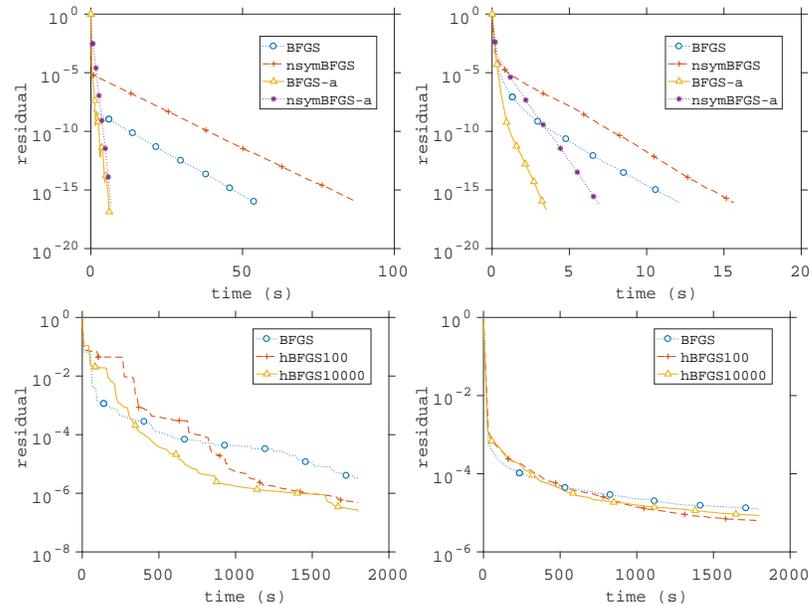

    \centering
\begin{minipage}{0.35\textwidth}
  \centering
\includegraphics[width =  \textwidth ]{convenientrandomjj=4-time}
\end{minipage}%
\begin{minipage}{0.35\textwidth}
  \centering
\includegraphics[width =  \textwidth ]{gaussrandomjj=1-time}
\end{minipage}%
\\
\begin{minipage}{0.35\textwidth}
  \centering
\includegraphics[width =  \textwidth ]{coordheur2epsilon-normalized-time}
\end{minipage}%
\begin{minipage}{0.35\textwidth}
  \centering
\includegraphics[width =  \textwidth ]{convenientheur2svhn-time}
\end{minipage}%
       \caption{Accelerated matrix inversion on synthetic data. From left to right: (i) Eigenvalues of $A \in \R^{100\times 100}$ are $1,10^3, 10^3, \dots,10^3$ and coordinate sketches with convenient probabilities are used. (ii) Eigenvalues of $A \in \R^{100\times 100}$  are $1,2,\dots,n$ and Gaussian sketches are used. Label ``nsym'' indicates non-enforcing symmetry and  ``-a'' indicates acceleration. (iii)   Epsilon dataset ($n=2000$), coordinate sketches with uniform probabilities. (iv) SVHN dataset ($n=3072$), coordinate sketches with convenient probabilities. Label ``h'' indicates that $\lambda_{\min}$ was not precomputed, but $\theta$ was chosen as described in the text.}
\label{fig:ami_artificial_paper}
\end{figure}

The experiments suggest that once the parameters $\theta,\nu$ are estimated exactly, we get a speedup comparing to the nonaccelerated method; and the amount of speedup depends on the structure of $A$ and the sketching strategy. We observe from Figure~\ref{fig:ami_artificial_paper} that we gain a great speedup for ill conditioned problems once the eigenvalues are concentrated around the largest eigenvalue. We also observe from Figure~\ref{fig:ami_artificial_paper} that enforcing symmetry combines well with $\theta, \nu$ computed by \eqref{eq:ami_munu_conv_paper}, which does not consider the symmetrya. On top of that, choice of $\theta, \nu$ per \eqref{eq:ami_munu_conv_paper} seems to be robust to different sketching strategies, and in worst case performs as fast as the nonaccelerated algorithm.

\subsubsection{The second experiment: well understood artificial data \label{sec:ami_bl}}

Let us consider inverting the matrix $\mA = \eta \mI + \beta {\bf 1 1^\top}$ for $\eta>0$ and $\beta\geq -\frac{\eta}{n}$ so as in this case we have control over both $\theta$ and $\nu$. This artificial example was considered in \cite{tu2017breaking} for solving linear systems. In particular, we show that for coordinate sketches with convenient probabilities (which is indeed the same as uniform probabilities in this example), we have
\begin{eqnarray*}
\theta^P &\eqdef& \lambda_{\min}(\E{\mP}) =  \frac{\min \left(\eta, \eta+n\beta \right)}{n(\eta+\beta)},\\
\nu^P &\eqdef& \lambda_{\max}\left(\E{ \E{\mP}^{-\frac12 } \mP \E{\mP}^{-1}\mP \E{\mP}^{-\frac12}}\right)  = n. \\
\end{eqnarray*}

Due to the fact that we do not have a theoretical justification of $\theta,\nu$ for $n>2$ when enforcing symmetry, we set $\theta=\theta^P$ and $\nu=\nu^P$ for Gaussian sketches as well.

\begin{figure}[H]
    \centering
\begin{minipage}{0.40\textwidth}
  \centering
\includegraphics[width =  \textwidth]{breakingconvenient0-1-time}
\end{minipage}%
\begin{minipage}{0.40\textwidth}
  \centering
\includegraphics[width =  \textwidth ]{breakinggauss0-1-time}
\end{minipage}%
    \caption{Accelerated matrix inversion on synthetic data. Parameter choice: $\eta=1+10^{-1}, \beta=-n^{-1}, n=100$. From left to right we have: Coordinate sketch with uniform (convenient) probabilities and Gaussian sketch respectively. 
}\label{fig:ami_bl_ex}
\end{figure}

\begin{figure}[H]
    \centering
\begin{minipage}{0.40\textwidth}
  \centering
\includegraphics[width =  \textwidth]{breakingconvenient0-001-time}
\end{minipage}%
\begin{minipage}{0.40\textwidth}
  \centering
\includegraphics[width =  \textwidth ]{breakinggauss0-001-time}
\end{minipage}%
    \caption{Accelerated matrix inversion on synthetic data. Parameter choice: $\eta=1+10^{-3}, \beta=-n^{-1}, n=100$. From left to right we have: Coordinate sketch with uniform (convenient) probabilities and Gaussian sketch respectively. 
}\label{fig:ami_bl_ex}
\end{figure}

\begin{figure}[H]
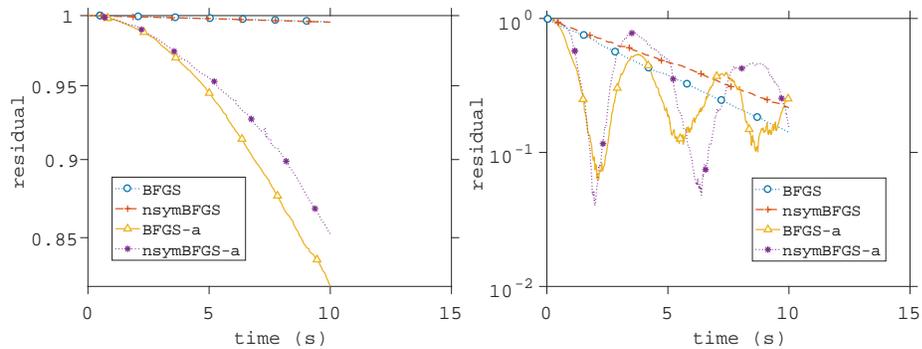

    \centering
\begin{minipage}{0.40\textwidth}
  \centering
\includegraphics[width =  \textwidth]{breakingconvenient1e-05-time}
\end{minipage}%
\begin{minipage}{0.40\textwidth}
  \centering
\includegraphics[width =  \textwidth ]{breakinggauss1e-05-time}
\end{minipage}%
    \caption{Accelerated matrix inversion on synthetic data. Parameter choice: $\eta=1+10^{-5}, \beta=-n^{-1}, n=100$. From left to right we have: Coordinate sketch with uniform (convenient) probabilities and Gaussian sketch, respectively. 
}\label{fig:ami_bl_ex}
\end{figure}

As expected from the theory, as the matrix to be inverted becomes more ill conditioned, the accelerated method performs significantly better compared to the nonaccelerated method for coordinate sketches. In fact, an arbitrary speedup can be obtained by setting $\beta=-n^{-1}$ and $\eta \rightarrow 1$ for the coordinate sketches setup. On the other hand, Gaussian sketches report the slowing of the algorithm, most likely caused by the fact that the theoretical parameters $\theta, \nu$ for Gaussian sketches with enforced symmetry are different to $\theta^P, \nu^P$, which are estimated for coordinate sketches without enforced symmetry. In the case of coordinate sketches with symmetry enforced, we suspect a great speedup even though the parameters $\theta, \nu$ were set to $\theta^P, \nu^P$.

\subsubsection{The third experiment: more complex artificial data}

We randomly generate an orthonormal matrix $\mU$, choose diagonal matrix $\mD$, and set $\mA=\mU\mD\mU^\top$. Clearly, diagonal elements of $\mD$ are eigenvalues of $\mA$. We set them in the following way:

\begin{itemize}
\item Uniform grid. The eigenvalues are set to $1,2,\dots,n$.
\item One small, the rest larger. The smallest eigenvalue is $1$, remaining eigenvalues are all $10$ in the first example, all $100$ in the second example and all $1000$ in the third example in this category.
\item One large, the rest small. The largest eigenvalue is $10^4$, the remaining eigenvalues are all $1$.
\end{itemize}

Firstly, consider coordinate sketches with convenient probabilities. Notice that we can easily estimate $\nu^P, \theta^P$ due to the results from Section~\ref{sec:ami_convenient_munu} since we have control of $\lambda_{\min}(\mA)$ and therefore also of $\theta$. Therefore, we set $\theta=\theta^P=\min \mD_{i,i}$ and $\nu=\nu^P$ for Algorithm~\ref{alg:ami_qn}. Then, we consider coordinate sketches with uniform probabilities and Gaussian sketches. In both cases, we set the parameters $\theta,\nu$ as for coordinate sketches with convenient probabilities.

\begin{figure}[H]
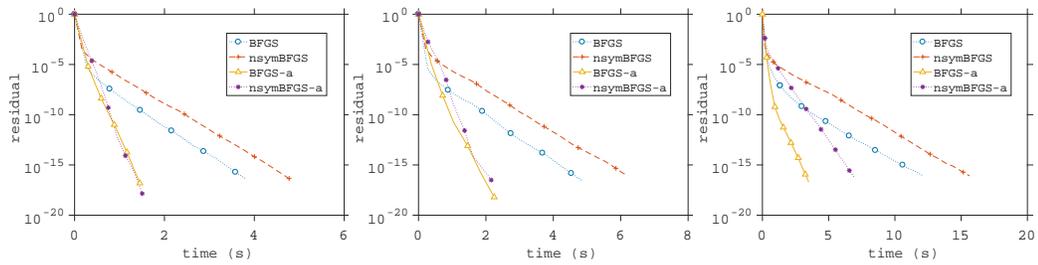

    \centering
\begin{minipage}{0.30\textwidth}
  \centering
\includegraphics[width =  \textwidth ]{convenientrandomjj=1-time}
\end{minipage}%
\begin{minipage}{0.30\textwidth}
  \centering
\includegraphics[width =  \textwidth ]{coordrandomjj=1-time}
\end{minipage}%
\begin{minipage}{0.30\textwidth}
  \centering
\includegraphics[width =  \textwidth ]{gaussrandomjj=1-time}
\end{minipage}%
    \caption{   Eigenvalues set to $1,2,3,\dots n$. From left to right we have: Coordinate sketch with convenient probabilities, coordinate sketch with uniform probabilities and Gaussian sketch respectively. 
}
\label{fig:ami_rand_conv_lin}
\end{figure}

\begin{figure}[H]
    \centering
\begin{minipage}{0.30\textwidth}
  \centering
\includegraphics[width =  \textwidth ]{convenientrandomjj=2-time}
\end{minipage}%
\begin{minipage}{0.30\textwidth}
  \centering
\includegraphics[width =  \textwidth ]{coordrandomjj=2-time}
\end{minipage}%
\begin{minipage}{0.30\textwidth}
  \centering
\includegraphics[width =  \textwidth ]{gaussrandomjj=2-time}
\end{minipage}%
    \caption{   Eigenvalues set to $1,10,10,\dots 10$. From left to right we have: Coordinate sketch with convenient probabilities, coordinate sketch with uniform probabilities and Gaussian sketch respectively. 
}\label{fig:ami_rand_conv_10}
\end{figure}

\begin{figure}[H]
    \centering
\begin{minipage}{0.30\textwidth}
  \centering
\includegraphics[width =  \textwidth ]{convenientrandomjj=3-time}
\end{minipage}%
\begin{minipage}{0.30\textwidth}
  \centering
\includegraphics[width =  \textwidth ]{coordrandomjj=3-time}
\end{minipage}%
\begin{minipage}{0.30\textwidth}
  \centering
\includegraphics[width =  \textwidth ]{gaussrandomjj=3-time}
\end{minipage}%
    \caption{  Accelerated matrix inversion on synthetic data.  Eigenvalues set to $1,100,100,\dots 100$. From left to right we have: Coordinate sketch with convenient probabilities, coordinate sketch with uniform probabilities and Gaussian sketch respectively. 
}\label{fig:ami_rand_conv_100}
\end{figure}

\begin{figure}[H]
    \centering
\begin{minipage}{0.30\textwidth}
  \centering
\includegraphics[width =  \textwidth ]{convenientrandomjj=4-time}
\end{minipage}%
\begin{minipage}{0.30\textwidth}
  \centering
\includegraphics[width =  \textwidth ]{coordrandomjj=4-time}
\end{minipage}%
\begin{minipage}{0.30\textwidth}
  \centering
\includegraphics[width =  \textwidth ]{gaussrandomjj=4-time}
\end{minipage}%
    \caption{ Accelerated matrix inversion on synthetic data.  Eigenvalues set to $1$, $1000$, $1000,$ $\dots,$ $1000$. From left to right we have: Coordinate sketch with convenient probabilities, coordinate sketch with uniform probabilities and Gaussian sketch respectively. 
}\label{fig:ami_rand_conv_1000}
\end{figure}

\begin{figure}[H]
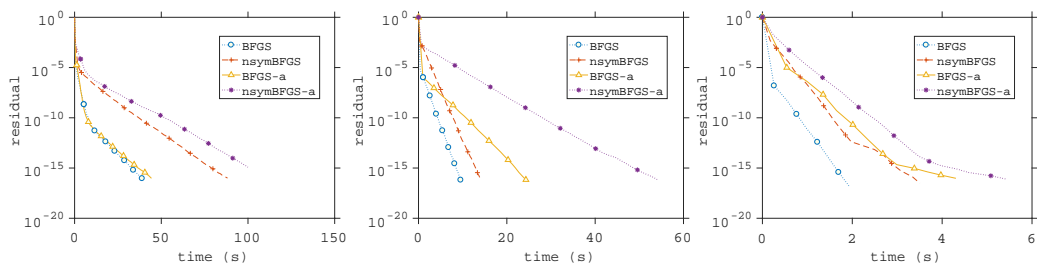

    \centering
\begin{minipage}{0.30\textwidth}
  \centering
\includegraphics[width =  \textwidth ]{convenientrandomjj=5-time}
\end{minipage}%
\begin{minipage}{0.30\textwidth}
  \centering
\includegraphics[width =  \textwidth ]{coordrandomjj=5-time}
\end{minipage}%
\begin{minipage}{0.30\textwidth}
  \centering
\includegraphics[width =  \textwidth ]{gaussrandomjj=5-time}
\end{minipage}%
    \caption{  Accelerated matrix inversion on synthetic data. Eigenvalues set to $10000$, $1$,$1$, $\dots$, $1$. From left to right we have: Coordinate sketch with convenient probabilities, coordinate sketch with uniform probabilities and Gaussian sketch respectively. 
}\label{fig:ami_rand_conv_110000}
\end{figure}

The numerical experiments in this section indicate that one might choose $\theta ,\nu$ as per Section~\ref{sec:ami_convenient_munu}. In other words, one might pretend to be in the setting when symmetry is not enforced and coordinate sketches with convenient probabilities are used. In fact, the practical speedup coming from the acceleration depends very strongly on the structure of matrix $\mA$.
 Another message to be delivered is that both preserving symmetry and acceleration yield a better convergence and they combine together well.

We also consider a problem where we pretend to not have access to $\lambda_{\min}(\mA)$, therefore we cannot choose $\theta=\theta^P$. Instead, we naively choose $\theta=\frac{1}{100\nu}$ and $\theta=\frac{1}{10000\nu}$.

\begin{figure}[H]
    \centering
\begin{minipage}{0.30\textwidth}
  \centering
\includegraphics[width =  \textwidth ]{convenientrandom-heurjj=1-time}
\end{minipage}%
\begin{minipage}{0.30\textwidth}
  \centering
\includegraphics[width =  \textwidth ]{coordrandom-heurjj=1-time}
\end{minipage}%
\begin{minipage}{0.30\textwidth}
  \centering
\includegraphics[width =  \textwidth ]{gaussrandom-heurjj=1-time}
\end{minipage}%
    \caption{ Accelerated matrix inversion on synthetic data.  Eigenvalues set to $1$, $2,$ $\dots$, $n$. From left to right we have: Coordinate sketch with convenient probabilities, coordinate sketch with uniform probabilities and Gaussian sketch respectively. 
}\label{fig:ami_rand_conv_110000}
\end{figure}

\begin{figure}[H]
    \centering
\begin{minipage}{0.30\textwidth}
  \centering
\includegraphics[width =  \textwidth ]{convenientrandom-heurjj=2-time}
\end{minipage}%
\begin{minipage}{0.30\textwidth}
  \centering
\includegraphics[width =  \textwidth ]{coordrandom-heurjj=2-time}
\end{minipage}%
\begin{minipage}{0.30\textwidth}
  \centering
\includegraphics[width =  \textwidth ]{gaussrandom-heurjj=2-time}
\end{minipage}%
    \caption{Accelerated matrix inversion on synthetic data. Eigenvalues set to $1$, $10$, $10,$ $\dots$, $ 10$. Coordinate sketch with convenient probabilities, coordinate sketch with uniform probabilities and Gaussian sketch respectively. 
}\label{fig:ami_rand_conv_110000}
\end{figure}

\begin{figure}[H]
    \centering
\begin{minipage}{0.30\textwidth}
  \centering
\includegraphics[width =  \textwidth ]{convenientrandom-heurjj=3-time}
\end{minipage}%
\begin{minipage}{0.30\textwidth}
  \centering
\includegraphics[width =  \textwidth ]{coordrandom-heurjj=3-time}
\end{minipage}%
\begin{minipage}{0.30\textwidth}
  \centering
\includegraphics[width =  \textwidth ]{gaussrandom-heurjj=3-time}
\end{minipage}%
    \caption{Accelerated matrix inversion on synthetic data. Eigenvalues set to $1$, $100$, $100$, $\dots$, $100$. From left to right we have: Coordinate sketch with convenient probabilities, coordinate sketch with uniform probabilities and Gaussian sketch respectively. 
}\label{fig:ami_rand_conv_110000}
\end{figure}

\begin{figure}[H]
    \centering
\begin{minipage}{0.30\textwidth}
  \centering
\includegraphics[width =  \textwidth ]{convenientrandom-heurjj=4-time}
\end{minipage}%
\begin{minipage}{0.30\textwidth}
  \centering
\includegraphics[width =  \textwidth ]{coordrandom-heurjj=4-time}
\end{minipage}%
\begin{minipage}{0.30\textwidth}
  \centering
\includegraphics[width =  \textwidth ]{gaussrandom-heurjj=4-time}
\end{minipage}%
    \caption{Accelerated matrix inversion on synthetic data. Eigenvalues set to $1$, $1000$, $1000$, $\dots$, $1000$. From left to right we have: Coordinate sketch with convenient probabilities, coordinate sketch with uniform probabilities and Gaussian sketch respectively. 
}\label{fig:ami_rand_conv_110000}
\end{figure}

\begin{figure}[H]
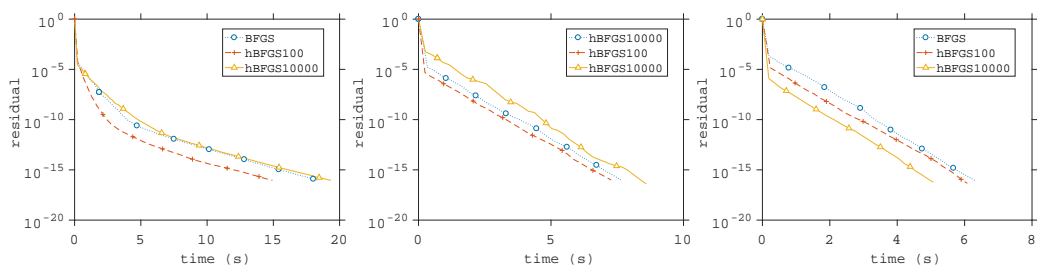

    \centering
\begin{minipage}{0.30\textwidth}
  \centering
\includegraphics[width =  \textwidth ]{convenientrandom-heurjj=5-time}
\end{minipage}%
\begin{minipage}{0.30\textwidth}
  \centering
\includegraphics[width =  \textwidth ]{coordrandom-heurjj=5-time}
\end{minipage}%
\begin{minipage}{0.30\textwidth}
  \centering
\includegraphics[width =  \textwidth ]{gaussrandom-heurjj=5-time}
\end{minipage}%
    \caption{Accelerated matrix inversion on synthetic data. Eigenvalues set to $10000$, $1$, $1$, $\dots$, $1$. From left to right we have: Coordinate sketch with convenient probabilities, coordinate sketch with uniform probabilities and Gaussian sketch respectively. 
}\label{fig:ami_rand_conv_110000}
\end{figure}

Notice that once the acceleration parameters are not set exactly (but they are still reasonable), we observe that the performance of the accelerated algorithm is essentially the same as the performance of the nonaccelerated algorithm. We have observed the similar behavior when setting $\theta=\theta^P$ for Gaussian sketches.

\subsubsection{The fifth experiment: LibSVM data}

Next we investigate if the accelerated {\tt BFGS} update improves upon the standard {\tt BFGS} update when applied to the Hessian $\nabla^2 f(x)$  of ridge regression problems
of the form
\begin{equation}\label{eq:ami_ridgeMatrix}
\min_{x\in \R^d}f(x)\eqdef \frac{1}{2}\norm{\mA x-b}_2^2 + \frac{\lambda}{2} \norm{x}_2^2,\quad  \quad\nabla^2 f(x) = \mA^\top \mA+\lambda \mI,
\end{equation}
using data from LibSVM~\cite{chang2011libsvm}. Datapoints (rows of $\mA$) were normalized such that $\|\mA_{i:}\|^2=1$ for all $i$ and the regularization parameter was chosen as $\lambda=\frac{1}{m}$.

First, we run the experiments on smaller problems when parameters $\theta$, $\nu$ are precomputed for coordinate sketches with convenient probabilities \eqref{eq:ami_munu_conv_paper}.

\begin{figure}[H]
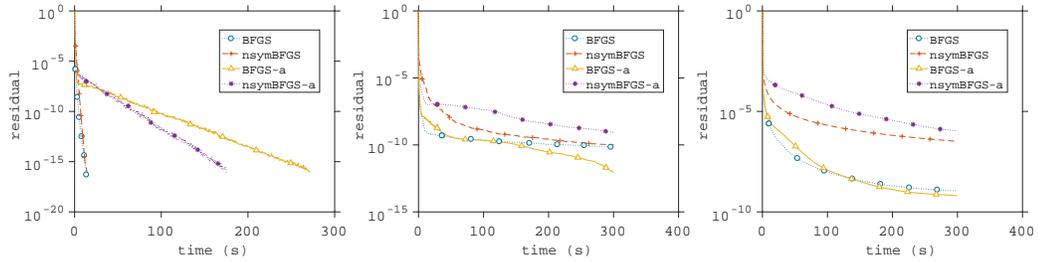

    \centering
\begin{minipage}{0.30\textwidth}
  \centering
\includegraphics[width =  \textwidth ]{convenientaloi-time}
\end{minipage}%
\begin{minipage}{0.30\textwidth}
  \centering
\includegraphics[width =  \textwidth ]{coordaloi-time}
\end{minipage}%
\begin{minipage}{0.30\textwidth}
  \centering
\includegraphics[width =  \textwidth ]{gaussaloi-time}
\end{minipage}%
    \caption{ Accelerated matrix inversion on real data.  Dataset aloi: $n=128$. From left to right we have: Coordinate sketch with convenient probabilities, coordinate sketch with uniform probabilities and Gaussian sketch respectively. 
}\label{fig:ami_aloig}
\end{figure}

\begin{figure}[H]
    \centering
\begin{minipage}{0.30\textwidth}
  \centering
\includegraphics[width =  \textwidth ]{convenientw1a-time}
\end{minipage}%
\begin{minipage}{0.30\textwidth}
  \centering
\includegraphics[width =  \textwidth ]{coordw1a-time}
\end{minipage}%
\begin{minipage}{0.30\textwidth}
  \centering
\includegraphics[width =  \textwidth ]{gaussw1a-time}
\end{minipage}%
    \caption{ Accelerated matrix inversion on real data.  Dataset w1a: $n=300$. From left to right we have: Coordinate sketch with convenient probabilities, coordinate sketch with uniform probabilities and Gaussian sketch respectively. 
}\label{fig:ami_w1a}
\end{figure}

\begin{figure}[H]
    \centering
\begin{minipage}{0.30\textwidth}
  \centering
\includegraphics[width =  \textwidth ]{convenientw2a-time}
\end{minipage}%
\begin{minipage}{0.30\textwidth}
  \centering
\includegraphics[width =  \textwidth ]{coordw2a-time}
\end{minipage}%
\begin{minipage}{0.30\textwidth}
  \centering
\includegraphics[width =  \textwidth ]{gaussw2a-time}
\end{minipage}%
    \caption{ Accelerated matrix inversion on real data.  Dataset w2a: $n=300$. From left to right we have: Coordinate sketch with convenient probabilities, coordinate sketch with uniform probabilities and Gaussian sketch respectively. 
}\label{fig:ami_w2a}
\end{figure}

\begin{figure}[H]
    \centering
\begin{minipage}{0.30\textwidth}
  \centering
\includegraphics[width =  \textwidth ]{convenientmushrooms-time}
\end{minipage}%
\begin{minipage}{0.30\textwidth}
  \centering
\includegraphics[width =  \textwidth ]{coordmushrooms-time}
\end{minipage}%
\begin{minipage}{0.30\textwidth}
  \centering
\includegraphics[width =  \textwidth ]{gaussmushrooms-time}
\end{minipage}%
    \caption{ Accelerated matrix inversion on real data.  Dataset mushrooms: $n=112$. From left to right we have: Coordinate sketch with convenient probabilities, coordinate sketch with uniform probabilities and Gaussian sketch respectively. 
}\label{fig:ami_mushroomsxx}
\end{figure}

\begin{figure}[H]
\centering
\begin{minipage}{0.30\textwidth}
  \centering
\includegraphics[width =  \textwidth ]{convenientprotein-time}
\end{minipage}%
\begin{minipage}{0.30\textwidth}
  \centering
\includegraphics[width =  \textwidth ]{coordprotein-time}
\end{minipage}%
\begin{minipage}{0.30\textwidth}
  \centering
\includegraphics[width =  \textwidth ]{gaussprotein-time}
\end{minipage}%
    \caption{  Accelerated matrix inversion on real data. Dataset protein: $n=357$. From left to right we have: Coordinate sketch with convenient probabilities, coordinate sketch with uniform probabilities and Gaussian sketch respectively. 
}\label{fig:ami_protein}
\end{figure}

\begin{figure}[H]
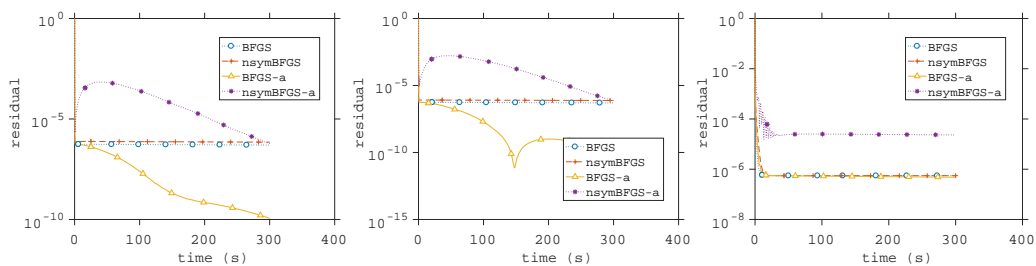

    \centering
\begin{minipage}{0.30\textwidth}
  \centering
\includegraphics[width =  \textwidth ]{convenientphishing-time}
\end{minipage}%
\begin{minipage}{0.30\textwidth}
  \centering
\includegraphics[width =  \textwidth ]{coordphishing-time}
\end{minipage}%
\begin{minipage}{0.30\textwidth}
  \centering
\includegraphics[width =  \textwidth ]{gaussphishing-time}
\end{minipage}%
    \caption{  Accelerated matrix inversion on real data. Dataset phishing: $n=68$. From left to right we have: Coordinate sketch with convenient probabilities, coordinate sketch with uniform probabilities and Gaussian sketch respectively. 
}\label{fig:ami_splice-scale}
\end{figure}

In the vast majority of examples, the accelerated method performed significantly better than the nonaccelerated method for coordinate sketches (with both convenient and uniform probabilities), however the methods were comparable for Gaussian sketches. We believe that this is due to the fact that choice of parameters as per~\eqref{eq:ami_munu_conv_paper} is close to the optimal parameters for coordinate sketches, and further for Gaussian sketches. However, the experiments on coordinate sketches indicates that for some classes of problems, accelerated algorithms with finely tuned parameters bring a great speedup compared to nonaccelerated ones. 

We also consider a problem where we do not compute $\lambda_{\min}(\mA)$, and therefore we cannot choose $\theta=\theta^P$ in \eqref{eq:ami_munu_conv_paper}. Instead, we choose $\theta=\frac{1}{100\nu}$ and $\theta=\frac{1}{10000\nu}$.

\begin{figure}[H]
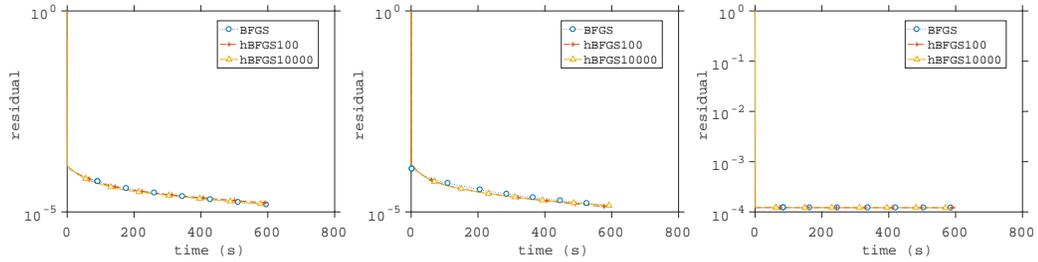

    \centering
\begin{minipage}{0.30\textwidth}
  \centering
\includegraphics[width =  \textwidth ]{convenientheur2madelon-time}
\end{minipage}%
\begin{minipage}{0.30\textwidth}
  \centering
\includegraphics[width =  \textwidth ]{coordheur2madelon-time}
\end{minipage}%
\begin{minipage}{0.30\textwidth}
  \centering
\includegraphics[width =  \textwidth ]{gaussheur2madelon-time}
\end{minipage}%
    \caption{ Accelerated matrix inversion on real data.  Dataset madelon: $n=500$. From left to right we have: Coordinate sketch with convenient probabilities, coordinate sketch with uniform probabilities and Gaussian sketch respectively. 
}\label{fig:ami_splice-scale}
\end{figure}

\begin{figure}[H]
    \centering
\begin{minipage}{0.30\textwidth}
  \centering
\includegraphics[width =  \textwidth ]{convenientheur2epsilon-normalized-time}
\end{minipage}%
\begin{minipage}{0.30\textwidth}
  \centering
\includegraphics[width =  \textwidth ]{coordheur2epsilon-normalized-time}
\end{minipage}%
\begin{minipage}{0.30\textwidth}
  \centering
\includegraphics[width =  \textwidth ]{gaussheur2epsilon-normalized-time}
\end{minipage}%
    \caption{ Accelerated matrix inversion on real data.  Dataset epsilon: $n=2000$. From left to right we have: Coordinate sketch with convenient probabilities, coordinate sketch with uniform probabilities and Gaussian sketch respectively. 
}\label{fig:ami_splice-scale}
\end{figure}

\begin{figure}[H]
    \centering
\begin{minipage}{0.30\textwidth}
  \centering
\includegraphics[width =  \textwidth ]{convenientheur2svhn-time}
\end{minipage}%
\begin{minipage}{0.30\textwidth}
  \centering
\includegraphics[width =  \textwidth ]{coordheur2svhn-time}
\end{minipage}%
\begin{minipage}{0.30\textwidth}
  \centering
\includegraphics[width =  \textwidth ]{gaussheur2svhn-time}
\end{minipage}%
    \caption{ Accelerated matrix inversion on real data.  Dataset svhn: $n=3072$. From left to right we have: Coordinate sketch with convenient probabilities, coordinate sketch with uniform probabilities and Gaussian sketch respectively. 
}\label{fig:ami_splice-scale}
\end{figure}

\begin{figure}[H]
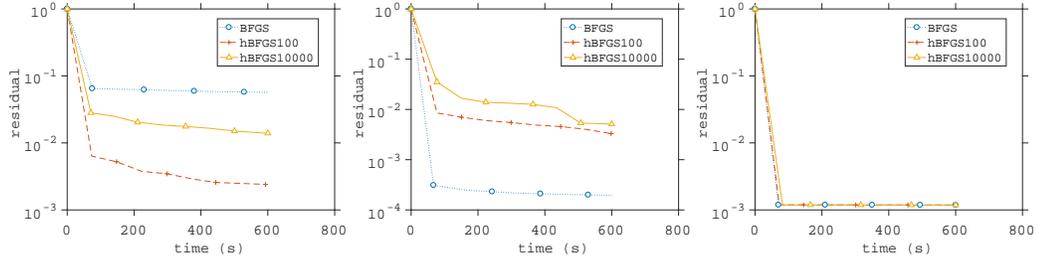

    \centering
\begin{minipage}{0.30\textwidth}
  \centering
\includegraphics[width =  \textwidth ]{convenientheur2gisette-scale-time}
\end{minipage}%
\begin{minipage}{0.30\textwidth}
  \centering
\includegraphics[width =  \textwidth ]{coordheur2gisette-scale-time}
\end{minipage}%
\begin{minipage}{0.30\textwidth}
  \centering
\includegraphics[width =  \textwidth ]{gaussheur2gisette-scale-time}
\end{minipage}%
    \caption{  Accelerated matrix inversion on real data. Dataset gisette: $n=5000$. From left to right we have: Coordinate sketch with convenient probabilities, coordinate sketch with uniform probabilities and Gaussian sketch respectively. 
}\label{fig:ami_splice-scale}
\end{figure}

Notice that once the acceleration parameters are not set exactly (but they are still reasonable), we observe that the performance of the accelerated algorithm is essentially the same as the performance of the nonaccelerated algorithm, which is essentially the same conclusion as for artificially generated examples.

\subsubsection{The fourth experiment: sensitivity to the acceleration parameters}
Here we investigate the sensitivity of the accelerated {\tt BFGS} to the parameters $\theta$ and $\nu$. First we compute $\nu^P, \theta^P$ and from this we extract the following exponential grids: $\theta_i =2^{i-4}\theta$ and $\nu_i=5^{i-4}\nu$ for $i=1,2,\dots 7$. To gauge the gain is using acceleration with a particular $(\theta, \nu)$ pair, we run the accelerated algorithm for a fixed time then store the error of the final iterate. We then compute average per iteration decrease and divide it by average per iteration decrease of nonaccelerated algorithm. Thus if the resulting difference is less than one, then the accelerated algorithm was faster to nonaccelerated. 

In the plots below, $n=200$ was chosen. We focused on 2 problems described in the previous section---when the eigenvalues are uniformly distributed and when the the largest eigenvalue have multiplicity $n-1$. 

\begin{figure}[H]
    \centering
\begin{minipage}{0.33\textwidth}
  \centering
\includegraphics[width =  \textwidth]{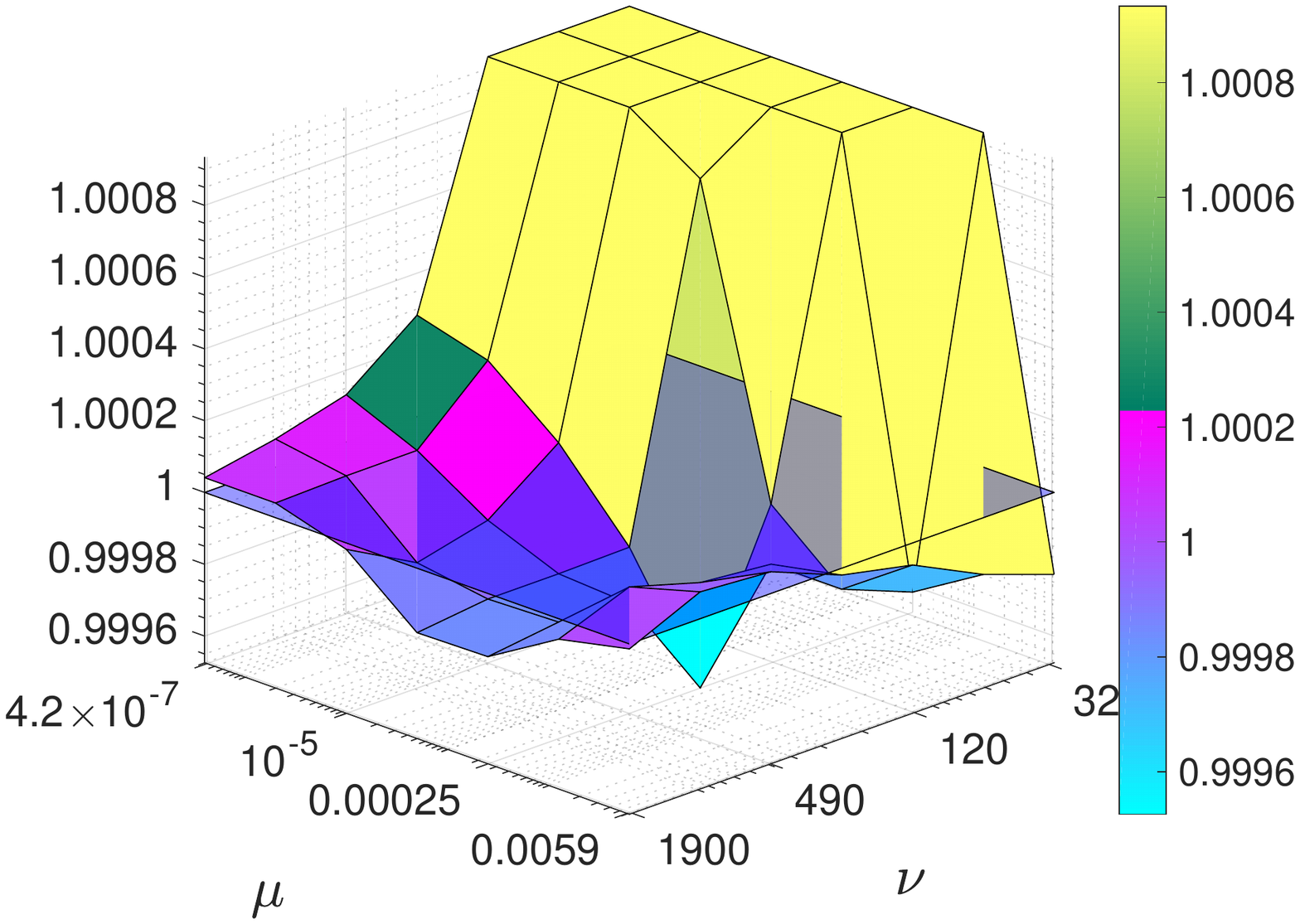}
\end{minipage}%
\begin{minipage}{0.33\textwidth}
  \centering
\includegraphics[width =  \textwidth ]{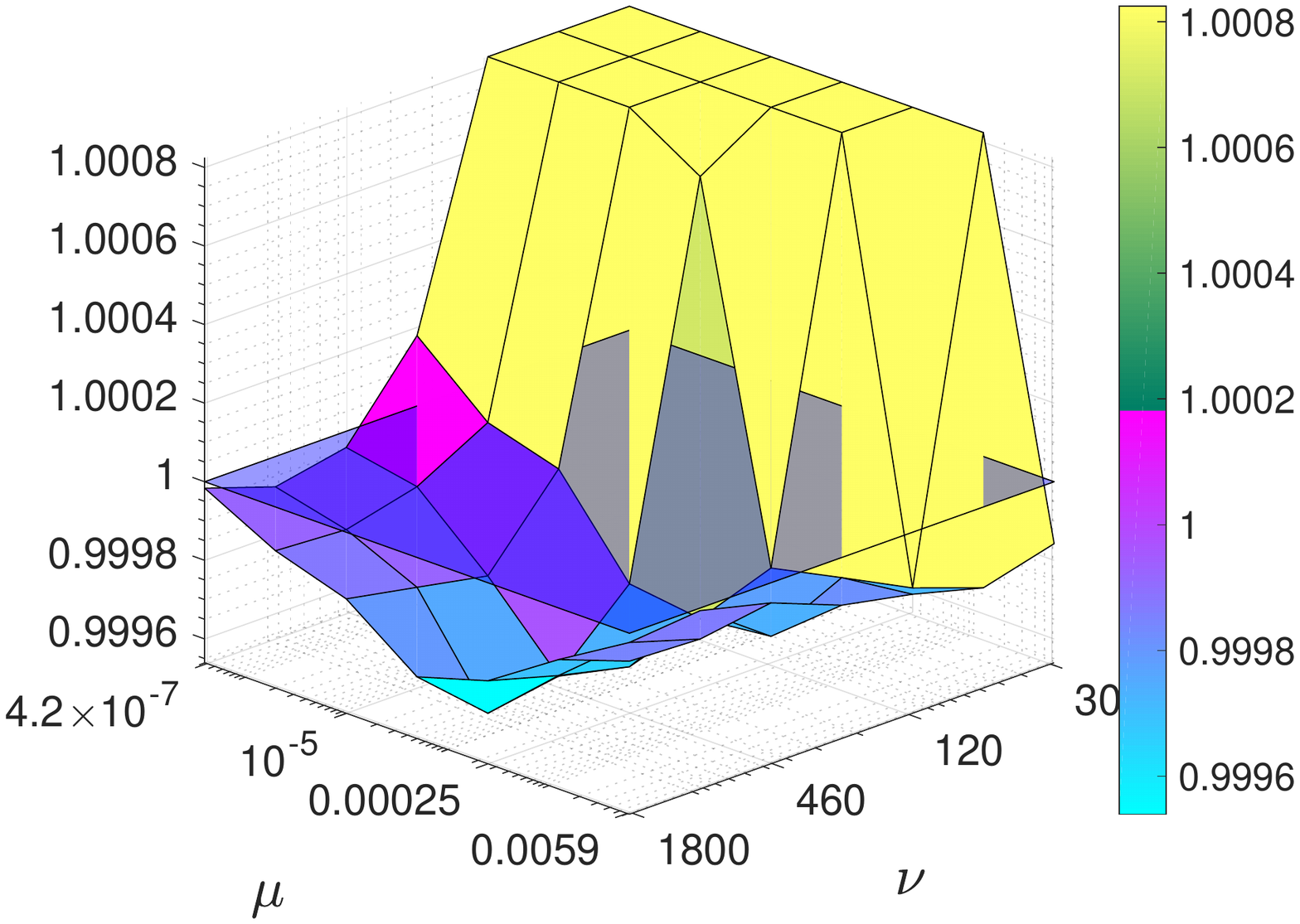}
\end{minipage}%
\begin{minipage}{0.33\textwidth}
  \centering
\includegraphics[width =  \textwidth ]{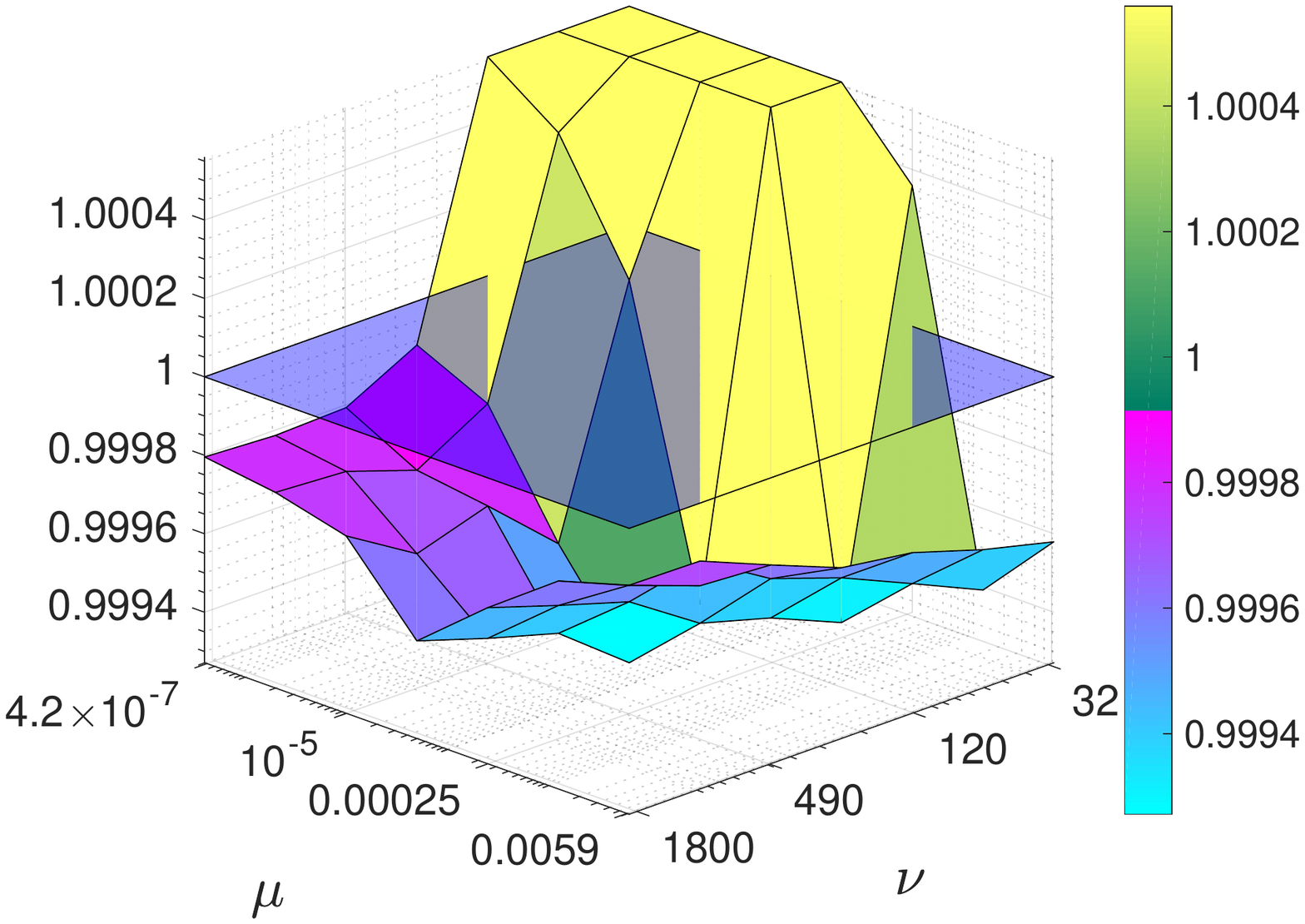}
\end{minipage}%
    \caption{Accelerated matrix inversion on synthetic data. Sensitivity to acceleration parameters. Eigenvalues of $A$ are set to $1,2\dots,n$. From left to right we have: Coordinate sketches with convenient probabilities, coordiante sketches with uniform probabilities and Gaussian sketches. Choice of parameters as per \eqref{eq:ami_munu_conv_paper} in the middle of plots. Each instance was run for 5 seconds. }
\label{fig:ami_heat_gauss_rand}
\end{figure}

\begin{figure}[H]
    \centering
\begin{minipage}{0.33\textwidth}
  \centering
\includegraphics[width =  \textwidth]{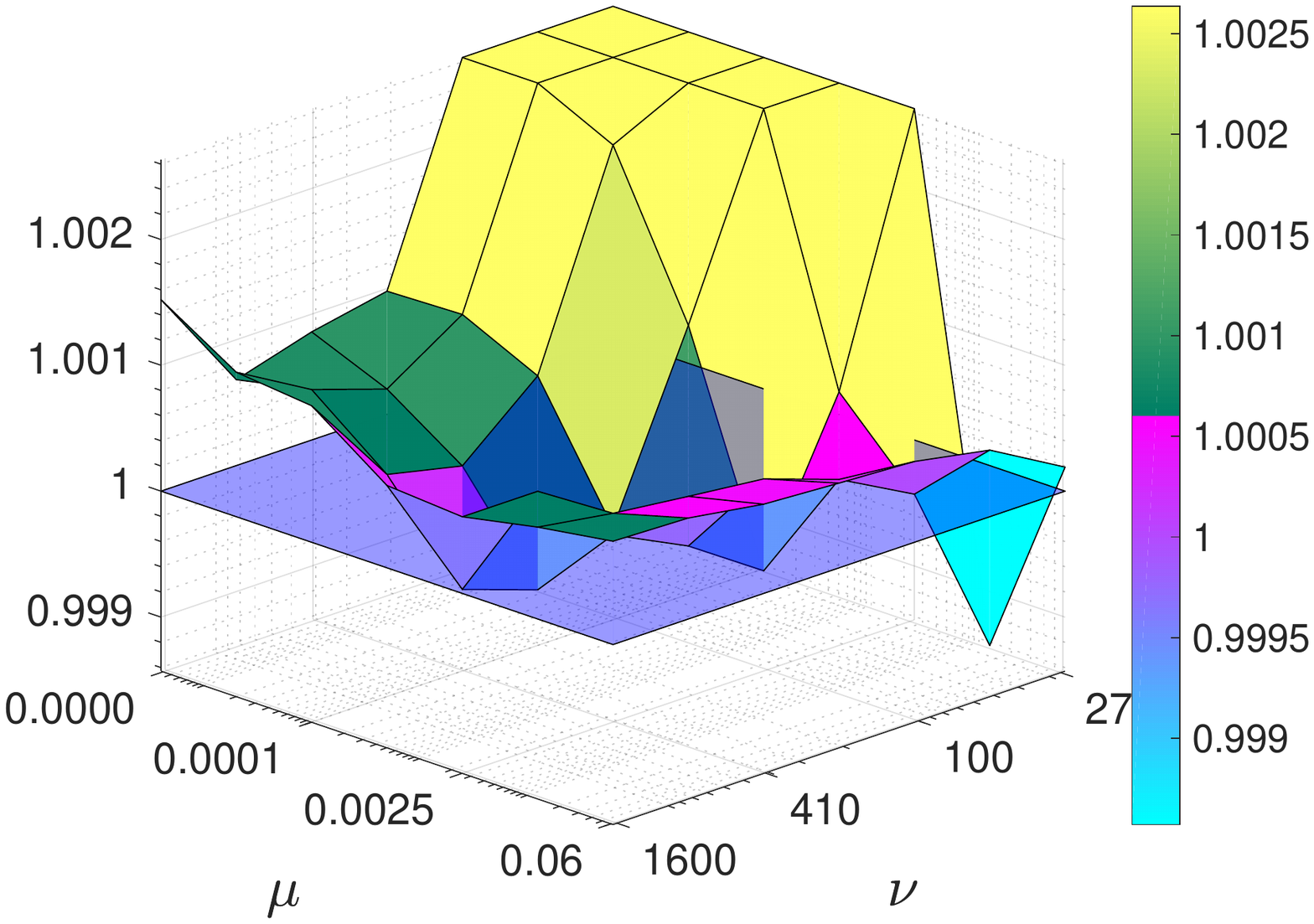}
\end{minipage}%
\begin{minipage}{0.33\textwidth}
  \centering
\includegraphics[width =  \textwidth ]{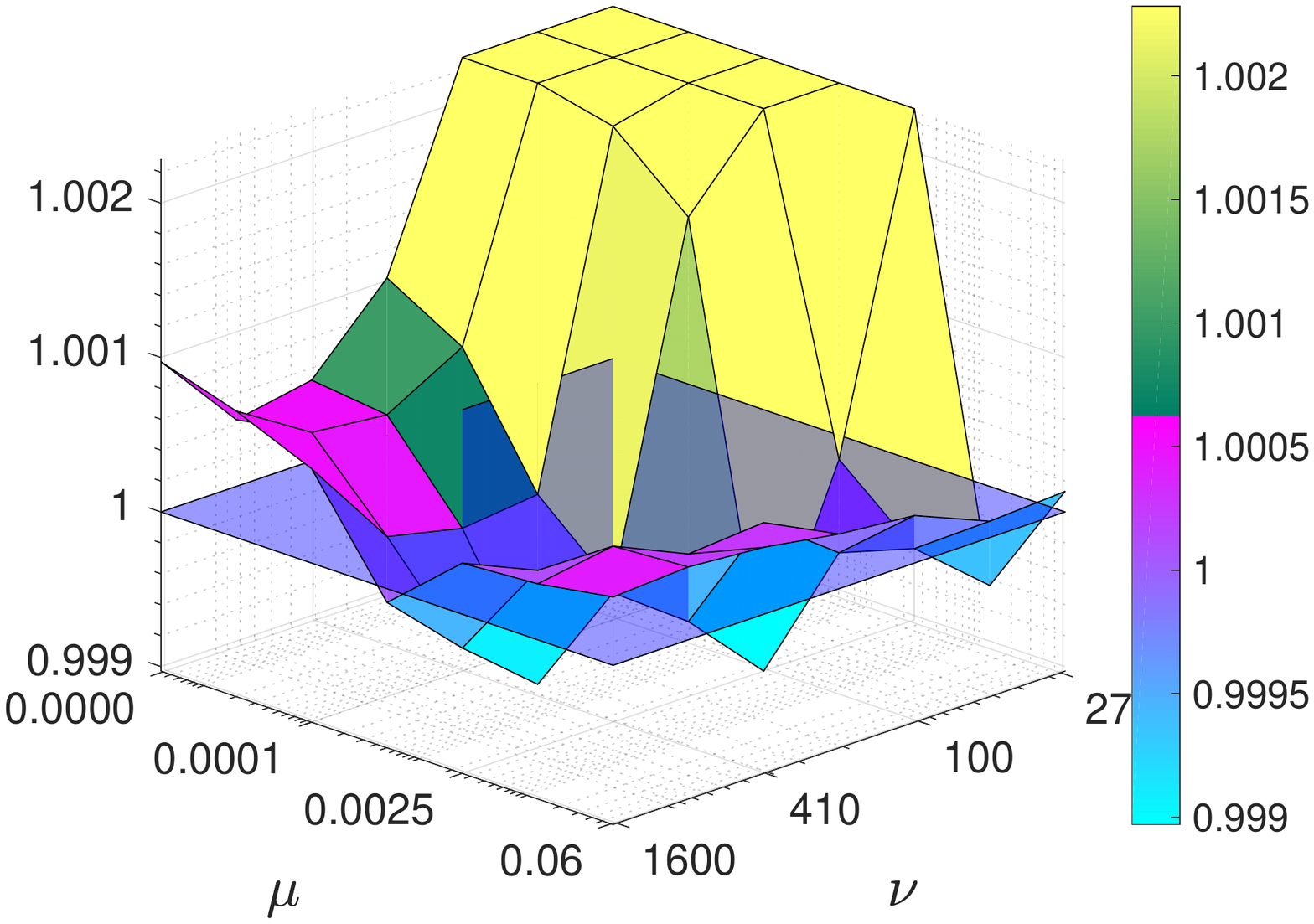}
\end{minipage}%
\begin{minipage}{0.33\textwidth}
  \centering
\includegraphics[width =  \textwidth ]{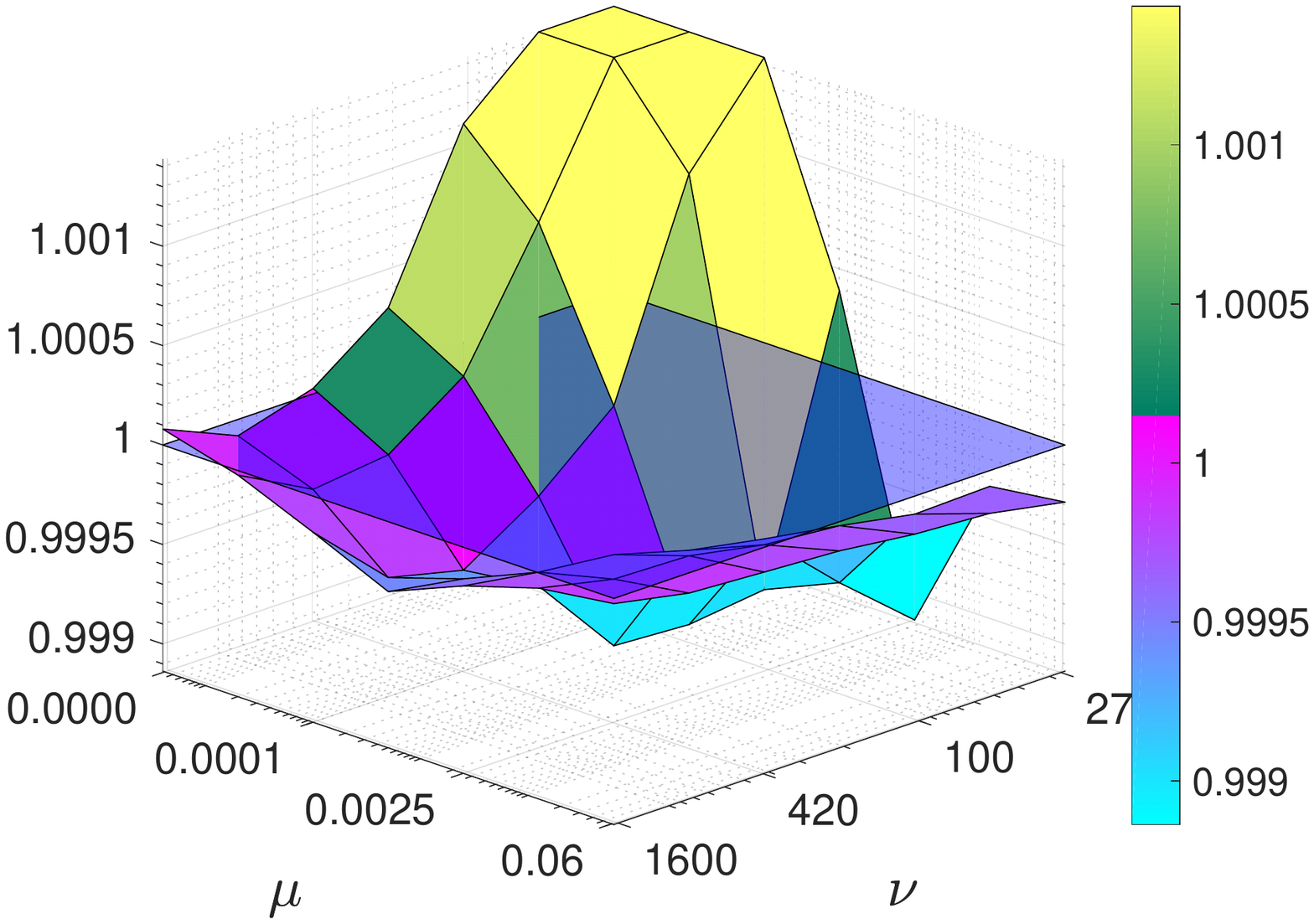}
\end{minipage}%
    \caption{Accelerated matrix inversion on synthetic data. Sensitivity to acceleration parameters. Eigenvalues of $\mA$ are set to $1,10, 10,\dots,10 $. From left to right we have: Coordinate sketches with convenient probabilities, coordiante sketches with uniform probabilities and Gaussian sketches. Choice of parameters as per \eqref{eq:ami_munu_conv_paper} in the middle of plots. Each instance was run for 2 seconds.}
    \label{fig:ami_heat_gauss_rand}
\end{figure}

\begin{figure}[H]
    \centering
\begin{minipage}{0.33\textwidth}
  \centering
\includegraphics[width =  \textwidth]{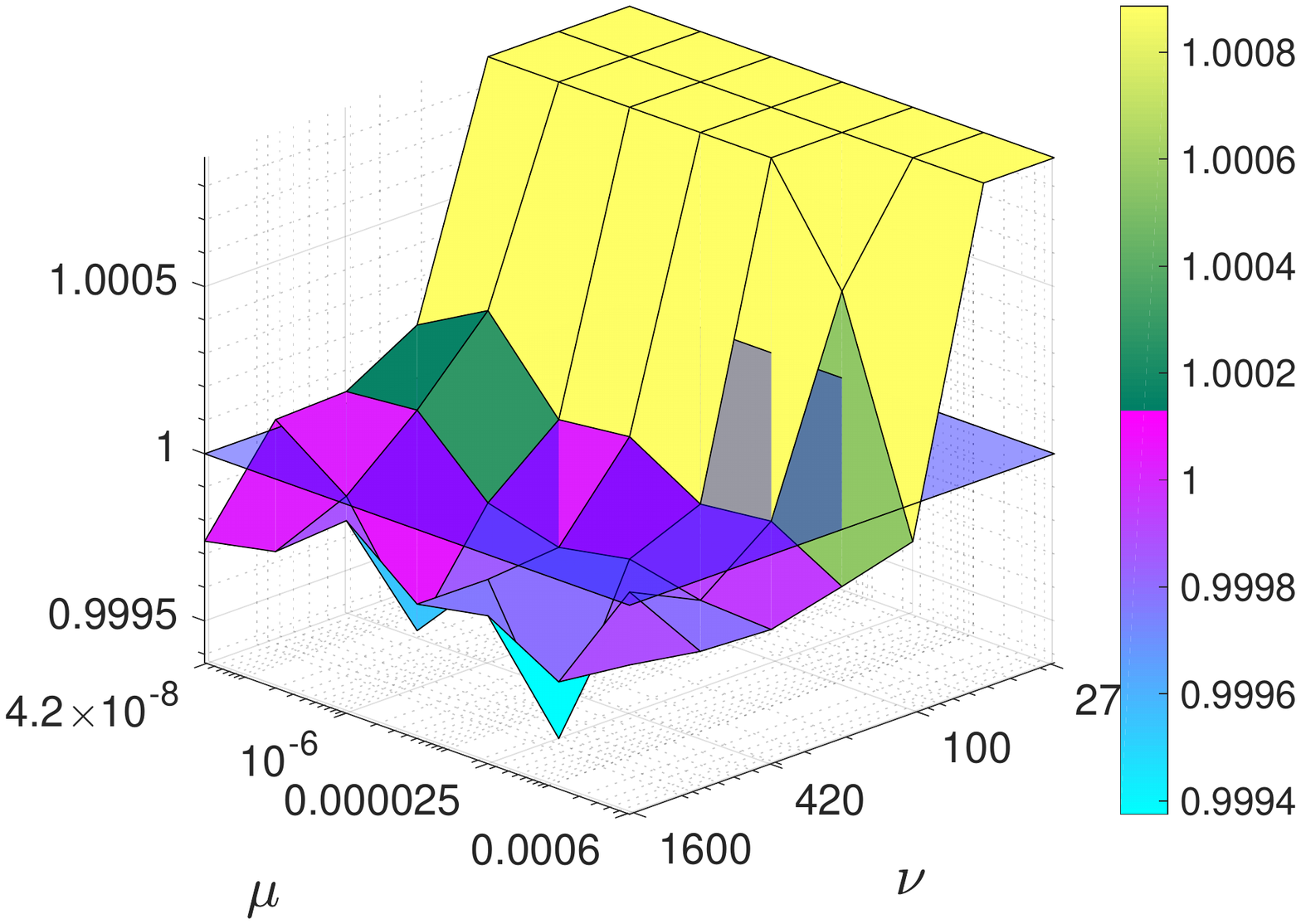}
\end{minipage}%
\begin{minipage}{0.33\textwidth}
  \centering
\includegraphics[width =  \textwidth ]{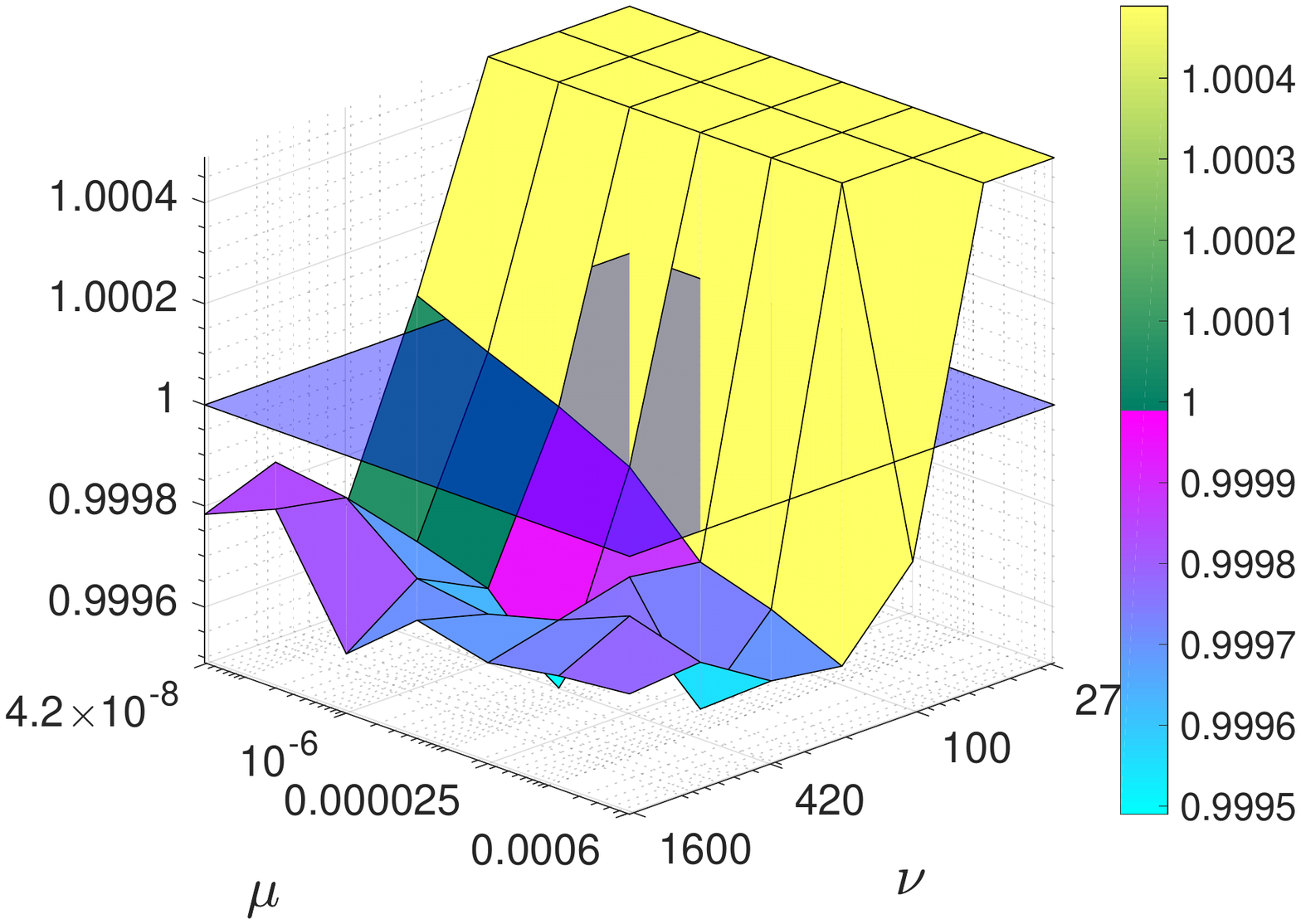}
\end{minipage}%
\begin{minipage}{0.33\textwidth}
  \centering
\includegraphics[width =  \textwidth ]{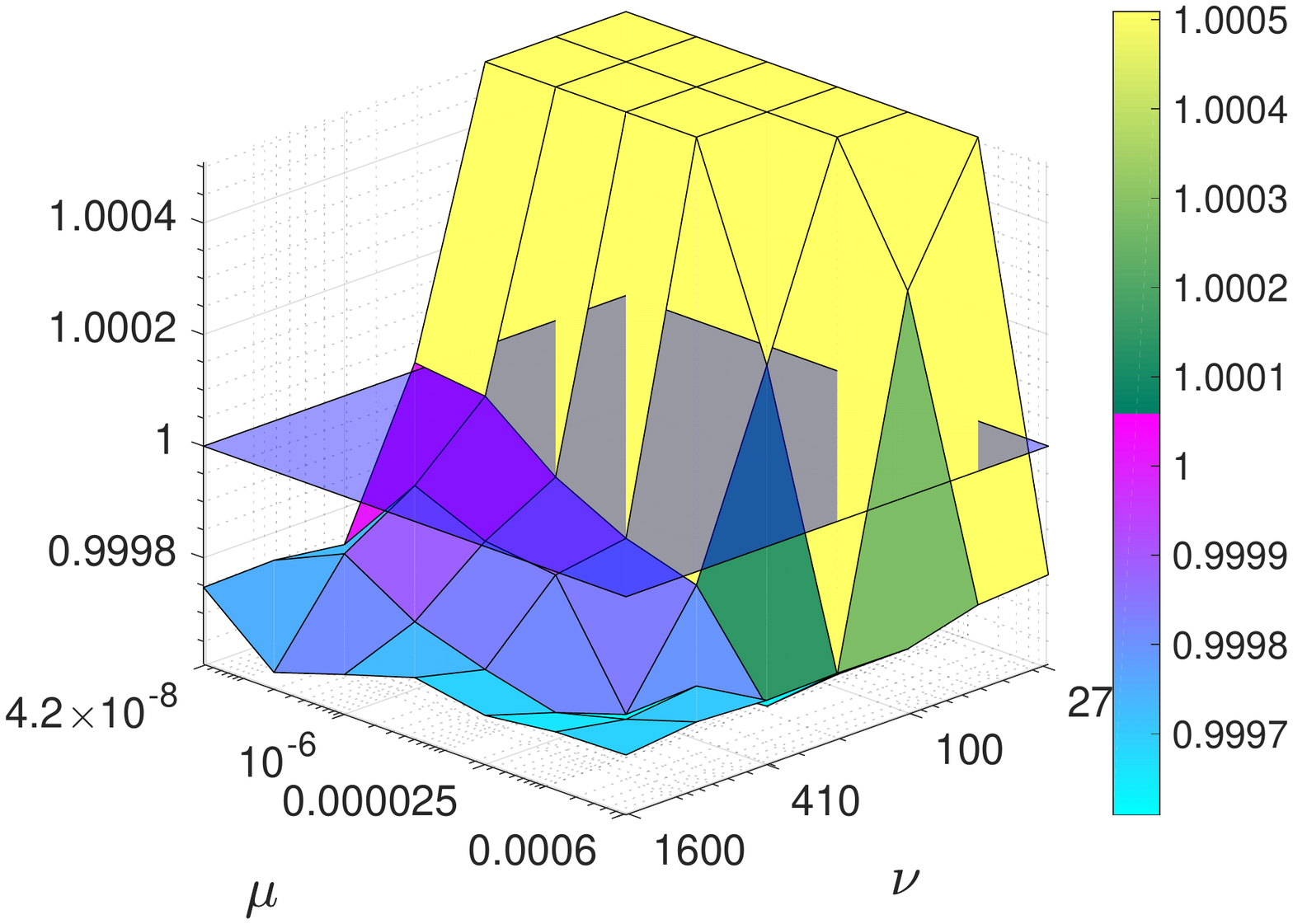}
\end{minipage}%
       \caption{Accelerated matrix inversion on synthetic data. Sensitivity to acceleration parameters. Eigenvalues of $\mA$ are set to $1,1000, 1000, \dots,1000$. From left to right we have: Coordinate sketches with convenient probabilities, coordiante sketches with uniform probabilities and Gaussian sketches. Choice of parameters as per \eqref{eq:ami_munu_conv_paper} in the middle of plots. Each instance was run for 10 seconds.}
\label{fig:ami_heat_gauss_rand}
\end{figure}

The crucial aspect to make the accelerated algorithm to converge is to set $\nu$ large enough.  In fact, combination of both small $\nu$ and small $\theta$ leads almost always to non-convergent algorithm. On the other hand, it seems that once $\nu$ is chosen correctly, big enough $\theta$ leads to fast convergence. This indicates how to compute $\theta$ in practice (recall that computing $\nu$ is feasible)---one needs just to choose it small enough (definitely smaller than $\frac{1}{\nu}$).

\subsection{{\tt BFGS} optimization method}
\label{sec:ami_accBFGSexp}
 We test Algorithm~\ref{alg:ami_bfgs_opt} on several logistic regression problems using data from  LibSVM~\cite{chang2011libsvm}.  In all our tests we centered and normalized the data, included a bias term (a linear intercept), and choose the regularization parameter as $\lambda =1/m$, where $m$ is the number of data points. To keep things as simple as possible, we also used a fixed stepsize which 
was determined using grid search.  Since our theory regarding the choice for the parameters $\theta$ and $\nu$ does not apply in this setting, we simply probed the space of parameters manually and reported the best found result, see Figure~\ref{fig:ami_australian}. In the legend  we use \texttt{BFGS}-a-$\theta$-$\nu$ to denote the accelerated {\tt BFGS} method (Algorithm~\ref{alg:ami_bfgs_opt}) with parameters $\theta$ and $\nu$.


\begin{figure}[!h]
    \centering
\begin{minipage}{0.35 \textwidth}
\includegraphics[width =  \textwidth ]{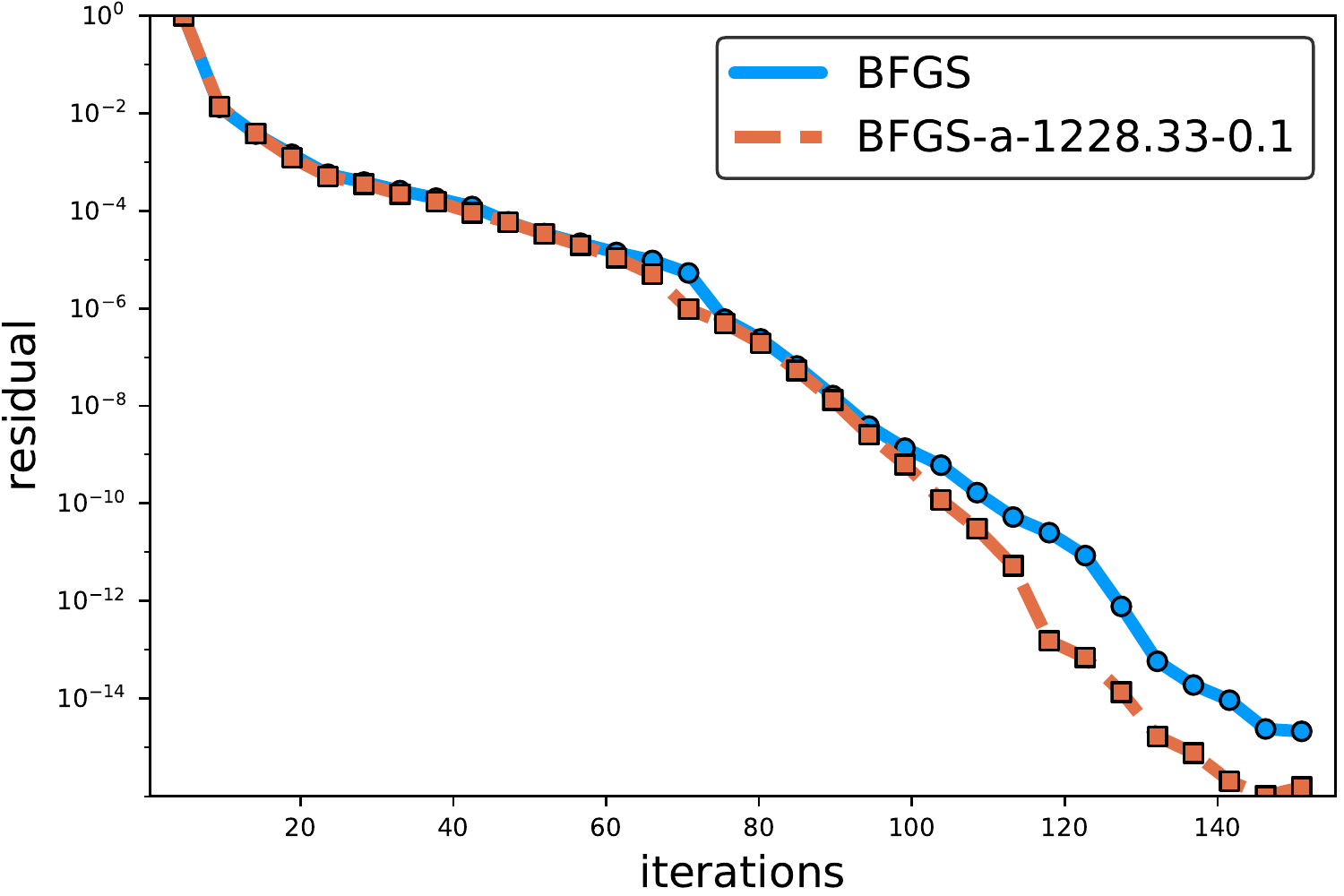}
\end{minipage}%
\begin{minipage}{0.32 \textwidth}
\includegraphics[width =  \textwidth ]{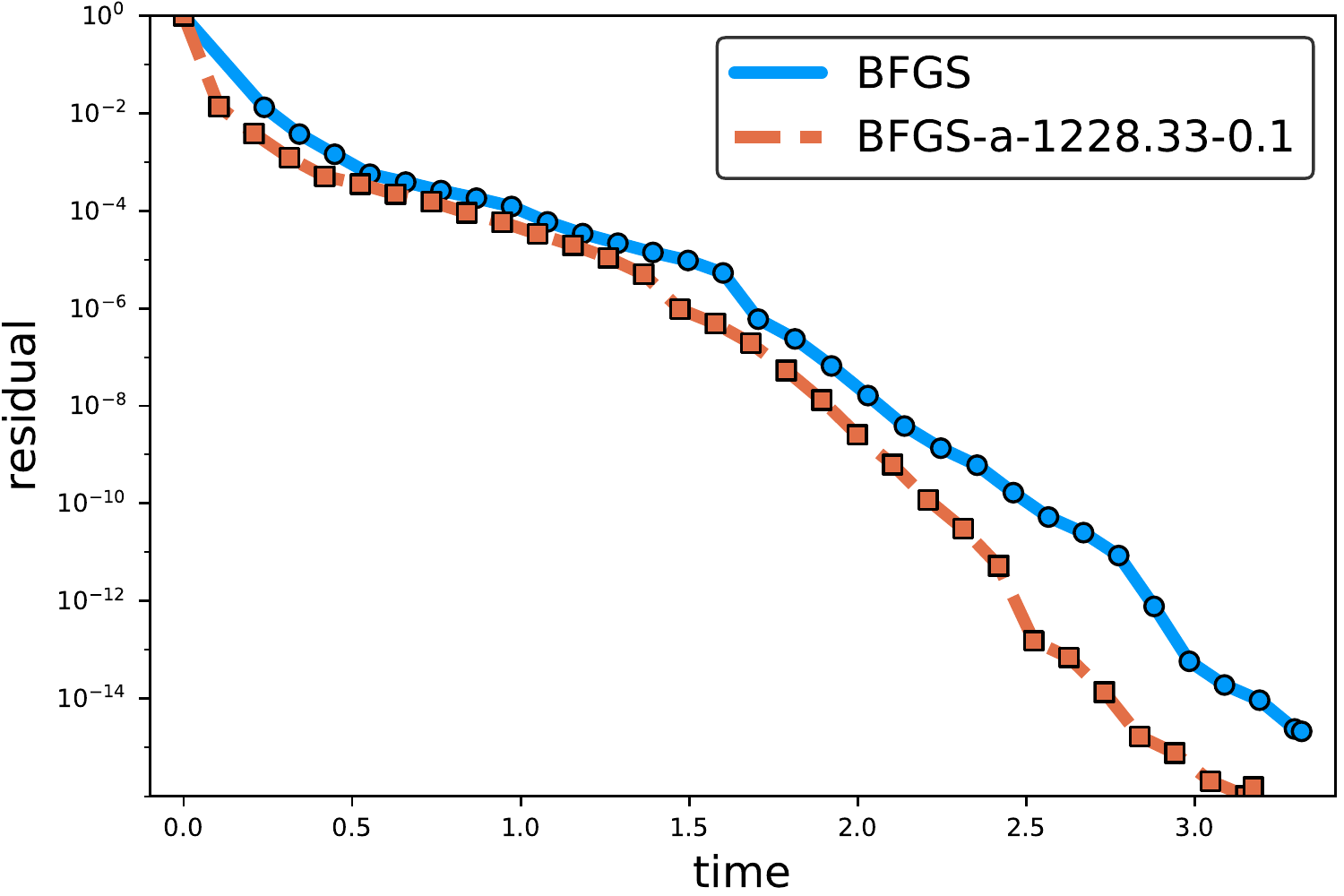}
\end{minipage}%
\\
\begin{minipage}{0.35 \textwidth}
\includegraphics[width =  \textwidth ]{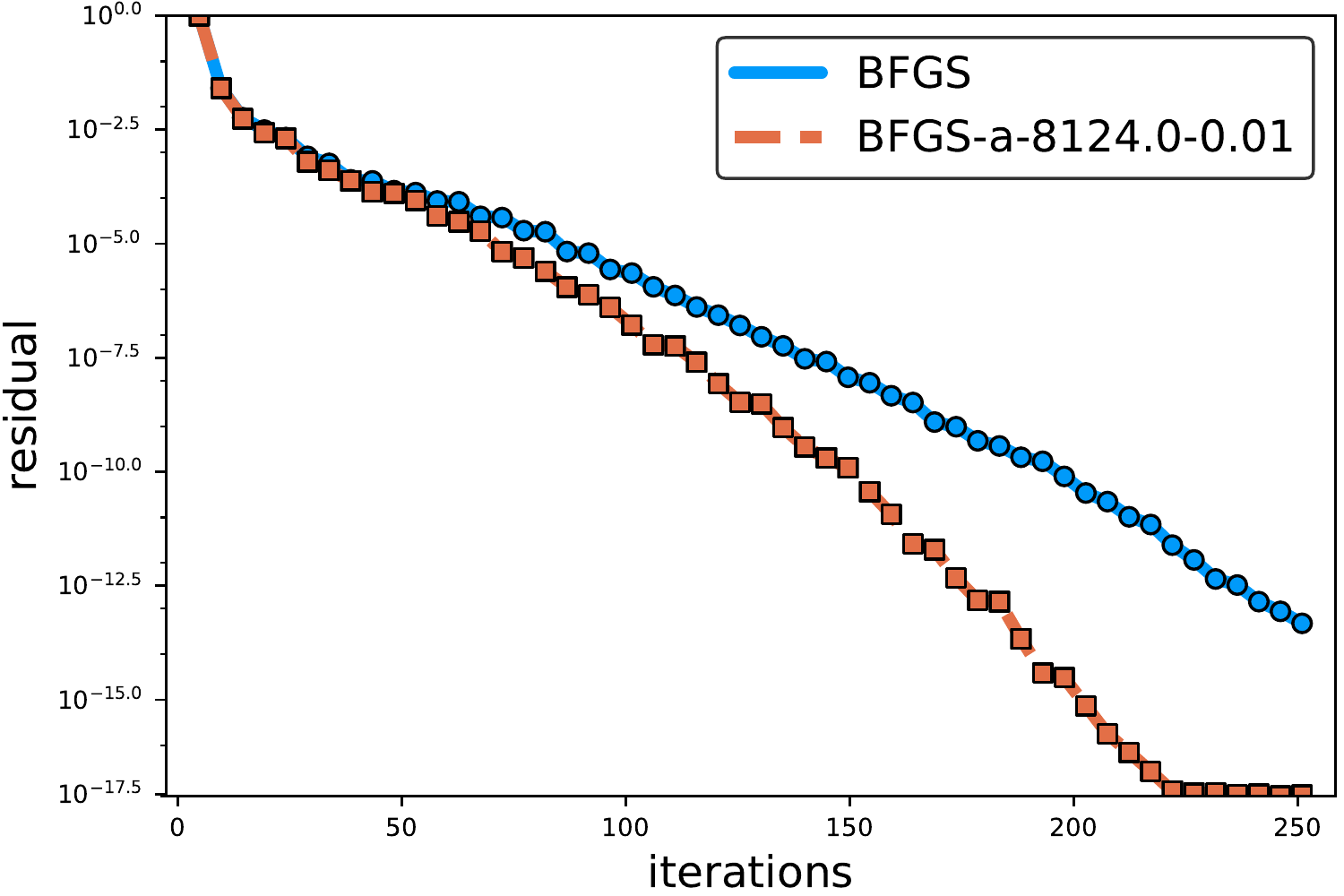}
\end{minipage}
\begin{minipage}{0.35 \textwidth}
\includegraphics[width =  \textwidth ]{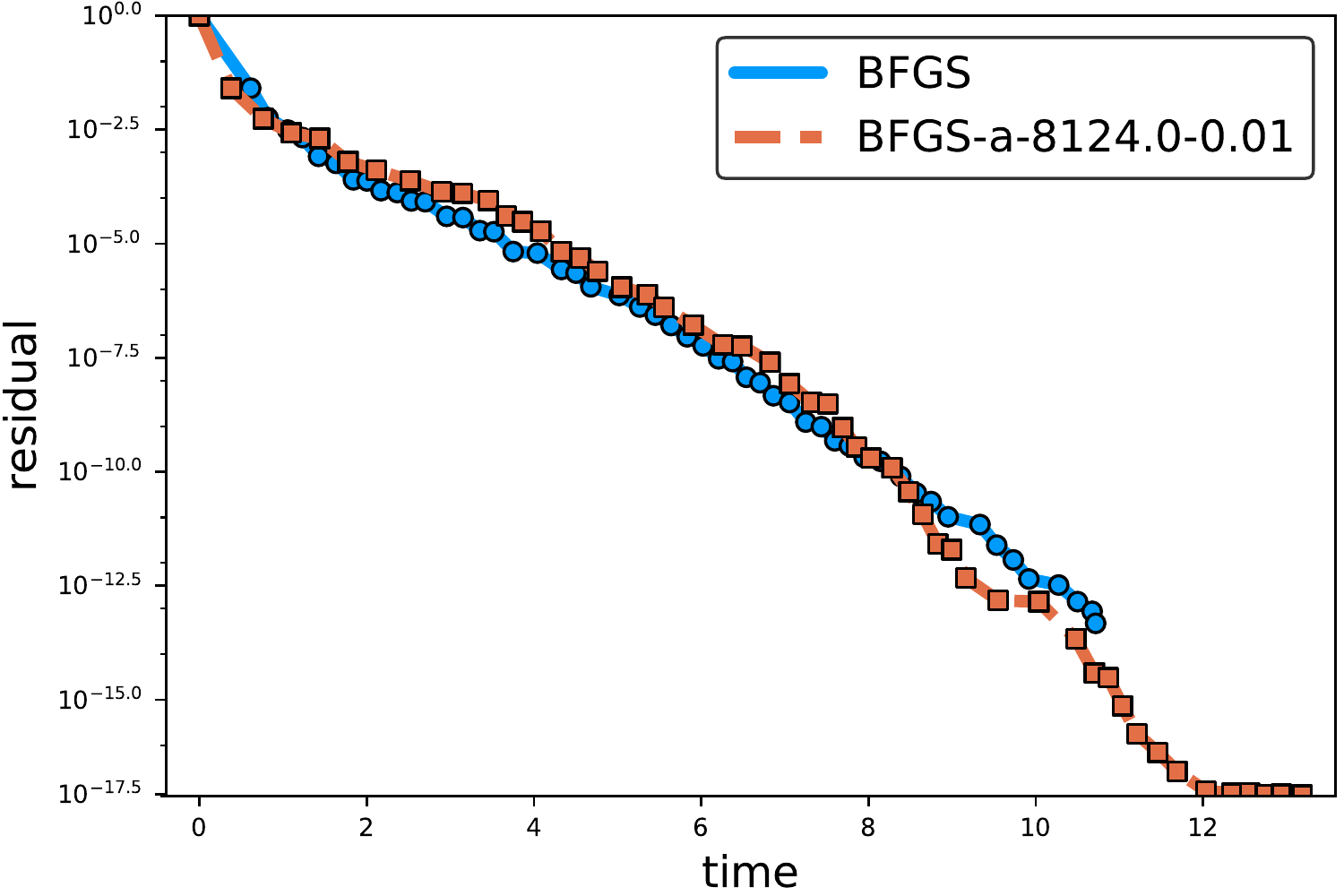}
\end{minipage}
\\
\begin{minipage}{0.35	 \textwidth}
\includegraphics[width =  \textwidth ]{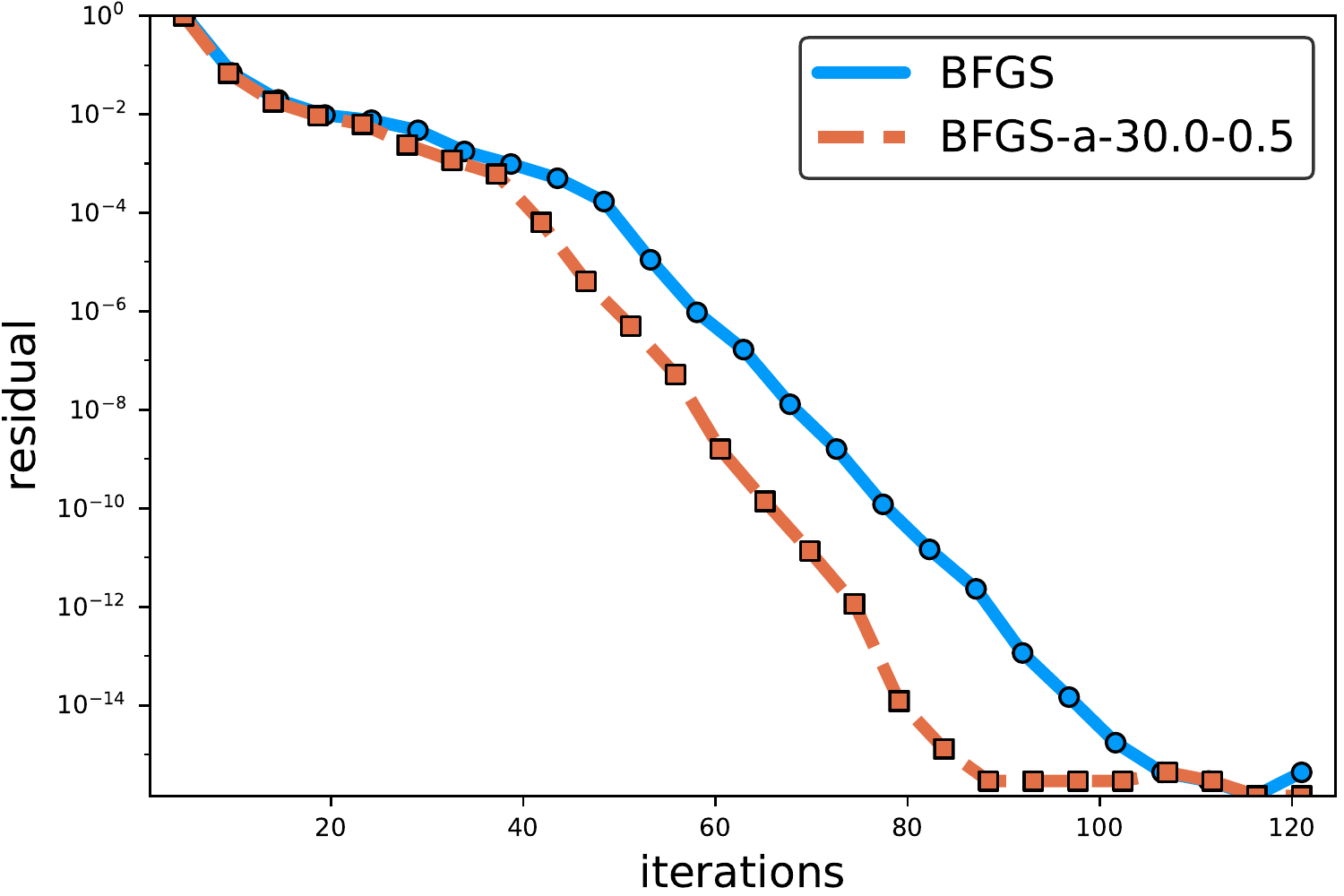}	
\end{minipage}%
\begin{minipage}{0.35	 \textwidth}
\includegraphics[width =  \textwidth ]{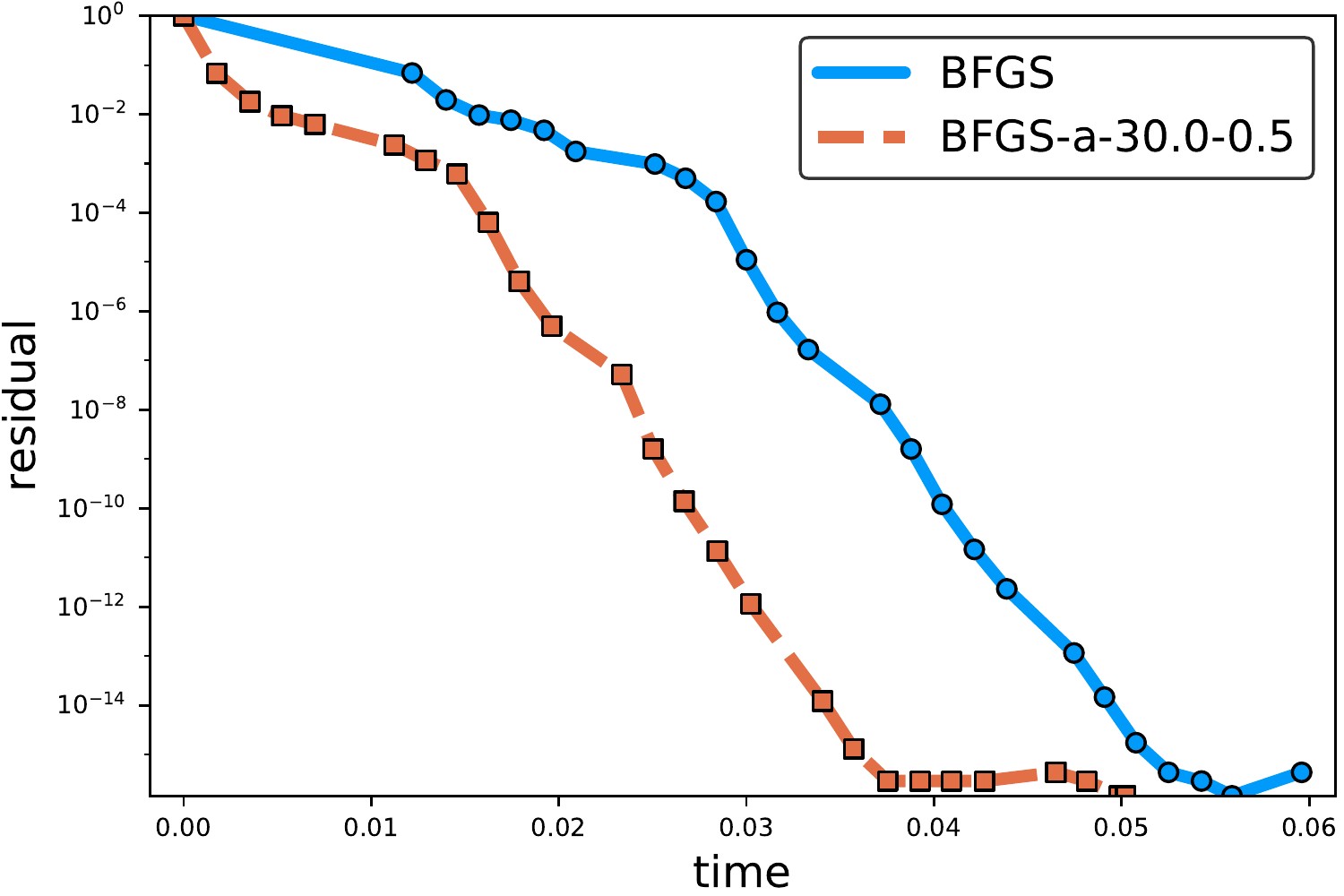}	
\end{minipage}%
\\
\begin{minipage}{0.35 \textwidth}
\includegraphics[width =  \textwidth ]{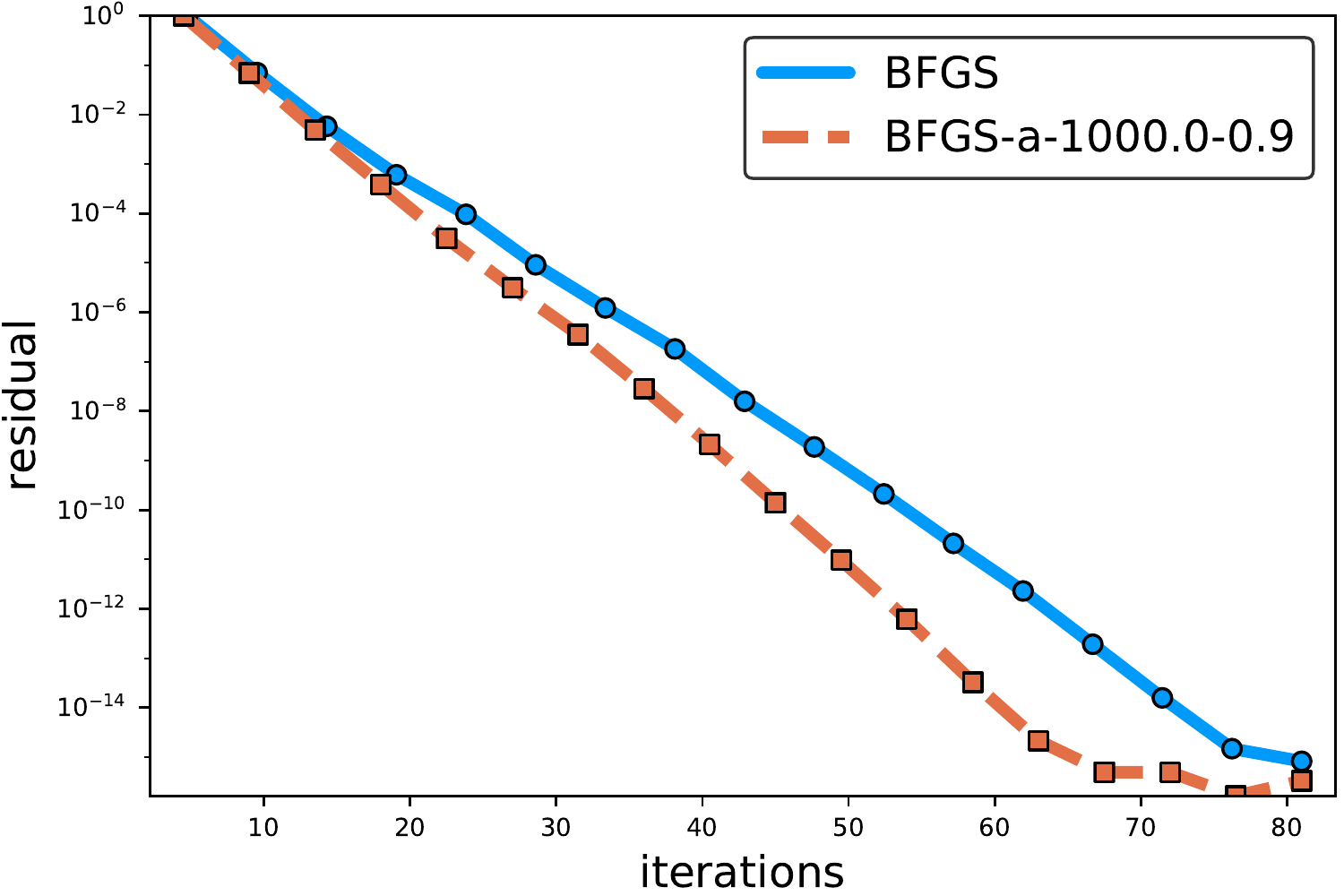}
\end{minipage}%
\begin{minipage}{0.35 \textwidth}
\includegraphics[width =  \textwidth ]{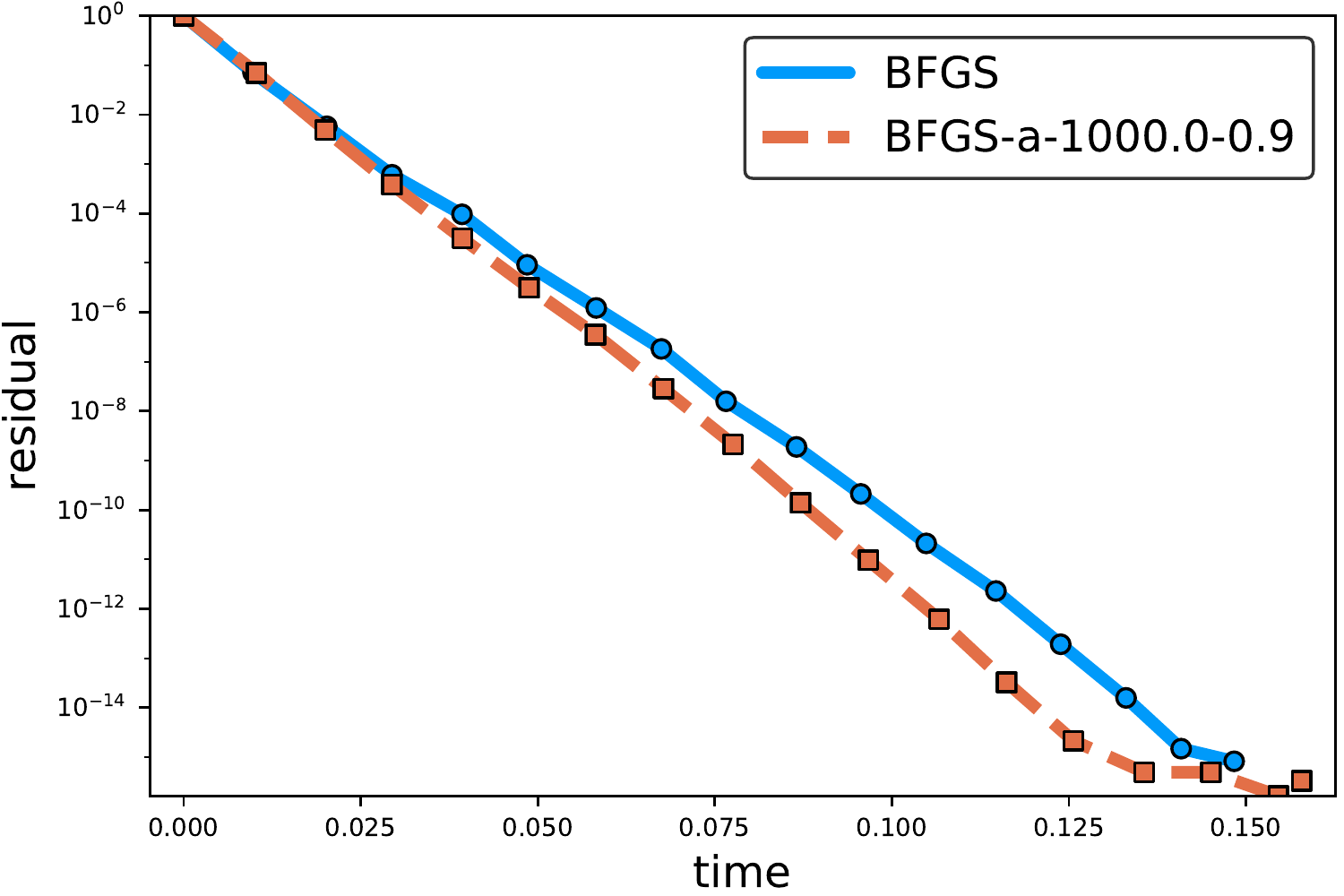}
\end{minipage}%
\caption{Algorithm~\ref{alg:ami_bfgs_opt} ({\tt BFGS} with accelerated matrix inversion quasi-Newton update) vs standard {\tt BFGS}. Left column: time, right column: iteration. From top to bottom:  \texttt{phishing}, \texttt{mushrooms}, \texttt{australian} and \texttt{splice} dataset.}
\label{fig:ami_australian}
\end{figure}

On all four datasets, our method outperforms the classic {\tt BFGS} method, indicating that replacing classic {\tt BFGS} update rules for learning the inverse Hessian by our new accelerated rules can be beneficial in practice.

Much like the \texttt{phishing} problem in Figure~\ref{fig:ami_australian}, the problems \texttt{madelon}, \texttt{covtype} and \texttt{a9a} in Figure~\ref{fig:ami_bfgs_opt_libsvm} did not benefit that much from acceleration.

Indeed, we found in our experiments that even when choosing extreme values of $\theta$ and $\nu$, the generated inverse Hessian would not significantly deviate from the estimate that one would obtain using the standard {\tt BFGS} update. Thus on these two problems there is apparently little room for improvement by using acceleration.    

\begin{figure}[H]
	\centering
	\begin{minipage}{0.3\textwidth}
		\centering
		\includegraphics[width = \textwidth ]{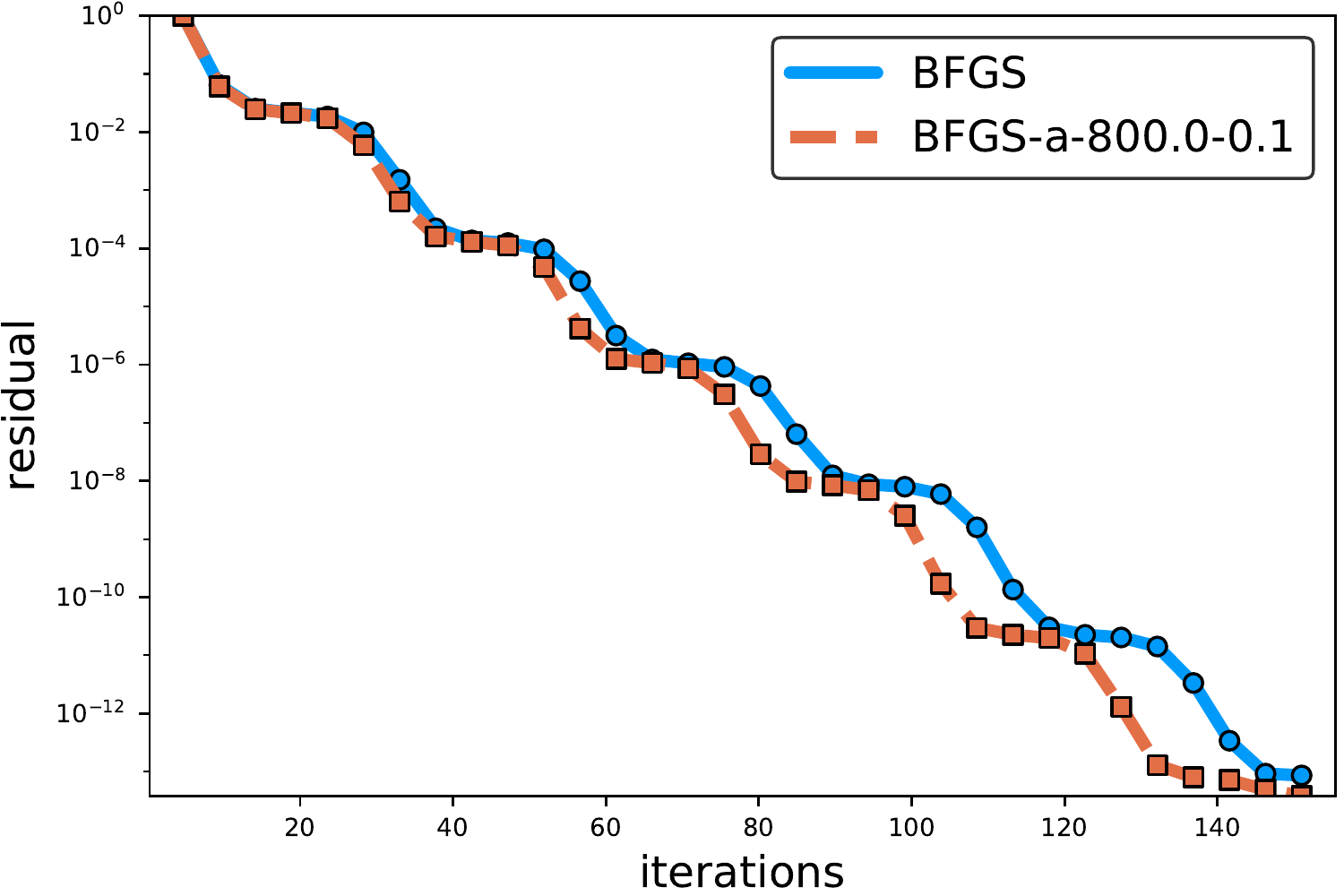}
	\end{minipage}
	\begin{minipage}{0.3\textwidth}
		\centering
		\includegraphics[width=\textwidth]{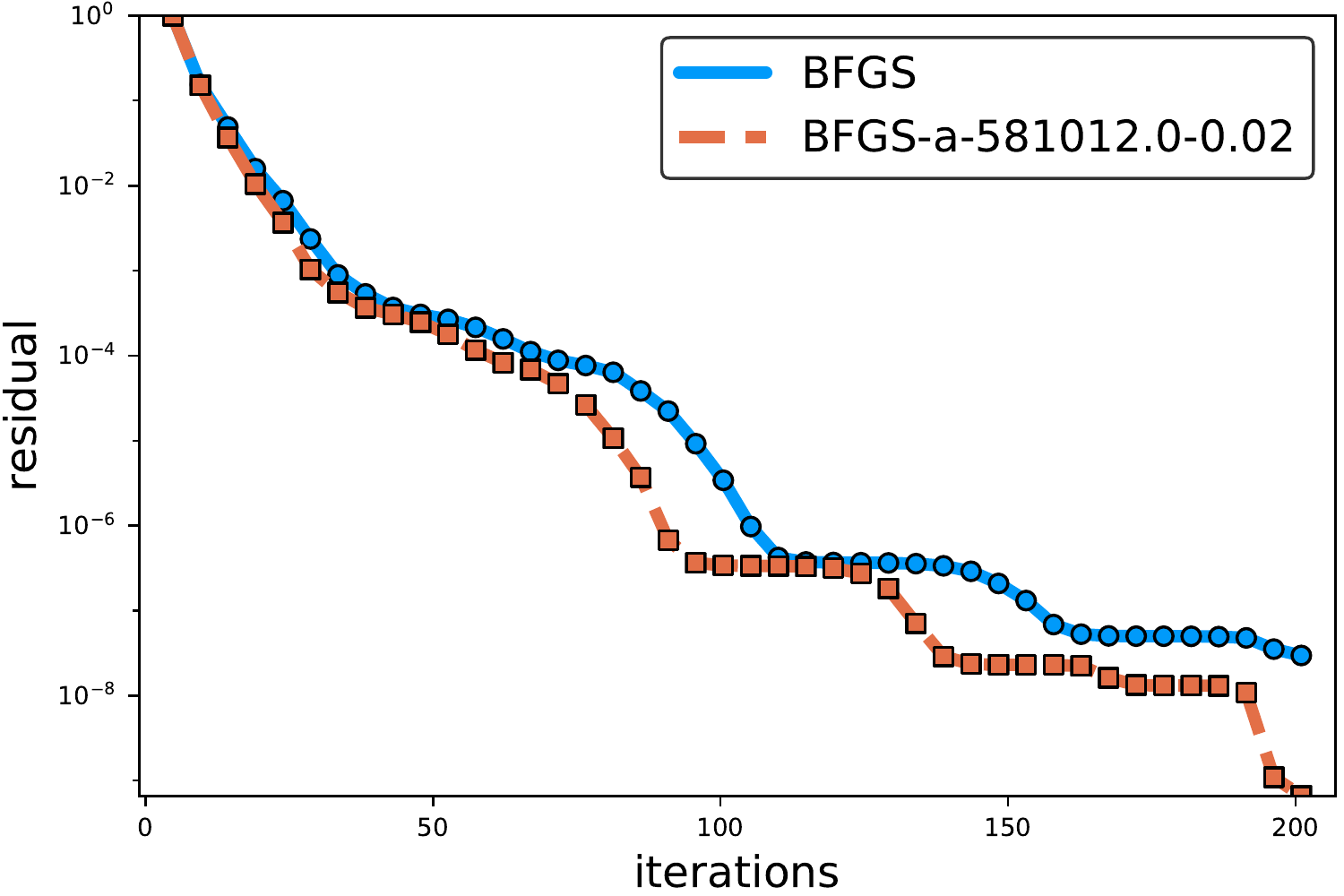}
	\end{minipage}%
		\begin{minipage}{0.3\textwidth}
			\centering
			\includegraphics[width=\textwidth]{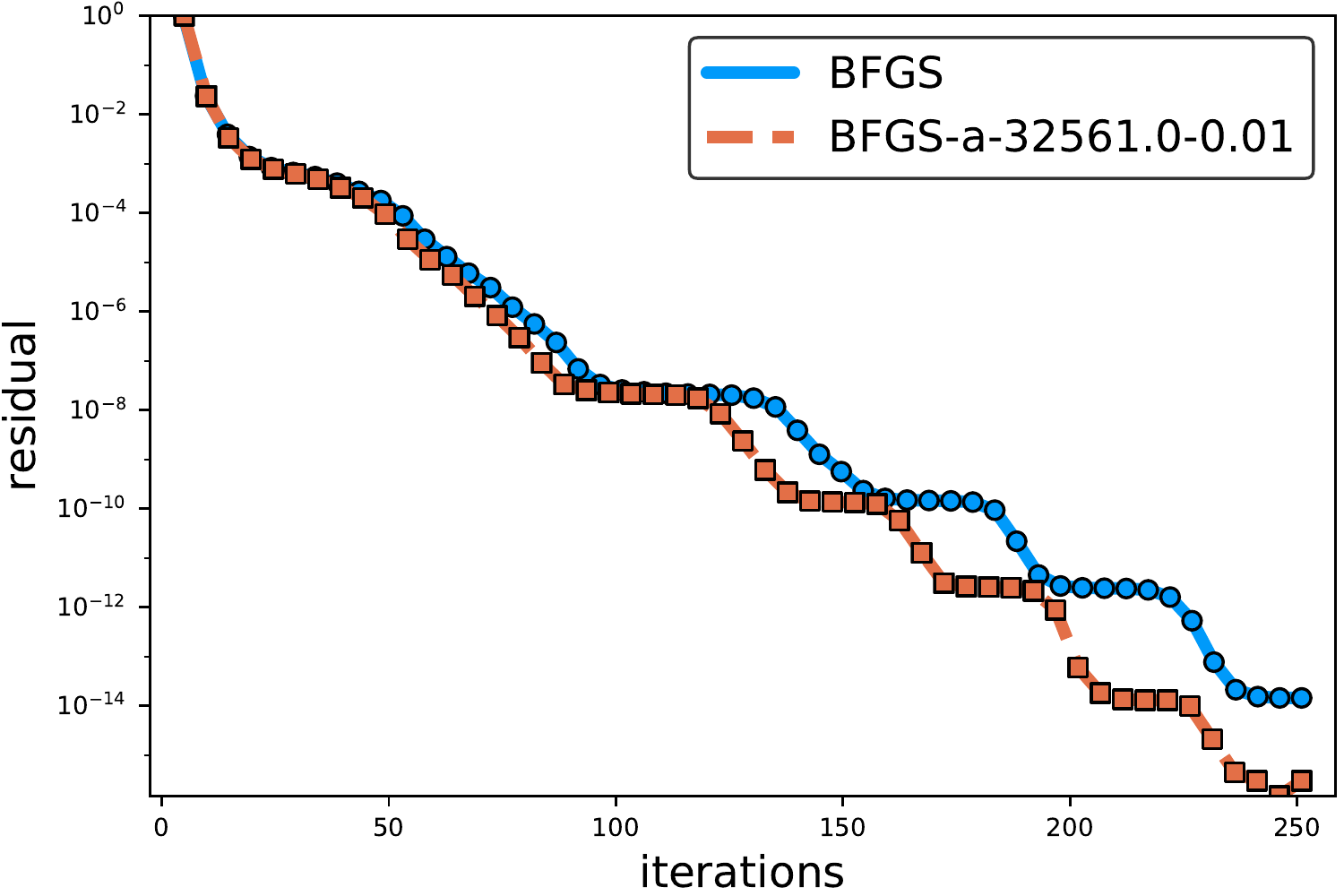}	
		\end{minipage}%
	\caption{Accelerated {\tt BFGS} applied on real data. Left to right:  \texttt{madelon},\texttt{covtype},\texttt{a9a} }
		\label{fig:ami_bfgs_opt_libsvm}
\end{figure}

\section{Conclusion}
In this chapter, we developed an accelerated sketch-and-project method for solving linear systems in Euclidean spaces. The method was applied to invert positive definite matrices, while keeping their symmetric structure. Our accelerated matrix inversion algorithm was then incorporated into an optimization framework to develop both accelerated stochastic and deterministic {\tt BFGS}, which to the best of our knowledge, are {\em the first  accelerated quasi-Newton updates.} 

We show that under a careful choice of the  parameters of the method, and depending on the problem structure and conditioning, acceleration might result into significant speedups both for the matrix inversion problem and for the stochastic {\tt BFGS} algorithm. We confirm experimentally that our accelerated methods can lead to speed-ups when compared to the classical {\tt BFGS} algorithm.

As a future line of research, it might be interesting to study the accelerated {\tt BFGS} algorithm (either deterministic or stochastic)  further, and provide a convergence analysis on a suitable class of functions. Another interesting area of research might be to combine accelerated {\tt BFGS} with limited memory \cite{liu1989limited} or engineer the method so that it can efficiently compete with first-order algorithms for some empirical risk minimization problems, such as, for example \cite{sbfgs}. 

As we show in this work, {\em Nesterov's acceleration can be applied to quasi-Newton updates}. We believe this is a surprising fact, as quasi-Newton updates have not been understood as optimization algorithms, which prevented the idea of applying acceleration in this context.

Since since second-order methods are becoming more and more ubiquitous in machine learning and data science, we hope that our work will motivate further advances at the frontiers of big data optimization.

Note that this whole chapter was devoted to accelerating and algorithm for solving linear systems, and applying the obtained knowledge to obtain an algorithm to solve general convex optimization~\eqref{eq:ami_opt_main}. In the next chapter, we will tackle problem~\eqref{eq:ami_opt_main} directly, developping a specific accelerated stochastic algorithm.



\chapter{Concluding Remarks}

\section{Summary}

In this work, we developed a number of stochastic iterative optimization algorithms with primary focus on solving supervised machine learning problems cast in the  form of regularized empirical risk minimization problems of the form
 \begin{equation}\label{eq:finitesum_conclusion}
\min_{x\in \R^d} \left\{  F(x) \eqdef \underbrace{\frac1n \sum_{i=1}^n f_i(x)}_{\eqdef f(x)} + \psi(x) \right\}.
 \end{equation}

Each chapter of the thesis introduced a state-of-the-art approach for solving~\eqref{eq:finitesum_conclusion} under further assumptions on both the problem structure and the oracle model.

In Chapter~\ref{acd}, we considered an instance of~\eqref{eq:finitesum_conclusion} in the regime where  the  dimension $d$ is very large, and in the simplified setting with $n=1$ and $\psi\equiv 0$. In high-dimensional  optimization, variants of coordinate descent methods reign supreme.  In this chapter, we proposed an accelerated coordinate descent ({\tt ACD}) method, explicitly allowing for {\em arbitrary} sampling. This generalized previous results which considered sampling of a single coordinate only. We further designed a novel non-uniform minibatch sampling strategy can provably outperform uniform minibatch sampling, which is the first result of its kind in the literature. The mentioned sampling can be applied in many contexts beyond coordinate descent methods, as demonstrated in Chapters~\ref{jacsketch} and~\ref{sigmak} where we applied this sampling in the context of stochastic gradient descent methods. 

In Chapter~\ref{sega}, we considered problem \eqref{eq:finitesum_conclusion} in the $n=1$ setting, also assuming that $d$ is very large. As mentioned above, in this regime, randomized {\tt CD} methods are the state of the art. However, as they may not converge otherwise, these methods are always studied either in the nonregularized regime, or with a separable regularizer. With a goal to lift this fundamental limitation, in this chapter we designed new variants of {\tt CD} methods which can provably work with any, even nonseparable, regularizer $\psi$. We further consider a more general subspace gradient oracle: we allow our method access to  gradients of $f$ over a random subspace of $\R^d$ with a very general notion of  randomness going further than that used in randomized {\tt CD} methods. Our algorithm, {\tt SEGA}, can be interpreted as a variance-reduced {\tt CD} method, with a novel variance reduction strategy aimed at removing the adverse effect the (nonseparable) regularizer has on classical {\tt CD} methods. {\tt SEGA} is  the first {\tt CD}-type method capable of converging even with nonseparable regularizers. Moreover, when specialized to random subspaces spanned coordinate vectors, our convergence results for {\tt SEGA} match the state-of-the-art convergence rates of {\tt CD} methods, up to a small constant factor (less than 10). For instance, {\tt SEGA} can be coupled with both importance sampling and Nesterov's acceleration, achieving rates similar to those of the {\tt ACD} method developed in Chapter~\ref{acd}.

In Chapter~\ref{99}, we considered a distributed optimization problem with $n$ workers and a centralized parameter server of the form~\eqref{eq:finitesum_conclusion}, with $f_i$ interpreted as the loss of model $x$ on the data stored on the $i$th worker.  We proposed a generic technique for reducing worker$\rightarrow$server communication in several popular iterative methods for solving \eqref{eq:finitesum_conclusion} by a factor up to $n$, without hurting the convergence rate by more than a small constant. The key idea is based on combining the underlying iterative solver with a novel independent coordinate descent / random sparsification mechanism. We demonstrated that our approach can be incorporated on top of distributed implementations of various algorithms such as {\tt GD}, {\tt SGD}, or {\tt SAGA}. 

In Chapter~\ref{jacsketch}, we proposed and analyzed a remarkably general randomized algorithm for solving \eqref{eq:finitesum_conclusion}---Generalized Jacobian Sketching method ({\tt GJS})---capable of reducing variance coming both from subsampling the data and parameters.  {\tt GJS} is the first variance reduced method with this property.  This was enabled by studying gradient estimators arising from arbitrary sketching operators applied to the Jacobian matrix of the mapping $x\mapsto (f_1(x),\dots, f_n(x))$.  {\tt GJS} recovers many well-known and recently developed algorithms as a special case, including {\tt SAGA}~\cite{saga, qian2019saga}, {\tt LSVRG}~\cite{hofmann2015variance, kovalev2019don},  {\tt SEGA} from Chapter~\ref{sega} and {\tt ISEGA} from Chapter~\ref{99}. Remarkably, our convergence theory for {\tt GJS} either recovers the best-known convergence rate in each special case, or improves upon the current best rates. Besides the unification and improvement upon the well-established algorithms, {\tt GJS} can specialize to a large number of new specific methods with intriguing properties. 

In Chapter~\ref{sigmak}, we go even one step further in our unification efforts.  In particular,  we proposed a framework capable of analyzing both variance-reduced  and non-variance reduced variants of {\tt SGD} at the same time. This is the first framework with this property. Further, our framework is capable of accurate modeling  standard, parallel and distributed {\tt SGD} methods, with gradient estimator formed through various mechanisms, including subsampling (e.g., minibatching and importance sampling) and compression (e.g., sparsification and quantization) and their combination. The development of our framework was motivated by the need establish a general theory capable of taming large swaths of the almost impenetrable wilderness of {\tt SGD} methods, while at the same time facilitating faster development of new variants. Our framework includes {\tt GJS} as a special case, as well as many other methods.

In Chapter~\ref{asvrcd}, we establish a novel and fundamental link between the novel variance reduced {\tt CD} methods first developed in this thesis, and the world of {\tt SGD} methods for finite-sum optimization. In particular, we were able to show that our {\tt SEGA} (resp.\ {\tt SVRCD}\footnote{{\tt SVRCD} is an algorithm similar to {\tt SEGA}, proposed in Chapter~\ref{jacsketch}.}) method reduces to the well-known {\tt SAGA} (resp.\ {\tt LSVRG}) method when applied to a carefully constructed problem with a special regularizer. Moreover, our general theory is able to recover the best best-known oracle complexity for these methods. We have also tightened the analysis of both {\tt SEGA} and {\tt SVRCD}, so that these methods are capable of exploiting the structure of the regularizer $\psi$ for faster convergence.   Further,  we succeeded in incorporating Nesterov's momentum into the {\tt SVRCD} algorithm. The resulting method---{\tt ASVRCD}---is the first accelerated variance reduced coordinate descent method. Moreover, in a special case, {\tt ASVRCD} reduces to a (variant of) the celebrated Katyusha algorithm~\cite{allen2017katyusha,  l-svrg-as}, thus achieving the optimal rate for  finite-sum problems.

In Chapter~\ref{local}, we introduced a novel optimization formulation of federated learning aimed at allowing device-specific personalization. Unlike the standard federated learning formulation which seeks to find a single global model to be used on all devices, we constructed a separate loss for each device, with an extra penalty that ensures the local models do not deviate too much from each other. We first applied standard {\tt SGD} to our two-sum formulation, and observed that the resulting method is a novel variant of the celebrated local gradient descent method. However, in sharp contrast with all preceding formulations of federated learning, we showed that local methods need fewer communication rounds if more personalization is desired. This is the first result in the literature suggesting that local {\tt GD} methods can lead to communication complexity benefits. However, we went much beyond this in the chapter, proposing a number of new methods capable of working with a regularizer, capable of achieving variance reduction and partial participation.

In Chapter~\ref{sscn}, we proposed a second-order subspace descent algorithm ({\tt SSCN}) designed to solve~\eqref{eq:finitesum_conclusion} with $n=1$, separable $\psi$, and large $d$. We proved that {\tt SSCN}  enjoys a global convergence rate and a fast local convergence rate. While our global result interpolates between the rate of {\tt CD} and the rate of the cubically regularized Newton method of Polyak and Nesterov, the local rate is identical to the convergence of stochastic subspace descent applied to minimizing quadratic function $\frac12 (x-x^*)^\top \nabla^2 f(x^*) (x-x^*)$, which we find remarkable.

In Chapter~\ref{ami} we developed an accelerated stochastic algorithm for solving linear systems in Euclidean spaces. Our method can be specialized to obtain a large class of accelerated algorithms designed to invert a positive definite matrix. In particular, we were able to accelerate the subroutine employed in stochastic quasi-Newton methods which updates the inverse Hessian estimator. Despite more than half a century of research into quasi-Newton methods, we have developed the first provably accelerated quasi-Newton matrix update formula.

\section{Future Research Work}

In this section, we outline a few challenges that remain open problems to be addressed in the future. 

\begin{itemize}

\item {\bf Unified framework for accelerated stochastic algorithms.} We believe that an  accelerated variant of the {\tt GJS} algorithm and/or a general analysis of accelerated stochastic algorithms analogous to Chapter~\ref{sigmak} would be of immense value. Such  results would immediately lead to countless optimal algorithms in terms of  oracle complexity. 

\item {\bf Understanding Nesterov's acceleration.} Throughout this thesis, we incorporated some form of Nesterov's acceleration into several algorithms. The effect of the acceleration mechanism, however, slightly varied from chapter to chapter.\footnote{This is consistent with the related literature. } This can be well demonstrated comparing the effect of the sampling on {\tt ACD} and {\tt ASVRCD}. Specifically, the optimal sampling for {\tt SVRCD} and {\tt ASVRCD} are identical, while the optimal sampling for {\tt CD} and {\tt ACD} are vastly different. The effect of the acceleration on the sketch-and-project algorithm studied in Chapter~\ref{ami} is even more complex. We are certain that a better understanding of randomized algorithms with acceleration would enable the development of a broader range of methods, such as accelerated {\tt SSCN}, for example.

\item {\bf {\tt SSCN} with non-separable $\psi$.} In Chapter~\ref{sega}, we discovered a mechanism  allowing {\tt CD} algorithms to deal with a non-separable regularizer $\psi$. Can  a similar result be established for second-order algorithms? We believe that control variates---a tool from statistics widely used throughout this work---might help to resolve the issue. However, to prove a tight convergence rate of such a method is highly non-trivial, especially since the literature on variance-reduced second-order methods is very limited at the moment.

\end{itemize}

\begin{onehalfspacing}
\renewcommand*\bibname{\centerline{REFERENCES}} 
\addcontentsline{toc}{chapter}{References}
\newcommand{\BIBdecl}{\setlength{\itemsep}{0pt}}

\bibliographystyle{plain}

\bibliography{sorted}

\end{onehalfspacing}

\appendix

\newpage

\begingroup
\let\clearpage\relax
\begin{center}
\vspace*{2\baselineskip}
{ \textbf{{\large APPENDICES}}} 
\addcontentsline{toc}{chapter}{Appendices} 
\end{center}
\endgroup

\chapter{ Table of Frequently Used Notation~\label{sec:table}}

\renewcommand{\EE}{\mathbb{E}}

{
\footnotesize

 \begin{longtable}{| p{.20\textwidth} | p{.80\textwidth}| } 
\hline
\multicolumn{2}{|c|}{{\bf \normalsize For all chapters} }\\
\hline
\multicolumn{2}{|c|}{ Basic }\\
\hline
$\E{\cdot}$,  $\Probbb{\cdot}$  & Expectation / Probability \\
$\langle \cdot ,\cdot \rangle$,  $\| \cdot \|$ & Standard inner product and norm in $\R^d$: $\langle x ,y \rangle =\sum_{i=1}^d x_iy_i$; $\| x \| =\sqrt{\langle x ,x \rangle} $ \\
$\langle \cdot ,\cdot \rangle_{\mB}$,  $\| \cdot \|_{\mB}$ & Weighted inner product and norm in $\R^d$: $\langle x ,y \rangle =x^\top \mB y$; $\| x \| =\sqrt{\langle x ,x \rangle_\mB} $ \\
$\lambda_{\max} (\cdot), \lambda_{\min}(\cdot)$ & Maximal eigenvalue / minimal eigenvalue   \\
$\nabla h(x)$& Gradient of a differentiable function $h$ \\
$\nabla_i h(x)$& $i$th partial derivative of a differentiable function $h$  \\
$\nabla^2 h(x)$& Hessian of a twice differentiable function $h$  \\
\hline
\multicolumn{2}{|c|}{ Objective }\\
\hline
  $d$ & Dimensionality of space $x\in \RR^d$  \\
 $F: \R^d \rightarrow \R$ & Objective function \\
  $f: \R^d \rightarrow \R$ & Smooth part of the objective, often of finite sum structure ($=\frac1n \sum_{i=1}^n f_i(x)$) \\
    $f_j: \R^d \rightarrow \R$ & Differentiable convex function ($1\leq j \leq n$) \\
 $\psi: \R^d \rightarrow \R \cup \{ +\infty \}$ & Non-smooth part of the objective\\
 $x^*$ & Global optimum of \eqref{eq:finitesum}\\
 $F^{*}$ & $ \eqdef F(x^{*})$, the optimum value of the objective \\
$L$, $\mM$ & Smoothness constant/smoothness matrix of $f$ \\
$\mu$ & Strong convexity of $f$  \\
  \hline
 \multicolumn{2}{|c|}{ Linear operators $\R^{d\times n} \to \R^{d\times n}$} \\
 \hline
  $\cA$ & A generic  linear operator \\ 
  $\cA^*$ & Adjoint of $\cA$: $\langle \cA \mX, \mY \rangle \equiv \langle  \mX, \cA^*\mY \rangle $ for all $\mX,\mY\in \R^{d\times n}$\\ 
   $\cA^\dagger$ & Moore Penrose pseudoinverse of $\cA$ \\ 
   $\Range{\cA}$ & Image (range space) of $\cA$: $\Range{\cA} \eqdef \{\cA \mX \;:\; \mX\in \R^{d\times n}\}$ \\
   $\Range{\cA}^\perp$ & Orthogonal complement of $\Range{\cA}$ \\   
      $\Null{\cA}$ & Kernel (null space) of $\cA$: $\Null{\cA} \eqdef \{ \mX \in \R^{d\times n} \;:\; \cA \mX = 0\}$ \\
 $\cI$ & Identity operator: $\cI \mX \equiv \mX$\\ 
 \hline
\multicolumn{2}{|c|}{ Other }\\
\hline
$e_i$ & $i$th vector from the standard basis  \\
$\mI$ & Identity matrix  \\
  $\prox_{\alpha \psi}(x)$ & Proximal operator of $\psi$: $\prox_{\alpha \psi}(x) \eqdef \argmin_{u\in \R^d} \{ \alpha\psi(u) + \frac{1}{2}\|u-x\|^2\}$  \\
  $\langle \mX, \mY \rangle$ & Trace inner product of matrices $\mX$ and $\mY$: $\langle \mX, \mY \rangle \eqdef \Tr{\mX^\top \mY}$ \\
  $\|\mX\| $ & Frobenius norm of matrix $\mX$: $\|\mX\| = \langle \mX, \mX\rangle^{\frac{1}{2}} $\\
  $\mX \circ \mY$ & Hadamard product: $(\mX \circ \mY)_{ij} = \mX_{ij} \mY_{ij}$  \\
    $\mX \otimes \mY$ & Kronecker product  \\
    $ \Tr{\cdot}$ & Trace \\
  $\diag(x)$ & Diagonal matrix with vector $x$ on the diagonal \\
     $[n]$ & Set $\{1,2,\dots,n\}$ \\
          $D_h(x,y)$ & Bregamn distance $D_h(x,y) \eqdef h(x) -h(y)-\langle \nabla h(y),x-y \rangle$ \\
\hline
\caption{Summary of frequently used notation.}\label{tbl:notation}
\end{longtable}

\begin{longtable}{| p{.20\textwidth} | p{.80\textwidth}| } 
\hline
\multicolumn{2}{|c|}{{\bf \normalsize Chapter~\ref{sega}} }\\
\hline
\multicolumn{2}{|c|}{{ Basic} }\\
\hline
  $\cD$ & Distribution over sketch matrices $\mS$ \\
 $\mS$ & Sketch matrix from \eqref{eq:sega_sketch-n-project}\\
  $b$ & Random variable such that $\mS\in \R^{n\times b}$  \\
  $\zeta(\mS, x)$ & Sketched gradient at $x$ from \eqref{eq:sega_sketched_grad}\\
  $\mZ$ &  $\mS \left(\mS^\top  \mS\right)^\dagger\mS^\top$  \\
  $\theta$ & Random variable for which $\E{\theta  \mZ} = \mI$ from \eqref{eq:sega_unbiased}\\
  $\mC$ & $ \E{ \theta^2 \mZ }$ from Theorem~\ref{thm:sega_main}\\
 $h, g$ & Biased and unbiased gradient estimators from \eqref{eq:sega_h^{k+1}}, \eqref{eq:sega_g^k}\\
    $\Lgen$ & Lyapunov function from Theorem~\ref{thm:sega_main}, \\ 
     $\sigma$ & Parameter for Lyapunov function from Theorems~\ref{thm:sega_main}, \ref{t:imp_dacc}\\ 
 \hline
 \multicolumn{2}{|c|}{{ Extra Notation for Section~\ref{sec:sega_CD} }}\\
 \hline
  $p$, $\Probmat$ & Probability vector and matrix \\
    $v$ & vector of ESO parameters from \eqref{eq:sega_ESO}\\
  $\mPdiag,\mVdiag$ & $\diag(p),\diag(v)$  \\
  $\gamma$& $\alpha - \alpha^2\max_{i}\{\frac{v_{i}}{p_{i}}\}-\sigma$ from Theorem~\ref{t:imp_dacc} \\
$y,z$ & Extra sequences of iterates for \texttt{ASEGA}  \\
$\tau,\beta$ & Parameters for \texttt{ASEGA}  \\
  $\Lnacc, \Lacc$ & Lyapunov functions from Theorems~\ref{t:imp_dacc}, \ref{t:imp_acc}\\ 
   $ \TD(v,p)$ & $\max_i \frac{\sqrt{v_i}}{p_i}$  \\
\hline
\caption{Summary of frequently used notation specific to Chapter~\ref{sega}.}\label{tbl:notation_sega}
\end{longtable}

\begin{longtable}{| p{.20\textwidth} | p{.80\textwidth}| } 
\hline
\multicolumn{2}{|c|}{{\bf \normalsize Chapter~\ref{99}} }\\
\hline
\multicolumn{2}{|c|}{{ General} }\\
\hline
  $n$ & Number of parallel workers/machines  \\
  $\tau$ & Ratio of coordinate blocks to be sampled by each machine \\
   $m$ & Number of coordinate blocks \\
   $f_i$ & Part of the objective owned by machine $i$ from \eqref{eq:99_problem}\\
   $L$ & Each $f_i$ is $L$ smooth (Assumption~\ref{as:99_smooth_sc})\\
   $U_i^t$ & Subset of blocks sampled at iteration $t$ and worker $i$  \\
    $g$ & Unbiased gradient estimator  \\
     \hline
     \multicolumn{2}{|c|}{{ {\tt ISAGA}} }\\
 \hline
  $\mJ_j$ & Delayed estimate of $j$th gradient from \eqref{eq:99_saga_alpha_dist},  \eqref{eq:99_saga_alpha_sm}  \\
  $N$ & Finite sum size for shared data problem from \eqref{eq:99_problem_saga_sm} \\
  $l$ & Number of datapoints per machine in distributed setup  from \eqref{eq:99_problem_saga_dist} \\
  $\cL^t $ & Lyapunov function from \eqref{eq:99_saga_lyapunov}\\
 \hline
\multicolumn{2}{|c|}{{ {\tt ISGD} } }\\
 \hline
  $g_i^t$ & Unbiased stochastic gradient; $\EE [g_i^t] = \nabla f_i(x^t)$     \\
  $\sigma^2$ & An upper bound on the variance of stochastic gradients from Assumption \ref{as:99_bounded_noise} \\
   \hline
\multicolumn{2}{|c|}{{ {\tt ISEGA}} }\\
 \hline
  $h_i^t$ & Sequence  of biased estimators for $\nabla f_i(x^t)$  from \eqref{eq:99_sega_h} \\
  $g_i^t$ & Sequence of unbiased estimators for $\nabla f_i(x^t)$  from \eqref{eq:99_sega_g} \\
  $\Lgen^t$ & Lyapunov function from Theorem~\ref{thm:99_sega} \\
\hline
\caption{Summary of frequently used notation specific to Chapter~\ref{99}.}\label{tbl:notation_99}
\end{longtable}
\newpage
\begin{longtable}{| p{.20\textwidth} | p{.80\textwidth}| } 
\hline
\multicolumn{2}{|c|}{{\bf \normalsize Chapter~\ref{jacsketch}} }\\
\hline
 \multicolumn{2}{|c|}{ Sets} \\
 \hline
  $R$ &  random subset (``sampling'') of $[n]$ \\
   $R^k$ &  random subset (``sampling'') of $[n]$ drawn at iteration $k$ \\
  $L$ & a random subset (``sampling'') of $[d]$ \\
   $L^k$ &  random subset (``sampling'') of $[d]$ drawn at iteration $k$ \\ 
  $\pRj$ & probability that $j\in R$ \\
    $\pLi$ & probability that $i\in L$ \\   
 \hline
  \multicolumn{2}{|c|}{ Spaces $\R^n$ and $\R^d$} \\
   \hline
  $\eR \in \R^n$ & vector of all ones in $\R^n$ \\
  $\eL \in \R^d$ & vector of all ones in $\R^d$ \\
    $\eRj \in \R^n$ & $j$th standard unit basis vector in $\R^n$ \\
    $\eLi \in \R^n$ & $i$th standard unit basis vector in $\R^d$ \\   
    $x^k \in \R^d$  & the $k$th iterate produced by {\tt GJS}\\
    $\pL \in \R^d$ & the vector $(\pL_1, \dots,  \pL_d)$ \\
        $\ptL \in \R^d$ & the vector $(\ptL_1, \dots,  \ptL_d)$ \\
    $\pR \in \R^n$ & the vector $(\pR_1, \dots,  \pR_n)$ \\
      $\ptR \in \R^n$ & the vector $(\ptR_1, \dots,  \ptR_n)$ \\
     $\qR \in \R^n$ & the vector $(\qR_1, \dots,  \qR_n)$ \\
         $\qtR \in \R^n$ & the vector $(\qtR_1, \dots,  \qtR_n)$  \\
        $v \in \R^n$ & any vector for which~\eqref{eq:gjs_ESO_saga} holds  \\
          $x^{-1}$ &  elementwise inverse of $x$ \\
     $g^k$ & estimator of the gradient $\nabla f(x^k)$ produced by {\tt GJS} \\
  \hline    
 \multicolumn{2}{|c|}{ Matrices in $\R^{d\times d}$, $\R^{d\times n}$ and $\R^{n\times n}$} \\
 \hline 
 $\mI_d \in \R^{d\times d}$ & $d\times d$ identity matrix \\
  $\mI_n \in \R^{n\times n}$ & $n\times n$ identity matrix \\
  $\mG(x) \in \R^{d\times n}$ & the Jacobian matrix, i.e., $\mG(x) = [\nabla f_1(x), \dots, \nabla f_n(x)]$ \\
  $\mJ^k \in \R^{d\times n}$ & estimator of the Jacobian produced by {\tt GJS} \\
 $\mM_j \in \R^{d\times d}$ & smoothness matrix of $f_j$ (if $\mM_j = m^j \mI_d$, then this specializes to $m^j$-smoothness)\\
 $\mR \in \R^{n\times n}$ & a random matrix we use to multiply $\mJ$ or $\mG$ from the right \\
 $\mR_{R} \in \R^{n\times n}$ & the random matrix  $\mR_R \eqdef \sum_{j\in R} \eRj \eRj^\top$  \\
 $\mL  \in \R^{d\times d}$ & a random matrix we use to multiply $\mJ$ or $\mG$ from the left \\
  $\mL_{L} \in \R^{d\times d}$ & the random matrix  $\mL_L \eqdef \sum_{i\in L} \eLi \eLi^\top$  \\
    $\PR \in \R^{n\times n}$ & Matrix defined by $\PR_{jj'}=\Prob{j\in R, j'\in R}$\\
        $\PtR \in \R^{n\times n}$ & Matrix defined by $\PtR_{jj'}=\Prob{j\in R_\tR, j'\in R_\tR}$\\
  \hline
 \multicolumn{2}{|c|}{ Linear operators $\R^{d\times n} \to \R^{d\times n}$} \\
 \hline
  $\cU$ & any unbiased operator: $\E{\cU \mX} \equiv \mX$, i.e., $\E{\cU} \equiv \cI$\\ 
 $\cS$ &  any random projection operator \\
 $\cM$  & operator defined via $(\cM \mX)_{:j} = \mM_j \mX_{:j}$ \\
 $\cB$   & (a technical) operator used to define the Lyapunov function \eqref{eq:gjs_Lyapunov}   \\
  $\cR$   & (a technical) operator such that $\mJ^k - \mG(x^*)\in \Range{\cR}$  \\
 \hline  
 \multicolumn{2}{|c|}{ Miscellaneous} \\
 \hline
  $\Gamma$  & Random operator  $\Gamma : \R^{d\times n} \to \R^d$ defined by $\Gamma \mX = \cU \mX \eR$ \\
 \hline  
 \caption{Summary of frequently used notation specific to Chapter~\ref{jacsketch}.}\label{tbl:notation_jacsketch}
\end{longtable}
\newpage
\begin{longtable}{| p{.20\textwidth} | p{.80\textwidth}| } 
  \hline  
 \multicolumn{2}{|c|}{{\bf \normalsize Chapter~\ref{sscn}} }\\
 \hline
\multicolumn{2}{|c|}{{ From main body of the chapter} }\\
\hline
  $\mS \in \R^{d, \tau(\mS)}$ & Random matrix sampled from distribution $\cD$ from \eqref{eq:sscn_update_general}\\
    $S$ & Random  subset of $\{1,\dots ,d \}$ from \eqref{eq:sscn_update_general}\\
$M_{\mS}$ & Lipschitz constant of $\nabla^2 f(x)$ on the range of $\mS$ from\eqref{eq:sscn_MS_def} \\
$M$ & Lipschitz constant of $\nabla^2 f(x)$ on $\R^d$; $M = M_{ \mI^d}$  \\
$L$  & Lipschitz constant of $\nabla f(x)$ on $\R^d$  \\
$\mA_{\mS}$& $\eqdef \mS^\top \mA \mS \in \R^{\tau(\mS) \times \tau(\mS)}$, for 
a given matrix $\mA \in \R^{d\times d}$   \\
$\nabla_{\mS} f(x)$& $ \eqdef \mS^\top \nabla f(x)$  \\
$\nabla^2_{\mS} f(x)$& $ \eqdef (\nabla^2 f(x))_{\mS} = \mS^\top \nabla^2 f(x) \mS$  \\
$\mH_{\mS}(x)$&  $\eqdef  \nabla^2_{\mS} f(x) +\sqrt{ \frac{M_{\mS}}{2}} \| \nabla_{\mS} f(x)\|^{\frac12} \mI^{\tau(\mS)} $ from Lemma \ref{lem:sscn_decrease} \\
$\zeta$ & $   \eqdef \lambda_{\min} \left( \left(\nabla^2 f(x^*)\right)^{\frac12}  \E{\mS (\nabla^2_{\mS} f(x^*) )^{-1} \mS^\top}  \left(\nabla^2 f(x^*)\right)^{\frac12} \right)$ from \eqref{eq:sscn_sc_generalized}\\
$\mZ$& $ \eqdef \mS \left(\mS^\top \mS\right)^{-1} \mS^\top$, the projection onto range of $\mS$ from Section \ref{sec:sscn_setup} \\
$R$ & $\eqdef\sup\limits_{x \in \R^d} \Bigl\{  \|x - x^{*} \| \; : \; F(x) \leq F(x^{0})  \Bigr\}$  from \eqref{Rdef} \\ 
 $\lambda_f(x) $ & $ \eqdef \left(\nabla f(x)^\top \left(\nabla^2 f(x)\right)^{-1}\nabla f(x) \right)^\frac12$, Newton decrement from \eqref{eq:sscn_newton_decrement} \\
 $\level $& $ \eqdef \{x ; f(x)\leq f(x^0)\}$, sublevel set  \\
\hline
\caption{Summary of frequently used notation specific to Chapter~\ref{sscn}.}\label{tbl:notation_sscn}
\end{longtable}

}

\chapter{ Appendix for Chapter \ref{acd}}
\label{acd_appendix}

\graphicspath{{ACD/experiments/images/}}

\section{Proof of Theorem~\ref{th:acd}}
Before starting the proof, we mention that the proof technique we use is inspired by the work \cite{allen2014linear,allen2016even}, which takes the advantage of the coupling of gradient descent with mirror descent, resulting in a relatively simple proof.

\subsection{Proof of inequality \eqref{eq:acd_998dgff}} \label{subsec:sigma_w<=1}

By comparing \eqref{eq:acd_sc} and \eqref{eq:acd_M-smooth-intro} for $h=e_i$, we get  $\mu_w w_i \leq \mM_{ii}$, and the first inequality in \eqref{eq:acd_998dgff} follows. Using \eqref{eq:acd_v_def} it follows that  $e_i^\top (\mP \circ \mM) e_i \preceq e_i^\top \Diag{p\circ v} e_i$, which in turn implies $\mM_{ii} \leq v_i$ and the second inequality in \eqref{eq:acd_998dgff} follows.

\subsection{Descent lemma}

The following lemma is a consequence of $\mM$-smoothness of $f$, and ESO inequality~\eqref{eq:acd_v_def}.

\begin{lemma} Under the assumptions of Theorem~\ref{th:acd}, for all $k\geq 0$ we have the bound
\begin{equation}
\label{eq:acd_eso_inq}
f(x^{k+1})-\E{f(y^{k+1})\,|\, x^{k+1}} \geq  \frac12 \| \nabla f(x^{k+1})\|^2_{v^{-1}\circ p}.
\end{equation}
\end{lemma}
\begin{proof}
We have
\begin{eqnarray*}
 && \E{f(y^{k+1}) \;|\; x^{k+1}} \\
 &\stackrel{\eqref{eq:acd_y_update}}{=}&\E{f\left(x^{k+1}-\sum_{i\in S^k} \frac{1}{v_i} \nabla_if(x^{k+1})e_i\right)  \;|\; x^{k+1}}
\\
&\stackrel{\eqref{eq:acd_M-smooth-intro}}{\leq}&
f(x^{k+1}) - \|\nabla f(x^{k+1}) \|^2_{v^{-1}\circ p}+ \frac12\E{\left\| \sum_{i\in S^k} \frac{1}{v_i} \nabla_if(x^{k+1}) e_i \right\|^2_{\mM} \;|\; x^{k+1}}
\\
&\stackrel{\eqref{eq:acd_v_def}}{\leq}&
f(x^{k+1}) - \|\nabla f(x^{k+1}) \|^2_{v^{-1}\circ p}+ \frac12 \left\| \nabla f(x^{k+1})  \right\|^2_{v^{-1}\circ p}.
\end{eqnarray*}
\end{proof}

\subsection{Key technical inequality}

We first establish a  lemma which will play a key part in the analysis.

\begin{lemma}\label{l:mirror_acd}
For every $u$ we have
\begin{eqnarray*}
&&
\eta \sum_{i\in S^k} \left\langle \frac{1}{p_i}\nabla_i f(x^{k+1}) e_i, z^{k+1}-u  \right\rangle -\frac{\eta\mu_w}{2}\|x^{k+1}-u \|^2_w \\
&&
\qquad \qquad
\leq -\frac12\| z^k-z^{k+1}\|^2_w+\frac{1}{2}\|z^k-u \|_w^2-\frac{1+\eta \mu_w}{2}\| z^{k+1}-u\|^2_{w}.
\end{eqnarray*}
\end{lemma}

\begin{proof}
The proof is a direct generalization of the proof of analogous lemma of \cite{allen2016even}. We include it for completeness.
Notice that \eqref{eq:acd_z_update} is equivalent to
\[
z^{k+1}=\argmin_z h^k(z)\eqdef \argmin_z \frac{1}{2}\| z-z^k\|^2_{w}+\eta \sum_{i\in S^k} \langle \frac{1}{p_i}\nabla_i f(x^{k+1}) e_i, z \rangle   +\frac{\eta \mu_w}{2}\|z-x^{k+1} \|_{w}^2.
\]
Therefore, we have for every $u$
\begin{eqnarray}
\nonumber
0&=&\langle \nabla h^k(z^{k+1}),z^{k+1}-u \rangle_{w} \\
&=&
\langle z^{k+1}-z^k, z^{k+1}-u\rangle_{w} +\eta \sum_{i\in S^k} \langle \frac{1}{p_i}\nabla_i f(x^{k+1}) e_i, z^{k+1}-u  \rangle  \nonumber
\\
&& \qquad \qquad  +\eta \mu_w \langle z^{k+1}-x^{k+1}, z^{k+1}-u\rangle_{w}. \label{eq:acd_zk_plus_1_optimal}
\end{eqnarray}
Next, by generalized Pythagorean theorem we have
\begin{equation}\label{eq:acd_pyt_z}
\langle z^{k+1}-z^k,z^{k+1}-u \rangle_w =\frac12 \|z^k-z^{k+1} \|^2_w-\frac12 \|z^k-u\|^2_w+\frac12 \|u-z^{k+1} \|^2_w
\end{equation}
and
\begin{equation}\label{eq:acd_pyt_x}
\langle z^{k+1}-x^{k+1},z^{k+1}-u \rangle_w =\frac12 \|x^{k+1}-z^{k+1} \|^2_w-\frac12 \|x^{k+1}-u\|^2_w+\frac12 \|u-z^{k+1} \|^2_w.
\end{equation}
It remains to put~\eqref{eq:acd_pyt_z} and~\eqref{eq:acd_pyt_x} into~\eqref{eq:acd_zk_plus_1_optimal}.
\end{proof}

\subsection{Proof of the theorem}

To mitigate notational burden, consider all expectations in this proof to be taken with respect to the choice of the random subset of coordinates $S^k$.
Using Lemma~\ref{l:mirror_acd} we have
\begin{eqnarray*}
&&
\eta \sum_{i\in S^k} \langle \frac{1}{p_i}\nabla_i f(x^{k+1}) e_i , z^{k}-u  \rangle -\frac{\eta\mu_w}{2}\|x^{k+1}-u \|^2_w 
\\
&
\leq &
\eta \sum_{i\in S^k} \langle \frac{1}{p_i}\nabla_i f(x^{k+1}) e_i, z^{k}-z^{k+1}  \rangle
 -\frac12\| z^k-z^{k+1}\|^2_w+\frac{1}{2}\|z^k-u \|_w^2 \\
 && \qquad 
 -\frac{1+\eta \mu_w}{2}\| z^{k+1}-u\|^2_{w} 
\\
&\leq &
\frac{\eta^2}{2}\| \sum_{i\in S^k} \frac{1}{p_i}   \nabla_i f(x^{k+1})e_i \|^2_{w^{-1}}+\frac{1}{2}\|z^k-u \|_w^2
-\frac{1+\eta \mu_w}{2}\| z^{k+1}-u\|^2_{w} 
\\
&= &
\frac{\eta^2}{2} \| \sum_{i\in S^k}  \nabla_i f(x^{k+1})e_i \|^2_{w^{-1}\circ p^{-2}}+\frac{1}{2}\|z^k-u \|_w^2
-\frac{1+\eta \mu_w}{2}\| z^{k+1}-u\|^2_{w}  .
\end{eqnarray*}

Taking the expectation over the choice of $S^k$, we get
\begin{eqnarray}\nonumber
&&
\eta \langle \nabla f(x^{k+1}), z^{k}-u  \rangle -\frac{\eta\mu_w}{2}\|x^{k+1}-u \|^2_w 
\\ \nonumber
& \leq &
\frac{\eta^2}{2} \|  \nabla f(x^{k+1}) \|^2_{w^{-1}\circ p^{-1}}+\frac{1}{2}\|z^k-u \|_w^2
-\frac{1+\eta \mu_w}{2}\E{\| z^{k+1}-u\|^2_{w} }
\\ \nonumber
& 
\stackrel{\eqref{eq:acd_w_def}}{=} &
\frac{\eta^2}{2} \|  \nabla f(x^{k+1}) \|^2_{v^{-1}\circ p}+\frac{1}{2}\|z^k-u \|_w^2
-\frac{1+\eta \mu_w}{2}\E{\| z^{k+1}-u\|^2_{w} }
\\
\nonumber & 
\stackrel{\eqref{eq:acd_eso_inq}}{\leq} &
\eta^2\left( f(x^{k+1})-\E{f(y^{k+1}) }\right)+\frac{1}{2}\|z^k-u \|_w^2
-\frac{1+\eta \mu_w}{2}\E{\| z^{k+1}-u\|^2_{w} }. \label{eq:acd_acd_proof_almost}
\end{eqnarray}

Next, we have the following bounds
\begin{eqnarray*}
\eta \left( f(x^{k+1})- f(x^*)\right)
& \stackrel{\eqref{eq:acd_sc}}{\leq} &
\eta \langle \nabla f(x^{k+1}),x^{k+1}-x^*\rangle -\frac{\eta \mu_w}{2} \|x^*-x^{k+1} \|^2_w
\\
& = & \eta \langle \nabla f(x^{k+1}),x^{k+1}-z^{k}\rangle +\eta \langle \nabla f(x^{k+1}),z^{k}-x^*\rangle \\
&& \qquad  -\frac{\eta \mu_w}{2} \|x^*-x^{k+1} \|^2_w
\\
&  \stackrel{\eqref{eq:acd_x_update_acd}}{=} &
\frac{(1-\theta)\eta}{\theta} \langle \nabla f(x^{k+1}),y^k-x^{k+1}\rangle +\eta \langle \nabla f(x^{k+1}),z^{k}-x^*\rangle \\
&& \qquad  -\frac{\eta \mu_w}{2} \|x^*-x^{k+1} \|^2_w
\\
& \stackrel{\eqref{eq:acd_acd_proof_almost}}{\leq} &
\frac{(1-\theta)\eta}{\theta}\left( f(y^k)-f(x^{k+1})\right)+\eta^2\left(f(x^{k+1}) -\E{f(y^{k+1})}\right)
\\
&& \qquad +
\frac{1}{2}\|z^k-x^* \|_w^2
-\frac{1+\eta \mu_w}{2}\E{\| z^{k+1}-x^*\|^2_{w} }.
\end{eqnarray*}
Choosing $\eta=\frac1\theta$  and rearranging the above we obtain

\begin{eqnarray*}
&& \frac{1}{\theta^2}\left(\E{f(y^{k+1})}-f(x^*)\right)+\frac{1+\frac{\mu_w}{\theta}}{2}\E{\| z^{k+1}-x^*\|^2_{w} } \\
&& \qquad \qquad \leq
\frac{(1-\theta)}{\theta^2}\left( f(y^k)-f(x^*)\right)+
\frac{1}{2}\|z^k-x^* \|_w^2
\end{eqnarray*}

Finally, setting $\theta$ such that $1+\frac{\mu_w}{\theta}=\frac{1}{1-\theta}$, which coincides with~\eqref{eq:acd_tau_def_acd}, we get
\[
\E{P^{k+1}}
\leq
(1-\theta)P^k,
\]
as desired.

\section{Better rates for minibatch \texttt{CD} (without acceleration) } \label{sec:acd_cd_imp}

In this section we establish better rates for minibatch \texttt{CD} method than the current state of the art. Our starting point is the following complexity theorem.

\begin{theorem} Choose any proper sampling and let $\mP$ be its probability matrix and $p$ its probability vector. Let \[c(S,\mM) \eqdef \lambda_{\max}(\mP''\circ \mM),\] where $\mP'' \eqdef \mD^{-1} \mP \mD^{-1}$ and $\mD \eqdef \Diag{p}$. Then the vector $v$ defined by $v_i = c(S,\mM) p_i$ satisfies  the ESO inequality \eqref{eq:acd_v_def}. Moreover, if we run the non-accelerated \texttt{CD} method \eqref{eq:acd_Parallel-CD-intro}  with this sampling and stepsizes $\alpha_i = \frac{1}{c(S,\mM) p_i}$, then the iteration complexity of the method is
\begin{equation}\label{eq:acd_bis798f98gud}\frac{c(S,\mM)}{\mu} \log \frac{1}{\epsilon}.\end{equation}
\end{theorem}
\begin{proof} Let $v_i = c p_i$ for all $i$. The ESO inequality holds for this choice of $v$ if $\mP \circ \mM \preceq c \mD^2$. This is equivalent to
Since $\mD^{-1}(\mP \circ \mM)\mD^{-1} = \mP'' \circ \mM$, the above inequality is equivalent to $\mP'' \circ \mM \preceq c \mI$, which is equivalent to $c\geq \lambda_{\max}(\mP'' \circ \mM)$. So, choosing $c=c(S,\mM)$ works.  Plugging this choice of $v$ into the complexity result \eqref{eq:acd_NSync} gives \eqref{eq:acd_bis798f98gud}.
\end{proof}

\subsection{Two uniform samplings and one new importance sampling}

In the next theorem we compute now consider several special samplings. All of them choose in expectation a minibatch of size $\tau$ and are hence directly comparable.

\begin{theorem} \label{thm:acd_809h0s9s} The following statements hold:
\begin{itemize}
\item[(i)] Let $S_1$ be the $\tau$--nice sampling. Then 
\begin{equation}\label{eq:acd_c1}c_1\eqdef c(S_1,\mM) = \frac{d}{\tau} \lambda_{\max}\left( \frac{\tau-1}{d-1}\mM + \frac{d-\tau}{d-1} \Diag{\mM}\right).\end{equation}
\item[(ii)] Let $S_2$ be the independent uniform sampling with minibatch size $\tau$. That is, for all  $i$ we independently decide whether $i\in S$, and do so by picking $i$ with probability $p_i = \frac{\tau}{d}$. Then
\begin{equation}\label{eq:acd_c2}c_2\eqdef c(S_2,\mM) = \lambda_{\max}\left(\mM + \frac{d-\tau}{\tau}\Diag{\mM}\right).\end{equation}
\item[(iii)] Let $S_3$ be an independent sampling where we choose $p_i \propto \frac{\mM_{ii}}{\delta + \mM_{ii}}$ where $\delta>0$ is chosen so that $\sum_i p_i = \tau$.  Then
\begin{equation}\label{eq:acd_c3} c_3\eqdef c(S_3,\mM) = \lambda_{\max}(\mM) + \delta.\end{equation}
Moreover,
\begin{equation}\label{eq:acd_bi7gdv98dg} \delta \leq \frac{\trace{\mM}}{\tau}.\end{equation}
\end{itemize}
\end{theorem}
\begin{proof}
We will deal with each case separately:
\begin{enumerate}
\item[(i)] The probability matrix of $S_1$ is $\mP = \frac{\tau}{d}\left(\beta \mE + (1-\beta) \mI  \right),$ where $\beta = \frac{\tau-1}{d-1}$, and $\mD = \frac{\tau}{d}\mI$. Hence,
\begin{eqnarray*}\mP'' \circ \mM &=& (\mD^{-1} \mP \mD^{-1})\circ \mM  \\
&=& \frac{\tau}{d}\left( \beta \mD^{-1}\mE \mD^{-1} + (1-\beta) \mD^{-2} \right) \circ \mM\\
&=& \frac{\tau}{d}\left( \frac{\tau-1}{d-1}\mE + \frac{d-\tau}{d-1} \mI \right) \circ \mM \\
&=& \frac{\tau}{d}\left( \frac{\tau-1}{d-1}\mM + \frac{d-\tau}{d-1} \Diag{\mM}\right).
\end{eqnarray*}
\item[(ii)] The probability matrix of $S_2$ is $\mP = \frac{\tau}{d}\left(\frac{\tau}{d}\mE + (1-\frac{\tau}{d}) \mI  \right)$, and $\mD = \frac{\tau}{d}\mI$. Hence,
\begin{eqnarray*}\mP'' \circ \mM &=& (\mD^{-1} \mP \mD^{-1})\circ \mM  \\
&=&  \frac{\tau}{d}\left(\frac{\tau}{d}\mD^{-1}\mE \mD^{-1}+ \left(1-\frac{\tau}{d}\right) \mD^{-2}  \right) \circ \mM\\
&=& \left( \mE + \frac{d-\tau}{\tau} \mI \right) \circ \mM \\
&=& \mM + \frac{d-\tau}{\tau} \Diag{\mM}.
\end{eqnarray*}
\item[(iii)] The probability matrix of $S_3$ is $\mP = pp^\top + \mD - \mD^2$. Therefore,
\begin{eqnarray*}\mP'' \circ \mM &=& (\mD^{-1} \mP \mD^{-1})\circ \mM  \\
&=&   \left(\mD^{-1} pp^\top \mD^{-1}+ \mD^{-1} - \mI \right) \circ \mM\\
&=&   \left(\mE + \mD^{-1} - \mI \right) \circ \mM\\
&=&   \left(\mE + \delta(\Diag{\mM})^{-1} \right) \circ \mM\\
&=& \mM+ \delta \mI.
\end{eqnarray*}
To establish the bound on $\delta$, it suffices to note that
\[\tau = \sum_i p_i = \sum_{i} \frac{\mM_{ii}}{\delta + \mM_{ii}} \leq \sum_{i} \frac{\mM_{ii}}{\delta } = \frac{\trace{\mM}}{\delta}.\]
\end{enumerate}
\end{proof}

\subsection{Comparing the samplings}

In the next result we show that sampling $S_3$ is at most twice worse than $S_2$, which is at most twice worse than $S_1$. Note that $S_1$ is uniform; and it is the standard minibatch sampling used in the literature and applications. Our novel sampling $S_3$ is {\em non-uniform}, and is at most four times worse than $S_1$ in the worst case. However, {\em it can be substantially better}, as we shall show later by giving an example.

\begin{theorem} \label{thm:acd_nonacc_comp}The leading complexity terms $c_1,c_2$, and $c_3$ of \texttt{CD} (Algorithm~\eqref{eq:acd_Parallel-CD-intro}) with samplings $S_1, S_2$, and $S_3$, respectively, defined in Theorem~\ref{thm:acd_809h0s9s}, compare as follows:
\begin{itemize}
    \item[(i)] $c_3 \leq \frac{2d-\tau}{d-\tau} c_2$
    \item[(ii)] $c_2 \leq \frac{(d-1)\tau}{d(\tau-1)} c_1 \leq 2 c_1$
\end{itemize}
\end{theorem}
\begin{proof}
We have:
\begin{enumerate}
\item[(i)] 
\begin{eqnarray*}
c_3 &\overset{\eqref{eq:acd_c3}}{=}&\lambda_{\max}(\mM) + \delta \\
&\leq & \lambda_{\max}\left( \mM + \frac{d-\tau}{\tau}\Diag{\mM}\right) + \delta \\
&\overset{\eqref{eq:acd_c2}}{=}& c_2 + \delta \\
&\overset{\eqref{eq:acd_bi7gdv98dg}}{\leq }&  c_2 + \frac{\trace{\mM}}{\tau}\\
&\leq & c_2 + \frac{d \max_i \mM_{ii}}{\tau} \\
&= & c_2 + \frac{d}{d-\tau}\frac{d-\tau}{\tau} \max_i \mM_{ii} \\
&= & c_2 + \frac{d}{d-\tau}\lambda_{\max}\left(\frac{d-\tau}{\tau} \Diag{\mM}\right)\\
&\leq & c_2 + \frac{d}{d-\tau}\lambda_{\max}\left(\mM+\frac{d-\tau}{\tau} \Diag{\mM}\right)\\
&\overset{\eqref{eq:acd_c2}}{=}&  \frac{2d-\tau}{d-\tau} c_2.
\end{eqnarray*}
\item[(ii)] 
\begin{eqnarray*} 
c_2 &\overset{\eqref{eq:acd_c2}}{=}& \lambda_{\max}\left(\mM + \frac{d-\tau}{\tau}\Diag{\mM}\right) \\
&=&  \lambda_{\max}\left( \frac{d(\tau-1)}{\tau(d -1)} \mM + \frac{d-\tau}{\tau}\Diag{\mM} + \left(1-\frac{d(\tau-1)}{\tau(d -1)}\right)\mM\right)\\
&\overset{(\dagger)}{\leq} &  \lambda_{\max}\left( \frac{d(\tau-1)}{\tau(d -1)} \mM + \frac{d-\tau}{\tau}\Diag{\mM}\right) + \lambda_{\max}\left(\left(1-\frac{d(\tau-1)}{\tau(d -1)}\right)\mM\right)\\
&\leq & \lambda_{\max}\left( \frac{d(\tau-1)}{\tau(d -1)} \mM + \frac{d(d -\tau)}{\tau(d -1)}\Diag{\mM}\right) +\frac{d-\tau}{(d -1)\tau} \lambda_{\max}\left(\mM\right)\\
&\overset{\eqref{eq:acd_c1}}{=}& c_1 + \frac{d-\tau}{(d -1)\tau} \lambda_{\max}\left(\mM\right)\\
&\overset{\eqref{eq:acd_c2}}{\leq} & c_1 +  \frac{d-\tau}{(d -1)\tau} c_2.
\end{eqnarray*}
The statement follows by reshuffling the final inequality. In  step $(\dagger)$ we have used subadditivity of the function $\mA \mapsto \lambda_{\max}(\mA)$.
\end{enumerate}
\end{proof}

The next simple example shows that sampling $S_3$ can be arbitrarily better than sampling $S_1$.

\begin{example}\label{ex:nonacc_imp_diff}
Consider $d\gg 1$, and choose any $\tau$ and 
\[
\mM\eqdef\begin{pmatrix}
n&0^\top \\ 0 & \mI 
\end{pmatrix}
\]
for $\mI\in \R^{(d -1)\times (d -1)}$. Then, it is easy to verify that $c_1\stackrel{\eqref{eq:acd_c1}}{=}\frac{d^2}{\tau}$ and $c_3\stackrel{\eqref{eq:acd_c3}+\eqref{eq:acd_bi7gdv98dg}}{\leq} n+\frac{2n-1}{\tau}=\cO(\frac{d}{\tau})$. Thus, convergence rate of \texttt{CD} with $S_3$ sampling can be up to $\cO(d )$ times better than convergence rate of \texttt{CD} with $\tau$--nice sampling. 
\end{example}

\begin{remark}\label{rem:proport}
Looking only at diagonal elments of $\mM$, an intuition tells us that one should sample a coordinate corresponding to larger diagonal entry of $\mM$ with higher probability. However, this might lead to worse convergence, comparing to $\tau$--nice sampling. Therefore the results we provide in this section cannot be qualitatively better, i.e. there are examples of smoothness matrix, for which assigning bigger probability to bigger diagonal elements leads to worse rate. It is an easy exercise to verify that for $\mM\in \R^{10\times 10}$ such that
\[
\mM\eqdef \begin{pmatrix}
2 & 0^\top \\ 0& 11^\top
\end{pmatrix}, 
\] 
and $\tau\geq 2$ we have $c(S_{\text{nice}},\mM)\leq c(S',\mM)$ for any $S'$ satisfying $p(S')_i\geq p(S')_j$ if and only if $\mM_{ii}\geq \mM_{jj}$. 
\end{remark}

\section{Proofs for Section~\ref{sec:acd_import}}
\subsection{Proof of Theorem~\ref{thm:acd_LB}} \label{app:inequality_lower_bound}

We start with a lemma which allows us to focus on ESO parameters $v_i$ which are proportional to the squares of the probabilities $p_i$.

\begin{lemma} \label{lem:acd_b7f8gf8} Assume that the ESO inequality  \eqref{eq:acd_v_def}  holds. Let $j=\arg\max_i \frac{v_i}{p_i^2}$, $c = \frac{v_j}{p_j^2}$ and $v' = c p^2$ (i.e.,  $v'_i = c p_i^2$ for all $i$). Then the following statements hold:
\begin{enumerate}
\item[(i)] $v'\geq v$.
\item[(ii)] ESO inequality \eqref{eq:acd_v_def}  holds for $v'$ also.
\item[(iii)]  Assuming $f$ is $\mu$--convex, Theorem~\ref{th:acd} holds if we replace $v$ by $v'$, and the rate \eqref{eq:acd_ineq-minibatch-speedup} is unchanged if we replace $v$ by $v'$. 
\end{enumerate}
\end{lemma}

\begin{proof}
\begin{enumerate}
\item[(i)] $v'_i = c p_i^2 = \frac{v_j}{p_j^2} p_i^2 = \left( \frac{v_j}{p_j^2} \frac{p_i^2}{v_i}\right) v_i \geq v_i$.
\item[(ii)] This follows directly from (i).
\item[(iii)]  Theorem~\ref{th:acd} holds with $v$ replaced by $v'$ because ESO holds. To show that the rates are unchanged first note that
$\max_i \frac{v_i}{p_i^2} = \frac{v_j}{p_j^2} = c$. On the other hand, by construction, we have $c= \frac{v'_i}{p_i^2}$ for all $i$. So, in particular, $c = \max_i \frac{v'_i}{p_i^2}$.
\end{enumerate}
\end{proof}

In view of the above lemma, we can assume without loss of generality that $v=cp^2$. Hence, the rate in \eqref{eq:acd_ineq-minibatch-speedup}  can be written in the form
\begin{equation}\label{eq:acd_b87gs98scc}\sqrt{\max_i \frac{v_i}{p_i^2 \mu}}  = \sqrt{\frac{c}{\mu}}.\end{equation}

In what follows, we will establish a lower  bound on $c$, which will  lead to the lower bound on the rate expressed as inequality~\eqref{eq:acd_ineq-minibatch-speedup}. As a starting point, note that directly from \eqref{eq:acd_v_def} we get the bound
\begin{equation}\label{eq:acd_nb87d90sh}\mP \circ \mM \preceq  \Diag{p\circ v} = c \Diag{p^3}.\end{equation}
Let $\mD_1=\Diag{p}^{-1/2}$ and $\mD_2 = \Diag{p}^{-1}$. From \eqref{eq:acd_nb87d90sh} we get
$ \mD_1 \mD_2 (\mP\circ \mM) \mD_2 \mD_1 \preceq  c \mI$ and hence
\begin{equation} \label{eq:acd_bu987g98f} c \geq c(S,\mM) \eqdef \lambda_{\max}(\mD_1 \mD_2 (\mP\circ \mM) \mD_2 \mD_1).\end{equation}

At this point, the following identity will be useful.

\begin{lemma}\label{lem:acd_bgfd78gd8} Let $\mA,\mB,\mD_1,\mD_2\in \R^{d\times d}$, with $\mD_1,\mD_2$ being diagonal. Then
\begin{equation} \label{eq:acd_hadamard_diag}
\mD_1(\mA \circ \mB) \mD_2 = (\mD_1\mA \mD_2) \circ \mB =   \mA \circ(\mD_1\mB \mD_2).
\end{equation}
\end{lemma}
\begin{proof}
The proof is straightforward, and hence we do not include it. The identity is formulated as an exercise in \cite{FuzhenZhang}.
\end{proof}

Repeatedly applying Lemma~\ref{lem:acd_bgfd78gd8}, we get
\[\mD_1 \mD_2 (\mP\circ \mM) \mD_2 \mD_1 = \underbrace{(\mD_1 \mP \mD_1)}_{\mP'} \circ \underbrace{(\mD_2 \mM \mD_2)}_{\mM'}.\]
Plugging this back into \eqref{eq:acd_bu987g98f},  and since $\mP'_{ii}=1$ for all $i$, we get the bound
\begin{eqnarray} \label{eq:acd_ bi9s8g9snnn} c &\geq& c(S,\mM) =  \lambda_{\max}(\mP' \circ \mM') \geq  \max_i \; (\mP'\circ \mM')_{ii} =  \max_i \mP'_{ii} \mM'_{ii}  = \max_i \mM'_{ii} \notag \\
&=& \max_i \frac{\mM_{ii}}{p_i^2}   \geq   \frac{\left(\sum_{i=1}^d \mM_{ii}^{1/2}\right)^2}{\tau^2}.\label{eq:acd_bdg8db8d}\end{eqnarray}
The last inequality follows by observing that the optimal solution of the optimization problem
\[\min_p \left\{\max_i \frac{\mM_{ii}}{p_i^2} \;|\; p_1,\dots,p_d > 0, \; \sum_i p_i =\tau\right\}\]
is $p_i =\tau \frac{\mM_{ii}^{1/2}}{\sum_j \mM_{jj}^{1/2}} $. Inequality~\eqref{eq:acd_ineq-minibatch-speedup} now follows by substituting  
the lower bound on $c$ obtained in \eqref{eq:acd_bdg8db8d} into 
\eqref{eq:acd_b87gs98scc}.

\subsection{Proof of Lemma~\ref{thm:acd_special-ESO-result}}
\begin{eqnarray*}
 \Diag{p_1 v_1, \dots, p_d v_d} &=& c(S,\mM)\Diag{p_1^3,\dots ,p_d^3}\\
 &=&
c(S,\mM) \mD^{3} \\
 &=&
\lambda_{\max}\left(\left( \mD^{-1/2} \mP \mD^{-1/2}\right) \circ \left( \mD^{-1}\mM \mD^{-1}\right)\right)\mD^{3}
\\
& \succeq & 
\mD^{\frac32}\left(\left( \mD^{-1/2} \mP \mD^{-1/2}\right) \circ \left( \mD^{-1}\mM \mD^{-1}\right)\right)\mD^{\frac32}
\\
&\stackrel{\eqref{eq:acd_hadamard_diag}}{=}&
\mP\circ \mM.
\end{eqnarray*}
The last inequality came from the fact that $\mD$ is diagonal. 

\subsection{Bound on $c(S_1,\mM)$}

\begin{lemma}\label{lem:acd_tau-nice-2nd-derivation}
$c(S_1,\mM) \leq \frac{d^2}{\tau^2} ((1-\beta)\max_i \mM_{ii} + \beta L).$
\end{lemma}
\begin{proof} Recall that the probability matrix of $S_1$ is $\mP = \frac{\tau}{d}\left((1-\beta) \mI + \beta \mE\right)$. Since $p_i=\frac{\tau}{d}$ and $\mM \preceq L \mI$, we have
\begin{eqnarray*}
c(S_1,\mM) &=& \lambda_{\max}\left(\mP' \circ \mM' \right) \\
&=& \lambda_{\max}\left( (\mD^{-1/2}\mP \mD^{-1/2}) \circ (\mD^{-1}\mM \mD^{-1}) \right) \\ 
&=& \lambda_{\max} \left( \frac{\tau}{d}\left((1-\beta) \mD^{-1} + \beta \mD^{-1/2}\mE \mD^{-1/2}\right) \circ \mD^{-1} \mM \mD^{-1} \right)\\
&=&\frac{\tau}{d} \lambda_{\max} \left( \left((1-\beta) \mD^{-1} + \beta \mD^{-1/2}\mE \mD^{-1/2}\right) \circ \mD^{-1} \mM \mD^{-1} \right)\\
&=& \frac{\tau}{d}  \lambda_{\max}  \left((1-\beta) \Diag{\mM_{ii}/p_i^3} + \beta \mD^{-3/2}\mM \mD^{-3/2}\right) \\
&\preceq  & \frac{\tau}{d}  \lambda_{\max}  \left((1-\beta) \Diag{\mM_{ii}/p_i^3} + \beta L \mD^{-3} \right) \\
&=& \frac{\tau}{d}  \lambda_{\max}  \left((1-\beta) \frac{d^3}{\tau^3}\max_i \mM_{ii}+ \beta L \frac{d^3}{\tau^3} \right) \\
&=& \frac{d^2}{\tau^2} \left((1-\beta) \max_i \mM_{ii} + \beta L \right).
\end{eqnarray*}
\end{proof}

\subsection{Proof of Theorem~\ref{thm:acd_comparison}}
For the purpose of this proof, let $S_2$ be the independent uniform sampling with minibatch size $\tau$. That is, for all  $i$ we independently decide whether $i\in S$, and do so by picking $i$ with probability $p_i = \frac{\tau}{d}$. Recall that $S_3$ is the independent importance sampling.

For simplicity, let  $\mP_i$ be the probability matrix of sampling $S_i$,  $\mD_i\eqdef \Diag{\mP_i}$, and $\mM'_i\eqdef\mD^{-1/2}_i\mM \mD^{-1/2}_i$, for $i=1,3$. Next, we have
\begin{eqnarray}
c(S_i,\mM)&=& \lambda_{\max}\left(\left( \mD_i^{-1/2} \mP_i \mD_i^{-1/2}\right) \circ \left( \mD_i^{-1}\mM \mD_i^{-1}\right)\right)
\nonumber
\\
&\stackrel{\eqref{eq:acd_hadamard_diag}}{=}& \lambda_{\max}\left(\left( \mD_i^{-1} \mP_i \mD_i^{-1}\right) \circ \left( \mD_i^{-1/2}\mM \mD_i^{-1/2}\right)\right)
\nonumber
\\
\nonumber
&=& \lambda_{\max}\left(\left( \mE + \mD_i^{-1}-\mI  \right) \circ \mM_i'\right)
\\
&=& \lambda_{\max}\left(\mM_i' + \Diag{\mM_i'}\circ(\mD_i^{-1}-\mI ) \right),\label{eq:acd_cs_simple}
\end{eqnarray}
where the third identity holds since both $S_i$ is an independent sampling, which means that $\left( \mD_i^{-1} \mP_i \mD_i^{-1}\right)_{kl} = \frac{p_{kl}}{p_{k}p_l}$, where $p=\Diag{\mD_i}$.

Denote $c_i\eqdef c(S_i,\mM)$. Thus for $S_2$ we have
\begin{equation}
c_2 =\frac{\tau}{d}\lambda_{\max}\left( \mM +\frac{d-\tau}{\tau} \Diag{\mM} \right).\label{eq:acd_c2_acc}
\end{equation}

Let us now establish a technical lemma.

\begin{lemma}\label{lem:acd_import_main}
\begin{equation}\label{eq:acd_p_choice_ineq}
 \lambda_{\max}\left(\mM'_3 +\diag(\mM'_3)\circ(\mD^{-1}_3-\mI ) \right)
\leq
 \frac{2d-\tau}{d-\tau}\lambda_{\max}\left(\mM'_3 +\frac{d-\tau}{\tau} \Diag{\mM'_3}\right) .
\end{equation}
\end{lemma}
\begin{proof}
The statement follows immediately repeating the steps of the proof of (i) from Theorem~\ref{thm:acd_nonacc_comp} using the fact that for sampling $S_3$ we have $p_i/\mM_{ii}\propto p_i^{-1}-1$. 
\end{proof}

We can now proceed with comparing $c_2$ to $c_3$.
\begin{eqnarray}
c_3 &=&  \lambda_{\max}\left(\mM'_3 + \Diag{\mM'_3}\circ(\mD^{-1}_3-\mI ) \right) \nonumber \\
&\stackrel{\eqref{eq:acd_p_choice_ineq}}{\leq}  &
 \frac{2d-\tau}{d-\tau}\lambda_{\max}\left(\mM'_3 +\frac{d-\tau}{\tau} \Diag{\mM'_3}\right) \nonumber \\
 &\stackrel{(*)}{\leq}&
\frac{2d-\tau}{d-\tau}
\lambda_{\max}\left(d \Diag{\mM'_3} +\frac{d-\tau}{\tau} \Diag{\mM'_3}\right) \nonumber\\
&=&
\frac{2d-\tau}{d-\tau} \frac{d\tau+d-\tau}{\tau}\lambda_{\max}\left( \Diag{\mM'_3}\right) \nonumber\\
&\stackrel{(**)}{\leq}&
\frac{2d-\tau}{d-\tau} \frac{d\tau+d-\tau}{\tau} \frac{\tau}{d}\lambda_{\max}\left( \Diag{\mM}\right) \nonumber\\
&\leq&
\frac{2d-\tau}{d-\tau} \frac{d\tau+d-\tau}{\tau} \frac{\tau}{d}\frac{\tau}{d-\tau}\lambda_{\max}\left(\mM+\frac{d-\tau}{\tau} \Diag{\mM}\right)  \nonumber \\
&\stackrel{\eqref{eq:acd_cs_simple}}{=}&
\frac{2d-\tau}{d-\tau} \frac{d\tau+d-\tau}{\tau} \frac{\tau}{d-\tau} c_2 . \label{eq:acd_c3c2cmp}
\end{eqnarray}

Above, inequality  $(*)$ holds since for any $d\times d$ matrix $\mQ\succ 0$ we have $\mQ\preceq d\Diag{\mQ}$ and inequality $(**)$ holds since $(\mD_3)_{ii}\geq(\mD_3)_{jj} $ if and only if $\mM_{ii}\geq \mM_{jj}$ due to choice of $p$.

Let us now compare to $c_2$ and $c_1$. We have
\begin{eqnarray}
c_1&=& \lambda_{\max}\left(\left( \mD_1^{-1/2} \mP_1 \mD_1^{-1/2}\right) \circ \left( \mD_1^{-1}\mM \mD_1^{-1}\right)\right)
\nonumber
\\
&\stackrel{\eqref{eq:acd_hadamard_diag}}{=}& \lambda_{\max}\left(\left( \mD_1^{-1} \mP_1 \mD_1^{-1}\right) \circ \left( \mD_1^{-1/2}\mM \mD_1^{-1/2}\right)\right)
\nonumber
\\
\nonumber
&=& \lambda_{\max}\left(\left( \frac{\tau-1}{d-1}\frac{d}{\tau} \mE + \frac{d}{\tau} \mI-\frac{\tau-1}{d-1}\frac{d}{\tau}\mI  \right) \circ \mM''_1\right)
\\
\nonumber
&=& \frac{d}{\tau}\lambda_{\max}\left(\frac{\tau-1}{d-1}  \mM''_1 + \frac{d-\tau}{d-1} \Diag{ \mM''_1}  \right)
\\ \label{eq:acd_c1_acc}
&=& \left(\frac{d}{\tau}\right)^2\lambda_{\max}\left(\frac{\tau-1}{d-1}  \mM + \frac{d-\tau}{d-1} \Diag{\mM}  \right).
\end{eqnarray}
As~\eqref{eq:acd_c2_acc} and~\eqref{eq:acd_c1_acc} are established, following the proof of (ii) from Theorem~\ref{thm:acd_nonacc_comp}, we arrive at
\begin{equation}\label{eq:acd_c1c2cmp}
c_2 \leq \frac{(d -1)\tau}{d(\tau-1)} c_1 \leq 2 c_1.
\end{equation}
It remains to combine~\eqref{eq:acd_c3c2cmp} and~\eqref{eq:acd_c1c2cmp} to establish~\eqref{eq:acd_thm_comparison}. 

An example with $c_3\approx \left(\frac{\tau}{d}\right)^2 c_2$ follows.
\begin{example}\label{ex:acc_imp_diff}
Consider $d\geq 1$, choose any $d\geq \tau \geq 1 $ and 
\[
\mM\eqdef\begin{pmatrix}
N &0^\top \\ 0 & \mI 
\end{pmatrix}
\]
for $\mI\in \R^{(d -1)\times (d -1)}$. Then, it is easy to verify that $c_1\stackrel{\eqref{eq:acd_c1_acc}}{=}\left(\frac{d}{\tau}\right)^2N$. 
Moreover, for large enough $N$ we have 
\[
p\approx \left(1,\frac{\tau-1}{d-1},\dots, \frac{\tau-1}{d-1}\right)^\top \qquad \Rightarrow \qquad
\mM'_3 \approx  \Diag{ N, \frac{d-1}{\tau-1},\dots,\frac{d-1}{\tau-1} }. 
\]
Therefore, using~\eqref{eq:acd_cs_simple} and again for large enough $N$, we get $c_3\approx N$. Thus, $c_3\approx \left(\frac{\tau}{d}\right)^2 c_2$. 
\end{example}

\chapter{  Appendix for Chapter \ref{sega}}
\label{sega_appendix}

\graphicspath{{SEGA/experiments/}}

\section{Proofs for Section~\ref{sec:sega_analysis}}
\begin{lemma}\label{lem:sega_relate}
Suppose that $f$ is twice differentiable. 
Assumption~\ref{ass:sega_M_smooth_inv} is equivalent to~\eqref{eq:acd_M-smooth-intro} for $\mmM=\mM^{-1}$. 
\end{lemma}
\begin{proof}
We first establish that Assumption~\ref{ass:sega_M_smooth_inv} implies~\eqref{eq:acd_M-smooth-intro}. 
Summing up~\eqref{eq:sega_M_smooth_inv} for $(x,y)$ and $(y,x)$ yields
\[
 \langle  \nabla f(x) - \nabla f(y) ,x - y \rangle \geq \|\nabla f(x) - \nabla f(y)\|_{\mmM}^2. 
\]
Using Cauchy Schwartz inequality we obtain
\[
 \| x-y\|_{\mQ^{-1}} \geq \|\nabla f(x) - \nabla f(y)\|_{\mmM}. 
\]
By the mean value theorem, there is $z \in [x,y]$ such that $\nabla f(x) - \nabla f(y) = \nabla^2 f(z) (x-y)$. Thus
\[
 \| x-y\|_{\mQ^{-1}} \geq \|x-y\|_{\nabla^2 f(z) \mmM \nabla^2 f(z) }. 
\]
The above is equivalent to
\[
\left(\nabla^2 f(z)\right)^{-\frac12}\mQ^{-1}\left(\nabla^2 f(z)\right)^{-\frac12} \succeq 
\left(\nabla^2 f(z)\right)^{\frac12} \mQ\left(\nabla^2 f(z)\right)^{\frac12}
\]
Note that for any $\mM'\succ 0$ we have $\mM' \succeq \mM^{-1}$ if and only if $\mM \succeq \mI$. Thus 
\[
\left(\nabla^2 f(z)\right)^{-\frac12}\mQ^{-1}\left(\nabla^2 f(z)\right)^{-\frac12}\succeq \mI,
\]
which is equivalent to $\mQ^{-1} \succeq \nabla^2 f(z)$.
To establish the other direction, denote $\phi(y)=f(y)-\langle \nabla f(x),y\rangle$. Clearly, $x$ is minimizer of $\phi$ and therefore we have
\[
\phi(x)\leq \phi(x-\mM^{-1} \nabla f(y)) \leq \phi(y)-\frac12 \|\nabla f(y) \|^2_{\mM^{-1}}, 
\]
which is exactly~\eqref{eq:sega_M_smooth_inv} for $\mmM=\mM^{-1}$.  

\end{proof}

\begin{lemma}
\label{lem:sega_zkbzk}
    For $\mZ_k \eqdef  \mS_k (\mS_k^\top \mS_k)^\dagger \mS_k^\top$, then
    \begin{align}
    \label{eq:sega_zkbzk}
        \mZ_k^\top  \mZ_k = \mZ_k.
    \end{align}
\end{lemma}
\begin{proof}
	It is a property of pseudo-inverse that for any matrices $\mA, \mI$ it holds $((\mA\mI)^\dagger)^\top = (\mI^\top\mA^\top)^\dagger$, so $\mZ_k^\top = \mZ_k$. Moreover, we also know for any $\mA$ that $\mA^\dagger \mA \mA^\dagger = \mA^\dagger$ and, thus,
\[
        \mZ_k^\top \mZ_k 
        =  \mS_k (\mS_k^\top \mS_k)^\dagger \mS_k^\top  \mS_k(\mS_k^\top  \mS_k)^\dagger \mS_k^\top= \mZ_k.
\]
\end{proof}

\subsection{Proof of Theorem~\ref{thm:sega_main}}
We first state two lemmas which will be crucial for the analysis. They characterize key properties of the gradient learning process~\eqref{eq:sega_h^{k+1}},~\eqref{eq:sega_g^k} and will be used later to bound expected distances of both $h^{k+1}$ and $g^k$ from $\nabla f(x^*)$. The proofs are provided in Appendix~\ref{sec:sega_prl1} and~\ref{sec:sega_prl2} respectively

\begin{lemma} \label{lem:sega_2} For all $v\in \R^d$ we have
\begin{equation}\label{eq:sega_h_decomp}
\E{\|h^{k+1} - v\| ^2} =  \|h^k - v\|_{\mI - \E{\mZ}}^2 + \|\nabla f(x^k) - v\|_{\E{\mZ}}^2.
\end{equation}
\end{lemma}

\begin{lemma}
\label{lem:sega_gk_general}
        Let $\mC\eqdef \E{\theta^2\mZ}$. Then for all $v\in \R^d$ we have
\[
        \E{\|g^k - v\| ^2} \le 2\|\nabla f(x^k) - v\|^2_\mC + 2\|h^k - v\|^2_{\mC-\mI}.
\]
\end{lemma}
For notational simplicity, it will be convenient to define Bregman divergence between $x$ and $y$: \[D_f(x,y)\eqdef f(x) - f(y) - \langle \nabla f(y)),x - y \rangle\]
We can now proceed with the proof of Theorem~\ref{thm:sega_main}. 
    Let us start with bounding the first term in the expression for $\Lgen^{k+1}$. From Lemma~\ref{lem:sega_gk_general} and strong convexity it follows that
    \begin{eqnarray*}
        \E{\|x^{k+1} - x^*\|^2 }
        &=& \E{\|\prox_{\alpha \psi}(x^k - \alpha g^k) - \prox_{\alpha \psi}(x^* - \alpha \nabla f(x^*))\|^2 }\\
        &\le& \E{\|x^k - \alpha g^k - (x^* - \alpha \nabla f(x^*))\|^2 }\\ 
        &=& \|x^k - x^*\|^2  - 2\alpha\E{ (g^k - \nabla f(x^*))^\top  (x^k - x^*)}\\
        &&\qquad + \alpha^2\E{\|g^k - \nabla f(x^*)\|^2 }\\
        &\le & \|x^k - x^*\|^2  - 2\alpha (\nabla f(x^k) - \nabla f(x^*))^\top  (x^k - x^*) \\
        &&\qquad + 2\alpha^2\|\nabla f(x^k) - \nabla f(x^*)\|^2_\mC + 2\alpha^2\|h^k - \nabla f(x^*)\|^2_{\mC-\mI}\\
        &\le &\|x^k - x^*\|^2  - \alpha\mu\|x^k - x^*\|^2   - 2\alpha D_f(x^k, x^*) \\
        &&\qquad + 2\alpha^2\|\nabla f(x^k) - \nabla f(x^*)\|^2_\mC + 2\alpha^2\|h^k - \nabla f(x^*)\|^2_{\mC-\mI}.
    \end{eqnarray*}
    Using Assumption~\ref{ass:sega_M_smooth_inv} we get
    \begin{eqnarray*}
        -2\alpha D_f(x^k, x^*) \le -\alpha \|\nabla f(x^k) - \nabla f(x^*)\|^2_\mmM.
    \end{eqnarray*}
    As for the second term in $\Lgen^{k+1}$, we have by Lemma~\ref{lem:sega_2}
    \begin{eqnarray*}
        \alpha\sigma\E{\|h^{k+1} - \nabla f(x^*)\|^2 } = \alpha\sigma\|h^k - \nabla f(x^*)\|^2_{\mI - \E{\mZ}}  + \alpha\sigma\|\nabla f(x^k) - \nabla f(x^*)\|^2_{\E{\mZ}}.
    \end{eqnarray*}
    
    Combining it into Lyapunov function $\Lgen^k$,
    \begin{eqnarray*}
        \Lgen^{k+1} 
        &\le & (1  - \alpha\mu)\|x^k - x^*\|^2  + \alpha\sigma\|h^k - \nabla f(x^*)\|^2_{\mI - \E{\mZ}} + 2\alpha^2\|h^k - \nabla f(x^*)\|^2_{\mC-\mI} \\
        &&  + \alpha\sigma\|\nabla f(x^k) - \nabla f(x^*)\|^2_{\E{\mZ}} + 2\alpha^2\|\nabla f(x^k) - \nabla f(x^*)\|^2_\mC -\alpha \|\nabla f(x^k) \\
        && - \nabla f(x^*)\|^2_\mmM.
    \end{eqnarray*}
    To see that this gives us the theorem's statement, consider first
    \begin{eqnarray*}
        \alpha\sigma \E{\mZ} +2\alpha^2 \mC - \alpha\mmM 
        = 2\alpha \left(\alpha\mC - \frac{1}{2}(\mmM - \sigma \E{\mZ})\right)
        \le 0,
    \end{eqnarray*}
    so we can drop norms related to $\nabla f(x^k) - \nabla f(x^*)$. Next, we have
    \begin{eqnarray*}
        \alpha\sigma (\mI - \E{\mZ}) + 2\alpha^2(\mC - \mI) 
        &=& \alpha\left( \alpha(2(\mC - \mI) + \sigma\mu\mI) - \E{\mZ}\right) + \sigma \alpha(1 - \alpha\mu)\mI \\
        &\le& \sigma\alpha (1 - \alpha\mu)\mI,
    \end{eqnarray*}
    which follows from our assumption on $\alpha$.

\subsection{Proof of Lemma~\ref{lem:sega_2} \label{sec:sega_prl1}}

\begin{proof}
Keeping in mind that $\mZ_k^\top = \mZ_k$, we first write
\begin{eqnarray*}
    \E{\| h^{k+1}-v\| ^2} 
    &\overset{\eqref{eq:sega_988fgf}}{=}& \E{ \left\|h^k + \mZ_k (\nabla f(x^k) - h^k) - v \right\| ^2 }\\
 &=&  \E{ \left\| \left(\mI - \mZ_k \right)(h^k - v) + \mZ_k (\nabla f(x^k)-v)  \right\| ^2 }\\
 &=& \E{ \left\| \left(\mI - \mZ_k \right)(h^k-v) \right\| ^2} + \E{\left\|\mZ_k  ( \nabla f(x^k) - v)  \right\| ^2 }\\
 && \qquad + 2 (h^k-v)^\top\E{\left(\mI - \mZ_k \right)^\top  \mZ_k}  (\nabla f(x^k) - v) \\
 &=& (h^k-v)^\top \E{ \left(\mI - \mZ_k \right)^\top  \left(\mI - \mZ_k \right)} (h^k-v) \\
 && \qquad + (\nabla f(x^k)-v)^\top \E{\mZ_k\mZ_k} (\nabla f(x^k) -v)\\
 && \qquad + 2 (h^k-v)^\top\E{\mZ_k - \mZ_k\mZ_k}  (\nabla f(x^k) - v).
 \end{eqnarray*}
By Lemma \ref{lem:sega_zkbzk} we have $\mZ_k\mZ_k=\mZ_k$, so the last term in the expression above is equal to 0. As for the other two, expanding the matrix factor in the first term leads to
\begin{eqnarray*}
     \E{ \left(\mI - \mZ_k \right)^\top   \left(\mI - \mZ_k \right)}
     &=& \E{ \left(\mI - \mZ_k \right)   \left(\mI - \mZ_k \right)}\\
     &=&  \E{ \mI - \mZ_k\mI - \mI\mZ_k + \mZ_k  \mZ_k}\\
     &=& \mI - \E{\mZ_k}.
\end{eqnarray*}
We, thereby, have derived
 \begin{eqnarray*}
 \E{\| h^{k+1}-v\| ^2} 
 &=&(h^k-v)^\top \left(\mI - \E{\mZ_k}\right) (h^k - v) \\
 && \quad + (\nabla f(x^k)-v)^\top \E{\mZ_k\mZ_k} (\nabla f(x^k) - v)\\
 &=& \|h^k - v\|_{\mI - \E{\mZ}}^2 + \|\nabla f(x^k) - v\|_{\E{\mZ}}^2.
\end{eqnarray*}

\end{proof}

\subsection{Proof of Lemma~\ref{lem:sega_gk_general}\label{sec:sega_prl2}}
\begin{proof}
    Throughout this proof, we will use without any mention that $\mZ_k^\top = \mZ_k$.
    
    Writing $g^k - v = a + b$, where $a\eqdef (\mI - \theta_k\mZ_k)(h^k - v)$ and $b\eqdef \theta_k \mZ_k (\nabla f(x^k) - v)$, we get $\|g^k\| ^2\le 2(\|a\| ^2 + \|b\| ^2)$. Using Lemma~\ref{lem:sega_zkbzk} and the definition of $\theta_k$ yields
    \begin{eqnarray*}
        \E{\|a\| ^2} 
        &=& \E{\|\left(\mI - \theta_k\mZ_k\right)(h^k - v)\| ^2}\\
        &=& (h^k - v)^\top \E{\left( \mI - \theta_k\mZ_k\right)\left( \mI - \theta_k\mZ_k\right)} (h^k - v)\\
        &=&  (h^k - v)^\top \E{\left( \mI - \theta_k\mZ_k - \theta_k\mZ_k + \theta_k^2\mZ_k\mI\mZ_k \right)} (h^k - v)\\
        &=&(h^k - v)^\top \E{\left( \mI - 2\mI + \theta_k^2\mZ_k\right)} (h^k - v)\\
        &=& \|h^k - v\|^2_{\E{\theta^2\mZ} - \mI}.
    \end{eqnarray*}
    Similarly, the second term in the upper bound on $g^k$ can be rewritten as
    \begin{eqnarray*}
        \E{\|b\| ^2}
        &=& \E{\|\theta_k \mZ_k (\nabla f(x^k) - v)\|^2 }\\
        &=& (\nabla f(x^k) - v)^\top \E{\theta_k^2\mZ_k \mZ_k} (\nabla f(x^k) - v)\\
        &=& \|\nabla f(x^k) - v\|^2_{\mC}.
    \end{eqnarray*}
    Combining the pieces, we get the claim. 
\end{proof}

\section{Proofs for Section~\ref{sec:sega_CD} \label{sec:sega_proofs_CD}}
\subsection{Technical lemmas\label{sec:sega_techlem}}
We first start with an analogue of Lemma~\ref{lem:sega_gk_general} allowing for a norm different from $\|\cdot\| $. We remark that 
matrix $\mQ'$ in the lemma is  not to be confused with the smoothness matrix $\mQ$ from Assumption~\ref{ass:sega_M_smooth_inv}.

\begin{lemma} \label{lem:sega_1xxxx} Let $\mQ'\succ 0$. The variance of $g^k$ as an estimator of $\nabla f(x^k)$ can be bounded as follows:
\begin{equation} \label{eq:sega_almost_eso_g}
\frac12\E{\|g^k\|_{\mQ'}^2} \leq  \|h^k \|_{\mPdiag^{-1} (\Probmat \circ \mQ')\mPdiag^{-1}- \mQ'}^2 +  \|\nabla f(x^k)\|^2_{\mPdiag^{-1} (\Probmat \circ \mQ')\mPdiag^{-1}}.
\end{equation}
\end{lemma}

\begin{proof}  
Denote $\mS_{k}$ to be a matrix with columns $e_i$ for $i\in \Range{\mS_{k}}$.  
We first write \[g^k = \underbrace{h^k - \mPdiag^{-1}  \mS_{k} \mS_{k}^\top h^k}_{a} +  \underbrace{ \mPdiag^{-1}    \mS_{k}\mS_{k}^\top \nabla f(x^k)}_{b}.\]
Let us bound the expectation of each term individually. The first term is equal to
\begin{eqnarray*}
\E{\|a\|_{\mQ'}^2} &=&  
\E{\left\| \left(\mI -  \mPdiag^{-1}  \mS_{k} \mS_{k}^\top\right) h^k  \right\|_{\mQ'}^2}
\\
&=&(h^k)^\top \E{ \left(\mI -\mPdiag^{-1}  \mS_{k} \mS_{k}^\top \right)^\top \mQ' \left(\mI - \mPdiag^{-1}  \mS_{k} \mS_{k}^\top \right)} h^k
\\
&=&
(h^k)^\top  \E{\left(\mQ' -\mPdiag^{-1}  \mS_{k} \mS_{k}^\top \mQ'- \mQ'  \mS_{k} \mS_{k}^\top\mPdiag^{-1}\right)}h^k  
\\
&&
\qquad \qquad  +(h^k)^\top  \E{ \left(\mPdiag^{-1}  \mS_{k} \mS_{k}^\top \mQ'  \mS_{k} \mS_{k}^\top \mPdiag^{-1} \right)}h^k 
\\
&=& 
(h^k)^\top \left( \mPdiag^{-1} (\Probmat \circ \mQ')\mPdiag^{-1}- \mQ'  \right)h^k.
\end{eqnarray*}
The second term can be bounded as
\begin{eqnarray*}
\E{\|b\|_{\mQ'}^2} &=&  \E{\left\|\mPdiag^{-1} \mS_{k}^\top \nabla f(x^k)  \mS_{k} \right\|_{\mQ'}^2}=
\E{ \|\nabla f(x^k) \|^2_{\mPdiag^{-1}  \mS_{k} \mS_{k}^\top \mQ'  \mS_{k} \mS_{k}^\top\mPdiag^{-1} }}
\\
&=&
\|\nabla f(x^k) \|^2_{\mPdiag^{-1} (\Probmat \circ \mQ')\mPdiag^{-1}} .
\end{eqnarray*}
 It remains to combine the two bounds.  
\end{proof}

We also state the analogue of Lemma~\ref{lem:sega_2}, which allows for a different norm as well. 
\begin{lemma} \label{lem:sega_2xxxxxx} For all diagonal $\mD\succ 0$ we have
\begin{equation}\label{eq:sega_h_dec_general}
\E{\|h^{k+1} \|_{\mD}^2} = \|h^k \|_{\mD-\mPdiag \mD}^2 + \|\nabla f(x^k) \|_{\mPdiag \mD}^2.
\end{equation}
\end{lemma}

\begin{proof}  
Denote $\mS_{k}$ to be a matrix with columns $e_i$ for $i\in \mS_{k}$.  
We first write \[h^{k+1} = h^k - \mS_{k} \mS_{k}^\top  h^k + \mS_{k} \mS_{k}^\top  \nabla f(x^k)  .\]
Therefore
\begin{eqnarray*}
\E{\|h^{k+1}\|_{\mD}^2} &=& \E{\left\|(\mI- \mS_{k}\mS_{k}^\top )h^k + \mS_{k} \mS_{k}^\top  \nabla f(x^k) \right\|_{\mD}^2}
\\
&=& 
\E{\left\|(\mI- \mS_{k}\mS_{k}^\top )h^k\right\|_{\mD}^2} +\E{\left\| \mS_{k} \mS_{k}^\top  \nabla f(x^k) \right\|_{\mD}^2} 
\\
&& \qquad 
+ 2\E{{h^k}^\top (\mI- \mS_{k}\mS_{k}^\top ) \mD  \mS_{k} \mS_{k}^\top  \nabla f(x^k)} 
\\
&=&
\|h^k \|_{\mD-\mPdiag \mD}^2 + \|\nabla f(x^k) \|_{\mPdiag \mD}^2.
\end{eqnarray*}
  
\end{proof}

\subsection{Proof of Theorem~\ref{t:imp_dacc}}
\begin{proof}
Throughout the proof, we will use the following Lyapunov function: 
\begin{equation*}
	\Lnacc^k \eqdef f(x^{k})-f(x^*)+ \sigma \|h^{k} \|^2_{\mP^{-1}}.
\end{equation*}
Following similar steps to what we did before, we obtain
\begin{eqnarray*}
\E{\Lnacc^{k+1}}
&\stackrel{\eqref{eq:acd_M-smooth-intro}}{\leq}& 
f(x^{k})-f(x^*)+\alpha \E{\langle\nabla f(x^k), g^k\rangle}+\frac{\alpha^2}{2} \E{\|g^k\|_\mM^2} +\sigma\E{\|h^{k+1} \|^2_{\mPdiag^{-1}}}
\\
&=&
f(x^{k})-f(x^*)-\alpha \|\nabla f(x^k) \|_2^2+\frac{\alpha^2}{2} \E{\|g^k\|_\mM^2} +\sigma\E{\|h^{k+1} \|^2_{\mPdiag^{-1}}}
\\
&\stackrel{\eqref{eq:sega_almost_eso_g}}{\leq} &
f(x^{k})-f(x^*)-\alpha \|\nabla f(x^k)\|_2^2
+\alpha^2  \|\nabla f(x^k) \|^2_{\mPdiag^{-1} (\Probmat \circ \mM)\mPdiag^{-1}} 
\\
&& \qquad
+ \alpha^2\|h^k\|^2_{\mPdiag^{-1} (\Probmat \circ \mM)\mPdiag^{-1}-\mM} +\sigma\E{\|h^{k+1} \|^2_{\mPdiag^{-1}}}.
\end{eqnarray*}
This is the place where the ESO~assumption comes into play. By applying it to the right-hand side of the bound above, we obtain
\begin{eqnarray*}
\E{\Lnacc^{k+1}} 
&\stackrel{\eqref{eq:sega_ESO}}{\leq} &
f(x^{k})-f(x^*)-\alpha \|\nabla f(x^k)\|_2^2
+\alpha^2  \|\nabla f(x^k) \|^2_{\mVdiag\mPdiag^{-1}} + \alpha^2\|h^k\|^2_{\mVdiag\mPdiag^{-1}-\mM}
\\
&& \qquad
+\sigma\E{\|h^{k+1} \|^2_{\mPdiag^{-1}}}
\\
&\stackrel{\eqref{eq:sega_h_dec_general}}{=} &
f(x^{k})-f(x^*)-\alpha \|\nabla f(x^k)\|_2^2
+\alpha^2  \|\nabla f(x^k) \|^2_{\mVdiag\mPdiag^{-1}}  + \alpha^2\|h^k\|^2_{\mVdiag\mPdiag^{-1}-\mM}
\\
&& \qquad
+\sigma\| \nabla f(x^k)\|_2^2+\sigma\|h^k \|^2_{\mPdiag^{-1}-\mI}
\\
&= &
f(x^{k})-f(x^*)-\left(\alpha - \alpha^2\max_{i}\frac{v_{i}}{p_{i}}- \sigma \right)\|\nabla f(x^k)\|_2^2
\\
&& \qquad
+  \|h^k\|^2_{\alpha^2(\mVdiag\mPdiag^{-1}-\mM)+\sigma(\mPdiag^{-1}-\mI)}.
\end{eqnarray*}
Due to Polyak-\L{}ojasiewicz inequality, we can further upper bound the last expression by
\begin{eqnarray*}
\left(1-\left(\alpha - \alpha^2\max_{i}\frac{v_{i}}{p_{i}}- \sigma \right)\mu \right)(f(x^{k})-f(x^*))
+  \|h^k\|^2_{\alpha^2(\mVdiag\mP^{-1}-\mM)+\sigma(\mP^{-1}-\mI)}.
\end{eqnarray*}
To finish the proof, it remains to use~\eqref{eq:sega_assumption}.  
\end{proof}
\subsection{Proof of Corollary~\ref{cor:sega_imp_dacc}}
The claim was obtained by choosing carefully $\alpha$ and~$\sigma$ using numerical grid search. Note that by strong convexity we have  $\mI\succeq \mu \diag(\mM)^{-1}$, so we can satisfy assumption~\eqref{eq:sega_assumption}. Then, the claim follows immediately noticing that we can also set $\mVdiag=\diag(\mM)$ while maintaining
\[
\left(\alpha - \alpha^2\max_{i}\frac{\mM_{ii}}{p_{i}}- \sigma \right) \geq \frac{0.117 }{\tracee(\mM)}.
\]

\subsection{Accelerated \texttt{SEGA} with arbitrary sampling\label{sec:sega_acc_thm}} 
Before establishing the main theorem, we first state two technical lemmas which will be crucial for the analysis. First one, Lemma~\ref{lem:sega_mirror} provides a key inequality following from~\eqref{eq:sega_z_update}. The second one, Lemma~\ref{lem:sega_acc_grad}, analyzes update~\eqref{eq:sega_y_update} and was technically established throughout the proof of Theorem~\ref{t:imp_dacc}. We include a proof of lemmas in Appendix~\ref{app:mirror} and~\ref{app:grad} respectively. 

\begin{lemma}\label{lem:sega_mirror}
For every $u\in \R^d$ we have
\begin{eqnarray}
&& \beta\langle \nabla f(x^{k+1}), z^{k}-u  \rangle -\frac{\beta \mu}{2}\|x^{k+1}-u \|_2^2\nonumber
\\
&& \qquad  \qquad \leq 
\beta^2\frac12\E{\|g^k\|_2^2}+\frac{1}{2}\|z^k-u\|_2^2
-\frac{1+\beta \mu}{2}\E{\| z^{k+1}-u\|_2^2} \label{eq:sega_md_imp}
\end{eqnarray}
 \end{lemma}

 \begin{lemma} \label{lem:sega_acc_grad} Letting $ \TD(v,p)\eqdef\max_i \frac{\sqrt{v_i}}{p_i}$, we have
\begin{eqnarray}
f(x^{k+1})-\E{f(y^{k+1})} +\|h^k\|^2_{\alpha^2(\mVdiag\mPdiag^{-3}-\mPdiag^{-1}\mM\mPdiag^{-1})} 
\geq  \left(\alpha - \alpha^2\TD(v,p)^2 \right)
 \| \nabla f(x^k)\|^2_{\mPdiag^{-1}}.  \label{eq:sega_gd}
\end{eqnarray}
 \end{lemma}
 
Now we state the main theorem of Section~\ref{s:acc}, providing a convergence rate of \texttt{ASEGA} (Algorithm~\ref{alg:sega_acc}) for arbitrary minibatch sampling.
As we mentioned, the convergence rate is, up to a constant factor, same as state-of-the-art minibatch accelerated coordinate descent~\cite{hanzely2018accelerated}. 

  \begin{theorem}\label{t:imp_acc}
Assume $\mM$-smoothness and $\mu$-strong convexity and that $v$ satisfies~\eqref{eq:sega_ESO}. 
 Denote 
 \[
 \Lacc^{k} \eqdef\frac{2}{75} \frac{  \TD(v,p)^{-2}}{ \tau^2 } \left(\E{f(y^{k})}-f(x^{*})\right)
 + \frac{1+\beta\mu}{2}\E{\| z^{k}-x^*\|_2^2}+
 \sigma \E{\|h^{k} \|^2_{\mPdiag^{-2}}}
 \] 
 and choose 
\begin{eqnarray}
c_1&=&\max \left( 1,  \TD(v,p)^{-1} \frac{\sqrt{\mu}}{\min_i p_i} \right), \label{eq:sega_c4_choice}
\\
\alpha &=&   \frac{1}{5 \TD(v,p)^{2}}, \label{eq:sega_alfa_choice}
\\
\beta &=& \frac{2}{75 \tau  \TD(v,p)^{2}}, \label{eq:sega_beta_choice}
\\
\sigma &=& 5\beta^2,  \label{eq:sega_sigma_choice}
\\
\tau &=&\frac{\sqrt{\frac{4}{9\cdot 5^4}  \TD(v,p)^{-4} \mu^2+\frac{8}{75}  \TD(v,p)^{-2} \mu }-\frac{2}{75} \TD(v,p)^{-2} \mu}{2} .  \label{eq:sega_tau_choice}
\end{eqnarray}
Then, we have
\[
\E{\Lacc^{k}}\leq \left(1-c_1^{-1}\tau \right)^k\Lacc^0.
\]
\end{theorem}

\begin{proof}  The proof technique is inspired by Allen-Zhu and Orecchia \cite{allen2014linear}. First of all, let us see what strong convexity of $f$ gives us:
\begin{eqnarray*}
\beta \left( f(x^{k+1})- f(x^*)\right)
\leq
\beta \langle \nabla f(x^{k+1}),x^{k+1}-x^*\rangle -\frac{\beta \mu}{2} \|x^*-x^{k+1}\|_2^2.
\end{eqnarray*}
Thus, we are interested in finding an upper bound for the scalar product that appeared above. We have
 \begin{eqnarray*}
&& \beta\langle \nabla f(x^{k+1}), z^{k}-u  \rangle -\frac{\beta \mu}{2}\|x^{k+1}-u \|_2^2+ \sigma \E{\|h^{k+1} \|^2_{\mPdiag^{-2}}}
\\
&& \qquad  \qquad 
\stackrel{\eqref{eq:sega_md_imp}}{\leq} 
\beta^2\frac12\E{\|g^k\|_2^2}+\frac{1}{2}\|z^k-u\|_2^2
-\frac{1+\beta \mu}{2}\E{\| z^{k+1}-u\|_2^2} +  \sigma \E{\|h^{k+1} \|^2_{\mPdiag^{-2}}}.
\end{eqnarray*}
Using the Lemmas introduced above, we can upper bound the norms of $g^k$ and $h^{k+1}$ by using norms of $h^k$ and $\nabla f(x^k)$ to get the following:
\begin{eqnarray*}
&&\beta^2\frac12\E{\|g^k\|_2^2} +  \sigma \E{\|h^{k+1} \|^2_{\mPdiag^{-2}}}\\
&&\qquad\qquad\stackrel{\eqref{eq:sega_h_dec_general}}{\leq}
\beta^2\frac12\E{\|g^k\|_2^2}
+\sigma\|h^{k} \|^2_{\mPdiag^{-2}-\mPdiag^{-1}} +\sigma  \|\nabla f(x^k)\|^2_{\mPdiag^{-1}}
\\
&&\qquad\qquad
\stackrel{\eqref{eq:sega_almost_eso_g}}{\leq}
\beta^2\|h^k \|_{\mPdiag^{-1}- \mI}^2 + \beta^2 \|\nabla f(x^k)\|^2_{\mPdiag^{-1}}
+\sigma\|h^{k} \|^2_{\mPdiag^{-2}-\mPdiag^{-1}} +\sigma  \|\nabla f(x^k)\|^2_{\mPdiag^{-1}}.
\end{eqnarray*}
Now, let us get rid of $\nabla f(x^k)$ by using the gradients property from Lemma~\ref{lem:sega_acc_grad}:
\begin{eqnarray*}
&&\beta^2\frac12\E{\|g^k\|_2^2} +  \sigma \E{\|h^{k+1} \|^2_{\mPdiag^{-2}}}\\
&&\qquad\qquad
\stackrel{\eqref{eq:sega_gd}}{\leq} 
\beta^2\|h^k \|_{\mPdiag^{-1}- \mI}^2+ \left( \beta^2 +\sigma \right) \frac{f(x^{k+1})-f(y^{k+1}) +\|h^k\|^2_{\alpha^2(\mVdiag\mPdiag^{-3}-\mPdiag^{-1}\mM\mPdiag^{-1})}}{ \alpha - \alpha^2\TD(v,p)^2} 
\\
&& \qquad\qquad \qquad+\sigma\|h^{k} \|^2_{\mPdiag^{-2}-\mPdiag^{-1}}
\\
&& \qquad\qquad
=
\|h^k \|^2_{\beta^2(\mPdiag^{-1}- \mI)+ \frac{(\beta^2+\sigma) \alpha^2}{ \alpha - \alpha^2\TD(v,p)^2} (\mVdiag\mPdiag^{-3}-\mPdiag^{-1}\mM\mPdiag^{-1}) +\sigma (\mPdiag^{-2}-\mPdiag^{-1})}
\\
&&\qquad\qquad\qquad
+  \frac{\beta^2 + \sigma }{ \alpha - \alpha^2\TD(v,p)^2}( f(x^{k+1})-\E{f(y^{k+1})}) 
\\
&& \qquad\qquad
\leq
\|h^k \|^2_{\beta^2\mPdiag^{-1}+ \frac{(\beta^2+\sigma) \alpha^2}{ \alpha - \alpha^2\TD(v,p)^2} \mVdiag\mPdiag^{-3} +\sigma (\mPdiag^{-2}-\mPdiag^{-1})}
\\
&&\qquad\qquad\qquad
+  \frac{\beta^2 + \sigma }{ \alpha - \alpha^2\TD(v,p)^2}( f(x^{k+1})-\E{f(y^{k+1})}) .
\end{eqnarray*}
Plugging this into the bound with which we started the proof, we deduce
\begin{eqnarray*}
&& \beta\langle \nabla f(x^{k+1}), z^{k}-u  \rangle -\frac{\beta \mu}{2}\|x^{k+1}-u \|_2^2+ \sigma \E{\|h^{k+1} \|^2_{\mPdiag^{-2}}}\\
&& \qquad  \qquad 
\leq
\|h^k \|^2_{\beta^2\mPdiag^{-1}+ \frac{(\beta^2+\sigma) \alpha^2}{ \alpha - \alpha^2\TD(v,p)^2} \mVdiag\mPdiag^{-3} +\sigma (\mPdiag^{-2}-\mPdiag^{-1})}
\\
&& \qquad \qquad \qquad 
+  \frac{\beta^2 + \sigma }{ \alpha - \alpha^2\TD(v,p)^2}( f(x^{k+1})-\E{f(y^{k+1})}) 
+\frac{1}{2}\|z^k-u \|_2^2
\\
&& \qquad \qquad \qquad
-\frac{1+\beta \mu}{2}\E{\| z^{k+1}-u\|_2^2 }.
\end{eqnarray*} 
Recalling our first step, we get with a few rearrangements
\begin{eqnarray*}
& & \beta \left( f(x^{k+1})- f(x^*)\right)
 \\
& & \qquad \qquad 
\leq
\beta \langle \nabla f(x^{k+1}),x^{k+1}-x^*\rangle -\frac{\beta \mu}{2} \|x^*-x^{k+1}\|_2^2
\\
& & \qquad \qquad 
=
\beta \langle \nabla f(x^{k+1}),x^{k+1}-z^{k}\rangle +\beta \langle \nabla f(x^{k+1)},z^{k}-x^*\rangle -\frac{\beta \mu}{2} \|x^*-x^{k+1}\|_2^2
\\
& & \qquad \qquad 
=
\frac{(1-\tau)\beta}{\tau} \langle \nabla f(x^{k+1}),y^k-x^{k+1}\rangle +\beta \langle \nabla f(x^{k+1}),z^{k}-x^*\rangle -\frac{\beta \mu}{2} \|x^*-x^{k+1}\|_2^2
\\
& &\qquad \qquad  
\leq 
\frac{(1-\tau)\beta}{\tau}\left( f(y^k)-f(x^{k+1})\right)
+
\|h^k \|^2_{\beta^2\mPdiag^{-1}+ \frac{(\beta^2+\sigma) \alpha^2}{ \alpha - \alpha^2\TD(v,p)^2} \mVdiag\mPdiag^{-3} +\sigma (\mPdiag^{-2}-\mPdiag^{-1})}
\\
& & \qquad \qquad\qquad+
  \frac{\beta^2 + \sigma }{ \alpha - \alpha^2\TD(v,p)^2}( f(x^{k+1})-\E{f(y^{k+1})}) 
   +
\frac{1}{2}\|z^k-x^*\|_2^2
\\
& & \qquad \qquad\qquad
-\frac{1+\beta\mu}{2}\E{\| z^{k+1}-x^*\|_2^2} - \sigma \E{\|h^{k+1} \|^2_{\mPdiag^{-2}}}.
\end{eqnarray*}
Let us choose $\sigma$, $\beta$ such that for some constant $c_2$ (which we choose at the end) we have
 \[
c_2\sigma=\beta^2, \qquad \beta=\frac{\alpha - \alpha^2\TD(v,p)^2}{(1+c_2^{-1}) \tau }.
 \]
 Consequently, we have
 \begin{eqnarray*}
&&
\frac{\alpha - \alpha^2\TD(v,p)^2}{(1+c_2^{-1}) \tau^2 } \left(\E{f(y^{k+1})}-f(x^{*})\right)
 + \frac{1+\beta\mu}{2}\E{\| z^{k+1}-x^*\|_2^2}+
 \sigma \E{\|h^{k+1} \|^2_{\mPdiag^{-2}}}
 \\
  &&
 \qquad \qquad \leq 
(1-\tau)\frac{\alpha - \alpha^2\TD(v,p)^2}{(1+c_2^{-1}) \tau^2 }\left( f(y^k)-f(x^{*})\right)+
  \frac{1}{2}\|z^k-x^*\|_2^2
  \\
&& 
\qquad \qquad \qquad+
  \|h^k \|^2_{\left(\mPdiag^{-1} -(1-c_2)\mI+  \frac{(1+c_2)  \alpha^2}{ \alpha - \alpha^2\TD(v,p)^2} \mVdiag\mPdiag^{-2} \right) \sigma \mPdiag^{-1}} .
 \end{eqnarray*}
 Let us make a particular choice of $\alpha$, so that for some constant $c_3$ (which we choose at the end) we can obtain the equations below:
\begin{eqnarray*}
  \alpha = \frac{1}{c_3\TD(v,p)^2} 
  \quad
  \Rightarrow 
  \quad
  \alpha - \alpha^2\TD(v,p)^2
 &=&
   \frac{c_3-1}{c_3^2}\TD(v,p)^{-2},
   \\
  \frac{\alpha^2}{\alpha - \alpha^2\TD(v,p)^2}
  &=&
  \frac{1}{(c_3-1)\TD(v,p)^2}
  .
\end{eqnarray*}
 Thus
  \begin{eqnarray*}
&&
\frac{ \frac{c_3-1}{c_3^2}\TD(v,p)^{-2}}{(1+c_2^{-1}) \tau^2 } \left(\E{f(y^{k+1})}-f(x^{*})\right)
 + \frac{1+\beta\mu}{2}\E{\| z^{k+1}-x^*\|_2^2}+
 \sigma \E{\|h^{k+1} \|^2_{\mPdiag^{-2}}}
 \\
  &&
 \qquad \qquad \leq 
(1-\tau)\frac{ \frac{c_3-1}{c_3^2}\TD(v,p)^{-2}}{(1+c_2^{-1}) \tau^2 }\left( f(y^k)-f(x^{*})\right)+
  \frac{1}{2}\|z^k-x^*\|_2^2
  \\
&& 
\qquad \qquad \qquad+
  \|h^k \|^2_{\left(\mPdiag^{-1} -(1-c_2)\mI+  \frac{(1+c_2) }{ (c_3-1)\TD(v,p)^{2}} \mVdiag\mPdiag^{-2} \right) \sigma \mPdiag^{-1}}.
 \end{eqnarray*}
 Using the definition of $\TD(v,p)$, one can see that the above gives

  \begin{eqnarray*}
&&
\frac{ \frac{c_3-1}{c_3^2} \TD(v,p)^{-2}}{(1+c_2^{-1}) \tau^2 } \left(\E{f(y^{k+1})}-f(x^{*})\right)
 + \frac{1+\beta\mu}{2}\E{\| z^{k+1}-x^*\|_2^2}+
 \sigma \E{\|h^{k+1} \|^2_{\mPdiag^{-2}}}
 \\
  &&
 \qquad \qquad \leq 
(1-\tau)\frac{ \frac{c_3-1}{c_3^2} \TD(v,p)^{-2}}{(1+c_2^{-1}) \tau^2 }\left( f(y^k)-f(x^{*})\right)+
  \frac{1}{2}\|z^k-x^*\|_2^2 \\
    &&
 \qquad \qquad \qquad
+
  \|h^k \|^2_{\left(\mPdiag^{-1} -(1-c_2)\mI+  \frac{1+c_2 }{ c_3-1} \mI \right) \sigma \mPdiag^{-1}}.
 \end{eqnarray*}
 To get the convergence rate, we shall establish
 \begin{equation}\label{eq:sega_h_suff_acc_imp}
 \left(1-c_2-  \frac{1+c_2 }{ c_3-1}\right) c_1 \mI \succeq \tau\mPdiag^{-1}
 \end{equation}
 and
 \begin{equation}\label{eq:sega_acc_tau_implicit}
1+\beta \mu \geq \frac{1}{1-\tau}.
 \end{equation}
 To this end, let us recall that
 \[
 \beta=\frac{c_3-1}{c^2_2}  \TD(v,p)^{-2}\tau^{-1} \frac{1}{1+c_2^{-1}}.
 \]
Now we would like to set equality in~\eqref{eq:sega_acc_tau_implicit}, which yields
\[
0=\tau^2+\frac{c_3-1}{c^2_2} \TD(v,p)^{-2}\frac{1}{1+c_2^{-1}}  \mu  \tau- \frac{c_3-1}{c^2_2}  \TD(v,p)^{-2}\frac{1}{1+c_2^{-1}} \mu =0.
\]
 This, in turn, implies
 \begin{eqnarray*}
\tau
&=&
\frac{\sqrt{\left(\frac{c_3-1}{c^2_2}\right)^2  \TD(v,p)^{-4}\frac{1}{\left(1+c_2^{-1}\right)^2} \mu^2+4\frac{c_3-1}{c^2_2}  \TD(v,p)^{-2} \frac{1}{1+c_2^{-1}}\mu }-\frac{c_3-1}{c^2_2}  \TD(v,p)^{-2} \frac{1}{1+c_2^{-1}}\mu}{2}
\\
&=&
\cO\left(\sqrt{\frac{c_3-1}{c^2_2}}\frac{1}{\sqrt{1+c_2^{-1}}} \TD(v,p)^{-1}\sqrt{\mu} \right).
  \end{eqnarray*}
 Notice that for any $c\leq 1$ we have $\frac{\sqrt{c^2+4c}-c}{2}\leq \sqrt{c}$ and therefore
\begin{equation}\label{eq:sega_tau_bound}
\tau \leq \sqrt{\frac{c_3-1}{c^2_2}}  \TD(v,p)^{-1} \frac{1}{\sqrt{1+c_2^{-1}}} \sqrt{\mu} .
\end{equation}
 Using this inequality and a particular choice of constants, we can upper bound $\mP^{-1}$ by a matrix proportional to identity as shown below:
 \begin{eqnarray*}
 \tau \mPdiag^{-1} &\stackrel{\eqref{eq:sega_tau_bound}}{\preceq}&
 \sqrt{\frac{c_3-1}{c^2_2}} \TD(v,p)^{-1} \frac{1}{\sqrt{1+c_2^{-1}}} \sqrt{\mu} \mPdiag^{-1} 
 \\
 &\preceq &
   \sqrt{\frac{c_3-1}{c^2_2}} \TD(v,p)^{-1}\frac{1}{\sqrt{1+c_2^{-1}}} \frac{\sqrt{\mu}}{\min_i p_i}  \mI
 \\
 &\stackrel{\eqref{eq:sega_c4_choice}}{\preceq} &
   \sqrt{\frac{c_3-1}{c^2_2}}\frac{1}{\sqrt{1+c_2^{-1}}} c_1 \mI
 \\
 &\stackrel{(*)}{\preceq} & 
  \left(1-c_2-  \frac{1+c_2 }{ c_3-1}\right) c_1 \mI ,
 \end{eqnarray*}
 which is exactly \eqref{eq:sega_h_suff_acc_imp}. Above, $(*)$ holds for choice $c_3=5$ and $c_2=\frac{1}{5}$. It remains to verify that~\eqref{eq:sega_alfa_choice}, \eqref{eq:sega_beta_choice}, \eqref{eq:sega_sigma_choice} and~\eqref{eq:sega_tau_choice} indeed correspond to our derivations.
  
\end{proof}

We also mention, without a proof, that acceleration parameters can be chosen in general such that $c_1$ can be lower bounded by constant and therefore the rate from Theorem~\ref{t:imp_acc} coincides with the rate from Table~\ref{tab:CDcmp}. Corollary~\ref{cor:sega_acc_imp} is in fact a weaker result of that type. 
 
 \subsubsection{Proof of Corollary~\ref{cor:sega_acc_imp}}
It suffices to verify that one can choose $v=\diag(\mM)$ in~\eqref{eq:sega_ESO} and that due to $p_i\propto \sqrt{\mM_{ii}}$ we have $c_1= 1$.

\subsection{Proof of Lemma~\ref{lem:sega_mirror}\label{app:mirror}} 
 \begin{proof}
 Firstly \eqref{eq:sega_z_update}, is equivalent to
\[
z^{k+1}=\argmin_z \psi^k(z)\eqdef \frac{1}{2}\| z-z^k\|_2^2+ \beta \langle g^k, z \rangle   +\frac{\beta \mu}{2}\|z-x^{k+1}\|_2^2.
\]
Therefore, we have for every $u$
\begin{align}
\nonumber
0&=\langle \nabla \psi^k(z^{k+1}),z^{k+1}-u \rangle \\
&=
\langle z^{k+1}-z^k, z^{k+1}-u\rangle +\beta \langle g^k, z^{k+1}-u  \rangle +\beta \mu \langle z^{k+1}-x^{k+1}, z^{k+1}-u\rangle. \label{eq:sega_zk_plus_1_optimal}
\end{align}
Next, by generalized Pythagorean theorem we have
\begin{equation}\label{eq:sega_pyt_z}
\langle z^{k+1}-z^k,z^{k+1}-u \rangle =\frac12 \|z^k-z^{k+1}\|_2^2-\frac12 \|z^k-u\|_2^2+\frac12 \|u-z^{k+1}\|_2^2
\end{equation}
and
\begin{equation}\label{eq:sega_pyt_x}
\langle z^{k+1}-x^{k+1},z^{k+1}-u \rangle =\frac12 \|x^{k+1}-z^{k+1}\|_2^2-\frac12 \|x^{k+1}-u\|_2^2+\frac12 \|u-z^{k+1}\|_2^2.
\end{equation}
Plugging~\eqref{eq:sega_pyt_z} and~\eqref{eq:sega_pyt_x} into~\eqref{eq:sega_zk_plus_1_optimal} we obtain
\begin{eqnarray*}
&&
\beta \langle g^k, z^{k}-u  \rangle -\frac{\beta \mu}{2}\|x^{k+1}-u\|_2^2
\\
&& \qquad \qquad
\leq 
\beta  \langle g^k, z^{k}-z^{k+1}  \rangle
 -\frac12\| z^k-z^{k+1}\|_2^2+\frac{1}{2}\|z^k-u\|_2^2-\frac{1+\beta \mu}{2}\| z^{k+1}-u\|_2^2
\\
&& \qquad \qquad
\stackrel{(*)}{\leq} 
\frac{\beta^2}{2}\| g^k \|_2^2+\frac{1}{2}\|z^k-u\|_2^2
-\frac{1+\beta \mu}{2}\| z^{k+1}-u\|_2^2.
\end{eqnarray*}
The step marked by $(*)$ holds due to Cauchy-Schwartz inequality. 
It remains to take the expectation conditioned on $x^{k+1}$ and use~\eqref{eq:sega_unbiased_estimator}.

\end{proof}
 
\subsection{Proof of Lemma~\ref{lem:sega_acc_grad}\label{app:grad}}
 \begin{proof}
 The shortest, although not the most intuitive, way to write the proof is to put matrix factor into norms. Apart from this trick, the proof is quite simple consists of applying smoothness followed by ESO:
\begin{eqnarray*}
\E{f(y^{k+1})} - f(x^{k+1})
&\stackrel{\eqref{eq:acd_M-smooth-intro}}{\leq}& 
-\alpha \E{\langle\nabla f(x^k), \mPdiag^{-1} g^k\rangle}+\frac{\alpha^2}{2} \E{\|\mPdiag^{-1} g^k\|_\mM^2} 
\\
&=&
-\alpha \|\nabla f(x^k)\|^2_{\mPdiag^{-1}}  +\frac{\alpha^2}{2} \E{\| g^k\|_{\mPdiag^{-1}\mM\mPdiag^{-1}}} 
\\
&\stackrel{\eqref{eq:sega_almost_eso_g}}{\leq} &
-\alpha \|\nabla f(x^k)\|^2_{\mPdiag^{-1}}  +\alpha^2\|\nabla f(x^k) \|^2_{\mPdiag^{-1}(\Probmat \circ \mPdiag^{-1}\mM\mPdiag^{-1})\mPdiag^{-1}}
\\
&& \qquad 
+\alpha^2 \| h^k\|^2_{\mPdiag^{-1}(\Probmat \circ \mPdiag^{-1}\mM\mPdiag^{-1})\mPdiag^{-1}-\mPdiag^{-1}\mM\mPdiag^{-1}} 
\\
&= &
-\alpha \|\nabla f(x^k)\|^2_{\mPdiag^{-1}}  +\alpha^2\|\nabla f(x^k) \|^2_{\mPdiag^{-2}(\Probmat \circ\mM)\mPdiag^{-2}}
\\
&& \qquad 
+\alpha^2 \| h^k\|^2_{\mPdiag^{-2}(\Probmat \circ \mM)\mPdiag^{-2}-\mPdiag^{-1}\mM\mPdiag^{-1}} 
\\
&\stackrel{\eqref{eq:sega_ESO}}{\leq} &
-\alpha \|\nabla f(x^k)\|^2_{\mPdiag^{-1}}  +\alpha^2\|\nabla f(x^k) \|^2_{\mVdiag \mPdiag^{-3}}
\\
&& \qquad 
+\alpha^2 \| h^k\|^2_{\mVdiag \mPdiag^{-3}-\mPdiag^{-1}\mM\mPdiag^{-1}} 
\\
&\leq&
-\left(\alpha - \alpha^2\max_{i}\frac{v_i}{p_{i}^2} \right)\|f(x^k)\|^2_{\mPdiag^{-1}}  
+\alpha^2 \| h^k\|^2_{\mVdiag\mPdiag^{-3}-\mPdiag^{-1}\mM\mPdiag^{-1}} .
\end{eqnarray*}

\end{proof}

\section{Subspace \texttt{SEGA}: a more aggressive approach \label{sec:sega_subSEGA}}

In this section we describe a {\em more aggressive} variant of \texttt{SEGA}, one that exploits the fact that the gradients of $f$ lie in a lower dimensional subspace if this is indeed the case.

In particular, assume that $F(x) = f(x) + \psi(x)$ and \[f(x) = \phi(\mA x),\] where $\mA\in \R^{m\times d}$.\footnote{Strong convexity is not compatible with the assumption that $\mA$ does not have full rank, so a different type of analysis using Polyak-\L{}ojasiewicz inequality is required to give a formal justification. However, we proceed with the analysis anyway to build the intuition why this approach leads to better rates.} Note that $\nabla f(x)$ lies in $\Range{\mA^\top}$. There are situations where the dimension of $\Range{\mA^\top}$ is much smaller than $n$. For instance, this happens when $m\ll d$. However, standard coordinate descent methods still move around in directions $e_i\in \R^d$ for all $i$. We can modify the gradient sketch method to force our gradient estimate to lie in $\Range{\mA^\top}$, hoping that this will lead to faster convergence.

\subsection{The algorithm}

Let $x^k$ be the current iterate, and let  $h^k$ be the current estimate of the gradient of $f$. Assume that the sketch $\mS_k^\top \nabla f(x^k)$ is available. We can now define $h^{k+1}$ through the following modified sketch-and-project process:
\begin{eqnarray}
h^{k+1} &=& \arg \min_{h\in \R^{d}} \| h -  h^k\| ^2 \notag \\
&& \text{subject to} \quad \mS_k^\top h = \mS_k^\top \nabla f(x^k), \label{eq:sega_sketch-n-project2B}\\
&& \phantom{subject to} \quad h \in \Range{\mA^\top}.\notag
\end{eqnarray}

Standard arguments reveal that the closed-form solution of \eqref{eq:sega_sketch-n-project2B} is
\begin{equation}h^{k+1} = \mH\ \left(h^k - \mS_k(\mS_k^\top\mH\mS_k)^\dagger \mS_k^\top(\mH h^k - \nabla f(x^k)) \right), \label{eq:sega_h^{k+1}3}\end{equation}
where  
\begin{align}\label{eq:sega_H}
\mH \eqdef \mA^\top (\mA \mA^\top)^\dagger \mA
\end{align}
 is the projector onto $\Range{\mA^\top}$.
A quick sanity check reveals that this gives the same formula as \eqref{eq:sega_h^{k+1}} in the case where $\Range{\mA^\top} = \R^d$. We can also write
\begin{equation} \label{eq:sega_general_update_of_h}h^{k+1} = \mH h^k -   \mH\mZ_k(\mH h^k - \nabla f(x^k)) = \left(\mI - \mH\mZ_k\right)\mH h^k + \mH\mZ_k\nabla f(x^k),\end{equation}
where 
\begin{align}\label{eq:sega_Z_k}
\mZ_k \eqdef \mS_k(\mS_k^\top\mH\mS_k)^\dagger \mS_k^\top.
\end{align}
 Assume that $\theta_k$ is chosen in such a way that
 \begin{equation*}
     \E{\theta_k\mZ_k} = \mI.
 \end{equation*}
 Then, the following estimate of $\nabla f(x^k)$
 \begin{align}\label{eq:sega_g^k_agressive}
     g^k\eqdef \mH h^k + \theta_k\mH\mZ_k (\nabla f(x^k) - \mH h^k)
 \end{align}
 is unbiased, i.e.\ $\E{g^k} = \nabla f(x^k)$. After evaluating $g^k$, we perform the same step as in \texttt{SEGA}: \[x^{k+1} = \prox_{\alpha \psi}(x^k - \alpha g^k). \]
By inspecting \eqref{eq:sega_sketch-n-project2B}, \eqref{eq:sega_H} and \eqref{eq:sega_g^k_agressive}, we get the following simple observation.
\begin{lemma} \label{lem:sega_range_preserve}If $h^0\in \Range{\mA^\top}$, then $h^k, g^k \in \Range{\mA^\top}$ for all $k$.
\end{lemma}

Consequently, if $h^0 \in \Range{\mA^\top}$, \eqref{eq:sega_h^{k+1}3}  simplifies to
\begin{equation} \label{eq:sega_b98g9fd2}
h^{k+1} =  h^k - \mH\mS_k(\mS_k^\top\mH\mS_k)^\dagger \mS_k^\top(h^k - \nabla f(x^k)) 
\end{equation}
and \eqref{eq:sega_g^k_agressive} simplifies to
\begin{equation} \label{eq:sega_g^k2-xx}
 g^k \eqdef h^k + \theta_k\mH\mZ_k (\nabla f(x^k) - h^k).\end{equation}

\begin{example}[Coordinate sketch]\label{ex:coord_setup2}
Consider $\cD$ given by $\mS=e_i $ with probability $p_i>0$. Then we can choose the bias-correcting random variable as $\theta = \theta(s) = \frac{w_i}{p_i} $, where $w_i \eqdef \|\mH e_i\|_2^2 = e_i^\top \mH e_i$. Indeed, with this choice, \eqref{eq:sega_unbiased} is satisfied. For simplicity, further choose $p_i = 1/n$ for all $i$. We then have
\begin{equation} \label{eq:sega_b98g9fd2-coord} h^{k+1} = h^k -   \frac{e_{i}^\top  h^k - e_{i}^\top \nabla f(x^k)}{w_{i}} \mH e_{i} = \left(\mI - \frac{\mH e_{i} e_{i}^\top}{w_{i}}\right)  h^k + \frac{ \mH e_{i} e_{i}^\top}{w_{i}}\nabla f(x^k)\end{equation}
and \eqref{eq:sega_g^k2-xx} simplifies to
\begin{equation} \label{eq:sega_g^k2-xx-coord} g^k \eqdef (1-\theta_k) h^k + \theta_k h^{k+1} = h^k + n \mH e_{i} e_{i}^\top \left( \nabla f(x^k) -  h^k \right).\end{equation}
\end{example}

\subsection{Lemmas}
All theory provided in this subsection is, in fact, a straightforward generalization of our non-subspace results. The reader can recognize similarities in both statements and proofs with that of previous sections.
\begin{lemma}\label{lem:sega_properties_of_Z_and_H}
    Define $\mZ_k$ and $\mH$ as in equations~\eqref{eq:sega_Z_k} and~\eqref{eq:sega_H}. Then $\mZ_k$ is symmetric, $\mZ_k \mH \mZ_k = \mZ_k$, $\mH^2=\mH$ and $\mH = \mH^\top $.
\end{lemma}
\begin{proof}
    The symmetry of $\mZ_k$ follows from its definition. The second statement is a corollary of the equations $((\mA_1\mA_2)^\dagger)^\top = (\mA_2^\top\mA_1^\top)^\dagger$ and $\mA_1^\dagger \mA_1 \mA_1^\dagger = \mA_1^\dagger$, which are true for any matrices $\mA_1, \mA_2$. 
Finally, the last two rules follow directly from the definition of $\mH$ and the property $\mA_1^\dagger \mA_1 \mA_1^\dagger = \mA_1^\dagger$. 
\end{proof}

\begin{lemma}\label{lem:sega_bug79dv987gs}
    Assume $h^k\in \Range{\mA^\top}$. Then \[\E{\| h^{k+1} - v\| ^2} = \|h^k - v\|_{\mI - \E{\mZ}}^2 + \|\nabla f(x^k) - v\|_{\E{\mZ}}^2\] for any vector $v\in \Range{\mA^\top}$.
\end{lemma}
\begin{proof}
    By Lemma \ref{lem:sega_properties_of_Z_and_H} we can rewrite $\mH$ as $\mH^\top$, so
\begin{eqnarray}
    \E{\| h^{k+1}-v\| ^2} 
    &\overset{\eqref{eq:sega_general_update_of_h}}{=}& \E{ \left\|h^k - \mH\mZ_k (h^k - \nabla f(x^k)) - v \right\| ^2 } \nonumber\\
 &=&  \E{ \left\| \left(\mI -\mH \mZ_k \right)(h^k - v) + \mH\mZ_k (\nabla f(x^k)-v)  \right\| ^2 }\nonumber\\
 &=&  \E{ \left\| \left(\mI - \mH^\top\mZ_k \right)(h^k - v) + \mH\mZ_k(\nabla f(x^k)-v)  \right\| ^2 }\nonumber\\
 &=& \E{ \left\| \left(\mI - \mH^\top\mZ_k \right)(h^k-v) \right\| ^2} + \E{\left\|\mH\mZ_k  ( \nabla f(x^k) - v)  \right\| ^2 }\nonumber\\
 && \quad + 2 (h^k-v)^\top\E{\left(\mI - \mH^\top\mZ_k \right)^\top \mH\mZ_k}  (\nabla f(x^k) - v) \nonumber\\
 &=& (h^k-v)^\top \E{ \left(\mI - \mH^\top\mZ_k \right)^\top   \left(\mI - \mH\mZ_k \right)} (h^k-v) \nonumber\\
 && \quad + (\nabla f(x^k)-v)^\top \E{\mZ_k\mH^\top\mH\mZ_k} (\nabla f(x^k) -v)\nonumber\\
 && \quad + 2 (h^k-v)^\top\E{ \mH\mZ_k - \mZ_k\mH\mH\mZ_k}  (\nabla f(x^k) - v). \label{eq:sega_last_term_in_h_minus_v}
 \end{eqnarray}
By Lemma~\ref{lem:sega_properties_of_Z_and_H} we have 
\begin{align*}
	\mZ_k\mH\mH\mZ_k = \mZ_k\mH\mZ_k = \mZ_k,
\end{align*}
so the last term in~\eqref{eq:sega_last_term_in_h_minus_v} is equal to 0. As for the other two, expanding the matrix factor in the first term leads to
\begin{eqnarray*}
      \left(\mI - \mH^\top\mZ_k \right)^\top   \left(\mI - \mH\mZ_k \right)
     &=& \left(\mI - \mZ_k \mH \right)   \left(\mI - \mH\mZ_k\right)\\
     &=& \mI - \mZ_k\mH - \mH^\top\mZ_k + \mZ_k\mH  \mH\mZ_k\\
     &=& \mI - \mZ_k\mH - \mH^\top\mZ_k + \mZ_k.
\end{eqnarray*}
 Let us mention that $\mH(h^k - v) = h^k - v$ and $(h^k - v)^\top\mH^\top=(h^k - v)^\top$ as both vectors $h^k$ and $v$ belong to $\Range{\mA^\top}$. Therefore,
 \begin{align*}
     (h^k-v)^\top \E{ \mI - \mZ_k\mH - \mH^\top\mZ_k + \mZ_k} (h^k-v) = (h^k-v)^\top \left(\mI - \E{\mZ_k}\right) (h^k-v).
 \end{align*}
 It remains to consider
 \begin{align*}
 	\E{\mZ_k\mH^\top\mH\mZ_k} = \E{\mZ_k\mH\mH\mZ_k} = \E{\mZ_k}.
 \end{align*}
We, thereby, have derived
 \begin{eqnarray*}
 \E{\| h^{k+1}-v\| ^2} 
 &=&(h^k-v)^\top \left( \mI - \E{\mZ_k}\right) (h^k - v)\\
 &&\quad + (\nabla f(x^k)-v)^\top \E{\mZ_k\mZ_k} (\nabla f(x^k) - v)\\
 &=& \|h^k - v\|_{ \mI - \E{\mZ_k}}^2 + \|\nabla f(x^k) - v\|_{\E{\mZ}}^2.
\end{eqnarray*}

 \end{proof}

\begin{lemma}\label{lem:sega_db98gdf9jf}
    Suppose $h^k\in\Range{\mA^\top}$ and $g^k$ is defined by~\eqref{eq:sega_g^k_agressive}. Then
    \begin{align}\label{eq:sega_lemma2_agressive}
        \E{\|g^k - v\|^2 } \le \|h^k - v\|^2_{\mC - \mI} + \|\nabla f(x^k) -v\|^2_{\mC}
    \end{align}
     for any $v\in\Range{\mA^\top}$, where
     \begin{align}\label{eq:sega_mC}
       \mC\eqdef \E{\theta^2\mZ}.
     \end{align}
\end{lemma}
\begin{proof}
    Writing $g^k - v = a + b$, where $a\eqdef (\mI - \theta_k\mH\mZ_k)(h^k - v)$ and $b\eqdef \theta_k \mH\mZ_k (\nabla f(x^k) - v)$, we get $\|g^k\| ^2\le 2(\|a\| ^2 + \|b\| ^2)$.  By definition of $\theta_k$,
    \begin{eqnarray*}
        \E{\|a\| ^2} 
        &=& \E{\|\left(\mI - \theta_k\mH\mZ_k\right)(h^k - v)\| ^2}\\
        &=& (h^k - v)^\top \E{\left( \mI - \theta_k\mZ_k\mH\right)\left( \mI - \theta_k\mH\mZ_k\right)} (h^k - v)\\
        &=&  \|h^k - v\|^2_{\E{\left( \mI - \theta_k\mZ_k\mH\mI - \theta_k\mH\mZ_k + \theta_k^2\mZ_k\mH\mI\mH\mZ_k \right)} }.
    \end{eqnarray*}
According to Lemma~\ref{lem:sega_properties_of_Z_and_H},  $\mH=\mH$ and $\mZ_k\mH\mZ_k=\mZ_k$, so
    \begin{eqnarray*}
    \E{\|a\| ^2} 
        &=&(h^k - v)^\top \E{\left( \mI - \theta_k\mZ_k\mH - \theta_k\mH^\top\mZ_k + \theta_k^2\mZ_k\right)} (h^k - v)\\
        &=& \|h^k - v\|^2_{\E{\theta^2\mZ} - \mI},
    \end{eqnarray*}
    where in the last step we used the assumption that $h^k$ and $v$ are from $\Range{\mA^\top}$ and $\mH$ is the projector operator onto $\Range{\mA^\top}$.
    
    Similarly, the second term in the upper bound on $g^k$ can be rewritten as
    \begin{eqnarray*}
        \E{\|b\| ^2}
        &=& \E{\|\theta_k \mH\mZ_k (\nabla f(x^k) - v)\|^2 }\\
        &=& (\nabla f(x^k) - v)^\top \E{\theta_k^2\mZ_k \mH^\top\mH\mZ_k} (\nabla f(x^k) - v)\\
        &=& \|\nabla f(x^k) - v\|^2_{\E{\theta_k^2\mZ_k}}.
    \end{eqnarray*}
    Combining the pieces, we get the claim. 
\end{proof}

\subsection{Main result}

The main result of this section is:

\begin{theorem}\label{thm:sega_main_agressive}
    Assume that $f$ is $\mmM$-smooth, $\mu$-strongly convex, and that $\alpha>0$ is such that
\begin{equation}
        \alpha\left(2(\mC - \mI) +\sigma \mu \mI\right) \le \sigma\E{\mZ},\qquad \alpha \mC \le \frac{1}{2}\left(\mmM - \sigma \E{\mZ}\right). \label{eq:sega_general_bound_on_stepsize_agressive}        
    \end{equation}
    If we define $\Lgen^k \eqdef \|x^k - x^*\|^2  + \sigma \alpha \|h^k - \nabla f(x^k)\|^2 $, then
$
        \E{\Lgen^{k}} \le (1 - \alpha\mu)^k \Lgen^0.
$
\end{theorem}
\begin{proof}
Having established Lemmas~\ref{lem:sega_properties_of_Z_and_H}, \ref{lem:sega_bug79dv987gs} and
\ref{lem:sega_db98gdf9jf}, the proof follows the same steps as the proof of  Theorem~\ref{thm:sega_main}.  
\end{proof}

\subsection{The conclusion of subspace \texttt{SEGA}}
Let us recall that $g^k = h^k + \theta_k\mZ_k(\nabla f(x^k) - h^k)$. A careful examination shows that when we reduce $\theta_k$ from $\cO(n)$ to $\cO(d)$, we put more trust in the value of $h^k$ with the benefit of reducing the variance of $g^k$. This insight points out that a practical implementation of the algorithm may exploit the fact that $h^k$ learns the gradient of $f$ by using smaller $\theta_k$.

It is also worth noting that \texttt{SEGA} is a stationary point algorithm regardless of the value of $\theta_k$. Indeed, if one has $x^k = x^*$ and $h^k = \nabla f(x^*)$, then $g^k = \nabla f(x^*)$ for any $\theta_k$. Therefore, once we get a reasonable $h^k$, it is well grounded to choose $g^k$ to be closer to $h^k$. This argument is also supported by our experiments.

Finally, the ability to take bigger stepsizes is also of high interest. One can think of extending other methods in this direction, especially if interested in applications with a small rank of matrix~$\mA$.

\section{Simplified analysis of \texttt{SEGA}  \label{sec:sega_simple_SEGA}}

In this section we consider the setup from Example~\ref{ex:coord_setup} with uniform probabilities: $p_i=1/d$ for all $i$. We now state the main complexity result.

\begin{theorem} \label{thm:sega_simple} Choose $\cD$ to be the uniform distribution over unit basis vectors in $\R^d$. For any $\sigma>0$ define \[\Lgen^k \eqdef \|x^k-x^*\|_2^2 + \sigma \alpha \|h^k\|_2^2,\] where $\{x^k,h^k\}_{k\geq 0}$ are the iterates of the  gradient sketch  method.  If the stepsize  satisfies 
\begin{equation}\label{eq:sega_alpha_bound} 0<\alpha \leq \min\left\{ \frac{1-\frac{L\sigma}{n}}{2Ld}, \frac{1}{d\left(\mu + \frac{2(d-1)}{\sigma}\right)} \right\},\end{equation}
then
$\E{\Lgen^{k+1}} \leq (1-\alpha \mu) \Lgen^{k}.$
This means that 
\[k \geq \frac{1}{\alpha \mu} \log \frac{1}{\epsilon} \quad \Rightarrow \quad \E{\Lgen^k} \leq \epsilon \Lgen^0.\]

\end{theorem}

In particular, if we let $\sigma = \frac{d}{2L}$, then  $\alpha = \frac{1}{(4L+\mu)d}$ satisfies \eqref{eq:sega_alpha_bound}, and we have the iteration complexity \[d\left(4 + \frac{1}{\kappa}\right)\kappa \log \frac{1}{\epsilon} = \tilde{\cO}(d\kappa),\] where $\kappa\eqdef  \frac{L}{\mu}$ is the condition number.

This is the same complexity as \texttt{NSync} \cite{nsync} under the same assumptions on $f$. \texttt{NSync} also needs just access to partial derivatives. However, \texttt{NSync} uses variable stepsizes, while \texttt{SEGA} can do the same with  {\em fixed} stepsizes. This is because \texttt{SEGA} {\em learns} the direction $g^k$ using  past information.

\subsection{Technical lemmas}

Since $f$ is $L$-smooth, we have
\begin{equation}\label{eq:sega_L-smooth_inequality}
\|\nabla f(x^k)\|_{2}^2 \leq 2L(f(x^k) - f(x^*)).
\end{equation}
On the other hand, by $\mu$-strong convexity of $f$ we have \begin{equation}\label{eq:sega_8998sgjfbif}f(x^*) \geq f(x^k) + \langle \nabla f(x^k), x^*-x^k \rangle + \frac{\mu}{2}\|x^*-x^k\|_{2}^2.\end{equation}

\begin{lemma} \label{lem:sega_1} The variance of $g^k$ as an estimator of $\nabla f(x^k)$ can be bounded as follows:
\begin{equation} \label{eq:sega_bu9808hf9}
\E{\|g^k\|_{2}^2} \leq 4Ln (f(x^k) - f(x^*)) + 2(d-1)\|h^k \|_{2}^2.
\end{equation}
\end{lemma}

\begin{proof}  
In view of \eqref{eq:sega_8h0h09ffs}, we first write \[g^k = \underbrace{h^k - \frac{1}{p_i}  e_{i}^\top  h^k e_{i}}_{a} +  \underbrace{ \frac{1}{p_i} e_{i}^\top \nabla f(x^k)  e_{i}}_{b},\]
and note that $p_i=1/n$ for all $i$. Let us bound the expectation of each term individually. The first term is equal to
\begin{eqnarray*}
\E{\|a\|_{2}^2} &=& \E{\left\|h^k - d e_{i}^\top  h^k e_{i}\right\|_{2}^2}\\
&=& \E{\left\| \left(\mI - d e_{i} e_{i}^\top \right) h^k  \right\|_{2}^2}\\
&=&(h^k)^\top \E{ \left(\mI - d e_{i} e_{i}^\top \right)^\top \left(\mI - d e_{i} e_{i}^\top \right)} h^k\\
&=& (d-1)\|h^k\|_{2}^2.
\end{eqnarray*}
The second term can be bounded as
\begin{eqnarray*}
\E{\|b\|_{2}^2} &=&  \E{\left\| d  e_{i}^\top \nabla f(x^k)  e_{i} \right\|_{2}^2}\\
&=& d^2 \sum_{i=1}^d \frac{1}{d} (e_i^\top \nabla f(x^k))^2 \\
&=& d \|\nabla f(x^k)\|_2^2 \\
&=& d \|\nabla f(x^k) - \nabla f(x^*)\|_2^2 \\
&\overset{\eqref{eq:sega_L-smooth_inequality}}{\leq} & 2Ld (f(x^k) - f(x^*)),
\end{eqnarray*}
where in the last step we used $L$-smoothness of $f$. It remains to combine the two bounds.

\end{proof}

\begin{lemma} \label{lem:sega_2_easy} For all $v\in \R^d$ we have
\begin{equation}\label{eq:sega_h_decomp_easy}
\E{\|h^{k+1} \|_2^2} = \left( 1-\frac{1}{d}\right) \|h^k \|_2^2 +\frac{1}{d} \|\nabla f(x^k) - v\|_2^2.
\end{equation}
\end{lemma}
\begin{proof}
We have
\begin{eqnarray*}
\E{\| h^{k+1}\|_2^2} &\overset{\eqref{eq:sega_988fgf}}{=}& \E{ \left\|h^k + e_{i_k}^\top (\nabla f(x^k) - h^k) e_{i_k} \right\|_2^2 }\\
 &=&  \E{ \left\| \left(\mI - e_{i_k} e_{i_k}^\top \right)h^k + e_{i_k} e_{i_k}^\top \nabla f(x^k)  \right\|_2^2 }\\
 &=& \E{ \left\| \left(\mI - e_{i_k} e_{i_k}^\top \right)h^k \right\|_2^2} + \E{\left\|e_{i_k} e_{i_k}^\top \nabla f(x^k)  \right\|_2^2 }\\
 &=& (h^k)^\top \E{\left(\mI - e_{i_k} e_{i_k}^\top \right)^\top \left(\mI - e_{i_k} e_{i_k}^\top \right)} h^k \\
 && \qquad +  (\nabla f(x^k))^\top \E{(e_{i_k} e_{i_k}^\top)^\top e_{i_k} e_{i_k}^\top} \nabla f(x^k)\\
 &=&(h^k)^\top \E{\mI - e_{i_k} e_{i_k}^\top } h^k + (\nabla f(x^k))^\top \E{e_{i_k} e_{i_k}^\top} \nabla f(x^k)\\
 &=& \left(1-\frac{1}{d}\right) \|h^k\|_2^2 + \frac{1}{d} \|\nabla f(x^k)\|_2^2.
\end{eqnarray*}

\end{proof}

\subsection{Proof of Theorem~\ref{thm:sega_simple}}

 We can now write
\begin{eqnarray*}
\E{\|x^{k+1}-x^*\|_{2}^2} &=& \E{\|x^k - \alpha g^k - x^*\|_{2}^2}\\
&=& \|x^k -x^*\|_2^2 + \alpha^2 \E{\|g^k\|_2^2} - 2\alpha \langle \E{g^k}, x^k - x^* \rangle \\
&\overset{\eqref{eq:sega_unbiased_estimator}}{=}& \|x^k -x^*\|_2^2 + \alpha^2 \E{\|g^k\|_2^2} - 2\alpha \langle \nabla f(x^k), x^k - x^* \rangle\\
&\overset{\eqref{eq:sega_8998sgjfbif}}{\leq} & (1-\alpha \mu)\|x^k -x^*\|_2^2  + \alpha^2 \E{\|g^k\|_2^2} - 2\alpha (f(x^k)-f(x^*)).
\end{eqnarray*}
Using Lemma~\ref{lem:sega_1}, we can further estimate 
\begin{eqnarray*}
\E{\|x^{k+1}-x^*\|_2^2} &\leq & (1-\alpha \mu)\|x^k -x^*\|_2^2  \\
&&\qquad  + 2\alpha(2Ld \alpha -1) (f(x^k) - f(x^*)) + 2(d-1)\alpha^2\|h^k \|_2^2 .
\end{eqnarray*}
Let us now add $\sigma \alpha \E{ \|h^{k+1}\|_2^2}$ to both sides of the last inequality. Recalling the definition of the Lyapunov function, and applying Lemma~\ref{lem:sega_2}, we get
\begin{eqnarray*}
\E{\Lgen^{k+1}} &\leq & (1-\alpha \mu)\|x^k -x^*\|_2^2  + 2\alpha(2Ld \alpha -1) (f(x^k) - f(x^*))  \\
&& \quad + 2(d-1)\alpha^2\|h^k \|_2^2 + \sigma \alpha \left(1-\frac{1}{d} \right) \|h^k\|_2^2 +  \frac{\sigma \alpha}{d}\|\nabla f(x^k)\|_2^2\\
&\overset{\eqref{eq:sega_L-smooth_inequality}}{\leq}& (1-\alpha \mu)\|x^k -x^*\|_2^2 +   2\alpha\underbrace{\left(2Ld \alpha + \frac{L\sigma}{d} -1\right)}_{\text{I}} (f(x^k) - f(x^*)) \\
&&\qquad  +\underbrace{ \left( 1-\frac{1}{d} + \frac{2(d-1)\alpha}{\sigma}\right)}_{\text{II}} \sigma \alpha \|h^k\|_2^2.
\end{eqnarray*}
Let us choose $\alpha$ so that $\text{I} \leq 0$ and $\text{II} \leq 1-\alpha \mu$. This leads to the bound \eqref{eq:sega_alpha_bound}. For any $\alpha > 0$ satisfying this bound we therefore have
$
\E{\Lgen^{k+1}} \leq  (1-\alpha \mu) \Lgen^k,
$
as desired. Lastly, as we have freedom to choose $\sigma$, let us pick it so as to maximize the upper bound on the stepsize.

\chapter{Appendix for Chapter~\ref{99}}
\label{99_appendix}

\renewcommand{\EE}{\mathbb{E}}

\graphicspath{{99/experiments/}}

\section{{\tt IBGD}: Bernoulli alternative to {\tt IBCD}}
As an alternative to computing a random block of partial derivatives of size $\tau m$, it is possible to compute the whole gradient with probability $\tau$, and attain the same complexity result. While this can be inserted in all algorithms we propose, we only present an alternative to {\acrshort{IBCD}}, which we call {\tt IBGD}.

\begin{algorithm}[h]
  \caption{Independent Bernoulli Gradient Descent ({\tt IBGD})}
  \label{alg:ibd}
\begin{algorithmic}[1]
\State{\bfseries Input: } {$x^0\in\RR^d$, probability of computing the whole gradient $\tau$, stepsize $\alpha$, \# of parallel units $n$}
  \For{$k=0,1,2,\dotsc$}
    \For{$i=1,\dotsc,n$ in parallel}
        \State Set $g^k_i = \begin{cases} \nabla f_i(x^k) & \text{with probability} \quad \tau \\  0 & \text{with probability} \quad 1-\tau  \end{cases}\quad $ independently 
        \State $x_i^{k+1} = x^k - \alpha  g^k_i$
    \EndFor
    \State $x^{k+1} = \frac{1}{n}\sumin x_i^{k+1}$
  \EndFor
\end{algorithmic}
\end{algorithm}

\begin{theorem}\label{th:ibd} 
Suppose that Assumptions~\ref{as:99_smooth_sc} holds and $\nabla f_i(x^*) = 0$ for all $i$. For Algorithm~\ref{alg:ibd} with $\alpha = \frac{n}{\tau n + 2(1 - \tau)}\frac{1}{2L}$ we have

\[
\EE\left[ \|x^{k} -x^*\|^2 \right] \leq \left(1-\frac{\mu}{2L} \frac{\tau n}{\tau n + 2(1-\tau)}\right)^k\|x^{0} -x^*\|^2.
\]
\end{theorem}

Note that {\tt IBGD} does not perform sparse updates to the server; it is either full (dense), or none. This resembles the most naive asynchronous setup -- where each iteration, a random subset of machines communicates with the server\footnote{In reality, there subset is not drawn from a fixed distribution}. Our findings thus show that we can expect perfect linear scaling for such unreal asynchronous setup. In the honest asynchronous setup, we shall still expect good parallel scaling once the sequence of machines that communicate with server somewhat resembles a fixed uniform distribution.

\section{Asynchronous {\tt ISGD}}\label{sec:99_asynch}
In this section we extend {\tt ISGD} algorithm to the asynchronous setup.  In particular, we revisit the method that was considered by~\cite{grishchenko2018asynchronous}, extend its convergence to stochastic oracle and show better dependency on quantization noise.

\begin{algorithm}[h]
  \caption{Asynchronous {\tt ISGD}}
  \label{alg:asynch_sgd}
\begin{algorithmic}[1]
\State{\bfseries Input: }{$x^0\in\RR^d$, partition of $\RR^d$ into $m$ blocks $u_1,\dotsc, u_m$, ratio of blocks to be sampled $\tau$, stepsize $\alpha$, \# parallel units $n$}
  \For{$k=0,1,2,\dotsc$}
  	\State Worker $i=i_t$ is making update
       \State $w^{t - d_i^k} = \avejn x_j^{t-d_i^k}$
    	   \State $z_i^{k+1} = \proxR(w^{t-d_i^k})$
        \State Sample independently and uniformly a subset of $\tau m$ blocks $U_i^k \subseteq \{u_1, \dotsc, u_m\}$
        \State Sample blocks of stochastic gradient $(g_i^k)_{U_i^k}$  such that $\EE [g_i^k \, |\, x^k] = \nabla f_i(z_i^{k+1})$
        \State $x_i^k = \dotsb = x_i^{t-d_i^k}$
        \State $x_i^{k+1} = x_i^k + \frac{1}{\tau + \frac{1}{n}}(z_i^{k+1} - \alpha  (g_i^k) - x_i^k)_{U_i^k}$
        \State Send $x_i^{k+1} - x_i^k$ and receive $w^{k+1} = w^k + (x_i^{k+1} - x_i^{t})$
  \EndFor
  \State {\bfseries Output:} $x^k = \proxR(w^k)$
\end{algorithmic}
\end{algorithm}

Let us denote the delay of worker $i$ at moment $t$ by $d_i^k$.
\begin{theorem}\label{th:asynch}
	Assume $f_1,\dotsc, f_n$ are $L$-smooth and $\mu$-strongly convex and let Assumption~\ref{as:99_bounded_noise_at_opt} be satisfied. Let us run Algorithm~\ref{alg:asynch_sgd} for $t$ iterations and assume that delays are bounded: $d_i^k\le M$ for any $i$ and $t$. If $\alpha\le \frac{1}{2L(\tau + \frac{2}{n})}$, then
	\begin{align*}
		\EE \|x^k - x^*\|^2
		&\le \left(1 - \tau\alpha\mu\right)^{\lfloor t/M\rfloor}C + 4\alpha\frac{\sigma^2}{n},
	\end{align*}
	where $C\eqdef \max_{i=1,\dotsc, n}\|x^0 - x_i^*\|^2$, $x_i^*\eqdef x^* - \tau\alpha \nabla f_i(x^*)$ and $\lfloor \cdot \rfloor$ is the floor operator.
\end{theorem}
Plugging $\alpha=\frac{1}{2L(\tau + \frac{2}{n})}$ gives complexity that will be significantly improving from increasing $\tau$ until $\tau = \frac{1}{n}$, and then only if $\tau$ jumps from $\frac{1}{n}$ to 1. In contrast, doubling $\tau$ from $\frac{2}{n}$ to $\frac{4}{n}$ would make little difference.

We note that if $\ell_1$ penalty is used, in practice $z_i^k$ should be rather computed on the parameter server side because it will sparsify the vector for communication back.

\section{Proofs for Section~\ref{sec:99_basic}}

\subsection{Key techniques}

The most important equality used many times to prove the results of this chapter is a simple decomposition of expected distances into the distance of expectation and variance:
\begin{align}
    \EE \|X - a\|^2 = \|\EE X - a\|^2 + \EE \|X - \EE X\|^2,\label{eq:99_variance_decompos}
\end{align}
where $X$ is any random vector with finite variance and $a$ is an arbitrary vector from $\RR^d$. 

As almost every algorithm we propose average all updates coming from workers, it will be useful to bound the expected distance of mean of $n$ random variables from the optimum. Lemma~\ref{lem:99_expected_distance} provides the result.  

\begin{lemma}\label{lem:99_expected_distance}
    Suppose that $x^{k+1} = \frac{1}{n}\sumin x_i^k$. Then, we have
    \begin{align*}
        \EE \|x^{k+1} - x^*\|^2
        \le & \left\|\frac{1}{n}\sumin \EE x_i^{k+1} - x^* \right\|^2 
         + \frac{1}{n^2}\sumin \EE \|x_i^{k+1} - \EE x_i^{k+1}\|^2.
    \end{align*}
\end{lemma}
\begin{proof}
    First of all, we have
    \begin{align*}
        \|x^{t + 1} - x^*\|
        = \left\|\frac{1}{n}\sumin x_i^k - (x^* - \alpha \nabla f(x^*)) \right\|.
    \end{align*}
    Now let us proceed to expectations. Note that for any random vector $X$ we have
    \begin{align*}
        \EE \|X\|^2  = \|\EE X\|^2 + \EE \|X - \EE X\|^2.
    \end{align*}
    Applying this to random vector $X \eqdef \frac{1}{n}\sumin x_i^k - x^* $, we get
    \begin{align*}
        \EE \left\|\frac{1}{n}\sumin x_i^k - x^* \right\|^2 = \left\|\EE \frac{1}{n}\sumin x_i^k - x^* \right\|^2 + \EE \left\|\frac{1}{n}\sumin x_i^k - \EE \frac{1}{n}\sumin x_i^k\right\|^2.
    \end{align*}
    In addition, in all minibatching schemes  $x_i^{k+1}$ are conditionally independent given $x^k$. Therefore, for the variance term we get 
    \begin{align}
    \EE \left \|\frac{1}{n}\sumin x_i^k - \EE \frac{1}{n}\sumin x_i^k \right \|^2 = \frac{1}{n^2}\sumin \EE \|x_i^{k+1} - \EE x_i^{k+1}\|^2. \label{eq:99_variance_decomposition}
    \end{align} 
    Plugging it into our previous bounds concludes the proof.
\end{proof}

\subsection{Proof of Theorem~\ref{th:cd}}
\begin{proof}
From Lemma~\ref{lem:99_sgd_recur}, using $\sigma=0$ and $\nabla f_i(x^*)=0$, we immediately obtain
\[
\EE\left[ \|x^{k+1} -x^*\|^2 \,\|\, x^k \right]\leq  (1-\mu \alpha \tau)\|x^{k} -x^*\|^2=\left(1-\frac{\mu}{2L} \frac{\tau n}{\tau n + 2(1-\tau)}\right)\|x^{k} -x^*\|^2 .
\]
It remains to apply the above inequality recursively. 
\end{proof}

\subsection{Proof of Theorem~\ref{th:ibd}}
\begin{proof}
    Clearly, 
    \begin{align*}
        \EE x_i^{k+1}
        = x^k - \alpha\tau\nabla f_i(x^k).
    \end{align*} 
        Let us now elaborate on the second moments.

    Thus,
\[
        \EE \|x_i^{k+1} - \EE x_i^{k+1}\|^2 
        = \alpha^2\EE \left[\| g_i^k - \tau \nabla f_i(x^k)\|^2 \right] 
        = \tau (1-\tau) \| \nabla f_i(x)\|^2.
\]
 Therefore, the conditional variance of $x^{k+1}$ variance is equal to
    \begin{align*}
        \EE \|x^{k+1} - \EE x^{k+1}\|^2 
        = \frac{1}{n^2}\sumin \EE \|x_i^{k+1} - x_i^{k+1}\|^2  = \frac{\tau (1-\tau)}{n^2}\sumin  \| \nabla f_i(x)\|^2.
    \end{align*}
Note that the above equality is exactly~\eqref{eq:99_x_moments} with $\sigma=0$. Thus, one can use Lemma~\ref{lem:99_sgd_recur} (with using $\sigma=0$ and $\nabla f_i(x^*)=0$) obtaining
\[
\EE\left[ \|x^{k+1} -x^*\|^2 \,\|\, x^k \right]\leq  (1-\mu \alpha \tau)\|x^{k} -x^*\|^2=\left(1-\frac{\mu}{2L} \frac{\tau n}{\tau n + 2(1-\tau)}\right)\|x^{k} -x^*\|^2 .
\]
It remains to apply the above inequality recursively. 
\end{proof}

\section{Missing parts from Sections~\ref{sec:99_saga} and~\ref{sec:99_saga_dist}}

\subsection{Useful lemmata}

Let us start with a variance bound, which will be useful for both Algorithm~\ref{alg:saga_dist} and Algorithm~\ref{alg:saga}. 
Define $\Phi(x) = \frac1\qwerty \sum_{i=1}^\qwerty \phi_i(x)$, and define $x^+ =x - \alpha (\nabla f_j(x)-\mJ_j+\bar{\mJ})_U$ for (uniformly) randomly chosen index $1\leq j \leq \qwerty$ and subset of blocks $U$ of size $\tau m $. Define also $\bar{\mJ} = \frac1\qwerty \sum_{i=1}^\qwerty \mJ_i$.

\begin{lemma}[Variance bound]\label{lem:99_saga_variance}
Assume $\phi$ is $\mu$-strongly convex and $\phi_j$ is $L$-smooth and convex for all $j$. Suppose that $x^* = \argmin \Phi(x)$. Then, for any $x$ we have
    \begin{align}\label{eq:99_saga_variance}
        \EE \|x^{+} - \EE x^{+}\|^2
        \le 2\alpha^2\tau\left(2L(\phi(x) - \phi(x^*) + \frac1\qwerty\sum_{j=1}^\qwerty \|\mJ_j - \nabla \phi_j(x^*)\|^2 \right).
    \end{align}
\end{lemma}
\begin{proof}
Since $x^{+} = x - \alpha (\nabla \phi_{j}(x) - \mJ_j + \overline \mJ)_{U}$ and $\EE x^{+} = x - \alpha \tau \nabla \phi(x)$, we get
    \begin{align*}
       &  \EE \|x^{+} - \EE x^{+}\|^2 \\
        &= \alpha^2 \EE \left\|\tau \nabla \phi(x) -(\nabla \phi_{j}(x) - \mJ_j + \overline \mJ)_{U} \right\|^2 \\
        &= \alpha^2 \EE  \|(\tau \nabla \phi(x) - (\nabla \phi_{j}(x) - \mJ_j + \overline \mJ))_{U}\|^2
         + \alpha^2 \EE \|\tau\nabla \phi(x)- (\tau \nabla \phi(x))_{U}\|^2 \\
        &= \alpha^2 \tau \EE \|\tau \nabla \phi(x) - (\nabla \phi_{j}(x) - \mJ_j + \overline \mJ)\|^2
         + \alpha^2(1 - \tau)\tau^2\|\nabla \phi(x)\|^2.
    \end{align*}
    We will leave the second term as is for now and obtain a bound for the first one. Note that the expression inside the norm is now biased: $\EE[\tau \nabla \phi(x) - (\nabla \phi_{j}(x) - \mJ_j + \overline \mJ)]=(\tau - 1)\nabla \phi(x)$. Therefore,
    {
    \footnotesize
    \begin{align*}
        \EE \|\tau \nabla \phi(x) - (\nabla \phi_{j}(x) - \mJ_j + \overline \mJ)\|^2
        = (1 - \tau)^2 \|\nabla \phi(x)\|^2 + \EE\|\nabla \phi(x) - (\nabla \phi_{j}(x) - \mJ_j + \overline \mJ)\|^2.
    \end{align*}
    }
    Now, since $\nabla \phi_{j}(x)$ and $\mJ_j$ are not independent, we shall decouple them using inequality $\|a + b\|^2 \le 2\|a\|^2 + 2\|b\|^2$. In particular,
    \begin{align*}
      &  \EE \|\nabla \phi(x) - (\nabla \phi_{j}(x) - \mJ_j + \overline \mJ)\|^2 \\
        &= \EE \|\nabla \phi(x) - \nabla \phi_{j}(x) + \nabla \phi_{j}(x^*) - \nabla \phi_{j}(x^*) + \mJ_j - \overline \mJ\|^2\\
        &\le 2\EE \|\nabla \phi(x) - \nabla \phi_{j}(x) + \nabla \phi_{j}(x^*)\|^2 + 2\EE\|\mJ_j - \nabla \phi_{j}(x^*) - \overline \mJ\|^2.
    \end{align*}
    Both terms can be simplified by expanding the squares. For the first one we have:
    \begin{align*}
        \EE \|\nabla \phi(x) - \nabla \phi_{j}(x) + \nabla \phi_{j}(x^*)\|^2
        &= \|\nabla \phi(x)\|^2 - 2\< \nabla \phi(x), \EE\left[\nabla \phi_{j}(x) - \nabla \phi_{j}(x^*) \right]> \\
        &\quad + \EE \|\nabla \phi_{j}(x) - \nabla \phi_{j}(x^*)\|^2 \\
        &= -\|\nabla \phi(x)\|^2 + \frac1\qwerty\sum_{j=1}^\qwerty \|\nabla \phi_j(x) - \nabla \phi_j(x^*)\|^2.
    \end{align*}
    Similarly,
    \begin{align*}
        \EE\|\mJ_j - \nabla \phi_{j}(x^*) - \overline \mJ\|^2
        &= \frac1\qwerty\sum_{j=1}^\qwerty\|\mJ_j - \nabla \phi_{j}(x^*)\|^2 - 2\EE\<\mJ_j - \nabla \phi_{j}(x^*), \overline \mJ> + \|\overline \mJ\|^2 \\
        &=\frac1\qwerty\sum_{j=1}^\qwerty\|\mJ_j - \nabla \phi_{j}(x^*)\|^2 - \|\overline \mJ\|^2 \\
        &\le \frac1\qwerty\sum_{j=1}^\qwerty\|\mJ_j - \nabla \phi_{j}(x^*)\|^2.
    \end{align*}
    Coming back to the first bound that we obtained for this lemma, we deduce
    \begin{align*}
       & \EE \|x^{+} - \EE x^{+}\|^2 \\
        &\le \alpha^2\tau\left((1 - \tau)^2\|\nabla \phi(x)\|^2 - 2\|\nabla \phi(x)\|^2 + \frac{2}{k}\sum_{j=1}^\qwerty\|\nabla \phi_j(x) - \nabla \phi_j(x^*)\|^2  \right)\\
        &\quad +\alpha^2\tau\frac{2}{k}\sum_{j=1}^\qwerty \|\mJ_j - \nabla \phi_j(x^*)\|^2 + \alpha^2(1 - \tau)\tau^2\|\nabla \phi(x)\|^2.
    \end{align*}
    The coefficient before $\|\nabla \phi(x)\|^2$ is equal to $\alpha^2\tau((1 - \tau)^2 - 2 + (1 - \tau)\tau) = \alpha^2\tau(1 - \tau - 2) < 0$, so we can drop this term. By smoothness of each $\phi_j$,
    \begin{eqnarray}\nonumber
        \frac1\qwerty\sum_{j=1}^\qwerty \|\nabla \phi_j(x) - \nabla \phi(x^*)\|^2
        &\le& \frac{2L}{\qwerty}\sum_{j=1}^\qwerty \left(\phi_j(x) - \phi(x^*) - \<\nabla \phi_j(x^*), x - x^*> \right) \\
        &=& 2L (\phi(x) - \phi(x^*)), \label{eq:99_saga_grad_sum}
    \end{eqnarray}
    where in the last step we used $\frac1\qwerty\sum_{j=1}^\qwerty\nabla \phi_j(x^*)=0$.
\end{proof}

\begin{lemma}\label{lem:99_saga_grads}
For {\tt ISAGA} with shared data we have (given the setting from Theorem~\ref{th:saga_shared})
    \begin{align*}
        \EE \left[\sum_{j=1}^N \|\mJ_j^{k+1} - \nabla f_j(x^*)\|^2 \, \Big| \, x^k \right]
        \le 2\tau Ln (f(x^k) - f(x^*)) + \left(1 - \frac{\tau n}{N}\right) \sum_{j=1}^N \|\mJ_j^{k} - \nabla f_j(x^*)\|^2.
    \end{align*}
    On the other hand, for distributed {\tt ISAGA} we have for all $i$ (given the setting from Theorem~\ref{th:saga_dist}):
        \begin{align*}
        \EE \left[\sum_{j=1}^{l} \|\mJ_{ij}^{k+1} - \nabla f_{ij}(x^*)\|^2 \, \Big| \, x^k  \right]
        \le 2\tau L (f_i(x^k) - f_i(x^*)) + \left(1 - \frac{\tau }{l}\right) \sum_{j=1}^{l} \|\mJ_{ij}^{k} - \nabla f_{ij}(x^*)\|^2.
    \end{align*}
\end{lemma}
\begin{proof}
Consider all expectations throughout this proof to be conditioned on $x^k$. 
Let $j^k$ be the function index used to obtain $x_i^{k+1}$ from $x^k$. Then we have $(\mJ_{j^k}^{k+1})_{U_i^k} = (\nabla f_{j^k}(x^k))_{U_i^k}$. In the rest of the blocks, $\mJ_{j^k}^{k+1}$ coincides with its previous value. This implies
    \begin{align}
        \EE \left[\|\mJ_{j^k}^{k+1} - \nabla f_{j^k}(x^*)\|^2\mid j^k \right]
        = \tau \|\nabla f_{j^k}(x^k) - \nabla f_{j^k}(x^*)\|^2 + (1 - \tau) \|\mJ_{j^k}^{k} - \nabla f_{j^k}(x^*)\|^2. \label{eq:99_saga_134278238213698}
    \end{align}
    Taking expectation with respect to sampling of $j^k$ we obtain for shared data setup:
    \begin{eqnarray*}
        && \EE \left[\sum_{j=1}^N \|\mJ_j^{k+1} - \nabla f_j(x^*)\|^2 \right]
        \\
        & =& 
        \sum_{j=1}^N \EE \left[ \|\mJ_j^{k+1} - \nabla f_j(x^*)\|^2 \right] \\ 
       & \stackrel{\eqref{eq:99_saga_134278238213698}}{=}& \sum_{j=1}^N \frac{n}{N} \left( \tau   \|\nabla f_j(x^k) - \nabla f_j(x^*)\|^2 + \left(1 - \tau \right)  \|\mJ_j^{k} - \nabla f_j(x^*)\|^2 \right) \\
       && \qquad +
 \sum_{j=1}^N \left(1-\frac{n}{N}\right)    \left(   \left(1 - \tau \right)  \|\mJ_j^{k} - \nabla f_j(x^*)\|^2 \right)\\
       & =& \tau \frac{n}{N} \sum_{j=1}^N \|\nabla f_j(x^k) - \nabla f_j(x^*)\|^2 + \left(1 - \frac{\tau n}{N}\right) \sum_{j=1}^N \|\mJ_j^{k} - \nabla f_j(x^*)\|^2.
    \end{eqnarray*}
    
Similarly, for distributed setup we get   
{
\footnotesize  
   \begin{align*}
        \EE \left[\sum_{j=1}^{l} \|\mJ_{ij}^{k+1} - \nabla f_{ij}(x^*)\|^2 \right]
        \le \tau \frac{1}{l} \sum_{j=1}^l \|\nabla f_{ij}(x^k) - \nabla f_{ij}(x^*)\|^2+ \left(1 - \frac{\tau }{l}\right) \sum_{j=1}^{l} \|\mJ_{ij}^{k} - \nabla f_{ij}(x^*)\|^2
    \end{align*}
    }
    Using~\eqref{eq:99_saga_grad_sum}, the first sum of right hand side can be bounded by $2LN(f(x^k) - f(x^*))$ or $2Ll(f(x^k) - f(x^*))$. 
\end{proof}
\subsection{Proof of Theorem~\ref{th:saga_dist}}
\begin{proof}
    First of all, let us verify that it indeed holds $c>0$ and $\rho\ge 0$. As $\alpha \le \frac{1}{L\left(\frac{3}{n} + \tau\right)}$, we have $c = \frac{1}{n}\left( \frac{1}{\alpha L} - \frac{1}{n} - \tau\right) \ge \frac{1}{n}\left( \frac{3}{n} + \tau - \frac{1}{n} - \tau\right) > 0$. Furthermore, $\alpha\mu \ge 0$, so to show $\rho\ge 0$, it is enough to mention that $$\frac{1}{l} - \frac{2}{n^2lc} = \frac{1}{l} - \frac{2}{n^2l\left(\frac{1}{\alpha L} - \frac{1}{n} - \tau \right)} \ge \frac{1}{l} - \frac{2}{n^2l\left(\frac{3}{n} + \tau - \frac{1}{n} - \tau \right)} = 0.$$
    
    Now we proceed to the proof of convergence. We are going to decompose the expected distance from $x^{k+1}$ to $x^*$ into its variance and the distance of expected iterates, so let us analyze them separately. The variance can be bounded as follows:
    \begin{eqnarray} 
            & &    \EE \left[\|x^{k+1} - \EE [x^{k+1} \mid x^k]\|^2\mid x^k \right]\nonumber \\
       & &   \stackrel{\eqref{eq:99_variance_decomposition}+\eqref{eq:99_saga_variance}}{\le} 
2\frac{\alpha^2\tau}{n}\left(2L(f(x) - f(x^*) + \frac{1}{ln}\sum_{j=1}^n\sum_{j=1}^l \|\mJ_{ij} - \nabla f_{ij}(x^*)\|^2 \right). \label{eq:99_saga_distrib_variance}
    \end{eqnarray}
    For the distance of the expected iterates we write
    \begin{align*}
        \|\EE [x^{k+1} \mid x^k] - x^*\|^2
        &= \|x^k - \alpha\tau \nabla f(x^k) - x^*\|^2 \\
        &\le (1 - \alpha\tau\mu)\|x^k - x^*\|^2 - 2\alpha\tau (f(x^k) - f(x^*)) \\
        & \qquad + 2\alpha^2\tau^2 L(f(x^k) - f(x^*)).
    \end{align*}
    As is usually done for {\tt SAGA}, we are going to prove convergence using a Lyapunov function. Namely, let us define
    \begin{align}\label{eq:99_saga_lyapunov}
        \cL^k \eqdef \EE \left[\|x^k - x^*\|^2 + c \alpha^2\sum_{i=1}^n\sum_{j=1}^l \|\mJ_{ij}^k - \nabla f_{ij}(x^*)\|^2 \right],
    \end{align}
    where $c =  \frac{1}{n}\left( \frac{1}{\alpha L} - \frac{1}{n} - \tau\right)$. Using Lemma~\ref{lem:99_saga_grads} together with the bounds above, we get
    \begin{align*}
        \cL^{k+1} 
        &\le \EE \left[(1 - \alpha\tau\mu)\|x^k - x^*\|^2 + \left(\frac{2\alpha^2\tau}{n^2l} + c\alpha^2\left(1 - \frac{\tau }{l}\right)\right)\sum_{i=1}^n\sum_{j=1}^l \|\mJ_{ij}^k - \nabla f_{ij}(x^*)\|^2 \right] \\
        &\quad + 2\alpha\tau\EE\Bigl[\Bigl(\underbrace{\alpha\tau L + \frac{\alpha L}{n} + c\alpha L n - 1}_{= 0\text{ by our choice of } c } \Bigr) (f(x^k) - f(x^*)) \Bigr].
    \end{align*}
    In fact, we chose $c$ exactly to make the last expression equal to zero.
    After dropping it, we reduce the bound to
    \begin{align*}
        \cL^{k+1} 
        \le (1 - \rho)\EE \left[ \|x^k - x^*\|^2 + c\alpha^2\sum_{i=1}^n\sum_{j=1}^l \|\mJ_{ij}^k - \nabla f_{ij}(x^*)\|^2 \right]
        = (1 - \rho)\cL^{k},
    \end{align*}
    where $\rho = \min\left\{\alpha\tau\mu, \frac{\tau }{l} - \frac{2\tau}{n^2lc} \right\}$. Note that $\EE \|x^k - x^*\|^2\le \cL^k \le (1 - \rho)^k\cL^0$ by induction, so we have the stated linear rate.
\end{proof}
\subsection{Proof of Theorem~\ref{th:saga_shared}}
\begin{proof}
    First of all, let us verify that it indeed holds $c>0$ and $\rho\ge 0$. As $\alpha \le \frac{1}{L\left(\frac{3}{n} + \tau\right)}$, we have $c = \frac{1}{n}\left( \frac{1}{\alpha L} - \frac{1}{n} - \tau\right) \ge \frac{1}{n}\left( \frac{3}{n} + \tau - \frac{1}{n} - \tau\right) > 0$. Furthermore, $\alpha\mu \ge 0$, so to show $\rho\ge 0$, it is enough to mention that $$\frac{n}{N} - \frac{2}{nNc} = \frac{n}{N} - \frac{2}{N\left(\frac{1}{\alpha L} - \frac{1}{n} - \tau \right)} \ge \frac{n}{N} - \frac{2}{N\left(\frac{3}{n} + \tau - \frac{1}{n} - \tau \right)} = 0.$$
    
    Now we proceed to the proof of convergence. We are going to decompose the expected distance from $x^{k+1}$ to $x^*$ into its variance and the distance of expected iterates, so let us analyze them separately. The variance term can be bounded as follows
    \begin{eqnarray} \label{eq:99_saga_shared_variance}
       && \EE \left[\|x^{k+1} - \EE [x^{k+1} \mid x^k]\|^2\mid x^k \right] \nonumber \\
        && \qquad \stackrel{\eqref{eq:99_variance_decomposition}+\eqref{eq:99_saga_variance}}{\le} 
        \frac{2\alpha^2\tau}{n}\left(2L(f(x^k) - f(x^*) + \frac{1}{N}\sum_{j=1}^N \|\mJ_j^k - \nabla f_j(x^*)\|^2 \right) 
    \end{eqnarray}
    For the distance of the expected iterates we write
    {
    \footnotesize
    \begin{align*}
        \|\EE [x^{k+1} \mid x^k] - x^*\|^2
        &= \|x^k - \alpha\tau \nabla f(x^k) - x^*\|^2 \\
        &\le (1 - \alpha\tau\mu)\|x^k - x^*\|^2 - 2\alpha\tau (f(x^k) - f(x^*)) + 2\alpha^2\tau^2 L(f(x^k) - f(x^*)).
    \end{align*}
    }
    As is usually done for {\tt SAGA}, we are going to prove convergence using a Lyapunov function. Namely, let us define
    \begin{align*}
        \cL^k \eqdef \EE \left[\|x^k - x^*\|^2 + c \alpha^2\sum_{j=1}^N \|\mJ_j^k - \nabla f_j(x^*)\|^2 \right],
    \end{align*}
    where $c =  \frac{1}{n}\left( \frac{1}{\alpha L} - \frac{1}{n} - \tau\right)$. Using Lemma~\ref{lem:99_saga_grads} together with the bounds above, we get
    \begin{align*}
        \cL^{k+1} 
        &\le \EE \left[(1 - \alpha\tau\mu)\|x^k - x^*\|^2 + \left(\frac{2\alpha^2\tau}{nN} + c\alpha^2\left(1 - \frac{\tau n}{N}\right)\right)\sum_{j=1}^N\|\mJ_j^k - \nabla f_j(x^*)\|^2 \right] \\
        &\quad + 2\alpha\tau\EE\Bigl[\Bigl(\underbrace{\alpha\tau L + \frac{\alpha L}{n} + c\alpha L n - 1}_{= 0\text{ by our choice of } c } \Bigr) (f(x^k) - f(x^*)) \Bigr].
    \end{align*}
    In fact, we chose $c$ exactly to make the last expression equal to zero.
    After dropping it, we reduce the bound to
    \begin{align*}
        \cL^{k+1} 
        \le (1 - \rho)\EE \left[ \|x^k - x^*\|^2 + c\alpha^2\sum_{j=1}^N\|\mJ_j^k - \nabla f_j(x^*)\|^2\right]
        = (1 - \rho)\cL^{k},
    \end{align*}
    where $\rho = \min\left\{\alpha\tau\mu, \frac{\tau n}{N} - \frac{2\tau}{nNc} \right\}$. Note that $\EE \|x^k - x^*\|^2\le \cL^k \le (1 - \rho)^k\cL^0$ by induction, so we have the stated linear rate.
\end{proof}

\section{Proofs for Section~\ref{sec:99_sgd}}
\subsection{Useful lemmas}

The next lemma is a key technical tool to analyze Algorithm~\ref{alg:sgd}. It provides a better expression for first and second moments of algorithm iterates.

\begin{lemma}[{\tt SGD} moments]\label{lem:99_moments}
    Consider the randomness of the update of Algorithm~\ref{alg:sgd} at moment $t$. The first moments of the generated iterates are simply $\EE x_i^{k+1} = x^k - \alpha\tau \nabla f_i(x^k)$ and $\EE x^{k+1} = x^k - \alpha\tau \nabla f(x^k)$, while their second moments are:
    \begin{align}
        \EE\left[ \|x_i^{k+1} - \EE x_i^{k+1}\|^2 \, | \, x^k\right]&= \alpha^2 \tau\Big(\left(1 - \tau \right) \|\nabla f_i(x^k)\|^2  + \EE \|g_i^k - \nabla f_i(x^k)\|^2 \Big), \label{eq:99_x_i_moments}\\
        \EE\left[\|x^{k+1} - \EE x^{k+1}\|^2 \, | \, x^k\right] &= \alpha^2 \frac{\tau}{n^2}\sumin \Big((1 - \tau)\|\nabla f_i(x^k)\|^2  + \EE \|g_i^k - \nabla f_i(x^k)\|^2 \Big) \label{eq:99_x_moments}.
    \end{align}
\end{lemma}

\begin{proof}
    Clearly, 
    \begin{align*}
        \EE x_i^{k+1}
        = x^k - \alpha \EE \left[(g_i^k)_{U_i^k} \right]
        = x^k - \alpha \EE \left[ \left( \nabla f_i(x^k)\right)_{U_i^k} \right]
        = x^k - \alpha\tau\nabla f_i(x^k)
    \end{align*}
    and, therefore, $\EE x^{k+1} = x^k - \alpha\tau\nabla f(x^k)$.
    Let us now elaborate on the second moments. Using the obtained formula for $\EE x_i^{k+1}$, we get $(x_i^{k+1} - \EE x_i^{k+1})_{U_i^k} = -\alpha (g_i^k - \tau \nabla f_i(x^k))_{U_i^k}$ and $(x_i^{k+1} - \EE x_i^{k+1})_{\bar{U}_i^k} = \alpha\tau\left( \nabla f_i(x^k)\right)_{\bar{U}_i^k}$ where $\bar{U}_i^k$ is a set of blocks not contained in $U_i^k$. Thus,
    \begin{align*}
        \EE \|x_i^{k+1} - \EE x_i^{k+1}\|^2 
        &= \alpha^2\EE \left[\| (g_i^k - \tau \nabla f_i(x^k))_{U_i^k}\|^2 + \tau^2\|\left( \nabla f_i(x^k)\right)_{\bar{U}_i^k}\|^2 \right] \\
        &= \alpha^2 \left(\tau\EE\| g_i^k - \tau \nabla f_i(x^k)\|^2  + \tau^2\left(1 - \tau\right)   \|\nabla f_i(x^k)\|^2\right).
    \end{align*}
    Note that $\EE g_i^k - \tau \nabla f_i(x^k) = (1 - \tau)\nabla f_i(x^k)$, so we can use decomposition~\eqref{eq:99_variance_decompos} to write $\EE\| g_i^k - \tau \nabla f_i(x^k)\|^2 = (1 - \tau)^2\|\nabla f_i(x^k)\|^2 + \EE \|g_i^k - \nabla f_i(x^k)\|^2$. This develops our previous statement into
    {
    \footnotesize
    \begin{align*}
        \EE \|x_i^{k+1} - \EE x_i^{k+1}\|^2 
        & = \alpha^2 \left(\tau\left((1 - \tau)^2\|\nabla f_i(x^k)\|^2 + \EE \|g_i^k - \nabla f_i(x^k)\|^2\right)  + \tau^2\left(1 - \tau\right)   \|\nabla f_i(x^k)\|^2\right)\\
        &= \alpha^2 \tau\left((1 - \tau)\|\nabla f_i(x^k)\|^2 + \EE \|g_i^k - \nabla f_i(x^k)\|^2 \right),
    \end{align*}
    }
    which coincides with what we wanted to prove for $x_i^{k+1}$. As for $x^{k+1}$, it is merely the average of independent random variables conditioned on $x^k$. Therefore, its variance is equal to
    \begin{align*}
        \EE \|x^{k+1} - \EE x^{k+1}\|^2 
        = \frac{1}{n^2}\sumin \EE \|x_i^{k+1} - x_i^{k+1}\|^2 .
    \end{align*}
    This concludes the proof.
\end{proof}

\begin{lemma}\label{lem:99_expected_squared_grad}
    Let $f_i$ be $L$-smooth and convex for all $i$. Then, 
    \begin{align}
        \frac{1}{n}\sumin \|\nabla f_i(x^k)\|^2
        \le 4L(f(x^k) - f(x^*)) + \frac{2}{n}\sumin \|\nabla f_i(x^*)\|^2. \label{eq:99_lemma3}
    \end{align}
\end{lemma}
\begin{proof}
    If $\nabla f_i(x^*)=0$ for all $i$, we can simply write $\|\nabla f_i(x^k)\|^2 = \|\nabla f_i(x^k) - \nabla f_i(x^*)\|^2 \leq 2L(f(x^k) - f(x^*) - \<\nabla f_i(x^*), x^k - x^*>) = 2L (f(x^k) - f(x^*))$. 
    Otherwise, we have to use inequality $\|a+b\|^2 \le 2\|a\|^2 + 2\|b\|^2$ with $a = \nabla f_i(x^k) - \nabla f_i(x^*)$ and $b = \nabla f_i(x^*)$. We get
    \begin{align*}
        \sumin \|\nabla f_i(x^k)\|^2
        &\le 2 \sumin \|\nabla f_i(x^k) - \nabla f_i(x^*)\|^2 + 2\sumin\|\nabla f_i(x^*)\|^2 \\
        &\le 4L \sumin (f_i(x^k) - f_i(x^*) - \<\nabla f_i(x^*), x^k - x^*>) + 2\sumin\|\nabla f_i(x^*)\|^2\\
        &= 4Ln (f(x^k) - f(x^*)) + 2\sumin\|\nabla f_i(x^*)\|^2.
    \end{align*}
\end{proof}

\begin{lemma}\label{lem:99_second_momentum_of_stoch_grad}
	Let $f= \EE f(\cdot ; \xi)$ be $\mu$-strongly and $f(\cdot; \xi)$ be $L$-smooth and convex almost surely. Then, for any $x$ and $y$
	\begin{align*}
		\EE \|\nabla f(x; \xi)\|^2
		\le 4L(f(x) - f(y) - \<\nabla f(y), x - y>) + 2\EE \|\nabla f(y; \xi)\|^2.
	\end{align*}
\end{lemma}
\begin{proof}
	The proof proceeds exactly the same way as that of Lemma~\ref{lem:99_expected_squared_grad}.
\end{proof}

\begin{lemma}\label{lem:99_sgd_recur}
  Suppose that Assumption~\ref{as:99_smooth_sc} holds. Then, if we have
  {
  \footnotesize
    \begin{align*}
       & 2\alpha\tau\left(1 - \alpha\tau L - \frac{2\alpha L(1 - \tau)}{n} \right)\EE[f(x^k) - f(x^*)] \\
       & \le (1 - \alpha\tau\mu)\EE\|x^k - x^*\|^2 - \EE \|x^{k+1} - x^*\|^2  + \alpha^2\frac{\tau}{n}\left(\sigma^2 + 2\frac{1 - \tau}{n}\sumin \|\nabla f_i(x^*)\|^2 \right).
    \end{align*}
    }
\end{lemma}

\begin{proof}
    Substituting Assumption~\ref{as:99_bounded_noise} into~\eqref{eq:99_x_moments}, we obtain
    \begin{align}
        \EE \left[\|x^{k+1} - \EE [x^{k+1} \mid x^k]\|^2 \mid x^k\right]
        &\le \alpha^2 \frac{\tau}{n^2}\sumin \left((1 - \tau)\|\nabla f_i(x^k)\|^2 + \sigma^2 \right).\label{eq:99_sgd_varxt}
    \end{align}
    We use it together with decomposition~\eqref{eq:99_variance_decompos} to write
    \begin{align*}
        &\EE \left[\|x^{k+1} - x^*\|^2 \mid x^k \right] \\
        &= \|\EE [x^{k+1} \mid x^k]  - x^*\|^2 + \EE\left[ \|x^{k+1} - \EE [x^{k+1} \mid x^k]\|^2\right]\\
        &\stackrel{\eqref{eq:99_sgd_varxt}}{\le} \|x^k -\alpha\tau \nabla f(x^k) - x^*\|^2 + \alpha^2 \frac{\tau}{n^2}\sumin \left((1 - \tau)\|\nabla f_i(x^k)\|^2 + \sigma^2 \right) \\
        & \stackrel{\eqref{eq:99_lemma3}}{\le} \|x^k -\alpha\tau \nabla f(x^k) - x^*\|^2 \\
        & \qquad + \alpha^2 \frac{\tau}{n}\left((1 - \tau)\left(4L(f(x^k) - f(x^*)) + \frac{2}{n}\sumin \|\nabla f_i(x^*)\|^2\right) + \sigma^2\right).
    \end{align*}
    Let us expand the first square:
    \begin{align*}
        & \|x^k -\alpha\tau \nabla f(x^k) - x^*\|^2 \\
        & \qquad  \qquad = \|x^k - x^*\|^2 - 2 \alpha\tau\< x^k - x^*, \nabla f(x^k)> + \alpha^2 \tau^2\|\nabla f(x^k)\|^2 \\
        & \qquad  \qquad\le \|x^k - x^*\|^2 - 2 \alpha\tau\< x^k - x^*, \nabla f(x^k)> + \alpha^2 \tau 2L(f(x^k) - f(x^*)).
    \end{align*}
    The scalar product gives
    \begin{align*}
         \<\nabla f(x^k), x^k - x^*> 
        \ge f(x^k) - f(x^*) + \frac{\mu}{2}\|x^k - x^*\|^2.
    \end{align*}
    Combining the produced bounds, we show that
    \begin{align*}
        & \EE\left[\|x^{k+1} - x^*\|^2  \mid x^k\right] \\
        &\le \left(1 - \alpha\tau\mu\right)\|x^k - x^*\|^2 + \left(2\alpha^2 \tau L - 2\alpha\tau + \alpha^2(1 - \tau)\frac{4\tau L}{n}\right) (f(x^k) - f(x^*))\\
        &\quad + \alpha^2\frac{\tau}{n}\left(2(1 - \tau)\frac{1}{n}\sumin \|\nabla f_i(x^*)\|^2 + \sigma^2 \right).
    \end{align*}
    This is equivalent to our claim.
\end{proof}

\subsection{Proof of Theorem~\ref{th:sgd}}
\begin{proof}
Only for the purpose of this proof, denote $\alpha_k \eqdef \alpha^k$ in order to not confuse superscript with power. From the choice of $\alpha_k$ we deduce that $2\alpha_k\tau\left(1 - \alpha_k\tau L - \frac{2\alpha_k L(1 - \tau)}{n} \right)\ge \alpha_k\tau$.  Therefore, the result of Lemma~\ref{lem:99_sgd_recur} simplifies to
    \begin{align} \label{eq:99_lem5_consequence}
        \EE[f(x^\qwerty) - f(x^*)]
        \le \frac{1}{\alpha_k\tau}(1 - \alpha_k\tau\mu)\EE \|x^\qwerty - x^*\|^2 - \frac{1}{\alpha_k\tau}\EE \|x^{\qwerty+1} - x^*\|^2 + \alpha_k\frac{E}{n},
    \end{align}
    where $E \eqdef \sigma^2 + (1 - \tau)\frac{2}{n}\sumin\|\nabla f_i(x^*)\|^2$. Dividing~\eqref{eq:99_lem5_consequence} by $\alpha_k$ and summing it for $k=0,\dots, t$ we obtain
     \begin{align*}
        \sum_{k=0}^t \frac{1}{\alpha_k} \EE[f(x^\qwerty) - f(x^*)]
        &\le \frac{1}{\alpha_0^2\tau}(1 - \alpha_0\tau\mu)\|x^0 - x^*\|^2 - \frac{1}{\alpha_k^2\tau}\EE \|x^{k+1} - x^*\|^2 +t \frac{E}{n} \\
        &\quad + \frac1\tau \sum_{k=1}^{t-1}\left( \frac{1}{\alpha_k^2} \left( 1-\alpha_k\tau\mu\right) - \frac{1}{\alpha_{k-1}^2}\right)\EE\|x^\qwerty - x^*\|^2.\\
    \end{align*}
    Next, notice that 
             \begin{align*}
\frac{1}{\alpha_{k}^2} - \frac{1}{\alpha_{k+1}^2} \left( 1-\alpha_{k+1}\tau\mu\right) &= 
 \frac{1}{\alpha_{k}^2} \left(1-  \frac{\alpha_{k}^2}{\alpha_{k+1}^2} \left( 1-\alpha_{k+1}\tau\mu\right)\right) \\
 &=
  \frac{1}{\alpha_{k}^2} \left(1- \left(1+ \frac{c}{a+ck} \right)^2  \left( 1-\frac{\tau\mu}{a+c(k+1)}\right)\right)
  \\
 &\stackrel{ (*)}{\geq} 
   \frac{1}{\alpha_{k}^2} \left(1- \left(1+ \frac{2.125\,c}{a+ck} \right) \left( 1-\frac{\tau\mu}{a+c(k+1)}\right)\right)
   \\
   &=
     \frac{1}{\alpha_{k}^2} \left(1- \left(1+ \frac{2.125}{4} \frac{1  }{\frac{a}{\tau \mu }+  \frac14 k} \right) \left( 1-\frac{1}{\frac{a}{\tau \mu }+ \frac14 k +\frac14}\right)\right) 
       \\
          &\stackrel{(**)}{\geq }0.
    \end{align*}
    Above $(*)$ holds since $\frac{c}{a+ck}\leq \frac18$ and $(1+\epsilon)^2\leq (1+2.125\epsilon)$ for $\epsilon \leq \frac18$. Next, inequality $(**)$ holds since function $\varphi(y) = (1+\frac{2.125}{4y})\left(1-\frac{1}{y+\frac14}\right)$ is upper bounded by 1 on $[0,\infty)$. Thus, we have

       \begin{align*}
        \sum_{\qwerty=0}^k \frac{1}{\alpha_k} \EE[f(x^\qwerty) - f(x^*)]
\le \frac{a^2}{\tau}\left(1 - \frac{\tau\mu}{a}\right)\|x^0 - x^*\|^2 +t \frac{E}{n}.
    \end{align*}
    
    All that remains is to mention that by Jensen's inequality $\EE f(\hat x^k) \le \frac{1}{(k+1)a+ \frac{c}{2}k(k+1)}\sum_{k=0}^k(a+ck)\EE f(x^\qwerty)=  \frac{1}{\sum_{k=0}^k\alpha_k^{-1}}\sum_{\qwerty=0}^k\alpha_k^{-1} \EE f(x^\qwerty) $.

\end{proof}

\subsection{Proof of Theorem~\ref{th:sgd_ncvx}}
It will be useful to establish a technical lemma first.
\begin{lemma}\label{lem:99_func_impr}
    Let $f$ be $L$-smooth and assume that $\frac{1}{n}\sumin \|\nabla f_i(x) - \nabla f(x)\|^2 \le \nu^2$ for all $x$. Then, considering only randomness from iteration $t$ of Algorithm~\ref{alg:sgd},
    \begin{align*}
        \EE f(x^{k+1})
        \le f(x^k) - \alpha\tau\left(1 - \frac{\alpha\tau L }{2} - \alpha L\left(1 - \tau\right)\frac{1}{n}\right)\|\nabla f(x^k)\|^2 + 
       \alpha^2L\tau\frac{\left(1 - \tau\right)\nu^2+\frac12\sigma^2}{n}.
    \end{align*}
\end{lemma}
\begin{proof}
    Using smoothness of $f$ and assuming $x^k$ is fixed, we write
    \begin{align*}
        \EE f(x^{k+1}) 
        &\le f(x^k) + \< \nabla f(x^k), \EE\ x^{k+1} - x^k> + \frac{L}{2}\EE\|x^{k+1} - x^k\|^2 \\
        &= f(x^k) - \alpha\tau \|\nabla f(x^k)\|^2 + \frac{L}{2}\EE\|x^{k+1} - x^k\|^2 .
    \end{align*}
    It holds 
    \begin{eqnarray*}
       &&  \EE \|x^{k+1} - x^k\|^2 \\
        &=& \left\|\EE x^{k+1} - x^k \right\|^2 + \EE\left\|x^{k+1} - \EE\left[x^{k+1} \mid x^k \right]\right\|^2
        \\
        &\overset{\eqref{eq:99_x_moments}}{=} &
        \alpha^2\tau^2  \left\|\nabla f(x^k) \right\|^2 + \alpha^2\tau\frac{1}{n^2}\sumin \left((1 - \tau) \|\nabla f_i(x^k)\|^2 +\EE \|g_i^k - \nabla f_i(x^k)\|^2\right) 
        \\
        &\overset{\text{As}.~\ref{as:99_bounded_noise}}{\leq}&
         \alpha^2\tau^2  \left\|\nabla f(x^k) \right\|^2 + \alpha^2\tau\frac{1}{n^2}\sumin \left((1 - \tau) \|\nabla f_i(x^k)\|^2 +\sigma^2\right).
    \end{eqnarray*}
    Using inequality $\|a+b\|^2 \le \|a\|^2 + \|b\|^2$ with $a= \nabla f_i(x^k) - \nabla f(x^k)$ and $b = \nabla f(x^k)$ yields
    \begin{align*}
        \frac{1}{n}\sumin \|\nabla f_i(x^k)\|^2
        \le \frac{2}{n}\sumin \|\nabla f_i(x^k) - \nabla f(x^k)\|^2 + 2 \|\nabla f(x^k)\|^2
        \le 2\nu^2 + 2 \|\nabla f(x^k)\|^2.
    \end{align*}
    Putting the pieces together, we prove the claim.
\end{proof}
We now proceed with Proof of Theorem~\ref{th:sgd_ncvx}.
    \begin{proof}
    Taking full expectation in Lemma~\ref{lem:99_func_impr} and telescoping this inequality from 0 to $t$, we obtain
    \begin{align*}
        0&\le \EE f(x^{k+1}) - f^* \\ & \le f(x^0) - f^* - \alpha\tau\left(1 - \frac{\alpha\tau L}{2} - \alpha L\left(1 - \tau\right)\frac{1}{n}\right)\sum_{k=0}^k\|\nabla f(x^\qwerty)\|^2  \\
        & \qquad \qquad + t\alpha^2 L\tau\frac{\left(1 - \tau\right)\nu^2+\frac12\sigma^2}{n}.
    \end{align*}
    Rearranging the gradients and dividing by the coefficient before it, we get the result.
    \end{proof}

\section{Missing parts from Section~\ref{sec:99_ABCDE}}

\subsection{Proof of Lemma~\ref{lem:99_stronggrowth}}
\begin{proof}

Let us first bound a variance of $\frac1\tau (g_i)_{U_i}$ -- an unbiased estimate of $\nabla f_i(x)$, as it will appear later in the derivations:
\begin{eqnarray}
& & \EE\left[ \left\|\frac{1}{\tau} (g_i)_{U_i} -\nabla f_i(x)\right \|^2 \right]
\\
&=&
\EE_g\left[\EE_U\left[ \left\|\frac{1}{\tau} (g_i)_{U_i} -\nabla f_i(x)\right \|^2 \right]\right]
\nonumber \\
&=&
\EE_g\left[  (1-\tau) \|\nabla f_i(x) \|^2+ \tau \left\|  \frac1\tau g_i- \nabla f_i(x)\right\|^2   
\nonumber\right] \\
&\stackrel{\eqref{eq:99_variance_decompos}}{=}&
(1-\tau) \|\nabla f_i(x) \|^2+ \tau \left\|  \left( \frac1\tau -1\right) \nabla f_i(x)\right\|^2 +\tau \EE_g\left[\left\|  \frac1\tau (g_i- \nabla f_i(x))\right\|^2 \right]
 \nonumber \\
&=&
(1-\tau) \|\nabla f_i(x) \|^2+ \tau \left( \frac1\tau -1\right)^2\left\|  \nabla f_i(x)\right\|^2+\frac1\tau \left\|  g_i- \nabla f_i(x)\right\|^2 \nonumber \\
&\stackrel{\eqref{eq:99_acc_sg_fi}}{\leq}&
(1-\tau) \|\nabla f_i(x) \|^2+ \tau \left( \frac1\tau -1\right)^2\left\|  \nabla f_i(x)\right\|^2 +\frac{\bar{\rho}}{\tau} \left\|  \nabla f_i(x)\right\|^2 + \frac{\bar{\sigma}^2}{\tau} \nonumber \\
&=&
\left(\frac1\tau-1+\frac{\bar{\rho}}{\tau} \right) \|\nabla f_i(x) \|^2  + \frac{\bar{\sigma}^2}{\tau}.
\label{eq:99_acc_gi_bound}
\end{eqnarray}

Next we proceed with bounding the second moment of gradient estimator: 
\begin{eqnarray*}
\EE\left[\|q\|^2\right] &=& 
\EE\left[ \left \|\frac{1}{n\tau} \sum_{i=1}^n(g_i)_{U_i} \right\|^2\right] 
\\
&\stackrel{\eqref{eq:99_variance_decompos}}{=}&
\|\nabla f(x) \|^2 + \EE\left[ \left\|\frac{1}{n\tau} \sum_{i=1}^n\left((g_i)_{U_i} -\nabla f_i(x)\right)\right\|^2 \right]
\\
&\stackrel{(*)}{=}&
\|\nabla f(x) \|^2 + \frac{1}{n^2} \sum_{i=1}^n\EE\left[ \left\|\frac{1}{\tau} (g_i)_{U_i} -\nabla f_i(x)\right \|^2 \right]
\\
&\stackrel{\eqref{eq:99_acc_gi_bound}}{\leq}&
\|\nabla f(x) \|^2 + \frac{1}{n^2} \sum_{i=1}^n\left( \left(\frac1\tau-1+\frac{\bar{\rho}}{\tau} \right) \|\nabla f_i(x) \|^2  + \frac{\bar{\sigma}^2}{\tau} 
 \right) 
 \\ 
  &=&
\|\nabla f(x) \|^2 + \frac{\bar{\sigma}^2}{n\tau} + \left(\frac1\tau-1+\frac{\bar{\rho}}{\tau} \right) \frac{1}{n^2} \sum_{i=1}^n \|\nabla f_i(x) \|^2
\\ 
  &\stackrel{\eqref{eq:99_acc_sg_f}}{\leq}&
\|\nabla f(x) \|^2 + \frac{\bar{\sigma}^2}{n\tau} + \left(\frac1\tau-1+\frac{\bar{\rho}}{\tau} \right) \frac{1}{n} \left(  \tilde{\rho}\| \nabla f(x)\|^2 + \tilde{\sigma}^2\right)
\\ 
  &=&
\left(1+ \frac{\tilde{\rho}}{n}   \left(\frac1\tau-1+\frac{\bar{\rho}}{\tau} \right) \right)\|\nabla f(x) \|^2 + \frac{\bar{\sigma}^2}{n\tau} + \frac{\tilde{\sigma}^2}{n}\left(\frac1\tau-1+\frac{\bar{\rho}}{\tau} \right) 
\\ 
  &\stackrel{\eqref{eq:99_acc_rho}+ \eqref{eq:99_acc_sigma}}{=}&
\hat{\rho}\|\nabla f(x) \|^2 + \frac{\bar{\sigma}^2}{n\tau} + \hat{\sigma}^2.
\end{eqnarray*}
Above, $(*)$ holds since $\frac{1}{\tau} (g_i)_{U_i} -\nabla f_i(x)$ is zero mean for all $i$ and $U_i, U_j$ are independent for $i\neq j$.

\end{proof}

\section{Proofs for Section~\ref{sec:99_sega}}
\subsection{Useful lemmata}
First, we mention a basic property of the proximal operator. 
\begin{proposition}
    \label{pr:prox_contraction}
    Let $R$ be a closed and convex function. Then for any $x, y\in\RR^d$
    \begin{align} \label{eq:99_prox_contraction}
        \|\proxR(x) - \proxR(y)\|
        \le \|x - y\|.
    \end{align}
\end{proposition}

The next lemma, taken from \cite[Lemma B.3]{sega}, gives a basic recurrence for the sequence $\{ h_i^k\}_{t=1}^\infty$ from {\tt ISEGA}. 
\begin{lemma} If $h_i^{k+1} \eqdef h_i^k + \tau (g_i^k - h^k)$, where $g_i^k \eqdef h_i^k + \frac{1}{\tau} (\nabla f_i(x^k) - h_i^k)_{U_i^k}$, then
\begin{equation} \label{eq:99_sega_h_bound}
\EE\left[ \|h^{k+1}_i - \nabla f_i(x^*)\|^2 \right] =(1-\tau) \|h^k_i - \nabla f_i(x^*) \|^2+\tau \| \nabla f_i(x^k)-\nabla f_i(x^*)\|^2.
\end{equation}
\end{lemma}

We will also require a recurrent bound on sequence $\{ g^k\}_{t=1}^\infty$ from {\tt ISEGA}. 
\begin{lemma}
Consider any vectors $v_i$ and set $v\eqdef \frac1n \sum_{i=1}^n v_i$. Then, we have
\begin{equation} \label{eq:99_sega_g_bound}
\EE \left[ \| g^k - v \|^2\right] \leq 
\frac{2}{n^2} \sum_{i=1}^n \left( \left( \frac{1}{\tau}+ (n-1)\right) \| \nabla f_i(x^k) - v_i\|^2
+    \left( \frac{1}{\tau}-1\right)\|h^k_i - v_i\|^2\right) .
\end{equation}
\end{lemma}
\begin{proof}
   Writing $g^k - v = a + b$, where 
   \[a\eqdef \frac1n \sum_{i=1}^n \left( h_i^k - v_i - \tau^{-1}(h_i^k - v_i)_{U_i^k}  \right)\] and 
   \[b\eqdef \frac1n \sum_{i=1}^n \tau^{-1}  (\nabla f_i(x^k) - v_i)_{U_i^k} \]
 we get $\|g^k-v\|^2=\|a+b\|^2=\leq 2(\|a\|^2 + \|b\|^2)$. 
  
 Let us bound $\EE \left[\|b\|^2\right]$ using Young's inequality $2\<x, y> \le \|x\|^2 + \|y\|^2$:
 {
\footnotesize
\begin{eqnarray*}
&& \EE \left[\|b\|^2\right]
\\
&=& 
\frac{1}{n^2} \EE \left[ \< \sum_{i=1}^n  \tau^{-1}(\nabla f_i(x^k) - v_i)_{U_i^k} ,  \sum_{i=1}^n  \tau^{-1}  (\nabla f_i(x^k) - v_i)_{U_i^k} >\right]
\\
&=& 
\frac{1}{\tau^2 n^2} \EE \left[\sum_{i=1}^n \left\| (\nabla f_i(x^k) - v_i)_{U_i^k} \right\|^2  \right] 
 + 
\frac{2}{\tau^2n^2} \EE \left[\sum_{i\neq j} \<(\nabla f_i(x^k) - v_i)_{U_i^k}  , (\nabla f_j(x^k) - v_j)_{U_i^k} > \right]  
\\
&=& 
\frac{1}{\tau n^2}\sum_{i=1}^n \left\| \nabla f_i(x^k) - v_i \right\|^2
 + 
\frac{2}{n^2} \sum_{i\neq j}  \<\nabla f_i(x^k) - v_i, \nabla f_j(x^k) - v_j> 
\\
&\leq& 
\frac{1}{\tau n^2}\sum_{i=1}^n \left\| \nabla f_i(x^k) - v_i \right\|^2 
+
 \frac{1}{n^2} \sum_{i\neq j} \left( \| \nabla f_i(x^k) - v_i\|^2 +\| \nabla f_j(x^k) - v_j\|^2\right)
\\
&= & 
\frac{1}{n^2}\left( \frac1\tau+ n-1 \right)  \sum_{i=1}^n \| \nabla f_i(x^k) - v_i\|^2.
\end{eqnarray*}
}
Similarly we bound $\EE \left[\|a\|^2\right]$:
{
\footnotesize
\begin{eqnarray*}
\EE \left[\|a\|^2\right]
&=& 
\frac{1}{n^2} \EE \left[\<  \sum_{i=1}^n \left( h_i^k - v_i - \tau^{-1}(h_i^k - v_i)_{U_i^k}  \right),   \sum_{i=1}^n \left( h_i^k - v_i - \tau^{-1} (h_i^k - v_i)_{U_i^k}  \right)> \right]
\\
&=& 
\frac{1}{n^2} \EE \left[\sum_{i=1}^n \< \left( h_i^k - v_i - \tau^{-1} (h_i^k - v_i)_{U_i^k} \right),\left( h_i^k - v_i - \tau^{-1}(h_i^k - v_i)_{U_i^k}  \right)> \right] 
\\
&& \qquad + 
\frac{2}{n^2}\EE \left[\sum_{i\neq j}\< \left( h_i^k - v_i - \tau^{-1} (h_i^k - v_i)_{U_i^k}  \right),\left( h_j^k - v_j - \tau^{-1} (h_j^k - v_j)_{U_i^k}  \right)> \right]  
\\
&=& 
\frac{\tau^{-1}-1}{n^2}\sum_{i=1}^n  \|h^k_i - v_i\|^2.
\end{eqnarray*}
}
It remains to combine the above results.
\end{proof}
\subsection{Proof of Theorem~\ref{thm:99_sega}}
\begin{proof}
For convenience, denote $g^k \eqdef \frac1n \sum_{i=1}^n g_i^k$. It holds
{
\footnotesize
\begin{eqnarray}
\nonumber
&& \EE[\|x^{k+1} - x^*\|^2] \\
\nonumber
&=& \EE\left[\|\prox_{\alpha \psi}(x^k - \alpha g^k) - \prox_{\alpha \psi}(x^* - \alpha \nabla f(x^*))\|^2\right]\\ \nonumber
&\stackrel{\eqref{eq:99_prox_contraction}}{\leq}& \EE\left[\|x^k - \alpha g^k - (x^* - \alpha \nabla f(x^*))\|^2\right]\\ \nonumber
&=& \|x^k - x^*\|^2 - 2\alpha \<\nabla f(x^k) - \nabla f(x^*),  x^k - x^*>+ \alpha^2\EE\left[ \|g^k - \nabla f(x^*)\|^2\right] 
\\ \nonumber
&\stackrel{\eqref{eq:99_sega_g_bound}}{\leq}&
\|x^k - x^*\|^2 - 2\alpha \<\nabla f(x^k) - \nabla f(x^*),x^k - x^*>
\\ \nonumber &&  +
\alpha^2 \frac{2}{n^2} \sum_{i=1}^n \left( \left( \frac{1}{\tau}+ (n-1)\right) \| \nabla f_i(x^k) - \nabla f_i(x^*)\|^2
+    \left( \frac{1}{\tau}-1\right)\|h^k_i - \nabla f_i(x^*)\|^2\right)
\\ \nonumber
&\leq&
\|x^k - x^*\|^2 -  \alpha \mu \|x^k-x^* \|^2  -2\alpha D_f(x^k,x^*)
\\ &&  +
\frac{2}{n^2} \sum_{i=1}^n \left( \left( \frac{1}{\tau}+ (n-1)\right) \| \nabla f_i(x^k) - \nabla f_i(x^*)\|^2
+    \left( \frac{1}{\tau}-1\right)\|h^k_i - \nabla f_i(x^*)\|^2\right).
\label{eq:99_sega_first}
\end{eqnarray}
}
Moreover, we have from smoothness and convexity of $f_i$
\begin{equation}\label{eq:99_sega_smoothness}
-2D_{f_i}(x^k,x^*)\leq -\frac{1}{L} \|\nabla f_i(x^k)-\nabla f_i(x^*) \|^2. 
\end{equation}
Combining the above, for any $\seganu \geq 0$ (which we choose later) we get
{
\footnotesize
\begin{eqnarray*}
&&\EE[\|x^{k+1} - x^*\|^2] + \alpha \seganu  \frac1n\sum_{i=1}^n \EE\left[ \|h^{k+1}_i - \nabla f_i(x^*)\|^2 \right]  \\
&& \quad 
\stackrel{\eqref{eq:99_sega_first}+\eqref{eq:99_sega_h_bound}}{\leq}
\|x^k - x^*\|^2 -  \alpha \mu \|x^k-x^* \|^2  -2\alpha D_f(x^k,x^*)
\\ && \qquad \qquad +
\alpha^2 \frac{2}{n^2} \sum_{i=1}^n \left( \left( \frac{1}{\tau}+ (n-1)\right) \| \nabla f_i(x^k) - \nabla f_i(x^*)\|^2
+    \left( \frac{1}{\tau}-1\right)\|h^k_i - \nabla f_i(x^*)\|^2\right)
\\ && \qquad \qquad + 
 \alpha \seganu  \frac1n\sum_{i=1}^n \left((1-\tau) \|h^k_i - \nabla f_i(x^*) \|^2+\tau \| \nabla f_i(x^k)-\nabla f_i(x^*)\|^2\right)\\
 && \qquad 
 \stackrel{\eqref{eq:99_sega_smoothness}}{\leq}
\|x^k - x^*\|^2 -  \alpha \mu \|x^k-x^* \|^2  - \frac{\alpha}{nL} \sum_{i=1}^n \|\nabla f_i(x^k)-\nabla f_i(x^*) \|^2
\\ && \qquad \qquad +
\alpha^2\frac{2}{n^2} \sum_{i=1}^n \left( \left( \frac{1}{\tau}+ (n-1)\right) \| \nabla f_i(x^k) - \nabla f_i(x^*)\|^2
+    \left( \frac{1}{\tau}-1\right)\|h^k_i - \nabla f_i(x^*)\|^2\right)
\\ && \qquad \qquad + 
 \alpha \seganu  \frac1n\sum_{i=1}^n \left((1-\tau) \|h^k_i - \nabla f_i(x^*) \|^2+\tau \| \nabla f_i(x^k)-\nabla f_i(x^*)\|^2\right)\\
 && \qquad 
 =
 (1-\alpha \mu)\|x^k - x^*\|^2 + \left(\seganu\tau + \frac{2\alpha}{n}\left(\frac1\tau+n-1\right)-\frac{1}{L}\right)\frac{\alpha}{n}\sum_{i=1}^n\| \nabla f_i(x^k)-\nabla f_i(x^*)\|^2
 \\ && \qquad \qquad +
  \left( \frac{2\alpha}{n}\left( \frac1\tau -1\right)+ \seganu(1-\tau) \right) \frac{\alpha}{n}\sum_{i=1}^n \|h^k_i - \nabla f_i(x^*) \|^2.
\end{eqnarray*}
}
To get rid of gradient differences in this bound, we want to obtain $\frac{1}{L}\geq \frac{2\alpha}{n}(\frac{1}{\tau}+n-1)+ \seganu\tau$, which, in turn, is satisfied if
\begin{eqnarray*}
\alpha &=& {\cal O}\left(\frac{1+\tau n}{L}\right),\\
\omega &=& {\cal O}\left(\frac{1}{L\tau}\right).
\end{eqnarray*} 
Next, we want to prove contraction with factor $(1 - \alpha\mu)$ in terms of $\|h_i^k-\nabla f_i(x^*)\|^2$, so we require
\[
(1-\alpha\mu ) \nu\geq \omega (1-\tau)+\frac{2\alpha}{n}\left(\frac1\tau-1\right)
\]
we shall choose $\alpha$ such that the following two properties hold:
\begin{eqnarray*}
\alpha &=& {\cal O}\left(\frac{\tau}{\mu}\right),\\
\alpha &=& {\cal O}\left(\frac{n\tau^2\omega}{1-\tau}\right) = {\cal O}\left(\frac{n\tau}{(1-\tau)L}\right)\geq  {\cal O}\left(\frac{n\tau}{L}\right).
\end{eqnarray*} 
In particular, the choice $\omega = \frac{1}{2L\tau}$ and $\alpha= \min\left( \frac{1}{4L\left( 1+\frac{1}{n\tau }\right)}, \frac{1}{\frac{\mu}{\tau}+ \frac{4L}{n\tau}}\right) $ works. 
\end{proof}
\section{Proofs for Section~\ref{sec:99_asynch}}
One way to analyze a delayed algorithm is to define sequence of epoch start moments $T_0, T_1, \dotsc$ such that $T_0=0$ and $T_{k+1} = \min\{t: t - \max_{i=1,\dotsc,n} d_i^k \ge T_k\}$. In case delays are bounded uniformly, i.e.\ for some number $M$ it holds $d_i^k \le M$ for all $i$ and $t$, one can show by induction~\cite{mishchenko2018distributed} that $T_k\le Mk$.

In addition, we define for every $i$ sequence
\begin{align*}
	z_i^k = x^{t - d_i^k}.
\end{align*}
For notational simplicity, we will assume that if worker $i$ does not perform an update at iteration $t$, then all related vectors increase their counter without changing their value,  i.e.\ $g_i^{k+1} = g_i^k$, $U_i^{k+1}=U_i^k$, $z_i^{k+1} = z_i^k$ and $x_i^{k+1} = x_i^k$. Then, we can write a simple identity for $x_i^k$ that holds for any $i$ and $k$,
\begin{align}
	x_i^k
	= x^{t - d_i^k}  - \alpha (g_i^{k})_{U_i^{k}}
	= z_i^k  - \alpha (g_i^{k})_{U_i^{k}}. \label{eq:99_delayed_recurrence}
\end{align}
\subsection{Useful lemmata}
\begin{lemma}\label{lem:99_asynch}
	Let Assumption~\ref{as:99_bounded_noise_at_opt} be satisfied and assume without loss of generality that $d_1^k<\dotsc< d_n^k$. Then, for any $i$
	\begin{align}
		\EE\|x_i^k - \EE[x_i^k \mid z_i^k, x_{i+1}^k, \dotsc, x_n^k]\|^2
		&\le 4\alpha^2\tau\EE\left[\sigma^2 + 2L(f_i(z_i^k) - f_i(x^*) - \<\nabla f_i(x^*), z_i^k -x^*>)\right].\label{eq:99_conditioned_variance}
	\end{align}
\end{lemma}
\begin{proof}
	Denote by $\cF_i^k$ the sigma-algebra generated by $z_i^k, x_{i+1}^k, \dotsc, x_n^k$. Then, \[\EE\left[\cdot \mid z_i^k, x_{i+1}^k, \dotsc, x_n^k\right] = \EE\left[\cdot \mid \cF_i^k\right].\]

	Since $d_1^k<\dotsc< d_n^k$, $x_i^k$ is independent of the randomness in $x_1^k, \dotsc, x_{i-1}^k$ as those vectors were obtained after $x_i^k$. Recall that
	\begin{align*}
		x_i^k 
		\overset{\eqref{eq:99_delayed_recurrence}}{=} z_i^k  - \alpha (g_i^{k})_{U_i^{k}}
	\end{align*}
	and denote $\tilde x_i^k \eqdef z_i^{k} - \nabla f_i(z_i^k)$. Clearly, by uniform sampling of the blocks $\EE[ x_i^k\mid \cF_i^k] = z_i^{k} - \tau\alpha\EE[g_i^k\mid \cF_i^k] = z_i^k - \alpha\tau \nabla f_i(z_i^k)$. Thus,
	\begin{align*}
		& \EE\|x_i^k - \EE[x_i^k \mid \cF_i^k]\|^2 \\
		&= \alpha^2 \EE\|(g_i^k)_{U_i^k} - \tau \nabla f_i(z_i^k)\|^2\\
		&= (1 - \tau)\alpha^2\EE\|\tau \nabla f_i(z_i^k)\|^2 +\tau \alpha^2\EE \| g_i^k - \tau \nabla f_i(z_i^k)\|^2 \\
		&= (1 - \tau)\alpha^2\tau^2\EE\|\nabla f_i(z_i^k)\|^2 +\tau \alpha^2\EE \left[\| \nabla f_i(z_i^k) - \tau \nabla f_i(z_i^k)\|^2 + \|g_i^k - \nabla f_i(z_i^k)\|^2 \right]\\
		&\le \tau \alpha^2\EE\left[\|\nabla f_i(z_i^k)\|^2 + \|g_i^k - \nabla f_i(z_i^k)\|^2\right] \\
		&\le\tau \alpha^2\EE\left[\|\nabla f_i(z_i^k)\|^2 + 2\sigma^2 + 4L(f_i(z_i^k) - f_i(x^*) - \<\nabla f_i(x^*), z_i^k - x^*>)\right] .
	\end{align*}
	In addition,
	\begin{align*}
		\|\nabla f_i(z_i^k)\|^2
		& \le 2\|\nabla f_i(z_i^k) - \nabla f_i(x^*)\|^2 + 2 \|\nabla f_i(x^*)\|^2 \\
		& \le 4L(f_i(z_i^k) - f_i(x^*) - \<\nabla f_i(x^*), z_i^k - x^*>) + 2 \sigma^2.
	\end{align*}
\end{proof}
We will use in the proof of Theorem~\ref{th:asynch} Jensen's inequality for a set of vectors $a_1,\dotsc, a_n\in\RR^d$ in the form
\begin{align*}
	\left\|\avein a_i \right\|^2
	\le \avein \|a_i\|^2.
\end{align*}
\begin{lemma}\label{lem:99_asynch_contraction}
	Assume that $f_i$ is $L$-smooth and $\mu$-strongly convex. If $\tilde x_i^k \eqdef z_i^k - \tau\alpha \nabla f_i(z_i^k)$ and $x_i^* \eqdef x^* - \tau\alpha \nabla f_i(x^*)$, we have
	\begin{align*}
		\|\tilde x_i^k - x_i^*\|^2
		\le (1 - \tau\alpha\mu)\|z_i^k - x^*\|^2 - 2\alpha\tau(1 - \tau\alpha L)(f_i(z_i^k) - f_i(x^*) - \<\nabla f_i(z_i^k), z_i^k - x^*>).
	\end{align*}
\end{lemma}
\begin{proof}
	It holds
	\begin{align*}
		 \|\tilde x_i^k - x_i^*\|^2
		= \|z_i^k - x^*\|^2 - 2\alpha\tau\<\nabla f_i(z_i^k) - \nabla f_i(x^*), z_i^k - x^*> + \alpha^2\tau^2\|\nabla f_i(z_i^k) - \nabla f_i(x^*)\|^2.
	\end{align*}
	Moreover, by strong convexity and smoothness of $f_i$ (see e.g.~\cite{nesterov2018lectures})
	\begin{align*}
		2\<\nabla f_i(z_i^k) - \nabla f_i(x^*), z_i^k - x^*>
		\ge \mu \|z_i^k - x^*\|^2 + 2(f_i(z_i^k) - f_i(x^*) - \<\nabla f_i(z_i^k), z_i^k - x^*>).
	\end{align*}
	On the other hand, convexity and smoothness of $f_i$ together imply
	\begin{align*}
		\|\nabla f_i(z_i^k) - \nabla f_i(x^*)\|^2
		\le 2L(f_i(z_i^k) - f_i(x^*) - \<\nabla f_i(z_i^k), z_i^k - x^*>).
	\end{align*}
	Consequently,
	\begin{align*}
		 \|\tilde x_i^k - x_i^*\|^2
		\le (1 - \tau\alpha\mu)\|z_i^k - x^*\|^2 - 2\tau\alpha(1 - \tau\alpha L)(f_i(z_i^k) - f_i(x^*) - \<\nabla f_i(z_i^k), z_i^k - x^*>).
	\end{align*}
\end{proof}

\subsection{Proof of Theorem~\ref{th:asynch}}
We are going to prove a more general result that does not need uniform boundedness of delays over time. Theorem~\ref{th:asynch} will follow as a special case of the more general theorem.
\begin{theorem}\label{th:asynch_epoch}
	Assume that every $f_i$ is $L$-smooth and $\mu$-strongly convex and also that the gradients noise has bounded variance at $x^*$ as in Assumption~\ref{as:99_bounded_noise_at_opt}. If also $\alpha \le \frac{1}{2L(\tau + \frac{2}{n})}$, then for any $t\in[T_k, T_{k+1})$
	\begin{align*}
		\EE\|x^k - x^*\|^2
		\le \left(1 - \tau\alpha\mu\right)^\qwerty\max_{i=1,\dotsc, n}\|x^0 - x_i^*\|^2 + 4\alpha\frac{\sigma^2}{\mu n}.
	\end{align*}
\end{theorem}
\begin{proof}
	Recall that we use in the Algorithm $w^k = \avein x_i^k$ and that $x^k =\proxR(w^k)$. Next, by non-expansiveness of the proximal operator it holds for all $t$
	\begin{align*}
		\|x^k - x^*\|^2
		&= \|\proxR(w^k) - \proxR(x^* - \alpha \nabla f(x^*)\|^2\\
		&\le \|w^k - (x^* - \alpha \nabla f(x^*))\|^2.
	\end{align*}
	Denote for simplicity $w^* \eqdef x^* - \alpha \nabla f(x^*)$. Then, we have shown $\|x^k - x^*\|^2 \le \|w^k - w^*\|^2$.
	
	Fix any $t$ and assume without loss of generality that $d_1^k < d_2^k < \dotsb < d_n^k$. Then, using the tower property of expectation
	\begin{align*}
		\EE \|w^k - w^*\|^2
		= \EE\left[\EE\left[\|w^k - w^*\|^2 \mid x_2^k, \dotsc, x_n^k\right] \right].
	\end{align*}
	At the same time, conditioned on $z_1^k, x_2^k, \dotsc, x_n^k$ the only randomness in $w^k$ is from $x_1^k$, so
	\begin{align*}
		\EE\left[\|w^k - w^*\|^2 \mid z_1^k, x_2^k, \dotsc, x_n^k\right]
		&= \|\EE\left[w^k - w^* \mid z_1^k, x_2^k, \dotsc, x_n^k\right]\|^2 \\
		& \qquad + \frac{1}{n^2}\EE\|x_1^k - \EE[x_1^k \mid z_1^k,x_2^k, \dotsc, x_n^k]\|^2.
	\end{align*}
	By continuing unrolling the first term in the right-hand side we arrive at
	\begin{align*}
		\EE \|w^k - w^*\|^2
		&\le \|\EE w^k - w^*\|^2 + \frac{1}{n^2}\sum_{i=1}^n \EE\|x_i^k - \EE[x_i^k \mid z_i^k, x_{i+1}^k, \dotsc, x_n^k]\|^2\\
		&= \|\EE w^k - w^*\|^2 + \frac{1}{n^2}\sum_{i=1}^n \EE\|x_i^k - \EE[x_i^k \mid z_i^k, x_{i+1}^k, \dotsc, x_n^k]\|^2\\
		&\overset{\eqref{eq:99_conditioned_variance}}{\le } \|\EE w^k - w^*\|^2 + 4\tau\alpha^2\frac{\sigma^2}{n}  \\
		& \qquad +\frac{8\tau \alpha^2L}{n^2}\sumin(f_i(z_i^k) - f_i(x^*) - \<\nabla f_i(x^*), z_i^k -x^*>).
	\end{align*}
	Moreover, by Jensen's inequality
	\begin{align*}
		\|\EE w^k - w^*\|^2
		&= \left\| \avein \EE[x_i^k - x_i^*]\right\|^2 \\
		&\le \avein\left\|\EE [x_i^k - x_i^*]\right\|^2 \\
		&= \avein\left\|\EE [\EE[x_i^k - x_i^*\mid \cF_i^k]\right\|^2 \\
		&\le \avein\EE\left\|\EE[x_i^k - x_i^*\mid \cF_i^k]\right\|^2.
	\end{align*}
	Combining it with our older results, we get
	{
\footnotesize
	\begin{align*}
		\EE\|w^k - w^*\|^2
		\le \avein\EE\left\|\tilde x_i^k - x_i^*\right\|^2 + 4\tau\alpha^2\frac{\sigma^2}{n}  +\frac{8\tau \alpha^2L}{n^2}\sumin(f_i(z_i^k) - f_i(x^*) - \<\nabla f_i(x^*), z_i^k -x^*>).
	\end{align*}
	}
	Let us apply Lemma~\ref{lem:99_asynch_contraction} to verify that
	{
\footnotesize
	\begin{align*}
		&\avein \|\tilde x_i^k  - x_i^*\|^2 + \frac{8\tau\alpha^2 L}{n^2}\sumin (f_i(z_i^k) - f_i(x^*) - \<\nabla f_i(x^*), z_i^k - x^*>) \\
		& \le (1 - \tau\alpha\mu)\avein \|z_i^k - x^*\|^2 \\
		& \qquad - 2\tau\alpha\underbrace{\left(1 - 2\tau\alpha L - \frac{4 \alpha L}{n}\right)}_{\ge 0}\avein (f_i(z_i^k) - f_i(x^*) - \<\nabla f_i(x^*), z_i^k - x^*>) \\
		& \le (1 - \tau\alpha\mu)\avein \|z_i^k - x^*\|^2.
	\end{align*}
	}
	Since $z_i^k =  x^{t-d_i^k}$, we have proved
	\begin{align*}
		\EE\|x^k - x^*\|^2
		&\le (1 - \tau\alpha\mu)\avein \EE\|z_i^k-x^*\|^2 + 4\tau\alpha^2\frac{\sigma^2}{n} \\
		&=  (1 - \tau\alpha\mu)\avein \EE\|x^{t-d_i^k}-x^*\|^2 + 4\tau\alpha^2\frac{\sigma^2}{n}\\
		&\le (1 - \tau\alpha\mu)\max_i \EE\|x^{t-d_i^k}-x^*\|^2 + 4\tau\alpha^2\frac{\sigma^2}{n}.
	\end{align*}
	If we define sequences \[\upsilon^k\eqdef \max_{i=1,\dotsc, n}\EE\|x^{t-d_i^k} - x^*\|^2\] and \[\Psi^\qwerty \eqdef \max_{t\in [T_k, T_{k+1})} \left\{\max\{0, \upsilon^k -  4\alpha\frac{\sigma^2}{\mu n}\}\right\},\] it follows from the above that
	\begin{align*}
		\EE\|x^k - x^*\|^2 - 4\alpha\frac{\sigma^2}{\mu n}
		\le (1 - \tau\alpha\mu)(\max_i \EE\|x^{t-d_i^k}-x^*\|^2 - 4\alpha\frac{\sigma^2}{\mu n}).
	\end{align*}
	Therefore, if $\Psi^{k_0}=0$ for some $k_0$ then $\Psi^\qwerty=0$ for all $k\ge k_0$. Otherwise,
	 $\Psi^{\qwerty+1} \le (1 - \tau\alpha\mu)\Psi^\qwerty$ and for any $t\in[T_k, T_{k+1})$
	\begin{align*}
		\EE\|x^k - x^*\|^2
		\le \psi^{k}
		\le \Psi^\qwerty + 4\alpha\frac{\sigma^2}{\mu n}
		\le (1 - \tau\alpha\mu)^\qwerty \|x^0 - x^*\|^2 + 4\alpha\frac{\sigma^2}{\mu n}.
	\end{align*}
\end{proof}
\begin{proposition}[\cite{mishchenko2018distributed}]\label{pr:epoch}
	If delays are uniformly bounded over time, i.e.\ $d_i^k\le M$ for any $i$ and $t$, then $T_k\le Mk$.
\end{proposition}
Combining Theorem~\ref{th:asynch_epoch} and Proposition~\ref{pr:epoch} gives a lower bound on $k$ and implies Theorem~\ref{th:asynch}.

\chapter{ Appendix for Chapter \ref{jacsketch}}
\label{jacsketch_appendix}

\graphicspath{{jacsketch/images/}}

\section{Summary of complexity results}

We provide a comprehensive table for faster navigation through special cases and their iteration complexities. In particular, for each special case of {\tt GJS}, we provide the leading complexity term  (i.e., a $\log \frac{1}{\varepsilon}$ factor is omitted in all results) and a reference to the corresponding corollary where this result is established. We also indicate how the operator $\cB$ appearing in the Lyapunov function is picked (this is not needed to run the method; it is only used in the analysis). All details can be found later in the Appendix.

\begin{table}[!h]
\begin{center}
\footnotesize
\begin{tabular}{|c|c|c|c|c|}
\hline
 \multicolumn{2}{|c|}{\bf Algorithm} & \multicolumn{3}{|c|}{\bf Theory} \\
\hline
 \#  &  Name  &  Cor.\ of Thm~\ref{thm:gjs_main} &  $\cB \mX $ & Leading complexity term (i.e., $\log \frac{1}{\varepsilon}$ factor omitted) \\
\hline
\hline
 \ref{alg:gjs_SAGA} &  {\tt SAGA} &  Corollary~\ref{cor:gjs_saga} & $\beta \mX$ & $ \nR + \frac{4 m}{\mu}$ \\
\hline
\ref{alg:gjs_SAGA_AS_ESO} & {\tt SAGA}  & Corollary~\ref{cor:gjs_saga_as2}  & $\mX\diag(b)$& $\max  \limits_j  \left(   \frac{1}{\pRj} +   \frac{1}{  \pRj} \frac{4 v_j }{\mu n}  \right)$ \\
\hline
   \ref{alg:gjs_SEGA} & {\tt SEGA} &  Corollary~\ref{cor:gjs_sega} &  $\beta \mX$ & $ \dL+  \dL        \frac{4  m}{\mu}  $ \\
\hline
 \ref{alg:gjs_SEGAAS} & {\tt SEGA} &  Corollary~\ref{cor:gjs_sega_is_11} &  $\diag(b)\mX$ & $  \max  \limits_i\left(    \frac{1}{\pLi} +  \frac{1}{\pLi} \frac{4 m_i}{\mu} \right) $ \\
\hline
  \ref{alg:gjs_SVRCD} & {\tt SVRCD} & Corollary~\ref{cor:gjs_svrcd} & $\beta \mX$& $\frac{1}{\probx} + \max  \limits_i \frac{1}{\pLi}\frac{4 m_i}{\mu}   $ \\
\hline
  \ref{alg:gjs_SGD_AS} & {\tt SGD-star} & Corollary~\ref{cor:gjs_sgd}  & $0$&  $\max  \limits_j  \frac{1}{ \pRj}\frac{v_j}{\mu  n} $ \\
\hline
  \ref{alg:gjs_LSVRG-AS} & {\tt LSVRG}  & Corollary~\ref{cor:gjs_lsvrg_as} & $\beta \mX$& $\frac{1}{ \probx } +  \max  \limits_j \frac{1}{ \pRj} \frac{4v_j }{\mu n }  $ \\
\hline
\ref{alg:gjs_B2} & {\tt B2} & Corollary~\ref{cor:gjs_B2} &  $\beta \mX$ & $ \frac{1}{\probx} + \frac{1}{\proby}\frac{4 m}{\mu}  $ \\
\hline
  \ref{alg:gjs_invsvrg}  & {\tt LSVRG-inv}  & Corollary~\ref{cor:gjs_inverse_svrg} & $ \mX \diag(b)$ & $\max \limits_j  \frac{1}{\pRj}   + \frac{1}{\proby}\frac{4 m}{\mu}  $ \\
\hline
\ref{alg:gjs_B_sega} & {\tt SVRCD-inv} & Corollary~\ref{cor:gjs_SVRCD-inv} &  $\diag(b) \mX$&  $ \max \limits_i \frac{1}{\pLi}  + \frac{1}{\proby}\frac{4 m}{\mu}    $ \\
\hline
 \ref{alg:gjs_RL} & {\tt RL} &  Corollary~\ref{cor:gjs_RL} &  $\mX \diag(b)$ & $\max \limits_{i,j} \left( \frac{1}{\pRj} + \frac{1}{\pLi} \frac{4 m_i^j}{\mu }  \right)$  \\
 \hline
 \ref{alg:gjs_LR} & {\tt LR} &  Corollary~\ref{cor:gjs_LR} &  $\diag(b) $ \mX  &  $\max \limits_{i,j} \left(  \frac{1}{\pLi}  + \frac{1}{\pRj} \frac{4 v_j}{\mu }\right) $ \\
 \hline
\ref{alg:gjs_saega} & {\tt SAEGA} & Corollary~\ref{cor:gjs_saega} &  $\mB \circ \mX$& $ \max \limits_{i,j}   \left(\frac{1}{\pLi \qRj}  +  \frac{1}{\pLi \qRj} \frac{4 m^j_i}{\mu n }   \right)$ \\
 \hline
  \ref{alg:gjs_svrcdg} & {\tt SVRCDG} & Corollary~\ref{cor:gjs_svrcdg} &  $\beta  \mX$& $ \frac{1}{\probx} + \max \limits_{i,j}   \frac{1}{\pLi \qRj}  \frac{4 m^j_i}{\mu n}  $ \\
 \hline
 \ref{alg:gjs_isaega} & {\tt ISAEGA} & Corollary~\ref{cor:gjs_isaega}  & $\mB \circ \mX$& $  \max \limits_{j\in N_\tR, i,\tR}   \left(\frac{1}{\ptLi \qtRj}  + \left( 1+ \frac{1}{ n\ptLi \qtRj} \right) \frac{4 m_i^j}{\mu}  \right)$ \\
\hline
 \ref{alg:gjs_isega} & {\tt ISEGA} & Corollary~\ref{cor:gjs_isega}  & $\mB \circ \mX$& $  \max \limits_{j\in N_\tR, i,\tR}   \left(\frac{1}{\ptLi |\NRt|}  + \left( 1+ \frac{1}{ n\ptLi |\NRt|} \right) \frac{4 m_i^j}{\mu}  \right)$ \\
\hline
  \ref{alg:gjs_jacsketch} & {\tt JS}  & Corollary~\ref{cor:gjs_jacsketh} & $\beta\mX\mB$& $ \frac{ 4n^{-1}\ugly \mu^{-1} \lambda_{\max}\left(  \mB^\top \E{\mR} \mB  \right) + \lambda_{\max} \left(\mB^\top \mB\right) }{\lambda_{\min} \left(\mB^\top \E{\mR} \mB \right) }  $ \\
\hline
\end{tabular}
\end{center}
\caption{Iteration complexity of selected special cases of {\tt GJS} (Algorithm~\ref{alg:gjs_SketchJac}). Whenever $m$ appears in a result, we assume that $\mM_j = m \mI_d$ for all $j$ (i.e., $f_j$ is $m$-smooth).  Whenever $m_i$ appears in a result, we assume that $f$ is $\mM$-smooth with $\mM= \diag(m_1,\dots,m_d)$. Whenever $m_i^j$ appears in a result, we assume that $\mM_j=\diag(m_1^j,\dots,m_d^j)$. Quantities $\pLi$ for $i \in [d]$, $\pRj$ for $j\in [n]$, $\probx$ and $\proby$ are probabilities defining the algorithms. }
\label{tbl:gjs_all_special_cases_theory}
\end{table}

\section{Several lemmas}
 
\subsection{Existence lemma} 
\begin{lemma} \label{lem:gjs_existence}
Suppose that $\cX \in \Range{\cM}$. Denote $\piop(\mX)\eqdef \cU\mX\eR$. Suppose that $ \E{\left( \piop\cM^{\frac12} \right)^*\piop \cM^{\frac12} }$ exists and $\lambda_{\min}\left(\E{\cS}\right)>0$. Then, there are $\alpha>0$ and $\cB$ such that~\eqref{eq:gjs_small_step} and~\eqref{eq:gjs_small_step2} hold. Moreover, inequalities \eqref{eq:gjs_small_step}, \eqref{eq:gjs_small_step2} hold for $\alpha =0, \cB=0$ without any extra assumptions. 

\end{lemma}

\begin{proof}
Consider only $\alpha, \cB$ such that that $\alpha< \lambda_{\min}\left(\E{\cS}\right) \mu^{-1}$, $\lambda_{\min}\left(\cB^*\cB\right)>0$, $\lambda_{\max}\left(\cB^*\cB\right)<\infty$.
Let $\mY  = {\cM^\dagger}^{\frac12} \mX$. Thus we have $  \E{ \norm{ \cU  \mX \eR }^2 } \leq   \|\mY\|^2\lambda_{\max}  \E{\left( \piop\cM^{\frac12} \right)^*\piop \cM^{\frac12} }$. 

 Thus 
\begin{eqnarray*}
  (1-\alpha \mu) \NORMG{ \cB\mY } - \NORMG{  \left(\cI - \E{\cS} \right)^{\frac12}\cB  \mY} 
  &=& \<(\cB\mY)^\top, ( \E{\cS} - \alpha\mu \cI) \cB\mY > \\
  &\geq& \left(\lambda_{\min}\left(\E{\cS}\right) - \alpha \mu\cI\right)\| \cB\mY\|^2 \\
  & \geq&   \left(\lambda_{\min}\left(\E{\cS}\right) - \alpha \mu\cI \right)\lambda_{\min}\left(\cB^*\cB\right)\| \mY\|^2\,.
\end{eqnarray*}
Therefore, to have~\eqref{eq:gjs_small_step}, it suffices to set
\[
\alpha \leq  \frac{ \lambda_{\min}\left(\E{\cS}\right) \lambda_{\min}\left(\cB^*\cB\right)}{  \mu \lambda_{\min}\left(\cB^*\cB\right) + \frac{2}{n^2}  \lambda_{\max}  \left[\E{\left( \piop\cM^{\frac12} \right)^*\piop \cM^{\frac12} }\right]}.
\]
Similarly, to satisfy~\eqref{eq:gjs_small_step2}, it suffices to have
\[
\frac{2\alpha}{n} \lambda_{\max}\left(\E{\left( \piop\cM^{\frac12} \right)^*\piop \cM^{\frac12} } \right) +  n\lambda_{\min}\left(\E{\cS}\right) \lambda_{\max}\left(\cB^*\cB\right) \leq 1.
\]
A valid choice to satisfy the above is for example $\alpha, \cB$ such that
\[
\lambda_{\max}\left(\cB^*\cB\right) \leq \frac{1}{2n\lambda_{\min} \left(\E{\cS}\right) }, \quad \alpha \leq  \frac{1}{\frac{1}{n} \lambda_{\max}\left(\E{\piop \left(\cM^{\frac12} \right)^*\piop \cM^{\frac12} }\right) }.
\]

\end{proof}
 
 \subsection{Smoothness lemmas}
 
Let $h:\R^d\to \R$ be a differentiable and convex function. The  Bregman distance of $x$ and $y$ with respect to $h$ is defined by
\begin{equation} \label{eq:gjs_b987gf98f}
D_h(x,y) \eqdef h(x) - h(y) - \<\nabla h(y),x-y >.
\end{equation}

 \begin{lemma}[Lemma~\ref{lem:sega_relate} from the appendix of Chapter~\ref{sega}]\label{lem:gjs_smooth22}
 Suppose that function $h:\R^d\to \R$ is convex and $\mM$-smooth, where $\mM\succeq 0$. Then 
\begin{equation}\label{eq:gjs_smooth_lemma}
D_h(x,y) \geq  \frac12 \norm{\nabla h(y)-\nabla h(x)}^2_{\mM^{\dagger}}, \quad \forall x,y \in \R^d.
\end{equation}
Further, 
\begin{equation} \label{eq:gjs_smooth_dotprod}
\<\nabla h(x) - \nabla h(y), x-y> \geq \| \nabla h(x) - \nabla h(y)\|^2_{\mM^{\dagger}}.
\end{equation}
 \end{lemma}

\begin{proof}
Fix $y$ and consider the function $\phi(x) \eqdef h(x) - \<\nabla h(y), x >$. Clearly, $\phi$ is $\mM$-smooth, and hence
\begin{equation}\label{eq:gjs_n98hf8gf}
\phi(x+d) \leq \phi(x) + \langle \nabla \phi(x),d\rangle + \frac{1}{2}\|d\|_{\mM}^2, \quad \forall x,d\in \R^d.
\end{equation}
 Moreover, since $h$ is convex, $\phi$ is convex, non-negative and is minimized at  $y$. Letting $t=\nabla h(x)-\nabla h(y)$, this implies that
\begin{eqnarray*}
\phi(y) &\leq & \phi \left(x - \mM^{\dagger}t \right) \\
&\overset{\eqref{eq:gjs_n98hf8gf}}{\leq} & \phi(x) -  \langle \nabla \phi(x),\mM^{\dagger}t\rangle + \frac{1}{2}\|\mM^{\dagger}t\|_{\mM}^2 \\
& = & \phi(x) -  \langle t ,\mM^{\dagger}t\rangle + \frac{1}{2}\|\mM^{\dagger}t\|_{\mM}^2 \\
&= & \phi(x) - \frac12\|t  \|^2_{\mM^{\dagger}},
\end{eqnarray*}
which is equivalent to~\eqref{eq:gjs_smooth_lemma}. In the last step we have used the identities $(\mM^\dagger)^\top = \left(\mM^\top \right)^\dagger = \mM^\dagger$ and $\mM^\dagger \mM \mM^\dagger = \mM^\dagger$.

To show~\eqref{eq:gjs_smooth_dotprod}, it suffices to sum inequality~\ref{eq:gjs_smooth_lemma} applied on vector pairs $(x,y)$ and $(y,x)$.

\end{proof}

\begin{lemma} \label{lem:gjs_smooth2} Let \eqref{eq:gjs_smooth_ass} hold. That is, assume that function $f_j$ are convex and $\mM_j$-smooth. Then
\begin{equation}\label{eq:gjs_smooth}
D_{f_j}(x,y) \geq \frac12 \norm{ \nabla f_j(x)-\nabla f_j(y) }^2_{\mM_j^{\dagger}}, \quad \forall x,y\in \R^d .
\end{equation}
If $x-y\in \Null{\mM_j}$, then 
\begin{enumerate} 
\item[(i)]  \begin{equation}\label{eq:gjs_linear_on_subspace} f_j(x) = f_j(y) + \langle \nabla f_j(y), x-y\rangle,\end{equation}
\item[(ii)]
\begin{equation} \label{eq:gjs_n98g8ff} \nabla f_j(x)-\nabla f_j(y) \in \Null{\mM_j},\end{equation}
\item[(iii)]  
\begin{equation} \label{eq:gjs_nb87sgb} \langle \nabla f_j(x) - \nabla f_j(y),x-y\rangle =0.\end{equation}
\end{enumerate}

If, in addition, $f_j$ is bounded below, then $\nabla f_j(x)  \in \Range{\mM_j}$ for all $x$.
\end{lemma}

\begin{proof} 
Inequality \eqref{eq:gjs_smooth} follows by applying Lemma~\ref{lem:gjs_smooth22} for $h=f_j$ and $\mM=\mM_j$. Identity \eqref{eq:gjs_linear_on_subspace} is a direct consequence of  \eqref{eq:gjs_smooth_ass}.  Combining \eqref{eq:gjs_smooth} and \eqref{eq:gjs_linear_on_subspace}, we get $0 \geq \frac12 \norm{ \nabla f_j(x)-\nabla f_j(y) }^2_{\mM_j^{\dagger}}$ , which implies that \begin{equation}\label{eq:gjs_nbui8gf7gdnjs}\nabla f_j(x)-\nabla f_j(y) \in \Null{\mM_j^\dagger} = \Null{\mM_j^\top} =\Null{\mM_j}  ,\end{equation} 
recovering \eqref{eq:gjs_n98g8ff}.  By adding two copies of  \eqref{eq:gjs_linear_on_subspace} (with the roles of $x$ and $y$ exchanged), we get \eqref{eq:gjs_nb87sgb}. Finally, if $f_j$ is bounded below, then in view of \eqref{eq:gjs_linear_on_subspace} there exists $c\in \R$ such that,
\[c \leq   \inf_{x \in  y +\Null{\mM_j}} f_j(x)  \overset{\eqref{eq:gjs_linear_on_subspace} }{=} \inf_{x \in  y +\Null{\mM_j}} f_j(y) +  \langle \nabla f_j(y), x-y\rangle  .\]
This implies that $\nabla f_j(y) \in \Range{\mM_j^\top} = \Range{\mM_j} $.


\end{proof}

\begin{lemma} \label{lem:gjs_smooth3}
Assume $f$ is twice continuously differentiable. Then $\mG(x)-\mG(y)\in \Range{\cM}$ for all $x,y\in \R^d$. 
\end{lemma}

\begin{proof}
For $\mG(x)-\mG(y)\in \Range{\cM}$, it suffices to show that $\nabla f_j(x) - \nabla f_j(y)\in \Range{\mM_j}$.
Without loss of generality, suppose that $f(z,w)$ (for $x  = [z,w]$) is such that $f(z,\cdot)$ is linear (for fixed $z$; from~\eqref{eq:gjs_linear_on_subspace}) and $f(\cdot,w)$ is $\mM'$ smooth for full rank $\mM'$. Note that

\[
0\preceq \nabla^2f(x) = \begin{pmatrix}
 \nabla_{ww}^2f(w,z) &  \nabla_{wz}^2f(w,z) \\
  \nabla_{zw}^2f(w,z) &  \nabla_{zz}^2f(w,z)
\end{pmatrix} = 
\begin{pmatrix}
 \nabla_{ww}^2f(w,z) &  \nabla_{wz}^2f(w,z) \\
  \nabla_{zw}^2f(w,z) &  0
\end{pmatrix} .
\]
Since every submatrix of the above must be positive definite, it is easy to see that we must have both $ \nabla_{wz}^2f(w,z)= 0$, $ \nabla_{zw}^2f(w,z) = 0$. This, however, means that $f(w,z)$ is separable in $z,w$. Therefore indeed  $\nabla f_j(x) - \nabla f_j(y)\in \Range{\mM_j}$ for all $x,y\in \R^d$ and all $j\in [n]$.

\end{proof}

\subsection{Projection lemma}

In the next lemma, we establish some basic properties of the interaction of the random projection matrices $\cS$ and $\cI - \cS$ with various matrices, operators, and norms.

\begin{lemma} \label{lem:gjs_nb98gd8fdx}
Let $\cS$ be a random projection operator and $\cA$ any deterministic linear operator commuting with $\cS$, i.e.,  $\cA \cS = \cS \cA$. Further, let $\mX,\mY \in \R^{d\times n}$ and define $\mZ = (\cI-\cS) \mX + \cS \mY$. Then
\begin{itemize}
\item[(i)] $\cA \mZ = (\cI-\cS) \cA \mX + \cS \cA \mY $,
\item[(ii)] $\norm{\cA \mZ}^2 = \norm{(\cI-\cS) \cA \mX}^2  + \norm{\cS \cA \mY}^2 $,
\item[(iii)] $\E{\norm{\cA \mZ}^2} = \norm{(\cI-\E{\cS})^{1/2} \cA \mX}^2  + \norm{\E{\cS}^{1/2} \cA \mY}^2 $, where the expectation is with respect to $\cS$.
\end{itemize}
\end{lemma}
\begin{proof} Part (i) follows by noting that $\cA$ commutes with $\cI-\cS$. Part (ii) follows from (i) by expanding the square, and noticing that  $(\cI-\cS)\cS = 0$. Part (iii) follows from (ii) after using the definition of the Frobenius norm, i.e., $\|\mM\|^2 =\Tr{\mM^\top \mM}$, the identities $(\cI -\cS)^2  = \cI -\cS$, $\cS^2 =\cS$,  and taking expectation on both sides. 
\end{proof}

\subsection{Decomposition lemma}

In the next lemma, we give a  bound on the expected squared distance of the gradient estimator $g^k$ from $\nabla f(x^*)$.

\begin{lemma}\label{lem:gjs_g_lemma}
For all $k\geq 0$ we have
\begin{equation}\label{eq:gjs_g_lemma}
\E{ \norm{ g^k - \nabla f(x^*)}^2} \leq    \frac{2}{n^2} \E{  \norm{ \cU \left(\mG(x^k) - \mG(x^*) \right) \eR}^2  } + \frac{2}{n^2}  \E{ \norm{ \cU \left(\mJ^k - \mG(x^*) \right) \eR}^2}.
\end{equation}
\end{lemma}
\begin{proof}
In view of \eqref{eq:gjs_ni98hffs} and since $\nabla f(x^*) = \frac{1}{n}\mG(x^*) \eR$, we have
\begin{equation}\label{eq:gjs_nb87fvdbs8s}
g^k-\nabla f(x^*) = \underbrace{\frac1n \cU \left(\mG(x^k)- \mG(x^*) \right) \eR}_{a}  +  \underbrace{\frac1n \left(\mJ^k - \mG(x^*) \right) \eR - \frac1n \cU \left( \mJ^k - \mG(x^*) \right) \eR}_{b} .
\end{equation}
Applying the bound $\norm{a+b}^2 \leq 2\norm{a}^2 + 2\norm{b}^2$ to \eqref{eq:gjs_nb87fvdbs8s} and taking expectations, we get
\begin{eqnarray*}
\E{ \norm{ g^k -\nabla f(x^*) }^2} &\leq &
 \E{ \frac{2}{n^2} \norm{  \cU \left(\mG(x^k)- \mG(x^*) \right) \eR  }^2} \\
 && \qquad + \E{ \frac{2}{n^2} \norm{   \left(\mJ^k-\mG(x^*) \right) \eR -   \cU \left(\mJ^k - \mG(x^*) \right) \eR }^2  }\\
&=&
 \frac{2}{n^2} \E{  \norm{ \cU \left(\mG(x^k) - \mG(x^*) \right) \eR }^2  } \\
 && \qquad + \frac{2}{n^2} \E{  \norm{ \left(\cI- \cU \right) \left(\mJ^k - \mG(x^*) \right) \eR }^2 }.
\end{eqnarray*}

It remains to note that
\begin{eqnarray*}
\E{  \norm{ \left(\cI-  \cU) (\mJ^k - \mG(x^*) \right) \eR }^2 }
 &=&
\E{ \norm{ \cU \left(\mJ^k - \mG(x^*) \right) \eR  }^2 } - \norm{ \left( \mJ^k - \mG(x^*) \right) \eR }^2
  \\
  &\leq& 
  \E{ \norm{ \cU \left(\mJ^k - \mG(x^*) \right) \eR }^2 } .
\end{eqnarray*}

\end{proof}

\section{Proof of Theorem~\ref{thm:gjs_main}}

For simplicity of notation, in this proof, all expectations are conditional on $x^k$, i.e., the expectation is taken with respect to the randomness of $g^k$.

Since
\begin{equation}\label{eq:gjs_prox_opt}
x^* =  \prox_{\alpha \psi} (x^* - \alpha \nabla f(x^*)),
\end{equation}
and since the prox operator is non-expansive, we have
\begin{eqnarray}
\E{\norm{x^{k+1} -x^*}^2 } &\overset{ \eqref{eq:gjs_prox_opt}}{=} &
 \E{\norm{ \prox_{\alpha \psi} (x^k-\alpha g^k) -  \prox_{\alpha \psi} (x^*-\alpha \nabla f(x^*))  }^2}  
 \notag \\
 &\leq &
 \E{\norm{x^k-  x^* -\alpha(   g^{k} - \nabla f(x^*) ) }^2}  
 \notag\\
& \overset{\eqref{eq:gjs_unbiased_xx}}{=}& 
 \norm{x^k  -x^*}^2 -2\alpha \<  \nabla f(x^k)- \nabla f(x^*) , x^k  -x^*> \notag \\ 
 && \qquad  + \alpha^2\E{\norm{ g^{k}- \nabla f(x^*) }^2}
 \notag  \\ 
&\overset{\eqref{eq:gjs_strconv3}+ \eqref{eq:gjs_b987gf98f}}{\leq} & 
   (1-\alpha\mu)\norm{x^k  -x^*}^2 +\alpha^2\E{\norm{  g^{k} -\nabla f(x^*) }^2}  \notag \\
   && \qquad -2\alpha D_f(x^k,x^*).
 \label{eq:gjs_convstepsub1XX}
\end{eqnarray}

Since $f(x)=\frac{1}{n}\sum_{j=1}^n f_j(x)$, in view of \eqref{eq:gjs_b987gf98f} and \eqref{eq:gjs_smooth} we have
\begin{eqnarray}
D_{f}(x^k,x^*) \quad \overset{\eqref{eq:gjs_b987gf98f}}{=} \quad  \frac{1}{n}\sum_{j=1}^n D_{f_j}(x^k,x^*) & \overset{\eqref{eq:gjs_smooth}}{\geq} &\frac{1}{2n} \sum_{j=1}^n \norm{\nabla f_j(x^k) -\nabla f_j(x^*)  }^2_{{\mM_j^{\dagger}} }  \notag\\
&=&  
 \frac{1}{2n}  \left\|{\cM^\dagger}^{\frac12}\left(\mG(x^k) -\mG(x^*) \right) \right\|^2. \label{eq:gjs_nb98gd8ff}
\end{eqnarray}

By combining \eqref{eq:gjs_convstepsub1XX} and \eqref{eq:gjs_nb98gd8ff}, we get
\begin{eqnarray*}
\E{\norm{x^{k+1} -x^*}^2 } & \leq &
   (1-\alpha\mu)\norm{x^k  -x^*}^2 +\alpha^2\E{\norm{g^{k} -\nabla f(x^*) }^2}  \\
   && \qquad -\frac{\alpha}{n}  \norm{ {\cM^\dagger}^{\frac12} \left( \mG(x^k) -\mG(x^*) \right) }^2.
\end{eqnarray*}

Next, applying Lemma~\ref{lem:gjs_g_lemma} leads to the estimate
\begin{eqnarray}
\E{\norm{x^{k+1} -x^*}^2 } &\leq &
   (1-\alpha\mu)\norm{x^k  -x^*}^2 -\frac{\alpha}{n} \norm{ {\cM^\dagger}^{\frac12} \left(\mG(x^k) -\mG(x^*) \right) }^2 \notag 
   \\
   && \qquad  +  \frac{2\alpha^2}{n^2} \E{ \norm{ \cU \left( \mG(x^k) - \mG(x^*) \right) \eR }^2  } \notag \\
   && \qquad + \frac{2\alpha^2}{n^2}  \E{ \norm{ \cU \left( \mJ^k - \mG(x^*) \right) \eR  }^2}   \label{eq:gjs_48u34719841234} .
\end{eqnarray}

In view of~\eqref{eq:gjs_nio9h8fbds79kjh}, we have $\mJ^{k+1} = (\cI-\cS)\mJ^k + \cS \mG(x^k)$, whence
\begin{equation} \label{eq:gjs_h98gf9hh89dsd}
\underbrace{\mJ^{k+1} - \mG(x^*)}_{\mZ} = (\cI-\cS) \underbrace{(\mJ^k -\mG(x^*))}_{\mX} + \cS \underbrace{(\mG(x^k) - \mG(x^*))}_{\mY}.
\end{equation}
Since, by assumption,  both $\cB$ and ${\cM^\dagger}^{\frac12}$ commute with $\cS$, so does their composition $\cA \eqdef \cB {\cM^\dagger}^{\frac12}$. Applying Lemma~\ref{lem:gjs_nb98gd8fdx}, we get
 \begin{eqnarray}\label{eq:gjs_J_jac_bound}\nonumber
\E{ \NORMG{\cB {\cM^\dagger}^{\frac12} \left(\mJ^{k+1}-\mG(x^*)  \right) }}  &=&  \NORMG{ (\cI - \E{\cS})^{\frac12}  \cB {\cM^\dagger}^{\frac12} \left(\mJ^k-\mG(x^*) \right) }  \\
&& +  \NORMG{\E{\cS}^{\frac12}  \cB {\cM^\dagger}^{\frac12} \left(\mG(x^k)-\mG(x^*) \right) } .
\end{eqnarray}

Adding $\alpha$-multiple of~\eqref{eq:gjs_J_jac_bound}  to~\eqref{eq:gjs_48u34719841234} yields
\begin{eqnarray*}
&&\E{\norm{x^{k+1} -x^*}^2 } + \alpha\E{ \NORMG{\cB {\cM^\dagger}^{\frac12} \left(\mJ^{k+1}-\mG(x^*)\right)}}
\\
&& \qquad  \leq 
   (1-\alpha\mu)\norm{x^k  -x^*}^2 + \frac{2\alpha^2}{n^2} \E{  \norm{ \cU \left(\mG(x^k) - \mG(x^*) \right) \eR}^2  }    \\
   && \qquad \qquad  
 + \frac{2\alpha^2}{n^2}  \E{ \norm{ \cU \left(\mJ^k - \mG(x^*) \right) \eR }^2 } +   \alpha\NORMG{ (\cI - \E{\cS})^{\frac12} \cB {\cM^\dagger}^{\frac12} \left(\mJ^k-\mG(x^*) \right)} 
 \\
    && \qquad \qquad   + \alpha \NORMG{\E{\cS}^{\frac12}  \cB {\cM^\dagger}^{\frac12}\left(\mG(x^k)-\mG(x^*) \right)}-\frac{\alpha}{n} \norm{ {\cM^\dagger}^{\frac12}\left(\mG(x^k) -\mG(x^*) \right) }^2
   \\
&& \qquad  \stackrel{\eqref{eq:gjs_small_step}}{\leq} 
(1-\alpha\mu)\norm{x^k  -x^*}^2 + (1-\alpha\mu)\alpha \NORMG{\cB {\cM^\dagger}^{\frac12} \left(\mJ^{k}-\mG(x^*) \right)} \\
&& \qquad \qquad  
+\frac{2\alpha^2}{n^2} \E{  \norm{ \cU \left(\mG(x^k) - \mG(x^*) \right) \eR }^2  } +  \alpha \NORMG{\E{\cS}^{\frac12}  \cB {\cM^\dagger}^{\frac12} \left(\mG(x^k)-\mG(x^*) \right)}  \\
    && \qquad \qquad  
 -\frac{\alpha}{n} \left\| {\cM^\dagger}^{\frac12} \left(\mG(x^k) -\mG(x^*) \right) \right\|^2
   \\
&& \qquad  \stackrel{\eqref{eq:gjs_small_step2}}{\leq} 
(1-\alpha\mu) \left( \norm{x^k  -x^*}^2 +\alpha \NORMG{\cB {\cM^\dagger}^{\frac12} \left(\mJ^{k}-\mG(x^*) \right)} \right).
\end{eqnarray*}
Above, we have used~\eqref{eq:gjs_small_step} with $\mX=  \mJ^k-\mG(x^*) $ and~\eqref{eq:gjs_small_step2} with $\mX=  \mG(x^k) -\mG(x^*) $.

\section{Special cases: {\tt SAGA}-like methods}

\subsection{Basic variant of {\tt SAGA}~\cite{saga} \label{sec:gjs_saga_basic}}

Suppose that for all $j$, $f_j$ is $m$-smooth (i.e., $\mM_j = m \mI_d$). To recover basic SAGA~\cite{saga}, consider the following choice of random operators $\cS,  \cU$: 
\[
(\forall j) \text{ with probability } \frac1n :  \quad \cS \mX =\mX \eRj \eRj^\top \quad \text{and} \quad   \cU \mX =\mX n \eRj \eRj^\top.
 \]

The resulting algorithm is stated as Algorithm~\ref{alg:gjs_SAGA}. Further, as a direct consequence of Theorem~\ref{thm:gjs_main}, convergence rate of {\tt SAGA} (Algorithm~\ref{alg:gjs_SAGA}) is presented in Corollary~\ref{cor:gjs_saga}.

\begin{algorithm}[h]
    \caption{{\tt SAGA} \cite{saga}}
    \label{alg:gjs_SAGA}
    \begin{algorithmic}
        \Require learning rate $\alpha>0$, starting point $x^0\in\R^d$
        \State Set $\psi_j^0 = x^0$ for each $j\in \{1,2,\dots,n\}$
        \For{ $k=0,1,2,\ldots$ }
        \State{Sample  $j \in [n]$ uniformly at random}
        \State{Set $\phi_j^{k+1} = x^k$ and $\phi_i^{k+1} = \phi_i^{k}$ for $i\neq j$}
        \State{$g^k = \nabla f_j(\phi_j^{k+1}) - \nabla f_j(\phi_j^k) + \frac{1}{n}\sum\limits_{i=1}^n\nabla f_i(\phi_i^k)$}
        \State{$x^{k+1} =  \prox_{\alpha \psi} ( x^k - \alpha g^k)$}
        \EndFor
    \end{algorithmic}
\end{algorithm}

\begin{corollary}[Convergence rate of {\tt SAGA}]\label{cor:gjs_saga}  Let $\alpha =  \frac{1}{4m + \mu n}$. Then, iteration complexity of Algorithm~\ref{alg:gjs_SAGA} (proximal {\tt SAGA}) is $\left( 4\frac{m}{\mu}+ n \right) \log\frac1\epsilon$. 
\end{corollary}

\subsection{{\tt SAGA} with arbitrary sampling \label{sec:gjs_saga_as}}

In contrast to Section~\ref{sec:gjs_saga_basic}, here we use the general matrix smoothness assumption, i.e., that $f_j$ is $\mM_j$ smooth. We  recover results from~\cite{qian2019saga}. Denote $\pR$ to be probability vector, i.e., $\pR_i = \Probbb{i\in R}$ where $R$ is a random subset of $[n]$.

We shall consider the following choice of random operators $\cS,  \cU$: 
\[
(\forall R) \text{ with probability } \pRR :  \quad \cS \mX = \mX \sum_{j\in R}  \eRj \eRj^\top \quad \text{and} \quad   \cU \mX =\mX \sum_{j\in R} \frac{1}{\pRj} \eRj \eRj^\top .
 \]

The resulting algorithm is stated as Algorithm~\ref{alg:gjs_SAGA_AS_ESO}.
\begin{algorithm}[h]
    \caption{{\tt SAGA} with arbitrary sampling (a variant of \cite{qian2019saga})}
    \label{alg:gjs_SAGA_AS_ESO}
    \begin{algorithmic}
        \Require learning rate $\alpha>0$, starting point $x^0\in\R^d$, random sampling $R \subseteq \{1,2,\dots,n\}$
        \State Set $\phi_j^0 = x^0$ for each $j\in [n]$
        \For{ $k=0,1,2,\ldots$ }
        \State{Sample random  $R^k \subseteq \{1,2,\dots,n\}$}
        \State{Set $\phi_j^{k+1} = \begin{cases} x^k & j\in R^k\\  \phi_j^{k} &  j\not \in R^k \end{cases}$}
        \State{$g^k = \frac{1}{n}\sum\limits_{j=1}^n\nabla f_j(\phi_j^k) +\sum \limits_{j\in R^k} \frac{1}{n \pRj} \left( \nabla f_j(\phi_j^{k+1}) - \nabla f_j(\phi_j^k) \right)$}
        \State{$x^{k+1} =  \prox_{\alpha \psi} (x^k - \alpha g^k)$}
        \EndFor
    \end{algorithmic}
\end{algorithm}

In order to give tight rates under $\mM$-smoothness, we need to do a  bit more work. First, let $v\in \R^n$ be a vector for which the following inequality {\em expected separable overapproximation} inequality holds
\begin{equation} \label{eq:gjs_ESO_saga}
\E{\left\|\sum_{j \in R} \mM^{\frac12}_{j} h_{j}\right\|^{2}} \leq \sum_{j=1}^{n} \pRj v_{j}\left\|h_{j}\right\|^{2}, \qquad \forall h_1,\dots,h_n \in \R^{d} .
\end{equation}
Since the function on the left is a quadratic in $h = (h_1,\dots,h_n) \in \R^{nd}$, this inequality is satisfied for large enough values of $v_j$.
A variant of~\eqref{eq:gjs_ESO_saga} was used to obtain the best known rates for coordinate descent with arbitrary sampling~\cite{qu2016coordinate1, qu2016coordinate2}.

Further, we shall consider the following assumption:
\begin{assumption}\label{as:gjs_pseudoinverse}
Suppose that  for all $k$ 
\begin{equation}\label{eq:gjs_g_in_range}
\mG(x^k) - \mG(x^*) = \cM^{\dagger} \cM \left(\mG(x^k) - \mG(x^*)\right)
\end{equation}
 and 
 \begin{equation}\label{eq:gjs_j_in_range}
\mJ^k - \mG(x^*) = \cM^{\dagger} \cM \left(\mJ^k - \mG(x^*)\right).
\end{equation}
\end{assumption}
The assumption, although in a slightly less general form, was demonstrated to obtain tightest complexity results for {\tt SAGA}~\cite{qian2019saga}. Note that if for each $j$, $f_j$ corresponds to loss function of a linear model, then~\eqref{eq:gjs_g_in_range} and~\eqref{eq:gjs_j_in_range} follow for free. Further, Lemmas~\ref{lem:gjs_smooth2} and~\ref{lem:gjs_smooth3} give some easy-to-interpret sufficient  sufficient conditions, such as lower boundedness of all functions $f_j$ (which happens for any loss function), or twice differentiability of all functions $f_j$.

\begin{corollary}[Convergence rate of {\tt SAGA}]\label{cor:gjs_saga_as2}  Let $ \alpha  = \min_j \frac{n \pRj}{4v_j + n\mu}$. Then the iteration complexity of Algorithm~\ref{alg:gjs_SAGA_AS_ESO} is $\max_j \left(  \frac{4v_j  + n\mu }{n \mu \pRj}\right) \log\frac1\epsilon$. 
\end{corollary}

\begin{remark}
Corollary~\ref{cor:gjs_saga_as2} is slightly more general than Theorem 4.6 from~\cite{qian2019saga} does not explicitly require linear models and $\mM$ smoothness implied by the linearity.
\end{remark}

\section{Special cases: {\tt SEGA}-like methods}

Let $n=1$. Note that now operators $\cS$ and $\cU$ act on $d\times n$ matrices, i.e., on vectors in $\R^d$. To simplify notation, instead of $\mX\in \R^{d\times n}$ we will write $x = (x_1,\dots,x_d) \in \R^{d}$.

\subsection{Basic variant of {\tt SEGA}~\cite{sega}  \label{sec:gjs_sega}}

Suppose that  $f$ is $m$-smooth (i.e., $\mM_1=m\mI_d$) with $m>0$. To recover  basic SEGA from~\cite{sega}, consider the following choice of random operators $\cS$ and $\cU$: 
\[
(\forall i) \text{ with probability } \frac1d :  \quad \cS x = \eLi \eLi^\top x = x_i \eLi \quad \text{and} \quad   \cU x = d \eLi \eLi^\top x =  dx_i \eLi.
 \]
The resulting algorithm is stated as Algorithm~\ref{alg:gjs_SEGA}.

\begin{algorithm}[h]
    \caption{{\tt SEGA} \cite{sega}}
    \label{alg:gjs_SEGA}
    \begin{algorithmic}
        \Require Stepsize $\alpha>0$, starting point $x^0\in\R^d$
        \State Set $h^0 = 0$
        \For{ $k=0,1,2,\ldots$ }
        \State{Sample  $i\in \{1,2,\dots d \}$ uniformly at random}
        \State{Set $h^{k+1} = h^{k} +( \nabla_i f(x^k) - h^{k}_i)\eLi$}
        \State{$g^k = h^k+ d (\nabla_i f(x^k) - h_i^k) \eLi$}
        \State{$x^{k+1} =  \prox_{\alpha \psi} (x^k - \alpha g^k)$}
        \EndFor
    \end{algorithmic}
\end{algorithm}

\begin{corollary}[Convergence rate of SEGA]\label{cor:gjs_sega}  Let $\alpha =  \frac{1}{4md+ \mu d}$. Then the iteration complexity of Algorithm~\ref{alg:gjs_SEGA} is $\left( 4\frac{md}{\mu}+d \right) \log\frac1\epsilon$. 
\end{corollary}

\subsection{{\tt SEGA} with arbitrary sampling \label{sec:gjs_sega_is_v1}}
Consider a more general setup to that in Section~\ref{sec:gjs_sega} and let us allow the smoothness matrix to be an arbitrary diagonal (positive semidefinite) matrix: $\mM = \diag(m_1,\dots,m_d)$ with $m_1,\dots,m_d>0$. In this regime, we will establish a convergence rate for an arbitrary sampling strategy, and then use this to develop importance sampling. 

Let $\pL\in \R^d$ be a probability vector with entries $\pLi = \Probbb{i\in L}$. Consider the following choice of random operators $\cS$ and  $\cU$: 
\begin{equation} \label{eq:gjs_sega_choice}
(\forall L) \text{ with prob. } \pLL :  \; \cS x = \sum_{i\in L} \eLi \eLi^\top x = \sum_{i\in L} x_i \eLi  \quad \text{and} \quad   \cU x =  \sum_{i\in L} \frac{1}{\pLi} \eLi \eLi^\top x = \sum_{i\in L} \frac{x_i}{\pLi} \eLi .
\end{equation}
The resulting algorithm is stated as Algorithm~\ref{alg:gjs_SEGAAS}.

\begin{algorithm}[h]
    \caption{{\tt SEGA} with arbitrary sampling}
    \label{alg:gjs_SEGAAS}
    \begin{algorithmic}
        \Require Stepsize $\alpha>0$, starting point $x^0\in\R^d$, random sampling $L\subseteq \{1,2,\dots,d\}$
        \State Set $h^0 = 0$
        \For{ $k=0,1,2,\ldots$ }
        \State{Sample random  $L^k \subseteq \{1,2,\dots,d\}$}
        \State{Set $h^{k+1} = h^{k} +\sum \limits_{i\in L^k}( \nabla_i f(x^k) - h^{k}_i)\eLi$}
        \State{$g^k = h^k+\sum \limits_{i\in L^k}  \frac{1}{\pLi}(\nabla_i f(x^k) - h_i^k)\eLi$}
        \State{$x^{k+1} =  \prox_{\alpha \psi} (x^k - \alpha g^k)$}
        \EndFor
    \end{algorithmic}
\end{algorithm}

\begin{corollary}[Convergence rate of {\tt SEGA}]\label{cor:gjs_sega_is_11}  Iteration complexity of Algorithm~\ref{alg:gjs_SEGAAS} with $\alpha  = \min_i \frac{\pLi}{4 m_i+ \mu}$ is $  \max_i\left(\frac{4 m_i +\mu}{\pLi\mu} \right)\log\frac1\epsilon$. 
\end{corollary}
Corollary~\ref{cor:gjs_sega_is_11} indicates an up to constant factor optimal choice $\pLi \propto m_i$, which yields, up to a constant factor, $\frac{\sum_{i=1}^d m_i}{\mu}\log\frac1\epsilon$ complexity.  
In the applications where $m$ is not unique\footnote{For example when a general matrix smoothness holds; one has to upper bound it by a diagonal matrix in order to comply with the assumptions of the section. In such case, there is an infinite array of possible choices of $m$.}, it is the best to choose one which minimizes $m^\top \eL$.

\begin{remark}
Note that if $\pLi=1$ for all $i$ (i.e., if $\cU=\cI$), we recover proximal gradient descent as a special case.
\end{remark}

\subsection{{\tt SVRCD} with arbitrary sampling \label{sec:gjs_svrcd_is2}}

As as a particular special case of Algorithm~\ref{alg:gjs_SketchJac} we get a new method, which  we call {\em Stochastic Variance Reduced Coordinate Descent} ({\tt SVRCD}). The algorithm is similar to {\tt SEGA}. The main difference is that {\tt SVRCD} does not {\em update} a subset $L$ of coordinates of vector $h^k$ each iteration. Instead, with probability $\probx$, it {\em sets} $h^k$ to $\nabla f(x^k)$.  

We choose $\cS$ and $\cU$ via
\[
  \cS \mX = 
  \begin{cases}
    0  & \text{w.p.}\quad 1- \probx \\
    \mX              & \text{w.p.}\quad  \probx
\end{cases} 
\quad \text{and} \quad   (\forall L)\quad \text{w.p.}\quad \pLL : \,  \cU \mX =  \sum_{i\in L} \frac{1}{\pLi} \eLi \eLi^\top \mX,
 \]
 where again $\pLi = \Probbb{i\in L}$. The randomness of $\cS$ is independent from the randomness of $\cU$ (which comes from the randomness of $L$). The resulting algorithm is stated as Algorithm~\ref{alg:gjs_SVRCD}.

\begin{algorithm}[h]
    \caption{{\tt SVRCD} {\bf [NEW METHOD]}}
    \label{alg:gjs_SVRCD}
    \begin{algorithmic}
        \Require starting point $x^0\in\R^d$, random sampling $L \subseteq \{1,2,\dots,d\}$, probability $\probx$, stepsize $\alpha>0$
        \State Set $h^0 = 0$
        \For{ $k=0,1,2,\ldots$ }
        \State{Sample random  $L^k \subseteq \{1,2,\dots,d\}$}
        \State{$g^k = h^k+\sum \limits_{i\in L^k}  \frac{1}{\pLi}(\nabla_i f(x^k) - h_i^k) \eLi$}
        \State{$x^{k+1} =  \prox_{\alpha \psi} (x^k - \alpha g^k)$}
        \State{Set $h^{k+1} =   \begin{cases}
    h^k  & \text{with probability} \quad 1- \probx \\
    \nabla f(x^k)              & \text{with probability} \quad \probx
\end{cases} $}
        \EndFor
    \end{algorithmic}
\end{algorithm}

As in Section~\ref{sec:gjs_sega_is_v1}, we shall assume that $f$ is $\mM= \diag(m_1, \dots, m_d)$- smooth.

\begin{corollary}\label{cor:gjs_svrcd}
The iteration complexity of Algorithm~\ref{alg:gjs_SVRCD} with $\alpha  = \min_i \frac{1}{4m_i / \pLi + \mu / \probx}$ is $$ \left(\frac{1}{\probx}  +  \max_i \frac{4m_i}{\pLi\mu}\right) \log\frac1\epsilon .$$
\end{corollary}

Corollary~\ref{cor:gjs_svrcd} indicates optimal choice $p \propto m$.

\begin{remark}
If $\pLi=1$ for all $i$ and $\probx=1$, we recover proximal gradient descent as a special case.
\end{remark}

\section{Special cases: {\tt SGD-star}} \label{sec:gjs_SGD-AS-star}

Suppose that $\mG(x^*) $ is known. We will show that shifted a version of {\tt SGD-AS} converges with linear rate in such case.  Let $\mJ^0 = \mG(x^*)$. Consider the following choice of random operators $\cS$,  $\cU$: 
\[
\cS \mX = 0  \quad \text{and} \quad    (\forall R) \text{ with probability } \pRR :  \quad \cU \mX =  \mX  \sum_{j\in R} \frac{1}{\pRj} \eRj \eRj^\top.
 \]
The resulting algorithm is stated as Algorithm~\ref{alg:gjs_SGD_AS}, which is in fact arbitrary sampling version of {\tt SGD-star} from~\cite{sigma_k}.

\begin{algorithm}[h]
    \caption{{\tt SGD-star} \cite{sigma_k} }
    \label{alg:gjs_SGD_AS}
    \begin{algorithmic}
        \Require learning rate $\alpha>0$, starting point $x^0\in\R^d$,  random sampling  $R \subseteq \{1,2,\dots,n\}$
        \For{ $k=0,1,2,\ldots$ }
        \State{Sample random  $R^k \subseteq \{1,2,\dots,n\}$ }
        \State{$g^k = \frac{1}{n}\mG(x^*) \eR +\sum \limits_{j\in R^k} \frac{1}{n \pRj} \left( \nabla f_j(x^k) -\nabla f_j(x^*) \right)$}
        \State{$x^{k+1} =  \prox_{\alpha \psi} (x^k - \alpha g^k)$}
        \EndFor
    \end{algorithmic}
\end{algorithm}

\begin{corollary}[Convergence rate of {\tt SGD-AS-star}]\label{cor:gjs_sgd}  Suppose that $f_j$ is $\mM_j$-smooth for all $j$ and suppose that $v$ satisfies~\eqref{eq:gjs_ESO_saga}. Let $\alpha =n \min_j \frac{\pRj}{v_j} $. Then, the iteration complexity of Algorithm~\ref{alg:gjs_SGD_AS} is $$ \max_j \left( \frac{v_j}{n \pRj \mu} \right) \log\frac1\epsilon .$$
\end{corollary}

\begin{remark}
In overparameterized models, one has  $\mG(x^*) =0 $. In such a case, Algorithm~\ref{alg:gjs_SGD_AS} becomes {\tt SGD-AS} \cite{pmlr-v97-qian19b}, and we recover its tight convergence rate.
\end{remark}

\section{Special cases: loopless {\tt SVRG} with arbitrary sampling ({\tt LSVRG})} \label{sec:gjs_LSVRG-AS}

In this section we extend Loopless {\tt SVRG} (i.e., {\tt LSVRG}) from~\cite{hofmann2015variance, kovalev2019don} to arbitrary sampling. The main difference to {\tt SAGA} is that {\tt LSVRG} does not update $\mJ^k$ at all with probability $1-\probx$. However, with probability $1-\probx$, it sets $\mJ^k$ to $\mG(x^k)$.  Define $\cS$ and $\cU$ as follows:
\[
  \cS \mX = 
  \begin{cases}
    0  & \text{w.p.}\quad 1- \probx \\
    \mX              & \text{w.p.}\quad  \probx
\end{cases} 
\quad \text{and} \quad   (\forall R) \text{ with probability }\; \pRR : \,  \cU \mX =  \mX  \sum_{i\in R} \frac{1}{\pRj} \eRj \eRj^\top,
 \]
 where $\pRj = \Probbb{j\in R}$.

The resulting algorithm is stated as Algorithm~\ref{alg:gjs_LSVRG-AS}.

\begin{algorithm}[h]
    \caption{{\tt LSVRG} ({\tt LSVRG} \cite{hofmann2015variance, kovalev2019don} with arbitrary sampling) {\bf [NEW METHOD]}}
    \label{alg:gjs_LSVRG-AS}
    \begin{algorithmic}
        \Require learning rate $\alpha>0$, starting point $x^0\in\R^d$, random sampling $R\subseteq \{1,2,\dots,n\}$
        \State Set $\phi = x^0$
        \For{ $k=0,1,2,\ldots$ }
        \State{Sample a random subset $R^k \subseteq \{1,2,\dots n \}$ }
        \State{$g^k = \frac{1}{n}\sum\limits_{j=1}^n \nabla f_j(\phi^k) +\sum \limits_{j\in R^k} \frac{1}{n \pRj} \left( \nabla f_j(x^{k}) - \nabla f_j(\phi^k) \right)$}
        \State{$x^{k+1} =  \prox_{\alpha \psi} (x^k - \alpha g^k)$}
        \State{Set $\phi^{k+1} = \begin{cases} x^k & \text{with probability} \quad  \probx \\ \phi^{k} & \text{with probability} \quad  1- \probx \end{cases}$}
        \EndFor
    \end{algorithmic}
\end{algorithm}

In order to give tight rates under $\mM$-smoothness, we shall consider ESO assumption~\eqref{eq:gjs_ESO_saga} and Assumption~\ref{as:gjs_pseudoinverse} (same as for {\tt SAGA-AS}).

The next corollary shows the convergence result. 
\begin{corollary}[Convergence rate of {\tt LSVRG}]\label{cor:gjs_lsvrg_as}  Let $ \alpha  = \min_j \frac{n}{4 \frac{v_j}{ \pRj} + \frac{ \mu n}{\probx}}$. Then, the iteration complexity of Algorithm~\ref{alg:gjs_LSVRG-AS} is $$\max_j \left(  4 \frac{v_j}{n \mu \pRj }  + \frac{1}{ \probx } \right) \log\frac1\epsilon .$$
\end{corollary}

\begin{remark}
One can consider a slightly more general setting with \[\cS \mX = 
  \begin{cases}
    0  & \text{w.p.}\quad 1- \probx \\
    \mX   \sum_{i\in R'} \eRj\eRj^\top    & \text{w.p.}\quad  \probx \end{cases},\] where distribution of $R'\subseteq [n] $ is arbitrary. Clearly, such methods is a special case of Algorithm~\ref{alg:gjs_SketchJac}, and setting $R' =[n]$ with probability 1, {\tt LSVRG} is obtained. However, in a general form, such algorithm resembles {\tt SCSG}~\cite{lei2017less}. However, unlike {\tt SCSG}, the described method converges linearly, thus is superior to {\tt SCSG}.
\end{remark}

\section{Special cases: methods with Bernoulli $\cU$}
Throughout this section, we will suppose that $ \mM_j = m \mI_d$ for all $j$. This is sufficient to establish strong results. Indeed, Bernoulli $\cU$ does not allow for an efficient importance sampling and hence one can't develop arbitrary sampling results similar to those in Section~\ref{sec:gjs_saga_as} or Section~\ref{sec:gjs_sega_is_v1}.  

\subsection{{\tt B2} (Bernoulli $\cS$)} \label{sec:gjs_B2}
Let $n=1$. Note that now operators $\cS$ and $\cU$ act on $d\times n$ matrices, i.e., on vectors in $\R^d$. To simplify notation, instead of $\mX\in \R^{d\times n}$ we will write $x = (x_1,\dots,x_d) \in \R^{d}$. Given  probabilities $0< \probx, \proby \leq 1$, let both $\cS$ and $ \cU$ be  Bernoulli (i.e., scaling) sketches: 
\[
\cS x= 
  \begin{cases}
    0  & \text{w.p.}\quad 1- \probx \\
    x              & \text{w.p.}\quad  \probx
\end{cases} 
\quad \text{and} \quad
\cU x = 
  \begin{cases}
    0  & \text{w.p.}\quad 1- \proby \\
   \frac{1}{\proby}x             & \text{w.p.}\quad  \proby
\end{cases} .
 \]
The resulting algorithm is stated as Algorithm~\ref{alg:gjs_B2}.

\begin{algorithm}[h]
    \caption{{\tt B2} {\bf [NEW METHOD]}}
    \label{alg:gjs_B2}
    \begin{algorithmic}
        \Require learning rate $\alpha>0$, starting point $x^0\in\R^d$, probabilities $\proby \in (0,1]$ and $\probx \in (0,1]$
        \State Set $\phi = x^0$
        \For{ $k=0,1,2,\ldots$ }
        \State{$g^k =   \begin{cases}
    \nabla f(\phi^k) & \text{with probability} \quad 1- \proby \\
   \frac{1}{\proby}\nabla f(x^k) - \left(\frac{1}{\proby}-1\right) \nabla f(\phi^k)             & \text{with probability} \quad  \proby
\end{cases} $}
        \State{$x^{k+1} =  \prox_{\alpha \psi} (x^k - \alpha g^k)$}
        \State{Set $\phi^{k+1} = \begin{cases} x^k & \text{with probability} \quad \probx\\ \phi^{k} &  \text{with probability}  \quad 1 -\probx \end{cases}$}
        \EndFor
    \end{algorithmic}
\end{algorithm}

\begin{corollary}[Convergence rate {\tt B2}]\label{cor:gjs_B2}  Suppose that $f$ is $m$-smooth. Let $\alpha = \frac{1}{4\frac{m}{ \proby} + \frac{\mu}{ \probx} }$. Then, the iteration complexity of Algorithm~\ref{alg:gjs_B2} is $$\left( 4 \frac{m}{\mu \proby}+ \frac{1}{\probx} \right) \log\frac1\epsilon .$$ 
\end{corollary}

\begin{remark}
It is possible to choose correlated $\cS$ and $\cU$ without any sacrifice in the rate. 
\end{remark}

 \subsection{{\tt LSVRG-inv} (right $\cS$)} \label{sec:gjs_SVRG-1}

Given a probability scalar $0< \proby \leq 1$, consider choosing operators $\cS$ and $\cU$ as follows: 
\[
\cS \mX = \mX \sum_{j\in R} \eRj \eRj^\top \quad \text{w.p.}\quad \pRR
\quad \quad\text{and} \quad\quad
\cU \mX = 
  \begin{cases}
    0  & \text{w.p.}\quad 1- \proby \\
   \frac{1}{\proby} \mX              & \text{w.p.}\quad  \proby.
\end{cases} 
 \]
The resulting algorithm is stated as Algorithm~\ref{alg:gjs_invsvrg}.

\begin{algorithm}[h]
    \caption{{\tt LSVRG-inv} {\bf [NEW METHOD]}}
    \label{alg:gjs_invsvrg}
    \begin{algorithmic}
        \Require starting point $x^0\in\R^d$, random sampling $R\subseteq \{1,2,\dots,n\}$, probability $\proby \in (0,1]$ , learning rate $\alpha>0$
        \State Set $\phi_j^0 = x^0$ for $j=1,2,\dots,n$
        \For{ $k=0,1,2,\ldots$ }
        \State{$g^k =   \begin{cases}
   \frac1n\sum \limits_{j=1}^n \nabla f_j(\phi_j^k) & \text{with probability} \quad 1- \proby \\
   \frac{1}{\proby}\nabla f(x^k) - \left(\frac{1}{\proby}-1 \right)  \frac1n\sum \limits_{j=1}^n \nabla f_j(\phi_j^k)            & \text{with probability} \quad  \proby
\end{cases} $}
        \State{$x^{k+1} =  \prox_{\alpha \psi} (x^k - \alpha g^k)$}
    \State{Sample a random subset $R^k \subseteq \{1,2,\dots n \}$ }
        \State{Set $\phi_j^{k+1} = \begin{cases} x^k &  \quad j\in R^k\\ \phi_j^{k} &  \quad j \notin R^k \end{cases}$}        
        \EndFor
    \end{algorithmic}
\end{algorithm}

\begin{corollary}[Convergence rate of {\tt LSVRG-inv}]\label{cor:gjs_inverse_svrg}  Suppose that each $f_i$ is $m$-smooth. Let $\alpha = \min_j \frac{1}{4\frac{m}{\proby} + \frac{\mu}{\pRj}}$. Then, the iteration complexity of Algorithm~\ref{alg:gjs_invsvrg} is $$\max_j \left( 4\frac{m}{\mu \proby}+ \frac{1}{\pRj} \right) \log\frac1\epsilon .$$
\end{corollary}

 \subsection{{\tt SVRCD-inv} (left $\cS$)} \label{sec:gjs_SVRCD-inv}

Let $n=1$. Note that now operators $\cS$ and $\cU$ act on $d\times n$ matrices, i.e., on vectors in $\R^d$. To simplify notation, instead of $\mX\in \R^{d\times n}$ we will write $x = (x_1,\dots,x_d) \in \R^{d}$.

Consider again setup where $n=1$. Choose operators $\cS$ and $\cU$ as follows: 
\[
\cS x= \sum_{i\in L} \eLi \eLi^\top x \quad \text{w.p.}\quad  \pLL
\quad \text{and} \quad
\cU x = 
  \begin{cases}
    0  & \text{w.p.}\quad 1- \proby \\
   \frac{1}{\proby} x              & \text{w.p.}\quad  \proby \;.
\end{cases} 
 \]
For convenience, let $\pL$ be the probability vector defined as: $\pLi = \Probbb{i\in L}$.

The resulting algorithm is stated as Algorithm~\ref{alg:gjs_B_sega}.

\begin{algorithm}[h]
    \caption{{\tt SVRCD-inv} {\bf [NEW METHOD]}}
    \label{alg:gjs_B_sega}
    \begin{algorithmic}
        \Require  starting point $x^0\in\R^d$, random sampling $L\subseteq \{1,2,\dots,d\}$, probability $\proby \in (0,1]$, learning rate $\alpha>0$
        \State Choose $h^0 \in \R^d$
        \For{ $k=0,1,2,\ldots$ }
        \State{$g^k =   \begin{cases}
    h^k & \text{with probability} \quad 1- \proby \\
   \frac{1}{\proby}\nabla f(x^k) - \left(\frac{1}{\proby}-1 \right)h^k             & \text{with probability} \quad  \proby
\end{cases} $}
        \State{$x^{k+1} =  \prox_{\alpha \psi} (x^k - \alpha g^k)$}
            \State{Sample a random subset $L^k \subseteq \{1,2,\dots d \}$ }
        \State{Set $h^{k+1} = h^k + \sum \limits_{i\in L^k} (\nabla_i f(x^k)  - h^k_i )\eLi$}
        \EndFor
    \end{algorithmic}
\end{algorithm}

\begin{corollary}[Convergence rate of {\tt SVRCD-inv}]\label{cor:gjs_SVRCD-inv}  Suppose that each $f_j$ is $m$-smooth. Let $\alpha = \min_i \frac{1}{4 \frac{m}{\proby} + \frac{\mu}{ \pLi}}$. Then, the iteration complexity of Algorithm~\ref{alg:gjs_B_sega} is $$\max_i \left( 4\frac{m}{\mu \proby}+ \frac{1}{\pLi} \right) \log\frac1\epsilon.$$
\end{corollary}

\section{Special cases: combination of left and right sketches}

 \subsection{{\tt RL} (right sampling $\cS$, left unbiased sampling $\cU$)}  \label{sec:gjs_RL}

Consider choosing $\cS$ and $\cU$ as follows: 
\[
\cS \mX= \mX \sum_{j\in R} \eRj \eRj^\top \quad \text{w.p.} \quad \pR_{R}
\quad \text{and} \quad
\cU \mX= \sum_{i\in L} \frac{1}{\pLi} \eLi  \eLi^\top \mX  \quad \text{w.p.} \quad \pL_{L} \;.
 \]
 
The resulting algorithm is stated as Algorithm~\ref{alg:gjs_RL}.

\begin{algorithm}[h]
    \caption{{\tt RL} {\bf [NEW METHOD]}}
    \label{alg:gjs_RL}
    \begin{algorithmic}
        \Require starting point $x^0\in\R^d$, random sampling $L\subseteq \{1,2,\dots,d\}$,  random sampling $R\subseteq \{1,2,\dots,n\}$, learning rate $\alpha>0$
        \State Set $\phi_j^0 = x^0$ for each $j$
        \For{ $k=0,1,2,\ldots$ }
        \State{Sample random $R^k \subseteq \{1,2,\dots,n\}$ }
        \State{Set $\phi_j^{k+1} = \begin{cases} x^k & \quad j\in R^k\\ \phi_j^{k} & \quad j\not \in R^k \end{cases}$}
        \State{Sample random $L^k \subseteq \{1,2,\dots,d\}$}
        \State{$g^k = \frac{1}{n}\sum\limits_{j=1}^n\nabla f_j(\phi_j^k) +\sum \limits_{i\in L^k} \frac{1}{\pLi}\left( \nabla_i f(x^k) -\frac{1}{n}\sum\limits_{j=1}^n\nabla_i f_j(\phi_j^k) \right)\eLi$}
        \State{$x^{k+1} =  \prox_{\alpha \psi} (x^k - \alpha g^k)$}
        \EndFor
    \end{algorithmic}
\end{algorithm}

\begin{corollary}[Convergence rate of {\tt RL}]\label{cor:gjs_RL}  Suppose that each $f_j$ is $\diag(m^j)$-smooth, where $m^j\in \R^d$ and $\diag(m^j)\succ 0$. Let $\alpha =\min_{i,j} \left( 4 \frac{m_i^j }{\pLi} + \frac{\mu}{ \pRj} \right)^{-1}$. Then, the iteration complexity of Algorithm~\ref{alg:gjs_RL} is $$\max_{i,j} \left( 4\frac{m_i^j}{\mu \pLi} +  \frac{1}{\pRj} \right) \log\frac1\epsilon.$$
\end{corollary}

 \subsection{{\tt LR} (left sampling $\cS$, right unbiased sampling $\cU$)}  \label{sec:gjs_LR}
 
 Consider choosing $\cS$ and $\cU$ as follows: 
\[
\cS \mX=  \sum_{i \in L} \eLi \eLi^\top  \mX\quad \text{w.p.} \quad \pL_{L}
\quad \text{and} \quad
\cU \mX=  \mX  \sum_{j\in R} \frac{1}{\pRj} \eRj  \eRj^\top \quad \text{w.p.} \quad \pR_{R} \;.
 \]

The resulting algorithm is stated as Algorithm~\ref{alg:gjs_LR}.

\begin{algorithm}[h]
    \caption{{\tt LR} {\bf [NEW METHOD]}}
    \label{alg:gjs_LR}
    \begin{algorithmic}
        \Require starting point $x^0\in\R^d$, random sampling $L\subseteq \{1,2,\dots,d\}$,  random sampling $R\subseteq \{1,2,\dots,n\}$, learning rate $\alpha>0$
        \State Set $h^0 = x^0$ for each $j$
        \For{ $k=0,1,2,\ldots$ }
        \State{Sample random $L^k \subseteq \{1,2,\dots,d\}$ }
        \State{Set $h^{k+1} = h^{k} +\sum \limits_{i\in L^k} ( \nabla_i f(x^k) - h^{k}_i)\eLi$}
        \State{Sample random $R^k \subseteq \{1,2,\dots,n\}$ }
        \State{$g^k = \nabla f(h^k) +\sum \limits_{j\in R^k} \frac{1}{n\pRj}\left( \nabla f_j(x^k) - \nabla f_j(h^k) \right)$}
        \State{$x^{k+1} =  \prox_{\alpha \psi} (x^k - \alpha g^k)$}
        \EndFor
    \end{algorithmic}
\end{algorithm}

\begin{corollary}[Convergence rate of {\tt LR}]\label{cor:gjs_LR}  Suppose that each $f_j$ is $\mM_j$-smooth, and suppose that $v \in \R^{n}$ is such that~\eqref{eq:gjs_ESO_saga} holds. Let $\alpha = \min_{i,j} \frac{1}{4 v_j\pRj^{-1} + \mu\pLi^{-1}}$. Then, the iteration complexity of Algorithm~\ref{alg:gjs_LR} is $$\max_{i,j} \left( 4\frac{v_i}{\mu \pRj}+ \frac{1}{\pLi} \right) \log\frac1\epsilon.$$
\end{corollary}

\section{Special cases: joint left and right sketches}

\subsection{{\tt SAEGA} \label{sec:gjs_SAEGA}}

Another  new special case of Algorithm~\ref{alg:gjs_SketchJac} we propose is {\tt SAEGA} (the name comes from the combination of names {\tt SAGA} and {\tt SEGA}). In {\tt SAEGA}, both $\cS$ and $\cU$ are fully correlated and consist of right and left sketch. However, the mentioned right and left sketches are independent. In particular, we have
\[
\cS \mX = \mX_{LR}  = \left( \sum_{i\in L} \eLi \eLi^\top \right)\mX  \left( \sum_{j\in R} \eRj \eRj^\top \right), 
\]
where $L\subset[d]$, and $R\subset[n]$  are independent random sets. Next, $\cU$ is chosen as 
\[
\cU \mX = \cS \left( \left( {\pL}^{-1} \left(\pR^{-1}\right)^\top\right) \circ \mX\right),
\] 
where $\pLi = \Probbb{i\in  L}$ and $\pRj =  \Probbb{j\in  R}$. The resulting algorithm is stated as Algorithm~\ref{alg:gjs_saega}. 

\begin{algorithm}[h]
  \caption{{\tt SAEGA} {\bf [NEW METHOD]}}
  \label{alg:gjs_saega}
\begin{algorithmic}
\State{\bfseries Input: }{$x^0\in\R^d$, random sampling $L\subseteq \{1,2,\dots,d\}$, random sampling  $R \subseteq \{1,2,\dots,n\}$, stepsize $\alpha$ }
\State $\mJ^0  = 0$
  \For{$k=0,1,2,\dotsc$}
         \State Sample random $ L^k \subseteq \{1,2,\dots,d\}$ and $ R^k \subseteq \{1,2,\dots,n\}$
        \State Compute $\nabla_i f_{j}(x^k)$ for all $i\in L^k$ and $j \in R^k $ 
        \State $\mJ^{k+1}_{ij} = \begin{cases} \nabla_{i} f_{j}(x^k) &   i\in L^k \text{ and } j\in  R^k \\  \mJ^k_{ij} & \text{otherwise} \end{cases}$
    \State $g^k = \left(\mJ^k + \left( {\pL}^{-1} \left(\pR^{-1}\right)^\top\right)\circ(\mJ^{k+1}- \mJ^k)\right) \eR$
    \State $x^{k+1} = \prox_{\alpha \psi} (x^k - \alpha g^k)$
  \EndFor
\end{algorithmic}
\end{algorithm}

Suppose that for all $j\in [n]$, $\mM_j = \diag(m^j)\succ 0$ is diagonal matrix\footnote{A block diagonal matrix $\mM_j$ with blocks such that $\cM\cP_\cS = \cP_\cS \cM $ would work as well}.  Let $\PR \in \R^{n\times n}$ be the probability matrix with respect to $R$-sampling , i.e., $\PR_{ jj'} = \Probbb{j\in R, j' \in R}$. 

\begin{corollary}\label{cor:gjs_saega}  Consider any (elementwise) positive vector $\qR$ such that $$ \diag(\pR)^{-1} \PR \diag(\pR)^{-1}  \preceq  \diag(\qR)^{-1}.$$  Let $\alpha =  \min_{i,j} \frac{n\pLi \qRj}{4m^j_i + n\mu}$. Then, the iteration complexity of Algorithm~\ref{alg:gjs_saega} is $$ \max_{i,j}   \left( 4\frac{m^j_i}{\mu n \pLi \qRj}  + \frac{1}{\pLi} \frac{1}{\qRj}\right)
\log\frac1\epsilon.$$
\end{corollary}

\subsection{{\tt SVRCDG} \label{sec:gjs_SVRCDG}} 


Next new special case of Algorithm~\ref{alg:gjs_SketchJac} we propose is {\tt SVRCDG}. {\tt SVRCDG} uses the same random operator $\cU$ as {\tt SAEGA}. The difference to {\tt SAEGA} lies in operator $\cS$ which is Bernoulli random variable:
\[
\cS \mX = 
\begin{cases}
0 &   \text{w.p.}\quad 1-\probx \\
\mX & \text{w.p.}\quad \probx \\
\end{cases},\qquad  
\cU \mX = \mI_{L:}\left( \left( {\pL}^{-1} \left(\pR^{-1}\right)^\top\right) \circ \mX\right)\mI_{:R},
\] 
where $L\subseteq [d]$, and $R\subseteq [n]$ are independent random sets and  $\pLi = \Probbb{i\in  L}$ and $\pRj =  \Probbb{j\in  R}$.

 The resulting algorithm is stated as Algorithm~\ref{alg:gjs_svrcdg}. 

\begin{algorithm}[h]
  \caption{{\tt SVRCDG} {\bf [NEW METHOD]}}
  \label{alg:gjs_svrcdg}
\begin{algorithmic}
\State{\bfseries Input: }{$x^0\in\R^d$, random sampling $L\subseteq \{1,2,\dots,d\}$, random sampling  $R \subseteq \{1,2,\dots,n\}$, stepsize $\alpha$, probability $\probx$ }
\State $\mJ^0  = 0$
  \For{$k=0,1,2,\dotsc$}
         \State Sample random $ L^k \subseteq \{1,2,\dots,d\}$ and $ R^k \subseteq \{1,2,\dots,n\}$
        \State Observe $\nabla_i f_{j}(x^k)$ for all $i\in L^k$ and $j \in R^k$ 
    \State $g^k = \left(\mJ^k + \left( {\pL}^{-1}\left(\pR^{-1}\right)^\top\right)\circ \left(\mI_{L^k:} \left(\mG(x^k)- \mJ^k \right) \mI_{:R^k} \right)\right) \eR$
    \State $x^{k+1} = \prox_{\alpha \psi} (x^k - \alpha g^k)$
                \State  $\mJ^{k+1} = \begin{cases} \mG(x^k) & \text{with probability} \quad\probx \\
                 \mJ^k & \text{with probability} \quad 1-\probx \end{cases}$
  \EndFor
\end{algorithmic}
\end{algorithm}

Suppose that for all $j$, $\mM_j = \diag(m^j)$ is diagonal matrix\footnote{Block diagonal $\mM_j$ with blocks such that $\cM\cS = \cS \cM $ would work as well}.  For notational simplicity, denote $\mM'\in \R^{d\times n}$ to be the matrix with $j$th column equal to $m_j$. Let $\PR \in \R^{n\times n}$ be the probability matrix with respect to $R$ - sampling , i.e., $\PR_{ jj'} = \Probbb{j\in R, j' \in R}$. 

\begin{corollary}\label{cor:gjs_svrcdg}   Consider any (elementwise) positive vector $\qR$ such that $$ \diag(\pR)^{-1} \PR \diag(\pR)^{-1}  \preceq  \diag(\qR)^{-1}.$$  Let $\alpha =  \min_{i,j} \frac{1}{4\frac{m^j_i}{\pLi \qRj n} + \frac{1}{\probx}\mu}$. Then, the iteration complexity of Algorithm~\ref{alg:gjs_svrcdg} is $$  \max_{i,j}   \left( 4\frac{m^j_i}{\mu n \pLi \qRj}  + \frac{1}{\probx}\right)
\log\frac1\epsilon.$$


\end{corollary}

\subsection{{\tt ISAEGA}  (with distributed data) \label{sec:gjs_ISAEGA}}

In this section, we consider a distributed setting from~\cite{mishchenko201999}. In particular,~\cite{mishchenko201999} proposed a strategy of running coordinate descent on top of various optimization algorithms such as {\tt GD}, {\tt SGD} or {\tt SAGA}, while keeping the convergence rate of the original method. This allows for sparse communication from workers to master. 

However, {\tt ISAGA} (distributed {\tt SAGA} with {\tt RCD} on top of it), as proposed, assumes zero gradients at the optimum which only holds for overparameterized models. It was stated as an open question whether it is possible to derive {\tt SEGA} on top of it such that the mentioned assumption can be dropped. We answer this question positively, proposing {\tt ISAEGA} (Algorithm~\ref{alg:gjs_isaega}). Next, algorithms proposed in~\cite{mishchenko201999} only allow for uniform sampling under simple smoothness. In contrast, we develop an arbitrary sampling strategy for general matrix smoothness\footnote{We do so only for {\tt ISAEGA}. However, our framework allows obtaining arbitrary sampling results for {\tt ISAGA}, {\tt ISEGA} and {\tt ISGD} (with no variance at optimum) as well. We omit it for space limitations}. 

Assume that we have $\TR$ parallel units, each owning set of indices $\NRt$ (for $1\leq \tR\leq \TR$). Next, consider distributions $\cDR_\tR$ over subsets of $\NRt$ and distributions $\cDL_{\tR}$ over subsets coordinates $[d]$ for each machine. Each iteration we sample $R_\tR \sim\cDR_{\tR}, L_\tR \sim\cDL_{\tR}$ (for $1\leq \tR\leq \TR$) and observe the corresponding part of Jacobian $\mJ^k_{\cap_{\tR} (L_\tR,R_\tR)}$. Thus the corresponding random Jacobian sketch becomes 
\[
\cS \mX = \mX_{\cap_\tR (L_\tR,R_\tR)}  = \sum_{\tR=1}^\TR  \left( \sum_{i\in L_\tR} \eLi \eLi^\top \right)\mX_{:\NRt} \left( \sum_{j\in R_\tR} \eRj \eRj^\top \right).
\]

 Next, for each $1\leq \tR \leq \TR$ consider vector $\ptL \in \R^d$, $\ptR\in \R^{|\NRt|}$ such that $\Probbb{ i\in L_\tR} =\ptLi$ and $\Probbb{ j\in R_\tR} =\ptRj$. Given the notation, random operator $\cU$ is chosen as 

\[
\cU \mX = \sum_{\tR=1}^\TR   \left( \left(\ptL\right)^{-1} \left(\left(\ptR\right)^{-1} \right)^\top\right) \circ \left( \left( \sum_{i\in L_\tR} \eLi \eLi^\top \right)\mX_{:\NRt} \left( \sum_{j\in R_\tR} \eRj \eRj^\top \right)\right). 
\] 
The resulting algorithm is stated as Algorithm~\ref{alg:gjs_isaega}. 

\begin{algorithm}[h]
  \caption{{\tt ISAEGA} {\bf [NEW METHOD]}}
  \label{alg:gjs_isaega}
\begin{algorithmic}
\State{\bfseries Input: }{$x^0\in\R^d$, \# parallel units $\TR$, each owning set of indices $N_\tR$ (for $1\leq \tR\leq \TR$), distributions $\cDR_t$ over subsets of $ \NRt$, distributions $\cDL_{\tR}$ over subsets coordinates $[d]$, stepsize $\alpha$ }
\State $\mJ^0  = 0$
  \For{$k=0,1,2,\dotsc$}
    \For{$\tR=1,\dotsc,\TR$ in parallel}
         \State Sample $ R_\tR \sim \cDR_\tR$; $R_\tR\subseteq  \NRt$ (independently on each machine)
        \State Sample $ L_\tR \sim \cDL_\tR$; $L_\tR\subseteq [d]$  (independently on each machine)
        \State Observe $\nabla_{L_\tR} f_{j}(x^k)$ for $j \in R_\tR $ 
        \State For $i\in [d], j\in \NRt$ set $\mJ^{k+1}_{i,j} = \begin{cases}
        \nabla_{i} f_{j}(x^k) & \text{if} \quad i\in [d], j\in R_\tR, i\in L_{\tR} \\
        \mJ^{k}_{i,j} & \text{otherwise}
        \end{cases} $
        \State Send $\mJ^{k+1}_{:\NRt }- \mJ^k_{:\NRt }$ to master \Comment{Sparse; low communication}
    \EndFor
    \State $g^k = \left(\mJ^k + \sum \limits_{\tR=1}^\TR   \left( {{\ptL}^{-1}} {{\ptR}^{-1}}^\top\right) \circ \left( \left( \sum_{i\in L_\tR} \eLi \eLi^\top \right)\left(\mJ^{k+1}-\mJ^k   \right)_{:\NRt} \left( \sum_{j\in R_\tR} \eRj \eRj^\top \right)\right) \right)\eR$
    \State $x^{k+1} = \prox_{\alpha \psi} (x^k - \alpha g^k)$
  \EndFor
\end{algorithmic}
\end{algorithm}

Suppose that for all $1\leq j\leq n$, $\mM_j = \diag(m^j)$ is diagonal matrix\footnote{block diagonal $\mM_j$ with blocks such that $\cM\cS = \cS \cM $ would work as well}. Let $\PtR \in \R^{ \| \NRt \|\times  \| \NRt \|}$ be the probability matrix with respect to $R_{\tR}$ - sampling , i.e., $\PtR_{jj'} = \Probbb{j\in R_{\tR}, j' \in R_{\tR}}$.

\begin{corollary}\label{cor:gjs_isaega} 
For all $\tR$ consider any (elementwise) positive vector $\qtR$ such that $ \diag(\ptR)^{-1} \PtR \diag(\ptR)^{-1}  \preceq  \diag(\qtR)^{-1}$. Let $\alpha =  \min_{j\in \NRt, i,\tR} \frac{1}{4 m^j_i \left( 1+ \frac{1}{n\ptLi \qtRj}\right) + \frac{\mu}{\ptLi }\qtRj}$. Then, iteration complexity of Algorithm~\ref{alg:gjs_isaega} is $  \max_{j\in \NRt, i,\tR}  \left( 4\frac{m^j_i}{\mu}\left( 1+ \frac{1}{ n\ptLi\qtRj} \right)  +  \frac{1}{\ptLi\qtRj} \right)
\log\frac1\epsilon$. 
\end{corollary}
Thus, for all $j$, it does not make sense to increase sampling size beyond point where $\ptLi\qtRj \geq \frac1n$as the convergence speed would not increase significantly\footnote{For indices $i,j,\tR$ which maximize the rate from Corollary~\ref{cor:gjs_isaega}.} .

\begin{remark}
In special case when $R_\tR = \NRt$ always, {\tt ISAEGA} becomes {\tt ISEGA} from~\cite{mishchenko201999}. However~\cite{mishchenko201999} assumes that $|\NRt|$ is constant in $\tR$ and $L_{\tR} = \eLi$ with probability $\frac1d$. Thus, even special case of Corollary~\ref{cor:gjs_isaega} generalizes results on {\tt ISEGA} from~\cite{mishchenko201999}. For completeness, we state {\tt ISEGA} as Algorithm~\ref{alg:gjs_isega} and Corollary~\ref{cor:gjs_isega} provides its iteration complexity. 
\end{remark}

\begin{algorithm}[h]
  \caption{{\tt ISEGA} ({\tt ISEGA}~\cite{mishchenko201999} with arbitrary sampling) {\bf [NEW METHOD]}}
  \label{alg:gjs_isega}
\begin{algorithmic}
\State{\bfseries Input: }{$x^0\in\R^d$, \# parallel units $\TR$, each owning set of indices $N_\tR$ (for $1\leq \tR\leq \TR$), distributions $\cDL_{\tR}$ over subsets coordinates $[d]$, stepsize $\alpha$ }
\State $\mJ^0  = 0$
  \For{$k=0,1,2,\dotsc$}
    \For{$\tR=1,\dotsc,\TR$ in parallel}
        \State Sample $ L_\tR \sim \cDL_\tR$; $L_\tR\subseteq [d]$  (independently on each machine)
        \State Observe $\nabla_{L_\tR} f_{j}(x^k)$ for $j \in \NRt $ 
        \State For $i\in [d], j\in \NRt$ set $\mJ^{k+1}_{i,j} = \begin{cases}
        \nabla_{i} f_{j}(x^k) & \text{if} \quad i\in [d], j\in \NRt , i\in L_{\tR} \\
        \mJ^{k}_{i,j} & \text{otherwise}
        \end{cases} $
        \State Send $\mJ^{k+1}_{:\NRt }- \mJ^k_{:\NRt }$ to master \Comment{Sparse; low communication}
    \EndFor
    \State $g^k = \left(\mJ^k + \sum_{\tR=1}^\TR   \left( {{\ptL}^{-1}} {\eR}^\top\right) \circ \left( \left( \sum_{i\in L_\tR} \eLi \eLi^\top \right)\left(\mJ^{k+1}-\mJ^k   \right)_{:\NRt} \right) \right)\eR$
    \State $x^{k+1} = \prox_{\alpha \psi} (x^k - \alpha g^k)$
  \EndFor
\end{algorithmic}
\end{algorithm}

 \begin{corollary}\label{cor:gjs_isega} 
 Let $\alpha =  \min_{j\in \NRt, i,\tR} \frac{1}{4 m^j_i \left( 1+ \frac{1}{n\ptLi |\NRt |}\right) + \frac{\mu}{\ptLi }|\NRt |}$. Then, iteration complexity of Algorithm~\ref{alg:gjs_isaega} is $  \max_{j\in \NRt, i,\tR}  \left( 4\frac{m^j_i}{\mu}\left( 1+ \frac{1}{ n\ptLi|\NRt |} \right)  +  \frac{1}{\ptLi|\NRt |} \right)
\log\frac1\epsilon$. 
\end{corollary}

\section{Special cases: {\tt JacSketch} \label{sec:gjs_jacsketch}}
As next special case of {\tt GJS} (Algorithm~\ref{alg:gjs_SketchJac}) we present {\tt JacSketch} ({\tt JS}) motivated by~\cite{jacsketch}. The algorithm observes every iteration a single right sketch of the Jacobian and constructs operators $\cS, \cU$ in the following fashion:
\[
  \cS \mX = \mX \mR 
\quad \text{and} \quad   \cU \mX = \mX \mR \E{\mR}^{-1}
 \]
where $\mR\in \R^{n\times n}$ is random projection matrix.

\begin{algorithm}[!h]
\begin{algorithmic}[1]
\State \textbf{Parameters:} Stepsize $\alpha>0$, Distribution $\cD$ over random projector matrices $\mR\in \R^{n\times n}$
\State \textbf{Initialization:} Choose solution estimate $x^0 \in \R^d$ and Jacobian estimate $ \mJ^0\in \R^{d\times n}$ 
\For{$k =  0, 1, 2,\dots$}
\State Sample realization of $\mR\sim \cD$ perform sketches $\mG(x^k)\mR$ 
\State  $\mJ^{k+1} = \mJ^k - (\mJ^k - \mG(x^k)\mR)$ 
\State $g^k = \frac1n \mJ^k \eR + \frac1n  \left(\mG(x^k) -\mJ^k\right)\mR \E{\mR}^{-1}\eR$  
    \State $x^{k+1} =  \prox_{\alpha \psi}  (x^k - \alpha g^k)$ 
\EndFor
\end{algorithmic}
\caption{{\tt JS} ({\tt JacSketch})}
\label{alg:gjs_jacsketch}
\end{algorithm}

Note that Algorithm~\ref{alg:gjs_jacsketch} differs to what was proposed in~\cite{jacsketch} in the following points. 
\begin{itemize}
\item Approach from~\cite{jacsketch} uses a scalar random variable $\theta_\mR $ to set $\cU \mX = \theta_\mR \mX \mR$. Instead, we set $\E{\cU}= \mX \mR \E{\mR}^{-1}$. This tweak allows Algorithm~\ref{alg:gjs_SketchJac} to recover the tightest known analysis of {\tt SAGA} as a special case. Note that the approach from~\cite{jacsketch} only recovers tight rates for {\tt SAGA} under uniform sampling.
\item Unlike~\cite{jacsketch}, our setup allows for proximable regularizer, thus is more general.
\item Approach from~\cite{jacsketch} allows projections under a general weighted norm. Algorithm~\ref{alg:gjs_SketchJac} only allows for non-weighted norm; which is only done for the sake of simplicity as the chapter is already very notation-heavy. However, {\tt GJS} (Algorithm~\ref{alg:gjs_SketchJac}) is general enough to alow for an arbitrary weighted norm. 
\end{itemize} 

The next corollary shows the convergence result. 
\begin{corollary}[Convergence rate of {\tt JacSketch}]\label{cor:gjs_jacsketh} Suppose that operator $\cM$ is commutative with right multiplication by $\mR$ always. Consider any $\mB\in \R^{n\times n}$ which commutes with $\mR$ always.
Denote
\[
\mM^{\frac12} \eqdef  \begin{pmatrix}
\mM_1^{\frac12} & & \\
& \ddots &\\
& & \mM_n^{\frac12}
\end{pmatrix}\]
and \[ \ugly \eqdef 
\lambda_{\max}\left( 
{\mM^{\frac12}}^\top 
 \left(\E{ \mR \E{\mR}^{-1} \eR \eR^\top \E{\mR}^{-1} \mR} \otimes \mI_{d}\right)
\mM^{\frac12}
 \right),
\]
and let 
  \[
\alpha =\frac{\lambda_{\min} \left(\mB^\top \E{\mR} \mB \right) }{ 4n^{-1}\ugly \lambda_{\max}\left(  \mB^\top \E{\mR} \mB  \right) + \mu\lambda_{\max} \left(\mB^\top \mB\right) }.
\]
Then, the iteration complexity of Algorithm~\ref{alg:gjs_jacsketch} is 

\[
 \frac{ 4n^{-1}\ugly \mu^{-1} \lambda_{\max}\left(  \mB^\top \E{\mR} \mB  \right) + \lambda_{\max} \left(\mB^\top \mB\right) }{\lambda_{\min} \left(\mB^\top \E{\mR} \mB \right) } \log\frac1\epsilon.
 \]

\end{corollary}

 \section{Special cases: proofs}

In this section, we provide the proofs of all corollaries listed in previous sections.  For simplicity, we will use the following notation throughout this section: $\Gamma(\mX)=  \cU(\mX)e$.
 
\subsection{{ \tt SAGA} methods: proofs}
\subsubsection{Setup for Corollary~\ref{cor:gjs_saga}}
Note first that the choice of $\cS, \cU$ yields
\begin{eqnarray*}
\E{\cS(\mX)} &=& \frac1n \mX , \qquad \\
\E{\|\Gamma(\mX) \|^2}  &=& n^2  \E{\< \mX^\top,  \eRj \eRj^\top \eR \eR^\top \eRj \eRj^\top \mX^\top> }  = n\|\mX\|^2.
\end{eqnarray*}

Next, as we have no prior knowledge about $\mG(x^*)$, let $\cR\equiv \cI$; i.e. $\Range{ \cR} = \R^{d\times n}$. Lastly, consider $\cB$ operator to be a multiplication with constant $\beta$: $\cB (\mX)  = \beta \mX$. Thus, for \eqref{eq:gjs_small_step} we should have
\[
\frac{2\alpha}{n}m + \beta^2 \left(1-\frac1n\right) \leq (1-\alpha \mu) \beta^2
\]
and for~\eqref{eq:gjs_small_step2} we should have
\[
\frac{2\alpha}{n} m + \frac{\beta^2}{n} \leq \frac{1}{n}.
\]

It remains to notice that choices $\alpha = \frac{1}{4m + \mu n}$ and $\beta^2 = \frac{1}{2}$ are valid to satisfy the above bounds.

\subsubsection{Setup for Corollary~\ref{cor:gjs_saga_as2} \label{sec:gjs_cor:saga_as2}}
First note that $\E{\cS (\mX)} = \mX \diag (\pR)$. Next, due to~\eqref{eq:gjs_j_in_range},~\eqref{eq:gjs_g_in_range}, inequalities~\eqref{eq:gjs_small_step} and \eqref{eq:gjs_small_step2} with choice $\mY=  {\cM^\dagger}^{\frac12}\mX$ become respectively:
\begin{equation}\label{eq:gjs_linear_1}
\frac{2\alpha}{n^2} \E{\left \|  
 \sum_{j\in R} \pRj^{-1}
\mM_j^{\frac12} \mY_{:j}\right \|^2} + 
\left\|\left(\cI-\E{\cS}\right)^{\frac{1}{2}}\cB(\mY)\right\|^{2} \leq(1-\alpha \mu)\|\cB(\mY)\|^{2}
\end{equation}

\begin{equation}\label{eq:gjs_linear_2}
\frac{2\alpha}{n^2}  \E{\left \|  \sum_{j\in R} \pRj^{-1}
\mM_i^{\frac12} \mY_{:i}\right \|^2}  + 
\left\|\left(\E{\cS}\right)^{\frac{1}{2}}\cB(\mY)\right\|^{2} \leq \frac1n  \|\mY \|^2
\end{equation}

Note that 
\[
\E{\left \|  
 \sum_{j\in R} \pRj^{-1}
\mM_j^{\frac12} \mY_{:j}\right \|^2}  =
  \E{   \left \| \sum_{j\in R} 
\mM_j^{\frac12}(\pRj^{-1} \mY_{:j}) \right \|^2}  
\leq 
\sum_{j=1}^n \pRj^{-1}v_j \|\mY_{:j}\|^2,
\]
where we used ESO assumption~\eqref{eq:gjs_ESO_saga} in the last bound above. Next, choose $\cB$ to be right multiplication with $\diag(b)$. Thus, for~\eqref{eq:gjs_linear_1} it suffices to have for all $j\in [n]$
\[
\frac{2\alpha}{n^2} v_j \pRj^{-1} + b_j^2(1-\pRj) \leq b_j^2 (1-\alpha \mu) \qquad \Rightarrow \qquad
\frac{2\alpha}{n^2} v_j \pRj^{-1} + b_j^2\alpha \mu \leq b_j^2\pRj 
\]
For~\eqref{eq:gjs_linear_2} it suffices to have for all $j\in [n]$
\[
\frac{2\alpha}{n^2} v_j \pRj^{-1} + b_j^2 \pRj \leq \frac{1}{n}
\]
It remains to notice that the choice $b_j^2 = \frac{1}{2n\pRj}$ and $\alpha  = \min_j \frac{n \pRj}{4v_j + n\mu}$ is valid.

\subsection{{\tt SEGA} methods: proofs}

\subsubsection{Setup for Corollary~\ref{cor:gjs_sega}}

 Note that 
\begin{eqnarray*}
\E{\cS x} &=& \frac{1}{d} x, \\
\E{\|\Gamma(x) \|^2}  &=& d^2  \E{\< x,  \eLi \eLi^\top \eLi \eLi^\top x> }  = d\|x\|^2.
\end{eqnarray*}

Next, choose operator $\cB$  to be constant; in particular $\cB x = \beta x$. Thus to satisfy~\eqref{eq:gjs_small_step} it suffices to have
\[
2\alpha d m  +  \beta^2\left(1-\frac1d \right  ) \leq \beta^2(1-\alpha \mu ) \qquad 
\Rightarrow \qquad 
2\alpha d m + \alpha \mu \beta^2 \leq \frac{\beta^2}{d}.
\]
To satisfy~\eqref{eq:gjs_small_step2}, it suffices to have
\[
2\alpha d  m+ \frac{\beta^2}{d} \leq 1.
\]
It remains to notice that $\beta^2 = \frac{d}{2}$ and $\alpha= \frac{1}{4md+ \mu d}$ satisfies the above conditions.

\subsubsection{Setup for Corollary~\ref{cor:gjs_sega_is_11}}
Note that 
$
\E{\cS (x)} = \diag(\pL) x
$
and
\begin{equation*}
\E{\left\|\Gamma(x)\right\|^{2}} = \| x\|^2_{ \E{\sum_{i\in L} \frac{1}{\pLi} \eLi \eLi^\top \sum_{i\in L} \frac{1}{\pLi} \eLi \eLi^\top } }   =  \| x\|^2_{\pL^{-1}}.
\end{equation*}

Let us consider $\cB$ to be the operator corresponding to left multiplication with matrix $\diag(b)$: $\cB(x) = \diag(b)x$. Thus, for~\eqref{eq:gjs_small_step} it suffices to have for all $i$
\[
2\alpha m_i \pLi^{-1} + b_i^2(1-\pLi) \leq b_i^2 (1-\alpha \mu) \qquad \Rightarrow \qquad
2\alpha m_i \pLi^{-1} + b_i^2\alpha \mu \leq b_i^2\pLi .
\]

For~\eqref{eq:gjs_small_step2} it suffices to have for all $i$
\[
2\alpha m_i \pLi^{-1} + b_i^2 \pLi \leq 1 
\]
It remains to notice that choice $b_i^2 = \frac{1}{2\pLi}$ and $\alpha  = \min_i \frac{\pLi}{4m_i+ \mu}$ is valid.

\subsubsection{Setup for Corollary~\ref{cor:gjs_svrcd}}

Note that 
$
\E{\cS (x)} = \probx x
$
and
\begin{equation*}
\E{\left\|\Gamma(x)\right\|^{2}} = \| x\|^2_{ \E{\sum_{i\in L} \frac{1}{\pLi} \eLi \eLi^\top \sum_{i\in L} \frac{1}{\pLi} \eLi \eLi^\top } }   =  \| x\|^2_{\pL^{-1}}.
\end{equation*}

Let us consider $\cB$ to be the operator corresponding to scalar multiplication with $\beta$. Thus, for~\eqref{eq:gjs_small_step} it suffices to have for all $i$
\[
2\alpha w_i \pLi^{-1} + \beta^2(1-\probx) \leq \beta^2 (1-\alpha \mu) \qquad \Rightarrow \qquad
2\alpha w_i  \pLi^{-1} +\beta^2\alpha \mu \leq\beta^2\probx.
\]
For~\eqref{eq:gjs_small_step2} it suffices to have for all $i$,
\[
2\alpha w_i  \pLi^{-1} + \beta^2 \probx \leq 1 .
\]
It remains to notice that choice $\beta^2 = \frac{1}{2\probx}$ and $\alpha  = \min_i \frac{1}{4w_i\pLi^{-1}+ \mu\probx^{-1}}$ is valid.

\subsection{Setup for Corollary~\ref{cor:gjs_sgd}}
Choose $\cB$ to be operator which maps everything into 0. On top of that, by construction we have $\cR = 0$ and thus~\eqref{eq:gjs_small_step} is satisfied for free. Moreover, from~\eqref{eq:gjs_ESO_saga} we have (following the steps from Section~\ref{sec:gjs_cor:saga_as2}):
\[
\E{\|\Gamma(\cM^{\frac12}(\mX))\|^2}\leq \sum_{j=1}^n p_j^{-1}v_j \| \mX_{:j} \|^2.
\]
Further, due to~\eqref{eq:gjs_g_in_range} and~\eqref{eq:gjs_j_in_range}, to satisfy \eqref{eq:gjs_small_step2} we shall have 
\[
\frac{2\alpha}{n^2} \sum_{j=1}^n p_j^{-1}v_j \| \mY_{:j} \|^2 \leq \frac1n \| \mY\|^2,
\]
which simplifies to
\[
\frac{2\alpha}{n} \frac{v_j}{p_j}\leq 1
\]
and thus it suffices to choose $\alpha  =  \frac{n}{2} \min_j \frac{p_j}{v_j} $.

\begin{remark}
Factor 2 can be omitted since for Lemma~\ref{lem:gjs_g_lemma}, the second factor is 0 and thus we no longer need the Jensen's inequality.
\end{remark}

\subsection{Setup for Corollary~\ref{cor:gjs_lsvrg_as}}

First note that $\E{\cS (\mX)} = \probx \mX $. Next, due to~\eqref{eq:gjs_j_in_range},~\eqref{eq:gjs_g_in_range}, inequalities~\eqref{eq:gjs_small_step} and \eqref{eq:gjs_small_step2} with choice $\mY=  {\cM^\dagger}^{\frac12}\mX$ become respectively:
\begin{equation}\label{eq:gjs_linear_1_svrg}
\frac{2\alpha}{n^2} \E{\left \|  
 \sum_{j\in R} \pRj^{-1}
\mM_j^{\frac12} \mY_{:j}\right \|^2} + 
\left\|\left(\cI-\E{\cS}\right)^{\frac{1}{2}}\cB(\mY)\right\|^{2} \leq(1-\alpha \mu)\|\cB(\mY)\|^{2},
\end{equation}
and
\begin{equation}\label{eq:gjs_linear_2_svrg}
\frac{2\alpha}{n^2}  \E{\left \|  \sum_{j\in R} \pRj^{-1}
\mM_i^{\frac12} \mY_{:i}\right \|^2}  + 
\left\|\left(\E{\cS}\right)^{\frac{1}{2}}\cB(\mY)\right\|^{2} \leq \frac1n  \|\mY \|^2.
\end{equation}

Note next that
\[
\E{\left \|  
 \sum_{j\in R} \pRj^{-1}
\mM_j^{\frac12} \mY_{:j}\right \|^2}  =
  \E{   \left \| \sum_{j\in R} 
\mM_j^{\frac12}(\pRj^{-1} \mY_{:j}) \right \|^2}  
\leq 
\sum_{j=1}^n \pRj^{-1}v_j \|\mY_{:j}\|^2,
\]
where we used ESO assumption~\eqref{eq:gjs_ESO_saga} in the last bound above. Next, choose $\cB$ to be multiplication with scalar $\beta$. Thus, for~\eqref{eq:gjs_linear_1_svrg} it suffices to have for all $j\in [n]$
\[
\frac{2\alpha}{n^2} v_j \pRj^{-1} + \beta^2(1-\probx) \leq \beta^2 (1-\alpha \mu) \qquad \Rightarrow \qquad
\frac{2\alpha}{n^2} v_j \pRj^{-1} + \beta^2\alpha \mu \leq \beta^2\probx
\]
For~\eqref{eq:gjs_linear_2_svrg} it suffices to have for all $j\in [n]$
\[
\frac{2\alpha}{n^2} v_j \pRj^{-1} + \beta^2 \probx \leq \frac{1}{n}
\]
It remains to notice that choice $\beta^2 = \frac{1}{2n\probx}$ and $\alpha  = \min_j \frac{n }{4v_j\pRj^{-1} + n\mu \probx^{-1}}$ is valid.

\subsection{Methods with Bernoulli $\cU$: proofs}

\subsubsection{Setup for Corollary~\ref{cor:gjs_B2}}

Note first that the choice of $\cS, \cU$ yield
\begin{eqnarray*}
\E{\cS(x)} &=& \probx x,  \\
\E{\|\Gamma(x) \|^2}  &=& \E{\| \cU x\|^2 }  = \proby^{-1}\|x \|^2.
\end{eqnarray*}

Next, consider $\cB$ operator to be a multiplication with a constant $ b$.

Thus for \eqref{eq:gjs_small_step} we should have

\[
2\alpha\proby^{-1}L + b^2 \left(1-\probx\right) \leq (1-\alpha \mu) b^2
\]
and for~\eqref{eq:gjs_small_step2} we should have
\[
2\alpha \proby^{-1} L + \probx b^2 \leq 1.
\]

It remains to notice that choices $\alpha = \frac{1}{4\proby^{-1}L + \mu \probx^{-1}}$ and $b^2 = \frac{1}{2 \probx}$ are valid to satisfy the above bounds.

\subsubsection{Setup for Corollary~\ref{cor:gjs_inverse_svrg}}

Note first that the choice of $\cS, \cU$ yields
\begin{eqnarray*}
\E{\cS(\mX)} &=&  \mX \diag(\pR) \\
\E{\|\Gamma(\mX) \|^2}  &=& \E{\| \cU(\mX)e\|^2 }  = \proby^{-1}\|\mX e\|^2 \leq  \proby^{-1}n \|\mX \|^2
\end{eqnarray*}

Next, as we have no prior knowledge about $\mG(x^*)$, consider $\cR$ to be identity operator; i.e. $\Range{ \cR} = \R^{d\times n}$. Lastly, consider $\cB$ operator to be a right multiplication with $ \diag(b)$.

Thus for \eqref{eq:gjs_small_step} we should have

\[
\forall j:  \quad 
\frac{2\alpha}{n}\proby^{-1}m +\alpha \mu b^2_j  \leq b^2_j \pRj
\]
and for~\eqref{eq:gjs_small_step2} we should have
\[
\forall j: \quad  \frac{2\alpha}{n} \proby^{-1} m + \pRj b^2_j \leq \frac{1}{n}
\]

It remains to notice that choices $\alpha = \min_j \frac{1}{4\proby^{-1}m + \mu \pRj^{-1}}$ and $b^2_j = \frac{1}{2n \pRj}$ are valid to satisfy the above bounds.

\subsubsection{Setup for Corollary~\ref{cor:gjs_SVRCD-inv}}

Note first that the choice of $\cS, \cU$ yields
\begin{eqnarray*}
\E{\cS(x)} &=& \pL \circ x,  \\
\E{\|\Gamma(x) \|^2}  &=& \E{\| \cU(x)\|^2 }  = \proby^{-1}\|x \|^2 .
\end{eqnarray*}

Next, as we have no prior knowledge about $\mG(x^*)$, consider $\cR$ to be identity operator; i.e. $\Range{ \cR} = \R^{d\times n}$. Lastly, consider $\cB$ operator to be left multiplication with matrix $\diag( b)$.

Thus, for \eqref{eq:gjs_small_step} we should have
$2\alpha\proby^{-1}m + b^2_i \alpha \mu \leq b^2_i \pLi
$ for all $i$,
and for~\eqref{eq:gjs_small_step2}, we should have
$
2\alpha\proby^{-1}m + \pLi b^2_i \leq 1
.$

It remains to notice that the choices $\alpha = \min_i \frac{1}{4\proby^{-1}m + \mu \pLi^{-1}}$ and $b^2_i= \frac{1}{2 \pLi^{-1}}$ are valid to satisfy the above bounds.

\subsection{Combination of left and right sketches: proofs}

\subsubsection{Setup for Corollary~\ref{cor:gjs_RL}}
Note first that the choice of $\cS, \cU$ yields $\E{\cS(\mX)} =\mX \diag(\pR)$ and
\begin{eqnarray*}
\E{\|\piop(\mX) \|^2}  &=& \| \cM^{\frac12}(\mX)\eR \|^2_{\diag({\pL}^{-1})} \leq n\sum_{j=1}^n \| \mM_j \mX_{:j}\|^2_{\diag({\pL}^{-1})} = n \sum_{j=1}^n \|  \mX_{:j}\|^2_{ \diag(m^j\circ {\pL}^{-1})}.
\end{eqnarray*}

Let $\cB$ be right multiplication by $\diag(b)$.
Thus for \eqref{eq:gjs_small_step} we should have

\[
\forall i,j: \quad 2\frac{\alpha}{n}{m_i}^j {\pLi}^{-1} + b^2_j \alpha \mu \leq b^2_j \pRj
\]
and for~\eqref{eq:gjs_small_step2} we should have
\[
\forall i,j: \quad 2\frac{\alpha}{n}{m_i}^j {\pLi}^{-1}+ \pRj b^2_j \leq \frac1n.
\]

It remains to notice that choices $\alpha = \min_{i,j} \frac{1}{4{m_i}^j {\pLi}^{-1} + \mu \pRj^{-1}}$ and $b^2_j=\frac{1}{ 2\pRj n}$ are valid to satisfy the above bounds. 

\subsubsection{Setup for Corollary~\ref{cor:gjs_LR}}
Note first that the choice of $\cS, \cU$ yields

\begin{eqnarray*}
\E{\cS(\mX)} &=&    \diag(\pL)  \mX  \\
\E{\|\Gamma(\mX) \|^2}  &\leq & \sum_{j=1}^n \pRj^{-1} v_j \|\mX_{:j}\|^2
\end{eqnarray*}
The second inequality is a direct consequence of ESO (which is shown is Section~\ref{sec:gjs_cor:saga_as2}).

Let $\cB$ be left multiplication by $\diag(b)$. Thus for \eqref{eq:gjs_small_step} we should have

\[
\forall i,j: \quad 2\frac{\alpha}{n}v_j {\pRj}^{-1} + b^2_i \alpha \mu \leq b^2_i \pLi
\]
and for~\eqref{eq:gjs_small_step2} we should have
\[
\forall i,j: \quad  2\frac{\alpha}{n}v_j{\pRj}^{-1} + \pLi b^2_i \leq \frac1n .
\]

It remains to notice that choices $\alpha = \min_{i,j} \frac{1}{4{v_j} {\pRj}^{-1} + \mu \pLi^{-1}}$ and $b^2_i = \frac{1}{ 2\pLi n}$ are valid to satisfy the above bounds.

\subsection{Joint sketches: proofs}

\subsubsection{Setup for Corollary~\ref{cor:gjs_saega} \label{sec:gjs_corsaega}}

For notational simplicity, denote ${\mM'}^{\frac12}\in \R^{d\times n}$ to be the matrix with $j$th column equal to (elementwise) square root of $m_j$. We have 
\[
\E{\cS(\mX)}  = \left( {\pL}{\pR}^\top\right)\circ \mX
\]
and 

\begin{eqnarray}
 \nonumber
 &&
  \E{\|  \piop(\cM^{\frac12 }\mX)\|^2} 
  \\ \nonumber
   &   = &
 \E{ \left\| \left( \left( {\pL}^{-1}{\pR^{-1}}^\top\right) \circ \left( \left( \sum_{i\in L} \eLi \eLi^\top \right) (\mM^{'\frac12}\circ\mX) \left( \sum_{j\in R} \eRj \eRj^\top \right) \right) \right)\eR \right\|^2} 
 \\ \nonumber
 &   = & 
  \E{ \left\|   \left( \left( \sum_{i\in L, j\in R} \eLi \eRj^\top \right) \circ\left( {\pL}^{-1}{\pR^{-1}}^\top\right) \circ \mM^{'\frac12}\circ\mX \right)\eR \right\|^2}
   \\ \nonumber
    &   = & 
  \E{ \left\|  \left( \sum_{i\in L} \eLi \eLi^\top \right)  \left( \left( {\pL}^{-1}{\pR^{-1}}^\top\right) \circ \mM^{'\frac12}\circ\mX\right)\eRR \right\|^2}
   \\ \nonumber
       &   = & 
  \E{ \left\|  \left(     \left( {\pL}^{-\frac12}{\pR^{-1}}^\top\right) \circ \mM^{'\frac12}\circ\mX \right)\eRR \right\|^2}
   \\ \nonumber
             &   = & 
  \E{  \tracee\left(   \left(     \left( {\pL}^{-\frac12}{\pR^{-1}}^\top\right) \circ \mM^{'\frac12}\circ\mX \right)\mI_{R,R} \left(     \left( {\pL}^{-\frac12}{\pR^{-1}}^\top\right)^\top \circ {\mM^{'\frac12}}^\top\circ\mX^\top \right) \right)}
   \\ \nonumber
             &   = & 
\tracee\left(   \left(     \left( {\pL}^{-\frac12}{\pR^{-1}}^\top\right) \circ \mM^{'\frac12}\circ\mX \right)\PR \left(     \left( {\pL}^{-\frac12}{\pR^{-1}}^\top\right)^\top \circ  {\mM^{'\frac12}}^\top\circ\mX^\top \right) \right)
   \\ \nonumber
                &   = & 
\tracee\left(   \left(     \left( {\pL}^{-\frac12}{\eR}^\top\right) \circ \mM^{'\frac12}\circ\mX \right) \diag(\pR)^{-1}\PR \diag(\pR)^{-1} \left(     \left( {\pL}^{-\frac12}{\eR}^\top\right)^\top \circ {\mM^{'\frac12}}^\top\circ\mX^\top \right) \right)
   \\ \nonumber
                   &   \leq & 
\tracee\left(   \left(     \left( {\pL}^{-\frac12}{\eR}^\top\right) \circ \mM^{'\frac12}\circ\mX \right)  \diag(\qR)^{-1} \left(     \left( {\pL}^{-\frac12}{\eR}^\top\right)^\top \circ {\mM^{'\frac12}}^\top\circ\mX^\top \right) \right)
   \\
                      &   = & 
\left\| \mX \circ \mM^{'\frac12} \circ  \left( {\pL}^{-\frac12}{\qR^{-\frac12}}^\top \right) \right\|^2.
 \label{eq:gjs_SAEGA_ESO}
\end{eqnarray}

Next, choose operator $\mB$ to be such that $\cB(\mX) \eqdef \mB \circ \mX$ for $\mB\in \R^{d\times n}$. 
Thus, for~\eqref{eq:gjs_small_step} and \eqref{eq:gjs_small_step2} we shall have respectively
\[
\forall i, j: \quad 
\frac{2\alpha}{n^2}\left( \frac{m^j_i}{\pLi \qRj}\right)+ \mB_{i,j}^2\alpha \mu \leq \mB_{ij}^2 \pLi \qRj
\]
and
\[
\forall i, j: \quad 
\frac{2\alpha}{n^2}\left( \frac{m^j_i}{\pLi \qRj}\right)+ \mB_{ij}^2 \pLi \qRj \leq \frac{1}{n}.
\]

It remains to choose $\mB_{i,j}^2= \frac{1}{2n \pLi \qRj}$ and $\alpha = \min_{i,j} \frac{n\pLi \qRj}{4m^j_i + n\mu}$.

\subsubsection{Setup for Corollary~\ref{cor:gjs_svrcdg}}

We have 
\[
\E{\cS(\mX)}  =\probx\mX.
\]

Next, choose operator $\mB$ to be such that $\cB(\mX) \eqdef \beta \circ \mX$ for scalar $\beta$ which would be specified soon. 
Proceeding with bound~\eqref{eq:gjs_SAEGA_ESO}, for~\eqref{eq:gjs_small_step} and \eqref{eq:gjs_small_step2} we shall have respectively
\[
\forall i, j: \quad 
\frac{2\alpha}{n^2}\left( \frac{m^j_i}{\pLi \qRj}\right)+ \beta^2\alpha \mu \leq \beta^2 \probx
\]
and
\[
\forall i, j: \quad 
\frac{2\alpha}{n^2}\left( \frac{m^j_i}{\pLi \qRj}\right)+ \beta^2 \probx \leq \frac{1}{n}.
\]

It remains to choose $\beta^2= \frac{1}{2n \probx}$ and $\alpha = \min_{i,j} \frac{1}{4\frac{m^j_i}{n\pLi \qRj} + \probx^{-1}\mu}$.

\subsubsection{Setup for Corollary~\ref{cor:gjs_isaega}}

For notational simplicity, denote $\mM'\in \R^{d\times n}$ to be a matrix with $j$th column equal to $m_j$.

Let $  \piop_\tR(\mX_{:\NRt}) =  \left( {{\ptL}^{-1}} {{\ptR}^{-1}}^\top\right) \circ \left( \left( \sum_{i\in L_\tR} \eLi \eLi^\top \right)\mX_{:\NRt} \left( \sum_{j\in R_\tR} \eRj \eRj^\top \right)\right) e_{\NRt}$. 
Thus 
\[
\E{\cS(\mX)}  = 
\sum_{\tR=1}^\TR
 \left( {{\ptL}} {{\ptR}}^\top\right) \circ\mX_{:\NRt}  
\]
and 
\begin{eqnarray}
\E{\|  \piop(\mX)\|^2}   &=& \E{ \left\|   \sum_{\tR=1}^\TR  \piop_\tR(\mX_{:\NRt}) \right \|^2 } 
\\ \nonumber
&=&
\E{ \left\| \sum_{\tR=1}^\TR  \piop_\tR(\mX_{:\NRt})  -  \E{\sum_{\tR=1}^\TR  \piop_\tR(\mX_{:\NRt})  } \right \|^2 } 
+ 
 \left\|   \E{\sum_{\tR=1}^\TR  \piop_\tR(\mX_{:\NRt})  } \right \|^2 
 \\ \nonumber
 &=&
 \E{ \left\| \sum_{\tR=1}^\TR\left(    \piop_\tR(\mX_{:\NRt}) - \mX_{:\NRt}e_{\NRt}  \right)\right \|^2 } 
+ 
 \left\|   \mX \eR \right \|^2 
 \\ \nonumber
  &=&
 \sum_{\tR=1}^\TR  \E{ \left\|     \piop_\tR(\mX_{:\NRt}) - \mX_{:\NRt}e_{\NRt}  \right \|^2 } 
+ 
 \left\|     \mX \eR \right \|^2 
 \\ \nonumber
   &\leq&
 \sum_{\tR=1}^\TR\E{ \left\|   \piop_\tR(\mX_{:\NRt})\right \|^2 } 
+ 
 \left\|    \mX \eR \right \|^2 
 \\ \label{eq:gjs_ISEAGA_ESO}
    &\leq&
 \sum_{\tR=1}^\TR \E{ \left\|   \piop_\tR(\mX_{:\NRt}) \right \|^2 } 
+ 
 n \left\|  \mX \right \|^2.
\end{eqnarray}

Using the bounds from Section~\ref{sec:gjs_corsaega} we further get
\begin{eqnarray*}
  \E{\|  \piop(\cM^{\frac12 }\mX)\|^2} 
  &  \stackrel{\eqref{eq:gjs_ISEAGA_ESO}+\eqref{eq:gjs_SAEGA_ESO}}{\leq}&
 \sum_{\tR=1}^\TR   \left\|   \mX_{:\NRt} \circ \left( {\ptL}^{-\frac12} {{\qtR}^{-\frac12}}^\top \right) \circ \mM'_{:\NRt} \right \|^2 
+ 
 n \left\| \mM' \circ \mX \right \|^2. \\
\end{eqnarray*}

Next, choose operator $\mB$ to be such that for any $\mX$: $\cB(\mX) \eqdef \mB \circ \mX$ where $\cB \in \R^{d\times n}$. 
Thus, for~\eqref{eq:gjs_small_step} and \eqref{eq:gjs_small_step2} we shall have respectively
\[
\forall i, \tR, j\in \NRt: \quad 
\frac{2\alpha}{n^2}m^j_i \left( \frac{1}{ \ptLi \qtRj}+ n\right)+ \mB_{i,j}^2\alpha \mu \leq \mB_{i,j}^2 \ptLi \qtRj
\]
and
\[
\forall i, \tR, j\in \NRt: \quad 
\frac{2\alpha}{n^2}m^j_i \left( \frac{1}{ \ptLi \qtRj}+ n\right)+ \mB_{i,j}^2\ptLi \qtRj \leq \frac{1}{n}.
\]

It remains to choose $\mB_{i,j}^2= \frac{1}{2n \ptL \qtRj }$ and $\alpha = \min_{j\in \NRt, i,\tR} \frac{1}{4 m^j_i \left( 1+ \frac{1}{n\ptL \qtRj}\right) + \frac{\mu}{\ptL \qtRj}}$.

\subsection{Setup for Corollary~\ref{cor:gjs_jacsketh}}
Let $x$ be column-wise vectorization of $\mX$. 
Note that 
\[
\Gamma(\cM^\frac12(\mX)) = \cM^\frac12(\mX) \mR\E{\mR}^{-1}\eR = \left( \eR^\top \E{\mR}^{-1} \mR \otimes \mI_{d}\right)\begin{pmatrix}
\mM_1^{\frac12} & & \\
& \ddots &\\
& & \mM_n^{\frac12}
\end{pmatrix}x. 
\]
Thus,  $\E{\left\| \Gamma(\cM^\frac12(\mX))\right\|^2 } \leq 
\|\mX\|^2 \ugly.$ Let $\cB(\mX) = \beta \mX \mB$. Thus, we have 
\begin{eqnarray*}
 && 
  (1-\alpha \mu) \NORMG{ \cB\mX } - \NORMG{  \left(\cI - \E{\cS} \right)^{\frac12}\cB  \mY} 
  \\
 &=& \beta^2 \tracee\left(\mX \mB^\top ( \E{\mR} - \alpha\mu\mI) \mB\mX^\top \right) \\
 &\leq &\beta^2 \lambda_{\min} \left(\mB^\top ( \E{\mR} - \alpha\mu\mI) \mB \right) \|\mX\|^2 \\
  &\leq &\beta^2 \left( \lambda_{\min} \left(\mB^\top \E{\mR} \mB \right) -  \alpha\mu\lambda_{\max} \left(\mB^\top \mB\right)\right) \|\mX\|^2.
\end{eqnarray*}

Further, 
\[
\NORMG{\left(\E{\cS}\right)^{\frac12}  \cB  {\cM^\dagger}^{\frac12}\mX  }   =  \beta^2\tracee \left( \mX \mB^\top \E{\mR} \mB \mX^\top\right) \leq  \beta^2\|\mX\|^2 \lambda_{\max}\left(  \mB^\top \E{\mR} \mB  \right).
\]

Using the derived bounds together with~\eqref{eq:gjs_j_in_range},~\eqref{eq:gjs_g_in_range}, for conditions~\eqref{eq:gjs_small_step} and \eqref{eq:gjs_small_step2} it suffices to have:
\begin{equation}\label{eq:gjs_linear_1_js}
\frac{2\alpha}{n^2}\ugly + \beta^2\alpha\mu\lambda_{\max} \left(\mB^\top \mB\right)
  \leq 
\beta^2\lambda_{\min} \left(\mB^\top \E{\mR} \mB \right)  ,
\end{equation}
and
\begin{equation}\label{eq:gjs_linear_2_js}
\frac{2\alpha}{n^2}\ugly + \beta^2 \lambda_{\max}\left(  \mB^\top \E{\mR} \mB  \right) \leq \frac1n.
\end{equation}
It remains to notice that choices $\beta^2 = \frac{1}{2n \lambda_{\max}\left(  \mB^\top \E{\mR} \mB  \right)}$ and 
\[
\alpha =\frac{\lambda_{\min} \left(\mB^\top \E{\mR} \mB \right) }{ 4n^{-1}\ugly \lambda_{\max}\left(  \mB^\top \E{\mR} \mB  \right) + \mu\lambda_{\max} \left(\mB^\top \mB\right) }
\]
are valid.

\section{Convergence under strong growth condition \label{sec:gjs_sg}}

In this section, we extend the result of Algorithm~\ref{alg:gjs_SketchJac} to the case when $F\eqdef f+\psi$ satisfies a strong growth condition instead of quasi strong convexity. Note that strong growth is weaker (more general) than quasi strong convexity~\cite{karimi2016linear}.

Suppose that $\cX^*$ is a set of minimizers of convex function $F$. Clearly, $\cX^*$ must be convex. Define $[x]^*$ to be a projection of $x$ onto $\cX^*$.
\begin{assumption}\label{as:gjs_strong_growth}

Suppose that $F$ satisfies strong growth, i.e. for every $x$: 
\begin{equation} \label{eq:gjs_strong_growth}
F(x)- F([x]^*) \geq \frac{\mu}{2}\|x-[x]^* \|^2.
\end{equation}
\end{assumption}

\subsection{Technical proposition and lemma}

In order to establish the convergence results, it will be useful to establish Proposition~\ref{prop:41717892449y8} and Lemma~\ref{lem:gjs_gminusnablafk}.

\begin{proposition} (\cite{proxsvrg, qian2019saga}) \label{prop:41717892449y8}
Let $f$ be $\mM$-smooth and suppose that~\eqref{eq:gjs_strong_growth} holds. 
Suppose that $x^{k+1} = x^k  - \alpha g^k$ where $\E{g^k} = \nabla f(x^k)$ and $\alpha \leq \frac{1}{3\lambda_{\max}(\mM)}$. Then
 $$
\mathbb{E}_{k}\left[\left\|x^{k+1}-\left[x^{k+1}\right]^{*}\right\|^{2}\right] \leq \frac{1}{1+\mu \alpha} \mathbb{E}_{k}\left[\left\|x^{k}-\left[x^{k}\right]^{*}\right\|^{2}\right]+\frac{2 \alpha^{2}}{1+\mu \alpha} \mathbb{E}_{k}\left[\left\|g^{k}-\nabla f\left(x^{k}\right)\right\|^{2}\right].
$$
\end{proposition}

\begin{lemma}\label{lem:gjs_gminusnablafk}
For any $x^* \in \cX^*$ we have
\begin{equation} \label{eq:gjs_gminusnablafk}
\E{\| g^k - \nabla f(x^k)\|^2} \leq \frac{2}{n^2} \E{\left \| \cU (\mG(x^*)-\mJ^k)\eR \right\|^2} +\frac{2}{n^2}\E{\left\| \cU(\mG(x^k)-\mG(x^*))\eR  \right\|^2}.
\end{equation}
\end{lemma}

\begin{proof}
\begin{eqnarray*}
&& \E{\| g^k - \nabla f(x^k)\|^2}  \\
&=& \E{\left \| \frac1n \mJ^k \eR - \frac1n\cU(\mG(x^k)-\mJ^k)\eR - \frac1n\mG(x^k)\eR \right\|^2} \\
&=& \frac{1}{n^2} \E{ \left \| ( \mJ^k - \mG(x^*)) \eR- \cU(\mG(x^*)-\mJ^k)\eR + \cU(\mG(x^k)-\mG(x^*))\eR + (\mG(x^*)- \mG(x^k))\eR \right \|^2} \\
&\leq & 
\frac{2}{n^2} \E{\left \| (\mJ^k  - \mG(x^*))\eR- \cU(\mG(x^*)-\mJ^k)\eR \right\|^2}\\
&& \qquad+\frac{2}{n^2}\E{\left\| \cU(\mG(x^k)-\mG(x^*))\eR + (\mG(x^*)- \mG(x^k))\eR \right\|^2}\\
&\leq &
\frac{2}{n^2} \E{\left \| \cU (\mG(x^*)-\mJ^k)\eR \right\|^2} +\frac{2}{n^2}\E{\left\| \cU(\mG(x^k)-\mG(x^*))\eR  \right\|^2}.
\end{eqnarray*}
\end{proof}

Lastly, it is necessary to assume the null space consistency of solution set $\cX^*$ under $\mM$ smothness. A similar assumption was considered in~\cite{qian2019saga}.
\begin{assumption} \label{as:gjs_null_consistency}
For any $x^*, y^* \in \cX$ we have 
\begin{equation}\label{eq:gjs_null_consistency}
{\cM^\dagger}^\frac12\mG(x^*) = {\cM^\dagger}^\frac12\mG(y^*).
\end{equation}
\end{assumption}

\subsection{Convergence proof}
We next state the convergence result of Algorithm~\ref{alg:gjs_SketchJac} under strong growth condition.

\begin{theorem} \label{thm:gjs_strong}

Suppose that~\eqref{eq:gjs_strong_growth} holds. Let $\cB$ be any linear operator commuting with $\cS$, and assume ${\cM^\dagger}^{\frac{1}{2}}$ commutes with $\cS$. Let $\cR$ be any linear operator for which $\cR(\mJ^k) = \cR(\mG(x^*))$ for every $k\geq 0$. Define the Lyapunov function $\Psi^k$ as per~\eqref{eq:gjs_Lyapunov} for any $x^*\in \cX^*$. Suppose that $\alpha\leq \frac{1}{\lambda_{\max}(\mM)}$ and $\cB$ are chosen so that 
\begin{eqnarray} \nonumber 
&&  \frac{2\alpha}{n^2} \left(  \frac{3+ \mu \alpha}{1+\mu\alpha}\right) \E{ \norm{ \cU  \mX \eR }^2 }   +   \NORMG{  \left(\cI - \E{\cS} \right)^{\frac12}\cB  {\cM^\dagger}^{\frac12} \mX } \\
& & \qquad \qquad \qquad \leq \left(1-\frac{\alpha \mu}{2+ 2\alpha \mu}\right)  \NORMG{ \cB {\cM^\dagger}^{\frac12}\mX }\label{eq:gjs_sg_small_step}
\end{eqnarray}
whenever  $ \mX\in \Range{\cR}^\perp$ and
\begin{eqnarray}
\frac{2\alpha}{n^2}\left(  \frac{3+ \mu \alpha}{1+\mu\alpha}\right) \E{  \norm{ \cU  \mX  \eR }^2  } 
+   \NORMG{\left(\E{\cS}\right)^{\frac12}  \cB  {\cM^\dagger}^{\frac12}\mX  }   &  \leq & \frac{1}{n} \norm{{\cM^\dagger}^{\frac12}\mX}^2
\label{eq:gjs_sg_small_step2}
\end{eqnarray}
for all $\mX\in \R^{d\times n}$. Then  for all $k\geq 0$, we have
\[ \E{\Psi^{k}}\leq \left(1-\frac{\alpha \mu}{2+ 2\alpha \mu}\right)^k \Psi^{0}.\]

\end{theorem}
\begin{proof}

Consider any $x^*\in \cX^*$. Due to non-expansiveness of the prox operator we have
\begin{eqnarray*}
\E{\norm{x^{k+1} -[x^{k+1}]^*}_2^2 } &\leq&
\E{\norm{x^{k+1} -[x^{k}]^*}_2^2 }\\
&\overset{ \eqref{eq:gjs_prox_opt}}{=} &
 \E{\norm{ \prox_{\alpha \psi} (x^k-\alpha g^k) -  \prox_{\alpha \psi} ([x^k]^*-\alpha \nabla f([x^k]^*))  }_2^2}  
 \\
 &\leq&
 \E{\norm{x^k - \alpha g^{k} -( [x^k]^* -\alpha \nabla f([x^k]^*) ) }_2^2}  
\\
& =& 
 \norm{x^k  -[x^k]^*}_2^2 -2\alpha \dotprod{\nabla f(x^k)- \nabla f([x^k]^*) , x^k  -[x^k]^*}   \\
 && \qquad + \alpha^2\E{\norm{g^{k}- \nabla f([x^k]^*) }_2^2}
 \\
& \stackrel{\eqref{eq:gjs_smooth_dotprod}}{\leq}& 
 \norm{x^k  -[x^k]^*}_2^2 -\frac{2\alpha}{n} \left\| {\cM^{\dagger}}^\frac12( \mG(x^k)- \mG([x^k]^*)) \right \|^2 \\
 && \qquad + \alpha^2\ED{\norm{g^{k}- \nabla f([x^k]^*) }_2^2}.
\end{eqnarray*}

Combining the above bound with Proposition~\ref{prop:41717892449y8} yields
\begin{eqnarray*}
&& \E{\norm{x^{k+1} -[x^{k+1}]^*}_2^2 }\\
&\leq & 
\left(\frac{1}{2+2\alpha \mu}+\frac12\right)\left\|x^{k}-\left[x^{k}\right]^{*}\right\|^{2}-\frac{\alpha}{n}\left\| {\cM^{\dagger}}^\frac12( \mG(x^k)- \mG([x^k]^*)) \right\|^2 \\
&& 
+\frac12\alpha^{2}\E{\left\|g^{k}-\nabla f\left(\left[x^{k}\right]^{*}\right)\right\|^{2}
}+\frac{ \alpha^{2} }{1+\mu \alpha} \E{\left\|g^{k}-\nabla f\left(x^{k}\right)\right\|^{2}}
\\
&\leq & 
\left(\frac{\alpha \mu +2}{2+2\alpha \mu}\right)\left\|x^{k}-\left[x^{k}\right]^{*}\right\|^{2}-\frac{\alpha}{n} \left\| {\cM^{\dagger}}^\frac12( \mG(x^k)- \mG([x^k]^*)) \right\|^2 
\\
&& 
+\frac12\alpha^{2}\E{\left\|g^{k}-\nabla f\left(\left[x^{k}\right]^{*}\right)\right\|^{2}}+\frac{ \alpha^{2} }{1+\mu \alpha}\E{\left\|g^{k}-\nabla f\left(x^{k}\right)\right\|^{2}}
\\
&\stackrel{\eqref{eq:gjs_gminusnablafk}}{\leq} & 
\left(\frac{\alpha \mu +2}{2+2\alpha \mu}\right)\left\|x^{k}-\left[x^{k}\right]^{*}\right\|^{2}-\frac{\alpha}{n}\left\| {\cM^{\dagger}}^\frac12( \mG(x^k)- \mG([x^k]^*))\right\|^2 
\\
&& 
+\frac{ 2\alpha^{2} }{n^2(1+\mu \alpha)} \left(\E{\left \| \cU(\mG(x^*)-\mJ^k)\eR \right\|^2}  +  \E{\left\| \cU(\mG(x^k)-\mG(x^*))\eR  \right\|^2} \right) 
\\
&&
+\frac12\alpha^{2}\E{\left\|g^{k}-\nabla f\left(\left[x^{k}\right]^{*}\right)\right\|^{2}}
\\
&\stackrel{\eqref{eq:gjs_g_lemma}}{\leq} & 
\left(\frac{\alpha \mu +2}{2+2\alpha \mu}\right)\left\|x^{k}-\left[x^{k}\right]^{*}\right\|^{2}-\frac{\alpha}{n} \left\| {\cM^{\dagger}}^\frac12 ( \mG(x^k)- \mG([x^k]^*)) \right\|^2 
\\
&& 
+\frac{\alpha^{2}}{n^2}\left( \frac{ 2 }{1+\mu \alpha}+1\right)\left(\E{\left \| \cU(\mG(x^*)-\mJ^k)\eR \right\|^2}  +  \E{\left\| \cU(\mG(x^k)-\mG(x^*))\eR  \right\|^2} \right)
\\
&\stackrel{\eqref{eq:gjs_null_consistency}}{\leq} & 
\left(\frac{\alpha \mu +2}{2+2\alpha \mu}\right)\left\|x^{k}-\left[x^{k}\right]^{*}\right\|^{2}-\frac{\alpha}{n} \left\| {\cM^{\dagger}}^\frac12 ( \mG(x^k)- \mG(x^*)) \right\|^2 
\\
&& 
+\frac{\alpha^{2}}{n^2}\left( \frac{ 2 }{1+\mu \alpha}+1\right)\left(\E{\left \| \cU(\mG(x^*)-\mJ^k)\eR \right\|^2}  +  \E{\left\| \cU(\mG(x^k)-\mG(x^*))\eR  \right\|^2} \right).
\end{eqnarray*}
 Since, by assumption,  both $\cB$ and ${\cM^\dagger}^{\frac12}$ commute with $\cS$, so does their composition $\cA \eqdef \cB {\cM^\dagger}^{\frac12}$. Applying Lemma~\ref{lem:gjs_nb98gd8fdx}, we get
 \begin{eqnarray}\label{eq:gjs_J_jac_bound9999}\nonumber
\E{ \NORMG{\cB {\cM^\dagger}^{\frac12} \left(\mJ^{k+1}-\mG(x^*)  \right) }}  &=&  \NORMG{ (\cI - \E{\cS})^{\frac12}  \cB {\cM^\dagger}^{\frac12} \left(\mJ^k-\mG(x^*) \right) }  \\
&& +  \NORMG{\E{\cS}^{\frac12}  \cB {\cM^\dagger}^{\frac12} \left(\mG(x^k)-\mG(x^*) \right) } .
\end{eqnarray}

 Adding $\alpha$ multiple of~\eqref{eq:gjs_J_jac_bound9999} to the previous bounds yields
 \begin{eqnarray*}
&& \E{\norm{x^{k+1} -[x^{k+1}]^*}_2^2 }  + \alpha\E{ \NORMG{\cB {\cM^\dagger}^{\frac12} \left(\mJ^{k+1}-\mG(x^*)  \right) }}    \\
&\leq & 
\left(1- \frac{\alpha \mu }{2+2\alpha \mu}\right)\left\|x^{k}-\left[x^{k}\right]^{*}\right\|^{2}-\frac{\alpha}{n} \left\| {\cM^{\dagger}}^\frac12 ( \mG(x^k)- \mG(x^*)) \right\|^2 
\\
&& 
+\frac{\alpha^{2}}{n^2}\left( \frac{ 3+ \mu\alpha }{1+\mu \alpha}\right)\left(\E{\left \| \cU(\mG(x^*)-\mJ^k) \eR\right\|^2}  +  \E{\left\| \cU(\mG(x^k)-\mG(x^*)) \eR \right\|^2} \right)
\\
&& 
+ \alpha \NORMG{ (\cI - \E{\cS})^{\frac12}  \cB {\cM^\dagger}^{\frac12} \left(\mJ^k-\mG(x^*) \right) }   + \alpha\NORMG{\E{\cS}^{\frac12}  \cB {\cM^\dagger}^{\frac12} \left(\mG(x^k)-\mG(x^*) \right) }
\\
&\stackrel{\eqref{eq:gjs_sg_small_step2}}{\leq} & 
\left(1- \frac{\alpha \mu }{2+2\alpha \mu}\right)\left\|x^{k}-\left[x^{k}\right]^{*}\right\|^{2}
+\frac{\alpha^{2}}{n^2}\left( \frac{ 3+ \mu\alpha }{1+\mu \alpha}\right)\E{\left \| \cU(\mG(x^*)-\mJ^k)\eR \right\|^2}  
\\
&& 
+ \alpha \NORMG{ (\cI - \E{\cS})^{\frac12}  \cB {\cM^\dagger}^{\frac12} \left(\mJ^k-\mG(x^*) \right) }
\\
&\stackrel{\eqref{eq:gjs_sg_small_step}}{\leq} & 
\left(1- \frac{\alpha \mu }{2+2\alpha \mu}\right) \left( \left\|x^{k}-\left[x^{k}\right]^{*}\right\|^{2}+  \alpha\NORMG{\cB {\cM^\dagger}^{\frac12} \left(\mJ^{k}-\mG(x^*)  \right) } \right).
\end{eqnarray*}
\end{proof}

\begin{remark}
Since $2+2\alpha \mu = \cO(1)$ and $\frac{3\mu \alpha }{1+ \mu\alpha} =\cO(1)$ the convergence rate under strong growth provided by Theorem~\ref{thm:gjs_strong} is of the same order as the convergence rate under quasi strong convexity (Theorem~\ref{thm:gjs_main}).  
\end{remark}


\chapter{Appendix for Chapter \ref{sigmak}}
\label{sigmak_appendix}

\graphicspath{{sigmak/images/}}

\section{Special cases}\label{sec:sk_special_cases}

\subsection{Proximal {\tt SGD} for stochastic optimization}\label{sec:sk_SGD}

\begin{algorithm}[h]
    \caption{{\tt SGD}}
    \label{alg:sk_sgd_prox}
    \begin{algorithmic}
        \Require learning rate $\alpha>0$, starting point $x^0\in\R^d$, distribution $\cD$ over $\xi $
        \For{ $k=0,1,2,\ldots$ }
        \State{Sample $\xi \sim \cD$}
        \State{$g^k = \nabla f_\xi (x^k)$}
        \State{$x^{k+1} = \proxR(x^k - \alpha g^k)$}
        \EndFor
    \end{algorithmic}
\end{algorithm}

We start with stating the problem, the assumptions on the objective and on the stochastic gradients for {\tt SGD} \cite{nguyen2018sgd}. Consider the expectation minimization problem
\begin{equation}\label{eq:sk_main_sgd}
    \min_{x\in\R^d} f(x) + \psi(x),\quad f(x) \eqdef \ED{f_\xi(x)}
\end{equation}
where $\xi \sim \cD$, $f_\xi(x)$ is differentiable and $L$-smooth almost surely in $\xi$.

 Lemma~\ref{lem:sk_lemma1_sgd} shows that the stochastic gradient $g^k = \nabla f_\xi(x^k)$ satisfies Assumption~\ref{as:sk_general_stoch_gradient}. The corresponding choice of parameters can be found in Table~\ref{tbl:sk_special_cases-parameters}. 

\begin{lemma}[Generalization of Lemmas 1,2 from~\cite{nguyen2018sgd}]\label{lem:sk_lemma1_sgd}
    Assume that $f_\xi(x)$ is convex in $x$ for every $\xi$. Then for every $x\in\R^d$
    \begin{equation}\label{eq:sk_lemma1_sgd}
        \ED{\norm{\nabla f_\xi(x)- \nabla f(x^*)}^2} \le 4L(D_f(x,x^*)) + 2\sigma^2,
    \end{equation}
    where $\sigma^2\eqdef \EEE_\xi\left[\norm{\nabla f_\xi(x^*)}^2\right]$. If further $f(x)$ is $\mu$-strongly convex with possibly non-convex $f_\xi$, then for every $x\in\R^d$
    \begin{equation}\label{eq:sk_lemma2_sgd}
        \ED{\norm{\nabla f_\xi(x) - \nabla f(x^*)}^2} \le 4L\kappa(D_f(x,x^*)) + 2\sigma^2,
    \end{equation}
    where $\kappa = \frac{L}{\mu}$.
\end{lemma}

\begin{corollary}\label{cor:sk_recover_sgd_rate}
    Assume that $f_\xi(x)$ is convex in $x$ for every $\xi$ and $f$ is $\mu$-strongly quasi-convex. Then {\tt SGD} with $\alpha \le \frac{1}{2L}$ satisfies
    \begin{equation}\label{eq:sk_recover_sgd_rate}
        \EEE\left[\norm{x^k - x^*}^2\right] \le (1-\alpha\mu)^k\norm{x^0-x^*}^2 + \frac{2\alpha \sigma^2}{\mu}.
    \end{equation}
    If we further assume that $f(x)$ is $\mu$-strongly convex with possibly non-convex $f_\xi(x)$,   {\tt SGD} with $\alpha \le \frac{1}{2L\kappa}$ satisfies~\eqref{eq:sk_recover_sgd_rate} as well. 
\end{corollary}
\begin{proof}
    It suffices to plug parameters from Table~\ref{tbl:sk_special_cases-parameters} into Theorem~\ref{thm:sk_main_gsgm}. 
\end{proof}

\subsubsection*{Proof of Lemma~\ref{lem:sk_lemma1_sgd}}
The proof is a direct generalization to the one from~\cite{nguyen2018sgd}.
Note that
\begin{eqnarray*}
&& \frac12 \ED{\norm{ \nabla f_\xi (x) - \nabla f(x^*)}^2} - \ED{\norm{\nabla f_\xi (x^*) - \nabla f(x^*)}^2} \\
&& \qquad \qquad=
\frac12 \ED{\norm{ \nabla f_\xi (x) - \nabla f(x^*)}^2 - \norm{\nabla f_\xi (x^*) - \nabla f(x^*)}^2} \\
&& \qquad \qquad\overset{\eqref{eq:sk_1/2a_minus_b}}{\le}
 \ED{\norm{ \nabla f_\xi (x) - \nabla f_{\xi}(x^*)}^2}  \\
&& \qquad \qquad\leq
2LD_f(x,x^*).
\end{eqnarray*}
It remains to rearrange the above to get~\eqref{eq:sk_lemma1_sgd}. To obtain~\eqref{eq:sk_lemma2_sgd}, we shall proceed similarly:
\begin{eqnarray*}
&& \frac12 \ED{\norm{ \nabla f_\xi (x) - \nabla f(x^*)}^2} - \ED{\norm{ \nabla f_\xi (x^*) - \nabla f(x^*)}^2} \\
&& \qquad \qquad=
\frac12 \ED{\norm{ \nabla f_\xi (x) - \nabla f(x^*)}^2 - \norm{ \nabla f_\xi (x^*) - \nabla f(x^*)}^2} \\
&& \qquad \qquad\overset{\eqref{eq:sk_1/2a_minus_b}}{\le}
 \ED{\norm{\nabla f_\xi (x) - \nabla f_{\xi}(x^*)}^2}  \\
&& \qquad \qquad\leq
L^2\norm{ x - x^*}^2 \\ 
&& \qquad \qquad\leq
2\frac{L^2}{\mu}D_f(x,x^*).
\end{eqnarray*}
Again, it remains to rearrange the terms.

\subsection{{\tt SGD-SR}}\label{SGD-AS}

In this section, we recover convergence result of {\tt SGD} under expected smoothness property from~\cite{pmlr-v97-qian19b}. This setup allows obtaining tight convergence rates of {\tt SGD} under arbitrary stochastic reformulation of finite sum minimization\footnote{For technical details on how to exploit expected smoothness for specific reformulations, see~\cite{pmlr-v97-qian19b}}. 

The stochastic reformulation is a special instance of~\eqref{eq:sk_main_sgd}:
\begin{equation}\label{eq:sk_problem_sgd-as}
    \min\limits_{x\in\R^d}f(x) + \psi(x),\quad f(x) =\ED{f_\xi(x)}, \quad  f_\xi(x) \eqdef \frac1n \sum_{i=1}^n \xi_i f_i(x) 
\end{equation}
where $\xi$ is a random vector from distribution $\cD$ such that for all $i$: $\ED{\xi_i}=1$ and $f_i$ (for all $i$) is smooth, possibly non-convex function. We next state the expextes smoothness assumption. A specific instances of this assumption allows to get tight convergence rates of {\tt SGD}, which we recover in this section.

\begin{algorithm}[h]
    \caption{{\tt SGD-SR}}
    \label{alg:sk_sgdas}
    \begin{algorithmic}
        \Require learning rate $\alpha>0$, starting point $x^0\in\R^d$, distribution $\cD$ over $\xi \in\R^n$ such that $\ED{\xi}$ is vector of ones
        \For{ $k=0,1,2,\ldots$ }
        \State{Sample $\xi \sim \cD$}
        \State{$g^k = \nabla f_\xi (x^k)$}
        \State{$x^{k+1} = \proxR(x^k - \alpha g^k)$}
        \EndFor
    \end{algorithmic}
\end{algorithm}

\begin{assumption}[Expected smoothness]\label{as:sk_exp_smoothness_sgd-as}
    We say that $f$ is $\cL$-smooth in expectation with respect to distribution $\cD$ if there exists $\cL = \cL(f,\cD) > 0$ such that
    \begin{equation}\label{eq:sk_exp_smoothness_sgd-as}
        \ED{\norm{\nabla f_\xi(x) - \nabla f_\xi(x^*)}^2} \le 2\cL D_f(x,x^*),
    \end{equation}
    for all $x\in\R^d$. For simplicity, we will write $(f,\cD) \sim ES(\cL)$ to say that \eqref{eq:sk_exp_smoothness_sgd-as} holds.
\end{assumption}

Next, we present Lemma~\ref{lem:sk_exp_smoothness_grad_up_bound_sgd-as} which shows that choice of constants for Assumption~\ref{as:sk_general_stoch_gradient} from Table~\ref{tbl:sk_special_cases-parameters} is valid. 

\begin{lemma}[Generalization of Lemma~2.4, \cite{pmlr-v97-qian19b}]\label{lem:sk_exp_smoothness_grad_up_bound_sgd-as}
    If $(f,\cD)\sim ES(\cL)$, then
    \begin{equation}\label{eq:sk_exp_smoothness_grad_up_bound_sgd-as}
        \ED{\norm{\nabla f_\xi(x) - \nabla f(x^*)}^2} \le 4\cL D_f(x,x^*) + 2\sigma^2.
    \end{equation}
    where $\sigma^2 \eqdef \ED{\norm{\nabla f_\xi(x^*) - \nabla f(x^*	)}^2}$.
\end{lemma}

A direct consequence of Theorem~\ref{thm:sk_main_gsgm} in this setup is Corollary~\ref{cor:sk_recover_sgd-as_rate}.

\begin{corollary}\label{cor:sk_recover_sgd-as_rate}
    Assume that $f(x)$ is $\mu$-strongly quasi-convex and $(f,\cD)\sim ES(\cL)$. Then {\tt SGD-SR} with $\alpha^k \equiv\alpha \le \frac{1}{2\cL}$ satisfies
    \begin{equation}\label{eq:sk_recover_sgd-as_rate}
        \EEE\left[\norm{x^k - x^*}^2\right] \le (1-\alpha\mu)^k\norm{x^0-x^*}^2 + \frac{2\alpha\sigma^2}{\mu}.
    \end{equation}
\end{corollary}


\subsubsection*{Proof of Lemma~\ref{lem:sk_exp_smoothness_grad_up_bound_sgd-as}}
Here we present the generalization of the proof of Lemma~2.4 from \cite{pmlr-v97-qian19b} for the case when $\nabla f(x^*) \neq 0$. In this proof all expectations are conditioned on $x^k$.
\begin{eqnarray*}
    \EEE\left[\norm{\nabla f_{\xi} (x) - \nabla f(x^*)}^2\right] &=& \EEE\left[ \norm{ \nabla f_\xi(x) - \nabla f_{\xi}(x^*) + \nabla f_{\xi}(x^*) - \nabla f(x^*) }^2\right] \\
    &\overset{\eqref{eq:sk_a_b_norm_squared}}{\le}&  2 \EEE\left[\norm{ \nabla f_{\xi}(x) - \nabla f_{\xi}(x^*)}^2\right] 
 + 2  \EEE\left[\norm{\nabla f_{\xi}(x^*) - \nabla f(x^*)}^2\right] \notag\\ 
    &\overset{\eqref{eq:sk_exp_smoothness_sgd-as}}{\le}& 4 \cL D_f(x,x^*) + 2 \sigma^2.
\end{eqnarray*}

\subsection{{\tt SGD-MB}}\label{sec:sk_SGD-MB}
In this section, we present a specific practical formulation of~\eqref{eq:sk_problem_sgd-as} which was not considered in~\cite{pmlr-v97-qian19b}. The resulting algorithm (Algorithm~\ref{alg:sk_SGD-MB}) is novel; it was not considered in~\cite{pmlr-v97-qian19b} as a specific instance of {\tt SGD-SR}. The key idea behind {\tt SGD-MB} is constructing unbiased gradient estimate via with-replacement sampling.

Consider random variable $\nu \sim \cD$ such that 
\begin{equation}\label{eq:sk_nb87fg87f}
 \Prob(\nu = i) =p_i;  \qquad \sum_{i=1}^np_i=1.
\end{equation} Notice that if we define \begin{equation}
\label{eq:sk_reform_function}f'_i(x)\eqdef \frac{1}{n p_i} f_i(x), \qquad i=1,2,\dots,n,\end{equation}
then 
\begin{equation} \label{eq:sk_reform_problem}f(x) = \frac{1}{n} \sum_{i=1}^n f_i(x)  \overset{\eqref{eq:sk_reform_function}}{=} \sum_{i=1}^n p_i f'_i(x) \overset{\eqref{eq:sk_nb87fg87f}}{=} \ED{f'_\nu (x)}.\end{equation}
So, we have rewritten the finite sum problem \eqref{eq:sk_f_sum} into the {\em equivalent stochastic optimization problem}  
\begin{equation}\label{eq:sk_reform_opt} \min_{x\in \R^d} \ED{f'_\nu (x)}.\end{equation}

We are now ready to describe our method. At each iteration $k$ we sample $\nu^k_i,\dots,\nu^k_{\tau}\sim \cD$ independently ($1\leq \tau \leq n$), and define $g^k\eqdef  \frac{1}{\tau}\sum_{i=1}^\tau \nabla f'_{\nu^k_i}(x^k) $. Further, we use $g^k$ as a stochastic gradient, resulting in Algorithm~\ref{alg:sk_SGD-MB}.

\begin{algorithm}[h]
    \caption{{\tt SGD-MB}}
    \label{alg:sk_SGD-MB}
    \begin{algorithmic}
        \Require learning rate $\alpha>0$, starting point $x^0\in\R^d$, distribution $\cD$ over $\nu$ such that~\eqref{eq:sk_nb87fg87f} holds.
        \For{ $k=0,1,2,\ldots$ }
        \State{Sample $\nu^k_i,\dots,\nu^k_{\tau}\sim \cD$ independently}
        \State{$g^k =  \frac{1}{\tau}\sum_{i=1}^\tau \nabla f'_{\nu^k_i}(x^k) $}
        \State{$x^{k+1} = x^k - \alpha g^k$}
        \EndFor
    \end{algorithmic}
\end{algorithm}

To remain in full generality, consider the following Assumption.

\begin{assumption}\label{ass:sk_main_mb} There exists constants $A'>0$ and $D'\geq 0$ such that \begin{equation}\label{eq:sk_ABassumpt} \ED{ \norm{\nabla f'_{\nu}(x)}^2 }\leq 2 A' (f(x) - f(x^*)) + D'\end{equation} for all $x\in \R^d$.  
\end{assumption}

Note that it is sufficient to have convex and smooth $f_i$ in order to satisfy Assumption~\ref{ass:sk_main_mb}, as Lemma~\ref{lem:sk_b98h0f} states. 

\begin{lemma}\label{lem:sk_b98h0f}  Let $\sigma^2\eqdef \ED{\norm{\nabla f'_\nu (x^*)}^2}$. If $f_i$ are convex and $L_i$-smooth, then Assumption~\ref{ass:sk_main_mb} holds for $A'=2\cL$ and $D'=2\sigma^2$, where \begin{equation}\label{eq:sk_ineq90f9f}\cL \leq \max_i \frac{L_i}{n p_i}.\end{equation} 
If moreover $\nabla f_i(x^*)=0$ for all $i$,  then Assumption~\ref{ass:sk_main_mb} holds for $A'=\cL$ and $D'=0$. 
\end{lemma}

Next, Lemma~\ref{lem:sk_mb} states that Algorithm~\ref{alg:sk_SGD-MB} indeed satisfies Assumption~\ref{as:sk_general_stoch_gradient}. 

\begin{lemma} \label{lem:sk_mb}
Suppose that Assumption~\ref{ass:sk_main_mb} holds. Then $g^k$ is unbiased; i.e. $\ED{g^k} = \nabla f(x^k)$. Further, 
\begin{eqnarray*}
\ED{\norm{g^k}^2} \leq  \frac{2A' + 2L(\tau-1)}{\tau}(f(x^k) - f(x^*)) +  \frac{  D'}{\tau}.
\end{eqnarray*}
\end{lemma}
Thus, parameters from Table~\ref{tbl:sk_special_cases-parameters} are validated. As a direct consequence of Theorem~\ref{thm:sk_main_gsgm} we get Corollary~\ref{cor:sk_mb}.

\begin{corollary}\label{cor:sk_mb}
    As long as $0< \alpha \leq \frac{\tau}{A' + L(\tau-1)}$, we have
    \begin{equation}\label{eq:sk_mb_rate}
\EEE \norm{x^k-x^*}^2 \leq (1-\alpha \mu)^k \norm{x^0-x^*}^2 + \frac{\alpha D'}{\mu \tau}.    \end{equation}
\end{corollary}

\begin{remark} \label{rem:sgdmb}
For $\tau=1$, {\tt SGD-MB} is a special of the method from~\cite{pmlr-v97-qian19b}, Section~3.2. However, for $\tau>1$, this is a different method; the difference lies in the with-replacement sampling. Note that with-replacement trick allows for efficient and implementation of independent importance sampling~\footnote{Distribution of random sets $S$ for which random variables $i\in S$ and $j\in S$ are independent for $j\neq i$.} with complexity $\cO(\tau \log(n))$. In contrast, implementation of without-replacement importance sampling has complexity $\cO(n)$, which can be significantly more expensive to the cost of evaluating $\sum_{i\in S} \nabla f_i(x)$.
 \end{remark}

\subsubsection*{Proof of Lemma~\ref{lem:sk_mb}}
Notice first that
\begin{eqnarray*}
\ED{g^k}
&\overset{\eqref{eq:sk_reform_function}}{=}&
  \frac{1}{\tau}\sum_{i=1}^\tau\ED{\frac{1}{n p_{\nu^k_{i}}} \nabla f_{\nu^k_i}(x^k)}
  \\
&=& 
\ED{\frac{1}{n p_{\nu}} \nabla f_{\nu}(x^k)}
  \\
&\overset{\eqref{eq:sk_nb87fg87f}}{=}& 
\sum_{i=1}^n p_i \frac{1}{n p_i} \nabla f_i(x^k) 
\\
&=& \nabla f(x_k).
\end{eqnarray*}

So, $g^k$ is an unbiased estimator of the gradient $\nabla f(x^k)$. Next, 
\begin{eqnarray*}
\ED{\norm{g^k}^2}&=& \ED{\norm{\frac{1 }{\tau}\sum_{i=1}^\tau \nabla f'_{\nu^k_i}(x^k)}^2 }\\ 
&=& \frac{1}{\tau^2}\ED{\sum_{i=1}^\tau \norm{\nabla f'_{\nu^k_i}(x^k)}^2  +  2    \sum_{i< j} \< \nabla f'_{\nu^k_i}(x^k),  \nabla f'_{\nu^k_j}(x^k)>  } \\
&=& \frac{1}{\tau}\ED{ \norm{\nabla f'_{\nu}(x^k)}^2}+  \frac{2}{\tau^2}    \sum_{i< j} \< \ED{\nabla f'_{\nu^k_i}(x^k) },  \ED{\nabla f'_{\nu^k_j}(x^k) }>     \\
&=& \frac{1}{\tau} \ED{\norm{\nabla f'_{\nu}(x^k)}^2} + \frac{\tau - 1}{\tau} \norm{\nabla f(x^k)}^2 \\
&\overset{\eqref{eq:sk_ABassumpt} }{\leq}& \frac{2A' (f(x^k) - f(x^*)) + D' + 2L(\tau-1)(f(x^k) - f(x^*))}{\tau}.
\end{eqnarray*}

\subsubsection*{Proof of Lemma~\ref{lem:sk_b98h0f}}
Let $\cL = \cL(f,\cD)>0$ be any constant for which
\begin{equation}\label{eq:sk_ES}\EEE_{\xi\sim \cD} \norm{\nabla \phi_{\xi}(x) - \nabla \phi_{\xi}(x^*)}^2 \leq 2\cL (f(x) - f(x^*))\end{equation}
 holds for all $x\in \R^d$. This is the expected smoothness property (for a single item sampling) from \cite{pmlr-v97-qian19b}. It was shown in \cite[Proposition 3.7]{pmlr-v97-qian19b}  that \eqref{eq:sk_ES} holds, and that $ \cL$ satisfies \eqref{eq:sk_ineq90f9f}. The claim now follows by applying \cite[Lemma~2.4]{pmlr-v97-qian19b}. 

\subsection{{\tt SGD-star}}\label{sec:sk_SGD-star}
Consider problem~\eqref{eq:sk_problem_sgd-as}. Suppose that $\nabla f_i(x^*)$ is known for all $i$. In this section we present a novel algorithm~---~{\tt SGD-star}~---~which is {\tt SGD-SR} shifted by the stochastic gradient in the optimum. The method is presented under Expected Smoothness Assumption~\eqref{eq:sk_exp_smoothness_sgd-as}, obtaining general rates under arbitrary sampling.  The algorithm is presented as Algorithm~\ref{alg:sk_SGD-star}. 

\begin{algorithm}[h]
    \caption{{\tt SGD-star}}
    \label{alg:sk_SGD-star}
    \begin{algorithmic}
        \Require learning rate $\alpha>0$, starting point $x^0\in\R^d$, distribution $\cD$ over $\xi \in \R^n$ such that $\ED{\xi}$ is vector of ones
        \For{ $k=0,1,2,\ldots$ }
        \State{Sample $\xi \sim \cD$}
        \State{$g^k = \nabla f_\xi(x^k) - \nabla f_\xi(x^*) + \nabla f(x^*)$}
        \State{$x^{k+1} = \proxR(x^k - \alpha g^k)$}
        \EndFor
    \end{algorithmic}
\end{algorithm}
Suppose that $(f,\cD) \sim ES(\cL)$. 
Note next that {\tt SGD-star} is just  {\tt SGD-SR} applied on objective $D_f(x,x^*)$ instead of $f(x)$ when $\nabla f(x^*) = 0$. This careful design of the objective yields $(D_f(\cdot, x^*),\cD) \sim ES(\cL)$ and $\ED{\norm{\nabla_x D_{f_\xi}(x,x^*)}^2 \, \mid x=x^*} = 0 $, and thus Lemma~\eqref{lem:sk_exp_smoothness_grad_up_bound_sgd-as} becomes
\begin{lemma}[Lemma~2.4, \cite{pmlr-v97-qian19b}]\label{lem:sk_SGD-star}
    If $(f,\cD)\sim ES(\cL)$, then
    \begin{equation}\label{eq:sk_sgd-star_lemma}
        \ED{\norm{g^k -\nabla f(x^*)}^2} \le 4\cL D_f(x^k,x^*).
    \end{equation}
\end{lemma}

A direct consequence of Corollary (thus also a direct consequence of Theorem~\ref{thm:sk_main_gsgm}) in this setup is Corollary~\ref{cor:sk_SGD-star}.

\begin{corollary}\label{cor:sk_SGD-star}
    Suppose that $(f,\cD)\sim ES(\cL)$. Then {\tt SGD-star} with $\alpha = \frac{1}{2\cL}$ satisfies
    \begin{equation}\label{eq:sk_SGD-shift-xx}
        \EEE\left[\norm{x^k - x^*}^2\right] \le \left(1-\frac{\mu}{2\cL} \right)^k\norm{x^0-x^*}^2.
    \end{equation}
\end{corollary}


\begin{remark}
Note that results from this section are obtained by applying results from~\ref{SGD-AS}. Since Section~\ref{sec:sk_SGD-MB} presets a specific sampling algorithm for {\tt SGD-SR}, the results can be thus extended to {\tt SGD-star} as well.
\end{remark}

\subsubsection*{Proof of Lemma~\ref{lem:sk_SGD-star}}
In this proof all expectations are conditioned on $x^k$.
\begin{eqnarray*}
    \ED{\norm{g^k -\nabla f(x^*)}^2}  &=& \ED{\norm{\nabla f_\xi(x^k) -\nabla f_\xi(x^*)}^2} \\ 
    &\overset{\eqref{eq:sk_exp_smoothness_sgd-as}}{\le}& 4 \cL D_f(x^k,x^*).
\end{eqnarray*}

\subsection{{\tt SAGA}}\label{sec:sk_saga}
In this section we show that our approach is suitable for {\tt SAGA} \cite{saga} (see Algorithm~\ref{alg:sk_SAGA}). Consider the finite-sum minimization problem
\begin{equation}\label{eq:sk_main_l-svrg}
    f(x) = \frac{1}{n}\sum\limits_{i=1}^n f_i(x) + \psi(x),
\end{equation}
where $f_i$ is convex, $L$-smooth for each $i$ and $f$ is $\mu$-strongly convex. 

\begin{algorithm}[h]
    \caption{{\tt SAGA} \cite{saga}}
    \label{alg:sk_SAGA}
    \begin{algorithmic}
        \Require learning rate $\alpha>0$, starting point $x^0\in\R^d$
        \State Set $\phi_j^0 = x^0$ for each $j\in[n]$
        \For{ $k=0,1,2,\ldots$ }
        \State{Sample  $j \in [n]$ uniformly at random}
        \State{Set $\phi_j^{k+1} = x^k$ and $\phi_i^{k+1} = \phi_i^{k}$ for $i\neq j$}
        \State{$g^k = \nabla f_j(\phi_j^{k+1}) - \nabla f_j(\phi_j^k) + \frac{1}{n}\sum\limits_{i=1}^n\nabla f_i(\phi_i^k)$}
        \State{$x^{k+1} = \proxR\left(x^k - \alpha g^k\right)$}
        \EndFor
    \end{algorithmic}
\end{algorithm}

\begin{lemma}\label{lem:sk_stoch_grad_second_moment_saga}
    We have
    \begin{equation}\label{eq:sk_stoch_grad_second_moment_saga1}
        \EEE\left[\norm{g^k- \nabla f(x^*)}^2\mid x^k\right] \le 4LD_f(x^k,x^*)+ 2\sigma_k^2
    \end{equation}
    and 
        \begin{equation}\label{eq:sk_stoch_grad_second_moment_saga2}
        \EEE\left[\sigma_{k+1}^2\mid x^k\right] \le \left(1 - \frac{1}{n}\right)\sigma_k^2 + \frac{2L}{n}D_f(x^k,x^*),
    \end{equation}
    where $\sigma_k^2 = \frac{1}{n}\sum\limits_{i=1}^n\norm{\nabla f_i(\phi_i^k) - \nabla f_i(x^*)}^2$.
\end{lemma}

Clearly, Lemma~\ref{lem:sk_stoch_grad_second_moment_saga} shows that Algorithm~\ref{alg:sk_SAGA} satisfies Assumption~\ref{as:sk_general_stoch_gradient}; the corresponding parameter choice can be found in Table~\ref{tbl:sk_special_cases-parameters}. Thus, as a direct consequence of Theorem~\ref{thm:sk_main_gsgm} with $M=4n$ we obtain the next corollary.

\begin{corollary}\label{thm:sk_recover_saga_rate}
    {\tt SAGA} with $\alpha = \frac{1}{6L}$ satisfies
    \begin{equation}\label{eq:sk_recover_saga_rate_coro}
        \EEE V^k \le \left(1-\min\left\{\frac{\mu}{6L},\frac{1}{2n}\right\}\right)^kV^0.
    \end{equation}
\end{corollary}

\subsubsection*{Proof of Lemma~\ref{lem:sk_stoch_grad_second_moment_saga}}
Note that  Lemma~\ref{lem:sk_stoch_grad_second_moment_saga} is a special case of Lemmas 3,4 from~\cite{mishchenko201999} without prox term. We reprove it with prox for completeness.

    Let all expectations be conditioned on $x^k$ in this proof. Note that $L$-smoothness and convexity of $f_i$ implies
\begin{equation}\label{eq:sk_norm_diff_grads}
    \frac{1}{2L}\norm{\nabla f_i(x) - \nabla f_i(y)}^2 \le f_i(x) - f_i(y) - \<\nabla f_i(y),x-y>,\quad \forall x,y\in\R^d, i\in[n].
\end{equation}

 By definition of $g^k$ we have
    \begin{eqnarray*}
       && \EEE\left[\norm{g^k - \nabla f(x^*)}^2\right] \\
        &  =&
        \EEE\left[ \norm{\nabla f_j(\phi_j^{k+1}) - \nabla f_j(\phi_j^k) + \frac{1}{n}\sum\limits_{i=1}^n\nabla f_i(\phi_i^k) - \nabla f(x^*)}^2\right]
        \\
        &  =&
         \EEE\left[\norm{\nabla f_j(x^k) - \nabla f_j(x^*) + \nabla f_j(x^*) - \nabla f_j(\phi_j^k) + \frac{1}{n}\sum\limits_{i=1}^n\nabla f_i(\phi_i^k) - \nabla f(x^*) }^2 \right]
         \\
                &  \overset{\eqref{eq:sk_a_b_norm_squared}}{\le} &
        2\EEE\left[\norm{\nabla f_j(x^k) - \nabla f_j(x^*)}^2\mid x^k\right]
        \\
        && \qquad \quad + 2\EEE\left[\norm{\nabla f_j(x^*) - \nabla f_j(\phi_j^k) - \EEE\left[\nabla f_j(x^*) - \nabla f_j(\phi_j^k)\right] }^2\right]
        \\
              & \overset{\eqref{eq:sk_variance_decomposition}+\eqref{eq:sk_norm_diff_grads}}{\le}&
         \frac{4L}{n}\sum\limits_{i=1}^nD_{f_i}(x^k,x^*)+ 2\EEE\left[\norm{\nabla f_j(x^*) - \nabla f_j(\phi_j^k)}^2\mid x^k\right]\\
              &  = & 4LD_f(x^k,x^*) + 2\underbrace{\frac{1}{n}\sum\limits_{i=1}^n\norm{\nabla f_i(\phi_i^k) - \nabla f_i(x^*)}^2}_{\sigma_k^2}.
    \end{eqnarray*}
    
To proceed with~\eqref{eq:sk_stoch_grad_second_moment_saga2}, we have
    \begin{eqnarray*}
        \EEE\left[\sigma_{k+1}^2\right] &=& \frac{1}{n}\sum\limits_{i=1}^n\EEE\left[\norm{\nabla f_i(\phi_i^{k+1}) - \nabla f_i(x^*)}^2\right]\\
        &=& \frac{1}{n}\sum\limits_{i=1}^n\left(\frac{n-1}{n}\norm{\nabla f_i(\phi_i^k) - \nabla f_i(x^*)}^2 + \frac{1}{n}\norm{\nabla f_i(x^k) - \nabla f_i(x^*)}^2\right)\\
        &\overset{\eqref{eq:sk_norm_diff_grads}}{\le}& \left(1 - \frac{1}{n}\right)\frac{1}{n}\sum\limits_{i=1}^n\norm{\nabla f_i(\phi_i^k) - \nabla f_i(x^*)}^2\\
        &&\quad + \frac{2L}{n^2}\sum\limits_{i=1}^n D_{f_i}(x^k,x^*)\\
        &=& \left(1 - \frac{1}{n}\right)\sigma_k^2 + \frac{2L}{n}D_f(x^k,x^*).
    \end{eqnarray*}

\subsection{{\tt N-SAGA}} \label{N-SAGA}

\begin{algorithm}[h]
    \caption{Noisy {\tt SAGA} ({\tt N-SAGA})}
    \label{alg:sk_N-SAGA}
    \begin{algorithmic}
        \Require learning rate $\alpha>0$, starting point $x^0\in\R^d$
        \State Set $\psi_j^0 = x^0$ for each $j\in[0]$
        \For{ $k=0,1,2,\ldots$ }
        \State{Sample  $j \in [n]$ uniformly at random and $\zeta$}
        \State{Set $g_j^{k+1} = g_j(x^k,\xi)$ and $g_i^{k+1} = g_i^{k} $ for $i\neq j$}
        \State{$g^k =g_j(x^k,\xi) - g_j^k + \frac{1}{n}\sum\limits_{i=1}^ng_i^k$}
        \State{$x^{k+1} = \proxR(x^k - \alpha g^k)$}
        \EndFor
    \end{algorithmic}
\end{algorithm}

Note that it can in practice happen that instead of $\nabla f_i(x)$ one can query $g_i(x,\zeta)$ such that $\EEE_\xi g_i(\cdot,\xi)=  \nabla f_i(\cdot) $ and $\EEE_\xi \norm{g_i(\cdot,\xi)}^2 \leq  \sigma^2$. This leads to a variant of {\tt SAGA} which only uses noisy  estimates of the stochastic gradients $\nabla_i(\cdot)$. We call this variant {\tt N-SAGA} (see Algorithm~\ref{alg:sk_N-SAGA}).

\begin{lemma}\label{lem:sk_saga_inexact1}
    We have
    \begin{equation}\label{eq:sk_stoch_grad_second_moment_saga_noise1}
        \EEE\left[\norm{g^k - \nabla f(x^*)}^2\mid x^k\right] \le 4LD_f(x^k,x^*) + 2\sigma_k^2 + 2\sigma^2,
    \end{equation}
    and
        \begin{equation}\label{eq:sk_stoch_grad_second_moment_saga_noise1}
        \EEE\left[\sigma_{k+1}^2\mid x^k\right] \le \left(1 - \frac{1}{n}\right)\sigma_k^2 + \frac{2L}{n}D_f(x^k,x^*) + \frac{\sigma^2}{n},
    \end{equation}
    where $\sigma_k^2 \eqdef \frac{1}{n}\sum\limits_{i=1}^n\norm{ g_i^k - \nabla f_i(x^*)}^2$.
\end{lemma}

\begin{corollary}\label{cor:sk_N-SAGA}
Let $\alpha = \frac{1}{6L}$. Then, iterates of Algorithm~\ref{alg:sk_N-SAGA} satisfy
\[
\EEE{V^k}\leq \left( 1- \min\left( \frac{\mu}{6L}, \frac{1}{2n} \right)\right)^k V^0 + \frac{\sigma^2 }{ L \min(\mu, \frac{3 L }{n})}.
\]
\end{corollary}

Analogous results can be obtained for {\tt LSVRG}.

\subsubsection*{Proof of Lemma~\ref{lem:sk_saga_inexact1}}

    Let all expectations be conditioned on $x^k$. By definition of $g^k$ we have
    \begin{eqnarray*}
    && 
        \EEE\left[\norm{g^k - \nabla f(x^*)}^2\right]  \\
        &\le& 
        \EEE\left[\norm{g_j(x^k,\zeta)- g_j^k + \frac{1}{n}\sum\limits_{i=1}^n g_i^k - \nabla f(x^*)}^2\right]
        \\
        &=&
         \EEE\left[\norm{g_j(x^k,\zeta) - \nabla f_j(x^*) + \nabla f_j(x^*) - g_j^k + \frac{1}{n}\sum\limits_{i=1}^n g_i^k  - \nabla f(x^*)}^2\right]
         \\
        &\overset{\eqref{eq:sk_a_b_norm_squared}}{\le}& 
        2\EEE\left[\norm{g_j(x^k,\zeta) - \nabla f_j(x^*)}^2\right]
        \\
        &&\quad 
        + 2\EEE\left[\norm{\nabla f_j(x^*) -g_j^k - \EEE\left[\nabla f_j(x^*) - g_j^k\right] }^2\right]
        \\
        &\overset{\eqref{eq:sk_variance_decomposition}}{\le}&
         2\EEE\left[\norm{g_j(x^k,\zeta) - \nabla f_j(x^*)}^2\right] + 2\EEE\left[\norm{\nabla f_j(x^*) - g_j^k }^2\right]
         \\
        &=& 
        2\EEE\left[\norm{g_j(x^k,\zeta) - \nabla f_j(x^*)}^2\right] + 2\underbrace{\frac{1}{n}\sum\limits_{i=1}^n\norm{g_i^k - \nabla f_i(x^*)}^2}_{\sigma_k^2} 
        \\
        &\overset{\eqref{eq:sk_variance_decomposition}}{\le}& 
                 2\EEE\left[ \norm{ \nabla f_j(x^k) - \nabla f_j(x^*)}^2\right] + 2\sigma^2 + 2\sigma_k^2
                \\
                &\overset{\eqref{eq:sk_norm_diff_grads}}{\le}&
                4LD_f(x^k,x^*) + 2\sigma_k^2 + 2\sigma^2
    \end{eqnarray*}

    For the second inequality, we have
    \begin{eqnarray*}
        \EEE\left[\sigma_{k+1}^2\right] &=& \frac{1}{n}\sum\limits_{i=1}^n\EEE\left[\norm{g_i^{k+1} - \nabla f_i(x^*)}^2\right]
        \\
        &=& 
        \frac{1}{n}\sum\limits_{i=1}^n\left(\frac{n-1}{n}\norm{g_i^k  - \nabla f_i(x^*)}^2 + \frac{1}{n}\EEE\left[ \norm{ g_i(x^k, \zeta) - \nabla f_i(x^*)}^2\right]\right)
        \\
        &\leq&
                \frac{1}{n}\sum\limits_{i=1}^n\left(\frac{n-1}{n}\norm{g_i^k  - \nabla f_i(x^*)}^2 + \frac{1}{n} \norm{ \nabla f_i(x^k) - \nabla f_i(x^*)}^2 + \frac{\sigma^2}{n}\right)
        \\
        &\overset{\eqref{eq:sk_norm_diff_grads}}{\le}&
        \left(1 - \frac{1}{n}\right)\sigma_k^2 + \frac{2L}{n}D_f(x^k,x^*)+ \frac{\sigma^2}{n}.
    \end{eqnarray*}

\subsection{{\tt SEGA}}\label{sec:sk_sega}

\begin{algorithm}[h]
    \caption{{\tt SEGA} \cite{sega}}
    \label{alg:sk_SEGA}
    \begin{algorithmic}
        \Require learning rate $\alpha>0$, starting point $x^0\in\R^d$
        \State Set $h^0 = 0$
        \For{ $k=0,1,2,\ldots$ }
        \State{Sample  $j \in [d]$ uniformly at random}
        \State{Set $h^{k+1} = h^{k} +e_i( \nabla_i f(x^k) - h^{k}_i)$}
        \State{$g^k = de_i (\nabla_i f(x^k) - h_i^k) + h^k$}
        \State{$x^{k+1} = \proxR(x^k - \alpha g^k)$}
        \EndFor
    \end{algorithmic}
\end{algorithm}

We show that the framework recovers the simplest version of {\tt SEGA} (i.e., setup from Theorem D1 from~\cite{sega}) in the proximal setting\footnote{General version for arbitrary gradient sketches instead of partial derivatives can be recovered as well, however, we omit it for simplicity}. 

\begin{lemma} (Consequence of Lemmas A.3., A.4. from~\cite{sega})
We have
\[
\EEE\left[\norm{g^{k}-\nabla f(x^*) \mid x^k}^{2}\right] \leq 2d\norm{\nabla f\left(x^{k}\right)-\nabla f(x^*)}^{2}+2d\sigma_k^2
\]
and
\[
\EEE
\left[\sigma_{k+1}^2 \mid x^k \right]= \left(1-\frac1d\right)\sigma_k^2+\frac1d \norm{\nabla f\left(x^{k}\right)-\nabla f(x^*)}^{2},
\]
where $\sigma_k^2 \eqdef \norm{h^{k}-\nabla f(x^*)}^2$.
\end{lemma}

Given that we have from convexity and smoothness $\norm{\nabla f(x^{k})-\nabla f(x^*)}^{2} \leq 2L D_f(x^k,x^*)$, Assumption~\ref{as:sk_general_stoch_gradient} holds the parameter choice as per Table~\ref{tbl:sk_special_cases-parameters}. Setting further $M = 4d^2$, we get the next corollary.
\begin{corollary}\label{cor:sk_sega}
{\tt SEGA} with $\alpha =\frac{1}{6dL} $ satisfies
\[
\EEE V^k \leq \left( 1- \frac{\mu}{6dL}\right)^kV^0.
\] 
\end{corollary}

\subsection{{\tt N-SEGA}} \label{N-SEGA}

\begin{algorithm}[h]
    \caption{Noisy {\tt SEGA} ({\tt N-SEGA})}
    \label{alg:sk_N-SEGA}
    \begin{algorithmic}
        \Require learning rate $\alpha>0$, starting point $x^0\in\R^d$
        \State Set $h^0 = 0$
        \For{ $k=0,1,2,\ldots$ }
        \State{Sample  $i \in [d]$ uniformly at random and sample $\xi$}
        \State{Set $h^{k+1} = h^{k} +e_i( g_i(x,\xi)- h^{k}_i)$}
        \State{$g^k = de_i (g_i(x,\xi) - h_i^k) +  h^k$}
        \State{$x^{k+1} = x^k - \alpha g^k$}
        \EndFor
    \end{algorithmic}
\end{algorithm}

Here we assume that $g_i(x,\zeta)$ is a noisy estimate of the partial derivative $\nabla_i f(x)$ such that $\EEE_\zeta g_i(x,\zeta) = \nabla_i f(x)$ and $\EEE_\zeta | g_i(x,\zeta) - \nabla_i f(x)|^2 \leq \frac{\sigma^2}{d}$. 

\begin{lemma} \label{lem:sk_sega_noise}
The following inequalities hold:
\[
\EEE\left[\norm{g^{k}-\nabla f(x^*)}^{2}\right]
 \leq 
 4dLD_f(x^k,x^*)+2d\sigma_k^2 +2d\sigma^2,
\]
\[
\EEE
\left[\sigma_{k+1}^2\right] \leq \left(1-\frac1d\right)\sigma_k^2+\frac{2L}{d}D_f(x^k,x^*) + \frac{\sigma^2}{d}, 
\]
where $\sigma_k^2 = \norm{h^{k}-\nabla f(x^*)}^2$.
\end{lemma}

\begin{corollary}\label{cor:sk_N-SEGA}
Let $\alpha = \frac{1}{6Ld}$. Applying Theorem~\ref{thm:sk_main_gsgm} with $M = 4d^2$, iterates of Algorithm~\ref{alg:sk_N-SEGA} satisfy
\[
\EEE{V^k}\leq \left( 1- \frac{\mu}{6dL}\right)^k V^0 + \frac{\sigma^2 }{ L \mu}.
\]
\end{corollary}

\subsubsection*{Proof of Lemma~\ref{lem:sk_sega_noise}}
 
Let all expectations be conditioned on $x^k$. 
For the first bound, we write \[g^k  - \nabla f(x^*)= \underbrace{h^k - \nabla f(x^*)- d  h^k_i e_{i}+d\nabla_i f(x^*)e_i}_{a} +   \underbrace{ d g_i(x^k,\xi)e_i- d\nabla_i f(x^*)e_i}_{b}.\]
 Let us bound the expectation of each term individually. The first term can be bounded as
\begin{eqnarray*}
\EEE{\norm{a}^2} &=& \EEE{\norm{ \left(\mI - d e_{i} e_{i}^\top \right) (h^k - \nabla f(x^*))  }_{2}^2}\\
&=& (d-1)\norm{h^k- \nabla f(x^*)}^2\\
&\leq& d\norm{h^k- \nabla f(x^*)}^2.
\end{eqnarray*}
The second term can be bounded as
\begin{eqnarray*}
\EEE{\norm{b}^2} &=&  \EEE_i \EEE_{\xi}{\norm{ d g_i(x,\xi)e_i- d\nabla f_i(x^*)e_i}^2}\\
&=& 
 \EEE_{i} \EEE_{\xi} \norm{ d g_i(x^k,\xi)e_i- d\nabla_i f(x^k)e_i }^2 + \EEE_i\norm{ d\nabla_i f(x^k)e_i-d\nabla f_i(x^*)e_i}^2 
\\
&\leq & d\sigma^2 + d  \norm{ \nabla f(x^k)- \nabla f(x^*)}^2 \\
&\leq& d\sigma^2 + 2Ld D_f(x^k,x^*),
\end{eqnarray*}
where in the last step we used $L$-smoothness of $f$. It remains to combine the two bounds.

For the second bound, we have
\begin{eqnarray*}
\EEE{\norm{ h^{k+1}  - \nabla f(x^*)}^2} &=& \EEE{ \norm{h^k + g_i(x^k,\xi)e_i - h^k_i - \nabla f(x^*) }^2 }\\
 &=&  \EEE{ \norm{\left(\mI - e_{i} e_{i}^\top \right)h^k + g_i(x^k,\xi)e_i -\nabla f(x^*)  }^2 }\\
 &=& \EEE{\norm{ \left(\mI - e_{i} e_{i}^\top \right) ( h^k - \nabla f(x^*)) }^2} + \EEE{\norm{g_i(x^k,\xi)e_i - \nabla_i f(x^*)e_i }^2 }\\
  &=& \left(1-\frac{1}{d}\right) \norm{h^k - \nabla f(x^*)}^2 + \EEE{\norm{g_i(x^k,\xi)e_i - \nabla_i f(x^k)e_i }^2 }   \\
  && \qquad + \EEE{\norm{\nabla_i f(x^k)e_i  - \nabla_i f(x^*)e_i}^2 } \\
 &=& \left(1-\frac{1}{d}\right) \norm{h^k - \nabla f(x^*)}^2 + \frac{\sigma^2}{d} + \frac{1}{d} \norm{\nabla f(x^k) - \nabla f(x^*)}^2 \\
  &\leq&
   \left(1-\frac{1}{d}\right) \norm{h^k - \nabla f(x^*)}^2 + \frac{\sigma^2}{d} + \frac{2L}{d}D_f(x^k,x^*).
\end{eqnarray*}

\subsection{{\tt SVRG}} \label{sec:sk_svrg}

\begin{algorithm}[h]
    \caption{{\tt SVRG} \cite{svrg}}
    \label{alg:sk_SVRG}
    \begin{algorithmic}
        \Require learning rate $\alpha>0$, epoch length $m$, starting point $x^0\in\R^d$
        \State  $\phi = x^0$
        \For{ $s=0,1,2,\ldots$ }
        \For{ $k=0,1,2,\ldots, m-1$ }
        \State{Sample  $i \in \{1,\ldots, n\}$ uniformly at random}
        \State{$g^k = \nabla f_i(x^k) - \nabla f_i(\phi) + \nabla f(\phi)$}
        \State{$x^{k+1} = \proxR(x^k - \alpha g^k)$}
        \EndFor
        \State  $\phi = x^0 = \frac1m \sum_{k=1}^m x^k$
        \EndFor
    \end{algorithmic}
\end{algorithm}
Let $\sigma_k^2 \eqdef \frac{1}{n}\sum\limits_{i=1}^n\norm{\nabla f_i(\phi) - \nabla f_i(x^*)}^2$. We will show that Lemma~\ref{lem:sk_iter_dec} recovers per-epoch analysis of {\tt SVRG} in a special case.

\begin{lemma}\label{lem:sk_svrg_lemma_1}
For $k \mod m \neq 0$ we have    
        \begin{equation}\label{eq:sk_svrg1}
        \EEE\left[\norm{g^k -\nabla f(x^*)}^2\mid x^k\right] \le 4LD_f(x^k,x^*) + 2\sigma_k^2
    \end{equation}
    and
    \begin{equation}\label{eq:sk_svrg3}
        \EEE\left[\sigma_{k+1}^2\mid x^k\right] = \sigma_{k+1}^2 = \sigma_k^2.
    \end{equation}
\end{lemma}
\begin{proof}
The proof of~\eqref{eq:sk_svrg1} is identical to the proof of~\eqref{eq:sk_stoch_grad_second_moment_saga1}. Next,~\eqref{eq:sk_svrg3} holds since $\sigma_k$ does not depend on $k$. 
\end{proof}

Thus, Assumption~\ref{as:sk_general_stoch_gradient} holds with parameter choice as per Table~\ref{tbl:sk_special_cases-parameters} and Lemma~\ref{lem:sk_iter_dec} implies the next corollary.
\begin{corollary}\label{cor:sk_svrg}
\begin{equation}\label{eq:sk_svrg_iter_dec}
 \EEE\norm{x^{k+1}-x^*}^2 + \alpha (1-2\alpha L)\EEE D_{f}(x^k,x^*) \leq 
        (1-\alpha\mu)\EEE\norm{x^k - x^*}^2 +2\alpha^2\EEE\sigma_k^2.
\end{equation}
\end{corollary}

\subsubsection*{Recovering {\tt SVRG} rate}
Summing~\eqref{eq:sk_svrg_iter_dec} for $k=0, \dots, m-1$ using $\sigma_k = \sigma_0$ we arrive at
\begin{eqnarray*}
\EEE\norm{x^{m}-x^*}^2+ \sum_{k=1}^m \alpha (1-2\alpha L)\EEE D_{f}(x^k,x^*)
&\leq& (1-\alpha\mu)\EEE\norm{x^0 - x^*}^2 + 2m\alpha^2\EEE\sigma_0^2 \\
&\leq &2 \left( \mu^{-1} +  2m\alpha^2   L\right) D_f(x^0,x^*) \;.
\end{eqnarray*}

Since $D_f$ is convex in the first argument, we have 
\[
m \alpha (1-2\alpha L)D_{f}\left( \frac1m \sum_{k=1}^m x^k,x^*\right) \leq \norm{x^{m}-x^*}^2+ \sum_{k=1}^m \alpha (1-2\alpha L)D_{f}(x^k,x^*)  
\]
and thus
\[
D_{f}\left( \frac1m \sum_{k=1}^m x^k,x^*\right) \leq \frac{ 2 \left( \mu^{-1} +  2m\alpha^2   L\right)}
{m \alpha (1-2\alpha L)} D_f(x^0,x^*) ,
 \]
which recovers rate from Theorem 1 in~\cite{svrg}.

\subsection{{\tt LSVRG}}\label{sec:sk_L-SVRG}

In this section we show that our approach also covers {\tt LSVRG} analysis from \cite{hofmann2015variance, kovalev2019don} (see Algorithm~\ref{alg:sk_L-SVRG}) with a minor extension -- it allows for proximable regularizer $\psi$. Consider the finite-sum minimization problem
\begin{equation}\label{eq:sk_main_l-svrg}
    f(x) = \frac{1}{n}\sum\limits_{i=1}^n f_i(x) + \psi(x),
\end{equation}
where each $f_i$ convex and  $L$-smooth for each $i$ and $f$ is $\mu$-strongly convex.

\begin{algorithm}[h]
    \caption{{\tt LSVRG} (\cite{hofmann2015variance, kovalev2019don})}
    \label{alg:sk_L-SVRG}
    \begin{algorithmic}
        \Require learning rate $\alpha>0$, probability $p\in (0,1]$, starting point $x^0\in\R^d$
        \State $w^0 = x^0$
        \For{ $k=0,1,2,\ldots$ }
        \State{Sample  $i \in \{1,\ldots, n\}$ uniformly at random}
        \State{$g^k = \nabla f_i(x^k) - \nabla f_i(w^k) + \nabla f(w^k)$}
        \State{$x^{k+1} = x^k - \alpha g^k$}
        \State{$w^{k+1} = \begin{cases}
            x^{k}& \text{with probability } p\\
            w^k& \text{with probability } 1-p
            \end{cases}$
        }
        \EndFor
    \end{algorithmic}
\end{algorithm}
Note that the gradient estimator is again unbiased, i.e. $\EEE\left[g^k\mid x^k\right] = \nabla f(x^k)$. Next, Lemma~\ref{lem:sk_l-svrg} provides with the remaining constants for Assumption~\ref{as:sk_general_stoch_gradient}. The corresponding choice is stated in Table~\ref{tbl:sk_special_cases-parameters}.

\begin{lemma}[Lemma~4.2 and Lemma~4.3 from~\cite{kovalev2019don} extended to prox setup]\label{lem:sk_l-svrg}
    We have
    \begin{equation}\label{eq:sk_lemma4.2_l-svrg}
        \EEE\left[\norm{g^k - \nabla f(x^*)}^2\mid x^k\right] \le 4LD_f(x^k,x^*) + 2\sigma_k^2
    \end{equation}
    and 
        \begin{equation}\label{eq:sk_lemma4.3_l-svrg}
        \EEE\left[\sigma_{k+1}^2\mid x^k\right] \le (1-p)\sigma_k^2 + 2LpD_f(x^k, x^*),
    \end{equation}
    where $\sigma_k^2 \eqdef \frac{1}{n}\sum\limits_{i=1}^n\norm{\nabla f_i(w^k) - \nabla f_i(x^*)}^2$.
\end{lemma}

Next, applying Theorem~\ref{thm:sk_main_gsgm} on Algorithm~\ref{alg:sk_L-SVRG} with $M=\frac{4}{p}$ we get Corollary~\ref{cor:sk_recover_l-svrg_rate}. 

\begin{corollary}\label{cor:sk_recover_l-svrg_rate}
    {\tt LSVRG} with $\alpha = \frac{1}{6L}$ satisfies
    \begin{equation}\label{eq:sk_recover_l-svrg_rate}
        \EEE V^k \le \left(1-\min\left\{\frac{\mu}{6L}, \frac{p}{2}\right\}\right)^kV^0.
    \end{equation}
\end{corollary}

\subsubsection*{Proof of Lemma~\ref{lem:sk_l-svrg}}\label{sec:sk_proofs_l-svrg}

    Let all expectations be conditioned on $x^k$. Using definition of $g^k$
    \begin{eqnarray*}
      &&  \EEE\left[\norm{g^k - \nabla f(x^*)}^2\right] \\
        & \overset{\text{Alg.}~\ref{alg:sk_L-SVRG}}{=}& 
        \EEE\left[\norm{
            \nabla f_i(x^k) - \nabla f_i(x^*) + \nabla f_i(x^*) - \nabla f_i(w^k) + \nabla f(w^k)
        - \nabla f(x^*)}^2\right]
        \\
           & \overset{\eqref{eq:sk_a_b_norm_squared}}{\leq}&
        2\EEE\left[\norm{\nabla f_i(x^k) - \nabla f_i(x^*)}^2\right]\\
        && \qquad \quad+ 
        2\EEE\left[\norm{\nabla f_i(x^*) - \nabla f_i(w^k) - \EEE\left[\nabla f_i(x^*) - \nabla f_i(w^k)\mid x^k\right]}^2\right]\\
           & \overset{\eqref{eq:sk_norm_diff_grads},\eqref{eq:sk_variance_decomposition}}{\leq}&
        4LD_f(x^k,x^*) + 2 \EEE\left[\norm{\nabla f_i(w^k) - \nabla f_i(x^*)}^2\right] \\
           &=&
        4LD_f(x^k,x^*) + 2 \sigma_k^2.
    \end{eqnarray*}
    For the second bound, we shall have
        \begin{eqnarray*}
        \EEE\left[\sigma_{k+1}^2\right] &\overset{\text{Alg.}~\ref{alg:sk_L-SVRG}}{=}& (1-p) \sigma_k^2 +  \frac{p}{n} \sum\limits_{i=1}^{n} \norm{\nabla f_i(x^k) - \nabla f_i(x^*)}^2\\
        &\overset{\eqref{eq:sk_norm_diff_grads}}{\leq}& (1-p) \sigma_k^2 + 2Lp D_f(x^k,x^*).
    \end{eqnarray*}

\subsection{{\tt DIANA}}\label{sec:sk_diana}
In this section we consider a distributed setup where each function $f_i$ from~\eqref{eq:sk_f_sum} is owned by $i$th machine (thus, we have all together $n$ machines). 

We show that our approach covers the analysis of {\tt DIANA} from \cite{mishchenko2019distributed, horvath2019stochastic}. {\tt DIANA} is a specific algorithm for distributed optimization with { \em quantization} -- lossy compression of gradient updates, which reduces the communication between the server and workers\footnote{It is a well-known problem in distributed optimization that the communication between machines often takes more time than actual computation.}. 

In particular, {\tt DIANA} quantizes gradient differences instead of the actual gradients. This trick allows for the linear convergence to the optimum once the full gradients are evaluated on each machine, unlike other popular quantization methods such as {\tt QSGD} \cite{alistarh2017qsgd} or {\tt TernGrad} \cite{wen2017terngrad}. In this case, {\tt DIANA} behaves as variance reduced method -- it reduces a variance that was injected due to the quantization. However, {\tt DIANA} also allows for evaluation of stochastic gradients on each machine, as we shall further see.

First of all, we introduce the notion of quantization operator.

\begin{definition}[Quantization]\label{def:quantization}
    We say that $\hat \Delta$ is a \textit{quantization} of vector $\Delta\in\R^d$ and write $\hat \Delta \sim {\rm Q}(\Delta)$ if
    \begin{equation}\label{eq:sk_quantization}
        \EEE\hat\Delta = \Delta, \qquad \EEE\norm{\hat \Delta - \Delta}^2 \le \omega \norm{\Delta}^2
    \end{equation}
    for some $\omega > 0$.
\end{definition}

\begin{algorithm}[t]
   \caption{{\tt DIANA} \cite{mishchenko2019distributed, horvath2019stochastic}}
   \label{alg:sk_diana}
\begin{algorithmic}[1]
   \Require learning rates $\gamma>0$ and $\alpha>0$, initial vectors $x^0, h_1^0,\dotsc, h_n^0 \in \R^d$ and $h^0 = \frac{1}{n}\sum_{i=1}^n h_i^0$
   \For{$k=0,1,2,\dotsc$}
       \State Broadcast $x^{k}$ to all workers
        \For{$i=1,\dotsc,n$ in parallel}
            \State Sample $g^{k}_i$ such that $\EEE [g^k_i \;|\; x^k]  =\nabla f_i(x^k)$ 
            \State $\Delta^k_i = g^k_i - h^k_i$
            \State Sample $\hat \Delta^k_i \sim {\rm Q}(\Delta^k_i)$
            \State $h_i^{k+1} = h_i^k + \gamma \hat \Delta_i^k$
            \State $\hat g_i^k = h_i^k + \hat \Delta_i^k$
        \EndFor
       \State $\hat \Delta^k = \frac{1}{n}\sum_{i=1}^n \hat \Delta_i^k$
       \State $ g^k = \frac{1}{n}\sum_{i=1}^n \hat g_i^k = h^k + \hat \Delta^k $
        \State $x^{k+1} = \proxR\left(x^k - \alpha  g^k \right)$
        \State $h^{k+1}  = \frac{1}{n}\sum_{i=1}^n h_i^{k+1} = h^k + \gamma \hat \Delta^k$
   \EndFor
\end{algorithmic}
\end{algorithm}

The aforementioned method is applied to solve problem \eqref{eq:sk_problem_gen}+\eqref{eq:sk_f_sum} where each $f_i$ is convex and $L$-smooth and $f$ is $\mu$-strongly convex. 

\begin{lemma}[Lemma 1 and consequence of Lemma 2 from \cite{horvath2019stochastic}]\label{lem:sk_lemma1_diana}
    Suppose that $\gamma \le \frac{1}{1+\omega}$. For all iterations $k\ge 0$ of Algorithm~\ref{alg:sk_diana} it holds
{
\footnotesize   
    \begin{eqnarray}
        \EEE\left[ g^k\mid x^k\right] &=& \nabla f(x^k), \label{eq:sk_unbiased_diana}\\
        \EEE\left[\norm{g^k - \nabla f(x^*)}^2\mid x^k\right] &\le& \left(1+\frac{2\omega}{n}\right)\frac{1}{n}\sum\limits_{i=1}^n\norm{\nabla f_i(x^k) - \nabla f_i(x^*)}^2\notag\\
        &&\quad + \frac{2\omega\sigma_k^2}{n} + \frac{(1+\omega)\sigma^2}{n},\label{eq:sk_second_moment_diana}
        \\
    \EEE\left[\sigma_{k+1}^2 \mid x^k\right] &\le& (1-\gamma)\sigma_k^2 + \frac{\gamma}{n}\sum\limits_{i=1}^n\norm{\nabla f_i(x^k) - \nabla f_i(x^*)}^2 + \gamma \sigma^2.\label{eq:sk_h_i_k_sec_moment_diana}
    \end{eqnarray}
    }
    where $\sigma_k^2  = \frac{1}{n}\sum\limits_{i=1}^n\norm{h_i^k -  \nabla f_i(x^*)}^2$ and $\sigma^2$ is such that $\frac{1}{n}\sum\limits_{i=1}^n\EEE\left[\norm{g_i^k - \nabla f_i(x^k)}^2\mid x^k\right]\le \sigma^2$.
\end{lemma}
Bounding further $\frac1n \sum_{i=1}^n\norm{\nabla f_i(x^k) - \nabla f_i(x^*)}^2 \leq 2L D_{f}(x^k,x^*)$ in the above Lemma, we see that Assumption~\ref{as:sk_general_stoch_gradient} as per Table~\ref{tbl:sk_special_cases-parameters} is valid. Thus, as a special case of Theorem~\ref{thm:sk_main_gsgm}, we obtain the following corollary.

\begin{corollary}\label{cor:sk_main_diana}
    Assume that $f_i$ is convex and $L$-smooth for all $i\in[n]$ and $f$ is $\mu$ strongly convex, $\gamma \le \frac{1}{\omega+1}$, $\alpha \le \frac{1}{\left(1+\frac{2\omega}{n}\right)L + ML\gamma}$ where $M > \frac{2\omega}{n\gamma}$. Then the iterates of {\tt DIANA} satisfy
    \begin{equation}\label{eq:sk_convergence_diana}
        \EEE\left[V^k\right] \le \max\left\{(1-\alpha\mu)^k, \left(1 + \frac{2\omega}{nM} - \gamma\right)^k\right\}V^0 + \frac{\left(\frac{1+\omega}{n} + M\gamma\right)\sigma^2\alpha^2}{\min\left\{\alpha\mu, \gamma - \frac{2\omega}{nM}\right\}},
    \end{equation}
    where the Lyapunov function $V^k$ is defined by $V^k \eqdef \norm{x^k - x^*}^2 + M\alpha^2\sigma_k^2$. For the particular choice $\gamma = \frac{1}{\omega+1}$, $M = \frac{4\omega(\omega+1)}{n}$, $\alpha = \frac{1}{\left(1 + \frac{6\omega}{n}\right)L}$, then {\tt DIANA} converges to a solution neighborhood and the leading iteration complexity term is
    \begin{equation}\label{eq:sk_diana_leading_term}
        \max\left\{\frac{1}{\alpha\mu}, \frac{1}{\gamma - \frac{2\omega}{nM}}\right\} = \max\left\{\kappa + \kappa \frac{6\omega}{n}, 2(\omega+1)\right\},
    \end{equation}
    where $\kappa = \frac{L}{\mu}$.
\end{corollary}

As mentioned, once the full (deterministic) gradients are evaluated on each machine, {\tt DIANA} converges linearly to the exact optimum. In particular, in such case we have $\sigma^2 = 0$. Corollary~\ref{cor:sk_main_diana_special_case} states the result in the case when $n=1$, i.e. there is only a single node~\footnote{node = machine}. For completeness, we present the mentioned simple case of {\tt DIANA} as Algorithm~\ref{alg:sk_diana_case}.

\begin{algorithm}[t]
   \caption{{\tt DIANA}: 1 node $\&$ exact gradients \cite{mishchenko2019distributed, horvath2019stochastic}}
   \label{alg:sk_diana_case}
\begin{algorithmic}[1]
   \Require learning rates $\gamma>0$ and $\alpha>0$, initial vectors $x^0, h^0 \in \R^d$
   \For{$k=0,1,2,\dotsc$}
       \State $\Delta^k = \nabla f(x^k) - h^k$
       \State Sample $\hat \Delta^k \sim {\rm Q}(\Delta^k)$
       \State $h^{k+1} = h^k + \gamma \hat \Delta^k$
       \State $g^k = h^k + \hat \Delta^k$
       \State $x^{k+1} = \proxR\left(x^k - \alpha  g^k \right)$
   \EndFor
\end{algorithmic}
\end{algorithm}

\begin{corollary}\label{cor:sk_main_diana_special_case}
    Assume that $f_i$ is $\mu$-strongly convex and $L$-smooth for all $i\in[n]$, $\gamma \le \frac{1}{\omega+1}$, $\alpha \le \frac{1}{\left(1+2\omega\right)L + ML\gamma}$ where $M > \frac{2\omega}{\gamma}$. Then the stochastic gradient $\hat g^k$ and the objective function $f$ satisfy Assumption~\ref{as:sk_general_stoch_gradient} with $A = \left(1+2\omega\right)L, B = 2\omega, \sigma_k^2 = \norm{h^k - h^*}^2, \rho = \gamma, C = L\gamma, D_1 = 0, D_2 = 0$ and 
    \begin{equation}\label{eq:sk_convergence_diana_special_case}
        \EEE\left[V^k\right] \le \max\left\{(1-\alpha\mu)^k, \left(1 + \frac{2\omega}{M} - \gamma\right)^k\right\}V^0,
    \end{equation}
    where the Lyapunov function $V^k$ is defined by $V^k \eqdef \norm{x^k - x^*}^2 + M\alpha^2\sigma_k^2$. For the particular choice $\gamma = \frac{1}{\omega+1}$, $M = 4\omega(\omega+1)$, $\alpha = \frac{1}{\left(1 + 6\omega\right)L}$ the leading term in the iteration complexity bound is
    \begin{equation}\label{eq:sk_diana_leading_term_special_case}
        \max\left\{\frac{1}{\alpha\mu}, \frac{1}{\gamma - \frac{2\omega}{M}}\right\} = \max\left\{\kappa + 6\kappa\omega, 2(\omega+1)\right\},
    \end{equation}
    where $\kappa = \frac{L}{\mu}$.
\end{corollary}

\subsection{{\tt Q-SGD-SR}}\label{Q-SGD-AS}

In this section, we consider a quantized version of {\tt SGD-SR}. 

\begin{algorithm}[h]
    \caption{{\tt Q-SGD-SR}}
    \label{alg:sk_qsgdas}
    \begin{algorithmic}
        \Require learning rate $\alpha>0$, starting point $x^0\in\R^d$, distribution $\cD$ over $\xi \in\R^n$ such that $\ED{\xi}$ is vector of ones
        \For{ $k=0,1,2,\ldots$ }
        \State{Sample $\xi \sim \cD$}
        \State{$g^k \sim {\rm Q}(\nabla f_\xi (x^k))$}
        \State{$x^{k+1} = \proxR(x^k - \alpha g^k)$}
        \EndFor
    \end{algorithmic}
\end{algorithm}

\begin{lemma}[Generalization of Lemma~2.4, \cite{pmlr-v97-qian19b}]\label{lem:sk_exp_smoothness_grad_up_bound_q-sgd-as}
    If $(f,\cD)\sim ES(\cL)$, then
    \begin{equation}\label{eq:sk_exp_smoothness_grad_up_bound_sgd-as}
        \ED{\norm{g^k - \nabla f(x^*)}^2} \le 4\cL(1+\omega)D_f(x^k,x^*) + 2\sigma^2(1+\omega).
    \end{equation}
    where $\sigma^2 \eqdef \ED{\norm{\nabla f_\xi(x^*)}^2}$.
\end{lemma}

A direct consequence of Theorem~\ref{thm:sk_main_gsgm} in this setup is Corollary~\ref{cor:sk_recover_q-sgd-as_rate}. 

\begin{corollary}\label{cor:sk_recover_q-sgd-as_rate}
    Assume that $f(x)$ is $\mu$-strongly quasi-convex and $(f,\cD)\sim ES(\cL)$. Then {\tt Q-SGD-SR} with $\alpha^k \equiv\alpha \le \frac{1}{2(1+\omega)\cL}$ satisfies
    \begin{equation}\label{eq:sk_recover_q-sgd-as_rate}
        \EEE\left[\norm{x^k - x^*}^2\right] \le (1-\alpha\mu)^k\norm{x^0-x^*}^2 + \frac{2\alpha(1+\omega)\sigma^2}{\mu}.
    \end{equation}
\end{corollary}

\subsubsection*{Proof of Lemma~\ref{lem:sk_exp_smoothness_grad_up_bound_q-sgd-as}}
In this proof all expectations are conditioned on $x^k$. First of all, from Lemma~\ref{lem:sk_exp_smoothness_grad_up_bound_sgd-as} we have
\begin{eqnarray*}
    \ED{\norm{\nabla f_\xi(x^k) - \nabla f(x^*)}^2} \le 4\cL D_f(x^k,x^*) + 2\sigma^2.
\end{eqnarray*}
The remaining step is to understand how quantization of $\nabla f_\xi(x^k)$ changes the above inequality if we put $g^k\sim {\rm Q}(\nabla f_\xi(x^k))$ instead of $\nabla f_\xi(x^k)$. Let us denote mathematical expectation with respect randomness coming from quantization by $\EEE_Q\left[\cdot\right]$. Using tower property of mathematical expectation we get
\begin{eqnarray*}
    \EEE\left[\|g^k - \nabla f(x^*)\|^2\right] &=& \EEE_{\cD}\left[\EEE_Q\|g^k - \nabla f(x^*)\|^2\right]\\
    &\overset{\eqref{eq:sk_variance_decomposition}}{=}& \EEE\left[\|g^k - \nabla f_\xi(x^k)\|^2\right] + \EEE\left[\|\nabla f_\xi(x^k) - \nabla f(x^*)\|^2\right]\\
    &\overset{\eqref{eq:sk_exp_smoothness_grad_up_bound_sgd-as}}{\le}&  \EEE\left[\|g^k - \nabla f_\xi(x^k)\|^2\right] + 4\cL D_f(x^k,x^*) + 2\sigma^2.
\end{eqnarray*} 
Next, we estimate the first term in the last row of the previous inequality
\begin{eqnarray*}
    \EEE\left[\|g^k - \nabla f_\xi(x^k)\|^2\right] &\overset{\eqref{eq:sk_quantization}}{\le}& \omega\EEE\left[\|\nabla f_\xi(x^k)\|^2\right]\\
    &\overset{\eqref{eq:sk_a_b_norm_squared}}{\le}& 2\omega\EEE\left[\|\nabla f_\xi(x^k) - \nabla f_\xi(x^*)\|^2\right] + 2\omega\EEE\left[\|\nabla f_\xi(x^*)\|^2\right]\\
    &\le& 4\omega\cL D_f(x^k,x^*) + 2\omega\sigma^2.
\end{eqnarray*}
Putting all together we get the result.

\subsection{{\tt VR-DIANA}} \label{sec:sk_VR-DIANA}
Corollary~\ref{cor:sk_main_diana} shows that once each machine evaluates a stochastic gradient instead of the full gradient, {\tt DIANA} converges linearly only to a certain neighborhood. In contrast, {\tt VR-DIANA}~\cite{horvath2019stochastic} uses a variance reduction trick within each machine, which enables linear convergence to the exact solution. In this section, we show that our approach recovers {\tt VR-DIANA} as well. 

\begin{algorithm}[t]
   \caption{{\tt VR-DIANA} based on {\tt LSVRG} (Variant 1), {\tt SAGA} (Variant 2), \cite{horvath2019stochastic}}
   \label{alg:sk_vr-diana}
\begin{algorithmic}[1]
        \Require{learning rates $\gamma > 0$ and $\alpha > 0$, initial vectors $x^0, h_{1}^0, \dots, h_{n}^0$, $h^0 = \frac{1}{n}\sum_{i=1}^n h_i^0$}
        \For{$k = 0,1,2,\ldots$}
        \State Sample random 
            $
                u^k = \begin{cases}
                    1,& \text{with probability } \frac{1}{m}\\
                    0,& \text{with probability } 1 - \frac{1}{m}\\
                \end{cases}
            $ \Comment{only for Variant 1}
        \State Broadcast $x^k$, $u^k$ to all workers\;
            \For{$i = 1, \ldots, n$ in parallel} \Comment{Worker side}
            \State Pick random $j_i^k \sim_{\rm u.a.r.} [m]$\;
            \State $\mu_i^k = \frac{1}{m} \sum\limits_{j=1}^{m} \nabla f_{ij}(w_{ij}^k)$\label{ln:mu} \;
            \State $g_i^k = \nabla f_{ij_i^k}(x^k) - \nabla f_{ij_i^k}(w_{ij_i^k}^k) + \mu_i^k$\;
            \State $\hat{\Delta}_i^k = Q(g_i^k - h_i^k)$\;
            \State $h_i^{k+1} = h_i^k + \gamma \hat{\Delta}_i^k$\;
                \For{$j = 1, \ldots, m$}
                    \State
                    $
                    w_{ij}^{k+1} =
                    \begin{cases}
                        x^k, & \text{if } u^k = 1 \\
                        w_{ij}^k, &\text{if } u^k = 0\\
                    \end{cases}
                    $ \Comment{Variant 1 ({\tt LSVRG}): update epoch gradient if $u^k = 1$}
                    \State
                    $
                    w_{ij}^{k+1} =
                    \begin{cases}
                    x^k, & j = j_i^k\\
                    w_{ij}^k, & j \neq j_i^k\\
                    \end{cases}
                    $ \Comment{Variant 2 ({\tt SAGA}): update gradient table}
                \EndFor
            \EndFor
            \State $h^{k+1} \! = \! h^k \!+\! \frac{\gamma}{n} \displaystyle \sum_{i=1}^n \hat{\Delta}_i^k$ \Comment{Gather quantized updates} 
            \State $g^k = \frac{1}{n}\sum\limits_{i=1}^{n} (\hat{\Delta}_i^k + h_i^k)$\;
            \State $x^{k+1} = x^k - \alpha g^k$\;
        \EndFor

\end{algorithmic}  
\end{algorithm}

The aforementioned method is applied to solve problem \eqref{eq:sk_problem_gen}+\eqref{eq:sk_f_sum} where each $f_i$ is also of a finite sum structure, as in \eqref{eq:sk_f_i_sum}, with  each $f_{ij}(x)$ being convex and $L$-smooth, and $f_i(x)$ being $\mu$-strongly convex. Note that $\nabla f(x^*) = 0$ and, in particular, $D_f(x,x^*) = f(x) - f(x^*)$ since the problem is considered without regularization.

\begin{lemma}[Lemmas 3, 5, 6 and 7 from \cite{horvath2019stochastic}]\label{lemmas_vr_diana}
    Let $\gamma \le \frac{1}{\omega+1}$. Then for all iterates $k\ge 0$ of Algorithm~\ref{alg:sk_vr-diana} the following inequalities hold:
    {
    \footnotesize
    \begin{eqnarray}
        \EEE\left[g^k\mid x^k\right] &=& \nabla f(x^k),\label{eq:sk_unbiased_g_k_vr_diana}\\
        \EEE\left[H^{k+1}\mid x^k\right] &\le& \left(1-\gamma\right)H^k + \frac{2\gamma}{m}D^k + 8\gamma Ln\left(f(x^k) - f(x^*)\right),\label{eq:sk_H_k+1_bound_vr_diana}\\
        \EEE\left[D^{k+1}\mid x^k\right] &\le& \left(1 - \frac{1}{m}\right)D^k + 2Ln\left(f(x^k) - f(x^*)\right),\label{eq:sk_D_k+1_bound_vr_diana}\\
        \EEE\left[\norm{g^k}^2\mid x^k\right] &\le& 2L\left(1+\frac{4\omega + 2}{n}\right)\left(f(x^k)-f(x^*)\right) + \frac{2\omega}{n^2}\frac{D^k}{m} + \frac{2(\omega+1)}{n^2}H^k,\label{eq:sk_second_moment_g_k_vr_diana}
    \end{eqnarray}
    }
    where $H^k = \sum\limits_{i=1}^n\norm{h_i^k - \nabla f_i(x^*)}^2$ and $D^k = \sum\limits_{i=1}^n\sum\limits_{j=1}^m\norm{\nabla f_{ij}(w_{ij}^k) - \nabla f_{ij}(x^*)}^2$.
\end{lemma}

\begin{corollary}\label{cor:sk_vr_diana_meets_assumption}
    Let $\gamma \le \min\left\{\frac{1}{3m},\frac{1}{\omega+1}\right\}$. Then stochastic gradient $\hat g^k$ (Algorithm~\ref{alg:sk_vr-diana}) and the objective function $f$ satisfy Assumption~\ref{as:sk_general_stoch_gradient} with $A = \left(1+\frac{4\omega + 2}{n}\right)L, B = \frac{2(\omega+1)}{n}, \rho = \gamma, C = L\left(\frac{1}{m}+4\gamma\right), D_1 = 0, D_2 = 0$ and
    \[
        \sigma_k^2 = \frac{H^k}{n} + \frac{D^k}{nm} = \frac{1}{n}\sum\limits_{i=1}^n\norm{h_i^k - \nabla f_i(x^*)}^2 + \frac{1}{nm}\sum\limits_{i=1}^n\sum\limits_{j=1}^m\norm{\nabla f_{ij}(w_{ij}^k) - \nabla f_{ij}(x^*)}^2.
    \]
\end{corollary}
\begin{proof}
    Indeed, \eqref{eq:sk_general_stoch_grad_unbias} holds due to \eqref{eq:sk_unbiased_g_k_vr_diana}. Inequality \eqref{eq:sk_general_stoch_grad_second_moment} follows from \eqref{eq:sk_second_moment_g_k_vr_diana} with $A = \left(1+\frac{4\omega + 2}{n}\right)L, B = \frac{2(\omega+1)}{n}, D_1 = 0, \sigma_k^2 = \frac{H^k}{n} + \frac{D^k}{nm}$ if we take into account that $\frac{2\omega}{n^2}\frac{D^k}{m} + \frac{2(\omega+1)}{n^2}H^k \le \frac{2(\omega+1)}{n}\left(\frac{D^k}{nm} + \frac{H^k}{n}\right)$. Finally, summing inequalities \eqref{eq:sk_H_k+1_bound_vr_diana} and \eqref{eq:sk_D_k+1_bound_vr_diana} and using $\gamma\le\frac{1}{3m}$
     { 
\footnotesize
    \begin{eqnarray*}
        \EEE\left[\sigma_k^2\mid x^k\right] &=& \frac{1}{n}\EEE\left[H^{k+1}\mid x^k\right] + \frac{1}{nm}\EEE\left[D^{k+1}\mid x^k\right]\\
        &\overset{\eqref{eq:sk_H_k+1_bound_vr_diana}+\eqref{eq:sk_D_k+1_bound_vr_diana}}{\le}& \left(1-\gamma\right)\frac{H^k}{n} + \left(1+2\gamma-\frac{1}{m}\right)\frac{D^k}{nm} + 2L\left(\frac{1}{m}+4\gamma\right)\left(f(x^k)-f(x^*)\right)\\
        &\le& \left(1-\gamma\right)\sigma_k^2 + 2L\left(\frac{1}{m}+4\gamma\right)\left(f(x^k)-f(x^*)\right)
    \end{eqnarray*}
    }
    we get \eqref{eq:sk_gsg_sigma} with $\rho = \gamma, C = L\left(\frac{1}{m}+4\gamma\right), D_2 = 0$.
\end{proof}

\begin{corollary}\label{cor:sk_main_vr_diana}
    Assume that $f_i$ is $\mu$-strongly convex and $f_{ij}$ is convex and $L$-smooth for all $i\in[n], j\in[m]$, $\gamma \le \min\left\{\frac{1}{3m},\frac{1}{\omega+1}\right\}$, $\alpha \le \frac{1}{\left(1+\frac{4\omega + 2}{n}\right)L + ML\left(\frac{1}{m}+4\gamma\right)}$ where $M > \frac{2(\omega+1)}{n\gamma}$. Then the iterates of {\tt VR-DIANA} satisfy
    \begin{equation}\label{eq:sk_convergence_vr_diana}
        \EEE\left[V^k\right] \le \max\left\{(1-\alpha\mu)^k, \left(1 + \frac{2(\omega+1)}{nM} - \gamma\right)^k\right\}V^0,
    \end{equation}
    where the Lyapunov function $V^k$ is defined by $V^k \eqdef \norm{x^k - x^*}^2 + M\alpha^2\sigma_k^2$. Further, if we set  $\gamma = \min\left\{\frac{1}{3m},\frac{1}{\omega+1}\right\}$, $M = \frac{4(\omega+1)}{n\gamma}$, $\alpha = \frac{1}{\left(1 + \frac{20\omega+18}{n} + \frac{4\omega+4}{n\gamma m}\right)L}$, then to achieve precision $\EEE\left[\norm{x^k-x^*}^2\right] \le \varepsilon V^0$ {\tt VR-DIANA} needs \[\cO\left(\max\left\{\kappa+\kappa\frac{\omega+1}{n}+\kappa\frac{(\omega+1)\max\left\{m,\omega+1\right\}}{nm},m,\omega+1\right\}\log\frac{1}{\varepsilon}\right)\] iterations, where $\kappa = \frac{L}{\mu}$.
\end{corollary}
\begin{proof}
    Using Corollary~\ref{cor:sk_vr_diana_meets_assumption} we apply Theorem~\ref{thm:sk_main_gsgm} and get the result.
\end{proof}

\begin{remark}
{\tt VR-DIANA} can be easily extended to the proximal setup in our framework.
\end{remark}

\subsection{{\tt JacSketch}} \label{sec:sk_JacSketch}

In this section, we show that our approach covers the analysis of {\tt JacSketch} from \cite{jacsketch}. {\tt JacSketch} is a generalization of {\tt SAGA} in the following manner. {\tt SAGA} observes every iteration $\nabla f_i(x)$ for random index $i$ and uses it to build both stochastic gradient as well as the control variates on the stochastic gradient in order to progressively decrease variance. In contrast, {\tt JacSketch}  observes every iteration the random sketch of the Jacobian, which is again used to build both stochastic gradient as well as the control variates on the stochastic gradient.

For simplicity, we do not consider proximal setup, since~\cite{jacsketch} does not either.

We first introduce the necessary notation (same as in \cite{jacsketch}). Denote first the Jacobian the objective \begin{equation}\label{eq:sk_jac_def}\Jac(x) \eqdef [\nabla f_1(x), \ldots, \nabla f_n(x)] \in \R^{d\times n}.\end{equation} 
Every iteration of the method, a random sketch of Jacobian $\nabla F(x^k)\mS$ (where $\mS\sim \cD$) is observed. Then, the method builds a variable $\mJ^k$, which is the current Jacobian estimate, updated using so-called sketch and project iteration~\cite{gower2015randomized}:
\[
\mJ^{k+1}  = \mJ^k(\mI - \Proj_{\mS_k}) + \Jac(x^k)\Proj_{\mS_k},
\]

where $\Proj_\mS$ is a projection under $\mW$ norm\footnote{Weighted Frobenius norm of matrix $\mX\in\R^{n\times n}$ with a positive definite weight matrix $\mW\in \R^{n\times n}$ is defined as 
$\norm{\mX}_{\mW^{-1}} \eqdef \sqrt{\Tr{ \mX \mW^{-1} \mX^\top}}.$
} ($\mW\in \R^{n\times n}$ is some positive definite weight matrix) defined as
$\Proj_\mS \eqdef  \mS (\mS^\top \mW \mS)^{\dagger} \mS^\top \mW$\footnote{Symbol $\dagger$ stands for Moore-Penrose pseudoinverse.}.

Further, in order to construct unbiased stochastic gradient, an access to the random scalar $\theta_{\mS}$ such that
\begin{equation}\label{eq:sk_unbiased}
\ED{\theta_{\mS} \Proj_\mS} \ones  =  \ones,
\end{equation}
where $e$ is the vector of all ones.

Next, the simplest option for the choice of the stochastic gradient is $\nabla f_{\mS}(x)$ -- an unbiased estimate of $\nabla f$ directly constructed using $\mS,\theta_{\mS} $:
\begin{equation}
\label{eq:sk_stochgradplain}
\nabla f_{\mS}(x)   = \frac{\theta_{\mS}}{n}\Jac(x) \Proj_\mS \ones.
\end{equation}

However, one can build a smarter estimate $ \nabla f_{\mS,\mJ}(x) $ via control variates constructed from $\mJ$:
\begin{equation}\label{eq:sk_controlgradJ} 
 \nabla f_{\mS,\mJ}(x) = \frac{\theta_{\mS}}{n} (\Jac(x)-\mJ)   \Proj_\mS\ones  + \frac{1}{n} \mJ \ones.
\end{equation}
The resulting algorithm is stated as Algorithm~\ref{alg:sk_jacsketch}.

\begin{algorithm}
    \begin{algorithmic}[1]
        \Require $\left(\cD, \mW, \theta_{\mS} \right)$, $x^0\in \R^d$, Jacobian estimate $\mJ^0 \in \R^{d \times n}$, stepsize $\alpha>0$         
        
        \For {$k =  0, 1, 2, \dots$}
        \State Sample a fresh copy $\mS_k\sim \cD$

        \State $ \mJ^{k+1}  = \mJ^k(\mI - \Proj_{\mS_k}) + \Jac(x^k)\Proj_{\mS_k}$
         \label{ln:jacupdate}        
        \State $g^{k} =  \nabla f_{\mS_k, \mJ^k}(x^k)$       \label{ln:gradupdate}    
                    
        \State $x^{k+1} = x^k - \alpha g^{k}$ \label{ln:xupdate}    
         
        \EndFor
    \end{algorithmic}
    \caption{{\tt JacSketch} \cite{jacsketch}}
    \label{alg:sk_jacsketch}
\end{algorithm}

Next we present Lemma~\ref{lem:sk_lemmas39_310_jacsketch} which directly justifies the parameter choice from Table~\ref{tbl:sk_special_cases2}. 

\begin{lemma}[Lemmas 2.5, 3.9 and 3.10 from \cite{jacsketch}]\label{lem:sk_lemmas39_310_jacsketch}
    Suppose that there are constants $\cL_1, \cL_2>0$ such that 
    \begin{eqnarray*}
      \ED{ \norm{ \nabla f_{\mS}(x) - \nabla f_{\mS}(x^*)}_2^2 } &\leq& 2  \cL_1 (f(x)-f(x^*)), \qquad \forall x\in \R^d \\
      \ED{\norm{(\Jac(x)-\Jac(x^*)) \Proj_{\mS}  }_{\mW^{-1}}^2 }& \leq & 2\cL_2 (f(x) -f(x^*)), \qquad \forall x\in \R^d, \label{eq:sk_ES2}
\end{eqnarray*}    
    
     Then
    \begin{equation}\label{eq:sk_jacs_contraction}
     \ED{\norm{\mJ^{k+1} -\Jac(x^*)}_{\mW^{-1}}^2} \le (1-\lambda_{\min}) \norm{\mJ^{k}-\Jac(x^*)}_{\mW^{-1}}^2 +  2\cL_2(f(x^k) -f(x^*)),
     \end{equation}
     \begin{equation}\label{eq:sk_gradbndsubdeltaXX}
    \ED{\norm{g^k}_2^2 } \le 4 \cL_1  (f(x^k) -f(x^*)) +  2 \frac{\lambda_{\max}}{n^2} \norm{\mJ^{k} -\Jac(x^*)}_{\mW^{-1}}^2,
    \end{equation}
    where $\lambda_{\min} = \lambda_{\min}\left(\ED{\Proj_{\mS}}\right)$ and $\lambda_{\max} = \lambda_{\max}\left( \mW^{1/2}\left( \ED{\theta_{\mS}^2 \Proj_{\mS} \ones \ones^\top \Proj_{\mS}} -\ones \ones^\top\right) \mW^{1/2}\right)$. Further,     $
    \ED{ \nabla f_{\mS,\mJ}(x)} = \nabla f(x)$.
\end{lemma}

Thus, as a direct consequence of Theorem~\ref{thm:sk_main_gsgm}, we obtain the next corollary.

\begin{corollary}\label{thm:sk_main_jacsketch}
Consider the setup from Lemma~\ref{lem:sk_lemmas39_310_jacsketch}. Suppose that $f$ is $\mu$-strongly convex and choose $\alpha \le \min\left\{\frac{1}{\mu},\frac{1}{2\cL_1 + M\frac{\cL_2}{n}}\right\}$ where $M > \frac{2\lambda_{\max}}{n\lambda_{\min}}$. Then the iterates of {\tt JacSketch} satisfy
    \begin{equation}\label{eq:sk_convergence_jacsketch}
        \EEE\left[V^k\right] \le \max\left\{(1-\alpha\mu)^k, \left(1 + \frac{2\lambda_{\max}}{nM} - \lambda_{\min}\right)^k\right\}V^0.
    \end{equation}

\end{corollary}

\subsection{Interpolation between methods ~\label{sec:sk_interpol}}

Given that a set of stochastic gradients satisfy Assumption~\ref{as:sk_general_stoch_gradient}, we show that an any convex combination of the mentioned stochastic gradients satisfy Assumption~\ref{as:sk_general_stoch_gradient} as well.

\begin{lemma}\label{lem:sk_convex_comb}
    Assume that sequences of stochastic gradients $\{g_1^k\}_{k\ge 0}, \ldots, \{g_m^k\}_{k\ge 0}$ at the common iterates $\{x^k\}_{k\ge 0}$ satisfy the Assumption~\ref{as:sk_general_stoch_gradient} with parameters \[A(j),B(j),\{\sigma_k^2(j)\}_{k\ge 0}, C(j),\rho(j),D_1(j),D_2(j), \qquad j\in[m]\] respectively. Then for any vector $\tau = (\tau_1,\ldots,\tau_m)^\top$ such as $\sum\limits_{j=1}^m\tau_j = 1$ and $\tau_j \ge 0, j\in[m]$ stochastic gradient $g_\tau^k = \sum\limits_{j=1}^m\tau_j g_j^k$ satisfies the Assumption~\ref{as:sk_general_stoch_gradient} with parameters:
    \begin{eqnarray}
        A_\tau = \sum\limits_{j=1}^m\tau_j A(j),\quad B_\tau = 1,\quad \sigma_{\tau,k}^2 = \sum\limits_{j=1}^m B(j)\tau_j \sigma_k^2(j),\quad \rho_\tau = \min\limits_{j\in[m]}\rho(j),\notag\\
        C_\tau = \sum\limits_{j=1}^m\tau_j C(j) B(j),\quad D_{\tau,1} = \sum\limits_{j=1}^m\tau_j D_1(j),\quad D_{\tau,2} = \sum\limits_{j=1}^m\tau_j D_2(j) B(j).\label{eq:sk_conv_comb_params}
    \end{eqnarray}
    Furthermore, if stochastic gradients $g_1^k, \dots, g_m^k$ are independent for all $k$, Assumption~\ref{as:sk_general_stoch_gradient} is satisfied with parameters
        \begin{eqnarray}
        A_\tau = L+\sum\limits_{j=1}^m\tau_j^2 A(j),\quad B_\tau = 1,\quad \sigma_{\tau,k}^2 = \sum\limits_{j=1}^m B(j)\tau_j^2 \sigma_k^2(j),\quad \rho_\tau = \min\limits_{j\in[m]}\rho(j),\notag\\
        C_\tau = \sum\limits_{j=1}^m\tau_j^2 C(j) B(j),\quad D_{\tau,1} = \sum\limits_{j=1}^m\tau_j^2 D_1(j),\quad D_{\tau,2} = \sum\limits_{j=1}^m\tau_j^2 D_2(j) B(j).\label{eq:sk_conv_comb_params_indep}
    \end{eqnarray}
\end{lemma}

What is more, instead of taking convex combination one can choose stochastic gradient at random. Lemma~\ref{lem:sk_flipping_a_coin} provides the result. 
\begin{lemma}\label{lem:sk_flipping_a_coin}
    Assume that sequences of stochastic gradients $\{g_1^k\}_{k\ge 0}, \ldots, \{g_m^k\}_{k\ge 0}$ at the common iterates $\{x^k\}_{k\ge 0}$ satisfy the Assumption~\ref{as:sk_general_stoch_gradient} with parameters \[A(j),B(j),\{\sigma_k^2(j)\}_{k\ge 0}, C(j),\rho(j),D_1(j),D_2(j), \qquad j\in[m],\]  respectively. Then for any vector $\tau = (\tau_1,\ldots,\tau_m)^\top$ such as $\sum\limits_{j=1}^m\tau_j = 1$ and $\tau_j \ge 0, j\in[m]$ stochastic gradient $g_\tau^k$ which equals $g_j^k$ with probability $\tau_j$ satisfies the Assumption~\ref{as:sk_general_stoch_gradient} with parameters:
    \begin{eqnarray}
        A_\tau = \sum\limits_{j=1}^m\tau_j A(j),\quad B_\tau =1,\quad \sigma_{\tau,k}^2 = \sum\limits_{j=1}^m\tau_j B(j) \sigma_k^2(j),\quad \rho_\tau = \min\limits_{j\in[m]}\rho(j),\notag\\
        C_\tau = \sum\limits_{j=1}^m\tau_j  B(j)C(j),\quad D_{\tau,1} = \sum\limits_{j=1}^m\tau_j D_1(j),\quad D_{\tau,2} = \sum\limits_{j=1}^mB(j)\tau_j D_2(j).\label{eq:sk_flipping_a_coin}
    \end{eqnarray}
    Furthermore, if stochastic gradients $g_1^k, \dots, g_m^k$ are independent for all $k$, Assumption~\ref{as:sk_general_stoch_gradient} is satisfied with parameters
        \begin{eqnarray}
        A_\tau = L+\sum\limits_{j=1}^m\tau_j^2 A(j),\quad B_\tau = 1,\quad \sigma_{\tau,k}^2 = \sum\limits_{j=1}^m B(j)\tau_j^2 \sigma_k^2(j),\quad \rho_\tau = \min\limits_{j\in[m]}\rho(j),\notag\\
        C_\tau = \sum\limits_{j=1}^m\tau_j^2 C(j) B(j),\quad D_{\tau,1} = \sum\limits_{j=1}^m\tau_j^2 D_1(j),\quad D_{\tau,2} = \sum\limits_{j=1}^m\tau_j^2 D_2(j) B(j).\label{eq:sk_flipping_a_coin_indep}
    \end{eqnarray}
\end{lemma}

\begin{example}[{\tt $\tau$-L-SVRG}]
    Consider the following method~--- {\tt $\tau$-L-SVRG}~--- which interpolates between vanilla {\tt SGD} and {\tt LSVRG}.
    \begin{algorithm}[h]
    \caption{{\tt $\tau$-L-SVRG}}
    \label{alg:sk_tau-L-SVRG}
    \begin{algorithmic}
        \Require learning rate $\alpha>0$, probability $p\in (0,1]$, starting point $x^0\in\R^d$, convex combination parameter $\tau\in[0,1]$
        \State $w^0 = x^0$
        \For{ $k=0,1,2,\ldots$ }
        \State{Sample  $i \in \{1,\ldots, n\}$ uniformly at random}
        \State{$g^k_{{\tt LSVRG}} = \nabla f_i(x^k) - \nabla f_i(w^k) + \nabla f(w^k)$}
        \State{Sample  $j \in \{1,\ldots, n\}$ uniformly at random}
        \State{$g^k_{{\tt SGD}} = \nabla f_j(x^k)$}
        \State{$g^k = \tau g^k_{{\tt SGD}} + (1-\tau)g^k_{{\tt LSVRG}}$}        
        \State{$x^{k+1} = x^k - \alpha g^k$}
        \State{$w^{k+1} = \begin{cases}
            x^{k}& \text{with probability } p\\
            w^k& \text{with probability } 1-p
            \end{cases}$
        }
        \EndFor
    \end{algorithmic}
    \end{algorithm}
    When $\tau = 0$ the Algorithm~\ref{alg:sk_tau-L-SVRG} becomes {\tt LSVRG} and when $\tau = 1$ it is just {\tt SGD} with uniform sampling.    Notice that Lemmas~\ref{lem:sk_l-svrg} and~\ref{lem:sk_exp_smoothness_grad_up_bound_sgd-as} still hold as they does not depend on the update rule for $x^{k+1}$.
    
    Thus, sequences $\{g_{SGD}^k\}_{k\ge 0}$ and $\{g_{L-SVRG}^k\}_{k\ge 0}$ satisfy the conditions of Lemma~\ref{lem:sk_convex_comb} and, as a consequence, stochastic gradient $g^k$ from {\tt $\tau$-L-SVRG} meets the Assumption~\ref{as:sk_general_stoch_gradient} with the following parameters:
    \begin{eqnarray*}
        A_\tau =  L + 2\tau^2\cL + 2(1-\tau)^2L,\quad B_\tau = 1,\quad \sigma_{\tau,k}^2 = 2\frac{(1-\tau)^2}{n}\sum\limits_{i=1}^n\norm{\nabla f_i(w^k) - \nabla f_i(x^*)}^2,\notag\\
         \rho_\tau = p,\quad C_\tau = 2(1-\tau)^2Lp,\quad D_{\tau,1} = 2\tau^2\sigma^2,\quad D_{\tau,2} = 0.
    \end{eqnarray*}
\end{example}

\begin{remark}
Similar interpolation with the analogous analysis can be considered between {\tt SGD} and {\tt SAGA}, or {\tt SGD} and {\tt SVRG}. 
\end{remark}

\subsubsection*{Proof of Lemma~\ref{lem:sk_convex_comb}}

    Indeed, \eqref{eq:sk_general_stoch_grad_unbias} holds due to linearity of mathematical expectation. Next, summing inequalities \eqref{eq:sk_general_stoch_grad_second_moment} for $g_1^k,\ldots,g_m^k$ and using convexity of $\norm{\cdot}^2$ we get
    \begin{eqnarray*}
        \EEE\left[\norm{g_\tau^k -\nabla f(x^*)}^2\mid x^k\right] &\le& \sum\limits_{j=1}^m\tau_j\EEE\left[\norm{g_j^k-\nabla f(x^*)}^2\mid x^k\right]
        \\
        & \overset{\eqref{eq:sk_general_stoch_grad_second_moment}}{\le}& 2\sum\limits_{j=1}^m\tau_j A(j) D_f(x^k,x^*) + \sum\limits_{j=1}^mB(j)\tau_j\sigma_k^2(j) +     \sum\limits_{j=1}^m\tau_j D_1(j),
    \end{eqnarray*}    
    which implies \eqref{eq:sk_general_stoch_grad_second_moment} for $g_\tau^k$ with $A_\tau = \sum\limits_{j=1}^m\tau_j A(j), B_\tau =1, \sigma_{\tau,k}^2 = \sum\limits_{j=1}^m\tau_j B(j)\sigma_k^2(j), D_{\tau,1} = \sum\limits_{j=1}^m\tau_j D_1(j)$.
    Finally, summing \eqref{eq:sk_gsg_sigma} for $g_1^k,\ldots,g_m^k$ gives us
    \begin{equation*}
        \EEE\left[\sigma_{\tau,k+1}^2\mid\sigma_{\tau,k}^2\right] \overset{\eqref{eq:sk_gsg_sigma}}{\le} \left(1-\min\limits_{j\in[m]}\rho(j)\right)\sigma_{\tau,k}^2 + 2\sum\limits_{j=1}^m\tau_j B(j)C(j)D_f(x^k,x^*) + \sum\limits_{j=1}^m\tau_j B(j)D_2(j),
    \end{equation*}
    which is exactly \eqref{eq:sk_gsg_sigma} for $\sigma_{\tau,k}^2$ with $\rho =\min\limits_{j\in[m]}\rho(j), C_\tau = \sum\limits_{j=1}^m\tau_j C(j), D_{\tau,2} = \sum\limits_{j=1}^m\tau_j D_2(j)$.

Next, for independent gradients we have
\begin{eqnarray}
\nonumber
& &\EEE\left[\norm{g_\tau^k -\nabla f(x^*)}^2\mid x^k\right]  \\
&  & \qquad =
 \sum_{j=1}^m \tau_j^2 \EEE \left[\norm{g_j^k -\nabla f(x^*)}^2\mid x^k\right] + 2\sum_{i< j} \tau_i \tau_j\EEE\<g_j^k -\nabla f(x^*),g_i^k -\nabla f(x^*) >
 \nonumber
 \\
&  & \qquad =
 \sum_{j=1}^m \tau_j^2 \EEE \left[\norm{g_j^k -\nabla f(x^*)}^2\mid x^k\right] + 2\sum_{i< j}\tau_i \tau_j \norm{\nabla f(x^k) -\nabla f(x^*)}^2
  \nonumber
 \\
&  & \qquad \leq
 \sum_{j=1}^m \tau_j^2 \EEE \left[\norm{g_j^k -\nabla f(x^*)}^2\mid x^k\right] + \left(\sum_{ j=1}^m \tau_j\right)^2 \norm{\nabla f(x^k) -\nabla f(x^*)}^2
  \nonumber
  \\
&  & \qquad =
 \sum_{j=1}^m \tau_j^2 \EEE \left[\norm{g_j^k -\nabla f(x^*)}^2\mid x^k\right] +  \norm{\nabla f(x^k) -\nabla f(x^*)}^2
  \nonumber
   \\
&  & \qquad \leq
 \sum_{j=1}^m \tau_j^2 \EEE \left[\norm{g_j^k -\nabla f(x^*)}^2\mid x^k\right] +  2LD_f(x^k,x^*).
 \label{eq:sk_indp_bounding}
\end{eqnarray}
and further the bounds follow. 

\subsubsection*{Proof of Lemma~\ref{lem:sk_flipping_a_coin}}

    Indeed, \eqref{eq:sk_general_stoch_grad_unbias} holds due to linearity and tower property of mathematical expectation. Next, using tower property of mathematical expectation and inequalities \eqref{eq:sk_general_stoch_grad_second_moment} for $g_1^k,\ldots,g_m^k$ we get
     { 
\footnotesize
    \begin{eqnarray*}
        \EEE\left[\norm{g_\tau^k-\nabla f(x^*)}^2\mid x^k\right] &=& \EEE\left[\EEE_\tau\left[\norm{g_\tau^k-\nabla f(x^*)}^2\right]\mid x^k\right] = \sum\limits_{j=1}^m\tau_j\EEE\left[\norm{g_j^k-\nabla f(x^*)}^2\mid x^k\right]\\ &\overset{\eqref{eq:sk_general_stoch_grad_second_moment}}{\le}& 2\sum\limits_{j=1}^m\tau_j A(j) D_f(x^k,x^*) + \sum\limits_{j=1}^m B(j)\tau_j\sigma_k^2(j) +     \sum\limits_{j=1}^m\tau_j D_1(j),
    \end{eqnarray*}
    }
    which implies \eqref{eq:sk_general_stoch_grad_second_moment} for $g_\tau^k$ with $A_\tau = \sum\limits_{j=1}^m\tau_j A(j), B_\tau = 1, \sigma_{\tau,k}^2 = \sum\limits_{j=1}^m\tau_j B(j) \sigma_k^2(j), D_{\tau,1} = \sum\limits_{j=1}^m\tau_j D_1(j)$.
    Finally, summing \eqref{eq:sk_gsg_sigma} for $g_1^k,\ldots,g_m^k$ gives us
    \begin{equation*}
        \EEE\left[\sigma_{\tau,k+1}^2\mid\sigma_{\tau,k}^2\right] \overset{\eqref{eq:sk_gsg_sigma}}{\le} \left(1-\min\limits_{j\in[m]}\rho(j)\right)\sigma_{\tau,k}^2 + 2\sum\limits_{j=1}^m\tau_j B(j) C(j)D_f(x^k,x^*) + \sum\limits_{j=1}^m\tau_jB(j)D_2(j),
    \end{equation*}
    which is exactly \eqref{eq:sk_gsg_sigma} for $\sigma_{\tau,k}^2$ with $\rho =\min\limits_{j\in[m]}\rho(j), C_\tau = \sum\limits_{j=1}^m\tau_j B(j)C(j), D_{\tau,2} = \sum\limits_{j=1}^m\tau_j B(j)D_2(j)$.
To show~\eqref{eq:sk_flipping_a_coin_indep}, it suffices to combine above bounds with the trick~\eqref{eq:sk_indp_bounding}.

\begin{remark}
Recently, \cite{tran2019hybrid} demonstrated in that the convex combination of {\tt SGD} and {\tt SARAH}~\cite{nguyen2017sarah} performs very well on non-convex problems. 
\end{remark}

\section{Proofs for Section~\ref{sec:sk_main_res}}

\subsection{Basic facts and inequalities}\label{sec:sk_basic_inequalities}

For all $a,b\in\R^d$ and $\xi > 0$ the following inequalities holds:
\begin{equation}\label{eq:sk_fenchel}
    \langle a,b\rangle \le \frac{\norm{a}^2}{2\xi} + \frac{\xi\norm{b}^2}{2},
\end{equation}
\begin{equation}\label{eq:sk_a_b_norm_squared}
    \norm{a+b}^2 \le 2\norm{a}^2 + 2\norm{b}^2,
\end{equation}
and
\begin{equation}\label{eq:sk_1/2a_minus_b}
    \frac{1}{2}\norm{a}^2 - \norm{b}^2 \le \norm{a+b}^2.
\end{equation}

For a random vector $\xi \in \R^d$ and any $x\in \R^d$ the variance can be decomposed as
\begin{equation}\label{eq:sk_variance_decomposition}
    \EEE\left[\norm{\xi - \EEE\xi}^2\right] = \EEE\left[\norm{\xi-x}^2\right] -\norm{\EEE\xi - x}^2 \;.
\end{equation}

\subsection{Proof of Lemma~\ref{lem:sk_iter_dec}}

    We start with estimating the first term of the Lyapunov function. Let $r^k = x^k - x^*$. Then
    \begin{eqnarray*}
        \norm{r^{k+1}}^2 &=& \norm{ \prox_{\alpha \psi} (x^k- \alpha  g^k) - \prox_{\alpha \psi} ( x^* - \alpha\nabla f(x^*)) }^2 \\
         &\leq & 
          \norm{x^k- x^* - \alpha (  g^k - \nabla f(x^*)) }^2
          \\
          &= & 
           \norm{r^k}^2 - 2\alpha\langle r^k,g^k  - \nabla f(x^*)\rangle + \alpha^2\norm{ g^k  - \nabla f(x^*)}^2.
    \end{eqnarray*}
    Taking expectation conditioned on $x^k$ we get
     { 
\footnotesize
    \begin{eqnarray*}
        \EEE\left[\norm{r^{k+1}}^2\mid x^k\right] &=& \norm{r^k}^2 - 2\alpha\langle r^k,\nabla f(x^k)-\nabla f(x^*)\rangle + \alpha^2\EEE\left[\norm{ g^k - \nabla f(x^*)}^2\mid x^k\right]
        \\
        &\overset{\eqref{eq:sk_mu_strongly_quasi_convex}}{\le}& 
        (1-\alpha\mu)\norm{r^k}^2 - 2\alpha D_f(x^k,x^*) + \alpha^2\EEE\left[\norm{g^k - \nabla f(x^*)}^2\mid x^k\right]
        \\
        &\overset{\eqref{eq:sk_general_stoch_grad_unbias}+\eqref{eq:sk_general_stoch_grad_second_moment}}{\le}& 
        (1-\alpha\mu)\norm{r^k}^2 + 2\alpha\left(A\alpha - 1\right)D_f(x^k,x^*) + B\alpha^2\sigma_k^2 + \alpha^2 D_1.
    \end{eqnarray*}
        }
    Using this we estimate the full expectation of $V^{k+1}$ in the following way:
    \begin{eqnarray}
    && \EEE\norm{x^{k+1}-x^*}^2 + M\alpha^2\EEE\sigma_{k+1}^2\notag
        \\
        &\overset{\eqref{eq:sk_gsg_sigma}}{\le}& (
        1-\alpha\mu)\EEE\norm{x^k-x^*}^2 + 2\alpha\left(A\alpha - 1\right)D_f(x^k,x^*) + B\alpha^2\EEE\sigma_k^2\notag
        \\
        &&\quad 
        + (1-\rho)M\alpha^2\EEE\sigma_k^2  + 2CM\alpha^2\EEE\left[D_f(x^k,x^*)\right] + (D_1+MD_2)\alpha^2\notag\\
        &=& (1-\alpha\mu)\EEE\norm{x^k - x^*}^2 + \left(1 + \frac{B}{M} - \rho\right)M\alpha^2\EEE\sigma_k^2\notag
        \\
        &&\quad
         + 2\alpha\left(\alpha(A+CM)-1\right)\EEE\left[D_f(x^k,x^*)\right] + (D_1+MD_2)\alpha^2\notag \,.
        \label{eq:sk_gsgm_recurrence}
    \end{eqnarray}

It remains to rearrange the terms.

\subsection{Proof of Theorem~\ref{thm:sk_main_gsgm}}

Note first that due to~\eqref{eq:sk_gamma_condition_gsgm} we have $2\alpha\left(1-\alpha(A+CM)\right)\EEE D_f(x^k,x^*)>0$, thus we can omit the term.

    Unrolling the recurrence from Lemma~\ref{lem:sk_iter_dec} and using the Lyapunov function notation gives us
    \begin{eqnarray*}
        \EEE V^{k} &\le& \max\left\{(1-\alpha\mu)^k,\left(1+\frac{B}{M}-\rho\right)^k\right\}V^0\\
        &&\quad + (D_1+MD_2)\alpha^2\sum\limits_{l=0}^{k-1}\max\left\{(1-\alpha\mu)^l,\left(1+\frac{B}{M}-\rho\right)^l\right\}\\
        &\le& \max\left\{(1-\alpha\mu)^k,\left(1+\frac{B}{M}-\rho\right)^k\right\}V^0\\
        &&\quad + (D_1+MD_2)\alpha^2\sum\limits_{l=0}^{\infty}\max\left\{(1-\alpha\mu)^l,\left(1+\frac{B}{M}-\rho\right)^l\right\}\\
        &\le& \max\left\{(1-\alpha\mu)^k,\left(1+\frac{B}{M}-\rho\right)^k\right\} V^0 + \frac{(D_1+MD_2)\alpha^2}{\min\left\{\alpha\mu, \rho - \frac{B}{M}\right\}}.
    \end{eqnarray*}

\chapter{Appendix for Chapter \ref{asvrcd}}
\label{asvrcd_appendix}

\graphicspath{{asvrcd/images/}}

 \section{Missing lemmas and proofs: {\tt SAGA}/{\tt LSVRG} is a special case of {\tt SEGA}/{\tt SVRCD}} 
 \subsection{Proof of Lemma~\ref{lem:asvrcd_equivalent_objectives} \label{sec:asvrcd_equivalent_objectivesproof}}
 Let $\Popt' \eqdef \frac1n \ee \ee^\top \otimes \mI$ and denote $\BD(\ueWMC)\eqdef \blockdiag(\ueWMC_{1}, \dots, \ueWMC_n)$ for simplicity.  Now clearly $x^0 \in \range{\Popt'}$, while $\Popt'$ is a projection matrix such that $\Ind{}(x)<\infty $ if and only if $\Popt'x=x$. Consequently, $\Popt=\Popt'$. 
 Next, if $x,y \in \range{\Popt}$, there is $\xx, \yy \in \R^\dd$ such that $x = \Lift{\xx}, y = \Lift{\yy}$. Therefore we can write
\begin{eqnarray*}
f(x) &=& f(\Popt(x)) = \frac1n \sum_{j=1}^n \ff_j(\xx) \\
& \geq&  \frac1n \sum_{j=1}^n \ff_j(\yy) + \< \nabla \left( \frac1n \sum_{j=1}^n \ff_j(\yy)\right), \xx-\yy  >  + \frac{\mmu}{2} \| \xx-\yy \|^2
\\
&=&
 f(y)+ \< \nabla f(y), x-y  >  + \frac{\mmu}{2n} \| x-y\|^2.
\end{eqnarray*}

Similarly,
\begin{eqnarray*}
f(x) &=& \frac1n \sum_{j=1}^n \ff_j(\xx) \\
& \leq&  \frac1n \sum_{j=1}^n \ff_j(\yy) + \< \nabla \left( \frac1n \sum_{j=1}^n\ff_j(\yy)\right), \xx-\yy  >  + \sum_{j=1}^n \frac{1}{2n} \| \xx-\yy \|^2_{\ueWMC_j}
\\
&=&
 f(y)+ \< \nabla f(y), x-y  >  + \frac{1}{2n} \|x-y\|^2_{\BD(\ueWMC)}.
\end{eqnarray*}

Thus we conclude $\mu = \frac{ \mmu}{n}$ and $\mM = \frac1n\BD(\ueWMC)$.  Further, for any $h\in \R^d$, we have:
\begin{eqnarray*}
&& h^\top \mM^{\frac12} \E{  \sum_{i\in S}\pLi^{-1}\eLi \eLi^\top\Popt  \sum_{i\in S} \pLi^{-1}\eLi \eLi^\top}\mM^{\frac12}  {h}  
\\
&& \qquad  =
\frac1n \|\BD(\ueWMC)^{\frac12}  {h} \|^2 _{\E{ \left( \sum_{i\in \sS}\pp_i^{-1}\left(\sum_{j\in R_i}e_j e_j^\top\right)\right)\Popt  \left( \sum_{i\in \sS}\pp_i^{-1}\left(\sum_{j\in R_i}e_j e_j^\top\right)\right)}} 
\\
&& \qquad =
 \frac{1}{n}\E{\left\| \sum_{i\in \sS }\ueWMC^{\frac12}_i  \pp_i^{-1}h_{R_i}   \right\|^2}
 \\
&& \qquad \stackrel{\eqref{eq:asvrcd_ESO_saga}}{\leq }
 \frac{1}{n}\sum_{i=1}^{n} \pp_i \vv_{i}\left\|h_{R_{i}}\right\|^{2}
\end{eqnarray*}

and thus~\eqref{eq:asvrcd_ESO_sega_good} holds with $v= \frac1n \vv$ as desired.

 \subsection{Proof of Lemma~\ref{lem:asvrcd_saga_from_sega} \label{sec:asvrcd_sagasegaproof}}
Denote $\Vect{\cdot}$ to be the vectorization operator, i.e.,  operator which takes a matrix as an input, and returns a vector constructed by a column-wise stacking of the matrix columns. We will show both \begin{equation}\label{eq:asvrcd_hj_equivalence}
h^k = \frac1n\Vect{\mJ^k}
\end{equation} 
and~\eqref{eq:asvrcd_iterates_equivalence} using mathematical induction. Clearly, if $k=0$ both~\eqref{eq:asvrcd_hj_equivalence} and~\eqref{eq:asvrcd_iterates_equivalence} hold. Now, let us proceed with the second induction step.

{ \footnotesize
\begin{eqnarray}\nonumber
x^{k+1} 
&=& 
 \prox_{\alpha \psi}(x^k - \alpha g^k) = \argmin_{x\in \R^d}\,  \alpha \Ind{}(x) +  \alpha\ppsi(x_{R_1}) + \|x - (x^k - \alpha g^k )\|^2
 \\
 \nonumber
&=& 
\argmin_{x\in \R^d}\,   \alpha \Ind{}(x) + \alpha \ppsi(x_{R_1}) + \left\|x- x^k + \alpha \left( h^k+\sum \limits_{i\in S}  \frac{1}{\pLi}(\nabla_i f(x^k) - h_i^k)\eLi\right) \right  \|^2
 \\
 \nonumber
&=& 
\argmin_{x = \Popt x}\,    \alpha \ppsi(x_{R_1}) + \left\|x- x^k + \alpha \left( h^k+\sum \limits_{i\in S}  \frac{1}{\pLi}(\nabla_i f(x^k) - h_i^k)\eLi\right) \right  \|^2
 \\
 \nonumber
&=& 
\argmin_{x = \Popt x}\,    \alpha \ppsi(x_{R_1}) + \left\|x- x^k + \alpha \left( h^k+\sum \limits_{i\in S}  \frac{1}{\pLi}(\nabla_i f(x^k) - h_i^k)\eLi\right) \right  \|^2_{\Popt}
 \\
 \nonumber
&\stackrel{\eqref{eq:asvrcd_iterates_equivalence}}{=}& 
\Lift{\argmin_{\xx \in \R^\dd}\,    \alpha \ppsi(\xx) + \frac1n \left\|n\xx- n\xx^k + \alpha \left( \sum_{i=1}^n h^k_{R_{i}}+\sum \limits_{i\in \sS}  \frac{1}{\pp_i} \left(  \frac1n\nabla \ff_i(\xx^k) - h_{R_i }^k\right)\right)\right  \|^2}
 \\
 \nonumber
&\stackrel{\eqref{eq:asvrcd_hj_equivalence}}{=}& 
\Lift{\argmin_{\xx \in \R^\dd}\,    \alpha \ppsi(\xx) + \frac1n \left\|n\xx- n\xx^k + \alpha \left( \frac1n\mJ^k\ee + \frac1n\sum \limits_{i\in \sS}  \frac{1}{\pp_i} \left((\nabla\ff_i(\xx^k) -\mJ^k_{:,i})\right)\right) \right  \|^2}
 \\
 \nonumber
&=& 
\Lift{\argmin_{\xx \in \R^\dd}\,   \aalpha \ppsi(\xx) +  \left\|\xx- \xx^k + \aalpha \left( \frac1n\mJ^k\ee + \frac1n\sum \limits_{i\in \sS}  \frac{1}{\pp_i} \left((\nabla\ff_i(\xx^k) -\mJ^k_{:,i})\right)\right) \right  \|^2}
 \\
&=& 
\Lift{\xx^{k+1}}. \label{eq:asvrcd_sequence}
\end{eqnarray}
}

It remains to notice that since $x^{k+1} = \Lift{\xx^k}$, we have $h^{k+1}= \frac1n\Vect{\mJ^{k+1}}$ as desired.

 \section{Missing lemmas and proofs: {\tt ASVRCD}}

\subsection{Technical lemmas}
We first start with two key technical lemmas.

\begin{lemma}
	Suppose that \begin{equation}
	\eta \leq \frac{1}{2L}. \label{eq:asvrcd_mkdkmd}
	\end{equation}
	Then, for all $x \in \range{\Popt}$ the following inequality holds:
	\begin{align} \nonumber
		& \frac{1}{\eta} \E{\<x - x^k,x^k - y^{k+1}> } \\
		& \qquad \qquad \leq 
		\E{F(x) - F(y^{k+1}) - \frac{1}{4\eta}\norm{y^{k+1} - x^k}^2 + \frac{\eta}{2}\norm{g^k - \nabla f(x^k)}^2_\Popt}
		\nonumber \\
		&\qquad \qquad \qquad \qquad -D_f(x,x^k). \label{eq:asvrcd_keylemma_acc}
	\end{align}
\end{lemma}

\begin{proof}
	From the definition of $y^{k+1}$ we get
	\begin{equation*}
		y^{k+1} = x^k - \eta g^k - \eta \Delta,
	\end{equation*}
	where $\Delta \in \partial \psi(y^{k+1})$.
	Therefore,
	\begin{align}\nonumber
	&	\E{\frac{1}{\eta}\<x - x^k, x^k - y^{k+1}>}
	\\
     &	\qquad =
		\E{\<x - x^k, g^k + \Delta>}
		\nonumber \\ \nonumber
		&	\qquad=
		\<x - x^k, \nabla f(x^k)>
		+
		\E{\<x - y^{k+1}, \Delta> + \< y^{k+1} - x^k, \Delta>}\\
		&	\qquad \leq
		f(x) - f(x^k) - D_f(x,x^k)
		+
		\E{\psi(x) - \psi(y^{k+1})}
		+
		\E{\< y^{k+1} - x^k, \Delta>} \label{eq:asvrcd_dnjansdjkajksd}
	\end{align}
	Now, we use the fact that $f$ is $L$-smooth over the set where iterates live (i.e.,  over $\{ x^0 + \range{\Popt}\}$):
	\begin{eqnarray}\nonumber
			f(y^{k+1})& \leq& f(x^k) + \<\nabla f(x^k), y^{k+1} - x^k> + \frac{L}{2}\norm{y^{k+1} - x^k}^2 \\
		& = & f(x^k) + \<\Popt\nabla f(x^k), y^{k+1} - x^k> + \frac{L}{2}\norm{y^{k+1} - x^k}^2.\label{eq:asvrcd_dabhusdbhu}
	\end{eqnarray}
Thus, we have
 { 
\footnotesize
	\begin{eqnarray*}
&&	\E{\frac{1}{\eta}\<x - x^k, x^k - y^{k+1}>} \\
	&& \,\,\,\,\,  \stackrel{\eqref{eq:asvrcd_dnjansdjkajksd}+\eqref{eq:asvrcd_dabhusdbhu}}{\leq}
	\E{F(x) - F(y^{k+1})
	+
	\< y^{k+1} - x^k, \Popt(\Delta + \nabla f(x^k))>
	+
	\frac{L}{2}\norm{y^{k+1} - x^k}^2} \\
	&& \qquad \qquad - D_f(x,x^k)
	 \\
	&& \qquad =
	\E{F(x) - F(y^{k+1})
		+
		\< y^{k+1} - x^k, \Popt(\nabla f(x^k) - g^k)>
		-
		\frac{1}{\eta}\norm{y^{k+1} - x^k}^2}
\\
&& \qquad
 \qquad 		
		+\E{\frac{L}{2}\norm{y^{k+1} - x^k}^2}- D_f(x,x^k) \\
 	&& \qquad \leq
	\E{F(x) - F(y^{k+1})
		+
		\frac{\eta}{2}\norm{\nabla f(x^k) - g^k}^2_\Popt
		-
		\frac{1}{2\eta}\norm{y^{k+1} - x^k}^2
		+
		\frac{L}{2}\norm{y^{k+1} - x^k}}\\
&& \qquad
  \qquad - D_f(x,x^k)\\
 	&& \qquad \stackrel{\eqref{eq:asvrcd_mkdkmd}}{\leq}
		\E{F(x) - F(y^{k+1})
		-
		\frac{1}{4\eta}\norm{y^{k+1} - x^k}^2
		+
		\frac{\eta}{2}\norm{\nabla f(x^k) - g^k}^2_\Popt} - D_f(x,x^k),  \\
	\end{eqnarray*}
	}
	 	which concludes the proof. 
\end{proof}

\begin{lemma}
	Suppose, the following choice of parameters is used:
	\begin{equation*}
		\eta =  \frac14 \max\{\Lcac, L\}^{-1},\quad
		\gamma = \frac{1}{\max\{2\mu, 4\theta_1/\eta\}},\quad
		\beta = 1 - \gamma\mu,\quad
		\theta_2 = \frac{\Lcac}{2\max\{L, \Lcac\}}.
	\end{equation*}
	Then the following inequality holds:
	\begin{align}\nonumber
		&\E{\norm{z^{k+1} - x^*}^2 + \frac{2\gamma\beta}{\theta_1}\left[F(y^{k+1}) - F(x^*)\right]}\\
		& \qquad \leq
		\beta \norm{z^k - x^*}^2
		+
		\frac{2\gamma\beta\theta_2}{\theta_1}\left[F(w^k)  - F(x^*)\right]
		+
		\frac{2\gamma\beta(1-\theta_1-\theta_2)}{\theta_1}
		\left[F(y^k) - F(x^*)\right].\label{eq:asvrcd_nhivbhi}
	\end{align}
\end{lemma}
\begin{proof}
{\footnotesize
	\begin{align*}
		& \E{\norm{z^{k+1} - x^*}^2}
		\\
		& \qquad =
		\E{\norm{\beta z^k + (1-\beta)x^k - x^* + \frac{\gamma}{\eta}(y^{k+1} - x^k)}^2}
		\\
		& \qquad\leq
		\beta \norm{z^k - x^*}^2
		+
		(1-\beta)\norm{x^k - x^*}^2
		+
		\frac{\gamma^2}{\eta^2}\E{\norm{y^{k+1} - x^k}^2}
		\\
		& \qquad\qquad +
		\frac{2\gamma}{\eta}\E{\<y^{k+1} - x^k, \beta z^k + (1-\beta)x^k - x^*>}
		\\
		&\qquad=
		\beta \norm{z^k - x^*}^2
		+
		(1-\beta)\norm{x^k - x^*}^2
		+
		\frac{\gamma^2}{\eta^2}\E{\norm{y^{k+1} - x^k}^2}
		\\
		& \qquad\qquad +
		\frac{2\gamma}{\eta}\E{\<y^{k+1} - x^k, x^k - x^*>}
		+
		\frac{2\gamma\beta}{\eta}\E{\<y^{k+1} - x^k, z^k - x^k>}
		\\
		&\qquad=	
		\beta \norm{z^k - x^*}^2
		+
		(1-\beta)\norm{x^k - x^*}^2
		+
		\frac{\gamma^2}{\eta^2}\E{\norm{y^{k+1} - x^k}^2}
		+
		\frac{2\gamma}{\eta}\E{\<x^k - y^{k+1}, x^* - x^k>}
		\\
		& \qquad\qquad +
		\frac{2\gamma\beta\theta_2}{\eta\theta_1}\E{\<x^k - y^{k+1} , w^k - x^k>}
		+
		\frac{2\gamma\beta( 1  - \theta_1 -\theta_2 )}{\eta\theta_1}\E{\<x^k - y^{k+1} , y^k - x^k>}
		\\
		&\qquad\stackrel{\eqref{eq:asvrcd_keylemma_acc}}{\leq}
		\beta \norm{z^k - x^*}^2
		+
		(1-\beta)\norm{x^k - x^*}^2
		+
		\frac{\gamma^2}{\eta^2}\E{\norm{y^{k+1} - x^k}^2}
		\\
		&\qquad\qquad		+
		2\gamma \E{F(x^*) - F(y^{k+1}) - \frac{1}{4\eta}\norm{y^{k+1} - x^k}^2 - D_f(x^*,x^k) + \frac{\eta}{2}\norm{g^k-\nabla f(x^k)}^2_\Popt}
		\\
		& \qquad\qquad +
		\frac{2\gamma\beta\theta_2}{\theta_1}\E{F(w^k) - F(y^{k+1}) - \frac{1}{4\eta}\norm{y^{k+1} - x^k}^2 - D_f(w^k,x^k) + \frac{\eta}{2}\norm{g^k-\nabla f(x^k)}^2_\Popt}
		\\
		&\qquad\qquad +
		\frac{2\gamma\beta( 1  - \theta_1 -\theta_2 )}{\theta_1}\E{F(y^k) - F(y^{k+1}) - \frac{1}{4\eta}\norm{y^{k+1} - x^k}^2 + \frac{\eta}{2}\norm{g^k - \nabla f(x^k)}^2_\Popt}
		\\
		&\qquad\stackrel{\eqref{eq:asvrcd_sc}}{\leq}
		\beta \norm{z^k - x^*}^2
		+
		(1-\beta - \gamma\mu)\norm{x^k - x^*}^2
		+
		\frac{\gamma^2}{\eta^2}\E{\norm{y^{k+1} - x^k}^2}
		\\
		&\qquad\qquad 		+
		2\gamma\beta \E{F(x^*) - F(y^{k+1}) - \frac{1}{4\eta}\norm{y^{k+1} - x^k}^2} +  \eta\gamma\E{\norm{g^k-\nabla f(x^k)}^2_\Popt}
		\\
		&\qquad \qquad+
		\frac{2\gamma\beta\theta_2}{\theta_1}\E{F(w^k) - F(y^{k+1}) - \frac{1}{4\eta}\norm{y^{k+1} - x^k}^2 -  D_f(w^k,x^k)+ \frac{\eta}{2}\norm{g^k-\nabla f(x^k)}^2_\Popt}
		\\
		&\qquad \qquad +
		\frac{2\gamma\beta(1  - \theta_1 -\theta_2)}{\theta_1}\E{F(y^k) - F(y^{k+1}) - \frac{1}{4\eta}\norm{y^{k+1} - x^k}^2 + \frac{\eta}{2}\norm{g^k - \nabla f(x^k)}^2_\Popt}.
	\end{align*}
	}
		Using $\beta = 1 - \gamma\mu$ we get
		{ \footnotesize
	\begin{align*}
	\E{\norm{z^{k+1} - x^*}^2}
	&\leq
	\beta \norm{z^k - x^*}^2
	+
	\left[\frac{\gamma^2}{\eta^2} - \frac{\gamma\beta}{2\eta\theta_1}\right]\E{\norm{y^{k+1} - x^k}^2}
	+
	\frac{\eta\gamma}{\theta_1}\E{\norm{g^k - \nabla f(x^k)}^2_\Popt}
		\\
	&
	\qquad -
	\frac{2\gamma\beta\theta_2}{\theta_1}D_f(w^k,x^k)
	+
	2\gamma\beta \E{F(x^*) - F(y^{k+1})}\\
	&\qquad
	+
	\frac{2\gamma\beta\theta_2}{\theta_1}\E{F(w^k) - F(y^{k+1})}
	+
	\frac{2\gamma\beta(1-  \theta_1 -\theta_2 )}{\theta_1}\E{F(y^k) - F(y^{k+1})}.
	\end{align*}
	}
	Using stepsize $\gamma \leq \frac{\beta\eta}{2\theta_1}$ we get
	{ \footnotesize
	\begin{align*}
	\E{\norm{z^{k+1} - x^*}^2}
	&\leq
	\beta \norm{z^k - x^*}^2
	+
	\frac{\eta\gamma}{\theta_1}\E{\norm{g^k - \nabla f(x^k)}^2_\Popt}
	-
	\frac{2\gamma\beta\theta_2}{\theta_1}D_f(w^k,x^k)
		\\
	&
	\qquad
		+
	2\gamma\beta \E{F(x^*) - F(y^{k+1})}
  +
	\frac{2\gamma\beta\theta_2}{\theta_1}\E{F(w^k) - F(y^{k+1})}
	+
			\\
	&
	\qquad
	\frac{2\gamma\beta(1  - \theta_1 - \theta_2)}{\theta_1}\E{F(y^k) - F(y^{k+1})}.
	\end{align*}
	}
	Now, using the expected smoothness from inequality~\eqref{eq:asvrcd_exp:smooth}:
	\begin{equation}
		\E{\norm{g^k - \nabla f(x^k)}^2_\Popt} \leq 2\Lcac D_f(w^k,x^k)
	\end{equation}
	 and stepsize $\eta \leq \frac{\beta\theta_2}{\Lcac}$ we get
	 { \footnotesize
	 \begin{align*}
	 \E{\norm{z^{k+1} - x^*}^2}
	 &\leq
	 \beta \norm{z^k - x^*}^2
   +
	 \frac{2\Lcac\eta\gamma }{\theta_1} D_f(w^k,x^k)
	 -
	 \frac{2\gamma\beta\theta_2}{\theta_1} D_f(w^k,x^k)
	 +
	 2\gamma\beta \E{F(x^*) - F(y^{k+1})}
	 \\
	 &
	 \qquad 
	 +
	 \frac{2\gamma\beta\theta_2}{\theta_1}\E{F(w^k) - F(y^{k+1})}
	 +
	 \frac{2\gamma\beta(1 - \theta_1- \theta_2 )}{\theta_1}\E{F(y^k) - F(y^{k+1})}\\
	 &\leq
	  \beta \norm{z^k - x^*}^2
	  +
	  2\gamma\beta \E{F(x^*) - F(y^{k+1})}
	  +
	  \frac{2\gamma\beta\theta_2}{\theta_1}\E{F(w^k) - F(y^{k+1})}\\
	  & \qquad +
	  \frac{2\gamma\beta( 1 - \theta_1- \theta_2 )}{\theta_1}\E{F(y^k) - F(y^{k+1})}\\
	  &=
	  \beta \norm{z^k - x^*}^2
	  -
	  \frac{2\gamma\beta}{\theta_1}\E{F(y^{k+1}) - F(x^*)}
	  +
	  \frac{2\gamma\beta\theta_2}{\theta_1}\left[F(w^k)  - F(x^*)\right] \\
	  & \qquad +
	  \frac{2\gamma\beta( 1-\theta_1-\theta_2)}{\theta_1}
	  \left[F(y^k) - F(x^*)\right].
	 \end{align*}
	 }
	 It remains to rearrange the terms.
\end{proof}

\subsection{Proof of Theorem~\ref{thm:asvrcd_acc}} 
 
	One can easily show that
	\begin{equation}\label{eq:asvrcd_bfrbuf}
		\E{F(w^{k+1})} = \probx F(y^k) + (1- \probx)F(w^k).
	\end{equation}
	Using that, we obtain
	{\footnotesize
	\begin{eqnarray*}
		\E{\Psi^{k+1}}
		&\stackrel{\eqref{eq:asvrcd_nhivbhi}+\eqref{eq:asvrcd_bfrbuf}}{\leq}&
		\beta \norm{z^k - x^*}^2
		+
		\frac{2\gamma\beta\theta_2}{\theta_1}\left[F(w^k)  - F(x^*)\right]
		+
		\frac{2\gamma\beta( 1-\theta_1 -\theta_2)}{\theta_1}
		\left[F(y^k) - F(x^*)\right]\\
		&& \qquad +
		\frac{(2\theta_2 + \theta_1)\gamma\beta}{\theta_1\probx}\left[ \probx F(y^{k}) + (1- \probx)F(w^k) - F(x^*)\right]\\
		&=&
		\beta \norm{z^k - x^*}^2
		+
		\frac{2\gamma\beta(1- \theta_1/2)}{\theta_1}
		\left[F(y^k) - F(x^*)\right]
		\\
		&& \qquad 
		+
		\frac{(2\theta_2 + \theta_1)\gamma\beta}{\theta_1\probx}\left[1 -  \probx + \frac{2\probx\theta_2}{2\theta_2 + \theta_1}\right]\left[F(w^k) - F(x^*)\right]\\
		&\leq&
		\max\left\{1 - \frac{1}{\max\{2, 4\theta_1/(\eta\mu)\}}, 1 -  \frac{\theta_1}{2}, 1 - \frac{ \probx\theta_1}{2\max \{2\theta_2,\theta_1\}} \right\}\Psi^k\\
		&=&
		 \left[1 -  \max\left\{\frac{2}{\probx},\frac{4}{\theta_1}\max\left\{\frac{1}{2}, \frac{\theta_2}{\rho}\right\}, \frac{4\theta_1}{\eta\mu}    \right\}^{-1} \right]\Psi^k.
	\end{eqnarray*}
	}
	Using $\theta_1 = \min\left\{\frac{1}{2},\sqrt{\eta\mu \max\left\{\frac{1}{2}, \frac{\theta_2}{\rho}\right\}}\right\} $ we get
	\begin{align*}
		\E{\Psi^{k+1}} &\leq
		\left[1 -  \max\left\{\frac{2}{\probx},8\max\left\{\frac{1}{2}, \frac{\theta_2}{\rho}\right\}, 4\sqrt{\frac{\max\left\{\frac{1}{2}, \frac{\theta_2}{\rho}\right\}}{\eta\mu}} \right\}^{-1} \right]\Psi^k\\
		&\leq
		\left[1 -  \frac{1}{4}\max\left\{\frac{1}{\probx}, \sqrt{\frac{2\max\left\{L, \frac{\Lcac}{\rho}\right\}}{\mu}} \right\}^{-1} \right]\Psi^k,
	\end{align*}
as desired.

\subsection{Proof of Lemma~\ref{lem:asvrcd_exp_smooth_vs_ESO}}

To establish that that we can choose $\Lcac = \ccL$, it suffices to see
{\small
\begin{eqnarray*}
\E{\norm{g^k - \nabla f(x^k)}^2_\Popt}  &= & \E{\norm{  \sum \limits_{i\in S}  \frac{1}{\pLi}(\nabla_i f(x^k) - \nabla_i f(w^k))\ones_i + \nabla f(w^k) - \nabla f(x^k)}^2_\Popt}  \\
& \leq  & 
 \E{\norm{  \sum \limits_{i\in S}  \frac{1}{\pLi}(\nabla_i f(x^k) - \nabla_i f(w^k))\ones_i}^2_\Popt}   \\
   &\stackrel{\eqref{eq:asvrcd_ccLdef}}{\leq} & 
\ccL \norm{ \nabla f(x^k) - \nabla f(w^k) }^2_{\mM^{-1}}\\
     &\stackrel{\eqref{eq:asvrcd_smooth}}{\leq} & 
 2 \ccL D_f(w^k,x^k).
\end{eqnarray*}
}
Next, to establish $\ccL \geq L $, let $\mQ \eqdef   \sum_{i\in S} \frac{1}{\pLi} \eLi \eLi^\top\Popt$. Consequently, we get

\begin{eqnarray*}
\ccL &\stackrel{\eqref{eq:asvrcd_ccLdef}}{=}&
 \lambda_{\max}\left( \mM^{\frac12} \E{  \sum_{i\in S} \frac{1}{\pLi} \eLi \eLi^\top\Popt \sum_{i\in S} \frac{1}{\pLi} \eLi \eLi^\top} \mM^{\frac12} \right)
 \\
 &=&
  \lambda_{\max}\left( \mM^{\frac12} \E{  \sum_{i\in S} \frac{1}{\pLi} \eLi \eLi^\top\Popt^2 \sum_{i\in S} \frac{1}{\pLi} \eLi \eLi^\top} \mM^{\frac12} \right)
 \\
 & =&
  \lambda_{\max}\left( \mM^{\frac12} \E{  \mQ \mQ^\top} \mM^{\frac12} \right)
 \\
& \geq &
   \lambda_{\max}\left( \mM^{\frac12} \E{  \mQ}\E{ \mQ^\top} \mM^{\frac12} \right)
 \\
  &=&
  \lambda_{\max}\left( \mM^{\frac12}\Popt^2 \mM^{\frac12} \right)
   \\
  &=&
  \lambda_{\max}\left( \mM^{\frac12}\Popt \mM^{\frac12} \right)
     \\
  &=&
 L,
\end{eqnarray*}
as desired. 

\subsection{Proof of Lemma~\ref{lem:asvrcd_acc_example}}
Let us look first at $\Popt = \mI$. In such case, it is easy to see that 
\begin{align*}
\E{\norm{g^k - \nabla f(x^k)}^2_{\Popt}} &= \E{\norm{g^k - \nabla f(x^k)}^2} \\
 &\leq  \E{\norm{d ( \nabla_i f(x^k) - \nabla_i f(w^k) ) e_i }^2}   \\
&= d \norm{ \nabla f(x^k) - \nabla f(w^k) }^2 \\
&\leq  2d \lambda_{\max} \mM D_f(w^k,x^k),
\end{align*}
 i.e., we can choose $\Lcac = d \lambda_{\max} \mM$. Noting that $ \lambda_{\max} \mM\geq L$, the iteration complexity of Algorithm~\ref{alg:asvrcd_acc} is $\cO\left(d\sqrt{\frac{L}{\mu}} \log \frac1\epsilon\right)$. On the other hand, if $\Popt = \frac1d \ee^\top$, we have 
\begin{align*}
\E{\norm{g^k - \nabla f(x^k)}^2_{\Popt}} =&\E{\norm{g^k - \nabla f(x^k)}^2_{\frac1d ee^\top}} \\
 &  \leq  \E{\norm{d ( \nabla_i f(x^k) - \nabla_i f(w^k) ) e_i }^2_{\frac1d ee^\top}}   \\
&=  \norm{ \nabla f(x^k) - \nabla f(w^k) }^2 \\
&\leq  2 \lambda_{\max} \mM D_f(w^k,x^k),
\end{align*}
and therefore $\Lcac =  \lambda_{\max} \mM$, which yields $\cO\left(\sqrt{\frac{d \lambda_{\max} \mM}{\mu}} \log \frac1\epsilon\right)$ convergence rate.

 \section{Missing lemmas and proofs: L-Katyusha as a particular case of {\tt ASVRCD}}

\subsection{Proof of Lemma~\ref{lem:asvrcd_katyusha_from_asvrcd}}
Let us proceed by induction. We will show the following for all $k\geq 0$ we have 
\begin{eqnarray} \nonumber
\xx^k = x^k_{R_1}=\dots =x^k_{R_n}, &&  \yy^k= y^k_{R_1}=\dots =y^k_{R_n}, \\
  \zzz^k = z^k_{R_1}=\dots =z^k_{R_n} & \mathrm{and} & \ww^k = w^k_{R_1}=\dots =w^k_{R_n}.
\label{eq:asvrcd_induction_acc}
\end{eqnarray}

Clearly, for $k=0$, the above claim holds. Let us proceed with the second induction step and assume that~\eqref{eq:asvrcd_induction_acc} holds for some $k\geq 0$. First, the update rule on $\{x^k\}$ for {\tt ASVRCD} together with the update rule on $\{\xx^k\}$ yields 
\begin{equation}\label{eq:asvrcd_induction_acc_x}
\xx^{k+1} = x^{k+1}_{R_1}=\dots= x^{k+1}_{R_n}. 
\end{equation}

To show 
\begin{equation}\label{eq:asvrcd_induction_acc_y}
\yy^{k+1}= y^{k+1}_{R_1}=\dots =y^{k+1}_{R_n},
\end{equation}

we essentially repeat the proof of Lemma~\ref{lem:asvrcd_saga_from_sega}. In particular, it is sufficient to repeat the sequence of inequalities~\eqref{eq:asvrcd_sequence} where variables \[(x^{k+1}, \xx^{k+1}, h^k, \mJ^{k} \alpha, \aalpha)\] are replaced by \[(y^{k+1}, \yy^{k+1},\nabla f(w^k), [\nabla \ff_1(\ww^k), \dots,\nabla \ff_n(\ww^k) ], \eta, \eeta),\] respectively. 

Next, $ \zzz^{k+1}= z^{k+1}_{R_1}=\dots= z^{k+1}_{R_n}$ follows from~\eqref{eq:asvrcd_induction_acc},~\eqref{eq:asvrcd_induction_acc_x} and~\eqref{eq:asvrcd_induction_acc_y} together with the update rule (on $\{z^k\}$ and $\{\zzz^k\}$) of both algorithms and the fact that $\frac{\gamma}{\eta} = \frac{\ggamma}{\eeta}$.

To finish the proof of the algorithms equivalence, we shall notice that $ \ww^{k+1}= w^{k+1}_{R_1}=\dots = w^{k+1}_{R_n}$ follows from~\eqref{eq:asvrcd_induction_acc},~\eqref{eq:asvrcd_induction_acc_y} together with the update rule (on $\{w^k\}$ and $\{\ww^k\}$) of both algorithms.

To show $\Lcac = \frac{\LcLc}{n} $ it is sufficient to see
{\small
\begin{eqnarray*}
&& \E{\norm{g^k - \nabla f(x^k)}^2_\Popt} \\
&& \qquad =  
\E{  
\norm{\sum_{i\in \sS} p_{i}^{-1}\left( \sum_{j\in R_i} \left( \nabla_{j} f(x^k) - \nabla_{j} f(w^k)\right) e_j \right) - \left(\nabla f(x^k) - \nabla f(w^k)\right)}^2_\Popt
}
\\
&& \qquad =  
\E{  
\norm{\Popt \left(\sum_{i\in \sS} p_{i}^{-1}\left( \sum_{j\in R_i} \left( \nabla_{j} f(x^k) - \nabla_{j} f(w^k)\right) e_j \right)\right) - \Popt\left(\nabla f(x^k) - \nabla f(w^k)\right)}^2
}
\\
&& \qquad =  
\frac1n \E{  
\norm{ \left( \frac1n \sum_{i\in \sS} \pp_{i}^{-1} \left( \nabla \ff_i(\xx^k) - \nabla \ff_i(\ww^k)\right) \right) -  \left(\nabla \ff(\xx^k) - \nabla \ff(\ww^k)\right)}^2
}
\\
&& \qquad =  
\frac1n\E{  
\norm{ \ggggg^k-  \nabla \ff(\ww^k)}^2
}
\\
&& \qquad \leq 
 2\frac{\LcLc}{n} D_\ff(\ww^k,\xx^k)
 \\
&& \qquad =  
 2\frac{\LcLc}{n} \left( \frac1n\sum_{i=1}^n D_{\ff_i}(w^k_{R_i},x^k_{R_i})\right)
  \\
&& \qquad =  
 2\frac{\LcLc}{n} D_{f}(w^k,x^k).
\end{eqnarray*}
}

Lastly, if $x,y \in \range{\Popt}$, there is $\xx, \yy \in \R^\dd$ such that  $x = \Lift{\xx}, y = \Lift{\yy}$. Therefore we can write
\begin{eqnarray*}
f(x) &=& f(\Popt(x)) = \frac1n \sum_{j=1}^n \ff_j(\xx) \\
& \leq&  \frac1n\sum_{j=1}^n \ff_j(\yy) + \< \nabla \left( \frac1n \sum_{j=1}^n \ff_j(\yy)\right), \xx-\yy  >  + \frac{\LL}{2} \| \xx-\yy \|^2
\\
&=&
 f(y)+ \< \nabla f(y), x-y  >  + \frac{\LL}{2n} \| x-y\|^2,
\end{eqnarray*}
 and thus $L = \lambda_{\max} \left( \mM^{\frac12} \Popt \mM^{\frac12}\right)  \leq  n^{-1}\LL$.

 \section{Tighter rates for {\tt GJS} by exploiting prox \label{sec:asvrcd_analysis2} and proof of Theorem~\ref{thm:asvrcd_sega_as}}

In this section, we show that specific nonsmooth function $\psi$ might lead to faster convergence of variance reduced methods. We exploit the well-known fact that under some circumstances, a proximal operator might change the smoothness structure of the objective~\cite{gutman2019condition}. In particular, we consider {\tt GJS} from Chapter~\ref{jacsketch}. We generalize Theorem~\ref{thm:gjs_main} therein, which allows for a tighter rate if $\psi$ has a specific structure.

\begin{theorem}[Extension of Theorem~\ref{thm:gjs_main} from Chapter~\ref{jacsketch}] \label{thm:asvrcd_main2}

Define $f(x) \eqdef \frac1n \sum_{i=1}^n f_i(x)$. Let Assumption~\ref{ass:asvrcd_indicator} hold and suppose that  ${\cM^\dagger}^{\frac{1}{2}}$ commutes with $\cS$. Next, let $\alpha$ and $\cB$ are such that for every $\mX\in \R^{d\times n}$ we have
\begin{equation}\label{eq:asvrcd_small_step_v2}
 \frac{2\alpha}{n^2}  \E{ \norm{ \cU  \mX e }^2_{\Popt} }   +   \NORMG{  \left(\cI - \E{\cS} \right)^{\frac12}\cB  {\cM^\dagger} \mX }  \leq (1-\alpha \sigma') \NORMG{ \cB {\cM^\dagger}\mX },
\end{equation}
\begin{equation}
\frac{2\alpha}{n^2} \E{  \norm{ \cU  \mX  e }_{\Popt}^2  } 
+   \NORMG{\left(\E{\cS}\right)^{\frac12}  \cB  {\cM^\dagger} \mX  }    \leq \frac{1}{n} \norm{{\cM^\dagger}\mX}^2\label{eq:asvrcd_small_step2_v2}
\end{equation}
and $\cB$ commutes with $\cS$. Then  for all $k\geq 0$, we have
$$\E{\Psi^{k}}\leq \left( 1-\alpha\sigma'\right)^k \Psi^{0},$$ where
\begin{eqnarray*}\label{eq:asvrcd_Lyapunov}
\Psi^k & \eqdef & \norm{ x^k - x^* }^2 + \alpha \NORMG{ \cB {\cM^\dagger}^{\frac12} \left( \mJ^k - \mG(x^*)\right)}.
\end{eqnarray*}
\end{theorem}

\subsection{Towards the proof of Theorem~\ref{thm:asvrcd_main2}}

\begin{lemma} (Slight extension of Lemma~\ref{lem:gjs_g_lemma})\label{lem:asvrcd_g_lemma}
Let $\cU$ be random linear operator which is identity in expectation. Let $\mG(x)$ be Jacobian at $x$ and $g^k = \frac1n \cU( \mG(x) -\mJ^k)\eR   -\frac1n\mJ^k\eR$.
Then for any $\mQ \in \R^{d\times d}, \mQ\succeq 0$ and all $k\geq 0$ we have
\begin{equation}\label{eq:asvrcd_g_lemma}
\E{ \norm{ g^k - \nabla f(x^*)}_\mQ^2} \leq    \frac{2}{n^2} \E{  \norm{ \cU \left(\mG(x^k) - \mG(x^*) \right) \eR}_\mQ^2  } + \frac{2}{n^2}  \E{ \norm{ \cU \left(\mJ^k - \mG(x^*) \right) \eR}_\mQ^2}.
\end{equation}
\end{lemma}
\begin{proof}
Since $\nabla f(x^*) = \frac{1}{n}\mG(x^*) \eR$, we have
\begin{equation}\label{eq:asvrcd_nb87fvdbs8s}
g^k-\nabla f(x^*) = \underbrace{\frac1n \cU \left(\mG(x^k)- \mG(x^*) \right) \eR}_{a}  +  \underbrace{\frac1n \left(\mJ^k - \mG(x^*) \right) \eR - \frac1n \cU \left( \mJ^k - \mG(x^*) \right) \eR}_{b} .
\end{equation}
Applying the bound $\norm{a+b}_{\mQ}^2 \leq 2\norm{a}_{\mQ}^2 + 2\norm{b}_{\mQ}^2$ to \eqref{eq:asvrcd_nb87fvdbs8s} and taking expectations, we get
\begin{eqnarray*}
\E{ \norm{ g^k -\nabla f(x^*) }_\mQ^2} &\leq &
 \E{ \frac{2}{n^2} \norm{  \cU \left(\mG(x^k)- \mG(x^*) \right) \eR  }_\mQ^2} \\
 && \qquad + \E{ \frac{2}{n^2} \norm{   \left(\mJ^k-\mG(x^*) \right) \eR -   \cU \left(\mJ^k - \mG(x^*) \right) \eR }_\mQ^2  }\\
&=&
 \frac{2}{n^2} \E{  \norm{ \cU \left(\mG(x^k) - \mG(x^*) \right) \eR }_\mQ^2  } \\
 && \qquad + \frac{2}{n^2} \E{  \norm{ \left(\cI- \cU \right) \left(\mJ^k - \mG(x^*) \right) \eR }_\mQ^2 }.
\end{eqnarray*}

It remains to note that
\begin{eqnarray*}
\E{  \norm{ \left(\cI-  \cU) (\mJ^k - \mG(x^*) \right) \eR }_\mQ^2 }
 &=&
\E{ \norm{ \cU \left(\mJ^k - \mG(x^*) \right) \eR  }_{\mQ}^2 } - \norm{ \left( \mJ^k - \mG(x^*) \right) \eR }_\mQ^2
  \\
  &\leq& 
  \E{ \norm{ \cU \left(\mJ^k - \mG(x^*) \right) \eR }_\mQ^2 } .
\end{eqnarray*}

\end{proof}
 
 Next, we restate two lemmas from the appendix of Chapter~\ref{jacsketch} which we need to show the convergence. 
 \begin{lemma} (Chapter~\ref{jacsketch_appendix}, Lemma~\ref{lem:gjs_smooth2}) \label{lem:asvrcd_smooth_consequence} Assume that function $f_j$ are convex and $\mM_j$-smooth. Then
\begin{equation}\label{eq:asvrcd_smooth}
D_{f_j}(x,y) \geq \frac12 \norm{ \nabla f_j(x)-\nabla f_j(y) }^2_{\mM_j^{\dagger}}, \quad \forall x,y\in \R^d, 1\leq j \leq n.
\end{equation}
If $x-y\in \Null{\mM_j}$, then 
\begin{enumerate} 
\item[(i)]  \begin{equation}\label{eq:asvrcd_linear_on_subspace} f_j(x) = f_j(y) + \langle \nabla f_j(y), x-y\rangle,\end{equation}
\item[(ii)]
\begin{equation} \label{eq:asvrcd_n98g8ff} \nabla f_j(x)-\nabla f_j(y) \in \Null{\mM_j},\end{equation}
\item[(iii)]  
\begin{equation} \label{eq:asvrcd_nb87sgb} \langle \nabla f_j(x) - \nabla f_j(y),x-y\rangle =0.\end{equation}
\end{enumerate}

If, in addition, $f_j$ is bounded below, then $\nabla f_j(x)  \in \range{\mM_j}$ for all $x$.
\end{lemma}

\begin{lemma}(Chapter~\ref{jacsketch_appendix}, Lemma~\ref{lem:gjs_nb98gd8fdx}) \label{lem:asvrcd_nb98gd8fdx}
Let $\cS$ be a random projection operator and $\cA$ any deterministic linear operator commuting with $\cS$, i.e.,  $\cA \cS = \cS \cA$. Further, let $\mX,\mY \in \R^{d\times n}$ and define $\mZ = (\cI-\cS) \mX + \cS \mY$. Then
\begin{itemize}
\item[(i)] $\cA \mZ = (\cI-\cS) \cA \mX + \cS \cA \mY $,
\item[(ii)] $\norm{\cA \mZ}^2 = \norm{(\cI-\cS) \cA \mX}^2  + \norm{\cS \cA \mY}^2 $,
\item[(iii)] $\E{\norm{\cA \mZ}^2} = \norm{(\cI-\E{\cS})^{1/2} \cA \mX}^2  + \norm{\E{\cS}^{1/2} \cA \mY}^2 $, where the expectation is with respect to $\cS$.
\end{itemize}
\end{lemma}

\paragraph{Proof of Theorem~\ref{thm:asvrcd_main2}}
For simplicity of notation, in this proof, all expectations are conditional on $x^k$, i.e., the expectation is taken with respect to the randomness of $g^k$. First notice that
\begin{equation}\label{eq:asvrcd_unbiased_xx}
\E{g^k} = \nabla f(x^k).
\end{equation}

For any differentiable function $h$ let $D_h(x,y)$ to be Bregman distance with kernel $h$, i.e.,  $D_h(x,y)\eqdef h(x)-h(y) - \langle \nabla h(y),x-y\rangle$. Since
\begin{equation}\label{eq:asvrcd_prox_opt_xx}
x^* = \prox_{\alpha \psi}(x^* - \alpha \nabla f(x^*)),
\end{equation}
and since the prox operator is non-expansive, we have
\begin{eqnarray}
\E{\norm{x^{k+1} -x^*}^2 } &\overset{ \eqref{eq:asvrcd_prox_opt_xx}}{=} &
 \E{\norm{\prox_{\alpha \psi}(x^k-\alpha g^k) - \prox_{\alpha \psi}(x^*-\alpha \nabla f(x^*))  }^2}  
 \notag \\
 &\stackrel{\eqref{eq:asvrcd_prox_cont}+\eqref{eq:asvrcd_q_identity}}{\leq} &
 \E{\norm{x^k-  x^* -\alpha \Popt(   g^{k} - \nabla f(x^*) ) }^2}  
 \notag\\
& \overset{\eqref{eq:asvrcd_unbiased_xx}}{=}& 
 \norm{x^k  -x^*}^2 -2\alpha \<  \nabla f(x^k)- \nabla f(x^*) , x^k  -x^*> \notag \\ 
 && \qquad  + \alpha^2\E{\norm{ g^{k}- \nabla f(x^*) }^2_{\Popt}}
 \notag  \\ 
&\leq & 
   (1-\alpha\sigma')\norm{x^k  -x^*}^2 +\alpha^2\E{\norm{  g^{k} -\nabla f(x^*) }^2_{\Popt}}  \notag \\
   && \qquad -2\alpha D_f(x^k,x^*).
 \label{eq:asvrcd_convstepsub1XX_v2}
\end{eqnarray}

Since $f(x)=\frac{1}{n}\sum_{i=1}^n f_i(x)$, in view of \eqref{eq:asvrcd_smooth} we have
\begin{eqnarray}
D_{f}(x^k,x^*) \quad = \quad  \frac{1}{n}\sum_{i=1}^n D_{f_i}(x^k,x^*) & \overset{\eqref{eq:asvrcd_smooth}}{\geq} &\frac{1}{2n} \sum_{i=1}^n \norm{\nabla f_i(x^k) -\nabla f_i(x^*)  }^2_{{\mM_i^{\dagger}} }  \notag\\
&=&  
 \frac{1}{2n}  \left\|{\cM^\dagger}^{\frac12}\left(\mG(x^k) -\mG(x^*) \right) \right\|^2. \label{eq:asvrcd_nb98gd8ff_v2}
\end{eqnarray}

By combining \eqref{eq:asvrcd_convstepsub1XX_v2} and \eqref{eq:asvrcd_nb98gd8ff_v2}, we get
\begin{eqnarray*}
\E{\norm{x^{k+1} -x^*}^2 } & \leq &
   (1-\alpha\sigma')\norm{x^k  -x^*}^2 +\alpha^2\E{\norm{g^{k} -\nabla f(x^*) }_{\Popt}^2}  \\
   && \qquad -\frac{\alpha}{n} \left\|{\cM^\dagger}^{\frac12}( \mG(x^k) -\mG(x^*)) \right\|^2.
\end{eqnarray*}

Next, applying Lemma~\ref{lem:asvrcd_g_lemma} with $\mQ = \Popt$ leads to the estimate
\begin{eqnarray}
\E{\norm{x^{k+1} -x^*}^2 } &\leq &
   (1-\alpha\sigma')\norm{x^k  -x^*}^2 -\frac{\alpha}{n} \norm{ {\cM^\dagger}^{\frac12} \left(\mG(x^k) -\mG(x^*) \right) }^2 \notag 
   \\
   && \qquad  +  \frac{2\alpha^2}{n^2} \E{ \norm{ \cU  \left( \mG(x^k) - \mG(x^*) \right)e }_{\Popt}^2  } \notag \\
   && \qquad + \frac{2\alpha^2}{n^2}  \E{\left\|\cU \left( \mJ^k - \mG(x^*) \right)e \right\|_{\Popt}^2}   \label{eq:asvrcd_48u34719841234_v2} .
\end{eqnarray}

Since, by assumption,  both $\cB$ and ${\cM^\dagger}^{\frac12}$ commute with $\cS$, so does their composition $\cA \eqdef \cB {\cM^\dagger}^{\frac12}$. Applying Lemma~\ref{lem:asvrcd_nb98gd8fdx}, we get
 \begin{eqnarray}\label{eq:asvrcd_J_jac_bound}\nonumber
\E{ \NORMG{\cB {\cM^\dagger}^{\frac12} \left(\mJ^{k+1}-\mG(x^*)  \right) }}  &=&  \NORMG{ (\cI - \E{\cS})^{\frac12}  \cB {\cM^\dagger}^{\frac12} \left(\mJ^k-\mG(x^*) \right) }  \\
&& +  \NORMG{\E{\cS}^{\frac12}  \cB {\cM^\dagger}^{\frac12} \left(\mG(x^k)-\mG(x^*) \right) } .
\end{eqnarray}

Adding $\alpha$-multiple of~\eqref{eq:asvrcd_J_jac_bound} for $\cC = {\cM^\dagger}^{\frac12}$ to~\eqref{eq:asvrcd_48u34719841234_v2} yields
 { 
\footnotesize
\begin{eqnarray*}
&&\E{\norm{x^{k+1} -x^*}^2 } + \alpha\E{ \NORMG{\cB\left({\cM^\dagger}^{\frac12}(\mJ^{k+1}-\mG(x^*) )\right)}}
\\
&& \qquad  \leq 
   (1-\alpha\sigma')\norm{x^k  -x^*}^2 + \frac{2\alpha^2}{n^2} \E{  \norm{ \cU \left(\mG(x^k) - \mG(x^*) \right)e}_{\Popt}^2  } +   \\
   && \qquad \qquad  
 \frac{2\alpha^2}{n^2}  \E{\|\cU(\mJ^k - \mG(x^*))e \|_{\Popt}^2} +   \alpha\NORMG{ (\cI - \E{\cS})^{\frac12} \left(\cB\left({\cM^\dagger}^{\frac12}(\mJ^k-\mG(x^*))\right)\right)} 
 \\
    && \qquad \qquad   + \alpha \NORMG{\E{\cS}^{\frac12}  \left(\cB \left({\cM^\dagger}^{\frac12}(\mG(x^k)-\mG(x^*))\right)\right)}-\frac{\alpha}{n} \left\| {\cM^\dagger}^{\frac12}(\mG(x^k) -\mG(x^*)) \right\|^2
   \\
&& \qquad  \stackrel{\eqref{eq:asvrcd_small_step_v2}}{\leq} 
(1-\alpha\sigma')\norm{x^k  -x^*}^2 + (1-\alpha\sigma')\alpha \E{\NORMG{\cB\left({\cM^\dagger}^{\frac12} \left(\mJ^{k}-\mG(x^*) \right)\right)}} \\
&& \qquad \qquad  
+\frac{2\alpha^2}{n^2} \E{  \| \cU(\mG(x^k) - \mG(x^*))e\|_{\Popt}^2  } +  \alpha \NORMG{\E{\cS}^{\frac12}  \left(\cB \left({\cM^\dagger}^{\frac12}(\mG(x^k)-\mG(x^*))\right)\right)}  \\
    && \qquad \qquad  
 -\frac{\alpha}{n} \left\| {\cM^\dagger}^{\frac12}(\mG(x^k) -\mG(x^*)) \right\|^2
   \\
&& \qquad  \stackrel{\eqref{eq:asvrcd_small_step2_v2}}{\leq} 
(1-\alpha\sigma') \left( \norm{x^k  -x^*}^2 +\alpha\E{ \NORMG{\cB\left({\cM^\dagger}^{\frac12}(\mJ^{k}-\mG(x^*) )\right)}} \right).
\end{eqnarray*}
}
Above, we have used~\eqref{eq:asvrcd_small_step_v2} with $\mX=\mJ^k-\mG(x^*)$ and~\eqref{eq:asvrcd_small_step2_v2} with $\mX= \mG(x^k) -\mG(x^*)$.

 \subsection{Proof of Theorem~\ref{thm:asvrcd_sega_as}}
First, due to our choice of $\cS$ we have
$
\E{\cS (x)} = \diag(p) x
$
and at the same time $\cS$ and ${\cM^\dagger}^{\frac12}$ commute. Next,~\eqref{eq:asvrcd_ESO_sega_good} implies
\begin{equation*}
\E{\left\|\cU(\mM^{\frac12} x)\right\|^{2}_{\Popt}} = \| x \|^2_{\mM^{\frac12} \E{\sum_{i\in S}\pLi^{-1}\eLi \eLi^\top \Popt \sum_{i\in S}\pLi^{-1} \eLi \eLi^\top }\mM^{\frac12} }   \leq  \| x \|^2_{p^{-1}\circ w}.
\end{equation*}

In order to satisfy~\eqref{eq:asvrcd_small_step_v2} and \eqref{eq:asvrcd_small_step2_v2} it remains to have (we substituted $y = {\mM^\dagger}^{\frac12}x$):
\begin{equation} \label{eq:asvrcd_linear_sega_1}
2\alpha\|y \|^2_{p^{-1}\circ w}+ 
\left\|\left(\cI-\E{\cS}\right)^{\frac{1}{2}}\cB(y)\right\|^{2} \leq(1-\alpha \sigma)\|\cB(y)\|^{2},
\end{equation}

\begin{equation} \label{eq:asvrcd_linear_sega}
2\alpha\|y \|^2_{p^{-1}\circ w} + 
\left\|\left(\E{\cS}\right)^{\frac{1}{2}}\cB(y)\right\|^{2} \leq \|y\|^2.
\end{equation}

Let us consider $\cB$ to be the operator corresponding to the left multiplication with matrix $\diag(b)$. Thus for satisfy~\eqref{eq:asvrcd_small_step_v2} it suffices to have for all $i \in [d]$:
\[
2\alpha m_i \pLi^{-1} + b_i^2(1-\pLi) \leq b_i^2 (1-\alpha \sigma) \qquad \Rightarrow \qquad
2\alpha m_i \pLi^{-1} + b_i^2\alpha \sigma \leq b_i^2\pLi .
\]
For~\eqref{eq:asvrcd_small_step2_v2} it suffices to have for all $i\in [d]$
\[
2\alpha m_i \pLi^{-1} + b_i^2 \pLi \leq 1 .
\]
It remains to notice that choice $b_i^2 = \frac{1}{2\pLi}$ and $\alpha  = \min_i \frac{\pLi}{4m_i+ \sigma}$ is valid.


\chapter{Appendix for Chapter \ref{local}}
\label{local_appendix}

\graphicspath{{local/plots/}}

\section{Remaining algorithms}

\subsection{Local {\tt GD} with variance reduction \label{sec:local_lgd}}
In this section, we present variance reduced local gradient descent with partial aggregation. In particular, the proposed algorithm (Algorithm~\ref{alg:local_vr_loc_gd}) incorporates control variates to Algorithm~\ref{alg:local_L2GD}. Therefore, the proposed method can be seen as a special case of Algorithm~\ref{alg:local_L2SGD} with $m=1$. We thus  present it for pedagogical purposes only, as it might shed additional insights into our approach.  In particular, the update rule of proposed method will be $$x^{k+1} =x^k - \alpha g^k,$$ where
\[
g^k \eqdef  \begin{cases}
\pagg^{-1} (\lambda\nabla \Phi(x^k) - \TRloc^{-1}\mJpsi^k) + \TRloc^{-1}\mJf^k + \TRloc^{-1}\mJpsi^k& \text{with probability} \,\, \pagg \\
(1-\pagg)^{-1} (\nabla f(x^k) - \TRloc^{-1}\mJf^k) + \TRloc^{-1}\mJf^k + \TRloc^{-1}\mJpsi^k &\text{with probability} \,\, 1-\pagg 
\end{cases}
\]
for some control variates vectors $\mJf^k, \mJpsi^k \in \R^{nd}$. A quick check  gives 
\[
\E{g^k\, | \, x^k} = \nabla f(x^k) + \lambda \nabla \Phi(x^k) = \nabla F(x^k), 
\]
thus the direction we are taking is unbiased regardless of the value of control variates $ \mJf^k,\mJpsi^k$. 
The goal is to make control variates $\mJf^k,\mJpsi^k$ correlated\footnote{Specifically we aim to have $\Corrloc{ \mJf^k,\TRloc \nabla f(x^k) } \rightarrow 1$ and $\Corrloc{ \TRloc^{-1}\mJpsi^k,\lambda \nabla \Phi(x^k) } \rightarrow 1$ as $x^k \rightarrow x^{*}$.} with $\TRloc \nabla f(x^k)$ and $\TRloc \lambda \nabla \Phi(x^k)$. One possible solution to the problem is for $\mJf^k,\mJpsi^k$ to track most recently observed values of $\TRloc \nabla f(\cdot)$ and $\TRloc\lambda \nabla \Phi(\cdot)$, which corresponds to the following update rule
\[
\left( \mJpsi^{k+1}, \mJf^{k+1}\right)  = 
\begin{cases}
\left(\TRloc\lambda \nabla \Phi(x^k), \mJf^{k} \right)& \text{with probability } \, \pagg \\
\left( \mJpsi^{k}, \TRloc\nabla f(x^k) \right) &\text{with probability } \, 1-\pagg.
\end{cases}
\]

A specific, distributed implementation of the described method is presented as Algorithm~\ref{alg:local_vr_loc_gd}. The only communication between the devices takes place when the average model $\bar{x}^k$ is being computed (with probability $\pagg$), which is analogous to standard local {\tt SGD}. Therefore we aim to set $\pagg$ rather small.

\begin{algorithm}[h]
  \caption{ Variance reduced local gradient descent }
  \label{alg:local_vr_loc_gd}
\begin{algorithmic}
\State{\bfseries Input: }{$x^0_1 = \dots = x_{n}^0\in\R^{d}$,  stepsize $\alpha$, probability $p$  }
\State $\mJf^0_1 = \dots =  \mJf^0_{\TRloc}  = \mJpsi^0_1 = \dots =  \mJpsi^0_{\TRloc} = 0 \in \R^{d}$
  \For{$k=0,1,2,\dotsc$}
  \State $\xi = 1$ with probability $p$ and $0$ with probability $1-p$ 
  \If {$\xi$}
     \State All {\color{red}Devices} $i=1,\dots,n$:      
        \State \hskip .3cm Compute $\nabla f_{\TRloc}(x^k_{\TRloc})$ 
        \State \hskip .3cm $x^{k+1}_\TRloc =  x^k_t - \alpha\left(\TRloc^{-1}(1-\pagg)^{-1} \nabla f_{\TRloc}(x^k_{\TRloc}) - \TRloc^{-1}\frac{\pagg}{1-\pagg} \mJf^k_{\TRloc} + \TRloc^{-1}\mJpsi^k_{\TRloc} \right)$
         		\State \hskip .3cm Set $\mJf^{k+1}_{\TRloc}=  \nabla f_{\TRloc}(x^k_{\TRloc}) $, $\mJpsi^{k+1}_{\TRloc} = \mJpsi^{k}_{\TRloc} $
         		  \Else 
  \State {\color{blue}Master} computes the average $\bar{x}^k = \frac{1}{n}\sum_{i=1}^n x_i^k$
\State {\color{blue}Master} does for all $i=1,\dots,n$:        
 		\State \hskip .3cm Set $x_{\TRloc}^{k+1} =x_{\TRloc}^k  - \alpha \left( \frac{\lambda}{\TRloc\pagg } (x^k_{\TRloc} - \bar{x}^k) -(\pagg^{-1} -1) \TRloc^{-1}\mJpsi^k_{\TRloc} + \TRloc^{-1}\mJf_{\TRloc}^k\right) $ 
 		\State \hskip .3cm Set $\mJpsi^{k+1}_{\TRloc} =  \lambda (x^k_{\TRloc} - \bar{x}^k )$, $\mJf^{k+1}_{\TRloc} = \mJf^{k}_{\TRloc} $
  \EndIf
  \EndFor
\end{algorithmic}
\end{algorithm}

Note that Algorithm~\ref{alg:local_vr_loc_gd} is a particular special case of {\tt SAGA} with importance sampling~\cite{qian2019saga}; thus, we obtain convergence rate of the method for free. We state it as Theorem~\ref{thm:local_basic}.

\begin{theorem}\label{thm:local_basic}
Let Assumption~\ref{as:local_smooth_sc_main} hold. Set $\alpha = \TRloc \min \left\{P\frac{(1-\pagg)}{4L+ \mu}, \frac{\pagg}{4\lambda +  \mu}  \right \}$. Then, iteration complexity of Algorithm~\ref{alg:local_vr_loc_gd} is \[
\max \left\{ \frac{4L +\mu}{ \mu(1-\pagg)}, \frac{4\lambda +\mu}{\mu \pagg} \right \} \log \frac1\varepsilon.\] 
\end{theorem}
\begin{proof}
Clearly, $$ F(x) = f(x) + \lambda\Phi(x) = \frac12 \left( \underbrace{2f(x)}_{\eqdef \floc(x)} +  \underbrace{2\lambda \Phi(x)}_{\eqdef \Philoc(x)}\right).$$
Note that $\Philoc$ is $\frac{2\lambda}{\TRloc}$-smooth and $\floc$ is $\frac{2L}{\TRloc}$-smooth. At the same time, $F$ is $\frac{\mu}{\TRloc}$-strongly convex. Using convergence theorem of {\tt SAGA} with importance sampling from~\cite{qian2019saga, gazagnadou2019optimal}, we get 
\[
\E{F(x^k) + \frac{\alpha}{2}\Upsilon(\mJf^k, \mJpsi^k) } 
\leq 
\left(1- \alpha \frac{\mu}{\TRloc} \right)^k
\left( 
F(x^0) +\frac{\alpha}{2}\Upsilon(\mJf^0, \mJpsi^0)
\right),
\]
where $$\Upsilon(\mJf^k, \mJpsi^k)\eqdef \frac{4}{n^2} \sum_{\TRloc=1}^\TRloc \left( \| \mJpsi^k_\TRloc- \lambda( x_\TRloc^{*} - \bar{x}^{*})\|^2 +\| \mJf_\TRloc^k -   \nabla f_\TRloc(x^{*}_\TRloc) \|^2\right)$$ and \[\alpha =\TRloc \min \left \{ \frac{(1-\pagg)}{4L+ \mu}, \frac{\pagg}{4\lambda +  \mu}  \right \}, \] as desired. 
\end{proof}

\begin{corollary}\label{cor:local_lgd_optimal_p}
Iteration complexity of Algorithm~\ref{alg:local_vr_loc_gd} is minimized for $\pagg  = \frac{4\lambda +\mu}{4\lambda +4L +2\mu}$, which yields complexity $4\left( \frac{\lambda }{\mu} + \frac{L}{\mu} +\frac12\right) \log \frac1\varepsilon$. The communication complexity is minimized for any $\pagg \leq  \frac{4\lambda +\mu}{4\lambda +4L +2\mu}$, in which case the total number of communication rounds to reach $\varepsilon$-solution is $\left(\frac{4\lambda}{\mu}+1 \right)\log\frac{1}{\varepsilon}$. 
\end{corollary}

As a direct consequence of Corollary~\ref{cor:local_lgd_optimal_p} we see that the optimal choice of $p$ that minimizes both communication and number of iterations to reach $\varepsilon$ solution of problem~\eqref{eq:local_main} is $\pagg =  \frac{4\lambda +\mu}{4\lambda +4L +2\mu}$.

\begin{remark}
While both Algorithm~\ref{alg:local_vr_loc_gd} and Algorithm~\ref{alg:local_L2SGD} are a special case of {\tt SAGA}, the practical version of variance reduced local {\tt SGD} (presented in Section~\ref{sec:local_general_full}) is not. In particular, we wish to run the {\tt SVRG}-like method locally in order to avoid storing the full gradient table.\footnote{{\tt SAGA} does not require storing a full gradient table for problems with linear models by memorizing the residuals. However, in full generality, {\tt SVRG}-like methods are preferable.} Therefore, variance reduced local {\tt SGD} that will be proposed in Section~\ref{sec:local_general_full} is neither a special case of {\tt SAGA} nor a special case of {\tt SVRG} (or a variant of {\tt SVRG}). However, it is still a special case of a {\tt GJS} from Chapter~\ref{jacsketch}.
\end{remark}

As mentioned, Algorithm~\ref{alg:local_L2SGD} is a generalization of Algorithm~\ref{alg:local_vr_loc_gd} when the local subproblem is a finite sum. Note that Algorithm~\ref{alg:local_vr_loc_gd} constructs a control variates for both local subproblem and aggregation function $\Phi$ and constructs corresponding unbiased gradient estimator. In contrast, Algorithm~\ref{alg:local_L2SGD} constructs extra control variates within the local subproblem in order to reduce the variance of gradient estimator coming from the local subsampling. 

\subsection{Efficient implementation of {\tt L2SGD+}}
Here we present an efficient implementation of {\tt L2SGD+} as Algorithm~\ref{alg:local_L2SGD_efficient} so that we do not have to communicate control variates. As a consequence, Algorithm~\ref{alg:local_L2SGD_efficient} needs to communicate on average $p(1-p)k$ times per $k$ iterations, while each communication consists of sending only local models to the master and back.

\begin{algorithm}[!h]
  \caption{{\tt L2SGD+}: Loopless Local  {\tt SGD} with Variance Reduction (communication-efficient implementation)}
  \label{alg:local_L2SGD_efficient}
\begin{algorithmic}
\State{\bfseries Input: }{$x^0_1 = \dots = x_{n}^0=\xxloc\in\R^{d}$, stepsize $\alpha$, probability $p$ }
\State Initialize control variates $\mJf^0_{\TRloc}  = 0 \in \R^{d\times m}, \mJpsi^0_{\TRloc}  = 0 \in \R^{d}$ (for $\TRloc = 1,\dots, \TRloc$), initial coin toss $\xi^{-1}=0$
  \For{$k=0,1,2,\dotsc$}
  \State $\xi^k = 1$ with probability $p$ and $0$ with probability $1-p$ 
  \If {$\xi^k = 0$}
\State \hskip -0.3cm All {\color{red}Devices} $i=1,\dots,n$:      
  \If {$\xi^{k-1} = 1$}
  \State Receive $ \xbloc^k_i, c$ from  {\color{blue}Master} 
  \State Reconstruct $\bar{x}^k = \bar{x}^{k-c}$ using  $ \xbloc^k_i, x^{k-c}_i, c$
   \State Set $x^k = {\color{blue} x}^k -c\alpha \frac{1}{\TRloc m} \mJf_{\TRloc}^{k} \onesmloc$,  $ \mJf_{\TRloc}^{k} = \mJf_{\TRloc}^{k-c} $, $ \mJpsi_{\TRloc}^{k} = \lambda(x^{k-c}_i - \bar{x}^k)$, 
  \EndIf
          \State   Sample $j \in \{1,\dots, m\}$ (uniformly at random)
       \State   $g^k_\TRloc = \frac{1}{\TRloc(1-\pagg)} \left(\nabla \flocc_{i,j}(x^k_{\TRloc})- \left(\mJf^k_{\TRloc}\right)_{:,j}\right) + \frac{ \mJf^k_{\TRloc} \onesmloc}{\TRloc m}  +   \frac{\mJpsi^k_{\TRloc}}{\TRloc}  $ 
         \State   $x^{k+1}_\TRloc = x^{k}_\TRloc - \alpha g^k_\TRloc$
        \State    Set $(\mJ^{k+1}_i)_{:,j} = \nabla \flocc_{i,j}(x^k_\TRloc)$, $\mJpsi^{k+1}_{\TRloc}= \mJpsi^k_{\TRloc}$,
        \State \hskip 0.53cm  $(\mJ^{k+1}_i)_{:,l} = (\mJ^{k+1}_i)_{:,l} $ for all $l\neq j$
 \Else 
\State \hskip -0.3cm {\color{blue}Master} does for all $i=1,\dots,n$:       
   \If {$\xi^{k-1} = 0$}
        \State Set $c=0$
       \State Receive $x_i^k$ from {\color{red}Device} and set ${\color{blue} \bar{x}} = \frac{1}{n}\sum_{i=1}^n x_i^k$, $\xbloc^k_i = x_i^k$
   \EndIf
 		\State   Set $\xbloc_{\TRloc}^{k+1} = \xbloc_{\TRloc}^k  - \alpha \left( \frac{\lambda }{\TRloc \pagg} (\xbloc^k_{\TRloc} - {\color{blue} \bar{x}} ) -\frac{\pagg^{-1} -1}{ \TRloc} \lambda (\xxloc -{\color{blue} \bar{x}}  ) \right) $
 		\State Set $\xxloc = \xbloc^{k}_i$
 		\State Set $c=c+1$
  \EndIf
  \EndFor
\end{algorithmic}
\end{algorithm}

\subsection{Local {\tt SGD} with variance reduction -- general method \label{sec:local_general_full}}

In this section, we present a fully general variance reduced local {\tt SGD}. We consider a more general instance of~\eqref{eq:local_main} where each local objective includes a possibly nonsmooth regularizer, which admits a cheap evaluation of proximal operator. In particular, the objective becomes

\begin{equation}\label{eq:local_problem_thispaper}
\compactify  \min_{x\in \R^{d\TRloc}}  F(x)\eqdef  \underbrace{ \frac1\nlocc \sum \limits_{i=1}^{\TRloc} \underbrace{  \left(\sum_{j=1}^{\NRtloc}\flocc_{i,j}(x_{i})\right)}_{= \frac{\nlocc}{\TRloc} f_{\TRloc}(x) } }_{= f(x)}+ \lambda\underbrace{\frac{1}{2\TRloc} \sum_{i=1}^{\TRloc} \| x_{i} - \bar{x}\|^2}_{= \Phi(x)} + \underbrace{\sum \limits_{i=1}^{\TRloc}\psi_{i}(x_{i})}_{\eqdef \psi(x)},
\end{equation}
where $\NRtloc$ is the number of data points owned by client $i$ and $\nlocc = \sum_{i=1}^\TRloc \NRtloc$.

In order to squeeze a faster convergence rate from minibatch samplings, we will assume that $\flocc_{i,j}$ is smooth with respect to a matrix $\mM_{i,j}$ (instead of scalar $\Lloc_{i,j} = \lambda_{\max} \mM_{i,j}$).

\begin{assumption}\label{as:local_smooth_sc}
Suppose that $\flocc_{i,j}$ is $\mM_{i,j}$ smooth ($\mM_{i,j} \in \R^{d\times d}, \mM_{i,j}\succ 0 $) and convex  for $1\leq j\leq \NRtloc,1\leq i \leq \TRloc$, i.e., for all $x,y\in \R^d$ we have
\begin{equation} 
\compactify \flocc_{i,j}(y) + \< \nabla \flocc_{i,j}(y) ,x- y> \leq \flocc_{i,j}(x)\leq    \flocc_{i,j}(y) + \< \nabla \flocc_{i,j}(y) ,x- y> +\frac{1}{2} \norm{y - x}_{{\mM_{i,j}}}^2 .  \label{eq:local_smooth_ass}
\end{equation}
Furthermore, assume that $\psi_i$ is convex for $1\leq i \leq \TRloc$. 
\end{assumption}

Our method (Algotihm~\ref{alg:local_vr_loc2}) allows for arbitrary aggregation probability (same as Algorithms~\ref{alg:local_vr_loc_gd}, \ref{alg:local_L2SGD}), arbitrary sampling of clients (to model the inactive clients) and arbitrary structure/sampling of the local objectives (i.e., arbitrary size of local datasets, arbitrary smoothness structure of each local objective and arbitrary subsampling strategy of each client). Moreover, it allows for the {\tt SVRG}-like update rule of local control variates $\mJ^k$, which requires less storage given an efficient implementation.

To be specific, each device owns a distribution $\cDR_i$ over subsets of $\NRtloc$. When the aggregation is not performed (with probability $1-\pagg$), a subset of active devices $\CS$ is selected ($\CS$ follows arbitrary fixed distribution $\cDR$). Each of the active clients ($i \in \CS$) samples a subset of local  indices $S_i \sim\cDR_{i}$ and observe the corresponding part of local Jacobian $\mG_i(x^k)_{ (:,S_i)}$ (where $\mG_{i}(x^k)\eqdef [\nabla \flocc_{i,1}(x^k), \nabla \flocc_{i,2}(x^k), \dots \nabla \flocc_{i, \NRtloc}(x^k)$). When the aggregation is performed (with probability $\pagg$) we evaluate $\bar{x}^k$ and distribute it to each device; using which each device computes a corresponding component of $\lambda\nabla\Phi(x^k)$. Those are the key components in constructing the unbiased gradient estimator (without control variates).

It remains to construct control variates and unbiased gradient estimator. If the aggregation is done, we just simply replace the last column of the gradient table. If the aggregation is not done, we have two options -- either keep replacing the columns of the Jacobian table (in such case, we obtain a particular case of {\tt SAGA}~\cite{saga}) or do {\tt LSVRG}-like replacement~\cite{hofmann2015variance, kovalev2019don} (in such case, the algorithm is a particular case of {\tt GJS} from Chapter~\ref{jacsketch}, but is not a special case of neither {\tt SAGA} nor {\tt LSVRG}. Note that {\tt LSVRG}-like replacement is preferrable in practice due to a better memory efficiency (one does not need to store the whole gradient table) for the models other than linear.

In order to keep the gradient estimate unbiased, it will be convenient to define vector $\pRloc \in \R^{\NRtloc}$ such that for each $j\in  \{ 1,\dots, \NRtloc\}$ we have $\Probbb{ j\in S_i} =\pRlocjloc$.

Next, to give a tight rate for any given pair of smoothness structure and sampling strategy, we use a rather standard tool called \emph{Expected Separable Overapproximation (ESO)} assumption -- it provides us with smoothness parameters of the objective which ``account'' for the given sampling strategy. 

\begin{assumption}\label{ass:local_eso}
Suppose that there is $v_i\in \R^{\NRtloc}$ such for each client we have:
\begin{equation} \label{eq:local_ESO_saga}
\E{\left\|\sum_{j \in S_i} \mM^{\frac12}_{i,j} h_{i,j}\right\|^{2}} \leq \sum_{j=1}^{\NRtloc} \pRlocjloc v_{i,j}\left\|h_{i,j}\right\|^{2},
\end{equation}
for all $ 1\leq i \leq \TRloc$, $  h_{i,j} \in \R^{\NRtloc}$, and $ j\in  \{ 1,\dots, \NRtloc\}$.
\end{assumption}

 Lastly, denote $\ptg$ to be the probability that worker $i$ is active and $\onestloc \in \R^{\NRtloc}$ to be the vector of ones. 

The resulting algorithm is stated as Algorithm~\ref{alg:local_vr_loc2}.

\begin{algorithm}[!h]
  \caption{{\tt L2SGD++}: Loopless Local {\tt SGD} with Variance Reduction and Partial Participation}
  \label{alg:local_vr_loc2}
\begin{algorithmic}
\State{\bfseries Input: }{$x^0_1, \dots x_{\TRloc}^0\in\R^{d}$, \# parallel units $\TRloc$, each of them owns $\NRtloc$ data points (for $1\leq i\leq \TRloc$), distributions $\cDR_t$ over subsets of $\{ 1, \dots, \NRtloc\}$, distribution $\cDR$ over subsets of $\{1,2,\dots \TRloc\}$, aggregation probability $\pagg$, stepsize $\alpha$ }
\State $\mJf^0_{i}  = 0 \in \R^{d\times \NRtloc},\mJpsi^0_{i}  = 0 \in \R^{d}$ (for $i = 1,\dots, \TRloc$)
  \For{$k=0,1,2,\dotsc$}
  \State $\xi = 1$ with probability $p$ and $0$ with probability $1-p$ 
  \If {$\xi=0$}
 \State Sample $\CS\sim \cDR$
 \State All {\color{red}Devices} $i\in \CS$:      
         \State \hskip .3cm  Sample $ S_i \sim \cDR_i$; $S_i\subseteq \{ 1,\dots, \NRtloc\}$ (independently on each machine)
        \State \hskip .3cm  Observe $\nabla \flocc_{i, j}(x^k_{i})$ for all $j \in S_i $ 
        \State \hskip .3cm $g^k_i = \frac{1}{\nlocc(1-\pagg)\ptg} \left( \sum_{j \in S_i} \pRlocjloc^{-1} \left(\nabla \flocc_{i,j}(x^k_{i})- \left(\mJf^k_{i}\right)_{:,j}\right)  \right) + \frac1\nlocc\mJf^k_{i} \onestloc  + i^{-1}\mJpsi^k_{i}  $
        \State \hskip .3cm $x^{k+1}_i =  \proxt (x^{k}_i - \alpha g^k_i)$
        \State \hskip .3cm For all $j\in \{ 1,\dots, \NRtloc\}$ set {\footnotesize $\mJ^{k+1}_{:,j} =
        \begin{cases}
        \begin{cases}
        \nabla \flocc_{i,j}(x^k_i) & \text{if} \quad  j\in S_i \\
        \mJ^{k}_{:,j} & \text{otherwise}
        \end{cases} & \text{if} \quad \mathrm{{\tt SAGA}} \\
                \begin{cases}
        \nabla \flocc_{i,j}(x^k_i); & \text{w. p.} \quad \psvrg_i \\
        \mJ^{k}_{:,j} & \text{otherwise}
        \end{cases} & \text{if} \quad \mathrm{L-{\tt SVRG}}
         \end{cases} $}
         \State \hskip .3cm Set $\mJpsi^{k+1}_{i}=\mJpsi^k_{i}$
 \State All {\color{red}Devices} $i \not \in \CS$:      
        \State \hskip .3cm $g^k_i = \frac1\nlocc\mJf^k_{i} \onestloc  + i^{-1}\mJpsi^k_{i}  $
        \State \hskip .3cm $x^{k+1}_i =  \proxt (x^{k}_i - \alpha g^k_i)$
         \State \hskip .3cm Set $\mJ^{k+1}_{i} = \mJ^{k}_{i}, \mJpsi^{k+1}_{i}= \mJpsi^k_{i}$
 \Else 
\State {\color{blue}Master} computes the average $\bar{x}^k = \frac{1}{n}\sum_{i=1}^n x_i^k$
\State {\color{blue}Master} does for all $i=1,\dots,n$:            
          \State \hskip .3cm $g^k_i=  \pagg^{-1}  \lambda (x^k_{i} - \bar{x}^k) -(\pagg^{-1} -1) i^{-1}\mJpsi^k_{i} + \frac1\nlocc \mJf_{i}^k \onestloc$
 		\State \hskip .3cm  Set $x_{i}^{k+1} = \proxt\left( x_{i}^k  - \alpha g^k_i \right)$
 		\State \hskip .3cm Set $\mJpsi^{k+1}_{i} =  \lambda (x_{i}^k - \bar{x}^k )$, $\mJf^{k+1}_{i} = \mJf^{k}_{i} $
  \EndIf
  \EndFor
\end{algorithmic}
\end{algorithm}

Next, Theorems~\ref{thm:local_saga_sam} and~\ref{thm:local_svrg_sam} present convergence rate of Algorithm~\ref{alg:local_vr_loc2} ({\tt SAGA} and {\tt SVRG} variant, respectively).

\begin{theorem}\label{thm:local_saga_sam}  Suppose that Assumptions~\ref{as:local_smooth_sc} and~\ref{ass:local_eso} hold. 
Let 
\[\alpha = \min \left \{ \min_{j\in  \{ 1,\dots, \NRtloc\}, 1\leq i\leq \TRloc} \frac{\nlocc (1-\pagg) \pRlocjloc \ptg}{ 4 v_j +\nlocc\frac{\mu}{\TRloc}} , \frac{\TRloc\pagg}{4\lambda+ \mu} \right\}.\] 
Then the iteration complexity of Algorithm~\ref{alg:local_vr_loc2} ({\tt SAGA} option)  is 
\[\max \left\{ \max_{j\in  \{ 1,\dots, \NRtloc\}, 1\leq  i \leq \TRloc} \left(  \frac{4v_j \frac{\TRloc}{\nlocc} + \mu}{\mu(1-\pagg) \pRlocjloc\ptg}\right), \frac{4\lambda +\mu}{\pagg\mu}  \right \}\log\frac1\varepsilon.\] 
\end{theorem}

\begin{theorem}\label{thm:local_svrg_sam} Suppose that Assumptions~\ref{as:local_smooth_sc} and~\ref{ass:local_eso} hold. Let 
\[\alpha = \min \left \{ \min_{j\in  \{ 1,\dots, \NRtloc\}, 1\leq i\leq \TRloc} \frac{\nlocc (1-\pagg)\ptg }{ 4 \frac{v_j}{\pRlocjloc} + \nlocc\frac{\mu}{\TRloc}
 \psvrg_i^{-1}} , \frac{\pagg\TRloc}{4\lambda+ \mu} \right\}.\]
 Then the iteration complexity of Algorithm~\ref{alg:local_vr_loc2} ({\tt LSVRG} option)  is
 \[\max \left\{ \max_{j\in  \{ 1,\dots, \NRtloc\}, 1\leq i\leq \TRloc} \left(  \frac{4v_j\frac{\TRloc}{\nlocc\pRlocjloc}  + \mu
 \psvrg_i^{-1}}{\ptg \mu
 (1-\pagg) }\right) , \frac{4\lambda +\mu}{\pagg\mu}  \right \}\log\frac1\varepsilon.\]
\end{theorem}

\begin{remark}
Algotihm~\ref{alg:local_vr_loc_gd} is a special case of Algorithm~\ref{alg:local_L2SGD} which is in turn special case of Algorithm~\ref{alg:local_vr_loc2}. Similarly, 
Theorem~\ref{alg:local_vr_loc_gd} is a special case of Theorem~\ref{thm:local_saga_simple} which is again special case of Theorem~\ref{thm:local_saga_sam}.
\end{remark}

\subsection{Local stochastic algorithms}

In this section, we present two more algorithms -- Local {\tt SGD} with partial variance reduction (Algorithm~\ref{alg:local_lsgd_partial}) and Local {\tt SGD} without variance reduction (Algorithm~\ref{alg:local_lsgd_none}). While Algorithm~\ref{alg:local_lsgd_none} uses no control variates at all (thus is essentially Algorithm~\ref{alg:local_L2GD} where local gradient descent steps are replaced with local {\tt SGD} steps), Algorithm~\ref{alg:local_lsgd_partial} constructs control variates for $\Phi$ only, resulting in locally drifted {\tt SGD} algorithm (with the constant drift between each consecutive rounds of communication). While we do not present the convergence rates of the methods here, we shall notice they can be easily obtained using the framework from~\cite{sigma_k}.

\begin{algorithm}[!h]
  \caption{Loopless Local {\tt SGD} ({\tt L2SGD}) }
  \label{alg:local_lsgd_none}
\begin{algorithmic}
\State{\bfseries Input: }{$x^0_1 = \dots = x_{n}^0\in\R^{d}$, stepsize $\alpha$, probability $p$ }
  \For{$k=0,1,2,\dotsc$}
  \State $\xi = 1$ with probability $p$ and $0$ with probability $1-p$ 
  \If {$\xi=0$}
 \State All {\color{red}Devices} $i=1,\dots,n$:      
          \State \hskip .3cm Sample $j \in \{1,\dots, m\}$ (uniformly at random)
       \State \hskip .3cm $g^k_i = \frac{1}{\TRloc(1-\pagg)} \left(\nabla \flocc_{i,j}(x^k_{i}) \right)$ 
         \State \hskip .3cm $x^{k+1}_i = x^{k}_i - \alpha g^k_i$
 \Else 
  \State {\color{blue}Master} computes the average $\bar{x}^k = \frac{1}{n}\sum_{i=1}^n x_i^k$
\State {\color{blue}Master} does for all $i=1,\dots,n$:       
   \State \hskip .3cm $g^k_i=   \frac{\lambda }{\TRloc \pagg} (x^k_{i} - \bar{x}^k) $
 		\State  \hskip .3cm Set $x_{i}^{k+1} = x_{i}^k  - \alpha g^k_i $
    \EndIf
      \EndFor
\end{algorithmic}
\end{algorithm}

\begin{algorithm}[!h]
  \caption{Loopless Local {\tt SGD} with partial variance reduction (L2SGD2)}  
  \label{alg:local_lsgd_partial}
\begin{algorithmic}
\State{\bfseries Input: }{$x^0_1 = \dots = x_{n}^0\in\R^{d}$, stepsize $\alpha$, probability $p$ }
\State $\mJpsi^0_{i}  = 0 \in \R^{d}$ (for $i = 1,\dots, \TRloc$)
  \For{$k=0,1,2,\dotsc$}
  \State $\xi = 1$ with probability $p$ and $0$ with probability $1-p$ 
  \If {$\xi=0$}
 \State All {\color{red}Devices} $i=1,\dots,n$:      
          \State \hskip .3cm Sample $j \in \{1,\dots, m\}$ (uniformly at random)
       \State \hskip .3cm $g^k_i = \frac{1}{\TRloc(1-\pagg)} \left(\nabla \flocc_{i,j}(x^k_{i}) \right) +   \frac{1}{\TRloc}  \mJpsi^k_{i}$ 
         \State \hskip .3cm $x^{k+1}_i = x^{k}_i - \alpha g^k_i$
        \State  \hskip .3cm Set  $\mJpsi^{k+1}_{i}= \mJpsi^k_{i}$
 \Else 
  \State {\color{blue}Master} computes the average $\bar{x}^k = \frac{1}{n}\sum_{i=1}^n x_i^k$
\State {\color{blue}Master} does for all $i=1,\dots,n$:       
   \State \hskip .3cm $g^k_i=   \frac{\lambda }{\TRloc \pagg} (x^k_{i} - \bar{x}^k) - \frac{\pagg^{-1} -1}{ \TRloc}\mJpsi^k_{i} $
 		\State  \hskip .3cm Set $x_{i}^{k+1} = x_{i}^k  - \alpha g^k_i $
 		\State  \hskip .3cm Set $\mJpsi^{k+1}_{i} =  \lambda (x_{i}^k - \bar{x}^k )$
    \EndIf
      \EndFor
\end{algorithmic}
\end{algorithm}

\clearpage 
\section{Missing lemmas and proofs}

\subsection{Gradient and Hessian of $\Phi$} \label{sec:local_grad}
\begin{lemma}
Let $\mI$ be the $d\times d$ identity matrix and $\mI_n$ be $n\times n$ identity matrix. Then, we have
\[
\nabla^2 \Phi(x) =  \frac1n\left(\mI_n - \frac{1}{n}ee^\top \right) \otimes \mI 
\qquad 
\mathrm{and} 
\qquad
\nabla \Phi(x) = \frac1n \left( x- \begin{pmatrix} \bar{x}\\ \vdots \\ \bar{x}\\ \bar{x} \\\bar{x} \\ \vdots \\ \bar{x} \end{pmatrix} \right).
\]
Furthermore, $L_{\Phi} = \frac1n$.
\end{lemma}

\begin{proof}

Let $\mO$ the $d\times d$ zero matrix
and let \[\mQ_i \eqdef [\underbrace{\mO, \dots, \mO}_{i-1}, \mI, \underbrace{\mO, \dots, \mO}_{n-i}] \in \R^{d\times dn}\]
and $\mQ \eqdef [\mI, \dots, \mI] \in \R^{d\times dn}$. Note that $x_i= \mQ_i x$, and $\bar{x} = \frac{1}{n} \mQ x$. So,
\[\Phi(x) = \frac{1}{2n}\sum_{i=1}^n \norm{\mQ_ix -\frac{1}{n} \mQ x }^2 = \frac{1}{2n}\sum_{i=1}^n \norm{\left(\mQ_i -\frac{1}{n} \mQ\right) x }^2 .\]
The Hessian of $\Phi$ is
 \begin{eqnarray*}\nabla^2 \Phi(x) &=&
 \frac1n \sum_{i=1}^n \left(\mQ_i -\frac{1}{n} \mQ\right)^\top \left(\mQ_i -\frac{1}{n} \mQ\right)\\
&=&  \frac1n \sum_{i=1}^n \left(\mQ_i^\top \mQ_i  -\frac{1}{n}\mQ_i^\top \mQ  - \frac{1}{n}\mQ^\top \mQ_i +\frac{1}{n^2}\mQ^\top \mQ\right)\\
&=&  \frac1n \sum_{i=1}^n \mQ_i^\top \mQ_i  - \frac{1}{n}\sum_{i=1}^n\frac{1}{n}\mQ_i^\top \mQ  - \frac{1}{n}\sum_{i=1}^n \frac{1}{n}\mQ^\top \mQ_i + \frac{1}{n}\sum_{i=1}^n\frac{1}{n^2}\mQ^\top \mQ\\
&=&  \frac1n\sum_{i=1}^n \mQ_i^\top \mQ_i  - \frac{1}{n^2}\mQ^\top \mQ\\
&=&  \frac1n \begin{pmatrix} \left(1-\frac{1}{n}\right)\mI & -\frac{1}{n} \mI & -\frac{1}{n} \mI   & \cdots & -\frac{1}{n} \mI    \\
 -\frac{1}{n} \mI & \left(1-\frac{1}{n}\right)\mI & -\frac{1}{n} \mI   & \cdots &-\frac{1}{n} \mI    \\
  -\frac{1}{n} \mI & -\frac{1}{n} \mI  & \left(1-\frac{1}{n}\right)\mI   & \cdots & -\frac{1}{n} \mI    \\
   \vdots & \vdots   & \vdots   & &  \vdots     \\
     -\frac{1}{n} \mI & -\frac{1}{n} \mI  &  -\frac{1}{n} \mI  & \cdots &  \left(1-\frac{1}{n}\right)\mI 
 \end{pmatrix}\\
 &=& \frac1n \begin{pmatrix} \left(1-\frac{1}{n}\right)& -\frac{1}{n}  & -\frac{1}{n}  & \cdots & -\frac{1}{n}   \\
 -\frac{1}{n}  & \left(1-\frac{1}{n}\right)& -\frac{1}{n}   & \cdots &-\frac{1}{n}  \\
  -\frac{1}{n} & -\frac{1}{n}  & \left(1-\frac{1}{n}\right)   & \cdots & -\frac{1}{n}    \\
   \vdots & \vdots   & \vdots   & &  \vdots     \\
     -\frac{1}{n} & -\frac{1}{n}  &  -\frac{1}{n} & \cdots &  \left(1-\frac{1}{n}\right)
 \end{pmatrix} \otimes \mI\\
 &=& \frac1n\left(\mI_n - \frac{1}{n}ee^\top \right) \otimes \mI.
 \end{eqnarray*}
 Notice that $\mI_n - \frac{1}{n}ee^\top$ is a circulant matrix, with eigenvalues $1$ (multiplicity $n-1$) and $0$ (multiplicity 1). Since the eigenvalues of a Kronecker product of two matrices are the products of pairs of  eigenvalues of the these matrices, we have
\[\lambda_{\max} (\nabla^2 \Phi(x)) = \lambda_{\max} \left( \frac1n\left(\mI_n - \frac{1}{n}ee^\top \right) \otimes \mI \right)= \frac1n \lambda_{\max}\left(\mI_n - \frac{1}{n}ee^\top \right) =  \frac1n.\]
So, $L_{\Phi}= \frac1n$.

The gradient of $\Phi$ is given by 
 { 
\footnotesize
\begin{eqnarray*}\nabla \Phi(x) &=& 
 \frac1n\sum_{i=1}^n \left(\mQ_i -\frac{1}{n} \mQ\right)^\top \left(\mQ_i -\frac{1}{n} \mQ\right) x \\
&=& \frac1n \sum_{i=1}^n \left(\mQ_i^\top \mQ_i  -\frac{1}{n}\mQ_i^\top \mQ  - \frac{1}{n}\mQ^\top \mQ_i +\frac{1}{n^2}\mQ^\top \mQ\right)x\\
&=& \frac1n\sum_{i=1}^n  \left[\begin{pmatrix} 0\\ \vdots \\ 0\\ x_i \\0 \\ \vdots \\ 0 \end{pmatrix} - \begin{pmatrix} 0\\ \vdots \\ 0\\ \bar{x} \\0 \\ \vdots \\ 0 \end{pmatrix} -  \begin{pmatrix} x_i/n\\ \vdots \\ x_i/n\\ x_i/n \\x_i/n \\ \vdots \\ x_i/n \end{pmatrix} + \begin{pmatrix} \bar{x}/n\\ \vdots \\ \bar{x}/n\\ \bar{x}/n \\\bar{x}/n \\ \vdots \\ \bar{x}/n \end{pmatrix} \right]\\
&=& \frac1n \left( \sum_{i=1}^n  \begin{pmatrix} 0\\ \vdots \\ 0\\ x_i \\0 \\ \vdots \\ 0 \end{pmatrix} - \sum_{i=1}^n \begin{pmatrix} 0\\ \vdots \\ 0\\ \bar{x} \\0 \\ \vdots \\ 0 \end{pmatrix} - \sum_{i=1}^n  \begin{pmatrix} x_i/n\\ \vdots \\ x_i/n\\ x_i/n \\x_i/n \\ \vdots \\ x_i/n \end{pmatrix} + \sum_{i=1}^n \begin{pmatrix} \bar{x}/n\\ \vdots \\ \bar{x}/n\\ \bar{x}/n \\\bar{x}/n \\ \vdots \\ \bar{x}/n \end{pmatrix} \right) \\
&=&  \frac1n \left(x - \begin{pmatrix} \bar{x}\\ \vdots \\ \bar{x}\\ \bar{x} \\\bar{x} \\ \vdots \\ \bar{x} \end{pmatrix} -\begin{pmatrix} \bar{x}\\ \vdots \\ \bar{x}\\ \bar{x} \\\bar{x} \\ \vdots \\ \bar{x} \end{pmatrix}
+\begin{pmatrix} \bar{x}\\ \vdots \\ \bar{x}\\ \bar{x} \\\bar{x} \\ \vdots \\ \bar{x} \end{pmatrix} \right) \\
&=&   \frac1n \left( x- \begin{pmatrix} \bar{x}\\ \vdots \\ \bar{x}\\ \bar{x} \\\bar{x} \\ \vdots \\ \bar{x} \end{pmatrix} \right).
\end{eqnarray*}
}
\end{proof}
\subsection{Proof of Theorem~\ref{thm:local_penalty}}

For any $\lambda, \theta \geq 0$ we have
\begin{eqnarray}  f(x(\lambda)) + \lambda \Phi(x(\lambda)) &\leq &  f(x(\theta)) + \lambda \Phi(x(\theta)) \label{eq:local_first}\\
 f(x(\theta)) + \theta \Phi(x(\theta)) &\leq &  f(x(\lambda)) + \theta \Phi(x(\lambda)).\label{eq:local_second}
\end{eqnarray}
By adding  inequalities \eqref{eq:local_first} and \eqref{eq:local_second}, we get
\[ ( \theta - \lambda)  (\Phi(x(\lambda)) - \Phi(x(\theta))) \geq 0,\]
which means that $\Phi(x(\lambda))$ is decreasing in $\lambda$. Assume $\lambda \geq \theta$. From the \eqref{eq:local_second} we get
\[f(x(\lambda )) \geq f(x(\theta)) + \theta  (\Phi(x(\theta)) - \Phi(x(\lambda))) \geq  f(x(\theta)), \]
where the last inequality follows since $\theta\geq 0$ and since $\Phi(x(\theta)) \geq \Phi(x(\lambda)) $. So, $f(x(\lambda))$ is increasing.

Notice that since $\Phi$ is a non-negative function and since $x(\lambda)$ minimizes $F$ and $\Phi(x(\infty))=0$, we have
\[ f(x(0)) \leq f(x(\lambda)) \leq f(x(\lambda)) + \lambda \Phi(x(\lambda)) \leq f(x(\infty)) ,\]
which implies \eqref{eq:local_dissimilarity} and \eqref{eq:local_bifg9dd8}.

\subsection{Proof of Theorem~\ref{thm:local_characterization}}
The equation $\nabla F(x^{*}(\lambda))=0$ can be equivalently written as
\[\nabla f_i(x_i^{*}(\lambda)) +\lambda (x_i^{*}(\lambda) - \overline{x}^{*}(\lambda)) = 0, \qquad i=1,2,\dots,n,\]
which is identical to~\eqref{eq:local_step_from_mean}. Averaging these identities over $i$, we get
$ \overline{x}^{*}(\lambda) =  \overline{x}^{*}(\lambda) - \frac{1}{\lambda } \frac{1}{n}\sum_{i=1}^n \nabla f_i(x^{*}_i(\lambda))$, which implies
\[
\sum_{i=1}^n \nabla f_i(x^{*}_i(\lambda)) = 0.
\] 
Further, we have 
\[\Phi(x^{*}(\lambda)) = 
\frac{1}{2n}\sum_{i=1}^n \norm{x^{*}_i(\lambda)-\overline{x}^{*}(\lambda)}^2  = \frac{1}{2 n\lambda^2 } \sum_{i=1}^n \norm{\nabla f_i(x^{*}_i(\lambda))}^2 = \frac{1}{2 \lambda^2  } \norm{\nabla f(x^{*}(\lambda))}^2,\]
 as desired.

\subsection{Proof of Lemma~\ref{lem:local_exp_smoothnes-basic}}

We first have
\begin{eqnarray*}
\E{ \norm{g(x) - G(x^*)}^2} &=& (1-p) \norm{ \frac{\nabla f(x)}{1-p} - \frac{\nabla f(x^*)}{1-p}}^2 + p\norm{ \lambda \frac{\nabla \Phi(x)}{p} -  \lambda \frac{\nabla \Phi(x^*)}{p}}^2\\
&=& \frac{1}{1-p}  \norm{ \nabla f(x) - \nabla f(x^*)}^2 + \frac{\lambda^2}{p} \norm{\nabla \Phi(x)-\nabla \Phi(x^*)}^2 \\
&\leq & \frac{2L_f}{1-p} D_f(x, x^*) + \frac{2\lambda^2 L_{\Phi}}{p} D_{\Phi}(x, x^*)
\\
&= & \frac{2L}{n(1-p)} D_f(x, x^*) + \frac{2\lambda^2}{np} D_{\Phi}(x, x^*).
\end{eqnarray*}
Since $D_f + \lambda D_{\Phi} = D_F$ and $\nabla F(x^*)=0$, we can continue:
\begin{eqnarray*}
\E{ \norm{g(x) - G(x^*)}^2}&\leq & \frac{2}{n}\max\left\{\frac{L}{1-p}, \frac{\lambda }{p}\right\} D_F(x, x^*) \\
&=& \frac{2}{n}\max\left\{\frac{L}{1-p}, \frac{\lambda }{p}\right\} \left(F(x)-F(x^*)\right) .
\end{eqnarray*}

Next, note that 
 \begin{align}
\sigma^2  &= \frac{1}{n^2}  \sum_{i=1}^n \left( \frac{1}{1-p} \| \nabla f_i(x^{*}_i)\|^2 + \frac{\lambda^2}{p} \|x_i^{*} - \overline{x}^{*} \|^2 \right) 
\nonumber\\
&=  \frac{1}{1-p} \norm{\nabla f(x^*)}^2 + \frac{\lambda^2}{p}\norm{\nabla \Phi(x^*)}^2 \nonumber\\
&= (1-p) \norm{ \frac{\nabla f(x^*)}{1-p} }^2 + p\norm{ \frac{\lambda  \nabla \Phi(x^*)}{p}) }^2 \nonumber\\
& = \E{\norm{G(x^{*})}^2}. \label{eq:local_sigma2}
 \end{align}
 Therefore, we have
 { 
\footnotesize
 \begin{eqnarray*}
 \E{\norm{g(x)}^2}& \leq &\E{ \norm{g(x) - G(x^*)}^2}+ 2\E{\norm{G(x^{*})}^2} \\
 & \stackrel{\text{Lemma}~\ref{lem:local_exp_smoothnes-basic}+\eqref{eq:local_sigma2}}{\leq} & 4 \cL (F(x)-F(x^*))  + 2\sigma^2
 \end{eqnarray*}
 }
 as desired. 
 
 \subsection{Proof of Theorem~\ref{thm:local_local_gd}}
 
 Follows from Lemma~\ref{lem:local_exp_smoothnes-basic} by applying Theorem 3.1 from~\cite{pmlr-v97-qian19b}.

\subsection{Proof of Corollary~\ref{cor:local_optimalp}}
Firstly, to minimize the total number of iterations, it suffices to minimize $\cL$ which is achieved with $p^{*}=\frac{\lambda}{L + \lambda}$. Let us look at the communication.
Fix $\varepsilon>0$, choose $\alpha = \frac{1}{2\cL}$ and let $k = \frac{2n\cL}{ \mu} \log \frac{1}{\varepsilon}$, so that
\[\left(1- \frac{\mu}{2n\cL}\right)^k \leq \varepsilon.\]
The expected number of communications to achieve this goal is equal to
\begin{eqnarray*} {\rm Comm}_p & \eqdef & p(1-p)k \\
&=& p(1-p)\frac{2\max\left\{\frac{L}{1-p}, \frac{\lambda }{p}\right\}}{ \mu}\log \frac{1}{\varepsilon} \\
&=& \frac{2\max\left\{ p L, (1-p)\lambda \right\}}{\mu} \log \frac{1}{\varepsilon}.
\end{eqnarray*}

The quantity ${\rm Comm}_p$ is minimized by choosing any $p$ such that $pL = (1-p)\lambda$, i.e., for $p =    \frac{\lambda}{\lambda + L} =p^{*}$,  as desired. The optimal expected number of communications is therefore equal to
\[ {\rm Comm}_{p^{*}} = \frac{2 \lambda }{ \lambda + L} \frac{L}{\mu}\log \frac{1}{\varepsilon}.\]

\subsection{Proof of Corollary~\ref{cor:local_lsd_optimal_p}}
Firstly, to minimize the total number of iterations, it suffices to solve
\[
\min_p \max \left\{   \frac{4\Lloc  + \mu m}{ (1-\pagg)\mu}, \frac{4\lambda +\mu}{\pagg\mu}\right \},
\] which is achieved with $p = p^{*}=\frac{4\lambda +\mu}{4\Lloc +4\lambda +(m+1)\mu}$.  The expected number of communications to reach $\varepsilon$-solution is
\begin{eqnarray*} {\rm Comm}_p & = & p(1-p)\max \left\{   \frac{4\Lloc  + \mu m}{ (1-\pagg)\mu}, \frac{4\lambda +\mu}{\pagg\mu}\right \} \log \frac{1}{\varepsilon} \\
&=& \frac{\max\left\{ p (4\Lloc + \mu m), (1-p)(4\lambda + \mu) \right\}}{\mu} \log \frac{1}{\varepsilon}.
\end{eqnarray*}

Minimizing the above in $p$  yield $p = p^{*}=\frac{4\lambda +\mu}{4\Lloc +4\lambda +(m+1)\mu}$,  as desired. The optimal expected number of communications is therefore equal to
\[ {\rm Comm}_{p^{*}} = \frac{4 \lambda + \mu }{4\Lloc + 4 \lambda + (m+1)\mu} \left( 4\frac{\Lloc}{\mu}+m\right) \log \frac{1}{\varepsilon}.\]

\subsection{Proof of Theorems~\ref{thm:local_saga_simple},~\ref{thm:local_saga_sam}, and~\ref{thm:local_svrg_sam}}
Note first that Algorithm~\ref{alg:local_L2SGD} is a special case of Algorithm~\ref{alg:local_vr_loc2}, and Theorem~\ref{thm:local_saga_simple} immediately follows from Theorem~\ref{thm:local_saga_sam}. Therefore it suffices to show Theorems~\ref{thm:local_saga_sam}, and~\ref{thm:local_svrg_sam}. In order to do so, we will cast Algorithm~\ref{alg:local_vr_loc2} as a special case of {\tt GJS} (Algorithm~\ref{alg:gjs_SketchJac}). As a consequence, Theorem~\ref{thm:local_saga_sam} will be a special cases of Theorem~\ref{thm:gjs_main}.

\subsubsection{Variance reduced local {\tt SGD} as special case of {\tt GJS}}
Let $\maploc(i,j ) \eqdef j+\sum_{l=1}^{i-1} \NRtloc$
In order to case problem~\eqref{eq:local_problem_thispaper} as~\eqref{eq:gjs_problem_gjs}, denote $\nloc \eqdef \nlocc+1$, $\floc_{\maploc(i,j)}(x) \eqdef \frac{\nlocc+1}{\nlocc} \flocc_{i,j}(x_i)$ and $\floc_{\nloc} \eqdef (\nlocc+1)\Phi $. Therefore the objective ~\eqref{eq:local_problem_thispaper} becomes
\begin{equation}
\min_{x\in \R^{\nlocc d}} \Ug(x) \eqdef\frac{1}{\nloc} \sum_{j=1}^{\nloc} \floc_j(x) + \psi(x).
\end{equation}
Let $\vg \in \R^{\nloc-1}$ be such that $\vg_{\maploc(i,j)} = \frac{\nlocc+1}{\nlocc} v_{i,j}$
and as a consequence of~\eqref{eq:local_ESO_saga} we have 
\begin{equation}\label{eq:local_ESO2}
\E{\left\|\sum_{j \in S_i} \mM^{\frac12}_{i,j} h_{i,j}\right\|^{2}} \leq \sum_{j=1}^{\NRtloc} \pRlocjloc \vg_{\maploc(i,j)}\left\|h_{i,j}\right\|^{2}, \quad \forall \; 1\leq i \leq \TRloc, \, \forall  h_{i,j} \in \R^{d}, j\in \{ 1, \dots, \NRtloc\}
\end{equation} 
At the same time, $\Ug$ is $\mug \eqdef \frac{\mu}{\TRloc}$ strongly convex.  

\subsubsection{Proof of Theorem~\ref{thm:local_saga_sam} and Theorem~\ref{thm:local_svrg_sam}}
Let $\ones \in \R^d$ be  a vector of ones and $\ptRloc\in \R^{\nlocc}$ is such that $\ptRloc_j = \pRlocjloc$ if $j\in \{1, \dots, \NRtloc\}$, otherwise $\ptRloc_j=0$.
 Given the notation, random operator $\cU$ is chosen as 
\[
\cU \mX = 
  \begin{cases} 
 (1-\pagg)^{-1} \sum_{i=1}^\TRloc   \left( \ptg^{-1} \ones \left(\left(\ptRloc\right)^{-1} \right)^\top\right) \circ \left( \mX_{:\NRtloc} \left( \sum_{j\in S_i} \eRj \eRj^\top \right)\right) 
& \text{w.p.}\quad (1- \pagg)  \\
\pagg^{-1} \mX_{:,\nloc}  & \text{w.p.}\quad \pagg 
\end{cases} 
\] 

We next give two options on how to update Jacobian -- first one is {\tt SAGA}-like, second one is {\tt SVRG} like. 
\begin{eqnarray*}
\text{{\tt SAGA}-like:} & \quad 
(\cS \mX)_{:,\NRtloc} =& 
  \begin{cases}
\mX_{:,S_i}  =   \mX_{:\NRtloc} \left( \sum_{j\in S_i} \eRj \eRj^\top \right),  & \text{w.p.}\quad (1- \pagg) \ptg, \\
0 & \text{w.p.}\quad (1- \pagg) (1-\ptg) + \pagg  \\
\end{cases}  
\\
&
\quad 
(\cS \mX)_{:,\nloc} =& 
  \begin{cases}
 \mX_{:,\nloc}             &  \text{w.p.}\quad  \pagg \\
0             &  \text{w.p.}\quad 1- \pagg 
\end{cases} 
 \\
\text{{\tt SVRG}-like:}& \quad 
(\cS \mX)_{:,\NRtloc}  =&
  \begin{cases}
   \mX_{:\NRtloc} b_i ; \; b_i = \begin{cases}
1 &  \text{w.p.}\quad  \psvrg_i \\
0 &  \text{w.p.}\quad  1- \psvrg_i
\end{cases}  
& \text{w.p.}\quad  (1- \pagg)\ptg \\
0            &  \text{w.p.}\quad  (1- \pagg) (1-\ptg) + \pagg \;
\end{cases}  
 \\
& \quad 
(\cS \mX)_{:,\nloc}  =&
  \begin{cases}
  \mX_{:,\nloc} & \text{w.p.}\quad  \pagg \\
0             &  \text{w.p.}\quad  1-\pagg \;.
\end{cases}  
\end{eqnarray*}

We can now proceed with the proof of Theorem~\ref{thm:local_saga_sam} and Theorem~\ref{thm:local_svrg_sam}. As $\nabla f_i(x) - \nabla f_i(y)\in \Range{\mM_i}$, we must have
\begin{equation}\label{eq:local_g_in_range}
\mG(x^k) - \mG(x^{*}) = \cM^{\dagger} \cM \left(\mG(x^k) - \mG(x^{*})\right)
\end{equation}
 and 
 \begin{equation}\label{eq:local_j_in_range}
\mJ^k - \mG(x^{*}) = \cM^{\dagger} \cM \left(\mJ^k - \mG(x^{*})\right).
\end{equation}
Due to~\eqref{eq:local_j_in_range},~\eqref{eq:local_g_in_range}, inequalities~\eqref{eq:gjs_small_step} and \eqref{eq:gjs_small_step2} with choice $\mY=  {\cM^\dagger}^{\frac12}\mX$ become respectively:
{ 
\footnotesize
\begin{eqnarray}
\nonumber
 && \frac{2\alpha}{{\nloc}^2}\pagg^{-1} \|\mM_{{\nloc}}^{\frac12}\mY_{:,{\nloc}}\|^2 +  \frac{2\alpha^2}{{\nloc}^2}(1-\pagg)^{-1}\sum_{i=1}^\TRloc\E{\left \|  
 \ptg^{-1}\sum_{j\in S_i} \pRlocjloc^{-1}
\mM_{i,j}^{\frac12} \mY_{:j}\right \|^2}   + 
\left\|\left(\cI-\E{\cS}\right)^{\frac{1}{2}}\cB(\mY)\right\|^{2}  \\
&& \qquad \qquad \qquad\leq(1-\alpha \mug)\|\cB(\mY)\|^{2}
\label{eq:local_linear_1_svrg}
\end{eqnarray}
}
{ 
\footnotesize
\begin{equation}\label{eq:local_linear_2_svrg}
 \frac{2\alpha}{{\nloc}^2} \pagg^{-1} \|\mM_{{\nloc}}^{\frac12}\mY_{:,{\nloc}}\|^2 +  \frac{2\alpha^2}{{\nloc}^2}(1-\pagg)^{-1}\sum_{i=1}^\TRloc\E{\left \|  
 \ptg^{-1}\sum_{j\in S_i} \pRlocjloc^{-1}
\mM_{i,j}^{\frac12} \mY_{:j}\right \|^2}  + 
\left\|\left(\E{\cS}\right)^{\frac{1}{2}}\cB(\mY)\right\|^{2} \leq \frac1\nloc  \|\mY \|^2
\end{equation}
}

Above, we have used 
{ 
\footnotesize
\[
\E{\|\cU\mX e\|^2} = \E{\|\cU{\cM}^{\frac12} \mY e\|^2} =  \pagg^{-1} \|\mM_{n}^{\frac12}\mY_{:,n}\|^2 + (1-\pagg)^{-1}\sum_{i=1}^\TRloc\E{\left \|  
 \ptg^{-1}\sum_{j\in S_i} \pRlocjloc^{-1}
\mM_{i,j}^{\frac12} \mY_{:j}\right \|^2} .
\]
}
Note that $\E{\cS (\mX)} = \mX \cdot \diag\left((1- \pagg)({ \color{blue} p} \circ \pg), \pagg\right) $ where ${ \color{blue} p}\in \R^{\nloc-1}$ such that ${ \color{blue} p}_{\maploc(i,j)} = \pRlocjloc$. Using~\eqref{eq:local_ESO2}, setting $\cB$ to be right multiplication with $\diag(b)$ and noticing that $\lambda_{\max} \mM_{\nloc} =   \nloc\lambda$ it suffices to have 

\[
 \frac{2\alpha}{\nloc}\pagg^{-1} \lambda+ (1- \pagg)b_{\nloc}^2 \leq (1- \alpha \mug) b_{\nloc}^2
\]
\[
 \frac{2\alpha}{{\nloc}^2}(1-\pagg)^{-1} \pRlocjloc^{-1} \ptg^{-1} \vg_{\maploc(i,j)}+ (1-  (1- \pagg)\pRlocjloc \ptg)b_j^2 \leq (1- \alpha \mug) b_j^2 \qquad  \forall j\in \{1,\dots,  \NRtloc\}, i \leq \TRloc
\]

\[
 \frac{2\alpha}{\nloc}\pagg^{-1}\lambda+  \pagg b_{\nloc}^2 \leq  \frac1{\nloc}
\]
\[
 \frac{2\alpha}{{\nloc}^2} (1-\pagg)^{-1} \pRlocjloc^{-1}  \ptg^{-1} \vg_{\maploc(i,j)}+ (1- \pagg)\pRlocjloc\ptg b_j^2  \leq  \frac1\nloc  \qquad \forall j\in \{1,\dots,  \NRtloc\}, i \leq \TRloc
\]
for {\tt SAGA} case and 
\[
 \frac{2\alpha}{\nloc}\pagg^{-1} \lambda+ (1- \pagg)b_{\nloc}^2 \leq (1- \alpha \mug) b_{\nloc}^2
\]
\[
 \frac{2\alpha}{{\nloc}^2}(1-\pagg)^{-1} \pRlocjloc^{-1}  \ptg^{-1} \vg_{\maploc(i,j)}+ (1-  (1- \pagg)\psvrg_i \ptg)b_j^2 \leq (1- \alpha \mug) b_j^2 \qquad \forall j\in \{1,\dots,  \NRtloc\}, i \leq \TRloc
\]

\[
 \frac{2\alpha}{\nloc}\pagg^{-1}\lambda+  \pagg b_{\nloc}^2 \leq  \frac1{\nloc}
\]
\[
 \frac{2\alpha}{{\nloc}^2} (1-\pagg)^{-1} \pRlocjloc^{-1}  \ptg^{-1} \vg_{\maploc(i,j)}+ (1- \pagg)\psvrg_i\ptg b_j^2  \leq  \frac1\nloc  \qquad  \forall j\in \{1,\dots,  \NRtloc\}, i \leq \TRloc
\]
for {\tt LSVRG} case. 

It remains to notice that to satisfy the {\tt SAGA} case, it suffices to set $b_{\nloc}^2 = \frac{1}{2{\nloc} \pagg}, b_{\maploc(i,j)}^2 = \frac1{2{\nloc} (1-\pagg)\pRlocjloc\ptg }$ (for $j\in \{1, \dots, \NRtloc\}, i\leq \TRloc$) and $\alpha = \min \left \{ \min_{j\in \{1,\dots,  \NRtloc\}, 1\leq i\leq \TRloc} \frac{{\nloc} (1-\pagg) \pRlocjloc\ptg}{ 4 \vg_{\maploc(i,j)} + {\nloc}\mug} , \frac{\pagg}{4\lambda+ \mug} \right\}$.

To satisfy {\tt LSVRG} case, it remains to set $b_{\nloc}^2 = \frac{1}{2{\nloc} \pagg}, b_{\maploc(i,j)}^2 = \frac1{2{\nloc} (1-\pagg)\psvrg_i \ptg}$ (for $j\in \{1, \dots, \NRtloc\}, i\leq \TRloc$) and $\alpha = \min \left \{ \min_{j\in \{1,\dots,  \NRtloc\}, 1\leq i\leq \TRloc} \frac{{\nloc} (1-\pagg)\ptg }{ 4 \frac{\vg_{\maploc(i,j)}}{\pRlocjloc} + {\nloc}\mug \psvrg_i^{-1}} , \frac{\pagg}{4\lambda+ \mug} \right\}$.

The last step to establish is to recall that $\nloc = \nlocc+1,\vg_{\maploc(i,j)} = \frac{\nlocc+1}{\nlocc} v_{i,j}$ and $\mug = \frac{\mu}{\TRloc}$ and note that the iteration complexity is $\frac{1}{\alpha \mug} \log\frac{1}{\varepsilon} = \frac{\TRloc}{\alpha \mu} \log\frac{1}{\varepsilon} $.

\subsubsection{Proof of Theorem~\ref{thm:local_saga_simple}}
To obtain convergence rate of Theorem~\ref{thm:local_saga_simple}, it remains to use Theorem~\ref{thm:local_saga_sam} with $\ptg=1, \NRtloc =m$ ($\forall i\leq \TRloc$), where each machine samples (when the aggregation is not performed) individual data points with probability $\frac1m$ and thus $p_j=\frac1m$ (for all $j\leq \nlocc$). The last remaining thing is to realize that $v_j =\Lloc$ for all $j\leq \nlocc$.


\chapter{Appendix for Chapter \ref{sscn}}
\label{sscn_appendix}

\graphicspath{{sscn/plots/}}

\section{Missing lemmas and proofs from Section~\ref{sec:sscn_preliminaries}}

\subsection{Explicit update}
\begin{lemma}\label{lem:sscn_explicit_update}
Let $x^+ = \argmin_y \langle g', y-x \rangle + \frac{H'}{2} \|x-y\|^2+ \frac{M'}{6}\|x-y \|^3$, where $H',M'  >0$.
Then we have
\begin{equation}\label{eq:sscn_x_update_implicit}
    x^+= x- \frac{2g'}{ H' + \sqrt{{H'}^2 + 2M'\|g'\|}}
\end{equation}
\end{lemma}
\begin{proof}
By first-order optimality conditions we have $g'+ H'(x^+-x) + \frac{M'}{2}\|x^+-x\|(x^+-x) = 0$ which immediately yields
\begin{equation} \label{eq:sscn_lpdasadissio}
 x^+ = x - \frac{g'}{H' + \frac{M'}{2}\|x^+-x\|}. 
\end{equation}
Rearranging the terms and taking the norm we have $\frac{M'}{2}\|x^+-x \|^2 + H'\|x^+-x \| + \|g'\| = 0$. Solving the quadratic equation we arrive at
\[
\|x^+-x \|  = \frac{\sqrt{{H'}^2 + 2M'\|g'\|}-H'}{M'}.
\]
Plugging it back to~\eqref{eq:sscn_lpdasadissio}, we get~\eqref{eq:sscn_x_update_implicit}.
\end{proof}

\subsection{Proof of Lemma~\ref{lem:sscn_ub}}
First, note that 
{ 
\begin{eqnarray*}
&& D_f(x^+,x) - \frac{1}{2} (x^+-x)^\top \nabla^2 f(x)(x^+-x) \\
&&  = 
\int_{0}^{1} \langle \nabla f(x+t(x^+-x)) - f(x), x^+-x \rangle \,dt - \frac{1}{2} (x^+-x)^\top \nabla^2 f(x)(x^+-x)
\\
&& = 
\int_{0}^{1} \int_{0}^{1} \langle t \nabla^2 f(x+st(x^+-x)),x^+-x, x^+-x \rangle \, ds\,dt- \frac{1}{2} (x^+-x)^\top \nabla^2 f(x)(x^+-x)
\\
&& = 
\int_{0}^{1} \int_{0}^{1} \langle t \nabla^2 f (x+st(x^+-x)) - \nabla^2 f(x),x^+-x,x^+-x \rangle \, ds\,dt 
\\ 
&& = 
\int_{0}^{1} \int_{0}^{1}  \int_{0}^{1}  \langle t^2 s \nabla^3 f (x+rst(x^+-x)),x^+-x,x^+-x,x^+-x \rangle \, dr\, ds\,dt.
\end{eqnarray*}
}
Using~\eqref{eq:sscn_update_general} we get
\begin{eqnarray*}
&& |f(x^+) - f(x) + \langle \nabla f(x),\mS h \rangle + \frac12 h^\top   \nabla^2_{\mS}  f(x) h | 
\\
& \stackrel{\eqref{eq:sscn_update_general} }{=}&   \left| \int_{0}^{1} \int_{0}^{1}  \int_{0}^{1}  \langle t^2 s \nabla^3 f (x+rst\mS h),\mS h,\mS h, \mS h \rangle \, dr\, ds\,dt \right| \\
& \refLE{eq:sscn_MS_def} &  
\int_{0}^{1} \int_{0}^{1}  \int_{0}^{1}   t^2 s M_{\mS} \|h_{\mS}\|^3 \, dr\, ds\,dt \\
&= &
 \frac{M_{\mS}}{6} \| h_{\mS}\|^3.
\end{eqnarray*}

\subsection{Proof of Lemma~\ref{lem:sscn_sharpness}}

First, $M \geq M_{\mS} $ is trivial. At the same time $M=M_{\mS}$ if $\nabla^3 f(x)$ is identity tensor always (which is clearly feasible) -- thus the inequality is tight. 

To show sharpness of $M_{\mS} \geq \left(\frac{\tau}{d}\right)^{\frac32} M$, consider $f(x) = \frac16 (x^\top e)^3$. In this case, we have\footnote{By $[e] \in \R^{d\times d\times d}$ we mean third order product of vector $e$.} $\nabla^3 f(x) = [e]^3$  and $\mS=e_i$. In such case, $M = d^\frac32$ and $M_{\mS} = \tau^\frac32$. Note that $f$ is non-convex in this example. However, $f$ is convex on a set where $xe\geq0$, hwere the argument follows through.

\section{Proofs for Section~\ref{sec:sscn_global}}

\subsection{Proof of Lemma~\ref{lem:sscn_exp_size}\label{sec:sscn_exp_size_proof}}
Let $\Tr{\mA}$ be a trace of square matrix $\mA$. We have
\begin{eqnarray*}
\E{\tau(\mS)} &=& \E{ \Tr{\mI^{\tau(\mS)}}} = \E{\Tr{ \mS^\top \mS \left(\mS^\top \mS\right)^{-1} } } = \E{\Tr{ \mS \left(\mS^\top \mS\right)^{-1}  \mS^\top }}\\
&=&
 \Tr{\E{ \mS \left(\mS^\top \mS\right)^{-1}  \mS^\top }} \stackrel{\eqref{eq:sscn_uniform_sampling}}{=} \Tr{\frac{\tau}{d} \mI^d}  = \tau.
\end{eqnarray*}

\subsection{Proof of Lemma~\ref{lem:sscn_keylemma}}

For any $h' \in \R^d$ denote 
$$
\Omega_{\mS}(x; h')  \Def  f(x) + \la \nabla f(x), \mZ h'\ra
+ \frac{1}{2}\la \nabla^2 f(x)\mZ h', \mZ h' \ra
+ \frac{H}{6}\|\mZ h' \|^3 + \psi(x +\mZ h')
$$
and 
$$
 T_{\mS}(x^k) \Def \argmin_{h' \in \R^d} \Omega_{\mS}(x; h').
$$ 

Then, for any fixed $y \in \R^d$ we have
$$
F(x^{k + 1})  \refLE{eq:sscn_coordinate_ub} 
\Omega_{\mS}(x^k; T_{\mS}(x^k)) \; \leq \; \Omega_{\mS}(x^k; y - x^k).
$$
Therefore,
\begin{eqnarray*}
\E{ F(x^{k + 1}) \, | \, x^k}& \leq & \E{ \Omega_{\mS}(x^k; y - x^k) } \\
\\
& = & f(x^k) + \frac{\tau}{d} \la \nabla f(x^k), y - x^k \ra
+ \E{  \frac{1}{2}\la \mZ \nabla^2 f(x^k)\mZ (y - x^k), y - x^k \ra }\\
\\
& \; & \quad + \quad
\frac{M}{6}\E{ \| \mZ (y - x^k) \|^3 }
+ \E{\psi(x +\mZ (y-x^k))} .
\end{eqnarray*}
Let us get rid of the expectations above. Firstly, we have 
\begin{eqnarray*}
\E{\psi(x +\mZ (y-x^k))} 
&=&
 \E{\left \langle \psi'\left( \left( \mI^d - \mZ \right)x^k +\mZ y\right) ,e  \right \rangle } \\
 &= &
  \E{\left \langle \left( \mI^d - \mZ \right) \psi'\left( x^k\right) ,e  \right \rangle } + \E{\left \langle \mZ  \psi'\left( y\right) ,e  \right \rangle } 
  \\
  &=&
 \left(1 - \frac{\tau}{d}\right)\psi(x^k) + \frac{\tau}{d} \psi(y) .
\end{eqnarray*}
For the cubed norm it can be estimated as follows
$$
\E{ \| \mZ h' \|^3 } \leq  \|h'\| \cdot \E{\| \mZ h'\|^2 }
\; = \; \frac{\tau}{d}\|h'\|^3, \qquad \forall h' \in \R^d.
$$
Lastly, note that 
{ 
\footnotesize
\begin{eqnarray*}
&& \E{\mZ \nabla^2 f(x^k) \mZ} \\
& = & \E{\mZ \left(\nabla^2 f(x^k)\right)^{\frac12}} \E{ \left(\nabla^2 f(x^k)\right)^{\frac12}\mZ} 
\\
&& \qquad + \E{ \left(  \mZ \left(\nabla^2 f(x^k)\right)^{\frac12} -   \E{\mZ \left(\nabla^2 f(x^k)\right)^{\frac12}} \right)\left(  \mZ \left(\nabla^2 f(x^k)\right)^{\frac12} -   \E{\mZ \left(\nabla^2 f(x^k)\right)^{\frac12}} \right)^\top  } 
\\
&=& \frac{\tau^2}{d^2} \nabla^2 f(x^k) +\E{  \left( \mZ - \frac{\tau}{d} \mI^d \right)\nabla^2 f(x^k) \left( \mZ - \frac{\tau}{d} \mI^d \right)  } 
\\
&\leq & \frac{\tau^2}{d^2} \nabla^2 f(x^k) + L \E{  \left( \mZ - \frac{\tau}{d} \mI^d \right)^2  }
 \\
&= & \frac{\tau^2}{d^2} \nabla^2 f(x^k) + \frac{\tau(d-\tau)}{d^2}L\mI^d.
\end{eqnarray*}
}

Therefore, we conclude
\begin{eqnarray*}
\E{  F(x^{k + 1}) \, | \, x^k }& \leq &
f(x^k) + \frac{\tau}{d} \la \nabla f(x^k), y - x^k \ra 
+ \frac{\tau(d-\tau)}{d^2}\cdot \frac{L}{2} \| y - x^k \|^2
\\
\\
& \; & \quad + \quad \frac{\tau^2}{d^2}  \cdot \frac{1}{2} \la \nabla^2 f(x^k)(y - x^k), y - x^k \ra
+ \frac{\tau}{d} \cdot \frac{M}{6}\|y - x^k\|^3 \\
\\
& \; & \quad + \quad \frac{\tau}{d} \psi(y) + \left(1 - \frac{\tau}{d}\right) \psi(x^k). \\
\end{eqnarray*}

Finally, by convexity and from Lipschitz continuity of the Hessian~\eqref{eq:sscn_coordinate_ub},
we have the following upper estimate:
\begin{eqnarray*}
& \la \nabla f(x^k), y - x^k \ra + \frac{\tau }{d} \cdot \frac{1}{2} \la \nabla^2 f(x^k)(y - x^k), y - x^k \ra \\
\\
& \qquad \qquad \; = \;
\frac{d - \tau}{d } \la \nabla f(x^k), y - x^k \ra
\\ \\
& \qquad \qquad \qquad+ \frac{\tau }{d } \Bigl( 
\la \nabla f(x^k),y - x^k \ra + \frac{1}{2} \la \nabla^2 f(x^k)(y - x^k), y - x^k \ra
\Bigr) \\
\\
& \qquad \qquad \; \leq \; 
\frac{d - \tau}{d } \Bigl( f(y) - f(x^k) \Bigr)
+ \frac{\tau }{d } \Bigl( 
f(y) - f(x^k) + \frac{M}{6}\|y - x^k\|^3
\Bigr) \\
\\
& \qquad \qquad  \; \leq \; f(y) - f(x^k) + \frac{M}{6}\|y - x^k\|^3.
\end{eqnarray*}
which completes the proof. \qed

\subsection{Proof of Theorem~\ref{thm:sscn_global_weakly}}

Let us denote the following auxiliary sequences:
\begin{eqnarray*}
a_{k} & \Def & k^2, \qquad A_{k} \; \Def \; A_{0} + \sum\limits_{i = 1}^k a_i, \qquad k \geq 1,
\end{eqnarray*}
and
\begin{eqnarray*}
A_0 & \Def & \frac{4}{3}\left( \frac{d}{\tau} \right)^3.
\end{eqnarray*}
Then, we have an estimate
\begin{equation}
 \label{A_k_grows}
A_{k} \; = \; A_0 + \sum\limits_{i = 1}^k i^2 \geq 
A_0 + \int\limits_{0}^k x^2 dx \; = \; A_0 + \frac{k^3}{3}.
\end{equation}
Now, let us fix iteration counter $k \geq 0$ and set
\begin{eqnarray*}
\alpha_k & \Def & \frac{d}{\tau} \frac{a_{k + 1}}{A_{k + 1}}
\quad \Leftrightarrow \quad 1 - \frac{\tau}{d} \alpha_k \; = \;  \frac{A_k}{A_{k + 1}}.
\end{eqnarray*}
Note that we have $\alpha_k \leq 1$ by the choice of $A_0$, since it holds
\begin{eqnarray*}
\max\limits_{\tau \geq 0} \frac{\tau^2}{A_0 + \frac{\tau^3}{3}} & = & \frac{\tau}{d}.
\end{eqnarray*}
Let us plug $y \equiv \alpha_k x^{*} + (1 - \alpha_k) x^k$ into~\eqref{GlobalUpper}.
By convexity we obtain
\begin{eqnarray*}
&&
\E{ F(x^{k + 1}) \, | \, x^k } 
\\
& \leq & 
\Bigl(1 - \frac{\tau}{d} \Bigr) F(x^k) + \frac{\tau}{d}\alpha_k F^{*} + \frac{\tau}{d}(1 - \alpha_k) F(x^k) \\
\\
& \; & \quad + \quad 
\frac{\tau}{d}\biggl( \frac{d - \tau}{d } \frac{L \|x^k - x^{*}\|^2}{2}  \alpha_k^2
+ \frac{M \|x^k - x^{*}\|^3}{3} \alpha_k^3 \biggr) \\
\\
& = &
\frac{A_k}{A_{k + 1}} F(x^k) + \frac{a_{k + 1}}{A_{k + 1}} F^{*}
+ \frac{d}{\tau} \frac{d - \tau}{d }  \frac{L \|x^k - x^{*}\|^2}{2} \left( \frac{a_{k + 1}}{A_{k + 1}} \right)^2 \\ 
\\
&& \qquad 
+ \left(\frac{d}{\tau}\right)^2 \frac{M \|x^k - x^{*}\|^3}{3} \left( \frac{a_{k + 1}}{A_{k + 1}}  \right)^3 \\
\\
& \leq & 
\frac{A_k}{A_{k + 1}} F(x^k) + \frac{a_{k + 1}}{A_{k + 1}} F^{*}
+ \frac{d - \tau}{2\tau} LR^2 \left( \frac{a_{k + 1}}{A_{k + 1}} \right)^2
+ \left( \frac{d}{\tau} \right)^2 \frac{M R^3}{3} \left( \frac{a_{k + 1}}{A_{k + 1}}  \right)^3.
\end{eqnarray*}
Therefore, for the residual $\delta_k \Def \E{ F(x^k) }- F^{*}$ we have the following bound
\begin{eqnarray*}
A_{k + 1} \delta_{k + 1} & \leq & A_k \delta_k 
+ \frac{d - \tau}{2\tau}LR^2 \frac{a_{k + 1}^2}{A_{k + 1}}
+ \left( \frac{d}{\tau} \right)^2 \frac{M R^3}{3} \frac{a_{k + 1}^3}{A_{k + 1}^2}, \quad k \geq 0.
\end{eqnarray*}
Summing up these inequalities for different $k$, we obtain
\begin{eqnarray*}
A_k \delta_k & \leq & A_0 \delta_0 
+ \frac{d - \tau}{2\tau} L R^2 \sum\limits_{i = 1}^k \frac{a_i^2}{A_i}
+ \left( \frac{d}{\tau} \right)^2 \frac{M R^3}{3} \sum\limits_{i = 1}^k \frac{a_i^3}{A_i^2}, \qquad k \geq 1.
\end{eqnarray*}
To finish the proof it remains to notice that
\begin{eqnarray*}
\sum\limits_{i = 1}^k \frac{a_i^2}{A_i} & \refLE{A_k_grows} &
\sum\limits_{i = 1}^k \frac{i^4}{A_0 + \frac{1}{3} i^3 }
\; \leq \;
3 \sum\limits_{i = 1}^k i \; \leq \; 3 k^2,
\end{eqnarray*}
and
\begin{eqnarray*}
\sum\limits_{i = 1}^k \frac{a_i^3}{A_i^2} & \refLE{A_k_grows} &
\sum\limits_{i = 1}^k \frac{i^6}{(A_0 + \frac{1}{3}i^3 )^2} 
\; \leq \;
9 k.
\end{eqnarray*}

\qed

\subsection{Proof of Theorem~\ref{thm:sscn_global_strongly}}

Given that Assumption~\ref{as:sscn_sc} (strong convexity) is satisfied, the following inequality holds
\begin{eqnarray*}
\frac{\mu}{2}\|x - x^{*}\|^2 & \leq & F(x) - F^{*}, \qquad  \forall x \R^d,
\end{eqnarray*}
and thus we have a bound for the radius of level sets~\eqref{Rdef}:
\begin{eqnarray*}
R^2 & \leq & \frac{2}{\mu}(F(x^0) - F^{*}).
\end{eqnarray*}
Combining the above with~\eqref{GlobalConv} we obtain the following convergence estimate for $k\geq 1$:
\begin{eqnarray*}
\E{ F(x^k) - F^{*} } & \leq & 
\left(
\frac{d - \tau}{\tau} \cdot \frac{18 L}{\mu k}
+ \bigl( \frac{d}{\tau} \bigr)^2 \cdot
\frac{18 M R }{\mu k^2}
+ \frac{1}{1 + \frac{1}{4}\bigl( \frac{\tau}{d} k \bigr)^3}
\right) \cdot 
\bigl( F(x^0) - F^{*} \bigr).
\end{eqnarray*}
Therefore, we get the linear decrease of the expected residual
\begin{eqnarray*}
\E{ F(x^k) - F^{*}  }
& \leq & \frac{1}{2}
\bigl( F(x^0) - F^{*} \bigr),
\end{eqnarray*}
as soon as the following three bounds for $k$ are all reached:
\begin{enumerate}
	\item $\frac{d - \tau}{\tau} \cdot \frac{18L}{\mu k} \leq \frac{1}{6}
	\quad \Leftrightarrow \quad
	k \geq 108 \frac{d - \tau}{\tau} \cdot \frac{L}{\mu}$.
	
	\item $\bigl( \frac{d}{\tau} \bigr)^2 \cdot \frac{18 M R}{\mu k^2} \leq \frac{1}{6}
	\quad \Leftrightarrow \quad
	k \geq \frac{d}{\tau} \sqrt{ 108  \frac{M R}{\mu}   }.
	$
	
	\item $\frac{1}{1 + \frac{1}{4}\bigl( \frac{\tau}{d} k^3 \bigr)^3} \leq \frac{1}{6}
	\quad \Leftrightarrow \quad
	k \geq \frac{d}{\tau} 20^{1/3}$.
\end{enumerate}
\qed

\section{Proofs for Section~\ref{sec:sscn_local} \label{sec:sscn_local_proofs}}
\subsection{Several technical lemmas}
 It will be convenient to denote the Newton decrement as follows:
 
 \begin{equation} 
 \label{eq:sscn_newton_decrement}
 \lambda_f(x) \eqdef \left(\nabla f(x)^\top \left(\nabla^2 f(x)\right)^{-1}\nabla f(x) \right)^\frac12
 \end{equation}
 and a sublevel set of $x^0$ as $\level$; i.e. $\level \eqdef \{x ; f(x)\leq f(x^0)\} $.

\begin{lemma} (Local bounds)
Suppose that $x^0$ is such that 
\[f(x^0)-f(x^*) \leq \varrho^4 \frac{2\left(\min_{x\in \level}\lambda_{\min} \nabla^2_{\mS} f(x)\right)^4}{LM_{\mS}^2 \|S \|^2}\] for some $\varrho>0$. Then, we have 
\begin{equation}\label{eq:sscn_ndajbdhahbj}
\sqrt{ \frac{M_{\mS}}{2}} \| \mS^\top \nabla f(x^k)\|^{\frac12} I \preceq \varrho \nabla^2_{\mS}f(x^k).
\end{equation} Suppose further that $f(x^0) - f(x^*)\leq \varphi^2 \frac{ \mu\left(\lambda_{\min} \nabla^2_{\mS} f(x^*)\right)^2 }{2 M_{\mS}^2}$ for some $\varphi>0$. Then we have 
\begin{equation}\label{eq:sscn_dasjnlsdhjkd}
(1+\varphi)^{-1} \nabla^2_{\mS}f(x^*) \preceq \nabla^2_{\mS}f(x^k) \preceq (1+\varphi) \nabla^2_{\mS}f(x^*).
\end{equation}
Lastly, if $f(x^0)-f(x^*)\leq \omega^{-1}\left( \frac{2\mu^{\frac32}}{(1+\gamma^{-1})M} \right)$ where $\omega(y)\eqdef y - \log(1+y)$ and $\gamma>0$, we have
\begin{equation} \label{eq:sscn_dojasjodasj}
f(x^k) -f(x^*)\leq  \frac12 (1+\gamma) \lambda_f(x^k)^2.
\end{equation}
\end{lemma}

\begin{proof}
For the sake of simplicity, let $x = x^k$ and $\mS = \mS^k$ throughout this proof.
For the first part, we have
\begin{eqnarray*}
\sqrt{ \frac{M_{\mS}}{2}} \| \mS^\top \nabla f(x)\|^{\frac12} \mI^{\tau(\mS)}
&\leq& 
\sqrt{ \frac{M_{\mS}}{2}} \| \mS\|^\frac12 \| \nabla f(x)\|^{\frac12}  \mI^{\tau(\mS)} 
\\
& \leq& 
\sqrt{ \frac{M_{\mS}}{2}}  \| \mS\|^\frac12  2^{\frac14}L^{\frac14}\left(f(x^0)-f(x^*) \right)^\frac14  \mI^{\tau(\mS)}
\\
&\leq&
\varrho\min_{x\in \level}\lambda_{\min} \nabla^2_{\mS} f(x)   \mI^{\tau(\mS)}
\\
& \leq&
 \varrho \nabla^2_{\mS}f(x).
\end{eqnarray*}

For the second part, we have
\begin{eqnarray*}
\nabla^2_{\mS} f(x^k)-\nabla^2_{\mS}  f(x^*) &\preceq & 
M_{\mS} \|x^k-x^* \| \mI^{\tau(\mS)} 
\\
&\preceq&  M_{\mS}\sqrt{\frac{2 (f(x^k) - f(x^*))}{\mu}}  \mI^{\tau(\mS)}
\\
&\preceq&   M_{\mS}\sqrt{\frac{2 (f(x^0) - f(x^*))}{\mu}}  \mI^{\tau(\mS)}
\\
&\preceq & 
\varphi  \nabla^2_{\mS} f(x^*).
\end{eqnarray*}
Therefore, we can conclude that $\nabla^2_{\mS} f(x) \preceq (1+\varphi) \nabla^2_{\mS} f(x^*)$. Analogously we can show $\nabla^2_{\mS} f(x^*) \preceq (1+\varphi) \nabla^2_{\mS} f(x)$ and thus~\eqref{eq:sscn_dasjnlsdhjkd} follows.

Lastly, if $f(x^0)-f(x^*)\leq \omega\left( \frac{2\mu^{\frac32}}{(1+\gamma^{-1})M} \right)$, then due to~\cite{nesterov2018lectures} we have
\[
\omega\left(\lambda_f(x^k)\right)\leq  f(x^k)-  f(x^*) \leq f(x^0)-  f(x^*) \leq  \omega\left( \frac{2\mu^{\frac32}}{(1+\gamma^{-1})M} \right)
\]
and thus $\lambda_f(x^k) \leq \frac{2\mu^{\frac32}}{(1+\gamma^{-1})M}$. Now~\eqref{eq:sscn_dojasjodasj} follows from Lemma~\ref{lem:sscn_selfconc} and Lemma~\ref{lem:sscn_sc2}.
\end{proof}

 \begin{lemma}\label{lem:sscn_selfconc}
 Function $f$ is $\frac{M}{\mu^{\frac32}}$ self-concordant.
 \end{lemma}
 \begin{proof}
 \[
 \frac{M}{\mu^{\frac32}} \|u\|^3_{\nabla^2 f(x)}  \geq  M\|u\|^3 \geq  \nabla^3 f(x)[u,u,u]
 \]
 \end{proof}
 
\begin{lemma}\label{lem:sscn_sc2}
Consider any $\gamma\in \R^+$ and suppose that $f$ is $\varsigma$ self-concordant. Then if  $\lambda_f(x)  < \frac{2}{(1+\gamma^{-1})\varsigma}$ we have
\begin{equation}\label{eq:sscn_sc2}
f(x) - f(x^*)\leq \frac12 \left(1+\gamma\right)  \lambda_f(x)^2   
\end{equation}
\end{lemma}

\begin{proof}
Define $\omega_*(z)\eqdef -z-\ln(1-z)$.
Note first that, $h(x)\eqdef \frac{\varsigma^2}{4} f(x)$ is 2 self concordant~\cite{nesterov2018lectures}. As a consequence, if $\lambda_h(x)<1$ we have~\cite{nesterov2018lectures}
\[
h(x) - h(x^*) \leq  \omega_*( \lambda_h(x)  ).
\] 
If further $\lambda_h(x) \leq \frac{1}{1+\gamma^{-1}}$ due to Lemma~\ref{lem:sscn_snjodanj}, we get
\[
\omega_*( \lambda_h(x)  ) \leq \left(1+\gamma\right) \frac {\lambda_h(x) ^2}{2}.
\] 
As $\lambda_h(x) = \frac{\varsigma}{2}\lambda_f(x)$, we get~\eqref{eq:sscn_sc2}.
\end{proof}

\begin{lemma}\label{lem:sscn_snjodanj}
Let $c\in \R^+$ and $0\leq y \leq \frac{1}{1+c}$. Then we have $\omega_*( y  ) \leq \left(1+\frac{1}{c}\right) \frac {y ^2}{2}.$
\end{lemma}
\begin{proof}
Clearly $\omega_*( y  ) = \sum_{i=2}^\infty \frac{y^i}{i}$ and thus function $\left(1+\frac{1}{c}\right) \frac {y ^2}{2} - \omega_*( y  ) $ is non-increasing for $y\geq 0$. Therefore, it suffices to check verify $\left(1+\frac{1}{c}\right) \frac {1}{2(1+c)^2} - \omega_*( \frac{1}{1+c} ) \geq 0$, which is easy task for Mathematica, see Figure~\ref{fig:sscn_proof}.
\begin{figure}[H]
\centering
\includegraphics[scale=0.6]{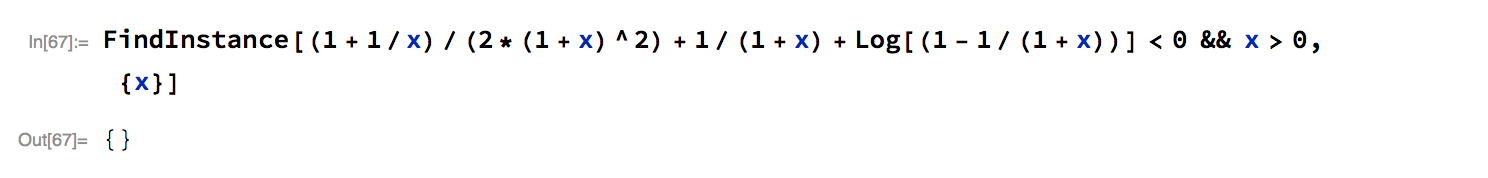}
\caption{Proof of $\left(1+\frac{1}{c}\right) \frac {1}{2(1+c)^2} - \omega_*( \frac{1}{1+c} ) \geq 0$ for all $c>0$.}
\label{fig:sscn_proof}
\end{figure}
\end{proof}

\subsection{Proof of Lemma~\ref{lem:sscn_decrease}}
Note that the update rule of {\tt SSCN} yields immediately (using first-order optimality conditions)

\begin{equation} \label{eq:sscn_grad_equality}
-\mS^\top \nabla f(x)= \left( \nabla^2_{\mS}f(x) + \frac{1}{2} M_{\mS}\|x^+-x\| \mI^{\tau(\mS)}\right)\left(x^+-x\right)
\end{equation}

and therefore

\begin{eqnarray}
\left \| \mS^\top \nabla f(x) \right \|^{\frac12}
 &=&
 \left(
 \left(x^+-x\right)^\top\left( \nabla^2_{\mS}f(x) + \frac{1}{2} M_{\mS}\|x^+-x\| I\right)^2\left(x^+-x\right) \right)^{\frac14}\nonumber
  \\
 & \geq & 
  \left(
 \left(x^+-x\right)^\top\left(  \frac{1}{2} M_{\mS}\|x^+-x\| I\right)^2\left(x^+-x\right) \right)^{\frac14}\nonumber
 \\ 
 &=&
 \sqrt{\frac{M_{\mS}}{2}}\|x^+-x\|.
 \label{eq:sscn_stupid_bound}
\end{eqnarray}

Furthermore, taking dot product of~\eqref{eq:sscn_grad_equality} with $(x^+-x)$ yields
\[
\left\langle \mS^\top \nabla f(x), x^+-x\right\rangle+\left\langle  \nabla^2_{\mS}f(x)\left(x^+-x\right), x^+-x\right\rangle+\frac{1}{2} M_{\mS} \|x^+-x\|^3=0
\]
and thus 

\[
\begin{aligned} f(x)-f(x^+) 
&\stackrel{\eqref{eq:sscn_coordinate_ub}}{\geq}
\left\langle  \mS^\top \nabla f(x), x-x^+\right\rangle-\frac{1}{2}\left\langle   \nabla^2_{\mS}f(x) (x^+-x), x^+-x\right\rangle-\frac{M_{\mS}}{6}  \|x^+-x\|^3
 \\ 
 &=
\frac{1}{2}\left\langle   \nabla^2_{\mS} f(x) (x^+-x), x^+-x\right\rangle+\frac{M_{\mS}}{3}  \|x^+-x\|^3 
\\
 &\stackrel{(*)}{\geq}
 \frac12
\left(x^+-x\right)^\top\left( \nabla^2_{\mS}f(x) + \frac{1}{2} M_{\mS}\|x^+-x\| I\right)\left(x^+-x\right)
\\
 &\stackrel{\eqref{eq:sscn_grad_equality}}{=}
 \frac12 \nabla f(x)^\top \mS \left( \nabla^2_{\mS}f(x) + \frac{1}{2} M_{\mS}\|x^+-x\| I\right)^{-1}  \mS^\top\nabla f(x)
 \\
 &\stackrel{\eqref{eq:sscn_stupid_bound}}{\geq }
 \frac12 \nabla f(x)^\top \mS \left( \nabla^2_{\mS}f(x) +\sqrt{ \frac{M_{\mS}}{2}} \| \mS^\top \nabla f(x)\|^{\frac12} I\right)^{-1}  \mS^\top \nabla f(x).
  \end{aligned}
\]
Above, in inequality $(*)$ we have used the fact that matrix $  \left( \nabla^2_{\mS}f(x) + \frac{1}{2} M_{\mS}\|x^+-x\| \mI^{\tau(\mS)}\right)$ is invertible since $f$ is strongly convex and thus $\nabla^2_{\mS}f(x) \succ 0$.

\subsection{Proof of Theorem~\ref{thm:sscn_local}}
First, suppose that $f(x^0)-f(x^*) \leq \varrho^4 \frac{2\left(\min_{x\in \level}\lambda_{\min} \nabla^2_{\mS} f(x)\right)^4}{LM_{\mS}^2\|\mS\|^{2}}$ for some $\varrho>0$. Using the fact that $\nabla^2_{\mS} f(x)$ is invertible ($\mS$ has full column rank and $\nabla^2 f(x) \succ 0$) we have
{ 
\footnotesize
\begin{eqnarray}
\E{ \frac12 \| \mS^\top \nabla f(x^k)\|^2_{ \left( \mH(x^k) \right)^{-1} } } 
& \stackrel{\eqref{eq:sscn_ndajbdhahbj}}{\geq}&
\E{  \frac12 \nabla f(x)^\top {\mS}\left((1+\varrho) \nabla^2_{\mS}f(x)\right)^{-1} \mS^\top \nabla f(x)}
\nonumber
\\
&=&
 \frac{1}{2(1+\varrho)} \nabla f(x)^\top \E{ {\mS}\left( \nabla^2_{\mS}f(x)\right)^{-1}{\mS}^\top }\nabla f(x).
 \label{eq:sscn_dnjasnjaaa}
\end{eqnarray}
}
If further $f(x^0) - f(x^*)\leq  \varphi^2\frac{ \mu \left(\lambda_{\min}\nabla^2_{\mS} f(x^*)\right)^2 }{2 M_{\mS}^2}$ for some $\varphi>0$ we get

\begin{eqnarray}
\nonumber
\E{ \frac12 \| \mS^\top \nabla f(x^k)\|^2_{ \left( \mH(x^k) \right)^{-1} } } 
&\stackrel{\eqref{eq:sscn_dnjasnjaaa}}{ \geq}&
  \frac{ \nabla f(x)^\top \E{ {\mS}\left( \nabla^2_{\mS}f(x)\right)^{-1}{\mS}^\top }\nabla f(x)}{2(1+\varrho)}
  \\
  \nonumber
  &\stackrel{\eqref{eq:sscn_dasjnlsdhjkd}}{ \geq}&
    \frac{ \nabla f(x)^\top \E{ {\mS}\left( \nabla^2_{\mS}f(x^*)\right)^{-1}{\mS}^\top }\nabla f(x) }{2(1+\varrho)(1+\varphi)}
  \\
  \nonumber
  &\stackrel{\eqref{eq:sscn_sc_generalized}}{ \geq}&
  \frac{\nabla f(x)^\top \left(  \zeta \left(  \nabla^2f(x^*)\right)^{-1}\right) \nabla f(x) }{2(1+\varrho)(1+\varphi)}
    \\
    \label{eq:sscn_cjdinusvsibu}
  &\stackrel{\eqref{eq:sscn_dasjnlsdhjkd}}{ \geq}&
    \frac{\zeta \lambda_f(x)^2}{2(1+\varrho)(1+\varphi)^2}
\end{eqnarray}

Lastly, if if $f(x^0)-f(x^*)\leq \omega^{-1}\left( \frac{2\mu^{\frac32}}{(1+\gamma^{-1})M} \right)$ where $\omega(y)\eqdef y - \log(1+y)$ and $\gamma>0$, we get

\begin{eqnarray*}
\E{ \frac12 \| \mS^\top \nabla f(x^k)\|^2_{ \left( \mH(x^k) \right)^{-1} } } 
  &\stackrel{\eqref{eq:sscn_cjdinusvsibu}}{ \geq}&
      \frac{\zeta \lambda_f(x)^2}{2(1+\varrho)(1+\varphi)^2}
      \\
  &\stackrel{\eqref{eq:sscn_dojasjodasj}}{ \geq}&
    \frac{\zeta (f(x)-f(x^*))}{(1+\varrho)(1+\varphi)^2(1+\gamma)}
\end{eqnarray*}

and thus~\eqref{eq:sscn_local_rate} follows. In particular for any $\varrho, \varphi, \gamma>0$, we can choose 
\[
\delta = \min \left\{  \varrho^4 \frac{2\left(\min_{x\in \level}\lambda_{\min} \nabla^2_S f(x)\right)^4}{LM_S^2},  \varphi^2\frac{ \mu \left(\lambda_{\min}\nabla^2_S f(x^*)\right)^2 }{2 M_S^2} ,\omega^{-1}\left( \frac{2\mu^{\frac32}}{(1+\gamma^{-1})M} \right) \right\}
\]
and 
\[
\varepsilon = 1 -\frac{1}{(1+\varrho)(1+\varphi)^2(1+\gamma)}.
\]
\qed

\chapter{Appendix for Chapter \ref{ami}}
\label{ami_appendix}

\graphicspath{{AMI/experiments/}}

\section{Proofs for Section~\ref{sec:ami_ami_paper} \label{sec:ami_proof_Euclidean}}

\subsection{Proof of Lemma~\ref{lem:ami_Z_bounds}}
First note that $\mZ$ is a self-adjoint positive operator and thus so is $\E{\mZ}.$  Consequently.
 \begin{eqnarray}
 \theta &\overset{\eqref{eq:ami_mu+nu}}{=}& \inf_{x \in \Range{\A^*}} \frac{\dotprod{\E{\mZ}x,x}}{\dotprod{x,x}} \nonumber \\
 &\overset{\eqref{eq:ami_exactness}}{=} &
 \inf_{x \in \Range{\E{\mZ}}} \frac{\dotprod{\E{\mZ}x,x}}{\dotprod{x,x}} \nonumber\\
 & \overset{\mbox{  Lemma~\ref{lem:ami_pseudo} item~\ref{it:pseudorange} }}{=}& \inf_{x \in \cX} \frac{\dotprod{\E{\mZ}\E{\mZ}^\dagger x,\E{\mZ}^\dagger x}}{\dotprod{\E{\mZ}^\dagger x,\E{\mZ}^\dagger x}} \nonumber\\
  & \overset{\mbox{  Lemma~\ref{lem:ami_pseudo} item~\ref{it:pseudoTTdagT} }}{=}& \inf_{x \in \cX} \frac{\dotprod{\E{\mZ}^\dagger x,x}}{\dotprod{\E{\mZ}^\dagger x,\E{\mZ}^\dagger x}} \nonumber\\
   & \overset{\mbox{  Lemma~\ref{lem:ami_squareroot} }}{=}& 
  \inf_{z \in \Range{(\E{\mZ}^\dagger)^{1/2}}} \frac{\dotprod{ z,z}}{\dotprod{\E{\mZ}^\dagger z, z}} \qquad (\mbox{set }z = (\E{\mZ}^\dagger)^{1/2}x)\nonumber  \\
  & \overset{\eqref{eq:ami_RangeGhalf} }{=}& 
  \frac{1}{\norm{\E{\mZ}^{\dagger}}}.\label{eq:ami_as9d8n923}
 \end{eqnarray} 
 
 For the bounds~\eqref{eq:ami_nubnds} we have that
\begin{eqnarray*}
 \nu &\overset{\eqref{eq:ami_mu+nu}}{=}& \sup_{x \in \Range{\A^*}} \frac{\E{\dotprod{\E{\mZ}^\dagger \mZ x,\mZ x}}}{\dotprod{\E{\mZ}x,x}}  \\
 &\leq & \sup_{x \in \Range{\A^*}} \frac{\norm{\E{\mZ}^\dagger} \E{\norm{\mZ x}_2^2}}{\dotprod{\E{\mZ}x,x}}\\
 & =& \norm{\E{\mZ}^\dagger}\\
 & \overset{\eqref{eq:ami_as9d8n923}}{\leq} & \frac{1}{\theta}.
 \end{eqnarray*} 
To bound $\nu$ from below we use that $\E{\mZ}^\dagger$ is self adjoint together with that the map $\mX \mapsto \dotprod{\mX\E{\mZ}^\dagger \mX x,x}$ is convex over the space of self-adjoint operators $\mX \in L(\cX)$ and for a fixed $x \in \cX$. Consequently by Jensen's inequality
\begin{equation}
\E{\dotprod{\mZ\E{\mZ}^\dagger \mZ x,x}} \geq \dotprod{\E{\mZ}\E{\mZ}^\dagger \E{\mZ}x,x}
\overset{ \mbox{  Lemma~\ref{lem:ami_pseudo} item~\ref{it:pseudoTTdagT} }}{=}\dotprod{\E{\mZ} x, x}.\label{eq:ami_ZEZZjen}
\end{equation}
Finally
\begin{eqnarray*}
\nu & \overset{\eqref{eq:ami_ZEZZjen}}{\geq} & \sup_{x \in \Range{\A^*}} \frac{\dotprod{\E{\mZ}x,x}}{\dotprod{\E{\mZ}x,x}}  = 1.
\end{eqnarray*}

Lastly, to show \eqref{eq:ami_nu_lower} we have
 \begin{eqnarray*}
 \Rank{\A^*} & \overset{\eqref{eq:ami_exactness}}{=}& \Rank{\E{\mZ}} \\
  & \overset{\mbox{  Lemma~\ref{lem:ami_projrank}+ Lemma~\ref{lem:ami_pseudo} (\emph{\ref{it:pseudoproj}}) }}{=}& \Tr{\E{\mZ}\E{\mZ}^\dagger} = \E{\Tr{\mZ \E{\mZ}^\dagger}}\\
 &= & \E{\Tr{\mZ \E{\mZ}^\dagger \mZ}} \\
 &\leq & \nu \E{\Tr{\mZ}} \overset{\mbox{  Lemma~\ref{lem:ami_projrank}}}{=} \nu   \E{\Rank{\mZ}},
 \end{eqnarray*}
where we used that  $\dotprod{ \E{\mZ \E{\mZ}^\dagger \mZ} u,u}  \leq  \nu \dotprod{\E{\mZ} u,u}$ for every $u \in \Range{\E{\mZ}}=\Range{\A^*}=\cX.$ \qed

 {\bf Proof} that $\mX \mapsto \dotprod{\mX\E{\mZ}^{\dagger}\mX x,x} = \norm{\mX x}_{\E{\mZ}^{\dagger}}^2$ is convex: Let $\mG= \E{\mZ}^{\dagger}$ then
 \begin{eqnarray*}
\norm{(\lambda \mX +(1-\lambda)\mY)x}_\mG^2 &= &
\lambda^2 \norm{\mX x}_\mG^2 +(1-\lambda)^2  \norm{\mY x}_\mG^2
+2\lambda(1-\lambda) \dotprod{x \mX \mG\mY ,x}\\
&= & -\lambda(1-\lambda)\norm{(\mX-\mY)x}_\mG^2 \\
& &+ \lambda  \norm{\mX x}_\mG^2+(1-\lambda)\norm{\mY x}_\mG^2\\
&\leq &\lambda  \norm{\mX x}_\mG^2+(1-\lambda)\norm{\mY x}_\mG^2 . \hskip5cm \qed 
\end{eqnarray*}

\subsection{Technical lemmas to prove Theorem~\ref{theo:conv}}

\begin{lemma} \label{lem:ami_invariant}For all $k \geq 0,$ the vectors $y_k-x_*, \, x_k-x_*$ and $v_k-x_*$ belong to $\Range{\A^*}.$
\end{lemma} 
\begin{proof}
 Note that $x_0 =y_0 =x_0$ and in view of~\eqref{eq:ami_xsol} we have $x_* \in x_0 + \Range{\A^*}.$ So $y_0 -x_* \in \Range{\A^*},$  $v_0 -x_* \in \Range{\A^*}$ and  $x_0 -x_* \in \Range{\A^*}.$ Assume by induction that  $y_k -x_* \in \Range{\A^*},$ $v_k -x_* \in \Range{\A^*}$ and  $x_k -x_* \in \Range{\A^*}.$ Since $g_k \in \Range{\A^*}$ and $x_{k+1} = y_k - g_k$ we have
 \[x_{k+1}-x_* = (y_k-x_*) -g_k \in \Range{\A^*}.\]
 Moreover, 
 \[v_{k+1}-x_* = \beta(v_k -x_*) + (1-\beta)(y_k-x*)-\gamma g_k \in \Range{\A^*}.\]
 Finally
  \[y_{k+1}-x_* = \eta v_{k+1}+(1-\eta)x_{k+1}-x_* = \eta (v_{k+1}-x_*)+(1-\eta)(x_{k+1}-x_*)\in \Range{\A^*}.\]
\end{proof}

\begin{lemma}
\begin{equation}\label{eq:ami_Enormnubnd}
\E{ \norm{\cZ_k(y_k -x_*)}_{\E{\cZ}^{\dagger}}^2 \, | \, y_k} \leq \nu \norm{y_k -x_*}_{\E{\cZ}}^2
\end{equation} 
\end{lemma}
\begin{proof}
Since $y_k -x_* \in \Range{\A^*}$ we have that
\begin{eqnarray*}
\E{ \norm{\cZ_k(y_k -x_*)}_{\E{\cZ}^{\dagger}}^2 \, | \, y_k} &=  & 
 \dotprod{\E{\cZ_k\E{\cZ}^{\dagger}\cZ_k}(y_k -x_*),(y_k -x_*)} \\
 & \overset{\eqref{eq:ami_mu+nu}}{\leq} &  \nu \dotprod{\E{\cZ}(y_k -x_*),(y_k -x_*)} \\
 &= & \nu \norm{y_k -x_*}_{\E{\cZ}}^2. 
\end{eqnarray*}
\end{proof}

\begin{lemma}
\label{lemma:threepoint}
\begin{equation}\label{eq:ami_ykxkident}
 \norm{y_k-x_*}_{\E{\cZ}}^2 = \norm{y_k-x_*}^2 -\E{\norm{x_{k+1}-x_*}^2 \, | \, y_k}\end{equation}
\end{lemma}
\begin{proof}
\begin{eqnarray*}
\E{\norm{x_{k+1}-x_*}^2 \, | \, y_k} &= & \E{\norm{(\cI-\cZ_k)( y_k -x_*)}^2 \, | \, y_k} \\
&= &\dotprod{(\cI-\E{\cZ})( y_k -x_*),y_k -x_*}\\
&= &\norm{y_k -x_*}^2 - \norm{y_k -x_*}^2_{\E{\cZ}}. 
\end{eqnarray*} 
\end{proof}
\subsection{Proof of Theorem \ref{theo:conv}}
\label{subsec:proof_main_thm}
Let $r_k \eqdef \norm{v_{k}-x_*}^2_{\E{\cZ}^\dagger}$. It follows that
\begin{eqnarray}
r_{k+1}^2 & =& \norm{v_{k+1}-x_*}^2_{\E{\cZ}^\dagger} \nonumber\\
& =& \norm{\beta v_k +(1-\beta)y_k -x_*- \gamma \cZ_k(y_k -x_*)}^2_{\E{\cZ}^\dagger}\nonumber\\
&=&\underbrace{\norm{\beta v_k +(1-\beta)y_k -x_*}^2_{\E{\cZ}^\dagger}}_{I} +\gamma^2\underbrace{\norm{\cZ_k(y_k -x_*)}^2_{\E{\cZ}^\dagger}}_{II} \nonumber \\
& &-2\gamma \underbrace{\dotprod{\beta (v_k-x_*) +(1-\beta)(y_k -x_*), \E{\cZ}^\dagger \cZ_k(y_k -x_*)}}_{III} \nonumber \\
& =& I + \gamma^2 II -2\gamma III.\label{eq:ami_proofstep1}
\end{eqnarray}

The first term can be upper bounded as follows
\begin{eqnarray}
I &=& \norm{\beta (v_k-x_*) +(1-\beta)(y_k -x_*)}^2_{\E{\cZ}^\dagger} \nonumber \\
& =&\beta^2 \norm{v_k-x_*}^2_{\E{\cZ}^\dagger}  + (1-\beta)^2 \norm{y_k-x_*}^2_{\E{\cZ}^\dagger} +2\beta(1-\beta) \dotprod{v_k-x_*,y_k -x_*}_{\E{\cZ}^\dagger}\nonumber \\
&\overset{\eqref{eq:ami_paral1}}{=} & \beta\norm{v_k-x_*}^2_{\E{\cZ}^\dagger}  + (1-\beta) \norm{y_k-x_*}^2_{\E{\cZ}^\dagger} -\beta(1-\beta) \norm{v_k-y_k}^2_{\E{\cZ}^\dagger} \nonumber \\
& \leq & \beta r_k^2 + (1-\beta) \norm{y_k-x_*}^2_{\E{\cZ}^\dagger},\label{eq:ami_Ibnded}
\end{eqnarray} 
where in the third equality we used a form of the parallelogram identity
\begin{equation}\label{eq:ami_paral1}
2\dotprod{u,v} = \norm{u}^2 + \norm{v}^2 - \norm{u-v}^2,
\end{equation}
with $u = v_k-x_*$ and $v = y_k-x_*.$

Taking expectation with to $\cS_k$ in the third term in~\eqref{eq:ami_proofstep1} gives
\begin{eqnarray}
\E{III \, | \, y_k, v_k, x_k} &= & \dotprod{\beta v_k +(1-\beta)y_k -x_*, \E{\cZ}^\dagger \E{\cZ}(y_k -x_*)} \nonumber \\
& =& \dotprod{\beta v_k +(1-\beta)y_k -x_*,y_k -x_*} \label{eq:ami_explainsteplat} \\
& =& \dotprod{\beta \left[\frac{1}{\eta}y_k - \frac{1-\eta}{\eta} x_k\right] +(1-\beta)y_k -x_*,y_k -x_*}\nonumber \\
&= & \dotprod{y_k - x_* +\beta \frac{1-\eta}{\eta}(y_k- x_k) ,y_k -x_*} \nonumber \\
& =& \norm{y_k - x_*}^2 + \beta \frac{1-\eta}{\eta}\dotprod{y_k- x_k,y_k-x_*}\nonumber \\
&=&\norm{y_k - x_*}^2 \nonumber \\
&&  - \beta \frac{1-\eta}{2\eta}\left(\norm{x_k-x_*}^2-\norm{y_k- x_k}^2-\norm{y_k- x_*}^2 \right)\label{eq:ami_IIIbnded}
\end{eqnarray}
where in the second equality~\eqref{eq:ami_explainsteplat} we used that $y_k -x_* \in \Range{\A^*} \overset{\eqref{eq:ami_exactness}}{ =} \Range{\E{\cZ}}$ together with a defining property of pseudoinverse operators $\E{\cZ}^\dagger\E{\cZ} w = w$ for all $w  \in \Range{\E{\cZ}}.$ In the last equality~\eqref{eq:ami_IIIbnded} we used yet again the identity~\eqref{eq:ami_paral1} with $u = y_k- x_k$ and $v = y_k-x_*.$

Plugging~\eqref{eq:ami_Ibnded} and~\eqref{eq:ami_IIIbnded} into~\eqref{eq:ami_proofstep1}
and taking conditional expectation gives
\begin{eqnarray}
\E{r_{k+1}^2 \, | \, y_k, v_k, x_k} & =& I + \gamma^2 \E{II \, | \, y_k}  -2 \gamma \E{III \, | \, y_k,v_k,x_k} \nonumber\\
& \overset{\eqref{eq:ami_Ibnded}+\eqref{eq:ami_IIIbnded}+\eqref{eq:ami_Enormnubnd}}{ =}& 
\beta r_k^2 + (1-\beta) \norm{y_k-x_*}^2_{\E{\cZ}^\dagger} + \gamma^2  \nu \norm{y_k-x_*}^2_{\E{\cZ}} \nonumber \\
& & -2\gamma \norm{y_k - x_*}^2 \nonumber \\
&& + \gamma \beta \frac{1-\eta}{\eta}\left(\norm{x_k-x_*}^2-\norm{y_k- x_k}^2-\norm{y_k- x_*}^2 \right)
\nonumber\\
& \overset{\eqref{eq:ami_ykxkident}+\eqref{eq:ami_nubnds}}{\leq } & 
\beta r_k^2 + \frac{1-\beta}{\theta} \norm{y_k-x_*}^2 + \gamma^2  \nu\norm{y_k-x_*}^2 \nonumber \\
&& -\gamma^2  \nu\E{\norm{x_{k+1}-x_*}^2 \, | \, y_k} -2\gamma\norm{y_k - x_*}^2\nonumber \\
& & + \left(   \beta \frac{1-\eta}{2\eta}\left(\norm{x_k-x_*}^2-\norm{y_k- x_*}^2\right)\right). \label{eq:ami_continue_omega_proof}
\end{eqnarray}
Therefore we have that
\begin{eqnarray}
&& \E{r_{k+1}^2  +\gamma^2  \nu\norm{x_{k+1}-x_*}^2\, | \, y_k, v_k, x_k }
\\
& \leq & 
\beta \left(r_k^2 + \underbrace{\gamma\frac{1-\eta}{\eta}}_{P_1}\norm{x_k-x_*}^2\right) \nonumber \\
& &+\left( \underbrace{\frac{1-\beta}{\theta} -2\gamma+\gamma^2  \nu - \beta \gamma\frac{1-\eta}{\eta}}_{P_2}\right)\norm{y_k-x_*}^2.    \nonumber 
\end{eqnarray}
To establish a recurrence, we need to choose the free parameters $\gamma, \eta$ and  $\beta$ so that $P_1 =\gamma^2  \nu$ and $P_2 =0.$ Furthermore we should try to set $\beta$ as small as possible so as to have a fast rate of convergence. 
Choosing $\beta  =1 - \sqrt{\frac{\theta}{\nu}},$ $\gamma = \sqrt{\frac{1}{\theta \nu}},$ $\eta = \frac{1}{1+\gamma\nu}$ gives $P_2 =0$, $\gamma^2 \nu = 1/\theta$ and 
 \begin{eqnarray*}
\E{r_{k+1}^2  +\frac{1}{\theta}\norm{x_{k+1}-x_*}^2\, | \, y_k, v_k, x_k } &\leq &  \left(1 - \sqrt{\frac{\theta}{\nu}} \right) \left(r_k^2 + \frac{1}{\theta}\norm{x_k-x_*}^2\right).
\end{eqnarray*}
Taking expectation and using the tower rules gives the result.\qed
%

\subsection{Changing norm \label{sec:ami_change_norm}}
Given an invertible positive self-adjoint $\cB \in L(\cX),$ suppose we want to find the least norm solution of~\eqref{eq:ami_primal} under the norm defined by
$\norm{x}_\cB  \eqdef \sqrt{\dotprod{\cB x,x}}$ as the metric in $\cX$. That is,
we want to solve
\begin{equation} \label{eq:ami_primalB_change_norm}
x^* \eqdef \arg\min_{x \in \cX} \frac{1}{2}\norm{x-x_0}_\cB^2, \quad \mbox{subject to} \quad \A x = b.
\end{equation}
By changing variables $x = \cB^{-1/2}z$ we have that the above is equivalent to solving
\begin{equation} \label{eq:ami_primalz}
z^* \eqdef \arg\min_{z \in \cX} \frac{1}{2}\norm{z-z_0}^2, \quad \mbox{subject to} \quad \A \cB^{-1/2} z = b,
\end{equation}
with $x^* = \cB^{-1/2}z^*$, and $B^{1/2}$ is the unique symmetric square root of $\cB$ (see Lemma~\ref{lem:ami_squareroot}). We can now apply Algorithm~\ref{alg:ami_SketchJac} to solve~\eqref{eq:ami_primalz} where $\A \cB^{-1/2}$ is  the system matrix.  Let $x_k$ and $v_k$ be the resulting iterates of applying Algorithm~\ref{alg:ami_SketchJac}.
To make explicit this change in the system matrix we define the matrix
\[
\cZ_\cB \eqdef \cB^{-1/2}\A^*\cS_k^*(\cS_k\A \cB^{-1}\A^*\cS_k^*)^{\dagger}\cS_k \A \cB^{-1/2},
\]
and the constants
\begin{equation} \label{eq:ami_muB}
\theta_B  \eqdef   \inf_{x \in \Range{\cB^{-1/2}\A^*}} \frac{\dotprod{\E{\cZ_\cB}x,x}}{\dotprod{x,x}}
\end{equation}
and
\begin{equation} \label{eq:ami_nuB}
\nu_B  \eqdef   \sup_{x \in \Range{\cB^{-1/2}\A^*}} \frac{\dotprod{\E{\cZ_\cB\E{\cZ_\cB}^\dagger \cZ_\cB}x,x}}{\dotprod{\E{\cZ_\cB}x,x}}.
\end{equation}

 Theorem~\ref{theo:conv} then guarantees that
 { \footnotesize
\[
\E{\norm{v_{k+1} -z_*}_{\E{\cZ_\cB }^\dagger}^2 +\frac{1}{\theta_B}\norm{x_{k+1} -z_*}^2 } \leq \left(1 - \sqrt{\frac{\theta_B}{\nu_B}} \right) \E{\norm{v_k -z_*}_{\E{\cZ_\cB }^\dagger}^2 + \frac{1}{\theta_B}\norm{x_k-z_*}^2}.
\]
}
 Reversing our change of variables $\bar{x}_k = \cB^{-1/2}x_k$ and $\bar{v}_k = \cB^{-1/2} v_k$ in the above displayed equation gives
  \begin{eqnarray}
&& \E{\norm{\bar{v}_{k+1} -x_*}_{\cB^{1/2}\E{\cZ_\cB}^\dagger \cB^{1/2}}^2 +\frac{1}{\theta_B}\norm{\bar{x}_{k+1} -x_*}_\cB^2 } \notag \\
&& \qquad \leq  \left(1 - \sqrt{\frac{\theta_B}{\nu_B}} \right) \E{\norm{\bar{v}_k -x_*}_{\cB^{1/2}\E{\cZ_\cB}^\dagger \cB^{1/2}}^2 + \frac{1}{\theta_B}\norm{\bar{x}_k-x_*}_\cB^2}.\label{eq:ami_convzvbar2}
\end{eqnarray} 
Thus we recover the same exact from the main theorem in~\cite{martinrichtarikaccell}, but in a much more general setting.

%

\section{Proof of Corollary~\ref{cor:ami_sss}}

Clearly,  $\mZ=\frac{1}{\mA_{i,i}}\mA^{\frac12}\mS\mS^\top \mA^{\frac12}$, and hence
$ \E{\mZ}=\frac{\mA}{\Tr{\mA}}$ and $ \theta^P=\frac{\lambda_{\min}(\mA)}{\Tr{\mA}}.$
After  simple algebraic manipulations we get
\[
\E{\E{\mZ}^{-\frac12 }\mZ \E{\mZ}^{-1} \mZ\E{\mZ}^{-\frac12 }}=\Tr{\mA}^2\E{ \frac{1}{\mA_{i,i}^2}\mS\mS^\top \mS\mS^\top }  =\Tr{\mA}\diag{\mA_{i,i}^{-1}},
\]
and therefore
$
\nu^P= \lambda_{\max} \E{\E{\mZ}^{-\frac12 }\mZ \E{\mZ}^{-1} \mZ\E{\mZ}^{-\frac12 }} = \frac{\Tr{\mA}}{\min_i \mA_{i,i}}.
$

 \section{Adding a stepsize \label{sec:ami_omega}}
 
 In this section we enrich Algorithm~\ref{alg:ami_SketchJac} with several {\em additional} parameters and study their effect on convergence of the resulting method.
 
 First, we consider an extension of Algorithm~\ref{alg:ami_SketchJac} to a variant which uses a {\em stepsize parameter}  $0<\omega<2$. That is, instead of performing the update  
 \begin{equation} \label{eq:ami_98gf98909} x_{k+1}=y_k - g_k,\end{equation}
 we perform the update
\begin{equation} \label{eq:ami_089d8h09ff} x_{k+1}=y_k-\omega g_k.\end{equation}  Parameters $\eta, \beta, \gamma$ are adjusted accordingly. The resulting method enjoys the rate
\[{\cal O}\left( \left(1-\sqrt{\frac{\nu}{\theta}\omega(2-\omega)}\right)^k \right),\]
recovering the rate from Theorem~\ref{theo:conv} as a special case for $\omega=1$.  The formal statement follows.

\begin{theorem}
\label{thm:ami_omega} 
Let $0 < \omega < 2$ be an arbitrary stepsize and define
\begin{align}
 \alpha \eqdef 2\omega - \omega^2 \geq 0\,.
\end{align}
Consider a modification of Algorithm~\ref{alg:ami_SketchJac} where instead of \eqref{eq:ami_98gf98909} we perform the update \eqref{eq:ami_089d8h09ff}. If we use the  parameters
\begin{align}
 \eta &= \frac{1}{1 + \gamma \nu} &
 \beta  &= 1 -\sqrt{\frac{\theta \alpha}{\nu}} & 
 \gamma &= \sqrt{\frac{\alpha}{\theta \nu}},
\end{align}
then the iterates $\{v_k,x_k\}_{k \geq 0}$ of Algorithm~\ref{alg:ami_SketchJac} satisfy
\begin{align*}
 &\E{\norm{v_{k} -x_*}_{\E{\cZ}^\dagger}^2 +\frac{1}{\theta}\norm{x_{k} -x_*}^2 } 
 \leq \left(1 - \sqrt{\frac{\theta \alpha}{\nu}} \right)^k \E{\norm{v_0 -x_*}_{\E{\cZ}^\dagger}^2 + \frac{1}{\theta}\norm{x_0-x_*}^2}.
\end{align*}
\end{theorem}
\begin{proof} See Appendix~\ref{sec:ami_89g8f9009jJJJ}.
\end{proof}
 
 \section{Allowing for  different $\eta$}\label{sec:ami_alpha}
 
 In this section we study how the choice of  the key parameter $\eta$ affects the convergence  rate.

This parameter determines how much the sequence $y_k = \eta v_k + (1-\eta) x_k$ resembles the sequence given by $x_k$ or by $v_k$. For instance, when $\eta = 0$, $y_k \equiv x_k$, i.e., we recover the steps of the non-accelerated method, and thus one would expect to obtain the same convergence rate as the non-accelerated method. Similar considerations hold in the other extreme, when $\eta \to 1$. We investigate this hypothesis, and especially discuss how $\beta$ and $\gamma$ must be chosen as a function of $\eta$ to ensure convergence. 
 
 The following statement is a generalization of Theorem~\ref{theo:conv}. For simplicity, we assume that the optional stepsize that was introduced in Theorem~\ref{thm:ami_omega} is set to one again, $\omega \equiv 1$.


\begin{theorem}
\label{thm:ami_family}
Let $0 < \eta < 1$ be fixed. Then the iterates $\{v_k,x_k\}_{k \geq 0}$ of Algorithm~\ref{alg:ami_SketchJac} 
with parameters
    \begin{align}
    \beta(s)  &= \frac{1+s - s \sqrt{\frac{\nu + 4\theta s - 2\nu s + \nu s^2}{\nu s^2}}}{2s}\,, &
    \gamma(s) &= \frac{1}{(1-s \beta(s))\nu}\,. \label{eq:ami_def_bg}
    \end{align}
where $\tau \eqdef \frac{1-\eta}{\eta}$ and $s \eqdef \frac{\tau}{\beta \gamma}$, 
satisfy
\begin{align*}
 &\E{\norm{v_{k} -x_*}_{\E{\cZ}^\dagger}^2 + \gamma \tau \norm{x_{k} -x_*}^2 } 
 \leq \rho^k \E{\norm{v_0 -x_*}_{\E{\cZ}^\dagger}^2 + \gamma \tau\norm{x_0-x_*}^2}.
\end{align*}
(or put differently):
\begin{align*}
 &\E{\norm{v_{k} -x_*}_{\E{\cZ}^\dagger}^2 + (1-\eta)\gamma \norm{x_{k} -x_*}^2 } 
 \leq \rho^k \E{\norm{v_0 -x_*}_{\E{\cZ}^\dagger}^2 + (1-\eta)\gamma \norm{x_0-x_*}^2}.
\end{align*}
where $\rho = \max\{\beta(s),s\beta(s)\}\leq 1$.
\end{theorem}

We can now exemplify a few special parameter settings.

\begin{example}
For $\eta = 1$, i.e., if $s \to 0$, we get the rate $\rho = 1-\frac{\theta}{\nu}$ with $\beta = 1-\frac{\theta}{\nu}$, $\gamma = \frac{1}{\nu}$.
\end{example}

\begin{example}
For $\eta \to 0$, i.e., in the limit $s \to \infty$, we get the rate  $\rho = 1-\frac{\theta}{\nu}$. 
\end{example}

\begin{example}
The rate $\rho$ is minimized for $s=1$, i.e.,  $\beta=1- \sqrt{\frac{\nu}{\theta}}$ and $\gamma = \sqrt{\frac{1}{\theta \nu}}$; recovering Theorem~\ref{theo:conv}.
\end{example}

The best case, in terms of convergence rate for both non-unit stepsize and a variable parameter choice happened to be the default parameter setup. The non-optimal parameter choice was studied in order to have theoretical guarantees for a wider class of parameters, as in practice one might be forced to rely on sub-optimal / inexact parameter choices.

\section{Proof of Theorem~\ref{thm:ami_omega}} \label{sec:ami_89g8f9009jJJJ}

The proof follows by slight modifications of the proof of Theorem~\ref{theo:conv}.

First we adapt Lemma~\ref{lemma:threepoint}. As we have $x_{k+1} - x_* = (\cI-\omega \cZ_k)(y_k - x_*)$ the following statement follows by the same arguments as in the proof of Lemma~\ref{lemma:threepoint}.
\begin{lemma}[Lemma~\ref{lemma:threepoint}']
\label{lemma:mod_lemma}
\begin{equation}
 \alpha \norm{y_k-x_*}_{\E{\cZ}}^2 = \norm{y_k-x_*}^2 -\E{\norm{x_{k+1}-x_*}^2 \, | \, y_k}
 \end{equation}
\end{lemma}
\begin{proof}
\begin{eqnarray*}
\E{\norm{x_{k+1}-x_*}^2 \, | \, y_k} &= & \E{\norm{(\cI-\cZ_k)( y_k -x_*)}^2 \, | \, y_k} \\
&= & \E{\dotprod{(\cI-\omega \cZ_k)( y_k -x_*),(\cI-\omega \cZ_k)y_k -x_*}}\\
&= &\norm{y_k -x_*}^2 - \alpha \norm{y_k -x_*}^2_{\E{\cZ}}. 
\end{eqnarray*} 
\end{proof}

We now follow the same steps as in proof of Theorem~\ref{theo:conv} in Section~\ref{subsec:proof_main_thm}.
We observe, that the first time Lemma~\ref{lemma:threepoint} is applied is in equation~\eqref{eq:ami_continue_omega_proof}. Using Lemma~\ref{lemma:mod_lemma} instead, gives

\begin{eqnarray}
&& \E{r_{k+1}^2 \, | \, y_k, v_k, x_k} 
\\
& \leq &  
\beta r_k^2 + \frac{1-\beta}{\theta} \norm{y_k-x_*}^2 + \frac{\gamma^2  \nu}{\alpha} \left( \norm{y_k-x_*}^2 -\E{\norm{x_{k+1}-x_*}^2 \, | \, y_k}\right) \nonumber \\
& & +2\gamma \left(-\norm{y_k - x_*}^2  + \beta \frac{1-\eta}{2\eta}\left(\norm{x_k-x_*}^2-\norm{y_k- x_*}^2\right)\right).
\end{eqnarray}
Therefore we have that
\begin{eqnarray}
&& \E{r_{k+1}^2  +\gamma^2  \nu\norm{x_{k+1}-x_*}^2\, | \, y_k, v_k, x_k }
\\
& \leq & 
\beta \left(r_k^2 + \underbrace{\gamma\frac{1-\eta}{\eta}}_{P_1'}\norm{x_k-x_*}^2\right) \nonumber \\
& &+\left( \underbrace{\frac{1-\beta}{\theta} -2\gamma+ \frac{\gamma^2  \nu}{\alpha} - \beta \gamma\frac{1-\eta}{\eta}}_{P_2'}\right)\norm{y_k-x_*}^2.    \nonumber 
\end{eqnarray}
Noting that $\frac{1-\eta}{\eta}= \gamma \nu$ and $\frac{\gamma^2 \nu}{\alpha} = \frac{\gamma(1-\eta)}{\alpha \eta} = \frac{1}{\theta}$, we observe $P_2' = 0$ and deduce the statement of Theorem~\ref{thm:ami_omega}.

\section{Proof of Theorem~\ref{thm:ami_family}} \label{sec:ami_j98**JJ*(}

It suffices to study equation~\eqref{eq:ami_continue_omega_proof}. We observe that for convergence the big bracket, $P_2$, should be negative,
\begin{align}
(1-\beta) \frac{1}{\theta} + \gamma^2 \nu - 2\gamma - \gamma \beta \frac{1-\eta}{\eta} \leq 0 \label{eq:ami_cond1}
\end{align} 
The convergence rate is then
\begin{align}
\rho \eqdef \max \left\{\beta, \frac{(1-\eta)\beta}{\eta \gamma \nu}\right\}\,. \label{eq:ami_rate}
\end{align}
or in the notation of Theorem~\ref{thm:ami_family}, $\rho = \max\{\beta,s\beta\}$.

This means, that in order to obtain the best convergence rate, we should therefore choose parameters $\beta$ and $\gamma$ such that $\beta$ is as small as possible. This observation is true regardless of the value of $s$ (which itself depends on $\gamma$).

With the notation $\tau=s \gamma \beta$, we reformulate~\eqref{eq:ami_cond1} to obtain
\begin{align}
\frac{1}{\theta} + \gamma^2 \nu - 2 \gamma \leq \beta \left( \frac{1}{\theta}  + s \gamma^2 \nu \right)
\end{align}
Thus we see, that $\beta$ cannot be chosen smaller than
\begin{align}
\beta^*(s,\gamma) = \frac{1 + \theta \gamma^2 \nu - 2 \theta \gamma}{1 + s \theta \gamma^2 \nu}
\end{align}
Minimizing this expression in $\gamma$ gives
\begin{align}
\beta^*(s) = \frac{1+s - s \sqrt{\frac{\nu + 4\theta s - 2\nu s + \nu s^2}{\nu s^2}}}{2s} \label{eq:ami_beta_star}
\end{align}
with $\gamma^*(s) = \frac{1}{(1-s \beta^*(s))\nu}$.

We further observe that this parameter setting indeed guarantees convergence, i.e. $\rho \leq 1$. From~\eqref{eq:ami_beta_star} we observe ($\nu > 0$, $s \geq 0$, $\theta \geq 0$):
\begin{align}
\beta^*(s) \leq  \frac{1+s - \sqrt{\frac{\nu - 2\nu s + \nu s^2}{\nu}}}{2s} = \frac{1+s - (s-1)}{2s} = \frac{1}{s}
\end{align}
Hence $ s \beta^*(s)  \leq 1$. On the other hand, $(1-s) \leq \sqrt{(1-s)^2 + \frac{4\theta s}{\nu}}$ and hence $(1+s) -\sqrt{(1-s)^2 + \frac{4\theta s}{\nu}} \leq 2s$, which shows $\beta^* (s) \leq 1$.

\section{Proofs and further comments on Section~\ref{sec:ami_asqn_pap}}

\subsection{Proof of Theorem~\ref{theo:qn}}

We perform a change of coordinates since it is easier to work with the standard Frobenius norm as opposed to the weighted Frobenius norm. Let $\hat{\mX} = \mA^{1/2}\mX \mA^{1/2}$  
so that~\eqref{eq:ami_primalqN} and~\eqref{eq:ami_qunac} become
\begin{equation} \label{eq:ami_primalqNH}
\hat{\mX}_* \eqdef    \arg\min \norm{\hat{\mX}}_{F}^2 \quad \mbox{subject to} \quad  \hat{\mX} = \mI, \quad \hat{\mX} = \hat{\mX}^\top,
\end{equation}
and
\begin{equation}\label{eq:ami_qunacI}
\hat{\mX}_{k+1}  =\mP+ \left(\mI-\mP\right) \hat{\mX}_{k}\left(\mI -\mP \right),
\end{equation}
respectively, where $\mP = \mA^{1/2}\mS(\mS^\top \mA\mS)^{-1}\mS^\top \mA^{1/2}.$  The linear operator that encodes the constaint in~\eqref{eq:ami_primalqNX} is given by
$\hat{\A}(\mX) = \left( \mX, \, \mX - \mX^\top\right)$ the adjoint of which is given by
$\hat{\A}^*(\mY_1,\mY_2) = \mY_1 + \mY_2-\mY_2^\top.$ Since
$\hat{\A}^*$ is clearly surjective, it follows that $\Range{\hat{\A}^*} = \R^{d \times d}$. 

Subtracting the identity matrix from both sides of~\eqref{eq:ami_qunacI} and using that $\mP$ is a projection matrix, we have that
\begin{equation}\label{eq:ami_qunacIres}
\hat{\mX}_{k+1} -\mI =\left(\mI-\mP\right) (\hat{\mX}_{k}-\mI)\left(\mI -\mP \right).
\end{equation}
To determine the $\mZ$ operator~\eqref{eq:ami_Z}, from~\eqref{eq:ami_IZprojres} and~\eqref{eq:ami_qunacIres} we know that
\[ \left(\mI-\mP\right) (\hat{\mX}_{k}-\mI)\left(\mI -\mP \right)= (\mI - \mZ) (\hat{\mX}_k -\mI).\]
Thus for every matrix $\mX \in \R^{d \times d}$ we have that
\begin{equation} \label{eq:ami_Zqn}
 \mZ(\mX) = \mX - \left(\mI-\mP\right) \mX\left(\mI -\mP \right) = \mX \mP + \mP\mX(\mI-\mP).\end{equation}
Denote column-wise vectorization of $\mX$ as $x$: $x\eqdef \Vect{\mX}$. To calculate a useful lower bound on~$\theta$, note that
\begin{eqnarray} \label{eq:ami_XZXlower}
 \Tr{\mX^\top \mZ(\mX)} &= &\Tr{\mX^\top \mX \mP} + \Tr{\mX^\top \mP\mX(\mI-\mP)}\nonumber 
 \\
 &= &
 x^\top \Vect{\mX \mP}+x^\top \Vect{ \mP\mX(\mI-\mP)}\nonumber 
 \\
  &= &
 x^\top (\mP\otimes \mI)x +x^\top((\mI-\mP)\otimes \mP)x \nonumber 
 \\
& \overset{\eqref{eq:ami_bigz}}{= } &
   x^\top \bigZ x ,
   \label{eq:ami_XZXeq}
\end{eqnarray} 
where we used that $\Tr{\mA^\top \mB}=\Vect{\mA}^\top \Vect{\mB}$ and $\Vect{ \mA\mX \mB}=(\mB^\top \otimes \mA) \Vect{x}$ holds for any $\mA,\mB,\mX$.

Consequently, $\theta$ is equal to
\[
\theta \overset{\eqref{eq:ami_nuqn}}{=}  \inf_{\mX \in \R^{n\times n}} \frac{\dotprod{\E{Z}\mX,\mX}_F}{\norm{\mX}_F^2}
 \overset{\eqref{eq:ami_XZXeq}}{= }  \inf_{x \in \R^{n^2\times n^2}} \frac{   x^\top \E{\bigZ} x  }{x^\top x} 
 = \lambda_{\min}(\E{\bigZ}).
\]
Notice that we have $2\lambda_{\min}(\E{\mP}) \geq\lambda_{\min}(\E{\bigZ})\geq \lambda_{\min}(\E{\mP})$ since $(\mP\otimes \mI)+(\mI\otimes \mP)\geq \bigZ\geq (\mP\otimes \mI)$.

In light of Algorithm~\ref{alg:ami_SketchJac}, the iterates of the accelerated version of~\eqref{eq:ami_qunacI} are given by
\begin{eqnarray}
 \hat{\mY}_k &=& \eta \hat{\mV}_k + (1-\eta) \hat{\mX}_k \nonumber\\
 \hat{\mG}_k & =& \mZ_k(\hat{\mY}_k - \mI ) \nonumber\\
\hat{\mX}_{k+1} &=& \hat{\mY}_k - \hat{\mG}_k \nonumber\\
\hat{\mV}_{k+1} &=& \beta \hat{\mV}_k +(1-\beta)\hat{\mY}_k - \gamma \hat{\mG}_k \label{eq:ami_qunacacc}\
\end{eqnarray}
where $\hat{\mY}_k,\hat{\mV}_{k},\hat{\mG} \in \R^{n\times n}.$ From Theorem~\ref{theo:conv} we have that $\hat{\mV}_k$ and $\hat{\mX}_k$ converge to the identity matrix according to
\begin{equation} \label{eq:ami_aas89hah}
 \E{\norm{\hat{\mV}_{k+1} -I}_{\E{\mZ}^\dagger}^2 +\frac{1}{\theta}\norm{\hat{\mX}_{k+1} -I}^2_{F} } \leq \left(1 - \sqrt{\frac{\theta}{\nu}} \right) \E{\norm{\hat{\mV}_k -\mI}_{\E{\mZ}^\dagger}^2 + \frac{1}{\theta}\norm{\hat{\mX}_k-\mI}^2_{F}},
 \end{equation}
where $\norm{\mX}_{\E{\mZ}^\dagger}^2 = \dotprod{\E{\mZ}^\dagger \mX, \mX}_F.$
Changing coordinates back to $ \hat{\mX}_k = \mA^{1/2}\mX_k \mA^{1/2}$ and defining
$\mY_k \eqdef \mA^{-1/2} \hat{\mY}_k \mA^{-1/2}$, $\mV_k \eqdef \mA^{-1/2} \hat{\mV}_k \mA^{-1/2}$ and $ \mG_k \eqdef A^{-1/2} \hat{\mG}_k \mA^{-1/2}$, we have that~\eqref{eq:ami_aas89hah} gives~\eqref{eq:ami_qnaccconv}. Furthermore, using the same coordinate change applied to the iterates~\eqref{eq:ami_qunacacc} gives Algorithm~\ref{alg:ami_qn}.

\subsection{Matrix inversion as linear system~\label{sec:ami_alternate}}
Denote $x=\Vect{\mX}$, i.e. $x$ is $d^2$ dimensional vector such that $x_{(n(i-1)+1):ni}=\mX_{:,i}$. Similarly, denote $e=\Vect{\mI}$. System \eqref{eq:ami_system} can be thus rewritten as 
\begin{equation}
 (\mI\otimes \mA)x=e.
\label{eq:ami_system_vector}
\end{equation}

Notice that all linear sketches of the original system $\mA \mX=\mI$ can be written as 
\begin{equation}\label{eq:ami_sketch_general}
\gS^\top  (\mI\otimes \mA)x=\gS^\top  e
\end{equation}
for a suitable $d^2\times d^2$ matrix $\gS$, therefore the setting is fairly general.

\subsubsection{Alternative proof of Theorem \ref{theo:qn}}
Let us now, for a purpose of this proof, consider sketch matrix $\gS$ to capture only sketching the original matrix system $\mA\mX=\mI$ by left multiplying by $\mS$, i.e. $\gS=(\mI\otimes \mS)$, as those are the considered sketches in the setting of Section~\ref{sec:ami_asqn_pap}. 

As we have 
\[
\Tr{\mB\mX^\top \mB\mX}=\Vect{\mB\mX \mB}^\top x=x^\top (\mB\otimes \mB) x, 
\]
weighted Frobenius norm of matrices is equivalent to a special weighted euclidean norm of vectors. Define also $\mC$ to be a matrix such that $\mC x=0$ if and only if $\mX=\mX^\top$. Therefore, \eqref{eq:ami_primalqNX} is equivalent to 

\begin{equation} \label{eq:ami_primalqNXvec}
x_{k+1}=   \arg\min \norm{x - x_k}_{\mA\otimes \mA}^2 \quad \mbox{subject to} \quad ( \mI\otimes \mS^\top) (\mI\otimes \mA) x = ( \mI\otimes \mS^\top) e, \quad \mC x = 0,
\end{equation}
which is a sketch-and-project method applied on the linear system, with update as per \eqref{eq:ami_qunac}:
\[
x^{k+1}=x^k-(\mH\otimes \mI)((\mI\otimes \mA) x -e ) -(\mI\otimes \mH)((\mI\otimes \mA) x -e )  + (\mH\mA\otimes \mH)  ((\mI\otimes \mA) x -e )
\]
for 
$\mH\eqdef \mS\left(\mS^\top \mA\mS\right)^{-1}\mS^\top.$
Using substitution $\hat{x}=(\mA^{\frac12}\otimes \mA^\frac12)x; \hat{\mS}=\mA^{\frac12}\mS$ and comparing to \eqref{eq:ami_IZprojres}, we get 
\[
\mZ=\mI\otimes \mI-(\mI-\mP)\otimes(\mI-\mP)
\]
for $\mP$ as defined inside the statement of Theorem~\ref{theo:qn}. 
Therefore, we have all necessary information to apply the results from \cite{martinrichtarikaccell}, recovering Theorem \ref{theo:qn}.

\section{Linear operators in Euclidean spaces}
Here we provide some technical lemmas and results for linear operators in Euclidean space, that we used in the main body of the chapter. Most of these results can be found in standard textbooks of analysis, such as~\cite{pedersen1996}. We give them here for completion.
  
Let $\cX, \cY, \cZ$ be Euclidean spaces, equipped with inner products. Formally, we should use a notation that distinguishes the inner product in each space. But instead we use $\dotprod{\cdot,\,\cdot}$ to denote the inner product on all spaces, as it will be easy to determine from which space the elements are in. That is, for $x_1,x_2 \in \cX$, we denote by $\dotprod{x_1,x_2}$ the inner product between $x_1$ and $x_2$ in $\cX.$

Let
\[\norm{T} \eqdef \sup_{\norm{x}\leq 1} \norm{T x},\]
denote the operator norm of $T$.  Let $0 \in L(\cX,\cY)$ denote the zero operator and $I \in L(\cX,\cY)$ the identity map.

\paragraph{The adjoint.}
Let $T^* \in L( \cY , \cX)$ denote the unique operator that satisfies
\[\dotprod{Tx,y} = \dotprod{x,T^*y},\]
for all $x \in \cX$ and $y \in \cY.$ We say that $T^*$ is the \emph{adjoint} of $T$. We say $T$ is \emph{self-adjoint} if $T =T^*.$
Since for all $x\in \cX$ and $s\in \cS$,
\[ \langle x , (ST)^*s \rangle = \langle ST  x , s \rangle_\cS = \langle T  x , S^* s \rangle_\cY =  \langle  x , T^* S^* s \rangle,\]
we have \[(ST)^* = T^* S^*.\]

\begin{lemma} \label{lem:ami_decomp}For $T \in L(\cX,\cY)$ we have that $\Range{T^*}^{\perp} = \Null{T}.$ Thus
	\begin{eqnarray}
	\cX &= & \Range{T^*} \oplus \Null{T}\\
	\cY &= & \Range{T} \oplus \Null{T^*}
	\end{eqnarray}
\end{lemma}
\begin{proof}
	See 3.2.6 in~\cite{pedersen1996}.
\end{proof}

\subsection{Positive operators}
We say that $G\in L(\cX)$ is positive if it is self-adjoint and if
$\dotprod{x,Gx} \geq 0$ for all $ x \in \cX$.
Let $(e_j)_{j=1}^{\infty} \in \cX$ be an orthonormal basis. The trace of $G$ is defined as
\begin{equation}\label{eq:ami_tracedef}
\Tr{G} \eqdef \sum_{j=1}^{\infty} \dotprod{G e_j, e_j}.
\end{equation}
The definition of trace is independent of the choice of basis due to the following lemma.
\begin{lemma} If $U$ is unitary and $G\geq 0$ then
$\Tr{UGU^*} = \Tr{G}.$
\end{lemma}
\begin{proof}
	See 3.4.3 and 3.4.4 in~\cite{pedersen1996}.
\end{proof}
\begin{lemma} \label{lem:ami_projrank} If $P \in L(\cX)$ is a projection matrix then
	$\Tr{P} = \dim(\Range{P}) = \Rank{P}.$
\end{lemma}
\begin{proof}  Let $d = \dim(\Range{P})$ which is possibly infinite.
	Given that $P$ is a projection we have that $\Range{P}$ is a closed subspace and thus there exists orthonormal basis  $(e_j)_{j=1}^{d} $ of $\Range{P}$. Consequently,
	$\Tr{P} \overset{\eqref{eq:ami_tracedef}}{=} \sum_{j=1}^{d} 1 = d= \dim(\Range{P}).$
\end{proof}

A \emph{square root} of an operator $G \in L(\cX)$ is an operator $R\in L(\cX) $ such that $R^2 =G.$
\begin{lemma} \label{lem:ami_squareroot}
	If $G: \cX \rightarrow \cX$ is positive, then there exists a unique positive square root of $G$ which we denote by $G^{1/2}.$ 
\end{lemma}
\begin{proof}
	See 3.2.11 in~\cite{pedersen1996}. 
\end{proof}

\begin{lemma} \label{lem:ami_decompTGT} For any $T \in L(\cX,\cY)$ and any $G \in L(\cY,\cY)$ that is positive and injective,
	\begin{equation} \label{eq:ami_NullTGT}
	\Null{T} =  \Null{T^*GT},
	\end{equation}
	and
	\begin{equation}\label{eq:ami_RangeTGT}
	\overline{\Range{T^*}} =  \overline{\Range{T^*GT}}.
	\end{equation}
\end{lemma}
\begin{proof}
	The inclusion $\Null{T} \subset \Null{T^*GT}$ is immediate. For the opposite inclusion, let $x \in \Null{T^*GT}.$ Since $G$ is positive we have by Lemma~\ref{lem:ami_squareroot} that there exists a square root with $G^{1/2}G^{1/2} =G.$ 
Therefore,
	$\dotprod{x,T^*GT x} = \dotprod{G^{1/2}Tx,G^{1/2}Tx} =0,$
	which implies that $G^{1/2}Tx =0$. Since $G$ is injective, it follows that $G^{1/2}$ is injective and thus $x \in \Null{T}$. Finally~\eqref{eq:ami_RangeTGT} follows by taking the orthogonal complements of~\eqref{eq:ami_NullTGT} and observing Lemma~\ref{lem:ami_decomp}. 
\end{proof}

As an immediate consequence of~\eqref{eq:ami_NullTGT} and~\eqref{eq:ami_RangeTGT} we have the following lemma.
\begin{corollary} For $G: \cX \rightarrow \cX$ positive we have that
	\begin{eqnarray}
	\Null{G^{1/2}} & =& \Null{G} \label{eq:ami_NullGhalf}\\
	\overline{\Range{G^{1/2}}} & =& \overline{\Range{G}}  \label{eq:ami_RangeGhalf}
	\end{eqnarray}
\end{corollary}

\subsection{Pseudoinverse}
 For a bounded linear operator $T$ define the pseudoinverse of $T$ as follows.
\begin{definition} \label{def:pseudoinverse} Let $T \in L(\cX,\cY)$ such that $\Range{T}$ is closed. $T^{\dagger}: \cY \rightarrow \cX$ is said to be the pseudoinverse if
\begin{enumerate}
\item $T^{\dagger} Tx = x $ for all $x \in \Range{T^* }.$
\item $T^{\dagger} x = 0$ for all $x \in \Null{T^* }.$
\item If $x \in \Null{T}$ and $y \in \Range{T^* }$ then $T^{\dagger}(x+y)=T^{\dagger}x +T^{\dagger}y.$
\end{enumerate}
\end{definition}
It follows directly from the definition (see~\cite{Desoer1963} for details) that $T^{\dagger}$ is a unique bounded linear operator. The following properties of pseudoinverse will be important.

\begin{lemma}[Properties of pseudoinverse] \label{lem:ami_pseudo} Let $T \in L(\cX,\cY)$ such that $\Range{T}$ is closed. It follows that 
\begin{enumerate}
\item $ T T^\dagger T = T$ \label{it:pseudoTTdagT}
\item $\Range{T^{\dagger}} = \Range{T^*}$ and $\Null{T^{\dagger}} = \Null{T^*}$\label{it:pseudorange}
\item $(T^*)^{\dagger} = (T^{\dagger})^*$ \label{it:pseudoadjoint}
\item If $T$ is self-adjoint and positive then $T^{\dagger}$ is self-adjoint and positive.  \label{it:pseudoadjpos}
\item   $T^{\dagger}T T^* = T^*$, that is, $T^{\dagger}T$ projects orthogonally onto $ \Range{T^* }$ and along $\Null{T}.$  \label{it:pseudoproj}
\item Consider the linear system $Tx=d$ where $d \in \Range{T}$. It follows that
\begin{equation}
\textstyle
\label{eq:ami_pseudoleastnorm}T^{\dagger}d = \arg \min_{x \in \cX} \frac{1}{2}\norm{x}^2 \quad \mbox{subject to} \quad Tx=d.
\end{equation}     \label{it:pseudoleastnorm}
\item $T^{\dagger} = T^* (TT^*)^{\dagger}$ \label{it:pseudoTadjTTadj}
\end{enumerate}
\end{lemma}
\begin{proof}
The proof of first five items can be found in~\cite{Desoer1963}. The proof of \eqref{it:pseudoleastnorm} is alternative characterization of the pseudoinverse and it can be established by using that $d \in \Range{T}$ together with item~{\it\ref{it:pseudoTTdagT}} thus $TT^{\dagger}d =d$. The proof then follows by using the orthogonal decomposition $\Range{T^*} \oplus \Null{T}$ to show that $T^{\dagger}d $ is indeed the minimum of~\eqref{eq:ami_pseudoleastnorm}. Finally item~\eqref{it:pseudoTadjTTadj} is a direct consequence of the previous items.
\end{proof}


\chapter{Accepted Papers}

\begin{quote}
\cite{sega} \bibentry{sega}.
\end{quote}

\begin{quote}
\cite{hanzely2018accelerated} \bibentry{hanzely2018accelerated}.
\end{quote}

\begin{quote}
\cite{gower2018accelerated} \bibentry{gower2018accelerated}.
\end{quote}

\begin{quote}
\cite{dutta2019nonconvex} \bibentry{dutta2019nonconvex}.
\end{quote}

\begin{quote}
\cite{sigma_k} \bibentry{sigma_k}.
\end{quote}  

\begin{quote}
\cite{mishchenko201999} \bibentry{mishchenko201999}.
\end{quote} 

\begin{quote}
\cite{dutta2019best} \bibentry{dutta2019best}.
\end{quote} 
                   
\begin{quote}
\cite{asvrcd} \bibentry{asvrcd}.
\end{quote}  

\begin{quote}
\cite{sscn} \bibentry{sscn}.
\end{quote}                  
              

\chapter{Submitted Papers}

\begin{quote}
\cite{abpg} \bibentry{abpg}.
\end{quote}

\begin{quote}
\cite{hanzely2019one} \bibentry{hanzely2019one}.
\end{quote}

\begin{quote}
\cite{local} \bibentry{local}.
\end{quote}  

\begin{quote}
\cite{fedoptimal} \bibentry{fedoptimal}.
\end{quote}

\end{document}